\UseRawInputEncoding
\documentclass[11pt,chapterheads]{book}
\usepackage{graphicx}
\usepackage{amsmath, amscd, amssymb, amsthm}
\usepackage[subnum]{cases}
\usepackage{changepage}
\usepackage{xcolor}
\usepackage{marginnote}
\usepackage[implicit=false, bookmarks=false]{hyperref}
\usepackage{cleveref}
\usepackage[all]{xy}
\usepackage{tabularx}
\usepackage{ltablex}
\usepackage{longtable}
\usepackage[T1]{fontenc}
\usepackage{mathpazo}
\usepackage{pdflscape}
\usepackage{enumitem}
\setlist[itemize]{leftmargin=+0.1in}

\usepackage[margin=1.1in]{geometry}

\newtheorem{theorem}{Theorem}[section]
\newtheorem{lemma}[theorem]{Lemma}

\newtheorem{proposition}[theorem]{Proposition}
\newtheorem{corollary}[theorem]{Corollary}
\newtheorem{conjecture}[theorem]{Conjecture}

\theoremstyle{definition}
\newtheorem{definition}[theorem]{Definition}
\newtheorem{remark}[theorem]{Remark}
\numberwithin{equation}{section}
\newtheorem{example}[theorem]{Example}
\newtheorem{assumption}[theorem]{Assumption}
\newtheorem{setting}[theorem]{Setting}

\usepackage{fancyhdr}
\begin{document}

\normalfont

\title{Geometric and Representation Theoretic Aspects of $p$-adic Motives}

\author{Xin Tong}
\date{}

\maketitle

\newpage

\section*{Abstract}
In this dissertation, we discuss mainly the corresponding geometric and representation theoretic aspects of relative $p$-adic Hodge theory and $p$-adic motives. To be more precise, we study the corresponding analytic geometry of the corresponding spaces over and attached to period rings in the relative $p$-adic Hodge theory, including derived topological de Rham complexes and derived topological logarithmic de Rham complexes after Bhatt, Gabber, Guo and Illusie which is in some sense equivalent to the derived prismatic cohomology of Bhatt-Scholze as shown in the work of Li-Liu, $\mathcal{O}\mathbb{B}_\mathrm{dR}$-sheaves after Scholze, $\varphi$-$\widetilde{C}_X$-sheaves and relative-$B$-pairs after Kedlaya-Liu, multidimensional rings after Carter-Kedlaya-Z\'abr\'adi and Pal-Z\'abr\'adi and many other possible general universal motivic rings or sheaves. Many contexts are expected to be sheafified, such as over Scholze's pro-\'etale sites of the considered analytic spaces by using perfectoids or the quasisyntomic sites by using quasiregular semiperfectoids as in the work of Bhatt-Morrow-Scholze and Bhatt-Scholze. The main motivation comes from the corresponding noncommutative Tamagawa Number conjectures after Burns-Flach-Fukaya-Kato, relative version of the generalized version of the period rings as in the work of Carter-Kedlaya-Z\'abr\'adi and Pal-Z\'abr\'adi, arithmetic families of the representations of fundamental groups in analytic geometry such as for analytification of the moduli stacks of algebraic curves after Reinecke, arithmetic families of general motivic structures in analytic geometry such as in the work of Andreatta-Brinon, Andreatta-Iovita, Berger, Bhatt-Morrow-Scholze, Bhatt-Scholze, Fargues-Fontaine, Fargues-Scholze, Fontaine and Kedlaya-Liu,  noncommutative analytic geometry and noncommutative deformation, derived noncommutative analytic geometry and derived noncommutative deformation, Langlands programs, analytic approach to algebraic topology and so on. Due to the natural though not functorial correspondence between the linear topology and the one induced by a Banach norm, we do not restrict ourselves to the functional analytic point of view when we take completion after Bambozzi-Ben-Bassat-Kremnizer, Clausen-Scholze, Gabber-Ramero, Huber, Kedlaya-Liu and Scholze.

\newpage
\section*{Acknowledgements}
We would like to thank Professor Kedlaya for his advice and direction along our study and research in the past years. All the results presented in this dissertation are influenced significantly by our discussion, which directly leads to the current presentation.\\

{\color{blue}\tableofcontents}

\newpage

\chapter{Preliminary and Introduction}

\newpage

\section{Preliminary}

This dissertation is a study of geometric and representation theoretic aspects of the corresponding $p$-adic motives. We have the following foundational materials.

\noindent  1. Noncommutative Motives: \cite{1NS},  \cite{1Ta};\\
\noindent 2. Noncommutative Harmonic Analysis, Noncommutative Microlocal Analysis and Pseudodifferential Analysis: \cite{1LVO};\\
\noindent 3. $\infty$-Categorical Topological Analysis: \cite{1BK}, \cite{1BBK}, \cite{1BBBK}, \cite{KKM}, \cite{BBM}, \cite{1CS1}, \cite{1CS2};\\
\noindent 4. Topological Division Rings and Topological Vector Spaces: \cite{1Bou};\\
\noindent 5. Topological Rings and Topological Modules: In general we need to consider very general topological rings not necessarily commutative, for instance see \cite{1AGM}, \cite{1CBCKSW}, \cite{1Hu1}, \cite{1Hu2},  \cite{1SW}, \cite{1TGI}, \cite{1U}, \cite{1W};\\
\noindent 6. Functional Analytic Rings and Functional Analytic Modules: In general we need to consider very general topological rings not necessarily commutative or archimedean, for instance see \cite{1AGM}, \cite{1BGR}, \cite{1KL1}, \cite{1KL2}, \cite{1U},   \cite{1W}. Much of the corresponding discussion in \cite{1KL1} and \cite{1KL2} works for noncommutative seminormed rings, also see \cite{1He}; \\
\noindent 7. Adic Rings:  For commutative see \cite{1CBCKSW}, \cite{1Hu1}, \cite{1Hu2},  \cite{1KL1}, \cite{1KL2}, \cite{1SW},  and for noncommutative see \cite{1FK};\\ 
\noindent 8. Distinguished Deformations of Rings:  \cite{1BMS1}, \cite{1BMS2}, \cite{1BS}, \cite{1CBCKSW}, \cite{1GR}, \cite{1KL1}, \cite{1KL2}, \cite{1Sch1}; \\
\noindent 9. Commutative Algebra: \cite{BourbakiAC}, \cite{1Lurie2}, \cite{1Lurie3}, \cite{1R}, \cite{1SP};\\
\noindent 10. Schemes: \cite{1EGAI}, \cite{1EGAII}, \cite{1EGAIII1}, \cite{1EGAIII2}, \cite{1EGAIV1},\\
 \cite{1EGAIV1}, \cite{1EGAIV1}, \cite{1EGAIV1}, \cite{1SGAI}, \cite{1SGAII}, \cite{1SGAIII1}, \cite{1SGAIII2}, \cite{1SGAIII3}, \cite{1SGAIV1}, \cite{1SGAIV2}, \cite{1SGAIV3}, \cite{1SGAIV.5}, \cite{1SGAV}, \cite{1SGAVI}, \cite{1SGAVII1},\\
  \cite{1SGAVII2}, \cite{1SP};\\
\noindent 11. Huber Spaces: \cite{1CBCKSW}, \cite{1Hu1}, \cite{1Hu2}, \cite{1KL1}, \cite{1KL2}, \cite{1Sch1}, \cite{1Sch2}, \cite{1Sch3}, \cite{1SW};\\
\noindent 12. Analytic Spaces: \cite{1BBBK}, \cite{1BK}, \cite{1BBK}, \cite{1CBCKSW}, \cite{1CS1}, \cite{1CS2}, \cite{1FS}, \cite{1Hu1}, \cite{1Hu2},  \cite{1KL1}, \cite{1KL2}, \cite{1Sch1}, \cite{1Sch2}, \cite{1Sch3}, \cite{1SW};\\
\noindent 13. Algebraic Topology: \cite{1E}, \cite{1M}, \cite{1MP}, \cite{1N};\\
\noindent 14. $\infty$-Categories and Their Models: \cite{1Bergner}, \cite{1Ci}, \cite{1J}, \cite{1Lurie1};\\
\noindent 15. Higher Algebra, Higher Toposes, Higher Geometries: \cite{1Lurie1}, \cite{1Lurie2}, \cite{1Lurie3},  \cite{1TV1}, \cite{1TV2}. Actually the corresponding $\infty$-categories of Banach ring spectra as in \cite{1BBBK} and the corresponding derived $I$-complete objects are really more relevant in the arithmetic geometry in our mind such as the corresponding objects in \cite{1BMS2}, \cite{1BS}, \cite{1BS2}, \cite{1Lurie3}, \cite{1NS}, \cite{1Po} and \cite{1Ye}\footnote{These already include many different types of $I$-adic completion in the derived sense, and the corresponding localization and completion in the algebraic topology and $\infty$-categorical theory. We would believe that the corresponding noncommutative consideration as these will be robust enough to deal with the main problems in Iwasawa theory and noncommutative derived analytic geometry.}. In fact on the $\infty$-categorical level, \cite{1BMS2} and \cite{1NS} considered the $I$-completion in the $\infty$-category theory as in \cite[Chapter 7.3]{1Lurie3}. The corresponding derived Nygaard-complete, derived Hodge-complete or other derived filtred complete objects are relevant while the corresponding derived Nygaard-incomplete, derived Hodge-incomplete or other derived filtred incomplete objects may be also relevant in certain situations.;\\
\noindent 16. EC, TMF, TAF, SAV, Sp-DG: \cite{BL1}, \cite{Lan2}, \cite{McCleary1}, \cite{DFHH1}, \cite{Lurie2}, \cite{Lurie3}, \cite{Lurie4}. This is the part for elliptic cohomology, topological modular forms, the topological automorphic forms, spectral abelian varieties, spectral $p$-divisible groups. We treat this part as a tiny degeneralization of the more general theory such as the chromatic homotopy, with restriction to derived abelian varieties and derived modular varieties. We highly recommend the reader to read Lurie's materials on elliptic cohomologies;\\
\noindent 17. Motives, Les Theories de Cohomologie d'apr\`es Weil, La Theorie de Cohomologie d'apr\`es Weil, Standard Conjectures: \cite{Motive1}-\cite{Motive33};\\
\noindent 18. Homotopty, Model Categories: \cite{Model1}, \cite{Model2}, \cite{Model3}, \cite{Model4}, \cite{Model5}, \cite{Model6}, \cite{Model7}, \cite{Model8}, \cite{Model9}, \cite{Model10}, \cite{Model11}, \cite{Model12}, \cite{Model13}, \cite{Model14}, \cite{Model15}, \cite{Model16}, \cite{Model17}, \cite{Model18}, \cite{Model19}, \cite{Model20}, \cite{Model21};\\
\noindent 19. K, THH, TCH, TP, TAQ, ArK: \cite{KTheory1}-\cite{KTheory17};\\  
\noindent 20. Motivic Homotopy, $\mathbb{A}^1$-Homotopy and $\mathbb{B}^1$-Homotopy: \cite{MotivicHomotopy1}, \cite{MotivicHomotopy2}, \cite{MotivicHomotopy3}, \cite{MotivicHomotopy4}, \cite{MotivicHomotopy5}, \cite{MotivicHomotopy6}, \cite{MotivicHomotopy7}, \cite{MotivicHomotopy8}, \cite{MotivicHomotopy9}, \cite{MotivicHomotopy10}, \cite{MotivicHomotopy11}.\\

\newpage
\subsection{Scholze's $v$-Spaces and Six Formalism}

Scholze's $v$-stacks happen over the category of perfectoid space over $\mathbb{F}_p$. We give the introduction following closely \cite{1Sch3} and \cite{1SW}. Our presentation is also following closely \cite{1Sch3}.

\begin{setting}
$v$-sheaves and $v$-stacks carry the topology which is called the $v$-topology which is finer than analytic topology, \'etale topology, pro-\'etale topology.
\end{setting}

\begin{example}
The two important situations are the following. First is the situation where the $v$-stack carries a basis of topological neighbourhoods consisting of perfectoids. The second situation is the key moduli of vector bundles over FF curves in \cite{1FS}.
\end{example}

\begin{definition}\mbox{\textbf{(Scholze \cite{1Sch3}, Analytic Prestacks)}}
Consider the following Grothendieck sites:
\begin{align}
\mathrm{Perfectoid}_{\mathbb{F}_p,\text{\'etale}},\mathrm{Perfectoid}_{\mathbb{F}_p,\text{pro\'etale}}, \mathrm{Perfectoid}_{\mathbb{F}_p,v}.
\end{align}
We just define a $(2,1)$-presheaf $\mathcal{F}$ over these sites to be any functor from 
\begin{align}
\mathrm{Perfectoid}_{\mathbb{F}_p,\text{\'etale}},\mathrm{Perfectoid}_{\mathbb{F}_p,\text{pro\'etale}}, \mathrm{Perfectoid}_{\mathbb{F}_p,v}.
\end{align}
to the groupoids. 

\end{definition}

\begin{definition}\mbox{\textbf{(Scholze \cite{1Sch3}, Analytic Stacks)}}
Consider the following Grothendieck sites:
\begin{align}
\mathrm{Perfectoid}_{\mathbb{F}_p,\text{\'etale}},\mathrm{Perfectoid}_{\mathbb{F}_p,\text{pro\'etale}}, \mathrm{Perfectoid}_{\mathbb{F}_p,v}.
\end{align}
We just define a $(2,1)$-sheaf $\mathcal{F}$ over these sites to be any functor from 
\begin{align}
\mathrm{Perfectoid}_{\mathbb{F}_p,\text{\'etale}},\mathrm{Perfectoid}_{\mathbb{F}_p,\text{pro\'etale}}, \mathrm{Perfectoid}_{\mathbb{F}_p,v}.
\end{align}
to the groupoids, which is further a stack in groupoids.

\end{definition}

\

\indent We then have the morphisms and sites of analytic stacks in certain situations.

\
\begin{definition}\mbox{\textbf{(Scholze \cite[Definition 1.20]{1Sch3}, Morphisms of Analytic Stacks)}}
The \'etale and quasi-pro-\'etale morphisms between small $v$-stacks are defined by using perfectoid coverings, and defining the corresponding \'etaleness and quasi-pro-\'etaleness after taking base changes along such perfectoid coverings.

\end{definition}

\begin{definition}\mbox{\textbf{(Scholze \cite[Definition 26.1]{1Sch3}, Sites of Analytic Stacks)}}
Let $X$ be a small $v$-stack, we have the $v$-site $X_v$\footnote{In order to study the cohomology of any small $v$-stacks, the approach taken by Scholze in \cite{1Sch3} and \cite{1FS} is not defining the particular \'etale or quasipro\'etale derived categories, instead the definition is rather not straightforward by looking at the desired subcategory of $v$-derived categories. The solidified derived category $D_\blacksquare$ in \cite{1FS} is also constructed in the same way.}. This then gives the $\infty$-category $D_{I}(X_v,\lambda)$ of derived $I$-complete objects in $\infty$-category of all $\lambda$-sheaves $D(X_v,\lambda)$ where $\lambda$ is a derived $I$-commutative algebra object.

\end{definition}

\

\begin{theorem} \mbox{\textbf{(Scholze \cite[Definition 26.1, Before the Remark 26.3]{1Sch3})}}
This $\infty$-category $D_{I}(X_v,\lambda)$ as well as the corresponding sub $\infty$-categories $D_{I,\text{qspro\'et}}(X_v,\lambda)$ and $D_{I,\text{\'et}}(X_v,\lambda)$ of the corresponding derived $I$-complete objects over the quasi-pro\'etale and \'etale sites of small $v$-stacks admit six formalism through: $\otimes, \mathrm{Hom}, f^!,f_!,f^*,f_*$ as in \cite[Definition 26.1, Before the Remark 26.3]{1Sch3}.\	
\end{theorem}

\indent Then as in \cite{1FS} for the corresponding solid $\infty$-category we can say the parallel things:

\begin{theorem} \mbox{\textbf{(Fargues-Scholze \cite[Chapter VII.2]{1FS})}}
This $\infty$-category $D_{\blacksquare}(X_v,\lambda)$ admits four formalism as in \cite[Chapter VII.2]{1FS}. Here the assumption on $\lambda$ will be definitely weaker solid condensed ring as assumed in \cite[Beginning of Chapter VII.2]{1FS}.	
\end{theorem}

\begin{remark}
In our situation, certainly what we have will be some derived $I$-completed Iwasawa modules instead of the usual Iwasawa Theoretic $I$-adic versions. This is definitely weaker, but things become robust and more well-defined. For instance if $\lambda$ is usual Iwasawa algebra, we have the following.	
\end{remark}

\begin{theorem} \mbox{\textbf{(Scholze \cite[Definition 26.1, Before the Remark 26.3]{1Sch3})}}
This $\infty$-category $D_{I}(X_v,\lambda)$ as well as the corresponding sub $\infty$-categories $D_{I,\text{qspro\'et}}(X_v,\lambda)$ and $D_{I,\text{\'et}}(X_v,\lambda)$ of the corresponding derived $I$-complete objects over the quasi-pro\'etale and \'etale sites of small $v$-stacks admit six formalism through: $\otimes, \mathrm{Hom}, f^!,f_!,f^*,f_*$ as in \cite[Definition 26.1, Before the Remark 26.3]{1Sch3}. Here we assume $\lambda$ be an Iwasawa algebra attached to some $\ell$-adic Lie group.	
\end{theorem}

\indent Then as in \cite{1FS} for the corresponding solid $\infty$-category we can say the parallel things:

\begin{theorem} \mbox{\textbf{(Fargues-Scholze \cite[Chapter VII.2]{1FS})}}
This $\infty$-category $D_{\blacksquare}(X_v,\lambda)$ admits four formalism as in \cite[Chapter VII.2]{1FS}. Here we assume $\lambda$ be an Iwasawa algebra attached to some $\ell$-adic Lie group.	
\end{theorem}

\newpage

\section{Introduction}

\subsection{The Motivation}

\indent Noncommutative Tamagawa number conjectures and noncommutative Iwasawa main conjectures are main topics in the field of arithmetic geometry. They give us the explicit approaches to study the motives after we take the realizations\footnote{One might want to consider motives over $\mathbb{Z}$ after some reasonable $p$-adic cohomology theories are established such as in \cite{1Ked3} and \cite{1Sch5}. Doing so might be relevant in some $p$-adic cohomological approach to Riemann Hypothesis over $\mathbb{Z}$ after  \cite{1Ked2} which is parallel to \cite{1De1}, \cite{1De2}, \cite{1Wei}.}. The associated Galois representations for instance could be used to define the corresponding $L$-functions, as well as local factors. Roughly speaking, the conjectures mention the following style  isomorphism up to certain factor\footnote{$\pi_1(*)$ is the profinite fundamental group of certain arithmetic scheme or analytic space.}:
\begin{align}
\mathrm{Determinant}_A(R{\Gamma}(\pi_1(*),p\mathrm{LieDef}_A(V)))\overset{}{\longrightarrow} \mathrm{Determinant}_A(0)
\end{align}
when taking some possibly trivial Iwasawa deformation along some $p$-adic Lie group tower and one could achieve the zero homotopy from the determinant of the cohomology in certain categories\footnote{As those categories in \cite{1De3}, \cite{1FK}, \cite{1MT}, \cite{1Wit3}.}, which should be closely related to special values of $L$-function and generalized functional equation.\\

\indent Our study was significantly inspired by the pictures of \cite{1BF1}, \cite{1BF2}, \cite{1FK}, \cite{1Na1}, \cite{1Wit1}, \cite{1Wit2}, \cite{1Wit3}\footnote{The commutative picture was already in \cite{1Ka1}, \cite{1Ka2}, \cite{1PR1}, \cite{1PR2}.}. On the other hand, work of Kedlaya-Liu \cite{1KL1}, \cite{1KL2} and work of Kedlaya-Pottharst \cite{1KP} also significantly inspired us to consider at least the geometrization or generalization of the work of \cite{1BF1}, \cite{1BF2}, \cite{1FK}, \cite{1Na1}, \cite{1Wit1}, \cite{1Wit2}, \cite{1Wit3}. In fact we have other important $p$-adic cohomology theories, especially in the integral setting we have the prismatic cohomology from Bhatt-Scholze \cite{1BS}. Therefore one might want to ask if we could have a well-posed prismatic Iwasawa deformation theory. To be more precise, whenever we have a well-defined $p$-adic cohomology theory, we should be able to consider some interesting Iwasawa deformation theory, possibly carrying some Banach coefficients as in \cite{1Na1}. It is actually very straightforward to first consider the integral picture, and it is very confusing to actually first understand the picture given certain Fr\'echet-Stein algebra in \cite{1ST} since the classical Iwasawa main conjecture happens over integral Iwasawa algebra $\Lambda$ as in \cite{1Iw}. \\

\subsection{The Results}

\indent In fact our study is just as described in the consideration above, namely we actually study some non-\'etale objects which will happen over certain period sheaves over analytic spaces carrying Frobenius actions. And then we can take the corresponding Banach deformation as well as the Iwasawa deformation. Note that the previous deformation will be very crucial such as in \cite{1KPX} and \cite{1Na1}. Meanwhile we can regard them as some certain sheaves over Fargues-Fontaine stacks in some deformed way after \cite{1FF},  \cite{1KL1} and \cite{1KL2}. And we need to work with noncommutative rings as in \cite{1BF1}, \cite{1BF2}, \cite{1FK}, \cite{1Wit1}, \cite{1Wit2}, \cite{1Wit3}, \cite{1Z}. These are reflected in the following:\\

\noindent 1. Hodge-Iwasawa Deformation: \cref{2corollary7.6}, \cref{2proposition2.3.10}, 
\cref{3propsition3.4.25}, \cref{3propsition3.4.26} and \cref{3proposition3.5.51}.\\ 
\noindent 2. Multidimensional Frobenius Modules: \cref{4theorem4.1.5} and \cref{4theorem4.1.6}.\\
\noindent 3. Hodge-Iwasawa Cohomology: \cref{proposition6.3}, \cref{6proposition6.4.26}, 
\cref{6proposition6.4.27}, \cref{6proposition6.4.28}, \cref{proposition4.42}, \cref{6proposition6.4.45}, \cref{proposition4.45}, \cref{6proposition6.4.51}. \\
\noindent 4. Noncommutative Hodge-Iwasawa Deformation: \cref{theorem3.10},  \cref{theorem3.11}, \cref{corollary3.12},  \cref{7theorem7.3.18}, \cref{7theorem7.3.21}, \cref{7theorem4.12}, \cref{proposition5.13}, \cref{8theorem8.1.1}, \cref{8theorem8.1.2}, \cref{8theorem8.1.3}, \cref{8theorem8.1.4}, \cref{8theorem8.1.5}, \cref{8theorem8.1.6}, \cref{8theorem8.1.7}, \cref{9theorem9.1.1}, \cref{9theorem9.1.2}.\\
\noindent 5. Derived Noncommutative Hodge-Iwasawa Deformation: \cref{10theorem10.1.1}, \cref{10theorem10.1.2}, \cref{10theorem10.1.3}, \cref{10theorem10.1.4}, \cref{10theorem10.5.10}.\\

\begin{remark}
Our goal is to find certain derived Iwasawa $\infty$-categories as in \cite{1De3}, \cite{1FK}, \cite{1MT}, \cite{1Wit3} in analytic geometry over certain deformed sheaves or rings. Note that this is very complicated and particularly far beyond Iwasawa deformation of regular motives. We believe Clausen-Scholze derived category $D(\mathrm{Mod}_{\Pi,\mathrm{condensed}})$ of condensed modules over centain condensed ring $\Pi$ in \cite{1CS1} could be more robust consideration. Here in our mind $\Pi$ should be some deformed period ring such as deformation $\widetilde{\Pi}^I_{R}\otimes^{\mathbb{L}\mathrm{Solidified}}\mathcal{A}$ of the Robba ring in \cite{1KL1} and \cite{1KL2}. Also the derived category $D_\mathrm{Solidified}(X)\subset D(X_v)$ of Fargues-Scholze in \cite{1FS} carrying condensed coefficients\footnote{As well as those $D_{\text{\'et}}(X)\subset D(X_v)$ and $D_{\text{qpro\'et}}(X)\subset D(X_v)$ in \cite{1Sch3}, with application in mind to seminormal rigid analytic spaces being regarded as small $v$-stacks.} and the work \cite{1BBK} could be more robust consideration as well. What we have achieved is literally some generalization of Kedlaya-Liu abelian categories of pseudocoherent sheaves in \cite{1KL2}. 
\end{remark}

\

\indent We now discuss some examples of the spaces and local charts in our study. As mention above, they at least will be some spaces over or attached to some period rings or sheaves, which sometimes are the corresponding local charts of some stacks such as the Fargues-Fontaine spaces or the spaces before taking the quotients by equivariance coming from the motivic structures.\\

%%%%%%%%%%%%%%%%%--------------------

\begin{example}\mbox{\textbf{($\infty$-Categorical Rings and $\infty$-Categorical Spaces)}}\\
\noindent 1. The rigid analytic affinoids and spaces in \cite{1Ta} are some key examples, which is actually some initial goal in our main consideration mentioned above inspired by \cite{1BF1}, \cite{1BF2}, \cite{1FK}, \cite{1Wit1}, \cite{1Wit2}, \cite{1Wit3}, \cite{1Z};\\
\noindent 2. The pseudorigid analytic affinoids and spaces over $\mathbb{Z}_p$ are also interesting to study, which is a geometric version of the arithmetic family in \cite{1Bel1} and \cite{1Bel2};\\
\noindent 3. We will consider the noncommutative $I$-adically complete rings as in \cite{1FK} in the Iwasawa consideration, as well as simplicial commutative rings which are derived $I$-adically complete rings over any interesting period rings, such as the prisms in \cite{1BS}, for instance the $\mathbb{A}_\mathrm{inf}(\mathcal{O}_{\mathbb{C}_p^\flat}), W(\mathbb{F}_p)[[u]], \mathbb{Z}_p[[q-1]]$ and Robba rings in \cite{1KL1} and \cite{1KL2}. One can consider for instance the topological rings over these period rings carrying the topology induced from the period rings. Tate adic Banach rings in \cite{1KL1} and \cite{1KL2} produce certain topological adic rings satisfying the open mapping property as in \cite{1CBCKSW}.\\

\end{example}

\begin{setting}\mbox{\textbf{($\infty$-Categorical Completion and $\infty$-Categorical Solidification)}} In fact one will have the chance to take many different types of completions, for instance the Hodge-completion or parallely the Nygaard-completion is happening in some filtered derived $\infty$-category which is very significant in the development of derived prismatic cohomology \cite{1BS}. Also we have the following foundations. First we have from \cite[Proposition 1.6]{1BBBK} the Banach completion on the Banach sets over some Banach ring $R$ which is the left adjoint functor of the inclusion:
\begin{align}
\mathrm{BanSets}_R\rightarrow \mathrm{NormSets}_R.
\end{align}
One may also consider the corresponding functor on the derived level between simplicial ind-Banach and ind-Normed sets:
\begin{align}
\mathrm{Simp}\mathrm{Ind}\mathrm{BanSets}_R\rightarrow \mathrm{Simp}\mathrm{Ind}\mathrm{NormSets}_R,
\end{align}
which may further be extended to the simplicial ind-Banach and ind-Normed abelian groups, rings and modules. And we have on the derived level the corresponding derived $I$-completion over some $\mathbb{E}_\infty$ ring $R$ where $I\subset \pi_0(R)$. When $I$ is finitely generated, if this happens in the derived $\infty$-category $\mathbb{D}_R$ as in \cite[Chapter 1 Notation]{1BS} and \cite[Chapter 3]{1BS2}, then this for instance could be regarded as the left adjoint functor of the inclusion:
\begin{align}
h\mathbb{D}_{R,I-\mathrm{com}}\rightarrow h\mathbb{D}_R,
\end{align}
when $R$ is classical as in \cite[Tag 091N, Section 15.90]{1SP}. Moreover we have the derived $I$-completion on the $\infty$-categorical level to get derived $I$-complete objects:
\begin{align}
\mathbb{D}_R\rightarrow  \mathbb{D}_{R,I-\mathrm{com}}. 
\end{align}
Also we have from \cite[Theorem 5.8]{1CS1} the derived solidification functor $(.)^{\mathbb{L}\mathrm{Solidified}}$ on the derived level which is the left adjoint of the following inclusion from derived category of solids to the derived category of condensed abelian groups:
\begin{align}
D_{\mathrm{condensed},\mathrm{ab},\mathrm{solid}}\rightarrow D_{\mathrm{condensed},\mathrm{ab}}.
\end{align}
And furthermore one can consider the corresponding functors on the derived $\infty$-categorical level (also see \cite[Chapter 1]{1Lurie2}):
\begin{align}
\mathbb{D}_{\mathrm{condensed},\mathrm{ab},\mathrm{solid}}\rightarrow \mathbb{D}_{\mathrm{condensed},\mathrm{ab}}.
\end{align}
Solids are abstractly characterized in \cite{1CS1} as those condensed abelian groups closed under limits, colimits and extensions. Therefore solidification is abstract characterization of completion. Solids may be very important in Iwasawa theory when we want to take solidified tensor products under some absolute topology. For instance we may take the product $\mathcal{F}\otimes^{\mathbb{L}\mathrm{Solidified}}P[1/I]^\wedge_p$ of some Fr\'echet algebra $\mathcal{F}$ such as in \cite{1ST} with certain ring $P[1/I]^\wedge_p$ (or any similar one) produced from some prism $(P,I)$ in \cite{1BS} under this framework. At least for instance in deformation theory happening over the Robba ring $\widetilde{\Pi}^I_{R}\otimes^{\mathbb{L}\mathrm{Solidified}}\mathcal{F}$ after Kedlaya-Liu \cite{1KL1} and \cite{1KL2}, sheafiness has been not a problem anymore after \cite{1BK} and \cite{1CS2}, which implies that one should look at ind-Banach sets instead of Banach sets and pro-\'etale solids instead of complete sets.
\end{setting}

\

\begin{remark}
Both the foundations in Bambozzi-Ben-Bassat-Kremnizer and Clausen-Scholze \cite{1BBBK}, \cite{1CS1} and \cite{1CS2} have noncommutative categories of noncommutative associative analytic rings and noncommutative associative Banach rings. Therefore this may allow one to study certain $\infty$-stacks in $\infty$-groupoids fibered over these categories under some descent consideration for instance after \cite{1KR1} and \cite{1KR2}. 
\end{remark}

\

\indent In chapter 2 and 3, we study the Frobenius modules over Robba rings carrying rigid analytic coefficients and Fr\'echet-Stein coefficients, in both equal-characteristic situation and mixed-characteristic situation. We call the theory Hodge-Iwasawa since the study of the Frobenius modules over the Robba rings and sheaves are significant in Galois deformation theory, deformation of representations of fundamental groups of analytic spaces and our generalizations of the picture in \cite{1BF1}, \cite{1BF2}, \cite{1FK}, \cite{1KP}, \cite{1Na1}, \cite{1Wit1}, \cite{1Wit2}, \cite{1Wit3}. We show the equivalence between categories of finite projective or pseudocoherent $(\varphi,\Gamma)$-modules over Robba rings with rigid analytic coefficients and Fr\'echet-Stein coefficients, which also are compared to sheaves over schematic and adic Fargues-Fontaine curves in some deformed sense. Especially when we are working over analytic fields, the picture is already interesting and significant enough in the Galois representation theory and Galois deformation theory.\\

\indent In chapter 4 and 5, we study multidimensional Robba rings and multidimensional $(\varphi,\Gamma)$-modules. This point of view of taking multidimensional analogs of the Robba rings and multidimensional $(\varphi,\Gamma)$-modules is actually motivated from some programs in making progress of the local Tamagawa number conjecture of Nakamura in \cite{1Na1} literally proposed in \cite{1PZ}. Following \cite{1CKZ} and \cite{1PZ}, we define the corresponding multidimensional Robba rings and multidimensional $(\varphi,\Gamma)$-modules by taking analytic function rings over $p$-adic rigid affinoids in rigid geometry. And we define the multidimensional $(\varphi,\Gamma)$-cohomologies, multidimensional $(\psi,\Gamma)$-cohomologies and multidimensional $\psi$-cohomologies. We carefully study the complexes of the multidimensional $(\varphi,\Gamma)$-cohomologies, multidimensional $(\psi,\Gamma)$-cohomologies and multidimensional $\psi$-cohomologies, and show that they are actually living in the derived category of the bounded perfect complexes. Chapter 4 mainly focuses on imperfect Robba rings, while in chapter 5 we define perfection of Robba rings in several variables and study the comparison of multidimensional $(\varphi,\Gamma)$-modules in certain situations carefully, which is literally following \cite{1KL1}, \cite{1KL2}, \cite{1KPX} and \cite{1Ked1}. \\

\indent In chapter 6, we apply the main results in chapter 3 to study the cohomologies and categories of relative $(\varphi,\Gamma)$-modules over Robba sheaves over certain analytic spaces. We mainly discuss three applications which are crucial in our project of Iwasawa deformation of motivic structures over some higher dimensional spaces. The first is the study of abelian property of the categories of relative $(\varphi,\Gamma)$-modules over Robba sheaves over rigid analytic spaces after \cite{1KL2}, which is important whenever one would like to construct some $K$-theoretic objects to formulate Iwasawa main conjectures as in \cite{1Wit1}, \cite{1Wit2}, \cite{1Wit3}. The second is the study of families of Riemann-Hilbert correspondences after \cite{1LZ} which is crucial in further application to the arithmetic geometry along our consideration. The last one is the consideration of the equivariant version of the Iwasawa main conjecture of Nakamura \cite{1Na2}.\\

\indent In chapter 7, 8 and 9, we consider generalization of our work of the generalization of the work of Kedlaya-Liu presented in chapter 2, 3 and 6, which is also motivated by our consideration in  chapter 4 and 5 after \cite{1CKZ}, \cite{1Na1}, \cite{1PZ} and \cite{1Z} in order to make further progress. We consider the deformation in possibly noncommutative general Banach rings such as perfectoid rings, preperfectoid rings, general quotient of the noncommutative free Tate rings and so on. Certainly one thing we have to deal with is the sheafiness of the deformed rings, which will produce some difficulty to apply Kedlaya-Liu's descent in \cite{1KL1}, \cite{1KL2} for vector bundles and stably-pseudocoherent sheaves. Even in the noncommutative coefficient situation we have not worked out a theory on the noncommutative analytic toposes, which implies there is no geometric method for us to study and apply. So finding new ideas is very important. In fact, we on the representation theoretic level have the result due to Kedlaya to have the descent for vector bundles. And one can in the commutative situation use Clausen-Scholze space \cite{1CS2} to achieve the similar result by embedding Huber spaces to Clausen-Scholze spaces in \cite{1CS2} and apply Clausen-Scholze descent. Also we could consider the derived analytic spaces from Bambozzi-Kremnizer in \cite{1BK}. In the noncommutative situation we generalize results in \cite{1KL1}, \cite{1KL2} and \cite{1CBCKSW} to deform the structure sheaves directly in analytic topology, \'etale topology and pro-\'etale but not $v$-topology, which allows us to compare certain stably-pseudocoherent sheaves and modules carrying Banach deformed coefficients even if they are noncommutative, which certainly provides possibility to make further progress in the study of the noncommutative situation in chapter 4 and 5. \\

\indent In chapter 10, we initiate the project on some topics on the geometric and representation theoretic aspects of period rings. In this first paper, we consider more general base spaces. To be more precise we discuss more general perfectoid rings. Distinguished deformation of rings is a generalization notion of the Fontaine-Wintenberger idempotent correspondence. For instance in \cite{1BS}, for any quasiregular semiperfectoid ring $A$ one can canonically associate a prism $(P_A,I_{A})$. This is a very general correspondence generalized from the notions in \cite{1BMS1} and \cite{1GR}. Therefore in our consideration we will consider the adic Banach rings in \cite{1KL1}, \cite{1KL2} which are not necessarily analytic in the sense of Kedlaya's AWS lecture notes in \cite{1CBCKSW}. Assuming certain perfectoidness after \cite{1KL1} and \cite{1KL2} we study the deformed Robba rings associated. Also after \cite{1KL1} and \cite{1KL2} we studied derived deformation of the Robba rings and the descent of finite projective module spectra over them, which we will believe has some application to conjectural derived eigenvarieties and derived Galois deformation for instance in \cite{1GV}. \\

\indent In chapter 11 and 12, we consider more widely the process of taking topological and functional analytic completions, in the derived sense coming from \cite{1B1}, \cite{1BBBK}, \cite{1BK},  \cite{1BMS2}, \cite{1BS}, \cite{1CS1}. Maybe this derived consideration will let us see the hidden $\infty$-categorical structures of the representations of the Iwasawa algebras in our Iwasawa-Prismatic theory, since we are taking the deformations in some derived sense. Besides the application in mind to the prismatic or derived de Rham Iwasawa theory, chapter 11 and chapter 12 initiate also the study of some relative $p$-adic motive and Hodge theory over general derived $I$-adic spaces after \cite{1B1}, \cite{1B2}, \cite{1BMS2}, \cite{1BS}, \cite{1CBCKSW}, \cite{1DLLZ1}, \cite{1DLLZ2}, \cite{1Dr1}, \cite{1GL}, \cite{1G1}, \cite{1Hu2}, \cite{1III1}, \cite{1III2}, \cite{1KL1}, \cite{1KL2}, \cite{1NS}, \cite{1O}, \cite{1Ro}, \cite{1Sch2} such as the pseudorigid analytic spaces and more general spaces carrying some derived $I$-adic topology, as well as the prelog simplicial commutative rings in \cite{1B1} carrying some derived $I$-adic topology as well.

\subsection{The Picture and The Future Consideration}

Let us try to discuss a little bit how one should think about our main motivation. We start from the following two settings of possible geometrizations of Iwasawa theory.\\

\noindent 1. \textbf{$\infty$-Categorical Iwasawa-Prismatic Theory}\footnote{This should be slightly more relevant in the geometrization of integral Iwasawa theory, but one can invert $p$ as well.}: After Bhatt-Lurie, Bhatt-Scholze and Drinfeld \cite{1BS}, \cite{1Dr}, \cite{1Sch4};\\
\noindent 2. \textbf{$\infty$-Categorical Hodge-Iwasawa Theory}\footnote{This should be slightly more relevant in the geometrization of rational Iwasawa theory, but one can not invert $p$ as well.}: After Kedlaya-Liu and Kedlaya-Pottharst \cite{1KL1}, \cite{1KL2}, \cite{1KP}.

\

\begin{setting} \textbf{($\infty$-Categorical Iwasawa-Prismatic Theory)}
Just like in \cite{1KP}, \cite{1Wit1}, \cite{1Wit2}, \cite{1Wit3}, and by using prismatic cohomology theory in \cite{1BS},   \cite{1Sch4} we can take a quasisyntomic formal ring $R$\footnote{As in \cite{1Sch4}, one can take a $p$-adic fornal scheme which is quasisyntomic. Certainly this is already a type of spaces which is interesting enough including a point situation and many situations where motivic comparisons could happen as in \cite{1BMS2}, \cite{1BS}.} over $\mathbb{Z}_p$ and have two sites of $\mathrm{Spf}R$. The first site is the corresponding prismatic site $(X_{\mathrm{prim}},\mathcal{O}_{X_{\mathrm{prim}}})$. The second site is the quasisyntomic site $(X_\mathrm{qsyn},\mathcal{O}_{X_{\mathrm{qsyn}}})$. Recall for the second one, for any quasiregular semiperfectoid affinoid $A$ in $X_\mathrm{qsyn}$ we have that one can canonically associate a prism $(P_A,I_A)$ to $A$ where we have $\mathcal{O}_{X_{\mathrm{qsyn}}}(A):=P_A$. Now by Bhatt-Scholze \cite{1Sch4} we have more sheaves from the structure sheaves here namely we have:
\begin{align}
\mathcal{O}_{X_{\mathrm{prim}}}[1/I_{\mathcal{O}_{X_{\mathrm{prim}}}}]^\wedge_p, \mathcal{O}_{X_{\mathrm{prim}}}[1/I_{\mathcal{O}_{X_{\mathrm{prim}}}}]^\wedge_p[1/p], \mathcal{O}_{X_{\mathrm{qsyn}}}[1/I_{\mathcal{O}_{X_{\mathrm{qsyn}}}}]^\wedge_p, \mathcal{O}_{X_{\mathrm{qsyn}}}[1/I_{\mathcal{O}_{X_{\mathrm{qsyn}}}}]^\wedge_p[1/p]. 
\end{align}
Carrying some integral Iwasawa algebra $\mathbb{Z}_p[[G]]$ for some compact $p$-adic Lie group, and after taking the derived completion\footnote{At this moment we are assuming that the derived completion is possible in our current setting. One may also consider the solidification of Clausen-Scholze. We want to mention that in the noncommutative setting there are many ways to do the completion in the derived sense, which is already subtle in the commutative setting.} we have:
\begin{align}
\mathcal{O}_{X_{\mathrm{prim}}}[1/I_{\mathcal{O}_{X_{\mathrm{prim}}}}]^\wedge_p{\otimes}^{\mathbb{L}\mathrm{com}}\mathbb{Z}_p[[G]],\\
 \mathcal{O}_{X_{\mathrm{prim}}}[1/I_{\mathcal{O}_{X_{\mathrm{prim}}}}]^\wedge_p{\otimes}^{\mathbb{L}\mathrm{com}}\mathbb{Z}_p[[G]][1/p],\\
  \mathcal{O}_{X_{\mathrm{qsyn}}}[1/I_{\mathcal{O}_{X_{\mathrm{qsyn}}}}]^\wedge_p{\otimes}^{\mathbb{L}\mathrm{com}}\mathbb{Z}_p[[G]], \\
  \mathcal{O}_{X_{\mathrm{qsyn}}}[1/I_{\mathcal{O}_{X_{\mathrm{qsyn}}}}]^\wedge_p{\otimes}^{\mathbb{L}\mathrm{com}}\mathbb{Z}_p[[G]][1/p]. 
\end{align}
Then one might want to ask if one can use such style deformation to establish the parallel story in \cite{1BF1}, \cite{1BF2}, \cite{1FK}, \cite{1KP}, \cite{1Wit1}, \cite{1Wit2}, \cite{1Wit3}. For instance taking the $J\subset \mathbb{Z}_p[[G]]$-adic quotient we have the Koszul complexes parametrized by such $J$:
\begin{align}
\mathrm{Kos}_J\mathcal{O}_{X_{\mathrm{prim}}}[1/I_{\mathcal{O}_{X_{\mathrm{prim}}}}]^\wedge_p{\otimes}\mathbb{Z}_p[[G]],\\
\mathrm{Kos}_J  \mathcal{O}_{X_{\mathrm{qsyn}}}[1/I_{\mathcal{O}_{X_{\mathrm{qsyn}}}}]^\wedge_p{\otimes}\mathbb{Z}_p[[G]]. \\ 
\end{align}
Then one may define the $\infty$-category of pro-systems of the Iwasawa complexes over these $\mathbb{E}_1$-rings, and consider the associated Waldhausen categories as in \cite{1Wit1}, \cite{1Wit2}, \cite{1Wit3}.

\end{setting}

\

\begin{setting} \textbf{($\infty$-Categorical Iwasawa-Prismatic Theory)}
Within the same framework we take $R$ to be $\mathcal{O}_K$ for some $p$-adic local field $K$. Bhatt-Scholze \cite{1BS}, \cite{1Sch4} showed that we have the category of Galois representations of $\mathbb{Z}_p$-coefficients of $\mathrm{Gal}_K$ is equivalent to the category of prismatic $F$-crystals over 
\begin{center}
$(X_\mathrm{prim},\mathcal{O}_{X_{\mathrm{prim}}}[1/I_{\mathcal{O}_{X_{\mathrm{prim}}}}]^\wedge_p)$, 
\end{center}
while the category of Galois representations of $\mathbb{Q}_p$-coefficients of $\mathrm{Gal}_K$ is equivalent to the category of prismatic $F$-crystals over $(X_\mathrm{prim},\mathcal{O}_{X_{\mathrm{prim}}}[1/I_{\mathcal{O}_{X_{\mathrm{prim}}}}]^\wedge_p[1/p])$. Then one could ask if we could consider some Iwasawa deformation through the some $p$-adic Lie quotient of $\mathrm{Gal}_K$ to establish the parallel story in \cite{1KP}, \cite{1Wit1}, \cite{1Wit2}, \cite{1Wit3}. Namely for any such $F$-crystal $M$ with associated representation $V$ over $\mathbb{Z}_p$ or $\mathbb{Q}_p$ we take the Iwasawa deformation $p\mathrm{DfLie}(M)$ by some Iwasawa $F$-crystal\footnote{How one should define this crystal will be determined by how one forms the completed tensor product in the Iwasawa deformation.} through the quotient from $G_K$ to some compact $p$-adic Lie group $G$, then we can ask if the following:
\begin{align}
R\Gamma(X_\mathrm{prim},p\mathrm{DfLie}(M)), R\Gamma(X_\mathrm{qsyn},p\mathrm{DfLie}(M))
\end{align}
recover the classical Iwasawa theory by using the Galois cohomology 
\begin{center}
$R\Gamma(\mathrm{Gal}_K,p\mathrm{DfLie}(V))$
\end{center}
of $p\mathrm{DfLie}(V)$, as well as the \'etale cohomology $R\Gamma(\mathrm{Spec}K,p\mathrm{DfLie}(\widetilde{V}))$ of the local system $p\mathrm{DfLie}(\widetilde{V})$ attached to $p\mathrm{DfLie}(V)$.
\end{setting}

\

\indent Beyond the somehow \'etale situations in the above picture, one could consider the corresponding category of prismatic crystals, which will be beyond the Galois representation theoretic consideration. Also one could regard these objects as certain sheaves over the prismatic stacks in \cite{1Dr}.\\ 

\indent In our study, we have the following picture. We will consider picture beyond \'etale situation, and we will study the Frobenius sheaves and Frobenius modules in the very general situation. And we will have the chance to regard the sheaves and modules with Frobenius actions as certain sheaves over Fargues-Fontaine stacks after \cite{1FF}, \cite{1KL1} and \cite{1KL2}. Actually we conjecture that the quasisyntomic descent and étale comparison results of Bhatt-Scholze \cite{1BS} will imply equivalence in some accurate sense beyond the vague similarity.\\

\chapter{$\infty$-Categorical Approaches to Hodge-Iwasawa Theory I}

\section{Motivation I}

\subsection{\text{Motivation I: Dememorization and Memorization}}
\begin{itemize}

\item<1-> Consider the cyclotomic tower $\{\mathbb{Q}_p(\zeta_{p^n})\}_n$ of $\mathbb{Q}_p$.

\item<2-> The infinite level of this tower is kind of special after the corresponding completion.

\item<3-> Over $\mathbb{Q}_p$, we could consider $\mathrm{Spa}(\mathbb{Q}_p,\mathfrak{o}_{\mathbb{Q}_p})_{\text{pro\'et}}$ due to Scholze \cite{Sch}, although the infinite level of the towers above participates in this topology but the corresponding pro-\'etale site forgets the corresponding cyclotomic tower while it is defined by using pro-systems of \'etale morphisms.

\item<4-> Work of Pottharst \cite{P1}, Kedlaya-Pottharst-Xiao \cite{KPX}, Kedlaya-Pottharst \cite{KP} implies one may see the corresponding cyclotomic tower back by considering the corresponding cyclotomic deformation as below.

\item<5-> One has the so-called $\psi$-cohomology originally dated back to Fontaine (see \cite[II.1.3]{CC}) attached to a $(\varphi,\Gamma)$-module $M$ (you could regarded this as a Galois representation):
\begin{align}
H_\psi(M)	
\end{align}
by using the operator $\psi$.

\item<6-> And we have the corresponding $(\varphi,\Gamma)$-module after Herr, but we consider the cyclotomic deformation as in Kedlaya-Pottharst-Xiao over the Robba ring $\mathcal{R}^\infty_{\mathbb{Q}_p}(\Gamma)$:
\begin{align}
H_{\varphi,\Gamma}(\mathbf{CycDef}(M)).	
\end{align}

\end{itemize}

\noindent{\text{Motivation I: Dememorization and Memorization}}
\begin{itemize}
\item<1-> This is defined by taking the corresponding external tensor product of $M$ with the corresponding module coming from the quotient $\Gamma$. This dates back to Pottharst on his analytic Iwasawa cohomology \cite{P1}.

\item<2-> Work of Kedlaya-Pottharst \cite{KP} observes that we can have the following sheaf version of the construction:
\begin{align}
H_\text{pro-\'etale}(\mathrm{Spa}(\mathbb{Q}_p,\mathfrak{o}_{\mathbb{Q}_p}),\mathbf{CycDef}(\widetilde{M})),	
\end{align}
which is defined by taking the corresponding external product of Kedlaya-Liu's sheaf $\widetilde{M}$ \cite{KL1} with the one defined by using the quotient $\Gamma$. 
\item<3-> The point is that we have the following comparison:
\begin{align}
H_{\psi}(M) \overset{\sim}{\rightarrow}	 H_{\varphi,\Gamma}(\mathbf{CycDef}(M))\overset{\sim}{\rightarrow}H_\text{pro-\'etale}(\mathrm{Spa}(\mathbb{Q}_p,\mathfrak{o}_{\mathbb{Q}_p}),\mathbf{CycDef}(\widetilde{M})).
\end{align}

\item<4-> Suppose $M(V)$ comes from a Galois representation $V$ of $G_{\mathbb{Q}_p}$ we even have the following comparison after Perrin-Riou \cite{PR}:
\begin{align}
&{\mathcal{O}_{\mathrm{Sp}\mathcal{R}^\infty_{\mathbb{Q}_p}(\Gamma)}}\widehat{\otimes}_\Lambda H_\mathrm{IW}(G_{\mathbb{Q}_p},V)\overset{\sim}{\rightarrow} \\
&\overset{\sim}{\rightarrow} H_{\psi}((M(V)) \overset{\sim}{\rightarrow}H_{\varphi,\Gamma}(\mathbf{CycDef}((M(V)))\overset{\sim}{\rightarrow}H_\text{pro-\'etale}(\mathrm{Spa}(\mathbb{Q}_p,\mathfrak{o}_{\mathbb{Q}_p}),\mathbf{CycDef}(\widetilde{(M(V)})).
\end{align}
\item<5-> Natural questions come:\\
I. How about the Lubin-Tate Iwasawa theory in Berger-Fourquaux-Schneider-Venjakob's work, observed by Kedlaya-Pottharst \cite{KP}, \cite{BF}, \cite{SV}.\\
II. How about higher dimensional toric towers and more general towers of rigid analytic spaces for instance.

\end{itemize}

\noindent{\text{Motivation I: Dememorization and Memorization}}
\begin{itemize}
\item<1-> These need us to generalize the corresponding framework to higher dimensional situation and more general deformed version. The problem is challenging, since we have some rigidized objects combined together. 
\item<2-> Rational coefficients are very complicated comparing to algebraic geometry, since sometimes we do not have the integral lattices over the \'etale sites. This is already a problem in the context of Kedlaya-Liu \cite{KL1}, \cite{KL2}. 
\item<3->It is not surprising much for us to consider generalizing the frame work of non-\'etale objects since even in the usual situations over a point work of Nakamura \cite{Nakamura1}, Kedlaya-Pottharst-Xiao \cite{KPX} and Kedlaya-Liu \cite{KL1}, \cite{KL2} implies that all kinds of families of Galois representations will be more conveniently studied by using $B$-pairs and $(\varphi,\Gamma)$-modules. 
\item<4-> If we only have some abelian group $G$ the corresponding deformation happens along the algebra $\mathbb{Q}_p[G]$ which gives rise to Galois representation of $\mathrm{Gal}(\overline{\mathbb{Q}}_p/\mathbb{Q}_p)$ with coefficient in $\mathbb{Q}_p[G]$ along the quotient $\mathrm{Gal}(\overline{\mathbb{Q}}_p/\mathbb{Q}_p)\rightarrow G$. One can then regard this as a sheave of module over some sheaf with deformed coefficient in $\mathbb{Q}_p[G]$. 
\item<5-> Note that we can also consider some deformation over an affinoid algebra in the rigid analytic geometry, which amounts to the $p$-adic families of special values. This is not available at once in archimedean functional analysis.	
\end{itemize}

\newpage

%\section{Motivation II}

\section{\text{Motivation II: Higher Dimensional modeling of the
 Weil Conjectures}}
\begin{itemize}
\item<1-> The corresponding equivariant consideration could be obviously generalized to the relative $p$-adic Hodge theory which is aimed at the study of the \'etale local systems over rigid analytic spaces.
\item<2-> This amounts to higher dimensional modeling of the generalized Weil conjecture after the work due to many people, to name a few Deligne \cite{De1}, \cite{De2} ($\ell$-adic \'etale sheaves), Kedlaya \cite{Ked1} ($p$-adic differential modules), Abe and Caro \cite{AC} ($p$-adic arithmetic $D$-modules) and so on.

\item<3-> The invariance comes from the quotient of \'etale fundamental groups of rigid analytic spaces, or the corresponding profinite fundamental groups of rigid analytic spaces.

\item<4-> (\text{Example}) One can consider the corresponding Fr\'echet-Stein algebras associated to the group $\mathbb{Z}_p\ltimes \mathbb{Z}_p^n$ which is Galois group (quotient of the corresponding profinite fundamental group) of a local chart of smooth proper rigid analytic spaces. Note that the top of this local chart in the smooth proper setting naturally participates in some nice topology.

\item<5-> (\text{Example}) One can consider the local systems in more general sense, for instance the locally constant sheaf $\underline{A}$	
attached to a topological ring, for instance an affinoid algebra in the rigid analytic geometry after Tate. This is somewhat special in the $p$-adic setting due to the fact the corresponding Hodge structures could achieve variation in $p$-adic rigid family.

\end{itemize}

\newpage

\section{What can be learnt from noncommutative Iwasawa Theory}

\subsection{\text{Integral \'Etale Noncommutative Iwasawa Theory}}

\begin{itemize}
\item<1-> Following some idea in the noncommutative Tamagawa Number conjecture after Fukaya-Kato \cite{FK} and the noncommutative Iwasawa theory over a scheme over finite field after Witte we would like to consider the following picture after Witte \cite{Wit1}.\\

\item<2-> Let $T$ be an adic ring in the sense of Fukaya-Kato \cite{FK}, which is a compact ring with two sided ideal $I$ such that we have each $T/I^n$ is finite for $n\geq 0$ and taking the inverse limit we recover the ring $T$ itself. This ring could be noncommutative, for instance the Iwasawa algebra attached to some $p$-adic Lie group.\\
\item<3-> (\text{Definition, after Witte \cite[Definition 5.4.1]{Wit1}}) Consider a rigid analytic space or a scheme $X/\mathbb{Q}_p$ separated and of finite type we consider the category $\mathbb{D}_\mathrm{perf}(X_\sharp,T)$ ($\sharp=\text{\'et},\text{pro\'et}$) which is the category of the inverse limit of perfect complexes of abelian sheaves of left modules over quotients of $T$ by open two-sided ideals of $T$ which are $DG$-flat, parametrized by open two-sided ideals of $T$.\\

\item<4-> (\text{Theorem, Witte \cite[Proposition 6.1.5]{Wit1}}) Let $p$ be a unit in $T$. The category defined above could be endowed with the structure of Waldhausen category\footnote{Strictly Speaking, these are the complicial biWaldhausen ones.} and the total direct image functor induces a well defined functor in the situation where $X$ is a scheme and $\sharp=\text{\'et}$:
\[
\xymatrix@R+0pc@C+2pc{
\mathbb{D}_\mathrm{perf}(X_\sharp,T)\ar[r]^{R\Gamma(X_\sharp,.)}\ar[r]\ar[r] &\mathbb{D}_\mathrm{perf}(T) 
}
\]	
which induces the corresponding map on the $K$-theory space:
\[
\xymatrix@R+0pc@C+2pc{
\mathbb{K}\mathbb{D}_\mathrm{perf}(X_\sharp,T)\ar[r]^{\mathbb{K}R\Gamma(X_\sharp,.)}\ar[r]\ar[r] &\mathbb{K}\mathbb{D}_\mathrm{perf}(T). 
}
\]
Then this map is homotopic to zero in some canonical way.

\end{itemize}

\newpage

\section{Introduction to the Interactions among Motives}

\subsection{\text{Equivariant relative $p$-adic Hodge Theory}}

\begin{itemize}
\item<1->Things discussed so far have motivated the corresponding equivariant relative $p$-adic Hodge Theory in the following sense. Witte \cite{Wit1} considered general framework of Grothendieck abelian categories, for instance one can consider the following categories:\\

%\item<2-> 1. The ind-category of all the arithmetic $D$-modules over a realizable scheme of finite type over a perfect field $k$ as those considered by Berthelot, Caro, Abe and etc (this is not covered in our work).\\
\item<1-> 1. The category of all the abelian sheaves over the \'etale or pro-\'etale sites of schemes of finite type over a field $k$ after Grothendieck, Scholze, Bhatt \cite{SGA4}, \cite{BS1} and etc;\\
\item<2-> 2. The category of all the abelian sheaves over the \'etale or pro-\'etale sites of adic spaces of finite type over a field $k$ after Huber, Scholze, Kedlaya-Liu \cite{Hu}, \cite{Sch}, \cite{KL1}, \cite{KL2};\\	
\item<3-> 3. The ind-category of the abelian category of the pseudocoherent Frobenius $\varphi$-sheaves over a rigid analytic space over a complete discrete valued field with perfect residue field $k$ after Kedlaya-Liu \cite{KL1}, \cite{KL2}.\\

\item<4-> 4. The category of abelian sheaves over the syntonic site  by covering of quasiregular semiperfect algebras, as in the work of Bhatt-Morrow-Scholze \cite{BMS}.  \\

\item<5-> One can naturally consider the corresponding $P$-objects throughout the categories listed above, where $P$ is noetherian for instance. For instance one can consider the third category and consider the corresponding local systems over $\underline{A}$ where $A$ is an affinoid algebra in rigid analytic geometry after Tate \cite{Ta1}, which are the $A$-objects in the corresponding category of all the abelian sheaves.

\end{itemize}

\noindent{\text{Equivariant relative $p$-adic Hodge Theory}}

\begin{itemize}

\item<1-> The corresponding $P$-objects are interesting, but in general are not that easy to study, especially we consider for instance those ring defined over $\mathbb{Q}_p$, let it alone if one would like to consider the categories of the complexes of such objects. 

\item<2-> We choose to consider the corresponding embedding of such objects into the categories of Frobenius sheaves with coefficients in $P$ after Kedlaya-Liu \cite{KL1}, \cite{KL2}. Again we expect everything will be more convenient to handle in the category of $(\varphi,\Gamma)$-modules.
	
\item<3-> Working over $R$ now a uniform Banach algebra with further structure of an adic ring over $\mathbb{F}_p$. And we assume that $R$ is perfect.
Let $\text{Robba}^\text{extended}_{I,R}$ be the Robba sheaves defined by Kedlaya-Liu \cite{KL1}, \cite{KL2}, with respect to some interval $I\subset (0,\infty)$, which are Fr\'echet completions of the ring of Witt vector of $R$ with respect to the Gauss norms induced from the norm on $R$. 

\item<4-> Taking suitable interval one can define the corresponding Robba rings $\text{Robba}^\text{extended}_{r,R}$, $\text{Robba}^\text{extended}_{\infty,R}$ and the corresponding full Robba ring $\text{Robba}^\text{extended}_{R}$.

\item<5-> We work in the category of Banach and ind-Fr\'echet spaces, which are commutative. Our generalization comes from those Banach reduced affinoid algebras $A$.

\end{itemize}

\noindent{\text{Equivariant relative $p$-adic Hodge Theory}}
\begin{itemize}

\item<1-> The $p$-adic functional analysis produces us some manageable structures within our study of relative $p$-adic Hodge theory, generalizing the original $p$-adic functional analytic framework of Kedlaya-Liu \cite{KL1}, \cite{KL2}.

\item<2-> Starting from Kedlaya-Liu's period rings, 
\begin{align}
&\text{Robba}^\text{extended}_{\infty,R},\text{Robba}^\text{extended}_{I,R},\text{Robba}^\text{extended}_{r,R},\text{Robba}^\text{extended}_{R},	\text{Robba}^\text{extended}_{{\mathrm{int},r},R},\\
&\text{Robba}^\text{extended}_{{\mathrm{int}},R},\text{Robba}^\text{extended}_{{\mathrm{bd},r},R},\text{Robba}^\text{extended}_{{\mathrm{bd}},R}
\end{align}
we can form the corresponding $A$-relative of the period rings:
\begin{align}
&\text{Robba}^\text{extended}_{\infty,R,A},\text{Robba}^\text{extended}_{I,R,A},\text{Robba}^\text{extended}_{r,R,A},\text{Robba}^\text{extended}_{R,A},	\text{Robba}^\text{extended}_{{\mathrm{int},r},R,A},\\
&\text{Robba}^\text{extended}_{\mathrm{int},R,A},\text{Robba}^\text{extended}_{{\mathrm{bd},r},R,A},\text{Robba}^\text{extended}_{{\mathrm{bd}},R,A}.	
\end{align} 
\item<3-> (\text{Remark}) There should be also many interesting contexts, for instance consider a finitely generated abelian group $G$, one can consider the group rings: 
\begin{align}
\text{Robba}^\text{extended}_{I,R}[G].	
\end{align}
\item<4-> And then consider the completion living inside the corresponding infinite direct sum Banach modules 
\begin{align}
\bigoplus\text{Robba}^\text{extended}_{I,R},	
\end{align}
over the corresponding period rings:
\begin{align}
\overline{\text{Robba}^\text{extended}_{I,R}[G]}.	
\end{align}
Then we take suitable intersection and union one can have possibly some interesting period rings $\overline{\text{Robba}^\text{extended}_{r,R}[G]}$ and $\overline{\text{Robba}^\text{extended}_{R}[G]}$.
\end{itemize}

\noindent{\text{Equivariant relative $p$-adic Hodge Theory}}
\begin{itemize}
\item<1-> The equivariant period rings in the situations we mentioned above carry relative Frobenius action $\varphi$ induced from the Witt vectors. 

\item<2-> They carry the corresponding Banach or (ind-)Fr\'echet spaces structures. So we can generalize the corresponding Kedlaya-Liu's construction to the following situations (here let $G$ be finite):

\item<3-> We can then consider the corresponding completed Frobenius modules over the rings in the equivariant setting. To be more precise over:
\begin{align}
\overline{\text{Robba}^\text{extended}_{R}[G]},\Omega_{\mathrm{int},R,A},\Omega_{R,A},\text{Robba}^\text{extended}_{R,A},\text{Robba}^\text{extended}_{\mathrm{bd},R,A}	
\end{align}
one considers the Frobenius modules finite locally free.

\item<4->  With the corresponding finite locally free models over
\begin{align}
\overline{\text{Robba}^\text{extended}_{r,R}[G]},\text{Robba}^\text{extended}_{r,R,A},\text{Robba}^\text{extended}_{{\mathrm{bd},r},R,A},	
\end{align} 
again carrying the corresponding semilinear Frobenius structures, where $r$ could be $\infty$.

\item<5-> One also consider families of Frobenius modules over 
\begin{align}
\overline{\text{Robba}^\text{extended}_{I,R}[G]},\text{Robba}^\text{extended}_{I,R,A},	
\end{align} 
in glueing fashion with obvious cocycle condition with respect to three intervals $I\subset J\subset K$. These are called the corresponding Frobenius bundles.
	
\end{itemize}

\newpage

\section{The Key Deformation}

\subsection{\text{Deformation of Schemes}}	

\begin{itemize}
\item<1-> One can consider the corresponding schemes attached to the above commutative rings, for instance
\begin{align}
\mathrm{Spec}\text{Robba}^\text{extended}_{r,R,A},\mathrm{Spec}\overline{\text{Robba}^\text{extended}_{r,R}[G]}.	
\end{align}
And consider the corresponding categories:
\begin{align}
\mathrm{Mod}(\mathcal{O}_{\mathrm{Spec}\text{Robba}^\text{extended}_{r,R,A}}),\mathrm{Mod}(\mathcal{O}_{\mathrm{Spec}\overline{\text{Robba}^\text{extended}_{r,R}[G]}}).	
\end{align}
 
\item<2-> These are very straightforward and even crucial especially when we consider 
\begin{align}
\varphi-\mathrm{Mod}(\mathcal{O}_{\mathrm{Spec}\text{Robba}^\text{extended}_{\infty,R,A}}),\varphi-\mathrm{Mod}(\mathcal{O}_{\mathrm{Spec}\overline{\text{Robba}^\text{extended}_{\infty,R}[G]}}),	
\end{align}
in some Frobenius equivariant way.

\item<3-> But on the other hand it is also very convenient to encode the Frobenius action inside the spaces themselves, which leads to Fargues-Fontaine Schemes as those in the work of Kedlaya-Liu \cite{KL1}, \cite{KL2}, \cite{FF}. 

\end{itemize}

\noindent{\text{Deformation of Schemes}}
\begin{itemize}
\item<1-> Roughly one takes the corresponding $\varphi=p^n$ equivariant elements in the full Robba ring, and putting them to be a commutative graded ring $\bigoplus P_{R,A,n}$, and then glueing them through the Proj construction by glueing subschemes taking the form of $\mathrm{Spec}P_{R,A}[1/f]_0$.
\item<2-> Roughly one takes the corresponding $\varphi=p^n$ equivariant elements in the full Robba ring $\overline{\text{Robba}^\text{extended}_{R}[G]}$, and putting them to be a commutative graded ring $\bigoplus P_{R,G,n}$, and then glueing them through the Proj construction by glueing subschemes taking the form of $\mathrm{Spec}P_{R,G}[1/f]_0$.

\item<3-> Therefore we have the natural functor:
\[
\xymatrix@R+0pc@C+0pc{
\mathrm{Mod}\mathcal{O}_{\mathrm{Proj}_{R,A}}\ar[r]\ar[r]\ar[r] &\mathrm{Mod}\mathcal{O}_{\mathrm{Spec}\text{Robba}^\text{extended}_{\infty,R,A}} , 
}
\]
defined by using the corresponding pullbacks.

\item<4-> \text{(Theorem, Tong \cite[Theorem 1.3]{T})} We have the following categories are equivalent (generalizing the work of \text{Kedlaya-Liu} \cite{KL1}, \text{Kedlaya-Pottharst} \cite{KP}): \\
I. The category of all the quasicoherent finite locally free sheaves over $\mathrm{Proj}\bigoplus P_{R,A,n}$;\\
%II. The category of all the $\varphi$-equivariant quasicoherent finite locally free sheaves over $\mathrm{Spec}\text{Robba}^\text{extended}^\infty_{R,A}$;\\
II. The category of all the Frobenius modules of the global sections of all the $\varphi$-equivariant quasicoherent finite locally free sheaves over $\mathrm{Spec}\text{Robba}^\text{extended}_{\infty,R,A}$;\\
III.  The category of all the Frobenius modules over $\text{Robba}^\text{extended}_{R,A}$;\\
IV. The category of all the Frobenius bundles over $\text{Robba}^\text{extended}_{R,A}$.

\item<5-> For the rings for general $G$, we expect one should also be able to establish some results parallel to this once the structures are more literally investigated.	We are also interested in the noncommutative coefficients as in Z\"ahringer's thesis \cite{Z}, but we need to use noncommutative topos.
\end{itemize}

\noindent{\text{Deformation of Schemes}}
\begin{itemize}
\item<1-> \text{(Theorem, Tong \cite[Proposition 3.16, Corollary 3.17]{T})} We have the following categories are equivalent (generalizing the work of \text{Kedlaya-Liu} \cite{KL1}, \text{Kedlaya-Pottharst} \cite{KP}): \\
I. The category of pro-systems of all the quasicoherent finite locally free sheaves over $\mathrm{Proj}\bigoplus P_{R,A_\infty,n}$;\\
%II. The category of all the $\varphi$-equivariant quasicoherent finite locally free sheaves over $\mathrm{Spec}\text{Robba}^\text{extended}^\infty_{R,A_\infty}$;\\
II. The category of pro-systems of all the Frobenius modules coming from the global sections of all the $\varphi$-equivariant quasicoherent finite locally free sheaves over $\mathrm{Spec}\text{Robba}^\text{extended}_{\infty,R,A_\infty}$;\\
III.  The category of pro-systems of all the Frobenius modules over $\text{Robba}^\text{extended}_{R,A_\infty}$;\\
IV. The category of pro-systems of all the Frobenius bundles over $\text{Robba}^\text{extended}_{R,A_\infty}$.	\\
Here $A_\infty$ is a Fr\'echet-Stein algebra attached to a compact $p$-adic Lie group such that the algebra is limit of (commutative) reduced affinoid algebras. And the finiteness is put on the infinite level of ind-scheme, actually one can also just put on each level. 

\end{itemize}

\noindent{\text{Deformation of Schemes}}

\begin{itemize}

\item<1-> \text{(Outline)} Following Kedlaya-Liu \cite{KL1}:\\
1. Construct the glueing process over the scheme $\mathrm{Spec}\text{Robba}^\text{extended}_{\infty,R,A}$;\\
2. The functors could be read off from the corresponding diagram above, namely one glues the resulting sheaves over each $\mathrm{Spec}\text{Robba}^\text{extended}_{\infty,R,A}[1/f]$ for each suitable element $f$ in the graded ring, then takes the corresponding global section;\\
3. Then from the last category back to the quasicoherent sheaves over the Fargues-Fontaine scheme we need to solve some Frobenius algebraic equation by $p$-adic analytic method to show that taking Frobenius invariance over each affine subspace is exact, where one uses Kedlaya-Liu's approach which could be dated back to Kedlaya's approach to slope filtration over extended Robba rings \cite{Ked2}.

\item<2-> Let us look back the functor:
\[
\xymatrix@R+0pc@C+0pc{
\mathrm{Mod}\mathcal{O}_{\mathrm{Proj}P_{R,A}}\ar[r]\ar[r]\ar[r] &\varphi-\mathrm{Mod}\mathcal{O}_{\mathrm{Spec}\text{Robba}^\text{extended}_{\infty,R,A}} \ar[r]\ar[r]\ar[r] &\varphi-\mathrm{Mod}\mathcal{O}_{\mathrm{Spec}\text{Robba}^\text{extended}_{R,A}} , 
}
\]
obviously one might want to generalize the picture above, which was also considered by Kedlaya-Liu in their original work \cite{KL1}, \cite{KL2}.
\item<3-> \text{(Theorem, Tong \cite[Theorem 1.4]{T})} We have the following categories are equivalent (generalizing the work of \text{Kedlaya-Liu} \cite{KL1}, \cite{KL2}, \text{Kedlaya-Pottharst} \cite{KP}): \\
I. The category of all the pseudocoherent sheaves over $\mathrm{Proj}\bigoplus P_{R,A,n}$;\\
%II. The category of all the $\varphi$-equivariant pseudocoherent sheaves over $\mathrm{Spec}\text{Robba}^\text{extended}^\infty_{R,A}$;\\
%III. The category of all the $\varphi$-equivariant pseudocoherent sheaves over $\mathrm{Spec}\text{Robba}^\text{extended}_{R,A}$;\\
%II. The category of all the Frobenius modules of the global sections of all the $\varphi$-equivariant pseudocoherent sheaves over $\mathrm{Spec}\text{Robba}^\text{extended}^\infty_{R,A}$;\\
%III. The category of all the Frobenius modules of the global sections of all the $\varphi$-equivariant pseudocoherent sheaves over $\mathrm{Spec}\text{Robba}^\text{extended}_{R,A}$;\\
II.  The category of all the pseudocoherent $\varphi$-equivariant modules over $\text{Robba}^\text{extended}_{R,A}$.
\end{itemize}

\newpage

\section{$K$-Theoretic Consideration}

\subsection{\text{The $K$-theory of Algebraic Relative Hodge-Iwasawa Modules}}
\begin{itemize}
\item<1-> Based on the study we did above, it should be very natural to consider more general pseudocoherent complexes in some higher categorical sense. Note that pseudocoherent objects were naturally emerging in SGA \cite{SGAVI} from some K-theoretic point of view. Also more importantly Hodge-Iwasawa theory to some extent will behave better if we forget the derived category, when we would like to study the K-theoretic aspects.

\item<2-> (\text{Definition}) Let $Ch\mathrm{Mod}\mathcal{O}_{\mathrm{Proj}P_{R}}$ denote the category of all the complexes of objects in $\mathrm{Mod}\mathcal{O}_{\mathrm{Proj}_{R}}$. 
\item<3-> (\text{Definition}) We now use the notations:
\begin{align}
D_{\mathrm{perf}}\mathrm{Proj}P_{R},D_{\mathrm{pseudo}}\mathrm{Proj}P_{R}
%,D_{\mathrm{perf}}\varphi-\mathrm{Mod}\mathcal{O}_{\mathrm{Spec}\text{Robba}^\text{extended}^\infty_{R}},D_{\mathrm{perf}}\varphi-\mathrm{Mod}\text{Robba}^\text{extended}^\infty_{R},D_{\mathrm{perf}}\varphi-\mathrm{Mod}\text{Robba}^\text{extended}_{R}
\end{align}
to denote the category of all the perfect and pseudocoherent complexes. 

\item<4-> (\text{Definition}) One also has the following subcategories:
\begin{align}
&D^{\mathrm{dg-flat}}_{\mathrm{perf}}\mathrm{Proj}P_{R},\\
&D^{\mathrm{str}}_{\mathrm{perf}}\mathrm{Proj}P_{R}.
%,D_{\mathrm{perf}}\varphi-\mathrm{Mod}\mathcal{O}_{\mathrm{Spec}\text{Robba}^\text{extended}^\infty_{R}},D_{\mathrm{perf}}\varphi-\mathrm{Mod}\text{Robba}^\text{extended}^\infty_{R},D_{\mathrm{perf}}\varphi-\mathrm{Mod}\text{Robba}^\text{extended}_{R}
\end{align}

\item<5-> (\text{Proposition, after Thomason-Trobaugh \cite{TT}}) These categories admit Waldhausen structure.

\item<6-> (\text{Question}) In the situation where $R=\widetilde{R}_\psi$ attached to the cyclotomic tower,  we would like to know if $D_{\mathrm{perf}}\mathrm{Proj}P_{R}$ and $D^{\mathrm{str}}_{\mathrm{perf}}\mathrm{Proj}P_{R}$ admit Waldhausen exact functors to $D_{\mathrm{perf}}(\mathbb{Q}_p)$ or $D^{\mathrm{str}}_{\mathrm{perf}}(\mathbb{Q}_p)$, which induce maps on the associated K-theory spaces.

\end{itemize}

\newpage

\section{Analytic $\infty$-Categorical Functional Analytic Hodge-Iwasawa Modules}

\subsection{$\infty$-Categorical Analytic Stacks and Descents I}

\noindent We now make the corresponding discussion after our previous work \cite{T2} on the homotopical functional analysis after many projects \cite{1BBBK}, \cite{1BBK}, \cite{BBM}, \cite{1BK} , \cite{1CS1}, \cite{1CS2}, \cite{KKM}. We choose to work over the Bambozzi-Kremnizer space \cite{1BK} attached to the corresponding Banach rings in our work after \cite{1BBBK}, \cite{1BBK}, \cite{BBM}, \cite{1BK}, \cite{KKM}. Note that what is happening is that attached to any Banach ring over $\mathbb{Q}_p$, say $B$, we attach a $(\infty,1)-$stack $\mathcal{X}(B)$ fibered over (in the sense of $\infty$-groupoid, and up to taking the corresponding opposite categories) after \cite{1BBBK}, \cite{1BBK}, \cite{BBM}, \cite{1BK}, \cite{KKM}:
\begin{align}
\mathrm{sComm}\mathrm{Simp}\mathrm{Ind}\mathrm{Ban}_{\mathbb{Q}_p},	
\end{align}
with 
\begin{align}
\mathrm{sComm}\mathrm{Simp}\mathrm{Ind}^m\mathrm{Ban}_{\mathbb{Q}_p}.	
\end{align}
associated with a $(\infty,1)$-ring object $\mathcal{O}_{\mathcal{X}(B)}$, such that we have the corresponding under the basic derived rational localization $\infty$-Grothendieck site
\begin{center}
 $(\mathcal{X}(B), \mathcal{O}_{\mathcal{X}(B),\mathrm{drl}})$ 
\end{center}
carrying the homotopical epimorphisms as the corresponding topology.

\begin{itemize}
\item<1-> By using this framework (certainly one can also consider \cite{1CS1} and \cite{1CS2} as the foundations, as in \cite{LBV}), we have the $\infty$-stack after Kedlaya-Liu \cite{KL1}, \cite{KL2}. Here in the following let $A$ be any Banach ring over $\mathbb{Q}_p$.

\item<2-> Generalizing Kedlaya-Liu's construction in \cite{KL1}, \cite{KL2} of the adic Fargues-Fontaine space we have a quotient (by using powers of the Frobenius operator) $X_{R,A}$ of the space 
\begin{align}
Y_{R,A}:=\bigcup_{0<s<r}\mathcal{X}(\text{Robba}^\text{extended}_{R,[s,r],A}).	
\end{align}

\item<3-> This is a locally ringed space $(X_{R,A},\mathcal{O}_{X_{R,A}})$, so one can consider the stable $\infty$-category $\mathrm{Ind}\mathrm{Banach}(\mathcal{O}_{X_{R,A}}) $ which is the $\infty$-category of all the $\mathcal{O}_{X_{R,A}}$-sheaves of inductive Banach modules over $X_{R,A}$. We have the parallel categories for $Y_{R,A}$, namely $\varphi\mathrm{Ind}\mathrm{Banach}(\mathcal{O}_{X_{R,A}})$ and so on. Here we will consider presheaves.
 
\item<4-> This is a locally ringed space $(X_{R,A},\mathcal{O}_{X_{R,A}})$, so one can consider the stable $\infty$-category $\mathrm{Ind}^m\mathrm{Banach}(\mathcal{O}_{X_{R,A}}) $ which is the $\infty$-category of all the $\mathcal{O}_{X_{R,A}}$-sheaves of inductive monomorphic Banach modules over $X_{R,A}$. We have the parallel categories for $Y_{R,A}$, namely $\varphi\mathrm{Ind}^m\mathrm{Banach}(\mathcal{O}_{X_{R,A}})$ and so on. Here we will consider presheaves.
   
\item<5-> In this context one can consider the $K$-theory as in the scheme situation by using the ideas and constructions from Blumberg-Gepner-Tabuada \cite{BGT}. Moreover we can study the Hodge Theory.

\item<6-> We expect that one can study among these big categories to find interesting relationships, since this should give us the right understanding of the $p$-adic Hodge theory. The corresponding pseudocoherent version comparison could be expected to be deduced as in Kedlaya-Liu's work \cite{KL1}, \cite{KL2}.

\begin{assumption}\label{assumtionpresheaves}
All the functors of modules or algebras below are presheaves.	
\end{assumption}

\item (\text{Proposition}) There is an equivalence between the $\infty$-categories of inductive Banach quasicoherent presheaves:
\[
\xymatrix@R+0pc@C+0pc{
\mathrm{Ind}\mathrm{Banach}(\mathcal{O}_{X_{R,A}})\ar[r]^{\mathrm{equi}}\ar[r]\ar[r] &\varphi\mathrm{Ind}\mathrm{Banach}(\mathcal{O}_{Y_{R,A}}).  
}
\]
\item (\text{Proposition}) There is an equivalence between the $\infty$-categories of monomorphic inductive Banach quasicoherent presheaves:
\[
\xymatrix@R+0pc@C+0pc{
\mathrm{Ind}^m\mathrm{Banach}(\mathcal{O}_{X_{R,A}})\ar[r]^{\mathrm{equi}}\ar[r]\ar[r] &\varphi\mathrm{Ind}^m\mathrm{Banach}(\mathcal{O}_{Y_{R,A}}).  
}
\]
\end{itemize}

\begin{itemize}

\item (\text{Proposition}) There is an equivalence between the $\infty$-categories of inductive Banach quasicoherent presheaves:
\[
\xymatrix@R+0pc@C+0pc{
\mathrm{Ind}\mathrm{Banach}(\mathcal{O}_{X_{R,A}})\ar[r]^{\mathrm{equi}}\ar[r]\ar[r] &\varphi\mathrm{Ind}\mathrm{Banach}(\mathcal{O}_{Y_{R,A}}).  
}
\]
\item (\text{Proposition}) There is an equivalence between the $\infty$-categories of monomorphic inductive Banach quasicoherent presheaves:
\[
\xymatrix@R+0pc@C+0pc{
\mathrm{Ind}^m\mathrm{Banach}(\mathcal{O}_{X_{R,A}})\ar[r]^{\mathrm{equi}}\ar[r]\ar[r] &\varphi\mathrm{Ind}^m\mathrm{Banach}(\mathcal{O}_{Y_{R,A}}).  
}
\]
\item (\text{Proposition}) There is an equivalence between the $\infty$-categories of inductive Banach quasicoherent commutative algebra $E_\infty$ objects:
\[
\xymatrix@R+0pc@C+0pc{
\mathrm{sComm}_\mathrm{simplicial}\mathrm{Ind}\mathrm{Banach}(\mathcal{O}_{X_{R,A}})\ar[r]^{\mathrm{equi}}\ar[r]\ar[r] &\mathrm{sComm}_\mathrm{simplicial}\varphi\mathrm{Ind}\mathrm{Banach}(\mathcal{O}_{Y_{R,A}}).  
}
\]
\item (\text{Proposition}) There is an equivalence between the $\infty$-categories of monomorphic inductive Banach quasicoherent commutative algebra $E_\infty$ objects:
\[
\xymatrix@R+0pc@C+0pc{
\mathrm{sComm}_\mathrm{simplicial}\mathrm{Ind}^m\mathrm{Banach}(\mathcal{O}_{X_{R,A}})\ar[r]^{\mathrm{equi}}\ar[r]\ar[r] &\mathrm{sComm}_\mathrm{simplicial}\varphi\mathrm{Ind}^m\mathrm{Banach}(\mathcal{O}_{Y_{R,A}}).  
}
\]

\item Then parallel as in \cite{LBV} we have the equivalence of the de Rham complex after \cite[Definition 5.9, Section 5.2.1]{KKM}:
\[
\xymatrix@R+0pc@C+0pc{
\mathrm{deRham}_{\mathrm{sComm}_\mathrm{simplicial}\mathrm{Ind}\mathrm{Banach}(\mathcal{O}_{X_{R,A}})\ar[r]^{\mathrm{equi}}}(-)\ar[r]\ar[r] &\mathrm{deRham}_{\mathrm{sComm}_\mathrm{simplicial}\varphi\mathrm{Ind}\mathrm{Banach}(\mathcal{O}_{Y_{R,A}})}(-), 
}
\]
\[
\xymatrix@R+0pc@C+0pc{
\mathrm{deRham}_{\mathrm{sComm}_\mathrm{simplicial}\mathrm{Ind}^m\mathrm{Banach}(\mathcal{O}_{X_{R,A}})\ar[r]^{\mathrm{equi}}}(-)\ar[r]\ar[r] &\mathrm{deRham}_{\mathrm{sComm}_\mathrm{simplicial}\varphi\mathrm{Ind}^m\mathrm{Banach}(\mathcal{O}_{Y_{R,A}})}(-). 
}
\]

\item Then we have the following equivalence of $K$-group $(\infty,1)$-spectrum from \cite{BGT}:
\[
\xymatrix@R+0pc@C+0pc{
\mathrm{K}^\mathrm{BGT}_{\mathrm{sComm}_\mathrm{simplicial}\mathrm{Ind}\mathrm{Banach}(\mathcal{O}_{X_{R,A}})\ar[r]^{\mathrm{equi}}}(-)\ar[r]\ar[r] &\mathrm{K}^\mathrm{BGT}_{\mathrm{sComm}_\mathrm{simplicial}\varphi\mathrm{Ind}\mathrm{Banach}(\mathcal{O}_{Y_{R,A}})}(-), 
}
\]
\[
\xymatrix@R+0pc@C+0pc{
\mathrm{K}^\mathrm{BGT}_{\mathrm{sComm}_\mathrm{simplicial}\mathrm{Ind}^m\mathrm{Banach}(\mathcal{O}_{X_{R,A}})\ar[r]^{\mathrm{equi}}}(-)\ar[r]\ar[r] &\mathrm{K}^\mathrm{BGT}_{\mathrm{sComm}_\mathrm{simplicial}\varphi\mathrm{Ind}^m\mathrm{Banach}(\mathcal{O}_{Y_{R,A}})}(-). 
}
\]
\end{itemize}

\noindent Now let $R=\mathbb{Q}_p(p^{1/p^\infty})^{\wedge\flat}$ and $R_k=\mathbb{Q}_p(p^{1/p^\infty})^{\wedge}\left<T_1^{\pm 1/p^{\infty}},...,T_k^{\pm 1/p^{\infty}}\right>^\flat$ we have the following Galois theoretic results with naturality along $f:\mathrm{Spa}(\mathbb{Q}_p(p^{1/p^\infty})^{\wedge}\left<T_1^{\pm 1/p^\infty},...,T_k^{\pm 1/p^\infty}\right>^\flat)\rightarrow \mathrm{Spa}(\mathbb{Q}_p(p^{1/p^\infty})^{\wedge\flat})$:

\begin{itemize}
\item (\text{Proposition}) There is an equivalence between the $\infty$-categories of inductive Banach quasicoherent presheaves:
\[
\xymatrix@R+6pc@C+0pc{
\mathrm{Ind}\mathrm{Banach}(\mathcal{O}_{X_{\mathbb{Q}_p(p^{1/p^\infty})^{\wedge}\left<T_1^{\pm 1/p^\infty},...,T_k^{\pm 1/p^\infty}\right>^\flat,A}})\ar[d]\ar[d]\ar[d]\ar[d] \ar[r]^{\mathrm{equi}}\ar[r]\ar[r] &\varphi\mathrm{Ind}\mathrm{Banach}(\mathcal{O}_{Y_{\mathbb{Q}_p(p^{1/p^\infty})^{\wedge}\left<T_1^{\pm 1/p^\infty},...,T_k^{\pm 1/p^\infty}\right>^\flat,A}}) \ar[d]\ar[d]\ar[d]\ar[d].\\
\mathrm{Ind}\mathrm{Banach}(\mathcal{O}_{X_{\mathbb{Q}_p(p^{1/p^\infty})^{\wedge\flat},A}})\ar[r]^{\mathrm{equi}}\ar[r]\ar[r] &\varphi\mathrm{Ind}\mathrm{Banach}(\mathcal{O}_{Y_{\mathbb{Q}_p(p^{1/p^\infty})^{\wedge\flat},A}}).\\ 
}
\]
\item (\text{Proposition}) There is an equivalence between the $\infty$-categories of monomorphic inductive Banach quasicoherent presheaves:
\[
\xymatrix@R+6pc@C+0pc{
\mathrm{Ind}^m\mathrm{Banach}(\mathcal{O}_{X_{R_k,A}})\ar[r]^{\mathrm{equi}}\ar[d]\ar[d]\ar[d]\ar[d]\ar[r]\ar[r] &\varphi\mathrm{Ind}^m\mathrm{Banach}(\mathcal{O}_{Y_{R_k,A}})\ar[d]\ar[d]\ar[d]\ar[d]\\
\mathrm{Ind}^m\mathrm{Banach}(\mathcal{O}_{X_{\mathbb{Q}_p(p^{1/p^\infty})^{\wedge\flat},A}})\ar[r]^{\mathrm{equi}}\ar[r]\ar[r] &\varphi\mathrm{Ind}^m\mathrm{Banach}(\mathcal{O}_{Y_{\mathbb{Q}_p(p^{1/p^\infty})^{\wedge\flat},A}}).\\  
}
\]
\item (\text{Proposition}) There is an equivalence between the $\infty$-categories of inductive Banach quasicoherent commutative algebra $E_\infty$ objects:
\[\displayindent=-0.1in
\xymatrix@R+6pc@C+0pc{
\mathrm{sComm}_\mathrm{simplicial}\mathrm{Ind}\mathrm{Banach}(\mathcal{O}_{X_{R_k,A}})\ar[d]\ar[d]\ar[d]\ar[d]\ar[r]^{\mathrm{equi}}\ar[r]\ar[r] &\mathrm{sComm}_\mathrm{simplicial}\varphi\mathrm{Ind}\mathrm{Banach}(\mathcal{O}_{Y_{R_k,A}})\ar[d]\ar[d]\ar[d]\ar[d]\\
\mathrm{sComm}_\mathrm{simplicial}\mathrm{Ind}\mathrm{Banach}(\mathcal{O}_{X_{\mathbb{Q}_p(p^{1/p^\infty})^{\wedge\flat},A}})\ar[r]^{\mathrm{equi}}\ar[r]\ar[r] &\mathrm{sComm}_\mathrm{simplicial}\varphi\mathrm{Ind}\mathrm{Banach}(\mathcal{O}_{Y_{\mathbb{Q}_p(p^{1/p^\infty})^{\wedge\flat},A}})  
}
\]
\item (\text{Proposition}) There is an equivalence between the $\infty$-categories of monomorphic inductive Banach quasicoherent commutative algebra $E_\infty$ objects:
\[\displayindent=-0.1in
\xymatrix@R+6pc@C+0pc{
\mathrm{sComm}_\mathrm{simplicial}\mathrm{Ind}^m\mathrm{Banach}(\mathcal{O}_{X_{R_k,A}})\ar[d]\ar[d]\ar[d]\ar[d]\ar[r]^{\mathrm{equi}}\ar[r]\ar[r] &\mathrm{sComm}_\mathrm{simplicial}\varphi\mathrm{Ind}^m\mathrm{Banach}(\mathcal{O}_{Y_{R_k,A}})\ar[d]\ar[d]\ar[d]\ar[d]\\
 \mathrm{sComm}_\mathrm{simplicial}\mathrm{Ind}^m\mathrm{Banach}(\mathcal{O}_{X_{\mathbb{Q}_p(p^{1/p^\infty})^{\wedge\flat},A}})\ar[r]^{\mathrm{equi}}\ar[r]\ar[r] &\mathrm{sComm}_\mathrm{simplicial}\varphi\mathrm{Ind}^m\mathrm{Banach}(\mathcal{O}_{Y_{\mathbb{Q}_p(p^{1/p^\infty})^{\wedge\flat},A}}) 
}
\]

\item Then parallel as in \cite{LBV} we have the equivalence of the de Rham complex after \cite[Definition 5.9, Section 5.2.1]{KKM}:
\[\displayindent=-0.2in
\xymatrix@R+6pc@C+0pc{
\mathrm{deRham}_{\mathrm{sComm}_\mathrm{simplicial}\mathrm{Ind}\mathrm{Banach}(\mathcal{O}_{X_{R_k,A}})\ar[r]^{\mathrm{equi}}}(-)\ar[d]\ar[d]\ar[d]\ar[d]\ar[r]\ar[r] &\mathrm{deRham}_{\mathrm{sComm}_\mathrm{simplicial}\varphi\mathrm{Ind}\mathrm{Banach}(\mathcal{O}_{Y_{R_k,A}})}(-)\ar[d]\ar[d]\ar[d]\ar[d]\\
\mathrm{deRham}_{\mathrm{sComm}_\mathrm{simplicial}\mathrm{Ind}\mathrm{Banach}(\mathcal{O}_{X_{\mathbb{Q}_p(p^{1/p^\infty})^{\wedge\flat},A}})\ar[r]^{\mathrm{equi}}}(-)\ar[r]\ar[r] &\mathrm{deRham}_{\mathrm{sComm}_\mathrm{simplicial}\varphi\mathrm{Ind}\mathrm{Banach}(\mathcal{O}_{Y_{\mathbb{Q}_p(p^{1/p^\infty})^{\wedge\flat},A}})}(-) 
}
\]
\[\displayindent=-0.2in
\xymatrix@R+6pc@C+0pc{
\mathrm{deRham}_{\mathrm{sComm}_\mathrm{simplicial}\mathrm{Ind}^m\mathrm{Banach}(\mathcal{O}_{X_{R_k,A}})\ar[r]^{\mathrm{equi}}}(-)\ar[d]\ar[d]\ar[d]\ar[d]\ar[r]\ar[r] &\mathrm{deRham}_{\mathrm{sComm}_\mathrm{simplicial}\varphi\mathrm{Ind}^m\mathrm{Banach}(\mathcal{O}_{Y_{R_k,A}})}(-)\ar[d]\ar[d]\ar[d]\ar[d]\\
\mathrm{deRham}_{\mathrm{sComm}_\mathrm{simplicial}\mathrm{Ind}^m\mathrm{Banach}(\mathcal{O}_{X_{\mathbb{Q}_p(p^{1/p^\infty})^{\wedge\flat},A}})\ar[r]^{\mathrm{equi}}}(-)\ar[r]\ar[r] &\mathrm{deRham}_{\mathrm{sComm}_\mathrm{simplicial}\varphi\mathrm{Ind}^m\mathrm{Banach}(\mathcal{O}_{Y_{\mathbb{Q}_p(p^{1/p^\infty})^{\wedge\flat},A}})}(-) 
}
\]

\item Then we have the following equivalence of $K$-group $(\infty,1)$-spectrum from \cite{BGT}:
\[
\xymatrix@R+6pc@C+0pc{
\mathrm{K}^\mathrm{BGT}_{\mathrm{sComm}_\mathrm{simplicial}\mathrm{Ind}\mathrm{Banach}(\mathcal{O}_{X_{R_k,A}})\ar[r]^{\mathrm{equi}}}(-)\ar[d]\ar[d]\ar[d]\ar[d]\ar[r]\ar[r] &\mathrm{K}^\mathrm{BGT}_{\mathrm{sComm}_\mathrm{simplicial}\varphi\mathrm{Ind}\mathrm{Banach}(\mathcal{O}_{Y_{R_k,A}})}(-)\ar[d]\ar[d]\ar[d]\ar[d]\\
\mathrm{K}^\mathrm{BGT}_{\mathrm{sComm}_\mathrm{simplicial}\mathrm{Ind}\mathrm{Banach}(\mathcal{O}_{X_{\mathbb{Q}_p(p^{1/p^\infty})^{\wedge\flat},A}})\ar[r]^{\mathrm{equi}}}(-)\ar[r]\ar[r] &\mathrm{K}^\mathrm{BGT}_{\mathrm{sComm}_\mathrm{simplicial}\varphi\mathrm{Ind}\mathrm{Banach}(\mathcal{O}_{Y_{\mathbb{Q}_p(p^{1/p^\infty})^{\wedge\flat},A}})}(-) 
}
\]
\[
\xymatrix@R+6pc@C+0pc{
\mathrm{K}^\mathrm{BGT}_{\mathrm{sComm}_\mathrm{simplicial}\mathrm{Ind}^m\mathrm{Banach}(\mathcal{O}_{X_{R_k,A}})\ar[r]^{\mathrm{equi}}}(-)\ar[d]\ar[d]\ar[d]\ar[d]\ar[r]\ar[r] &\mathrm{K}^\mathrm{BGT}_{\mathrm{sComm}_\mathrm{simplicial}\varphi\mathrm{Ind}^m\mathrm{Banach}(\mathcal{O}_{Y_{R_k,A}})}(-)\ar[d]\ar[d]\ar[d]\ar[d]\\
\mathrm{K}^\mathrm{BGT}_{\mathrm{sComm}_\mathrm{simplicial}\mathrm{Ind}^m\mathrm{Banach}(\mathcal{O}_{X_{\mathbb{Q}_p(p^{1/p^\infty})^{\wedge\flat},A}})\ar[r]^{\mathrm{equi}}}(-)\ar[r]\ar[r] &\mathrm{K}^\mathrm{BGT}_{\mathrm{sComm}_\mathrm{simplicial}\varphi\mathrm{Ind}^m\mathrm{Banach}(\mathcal{O}_{Y_{\mathbb{Q}_p(p^{1/p^\infty})^{\wedge\flat},A}})}(-) 
}
\]

\end{itemize}

\
\indent Then we consider further equivariance by considering the arithmetic profinite fundamental groups $\Gamma_{\mathbb{Q}_p}$ and $\mathrm{Gal}(\overline{\mathbb{Q}_p\left<T_1^{\pm 1},...,T_k^{\pm 1}\right>}/R_k)$ through the following diagram:

\[
\xymatrix@R+0pc@C+0pc{
\mathbb{Z}_p^k=\mathrm{Gal}(R_k/{\mathbb{Q}_p(p^{1/p^\infty})^\wedge\left<T_1^{\pm 1},...,T_k^{\pm 1}\right>}) \ar[r]\ar[r] \ar[r]\ar[r] &\Gamma_k:=\mathrm{Gal}(R_k/{\mathbb{Q}_p\left<T_1^{\pm 1},...,T_k^{\pm 1}\right>}) \ar[r] \ar[r]\ar[r] &\Gamma_{\mathbb{Q}_p}.
}
\]

\begin{itemize}
\item (\text{Proposition}) There is an equivalence between the $\infty$-categories of inductive Banach quasicoherent presheaves:
\[
\xymatrix@R+6pc@C+0pc{
\mathrm{Ind}\mathrm{Banach}_{\Gamma_k}(\mathcal{O}_{X_{\mathbb{Q}_p(p^{1/p^\infty})^{\wedge}\left<T_1^{\pm 1/p^\infty},...,T_k^{\pm 1/p^\infty}\right>^\flat,A}})\ar[d]\ar[d]\ar[d]\ar[d] \ar[r]^{\mathrm{equi}}\ar[r]\ar[r] &\varphi\mathrm{Ind}\mathrm{Banach}_{\Gamma_k}(\mathcal{O}_{Y_{\mathbb{Q}_p(p^{1/p^\infty})^{\wedge}\left<T_1^{\pm 1/p^\infty},...,T_k^{\pm 1/p^\infty}\right>^\flat,A}}) \ar[d]\ar[d]\ar[d]\ar[d].\\
\mathrm{Ind}\mathrm{Banach}(\mathcal{O}_{X_{\mathbb{Q}_p(p^{1/p^\infty})^{\wedge\flat},A}})\ar[r]^{\mathrm{equi}}\ar[r]\ar[r] &\varphi\mathrm{Ind}\mathrm{Banach}(\mathcal{O}_{Y_{\mathbb{Q}_p(p^{1/p^\infty})^{\wedge\flat},A}}).\\ 
}
\]
\item (\text{Proposition}) There is an equivalence between the $\infty$-categories of monomorphic inductive Banach quasicoherent presheaves:
\[
\xymatrix@R+6pc@C+0pc{
\mathrm{Ind}^m\mathrm{Banach}_{\Gamma_k}(\mathcal{O}_{X_{R_k,A}})\ar[r]^{\mathrm{equi}}\ar[d]\ar[d]\ar[d]\ar[d]\ar[r]\ar[r] &\varphi\mathrm{Ind}^m\mathrm{Banach}_{\Gamma_k}(\mathcal{O}_{Y_{R_k,A}})\ar[d]\ar[d]\ar[d]\ar[d]\\
\mathrm{Ind}^m\mathrm{Banach}_{\Gamma_0}(\mathcal{O}_{X_{\mathbb{Q}_p(p^{1/p^\infty})^{\wedge\flat},A}})\ar[r]^{\mathrm{equi}}\ar[r]\ar[r] &\varphi\mathrm{Ind}^m\mathrm{Banach}_{\Gamma_0}(\mathcal{O}_{Y_{\mathbb{Q}_p(p^{1/p^\infty})^{\wedge\flat},A}}).\\  
}
\]
\item (\text{Proposition}) There is an equivalence between the $\infty$-categories of inductive Banach quasicoherent commutative algebra $E_\infty$ objects:
\[
\xymatrix@R+6pc@C+0pc{
\mathrm{sComm}_\mathrm{simplicial}\mathrm{Ind}\mathrm{Banach}_{\Gamma_k}(\mathcal{O}_{X_{R_k,A}})\ar[d]\ar[d]\ar[d]\ar[d]\ar[r]^{\mathrm{equi}}\ar[r]\ar[r] &\mathrm{sComm}_\mathrm{simplicial}\varphi\mathrm{Ind}\mathrm{Banach}_{\Gamma_k}(\mathcal{O}_{Y_{R_k,A}})\ar[d]\ar[d]\ar[d]\ar[d]\\
\mathrm{sComm}_\mathrm{simplicial}\mathrm{Ind}\mathrm{Banach}_{\Gamma_0}(\mathcal{O}_{X_{\mathbb{Q}_p(p^{1/p^\infty})^{\wedge\flat},A}})\ar[r]^{\mathrm{equi}}\ar[r]\ar[r] &\mathrm{sComm}_\mathrm{simplicial}\varphi\mathrm{Ind}\mathrm{Banach}_{\Gamma_0}(\mathcal{O}_{Y_{\mathbb{Q}_p(p^{1/p^\infty})^{\wedge\flat},A}})  
}
\]
\item (\text{Proposition}) There is an equivalence between the $\infty$-categories of monomorphic inductive Banach quasicoherent commutative algebra $E_\infty$ objects:
\[
\xymatrix@R+6pc@C+0pc{
\mathrm{sComm}_\mathrm{simplicial}\mathrm{Ind}^m\mathrm{Banach}_{\Gamma_k}(\mathcal{O}_{X_{R_k,A}})\ar[d]\ar[d]\ar[d]\ar[d]\ar[r]^{\mathrm{equi}}\ar[r]\ar[r] &\mathrm{sComm}_\mathrm{simplicial}\varphi\mathrm{Ind}^m\mathrm{Banach}_{\Gamma_k}(\mathcal{O}_{Y_{R_k,A}})\ar[d]\ar[d]\ar[d]\ar[d]\\
 \mathrm{sComm}_\mathrm{simplicial}\mathrm{Ind}^m\mathrm{Banach}_{\Gamma_0}(\mathcal{O}_{X_{\mathbb{Q}_p(p^{1/p^\infty})^{\wedge\flat},A}})\ar[r]^{\mathrm{equi}}\ar[r]\ar[r] &\mathrm{sComm}_\mathrm{simplicial}\varphi\mathrm{Ind}^m\mathrm{Banach}_{\Gamma_0}(\mathcal{O}_{Y_{\mathbb{Q}_p(p^{1/p^\infty})^{\wedge\flat},A}}) 
}
\]

\item Then parallel as in \cite{LBV} we have the equivalence of the de Rham complex after \cite[Definition 5.9, Section 5.2.1]{KKM}:
\[\displayindent=-0.2in
\xymatrix@R+6pc@C+0pc{
\mathrm{deRham}_{\mathrm{sComm}_\mathrm{simplicial}\mathrm{Ind}\mathrm{Banach}_{\Gamma_k}(\mathcal{O}_{X_{R_k,A}})\ar[r]^{\mathrm{equi}}}(-)\ar[d]\ar[d]\ar[d]\ar[d]\ar[r]\ar[r] &\mathrm{deRham}_{\mathrm{sComm}_\mathrm{simplicial}\varphi\mathrm{Ind}\mathrm{Banach}_{\Gamma_k}(\mathcal{O}_{Y_{R_k,A}})}(-)\ar[d]\ar[d]\ar[d]\ar[d]\\
\mathrm{deRham}_{\mathrm{sComm}_\mathrm{simplicial}\mathrm{Ind}\mathrm{Banach}_{\Gamma_0}(\mathcal{O}_{X_{\mathbb{Q}_p(p^{1/p^\infty})^{\wedge\flat},A}})\ar[r]^{\mathrm{equi}}}(-)\ar[r]\ar[r] &\mathrm{deRham}_{\mathrm{sComm}_\mathrm{simplicial}\varphi\mathrm{Ind}\mathrm{Banach}_{\Gamma_0}(\mathcal{O}_{Y_{\mathbb{Q}_p(p^{1/p^\infty})^{\wedge\flat},A}})}(-) 
}
\]
\[\displayindent=-0.4in
\xymatrix@R+6pc@C+0pc{
\mathrm{deRham}_{\mathrm{sComm}_\mathrm{simplicial}\mathrm{Ind}^m\mathrm{Banach}_{\Gamma_k}(\mathcal{O}_{X_{R_k,A}})\ar[r]^{\mathrm{equi}}}(-)\ar[d]\ar[d]\ar[d]\ar[d]\ar[r]\ar[r] &\mathrm{deRham}_{\mathrm{sComm}_\mathrm{simplicial}\varphi\mathrm{Ind}^m\mathrm{Banach}_{\Gamma_k}(\mathcal{O}_{Y_{R_k,A}})}(-)\ar[d]\ar[d]\ar[d]\ar[d]\\
\mathrm{deRham}_{\mathrm{sComm}_\mathrm{simplicial}\mathrm{Ind}^m\mathrm{Banach}_{\Gamma_0}(\mathcal{O}_{X_{\mathbb{Q}_p(p^{1/p^\infty})^{\wedge\flat},A}})\ar[r]^{\mathrm{equi}}}(-)\ar[r]\ar[r] &\mathrm{deRham}_{\mathrm{sComm}_\mathrm{simplicial}\varphi\mathrm{Ind}^m\mathrm{Banach}_{\Gamma_0}(\mathcal{O}_{Y_{\mathbb{Q}_p(p^{1/p^\infty})^{\wedge\flat},A}})}(-) 
}
\]

\item Then we have the following equivalence of $K$-group $(\infty,1)$-spectrum from \cite{BGT}:
\[
\xymatrix@R+6pc@C+0pc{
\mathrm{K}^\mathrm{BGT}_{\mathrm{sComm}_\mathrm{simplicial}\mathrm{Ind}\mathrm{Banach}_{\Gamma_k}(\mathcal{O}_{X_{R_k,A}})\ar[r]^{\mathrm{equi}}}(-)\ar[d]\ar[d]\ar[d]\ar[d]\ar[r]\ar[r] &\mathrm{K}^\mathrm{BGT}_{\mathrm{sComm}_\mathrm{simplicial}\varphi\mathrm{Ind}\mathrm{Banach}_{\Gamma_k}(\mathcal{O}_{Y_{R_k,A}})}(-)\ar[d]\ar[d]\ar[d]\ar[d]\\
\mathrm{K}^\mathrm{BGT}_{\mathrm{sComm}_\mathrm{simplicial}\mathrm{Ind}\mathrm{Banach}_{\Gamma_0}(\mathcal{O}_{X_{\mathbb{Q}_p(p^{1/p^\infty})^{\wedge\flat},A}})\ar[r]^{\mathrm{equi}}}(-)\ar[r]\ar[r] &\mathrm{K}^\mathrm{BGT}_{\mathrm{sComm}_\mathrm{simplicial}\varphi\mathrm{Ind}\mathrm{Banach}_{\Gamma_0}(\mathcal{O}_{Y_{\mathbb{Q}_p(p^{1/p^\infty})^{\wedge\flat},A}})}(-) 
}
\]
\[
\xymatrix@R+6pc@C+0pc{
\mathrm{K}^\mathrm{BGT}_{\mathrm{sComm}_\mathrm{simplicial}\mathrm{Ind}^m\mathrm{Banach}_{\Gamma_k}(\mathcal{O}_{X_{R_k,A}})\ar[r]^{\mathrm{equi}}}(-)\ar[d]\ar[d]\ar[d]\ar[d]\ar[r]\ar[r] &\mathrm{K}^\mathrm{BGT}_{\mathrm{sComm}_\mathrm{simplicial}\varphi\mathrm{Ind}^m\mathrm{Banach}_{\Gamma_k}(\mathcal{O}_{Y_{R_k,A}})}(-)\ar[d]\ar[d]\ar[d]\ar[d]\\
\mathrm{K}^\mathrm{BGT}_{\mathrm{sComm}_\mathrm{simplicial}\mathrm{Ind}^m\mathrm{Banach}_{\Gamma_0}(\mathcal{O}_{X_{\mathbb{Q}_p(p^{1/p^\infty})^{\wedge\flat},A}})\ar[r]^{\mathrm{equi}}}(-)\ar[r]\ar[r] &\mathrm{K}^\mathrm{BGT}_{\mathrm{sComm}_\mathrm{simplicial}\varphi\mathrm{Ind}^m\mathrm{Banach}_{\Gamma_0}(\mathcal{O}_{Y_{\mathbb{Q}_p(p^{1/p^\infty})^{\wedge\flat},A}})}(-). 
}
\]

\end{itemize}

\

Furthermore we have the corresponding pro-\'etale version without the corresponding fundamental group equivariances.

\begin{itemize}
\item (\text{Proposition}) There is an equivalence between the $\infty$-categories of inductive Banach quasicoherent presheaves:
\[
\xymatrix@R+6pc@C+0pc{
\mathrm{Ind}\mathrm{Banach}(\mathcal{O}_{X_{\mathbb{Q}_p\left<T_1^{\pm 1},...,T_k^{\pm 1}\right>,\text{pro\'etale},A}})\ar[d]\ar[d]\ar[d]\ar[d] \ar[r]^{\mathrm{equi}}\ar[r]\ar[r] &\varphi\mathrm{Ind}\mathrm{Banach}(\mathcal{O}_{Y_{\mathbb{Q}_p\left<T_1^{\pm 1},...,T_k^{\pm 1}\right>,\text{pro\'etale},A}}) \ar[d]\ar[d]\ar[d]\ar[d].\\
\mathrm{Ind}\mathrm{Banach}(\mathcal{O}_{X_{\mathbb{Q}_p,\text{pro\'etale},A}})\ar[r]^{\mathrm{equi}}\ar[r]\ar[r] &\varphi\mathrm{Ind}\mathrm{Banach}(\mathcal{O}_{Y_{\mathbb{Q}_p,\text{pro\'etale},A}}).\\ 
}
\]
\item (\text{Proposition}) There is an equivalence between the $\infty$-categories of monomorphic inductive Banach quasicoherent presheaves:
\[
\xymatrix@R+6pc@C+0pc{
\mathrm{Ind}^m\mathrm{Banach}(\mathcal{O}_{X_{\mathbb{Q}_p\left<T_1^{\pm 1},...,T_k^{\pm 1}\right>,\text{pro\'etale},A}})\ar[r]^{\mathrm{equi}}\ar[d]\ar[d]\ar[d]\ar[d]\ar[r]\ar[r] &\varphi\mathrm{Ind}^m\mathrm{Banach}(\mathcal{O}_{Y_{\mathbb{Q}_p\left<T_1^{\pm 1},...,T_k^{\pm 1}\right>,\text{pro\'etale},A}})\ar[d]\ar[d]\ar[d]\ar[d]\\
\mathrm{Ind}^m\mathrm{Banach}(\mathcal{O}_{X_{\mathbb{Q}_p,\text{pro\'etale},A}})\ar[r]^{\mathrm{equi}}\ar[r]\ar[r] &\varphi\mathrm{Ind}^m\mathrm{Banach}(\mathcal{O}_{Y_{\mathbb{Q}_p,\text{pro\'etale},A}}).\\  
}
\]
\item (\text{Proposition}) There is an equivalence between the $\infty$-categories of inductive Banach quasicoherent commutative algebra $E_\infty$ objects:
\[\displayindent=-0.4in
\xymatrix@R+6pc@C+0pc{
\mathrm{sComm}_\mathrm{simplicial}\mathrm{Ind}\mathrm{Banach}(\mathcal{O}_{X_{\mathbb{Q}_p\left<T_1^{\pm 1},...,T_k^{\pm 1}\right>,\text{pro\'etale},A}})\ar[d]\ar[d]\ar[d]\ar[d]\ar[r]^{\mathrm{equi}}\ar[r]\ar[r] &\mathrm{sComm}_\mathrm{simplicial}\varphi\mathrm{Ind}\mathrm{Banach}(\mathcal{O}_{Y_{\mathbb{Q}_p\left<T_1^{\pm 1},...,T_k^{\pm 1}\right>,\text{pro\'etale},A}})\ar[d]\ar[d]\ar[d]\ar[d]\\
\mathrm{sComm}_\mathrm{simplicial}\mathrm{Ind}\mathrm{Banach}(\mathcal{O}_{X_{\mathbb{Q}_p,\text{pro\'etale},A}})\ar[r]^{\mathrm{equi}}\ar[r]\ar[r] &\mathrm{sComm}_\mathrm{simplicial}\varphi\mathrm{Ind}\mathrm{Banach}(\mathcal{O}_{Y_{\mathbb{Q}_p,\text{pro\'etale},A}})  
}
\]
\item (\text{Proposition}) There is an equivalence between the $\infty$-categories of monomorphic inductive Banach quasicoherent commutative algebra $E_\infty$ objects:
\[\displayindent=-0.7in
\xymatrix@R+6pc@C+0pc{
\mathrm{sComm}_\mathrm{simplicial}\mathrm{Ind}^m\mathrm{Banach}(\mathcal{O}_{X_{\mathbb{Q}_p\left<T_1^{\pm 1},...,T_k^{\pm 1}\right>,\text{pro\'etale},A}})\ar[d]\ar[d]\ar[d]\ar[d]\ar[r]^{\mathrm{equi}}\ar[r]\ar[r] &\mathrm{sComm}_\mathrm{simplicial}\varphi\mathrm{Ind}^m\mathrm{Banach}(\mathcal{O}_{Y_{\mathbb{Q}_p\left<T_1^{\pm 1},...,T_k^{\pm 1}\right>,\text{pro\'etale},A}})\ar[d]\ar[d]\ar[d]\ar[d]\\
 \mathrm{sComm}_\mathrm{simplicial}\mathrm{Ind}^m\mathrm{Banach}(\mathcal{O}_{X_{\mathbb{Q}_p,\text{pro\'etale},A}})\ar[r]^{\mathrm{equi}}\ar[r]\ar[r] &\mathrm{sComm}_\mathrm{simplicial}\varphi\mathrm{Ind}^m\mathrm{Banach}(\mathcal{O}_{Y_{\mathbb{Q}_p,\text{pro\'etale},A}}) 
}
\]

\item Then parallel as in \cite{LBV} we have the equivalence of the de Rham complex after \cite[Definition 5.9, Section 5.2.1]{KKM}:
\[\displayindent=-0.9in
\xymatrix@R+6pc@C+0pc{
\mathrm{deRham}_{\mathrm{sComm}_\mathrm{simplicial}\mathrm{Ind}\mathrm{Banach}(\mathcal{O}_{X_{\mathbb{Q}_p\left<T_1^{\pm 1},...,T_k^{\pm 1}\right>,\text{pro\'etale},A}})\ar[r]^{\mathrm{equi}}}(-)\ar[d]\ar[d]\ar[d]\ar[d]\ar[r]\ar[r] &\mathrm{deRham}_{\mathrm{sComm}_\mathrm{simplicial}\varphi\mathrm{Ind}\mathrm{Banach}(\mathcal{O}_{Y_{\mathbb{Q}_p\left<T_1^{\pm 1},...,T_k^{\pm 1}\right>,\text{pro\'etale},A}})}(-)\ar[d]\ar[d]\ar[d]\ar[d]\\
\mathrm{deRham}_{\mathrm{sComm}_\mathrm{simplicial}\mathrm{Ind}\mathrm{Banach}(\mathcal{O}_{X_{\mathbb{Q}_p,\text{pro\'etale},A}})\ar[r]^{\mathrm{equi}}}(-)\ar[r]\ar[r] &\mathrm{deRham}_{\mathrm{sComm}_\mathrm{simplicial}\varphi\mathrm{Ind}\mathrm{Banach}(\mathcal{O}_{Y_{\mathbb{Q}_p,\text{pro\'etale},A}})}(-) 
}
\]
\[\displayindent=-1.0in
\xymatrix@R+6pc@C+0pc{
\mathrm{deRham}_{\mathrm{sComm}_\mathrm{simplicial}\mathrm{Ind}^m\mathrm{Banach}(\mathcal{O}_{X_{\mathbb{Q}_p\left<T_1^{\pm 1},...,T_k^{\pm 1}\right>,\text{pro\'etale},A}})\ar[r]^{\mathrm{equi}}}(-)\ar[d]\ar[d]\ar[d]\ar[d]\ar[r]\ar[r] &\mathrm{deRham}_{\mathrm{sComm}_\mathrm{simplicial}\varphi\mathrm{Ind}^m\mathrm{Banach}(\mathcal{O}_{Y_{\mathbb{Q}_p\left<T_1^{\pm 1},...,T_k^{\pm 1}\right>,\text{pro\'etale},A}})}(-)\ar[d]\ar[d]\ar[d]\ar[d]\\
\mathrm{deRham}_{\mathrm{sComm}_\mathrm{simplicial}\mathrm{Ind}^m\mathrm{Banach}(\mathcal{O}_{X_{\mathbb{Q}_p,\text{pro\'etale},A}})\ar[r]^{\mathrm{equi}}}(-)\ar[r]\ar[r] &\mathrm{deRham}_{\mathrm{sComm}_\mathrm{simplicial}\varphi\mathrm{Ind}^m\mathrm{Banach}(\mathcal{O}_{Y_{\mathbb{Q}_p,\text{pro\'etale},A}})}(-) 
}
\]

\item Then we have the following equivalence of $K$-group $(\infty,1)$-spectrum from \cite{BGT}:
\[\displayindent=-0.4in
\xymatrix@R+6pc@C+0pc{
\mathrm{K}^\mathrm{BGT}_{\mathrm{sComm}_\mathrm{simplicial}\mathrm{Ind}\mathrm{Banach}(\mathcal{O}_{X_{\mathbb{Q}_p\left<T_1^{\pm 1},...,T_k^{\pm 1}\right>,\text{pro\'etale},A}})\ar[r]^{\mathrm{equi}}}(-)\ar[d]\ar[d]\ar[d]\ar[d]\ar[r]\ar[r] &\mathrm{K}^\mathrm{BGT}_{\mathrm{sComm}_\mathrm{simplicial}\varphi\mathrm{Ind}\mathrm{Banach}(\mathcal{O}_{Y_{\mathbb{Q}_p\left<T_1^{\pm 1},...,T_k^{\pm 1}\right>,\text{pro\'etale},A}})}(-)\ar[d]\ar[d]\ar[d]\ar[d]\\
\mathrm{K}^\mathrm{BGT}_{\mathrm{sComm}_\mathrm{simplicial}\mathrm{Ind}\mathrm{Banach}(\mathcal{O}_{X_{\mathbb{Q}_p,\text{pro\'etale},A}})\ar[r]^{\mathrm{equi}}}(-)\ar[r]\ar[r] &\mathrm{K}^\mathrm{BGT}_{\mathrm{sComm}_\mathrm{simplicial}\varphi\mathrm{Ind}\mathrm{Banach}(\mathcal{O}_{Y_{\mathbb{Q}_p,\text{pro\'etale},A}})}(-) 
}
\]
\[\displayindent=-0.4in
\xymatrix@R+6pc@C+0pc{
\mathrm{K}^\mathrm{BGT}_{\mathrm{sComm}_\mathrm{simplicial}\mathrm{Ind}^m\mathrm{Banach}(\mathcal{O}_{X_{\mathbb{Q}_p\left<T_1^{\pm 1},...,T_k^{\pm 1}\right>,\text{pro\'etale},A}})\ar[r]^{\mathrm{equi}}}(-)\ar[d]\ar[d]\ar[d]\ar[d]\ar[r]\ar[r] &\mathrm{K}^\mathrm{BGT}_{\mathrm{sComm}_\mathrm{simplicial}\varphi\mathrm{Ind}^m\mathrm{Banach}(\mathcal{O}_{Y_{\mathbb{Q}_p\left<T_1^{\pm 1},...,T_k^{\pm 1}\right>,\text{pro\'etale},A}})}(-)\ar[d]\ar[d]\ar[d]\ar[d]\\
\mathrm{K}^\mathrm{BGT}_{\mathrm{sComm}_\mathrm{simplicial}\mathrm{Ind}^m\mathrm{Banach}(\mathcal{O}_{X_{\mathbb{Q}_p,\text{pro\'etale},A}})\ar[r]^{\mathrm{equi}}}(-)\ar[r]\ar[r] &\mathrm{K}^\mathrm{BGT}_{\mathrm{sComm}_\mathrm{simplicial}\varphi\mathrm{Ind}^m\mathrm{Banach}(\mathcal{O}_{Y_{\mathbb{Q}_p,\text{pro\'etale},A}})}(-). 
}
\]
\\

\end{itemize}

\

\indent Now we consider \cite{1CS1} and \cite[Proposition 13.8, Theorem 14.9, Remark 14.10]{1CS2}\footnote{Note that we are motivated as well from \cite{LBV}.}, and study the corresponding solid perfect complexes, solid quasicoherent sheaves and solid vector bundles. Here we are going to use different formalism, therefore we will have different categories and functors. We use the notation $\circledcirc$ to denote any element of $\{\text{solid perfect complexes}, \text{solid quasicoherent sheaves}, \text{solid vector bundles}\}$ from \cite{1CS2} with the corresponding descent results of \cite[Proposition 13.8, Theorem 14.9, Remark 14.10]{1CS2}. Then we have the following:

\begin{itemize}
\item Generalizing Kedlaya-Liu's construction in \cite{KL1}, \cite{KL2} of the adic Fargues-Fontaine space we have a quotient (by using powers of the Frobenius operator) $X_{R,A}$ of the space by using \cite{1CS2}:
\begin{align}
Y_{R,A}:=\bigcup_{0<s<r}\mathcal{X}^\mathrm{CS}(\text{Robba}^\text{extended}_{R,[s,r],A}).	
\end{align}

\item This is a locally ringed space $(X_{R,A},\mathcal{O}_{X_{R,A}})$, so one can consider the stable $\infty$-category $\mathrm{Module}_{\text{CS},\mathrm{quasicoherent}}(\mathcal{O}_{X_{R,A}}) $ which is the $\infty$-category of all the $\mathcal{O}_{X_{R,A}}$-sheaves of solid modules over $X_{R,A}$. We have the parallel categories for $Y_{R,A}$, namely $\varphi\mathrm{Module}_{\text{CS},\mathrm{quasicoherent}}(\mathcal{O}_{X_{R,A}})$ and so on. Here we will consider sheaves.

\begin{assumption}\label{assumtionpresheaves}
All the functors of modules or algebras below are Clausen-Scholze sheaves \cite[Proposition 13.8, Theorem 14.9, Remark 14.10]{1CS2}.	
\end{assumption}

\item (\text{Proposition}) There is an equivalence between the $\infty$-categories of inductive solid sheaves:
\[
\xymatrix@R+0pc@C+0pc{
\mathrm{Module}_\circledcirc(\mathcal{O}_{X_{R,A}})\ar[r]^{\mathrm{equi}}\ar[r]\ar[r] &\varphi\mathrm{Module}_\circledcirc(\mathcal{O}_{Y_{R,A}}).  
}
\]
\end{itemize}

\begin{itemize}

\item (\text{Proposition}) There is an equivalence between the $\infty$-categories of inductive Banach quasicoherent commutative algebra $E_\infty$ objects:
\[
\xymatrix@R+0pc@C+0pc{
\mathrm{sComm}_\mathrm{simplicial}\mathrm{Module}_{\text{solidquasicoherentsheaves}}(\mathcal{O}_{X_{R,A}})\ar[r]^{\mathrm{equi}}\ar[r]\ar[r] &\mathrm{sComm}_\mathrm{simplicial}\varphi\mathrm{Module}_{\text{solidquasicoherentsheaves}}(\mathcal{O}_{Y_{R,A}}).  
}
\]

\item Then as in \cite{LBV} we have the equivalence of the de Rham complex after \cite[Definition 5.9, Section 5.2.1]{KKM}\footnote{Here $\circledcirc=\text{solidquasicoherentsheaves}$.}:
\[
\xymatrix@R+0pc@C+0pc{
\mathrm{deRham}_{\mathrm{sComm}_\mathrm{simplicial}\mathrm{Module}_\circledcirc(\mathcal{O}_{X_{R,A}})\ar[r]^{\mathrm{equi}}}(-)\ar[r]\ar[r] &\mathrm{deRham}_{\mathrm{sComm}_\mathrm{simplicial}\varphi\mathrm{Module}_\circledcirc(\mathcal{O}_{Y_{R,A}})}(-). 
}
\]

\item Then we have the following equivalence of $K$-group $(\infty,1)$-spectrum from \cite{BGT}\footnote{Here $\circledcirc=\text{solidquasicoherentsheaves}$.}:
\[
\xymatrix@R+0pc@C+0pc{
\mathrm{K}^\mathrm{BGT}_{\mathrm{sComm}_\mathrm{simplicial}\mathrm{Module}_\circledcirc(\mathcal{O}_{X_{R,A}})\ar[r]^{\mathrm{equi}}}(-)\ar[r]\ar[r] &\mathrm{K}^\mathrm{BGT}_{\mathrm{sComm}_\mathrm{simplicial}\varphi\mathrm{Module}_\circledcirc(\mathcal{O}_{Y_{R,A}})}(-). 
}
\]
\end{itemize}

\noindent Now let $R=\mathbb{Q}_p(p^{1/p^\infty})^{\wedge\flat}$ and $R_k=\mathbb{Q}_p(p^{1/p^\infty})^{\wedge}\left<T_1^{\pm 1/p^{\infty}},...,T_k^{\pm 1/p^{\infty}}\right>^\flat$ we have the following Galois theoretic results with naturality along $f:\mathrm{Spa}(\mathbb{Q}_p(p^{1/p^\infty})^{\wedge}\left<T_1^{\pm 1/p^\infty},...,T_k^{\pm 1/p^\infty}\right>^\flat)\rightarrow \mathrm{Spa}(\mathbb{Q}_p(p^{1/p^\infty})^{\wedge\flat})$:

\begin{itemize}
\item (\text{Proposition}) There is an equivalence between the $\infty$-categories of inductive Banach quasicoherent presheaves\footnote{Here $\circledcirc=\text{solidquasicoherentsheaves}$.}:
\[
\xymatrix@R+6pc@C+0pc{
\mathrm{Modules}_\circledcirc(\mathcal{O}_{X_{\mathbb{Q}_p(p^{1/p^\infty})^{\wedge}\left<T_1^{\pm 1/p^\infty},...,T_k^{\pm 1/p^\infty}\right>^\flat,A}})\ar[d]\ar[d]\ar[d]\ar[d] \ar[r]^{\mathrm{equi}}\ar[r]\ar[r] &\varphi\mathrm{Modules}_\circledcirc(\mathcal{O}_{Y_{\mathbb{Q}_p(p^{1/p^\infty})^{\wedge}\left<T_1^{\pm 1/p^\infty},...,T_k^{\pm 1/p^\infty}\right>^\flat,A}}) \ar[d]\ar[d]\ar[d]\ar[d].\\
\mathrm{Modules}_\circledcirc(\mathcal{O}_{X_{\mathbb{Q}_p(p^{1/p^\infty})^{\wedge\flat},A}})\ar[r]^{\mathrm{equi}}\ar[r]\ar[r] &\varphi\mathrm{Modules}_\circledcirc(\mathcal{O}_{Y_{\mathbb{Q}_p(p^{1/p^\infty})^{\wedge\flat},A}}).\\ 
}
\]
\item (\text{Proposition}) There is an equivalence between the $\infty$-categories of inductive Banach quasicoherent commutative algebra $E_\infty$ objects\footnote{Here $\circledcirc=\text{solidquasicoherentsheaves}$.}:
\[
\xymatrix@R+6pc@C+0pc{
\mathrm{sComm}_\mathrm{simplicial}\mathrm{Modules}_\circledcirc(\mathcal{O}_{X_{R_k,A}})\ar[d]\ar[d]\ar[d]\ar[d]\ar[r]^{\mathrm{equi}}\ar[r]\ar[r] &\mathrm{sComm}_\mathrm{simplicial}\varphi\mathrm{Modules}_\circledcirc(\mathcal{O}_{Y_{R_k,A}})\ar[d]\ar[d]\ar[d]\ar[d]\\
\mathrm{sComm}_\mathrm{simplicial}\mathrm{Modules}_\circledcirc(\mathcal{O}_{X_{\mathbb{Q}_p(p^{1/p^\infty})^{\wedge\flat},A}})\ar[r]^{\mathrm{equi}}\ar[r]\ar[r] &\mathrm{sComm}_\mathrm{simplicial}\varphi\mathrm{Modules}_\circledcirc(\mathcal{O}_{Y_{\mathbb{Q}_p(p^{1/p^\infty})^{\wedge\flat},A}}).  
}
\]
\item Then as in \cite{LBV} we have the equivalence of the de Rham complex after \cite[Definition 5.9, Section 5.2.1]{KKM}\footnote{Here $\circledcirc=\text{solidquasicoherentsheaves}$.}:
\[
\xymatrix@R+6pc@C+0pc{
\mathrm{deRham}_{\mathrm{sComm}_\mathrm{simplicial}\mathrm{Modules}_\circledcirc(\mathcal{O}_{X_{R_k,A}})\ar[r]^{\mathrm{equi}}}(-)\ar[d]\ar[d]\ar[d]\ar[d]\ar[r]\ar[r] &\mathrm{deRham}_{\mathrm{sComm}_\mathrm{simplicial}\varphi\mathrm{Modules}_\circledcirc(\mathcal{O}_{Y_{R_k,A}})}(-)\ar[d]\ar[d]\ar[d]\ar[d]\\
\mathrm{deRham}_{\mathrm{sComm}_\mathrm{simplicial}\mathrm{Modules}_\circledcirc(\mathcal{O}_{X_{\mathbb{Q}_p(p^{1/p^\infty})^{\wedge\flat},A}})\ar[r]^{\mathrm{equi}}}(-)\ar[r]\ar[r] &\mathrm{deRham}_{\mathrm{sComm}_\mathrm{simplicial}\varphi\mathrm{Modules}_\circledcirc(\mathcal{O}_{Y_{\mathbb{Q}_p(p^{1/p^\infty})^{\wedge\flat},A}})}(-). 
}
\]

\item Then we have the following equivalence of $K$-group $(\infty,1)$-spectrum from \cite{BGT}\footnote{Here $\circledcirc=\text{solidquasicoherentsheaves}$.}:
\[
\xymatrix@R+6pc@C+0pc{
\mathrm{K}^\mathrm{BGT}_{\mathrm{sComm}_\mathrm{simplicial}\mathrm{Modules}_\circledcirc(\mathcal{O}_{X_{R_k,A}})\ar[r]^{\mathrm{equi}}}(-)\ar[d]\ar[d]\ar[d]\ar[d]\ar[r]\ar[r] &\mathrm{K}^\mathrm{BGT}_{\mathrm{sComm}_\mathrm{simplicial}\varphi\mathrm{Modules}_\circledcirc(\mathcal{O}_{Y_{R_k,A}})}(-)\ar[d]\ar[d]\ar[d]\ar[d]\\
\mathrm{K}^\mathrm{BGT}_{\mathrm{sComm}_\mathrm{simplicial}\mathrm{Modules}_\circledcirc(\mathcal{O}_{X_{\mathbb{Q}_p(p^{1/p^\infty})^{\wedge\flat},A}})\ar[r]^{\mathrm{equi}}}(-)\ar[r]\ar[r] &\mathrm{K}^\mathrm{BGT}_{\mathrm{sComm}_\mathrm{simplicial}\varphi\mathrm{Modules}_\circledcirc(\mathcal{O}_{Y_{\mathbb{Q}_p(p^{1/p^\infty})^{\wedge\flat},A}})}(-).}
\]

\end{itemize}

\
\indent Then we consider further equivariance by considering the arithmetic profinite fundamental groups $\Gamma_{\mathbb{Q}_p}$ and $\mathrm{Gal}(\overline{\mathbb{Q}_p\left<T_1^{\pm 1},...,T_k^{\pm 1}\right>}/R_k)$ through the following diagram:

\[
\xymatrix@R+0pc@C+0pc{
\mathbb{Z}_p^k=\mathrm{Gal}(R_k/{\mathbb{Q}_p(p^{1/p^\infty})^\wedge\left<T_1^{\pm 1},...,T_k^{\pm 1}\right>}) \ar[r]\ar[r] \ar[r]\ar[r] &\Gamma_k:=\mathrm{Gal}(R_k/{\mathbb{Q}_p\left<T_1^{\pm 1},...,T_k^{\pm 1}\right>}) \ar[r] \ar[r]\ar[r] &\Gamma_{\mathbb{Q}_p}.
}
\]

\begin{itemize}
\item (\text{Proposition}) There is an equivalence between the $\infty$-categories of inductive Banach quasicoherent presheaves\footnote{Here $\circledcirc=\text{solidquasicoherentsheaves}$.}:
\[
\xymatrix@R+6pc@C+0pc{
\mathrm{Modules}_{\circledcirc,\Gamma_k}(\mathcal{O}_{X_{\mathbb{Q}_p(p^{1/p^\infty})^{\wedge}\left<T_1^{\pm 1/p^\infty},...,T_k^{\pm 1/p^\infty}\right>^\flat,A}})\ar[d]\ar[d]\ar[d]\ar[d] \ar[r]^{\mathrm{equi}}\ar[r]\ar[r] &\varphi\mathrm{Modules}_{\circledcirc,\Gamma_k}(\mathcal{O}_{Y_{\mathbb{Q}_p(p^{1/p^\infty})^{\wedge}\left<T_1^{\pm 1/p^\infty},...,T_k^{\pm 1/p^\infty}\right>^\flat,A}}) \ar[d]\ar[d]\ar[d]\ar[d].\\
\mathrm{Modules}_{\circledcirc,\Gamma_0}(\mathcal{O}_{X_{\mathbb{Q}_p(p^{1/p^\infty})^{\wedge\flat},A}})\ar[r]^{\mathrm{equi}}\ar[r]\ar[r] &\varphi\mathrm{Modules}_{\circledcirc,\Gamma_0}(\mathcal{O}_{Y_{\mathbb{Q}_p(p^{1/p^\infty})^{\wedge\flat},A}}).\\ 
}
\]

\item (\text{Proposition}) There is an equivalence between the $\infty$-categories of inductive Banach quasicoherent commutative algebra $E_\infty$ objects\footnote{Here $\circledcirc=\text{solidquasicoherentsheaves}$.}:
\[
\xymatrix@R+6pc@C+0pc{
\mathrm{sComm}_\mathrm{simplicial}\mathrm{Modules}_{\circledcirc,\Gamma_k}(\mathcal{O}_{X_{R_k,A}})\ar[d]\ar[d]\ar[d]\ar[d]\ar[r]^{\mathrm{equi}}\ar[r]\ar[r] &\mathrm{sComm}_\mathrm{simplicial}\varphi\mathrm{Modules}_{\circledcirc,\Gamma_k}(\mathcal{O}_{Y_{R_k,A}})\ar[d]\ar[d]\ar[d]\ar[d]\\
\mathrm{sComm}_\mathrm{simplicial}\mathrm{Modules}_{\circledcirc,\Gamma_0}(\mathcal{O}_{X_{\mathbb{Q}_p(p^{1/p^\infty})^{\wedge\flat},A}})\ar[r]^{\mathrm{equi}}\ar[r]\ar[r] &\mathrm{sComm}_\mathrm{simplicial}\varphi\mathrm{Modules}_{\circledcirc,\Gamma_0}(\mathcal{O}_{Y_{\mathbb{Q}_p(p^{1/p^\infty})^{\wedge\flat},A}}).  
}
\]

\item Then as in \cite{LBV} we have the equivalence of the de Rham complex after \cite[Definition 5.9, Section 5.2.1]{KKM}\footnote{Here $\circledcirc=\text{solidquasicoherentsheaves}$.}:
\[\displayindent=-0.4in
\xymatrix@R+6pc@C+0pc{
\mathrm{deRham}_{\mathrm{sComm}_\mathrm{simplicial}\mathrm{Modules}_{\circledcirc,\Gamma_k}(\mathcal{O}_{X_{R_k,A}})\ar[r]^{\mathrm{equi}}}(-)\ar[d]\ar[d]\ar[d]\ar[d]\ar[r]\ar[r] &\mathrm{deRham}_{\mathrm{sComm}_\mathrm{simplicial}\varphi\mathrm{Modules}_{\circledcirc,\Gamma_k}(\mathcal{O}_{Y_{R_k,A}})}(-)\ar[d]\ar[d]\ar[d]\ar[d]\\
\mathrm{deRham}_{\mathrm{sComm}_\mathrm{simplicial}\mathrm{Modules}_{\circledcirc,\Gamma_0}(\mathcal{O}_{X_{\mathbb{Q}_p(p^{1/p^\infty})^{\wedge\flat},A}})\ar[r]^{\mathrm{equi}}}(-)\ar[r]\ar[r] &\mathrm{deRham}_{\mathrm{sComm}_\mathrm{simplicial}\varphi\mathrm{Modules}_{\circledcirc,\Gamma_0}(\mathcal{O}_{Y_{\mathbb{Q}_p(p^{1/p^\infty})^{\wedge\flat},A}})}(-). 
}
\]

\item Then we have the following equivalence of $K$-group $(\infty,1)$-spectrum from \cite{BGT}\footnote{Here $\circledcirc=\text{solidquasicoherentsheaves}$.}:
\[
\xymatrix@R+6pc@C+0pc{
\mathrm{K}^\mathrm{BGT}_{\mathrm{sComm}_\mathrm{simplicial}\mathrm{Modules}_{\circledcirc,\Gamma_k}(\mathcal{O}_{X_{R_k,A}})\ar[r]^{\mathrm{equi}}}(-)\ar[d]\ar[d]\ar[d]\ar[d]\ar[r]\ar[r] &\mathrm{K}^\mathrm{BGT}_{\mathrm{sComm}_\mathrm{simplicial}\varphi\mathrm{Modules}_{\circledcirc,\Gamma_k}(\mathcal{O}_{Y_{R_k,A}})}(-)\ar[d]\ar[d]\ar[d]\ar[d]\\
\mathrm{K}^\mathrm{BGT}_{\mathrm{sComm}_\mathrm{simplicial}\mathrm{Modules}_{\circledcirc,\Gamma_0}(\mathcal{O}_{X_{\mathbb{Q}_p(p^{1/p^\infty})^{\wedge\flat},A}})\ar[r]^{\mathrm{equi}}}(-)\ar[r]\ar[r] &\mathrm{K}^\mathrm{BGT}_{\mathrm{sComm}_\mathrm{simplicial}\varphi\mathrm{Modules}_{\circledcirc,\Gamma_0}(\mathcal{O}_{Y_{\mathbb{Q}_p(p^{1/p^\infty})^{\wedge\flat},A}})}(-).
}
\]

\end{itemize}

\

Furthermore we have the corresponding pro-\'etale version without the corresponding fundamental group equivariances.

\begin{itemize}
\item (\text{Proposition}) There is an equivalence between the $\infty$-categories of inductive Banach quasicoherent presheaves\footnote{Here $\circledcirc=\text{solidquasicoherentsheaves}$.}:
\[
\xymatrix@R+6pc@C+0pc{
\mathrm{Modules}_\circledcirc(\mathcal{O}_{X_{\mathbb{Q}_p\left<T_1^{\pm 1},...,T_k^{\pm 1}\right>,\text{pro\'etale},A}})\ar[d]\ar[d]\ar[d]\ar[d] \ar[r]^{\mathrm{equi}}\ar[r]\ar[r] &\varphi\mathrm{Modules}_\circledcirc(\mathcal{O}_{Y_{\mathbb{Q}_p\left<T_1^{\pm 1},...,T_k^{\pm 1}\right>,\text{pro\'etale},A}}) \ar[d]\ar[d]\ar[d]\ar[d].\\
\mathrm{Modules}_\circledcirc(\mathcal{O}_{X_{\mathbb{Q}_p,\text{pro\'etale},A}})\ar[r]^{\mathrm{equi}}\ar[r]\ar[r] &\varphi\mathrm{Modules}_\circledcirc(\mathcal{O}_{Y_{\mathbb{Q}_p,\text{pro\'etale},A}}).\\ 
}
\]

\item (\text{Proposition}) There is an equivalence between the $\infty$-categories of inductive Banach quasicoherent commutative algebra $E_\infty$ objects\footnote{Here $\circledcirc=\text{solidquasicoherentsheaves}$.}:
\[\displayindent=-0.4in
\xymatrix@R+6pc@C+0pc{
\mathrm{sComm}_\mathrm{simplicial}\mathrm{Modules}_\circledcirc(\mathcal{O}_{X_{\mathbb{Q}_p\left<T_1^{\pm 1},...,T_k^{\pm 1}\right>,\text{pro\'etale},A}})\ar[d]\ar[d]\ar[d]\ar[d]\ar[r]^{\mathrm{equi}}\ar[r]\ar[r] &\mathrm{sComm}_\mathrm{simplicial}\varphi\mathrm{Modules}_\circledcirc(\mathcal{O}_{Y_{\mathbb{Q}_p\left<T_1^{\pm 1},...,T_k^{\pm 1}\right>,\text{pro\'etale},A}})\ar[d]\ar[d]\ar[d]\ar[d]\\
\mathrm{sComm}_\mathrm{simplicial}\mathrm{Modules}_\circledcirc(\mathcal{O}_{X_{\mathbb{Q}_p,\text{pro\'etale},A}})\ar[r]^{\mathrm{equi}}\ar[r]\ar[r] &\mathrm{sComm}_\mathrm{simplicial}\varphi\mathrm{Modules}_\circledcirc(\mathcal{O}_{Y_{\mathbb{Q}_p,\text{pro\'etale},A}}).  
}
\]

\item Then as in \cite{LBV} we have the equivalence of the de Rham complex after \cite[Definition 5.9, Section 5.2.1]{KKM}\footnote{Here $\circledcirc=\text{solidquasicoherentsheaves}$.}:
\[\displayindent=-0.9in
\xymatrix@R+6pc@C+0pc{
\mathrm{deRham}_{\mathrm{sComm}_\mathrm{simplicial}\mathrm{Modules}_\circledcirc(\mathcal{O}_{X_{\mathbb{Q}_p\left<T_1^{\pm 1},...,T_k^{\pm 1}\right>,\text{pro\'etale},A}})\ar[r]^{\mathrm{equi}}}(-)\ar[d]\ar[d]\ar[d]\ar[d]\ar[r]\ar[r] &\mathrm{deRham}_{\mathrm{sComm}_\mathrm{simplicial}\varphi\mathrm{Modules}_\circledcirc(\mathcal{O}_{Y_{\mathbb{Q}_p\left<T_1^{\pm 1},...,T_k^{\pm 1}\right>,\text{pro\'etale},A}})}(-)\ar[d]\ar[d]\ar[d]\ar[d]\\
\mathrm{deRham}_{\mathrm{sComm}_\mathrm{simplicial}\mathrm{Modules}_\circledcirc(\mathcal{O}_{X_{\mathbb{Q}_p,\text{pro\'etale},A}})\ar[r]^{\mathrm{equi}}}(-)\ar[r]\ar[r] &\mathrm{deRham}_{\mathrm{sComm}_\mathrm{simplicial}\varphi\mathrm{Modules}_\circledcirc(\mathcal{O}_{Y_{\mathbb{Q}_p,\text{pro\'etale},A}})}(-). 
}
\]

\item Then we have the following equivalence of $K$-group $(\infty,1)$-spectrum from \cite{BGT}\footnote{Here $\circledcirc=\text{solidquasicoherentsheaves}$.}:
\[\displayindent=-0.4in
\xymatrix@R+6pc@C+0pc{
\mathrm{K}^\mathrm{BGT}_{\mathrm{sComm}_\mathrm{simplicial}\mathrm{Modules}_\circledcirc(\mathcal{O}_{X_{\mathbb{Q}_p\left<T_1^{\pm 1},...,T_k^{\pm 1}\right>,\text{pro\'etale},A}})\ar[r]^{\mathrm{equi}}}(-)\ar[d]\ar[d]\ar[d]\ar[d]\ar[r]\ar[r] &\mathrm{K}^\mathrm{BGT}_{\mathrm{sComm}_\mathrm{simplicial}\varphi\mathrm{Modules}_\circledcirc(\mathcal{O}_{Y_{\mathbb{Q}_p\left<T_1^{\pm 1},...,T_k^{\pm 1}\right>,\text{pro\'etale},A}})}(-)\ar[d]\ar[d]\ar[d]\ar[d]\\
\mathrm{K}^\mathrm{BGT}_{\mathrm{sComm}_\mathrm{simplicial}\mathrm{Modules}_\circledcirc(\mathcal{O}_{X_{\mathbb{Q}_p,\text{pro\'etale},A}})\ar[r]^{\mathrm{equi}}}(-)\ar[r]\ar[r] &\mathrm{K}^\mathrm{BGT}_{\mathrm{sComm}_\mathrm{simplicial}\varphi\mathrm{Modules}_\circledcirc(\mathcal{O}_{Y_{\mathbb{Q}_p,\text{pro\'etale},A}})}(-). 
}
\]
	
\end{itemize}

\newpage
\subsection{$\infty$-Categorical Analytic Stacks and Descents II}

\indent Then by the corresponding \v{C}ech $\infty$-descent in \cite[Section 9.3]{KKM} and \cite{BBM} we have the following objects by directly taking the corresponding \v{C}ech $\infty$-descent. In the following the right had of each row in each diagram will be the corresponding quasicoherent Robba bundles over the Robba ring carrying the corresponding action from the Frobenius or the fundamental groups, defined by directly applying \cite[Section 9.3]{KKM} and \cite{BBM}. We then have the following global section functors:

\begin{itemize}

\item (\text{Proposition}) There is a functor (global section) between the $\infty$-categories of inductive Banach quasicoherent presheaves:
\[
\xymatrix@R+0pc@C+0pc{
\mathrm{Ind}\mathrm{Banach}(\mathcal{O}_{X_{R,A}})\ar[r]^{\mathrm{global}}\ar[r]\ar[r] &\varphi\mathrm{Ind}\mathrm{Banach}(\{\mathrm{Robba}^\mathrm{extended}_{{R,A,I}}\}_I).  
}
\]
\item (\text{Proposition}) There is a functor (global section) between the $\infty$-categories of monomorphic inductive Banach quasicoherent presheaves:
\[
\xymatrix@R+0pc@C+0pc{
\mathrm{Ind}^m\mathrm{Banach}(\mathcal{O}_{X_{R,A}})\ar[r]^{\mathrm{global}}\ar[r]\ar[r] &\varphi\mathrm{Ind}^m\mathrm{Banach}(\{\mathrm{Robba}^\mathrm{extended}_{{R,A,I}}\}_I).  
}
\]

\item (\text{Proposition}) There is a functor (global section) between the $\infty$-categories of inductive Banach quasicoherent presheaves:
\[
\xymatrix@R+0pc@C+0pc{
\mathrm{Ind}\mathrm{Banach}(\mathcal{O}_{X_{R,A}})\ar[r]^{\mathrm{global}}\ar[r]\ar[r] &\varphi\mathrm{Ind}\mathrm{Banach}(\{\mathrm{Robba}^\mathrm{extended}_{{R,A,I}}\}_I).  
}
\]
\item (\text{Proposition}) There is a functor (global section) between the $\infty$-categories of monomorphic inductive Banach quasicoherent presheaves:
\[
\xymatrix@R+0pc@C+0pc{
\mathrm{Ind}^m\mathrm{Banach}(\mathcal{O}_{X_{R,A}})\ar[r]^{\mathrm{global}}\ar[r]\ar[r] &\varphi\mathrm{Ind}^m\mathrm{Banach}(\{\mathrm{Robba}^\mathrm{extended}_{{R,A,I}}\}_I).  
}
\]
\item (\text{Proposition}) There is a functor (global section) between the $\infty$-categories of inductive Banach quasicoherent commutative algebra $E_\infty$ objects:
\[
\xymatrix@R+0pc@C+0pc{
\mathrm{sComm}_\mathrm{simplicial}\mathrm{Ind}\mathrm{Banach}(\mathcal{O}_{X_{R,A}})\ar[r]^{\mathrm{global}}\ar[r]\ar[r] &\mathrm{sComm}_\mathrm{simplicial}\varphi\mathrm{Ind}\mathrm{Banach}(\{\mathrm{Robba}^\mathrm{extended}_{{R,A,I}}\}_I).  
}
\]
\item (\text{Proposition}) There is a functor (global section) between the $\infty$-categories of monomorphic inductive Banach quasicoherent commutative algebra $E_\infty$ objects:
\[
\xymatrix@R+0pc@C+0pc{
\mathrm{sComm}_\mathrm{simplicial}\mathrm{Ind}^m\mathrm{Banach}(\mathcal{O}_{X_{R,A}})\ar[r]^{\mathrm{global}}\ar[r]\ar[r] &\mathrm{sComm}_\mathrm{simplicial}\varphi\mathrm{Ind}^m\mathrm{Banach}(\{\mathrm{Robba}^\mathrm{extended}_{{R,A,I}}\}_I).  
}
\]

\item Then parallel as in \cite{LBV} we have a functor (global section) of the de Rham complex after \cite[Definition 5.9, Section 5.2.1]{KKM}:
\[
\xymatrix@R+0pc@C+0pc{
\mathrm{deRham}_{\mathrm{sComm}_\mathrm{simplicial}\mathrm{Ind}\mathrm{Banach}(\mathcal{O}_{X_{R,A}})\ar[r]^{\mathrm{global}}}(-)\ar[r]\ar[r] &\mathrm{deRham}_{\mathrm{sComm}_\mathrm{simplicial}\varphi\mathrm{Ind}\mathrm{Banach}(\{\mathrm{Robba}^\mathrm{extended}_{{R,A,I}}\}_I)}(-), 
}
\]
\[
\xymatrix@R+0pc@C+0pc{
\mathrm{deRham}_{\mathrm{sComm}_\mathrm{simplicial}\mathrm{Ind}^m\mathrm{Banach}(\mathcal{O}_{X_{R,A}})\ar[r]^{\mathrm{global}}}(-)\ar[r]\ar[r] &\mathrm{deRham}_{\mathrm{sComm}_\mathrm{simplicial}\varphi\mathrm{Ind}^m\mathrm{Banach}(\{\mathrm{Robba}^\mathrm{extended}_{{R,A,I}}\}_I)}(-). 
}
\]

\item Then we have the following a functor (global section) of $K$-group $(\infty,1)$-spectrum from \cite{BGT}:
\[
\xymatrix@R+0pc@C+0pc{
\mathrm{K}^\mathrm{BGT}_{\mathrm{sComm}_\mathrm{simplicial}\mathrm{Ind}\mathrm{Banach}(\mathcal{O}_{X_{R,A}})\ar[r]^{\mathrm{global}}}(-)\ar[r]\ar[r] &\mathrm{K}^\mathrm{BGT}_{\mathrm{sComm}_\mathrm{simplicial}\varphi\mathrm{Ind}\mathrm{Banach}(\{\mathrm{Robba}^\mathrm{extended}_{{R,A,I}}\}_I)}(-), 
}
\]
\[
\xymatrix@R+0pc@C+0pc{
\mathrm{K}^\mathrm{BGT}_{\mathrm{sComm}_\mathrm{simplicial}\mathrm{Ind}^m\mathrm{Banach}(\mathcal{O}_{X_{R,A}})\ar[r]^{\mathrm{global}}}(-)\ar[r]\ar[r] &\mathrm{K}^\mathrm{BGT}_{\mathrm{sComm}_\mathrm{simplicial}\varphi\mathrm{Ind}^m\mathrm{Banach}(\{\mathrm{Robba}^\mathrm{extended}_{{R,A,I}}\}_I)}(-). 
}
\]
\end{itemize}

\noindent Now let $R=\mathbb{Q}_p(p^{1/p^\infty})^{\wedge\flat}$ and $R_k=\mathbb{Q}_p(p^{1/p^\infty})^{\wedge}\left<T_1^{\pm 1/p^{\infty}},...,T_k^{\pm 1/p^{\infty}}\right>^\flat$ we have the following Galois theoretic results with naturality along $f:\mathrm{Spa}(\mathbb{Q}_p(p^{1/p^\infty})^{\wedge}\left<T_1^{\pm 1/p^\infty},...,T_k^{\pm 1/p^\infty}\right>^\flat)\rightarrow \mathrm{Spa}(\mathbb{Q}_p(p^{1/p^\infty})^{\wedge\flat})$:

\begin{itemize}
\item (\text{Proposition}) There is a functor (global section) between the $\infty$-categories of inductive Banach quasicoherent presheaves:
\[
\xymatrix@R+6pc@C+0pc{
\mathrm{Ind}\mathrm{Banach}(\mathcal{O}_{X_{\mathbb{Q}_p(p^{1/p^\infty})^{\wedge}\left<T_1^{\pm 1/p^\infty},...,T_k^{\pm 1/p^\infty}\right>^\flat,A}})\ar[d]\ar[d]\ar[d]\ar[d] \ar[r]^{\mathrm{global}}\ar[r]\ar[r] &\varphi\mathrm{Ind}\mathrm{Banach}(\{\mathrm{Robba}^\mathrm{extended}_{{R_k,A,I}}\}_I) \ar[d]\ar[d]\ar[d]\ar[d].\\
\mathrm{Ind}\mathrm{Banach}(\mathcal{O}_{X_{\mathbb{Q}_p(p^{1/p^\infty})^{\wedge\flat},A}})\ar[r]^{\mathrm{global}}\ar[r]\ar[r] &\varphi\mathrm{Ind}\mathrm{Banach}(\{\mathrm{Robba}^\mathrm{extended}_{{R_0,A,I}}\}_I).\\ 
}
\]
\item (\text{Proposition}) There is a functor (global section) between the $\infty$-categories of monomorphic inductive Banach quasicoherent presheaves:
\[
\xymatrix@R+6pc@C+0pc{
\mathrm{Ind}^m\mathrm{Banach}(\mathcal{O}_{X_{R_k,A}})\ar[r]^{\mathrm{global}}\ar[d]\ar[d]\ar[d]\ar[d]\ar[r]\ar[r] &\varphi\mathrm{Ind}^m\mathrm{Banach}(\{\mathrm{Robba}^\mathrm{extended}_{{R_k,A,I}}\}_I)\ar[d]\ar[d]\ar[d]\ar[d]\\
\mathrm{Ind}^m\mathrm{Banach}(\mathcal{O}_{X_{\mathbb{Q}_p(p^{1/p^\infty})^{\wedge\flat},A}})\ar[r]^{\mathrm{global}}\ar[r]\ar[r] &\varphi\mathrm{Ind}^m\mathrm{Banach}(\{\mathrm{Robba}^\mathrm{extended}_{{R_0,A,I}}\}_I).\\  
}
\]
\item (\text{Proposition}) There is a functor (global section) between the $\infty$-categories of inductive Banach quasicoherent commutative algebra $E_\infty$ objects:
\[
\xymatrix@R+6pc@C+0pc{
\mathrm{sComm}_\mathrm{simplicial}\mathrm{Ind}\mathrm{Banach}(\mathcal{O}_{X_{R_k,A}})\ar[d]\ar[d]\ar[d]\ar[d]\ar[r]^{\mathrm{global}}\ar[r]\ar[r] &\mathrm{sComm}_\mathrm{simplicial}\varphi\mathrm{Ind}\mathrm{Banach}(\{\mathrm{Robba}^\mathrm{extended}_{{R_k,A,I}}\}_I)\ar[d]\ar[d]\ar[d]\ar[d]\\
\mathrm{sComm}_\mathrm{simplicial}\mathrm{Ind}\mathrm{Banach}(\mathcal{O}_{X_{\mathbb{Q}_p(p^{1/p^\infty})^{\wedge\flat},A}})\ar[r]^{\mathrm{global}}\ar[r]\ar[r] &\mathrm{sComm}_\mathrm{simplicial}\varphi\mathrm{Ind}\mathrm{Banach}(\{\mathrm{Robba}^\mathrm{extended}_{{R_0,A,I}}\}_I).  
}
\]
\item (\text{Proposition}) There is a functor (global section) between the $\infty$-categories of monomorphic inductive Banach quasicoherent commutative algebra $E_\infty$ objects:
\[
\xymatrix@R+6pc@C+0pc{
\mathrm{sComm}_\mathrm{simplicial}\mathrm{Ind}^m\mathrm{Banach}(\mathcal{O}_{X_{R_k,A}})\ar[d]\ar[d]\ar[d]\ar[d]\ar[r]^{\mathrm{global}}\ar[r]\ar[r] &\mathrm{sComm}_\mathrm{simplicial}\varphi\mathrm{Ind}^m\mathrm{Banach}(\{\mathrm{Robba}^\mathrm{extended}_{{R_k,A,I}}\}_I)\ar[d]\ar[d]\ar[d]\ar[d]\\
 \mathrm{sComm}_\mathrm{simplicial}\mathrm{Ind}^m\mathrm{Banach}(\mathcal{O}_{X_{\mathbb{Q}_p(p^{1/p^\infty})^{\wedge\flat},A}})\ar[r]^{\mathrm{global}}\ar[r]\ar[r] &\mathrm{sComm}_\mathrm{simplicial}\varphi\mathrm{Ind}^m\mathrm{Banach}(\{\mathrm{Robba}^\mathrm{extended}_{{R_0,A,I}}\}_I). 
}
\]

\item Then parallel as in \cite{LBV} we have a functor (global section) of the de Rham complex after \cite[Definition 5.9, Section 5.2.1]{KKM}:
\[\displayindent=-0.4in
\xymatrix@R+6pc@C+0pc{
\mathrm{deRham}_{\mathrm{sComm}_\mathrm{simplicial}\mathrm{Ind}\mathrm{Banach}(\mathcal{O}_{X_{R_k,A}})\ar[r]^{\mathrm{global}}}(-)\ar[d]\ar[d]\ar[d]\ar[d]\ar[r]\ar[r] &\mathrm{deRham}_{\mathrm{sComm}_\mathrm{simplicial}\varphi\mathrm{Ind}\mathrm{Banach}(\{\mathrm{Robba}^\mathrm{extended}_{{R_k,A,I}}\}_I)}(-)\ar[d]\ar[d]\ar[d]\ar[d]\\
\mathrm{deRham}_{\mathrm{sComm}_\mathrm{simplicial}\mathrm{Ind}\mathrm{Banach}(\mathcal{O}_{X_{\mathbb{Q}_p(p^{1/p^\infty})^{\wedge\flat},A}})\ar[r]^{\mathrm{global}}}(-)\ar[r]\ar[r] &\mathrm{deRham}_{\mathrm{sComm}_\mathrm{simplicial}\varphi\mathrm{Ind}\mathrm{Banach}(\{\mathrm{Robba}^\mathrm{extended}_{{R_0,A,I}}\}_I)}(-), 
}
\]
\[\displayindent=-0.4in
\xymatrix@R+6pc@C+0pc{
\mathrm{deRham}_{\mathrm{sComm}_\mathrm{simplicial}\mathrm{Ind}^m\mathrm{Banach}(\mathcal{O}_{X_{R_k,A}})\ar[r]^{\mathrm{global}}}(-)\ar[d]\ar[d]\ar[d]\ar[d]\ar[r]\ar[r] &\mathrm{deRham}_{\mathrm{sComm}_\mathrm{simplicial}\varphi\mathrm{Ind}^m\mathrm{Banach}(\{\mathrm{Robba}^\mathrm{extended}_{{R_k,A,I}}\}_I)}(-)\ar[d]\ar[d]\ar[d]\ar[d]\\
\mathrm{deRham}_{\mathrm{sComm}_\mathrm{simplicial}\mathrm{Ind}^m\mathrm{Banach}(\mathcal{O}_{X_{\mathbb{Q}_p(p^{1/p^\infty})^{\wedge\flat},A}})\ar[r]^{\mathrm{global}}}(-)\ar[r]\ar[r] &\mathrm{deRham}_{\mathrm{sComm}_\mathrm{simplicial}\varphi\mathrm{Ind}^m\mathrm{Banach}(\{\mathrm{Robba}^\mathrm{extended}_{{R_0,A,I}}\}_I)}(-). 
}
\]

\item Then we have the following a functor (global section) of $K$-group $(\infty,1)$-spectrum from \cite{BGT}:
\[
\xymatrix@R+6pc@C+0pc{
\mathrm{K}^\mathrm{BGT}_{\mathrm{sComm}_\mathrm{simplicial}\mathrm{Ind}\mathrm{Banach}(\mathcal{O}_{X_{R_k,A}})\ar[r]^{\mathrm{global}}}(-)\ar[d]\ar[d]\ar[d]\ar[d]\ar[r]\ar[r] &\mathrm{K}^\mathrm{BGT}_{\mathrm{sComm}_\mathrm{simplicial}\varphi\mathrm{Ind}\mathrm{Banach}(\{\mathrm{Robba}^\mathrm{extended}_{{R_k,A,I}}\}_I)}(-)\ar[d]\ar[d]\ar[d]\ar[d]\\
\mathrm{K}^\mathrm{BGT}_{\mathrm{sComm}_\mathrm{simplicial}\mathrm{Ind}\mathrm{Banach}(\mathcal{O}_{X_{\mathbb{Q}_p(p^{1/p^\infty})^{\wedge\flat},A}})\ar[r]^{\mathrm{global}}}(-)\ar[r]\ar[r] &\mathrm{K}^\mathrm{BGT}_{\mathrm{sComm}_\mathrm{simplicial}\varphi\mathrm{Ind}\mathrm{Banach}(\{\mathrm{Robba}^\mathrm{extended}_{{R_0,A,I}}\}_I)}(-), 
}
\]
\[
\xymatrix@R+6pc@C+0pc{
\mathrm{K}^\mathrm{BGT}_{\mathrm{sComm}_\mathrm{simplicial}\mathrm{Ind}^m\mathrm{Banach}(\mathcal{O}_{X_{R_k,A}})\ar[r]^{\mathrm{global}}}(-)\ar[d]\ar[d]\ar[d]\ar[d]\ar[r]\ar[r] &\mathrm{K}^\mathrm{BGT}_{\mathrm{sComm}_\mathrm{simplicial}\varphi\mathrm{Ind}^m\mathrm{Banach}(\{\mathrm{Robba}^\mathrm{extended}_{{R_k,A,I}}\}_I)}(-)\ar[d]\ar[d]\ar[d]\ar[d]\\
\mathrm{K}^\mathrm{BGT}_{\mathrm{sComm}_\mathrm{simplicial}\mathrm{Ind}^m\mathrm{Banach}(\mathcal{O}_{X_{\mathbb{Q}_p(p^{1/p^\infty})^{\wedge\flat},A}})\ar[r]^{\mathrm{global}}}(-)\ar[r]\ar[r] &\mathrm{K}^\mathrm{BGT}_{\mathrm{sComm}_\mathrm{simplicial}\varphi\mathrm{Ind}^m\mathrm{Banach}(\{\mathrm{Robba}^\mathrm{extended}_{{R_0,A,I}}\}_I)}(-). 
}
\]

\end{itemize}

\
\indent Then we consider further equivariance by considering the arithmetic profinite fundamental groups $\Gamma_{\mathbb{Q}_p}$ and $\mathrm{Gal}(\overline{\mathbb{Q}_p\left<T_1^{\pm 1},...,T_k^{\pm 1}\right>}/R_k)$ through the following diagram:

\[
\xymatrix@R+0pc@C+0pc{
\mathbb{Z}_p^k=\mathrm{Gal}(R_k/{\mathbb{Q}_p(p^{1/p^\infty})^\wedge\left<T_1^{\pm 1},...,T_k^{\pm 1}\right>}) \ar[r]\ar[r] \ar[r]\ar[r] &\Gamma_k:=\mathrm{Gal}(R_k/{\mathbb{Q}_p\left<T_1^{\pm 1},...,T_k^{\pm 1}\right>}) \ar[r] \ar[r]\ar[r] &\Gamma_{\mathbb{Q}_p}.
}
\]

\begin{itemize}
\item (\text{Proposition}) There is a functor (global section) between the $\infty$-categories of inductive Banach quasicoherent presheaves:
\[
\xymatrix@R+6pc@C+0pc{
\mathrm{Ind}\mathrm{Banach}_{\Gamma_k}(\mathcal{O}_{X_{\mathbb{Q}_p(p^{1/p^\infty})^{\wedge}\left<T_1^{\pm 1/p^\infty},...,T_k^{\pm 1/p^\infty}\right>^\flat,A}})\ar[d]\ar[d]\ar[d]\ar[d] \ar[r]^{\mathrm{global}}\ar[r]\ar[r] &\varphi\mathrm{Ind}\mathrm{Banach}_{\Gamma_k}(\{\mathrm{Robba}^\mathrm{extended}_{{R_k,A,I}}\}_I) \ar[d]\ar[d]\ar[d]\ar[d].\\
\mathrm{Ind}\mathrm{Banach}(\mathcal{O}_{X_{\mathbb{Q}_p(p^{1/p^\infty})^{\wedge\flat},A}})\ar[r]^{\mathrm{global}}\ar[r]\ar[r] &\varphi\mathrm{Ind}\mathrm{Banach}(\{\mathrm{Robba}^\mathrm{extended}_{{R_0,A,I}}\}_I).\\ 
}
\]
\item (\text{Proposition}) There is a functor (global section) between the $\infty$-categories of monomorphic inductive Banach quasicoherent presheaves:
\[
\xymatrix@R+6pc@C+0pc{
\mathrm{Ind}^m\mathrm{Banach}_{\Gamma_k}(\mathcal{O}_{X_{R_k,A}})\ar[r]^{\mathrm{global}}\ar[d]\ar[d]\ar[d]\ar[d]\ar[r]\ar[r] &\varphi\mathrm{Ind}^m\mathrm{Banach}_{\Gamma_k}(\{\mathrm{Robba}^\mathrm{extended}_{{R_k,A,I}}\}_I)\ar[d]\ar[d]\ar[d]\ar[d]\\
\mathrm{Ind}^m\mathrm{Banach}_{\Gamma_0}(\mathcal{O}_{X_{\mathbb{Q}_p(p^{1/p^\infty})^{\wedge\flat},A}})\ar[r]^{\mathrm{global}}\ar[r]\ar[r] &\varphi\mathrm{Ind}^m\mathrm{Banach}_{\Gamma_0}(\{\mathrm{Robba}^\mathrm{extended}_{{R_0,A,I}}\}_I).\\  
}
\]
\item (\text{Proposition}) There is a functor (global section) between the $\infty$-categories of inductive Banach quasicoherent commutative algebra $E_\infty$ objects:
\[\displayindent=-0.4in
\xymatrix@R+6pc@C+0pc{
\mathrm{sComm}_\mathrm{simplicial}\mathrm{Ind}\mathrm{Banach}_{\Gamma_k}(\mathcal{O}_{X_{R_k,A}})\ar[d]\ar[d]\ar[d]\ar[d]\ar[r]^{\mathrm{global}}\ar[r]\ar[r] &\mathrm{sComm}_\mathrm{simplicial}\varphi\mathrm{Ind}\mathrm{Banach}_{\Gamma_k}(\{\mathrm{Robba}^\mathrm{extended}_{{R_k,A,I}}\}_I)\ar[d]\ar[d]\ar[d]\ar[d]\\
\mathrm{sComm}_\mathrm{simplicial}\mathrm{Ind}\mathrm{Banach}_{\Gamma_0}(\mathcal{O}_{X_{\mathbb{Q}_p(p^{1/p^\infty})^{\wedge\flat},A}})\ar[r]^{\mathrm{global}}\ar[r]\ar[r] &\mathrm{sComm}_\mathrm{simplicial}\varphi\mathrm{Ind}\mathrm{Banach}_{\Gamma_0}(\{\mathrm{Robba}^\mathrm{extended}_{{R_0,A,I}}\}_I).  
}
\]
\item (\text{Proposition}) There is a functor (global section) between the $\infty$-categories of monomorphic inductive Banach quasicoherent commutative algebra $E_\infty$ objects:
\[\displayindent=-0.4in
\xymatrix@R+6pc@C+0pc{
\mathrm{sComm}_\mathrm{simplicial}\mathrm{Ind}^m\mathrm{Banach}_{\Gamma_k}(\mathcal{O}_{X_{R_k,A}})\ar[d]\ar[d]\ar[d]\ar[d]\ar[r]^{\mathrm{equi}}\ar[r]\ar[r] &\mathrm{sComm}_\mathrm{simplicial}\varphi\mathrm{Ind}^m\mathrm{Banach}_{\Gamma_k}(\{\mathrm{Robba}^\mathrm{extended}_{{R_k,A,I}}\}_I)\ar[d]\ar[d]\ar[d]\ar[d]\\
 \mathrm{sComm}_\mathrm{simplicial}\mathrm{Ind}^m\mathrm{Banach}_{\Gamma_0}(\mathcal{O}_{X_{\mathbb{Q}_p(p^{1/p^\infty})^{\wedge\flat},A}})\ar[r]^{\mathrm{equi}}\ar[r]\ar[r] &\mathrm{sComm}_\mathrm{simplicial}\varphi\mathrm{Ind}^m\mathrm{Banach}_{\Gamma_0}(\{\mathrm{Robba}^\mathrm{extended}_{{R_0,A,I}}\}_I). 
}
\]

\item Then parallel as in \cite{LBV} we have a functor (global section) of the de Rham complex after \cite[Definition 5.9, Section 5.2.1]{KKM}:
\[\displayindent=-0.4in
\xymatrix@R+6pc@C+0pc{
\mathrm{deRham}_{\mathrm{sComm}_\mathrm{simplicial}\mathrm{Ind}\mathrm{Banach}_{\Gamma_k}(\mathcal{O}_{X_{R_k,A}})\ar[r]^{\mathrm{global}}}(-)\ar[d]\ar[d]\ar[d]\ar[d]\ar[r]\ar[r] &\mathrm{deRham}_{\mathrm{sComm}_\mathrm{simplicial}\varphi\mathrm{Ind}\mathrm{Banach}_{\Gamma_k}(\{\mathrm{Robba}^\mathrm{extended}_{{R_k,A,I}}\}_I)}(-)\ar[d]\ar[d]\ar[d]\ar[d]\\
\mathrm{deRham}_{\mathrm{sComm}_\mathrm{simplicial}\mathrm{Ind}\mathrm{Banach}_{\Gamma_0}(\mathcal{O}_{X_{\mathbb{Q}_p(p^{1/p^\infty})^{\wedge\flat},A}})\ar[r]^{\mathrm{global}}}(-)\ar[r]\ar[r] &\mathrm{deRham}_{\mathrm{sComm}_\mathrm{simplicial}\varphi\mathrm{Ind}\mathrm{Banach}_{\Gamma_0}(\{\mathrm{Robba}^\mathrm{extended}_{{R_0,A,I}}\}_I)}(-), 
}
\]
\[\displayindent=-0.4in
\xymatrix@R+6pc@C+0pc{
\mathrm{deRham}_{\mathrm{sComm}_\mathrm{simplicial}\mathrm{Ind}^m\mathrm{Banach}_{\Gamma_k}(\mathcal{O}_{X_{R_k,A}})\ar[r]^{\mathrm{global}}}(-)\ar[d]\ar[d]\ar[d]\ar[d]\ar[r]\ar[r] &\mathrm{deRham}_{\mathrm{sComm}_\mathrm{simplicial}\varphi\mathrm{Ind}^m\mathrm{Banach}_{\Gamma_k}(\{\mathrm{Robba}^\mathrm{extended}_{{R_k,A,I}}\}_I)}(-)\ar[d]\ar[d]\ar[d]\ar[d]\\
\mathrm{deRham}_{\mathrm{sComm}_\mathrm{simplicial}\mathrm{Ind}^m\mathrm{Banach}_{\Gamma_0}(\mathcal{O}_{X_{\mathbb{Q}_p(p^{1/p^\infty})^{\wedge\flat},A}})\ar[r]^{\mathrm{global}}}(-)\ar[r]\ar[r] &\mathrm{deRham}_{\mathrm{sComm}_\mathrm{simplicial}\varphi\mathrm{Ind}^m\mathrm{Banach}_{\Gamma_0}(\{\mathrm{Robba}^\mathrm{extended}_{{R_0,A,I}}\}_I)}(-). 
}
\]

\item Then we have the following a functor (global section) of $K$-group $(\infty,1)$-spectrum from \cite{BGT}:
\[
\xymatrix@R+6pc@C+0pc{
\mathrm{K}^\mathrm{BGT}_{\mathrm{sComm}_\mathrm{simplicial}\mathrm{Ind}\mathrm{Banach}_{\Gamma_k}(\mathcal{O}_{X_{R_k,A}})\ar[r]^{\mathrm{global}}}(-)\ar[d]\ar[d]\ar[d]\ar[d]\ar[r]\ar[r] &\mathrm{K}^\mathrm{BGT}_{\mathrm{sComm}_\mathrm{simplicial}\varphi\mathrm{Ind}\mathrm{Banach}_{\Gamma_k}(\{\mathrm{Robba}^\mathrm{extended}_{{R_k,A,I}}\}_I)}(-)\ar[d]\ar[d]\ar[d]\ar[d]\\
\mathrm{K}^\mathrm{BGT}_{\mathrm{sComm}_\mathrm{simplicial}\mathrm{Ind}\mathrm{Banach}_{\Gamma_0}(\mathcal{O}_{X_{\mathbb{Q}_p(p^{1/p^\infty})^{\wedge\flat},A}})\ar[r]^{\mathrm{global}}}(-)\ar[r]\ar[r] &\mathrm{K}^\mathrm{BGT}_{\mathrm{sComm}_\mathrm{simplicial}\varphi\mathrm{Ind}\mathrm{Banach}_{\Gamma_0}(\{\mathrm{Robba}^\mathrm{extended}_{{R_0,A,I}}\}_I)}(-), 
}
\]
\[
\xymatrix@R+6pc@C+0pc{
\mathrm{K}^\mathrm{BGT}_{\mathrm{sComm}_\mathrm{simplicial}\mathrm{Ind}^m\mathrm{Banach}_{\Gamma_k}(\mathcal{O}_{X_{R_k,A}})\ar[r]^{\mathrm{global}}}(-)\ar[d]\ar[d]\ar[d]\ar[d]\ar[r]\ar[r] &\mathrm{K}^\mathrm{BGT}_{\mathrm{sComm}_\mathrm{simplicial}\varphi\mathrm{Ind}^m\mathrm{Banach}_{\Gamma_k}(\{\mathrm{Robba}^\mathrm{extended}_{{R_k,A,I}}\}_I)}(-)\ar[d]\ar[d]\ar[d]\ar[d]\\
\mathrm{K}^\mathrm{BGT}_{\mathrm{sComm}_\mathrm{simplicial}\mathrm{Ind}^m\mathrm{Banach}_{\Gamma_0}(\mathcal{O}_{X_{\mathbb{Q}_p(p^{1/p^\infty})^{\wedge\flat},A}})\ar[r]^{\mathrm{global}}}(-)\ar[r]\ar[r] &\mathrm{K}^\mathrm{BGT}_{\mathrm{sComm}_\mathrm{simplicial}\varphi\mathrm{Ind}^m\mathrm{Banach}_{\Gamma_0}(\{\mathrm{Robba}^\mathrm{extended}_{{R_0,A,I}}\}_I)}(-). 
}
\]

\end{itemize}

\

\

\begin{remark}
\noindent We can certainly consider the quasicoherent sheaves in \cite[Lemma 7.11, Remark 7.12]{1BBK}, therefore all the quasicoherent presheaves and modules will be those satisfying \cite[Lemma 7.11, Remark 7.12]{1BBK} if one would like to consider the the quasicoherent sheaves. That being all as this said, we would believe that the big quasicoherent presheaves are automatically quasicoherent sheaves (namely satisfying the corresponding \v{C}ech $\infty$-descent as in \cite[Section 9.3]{KKM} and \cite[Lemma 7.11, Remark 7.12]{1BBK}) and the corresponding global section functors are automatically equivalence of $\infty$-categories. 
\end{remark}

\

\indent In Clausen-Scholze formalism we have the following:

\begin{itemize}

\item (\text{Proposition}) There is a functor (global section) between the $\infty$-categories of inductive Banach quasicoherent sheaves:
\[
\xymatrix@R+0pc@C+0pc{
\mathrm{Module}_\circledcirc(\mathcal{O}_{X_{R,A}})\ar[r]^{\mathrm{global}}\ar[r]\ar[r] &\varphi\mathrm{Module}_\circledcirc(\{\mathrm{Robba}^\mathrm{extended}_{{R,A,I}}\}_I).  
}
\]

\item (\text{Proposition}) There is a functor (global section) between the $\infty$-categories of inductive Banach quasicoherent sheaves:
\[
\xymatrix@R+0pc@C+0pc{
\mathrm{Module}_\circledcirc(\mathcal{O}_{X_{R,A}})\ar[r]^{\mathrm{global}}\ar[r]\ar[r] &\varphi\mathrm{Module}_\circledcirc(\{\mathrm{Robba}^\mathrm{extended}_{{R,A,I}}\}_I).  
}
\]

\item (\text{Proposition}) There is a functor (global section) between the $\infty$-categories of inductive Banach quasicoherent commutative algebra $E_\infty$ objects\footnote{Here $\circledcirc=\text{solidquasicoherentsheaves}$.}:
\[
\xymatrix@R+0pc@C+0pc{
\mathrm{sComm}_\mathrm{simplicial}\mathrm{Module}_\circledcirc(\mathcal{O}_{X_{R,A}})\ar[r]^{\mathrm{global}}\ar[r]\ar[r] &\mathrm{sComm}_\mathrm{simplicial}\varphi\mathrm{Module}_\circledcirc(\{\mathrm{Robba}^\mathrm{extended}_{{R,A,I}}\}_I).  
}
\]

\item Then as in \cite{LBV} we have a functor (global section) of the de Rham complex after \cite[Definition 5.9, Section 5.2.1]{KKM}\footnote{Here $\circledcirc=\text{solidquasicoherentsheaves}$.}:
\[
\xymatrix@R+0pc@C+0pc{
\mathrm{deRham}_{\mathrm{sComm}_\mathrm{simplicial}\mathrm{Module}_\circledcirc(\mathcal{O}_{X_{R,A}})\ar[r]^{\mathrm{global}}}(-)\ar[r]\ar[r] &\mathrm{deRham}_{\mathrm{sComm}_\mathrm{simplicial}\varphi\mathrm{Module}_\circledcirc(\{\mathrm{Robba}^\mathrm{extended}_{{R,A,I}}\}_I)}(-). 
}
\]

\item Then we have the following a functor (global section) of $K$-group $(\infty,1)$-spectrum from \cite{BGT}\footnote{Here $\circledcirc=\text{solidquasicoherentsheaves}$.}:
\[
\xymatrix@R+0pc@C+0pc{
\mathrm{K}^\mathrm{BGT}_{\mathrm{sComm}_\mathrm{simplicial}\mathrm{Module}_\circledcirc(\mathcal{O}_{X_{R,A}})\ar[r]^{\mathrm{global}}}(-)\ar[r]\ar[r] &\mathrm{K}^\mathrm{BGT}_{\mathrm{sComm}_\mathrm{simplicial}\varphi\mathrm{Module}_\circledcirc(\{\mathrm{Robba}^\mathrm{extended}_{{R,A,I}}\}_I)}(-). 
}
\]
\end{itemize}

\noindent Now let $R=\mathbb{Q}_p(p^{1/p^\infty})^{\wedge\flat}$ and $R_k=\mathbb{Q}_p(p^{1/p^\infty})^{\wedge}\left<T_1^{\pm 1/p^{\infty}},...,T_k^{\pm 1/p^{\infty}}\right>^\flat$ we have the following Galois theoretic results with naturality along $f:\mathrm{Spa}(\mathbb{Q}_p(p^{1/p^\infty})^{\wedge}\left<T_1^{\pm 1/p^\infty},...,T_k^{\pm 1/p^\infty}\right>^\flat)\rightarrow \mathrm{Spa}(\mathbb{Q}_p(p^{1/p^\infty})^{\wedge\flat})$:

\begin{itemize}
\item (\text{Proposition}) There is a functor (global section) between the $\infty$-categories of inductive Banach quasicoherent sheaves\footnote{Here $\circledcirc=\text{solidquasicoherentsheaves}$.}:
\[
\xymatrix@R+6pc@C+0pc{
\mathrm{Module}_\circledcirc(\mathcal{O}_{X_{\mathbb{Q}_p(p^{1/p^\infty})^{\wedge}\left<T_1^{\pm 1/p^\infty},...,T_k^{\pm 1/p^\infty}\right>^\flat,A}})\ar[d]\ar[d]\ar[d]\ar[d] \ar[r]^{\mathrm{global}}\ar[r]\ar[r] &\varphi\mathrm{Module}_\circledcirc(\{\mathrm{Robba}^\mathrm{extended}_{{R_k,A,I}}\}_I) \ar[d]\ar[d]\ar[d]\ar[d].\\
\mathrm{Module}_\circledcirc(\mathcal{O}_{X_{\mathbb{Q}_p(p^{1/p^\infty})^{\wedge\flat},A}})\ar[r]^{\mathrm{global}}\ar[r]\ar[r] &\varphi\mathrm{Module}_\circledcirc(\{\mathrm{Robba}^\mathrm{extended}_{{R_0,A,I}}\}_I).\\ 
}
\]

\item (\text{Proposition}) There is a functor (global section) between the $\infty$-categories of inductive Banach quasicoherent commutative algebra $E_\infty$ objects\footnote{Here $\circledcirc=\text{solidquasicoherentsheaves}$.}:
\[
\xymatrix@R+6pc@C+0pc{
\mathrm{sComm}_\mathrm{simplicial}\mathrm{Module}_\circledcirc(\mathcal{O}_{X_{R_k,A}})\ar[d]\ar[d]\ar[d]\ar[d]\ar[r]^{\mathrm{global}}\ar[r]\ar[r] &\mathrm{sComm}_\mathrm{simplicial}\varphi\mathrm{Module}_\circledcirc(\{\mathrm{Robba}^\mathrm{extended}_{{R_k,A,I}}\}_I)\ar[d]\ar[d]\ar[d]\ar[d]\\
\mathrm{sComm}_\mathrm{simplicial}\mathrm{Module}_\circledcirc(\mathcal{O}_{X_{\mathbb{Q}_p(p^{1/p^\infty})^{\wedge\flat},A}})\ar[r]^{\mathrm{global}}\ar[r]\ar[r] &\mathrm{sComm}_\mathrm{simplicial}\varphi\mathrm{Module}_\circledcirc(\{\mathrm{Robba}^\mathrm{extended}_{{R_0,A,I}}\}_I).  
}
\]

\item Then as in \cite{LBV} we have a functor (global section) of the de Rham complex after \cite[Definition 5.9, Section 5.2.1]{KKM}\footnote{Here $\circledcirc=\text{solidquasicoherentsheaves}$.}:
\[
\xymatrix@R+6pc@C+0pc{
\mathrm{deRham}_{\mathrm{sComm}_\mathrm{simplicial}\mathrm{Module}_\circledcirc(\mathcal{O}_{X_{R_k,A}})\ar[r]^{\mathrm{global}}}(-)\ar[d]\ar[d]\ar[d]\ar[d]\ar[r]\ar[r] &\mathrm{deRham}_{\mathrm{sComm}_\mathrm{simplicial}\varphi\mathrm{Module}_\circledcirc(\{\mathrm{Robba}^\mathrm{extended}_{{R_k,A,I}}\}_I)}(-)\ar[d]\ar[d]\ar[d]\ar[d]\\
\mathrm{deRham}_{\mathrm{sComm}_\mathrm{simplicial}\mathrm{Module}_\circledcirc(\mathcal{O}_{X_{\mathbb{Q}_p(p^{1/p^\infty})^{\wedge\flat},A}})\ar[r]^{\mathrm{global}}}(-)\ar[r]\ar[r] &\mathrm{deRham}_{\mathrm{sComm}_\mathrm{simplicial}\varphi\mathrm{Module}_\circledcirc(\{\mathrm{Robba}^\mathrm{extended}_{{R_0,A,I}}\}_I)}(-). 
}
\]

\item Then we have the following a functor (global section) of $K$-group $(\infty,1)$-spectrum from \cite{BGT}\footnote{Here $\circledcirc=\text{solidquasicoherentsheaves}$.}:
\[
\xymatrix@R+6pc@C+0pc{
\mathrm{K}^\mathrm{BGT}_{\mathrm{sComm}_\mathrm{simplicial}\mathrm{Module}_\circledcirc(\mathcal{O}_{X_{R_k,A}})\ar[r]^{\mathrm{global}}}(-)\ar[d]\ar[d]\ar[d]\ar[d]\ar[r]\ar[r] &\mathrm{K}^\mathrm{BGT}_{\mathrm{sComm}_\mathrm{simplicial}\varphi\mathrm{Module}_\circledcirc(\{\mathrm{Robba}^\mathrm{extended}_{{R_k,A,I}}\}_I)}(-)\ar[d]\ar[d]\ar[d]\ar[d]\\
\mathrm{K}^\mathrm{BGT}_{\mathrm{sComm}_\mathrm{simplicial}\mathrm{Module}_\circledcirc(\mathcal{O}_{X_{\mathbb{Q}_p(p^{1/p^\infty})^{\wedge\flat},A}})\ar[r]^{\mathrm{global}}}(-)\ar[r]\ar[r] &\mathrm{K}^\mathrm{BGT}_{\mathrm{sComm}_\mathrm{simplicial}\varphi\mathrm{Module}_\circledcirc(\{\mathrm{Robba}^\mathrm{extended}_{{R_0,A,I}}\}_I)}(-). 
}
\]

\end{itemize}

\
\indent Then we consider further equivariance by considering the arithmetic profinite fundamental groups $\Gamma_{\mathbb{Q}_p}$ and $\mathrm{Gal}(\overline{\mathbb{Q}_p\left<T_1^{\pm 1},...,T_k^{\pm 1}\right>}/R_k)$ through the following diagram:

\[
\xymatrix@R+0pc@C+0pc{
\mathbb{Z}_p^k=\mathrm{Gal}(R_k/{\mathbb{Q}_p(p^{1/p^\infty})^\wedge\left<T_1^{\pm 1},...,T_k^{\pm 1}\right>}) \ar[r]\ar[r] \ar[r]\ar[r] &\Gamma_k:=\mathrm{Gal}(R_k/{\mathbb{Q}_p\left<T_1^{\pm 1},...,T_k^{\pm 1}\right>}) \ar[r] \ar[r]\ar[r] &\Gamma_{\mathbb{Q}_p}.
}
\]

\begin{itemize}
\item (\text{Proposition}) There is a functor (global section) between the $\infty$-categories of inductive Banach quasicoherent sheaves\footnote{Here $\circledcirc=\text{solidquasicoherentsheaves}$.}:
\[
\xymatrix@R+6pc@C+0pc{
{\mathrm{Module}_\circledcirc}_{\Gamma_k}(\mathcal{O}_{X_{\mathbb{Q}_p(p^{1/p^\infty})^{\wedge}\left<T_1^{\pm 1/p^\infty},...,T_k^{\pm 1/p^\infty}\right>^\flat,A}})\ar[d]\ar[d]\ar[d]\ar[d] \ar[r]^{\mathrm{global}}\ar[r]\ar[r] &\varphi{\mathrm{Module}_\circledcirc}_{\Gamma_k}(\{\mathrm{Robba}^\mathrm{extended}_{{R_k,A,I}}\}_I) \ar[d]\ar[d]\ar[d]\ar[d].\\
{\mathrm{Module}_\circledcirc}(\mathcal{O}_{X_{\mathbb{Q}_p(p^{1/p^\infty})^{\wedge\flat},A}})\ar[r]^{\mathrm{global}}\ar[r]\ar[r] &\varphi{\mathrm{Module}_\circledcirc}(\{\mathrm{Robba}^\mathrm{extended}_{{R_0,A,I}}\}_I).\\ 
}
\]

\item (\text{Proposition}) There is a functor (global section) between the $\infty$-categories of inductive Banach quasicoherent commutative algebra $E_\infty$ objects\footnote{Here $\circledcirc=\text{solidquasicoherentsheaves}$.}:
\[
\xymatrix@R+6pc@C+0pc{
\mathrm{sComm}_\mathrm{simplicial}{\mathrm{Module}_\circledcirc}_{\Gamma_k}(\mathcal{O}_{X_{R_k,A}})\ar[d]\ar[d]\ar[d]\ar[d]\ar[r]^{\mathrm{global}}\ar[r]\ar[r] &\mathrm{sComm}_\mathrm{simplicial}\varphi{\mathrm{Module}_\circledcirc}(\{\mathrm{Robba}^\mathrm{extended}_{{R_k,A,I}}\}_I)\ar[d]\ar[d]\ar[d]\ar[d]\\
\mathrm{sComm}_\mathrm{simplicial}{\mathrm{Module}_\circledcirc}_{\Gamma_0}(\mathcal{O}_{X_{\mathbb{Q}_p(p^{1/p^\infty})^{\wedge\flat},A}})\ar[r]^{\mathrm{global}}\ar[r]\ar[r] &\mathrm{sComm}_\mathrm{simplicial}\varphi{\mathrm{Modules}_\circledcirc}_{\Gamma_0}(\{\mathrm{Robba}^\mathrm{extended}_{{R_0,A,I}}\}_I).  
}
\]

\item Then as in \cite{LBV} we have a functor (global section) of the de Rham complex after \cite[Definition 5.9, Section 5.2.1]{KKM}\footnote{Here $\circledcirc=\text{solidquasicoherentsheaves}$.}:
\[\displayindent=-0.4in
\xymatrix@R+6pc@C+0pc{
\mathrm{deRham}_{\mathrm{sComm}_\mathrm{simplicial}{\mathrm{Modules}_\circledcirc}_{\Gamma_k}(\mathcal{O}_{X_{R_k,A}})\ar[r]^{\mathrm{global}}}(-)\ar[d]\ar[d]\ar[d]\ar[d]\ar[r]\ar[r] &\mathrm{deRham}_{\mathrm{sComm}_\mathrm{simplicial}\varphi{\mathrm{Modules}_\circledcirc}_{\Gamma_k}(\{\mathrm{Robba}^\mathrm{extended}_{{R_k,A,I}}\}_I)}(-)\ar[d]\ar[d]\ar[d]\ar[d]\\
\mathrm{deRham}_{\mathrm{sComm}_\mathrm{simplicial}{\mathrm{Modules}_\circledcirc}_{\Gamma_0}(\mathcal{O}_{X_{\mathbb{Q}_p(p^{1/p^\infty})^{\wedge\flat},A}})\ar[r]^{\mathrm{global}}}(-)\ar[r]\ar[r] &\mathrm{deRham}_{\mathrm{sComm}_\mathrm{simplicial}\varphi{\mathrm{Modules}_\circledcirc}_{\Gamma_0}(\{\mathrm{Robba}^\mathrm{extended}_{{R_0,A,I}}\}_I)}(-). 
}
\]

\item Then we have the following a functor (global section) of $K$-group $(\infty,1)$-spectrum from \cite{BGT}\footnote{Here $\circledcirc=\text{solidquasicoherentsheaves}$.}:
\[
\xymatrix@R+6pc@C+0pc{
\mathrm{K}^\mathrm{BGT}_{\mathrm{sComm}_\mathrm{simplicial}{\mathrm{Modules}_\circledcirc}_{\Gamma_k}(\mathcal{O}_{X_{R_k,A}})\ar[r]^{\mathrm{global}}}(-)\ar[d]\ar[d]\ar[d]\ar[d]\ar[r]\ar[r] &\mathrm{K}^\mathrm{BGT}_{\mathrm{sComm}_\mathrm{simplicial}\varphi{\mathrm{Modules}_\circledcirc}_{\Gamma_k}(\{\mathrm{Robba}^\mathrm{extended}_{{R_k,A,I}}\}_I)}(-)\ar[d]\ar[d]\ar[d]\ar[d]\\
\mathrm{K}^\mathrm{BGT}_{\mathrm{sComm}_\mathrm{simplicial}{\mathrm{Modules}_\circledcirc}_{\Gamma_0}(\mathcal{O}_{X_{\mathbb{Q}_p(p^{1/p^\infty})^{\wedge\flat},A}})\ar[r]^{\mathrm{global}}}(-)\ar[r]\ar[r] &\mathrm{K}^\mathrm{BGT}_{\mathrm{sComm}_\mathrm{simplicial}\varphi{\mathrm{Modules}_\circledcirc}_{\Gamma_0}(\{\mathrm{Robba}^\mathrm{extended}_{{R_0,A,I}}\}_I)}(-). 
}
\]
	
\end{itemize}

\begin{proposition}
All the global functors from \cite[Proposition 13.8, Theorem 14.9, Remark 14.10]{1CS2} achieve the equivalences on both sides.	
\end{proposition}

\newpage
\subsection{$\infty$-Categorical Analytic Stacks and Descents III}

\indent In the following the right had of each row in each diagram will be the corresponding quasicoherent Robba bundles over the Robba ring carrying the corresponding action from the Frobenius or the fundamental groups, defined by directly applying \cite[Section 9.3]{KKM} and \cite{BBM}. We now let $\mathcal{A}$ be any commutative algebra objects in the corresponding $\infty$-toposes over ind-Banach commutative algebra objects over $\mathbb{Q}_p$ or the corresponding born\'e commutative algebra objects over $\mathbb{Q}_p$ carrying the Grothendieck topology defined by essentially the corresponding monomorphism homotopy in the opposite category. Then we promote the construction to the corresponding $\infty$-stack over the same $\infty$-categories of affinoids. Let $\mathcal{A}$ vary in the category of all the Banach algebras over $\mathbb{Q}_p$ we have the following.

\begin{itemize}

\item (\text{Proposition}) There is a functor (global section) between the $\infty$-prestacks of inductive Banach quasicoherent presheaves:
\[
\xymatrix@R+0pc@C+0pc{
\mathrm{Ind}\mathrm{Banach}(\mathcal{O}_{X_{R,-}})\ar[r]^{\mathrm{global}}\ar[r]\ar[r] &\varphi\mathrm{Ind}\mathrm{Banach}(\{\mathrm{Robba}^\mathrm{extended}_{{R,-,I}}\}_I).  
}
\]
\item (\text{Proposition}) There is a functor (global section) between the $\infty$-prestacks of monomorphic inductive Banach quasicoherent presheaves:
\[
\xymatrix@R+0pc@C+0pc{
\mathrm{Ind}^m\mathrm{Banach}(\mathcal{O}_{X_{R,-}})\ar[r]^{\mathrm{global}}\ar[r]\ar[r] &\varphi\mathrm{Ind}^m\mathrm{Banach}(\{\mathrm{Robba}^\mathrm{extended}_{{R,-,I}}\}_I).  
}
\]

\item (\text{Proposition}) There is a functor (global section) between the $\infty$-prestacks of inductive Banach quasicoherent presheaves:
\[
\xymatrix@R+0pc@C+0pc{
\mathrm{Ind}\mathrm{Banach}(\mathcal{O}_{X_{R,-}})\ar[r]^{\mathrm{global}}\ar[r]\ar[r] &\varphi\mathrm{Ind}\mathrm{Banach}(\{\mathrm{Robba}^\mathrm{extended}_{{R,-,I}}\}_I).  
}
\]
\item (\text{Proposition}) There is a functor (global section) between the $\infty$-stacks of monomorphic inductive Banach quasicoherent presheaves:
\[
\xymatrix@R+0pc@C+0pc{
\mathrm{Ind}^m\mathrm{Banach}(\mathcal{O}_{X_{R,-}})\ar[r]^{\mathrm{global}}\ar[r]\ar[r] &\varphi\mathrm{Ind}^m\mathrm{Banach}(\{\mathrm{Robba}^\mathrm{extended}_{{R,-,I}}\}_I).  
}
\]
\item (\text{Proposition}) There is a functor (global section) between the $\infty$-prestacks of inductive Banach quasicoherent commutative algebra $E_\infty$ objects:
\[
\xymatrix@R+0pc@C+0pc{
\mathrm{sComm}_\mathrm{simplicial}\mathrm{Ind}\mathrm{Banach}(\mathcal{O}_{X_{R,-}})\ar[r]^{\mathrm{global}}\ar[r]\ar[r] &\mathrm{sComm}_\mathrm{simplicial}\varphi\mathrm{Ind}\mathrm{Banach}(\{\mathrm{Robba}^\mathrm{extended}_{{R,-,I}}\}_I).  
}
\]
\item (\text{Proposition}) There is a functor (global section) between the $\infty$-prestacks of monomorphic inductive Banach quasicoherent commutative algebra $E_\infty$ objects:
\[
\xymatrix@R+0pc@C+0pc{
\mathrm{sComm}_\mathrm{simplicial}\mathrm{Ind}^m\mathrm{Banach}(\mathcal{O}_{X_{R,-}})\ar[r]^{\mathrm{global}}\ar[r]\ar[r] &\mathrm{sComm}_\mathrm{simplicial}\varphi\mathrm{Ind}^m\mathrm{Banach}(\{\mathrm{Robba}^\mathrm{extended}_{{R,-,I}}\}_I).  
}
\]

\item Then parallel as in \cite{LBV} we have a functor (global section) of the de Rham complex after \cite[Definition 5.9, Section 5.2.1]{KKM}:
\[
\xymatrix@R+0pc@C+0pc{
\mathrm{deRham}_{\mathrm{sComm}_\mathrm{simplicial}\mathrm{Ind}\mathrm{Banach}(\mathcal{O}_{X_{R,-}})\ar[r]^{\mathrm{global}}}(-)\ar[r]\ar[r] &\mathrm{deRham}_{\mathrm{sComm}_\mathrm{simplicial}\varphi\mathrm{Ind}\mathrm{Banach}(\{\mathrm{Robba}^\mathrm{extended}_{{R,-,I}}\}_I)}(-), 
}
\]
\[
\xymatrix@R+0pc@C+0pc{
\mathrm{deRham}_{\mathrm{sComm}_\mathrm{simplicial}\mathrm{Ind}^m\mathrm{Banach}(\mathcal{O}_{X_{R,-}})\ar[r]^{\mathrm{global}}}(-)\ar[r]\ar[r] &\mathrm{deRham}_{\mathrm{sComm}_\mathrm{simplicial}\varphi\mathrm{Ind}^m\mathrm{Banach}(\{\mathrm{Robba}^\mathrm{extended}_{{R,-,I}}\}_I)}(-). 
}
\]

\item Then we have the following a functor (global section) of $K$-group $(\infty,1)$-spectrum from \cite{BGT}:
\[
\xymatrix@R+0pc@C+0pc{
\mathrm{K}^\mathrm{BGT}_{\mathrm{sComm}_\mathrm{simplicial}\mathrm{Ind}\mathrm{Banach}(\mathcal{O}_{X_{R,-}})\ar[r]^{\mathrm{global}}}(-)\ar[r]\ar[r] &\mathrm{K}^\mathrm{BGT}_{\mathrm{sComm}_\mathrm{simplicial}\varphi\mathrm{Ind}\mathrm{Banach}(\{\mathrm{Robba}^\mathrm{extended}_{{R,-,I}}\}_I)}(-), 
}
\]
\[
\xymatrix@R+0pc@C+0pc{
\mathrm{K}^\mathrm{BGT}_{\mathrm{sComm}_\mathrm{simplicial}\mathrm{Ind}^m\mathrm{Banach}(\mathcal{O}_{X_{R,-}})\ar[r]^{\mathrm{global}}}(-)\ar[r]\ar[r] &\mathrm{K}^\mathrm{BGT}_{\mathrm{sComm}_\mathrm{simplicial}\varphi\mathrm{Ind}^m\mathrm{Banach}(\{\mathrm{Robba}^\mathrm{extended}_{{R,-,I}}\}_I)}(-). 
}
\]
\end{itemize}

\noindent Now let $R=\mathbb{Q}_p(p^{1/p^\infty})^{\wedge\flat}$ and $R_k=\mathbb{Q}_p(p^{1/p^\infty})^{\wedge}\left<T_1^{\pm 1/p^{\infty}},...,T_k^{\pm 1/p^{\infty}}\right>^\flat$ we have the following Galois theoretic results with naturality along $f:\mathrm{Spa}(\mathbb{Q}_p(p^{1/p^\infty})^{\wedge}\left<T_1^{\pm 1/p^\infty},...,T_k^{\pm 1/p^\infty}\right>^\flat)\rightarrow \mathrm{Spa}(\mathbb{Q}_p(p^{1/p^\infty})^{\wedge\flat})$:

\begin{itemize}
\item (\text{Proposition}) There is a functor (global section) between the $\infty$-prestacks of inductive Banach quasicoherent presheaves:
\[
\xymatrix@R+6pc@C+0pc{
\mathrm{Ind}\mathrm{Banach}(\mathcal{O}_{X_{\mathbb{Q}_p(p^{1/p^\infty})^{\wedge}\left<T_1^{\pm 1/p^\infty},...,T_k^{\pm 1/p^\infty}\right>^\flat,-}})\ar[d]\ar[d]\ar[d]\ar[d] \ar[r]^{\mathrm{global}}\ar[r]\ar[r] &\varphi\mathrm{Ind}\mathrm{Banach}(\{\mathrm{Robba}^\mathrm{extended}_{{R_k,-,I}}\}_I) \ar[d]\ar[d]\ar[d]\ar[d].\\
\mathrm{Ind}\mathrm{Banach}(\mathcal{O}_{X_{\mathbb{Q}_p(p^{1/p^\infty})^{\wedge\flat},-}})\ar[r]^{\mathrm{global}}\ar[r]\ar[r] &\varphi\mathrm{Ind}\mathrm{Banach}(\{\mathrm{Robba}^\mathrm{extended}_{{R_0,-,I}}\}_I).\\ 
}
\]
\item (\text{Proposition}) There is a functor (global section) between the $\infty$-prestacks of monomorphic inductive Banach quasicoherent presheaves:
\[
\xymatrix@R+6pc@C+0pc{
\mathrm{Ind}^m\mathrm{Banach}(\mathcal{O}_{X_{R_k,-}})\ar[r]^{\mathrm{global}}\ar[d]\ar[d]\ar[d]\ar[d]\ar[r]\ar[r] &\varphi\mathrm{Ind}^m\mathrm{Banach}(\{\mathrm{Robba}^\mathrm{extended}_{{R_k,-,I}}\}_I)\ar[d]\ar[d]\ar[d]\ar[d]\\
\mathrm{Ind}^m\mathrm{Banach}(\mathcal{O}_{X_{\mathbb{Q}_p(p^{1/p^\infty})^{\wedge\flat},-}})\ar[r]^{\mathrm{global}}\ar[r]\ar[r] &\varphi\mathrm{Ind}^m\mathrm{Banach}(\{\mathrm{Robba}^\mathrm{extended}_{{R_0,-,I}}\}_I).\\  
}
\]
\item (\text{Proposition}) There is a functor (global section) between the $\infty$-prestacks of inductive Banach quasicoherent commutative algebra $E_\infty$ objects:
\[
\xymatrix@R+6pc@C+0pc{
\mathrm{sComm}_\mathrm{simplicial}\mathrm{Ind}\mathrm{Banach}(\mathcal{O}_{X_{R_k,-}})\ar[d]\ar[d]\ar[d]\ar[d]\ar[r]^{\mathrm{global}}\ar[r]\ar[r] &\mathrm{sComm}_\mathrm{simplicial}\varphi\mathrm{Ind}\mathrm{Banach}(\{\mathrm{Robba}^\mathrm{extended}_{{R_k,-,I}}\}_I)\ar[d]\ar[d]\ar[d]\ar[d]\\
\mathrm{sComm}_\mathrm{simplicial}\mathrm{Ind}\mathrm{Banach}(\mathcal{O}_{X_{\mathbb{Q}_p(p^{1/p^\infty})^{\wedge\flat},-}})\ar[r]^{\mathrm{global}}\ar[r]\ar[r] &\mathrm{sComm}_\mathrm{simplicial}\varphi\mathrm{Ind}\mathrm{Banach}(\{\mathrm{Robba}^\mathrm{extended}_{{R_0,-,I}}\}_I).  
}
\]
\item (\text{Proposition}) There is a functor (global section) between the $\infty$-prestacks of monomorphic inductive Banach quasicoherent commutative algebra $E_\infty$ objects:
\[
\xymatrix@R+6pc@C+0pc{
\mathrm{sComm}_\mathrm{simplicial}\mathrm{Ind}^m\mathrm{Banach}(\mathcal{O}_{X_{R_k,-}})\ar[d]\ar[d]\ar[d]\ar[d]\ar[r]^{\mathrm{global}}\ar[r]\ar[r] &\mathrm{sComm}_\mathrm{simplicial}\varphi\mathrm{Ind}^m\mathrm{Banach}(\{\mathrm{Robba}^\mathrm{extended}_{{R_k,-,I}}\}_I)\ar[d]\ar[d]\ar[d]\ar[d]\\
 \mathrm{sComm}_\mathrm{simplicial}\mathrm{Ind}^m\mathrm{Banach}(\mathcal{O}_{X_{\mathbb{Q}_p(p^{1/p^\infty})^{\wedge\flat},-}})\ar[r]^{\mathrm{global}}\ar[r]\ar[r] &\mathrm{sComm}_\mathrm{simplicial}\varphi\mathrm{Ind}^m\mathrm{Banach}(\{\mathrm{Robba}^\mathrm{extended}_{{R_0,-,I}}\}_I).
}
\]

\item Then parallel as in \cite{LBV} we have a functor (global section) of the de Rham complex after \cite[Definition 5.9, Section 5.2.1]{KKM}:
\[\displayindent=-0.4in
\xymatrix@R+6pc@C+0pc{
\mathrm{deRham}_{\mathrm{sComm}_\mathrm{simplicial}\mathrm{Ind}\mathrm{Banach}(\mathcal{O}_{X_{R_k,-}})\ar[r]^{\mathrm{global}}}(-)\ar[d]\ar[d]\ar[d]\ar[d]\ar[r]\ar[r] &\mathrm{deRham}_{\mathrm{sComm}_\mathrm{simplicial}\varphi\mathrm{Ind}\mathrm{Banach}(\{\mathrm{Robba}^\mathrm{extended}_{{R_k,-,I}}\}_I)}(-)\ar[d]\ar[d]\ar[d]\ar[d]\\
\mathrm{deRham}_{\mathrm{sComm}_\mathrm{simplicial}\mathrm{Ind}\mathrm{Banach}(\mathcal{O}_{X_{\mathbb{Q}_p(p^{1/p^\infty})^{\wedge\flat},-}})\ar[r]^{\mathrm{global}}}(-)\ar[r]\ar[r] &\mathrm{deRham}_{\mathrm{sComm}_\mathrm{simplicial}\varphi\mathrm{Ind}\mathrm{Banach}(\{\mathrm{Robba}^\mathrm{extended}_{{R_0,-,I}}\}_I)}(-), 
}
\]
\[\displayindent=-0.4in
\xymatrix@R+6pc@C+0pc{
\mathrm{deRham}_{\mathrm{sComm}_\mathrm{simplicial}\mathrm{Ind}^m\mathrm{Banach}(\mathcal{O}_{X_{R_k,-}})\ar[r]^{\mathrm{global}}}(-)\ar[d]\ar[d]\ar[d]\ar[d]\ar[r]\ar[r] &\mathrm{deRham}_{\mathrm{sComm}_\mathrm{simplicial}\varphi\mathrm{Ind}^m\mathrm{Banach}(\{\mathrm{Robba}^\mathrm{extended}_{{R_k,-,I}}\}_I)}(-)\ar[d]\ar[d]\ar[d]\ar[d]\\
\mathrm{deRham}_{\mathrm{sComm}_\mathrm{simplicial}\mathrm{Ind}^m\mathrm{Banach}(\mathcal{O}_{X_{\mathbb{Q}_p(p^{1/p^\infty})^{\wedge\flat},-}})\ar[r]^{\mathrm{global}}}(-)\ar[r]\ar[r] &\mathrm{deRham}_{\mathrm{sComm}_\mathrm{simplicial}\varphi\mathrm{Ind}^m\mathrm{Banach}(\{\mathrm{Robba}^\mathrm{extended}_{{R_0,-,I}}\}_I)}(-). 
}
\]

\item Then we have the following a functor (global section) of $K$-group $(\infty,1)$-spectrum from \cite{BGT}:
\[
\xymatrix@R+6pc@C+0pc{
\mathrm{K}^\mathrm{BGT}_{\mathrm{sComm}_\mathrm{simplicial}\mathrm{Ind}\mathrm{Banach}(\mathcal{O}_{X_{R_k,-}})\ar[r]^{\mathrm{global}}}(-)\ar[d]\ar[d]\ar[d]\ar[d]\ar[r]\ar[r] &\mathrm{K}^\mathrm{BGT}_{\mathrm{sComm}_\mathrm{simplicial}\varphi\mathrm{Ind}\mathrm{Banach}(\{\mathrm{Robba}^\mathrm{extended}_{{R_k,-,I}}\}_I)}(-)\ar[d]\ar[d]\ar[d]\ar[d]\\
\mathrm{K}^\mathrm{BGT}_{\mathrm{sComm}_\mathrm{simplicial}\mathrm{Ind}\mathrm{Banach}(\mathcal{O}_{X_{\mathbb{Q}_p(p^{1/p^\infty})^{\wedge\flat},-}})\ar[r]^{\mathrm{global}}}(-)\ar[r]\ar[r] &\mathrm{K}^\mathrm{BGT}_{\mathrm{sComm}_\mathrm{simplicial}\varphi\mathrm{Ind}\mathrm{Banach}(\{\mathrm{Robba}^\mathrm{extended}_{{R_0,-,I}}\}_I)}(-), 
}
\]
\[
\xymatrix@R+6pc@C+0pc{
\mathrm{K}^\mathrm{BGT}_{\mathrm{sComm}_\mathrm{simplicial}\mathrm{Ind}^m\mathrm{Banach}(\mathcal{O}_{X_{R_k,-}})\ar[r]^{\mathrm{global}}}(-)\ar[d]\ar[d]\ar[d]\ar[d]\ar[r]\ar[r] &\mathrm{K}^\mathrm{BGT}_{\mathrm{sComm}_\mathrm{simplicial}\varphi\mathrm{Ind}^m\mathrm{Banach}(\{\mathrm{Robba}^\mathrm{extended}_{{R_k,-,I}}\}_I)}(-)\ar[d]\ar[d]\ar[d]\ar[d]\\
\mathrm{K}^\mathrm{BGT}_{\mathrm{sComm}_\mathrm{simplicial}\mathrm{Ind}^m\mathrm{Banach}(\mathcal{O}_{X_{\mathbb{Q}_p(p^{1/p^\infty})^{\wedge\flat},-}})\ar[r]^{\mathrm{global}}}(-)\ar[r]\ar[r] &\mathrm{K}^\mathrm{BGT}_{\mathrm{sComm}_\mathrm{simplicial}\varphi\mathrm{Ind}^m\mathrm{Banach}(\{\mathrm{Robba}^\mathrm{extended}_{{R_0,-,I}}\}_I)}(-). 
}
\]

\end{itemize}

\
\indent Then we consider further equivariance by considering the arithmetic profinite fundamental groups $\Gamma_{\mathbb{Q}_p}$ and $\mathrm{Gal}(\overline{{Q}_p\left<T_1^{\pm 1},...,T_k^{\pm 1}\right>}/R_k)$ through the following diagram:

\[
\xymatrix@R+0pc@C+0pc{
\mathbb{Z}_p^k=\mathrm{Gal}(R_k/{\mathbb{Q}_p(p^{1/p^\infty})^\wedge\left<T_1^{\pm 1},...,T_k^{\pm 1}\right>}) \ar[r]\ar[r] \ar[r]\ar[r] &\Gamma_k:=\mathrm{Gal}(R_k/{\mathbb{Q}_p\left<T_1^{\pm 1},...,T_k^{\pm 1}\right>}) \ar[r] \ar[r]\ar[r] &\Gamma_{\mathbb{Q}_p}.
}
\]

\begin{itemize}
\item (\text{Proposition}) There is a functor (global section) between the $\infty$-prestacks of inductive Banach quasicoherent presheaves:
\[
\xymatrix@R+6pc@C+0pc{
\mathrm{Ind}\mathrm{Banach}_{\Gamma_k}(\mathcal{O}_{X_{\mathbb{Q}_p(p^{1/p^\infty})^{\wedge}\left<T_1^{\pm 1/p^\infty},...,T_k^{\pm 1/p^\infty}\right>^\flat,-}})\ar[d]\ar[d]\ar[d]\ar[d] \ar[r]^{\mathrm{global}}\ar[r]\ar[r] &\varphi\mathrm{Ind}\mathrm{Banach}_{\Gamma_k}(\{\mathrm{Robba}^\mathrm{extended}_{{R_k,-,I}}\}_I) \ar[d]\ar[d]\ar[d]\ar[d].\\
\mathrm{Ind}\mathrm{Banach}(\mathcal{O}_{X_{\mathbb{Q}_p(p^{1/p^\infty})^{\wedge\flat},-}})\ar[r]^{\mathrm{global}}\ar[r]\ar[r] &\varphi\mathrm{Ind}\mathrm{Banach}(\{\mathrm{Robba}^\mathrm{extended}_{{R_0,-,I}}\}_I).\\ 
}
\]
\item (\text{Proposition}) There is a functor (global section) between the $\infty$-prestacks of monomorphic inductive Banach quasicoherent presheaves:
\[
\xymatrix@R+6pc@C+0pc{
\mathrm{Ind}^m\mathrm{Banach}_{\Gamma_k}(\mathcal{O}_{X_{R_k,-}})\ar[r]^{\mathrm{global}}\ar[d]\ar[d]\ar[d]\ar[d]\ar[r]\ar[r] &\varphi\mathrm{Ind}^m\mathrm{Banach}_{\Gamma_k}(\{\mathrm{Robba}^\mathrm{extended}_{{R_k,-,I}}\}_I)\ar[d]\ar[d]\ar[d]\ar[d]\\
\mathrm{Ind}^m\mathrm{Banach}_{\Gamma_0}(\mathcal{O}_{X_{\mathbb{Q}_p(p^{1/p^\infty})^{\wedge\flat},-}})\ar[r]^{\mathrm{global}}\ar[r]\ar[r] &\varphi\mathrm{Ind}^m\mathrm{Banach}_{\Gamma_0}(\{\mathrm{Robba}^\mathrm{extended}_{{R_0,-,I}}\}_I).\\  
}
\]
\item (\text{Proposition}) There is a functor (global section) between the $\infty$-stacks of inductive Banach quasicoherent commutative algebra $E_\infty$ objects:
\[
\xymatrix@R+6pc@C+0pc{
\mathrm{sComm}_\mathrm{simplicial}\mathrm{Ind}\mathrm{Banach}_{\Gamma_k}(\mathcal{O}_{X_{R_k,-}})\ar[d]\ar[d]\ar[d]\ar[d]\ar[r]^{\mathrm{global}}\ar[r]\ar[r] &\mathrm{sComm}_\mathrm{simplicial}\varphi\mathrm{Ind}\mathrm{Banach}_{\Gamma_k}(\{\mathrm{Robba}^\mathrm{extended}_{{R_k,-,I}}\}_I)\ar[d]\ar[d]\ar[d]\ar[d]\\
\mathrm{sComm}_\mathrm{simplicial}\mathrm{Ind}\mathrm{Banach}_{\Gamma_0}(\mathcal{O}_{X_{\mathbb{Q}_p(p^{1/p^\infty})^{\wedge\flat},-}})\ar[r]^{\mathrm{global}}\ar[r]\ar[r] &\mathrm{sComm}_\mathrm{simplicial}\varphi\mathrm{Ind}\mathrm{Banach}_{\Gamma_0}(\{\mathrm{Robba}^\mathrm{extended}_{{R_0,-,I}}\}_I).  
}
\]
\item (\text{Proposition}) There is a functor (global section) between the $\infty$-prestacks of monomorphic inductive Banach quasicoherent commutative algebra $E_\infty$ objects:
\[\displayindent=-0.4in
\xymatrix@R+6pc@C+0pc{
\mathrm{sComm}_\mathrm{simplicial}\mathrm{Ind}^m\mathrm{Banach}_{\Gamma_k}(\mathcal{O}_{X_{R_k,-}})\ar[d]\ar[d]\ar[d]\ar[d]\ar[r]^{\mathrm{global}}\ar[r]\ar[r] &\mathrm{sComm}_\mathrm{simplicial}\varphi\mathrm{Ind}^m\mathrm{Banach}_{\Gamma_k}(\{\mathrm{Robba}^\mathrm{extended}_{{R_k,-,I}}\}_I)\ar[d]\ar[d]\ar[d]\ar[d]\\
 \mathrm{sComm}_\mathrm{simplicial}\mathrm{Ind}^m\mathrm{Banach}_{\Gamma_0}(\mathcal{O}_{X_{\mathbb{Q}_p(p^{1/p^\infty})^{\wedge\flat},-}})\ar[r]^{\mathrm{global}}\ar[r]\ar[r] &\mathrm{sComm}_\mathrm{simplicial}\varphi\mathrm{Ind}^m\mathrm{Banach}_{\Gamma_0}(\{\mathrm{Robba}^\mathrm{extended}_{{R_0,-,I}}\}_I). 
}
\]

\item Then parallel as in \cite{LBV} we have a functor (global section) of the de Rham complex after \cite[Definition 5.9, Section 5.2.1]{KKM}:
\[\displayindent=-0.4in
\xymatrix@R+6pc@C+0pc{
\mathrm{deRham}_{\mathrm{sComm}_\mathrm{simplicial}\mathrm{Ind}\mathrm{Banach}_{\Gamma_k}(\mathcal{O}_{X_{R_k,-}})\ar[r]^{\mathrm{global}}}(-)\ar[d]\ar[d]\ar[d]\ar[d]\ar[r]\ar[r] &\mathrm{deRham}_{\mathrm{sComm}_\mathrm{simplicial}\varphi\mathrm{Ind}\mathrm{Banach}_{\Gamma_k}(\{\mathrm{Robba}^\mathrm{extended}_{{R_k,-,I}}\}_I)}(-)\ar[d]\ar[d]\ar[d]\ar[d]\\
\mathrm{deRham}_{\mathrm{sComm}_\mathrm{simplicial}\mathrm{Ind}\mathrm{Banach}_{\Gamma_0}(\mathcal{O}_{X_{\mathbb{Q}_p(p^{1/p^\infty})^{\wedge\flat},-}})\ar[r]^{\mathrm{global}}}(-)\ar[r]\ar[r] &\mathrm{deRham}_{\mathrm{sComm}_\mathrm{simplicial}\varphi\mathrm{Ind}\mathrm{Banach}_{\Gamma_0}(\{\mathrm{Robba}^\mathrm{extended}_{{R_0,-,I}}\}_I)}(-), 
}
\]
\[\displayindent=-0.4in
\xymatrix@R+6pc@C+0pc{
\mathrm{deRham}_{\mathrm{sComm}_\mathrm{simplicial}\mathrm{Ind}^m\mathrm{Banach}_{\Gamma_k}(\mathcal{O}_{X_{R_k,-}})\ar[r]^{\mathrm{global}}}(-)\ar[d]\ar[d]\ar[d]\ar[d]\ar[r]\ar[r] &\mathrm{deRham}_{\mathrm{sComm}_\mathrm{simplicial}\varphi\mathrm{Ind}^m\mathrm{Banach}_{\Gamma_k}(\{\mathrm{Robba}^\mathrm{extended}_{{R_k,-,I}}\}_I)}(-)\ar[d]\ar[d]\ar[d]\ar[d]\\
\mathrm{deRham}_{\mathrm{sComm}_\mathrm{simplicial}\mathrm{Ind}^m\mathrm{Banach}_{\Gamma_0}(\mathcal{O}_{X_{\mathbb{Q}_p(p^{1/p^\infty})^{\wedge\flat},-}})\ar[r]^{\mathrm{global}}}(-)\ar[r]\ar[r] &\mathrm{deRham}_{\mathrm{sComm}_\mathrm{simplicial}\varphi\mathrm{Ind}^m\mathrm{Banach}_{\Gamma_0}(\{\mathrm{Robba}^\mathrm{extended}_{{R_0,-,I}}\}_I)}(-). 
}
\]

\item Then we have the following a functor (global section) of $K$-group $(\infty,1)$-spectrum from \cite{BGT}:
\[
\xymatrix@R+6pc@C+0pc{
\mathrm{K}^\mathrm{BGT}_{\mathrm{sComm}_\mathrm{simplicial}\mathrm{Ind}\mathrm{Banach}_{\Gamma_k}(\mathcal{O}_{X_{R_k,-}})\ar[r]^{\mathrm{global}}}(-)\ar[d]\ar[d]\ar[d]\ar[d]\ar[r]\ar[r] &\mathrm{K}^\mathrm{BGT}_{\mathrm{sComm}_\mathrm{simplicial}\varphi\mathrm{Ind}\mathrm{Banach}_{\Gamma_k}(\{\mathrm{Robba}^\mathrm{extended}_{{R_k,-,I}}\}_I)}(-)\ar[d]\ar[d]\ar[d]\ar[d]\\
\mathrm{K}^\mathrm{BGT}_{\mathrm{sComm}_\mathrm{simplicial}\mathrm{Ind}\mathrm{Banach}_{\Gamma_0}(\mathcal{O}_{X_{\mathbb{Q}_p(p^{1/p^\infty})^{\wedge\flat},-}})\ar[r]^{\mathrm{global}}}(-)\ar[r]\ar[r] &\mathrm{K}^\mathrm{BGT}_{\mathrm{sComm}_\mathrm{simplicial}\varphi\mathrm{Ind}\mathrm{Banach}_{\Gamma_0}(\{\mathrm{Robba}^\mathrm{extended}_{{R_0,-,I}}\}_I)}(-), 
}
\]
\[
\xymatrix@R+6pc@C+0pc{
\mathrm{K}^\mathrm{BGT}_{\mathrm{sComm}_\mathrm{simplicial}\mathrm{Ind}^m\mathrm{Banach}_{\Gamma_k}(\mathcal{O}_{X_{R_k,-}})\ar[r]^{\mathrm{global}}}(-)\ar[d]\ar[d]\ar[d]\ar[d]\ar[r]\ar[r] &\mathrm{K}^\mathrm{BGT}_{\mathrm{sComm}_\mathrm{simplicial}\varphi\mathrm{Ind}^m\mathrm{Banach}_{\Gamma_k}(\{\mathrm{Robba}^\mathrm{extended}_{{R_k,-,I}}\}_I)}(-)\ar[d]\ar[d]\ar[d]\ar[d]\\
\mathrm{K}^\mathrm{BGT}_{\mathrm{sComm}_\mathrm{simplicial}\mathrm{Ind}^m\mathrm{Banach}_{\Gamma_0}(\mathcal{O}_{X_{\mathbb{Q}_p(p^{1/p^\infty})^{\wedge\flat},-}})\ar[r]^{\mathrm{global}}}(-)\ar[r]\ar[r] &\mathrm{K}^\mathrm{BGT}_{\mathrm{sComm}_\mathrm{simplicial}\varphi\mathrm{Ind}^m\mathrm{Banach}_{\Gamma_0}(\{\mathrm{Robba}^\mathrm{extended}_{{R_0,-,I}}\}_I)}(-). 
}
\]

\end{itemize}

\

\begin{remark}
\noindent We can certainly consider the quasicoherent sheaves in \cite[Lemma 7.11, Remark 7.12]{1BBK}, therefore all the quasicoherent presheaves and modules will be those satisfying \cite[Lemma 7.11, Remark 7.12]{1BBK} if one would like to consider the the quasicoherent sheaves. That being all as this said, we would believe that the big quasicoherent presheaves are automatically quasicoherent sheaves (namely satisfying the corresponding \v{C}ech $\infty$-descent as in \cite[Section 9.3]{KKM} and \cite[Lemma 7.11, Remark 7.12]{1BBK}) and the corresponding global section functors are automatically equivalence of $\infty$-categories. \\
\end{remark}

\

\indent In Clausen-Scholze formalism we have the following:

\begin{itemize}

\item (\text{Proposition}) There is a functor (global section) between the $\infty$-prestacks of inductive Banach quasicoherent sheaves:
\[
\xymatrix@R+0pc@C+0pc{
{\mathrm{Modules}_\circledcirc}(\mathcal{O}_{X_{R,-}})\ar[r]^{\mathrm{global}}\ar[r]\ar[r] &\varphi{\mathrm{Modules}_\circledcirc}(\{\mathrm{Robba}^\mathrm{extended}_{{R,-,I}}\}_I).  
}
\]

\item (\text{Proposition}) There is a functor (global section) between the $\infty$-prestacks of inductive Banach quasicoherent sheaves:
\[
\xymatrix@R+0pc@C+0pc{
{\mathrm{Modules}_\circledcirc}(\mathcal{O}_{X_{R,-}})\ar[r]^{\mathrm{global}}\ar[r]\ar[r] &\varphi{\mathrm{Modules}_\circledcirc}(\{\mathrm{Robba}^\mathrm{extended}_{{R,-,I}}\}_I).  
}
\]

\item (\text{Proposition}) There is a functor (global section) between the $\infty$-prestacks of inductive Banach quasicoherent commutative algebra $E_\infty$ objects\footnote{Here $\circledcirc=\text{solidquasicoherentsheaves}$.}:
\[
\xymatrix@R+0pc@C+0pc{
\mathrm{sComm}_\mathrm{simplicial}{\mathrm{Modules}_\circledcirc}(\mathcal{O}_{X_{R,-}})\ar[r]^{\mathrm{global}}\ar[r]\ar[r] &\mathrm{sComm}_\mathrm{simplicial}\varphi{\mathrm{Modules}_\circledcirc}(\{\mathrm{Robba}^\mathrm{extended}_{{R,-,I}}\}_I).  
}
\]

\item Then as in \cite{LBV} we have a functor (global section) of the de Rham complex after \cite[Definition 5.9, Section 5.2.1]{KKM}\footnote{Here $\circledcirc=\text{solidquasicoherentsheaves}$.}:
\[
\xymatrix@R+0pc@C+0pc{
\mathrm{deRham}_{\mathrm{sComm}_\mathrm{simplicial}{\mathrm{Modules}_\circledcirc}(\mathcal{O}_{X_{R,-}})\ar[r]^{\mathrm{global}}}(-)\ar[r]\ar[r] &\mathrm{deRham}_{\mathrm{sComm}_\mathrm{simplicial}\varphi{\mathrm{Modules}_\circledcirc}(\{\mathrm{Robba}^\mathrm{extended}_{{R,-,I}}\}_I)}(-). 
}
\]

\item Then we have the following a functor (global section) of $K$-group $(\infty,1)$-spectrum from \cite{BGT}\footnote{Here $\circledcirc=\text{solidquasicoherentsheaves}$.}:
\[
\xymatrix@R+0pc@C+0pc{
\mathrm{K}^\mathrm{BGT}_{\mathrm{sComm}_\mathrm{simplicial}{\mathrm{Modules}_\circledcirc}(\mathcal{O}_{X_{R,-}})\ar[r]^{\mathrm{global}}}(-)\ar[r]\ar[r] &\mathrm{K}^\mathrm{BGT}_{\mathrm{sComm}_\mathrm{simplicial}\varphi{\mathrm{Modules}_\circledcirc}(\{\mathrm{Robba}^\mathrm{extended}_{{R,-,I}}\}_I)}(-). 
}
\]
\end{itemize}

\noindent Now let $R=\mathbb{Q}_p(p^{1/p^\infty})^{\wedge\flat}$ and $R_k=\mathbb{Q}_p(p^{1/p^\infty})^{\wedge}\left<T_1^{\pm 1/p^{\infty}},...,T_k^{\pm 1/p^{\infty}}\right>^\flat$ we have the following Galois theoretic results with naturality along $f:\mathrm{Spa}(\mathbb{Q}_p(p^{1/p^\infty})^{\wedge}\left<T_1^{\pm 1/p^\infty},...,T_k^{\pm 1/p^\infty}\right>^\flat)\rightarrow \mathrm{Spa}(\mathbb{Q}_p(p^{1/p^\infty})^{\wedge\flat})$:

\begin{itemize}
\item (\text{Proposition}) There is a functor (global section) between the $\infty$-prestacks of inductive Banach quasicoherent sheaves\footnote{Here $\circledcirc=\text{solidquasicoherentsheaves}$.}:
\[
\xymatrix@R+6pc@C+0pc{
{\mathrm{Modules}_\circledcirc}(\mathcal{O}_{X_{\mathbb{Q}_p(p^{1/p^\infty})^{\wedge}\left<T_1^{\pm 1/p^\infty},...,T_k^{\pm 1/p^\infty}\right>^\flat,-}})\ar[d]\ar[d]\ar[d]\ar[d] \ar[r]^{\mathrm{global}}\ar[r]\ar[r] &\varphi{\mathrm{Modules}_\circledcirc}(\{\mathrm{Robba}^\mathrm{extended}_{{R_k,-,I}}\}_I) \ar[d]\ar[d]\ar[d]\ar[d].\\
{\mathrm{Modules}_\circledcirc}(\mathcal{O}_{X_{\mathbb{Q}_p(p^{1/p^\infty})^{\wedge\flat},-}})\ar[r]^{\mathrm{global}}\ar[r]\ar[r] &\varphi{\mathrm{Modules}_\circledcirc}(\{\mathrm{Robba}^\mathrm{extended}_{{R_0,-,I}}\}_I).\\ 
}
\]
\item (\text{Proposition}) There is a functor (global section) between the $\infty$-prestacks of inductive Banach quasicoherent commutative algebra $E_\infty$ objects\footnote{Here $\circledcirc=\text{solidquasicoherentsheaves}$.}:
\[
\xymatrix@R+6pc@C+0pc{
\mathrm{sComm}_\mathrm{simplicial}{\mathrm{Modules}_\circledcirc}(\mathcal{O}_{X_{R_k,-}})\ar[d]\ar[d]\ar[d]\ar[d]\ar[r]^{\mathrm{global}}\ar[r]\ar[r] &\mathrm{sComm}_\mathrm{simplicial}\varphi{\mathrm{Modules}_\circledcirc}(\{\mathrm{Robba}^\mathrm{extended}_{{R_k,-,I}}\}_I)\ar[d]\ar[d]\ar[d]\ar[d]\\
\mathrm{sComm}_\mathrm{simplicial}{\mathrm{Modules}_\circledcirc}(\mathcal{O}_{X_{\mathbb{Q}_p(p^{1/p^\infty})^{\wedge\flat},-}})\ar[r]^{\mathrm{global}}\ar[r]\ar[r] &\mathrm{sComm}_\mathrm{simplicial}\varphi{\mathrm{Modules}_\circledcirc}(\{\mathrm{Robba}^\mathrm{extended}_{{R_0,-,I}}\}_I).  
}
\]

\item Then as in \cite{LBV} we have a functor (global section) of the de Rham complex after \cite[Definition 5.9, Section 5.2.1]{KKM}\footnote{Here $\circledcirc=\text{solidquasicoherentsheaves}$.}:
\[\displayindent=-0.4in
\xymatrix@R+6pc@C+0pc{
\mathrm{deRham}_{\mathrm{sComm}_\mathrm{simplicial}{\mathrm{Modules}_\circledcirc}(\mathcal{O}_{X_{R_k,-}})\ar[r]^{\mathrm{global}}}(-)\ar[d]\ar[d]\ar[d]\ar[d]\ar[r]\ar[r] &\mathrm{deRham}_{\mathrm{sComm}_\mathrm{simplicial}\varphi{\mathrm{Modules}_\circledcirc}(\{\mathrm{Robba}^\mathrm{extended}_{{R_k,-,I}}\}_I)}(-)\ar[d]\ar[d]\ar[d]\ar[d]\\
\mathrm{deRham}_{\mathrm{sComm}_\mathrm{simplicial}{\mathrm{Modules}_\circledcirc}(\mathcal{O}_{X_{\mathbb{Q}_p(p^{1/p^\infty})^{\wedge\flat},-}})\ar[r]^{\mathrm{global}}}(-)\ar[r]\ar[r] &\mathrm{deRham}_{\mathrm{sComm}_\mathrm{simplicial}\varphi{\mathrm{Modules}_\circledcirc}(\{\mathrm{Robba}^\mathrm{extended}_{{R_0,-,I}}\}_I)}(-). 
}
\]

\item Then we have the following a functor (global section) of $K$-group $(\infty,1)$-spectrum from \cite{BGT}\footnote{Here $\circledcirc=\text{solidquasicoherentsheaves}$.}:
\[
\xymatrix@R+6pc@C+0pc{
\mathrm{K}^\mathrm{BGT}_{\mathrm{sComm}_\mathrm{simplicial}{\mathrm{Modules}_\circledcirc}(\mathcal{O}_{X_{R_k,-}})\ar[r]^{\mathrm{global}}}(-)\ar[d]\ar[d]\ar[d]\ar[d]\ar[r]\ar[r] &\mathrm{K}^\mathrm{BGT}_{\mathrm{sComm}_\mathrm{simplicial}\varphi{\mathrm{Modules}_\circledcirc}(\{\mathrm{Robba}^\mathrm{extended}_{{R_k,-,I}}\}_I)}(-)\ar[d]\ar[d]\ar[d]\ar[d]\\
\mathrm{K}^\mathrm{BGT}_{\mathrm{sComm}_\mathrm{simplicial}{\mathrm{Modules}_\circledcirc}(\mathcal{O}_{X_{\mathbb{Q}_p(p^{1/p^\infty})^{\wedge\flat},-}})\ar[r]^{\mathrm{global}}}(-)\ar[r]\ar[r] &\mathrm{K}^\mathrm{BGT}_{\mathrm{sComm}_\mathrm{simplicial}\varphi{\mathrm{Modules}_\circledcirc}(\{\mathrm{Robba}^\mathrm{extended}_{{R_0,-,I}}\}_I)}(-). 
}
\]

\end{itemize}

\
\indent Then we consider further equivariance by considering the arithmetic profinite fundamental groups $\Gamma_{\mathbb{Q}_p}$ and $\mathrm{Gal}(\overline{{Q}_p\left<T_1^{\pm 1},...,T_k^{\pm 1}\right>}/R_k)$ through the following diagram:

\[
\xymatrix@R+0pc@C+0pc{
\mathbb{Z}_p^k=\mathrm{Gal}(R_k/{\mathbb{Q}_p(p^{1/p^\infty})^\wedge\left<T_1^{\pm 1},...,T_k^{\pm 1}\right>}) \ar[r]\ar[r] \ar[r]\ar[r] &\Gamma_k:=\mathrm{Gal}(R_k/{\mathbb{Q}_p\left<T_1^{\pm 1},...,T_k^{\pm 1}\right>}) \ar[r] \ar[r]\ar[r] &\Gamma_{\mathbb{Q}_p}.
}
\]

\begin{itemize}
\item (\text{Proposition}) There is a functor (global section) between the $\infty$-prestacks of inductive Banach quasicoherent sheaves\footnote{Here $\circledcirc=\text{solidquasicoherentsheaves}$.}:
\[
\xymatrix@R+6pc@C+0pc{
{\mathrm{Modules}_\circledcirc}_{\Gamma_k}(\mathcal{O}_{X_{\mathbb{Q}_p(p^{1/p^\infty})^{\wedge}\left<T_1^{\pm 1/p^\infty},...,T_k^{\pm 1/p^\infty}\right>^\flat,-}})\ar[d]\ar[d]\ar[d]\ar[d] \ar[r]^{\mathrm{global}}\ar[r]\ar[r] &\varphi{\mathrm{Modules}_\circledcirc}_{\Gamma_k}(\{\mathrm{Robba}^\mathrm{extended}_{{R_k,-,I}}\}_I) \ar[d]\ar[d]\ar[d]\ar[d].\\
{\mathrm{Modules}_\circledcirc}(\mathcal{O}_{X_{\mathbb{Q}_p(p^{1/p^\infty})^{\wedge\flat},-}})\ar[r]^{\mathrm{global}}\ar[r]\ar[r] &\varphi{\mathrm{Modules}_\circledcirc}(\{\mathrm{Robba}^\mathrm{extended}_{{R_0,-,I}}\}_I).\\ 
}
\]

\item (\text{Proposition}) There is a functor (global section) between the $\infty$-stacks of inductive Banach quasicoherent commutative algebra $E_\infty$ objects\footnote{Here $\circledcirc=\text{solidquasicoherentsheaves}$.}:
\[
\xymatrix@R+6pc@C+0pc{
\mathrm{sComm}_\mathrm{simplicial}{\mathrm{Modules}_\circledcirc}_{\Gamma_k}(\mathcal{O}_{X_{R_k,-}})\ar[d]\ar[d]\ar[d]\ar[d]\ar[r]^{\mathrm{global}}\ar[r]\ar[r] &\mathrm{sComm}_\mathrm{simplicial}\varphi{\mathrm{Modules}_\circledcirc}_{\Gamma_k}(\{\mathrm{Robba}^\mathrm{extended}_{{R_k,-,I}}\}_I)\ar[d]\ar[d]\ar[d]\ar[d]\\
\mathrm{sComm}_\mathrm{simplicial}{\mathrm{Modules}_\circledcirc}_{\Gamma_0}(\mathcal{O}_{X_{\mathbb{Q}_p(p^{1/p^\infty})^{\wedge\flat},-}})\ar[r]^{\mathrm{global}}\ar[r]\ar[r] &\mathrm{sComm}_\mathrm{simplicial}\varphi{\mathrm{Modules}_\circledcirc}_{\Gamma_0}(\{\mathrm{Robba}^\mathrm{extended}_{{R_0,-,I}}\}_I).  
}
\]

\item Then as in \cite{LBV} we have a functor (global section) of the de Rham complex after \cite[Definition 5.9, Section 5.2.1]{KKM}\footnote{Here $\circledcirc=\text{solidquasicoherentsheaves}$.}:
\[\displayindent=-0.4in
\xymatrix@R+6pc@C+0pc{
\mathrm{deRham}_{\mathrm{sComm}_\mathrm{simplicial}{\mathrm{Modules}_\circledcirc}_{\Gamma_k}(\mathcal{O}_{X_{R_k,-}})\ar[r]^{\mathrm{global}}}(-)\ar[d]\ar[d]\ar[d]\ar[d]\ar[r]\ar[r] &\mathrm{deRham}_{\mathrm{sComm}_\mathrm{simplicial}\varphi{\mathrm{Modules}_\circledcirc}_{\Gamma_k}(\{\mathrm{Robba}^\mathrm{extended}_{{R_k,-,I}}\}_I)}(-)\ar[d]\ar[d]\ar[d]\ar[d]\\
\mathrm{deRham}_{\mathrm{sComm}_\mathrm{simplicial}{\mathrm{Modules}_\circledcirc}_{\Gamma_0}(\mathcal{O}_{X_{\mathbb{Q}_p(p^{1/p^\infty})^{\wedge\flat},-}})\ar[r]^{\mathrm{global}}}(-)\ar[r]\ar[r] &\mathrm{deRham}_{\mathrm{sComm}_\mathrm{simplicial}\varphi{\mathrm{Modules}_\circledcirc}_{\Gamma_0}(\{\mathrm{Robba}^\mathrm{extended}_{{R_0,-,I}}\}_I)}(-). 
}
\]

\item Then we have the following a functor (global section) of $K$-group $(\infty,1)$-spectrum from \cite{BGT}\footnote{Here $\circledcirc=\text{solidquasicoherentsheaves}$.}:
\[
\xymatrix@R+6pc@C+0pc{
\mathrm{K}^\mathrm{BGT}_{\mathrm{sComm}_\mathrm{simplicial}{\mathrm{Modules}_\circledcirc}_{\Gamma_k}(\mathcal{O}_{X_{R_k,-}})\ar[r]^{\mathrm{global}}}(-)\ar[d]\ar[d]\ar[d]\ar[d]\ar[r]\ar[r] &\mathrm{K}^\mathrm{BGT}_{\mathrm{sComm}_\mathrm{simplicial}\varphi{\mathrm{Modules}_\circledcirc}_{\Gamma_k}(\{\mathrm{Robba}^\mathrm{extended}_{{R_k,-,I}}\}_I)}(-)\ar[d]\ar[d]\ar[d]\ar[d]\\
\mathrm{K}^\mathrm{BGT}_{\mathrm{sComm}_\mathrm{simplicial}{\mathrm{Modules}_\circledcirc}_{\Gamma_0}(\mathcal{O}_{X_{\mathbb{Q}_p(p^{1/p^\infty})^{\wedge\flat},-}})\ar[r]^{\mathrm{global}}}(-)\ar[r]\ar[r] &\mathrm{K}^\mathrm{BGT}_{\mathrm{sComm}_\mathrm{simplicial}\varphi{\mathrm{Modules}_\circledcirc}_{\Gamma_0}(\{\mathrm{Robba}^\mathrm{extended}_{{R_0,-,I}}\}_I)}(-). 
}
\]

\end{itemize}

\begin{proposition}
All the global functors from \cite[Proposition 13.8, Theorem 14.9, Remark 14.10]{1CS2} achieve the equivalences on both sides.	
\end{proposition}

\newpage
\subsection{$\infty$-Categorical Analytic Stacks and Descents IV}

\indent In the following the right had of each row in each diagram will be the corresponding quasicoherent Robba bundles over the Robba ring carrying the corresponding action from the Frobenius or the fundamental groups, defined by directly applying \cite[Section 9.3]{KKM} and \cite{BBM}. We now let $\mathcal{A}$ be any commutative algebra objects in the corresponding $\infty$-toposes over ind-Banach commutative algebra objects over $\mathbb{Q}_p$ or the corresponding born\'e commutative algebra objects over $\mathbb{Q}_p$ carrying the Grothendieck topology defined by essentially the corresponding monomorphism homotopy in the opposite category. Then we promote the construction to the corresponding $\infty$-stack over the same $\infty$-categories of affinoids. We now take the corresponding colimit through all the $(\infty,1)$-categories. Therefore all the corresponding $(\infty,1)$-functors into $(\infty,1)$-categories or $(\infty,1)$-groupoids are from the homotopy closure of $\mathbb{Q}_p\left<C_1,...,C_\ell\right>$ $\ell=1,2,...$ in $\mathrm{sComm}\mathrm{Ind}\mathrm{Banach}_{\mathbb{Q}_p}$ or $\mathbb{Q}_p\left<C_1,...,C_\ell\right>$ $\ell=1,2,...$ in $\mathrm{sComm}\mathrm{Ind}^m\mathrm{Banach}_{\mathbb{Q}_p}$ as in \cite[Section 4.2]{BBM}:
\begin{align}
\mathrm{Ind}^{\mathbb{Q}_p\left<C_1,...,C_\ell\right>,\ell=1,2,...}\mathrm{sComm}\mathrm{Ind}\mathrm{Banach}_{\mathbb{Q}_p},\\
\mathrm{Ind}^{\mathbb{Q}_p\left<C_1,...,C_\ell\right>,\ell=1,2,...}\mathrm{sComm}\mathrm{Ind}\mathrm{Banach}_{\mathbb{Q}_p}	.
\end{align}

\begin{itemize}

\item (\text{Proposition}) There is a functor (global section) between the $\infty$-prestacks of inductive Banach quasicoherent presheaves:
\[
\xymatrix@R+0pc@C+0pc{
\mathrm{Ind}\mathrm{Banach}(\mathcal{O}_{X_{R,-}})\ar[r]^{\mathrm{global}}\ar[r]\ar[r] &\varphi\mathrm{Ind}\mathrm{Banach}(\{\mathrm{Robba}^\mathrm{extended}_{{R,-,I}}\}_I).  
}
\]
The definition is given by the following:
\[
\xymatrix@R+0pc@C+0pc{
\mathrm{homotopycolimit}_i(\mathrm{Ind}\mathrm{Banach}(\mathcal{O}_{X_{R,-}})\ar[r]^{\mathrm{global}}\ar[r]\ar[r] &\varphi\mathrm{Ind}\mathrm{Banach}(\{\mathrm{Robba}^\mathrm{extended}_{{R,-,I}}\}_I))(\mathcal{O}_i),  
}
\]
each $\mathcal{O}_i$ is just as $\mathbb{Q}_p\left<C_1,...,C_\ell\right>,\ell=1,2,...$.
\item (\text{Proposition}) There is a functor (global section) between the $\infty$-prestacks of monomorphic inductive Banach quasicoherent presheaves:
\[
\xymatrix@R+0pc@C+0pc{
\mathrm{Ind}^m\mathrm{Banach}(\mathcal{O}_{X_{R,-}})\ar[r]^{\mathrm{global}}\ar[r]\ar[r] &\varphi\mathrm{Ind}^m\mathrm{Banach}(\{\mathrm{Robba}^\mathrm{extended}_{{R,-,I}}\}_I).  
}
\]
The definition is given by the following:
\[
\xymatrix@R+0pc@C+0pc{
\mathrm{homotopycolimit}_i(\mathrm{Ind}^m\mathrm{Banach}(\mathcal{O}_{X_{R,-}})\ar[r]^{\mathrm{global}}\ar[r]\ar[r] &\varphi\mathrm{Ind}^m\mathrm{Banach}(\{\mathrm{Robba}^\mathrm{extended}_{{R,-,I}}\}_I))(\mathcal{O}_i),  
}
\]
each $\mathcal{O}_i$ is just as $\mathbb{Q}_p\left<C_1,...,C_\ell\right>,\ell=1,2,...$.

\item (\text{Proposition}) There is a functor (global section) between the $\infty$-prestacks of inductive Banach quasicoherent presheaves:
\[
\xymatrix@R+0pc@C+0pc{
\mathrm{Ind}\mathrm{Banach}(\mathcal{O}_{X_{R,-}})\ar[r]^{\mathrm{global}}\ar[r]\ar[r] &\varphi\mathrm{Ind}\mathrm{Banach}(\{\mathrm{Robba}^\mathrm{extended}_{{R,-,I}}\}_I).  
}
\]
The definition is given by the following:
\[
\xymatrix@R+0pc@C+0pc{
\mathrm{homotopycolimit}_i(\mathrm{Ind}\mathrm{Banach}(\mathcal{O}_{X_{R,-}})\ar[r]^{\mathrm{global}}\ar[r]\ar[r] &\varphi\mathrm{Ind}\mathrm{Banach}(\{\mathrm{Robba}^\mathrm{extended}_{{R,-,I}}\}_I))(\mathcal{O}_i),  
}
\]
each $\mathcal{O}_i$ is just as $\mathbb{Q}_p\left<C_1,...,C_\ell\right>,\ell=1,2,...$.
\item (\text{Proposition}) There is a functor (global section) between the $\infty$-stacks of monomorphic inductive Banach quasicoherent presheaves:
\[
\xymatrix@R+0pc@C+0pc{
\mathrm{Ind}^m\mathrm{Banach}(\mathcal{O}_{X_{R,-}})\ar[r]^{\mathrm{global}}\ar[r]\ar[r] &\varphi\mathrm{Ind}^m\mathrm{Banach}(\{\mathrm{Robba}^\mathrm{extended}_{{R,-,I}}\}_I).  
}
\]
The definition is given by the following:
\[
\xymatrix@R+0pc@C+0pc{
\mathrm{homotopycolimit}_i(\mathrm{Ind}^m\mathrm{Banach}(\mathcal{O}_{X_{R,-}})\ar[r]^{\mathrm{global}}\ar[r]\ar[r] &\varphi\mathrm{Ind}^m\mathrm{Banach}(\{\mathrm{Robba}^\mathrm{extended}_{{R,-,I}}\}_I))(\mathcal{O}_i),  
}
\]
each $\mathcal{O}_i$ is just as $\mathbb{Q}_p\left<C_1,...,C_\ell\right>,\ell=1,2,...$.
\item (\text{Proposition}) There is a functor (global section) between the $\infty$-prestacks of inductive Banach quasicoherent commutative algebra $E_\infty$ objects:
\[
\xymatrix@R+0pc@C+0pc{
\mathrm{sComm}_\mathrm{simplicial}\mathrm{Ind}\mathrm{Banach}(\mathcal{O}_{X_{R,-}})\ar[r]^{\mathrm{global}}\ar[r]\ar[r] &\mathrm{sComm}_\mathrm{simplicial}\varphi\mathrm{Ind}\mathrm{Banach}(\{\mathrm{Robba}^\mathrm{extended}_{{R,-,I}}\}_I).  
}
\]
The definition is given by the following:
\[\displayindent=-0.4in
\xymatrix@R+0pc@C+0pc{
\mathrm{homotopycolimit}_i(\mathrm{sComm}_\mathrm{simplicial}\mathrm{Ind}\mathrm{Banach}(\mathcal{O}_{X_{R,-}})\ar[r]^{\mathrm{global}}\ar[r]\ar[r] &\mathrm{sComm}_\mathrm{simplicial}\varphi\mathrm{Ind}\mathrm{Banach}(\{\mathrm{Robba}^\mathrm{extended}_{{R,-,I}}\}_I))(\mathcal{O}_i),  
}
\]
each $\mathcal{O}_i$ is just as $\mathbb{Q}_p\left<C_1,...,C_\ell\right>,\ell=1,2,...$.
\item (\text{Proposition}) There is a functor (global section) between the $\infty$-prestacks of monomorphic inductive Banach quasicoherent commutative algebra $E_\infty$ objects:
\[
\xymatrix@R+0pc@C+0pc{
\mathrm{sComm}_\mathrm{simplicial}\mathrm{Ind}^m\mathrm{Banach}(\mathcal{O}_{X_{R,-}})\ar[r]^{\mathrm{global}}\ar[r]\ar[r] &\mathrm{sComm}_\mathrm{simplicial}\varphi\mathrm{Ind}^m\mathrm{Banach}(\{\mathrm{Robba}^\mathrm{extended}_{{R,-,I}}\}_I).  
}
\]
The definition is given by the following:
\[
\xymatrix@R+0pc@C+0pc{
\mathrm{homotopycolimit}_i(\mathrm{sComm}_\mathrm{simplicial}\mathrm{Ind}^m\mathrm{Banach}(\mathcal{O}_{X_{R,-}})\ar[r]^{\mathrm{global}}\ar[r]\ar[r] &\mathrm{sComm}_\mathrm{simplicial}\varphi\mathrm{Ind}^m\mathrm{Banach}(\{\mathrm{Robba}^\mathrm{extended}_{{R,-,I}}\}_I))(\mathcal{O}_i),  
}
\]
each $\mathcal{O}_i$ is just as $\mathbb{Q}_p\left<C_1,...,C_\ell\right>,\ell=1,2,...$.

\item Then parallel as in \cite{LBV} we have a functor (global section ) of the de Rham complex after \cite[Definition 5.9, Section 5.2.1]{KKM}:
\[
\xymatrix@R+0pc@C+0pc{
\mathrm{deRham}_{\mathrm{sComm}_\mathrm{simplicial}\mathrm{Ind}\mathrm{Banach}(\mathcal{O}_{X_{R,-}})\ar[r]^{\mathrm{global}}}(-)\ar[r]\ar[r] &\mathrm{deRham}_{\mathrm{sComm}_\mathrm{simplicial}\varphi\mathrm{Ind}\mathrm{Banach}(\{\mathrm{Robba}^\mathrm{extended}_{{R,-,I}}\}_I)}(-), 
}
\]
\[
\xymatrix@R+0pc@C+0pc{
\mathrm{deRham}_{\mathrm{sComm}_\mathrm{simplicial}\mathrm{Ind}^m\mathrm{Banach}(\mathcal{O}_{X_{R,-}})\ar[r]^{\mathrm{global}}}(-)\ar[r]\ar[r] &\mathrm{deRham}_{\mathrm{sComm}_\mathrm{simplicial}\varphi\mathrm{Ind}^m\mathrm{Banach}(\{\mathrm{Robba}^\mathrm{extended}_{{R,-,I}}\}_I)}(-). 
}
\]
The definition is given by the following:
\[
\xymatrix@R+0pc@C+0pc{
\mathrm{homotopycolimit}_i\\
(\mathrm{deRham}_{\mathrm{sComm}_\mathrm{simplicial}\mathrm{Ind}\mathrm{Banach}(\mathcal{O}_{X_{R,-}})\ar[r]^{\mathrm{global}}}(-)\ar[r]\ar[r] &\mathrm{deRham}_{\mathrm{sComm}_\mathrm{simplicial}\varphi\mathrm{Ind}\mathrm{Banach}(\{\mathrm{Robba}^\mathrm{extended}_{{R,-,I}}\}_I)}(-))(\mathcal{O}_i),  
}
\]
\[
\xymatrix@R+0pc@C+0pc{
\mathrm{homotopycolimit}_i\\
(\mathrm{deRham}_{\mathrm{sComm}_\mathrm{simplicial}\mathrm{Ind}^m\mathrm{Banach}(\mathcal{O}_{X_{R,-}})\ar[r]^{\mathrm{global}}}(-)\ar[r]\ar[r] &\mathrm{deRham}_{\mathrm{sComm}_\mathrm{simplicial}\varphi\mathrm{Ind}^m\mathrm{Banach}(\{\mathrm{Robba}^\mathrm{extended}_{{R,-,I}}\}_I)}(-))(\mathcal{O}_i),  
}
\]
each $\mathcal{O}_i$ is just as $\mathbb{Q}_p\left<C_1,...,C_\ell\right>,\ell=1,2,...$.\item Then we have the following a functor (global section) of $K$-group $(\infty,1)$-spectrum from \cite{BGT}:
\[
\xymatrix@R+0pc@C+0pc{
\mathrm{K}^\mathrm{BGT}_{\mathrm{sComm}_\mathrm{simplicial}\mathrm{Ind}\mathrm{Banach}(\mathcal{O}_{X_{R,-}})\ar[r]^{\mathrm{global}}}(-)\ar[r]\ar[r] &\mathrm{K}^\mathrm{BGT}_{\mathrm{sComm}_\mathrm{simplicial}\varphi\mathrm{Ind}\mathrm{Banach}(\{\mathrm{Robba}^\mathrm{extended}_{{R,-,I}}\}_I)}(-), 
}
\]
\[
\xymatrix@R+0pc@C+0pc{
\mathrm{K}^\mathrm{BGT}_{\mathrm{sComm}_\mathrm{simplicial}\mathrm{Ind}^m\mathrm{Banach}(\mathcal{O}_{X_{R,-}})\ar[r]^{\mathrm{global}}}(-)\ar[r]\ar[r] &\mathrm{K}^\mathrm{BGT}_{\mathrm{sComm}_\mathrm{simplicial}\varphi\mathrm{Ind}^m\mathrm{Banach}(\{\mathrm{Robba}^\mathrm{extended}_{{R,-,I}}\}_I)}(-). 
}
\]
The definition is given by the following:
\[\displayindent=-0.4in
\xymatrix@R+0pc@C+0pc{
\mathrm{homotopycolimit}_i(\mathrm{K}^\mathrm{BGT}_{\mathrm{sComm}_\mathrm{simplicial}\mathrm{Ind}\mathrm{Banach}(\mathcal{O}_{X_{R,-}})\ar[r]^{\mathrm{global}}}(-)\ar[r]\ar[r] &\mathrm{K}^\mathrm{BGT}_{\mathrm{sComm}_\mathrm{simplicial}\varphi\mathrm{Ind}\mathrm{Banach}(\{\mathrm{Robba}^\mathrm{extended}_{{R,-,I}}\}_I)}(-))(\mathcal{O}_i),  
}
\]
\[\displayindent=-0.4in
\xymatrix@R+0pc@C+0pc{
\mathrm{homotopycolimit}_i(\mathrm{K}^\mathrm{BGT}_{\mathrm{sComm}_\mathrm{simplicial}\mathrm{Ind}^m\mathrm{Banach}(\mathcal{O}_{X_{R,-}})\ar[r]^{\mathrm{global}}}(-)\ar[r]\ar[r] &\mathrm{K}^\mathrm{BGT}_{\mathrm{sComm}_\mathrm{simplicial}\varphi\mathrm{Ind}^m\mathrm{Banach}(\{\mathrm{Robba}^\mathrm{extended}_{{R,-,I}}\}_I)}(-))(\mathcal{O}_i),  
}
\]
each $\mathcal{O}_i$ is just as $\mathbb{Q}_p\left<C_1,...,C_\ell\right>,\ell=1,2,...$.
\end{itemize}

\noindent Now let $R=\mathbb{Q}_p(p^{1/p^\infty})^{\wedge\flat}$ and $R_k=\mathbb{Q}_p(p^{1/p^\infty})^{\wedge}\left<T_1^{\pm 1/p^{\infty}},...,T_k^{\pm 1/p^{\infty}}\right>^\flat$ we have the following Galois theoretic results with naturality along $f:\mathrm{Spa}(\mathbb{Q}_p(p^{1/p^\infty})^{\wedge}\left<T_1^{\pm 1/p^\infty},...,T_k^{\pm 1/p^\infty}\right>^\flat)\rightarrow \mathrm{Spa}(\mathbb{Q}_p(p^{1/p^\infty})^{\wedge\flat})$:

\begin{itemize}
\item (\text{Proposition}) There is a functor (global section) between the $\infty$-prestacks of inductive Banach quasicoherent presheaves:
\[
\xymatrix@R+6pc@C+0pc{
\mathrm{Ind}\mathrm{Banach}(\mathcal{O}_{X_{\mathbb{Q}_p(p^{1/p^\infty})^{\wedge}\left<T_1^{\pm 1/p^\infty},...,T_k^{\pm 1/p^\infty}\right>^\flat,-}})\ar[d]\ar[d]\ar[d]\ar[d] \ar[r]^{\mathrm{global}}\ar[r]\ar[r] &\varphi\mathrm{Ind}\mathrm{Banach}(\{\mathrm{Robba}^\mathrm{extended}_{{R_k,-,I}}\}_I) \ar[d]\ar[d]\ar[d]\ar[d].\\
\mathrm{Ind}\mathrm{Banach}(\mathcal{O}_{X_{\mathbb{Q}_p(p^{1/p^\infty})^{\wedge\flat},-}})\ar[r]^{\mathrm{global}}\ar[r]\ar[r] &\varphi\mathrm{Ind}\mathrm{Banach}(\{\mathrm{Robba}^\mathrm{extended}_{{R_0,-,I}}\}_I).\\ 
}
\]
The definition is given by the following:
\[
\xymatrix@R+0pc@C+0pc{
\mathrm{homotopycolimit}_i(\square)(\mathcal{O}_i),  
}
\]
each $\mathcal{O}_i$ is just as $\mathbb{Q}_p\left<C_1,...,C_\ell\right>,\ell=1,2,...$ and $\square$ is the relative diagram of $\infty$-functors.
\item (\text{Proposition}) There is a functor (global section) between the $\infty$-prestacks of monomorphic inductive Banach quasicoherent presheaves:
\[
\xymatrix@R+6pc@C+0pc{
\mathrm{Ind}^m\mathrm{Banach}(\mathcal{O}_{X_{R_k,-}})\ar[r]^{\mathrm{global}}\ar[d]\ar[d]\ar[d]\ar[d]\ar[r]\ar[r] &\varphi\mathrm{Ind}^m\mathrm{Banach}(\{\mathrm{Robba}^\mathrm{extended}_{{R_k,-,I}}\}_I)\ar[d]\ar[d]\ar[d]\ar[d]\\
\mathrm{Ind}^m\mathrm{Banach}(\mathcal{O}_{X_{\mathbb{Q}_p(p^{1/p^\infty})^{\wedge\flat},-}})\ar[r]^{\mathrm{global}}\ar[r]\ar[r] &\varphi\mathrm{Ind}^m\mathrm{Banach}(\{\mathrm{Robba}^\mathrm{extended}_{{R_0,-,I}}\}_I).\\  
}
\]
The definition is given by the following:
\[
\xymatrix@R+0pc@C+0pc{
\mathrm{homotopycolimit}_i(\square)(\mathcal{O}_i),  
}
\]
each $\mathcal{O}_i$ is just as $\mathbb{Q}_p\left<C_1,...,C_\ell\right>,\ell=1,2,...$ and $\square$ is the relative diagram of $\infty$-functors.

\item (\text{Proposition}) There is a functor (global section) between the $\infty$-prestacks of inductive Banach quasicoherent commutative algebra $E_\infty$ objects:
\[
\xymatrix@R+6pc@C+0pc{
\mathrm{sComm}_\mathrm{simplicial}\mathrm{Ind}\mathrm{Banach}(\mathcal{O}_{X_{R_k,-}})\ar[d]\ar[d]\ar[d]\ar[d]\ar[r]^{\mathrm{global}}\ar[r]\ar[r] &\mathrm{sComm}_\mathrm{simplicial}\varphi\mathrm{Ind}\mathrm{Banach}(\{\mathrm{Robba}^\mathrm{extended}_{{R_k,-,I}}\}_I)\ar[d]\ar[d]\ar[d]\ar[d]\\
\mathrm{sComm}_\mathrm{simplicial}\mathrm{Ind}\mathrm{Banach}(\mathcal{O}_{X_{\mathbb{Q}_p(p^{1/p^\infty})^{\wedge\flat},-}})\ar[r]^{\mathrm{global}}\ar[r]\ar[r] &\mathrm{sComm}_\mathrm{simplicial}\varphi\mathrm{Ind}\mathrm{Banach}(\{\mathrm{Robba}^\mathrm{extended}_{{R_0,-,I}}\}_I).  
}
\]
The definition is given by the following:
\[
\xymatrix@R+0pc@C+0pc{
\mathrm{homotopycolimit}_i(\square)(\mathcal{O}_i),  
}
\]
each $\mathcal{O}_i$ is just as $\mathbb{Q}_p\left<C_1,...,C_\ell\right>,\ell=1,2,...$ and $\square$ is the relative diagram of $\infty$-functors.

\item (\text{Proposition}) There is a functor (global section) between the $\infty$-prestacks of monomorphic inductive Banach quasicoherent commutative algebra $E_\infty$ objects:
\[
\xymatrix@R+6pc@C+0pc{
\mathrm{sComm}_\mathrm{simplicial}\mathrm{Ind}^m\mathrm{Banach}(\mathcal{O}_{X_{R_k,-}})\ar[d]\ar[d]\ar[d]\ar[d]\ar[r]^{\mathrm{global}}\ar[r]\ar[r] &\mathrm{sComm}_\mathrm{simplicial}\varphi\mathrm{Ind}^m\mathrm{Banach}(\{\mathrm{Robba}^\mathrm{extended}_{{R_k,-,I}}\}_I)\ar[d]\ar[d]\ar[d]\ar[d]\\
 \mathrm{sComm}_\mathrm{simplicial}\mathrm{Ind}^m\mathrm{Banach}(\mathcal{O}_{X_{\mathbb{Q}_p(p^{1/p^\infty})^{\wedge\flat},-}})\ar[r]^{\mathrm{global}}\ar[r]\ar[r] &\mathrm{sComm}_\mathrm{simplicial}\varphi\mathrm{Ind}^m\mathrm{Banach}(\{\mathrm{Robba}^\mathrm{extended}_{{R_0,-,I}}\}_I).
}
\]
The definition is given by the following:
\[
\xymatrix@R+0pc@C+0pc{
\mathrm{homotopycolimit}_i(\square)(\mathcal{O}_i),  
}
\]
each $\mathcal{O}_i$ is just as $\mathbb{Q}_p\left<C_1,...,C_\ell\right>,\ell=1,2,...$ and $\square$ is the relative diagram of $\infty$-functors.

\item Then parallel as in \cite{LBV} we have a functor (global section) of the de Rham complex after \cite[Definition 5.9, Section 5.2.1]{KKM}:
\[\displayindent=-0.4in
\xymatrix@R+6pc@C+0pc{
\mathrm{deRham}_{\mathrm{sComm}_\mathrm{simplicial}\mathrm{Ind}\mathrm{Banach}(\mathcal{O}_{X_{R_k,-}})\ar[r]^{\mathrm{global}}}(-)\ar[d]\ar[d]\ar[d]\ar[d]\ar[r]\ar[r] &\mathrm{deRham}_{\mathrm{sComm}_\mathrm{simplicial}\varphi\mathrm{Ind}\mathrm{Banach}(\{\mathrm{Robba}^\mathrm{extended}_{{R_k,-,I}}\}_I)}(-)\ar[d]\ar[d]\ar[d]\ar[d]\\
\mathrm{deRham}_{\mathrm{sComm}_\mathrm{simplicial}\mathrm{Ind}\mathrm{Banach}(\mathcal{O}_{X_{\mathbb{Q}_p(p^{1/p^\infty})^{\wedge\flat},-}})\ar[r]^{\mathrm{global}}}(-)\ar[r]\ar[r] &\mathrm{deRham}_{\mathrm{sComm}_\mathrm{simplicial}\varphi\mathrm{Ind}\mathrm{Banach}(\{\mathrm{Robba}^\mathrm{extended}_{{R_0,-,I}}\}_I)}(-), 
}
\]
\[\displayindent=-0.4in
\xymatrix@R+6pc@C+0pc{
\mathrm{deRham}_{\mathrm{sComm}_\mathrm{simplicial}\mathrm{Ind}^m\mathrm{Banach}(\mathcal{O}_{X_{R_k,-}})\ar[r]^{\mathrm{global}}}(-)\ar[d]\ar[d]\ar[d]\ar[d]\ar[r]\ar[r] &\mathrm{deRham}_{\mathrm{sComm}_\mathrm{simplicial}\varphi\mathrm{Ind}^m\mathrm{Banach}(\{\mathrm{Robba}^\mathrm{extended}_{{R_k,-,I}}\}_I)}(-)\ar[d]\ar[d]\ar[d]\ar[d]\\
\mathrm{deRham}_{\mathrm{sComm}_\mathrm{simplicial}\mathrm{Ind}^m\mathrm{Banach}(\mathcal{O}_{X_{\mathbb{Q}_p(p^{1/p^\infty})^{\wedge\flat},-}})\ar[r]^{\mathrm{global}}}(-)\ar[r]\ar[r] &\mathrm{deRham}_{\mathrm{sComm}_\mathrm{simplicial}\varphi\mathrm{Ind}^m\mathrm{Banach}(\{\mathrm{Robba}^\mathrm{extended}_{{R_0,-,I}}\}_I)}(-). 
}
\]

\item Then we have the following a functor (global section) of $K$-group $(\infty,1)$-spectrum from \cite{BGT}:
\[
\xymatrix@R+6pc@C+0pc{
\mathrm{K}^\mathrm{BGT}_{\mathrm{sComm}_\mathrm{simplicial}\mathrm{Ind}\mathrm{Banach}(\mathcal{O}_{X_{R_k,-}})\ar[r]^{\mathrm{global}}}(-)\ar[d]\ar[d]\ar[d]\ar[d]\ar[r]\ar[r] &\mathrm{K}^\mathrm{BGT}_{\mathrm{sComm}_\mathrm{simplicial}\varphi\mathrm{Ind}\mathrm{Banach}(\{\mathrm{Robba}^\mathrm{extended}_{{R_k,-,I}}\}_I)}(-)\ar[d]\ar[d]\ar[d]\ar[d]\\
\mathrm{K}^\mathrm{BGT}_{\mathrm{sComm}_\mathrm{simplicial}\mathrm{Ind}\mathrm{Banach}(\mathcal{O}_{X_{\mathbb{Q}_p(p^{1/p^\infty})^{\wedge\flat},-}})\ar[r]^{\mathrm{global}}}(-)\ar[r]\ar[r] &\mathrm{K}^\mathrm{BGT}_{\mathrm{sComm}_\mathrm{simplicial}\varphi\mathrm{Ind}\mathrm{Banach}(\{\mathrm{Robba}^\mathrm{extended}_{{R_0,-,I}}\}_I)}(-), 
}
\]
\[
\xymatrix@R+6pc@C+0pc{
\mathrm{K}^\mathrm{BGT}_{\mathrm{sComm}_\mathrm{simplicial}\mathrm{Ind}^m\mathrm{Banach}(\mathcal{O}_{X_{R_k,-}})\ar[r]^{\mathrm{global}}}(-)\ar[d]\ar[d]\ar[d]\ar[d]\ar[r]\ar[r] &\mathrm{K}^\mathrm{BGT}_{\mathrm{sComm}_\mathrm{simplicial}\varphi\mathrm{Ind}^m\mathrm{Banach}(\{\mathrm{Robba}^\mathrm{extended}_{{R_k,-,I}}\}_I)}(-)\ar[d]\ar[d]\ar[d]\ar[d]\\
\mathrm{K}^\mathrm{BGT}_{\mathrm{sComm}_\mathrm{simplicial}\mathrm{Ind}^m\mathrm{Banach}(\mathcal{O}_{X_{\mathbb{Q}_p(p^{1/p^\infty})^{\wedge\flat},-}})\ar[r]^{\mathrm{global}}}(-)\ar[r]\ar[r] &\mathrm{K}^\mathrm{BGT}_{\mathrm{sComm}_\mathrm{simplicial}\varphi\mathrm{Ind}^m\mathrm{Banach}(\{\mathrm{Robba}^\mathrm{extended}_{{R_0,-,I}}\}_I)}(-). 
}
\]
The definition is given by the following:
\[
\xymatrix@R+0pc@C+0pc{
\mathrm{homotopycolimit}_i(\square)(\mathcal{O}_i),  
}
\]
each $\mathcal{O}_i$ is just as $\mathbb{Q}_p\left<C_1,...,C_\ell\right>,\ell=1,2,...$ and $\square$ is the relative diagram of $\infty$-functors.

\end{itemize}

\
\indent Then we consider further equivariance by considering the arithmetic profinite fundamental groups $\Gamma_{\mathbb{Q}_p}$ and $\mathrm{Gal}(\overline{{Q}_p\left<T_1^{\pm 1},...,T_k^{\pm 1}\right>}/R_k)$ through the following diagram:

\[
\xymatrix@R+0pc@C+0pc{
\mathbb{Z}_p^k=\mathrm{Gal}(R_k/{\mathbb{Q}_p(p^{1/p^\infty})^\wedge\left<T_1^{\pm 1},...,T_k^{\pm 1}\right>}) \ar[r]\ar[r] \ar[r]\ar[r] &\Gamma_k:=\mathrm{Gal}(R_k/{\mathbb{Q}_p\left<T_1^{\pm 1},...,T_k^{\pm 1}\right>}) \ar[r] \ar[r]\ar[r] &\Gamma_{\mathbb{Q}_p}.
}
\]

\begin{itemize}
\item (\text{Proposition}) There is a functor (global section) between the $\infty$-prestacks of inductive Banach quasicoherent presheaves:
\[
\xymatrix@R+6pc@C+0pc{
\mathrm{Ind}\mathrm{Banach}_{\Gamma_k}(\mathcal{O}_{X_{\mathbb{Q}_p(p^{1/p^\infty})^{\wedge}\left<T_1^{\pm 1/p^\infty},...,T_k^{\pm 1/p^\infty}\right>^\flat,-}})\ar[d]\ar[d]\ar[d]\ar[d] \ar[r]^{\mathrm{global}}\ar[r]\ar[r] &\varphi\mathrm{Ind}\mathrm{Banach}_{\Gamma_k}(\{\mathrm{Robba}^\mathrm{extended}_{{R_k,-,I}}\}_I) \ar[d]\ar[d]\ar[d]\ar[d].\\
\mathrm{Ind}\mathrm{Banach}(\mathcal{O}_{X_{\mathbb{Q}_p(p^{1/p^\infty})^{\wedge\flat},-}})\ar[r]^{\mathrm{global}}\ar[r]\ar[r] &\varphi\mathrm{Ind}\mathrm{Banach}(\{\mathrm{Robba}^\mathrm{extended}_{{R_0,-,I}}\}_I).\\ 
}
\]
The definition is given by the following:
\[
\xymatrix@R+0pc@C+0pc{
\mathrm{homotopycolimit}_i(\square)(\mathcal{O}_i),  
}
\]
each $\mathcal{O}_i$ is just as $\mathbb{Q}_p\left<C_1,...,C_\ell\right>,\ell=1,2,...$ and $\square$ is the relative diagram of $\infty$-functors.

\item (\text{Proposition}) There is a functor (global section) between the $\infty$-prestacks of monomorphic inductive Banach quasicoherent presheaves:
\[
\xymatrix@R+6pc@C+0pc{
\mathrm{Ind}^m\mathrm{Banach}_{\Gamma_k}(\mathcal{O}_{X_{R_k,-}})\ar[r]^{\mathrm{global}}\ar[d]\ar[d]\ar[d]\ar[d]\ar[r]\ar[r] &\varphi\mathrm{Ind}^m\mathrm{Banach}_{\Gamma_k}(\{\mathrm{Robba}^\mathrm{extended}_{{R_k,-,I}}\}_I)\ar[d]\ar[d]\ar[d]\ar[d]\\
\mathrm{Ind}^m\mathrm{Banach}_{\Gamma_0}(\mathcal{O}_{X_{\mathbb{Q}_p(p^{1/p^\infty})^{\wedge\flat},-}})\ar[r]^{\mathrm{global}}\ar[r]\ar[r] &\varphi\mathrm{Ind}^m\mathrm{Banach}_{\Gamma_0}(\{\mathrm{Robba}^\mathrm{extended}_{{R_0,-,I}}\}_I).\\  
}
\]
The definition is given by the following:
\[
\xymatrix@R+0pc@C+0pc{
\mathrm{homotopycolimit}_i(\square)(\mathcal{O}_i),  
}
\]
each $\mathcal{O}_i$ is just as $\mathbb{Q}_p\left<C_1,...,C_\ell\right>,\ell=1,2,...$ and $\square$ is the relative diagram of $\infty$-functors.

\item (\text{Proposition}) There is a functor (global section) between the $\infty$-stacks of inductive Banach quasicoherent commutative algebra $E_\infty$ objects:
\[
\xymatrix@R+6pc@C+0pc{
\mathrm{sComm}_\mathrm{simplicial}\mathrm{Ind}\mathrm{Banach}_{\Gamma_k}(\mathcal{O}_{X_{R_k,-}})\ar[d]\ar[d]\ar[d]\ar[d]\ar[r]^{\mathrm{global}}\ar[r]\ar[r] &\mathrm{sComm}_\mathrm{simplicial}\varphi\mathrm{Ind}\mathrm{Banach}_{\Gamma_k}(\{\mathrm{Robba}^\mathrm{extended}_{{R_k,-,I}}\}_I)\ar[d]\ar[d]\ar[d]\ar[d]\\
\mathrm{sComm}_\mathrm{simplicial}\mathrm{Ind}\mathrm{Banach}_{\Gamma_0}(\mathcal{O}_{X_{\mathbb{Q}_p(p^{1/p^\infty})^{\wedge\flat},-}})\ar[r]^{\mathrm{global}}\ar[r]\ar[r] &\mathrm{sComm}_\mathrm{simplicial}\varphi\mathrm{Ind}\mathrm{Banach}_{\Gamma_0}(\{\mathrm{Robba}^\mathrm{extended}_{{R_0,-,I}}\}_I).  
}
\]
The definition is given by the following:
\[
\xymatrix@R+0pc@C+0pc{
\mathrm{homotopycolimit}_i(\square)(\mathcal{O}_i),  
}
\]
each $\mathcal{O}_i$ is just as $\mathbb{Q}_p\left<C_1,...,C_\ell\right>,\ell=1,2,...$ and $\square$ is the relative diagram of $\infty$-functors.

\item (\text{Proposition}) There is a functor (global section) between the $\infty$-prestacks of monomorphic inductive Banach quasicoherent commutative algebra $E_\infty$ objects:
\[\displayindent=-0.4in
\xymatrix@R+6pc@C+0pc{
\mathrm{sComm}_\mathrm{simplicial}\mathrm{Ind}^m\mathrm{Banach}_{\Gamma_k}(\mathcal{O}_{X_{R_k,-}})\ar[d]\ar[d]\ar[d]\ar[d]\ar[r]^{\mathrm{global}}\ar[r]\ar[r] &\mathrm{sComm}_\mathrm{simplicial}\varphi\mathrm{Ind}^m\mathrm{Banach}_{\Gamma_k}(\{\mathrm{Robba}^\mathrm{extended}_{{R_k,-,I}}\}_I)\ar[d]\ar[d]\ar[d]\ar[d]\\
 \mathrm{sComm}_\mathrm{simplicial}\mathrm{Ind}^m\mathrm{Banach}_{\Gamma_0}(\mathcal{O}_{X_{\mathbb{Q}_p(p^{1/p^\infty})^{\wedge\flat},-}})\ar[r]^{\mathrm{global}}\ar[r]\ar[r] &\mathrm{sComm}_\mathrm{simplicial}\varphi\mathrm{Ind}^m\mathrm{Banach}_{\Gamma_0}(\{\mathrm{Robba}^\mathrm{extended}_{{R_0,-,I}}\}_I). 
}
\]
The definition is given by the following:
\[
\xymatrix@R+0pc@C+0pc{
\mathrm{homotopycolimit}_i(\square)(\mathcal{O}_i),  
}
\]
each $\mathcal{O}_i$ is just as $\mathbb{Q}_p\left<C_1,...,C_\ell\right>,\ell=1,2,...$ and $\square$ is the relative diagram of $\infty$-functors.

\item Then parallel as in \cite{LBV} we have a functor (global section) of the de Rham complex after \cite[Definition 5.9, Section 5.2.1]{KKM}:
\[\displayindent=-0.4in
\xymatrix@R+6pc@C+0pc{
\mathrm{deRham}_{\mathrm{sComm}_\mathrm{simplicial}\mathrm{Ind}\mathrm{Banach}_{\Gamma_k}(\mathcal{O}_{X_{R_k,-}})\ar[r]^{\mathrm{global}}}(-)\ar[d]\ar[d]\ar[d]\ar[d]\ar[r]\ar[r] &\mathrm{deRham}_{\mathrm{sComm}_\mathrm{simplicial}\varphi\mathrm{Ind}\mathrm{Banach}_{\Gamma_k}(\{\mathrm{Robba}^\mathrm{extended}_{{R_k,-,I}}\}_I)}(-)\ar[d]\ar[d]\ar[d]\ar[d]\\
\mathrm{deRham}_{\mathrm{sComm}_\mathrm{simplicial}\mathrm{Ind}\mathrm{Banach}_{\Gamma_0}(\mathcal{O}_{X_{\mathbb{Q}_p(p^{1/p^\infty})^{\wedge\flat},-}})\ar[r]^{\mathrm{global}}}(-)\ar[r]\ar[r] &\mathrm{deRham}_{\mathrm{sComm}_\mathrm{simplicial}\varphi\mathrm{Ind}\mathrm{Banach}_{\Gamma_0}(\{\mathrm{Robba}^\mathrm{extended}_{{R_0,-,I}}\}_I)}(-), 
}
\]
\[\displayindent=-0.4in
\xymatrix@R+6pc@C+0pc{
\mathrm{deRham}_{\mathrm{sComm}_\mathrm{simplicial}\mathrm{Ind}^m\mathrm{Banach}_{\Gamma_k}(\mathcal{O}_{X_{R_k,-}})\ar[r]^{\mathrm{global}}}(-)\ar[d]\ar[d]\ar[d]\ar[d]\ar[r]\ar[r] &\mathrm{deRham}_{\mathrm{sComm}_\mathrm{simplicial}\varphi\mathrm{Ind}^m\mathrm{Banach}_{\Gamma_k}(\{\mathrm{Robba}^\mathrm{extended}_{{R_k,-,I}}\}_I)}(-)\ar[d]\ar[d]\ar[d]\ar[d]\\
\mathrm{deRham}_{\mathrm{sComm}_\mathrm{simplicial}\mathrm{Ind}^m\mathrm{Banach}_{\Gamma_0}(\mathcal{O}_{X_{\mathbb{Q}_p(p^{1/p^\infty})^{\wedge\flat},-}})\ar[r]^{\mathrm{global}}}(-)\ar[r]\ar[r] &\mathrm{deRham}_{\mathrm{sComm}_\mathrm{simplicial}\varphi\mathrm{Ind}^m\mathrm{Banach}_{\Gamma_0}(\{\mathrm{Robba}^\mathrm{extended}_{{R_0,-,I}}\}_I)}(-). 
}
\]
The definition is given by the following:
\[
\xymatrix@R+0pc@C+0pc{
\mathrm{homotopycolimit}_i(\square)(\mathcal{O}_i),  
}
\]
each $\mathcal{O}_i$ is just as $\mathbb{Q}_p\left<C_1,...,C_\ell\right>,\ell=1,2,...$ and $\square$ is the relative diagram of $\infty$-functors.

\item Then we have the following a functor (global section) of $K$-group $(\infty,1)$-spectrum from \cite{BGT}:
\[
\xymatrix@R+6pc@C+0pc{
\mathrm{K}^\mathrm{BGT}_{\mathrm{sComm}_\mathrm{simplicial}\mathrm{Ind}\mathrm{Banach}_{\Gamma_k}(\mathcal{O}_{X_{R_k,-}})\ar[r]^{\mathrm{global}}}(-)\ar[d]\ar[d]\ar[d]\ar[d]\ar[r]\ar[r] &\mathrm{K}^\mathrm{BGT}_{\mathrm{sComm}_\mathrm{simplicial}\varphi\mathrm{Ind}\mathrm{Banach}_{\Gamma_k}(\{\mathrm{Robba}^\mathrm{extended}_{{R_k,-,I}}\}_I)}(-)\ar[d]\ar[d]\ar[d]\ar[d]\\
\mathrm{K}^\mathrm{BGT}_{\mathrm{sComm}_\mathrm{simplicial}\mathrm{Ind}\mathrm{Banach}_{\Gamma_0}(\mathcal{O}_{X_{\mathbb{Q}_p(p^{1/p^\infty})^{\wedge\flat},-}})\ar[r]^{\mathrm{global}}}(-)\ar[r]\ar[r] &\mathrm{K}^\mathrm{BGT}_{\mathrm{sComm}_\mathrm{simplicial}\varphi\mathrm{Ind}\mathrm{Banach}_{\Gamma_0}(\{\mathrm{Robba}^\mathrm{extended}_{{R_0,-,I}}\}_I)}(-), 
}
\]
\[
\xymatrix@R+6pc@C+0pc{
\mathrm{K}^\mathrm{BGT}_{\mathrm{sComm}_\mathrm{simplicial}\mathrm{Ind}^m\mathrm{Banach}_{\Gamma_k}(\mathcal{O}_{X_{R_k,-}})\ar[r]^{\mathrm{global}}}(-)\ar[d]\ar[d]\ar[d]\ar[d]\ar[r]\ar[r] &\mathrm{K}^\mathrm{BGT}_{\mathrm{sComm}_\mathrm{simplicial}\varphi\mathrm{Ind}^m\mathrm{Banach}_{\Gamma_k}(\{\mathrm{Robba}^\mathrm{extended}_{{R_k,-,I}}\}_I)}(-)\ar[d]\ar[d]\ar[d]\ar[d]\\
\mathrm{K}^\mathrm{BGT}_{\mathrm{sComm}_\mathrm{simplicial}\mathrm{Ind}^m\mathrm{Banach}_{\Gamma_0}(\mathcal{O}_{X_{\mathbb{Q}_p(p^{1/p^\infty})^{\wedge\flat},-}})\ar[r]^{\mathrm{global}}}(-)\ar[r]\ar[r] &\mathrm{K}^\mathrm{BGT}_{\mathrm{sComm}_\mathrm{simplicial}\varphi\mathrm{Ind}^m\mathrm{Banach}_{\Gamma_0}(\{\mathrm{Robba}^\mathrm{extended}_{{R_0,-,I}}\}_I)}(-). 
}
\]
The definition is given by the following:
\[
\xymatrix@R+0pc@C+0pc{
\mathrm{homotopycolimit}_i(\square)(\mathcal{O}_i),  
}
\]
each $\mathcal{O}_i$ is just as $\mathbb{Q}_p\left<C_1,...,C_\ell\right>,\ell=1,2,...$ and $\square$ is the relative diagram of $\infty$-functors.

\end{itemize}

\

\begin{remark}
\noindent We can certainly consider the quasicoherent sheaves in \cite[Lemma 7.11, Remark 7.12]{1BBK}, therefore all the quasicoherent presheaves and modules will be those satisfying \cite[Lemma 7.11, Remark 7.12]{1BBK} if one would like to consider the the quasicoherent sheaves. That being all as this said, we would believe that the big quasicoherent presheaves are automatically quasicoherent sheaves (namely satisfying the corresponding \v{C}ech $\infty$-descent as in \cite[Section 9.3]{KKM} and \cite[Lemma 7.11, Remark 7.12]{1BBK}) and the corresponding global section functors are automatically equivalence of $\infty$-categories. 
\end{remark}

\

\indent In Clausen-Scholze formalism we have the following\footnote{Certainly the homotopy colimit in the rings side will be within the condensed solid animated analytic rings from \cite{1CS2}.}:

\begin{itemize}
\item (\text{Proposition}) There is a functor (global section) between the $\infty$-prestacks of inductive Banach quasicoherent sheaves:
\[
\xymatrix@R+0pc@C+0pc{
{\mathrm{Modules}_\circledcirc}(\mathcal{O}_{X_{R,-}})\ar[r]^{\mathrm{global}}\ar[r]\ar[r] &\varphi{\mathrm{Modules}_\circledcirc}(\{\mathrm{Robba}^\mathrm{extended}_{{R,-,I}}\}_I).  
}
\]
The definition is given by the following:
\[
\xymatrix@R+0pc@C+0pc{
\mathrm{homotopycolimit}_i({\mathrm{Modules}_\circledcirc}(\mathcal{O}_{X_{R,-}})\ar[r]^{\mathrm{global}}\ar[r]\ar[r] &\varphi{\mathrm{Modules}_\circledcirc}(\{\mathrm{Robba}^\mathrm{extended}_{{R,-,I}}\}_I))(\mathcal{O}_i),  
}
\]
each $\mathcal{O}_i$ is just as $\mathbb{Q}_p\left<C_1,...,C_\ell\right>,\ell=1,2,...$.

\item (\text{Proposition}) There is a functor (global section) between the $\infty$-prestacks of inductive Banach quasicoherent sheaves:
\[
\xymatrix@R+0pc@C+0pc{
{\mathrm{Modules}_\circledcirc}(\mathcal{O}_{X_{R,-}})\ar[r]^{\mathrm{global}}\ar[r]\ar[r] &\varphi{\mathrm{Modules}_\circledcirc}(\{\mathrm{Robba}^\mathrm{extended}_{{R,-,I}}\}_I).  
}
\]
The definition is given by the following:
\[
\xymatrix@R+0pc@C+0pc{
\mathrm{homotopycolimit}_i({\mathrm{Modules}_\circledcirc}(\mathcal{O}_{X_{R,-}})\ar[r]^{\mathrm{global}}\ar[r]\ar[r] &\varphi{\mathrm{Modules}_\circledcirc}(\{\mathrm{Robba}^\mathrm{extended}_{{R,-,I}}\}_I))(\mathcal{O}_i),  
}
\]
each $\mathcal{O}_i$ is just as $\mathbb{Q}_p\left<C_1,...,C_\ell\right>,\ell=1,2,...$.

\item (\text{Proposition}) There is a functor (global section) between the $\infty$-prestacks of inductive Banach quasicoherent commutative algebra $E_\infty$ objects\footnote{Here $\circledcirc=\text{solidquasicoherentsheaves}$.}:
\[
\xymatrix@R+0pc@C+0pc{
\mathrm{sComm}_\mathrm{simplicial}{\mathrm{Modules}_\circledcirc}(\mathcal{O}_{X_{R,-}})\ar[r]^{\mathrm{global}}\ar[r]\ar[r] &\mathrm{sComm}_\mathrm{simplicial}\varphi{\mathrm{Modules}_\circledcirc}(\{\mathrm{Robba}^\mathrm{extended}_{{R,-,I}}\}_I).  
}
\]
The definition is given by the following:
\[
\xymatrix@R+0pc@C+0pc{
\mathrm{homotopycolimit}_i(\mathrm{sComm}_\mathrm{simplicial}{\mathrm{Modules}_\circledcirc}(\mathcal{O}_{X_{R,-}})\ar[r]^{\mathrm{global}}\ar[r]\ar[r] &\mathrm{sComm}_\mathrm{simplicial}\varphi{\mathrm{Modules}_\circledcirc}(\{\mathrm{Robba}^\mathrm{extended}_{{R,-,I}}\}_I))(\mathcal{O}_i),  
}
\]
each $\mathcal{O}_i$ is just as $\mathbb{Q}_p\left<C_1,...,C_\ell\right>,\ell=1,2,...$.
\item Then as in \cite{LBV} we have a functor (global section ) of the de Rham complex after \cite[Definition 5.9, Section 5.2.1]{KKM}\footnote{Here $\circledcirc=\text{solidquasicoherentsheaves}$.}:
\[
\xymatrix@R+0pc@C+0pc{
\mathrm{deRham}_{\mathrm{sComm}_\mathrm{simplicial}{\mathrm{Modules}_\circledcirc}(\mathcal{O}_{X_{R,-}})\ar[r]^{\mathrm{global}}}(-)\ar[r]\ar[r] &\mathrm{deRham}_{\mathrm{sComm}_\mathrm{simplicial}\varphi{\mathrm{Modules}_\circledcirc}(\{\mathrm{Robba}^\mathrm{extended}_{{R,-,I}}\}_I)}(-), 
}
\]
The definition is given by the following:
\[
\xymatrix@R+0pc@C+0pc{
\mathrm{homotopycolimit}_i\\
(\mathrm{deRham}_{\mathrm{sComm}_\mathrm{simplicial}{\mathrm{Modules}_\circledcirc}(\mathcal{O}_{X_{R,-}})\ar[r]^{\mathrm{global}}}(-)\ar[r]\ar[r] &\mathrm{deRham}_{\mathrm{sComm}_\mathrm{simplicial}\varphi{\mathrm{Modules}_\circledcirc}(\{\mathrm{Robba}^\mathrm{extended}_{{R,-,I}}\}_I)}(-))(\mathcal{O}_i),  
}
\]
each $\mathcal{O}_i$ is just as $\mathbb{Q}_p\left<C_1,...,C_\ell\right>,\ell=1,2,...$.\item Then we have the following a functor (global section) of $K$-group $(\infty,1)$-spectrum from \cite{BGT}\footnote{Here $\circledcirc=\text{solidquasicoherentsheaves}$.}:
\[
\xymatrix@R+0pc@C+0pc{
\mathrm{K}^\mathrm{BGT}_{\mathrm{sComm}_\mathrm{simplicial}{\mathrm{Modules}_\circledcirc}(\mathcal{O}_{X_{R,-}})\ar[r]^{\mathrm{global}}}(-)\ar[r]\ar[r] &\mathrm{K}^\mathrm{BGT}_{\mathrm{sComm}_\mathrm{simplicial}\varphi{\mathrm{Modules}_\circledcirc}(\{\mathrm{Robba}^\mathrm{extended}_{{R,-,I}}\}_I)}(-). 
}
\]
The definition is given by the following:
\[\displayindent=-0.4in
\xymatrix@R+0pc@C+0pc{
\mathrm{homotopycolimit}_i(\mathrm{K}^\mathrm{BGT}_{\mathrm{sComm}_\mathrm{simplicial}{\mathrm{Modules}_\circledcirc}(\mathcal{O}_{X_{R,-}})\ar[r]^{\mathrm{global}}}(-)\ar[r]\ar[r] &\mathrm{K}^\mathrm{BGT}_{\mathrm{sComm}_\mathrm{simplicial}\varphi{\mathrm{Modules}_\circledcirc}(\{\mathrm{Robba}^\mathrm{extended}_{{R,-,I}}\}_I)}(-))(\mathcal{O}_i),  
}
\]
each $\mathcal{O}_i$ is just as $\mathbb{Q}_p\left<C_1,...,C_\ell\right>,\ell=1,2,...$.
\end{itemize}

\noindent Now let $R=\mathbb{Q}_p(p^{1/p^\infty})^{\wedge\flat}$ and $R_k=\mathbb{Q}_p(p^{1/p^\infty})^{\wedge}\left<T_1^{\pm 1/p^{\infty}},...,T_k^{\pm 1/p^{\infty}}\right>^\flat$ we have the following Galois theoretic results with naturality along $f:\mathrm{Spa}(\mathbb{Q}_p(p^{1/p^\infty})^{\wedge}\left<T_1^{\pm 1/p^\infty},...,T_k^{\pm 1/p^\infty}\right>^\flat)\rightarrow \mathrm{Spa}(\mathbb{Q}_p(p^{1/p^\infty})^{\wedge\flat})$:

\begin{itemize}
\item (\text{Proposition}) There is a functor (global section) between the $\infty$-prestacks of inductive Banach quasicoherent sheaves\footnote{Here $\circledcirc=\text{solidquasicoherentsheaves}$.}:
\[
\xymatrix@R+6pc@C+0pc{
{\mathrm{Modules}_\circledcirc}(\mathcal{O}_{X_{\mathbb{Q}_p(p^{1/p^\infty})^{\wedge}\left<T_1^{\pm 1/p^\infty},...,T_k^{\pm 1/p^\infty}\right>^\flat,-}})\ar[d]\ar[d]\ar[d]\ar[d] \ar[r]^{\mathrm{global}}\ar[r]\ar[r] &\varphi{\mathrm{Modules}_\circledcirc}(\{\mathrm{Robba}^\mathrm{extended}_{{R_k,-,I}}\}_I) \ar[d]\ar[d]\ar[d]\ar[d].\\
{\mathrm{Modules}_\circledcirc}(\mathcal{O}_{X_{\mathbb{Q}_p(p^{1/p^\infty})^{\wedge\flat},-}})\ar[r]^{\mathrm{global}}\ar[r]\ar[r] &\varphi{\mathrm{Modules}_\circledcirc}(\{\mathrm{Robba}^\mathrm{extended}_{{R_0,-,I}}\}_I).\\ 
}
\]
The definition is given by the following:
\[
\xymatrix@R+0pc@C+0pc{
\mathrm{homotopycolimit}_i(\square)(\mathcal{O}_i),  
}
\]
each $\mathcal{O}_i$ is just as $\mathbb{Q}_p\left<C_1,...,C_\ell\right>,\ell=1,2,...$ and $\square$ is the relative diagram of $\infty$-functors.

\item (\text{Proposition}) There is a functor (global section) between the $\infty$-prestacks of inductive Banach quasicoherent commutative algebra $E_\infty$ objects\footnote{Here $\circledcirc=\text{solidquasicoherentsheaves}$.}:
\[
\xymatrix@R+6pc@C+0pc{
\mathrm{sComm}_\mathrm{simplicial}{\mathrm{Modules}_\circledcirc}(\mathcal{O}_{X_{R_k,-}})\ar[d]\ar[d]\ar[d]\ar[d]\ar[r]^{\mathrm{global}}\ar[r]\ar[r] &\mathrm{sComm}_\mathrm{simplicial}\varphi{\mathrm{Modules}_\circledcirc}(\{\mathrm{Robba}^\mathrm{extended}_{{R_k,-,I}}\}_I)\ar[d]\ar[d]\ar[d]\ar[d]\\
\mathrm{sComm}_\mathrm{simplicial}{\mathrm{Modules}_\circledcirc}(\mathcal{O}_{X_{\mathbb{Q}_p(p^{1/p^\infty})^{\wedge\flat},-}})\ar[r]^{\mathrm{global}}\ar[r]\ar[r] &\mathrm{sComm}_\mathrm{simplicial}\varphi{\mathrm{Modules}_\circledcirc}(\{\mathrm{Robba}^\mathrm{extended}_{{R_0,-,I}}\}_I).  
}
\]
The definition is given by the following:
\[
\xymatrix@R+0pc@C+0pc{
\mathrm{homotopycolimit}_i(\square)(\mathcal{O}_i),  
}
\]
each $\mathcal{O}_i$ is just as $\mathbb{Q}_p\left<C_1,...,C_\ell\right>,\ell=1,2,...$ and $\square$ is the relative diagram of $\infty$-functors.

\item Then as in \cite{LBV} we have a functor (global section) of the de Rham complex after \cite[Definition 5.9, Section 5.2.1]{KKM}\footnote{Here $\circledcirc=\text{solidquasicoherentsheaves}$.}:
\[\displayindent=-0.4in
\xymatrix@R+6pc@C+0pc{
\mathrm{deRham}_{\mathrm{sComm}_\mathrm{simplicial}{\mathrm{Modules}_\circledcirc}(\mathcal{O}_{X_{R_k,-}})\ar[r]^{\mathrm{global}}}(-)\ar[d]\ar[d]\ar[d]\ar[d]\ar[r]\ar[r] &\mathrm{deRham}_{\mathrm{sComm}_\mathrm{simplicial}\varphi{\mathrm{Modules}_\circledcirc}(\{\mathrm{Robba}^\mathrm{extended}_{{R_k,-,I}}\}_I)}(-)\ar[d]\ar[d]\ar[d]\ar[d]\\
\mathrm{deRham}_{\mathrm{sComm}_\mathrm{simplicial}{\mathrm{Modules}_\circledcirc}(\mathcal{O}_{X_{\mathbb{Q}_p(p^{1/p^\infty})^{\wedge\flat},-}})\ar[r]^{\mathrm{global}}}(-)\ar[r]\ar[r] &\mathrm{deRham}_{\mathrm{sComm}_\mathrm{simplicial}\varphi{\mathrm{Modules}_\circledcirc}(\{\mathrm{Robba}^\mathrm{extended}_{{R_0,-,I}}\}_I)}(-), 
}
\]

\item Then we have the following a functor (global section) of $K$-group $(\infty,1)$-spectrum from \cite{BGT}\footnote{Here $\circledcirc=\text{solidquasicoherentsheaves}$.}:
\[
\xymatrix@R+6pc@C+0pc{
\mathrm{K}^\mathrm{BGT}_{\mathrm{sComm}_\mathrm{simplicial}{\mathrm{Modules}_\circledcirc}(\mathcal{O}_{X_{R_k,-}})\ar[r]^{\mathrm{global}}}(-)\ar[d]\ar[d]\ar[d]\ar[d]\ar[r]\ar[r] &\mathrm{K}^\mathrm{BGT}_{\mathrm{sComm}_\mathrm{simplicial}\varphi{\mathrm{Modules}_\circledcirc}(\{\mathrm{Robba}^\mathrm{extended}_{{R_k,-,I}}\}_I)}(-)\ar[d]\ar[d]\ar[d]\ar[d]\\
\mathrm{K}^\mathrm{BGT}_{\mathrm{sComm}_\mathrm{simplicial}{\mathrm{Modules}_\circledcirc}(\mathcal{O}_{X_{\mathbb{Q}_p(p^{1/p^\infty})^{\wedge\flat},-}})\ar[r]^{\mathrm{global}}}(-)\ar[r]\ar[r] &\mathrm{K}^\mathrm{BGT}_{\mathrm{sComm}_\mathrm{simplicial}\varphi{\mathrm{Modules}_\circledcirc}(\{\mathrm{Robba}^\mathrm{extended}_{{R_0,-,I}}\}_I)}(-), 
}
\]

The definition is given by the following:
\[
\xymatrix@R+0pc@C+0pc{
\mathrm{homotopycolimit}_i(\square)(\mathcal{O}_i),  
}
\]
each $\mathcal{O}_i$ is just as $\mathbb{Q}_p\left<C_1,...,C_\ell\right>,\ell=1,2,...$ and $\square$ is the relative diagram of $\infty$-functors.

\end{itemize}

\
\indent Then we consider further equivariance by considering the arithmetic profinite fundamental groups $\Gamma_{\mathbb{Q}_p}$ and $\mathrm{Gal}(\overline{{Q}_p\left<T_1^{\pm 1},...,T_k^{\pm 1}\right>}/R_k)$ through the following diagram:

\[
\xymatrix@R+0pc@C+0pc{
\mathbb{Z}_p^k=\mathrm{Gal}(R_k/{\mathbb{Q}_p(p^{1/p^\infty})^\wedge\left<T_1^{\pm 1},...,T_k^{\pm 1}\right>}) \ar[r]\ar[r] \ar[r]\ar[r] &\Gamma_k:=\mathrm{Gal}(R_k/{\mathbb{Q}_p\left<T_1^{\pm 1},...,T_k^{\pm 1}\right>}) \ar[r] \ar[r]\ar[r] &\Gamma_{\mathbb{Q}_p}.
}
\]

\begin{itemize}
\item (\text{Proposition}) There is a functor (global section) between the $\infty$-prestacks of inductive Banach quasicoherent sheaves\footnote{Here $\circledcirc=\text{solidquasicoherentsheaves}$.}:
\[
\xymatrix@R+6pc@C+0pc{
{\mathrm{Modules}_\circledcirc}_{\Gamma_k}(\mathcal{O}_{X_{\mathbb{Q}_p(p^{1/p^\infty})^{\wedge}\left<T_1^{\pm 1/p^\infty},...,T_k^{\pm 1/p^\infty}\right>^\flat,-}})\ar[d]\ar[d]\ar[d]\ar[d] \ar[r]^{\mathrm{global}}\ar[r]\ar[r] &\varphi{\mathrm{Modules}_\circledcirc}_{\Gamma_k}(\{\mathrm{Robba}^\mathrm{extended}_{{R_k,-,I}}\}_I) \ar[d]\ar[d]\ar[d]\ar[d].\\
{\mathrm{Modules}_\circledcirc}(\mathcal{O}_{X_{\mathbb{Q}_p(p^{1/p^\infty})^{\wedge\flat},-}})\ar[r]^{\mathrm{global}}\ar[r]\ar[r] &\varphi{\mathrm{Modules}_\circledcirc}(\{\mathrm{Robba}^\mathrm{extended}_{{R_0,-,I}}\}_I).\\ 
}
\]
The definition is given by the following:
\[
\xymatrix@R+0pc@C+0pc{
\mathrm{homotopycolimit}_i(\square)(\mathcal{O}_i),  
}
\]
each $\mathcal{O}_i$ is just as $\mathbb{Q}_p\left<C_1,...,C_\ell\right>,\ell=1,2,...$ and $\square$ is the relative diagram of $\infty$-functors.

\item (\text{Proposition}) There is a functor (global section) between the $\infty$-stacks of inductive Banach quasicoherent commutative algebra $E_\infty$ objects\footnote{Here $\circledcirc=\text{solidquasicoherentsheaves}$.}:
\[
\xymatrix@R+6pc@C+0pc{
\mathrm{sComm}_\mathrm{simplicial}{\mathrm{Modules}_\circledcirc}_{\Gamma_k}(\mathcal{O}_{X_{R_k,-}})\ar[d]\ar[d]\ar[d]\ar[d]\ar[r]^{\mathrm{global}}\ar[r]\ar[r] &\mathrm{sComm}_\mathrm{simplicial}\varphi{\mathrm{Modules}_\circledcirc}_{\Gamma_k}(\{\mathrm{Robba}^\mathrm{extended}_{{R_k,-,I}}\}_I)\ar[d]\ar[d]\ar[d]\ar[d]\\
\mathrm{sComm}_\mathrm{simplicial}{\mathrm{Modules}_\circledcirc}_{\Gamma_0}(\mathcal{O}_{X_{\mathbb{Q}_p(p^{1/p^\infty})^{\wedge\flat},-}})\ar[r]^{\mathrm{global}}\ar[r]\ar[r] &\mathrm{sComm}_\mathrm{simplicial}\varphi{\mathrm{Modules}_\circledcirc}_{\Gamma_0}(\{\mathrm{Robba}^\mathrm{extended}_{{R_0,-,I}}\}_I).  
}
\]
The definition is given by the following:
\[
\xymatrix@R+0pc@C+0pc{
\mathrm{homotopycolimit}_i(\square)(\mathcal{O}_i),  
}
\]
each $\mathcal{O}_i$ is just as $\mathbb{Q}_p\left<C_1,...,C_\ell\right>,\ell=1,2,...$ and $\square$ is the relative diagram of $\infty$-functors.

\item Then as in \cite{LBV} we have a functor (global section) of the de Rham complex after \cite[Definition 5.9, Section 5.2.1]{KKM}\footnote{Here $\circledcirc=\text{solidquasicoherentsheaves}$.}:
\[\displayindent=-0.4in
\xymatrix@R+6pc@C+0pc{
\mathrm{deRham}_{\mathrm{sComm}_\mathrm{simplicial}{\mathrm{Modules}_\circledcirc}_{\Gamma_k}(\mathcal{O}_{X_{R_k,-}})\ar[r]^{\mathrm{global}}}(-)\ar[d]\ar[d]\ar[d]\ar[d]\ar[r]\ar[r] &\mathrm{deRham}_{\mathrm{sComm}_\mathrm{simplicial}\varphi{\mathrm{Modules}_\circledcirc}_{\Gamma_k}(\{\mathrm{Robba}^\mathrm{extended}_{{R_k,-,I}}\}_I)}(-)\ar[d]\ar[d]\ar[d]\ar[d]\\
\mathrm{deRham}_{\mathrm{sComm}_\mathrm{simplicial}{\mathrm{Modules}_\circledcirc}_{\Gamma_0}(\mathcal{O}_{X_{\mathbb{Q}_p(p^{1/p^\infty})^{\wedge\flat},-}})\ar[r]^{\mathrm{global}}}(-)\ar[r]\ar[r] &\mathrm{deRham}_{\mathrm{sComm}_\mathrm{simplicial}\varphi{\mathrm{Modules}_\circledcirc}_{\Gamma_0}(\{\mathrm{Robba}^\mathrm{extended}_{{R_0,-,I}}\}_I)}(-), 
}
\]

The definition is given by the following:
\[
\xymatrix@R+0pc@C+0pc{
\mathrm{homotopycolimit}_i(\square)(\mathcal{O}_i),  
}
\]
each $\mathcal{O}_i$ is just as $\mathbb{Q}_p\left<C_1,...,C_\ell\right>,\ell=1,2,...$ and $\square$ is the relative diagram of $\infty$-functors.

\item Then we have the following a functor (global section) of $K$-group $(\infty,1)$-spectrum from \cite{BGT}\footnote{Here $\circledcirc=\text{solidquasicoherentsheaves}$.}:
\[
\xymatrix@R+6pc@C+0pc{
\mathrm{K}^\mathrm{BGT}_{\mathrm{sComm}_\mathrm{simplicial}{\mathrm{Modules}_\circledcirc}_{\Gamma_k}(\mathcal{O}_{X_{R_k,-}})\ar[r]^{\mathrm{global}}}(-)\ar[d]\ar[d]\ar[d]\ar[d]\ar[r]\ar[r] &\mathrm{K}^\mathrm{BGT}_{\mathrm{sComm}_\mathrm{simplicial}\varphi{\mathrm{Modules}_\circledcirc}_{\Gamma_k}(\{\mathrm{Robba}^\mathrm{extended}_{{R_k,-,I}}\}_I)}(-)\ar[d]\ar[d]\ar[d]\ar[d]\\
\mathrm{K}^\mathrm{BGT}_{\mathrm{sComm}_\mathrm{simplicial}{\mathrm{Modules}_\circledcirc}_{\Gamma_0}(\mathcal{O}_{X_{\mathbb{Q}_p(p^{1/p^\infty})^{\wedge\flat},-}})\ar[r]^{\mathrm{global}}}(-)\ar[r]\ar[r] &\mathrm{K}^\mathrm{BGT}_{\mathrm{sComm}_\mathrm{simplicial}\varphi{\mathrm{Modules}_\circledcirc}_{\Gamma_0}(\{\mathrm{Robba}^\mathrm{extended}_{{R_0,-,I}}\}_I)}(-), 
}
\]

The definition is given by the following:
\[
\xymatrix@R+0pc@C+0pc{
\mathrm{homotopycolimit}_i(\square)(\mathcal{O}_i),  
}
\]
each $\mathcal{O}_i$ is just as $\mathbb{Q}_p\left<C_1,...,C_\ell\right>,\ell=1,2,...$ and $\square$ is the relative diagram of $\infty$-functors.\\

\end{itemize}

\begin{proposition}
All the global functors from \cite[Proposition 13.8, Theorem 14.9, Remark 14.10]{1CS2} achieve the equivalences on both sides.	\\
\end{proposition}

\newpage
\subsection{$\infty$-Categorical Analytic Stacks and Descents V}

Here we consider the corresponding archimedean picture, after \cite[Problem A.4, Kedlaya's Lecture]{1CBCKSW}. Recall for any algebraic variety $R$ over $\mathbb{R}$ this $X_R(\mathbb{C})$ is defined to be the corresponding quotient:
\begin{align}
X_{R}(\mathbb{C}):=R(\mathbb{C})\times P^1(\mathbb{C})/\varphi,\\
Y_R(\mathbb{C}):=R(\mathbb{C})\times P^1(\mathbb{C}).	
\end{align}
The Hodge structure is given by $\varphi$. We define the relative version by considering a further algebraic variety over $\mathbb{C}$, say $A$ as in the following:
\begin{align}
X_{R,A}(\mathbb{C}):=R(\mathbb{C})\times P^1(\mathbb{C})\times A(\mathbb{C})/\varphi,\\
Y_{R,A}(\mathbb{C}):=R(\mathbb{C})\times P^1(\mathbb{C})\times A(\mathbb{C}).	
\end{align}

Then by \cite{1BBK} and \cite{1CS2} we have the corresponding $\infty$-category of $\infty$-sheaves of simplicial ind-Banach quasicoherent modules which in our situation will be assumed to the modules in \cite{1BBK}, as well as the corresponding associated Clausen-Scholze spaces:
\begin{align}
X_{R}(\mathbb{C}):=R(\mathbb{C})\times P^1(\mathbb{C})^\blacksquare/\varphi,\\
Y_R(\mathbb{C}):=R(\mathbb{C})\times P^1(\mathbb{C})^\blacksquare.	
\end{align}
\begin{align}
X_{R,A}(\mathbb{C}):=R(\mathbb{C})\times P^1(\mathbb{C})\times A(\mathbb{C})^\blacksquare/\varphi,\\
Y_{R,A}(\mathbb{C}):=R(\mathbb{C})\times P^1(\mathbb{C})\times A(\mathbb{C})^\blacksquare,	
\end{align}
with the $\infty$-category of $\infty$-sheaves of simplicial liquid quasicoherent modules, liquid vector bundles and liquid perfect complexes, with further descent \cite[Proposition 13.8, Theorem 14.9, Remark 14.10]{1CS2}. We call the resulting global sections are the corresponding $c$-equivariant Hodge Modules. Then we have the following direct analogy:

\begin{itemize}
\item (\text{Proposition}) There is an equivalence between the $\infty$-categories of inductive Banach quasicoherent presheaves:
\[
\xymatrix@R+0pc@C+0pc{
\mathrm{Ind}\mathrm{Banach}(\mathcal{O}_{X_{R,A}})\ar[r]^{\mathrm{equi}}\ar[r]\ar[r] &\varphi\mathrm{Ind}\mathrm{Banach}(\mathcal{O}_{Y_{R,A}}).  
}
\]
\item (\text{Proposition}) There is an equivalence between the $\infty$-categories of monomorphic inductive Banach quasicoherent presheaves:
\[
\xymatrix@R+0pc@C+0pc{
\mathrm{Ind}^m\mathrm{Banach}(\mathcal{O}_{X_{R,A}})\ar[r]^{\mathrm{equi}}\ar[r]\ar[r] &\varphi\mathrm{Ind}^m\mathrm{Banach}(\mathcal{O}_{Y_{R,A}}).  
}
\]
\end{itemize}

\begin{itemize}

\item (\text{Proposition}) There is an equivalence between the $\infty$-categories of inductive Banach quasicoherent presheaves:
\[
\xymatrix@R+0pc@C+0pc{
\mathrm{Ind}\mathrm{Banach}(\mathcal{O}_{X_{R,A}})\ar[r]^{\mathrm{equi}}\ar[r]\ar[r] &\varphi\mathrm{Ind}\mathrm{Banach}(\mathcal{O}_{Y_{R,A}}).  
}
\]
\item (\text{Proposition}) There is an equivalence between the $\infty$-categories of monomorphic inductive Banach quasicoherent presheaves:
\[
\xymatrix@R+0pc@C+0pc{
\mathrm{Ind}^m\mathrm{Banach}(\mathcal{O}_{X_{R,A}})\ar[r]^{\mathrm{equi}}\ar[r]\ar[r] &\varphi\mathrm{Ind}^m\mathrm{Banach}(\mathcal{O}_{Y_{R,A}}).  
}
\]
\item (\text{Proposition}) There is an equivalence between the $\infty$-categories of inductive Banach quasicoherent commutative algebra $E_\infty$ objects:
\[
\xymatrix@R+0pc@C+0pc{
\mathrm{sComm}_\mathrm{simplicial}\mathrm{Ind}\mathrm{Banach}(\mathcal{O}_{X_{R,A}})\ar[r]^{\mathrm{equi}}\ar[r]\ar[r] &\mathrm{sComm}_\mathrm{simplicial}\varphi\mathrm{Ind}\mathrm{Banach}(\mathcal{O}_{Y_{R,A}}).  
}
\]
\item (\text{Proposition}) There is an equivalence between the $\infty$-categories of monomorphic inductive Banach quasicoherent commutative algebra $E_\infty$ objects:
\[
\xymatrix@R+0pc@C+0pc{
\mathrm{sComm}_\mathrm{simplicial}\mathrm{Ind}^m\mathrm{Banach}(\mathcal{O}_{X_{R,A}})\ar[r]^{\mathrm{equi}}\ar[r]\ar[r] &\mathrm{sComm}_\mathrm{simplicial}\varphi\mathrm{Ind}^m\mathrm{Banach}(\mathcal{O}_{Y_{R,A}}).  
}
\]

\item Then parallel as in \cite{LBV} we have the equivalence of the de Rham complex after \cite[Definition 5.9, Section 5.2.1]{KKM}:
\[
\xymatrix@R+0pc@C+0pc{
\mathrm{deRham}_{\mathrm{sComm}_\mathrm{simplicial}\mathrm{Ind}\mathrm{Banach}(\mathcal{O}_{X_{R,A}})\ar[r]^{\mathrm{equi}}}(-)\ar[r]\ar[r] &\mathrm{deRham}_{\mathrm{sComm}_\mathrm{simplicial}\varphi\mathrm{Ind}\mathrm{Banach}(\mathcal{O}_{Y_{R,A}})}(-), 
}
\]
\[
\xymatrix@R+0pc@C+0pc{
\mathrm{deRham}_{\mathrm{sComm}_\mathrm{simplicial}\mathrm{Ind}^m\mathrm{Banach}(\mathcal{O}_{X_{R,A}})\ar[r]^{\mathrm{equi}}}(-)\ar[r]\ar[r] &\mathrm{deRham}_{\mathrm{sComm}_\mathrm{simplicial}\varphi\mathrm{Ind}^m\mathrm{Banach}(\mathcal{O}_{Y_{R,A}})}(-). 
}
\]

\item Then we have the following equivalence of $K$-group $(\infty,1)$-spectrum from \cite{BGT}:
\[
\xymatrix@R+0pc@C+0pc{
\mathrm{K}^\mathrm{BGT}_{\mathrm{sComm}_\mathrm{simplicial}\mathrm{Ind}\mathrm{Banach}(\mathcal{O}_{X_{R,A}})\ar[r]^{\mathrm{equi}}}(-)\ar[r]\ar[r] &\mathrm{K}^\mathrm{BGT}_{\mathrm{sComm}_\mathrm{simplicial}\varphi\mathrm{Ind}\mathrm{Banach}(\mathcal{O}_{Y_{R,A}})}(-), 
}
\]
\[
\xymatrix@R+0pc@C+0pc{
\mathrm{K}^\mathrm{BGT}_{\mathrm{sComm}_\mathrm{simplicial}\mathrm{Ind}^m\mathrm{Banach}(\mathcal{O}_{X_{R,A}})\ar[r]^{\mathrm{equi}}}(-)\ar[r]\ar[r] &\mathrm{K}^\mathrm{BGT}_{\mathrm{sComm}_\mathrm{simplicial}\varphi\mathrm{Ind}^m\mathrm{Banach}(\mathcal{O}_{Y_{R,A}})}(-). 
}
\]
\end{itemize}

\begin{assumption}\label{assumtionpresheaves}
All the functors of modules or algebras below are Clausen-Scholze sheaves \cite[Proposition 13.8, Theorem 14.9, Remark 14.10]{1CS2}.	
\end{assumption}

\begin{itemize}
\item (\text{Proposition}) There is an equivalence between the $\infty$-categories of inductive liquid sheaves:
\[
\xymatrix@R+0pc@C+0pc{
\mathrm{Module}_\circledcirc(\mathcal{O}_{X_{R,A}})\ar[r]^{\mathrm{equi}}\ar[r]\ar[r] &\varphi\mathrm{Module}_\circledcirc(\mathcal{O}_{Y_{R,A}}).  
}
\]
\end{itemize}

\begin{itemize}

\item (\text{Proposition}) There is an equivalence between the $\infty$-categories of inductive Banach quasicoherent commutative algebra $E_\infty$ objects:
\[\displayindent=-0.4in
\xymatrix@R+0pc@C+0pc{
\mathrm{sComm}_\mathrm{simplicial}\mathrm{Module}_{\text{liquidquasicoherentsheaves}}(\mathcal{O}_{X_{R,A}})\ar[r]^{\mathrm{equi}}\ar[r]\ar[r] &\mathrm{sComm}_\mathrm{simplicial}\varphi\mathrm{Module}_{\text{liquidquasicoherentsheaves}}(\mathcal{O}_{Y_{R,A}}).  
}
\]

\item Then as in \cite{LBV} we have the equivalence of the de Rham complex after \cite[Definition 5.9, Section 5.2.1]{KKM}\footnote{Here $\circledcirc=\text{liquidquasicoherentsheaves}$.}:
\[
\xymatrix@R+0pc@C+0pc{
\mathrm{deRham}_{\mathrm{sComm}_\mathrm{simplicial}\mathrm{Module}_\circledcirc(\mathcal{O}_{X_{R,A}})\ar[r]^{\mathrm{equi}}}(-)\ar[r]\ar[r] &\mathrm{deRham}_{\mathrm{sComm}_\mathrm{simplicial}\varphi\mathrm{Module}_\circledcirc(\mathcal{O}_{Y_{R,A}})}(-). 
}
\]

\item Then we have the following equivalence of $K$-group $(\infty,1)$-spectrum from \cite{BGT}\footnote{Here $\circledcirc=\text{liquidquasicoherentsheaves}$.}:
\[
\xymatrix@R+0pc@C+0pc{
\mathrm{K}^\mathrm{BGT}_{\mathrm{sComm}_\mathrm{simplicial}\mathrm{Module}_\circledcirc(\mathcal{O}_{X_{R,A}})\ar[r]^{\mathrm{equi}}}(-)\ar[r]\ar[r] &\mathrm{K}^\mathrm{BGT}_{\mathrm{sComm}_\mathrm{simplicial}\varphi\mathrm{Module}_\circledcirc(\mathcal{O}_{Y_{R,A}})}(-). 
}
\]
\end{itemize}

\newpage

\chapter{$\infty$-Categorical Approaches to Hodge-Iwasawa Theory II}

\section{Introduction to the Interactions among Motives}

\subsection{\text{Equivariant relative $p$-adic Hodge Theory}}

\noindent{\text{Equivariant relative $p$-adic Hodge Theory}}

\begin{itemize}

\item<1-> The corresponding $P$-objects are interesting, but in general are not that easy to study, especially we consider for instance those ring defined over $\mathbb{Q}_p$, let it alone if one would like to consider the categories of the complexes of such objects. 

\item<2-> We choose to consider the corresponding embedding of such objects into the categories of Frobenius sheaves with coefficients in $P$ after Kedlaya-Liu \cite{KL1}, \cite{KL2}. Again we expect everything will be more convenient to handle in the category of $(\varphi,\Gamma)$-modules.
	
\item<3-> Working over $R$ now a uniform Banach algebra with further structure of an adic ring over $\mathbb{F}_p$. And we assume that $R$ is perfect.
Let $\mathrm{Robba}^{\mathrm{extended}}_{I,R}$ be the Robba sheaves defined by Kedlaya-Liu \cite{KL1}, \cite{KL2}, with respect to some interval $I\subset (0,\infty)$, which are Fr\'echet completions of the ring of Witt vector of $R$ with respect to the Gauss norms induced from the norm on $R$. Here we consider the following assumption:

\begin{assumption}
We now assume $R$ comes from the local chart of a rigid space over $\mathbb{Q}_p$. \footnote{This could be made more general, but at this moment let us be closer to classical $p$-adic Hodge Theory.} This will give us the chance to consider the following period rings from \cite{1Sch2} and \cite[Definition 8.6.5]{KL2}:
\begin{align}
B^+_{R,\mathrm{dR}},B_{R,\mathrm{dR}},\\
\mathcal{O}B^+_{R,\mathrm{dR}},\mathcal{O}B_{R,\mathrm{dR}}	
\end{align}
over:
\begin{align}
B^+_{\mathbb{Q}_p(p^{1/p^\infty})^{\wedge,\flat},\mathrm{dR}},B_{\mathbb{Q}_p(p^{1/p^\infty})^{\wedge,\flat},\mathrm{dR}},\\
\mathcal{O}B^+_{\mathbb{Q}_p(p^{1/p^\infty})^{\wedge,\flat},\mathrm{dR}},\mathcal{O}B_{\mathbb{Q}_p(p^{1/p^\infty})^{\wedge,\flat},\mathrm{dR}}	
\end{align}	
where the smaller rings contain element $t=\mathrm{log}([1+\overline{\pi}])$. As in \cite{BCM} we can take the corresponding self $q$-th power product. Then we have by taking the corresponding self $q$-th power product following:
\begin{align}
B^+_{R,\mathrm{dR},q},B_{R,\mathrm{dR},q},\\
\mathcal{O}B^+_{R,\mathrm{dR},q},\mathcal{O}B_{R,\mathrm{dR},q}	
\end{align}
over:
\begin{align}
B^+_{\mathbb{Q}_p(p^{1/p^\infty})^{\wedge,\flat},\mathrm{dR},q},B_{\mathbb{Q}_p(p^{1/p^\infty})^{\wedge,\flat},\mathrm{dR},q},\\
\mathcal{O}B^+_{\mathbb{Q}_p(p^{1/p^\infty})^{\wedge,\flat},\mathrm{dR},q},\mathcal{O}B_{\mathbb{Q}_p(p^{1/p^\infty})^{\wedge,\flat},\mathrm{dR},q}.	
\end{align}
Here we have the action from the product of arithmetic profinite fundamental groups and the product of Frobenius operators.
\end{assumption}

\item<4-> Following Carter-Kedlaya-Z\'abr\'adi and Pal-Z\'abr\'adi \cite{1CKZ} and \cite{1PZ}, taking suitable interval one can define the corresponding Robba rings $\mathrm{Robba}^{\mathrm{extended},q}_{r,R}$, $\mathrm{Robba}^{\mathrm{extended},q}_{\infty,R}$ and the corresponding full Robba ring $\mathrm{Robba}^{\mathrm{extended},q}_{R}$ by the corresponding self $q$-th power product. 

\item<5-> We work in the category of Banach and ind-Fr\'echet spaces, which are commutative. Our generalization comes from those Banach reduced affinoid algebras $A$ over $\mathbb{Q}_p$.

\end{itemize}

\noindent{\text{Equivariant relative $p$-adic Hodge Theory}}
\begin{itemize}

\item<1-> The $p$-adic functional analysis produces us some manageable structures within our study of relative $p$-adic Hodge theory, generalizing the original $p$-adic functional analytic framework of Kedlaya-Liu \cite{KL1}, \cite{KL2}.

\item<2-> Starting from Kedlaya-Liu's period rings after taking product\footnote{When we are talking about $q$-th power as in this chapter, the radius $r$ is then multiradius which are allowed to be different in different components, and the interval $I$ is then multiinterval which are allowed to be different in different components.}, 
\begin{align}
&\mathrm{Robba}^{\mathrm{extended},q}_{\infty,R},\mathrm{Robba}^{\mathrm{extended},q}_{I,R},\mathrm{Robba}^{\mathrm{extended},q}_{r,R},\mathrm{Robba}^{\mathrm{extended},q}_{R},	\mathrm{Robba}^{\mathrm{extended},q}_{{\mathrm{int},r},R},\\
&\mathrm{Robba}^{\mathrm{extended},q}_{{\mathrm{int}},R},\mathrm{Robba}^{\mathrm{extended},q}_{{\mathrm{bd},r},R},\mathrm{Robba}^{\mathrm{extended},q}_{{\mathrm{bd}},R}
\end{align}
we can form the corresponding $A$-relative of the period rings\footnote{Taking products over $\mathbb{Q}_p$.}:
\begin{align}
&\mathrm{Robba}^{\mathrm{extended},q}_{\infty,R,A},\mathrm{Robba}^{\mathrm{extended},q}_{I,R,A},\mathrm{Robba}^{\mathrm{extended},q}_{r,R,A},\mathrm{Robba}^{\mathrm{extended},q}_{R,A},	\mathrm{Robba}^{\mathrm{extended},q}_{{\mathrm{int},r},R,A},\\
&\mathrm{Robba}^{\mathrm{extended},q}_{\mathrm{int},R,A},\mathrm{Robba}^{\mathrm{extended},q}_{{\mathrm{bd},r},R,A},\mathrm{Robba}^{\mathrm{extended},q}_{{\mathrm{bd}},R,A}.	
\end{align} 
\item<3-> (\text{Remark}) There should be also many interesting contexts, for instance consider a finitely generated abelian group $G$, one can consider the group rings: 
\begin{align}
\mathrm{Robba}^{\mathrm{extended},q}_{I,R}[G].	
\end{align}
\item<4-> And then consider the completion living inside the corresponding infinite direct sum Banach modules 
\begin{align}
\bigoplus\mathrm{Robba}^{\mathrm{extended},q}_{I,R},	
\end{align}
over the corresponding period rings:
\begin{align}
\overline{\mathrm{Robba}^{\mathrm{extended},q}_{I,R}[G]}.	
\end{align}
Then we take suitable intersection and union one can have possibly some interesting period rings $\overline{\mathrm{Robba}^{\mathrm{extended},q}_{r,R}[G]}$ and $\overline{\mathrm{Robba}^{\mathrm{extended},q}_{R}[G]}$.
\end{itemize}

\noindent{\text{Equivariant relative $p$-adic Hodge Theory}}
\begin{itemize}
\item<1-> The equivariant period rings in the situations we mentioned above carry relative multi-Frobenius action $\varphi_q$ induced from the Witt vectors. 

\item<2-> They carry the corresponding Banach or (ind-)Fr\'echet spaces structures. So we can generalize the corresponding Kedlaya-Liu's construction to the following situations (here let $G$ be finite):

\item<3-> We can then consider the corresponding completed Frobenius modules over the rings in the equivariant setting. To be more precise over:
\begin{align}
\overline{\mathrm{Robba}^{\mathrm{extended},q}_{R,A}[G]},\Omega_{\mathrm{int},R,A},\Omega_{R,A},\mathrm{Robba}^{\mathrm{extended},q}_{R,A},\mathrm{Robba}^{\mathrm{extended},q}_{\mathrm{bd},R,A}	
\end{align}
one considers the Frobenius modules finite locally free.

\item<4->  With the corresponding finite locally free models over
\begin{align}
\overline{\mathrm{Robba}^{\mathrm{extended},q}_{r,R,A}[G]},\mathrm{Robba}^{\mathrm{extended},q}_{r,R,A},\mathrm{Robba}^{\mathrm{extended},q}_{{\mathrm{bd},r},R,A},	
\end{align} 
again carrying the corresponding semilinear Frobenius structures, where $r$ could be $\infty$.

\item<5-> One also consider families of Frobenius modules over 
\begin{align}
\overline{\mathrm{Robba}^{\mathrm{extended},q}_{I,R,A}[G]},\mathrm{Robba}^{\mathrm{extended},q}_{I,R,A},	
\end{align} 
in glueing fashion with obvious cocycle condition with respect to three multi-intervals $I\subset J\subset K$. These are called the corresponding Frobenius bundles.
	
\end{itemize}

\newpage

\section{Analytic $\infty$-Categorical Functional Analytic Hodge-Iwasawa Modules}

\subsection{$\infty$-Categorical Analytic Stacks and Descents I}

\noindent We now make the corresponding discussion after our previous work \cite{T2} on the homotopical functional analysis after many projects \cite{1BBBK}, \cite{1BBK}, \cite{BBM}, \cite{1BK} , \cite{1CS1}, \cite{1CS2}, \cite{KKM}. We choose to work over the Bambozzi-Kremnizer space \cite{1BK} attached to the corresponding Banach rings in our work after \cite{1BBBK}, \cite{1BBK}, \cite{BBM}, \cite{1BK}, \cite{KKM}. Note that what is happening is that attached to any Banach ring over $\mathbb{Q}_p$, say $B$, we attach a $(\infty,1)-$stack $\mathcal{X}(B)$ fibered over (in the sense of $\infty$-groupoid, and up to taking the corresponding opposite categories) after \cite{1BBBK}, \cite{1BBK}, \cite{BBM}, \cite{1BK}, \cite{KKM}:
\begin{align}
\mathrm{sComm}\mathrm{Simp}\mathrm{Ind}\mathrm{Ban}_{\mathbb{Q}_p},	
\end{align}
with 
\begin{align}
\mathrm{sComm}\mathrm{Simp}\mathrm{Ind}^m\mathrm{Ban}_{\mathbb{Q}_p}.	
\end{align}
associated with a $(\infty,1)$-ring object $\mathcal{O}_{\mathcal{X}(B)}$, such that we have the corresponding under the basic derived rational localization $\infty$-Grothendieck site
\begin{center}
 $(\mathcal{X}(B), \mathcal{O}_{\mathcal{X}(B),\mathrm{drl}})$ 
\end{center}
carrying the homotopical epimorphisms as the corresponding topology.

\begin{itemize}
\item By using this framework (certainly one can also consider \cite{1CS1} and \cite{1CS2} as the foundations, as in \cite{LBV}), we have the $\infty$-stack after Kedlaya-Liu \cite{KL1}, \cite{KL2}. Here in the following let $A$ be any Banach ring over $\mathbb{Q}_p$.

\item Generalizing Kedlaya-Liu's construction in \cite{KL1}, \cite{KL2} of the adic Fargues-Fontaine space we have a quotient (by using power $\varphi_q$ of the Frobenius operator) $X_{R,A,q}$ of the space 
\begin{align}
Y_{R,A,q}:=\bigcup_{I,\text{multi}}\mathcal{X}(\mathrm{Robba}^{\mathrm{extended},q}_{R,I,A}).	
\end{align}

\item This is a locally ringed space $(X_{R,A,q},\mathcal{O}_{X_{R,A,q}})$, so one can consider the stable $\infty$-category $\mathrm{Ind}\mathrm{Banach}(\mathcal{O}_{X_{R,A,q}}) $ which is the $\infty$-category of all the $\mathcal{O}_{X_{R,A,q}}$-sheaves of inductive Banach modules over $X_{R,A,q}$. We have the parallel categories for $Y_{R,A,q}$, namely $\varphi\mathrm{Ind}\mathrm{Banach}(\mathcal{O}_{X_{R,A,q}})$ and so on. Here we only consider presheaves.

\item This is a locally ringed space $(X_{R,A,q},\mathcal{O}_{X_{R,A,q}})$, so one can consider the stable $\infty$-category $\mathrm{Ind}^m\mathrm{Banach}(\mathcal{O}_{X_{R,A,q}}) $ which is the $\infty$-category of all the $\mathcal{O}_{X_{R,A,q}}$-sheaves of inductive monomorphic Banach modules over $X_{R,A,q}$. We have the parallel categories for $Y_{R,A,q}$, namely $\varphi\mathrm{Ind}^m\mathrm{Banach}(\mathcal{O}_{X_{R,A,q}})$ and so on. Here we only consider presheaves.

\begin{assumption}\label{assumtionpresheaves}
All the functors of modules or algebras below are presheaves.	
\end{assumption}

\item In this context one can consider the $K$-theory as in the scheme situation by using the ideas and constructions from Blumberg-Gepner-Tabuada \cite{BGT}. Moreover we can study the Hodge Theory.

\item We expect that one can study among these big categories to find interesting relationships, since this should give us the right understanding of the $p$-adic Hodge theory. The corresponding pseudocoherent version comparison could be expected to be deduced as in Kedlaya-Liu's work \cite{KL1}, \cite{KL2}.

\item (\text{Proposition}) There is an equivalence between the $\infty$-categories of inductive Banach quasicoherent presheaves:
\[
\xymatrix@R+0pc@C+0pc{
\mathrm{Ind}\mathrm{Banach}(\mathcal{O}_{X_{R,A,q}})\ar[r]^{\mathrm{equi}}\ar[r]\ar[r] &\varphi_q\mathrm{Ind}\mathrm{Banach}(\mathcal{O}_{Y_{R,A,q}}).  
}
\]
\item (\text{Proposition}) There is an equivalence between the $\infty$-categories of monomorphic inductive Banach quasicoherent presheaves:
\[
\xymatrix@R+0pc@C+0pc{
\mathrm{Ind}^m\mathrm{Banach}(\mathcal{O}_{X_{R,A,q}})\ar[r]^{\mathrm{equi}}\ar[r]\ar[r] &\varphi_q\mathrm{Ind}^m\mathrm{Banach}(\mathcal{O}_{Y_{R,A,q}}).  
}
\]
\item (\text{Proposition}) There is an equivalence between the $\infty$-categories of inductive Banach quasicoherent commutative algebra $E_\infty$ objects:
\[
\xymatrix@R+0pc@C+0pc{
\mathrm{sComm}_\mathrm{simplicial}\mathrm{Ind}\mathrm{Banach}(\mathcal{O}_{X_{R,A,q}})\ar[r]^{\mathrm{equi}}\ar[r]\ar[r] &\mathrm{sComm}_\mathrm{simplicial}\varphi_q\mathrm{Ind}\mathrm{Banach}(\mathcal{O}_{Y_{R,A,q}}).  
}
\]
\item (\text{Proposition}) There is an equivalence between the $\infty$-categories of monomorphic inductive Banach quasicoherent commutative algebra $E_\infty$ objects:
\[
\xymatrix@R+0pc@C+0pc{
\mathrm{sComm}_\mathrm{simplicial}\mathrm{Ind}^m\mathrm{Banach}(\mathcal{O}_{X_{R,A,q}})\ar[r]^{\mathrm{equi}}\ar[r]\ar[r] &\mathrm{sComm}_\mathrm{simplicial}\varphi_q\mathrm{Ind}^m\mathrm{Banach}(\mathcal{O}_{Y_{R,A,q}}).  
}
\]

\item Then parallel as in \cite{LBV} we have the equivalence of the de Rham complex after \cite[Definition 5.9, Section 5.2.1]{KKM}:
\[
\xymatrix@R+0pc@C+0pc{
\mathrm{deRham}_{\mathrm{sComm}_\mathrm{simplicial}\mathrm{Ind}\mathrm{Banach}(\mathcal{O}_{X_{R,A,q}})\ar[r]^{\mathrm{equi}}}(-)\ar[r]\ar[r] &\mathrm{deRham}_{\mathrm{sComm}_\mathrm{simplicial}\varphi_q\mathrm{Ind}\mathrm{Banach}(\mathcal{O}_{Y_{R,A,q}})}(-), 
}
\]
\[
\xymatrix@R+0pc@C+0pc{
\mathrm{deRham}_{\mathrm{sComm}_\mathrm{simplicial}\mathrm{Ind}^m\mathrm{Banach}(\mathcal{O}_{X_{R,A,q}})\ar[r]^{\mathrm{equi}}}(-)\ar[r]\ar[r] &\mathrm{deRham}_{\mathrm{sComm}_\mathrm{simplicial}\varphi_q\mathrm{Ind}^m\mathrm{Banach}(\mathcal{O}_{Y_{R,A,q}})}(-). 
}
\]

\item Then we have the following equivalence of $K$-group $(\infty,1)$-spectrum from \cite{BGT}:
\[
\xymatrix@R+0pc@C+0pc{
\mathrm{K}^\mathrm{BGT}_{\mathrm{sComm}_\mathrm{simplicial}\mathrm{Ind}\mathrm{Banach}(\mathcal{O}_{X_{R,A,q}})\ar[r]^{\mathrm{equi}}}(-)\ar[r]\ar[r] &\mathrm{K}^\mathrm{BGT}_{\mathrm{sComm}_\mathrm{simplicial}\varphi_q\mathrm{Ind}\mathrm{Banach}(\mathcal{O}_{Y_{R,A,q}})}(-), 
}
\]
\[
\xymatrix@R+0pc@C+0pc{
\mathrm{K}^\mathrm{BGT}_{\mathrm{sComm}_\mathrm{simplicial}\mathrm{Ind}^m\mathrm{Banach}(\mathcal{O}_{X_{R,A,q}})\ar[r]^{\mathrm{equi}}}(-)\ar[r]\ar[r] &\mathrm{K}^\mathrm{BGT}_{\mathrm{sComm}_\mathrm{simplicial}\varphi_q\mathrm{Ind}^m\mathrm{Banach}(\mathcal{O}_{Y_{R,A,q}})}(-). 
}
\]
\end{itemize}

\noindent Now let $R=\mathbb{Q}_p(p^{1/p^\infty})^{\wedge\flat}$ and $R_k=\mathbb{Q}_p(p^{1/p^\infty})^{\wedge}\left<T_1^{\pm 1/p^{\infty}},...,T_k^{\pm 1/p^{\infty}}\right>^\flat$ we have the following Galois theoretic results with naturality along $f:\mathrm{Spa}(\mathbb{Q}_p(p^{1/p^\infty})^{\wedge}\left<T_1^{\pm 1/p^\infty},...,T_k^{\pm 1/p^\infty}\right>^\flat)\rightarrow \mathrm{Spa}(\mathbb{Q}_p(p^{1/p^\infty})^{\wedge\flat})$:

\begin{itemize}
\item (\text{Proposition}) There is an equivalence between the $\infty$-categories of inductive Banach quasicoherent presheaves:
\[
\xymatrix@R+6pc@C+0pc{
\mathrm{Ind}\mathrm{Banach}(\mathcal{O}_{X_{\mathbb{Q}_p(p^{1/p^\infty})^{\wedge}\left<T_1^{\pm 1/p^\infty},...,T_k^{\pm 1/p^\infty}\right>^\flat,A,q}})\ar[d]\ar[d]\ar[d]\ar[d] \ar[r]^{\mathrm{equi}}\ar[r]\ar[r] &\varphi_q\mathrm{Ind}\mathrm{Banach}(\mathcal{O}_{Y_{\mathbb{Q}_p(p^{1/p^\infty})^{\wedge}\left<T_1^{\pm 1/p^\infty},...,T_k^{\pm 1/p^\infty}\right>^\flat,A,q}}) \ar[d]\ar[d]\ar[d]\ar[d].\\
\mathrm{Ind}\mathrm{Banach}(\mathcal{O}_{X_{\mathbb{Q}_p(p^{1/p^\infty})^{\wedge\flat},A,q}})\ar[r]^{\mathrm{equi}}\ar[r]\ar[r] &\varphi_q\mathrm{Ind}\mathrm{Banach}(\mathcal{O}_{Y_{\mathbb{Q}_p(p^{1/p^\infty})^{\wedge\flat},A,q}}).\\ 
}
\]
\item (\text{Proposition}) There is an equivalence between the $\infty$-categories of monomorphic inductive Banach quasicoherent presheaves:
\[
\xymatrix@R+6pc@C+0pc{
\mathrm{Ind}^m\mathrm{Banach}(\mathcal{O}_{X_{R_k,A,q}})\ar[r]^{\mathrm{equi}}\ar[d]\ar[d]\ar[d]\ar[d]\ar[r]\ar[r] &\varphi_q\mathrm{Ind}^m\mathrm{Banach}(\mathcal{O}_{Y_{R_k,A,q}})\ar[d]\ar[d]\ar[d]\ar[d]\\
\mathrm{Ind}^m\mathrm{Banach}(\mathcal{O}_{X_{\mathbb{Q}_p(p^{1/p^\infty})^{\wedge\flat},A,q}})\ar[r]^{\mathrm{equi}}\ar[r]\ar[r] &\varphi_q\mathrm{Ind}^m\mathrm{Banach}(\mathcal{O}_{Y_{\mathbb{Q}_p(p^{1/p^\infty})^{\wedge\flat},A,q}}).\\  
}
\]
\item (\text{Proposition}) There is an equivalence between the $\infty$-categories of inductive Banach quasicoherent commutative algebra $E_\infty$ objects:
\[
\xymatrix@R+6pc@C+0pc{
\mathrm{sComm}_\mathrm{simplicial}\mathrm{Ind}\mathrm{Banach}(\mathcal{O}_{X_{R_k,A,q}})\ar[d]\ar[d]\ar[d]\ar[d]\ar[r]^{\mathrm{equi}}\ar[r]\ar[r] &\mathrm{sComm}_\mathrm{simplicial}\varphi_q\mathrm{Ind}\mathrm{Banach}(\mathcal{O}_{Y_{R_k,A,q}})\ar[d]\ar[d]\ar[d]\ar[d]\\
\mathrm{sComm}_\mathrm{simplicial}\mathrm{Ind}\mathrm{Banach}(\mathcal{O}_{X_{\mathbb{Q}_p(p^{1/p^\infty})^{\wedge\flat},A,q}})\ar[r]^{\mathrm{equi}}\ar[r]\ar[r] &\mathrm{sComm}_\mathrm{simplicial}\varphi_q\mathrm{Ind}\mathrm{Banach}(\mathcal{O}_{Y_{\mathbb{Q}_p(p^{1/p^\infty})^{\wedge\flat},A,q}})  
}
\]
\item (\text{Proposition}) There is an equivalence between the $\infty$-categories of monomorphic inductive Banach quasicoherent commutative algebra $E_\infty$ objects:
\[
\xymatrix@R+6pc@C+0pc{
\mathrm{sComm}_\mathrm{simplicial}\mathrm{Ind}^m\mathrm{Banach}(\mathcal{O}_{X_{R_k,A,q}})\ar[d]\ar[d]\ar[d]\ar[d]\ar[r]^{\mathrm{equi}}\ar[r]\ar[r] &\mathrm{sComm}_\mathrm{simplicial}\varphi_q\mathrm{Ind}^m\mathrm{Banach}(\mathcal{O}_{Y_{R_k,A,q}})\ar[d]\ar[d]\ar[d]\ar[d]\\
 \mathrm{sComm}_\mathrm{simplicial}\mathrm{Ind}^m\mathrm{Banach}(\mathcal{O}_{X_{\mathbb{Q}_p(p^{1/p^\infty})^{\wedge\flat},A,q}})\ar[r]^{\mathrm{equi}}\ar[r]\ar[r] &\mathrm{sComm}_\mathrm{simplicial}\varphi_q\mathrm{Ind}^m\mathrm{Banach}(\mathcal{O}_{Y_{\mathbb{Q}_p(p^{1/p^\infty})^{\wedge\flat},A,q}}) 
}
\]

\item Then parallel as in \cite{LBV} we have the equivalence of the de Rham complex after \cite[Definition 5.9, Section 5.2.1]{KKM}:
\[\displayindent=-0.4in
\xymatrix@R+6pc@C+0pc{
\mathrm{deRham}_{\mathrm{sComm}_\mathrm{simplicial}\mathrm{Ind}\mathrm{Banach}(\mathcal{O}_{X_{R_k,A,q}})\ar[r]^{\mathrm{equi}}}(-)\ar[d]\ar[d]\ar[d]\ar[d]\ar[r]\ar[r] &\mathrm{deRham}_{\mathrm{sComm}_\mathrm{simplicial}\varphi_q\mathrm{Ind}\mathrm{Banach}(\mathcal{O}_{Y_{R_k,A,q}})}(-)\ar[d]\ar[d]\ar[d]\ar[d]\\
\mathrm{deRham}_{\mathrm{sComm}_\mathrm{simplicial}\mathrm{Ind}\mathrm{Banach}(\mathcal{O}_{X_{\mathbb{Q}_p(p^{1/p^\infty})^{\wedge\flat},A,q}})\ar[r]^{\mathrm{equi}}}(-)\ar[r]\ar[r] &\mathrm{deRham}_{\mathrm{sComm}_\mathrm{simplicial}\varphi_q\mathrm{Ind}\mathrm{Banach}(\mathcal{O}_{Y_{\mathbb{Q}_p(p^{1/p^\infty})^{\wedge\flat},A,q}})}(-) 
}
\]
\[\displayindent=-0.4in
\xymatrix@R+6pc@C+0pc{
\mathrm{deRham}_{\mathrm{sComm}_\mathrm{simplicial}\mathrm{Ind}^m\mathrm{Banach}(\mathcal{O}_{X_{R_k,A,q}})\ar[r]^{\mathrm{equi}}}(-)\ar[d]\ar[d]\ar[d]\ar[d]\ar[r]\ar[r] &\mathrm{deRham}_{\mathrm{sComm}_\mathrm{simplicial}\varphi_q\mathrm{Ind}^m\mathrm{Banach}(\mathcal{O}_{Y_{R_k,A,q}})}(-)\ar[d]\ar[d]\ar[d]\ar[d]\\
\mathrm{deRham}_{\mathrm{sComm}_\mathrm{simplicial}\mathrm{Ind}^m\mathrm{Banach}(\mathcal{O}_{X_{\mathbb{Q}_p(p^{1/p^\infty})^{\wedge\flat},A,q}})\ar[r]^{\mathrm{equi}}}(-)\ar[r]\ar[r] &\mathrm{deRham}_{\mathrm{sComm}_\mathrm{simplicial}\varphi_q\mathrm{Ind}^m\mathrm{Banach}(\mathcal{O}_{Y_{\mathbb{Q}_p(p^{1/p^\infty})^{\wedge\flat},A,q}})}(-) 
}
\]

\item Then we have the following equivalence of $K$-group $(\infty,1)$-spectrum from \cite{BGT}:
\[
\xymatrix@R+6pc@C+0pc{
\mathrm{K}^\mathrm{BGT}_{\mathrm{sComm}_\mathrm{simplicial}\mathrm{Ind}\mathrm{Banach}(\mathcal{O}_{X_{R_k,A,q}})\ar[r]^{\mathrm{equi}}}(-)\ar[d]\ar[d]\ar[d]\ar[d]\ar[r]\ar[r] &\mathrm{K}^\mathrm{BGT}_{\mathrm{sComm}_\mathrm{simplicial}\varphi_q\mathrm{Ind}\mathrm{Banach}(\mathcal{O}_{Y_{R_k,A,q}})}(-)\ar[d]\ar[d]\ar[d]\ar[d]\\
\mathrm{K}^\mathrm{BGT}_{\mathrm{sComm}_\mathrm{simplicial}\mathrm{Ind}\mathrm{Banach}(\mathcal{O}_{X_{\mathbb{Q}_p(p^{1/p^\infty})^{\wedge\flat},A,q}})\ar[r]^{\mathrm{equi}}}(-)\ar[r]\ar[r] &\mathrm{K}^\mathrm{BGT}_{\mathrm{sComm}_\mathrm{simplicial}\varphi_q\mathrm{Ind}\mathrm{Banach}(\mathcal{O}_{Y_{\mathbb{Q}_p(p^{1/p^\infty})^{\wedge\flat},A,q}})}(-) 
}
\]
\[
\xymatrix@R+6pc@C+0pc{
\mathrm{K}^\mathrm{BGT}_{\mathrm{sComm}_\mathrm{simplicial}\mathrm{Ind}^m\mathrm{Banach}(\mathcal{O}_{X_{R_k,A,q}})\ar[r]^{\mathrm{equi}}}(-)\ar[d]\ar[d]\ar[d]\ar[d]\ar[r]\ar[r] &\mathrm{K}^\mathrm{BGT}_{\mathrm{sComm}_\mathrm{simplicial}\varphi_q\mathrm{Ind}^m\mathrm{Banach}(\mathcal{O}_{Y_{R_k,A,q}})}(-)\ar[d]\ar[d]\ar[d]\ar[d]\\
\mathrm{K}^\mathrm{BGT}_{\mathrm{sComm}_\mathrm{simplicial}\mathrm{Ind}^m\mathrm{Banach}(\mathcal{O}_{X_{\mathbb{Q}_p(p^{1/p^\infty})^{\wedge\flat},A,q}})\ar[r]^{\mathrm{equi}}}(-)\ar[r]\ar[r] &\mathrm{K}^\mathrm{BGT}_{\mathrm{sComm}_\mathrm{simplicial}\varphi_q\mathrm{Ind}^m\mathrm{Banach}(\mathcal{O}_{Y_{\mathbb{Q}_p(p^{1/p^\infty})^{\wedge\flat},A,q}})}(-) 
}
\]

\end{itemize}

\indent Then we consider further equivariance by considering the arithmetic profinite fundamental group and actually its $q$-th power $\mathrm{Gal}(\overline{{Q}_p\left<T_1^{\pm 1},...,T_k^{\pm 1}\right>}/R_k)^{\times q}$ through the following diagram:\\

\[
\xymatrix@R+6pc@C+0pc{
\mathbb{Z}_p^k=\mathrm{Gal}(R_k/{\mathbb{Q}_p(p^{1/p^\infty})^\wedge\left<T_1^{\pm 1},...,T_k^{\pm 1}\right>})\ar[d]\ar[d]\ar[d]\ar[d] \ar[r]\ar[r] \ar[r]\ar[r] &\mathrm{Gal}(\overline{{Q}_p\left<T_1^{\pm 1},...,T_k^{\pm 1}\right>}/R_k) \ar[d]\ar[d]\ar[d] \ar[r]\ar[r] &\Gamma_{\mathbb{Q}_p} \ar[d]\ar[d]\ar[d]\ar[d]\\
(\mathbb{Z}_p^k=\mathrm{Gal}(R_k/{\mathbb{Q}_p(p^{1/p^\infty})^\wedge\left<T_1^{\pm 1},...,T_k^{\pm 1}\right>}))^{\times q} \ar[r]\ar[r] \ar[r]\ar[r] &\Gamma_k^{\times q}:=\mathrm{Gal}(R_k/{\mathbb{Q}_p\left<T_1^{\pm 1},...,T_k^{\pm 1}\right>})^{\times q} \ar[r] \ar[r]\ar[r] &\Gamma_{\mathbb{Q}_p}^{\times q}.
}
\]

\

We then have the correspond arithmetic profinite fundamental groups equivariant versions:

\begin{itemize}
\item (\text{Proposition}) There is an equivalence between the $\infty$-categories of inductive Banach quasicoherent presheaves:
\[
\xymatrix@R+6pc@C+0pc{
\mathrm{Ind}\mathrm{Banach}_{\Gamma_{k}^{\times q}}(\mathcal{O}_{X_{\mathbb{Q}_p(p^{1/p^\infty})^{\wedge}\left<T_1^{\pm 1/p^\infty},...,T_k^{\pm 1/p^\infty}\right>^\flat,A,q}})\ar[d]\ar[d]\ar[d]\ar[d] \ar[r]^{\mathrm{equi}}\ar[r]\ar[r] &\varphi_q\mathrm{Ind}\mathrm{Banach}_{\Gamma_{k}^{\times q}}(\mathcal{O}_{Y_{\mathbb{Q}_p(p^{1/p^\infty})^{\wedge}\left<T_1^{\pm 1/p^\infty},...,T_k^{\pm 1/p^\infty}\right>^\flat,A,q}}) \ar[d]\ar[d]\ar[d]\ar[d].\\
\mathrm{Ind}\mathrm{Banach}_{\Gamma_{0}^{\times q}}(\mathcal{O}_{X_{\mathbb{Q}_p(p^{1/p^\infty})^{\wedge\flat},A,q}})\ar[r]^{\mathrm{equi}}\ar[r]\ar[r] &\varphi_q\mathrm{Ind}\mathrm{Banach}_{\Gamma_{0}^{\times q}}(\mathcal{O}_{Y_{\mathbb{Q}_p(p^{1/p^\infty})^{\wedge\flat},A,q}}).\\ 
}
\]
\item (\text{Proposition}) There is an equivalence between the $\infty$-categories of monomorphic inductive Banach quasicoherent presheaves:
\[
\xymatrix@R+6pc@C+0pc{
\mathrm{Ind}^m\mathrm{Banach}_{\Gamma_{k}^{\times q}}(\mathcal{O}_{X_{R_k,A,q}})\ar[r]^{\mathrm{equi}}\ar[d]\ar[d]\ar[d]\ar[d]\ar[r]\ar[r] &\varphi_q\mathrm{Ind}^m\mathrm{Banach}_{\Gamma_{k}^{\times q}}(\mathcal{O}_{Y_{R_k,A,q}})\ar[d]\ar[d]\ar[d]\ar[d]\\
\mathrm{Ind}^m\mathrm{Banach}_{\Gamma_{0}^{\times q}}(\mathcal{O}_{X_{\mathbb{Q}_p(p^{1/p^\infty})^{\wedge\flat},A,q}})\ar[r]^{\mathrm{equi}}\ar[r]\ar[r] &\varphi_q\mathrm{Ind}^m\mathrm{Banach}_{\Gamma_{0}^{\times q}}(\mathcal{O}_{Y_{\mathbb{Q}_p(p^{1/p^\infty})^{\wedge\flat},A,q}}).\\  
}
\]	
\end{itemize}

\

\

\indent Now we consider \cite{1CS1} and \cite[Proposition 13.8, Theorem 14.9, Remark 14.10]{1CS2}\footnote{Note that we are motivated as well from \cite{LBV}.}, and study the corresponding solid perfect complexes, solid quasicoherent sheaves and solid vector bundles. Here we are going to use different formalism, therefore we will have different categories and functors. We use the notation $\circledcirc$ to denote any element of $\{\text{solid perfect complexes}, \text{solid quasicoherent sheaves}, \text{solid vector bundles}\}$ from \cite{1CS2} with the corresponding descent results of \cite[Proposition 13.8, Theorem 14.9, Remark 14.10]{1CS2}. Then we have the following:

\begin{itemize}
\item Generalizing Kedlaya-Liu's construction in \cite{KL1}, \cite{KL2} of the adic Fargues-Fontaine space we have a quotient (by using powers of the Frobenius operator) $X_{R,A,q}$ of the space by using \cite{1CS2}:
\begin{align}
Y_{R,A,q}:=\bigcup_{I,\mathrm{multi}}\mathcal{X}^\mathrm{CS}(\text{Robba}^\text{extended}_{R,I,A,q}).	
\end{align}

\item This is a locally ringed space $(X_{R,A,q},\mathcal{O}_{X_{R,A,q}})$, so one can consider the stable $\infty$-category $\mathrm{Module}_{\text{CS},\mathrm{quasicoherent}}(\mathcal{O}_{X_{R,A,q}}) $ which is the $\infty$-category of all the $\mathcal{O}_{X_{R,A,q}}$-sheaves of solid modules over $X_{R,A,q}$. We have the parallel categories for $Y_{R,A,q}$, namely $\varphi\mathrm{Module}_{\text{CS},\mathrm{quasicoherent}}(\mathcal{O}_{X_{R,A,q}})$ and so on. Here we will consider sheaves.

\begin{assumption}\label{assumtionpresheaves}
All the functors of modules or algebras below are Clausen-Scholze sheaves \cite[Proposition 13.8, Theorem 14.9, Remark 14.10]{1CS2}.	
\end{assumption}

\item (\text{Proposition}) There is an equivalence between the $\infty$-categories of inductive solid sheaves:
\[
\xymatrix@R+0pc@C+0pc{
\mathrm{Modules}_\circledcirc(\mathcal{O}_{X_{R,A,q}})\ar[r]^{\mathrm{equi}}\ar[r]\ar[r] &\varphi_q\mathrm{Modules}_\circledcirc(\mathcal{O}_{Y_{R,A,q}}).  
}
\]
\end{itemize}

\begin{itemize}

\item (\text{Proposition}) There is an equivalence between the $\infty$-categories of inductive Banach quasicoherent commutative algebra $E_\infty$ objects:
\[\displayindent=-1.0in
\xymatrix@R+0pc@C+0pc{
\mathrm{sComm}_\mathrm{simplicial}\mathrm{Modules}_{\text{solidquasicoherentsheaves}}(\mathcal{O}_{X_{R,A,q}})\ar[r]^{\mathrm{equi}}\ar[r]\ar[r] &\mathrm{sComm}_\mathrm{simplicial}\varphi_q\mathrm{Modules}_{\text{solidquasicoherentsheaves}}(\mathcal{O}_{Y_{R,A,q}}).  
}
\]

\item Then as in \cite{LBV} we have the equivalence of the de Rham complex after \cite[Definition 5.9, Section 5.2.1]{KKM}\footnote{Here $\circledcirc=\text{solidquasicoherentsheaves}$.}:
\[
\xymatrix@R+0pc@C+0pc{
\mathrm{deRham}_{\mathrm{sComm}_\mathrm{simplicial}\mathrm{Modules}_\circledcirc(\mathcal{O}_{X_{R,A,q}})\ar[r]^{\mathrm{equi}}}(-)\ar[r]\ar[r] &\mathrm{deRham}_{\mathrm{sComm}_\mathrm{simplicial}\varphi_q\mathrm{Modules}_\circledcirc(\mathcal{O}_{Y_{R,A,q}})}(-). 
}
\]

\item Then we have the following equivalence of $K$-group $(\infty,1)$-spectrum from \cite{BGT}\footnote{Here $\circledcirc=\text{solidquasicoherentsheaves}$.}:
\[
\xymatrix@R+0pc@C+0pc{
\mathrm{K}^\mathrm{BGT}_{\mathrm{sComm}_\mathrm{simplicial}\mathrm{Modules}_\circledcirc(\mathcal{O}_{X_{R,A,q}})\ar[r]^{\mathrm{equi}}}(-)\ar[r]\ar[r] &\mathrm{K}^\mathrm{BGT}_{\mathrm{sComm}_\mathrm{simplicial}\varphi_q\mathrm{Modules}_\circledcirc(\mathcal{O}_{Y_{R,A,q}})}(-). 
}
\]
\end{itemize}

\noindent Now let $R=\mathbb{Q}_p(p^{1/p^\infty})^{\wedge\flat}$ and $R_k=\mathbb{Q}_p(p^{1/p^\infty})^{\wedge}\left<T_1^{\pm 1/p^{\infty}},...,T_k^{\pm 1/p^{\infty}}\right>^\flat$ we have the following Galois theoretic results with naturality along $f:\mathrm{Spa}(\mathbb{Q}_p(p^{1/p^\infty})^{\wedge}\left<T_1^{\pm 1/p^\infty},...,T_k^{\pm 1/p^\infty}\right>^\flat)\rightarrow \mathrm{Spa}(\mathbb{Q}_p(p^{1/p^\infty})^{\wedge\flat})$:

\begin{itemize}
\item (\text{Proposition}) There is an equivalence between the $\infty$-categories of inductive Banach quasicoherent presheaves\footnote{Here $\circledcirc=\text{solidquasicoherentsheaves}$.}:
\[
\xymatrix@R+6pc@C+0pc{
\mathrm{Modules}_\circledcirc(\mathcal{O}_{X_{\mathbb{Q}_p(p^{1/p^\infty})^{\wedge}\left<T_1^{\pm 1/p^\infty},...,T_k^{\pm 1/p^\infty}\right>^\flat,A,q}})\ar[d]\ar[d]\ar[d]\ar[d] \ar[r]^{\mathrm{equi}}\ar[r]\ar[r] &\varphi_q\mathrm{Modules}_\circledcirc(\mathcal{O}_{Y_{\mathbb{Q}_p(p^{1/p^\infty})^{\wedge}\left<T_1^{\pm 1/p^\infty},...,T_k^{\pm 1/p^\infty}\right>^\flat,A,q}}) \ar[d]\ar[d]\ar[d]\ar[d].\\
\mathrm{Modules}_\circledcirc(\mathcal{O}_{X_{\mathbb{Q}_p(p^{1/p^\infty})^{\wedge\flat},A,q}})\ar[r]^{\mathrm{equi}}\ar[r]\ar[r] &\varphi_q\mathrm{Modules}_\circledcirc(\mathcal{O}_{Y_{\mathbb{Q}_p(p^{1/p^\infty})^{\wedge\flat},A,q}}).\\ 
}
\]
\item (\text{Proposition}) There is an equivalence between the $\infty$-categories of inductive Banach quasicoherent commutative algebra $E_\infty$ objects\footnote{Here $\circledcirc=\text{solidquasicoherentsheaves}$.}:
\[
\xymatrix@R+6pc@C+0pc{
\mathrm{sComm}_\mathrm{simplicial}\mathrm{Modules}_\circledcirc(\mathcal{O}_{X_{R_k,A,q}})\ar[d]\ar[d]\ar[d]\ar[d]\ar[r]^{\mathrm{equi}}\ar[r]\ar[r] &\mathrm{sComm}_\mathrm{simplicial}\varphi_q\mathrm{Modules}_\circledcirc(\mathcal{O}_{Y_{R_k,A,q}})\ar[d]\ar[d]\ar[d]\ar[d]\\
\mathrm{sComm}_\mathrm{simplicial}\mathrm{Modules}_\circledcirc(\mathcal{O}_{X_{\mathbb{Q}_p(p^{1/p^\infty})^{\wedge\flat},A,q}})\ar[r]^{\mathrm{equi}}\ar[r]\ar[r] &\mathrm{sComm}_\mathrm{simplicial}\varphi_q\mathrm{Modules}_\circledcirc(\mathcal{O}_{Y_{\mathbb{Q}_p(p^{1/p^\infty})^{\wedge\flat},A,q}}).  
}
\]
\item Then as in \cite{LBV} we have the equivalence of the de Rham complex after \cite[Definition 5.9, Section 5.2.1]{KKM}\footnote{Here $\circledcirc=\text{solidquasicoherentsheaves}$.}:
\[\displayindent=-0.4in
\xymatrix@R+6pc@C+0pc{
\mathrm{deRham}_{\mathrm{sComm}_\mathrm{simplicial}\mathrm{Modules}_\circledcirc(\mathcal{O}_{X_{R_k,A,q}})\ar[r]^{\mathrm{equi}}}(-)\ar[d]\ar[d]\ar[d]\ar[d]\ar[r]\ar[r] &\mathrm{deRham}_{\mathrm{sComm}_\mathrm{simplicial}\varphi_q\mathrm{Modules}_\circledcirc(\mathcal{O}_{Y_{R_k,A,q}})}(-)\ar[d]\ar[d]\ar[d]\ar[d]\\
\mathrm{deRham}_{\mathrm{sComm}_\mathrm{simplicial}\mathrm{Modules}_\circledcirc(\mathcal{O}_{X_{\mathbb{Q}_p(p^{1/p^\infty})^{\wedge\flat},A,q}})\ar[r]^{\mathrm{equi}}}(-)\ar[r]\ar[r] &\mathrm{deRham}_{\mathrm{sComm}_\mathrm{simplicial}\varphi_q\mathrm{Modules}_\circledcirc(\mathcal{O}_{Y_{\mathbb{Q}_p(p^{1/p^\infty})^{\wedge\flat},A,q}})}(-). 
}
\]

\item Then we have the following equivalence of $K$-group $(\infty,1)$-spectrum from \cite{BGT}\footnote{Here $\circledcirc=\text{solidquasicoherentsheaves}$.}:
\[
\xymatrix@R+6pc@C+0pc{
\mathrm{K}^\mathrm{BGT}_{\mathrm{sComm}_\mathrm{simplicial}\mathrm{Modules}_\circledcirc(\mathcal{O}_{X_{R_k,A,q}})\ar[r]^{\mathrm{equi}}}(-)\ar[d]\ar[d]\ar[d]\ar[d]\ar[r]\ar[r] &\mathrm{K}^\mathrm{BGT}_{\mathrm{sComm}_\mathrm{simplicial}\varphi_q\mathrm{Modules}_\circledcirc(\mathcal{O}_{Y_{R_k,A,q}})}(-)\ar[d]\ar[d]\ar[d]\ar[d]\\
\mathrm{K}^\mathrm{BGT}_{\mathrm{sComm}_\mathrm{simplicial}\mathrm{Modules}_\circledcirc(\mathcal{O}_{X_{\mathbb{Q}_p(p^{1/p^\infty})^{\wedge\flat},A,q}})\ar[r]^{\mathrm{equi}}}(-)\ar[r]\ar[r] &\mathrm{K}^\mathrm{BGT}_{\mathrm{sComm}_\mathrm{simplicial}\varphi_q\mathrm{Modules}_\circledcirc(\mathcal{O}_{Y_{\mathbb{Q}_p(p^{1/p^\infty})^{\wedge\flat},A,q}})}(-).}
\]

\end{itemize}

\
\indent Then we consider further equivariance by considering the arithmetic profinite fundamental groups $\Gamma_{\mathbb{Q}_p}$ and $\mathrm{Gal}(\overline{\mathbb{Q}_p\left<T_1^{\pm 1},...,T_k^{\pm 1}\right>}/R_k)$ through the following diagram:

\[
\xymatrix@R+0pc@C+0pc{
\mathbb{Z}_p^k=\mathrm{Gal}(R_k/{\mathbb{Q}_p(p^{1/p^\infty})^\wedge\left<T_1^{\pm 1},...,T_k^{\pm 1}\right>}) \ar[r]\ar[r] \ar[r]\ar[r] &\Gamma_k:=\mathrm{Gal}(R_k/{\mathbb{Q}_p\left<T_1^{\pm 1},...,T_k^{\pm 1}\right>}) \ar[r] \ar[r]\ar[r] &\Gamma_{\mathbb{Q}_p}.
}
\]

\begin{itemize}
\item (\text{Proposition}) There is an equivalence between the $\infty$-categories of inductive Banach quasicoherent presheaves\footnote{Here $\circledcirc=\text{solidquasicoherentsheaves}$.}:
\[
\xymatrix@R+6pc@C+0pc{
\mathrm{Modules}_{\circledcirc,\Gamma_k^{\times q}}(\mathcal{O}_{X_{\mathbb{Q}_p(p^{1/p^\infty})^{\wedge}\left<T_1^{\pm 1/p^\infty},...,T_k^{\pm 1/p^\infty}\right>^\flat,A,q}})\ar[d]\ar[d]\ar[d]\ar[d] \ar[r]^{\mathrm{equi}}\ar[r]\ar[r] &\varphi_q\mathrm{Modules}_{\circledcirc,\Gamma_k^{\times q}}(\mathcal{O}_{Y_{\mathbb{Q}_p(p^{1/p^\infty})^{\wedge}\left<T_1^{\pm 1/p^\infty},...,T_k^{\pm 1/p^\infty}\right>^\flat,A,q}}) \ar[d]\ar[d]\ar[d]\ar[d].\\
\mathrm{Modules}_{\circledcirc,\Gamma_0^{\times q}}(\mathcal{O}_{X_{\mathbb{Q}_p(p^{1/p^\infty})^{\wedge\flat},A,q}})\ar[r]^{\mathrm{equi}}\ar[r]\ar[r] &\varphi_q\mathrm{Modules}_{\circledcirc,\Gamma_0^{\times q}}(\mathcal{O}_{Y_{\mathbb{Q}_p(p^{1/p^\infty})^{\wedge\flat},A,q}}).\\ 
}
\]

\item (\text{Proposition}) There is an equivalence between the $\infty$-categories of inductive Banach quasicoherent commutative algebra $E_\infty$ objects\footnote{Here $\circledcirc=\text{solidquasicoherentsheaves}$.}:
\[
\xymatrix@R+6pc@C+0pc{
\mathrm{sComm}_\mathrm{simplicial}\mathrm{Modules}_{\circledcirc,\Gamma_k^{\times q}}(\mathcal{O}_{X_{R_k,A,q}})\ar[d]\ar[d]\ar[d]\ar[d]\ar[r]^{\mathrm{equi}}\ar[r]\ar[r] &\mathrm{sComm}_\mathrm{simplicial}\varphi_q\mathrm{Modules}_{\circledcirc,\Gamma_k^{\times q}}(\mathcal{O}_{Y_{R_k,A,q}})\ar[d]\ar[d]\ar[d]\ar[d]\\
\mathrm{sComm}_\mathrm{simplicial}\mathrm{Modules}_{\circledcirc,\Gamma_0^{\times q}}(\mathcal{O}_{X_{\mathbb{Q}_p(p^{1/p^\infty})^{\wedge\flat},A,q}})\ar[r]^{\mathrm{equi}}\ar[r]\ar[r] &\mathrm{sComm}_\mathrm{simplicial}\varphi_q\mathrm{Modules}_{\circledcirc,\Gamma_0^{\times q}}(\mathcal{O}_{Y_{\mathbb{Q}_p(p^{1/p^\infty})^{\wedge\flat},A,q}}).  
}
\]

\item Then as in \cite{LBV} we have the equivalence of the de Rham complex after \cite[Definition 5.9, Section 5.2.1]{KKM}\footnote{Here $\circledcirc=\text{solidquasicoherentsheaves}$.}:
\[\displayindent=-0.4in
\xymatrix@R+6pc@C+0pc{
\mathrm{deRham}_{\mathrm{sComm}_\mathrm{simplicial}\mathrm{Modules}_{\circledcirc,\Gamma_k^{\times q}}(\mathcal{O}_{X_{R_k,A,q}})\ar[r]^{\mathrm{equi}}}(-)\ar[d]\ar[d]\ar[d]\ar[d]\ar[r]\ar[r] &\mathrm{deRham}_{\mathrm{sComm}_\mathrm{simplicial}\varphi_q\mathrm{Modules}_{\circledcirc,\Gamma_k^{\times q}}(\mathcal{O}_{Y_{R_k,A,q}})}(-)\ar[d]\ar[d]\ar[d]\ar[d]\\
\mathrm{deRham}_{\mathrm{sComm}_\mathrm{simplicial}\mathrm{Modules}_{\circledcirc,\Gamma_0^{\times q}}(\mathcal{O}_{X_{\mathbb{Q}_p(p^{1/p^\infty})^{\wedge\flat},A,q}})\ar[r]^{\mathrm{equi}}}(-)\ar[r]\ar[r] &\mathrm{deRham}_{\mathrm{sComm}_\mathrm{simplicial}\varphi_q\mathrm{Modules}_{\circledcirc,\Gamma_0^{\times q}}(\mathcal{O}_{Y_{\mathbb{Q}_p(p^{1/p^\infty})^{\wedge\flat},A,q}})}(-). 
}
\]

\item Then we have the following equivalence of $K$-group $(\infty,1)$-spectrum from \cite{BGT}\footnote{Here $\circledcirc=\text{solidquasicoherentsheaves}$.}:
\[
\xymatrix@R+6pc@C+0pc{
\mathrm{K}^\mathrm{BGT}_{\mathrm{sComm}_\mathrm{simplicial}\mathrm{Modules}_{\circledcirc,\Gamma_k^{\times q}}(\mathcal{O}_{X_{R_k,A,q}})\ar[r]^{\mathrm{equi}}}(-)\ar[d]\ar[d]\ar[d]\ar[d]\ar[r]\ar[r] &\mathrm{K}^\mathrm{BGT}_{\mathrm{sComm}_\mathrm{simplicial}\varphi_q\mathrm{Modules}_{\circledcirc,\Gamma_k^{\times q}}(\mathcal{O}_{Y_{R_k,A,q}})}(-)\ar[d]\ar[d]\ar[d]\ar[d]\\
\mathrm{K}^\mathrm{BGT}_{\mathrm{sComm}_\mathrm{simplicial}\mathrm{Modules}_{\circledcirc,\Gamma_0^{\times q}}(\mathcal{O}_{X_{\mathbb{Q}_p(p^{1/p^\infty})^{\wedge\flat},A,q}})\ar[r]^{\mathrm{equi}}}(-)\ar[r]\ar[r] &\mathrm{K}^\mathrm{BGT}_{\mathrm{sComm}_\mathrm{simplicial}\varphi_q\mathrm{Modules}_{\circledcirc,\Gamma_0^{\times q}}(\mathcal{O}_{Y_{\mathbb{Q}_p(p^{1/p^\infty})^{\wedge\flat},A,q}})}(-).
}
\]

\end{itemize}

\newpage
\subsection{$\infty$-Categorical Analytic Stacks and Descents II}

As before, we have the following:

\begin{itemize}

\item (\text{Proposition}) There is a functor (global section) between the $\infty$-categories of inductive Banach quasicoherent presheaves:
\[
\xymatrix@R+0pc@C+0pc{
\mathrm{Ind}\mathrm{Banach}(\mathcal{O}_{X_{R,A,q}})\ar[r]^{\mathrm{global}}\ar[r]\ar[r] &\varphi_q\mathrm{Ind}\mathrm{Banach}(\{\mathrm{Robba}^\mathrm{extended}_{{R,A,q,I}}\}_I).  
}
\]
\item (\text{Proposition}) There is a functor (global section) between the $\infty$-categories of monomorphic inductive Banach quasicoherent presheaves:
\[
\xymatrix@R+0pc@C+0pc{
\mathrm{Ind}^m\mathrm{Banach}(\mathcal{O}_{X_{R,A,q}})\ar[r]^{\mathrm{global}}\ar[r]\ar[r] &\varphi_q\mathrm{Ind}^m\mathrm{Banach}(\{\mathrm{Robba}^\mathrm{extended}_{{R,A,q,I}}\}_I).  
}
\]

\item (\text{Proposition}) There is a functor (global section) between the $\infty$-categories of inductive Banach quasicoherent presheaves:
\[
\xymatrix@R+0pc@C+0pc{
\mathrm{Ind}\mathrm{Banach}(\mathcal{O}_{X_{R,A,q}})\ar[r]^{\mathrm{global}}\ar[r]\ar[r] &\varphi_q\mathrm{Ind}\mathrm{Banach}(\{\mathrm{Robba}^\mathrm{extended}_{{R,A,q,I}}\}_I).  
}
\]
\item (\text{Proposition}) There is a functor (global section) between the $\infty$-categories of monomorphic inductive Banach quasicoherent presheaves:
\[
\xymatrix@R+0pc@C+0pc{
\mathrm{Ind}^m\mathrm{Banach}(\mathcal{O}_{X_{R,A,q}})\ar[r]^{\mathrm{global}}\ar[r]\ar[r] &\varphi_q\mathrm{Ind}^m\mathrm{Banach}(\{\mathrm{Robba}^\mathrm{extended}_{{R,A,q,I}}\}_I).  
}
\]
\item (\text{Proposition}) There is a functor (global section) between the $\infty$-categories of inductive Banach quasicoherent commutative algebra $E_\infty$ objects:
\[
\xymatrix@R+0pc@C+0pc{
\mathrm{sComm}_\mathrm{simplicial}\mathrm{Ind}\mathrm{Banach}(\mathcal{O}_{X_{R,A,q}})\ar[r]^{\mathrm{global}}\ar[r]\ar[r] &\mathrm{sComm}_\mathrm{simplicial}\varphi_q\mathrm{Ind}\mathrm{Banach}(\{\mathrm{Robba}^\mathrm{extended}_{{R,A,q,I}}\}_I).  
}
\]
\item (\text{Proposition}) There is a functor (global section) between the $\infty$-categories of monomorphic inductive Banach quasicoherent commutative algebra $E_\infty$ objects:
\[
\xymatrix@R+0pc@C+0pc{
\mathrm{sComm}_\mathrm{simplicial}\mathrm{Ind}^m\mathrm{Banach}(\mathcal{O}_{X_{R,A,q}})\ar[r]^{\mathrm{global}}\ar[r]\ar[r] &\mathrm{sComm}_\mathrm{simplicial}\varphi_q\mathrm{Ind}^m\mathrm{Banach}(\{\mathrm{Robba}^\mathrm{extended}_{{R,A,q,I}}\}_I).  
}
\]

\item Then parallel as in \cite{LBV} we have a functor (global section) of the de Rham complex after \cite[Definition 5.9, Section 5.2.1]{KKM}:
\[
\xymatrix@R+0pc@C+0pc{
\mathrm{deRham}_{\mathrm{sComm}_\mathrm{simplicial}\mathrm{Ind}\mathrm{Banach}(\mathcal{O}_{X_{R,A,q}})\ar[r]^{\mathrm{global}}}(-)\ar[r]\ar[r] &\mathrm{deRham}_{\mathrm{sComm}_\mathrm{simplicial}\varphi_q\mathrm{Ind}\mathrm{Banach}(\{\mathrm{Robba}^\mathrm{extended}_{{R,A,q,I}}\}_I)}(-), 
}
\]
\[
\xymatrix@R+0pc@C+0pc{
\mathrm{deRham}_{\mathrm{sComm}_\mathrm{simplicial}\mathrm{Ind}^m\mathrm{Banach}(\mathcal{O}_{X_{R,A,q}})\ar[r]^{\mathrm{global}}}(-)\ar[r]\ar[r] &\mathrm{deRham}_{\mathrm{sComm}_\mathrm{simplicial}\varphi_q\mathrm{Ind}^m\mathrm{Banach}(\{\mathrm{Robba}^\mathrm{extended}_{{R,A,q,I}}\}_I)}(-). 
}
\]

\item Then we have the following a functor (global section) of $K$-group $(\infty,1)$-spectrum from \cite{BGT}:
\[
\xymatrix@R+0pc@C+0pc{
\mathrm{K}^\mathrm{BGT}_{\mathrm{sComm}_\mathrm{simplicial}\mathrm{Ind}\mathrm{Banach}(\mathcal{O}_{X_{R,A,q}})\ar[r]^{\mathrm{global}}}(-)\ar[r]\ar[r] &\mathrm{K}^\mathrm{BGT}_{\mathrm{sComm}_\mathrm{simplicial}\varphi_q\mathrm{Ind}\mathrm{Banach}(\{\mathrm{Robba}^\mathrm{extended}_{{R,A,q,I}}\}_I)}(-), 
}
\]
\[
\xymatrix@R+0pc@C+0pc{
\mathrm{K}^\mathrm{BGT}_{\mathrm{sComm}_\mathrm{simplicial}\mathrm{Ind}^m\mathrm{Banach}(\mathcal{O}_{X_{R,A,q}})\ar[r]^{\mathrm{global}}}(-)\ar[r]\ar[r] &\mathrm{K}^\mathrm{BGT}_{\mathrm{sComm}_\mathrm{simplicial}\varphi_q\mathrm{Ind}^m\mathrm{Banach}(\{\mathrm{Robba}^\mathrm{extended}_{{R,A,q,I}}\}_I)}(-). 
}
\]
\end{itemize}

\noindent Now let $R=\mathbb{Q}_p(p^{1/p^\infty})^{\wedge\flat}$ and $R_k=\mathbb{Q}_p(p^{1/p^\infty})^{\wedge}\left<T_1^{\pm 1/p^{\infty}},...,T_k^{\pm 1/p^{\infty}}\right>^\flat$ we have the following Galois theoretic results with naturality along $f:\mathrm{Spa}(\mathbb{Q}_p(p^{1/p^\infty})^{\wedge}\left<T_1^{\pm 1/p^\infty},...,T_k^{\pm 1/p^\infty}\right>^\flat)\rightarrow \mathrm{Spa}(\mathbb{Q}_p(p^{1/p^\infty})^{\wedge\flat})$:

\begin{itemize}
\item (\text{Proposition}) There is a functor (global section) between the $\infty$-categories of inductive Banach quasicoherent presheaves:
\[
\xymatrix@R+6pc@C+0pc{
\mathrm{Ind}\mathrm{Banach}(\mathcal{O}_{X_{\mathbb{Q}_p(p^{1/p^\infty})^{\wedge}\left<T_1^{\pm 1/p^\infty},...,T_k^{\pm 1/p^\infty}\right>^\flat,A,q}})\ar[d]\ar[d]\ar[d]\ar[d] \ar[r]^{\mathrm{global}}\ar[r]\ar[r] &\varphi_q\mathrm{Ind}\mathrm{Banach}(\{\mathrm{Robba}^\mathrm{extended}_{{R_k,A,q,I}}\}_I) \ar[d]\ar[d]\ar[d]\ar[d].\\
\mathrm{Ind}\mathrm{Banach}(\mathcal{O}_{X_{\mathbb{Q}_p(p^{1/p^\infty})^{\wedge\flat},A,q}})\ar[r]^{\mathrm{global}}\ar[r]\ar[r] &\varphi_q\mathrm{Ind}\mathrm{Banach}(\{\mathrm{Robba}^\mathrm{extended}_{{R_0,A,q,I}}\}_I).\\ 
}
\]
\item (\text{Proposition}) There is a functor (global section) between the $\infty$-categories of monomorphic inductive Banach quasicoherent presheaves:
\[
\xymatrix@R+6pc@C+0pc{
\mathrm{Ind}^m\mathrm{Banach}(\mathcal{O}_{X_{R_k,A,q}})\ar[r]^{\mathrm{global}}\ar[d]\ar[d]\ar[d]\ar[d]\ar[r]\ar[r] &\varphi_q\mathrm{Ind}^m\mathrm{Banach}(\{\mathrm{Robba}^\mathrm{extended}_{{R_k,A,q,I}}\}_I)\ar[d]\ar[d]\ar[d]\ar[d]\\
\mathrm{Ind}^m\mathrm{Banach}(\mathcal{O}_{X_{\mathbb{Q}_p(p^{1/p^\infty})^{\wedge\flat},A,q}})\ar[r]^{\mathrm{global}}\ar[r]\ar[r] &\varphi_q\mathrm{Ind}^m\mathrm{Banach}(\{\mathrm{Robba}^\mathrm{extended}_{{R_0,A,q,I}}\}_I).\\  
}
\]
\item (\text{Proposition}) There is a functor (global section) between the $\infty$-categories of inductive Banach quasicoherent commutative algebra $E_\infty$ objects:
\[
\xymatrix@R+6pc@C+0pc{
\mathrm{sComm}_\mathrm{simplicial}\mathrm{Ind}\mathrm{Banach}(\mathcal{O}_{X_{R_k,A,q}})\ar[d]\ar[d]\ar[d]\ar[d]\ar[r]^{\mathrm{global}}\ar[r]\ar[r] &\mathrm{sComm}_\mathrm{simplicial}\varphi_q\mathrm{Ind}\mathrm{Banach}(\{\mathrm{Robba}^\mathrm{extended}_{{R_k,A,q,I}}\}_I)\ar[d]\ar[d]\ar[d]\ar[d]\\
\mathrm{sComm}_\mathrm{simplicial}\mathrm{Ind}\mathrm{Banach}(\mathcal{O}_{X_{\mathbb{Q}_p(p^{1/p^\infty})^{\wedge\flat},A,q}})\ar[r]^{\mathrm{global}}\ar[r]\ar[r] &\mathrm{sComm}_\mathrm{simplicial}\varphi_q\mathrm{Ind}\mathrm{Banach}(\{\mathrm{Robba}^\mathrm{extended}_{{R_0,A,q,I}}\}_I).  
}
\]
\item (\text{Proposition}) There is a functor (global section) between the $\infty$-categories of monomorphic inductive Banach quasicoherent commutative algebra $E_\infty$ objects:
\[\displayindent=-0.4in
\xymatrix@R+6pc@C+0pc{
\mathrm{sComm}_\mathrm{simplicial}\mathrm{Ind}^m\mathrm{Banach}(\mathcal{O}_{X_{R_k,A,q}})\ar[d]\ar[d]\ar[d]\ar[d]\ar[r]^{\mathrm{global}}\ar[r]\ar[r] &\mathrm{sComm}_\mathrm{simplicial}\varphi_q\mathrm{Ind}^m\mathrm{Banach}(\{\mathrm{Robba}^\mathrm{extended}_{{R_k,A,q,I}}\}_I)\ar[d]\ar[d]\ar[d]\ar[d]\\
 \mathrm{sComm}_\mathrm{simplicial}\mathrm{Ind}^m\mathrm{Banach}(\mathcal{O}_{X_{\mathbb{Q}_p(p^{1/p^\infty})^{\wedge\flat},A,q}})\ar[r]^{\mathrm{global}}\ar[r]\ar[r] &\mathrm{sComm}_\mathrm{simplicial}\varphi_q\mathrm{Ind}^m\mathrm{Banach}(\{\mathrm{Robba}^\mathrm{extended}_{{R_0,A,q,I}}\}_I). 
}
\]

\item Then parallel as in \cite{LBV} we have a functor (global section) of the de Rham complex after \cite[Definition 5.9, Section 5.2.1]{KKM}:
\[\displayindent=-0.4in
\xymatrix@R+6pc@C+0pc{
\mathrm{deRham}_{\mathrm{sComm}_\mathrm{simplicial}\mathrm{Ind}\mathrm{Banach}(\mathcal{O}_{X_{R_k,A,q}})\ar[r]^{\mathrm{global}}}(-)\ar[d]\ar[d]\ar[d]\ar[d]\ar[r]\ar[r] &\mathrm{deRham}_{\mathrm{sComm}_\mathrm{simplicial}\varphi_q\mathrm{Ind}\mathrm{Banach}(\{\mathrm{Robba}^\mathrm{extended}_{{R_k,A,q,I}}\}_I)}(-)\ar[d]\ar[d]\ar[d]\ar[d]\\
\mathrm{deRham}_{\mathrm{sComm}_\mathrm{simplicial}\mathrm{Ind}\mathrm{Banach}(\mathcal{O}_{X_{\mathbb{Q}_p(p^{1/p^\infty})^{\wedge\flat},A,q}})\ar[r]^{\mathrm{global}}}(-)\ar[r]\ar[r] &\mathrm{deRham}_{\mathrm{sComm}_\mathrm{simplicial}\varphi_q\mathrm{Ind}\mathrm{Banach}(\{\mathrm{Robba}^\mathrm{extended}_{{R_0,A,q,I}}\}_I)}(-), 
}
\]
\[\displayindent=-0.4in
\xymatrix@R+6pc@C+0pc{
\mathrm{deRham}_{\mathrm{sComm}_\mathrm{simplicial}\mathrm{Ind}^m\mathrm{Banach}(\mathcal{O}_{X_{R_k,A,q}})\ar[r]^{\mathrm{global}}}(-)\ar[d]\ar[d]\ar[d]\ar[d]\ar[r]\ar[r] &\mathrm{deRham}_{\mathrm{sComm}_\mathrm{simplicial}\varphi_q\mathrm{Ind}^m\mathrm{Banach}(\{\mathrm{Robba}^\mathrm{extended}_{{R_k,A,q,I}}\}_I)}(-)\ar[d]\ar[d]\ar[d]\ar[d]\\
\mathrm{deRham}_{\mathrm{sComm}_\mathrm{simplicial}\mathrm{Ind}^m\mathrm{Banach}(\mathcal{O}_{X_{\mathbb{Q}_p(p^{1/p^\infty})^{\wedge\flat},A,q}})\ar[r]^{\mathrm{global}}}(-)\ar[r]\ar[r] &\mathrm{deRham}_{\mathrm{sComm}_\mathrm{simplicial}\varphi_q\mathrm{Ind}^m\mathrm{Banach}(\{\mathrm{Robba}^\mathrm{extended}_{{R_0,A,q,I}}\}_I)}(-). 
}
\]

\item Then we have the following a functor (global section) of $K$-group $(\infty,1)$-spectrum from \cite{BGT}:
\[
\xymatrix@R+6pc@C+0pc{
\mathrm{K}^\mathrm{BGT}_{\mathrm{sComm}_\mathrm{simplicial}\mathrm{Ind}\mathrm{Banach}(\mathcal{O}_{X_{R_k,A,q}})\ar[r]^{\mathrm{global}}}(-)\ar[d]\ar[d]\ar[d]\ar[d]\ar[r]\ar[r] &\mathrm{K}^\mathrm{BGT}_{\mathrm{sComm}_\mathrm{simplicial}\varphi_q\mathrm{Ind}\mathrm{Banach}(\{\mathrm{Robba}^\mathrm{extended}_{{R_k,A,q,I}}\}_I)}(-)\ar[d]\ar[d]\ar[d]\ar[d]\\
\mathrm{K}^\mathrm{BGT}_{\mathrm{sComm}_\mathrm{simplicial}\mathrm{Ind}\mathrm{Banach}(\mathcal{O}_{X_{\mathbb{Q}_p(p^{1/p^\infty})^{\wedge\flat},A,q}})\ar[r]^{\mathrm{global}}}(-)\ar[r]\ar[r] &\mathrm{K}^\mathrm{BGT}_{\mathrm{sComm}_\mathrm{simplicial}\varphi_q\mathrm{Ind}\mathrm{Banach}(\{\mathrm{Robba}^\mathrm{extended}_{{R_0,A,q,I}}\}_I)}(-), 
}
\]
\[
\xymatrix@R+6pc@C+0pc{
\mathrm{K}^\mathrm{BGT}_{\mathrm{sComm}_\mathrm{simplicial}\mathrm{Ind}^m\mathrm{Banach}(\mathcal{O}_{X_{R_k,A,q}})\ar[r]^{\mathrm{global}}}(-)\ar[d]\ar[d]\ar[d]\ar[d]\ar[r]\ar[r] &\mathrm{K}^\mathrm{BGT}_{\mathrm{sComm}_\mathrm{simplicial}\varphi_q\mathrm{Ind}^m\mathrm{Banach}(\{\mathrm{Robba}^\mathrm{extended}_{{R_k,A,q,I}}\}_I)}(-)\ar[d]\ar[d]\ar[d]\ar[d]\\
\mathrm{K}^\mathrm{BGT}_{\mathrm{sComm}_\mathrm{simplicial}\mathrm{Ind}^m\mathrm{Banach}(\mathcal{O}_{X_{\mathbb{Q}_p(p^{1/p^\infty})^{\wedge\flat},A,q}})\ar[r]^{\mathrm{global}}}(-)\ar[r]\ar[r] &\mathrm{K}^\mathrm{BGT}_{\mathrm{sComm}_\mathrm{simplicial}\varphi_q\mathrm{Ind}^m\mathrm{Banach}(\{\mathrm{Robba}^\mathrm{extended}_{{R_0,A,q,I}}\}_I)}(-). 
}
\]

\end{itemize}

\
\indent Then we consider further equivariance by considering the arithmetic profinite fundamental group and actually its $q$-th power $\mathrm{Gal}(\overline{{Q}_p\left<T_1^{\pm 1},...,T_k^{\pm 1}\right>}/R_k)^{\times q}$ through the following diagram\:

\[
\xymatrix@R+6pc@C+0pc{
\mathbb{Z}_p^k=\mathrm{Gal}(R_k/{\mathbb{Q}_p(p^{1/p^\infty})^\wedge\left<T_1^{\pm 1},...,T_k^{\pm 1}\right>})\ar[d]\ar[d]\ar[d]\ar[d] \ar[r]\ar[r] \ar[r]\ar[r] &\mathrm{Gal}(\overline{{Q}_p\left<T_1^{\pm 1},...,T_k^{\pm 1}\right>}/R_k) \ar[d]\ar[d]\ar[d] \ar[r]\ar[r] &\Gamma_{\mathbb{Q}_p} \ar[d]\ar[d]\ar[d]\ar[d]\\
(\mathbb{Z}_p^k=\mathrm{Gal}(R_k/{\mathbb{Q}_p(p^{1/p^\infty})^\wedge\left<T_1^{\pm 1},...,T_k^{\pm 1}\right>}))^{\times q} \ar[r]\ar[r] \ar[r]\ar[r] &\Gamma_k^{\times q}:=\mathrm{Gal}(R_k/{\mathbb{Q}_p\left<T_1^{\pm 1},...,T_k^{\pm 1}\right>})^{\times q} \ar[r] \ar[r]\ar[r] &\Gamma_{\mathbb{Q}_p}^{\times q}.
}
\]

\

We then have the correspond arithmetic profinite fundamental groups equivariant versions:
\begin{itemize}
\item (\text{Proposition}) There is a functor (global section) between the $\infty$-categories of inductive Banach quasicoherent presheaves:
\[
\xymatrix@R+6pc@C+0pc{
\mathrm{Ind}\mathrm{Banach}_{\Gamma_{k}^{\times q}}(\mathcal{O}_{X_{\mathbb{Q}_p(p^{1/p^\infty})^{\wedge}\left<T_1^{\pm 1/p^\infty},...,T_k^{\pm 1/p^\infty}\right>^\flat,A,q}})\ar[d]\ar[d]\ar[d]\ar[d] \ar[r]^{\mathrm{global}}\ar[r]\ar[r] &\varphi_q\mathrm{Ind}\mathrm{Banach}_{\Gamma_{k}^{\times q}}(\{\mathrm{Robba}^\mathrm{extended}_{{R_k,A,q,I}}\}_I) \ar[d]\ar[d]\ar[d]\ar[d].\\
\mathrm{Ind}\mathrm{Banach}_{\Gamma_{0}^{\times q}}(\mathcal{O}_{X_{\mathbb{Q}_p(p^{1/p^\infty})^{\wedge\flat},A,q}})\ar[r]^{\mathrm{global}}\ar[r]\ar[r] &\varphi_q\mathrm{Ind}\mathrm{Banach}_{\Gamma_{0}^{\times q}}(\{\mathrm{Robba}^\mathrm{extended}_{{R_0,A,q,I}}\}_I).\\ 
}
\]
\item (\text{Proposition}) There is a functor (global section) between the $\infty$-categories of monomorphic inductive Banach quasicoherent presheaves:
\[
\xymatrix@R+6pc@C+0pc{
\mathrm{Ind}^m\mathrm{Banach}_{\Gamma_{k}^{\times q}}(\mathcal{O}_{X_{R_k,A,q}})\ar[r]^{\mathrm{global}}\ar[d]\ar[d]\ar[d]\ar[d]\ar[r]\ar[r] &\varphi_q\mathrm{Ind}^m\mathrm{Banach}_{\Gamma_{k}^{\times q}}(\{\mathrm{Robba}^\mathrm{extended}_{{R_k,A,q,I}}\}_I)\ar[d]\ar[d]\ar[d]\ar[d]\\
\mathrm{Ind}^m\mathrm{Banach}_{\Gamma_{0}^{\times q}}(\mathcal{O}_{X_{\mathbb{Q}_p(p^{1/p^\infty})^{\wedge\flat},A,q}})\ar[r]^{\mathrm{global}}\ar[r]\ar[r] &\varphi_q\mathrm{Ind}^m\mathrm{Banach}_{\Gamma_{0}^{\times q}}(\{\mathrm{Robba}^\mathrm{extended}_{{R_0,A,q,I}}\}_I).\\  
}
\]
\item (\text{Proposition}) There is a functor (global section) between the $\infty$-categories of inductive Banach quasicoherent commutative algebra $E_\infty$ objects:
\[\displayindent=-0.4in
\xymatrix@R+6pc@C+0pc{
\mathrm{sComm}_\mathrm{simplicial}\mathrm{Ind}\mathrm{Banach}_{\Gamma_{k}^{\times q}}(\mathcal{O}_{X_{R_k,A,q}})\ar[d]\ar[d]\ar[d]\ar[d]\ar[r]^{\mathrm{global}}\ar[r]\ar[r] &\mathrm{sComm}_\mathrm{simplicial}\varphi_q\mathrm{Ind}\mathrm{Banach}_{\Gamma_{k}^{\times q}}(\{\mathrm{Robba}^\mathrm{extended}_{{R_k,A,q,I}}\}_I)\ar[d]\ar[d]\ar[d]\ar[d]\\
\mathrm{sComm}_\mathrm{simplicial}\mathrm{Ind}\mathrm{Banach}_{\Gamma_{0}^{\times q}}(\mathcal{O}_{X_{\mathbb{Q}_p(p^{1/p^\infty})^{\wedge\flat},A,q}})\ar[r]^{\mathrm{global}}\ar[r]\ar[r] &\mathrm{sComm}_\mathrm{simplicial}\varphi_q\mathrm{Ind}\mathrm{Banach}_{\Gamma_{0}^{\times q}}(\{\mathrm{Robba}^\mathrm{extended}_{{R_0,A,q,I}}\}_I).  
}
\]
\item (\text{Proposition}) There is a functor (global section) between the $\infty$-categories of monomorphic inductive Banach quasicoherent commutative algebra $E_\infty$ objects:
\[\displayindent=-0.4in
\xymatrix@R+6pc@C+0pc{
\mathrm{sComm}_\mathrm{simplicial}\mathrm{Ind}^m\mathrm{Banach}_{\Gamma_{k}^{\times q}}(\mathcal{O}_{X_{R_k,A,q}})\ar[d]\ar[d]\ar[d]\ar[d]\ar[r]^{\mathrm{global}}\ar[r]\ar[r] &\mathrm{sComm}_\mathrm{simplicial}\varphi_q\mathrm{Ind}^m\mathrm{Banach}_{\Gamma_{k}^{\times q}}(\{\mathrm{Robba}^\mathrm{extended}_{{R_k,A,q,I}}\}_I)\ar[d]\ar[d]\ar[d]\ar[d]\\
 \mathrm{sComm}_\mathrm{simplicial}\mathrm{Ind}^m\mathrm{Banach}_{\Gamma_{0}^{\times q}}(\mathcal{O}_{X_{\mathbb{Q}_p(p^{1/p^\infty})^{\wedge\flat},A,q}})\ar[r]^{\mathrm{global}}\ar[r]\ar[r] &\mathrm{sComm}_\mathrm{simplicial}\varphi_q\mathrm{Ind}^m\mathrm{Banach}_{\Gamma_{0}^{\times q}}(\{\mathrm{Robba}^\mathrm{extended}_{{R_0,A,q,I}}\}_I). 
}
\]

\item Then parallel as in \cite{LBV} we have a functor (global section) of the de Rham complex after \cite[Definition 5.9, Section 5.2.1]{KKM}:
\[\displayindent=-0.5in
\xymatrix@R+6pc@C+0pc{
\mathrm{deRham}_{\mathrm{sComm}_\mathrm{simplicial}\mathrm{Ind}\mathrm{Banach}_{\Gamma_{k}^{\times q}}(\mathcal{O}_{X_{R_k,A,q}})\ar[r]^{\mathrm{global}}}(-)\ar[d]\ar[d]\ar[d]\ar[d]\ar[r]\ar[r] &\mathrm{deRham}_{\mathrm{sComm}_\mathrm{simplicial}\varphi_q\mathrm{Ind}\mathrm{Banach}_{\Gamma_{k}^{\times q}}(\{\mathrm{Robba}^\mathrm{extended}_{{R_k,A,q,I}}\}_I)}(-)\ar[d]\ar[d]\ar[d]\ar[d]\\
\mathrm{deRham}_{\mathrm{sComm}_\mathrm{simplicial}\mathrm{Ind}\mathrm{Banach}_{\Gamma_{0}^{\times q}}(\mathcal{O}_{X_{\mathbb{Q}_p(p^{1/p^\infty})^{\wedge\flat},A,q}})\ar[r]^{\mathrm{global}}}(-)\ar[r]\ar[r] &\mathrm{deRham}_{\mathrm{sComm}_\mathrm{simplicial}\varphi_q\mathrm{Ind}\mathrm{Banach}_{\Gamma_{0}^{\times q}}(\{\mathrm{Robba}^\mathrm{extended}_{{R_0,A,q,I}}\}_I)}(-), 
}
\]
\[\displayindent=-0.7in
\xymatrix@R+6pc@C+0pc{
\mathrm{deRham}_{\mathrm{sComm}_\mathrm{simplicial}\mathrm{Ind}^m\mathrm{Banach}_{\Gamma_{k}^{\times q}}(\mathcal{O}_{X_{R_k,A,q}})\ar[r]^{\mathrm{global}}}(-)\ar[d]\ar[d]\ar[d]\ar[d]\ar[r]\ar[r] &\mathrm{deRham}_{\mathrm{sComm}_\mathrm{simplicial}\varphi_q\mathrm{Ind}^m\mathrm{Banach}_{\Gamma_{k}^{\times q}}(\{\mathrm{Robba}^\mathrm{extended}_{{R_k,A,q,I}}\}_I)}(-)\ar[d]\ar[d]\ar[d]\ar[d]\\
\mathrm{deRham}_{\mathrm{sComm}_\mathrm{simplicial}\mathrm{Ind}^m\mathrm{Banach}_{\Gamma_{0}^{\times q}}(\mathcal{O}_{X_{\mathbb{Q}_p(p^{1/p^\infty})^{\wedge\flat},A,q}})\ar[r]^{\mathrm{global}}}(-)\ar[r]\ar[r] &\mathrm{deRham}_{\mathrm{sComm}_\mathrm{simplicial}\varphi_q\mathrm{Ind}^m\mathrm{Banach}_{\Gamma_{0}^{\times q}}(\{\mathrm{Robba}^\mathrm{extended}_{{R_0,A,q,I}}\}_I)}(-). 
}
\]

\item Then we have the following a functor (global section) of $K$-group $(\infty,1)$-spectrum from \cite{BGT}:
\[
\xymatrix@R+6pc@C+0pc{
\mathrm{K}^\mathrm{BGT}_{\mathrm{sComm}_\mathrm{simplicial}\mathrm{Ind}\mathrm{Banach}_{\Gamma_{k}^{\times q}}(\mathcal{O}_{X_{R_k,A,q}})\ar[r]^{\mathrm{global}}}(-)\ar[d]\ar[d]\ar[d]\ar[d]\ar[r]\ar[r] &\mathrm{K}^\mathrm{BGT}_{\mathrm{sComm}_\mathrm{simplicial}\varphi_q\mathrm{Ind}\mathrm{Banach}_{\Gamma_{k}^{\times q}}(\{\mathrm{Robba}^\mathrm{extended}_{{R_k,A,q,I}}\}_I)}(-)\ar[d]\ar[d]\ar[d]\ar[d]\\
\mathrm{K}^\mathrm{BGT}_{\mathrm{sComm}_\mathrm{simplicial}\mathrm{Ind}\mathrm{Banach}_{\Gamma_{0}^{\times q}}(\mathcal{O}_{X_{\mathbb{Q}_p(p^{1/p^\infty})^{\wedge\flat},A,q}})\ar[r]^{\mathrm{global}}}(-)\ar[r]\ar[r] &\mathrm{K}^\mathrm{BGT}_{\mathrm{sComm}_\mathrm{simplicial}\varphi_q\mathrm{Ind}\mathrm{Banach}_{\Gamma_{0}^{\times q}}(\{\mathrm{Robba}^\mathrm{extended}_{{R_0,A,q,I}}\}_I)}(-), 
}
\]
\[
\xymatrix@R+6pc@C+0pc{
\mathrm{K}^\mathrm{BGT}_{\mathrm{sComm}_\mathrm{simplicial}\mathrm{Ind}^m\mathrm{Banach}_{\Gamma_{k}^{\times q}}(\mathcal{O}_{X_{R_k,A,q}})\ar[r]^{\mathrm{global}}}(-)\ar[d]\ar[d]\ar[d]\ar[d]\ar[r]\ar[r] &\mathrm{K}^\mathrm{BGT}_{\mathrm{sComm}_\mathrm{simplicial}\varphi_q\mathrm{Ind}^m\mathrm{Banach}_{\Gamma_{k}^{\times q}}(\{\mathrm{Robba}^\mathrm{extended}_{{R_k,A,q,I}}\}_I)}(-)\ar[d]\ar[d]\ar[d]\ar[d]\\
\mathrm{K}^\mathrm{BGT}_{\mathrm{sComm}_\mathrm{simplicial}\mathrm{Ind}^m\mathrm{Banach}_{\Gamma_{0}^{\times q}}(\mathcal{O}_{X_{\mathbb{Q}_p(p^{1/p^\infty})^{\wedge\flat},A,q}})\ar[r]^{\mathrm{global}}}(-)\ar[r]\ar[r] &\mathrm{K}^\mathrm{BGT}_{\mathrm{sComm}_\mathrm{simplicial}\varphi_q\mathrm{Ind}^m\mathrm{Banach}_{\Gamma_{0}^{\times q}}(\{\mathrm{Robba}^\mathrm{extended}_{{R_0,A,q,I}}\}_I)}(-). 
}
\]

\end{itemize}

\

\

\begin{remark}
\noindent We can certainly consider the quasicoherent sheaves in \cite[Lemma 7.11, Remark 7.12]{1BBK}, therefore all the quasicoherent presheaves and modules will be those satisfying \cite[Lemma 7.11, Remark 7.12]{1BBK} if one would like to consider the the quasicoherent sheaves. That being all as this said, we would believe that the big quasicoherent presheaves are automatically quasicoherent sheaves (namely satisfying the corresponding \v{C}ech $\infty$-descent as in \cite[Section 9.3]{KKM} and \cite[Lemma 7.11, Remark 7.12]{1BBK}) and the corresponding global section functors are automatically equivalence of $\infty$-categories. 
\end{remark}

\

\indent In Clausen-Scholze formalism we have the following:

\begin{itemize}

\item (\text{Proposition}) There is a functor (global section) between the $\infty$-categories of inductive Banach quasicoherent sheaves:
\[
\xymatrix@R+0pc@C+0pc{
\mathrm{Modules}_\circledcirc(\mathcal{O}_{X_{R,A,q}})\ar[r]^{\mathrm{global}}\ar[r]\ar[r] &\varphi_q\mathrm{Modules}_\circledcirc(\{\mathrm{Robba}^{\mathrm{extended},q}_{{R,A,I}}\}_I).  
}
\]

\item (\text{Proposition}) There is a functor (global section) between the $\infty$-categories of inductive Banach quasicoherent sheaves:
\[
\xymatrix@R+0pc@C+0pc{
\mathrm{Modules}_\circledcirc(\mathcal{O}_{X_{R,A,q}})\ar[r]^{\mathrm{global}}\ar[r]\ar[r] &\varphi_q\mathrm{Modules}_\circledcirc(\{\mathrm{Robba}^{\mathrm{extended},q}_{{R,A,I}}\}_I).  
}
\]

\item (\text{Proposition}) There is a functor (global section) between the $\infty$-categories of inductive Banach quasicoherent commutative algebra $E_\infty$ objects\footnote{Here $\circledcirc=\text{solidquasicoherentsheaves}$.}:
\[
\xymatrix@R+0pc@C+0pc{
\mathrm{sComm}_\mathrm{simplicial}\mathrm{Modules}_\circledcirc(\mathcal{O}_{X_{R,A,q}})\ar[r]^{\mathrm{global}}\ar[r]\ar[r] &\mathrm{sComm}_\mathrm{simplicial}\varphi_q\mathrm{Modules}_\circledcirc(\{\mathrm{Robba}^{\mathrm{extended},q}_{{R,A,I}}\}_I).  
}
\]

\item Then as in \cite{LBV} we have a functor (global section) of the de Rham complex after \cite[Definition 5.9, Section 5.2.1]{KKM}\footnote{Here $\circledcirc=\text{solidquasicoherentsheaves}$.}:
\[
\xymatrix@R+0pc@C+0pc{
\mathrm{deRham}_{\mathrm{sComm}_\mathrm{simplicial}\mathrm{Modules}_\circledcirc(\mathcal{O}_{X_{R,A,q}})\ar[r]^{\mathrm{global}}}(-)\ar[r]\ar[r] &\mathrm{deRham}_{\mathrm{sComm}_\mathrm{simplicial}\varphi_q\mathrm{Modules}_\circledcirc(\{\mathrm{Robba}^{\mathrm{extended},q}_{{R,A,I}}\}_I)}(-). 
}
\]

\item Then we have the following a functor (global section) of $K$-group $(\infty,1)$-spectrum from \cite{BGT}\footnote{Here $\circledcirc=\text{solidquasicoherentsheaves}$.}:
\[
\xymatrix@R+0pc@C+0pc{
\mathrm{K}^\mathrm{BGT}_{\mathrm{sComm}_\mathrm{simplicial}\mathrm{Module}_\circledcirc(\mathcal{O}_{X_{R,A,q}})\ar[r]^{\mathrm{global}}}(-)\ar[r]\ar[r] &\mathrm{K}^\mathrm{BGT}_{\mathrm{sComm}_\mathrm{simplicial}\varphi_q\mathrm{Module}_\circledcirc(\{\mathrm{Robba}^{\mathrm{extended},q}_{{R,A,I}}\}_I)}(-). 
}
\]
\end{itemize}

\noindent Now let $R=\mathbb{Q}_p(p^{1/p^\infty})^{\wedge\flat}$ and $R_k=\mathbb{Q}_p(p^{1/p^\infty})^{\wedge}\left<T_1^{\pm 1/p^{\infty}},...,T_k^{\pm 1/p^{\infty}}\right>^\flat$ we have the following Galois theoretic results with naturality along $f:\mathrm{Spa}(\mathbb{Q}_p(p^{1/p^\infty})^{\wedge}\left<T_1^{\pm 1/p^\infty},...,T_k^{\pm 1/p^\infty}\right>^\flat)\rightarrow \mathrm{Spa}(\mathbb{Q}_p(p^{1/p^\infty})^{\wedge\flat})$:

\begin{itemize}
\item (\text{Proposition}) There is a functor (global section) between the $\infty$-categories of inductive Banach quasicoherent sheaves\footnote{Here $\circledcirc=\text{solidquasicoherentsheaves}$.}:
\[
\xymatrix@R+6pc@C+0pc{
\mathrm{Module}_\circledcirc(\mathcal{O}_{X_{\mathbb{Q}_p(p^{1/p^\infty})^{\wedge}\left<T_1^{\pm 1/p^\infty},...,T_k^{\pm 1/p^\infty}\right>^\flat,A,q}})\ar[d]\ar[d]\ar[d]\ar[d] \ar[r]^{\mathrm{global}}\ar[r]\ar[r] &\varphi_q\mathrm{Module}_\circledcirc(\{\mathrm{Robba}^{\mathrm{extended},q}_{{R_k,A,I}}\}_I) \ar[d]\ar[d]\ar[d]\ar[d].\\
\mathrm{Module}_\circledcirc(\mathcal{O}_{X_{\mathbb{Q}_p(p^{1/p^\infty})^{\wedge\flat},A,q}})\ar[r]^{\mathrm{global}}\ar[r]\ar[r] &\varphi_q\mathrm{Module}_\circledcirc(\{\mathrm{Robba}^{\mathrm{extended},q}_{{R_0,A,I}}\}_I).\\ 
}
\]

\item (\text{Proposition}) There is a functor (global section) between the $\infty$-categories of inductive Banach quasicoherent commutative algebra $E_\infty$ objects\footnote{Here $\circledcirc=\text{solidquasicoherentsheaves}$.}:
\[
\xymatrix@R+6pc@C+0pc{
\mathrm{sComm}_\mathrm{simplicial}\mathrm{Module}_\circledcirc(\mathcal{O}_{X_{R_k,A,q}})\ar[d]\ar[d]\ar[d]\ar[d]\ar[r]^{\mathrm{global}}\ar[r]\ar[r] &\mathrm{sComm}_\mathrm{simplicial}\varphi_q\mathrm{Module}_\circledcirc(\{\mathrm{Robba}^{\mathrm{extended},q}_{{R_k,A,I}}\}_I)\ar[d]\ar[d]\ar[d]\ar[d]\\
\mathrm{sComm}_\mathrm{simplicial}\mathrm{Module}_\circledcirc(\mathcal{O}_{X_{\mathbb{Q}_p(p^{1/p^\infty})^{\wedge\flat},A,q}})\ar[r]^{\mathrm{global}}\ar[r]\ar[r] &\mathrm{sComm}_\mathrm{simplicial}\varphi_q\mathrm{Module}_\circledcirc(\{\mathrm{Robba}^{\mathrm{extended},q}_{{R_0,A,I}}\}_I).  
}
\]

\item Then as in \cite{LBV} we have a functor (global section) of the de Rham complex after \cite[Definition 5.9, Section 5.2.1]{KKM}\footnote{Here $\circledcirc=\text{solidquasicoherentsheaves}$.}:
\[\displayindent=-0.4in
\xymatrix@R+6pc@C+0pc{
\mathrm{deRham}_{\mathrm{sComm}_\mathrm{simplicial}\mathrm{Module}_\circledcirc(\mathcal{O}_{X_{R_k,A,q}})\ar[r]^{\mathrm{global}}}(-)\ar[d]\ar[d]\ar[d]\ar[d]\ar[r]\ar[r] &\mathrm{deRham}_{\mathrm{sComm}_\mathrm{simplicial}\varphi_q\mathrm{Module}_\circledcirc(\{\mathrm{Robba}^{\mathrm{extended},q}_{{R_k,A,I}}\}_I)}(-)\ar[d]\ar[d]\ar[d]\ar[d]\\
\mathrm{deRham}_{\mathrm{sComm}_\mathrm{simplicial}\mathrm{Module}_\circledcirc(\mathcal{O}_{X_{\mathbb{Q}_p(p^{1/p^\infty})^{\wedge\flat},A,q}})\ar[r]^{\mathrm{global}}}(-)\ar[r]\ar[r] &\mathrm{deRham}_{\mathrm{sComm}_\mathrm{simplicial}\varphi_q\mathrm{Module}_\circledcirc(\{\mathrm{Robba}^{\mathrm{extended},q}_{{R_0,A,I}}\}_I)}(-). 
}
\]

\item Then we have the following a functor (global section) of $K$-group $(\infty,1)$-spectrum from \cite{BGT}\footnote{Here $\circledcirc=\text{solidquasicoherentsheaves}$.}:
\[
\xymatrix@R+6pc@C+0pc{
\mathrm{K}^\mathrm{BGT}_{\mathrm{sComm}_\mathrm{simplicial}\mathrm{Module}_\circledcirc(\mathcal{O}_{X_{R_k,A,q}})\ar[r]^{\mathrm{global}}}(-)\ar[d]\ar[d]\ar[d]\ar[d]\ar[r]\ar[r] &\mathrm{K}^\mathrm{BGT}_{\mathrm{sComm}_\mathrm{simplicial}\varphi_q\mathrm{Module}_\circledcirc(\{\mathrm{Robba}^{\mathrm{extended},q}_{{R_k,A,I}}\}_I)}(-)\ar[d]\ar[d]\ar[d]\ar[d]\\
\mathrm{K}^\mathrm{BGT}_{\mathrm{sComm}_\mathrm{simplicial}\mathrm{Module}_\circledcirc(\mathcal{O}_{X_{\mathbb{Q}_p(p^{1/p^\infty})^{\wedge\flat},A,q}})\ar[r]^{\mathrm{global}}}(-)\ar[r]\ar[r] &\mathrm{K}^\mathrm{BGT}_{\mathrm{sComm}_\mathrm{simplicial}\varphi_q\mathrm{Module}_\circledcirc(\{\mathrm{Robba}^{\mathrm{extended},q}_{{R_0,A,I}}\}_I)}(-). 
}
\]

\end{itemize}

\
\indent Then we consider further equivariance by considering the arithmetic profinite fundamental groups $\Gamma_{\mathbb{Q}_p}$ and $\mathrm{Gal}(\overline{\mathbb{Q}_p\left<T_1^{\pm 1},...,T_k^{\pm 1}\right>}/R_k)$ through the following diagram:

\[
\xymatrix@R+0pc@C+0pc{
\mathbb{Z}_p^k=\mathrm{Gal}(R_k/{\mathbb{Q}_p(p^{1/p^\infty})^\wedge\left<T_1^{\pm 1},...,T_k^{\pm 1}\right>}) \ar[r]\ar[r] \ar[r]\ar[r] &\Gamma_k^{\times q}:=\mathrm{Gal}(R_k/{\mathbb{Q}_p\left<T_1^{\pm 1},...,T_k^{\pm 1}\right>}) \ar[r] \ar[r]\ar[r] &\Gamma_{\mathbb{Q}_p}.
}
\]

\begin{itemize}
\item (\text{Proposition}) There is a functor (global section) between the $\infty$-categories of inductive Banach quasicoherent sheaves\footnote{Here $\circledcirc=\text{solidquasicoherentsheaves}$.}:
\[
\xymatrix@R+6pc@C+0pc{
{\mathrm{Module}_\circledcirc}_{\Gamma_k^{\times q}}(\mathcal{O}_{X_{\mathbb{Q}_p(p^{1/p^\infty})^{\wedge}\left<T_1^{\pm 1/p^\infty},...,T_k^{\pm 1/p^\infty}\right>^\flat,A,q}})\ar[d]\ar[d]\ar[d]\ar[d] \ar[r]^{\mathrm{global}}\ar[r]\ar[r] &\varphi_q{\mathrm{Module}_\circledcirc}_{\Gamma_k^{\times q}}(\{\mathrm{Robba}^{\mathrm{extended},q}_{{R_k,A,I}}\}_I) \ar[d]\ar[d]\ar[d]\ar[d].\\
{\mathrm{Module}_\circledcirc}(\mathcal{O}_{X_{\mathbb{Q}_p(p^{1/p^\infty})^{\wedge\flat},A,q}})\ar[r]^{\mathrm{global}}\ar[r]\ar[r] &\varphi_q{\mathrm{Module}_\circledcirc}(\{\mathrm{Robba}^{\mathrm{extended},q}_{{R_0,A,I}}\}_I).\\ 
}
\]

\item (\text{Proposition}) There is a functor (global section) between the $\infty$-categories of inductive Banach quasicoherent commutative algebra $E_\infty$ objects\footnote{Here $\circledcirc=\text{solidquasicoherentsheaves}$.}:
\[\displayindent=-0.4in
\xymatrix@R+6pc@C+0pc{
\mathrm{sComm}_\mathrm{simplicial}{\mathrm{Module}_\circledcirc}_{\Gamma_k^{\times q}}(\mathcal{O}_{X_{R_k,A,q}})\ar[d]\ar[d]\ar[d]\ar[d]\ar[r]^{\mathrm{global}}\ar[r]\ar[r] &\mathrm{sComm}_\mathrm{simplicial}\varphi_q{\mathrm{Module}_\circledcirc}(\{\mathrm{Robba}^{\mathrm{extended},q}_{{R_k,A,I}}\}_I)\ar[d]\ar[d]\ar[d]\ar[d]\\
\mathrm{sComm}_\mathrm{simplicial}{\mathrm{Module}_\circledcirc}_{\Gamma_0^{\times q}}(\mathcal{O}_{X_{\mathbb{Q}_p(p^{1/p^\infty})^{\wedge\flat},A,q}})\ar[r]^{\mathrm{global}}\ar[r]\ar[r] &\mathrm{sComm}_\mathrm{simplicial}\varphi_q{\mathrm{Modules}_\circledcirc}_{\Gamma_0^{\times q}}(\{\mathrm{Robba}^{\mathrm{extended},q}_{{R_0,A,I}}\}_I).  
}
\]

\item Then as in \cite{LBV} we have a functor (global section) of the de Rham complex after \cite[Definition 5.9, Section 5.2.1]{KKM}\footnote{Here $\circledcirc=\text{solidquasicoherentsheaves}$.}:
\[\displayindent=-0.5in
\xymatrix@R+6pc@C+0pc{
\mathrm{deRham}_{\mathrm{sComm}_\mathrm{simplicial}{\mathrm{Modules}_\circledcirc}_{\Gamma_k^{\times q}}(\mathcal{O}_{X_{R_k,A,q}})\ar[r]^{\mathrm{global}}}(-)\ar[d]\ar[d]\ar[d]\ar[d]\ar[r]\ar[r] &\mathrm{deRham}_{\mathrm{sComm}_\mathrm{simplicial}\varphi_q{\mathrm{Modules}_\circledcirc}_{\Gamma_k^{\times q}}(\{\mathrm{Robba}^{\mathrm{extended},q}_{{R_k,A,I}}\}_I)}(-)\ar[d]\ar[d]\ar[d]\ar[d]\\
\mathrm{deRham}_{\mathrm{sComm}_\mathrm{simplicial}{\mathrm{Modules}_\circledcirc}_{\Gamma_0^{\times q}}(\mathcal{O}_{X_{\mathbb{Q}_p(p^{1/p^\infty})^{\wedge\flat},A,q}})\ar[r]^{\mathrm{global}}}(-)\ar[r]\ar[r] &\mathrm{deRham}_{\mathrm{sComm}_\mathrm{simplicial}\varphi_q{\mathrm{Modules}_\circledcirc}_{\Gamma_0^{\times q}}(\{\mathrm{Robba}^{\mathrm{extended},q}_{{R_0,A,I}}\}_I)}(-). 
}
\]

\item Then we have the following a functor (global section) of $K$-group $(\infty,1)$-spectrum from \cite{BGT}\footnote{Here $\circledcirc=\text{solidquasicoherentsheaves}$.}:
\[
\xymatrix@R+6pc@C+0pc{
\mathrm{K}^\mathrm{BGT}_{\mathrm{sComm}_\mathrm{simplicial}{\mathrm{Modules}_\circledcirc}_{\Gamma_k^{\times q}}(\mathcal{O}_{X_{R_k,A,q}})\ar[r]^{\mathrm{global}}}(-)\ar[d]\ar[d]\ar[d]\ar[d]\ar[r]\ar[r] &\mathrm{K}^\mathrm{BGT}_{\mathrm{sComm}_\mathrm{simplicial}\varphi_q{\mathrm{Modules}_\circledcirc}_{\Gamma_k^{\times q}}(\{\mathrm{Robba}^{\mathrm{extended},q}_{{R_k,A,I}}\}_I)}(-)\ar[d]\ar[d]\ar[d]\ar[d]\\
\mathrm{K}^\mathrm{BGT}_{\mathrm{sComm}_\mathrm{simplicial}{\mathrm{Modules}_\circledcirc}_{\Gamma_0^{\times q}}(\mathcal{O}_{X_{\mathbb{Q}_p(p^{1/p^\infty})^{\wedge\flat},A,q}})\ar[r]^{\mathrm{global}}}(-)\ar[r]\ar[r] &\mathrm{K}^\mathrm{BGT}_{\mathrm{sComm}_\mathrm{simplicial}\varphi_q{\mathrm{Modules}_\circledcirc}_{\Gamma_0^{\times q}}(\{\mathrm{Robba}^{\mathrm{extended},q}_{{R_0,A,I}}\}_I)}(-). 
}
\]
	
\end{itemize}

\begin{proposition}
All the global functors from \cite[Proposition 13.8, Theorem 14.9, Remark 14.10]{1CS2} achieve the equivalences on both sides.	
\end{proposition}

\newpage
\subsection{$\infty$-Categorical Analytic Stacks and Descents III}

\indent As before, we have the following. Let $\mathcal{A}$ vary in the category of all the Banach algebras over $\mathbb{Q}_p$ we have the following.

\begin{itemize}

\item (\text{Proposition}) There is a functor (global section) between the $\infty$-prestacks of inductive Banach quasicoherent presheaves:
\[
\xymatrix@R+0pc@C+0pc{
\mathrm{Ind}\mathrm{Banach}(\mathcal{O}_{X_{R,-,q}})\ar[r]^{\mathrm{global}}\ar[r]\ar[r] &\varphi_q\mathrm{Ind}\mathrm{Banach}(\{\mathrm{Robba}^\mathrm{extended}_{{R,-,I,q}}\}_I).  
}
\]
\item (\text{Proposition}) There is a functor (global section) between the $\infty$-prestacks of monomorphic inductive Banach quasicoherent presheaves:
\[
\xymatrix@R+0pc@C+0pc{
\mathrm{Ind}^m\mathrm{Banach}(\mathcal{O}_{X_{R,-,q}})\ar[r]^{\mathrm{global}}\ar[r]\ar[r] &\varphi_q\mathrm{Ind}^m\mathrm{Banach}(\{\mathrm{Robba}^\mathrm{extended}_{{R,-,I,q}}\}_I).  
}
\]

\item (\text{Proposition}) There is a functor (global section) between the $\infty$-prestacks of inductive Banach quasicoherent presheaves:
\[
\xymatrix@R+0pc@C+0pc{
\mathrm{Ind}\mathrm{Banach}(\mathcal{O}_{X_{R,-,q}})\ar[r]^{\mathrm{global}}\ar[r]\ar[r] &\varphi_q\mathrm{Ind}\mathrm{Banach}(\{\mathrm{Robba}^\mathrm{extended}_{{R,-,I,q}}\}_I).  
}
\]
\item (\text{Proposition}) There is a functor (global section) between the $\infty$-stacks of monomorphic inductive Banach quasicoherent presheaves:
\[
\xymatrix@R+0pc@C+0pc{
\mathrm{Ind}^m\mathrm{Banach}(\mathcal{O}_{X_{R,-,q}})\ar[r]^{\mathrm{global}}\ar[r]\ar[r] &\varphi_q\mathrm{Ind}^m\mathrm{Banach}(\{\mathrm{Robba}^\mathrm{extended}_{{R,-,I,q}}\}_I).  
}
\]
\item (\text{Proposition}) There is a functor (global section) between the $\infty$-prestacks of inductive Banach quasicoherent commutative algebra $E_\infty$ objects:
\[
\xymatrix@R+0pc@C+0pc{
\mathrm{sComm}_\mathrm{simplicial}\mathrm{Ind}\mathrm{Banach}(\mathcal{O}_{X_{R,-,q}})\ar[r]^{\mathrm{global}}\ar[r]\ar[r] &\mathrm{sComm}_\mathrm{simplicial}\varphi_q\mathrm{Ind}\mathrm{Banach}(\{\mathrm{Robba}^\mathrm{extended}_{{R,-,I,q}}\}_I).  
}
\]
\item (\text{Proposition}) There is a functor (global section) between the $\infty$-prestacks of monomorphic inductive Banach quasicoherent commutative algebra $E_\infty$ objects:
\[
\xymatrix@R+0pc@C+0pc{
\mathrm{sComm}_\mathrm{simplicial}\mathrm{Ind}^m\mathrm{Banach}(\mathcal{O}_{X_{R,-,q}})\ar[r]^{\mathrm{global}}\ar[r]\ar[r] &\mathrm{sComm}_\mathrm{simplicial}\varphi_q\mathrm{Ind}^m\mathrm{Banach}(\{\mathrm{Robba}^\mathrm{extended}_{{R,-,I,q}}\}_I).  
}
\]

\item Then parallel as in \cite{LBV} we have a functor (global section) of the de Rham complex after \cite[Definition 5.9, Section 5.2.1]{KKM}:
\[
\xymatrix@R+0pc@C+0pc{
\mathrm{deRham}_{\mathrm{sComm}_\mathrm{simplicial}\mathrm{Ind}\mathrm{Banach}(\mathcal{O}_{X_{R,-,q}})\ar[r]^{\mathrm{global}}}(-)\ar[r]\ar[r] &\mathrm{deRham}_{\mathrm{sComm}_\mathrm{simplicial}\varphi_q\mathrm{Ind}\mathrm{Banach}(\{\mathrm{Robba}^\mathrm{extended}_{{R,-,I,q}}\}_I)}(-), 
}
\]
\[
\xymatrix@R+0pc@C+0pc{
\mathrm{deRham}_{\mathrm{sComm}_\mathrm{simplicial}\mathrm{Ind}^m\mathrm{Banach}(\mathcal{O}_{X_{R,-,q}})\ar[r]^{\mathrm{global}}}(-)\ar[r]\ar[r] &\mathrm{deRham}_{\mathrm{sComm}_\mathrm{simplicial}\varphi_q\mathrm{Ind}^m\mathrm{Banach}(\{\mathrm{Robba}^\mathrm{extended}_{{R,-,I,q}}\}_I)}(-). 
}
\]

\item Then we have the following a functor (global section) of $K$-group $(\infty,1)$-spectrum from \cite{BGT}:
\[
\xymatrix@R+0pc@C+0pc{
\mathrm{K}^\mathrm{BGT}_{\mathrm{sComm}_\mathrm{simplicial}\mathrm{Ind}\mathrm{Banach}(\mathcal{O}_{X_{R,-,q}})\ar[r]^{\mathrm{global}}}(-)\ar[r]\ar[r] &\mathrm{K}^\mathrm{BGT}_{\mathrm{sComm}_\mathrm{simplicial}\varphi_q\mathrm{Ind}\mathrm{Banach}(\{\mathrm{Robba}^\mathrm{extended}_{{R,-,I,q}}\}_I)}(-), 
}
\]
\[
\xymatrix@R+0pc@C+0pc{
\mathrm{K}^\mathrm{BGT}_{\mathrm{sComm}_\mathrm{simplicial}\mathrm{Ind}^m\mathrm{Banach}(\mathcal{O}_{X_{R,-,q}})\ar[r]^{\mathrm{global}}}(-)\ar[r]\ar[r] &\mathrm{K}^\mathrm{BGT}_{\mathrm{sComm}_\mathrm{simplicial}\varphi_q\mathrm{Ind}^m\mathrm{Banach}(\{\mathrm{Robba}^\mathrm{extended}_{{R,-,I,q}}\}_I)}(-). 
}
\]
\end{itemize}

\noindent Now let $R=\mathbb{Q}_p(p^{1/p^\infty})^{\wedge\flat}$ and $R_k=\mathbb{Q}_p(p^{1/p^\infty})^{\wedge}\left<T_1^{\pm 1/p^{\infty}},...,T_k^{\pm 1/p^{\infty}}\right>^\flat$ we have the following Galois theoretic results with naturality along $f:\mathrm{Spa}(\mathbb{Q}_p(p^{1/p^\infty})^{\wedge}\left<T_1^{\pm 1/p^\infty},...,T_k^{\pm 1/p^\infty}\right>^\flat)\rightarrow \mathrm{Spa}(\mathbb{Q}_p(p^{1/p^\infty})^{\wedge\flat})$:

\begin{itemize}
\item (\text{Proposition}) There is a functor (global section) between the $\infty$-prestacks of inductive Banach quasicoherent presheaves:
\[
\xymatrix@R+6pc@C+0pc{
\mathrm{Ind}\mathrm{Banach}(\mathcal{O}_{X_{\mathbb{Q}_p(p^{1/p^\infty})^{\wedge}\left<T_1^{\pm 1/p^\infty},...,T_k^{\pm 1/p^\infty}\right>^\flat,-,q}})\ar[d]\ar[d]\ar[d]\ar[d] \ar[r]^{\mathrm{global}}\ar[r]\ar[r] &\varphi_q\mathrm{Ind}\mathrm{Banach}(\{\mathrm{Robba}^\mathrm{extended}_{{R_k,-,I,q}}\}_I) \ar[d]\ar[d]\ar[d]\ar[d].\\
\mathrm{Ind}\mathrm{Banach}(\mathcal{O}_{X_{\mathbb{Q}_p(p^{1/p^\infty})^{\wedge\flat},-,q}})\ar[r]^{\mathrm{global}}\ar[r]\ar[r] &\varphi_q\mathrm{Ind}\mathrm{Banach}(\{\mathrm{Robba}^\mathrm{extended}_{{R_0,-,I,q}}\}_I).\\ 
}
\]
\item (\text{Proposition}) There is a functor (global section) between the $\infty$-prestacks of monomorphic inductive Banach quasicoherent presheaves:
\[
\xymatrix@R+6pc@C+0pc{
\mathrm{Ind}^m\mathrm{Banach}(\mathcal{O}_{X_{R_k,-,q}})\ar[r]^{\mathrm{global}}\ar[d]\ar[d]\ar[d]\ar[d]\ar[r]\ar[r] &\varphi_q\mathrm{Ind}^m\mathrm{Banach}(\{\mathrm{Robba}^\mathrm{extended}_{{R_k,-,I,q}}\}_I)\ar[d]\ar[d]\ar[d]\ar[d]\\
\mathrm{Ind}^m\mathrm{Banach}(\mathcal{O}_{X_{\mathbb{Q}_p(p^{1/p^\infty})^{\wedge\flat},-,q}})\ar[r]^{\mathrm{global}}\ar[r]\ar[r] &\varphi_q\mathrm{Ind}^m\mathrm{Banach}(\{\mathrm{Robba}^\mathrm{extended}_{{R_0,-,I,q}}\}_I).\\  
}
\]
\item (\text{Proposition}) There is a functor (global section) between the $\infty$-prestacks of inductive Banach quasicoherent commutative algebra $E_\infty$ objects:
\[
\xymatrix@R+6pc@C+0pc{
\mathrm{sComm}_\mathrm{simplicial}\mathrm{Ind}\mathrm{Banach}(\mathcal{O}_{X_{R_k,-,q}})\ar[d]\ar[d]\ar[d]\ar[d]\ar[r]^{\mathrm{global}}\ar[r]\ar[r] &\mathrm{sComm}_\mathrm{simplicial}\varphi_q\mathrm{Ind}\mathrm{Banach}(\{\mathrm{Robba}^\mathrm{extended}_{{R_k,-,I,q}}\}_I)\ar[d]\ar[d]\ar[d]\ar[d]\\
\mathrm{sComm}_\mathrm{simplicial}\mathrm{Ind}\mathrm{Banach}(\mathcal{O}_{X_{\mathbb{Q}_p(p^{1/p^\infty})^{\wedge\flat},-,q}})\ar[r]^{\mathrm{global}}\ar[r]\ar[r] &\mathrm{sComm}_\mathrm{simplicial}\varphi_q\mathrm{Ind}\mathrm{Banach}(\{\mathrm{Robba}^\mathrm{extended}_{{R_0,-,I,q}}\}_I).  
}
\]
\item (\text{Proposition}) There is a functor (global section) between the $\infty$-prestacks of monomorphic inductive Banach quasicoherent commutative algebra $E_\infty$ objects:
\[\displayindent=-0.4in
\xymatrix@R+6pc@C+0pc{
\mathrm{sComm}_\mathrm{simplicial}\mathrm{Ind}^m\mathrm{Banach}(\mathcal{O}_{X_{R_k,-,q}})\ar[d]\ar[d]\ar[d]\ar[d]\ar[r]^{\mathrm{global}}\ar[r]\ar[r] &\mathrm{sComm}_\mathrm{simplicial}\varphi_q\mathrm{Ind}^m\mathrm{Banach}(\{\mathrm{Robba}^\mathrm{extended}_{{R_k,-,I,q}}\}_I)\ar[d]\ar[d]\ar[d]\ar[d]\\
 \mathrm{sComm}_\mathrm{simplicial}\mathrm{Ind}^m\mathrm{Banach}(\mathcal{O}_{X_{\mathbb{Q}_p(p^{1/p^\infty})^{\wedge\flat},-,q}})\ar[r]^{\mathrm{global}}\ar[r]\ar[r] &\mathrm{sComm}_\mathrm{simplicial}\varphi_q\mathrm{Ind}^m\mathrm{Banach}(\{\mathrm{Robba}^\mathrm{extended}_{{R_0,-,I,q}}\}_I).
}
\]

\item Then parallel as in \cite{LBV} we have a functor (global section) of the de Rham complex after \cite[Definition 5.9, Section 5.2.1]{KKM}:
\[\displayindent=-0.4in
\xymatrix@R+6pc@C+0pc{
\mathrm{deRham}_{\mathrm{sComm}_\mathrm{simplicial}\mathrm{Ind}\mathrm{Banach}(\mathcal{O}_{X_{R_k,-,q}})\ar[r]^{\mathrm{global}}}(-)\ar[d]\ar[d]\ar[d]\ar[d]\ar[r]\ar[r] &\mathrm{deRham}_{\mathrm{sComm}_\mathrm{simplicial}\varphi_q\mathrm{Ind}\mathrm{Banach}(\{\mathrm{Robba}^\mathrm{extended}_{{R_k,-,I,q}}\}_I)}(-)\ar[d]\ar[d]\ar[d]\ar[d]\\
\mathrm{deRham}_{\mathrm{sComm}_\mathrm{simplicial}\mathrm{Ind}\mathrm{Banach}(\mathcal{O}_{X_{\mathbb{Q}_p(p^{1/p^\infty})^{\wedge\flat},-,q}})\ar[r]^{\mathrm{global}}}(-)\ar[r]\ar[r] &\mathrm{deRham}_{\mathrm{sComm}_\mathrm{simplicial}\varphi_q\mathrm{Ind}\mathrm{Banach}(\{\mathrm{Robba}^\mathrm{extended}_{{R_0,-,I,q}}\}_I)}(-), 
}
\]
\[\displayindent=-0.4in
\xymatrix@R+6pc@C+0pc{
\mathrm{deRham}_{\mathrm{sComm}_\mathrm{simplicial}\mathrm{Ind}^m\mathrm{Banach}(\mathcal{O}_{X_{R_k,-,q}})\ar[r]^{\mathrm{global}}}(-)\ar[d]\ar[d]\ar[d]\ar[d]\ar[r]\ar[r] &\mathrm{deRham}_{\mathrm{sComm}_\mathrm{simplicial}\varphi_q\mathrm{Ind}^m\mathrm{Banach}(\{\mathrm{Robba}^\mathrm{extended}_{{R_k,-,I,q}}\}_I)}(-)\ar[d]\ar[d]\ar[d]\ar[d]\\
\mathrm{deRham}_{\mathrm{sComm}_\mathrm{simplicial}\mathrm{Ind}^m\mathrm{Banach}(\mathcal{O}_{X_{\mathbb{Q}_p(p^{1/p^\infty})^{\wedge\flat},-,q}})\ar[r]^{\mathrm{global}}}(-)\ar[r]\ar[r] &\mathrm{deRham}_{\mathrm{sComm}_\mathrm{simplicial}\varphi_q\mathrm{Ind}^m\mathrm{Banach}(\{\mathrm{Robba}^\mathrm{extended}_{{R_0,-,I,q}}\}_I)}(-). 
}
\]

\item Then we have the following a functor (global section) of $K$-group $(\infty,1)$-spectrum from \cite{BGT}:
\[
\xymatrix@R+6pc@C+0pc{
\mathrm{K}^\mathrm{BGT}_{\mathrm{sComm}_\mathrm{simplicial}\mathrm{Ind}\mathrm{Banach}(\mathcal{O}_{X_{R_k,-,q}})\ar[r]^{\mathrm{global}}}(-)\ar[d]\ar[d]\ar[d]\ar[d]\ar[r]\ar[r] &\mathrm{K}^\mathrm{BGT}_{\mathrm{sComm}_\mathrm{simplicial}\varphi_q\mathrm{Ind}\mathrm{Banach}(\{\mathrm{Robba}^\mathrm{extended}_{{R_k,-,I,q}}\}_I)}(-)\ar[d]\ar[d]\ar[d]\ar[d]\\
\mathrm{K}^\mathrm{BGT}_{\mathrm{sComm}_\mathrm{simplicial}\mathrm{Ind}\mathrm{Banach}(\mathcal{O}_{X_{\mathbb{Q}_p(p^{1/p^\infty})^{\wedge\flat},-,q}})\ar[r]^{\mathrm{global}}}(-)\ar[r]\ar[r] &\mathrm{K}^\mathrm{BGT}_{\mathrm{sComm}_\mathrm{simplicial}\varphi_q\mathrm{Ind}\mathrm{Banach}(\{\mathrm{Robba}^\mathrm{extended}_{{R_0,-,I,q}}\}_I)}(-), 
}
\]
\[
\xymatrix@R+6pc@C+0pc{
\mathrm{K}^\mathrm{BGT}_{\mathrm{sComm}_\mathrm{simplicial}\mathrm{Ind}^m\mathrm{Banach}(\mathcal{O}_{X_{R_k,-,q}})\ar[r]^{\mathrm{global}}}(-)\ar[d]\ar[d]\ar[d]\ar[d]\ar[r]\ar[r] &\mathrm{K}^\mathrm{BGT}_{\mathrm{sComm}_\mathrm{simplicial}\varphi_q\mathrm{Ind}^m\mathrm{Banach}(\{\mathrm{Robba}^\mathrm{extended}_{{R_k,-,I,q}}\}_I)}(-)\ar[d]\ar[d]\ar[d]\ar[d]\\
\mathrm{K}^\mathrm{BGT}_{\mathrm{sComm}_\mathrm{simplicial}\mathrm{Ind}^m\mathrm{Banach}(\mathcal{O}_{X_{\mathbb{Q}_p(p^{1/p^\infty})^{\wedge\flat},-,q}})\ar[r]^{\mathrm{global}}}(-)\ar[r]\ar[r] &\mathrm{K}^\mathrm{BGT}_{\mathrm{sComm}_\mathrm{simplicial}\varphi_q\mathrm{Ind}^m\mathrm{Banach}(\{\mathrm{Robba}^\mathrm{extended}_{{R_0,-,I,q}}\}_I)}(-). 
}
\]

\end{itemize}

\
\indent Then we consider further equivariance by considering the arithmetic profinite fundamental groups and actually its $q$-th power $\mathrm{Gal}(\overline{{Q}_p\left<T_1^{\pm 1},...,T_k^{\pm 1}\right>}/R_k)^{\times q}$ through the following diagram:\\

\[
\xymatrix@R+6pc@C+0pc{
\mathbb{Z}_p^k=\mathrm{Gal}(R_k/{\mathbb{Q}_p(p^{1/p^\infty})^\wedge\left<T_1^{\pm 1},...,T_k^{\pm 1}\right>})\ar[d]\ar[d]\ar[d]\ar[d] \ar[r]\ar[r] \ar[r]\ar[r] &\mathrm{Gal}(\overline{{Q}_p\left<T_1^{\pm 1},...,T_k^{\pm 1}\right>}/R_k) \ar[d]\ar[d]\ar[d] \ar[r]\ar[r] &\Gamma_{\mathbb{Q}_p} \ar[d]\ar[d]\ar[d]\ar[d]\\
(\mathbb{Z}_p^k=\mathrm{Gal}(R_k/{\mathbb{Q}_p(p^{1/p^\infty})^\wedge\left<T_1^{\pm 1},...,T_k^{\pm 1}\right>}))^{\times q} \ar[r]\ar[r] \ar[r]\ar[r] &\Gamma_k^{\times q}:=\mathrm{Gal}(R_k/{\mathbb{Q}_p\left<T_1^{\pm 1},...,T_k^{\pm 1}\right>})^{\times q} \ar[r] \ar[r]\ar[r] &\Gamma_{\mathbb{Q}_p}^{\times q}.
}
\]

\

We then have the correspond arithmetic profinite fundamental groups equivariant versions:
\begin{itemize}
\item (\text{Proposition}) There is a functor (global section) between the $\infty$-prestacks of inductive Banach quasicoherent presheaves:
\[
\xymatrix@R+6pc@C+0pc{
\mathrm{Ind}\mathrm{Banach}_{\Gamma_{k}^{\times q}}(\mathcal{O}_{X_{\mathbb{Q}_p(p^{1/p^\infty})^{\wedge}\left<T_1^{\pm 1/p^\infty},...,T_k^{\pm 1/p^\infty}\right>^\flat,-,q}})\ar[d]\ar[d]\ar[d]\ar[d] \ar[r]^{\mathrm{global}}\ar[r]\ar[r] &\varphi_q\mathrm{Ind}\mathrm{Banach}_{\Gamma_{k}^{\times q}}(\{\mathrm{Robba}^\mathrm{extended}_{{R_k,-,I,q}}\}_I) \ar[d]\ar[d]\ar[d]\ar[d].\\
\mathrm{Ind}\mathrm{Banach}(\mathcal{O}_{X_{\mathbb{Q}_p(p^{1/p^\infty})^{\wedge\flat},-,q}})\ar[r]^{\mathrm{global}}\ar[r]\ar[r] &\varphi_q\mathrm{Ind}\mathrm{Banach}(\{\mathrm{Robba}^\mathrm{extended}_{{R_0,-,I,q}}\}_I).\\ 
}
\]
\item (\text{Proposition}) There is a functor (global section) between the $\infty$-prestacks of monomorphic inductive Banach quasicoherent presheaves:
\[
\xymatrix@R+6pc@C+0pc{
\mathrm{Ind}^m\mathrm{Banach}_{\Gamma_{k}^{\times q}}(\mathcal{O}_{X_{R_k,-,q}})\ar[r]^{\mathrm{global}}\ar[d]\ar[d]\ar[d]\ar[d]\ar[r]\ar[r] &\varphi_q\mathrm{Ind}^m\mathrm{Banach}_{\Gamma_{k}^{\times q}}(\{\mathrm{Robba}^\mathrm{extended}_{{R_k,-,I,q}}\}_I)\ar[d]\ar[d]\ar[d]\ar[d]\\
\mathrm{Ind}^m\mathrm{Banach}_{\Gamma_{0}^{\times q}}(\mathcal{O}_{X_{\mathbb{Q}_p(p^{1/p^\infty})^{\wedge\flat},-,q}})\ar[r]^{\mathrm{global}}\ar[r]\ar[r] &\varphi_q\mathrm{Ind}^m\mathrm{Banach}_{\Gamma_{0}^{\times q}}(\{\mathrm{Robba}^\mathrm{extended}_{{R_0,-,I,q}}\}_I).\\  
}
\]
\item (\text{Proposition}) There is a functor (global section) between the $\infty$-stacks of inductive Banach quasicoherent commutative algebra $E_\infty$ objects:
\[\displayindent=-0.4in
\xymatrix@R+6pc@C+0pc{
\mathrm{sComm}_\mathrm{simplicial}\mathrm{Ind}\mathrm{Banach}_{\Gamma_{k}^{\times q}}(\mathcal{O}_{X_{R_k,-,q}})\ar[d]\ar[d]\ar[d]\ar[d]\ar[r]^{\mathrm{global}}\ar[r]\ar[r] &\mathrm{sComm}_\mathrm{simplicial}\varphi_q\mathrm{Ind}\mathrm{Banach}_{\Gamma_{k}^{\times q}}(\{\mathrm{Robba}^\mathrm{extended}_{{R_k,-,I,q}}\}_I)\ar[d]\ar[d]\ar[d]\ar[d]\\
\mathrm{sComm}_\mathrm{simplicial}\mathrm{Ind}\mathrm{Banach}_{\Gamma_{0}^{\times q}}(\mathcal{O}_{X_{\mathbb{Q}_p(p^{1/p^\infty})^{\wedge\flat},-,q}})\ar[r]^{\mathrm{global}}\ar[r]\ar[r] &\mathrm{sComm}_\mathrm{simplicial}\varphi_q\mathrm{Ind}\mathrm{Banach}_{\Gamma_{0}^{\times q}}(\{\mathrm{Robba}^\mathrm{extended}_{{R_0,-,I,q}}\}_I).  
}
\]
\item (\text{Proposition}) There is a functor (global section) between the $\infty$-prestacks of monomorphic inductive Banach quasicoherent commutative algebra $E_\infty$ objects:
\[\displayindent=-0.4in
\xymatrix@R+6pc@C+0pc{
\mathrm{sComm}_\mathrm{simplicial}\mathrm{Ind}^m\mathrm{Banach}_{\Gamma_{k}^{\times q}}(\mathcal{O}_{X_{R_k,-,q}})\ar[d]\ar[d]\ar[d]\ar[d]\ar[r]^{\mathrm{global}}\ar[r]\ar[r] &\mathrm{sComm}_\mathrm{simplicial}\varphi_q\mathrm{Ind}^m\mathrm{Banach}_{\Gamma_{k}^{\times q}}(\{\mathrm{Robba}^\mathrm{extended}_{{R_k,-,I,q}}\}_I)\ar[d]\ar[d]\ar[d]\ar[d]\\
 \mathrm{sComm}_\mathrm{simplicial}\mathrm{Ind}^m\mathrm{Banach}_{\Gamma_{0}^{\times q}}(\mathcal{O}_{X_{\mathbb{Q}_p(p^{1/p^\infty})^{\wedge\flat},-,q}})\ar[r]^{\mathrm{global}}\ar[r]\ar[r] &\mathrm{sComm}_\mathrm{simplicial}\varphi_q\mathrm{Ind}^m\mathrm{Banach}_{\Gamma_{0}^{\times q}}(\{\mathrm{Robba}^\mathrm{extended}_{{R_0,-,I,q}}\}_I). 
}
\]

\item Then parallel as in \cite{LBV} we have a functor (global section) of the de Rham complex after \cite[Definition 5.9, Section 5.2.1]{KKM}:
\[\displayindent=-0.4in
\xymatrix@R+6pc@C+0pc{
\mathrm{deRham}_{\mathrm{sComm}_\mathrm{simplicial}\mathrm{Ind}\mathrm{Banach}_{\Gamma_{k}^{\times q}}(\mathcal{O}_{X_{R_k,-,q}})\ar[r]^{\mathrm{global}}}(-)\ar[d]\ar[d]\ar[d]\ar[d]\ar[r]\ar[r] &\mathrm{deRham}_{\mathrm{sComm}_\mathrm{simplicial}\varphi_q\mathrm{Ind}\mathrm{Banach}_{\Gamma_{k}^{\times q}}(\{\mathrm{Robba}^\mathrm{extended}_{{R_k,-,I,q}}\}_I)}(-)\ar[d]\ar[d]\ar[d]\ar[d]\\
\mathrm{deRham}_{\mathrm{sComm}_\mathrm{simplicial}\mathrm{Ind}\mathrm{Banach}_{\Gamma_{0}^{\times q}}(\mathcal{O}_{X_{\mathbb{Q}_p(p^{1/p^\infty})^{\wedge\flat},-,q}})\ar[r]^{\mathrm{global}}}(-)\ar[r]\ar[r] &\mathrm{deRham}_{\mathrm{sComm}_\mathrm{simplicial}\varphi_q\mathrm{Ind}\mathrm{Banach}_{\Gamma_{0}^{\times q}}(\{\mathrm{Robba}^\mathrm{extended}_{{R_0,-,I,q}}\}_I)}(-), 
}
\]
\[\displayindent=-0.7in
\xymatrix@R+6pc@C+0pc{
\mathrm{deRham}_{\mathrm{sComm}_\mathrm{simplicial}\mathrm{Ind}^m\mathrm{Banach}_{\Gamma_{k}^{\times q}}(\mathcal{O}_{X_{R_k,-,q}})\ar[r]^{\mathrm{global}}}(-)\ar[d]\ar[d]\ar[d]\ar[d]\ar[r]\ar[r] &\mathrm{deRham}_{\mathrm{sComm}_\mathrm{simplicial}\varphi_q\mathrm{Ind}^m\mathrm{Banach}_{\Gamma_{k}^{\times q}}(\{\mathrm{Robba}^\mathrm{extended}_{{R_k,-,I,q}}\}_I)}(-)\ar[d]\ar[d]\ar[d]\ar[d]\\
\mathrm{deRham}_{\mathrm{sComm}_\mathrm{simplicial}\mathrm{Ind}^m\mathrm{Banach}_{\Gamma_{0}^{\times q}}(\mathcal{O}_{X_{\mathbb{Q}_p(p^{1/p^\infty})^{\wedge\flat},-,q}})\ar[r]^{\mathrm{global}}}(-)\ar[r]\ar[r] &\mathrm{deRham}_{\mathrm{sComm}_\mathrm{simplicial}\varphi_q\mathrm{Ind}^m\mathrm{Banach}_{\Gamma_{0}^{\times q}}(\{\mathrm{Robba}^\mathrm{extended}_{{R_0,-,I,q}}\}_I)}(-). 
}
\]

\item Then we have the following a functor (global section) of $K$-group $(\infty,1)$-spectrum from \cite{BGT}:
\[
\xymatrix@R+6pc@C+0pc{
\mathrm{K}^\mathrm{BGT}_{\mathrm{sComm}_\mathrm{simplicial}\mathrm{Ind}\mathrm{Banach}_{\Gamma_{k}^{\times q}}(\mathcal{O}_{X_{R_k,-,q}})\ar[r]^{\mathrm{global}}}(-)\ar[d]\ar[d]\ar[d]\ar[d]\ar[r]\ar[r] &\mathrm{K}^\mathrm{BGT}_{\mathrm{sComm}_\mathrm{simplicial}\varphi_q\mathrm{Ind}\mathrm{Banach}_{\Gamma_{k}^{\times q}}(\{\mathrm{Robba}^\mathrm{extended}_{{R_k,-,I,q}}\}_I)}(-)\ar[d]\ar[d]\ar[d]\ar[d]\\
\mathrm{K}^\mathrm{BGT}_{\mathrm{sComm}_\mathrm{simplicial}\mathrm{Ind}\mathrm{Banach}_{\Gamma_{0}^{\times q}}(\mathcal{O}_{X_{\mathbb{Q}_p(p^{1/p^\infty})^{\wedge\flat},-,q}})\ar[r]^{\mathrm{global}}}(-)\ar[r]\ar[r] &\mathrm{K}^\mathrm{BGT}_{\mathrm{sComm}_\mathrm{simplicial}\varphi_q\mathrm{Ind}\mathrm{Banach}_{\Gamma_{0}^{\times q}}(\{\mathrm{Robba}^\mathrm{extended}_{{R_0,-,I,q}}\}_I)}(-), 
}
\]
\[
\xymatrix@R+6pc@C+0pc{
\mathrm{K}^\mathrm{BGT}_{\mathrm{sComm}_\mathrm{simplicial}\mathrm{Ind}^m\mathrm{Banach}_{\Gamma_{k}^{\times q}}(\mathcal{O}_{X_{R_k,-,q}})\ar[r]^{\mathrm{global}}}(-)\ar[d]\ar[d]\ar[d]\ar[d]\ar[r]\ar[r] &\mathrm{K}^\mathrm{BGT}_{\mathrm{sComm}_\mathrm{simplicial}\varphi_q\mathrm{Ind}^m\mathrm{Banach}_{\Gamma_{k}^{\times q}}(\{\mathrm{Robba}^\mathrm{extended}_{{R_k,-,I,q}}\}_I)}(-)\ar[d]\ar[d]\ar[d]\ar[d]\\
\mathrm{K}^\mathrm{BGT}_{\mathrm{sComm}_\mathrm{simplicial}\mathrm{Ind}^m\mathrm{Banach}_{\Gamma_{0}^{\times q}}(\mathcal{O}_{X_{\mathbb{Q}_p(p^{1/p^\infty})^{\wedge\flat},-,q}})\ar[r]^{\mathrm{global}}}(-)\ar[r]\ar[r] &\mathrm{K}^\mathrm{BGT}_{\mathrm{sComm}_\mathrm{simplicial}\varphi_q\mathrm{Ind}^m\mathrm{Banach}_{\Gamma_{0}^{\times q}}(\{\mathrm{Robba}^\mathrm{extended}_{{R_0,-,I,q}}\}_I)}(-). 
}
\]

\end{itemize}

\

\begin{remark}
\noindent We can certainly consider the quasicoherent sheaves in \cite[Lemma 7.11, Remark 7.12]{1BBK}, therefore all the quasicoherent presheaves and modules will be those satisfying \cite[Lemma 7.11, Remark 7.12]{1BBK} if one would like to consider the the quasicoherent sheaves. That being all as this said, we would believe that the big quasicoherent presheaves are automatically quasicoherent sheaves (namely satisfying the corresponding \v{C}ech $\infty$-descent as in \cite[Section 9.3]{KKM} and \cite[Lemma 7.11, Remark 7.12]{1BBK}) and the corresponding global section functors are automatically equivalence of $\infty$-categories.\\ 
\end{remark}

\

\indent In Clausen-Scholze formalism we have the following:

\begin{itemize}

\item (\text{Proposition}) There is a functor (global section) between the $\infty$-prestacks of inductive Banach quasicoherent sheaves:
\[
\xymatrix@R+0pc@C+0pc{
{\mathrm{Modules}_\circledcirc}(\mathcal{O}_{X_{R,-,q}})\ar[r]^{\mathrm{global}}\ar[r]\ar[r] &\varphi_q{\mathrm{Modules}_\circledcirc}(\{\mathrm{Robba}^{\mathrm{extended},q}_{{R,-,I}}\}_I).  
}
\]

\item (\text{Proposition}) There is a functor (global section) between the $\infty$-prestacks of inductive Banach quasicoherent sheaves:
\[
\xymatrix@R+0pc@C+0pc{
{\mathrm{Modules}_\circledcirc}(\mathcal{O}_{X_{R,-,q}})\ar[r]^{\mathrm{global}}\ar[r]\ar[r] &\varphi_q{\mathrm{Modules}_\circledcirc}(\{\mathrm{Robba}^{\mathrm{extended},q}_{{R,-,I}}\}_I).  
}
\]

\item (\text{Proposition}) There is a functor (global section) between the $\infty$-prestacks of inductive Banach quasicoherent commutative algebra $E_\infty$ objects\footnote{Here $\circledcirc=\text{solidquasicoherentsheaves}$.}:
\[
\xymatrix@R+0pc@C+0pc{
\mathrm{sComm}_\mathrm{simplicial}{\mathrm{Modules}_\circledcirc}(\mathcal{O}_{X_{R,-,q}})\ar[r]^{\mathrm{global}}\ar[r]\ar[r] &\mathrm{sComm}_\mathrm{simplicial}\varphi_q{\mathrm{Modules}_\circledcirc}(\{\mathrm{Robba}^{\mathrm{extended},q}_{{R,-,I}}\}_I).  
}
\]

\item Then as in \cite{LBV} we have a functor (global section) of the de Rham complex after \cite[Definition 5.9, Section 5.2.1]{KKM}\footnote{Here $\circledcirc=\text{solidquasicoherentsheaves}$.}:
\[
\xymatrix@R+0pc@C+0pc{
\mathrm{deRham}_{\mathrm{sComm}_\mathrm{simplicial}{\mathrm{Modules}_\circledcirc}(\mathcal{O}_{X_{R,-,q}})\ar[r]^{\mathrm{global}}}(-)\ar[r]\ar[r] &\mathrm{deRham}_{\mathrm{sComm}_\mathrm{simplicial}\varphi_q{\mathrm{Modules}_\circledcirc}(\{\mathrm{Robba}^{\mathrm{extended},q}_{{R,-,I}}\}_I)}(-). 
}
\]

\item Then we have the following a functor (global section) of $K$-group $(\infty,1)$-spectrum from \cite{BGT}\footnote{Here $\circledcirc=\text{solidquasicoherentsheaves}$.}:
\[
\xymatrix@R+0pc@C+0pc{
\mathrm{K}^\mathrm{BGT}_{\mathrm{sComm}_\mathrm{simplicial}{\mathrm{Modules}_\circledcirc}(\mathcal{O}_{X_{R,-,q}})\ar[r]^{\mathrm{global}}}(-)\ar[r]\ar[r] &\mathrm{K}^\mathrm{BGT}_{\mathrm{sComm}_\mathrm{simplicial}\varphi_q{\mathrm{Modules}_\circledcirc}(\{\mathrm{Robba}^{\mathrm{extended},q}_{{R,-,I}}\}_I)}(-). 
}
\]
\end{itemize}

\noindent Now let $R=\mathbb{Q}_p(p^{1/p^\infty})^{\wedge\flat}$ and $R_k=\mathbb{Q}_p(p^{1/p^\infty})^{\wedge}\left<T_1^{\pm 1/p^{\infty}},...,T_k^{\pm 1/p^{\infty}}\right>^\flat$ we have the following Galois theoretic results with naturality along $f:\mathrm{Spa}(\mathbb{Q}_p(p^{1/p^\infty})^{\wedge}\left<T_1^{\pm 1/p^\infty},...,T_k^{\pm 1/p^\infty}\right>^\flat)\rightarrow \mathrm{Spa}(\mathbb{Q}_p(p^{1/p^\infty})^{\wedge\flat})$:

\begin{itemize}
\item (\text{Proposition}) There is a functor (global section) between the $\infty$-prestacks of inductive Banach quasicoherent sheaves\footnote{Here $\circledcirc=\text{solidquasicoherentsheaves}$.}:
\[
\xymatrix@R+6pc@C+0pc{
{\mathrm{Modules}_\circledcirc}(\mathcal{O}_{X_{\mathbb{Q}_p(p^{1/p^\infty})^{\wedge}\left<T_1^{\pm 1/p^\infty},...,T_k^{\pm 1/p^\infty}\right>^\flat,-,q}})\ar[d]\ar[d]\ar[d]\ar[d] \ar[r]^{\mathrm{global}}\ar[r]\ar[r] &\varphi_q{\mathrm{Modules}_\circledcirc}(\{\mathrm{Robba}^{\mathrm{extended},q}_{{R_k,-,I}}\}_I) \ar[d]\ar[d]\ar[d]\ar[d].\\
{\mathrm{Modules}_\circledcirc}(\mathcal{O}_{X_{\mathbb{Q}_p(p^{1/p^\infty})^{\wedge\flat},-,q}})\ar[r]^{\mathrm{global}}\ar[r]\ar[r] &\varphi_q{\mathrm{Modules}_\circledcirc}(\{\mathrm{Robba}^{\mathrm{extended},q}_{{R_0,-,I}}\}_I).\\ 
}
\]
\item (\text{Proposition}) There is a functor (global section) between the $\infty$-prestacks of inductive Banach quasicoherent commutative algebra $E_\infty$ objects\footnote{Here $\circledcirc=\text{solidquasicoherentsheaves}$.}:
\[
\xymatrix@R+6pc@C+0pc{
\mathrm{sComm}_\mathrm{simplicial}{\mathrm{Modules}_\circledcirc}(\mathcal{O}_{X_{R_k,-,q}})\ar[d]\ar[d]\ar[d]\ar[d]\ar[r]^{\mathrm{global}}\ar[r]\ar[r] &\mathrm{sComm}_\mathrm{simplicial}\varphi_q{\mathrm{Modules}_\circledcirc}(\{\mathrm{Robba}^{\mathrm{extended},q}_{{R_k,-,I}}\}_I)\ar[d]\ar[d]\ar[d]\ar[d]\\
\mathrm{sComm}_\mathrm{simplicial}{\mathrm{Modules}_\circledcirc}(\mathcal{O}_{X_{\mathbb{Q}_p(p^{1/p^\infty})^{\wedge\flat},-,q}})\ar[r]^{\mathrm{global}}\ar[r]\ar[r] &\mathrm{sComm}_\mathrm{simplicial}\varphi_q{\mathrm{Modules}_\circledcirc}(\{\mathrm{Robba}^{\mathrm{extended},q}_{{R_0,-,I}}\}_I).  
}
\]

\item Then as in \cite{LBV} we have a functor (global section) of the de Rham complex after \cite[Definition 5.9, Section 5.2.1]{KKM}\footnote{Here $\circledcirc=\text{solidquasicoherentsheaves}$.}:
\[\displayindent=-0.4in
\xymatrix@R+6pc@C+0pc{
\mathrm{deRham}_{\mathrm{sComm}_\mathrm{simplicial}{\mathrm{Modules}_\circledcirc}(\mathcal{O}_{X_{R_k,-,q}})\ar[r]^{\mathrm{global}}}(-)\ar[d]\ar[d]\ar[d]\ar[d]\ar[r]\ar[r] &\mathrm{deRham}_{\mathrm{sComm}_\mathrm{simplicial}\varphi_q{\mathrm{Modules}_\circledcirc}(\{\mathrm{Robba}^{\mathrm{extended},q}_{{R_k,-,I}}\}_I)}(-)\ar[d]\ar[d]\ar[d]\ar[d]\\
\mathrm{deRham}_{\mathrm{sComm}_\mathrm{simplicial}{\mathrm{Modules}_\circledcirc}(\mathcal{O}_{X_{\mathbb{Q}_p(p^{1/p^\infty})^{\wedge\flat},-,q}})\ar[r]^{\mathrm{global}}}(-)\ar[r]\ar[r] &\mathrm{deRham}_{\mathrm{sComm}_\mathrm{simplicial}\varphi_q{\mathrm{Modules}_\circledcirc}(\{\mathrm{Robba}^{\mathrm{extended},q}_{{R_0,-,I}}\}_I)}(-). 
}
\]

\item Then we have the following a functor (global section) of $K$-group $(\infty,1)$-spectrum from \cite{BGT}\footnote{Here $\circledcirc=\text{solidquasicoherentsheaves}$.}:
\[
\xymatrix@R+6pc@C+0pc{
\mathrm{K}^\mathrm{BGT}_{\mathrm{sComm}_\mathrm{simplicial}{\mathrm{Modules}_\circledcirc}(\mathcal{O}_{X_{R_k,-,q}})\ar[r]^{\mathrm{global}}}(-)\ar[d]\ar[d]\ar[d]\ar[d]\ar[r]\ar[r] &\mathrm{K}^\mathrm{BGT}_{\mathrm{sComm}_\mathrm{simplicial}\varphi_q{\mathrm{Modules}_\circledcirc}(\{\mathrm{Robba}^{\mathrm{extended},q}_{{R_k,-,I}}\}_I)}(-)\ar[d]\ar[d]\ar[d]\ar[d]\\
\mathrm{K}^\mathrm{BGT}_{\mathrm{sComm}_\mathrm{simplicial}{\mathrm{Modules}_\circledcirc}(\mathcal{O}_{X_{\mathbb{Q}_p(p^{1/p^\infty})^{\wedge\flat},-,q}})\ar[r]^{\mathrm{global}}}(-)\ar[r]\ar[r] &\mathrm{K}^\mathrm{BGT}_{\mathrm{sComm}_\mathrm{simplicial}\varphi_q{\mathrm{Modules}_\circledcirc}(\{\mathrm{Robba}^{\mathrm{extended},q}_{{R_0,-,I}}\}_I)}(-). 
}
\]

\end{itemize}

\
\indent Then we consider further equivariance by considering the arithmetic profinite fundamental groups $\Gamma_{\mathbb{Q}_p}$ and $\mathrm{Gal}(\overline{{Q}_p\left<T_1^{\pm 1},...,T_k^{\pm 1}\right>}/R_k)$ through the following diagram:

\[
\xymatrix@R+0pc@C+0pc{
\mathbb{Z}_p^k=\mathrm{Gal}(R_k/{\mathbb{Q}_p(p^{1/p^\infty})^\wedge\left<T_1^{\pm 1},...,T_k^{\pm 1}\right>}) \ar[r]\ar[r] \ar[r]\ar[r] &\Gamma_k:=\mathrm{Gal}(R_k/{\mathbb{Q}_p\left<T_1^{\pm 1},...,T_k^{\pm 1}\right>}) \ar[r] \ar[r]\ar[r] &\Gamma_{\mathbb{Q}_p}.
}
\]

\begin{itemize}
\item (\text{Proposition}) There is a functor (global section) between the $\infty$-prestacks of inductive Banach quasicoherent sheaves\footnote{Here $\circledcirc=\text{solidquasicoherentsheaves}$.}:
\[
\xymatrix@R+6pc@C+0pc{
{\mathrm{Modules}_\circledcirc}_{\Gamma_k^{\times q}}(\mathcal{O}_{X_{\mathbb{Q}_p(p^{1/p^\infty})^{\wedge}\left<T_1^{\pm 1/p^\infty},...,T_k^{\pm 1/p^\infty}\right>^\flat,-,q}})\ar[d]\ar[d]\ar[d]\ar[d] \ar[r]^{\mathrm{global}}\ar[r]\ar[r] &\varphi_q{\mathrm{Modules}_\circledcirc}_{\Gamma_k^{\times q}}(\{\mathrm{Robba}^{\mathrm{extended},q}_{{R_k,-,I}}\}_I) \ar[d]\ar[d]\ar[d]\ar[d].\\
{\mathrm{Modules}_\circledcirc}(\mathcal{O}_{X_{\mathbb{Q}_p(p^{1/p^\infty})^{\wedge\flat},-,q}})\ar[r]^{\mathrm{global}}\ar[r]\ar[r] &\varphi_q{\mathrm{Modules}_\circledcirc}(\{\mathrm{Robba}^{\mathrm{extended},q}_{{R_0,-,I}}\}_I).\\ 
}
\]

\item (\text{Proposition}) There is a functor (global section) between the $\infty$-stacks of inductive Banach quasicoherent commutative algebra $E_\infty$ objects\footnote{Here $\circledcirc=\text{solidquasicoherentsheaves}$.}:
\[\displayindent=-0.4in
\xymatrix@R+6pc@C+0pc{
\mathrm{sComm}_\mathrm{simplicial}{\mathrm{Modules}_\circledcirc}_{\Gamma_k^{\times q}}(\mathcal{O}_{X_{R_k,-,q}})\ar[d]\ar[d]\ar[d]\ar[d]\ar[r]^{\mathrm{global}}\ar[r]\ar[r] &\mathrm{sComm}_\mathrm{simplicial}\varphi_q{\mathrm{Modules}_\circledcirc}_{\Gamma_k^{\times q}}(\{\mathrm{Robba}^{\mathrm{extended},q}_{{R_k,-,I}}\}_I)\ar[d]\ar[d]\ar[d]\ar[d]\\
\mathrm{sComm}_\mathrm{simplicial}{\mathrm{Modules}_\circledcirc}_{\Gamma_0^{\times q}}(\mathcal{O}_{X_{\mathbb{Q}_p(p^{1/p^\infty})^{\wedge\flat},-,q}})\ar[r]^{\mathrm{global}}\ar[r]\ar[r] &\mathrm{sComm}_\mathrm{simplicial}\varphi_q{\mathrm{Modules}_\circledcirc}_{\Gamma_0^{\times q}}(\{\mathrm{Robba}^{\mathrm{extended},q}_{{R_0,-,I}}\}_I).  
}
\]

\item Then as in \cite{LBV} we have a functor (global section) of the de Rham complex after \cite[Definition 5.9, Section 5.2.1]{KKM}\footnote{Here $\circledcirc=\text{solidquasicoherentsheaves}$.}:
\[\displayindent=-0.7in
\xymatrix@R+6pc@C+0pc{
\mathrm{deRham}_{\mathrm{sComm}_\mathrm{simplicial}{\mathrm{Modules}_\circledcirc}_{\Gamma_k^{\times q}}(\mathcal{O}_{X_{R_k,-,q}})\ar[r]^{\mathrm{global}}}(-)\ar[d]\ar[d]\ar[d]\ar[d]\ar[r]\ar[r] &\mathrm{deRham}_{\mathrm{sComm}_\mathrm{simplicial}\varphi_q{\mathrm{Modules}_\circledcirc}_{\Gamma_k^{\times q}}(\{\mathrm{Robba}^{\mathrm{extended},q}_{{R_k,-,I}}\}_I)}(-)\ar[d]\ar[d]\ar[d]\ar[d]\\
\mathrm{deRham}_{\mathrm{sComm}_\mathrm{simplicial}{\mathrm{Modules}_\circledcirc}_{\Gamma_0^{\times q}}(\mathcal{O}_{X_{\mathbb{Q}_p(p^{1/p^\infty})^{\wedge\flat},-,q}})\ar[r]^{\mathrm{global}}}(-)\ar[r]\ar[r] &\mathrm{deRham}_{\mathrm{sComm}_\mathrm{simplicial}\varphi_q{\mathrm{Modules}_\circledcirc}_{\Gamma_0^{\times q}}(\{\mathrm{Robba}^{\mathrm{extended},q}_{{R_0,-,I}}\}_I)}(-). 
}
\]

\item Then we have the following a functor (global section) of $K$-group $(\infty,1)$-spectrum from \cite{BGT}\footnote{Here $\circledcirc=\text{solidquasicoherentsheaves}$.}:
\[
\xymatrix@R+6pc@C+0pc{
\mathrm{K}^\mathrm{BGT}_{\mathrm{sComm}_\mathrm{simplicial}{\mathrm{Modules}_\circledcirc}_{\Gamma_k^{\times q}}(\mathcal{O}_{X_{R_k,-,q}})\ar[r]^{\mathrm{global}}}(-)\ar[d]\ar[d]\ar[d]\ar[d]\ar[r]\ar[r] &\mathrm{K}^\mathrm{BGT}_{\mathrm{sComm}_\mathrm{simplicial}\varphi_q{\mathrm{Modules}_\circledcirc}_{\Gamma_k^{\times q}}(\{\mathrm{Robba}^{\mathrm{extended},q}_{{R_k,-,I}}\}_I)}(-)\ar[d]\ar[d]\ar[d]\ar[d]\\
\mathrm{K}^\mathrm{BGT}_{\mathrm{sComm}_\mathrm{simplicial}{\mathrm{Modules}_\circledcirc}_{\Gamma_0^{\times q}}(\mathcal{O}_{X_{\mathbb{Q}_p(p^{1/p^\infty})^{\wedge\flat},-,q}})\ar[r]^{\mathrm{global}}}(-)\ar[r]\ar[r] &\mathrm{K}^\mathrm{BGT}_{\mathrm{sComm}_\mathrm{simplicial}\varphi_q{\mathrm{Modules}_\circledcirc}_{\Gamma_0^{\times q}}(\{\mathrm{Robba}^{\mathrm{extended},q}_{{R_0,-,I}}\}_I)}(-). 
}
\]

\end{itemize}

\begin{proposition}
All the global functors from \cite[Proposition 13.8, Theorem 14.9, Remark 14.10]{1CS2} achieve the equivalences on both sides.	
\end{proposition}

\newpage
\subsection{$\infty$-Categorical Analytic Stacks and Descents IV}

\indent In the following the right had of each row in each diagram will be the corresponding quasicoherent Robba bundles over the Robba ring carrying the corresponding action from the Frobenius or the fundamental groups, defined by directly applying \cite[Section 9.3]{KKM} and \cite{BBM}. We now let $\mathcal{A}$ be any commutative algebra objects in the corresponding $\infty$-toposes over ind-Banach commutative algebra objects over $\mathbb{Q}_p$ or the corresponding born\'e commutative algebra objects over $\mathbb{Q}_p$ carrying the Grothendieck topology defined by essentially the corresponding monomorphism homotopy in the opposite category. Then we promote the construction to the corresponding $\infty$-stack over the same $\infty$-categories of affinoids. We now take the corresponding colimit through all the $(\infty,1)$-categories. Therefore all the corresponding $(\infty,1)$-functors into $(\infty,1)$-categories or $(\infty,1)$-groupoids are from the homotopy closure of $\mathbb{Q}_p\left<C_1,...,C_\ell\right>$ $\ell=1,q,...$ in $\mathrm{sComm}\mathrm{Ind}\mathrm{Banach}_{\mathbb{Q}_p}$ or $\mathbb{Q}_p\left<C_1,...,C_\ell\right>$ $\ell=1,q,...$ in $\mathrm{sComm}\mathrm{Ind}^m\mathrm{Banach}_{\mathbb{Q}_p}$ as in \cite[Section 4.2]{BBM}:
\begin{align}
\mathrm{Ind}^{\mathbb{Q}_p\left<C_1,...,C_\ell\right>,\ell=1,q,...}\mathrm{sComm}\mathrm{Ind}\mathrm{Banach}_{\mathbb{Q}_p},\\
\mathrm{Ind}^{\mathbb{Q}_p\left<C_1,...,C_\ell\right>,\ell=1,q,...}\mathrm{sComm}\mathrm{Ind}\mathrm{Banach}_{\mathbb{Q}_p}	.
\end{align}

\begin{itemize}

\item (\text{Proposition}) There is a functor (global section) between the $\infty$-prestacks of inductive Banach quasicoherent presheaves:
\[
\xymatrix@R+0pc@C+0pc{
\mathrm{Ind}\mathrm{Banach}(\mathcal{O}_{X_{R,-,q}})\ar[r]^{\mathrm{global}}\ar[r]\ar[r] &\varphi_q\mathrm{Ind}\mathrm{Banach}(\{\mathrm{Robba}^\mathrm{extended}_{{R,-,I,q}}\}_I).  
}
\]
The definition is given by the following:
\[
\xymatrix@R+0pc@C+0pc{
\mathrm{homotopycolimit}_i(\mathrm{Ind}\mathrm{Banach}(\mathcal{O}_{X_{R,-,q}})\ar[r]^{\mathrm{global}}\ar[r]\ar[r] &\varphi_q\mathrm{Ind}\mathrm{Banach}(\{\mathrm{Robba}^\mathrm{extended}_{{R,-,I,q}}\}_I))(\mathcal{O}_i),  
}
\]
each $\mathcal{O}_i$ is just as $\mathbb{Q}_p\left<C_1,...,C_\ell\right>,\ell=1,q,...$.
\item (\text{Proposition}) There is a functor (global section) between the $\infty$-prestacks of monomorphic inductive Banach quasicoherent presheaves:
\[
\xymatrix@R+0pc@C+0pc{
\mathrm{Ind}^m\mathrm{Banach}(\mathcal{O}_{X_{R,-,q}})\ar[r]^{\mathrm{global}}\ar[r]\ar[r] &\varphi_q\mathrm{Ind}^m\mathrm{Banach}(\{\mathrm{Robba}^\mathrm{extended}_{{R,-,I,q}}\}_I).  
}
\]
The definition is given by the following:
\[
\xymatrix@R+0pc@C+0pc{
\mathrm{homotopycolimit}_i(\mathrm{Ind}^m\mathrm{Banach}(\mathcal{O}_{X_{R,-,q}})\ar[r]^{\mathrm{global}}\ar[r]\ar[r] &\varphi_q\mathrm{Ind}^m\mathrm{Banach}(\{\mathrm{Robba}^\mathrm{extended}_{{R,-,I,q}}\}_I))(\mathcal{O}_i),  
}
\]
each $\mathcal{O}_i$ is just as $\mathbb{Q}_p\left<C_1,...,C_\ell\right>,\ell=1,q,...$.

\item (\text{Proposition}) There is a functor (global section) between the $\infty$-prestacks of inductive Banach quasicoherent presheaves:
\[
\xymatrix@R+0pc@C+0pc{
\mathrm{Ind}\mathrm{Banach}(\mathcal{O}_{X_{R,-,q}})\ar[r]^{\mathrm{global}}\ar[r]\ar[r] &\varphi_q\mathrm{Ind}\mathrm{Banach}(\{\mathrm{Robba}^\mathrm{extended}_{{R,-,I,q}}\}_I).  
}
\]
The definition is given by the following:
\[
\xymatrix@R+0pc@C+0pc{
\mathrm{homotopycolimit}_i(\mathrm{Ind}\mathrm{Banach}(\mathcal{O}_{X_{R,-,q}})\ar[r]^{\mathrm{global}}\ar[r]\ar[r] &\varphi_q\mathrm{Ind}\mathrm{Banach}(\{\mathrm{Robba}^\mathrm{extended}_{{R,-,I,q}}\}_I))(\mathcal{O}_i),  
}
\]
each $\mathcal{O}_i$ is just as $\mathbb{Q}_p\left<C_1,...,C_\ell\right>,\ell=1,q,...$.
\item (\text{Proposition}) There is a functor (global section) between the $\infty$-stacks of monomorphic inductive Banach quasicoherent presheaves:
\[
\xymatrix@R+0pc@C+0pc{
\mathrm{Ind}^m\mathrm{Banach}(\mathcal{O}_{X_{R,-,q}})\ar[r]^{\mathrm{global}}\ar[r]\ar[r] &\varphi_q\mathrm{Ind}^m\mathrm{Banach}(\{\mathrm{Robba}^\mathrm{extended}_{{R,-,I,q}}\}_I).  
}
\]
The definition is given by the following:
\[
\xymatrix@R+0pc@C+0pc{
\mathrm{homotopycolimit}_i(\mathrm{Ind}^m\mathrm{Banach}(\mathcal{O}_{X_{R,-,q}})\ar[r]^{\mathrm{global}}\ar[r]\ar[r] &\varphi_q\mathrm{Ind}^m\mathrm{Banach}(\{\mathrm{Robba}^\mathrm{extended}_{{R,-,I,q}}\}_I))(\mathcal{O}_i),  
}
\]
each $\mathcal{O}_i$ is just as $\mathbb{Q}_p\left<C_1,...,C_\ell\right>,\ell=1,q,...$.
\item (\text{Proposition}) There is a functor (global section) between the $\infty$-prestacks of inductive Banach quasicoherent commutative algebra $E_\infty$ objects:
\[
\xymatrix@R+0pc@C+0pc{
\mathrm{sComm}_\mathrm{simplicial}\mathrm{Ind}\mathrm{Banach}(\mathcal{O}_{X_{R,-,q}})\ar[r]^{\mathrm{global}}\ar[r]\ar[r] &\mathrm{sComm}_\mathrm{simplicial}\varphi_q\mathrm{Ind}\mathrm{Banach}(\{\mathrm{Robba}^\mathrm{extended}_{{R,-,I,q}}\}_I).  
}
\]
The definition is given by the following:
\[
\xymatrix@R+0pc@C+0pc{
\mathrm{homotopycolimit}_i(\mathrm{sComm}_\mathrm{simplicial}\mathrm{Ind}\mathrm{Banach}(\mathcal{O}_{X_{R,-,q}})\ar[r]^{\mathrm{global}}\ar[r]\ar[r] &\mathrm{sComm}_\mathrm{simplicial}\varphi_q\mathrm{Ind}\mathrm{Banach}(\{\mathrm{Robba}^\mathrm{extended}_{{R,-,I,q}}\}_I))(\mathcal{O}_i),  
}
\]
each $\mathcal{O}_i$ is just as $\mathbb{Q}_p\left<C_1,...,C_\ell\right>,\ell=1,q,...$.
\item (\text{Proposition}) There is a functor (global section) between the $\infty$-prestacks of monomorphic inductive Banach quasicoherent commutative algebra $E_\infty$ objects:
\[
\xymatrix@R+0pc@C+0pc{
\mathrm{sComm}_\mathrm{simplicial}\mathrm{Ind}^m\mathrm{Banach}(\mathcal{O}_{X_{R,-,q}})\ar[r]^{\mathrm{global}}\ar[r]\ar[r] &\mathrm{sComm}_\mathrm{simplicial}\varphi_q\mathrm{Ind}^m\mathrm{Banach}(\{\mathrm{Robba}^\mathrm{extended}_{{R,-,I,q}}\}_I).  
}
\]
The definition is given by the following:
\[
\xymatrix@R+0pc@C+0pc{
\mathrm{homotopycolimit}_i(\mathrm{sComm}_\mathrm{simplicial}\mathrm{Ind}^m\mathrm{Banach}(\mathcal{O}_{X_{R,-,q}})\ar[r]^{\mathrm{global}}\ar[r]\ar[r] &\mathrm{sComm}_\mathrm{simplicial}\varphi_q\mathrm{Ind}^m\mathrm{Banach}(\{\mathrm{Robba}^\mathrm{extended}_{{R,-,I,q}}\}_I))(\mathcal{O}_i),  
}
\]
each $\mathcal{O}_i$ is just as $\mathbb{Q}_p\left<C_1,...,C_\ell\right>,\ell=1,q,...$.

\item Then parallel as in \cite{LBV} we have a functor (global section ) of the de Rham complex after \cite[Definition 5.9, Section 5.2.1]{KKM}:
\[
\xymatrix@R+0pc@C+0pc{
\mathrm{deRham}_{\mathrm{sComm}_\mathrm{simplicial}\mathrm{Ind}\mathrm{Banach}(\mathcal{O}_{X_{R,-,q}})\ar[r]^{\mathrm{global}}}(-)\ar[r]\ar[r] &\mathrm{deRham}_{\mathrm{sComm}_\mathrm{simplicial}\varphi_q\mathrm{Ind}\mathrm{Banach}(\{\mathrm{Robba}^\mathrm{extended}_{{R,-,I,q}}\}_I)}(-), 
}
\]
\[
\xymatrix@R+0pc@C+0pc{
\mathrm{deRham}_{\mathrm{sComm}_\mathrm{simplicial}\mathrm{Ind}^m\mathrm{Banach}(\mathcal{O}_{X_{R,-,q}})\ar[r]^{\mathrm{global}}}(-)\ar[r]\ar[r] &\mathrm{deRham}_{\mathrm{sComm}_\mathrm{simplicial}\varphi_q\mathrm{Ind}^m\mathrm{Banach}(\{\mathrm{Robba}^\mathrm{extended}_{{R,-,I,q}}\}_I)}(-). 
}
\]
The definition is given by the following:
\[\displayindent=-0.1in
\xymatrix@R+0pc@C+0pc{
\mathrm{homotopycolimit}_i\\
(\mathrm{deRham}_{\mathrm{sComm}_\mathrm{simplicial}\mathrm{Ind}\mathrm{Banach}(\mathcal{O}_{X_{R,-,q}})\ar[r]^{\mathrm{global}}}(-)\ar[r]\ar[r] &\mathrm{deRham}_{\mathrm{sComm}_\mathrm{simplicial}\varphi_q\mathrm{Ind}\mathrm{Banach}(\{\mathrm{Robba}^\mathrm{extended}_{{R,-,I,q}}\}_I)}(-))(\mathcal{O}_i),  
}
\]
\[\displayindent=-0.2in
\xymatrix@R+0pc@C+0pc{
\mathrm{homotopycolimit}_i\\
(\mathrm{deRham}_{\mathrm{sComm}_\mathrm{simplicial}\mathrm{Ind}^m\mathrm{Banach}(\mathcal{O}_{X_{R,-,q}})\ar[r]^{\mathrm{global}}}(-)\ar[r]\ar[r] &\mathrm{deRham}_{\mathrm{sComm}_\mathrm{simplicial}\varphi_q\mathrm{Ind}^m\mathrm{Banach}(\{\mathrm{Robba}^\mathrm{extended}_{{R,-,I,q}}\}_I)}(-))(\mathcal{O}_i),  
}
\]
each $\mathcal{O}_i$ is just as $\mathbb{Q}_p\left<C_1,...,C_\ell\right>,\ell=1,q,...$.\item Then we have the following a functor (global section) of $K$-group $(\infty,1)$-spectrum from \cite{BGT}:
\[
\xymatrix@R+0pc@C+0pc{
\mathrm{K}^\mathrm{BGT}_{\mathrm{sComm}_\mathrm{simplicial}\mathrm{Ind}\mathrm{Banach}(\mathcal{O}_{X_{R,-,q}})\ar[r]^{\mathrm{global}}}(-)\ar[r]\ar[r] &\mathrm{K}^\mathrm{BGT}_{\mathrm{sComm}_\mathrm{simplicial}\varphi_q\mathrm{Ind}\mathrm{Banach}(\{\mathrm{Robba}^\mathrm{extended}_{{R,-,I,q}}\}_I)}(-), 
}
\]
\[
\xymatrix@R+0pc@C+0pc{
\mathrm{K}^\mathrm{BGT}_{\mathrm{sComm}_\mathrm{simplicial}\mathrm{Ind}^m\mathrm{Banach}(\mathcal{O}_{X_{R,-,q}})\ar[r]^{\mathrm{global}}}(-)\ar[r]\ar[r] &\mathrm{K}^\mathrm{BGT}_{\mathrm{sComm}_\mathrm{simplicial}\varphi_q\mathrm{Ind}^m\mathrm{Banach}(\{\mathrm{Robba}^\mathrm{extended}_{{R,-,I,q}}\}_I)}(-). 
}
\]
The definition is given by the following:
\[
\xymatrix@R+0pc@C+0pc{
\mathrm{homotopycolimit}_i(\mathrm{K}^\mathrm{BGT}_{\mathrm{sComm}_\mathrm{simplicial}\mathrm{Ind}\mathrm{Banach}(\mathcal{O}_{X_{R,-,q}})\ar[r]^{\mathrm{global}}}(-)\ar[r]\ar[r] &\mathrm{K}^\mathrm{BGT}_{\mathrm{sComm}_\mathrm{simplicial}\varphi_q\mathrm{Ind}\mathrm{Banach}(\{\mathrm{Robba}^\mathrm{extended}_{{R,-,I,q}}\}_I)}(-))(\mathcal{O}_i),  
}
\]
\[\displayindent=-0.5in
\xymatrix@R+0pc@C+0pc{
\mathrm{homotopycolimit}_i(\mathrm{K}^\mathrm{BGT}_{\mathrm{sComm}_\mathrm{simplicial}\mathrm{Ind}^m\mathrm{Banach}(\mathcal{O}_{X_{R,-,q}})\ar[r]^{\mathrm{global}}}(-)\ar[r]\ar[r] &\mathrm{K}^\mathrm{BGT}_{\mathrm{sComm}_\mathrm{simplicial}\varphi_q\mathrm{Ind}^m\mathrm{Banach}(\{\mathrm{Robba}^\mathrm{extended}_{{R,-,I,q}}\}_I)}(-))(\mathcal{O}_i),  
}
\]
each $\mathcal{O}_i$ is just as $\mathbb{Q}_p\left<C_1,...,C_\ell\right>,\ell=1,q,...$.
\end{itemize}

\noindent Now let $R=\mathbb{Q}_p(p^{1/p^\infty})^{\wedge\flat}$ and $R_k=\mathbb{Q}_p(p^{1/p^\infty})^{\wedge}\left<T_1^{\pm 1/p^{\infty}},...,T_k^{\pm 1/p^{\infty}}\right>^\flat$ we have the following Galois theoretic results with naturality along $f:\mathrm{Spa}(\mathbb{Q}_p(p^{1/p^\infty})^{\wedge}\left<T_1^{\pm 1/p^\infty},...,T_k^{\pm 1/p^\infty}\right>^\flat)\rightarrow \mathrm{Spa}(\mathbb{Q}_p(p^{1/p^\infty})^{\wedge\flat})$:

\begin{itemize}
\item (\text{Proposition}) There is a functor (global section) between the $\infty$-prestacks of inductive Banach quasicoherent presheaves:
\[
\xymatrix@R+6pc@C+0pc{
\mathrm{Ind}\mathrm{Banach}(\mathcal{O}_{X_{\mathbb{Q}_p(p^{1/p^\infty})^{\wedge}\left<T_1^{\pm 1/p^\infty},...,T_k^{\pm 1/p^\infty}\right>^\flat,-,q}})\ar[d]\ar[d]\ar[d]\ar[d] \ar[r]^{\mathrm{global}}\ar[r]\ar[r] &\varphi_q\mathrm{Ind}\mathrm{Banach}(\{\mathrm{Robba}^\mathrm{extended}_{{R_k,-,I,q}}\}_I) \ar[d]\ar[d]\ar[d]\ar[d].\\
\mathrm{Ind}\mathrm{Banach}(\mathcal{O}_{X_{\mathbb{Q}_p(p^{1/p^\infty})^{\wedge\flat},-,q}})\ar[r]^{\mathrm{global}}\ar[r]\ar[r] &\varphi_q\mathrm{Ind}\mathrm{Banach}(\{\mathrm{Robba}^\mathrm{extended}_{{R_0,-,I,q}}\}_I).\\ 
}
\]
The definition is given by the following:
\[
\xymatrix@R+0pc@C+0pc{
\mathrm{homotopycolimit}_i(\square)(\mathcal{O}_i),  
}
\]
each $\mathcal{O}_i$ is just as $\mathbb{Q}_p\left<C_1,...,C_\ell\right>,\ell=1,q,...$ and $\square$ is the relative diagram of $\infty$-functors.
\item (\text{Proposition}) There is a functor (global section) between the $\infty$-prestacks of monomorphic inductive Banach quasicoherent presheaves:
\[
\xymatrix@R+6pc@C+0pc{
\mathrm{Ind}^m\mathrm{Banach}(\mathcal{O}_{X_{R_k,-,q}})\ar[r]^{\mathrm{global}}\ar[d]\ar[d]\ar[d]\ar[d]\ar[r]\ar[r] &\varphi_q\mathrm{Ind}^m\mathrm{Banach}(\{\mathrm{Robba}^\mathrm{extended}_{{R_k,-,I,q}}\}_I)\ar[d]\ar[d]\ar[d]\ar[d]\\
\mathrm{Ind}^m\mathrm{Banach}(\mathcal{O}_{X_{\mathbb{Q}_p(p^{1/p^\infty})^{\wedge\flat},-,q}})\ar[r]^{\mathrm{global}}\ar[r]\ar[r] &\varphi_q\mathrm{Ind}^m\mathrm{Banach}(\{\mathrm{Robba}^\mathrm{extended}_{{R_0,-,I,q}}\}_I).\\  
}
\]
The definition is given by the following:
\[
\xymatrix@R+0pc@C+0pc{
\mathrm{homotopycolimit}_i(\square)(\mathcal{O}_i),  
}
\]
each $\mathcal{O}_i$ is just as $\mathbb{Q}_p\left<C_1,...,C_\ell\right>,\ell=1,q,...$ and $\square$ is the relative diagram of $\infty$-functors.

\item (\text{Proposition}) There is a functor (global section) between the $\infty$-prestacks of inductive Banach quasicoherent commutative algebra $E_\infty$ objects:
\[
\xymatrix@R+6pc@C+0pc{
\mathrm{sComm}_\mathrm{simplicial}\mathrm{Ind}\mathrm{Banach}(\mathcal{O}_{X_{R_k,-,q}})\ar[d]\ar[d]\ar[d]\ar[d]\ar[r]^{\mathrm{global}}\ar[r]\ar[r] &\mathrm{sComm}_\mathrm{simplicial}\varphi_q\mathrm{Ind}\mathrm{Banach}(\{\mathrm{Robba}^\mathrm{extended}_{{R_k,-,I,q}}\}_I)\ar[d]\ar[d]\ar[d]\ar[d]\\
\mathrm{sComm}_\mathrm{simplicial}\mathrm{Ind}\mathrm{Banach}(\mathcal{O}_{X_{\mathbb{Q}_p(p^{1/p^\infty})^{\wedge\flat},-,q}})\ar[r]^{\mathrm{global}}\ar[r]\ar[r] &\mathrm{sComm}_\mathrm{simplicial}\varphi_q\mathrm{Ind}\mathrm{Banach}(\{\mathrm{Robba}^\mathrm{extended}_{{R_0,-,I,q}}\}_I).  
}
\]
The definition is given by the following:
\[
\xymatrix@R+0pc@C+0pc{
\mathrm{homotopycolimit}_i(\square)(\mathcal{O}_i),  
}
\]
each $\mathcal{O}_i$ is just as $\mathbb{Q}_p\left<C_1,...,C_\ell\right>,\ell=1,q,...$ and $\square$ is the relative diagram of $\infty$-functors.

\item (\text{Proposition}) There is a functor (global section) between the $\infty$-prestacks of monomorphic inductive Banach quasicoherent commutative algebra $E_\infty$ objects:
\[\displayindent=-0.4in
\xymatrix@R+6pc@C+0pc{
\mathrm{sComm}_\mathrm{simplicial}\mathrm{Ind}^m\mathrm{Banach}(\mathcal{O}_{X_{R_k,-,q}})\ar[d]\ar[d]\ar[d]\ar[d]\ar[r]^{\mathrm{global}}\ar[r]\ar[r] &\mathrm{sComm}_\mathrm{simplicial}\varphi_q\mathrm{Ind}^m\mathrm{Banach}(\{\mathrm{Robba}^\mathrm{extended}_{{R_k,-,I,q}}\}_I)\ar[d]\ar[d]\ar[d]\ar[d]\\
 \mathrm{sComm}_\mathrm{simplicial}\mathrm{Ind}^m\mathrm{Banach}(\mathcal{O}_{X_{\mathbb{Q}_p(p^{1/p^\infty})^{\wedge\flat},-,q}})\ar[r]^{\mathrm{global}}\ar[r]\ar[r] &\mathrm{sComm}_\mathrm{simplicial}\varphi_q\mathrm{Ind}^m\mathrm{Banach}(\{\mathrm{Robba}^\mathrm{extended}_{{R_0,-,I,q}}\}_I).
}
\]
The definition is given by the following:
\[
\xymatrix@R+0pc@C+0pc{
\mathrm{homotopycolimit}_i(\square)(\mathcal{O}_i),  
}
\]
each $\mathcal{O}_i$ is just as $\mathbb{Q}_p\left<C_1,...,C_\ell\right>,\ell=1,q,...$ and $\square$ is the relative diagram of $\infty$-functors.

\item Then parallel as in \cite{LBV} we have a functor (global section) of the de Rham complex after \cite[Definition 5.9, Section 5.2.1]{KKM}:
\[\displayindent=-0.4in
\xymatrix@R+6pc@C+0pc{
\mathrm{deRham}_{\mathrm{sComm}_\mathrm{simplicial}\mathrm{Ind}\mathrm{Banach}(\mathcal{O}_{X_{R_k,-,q}})\ar[r]^{\mathrm{global}}}(-)\ar[d]\ar[d]\ar[d]\ar[d]\ar[r]\ar[r] &\mathrm{deRham}_{\mathrm{sComm}_\mathrm{simplicial}\varphi_q\mathrm{Ind}\mathrm{Banach}(\{\mathrm{Robba}^\mathrm{extended}_{{R_k,-,I,q}}\}_I)}(-)\ar[d]\ar[d]\ar[d]\ar[d]\\
\mathrm{deRham}_{\mathrm{sComm}_\mathrm{simplicial}\mathrm{Ind}\mathrm{Banach}(\mathcal{O}_{X_{\mathbb{Q}_p(p^{1/p^\infty})^{\wedge\flat},-,q}})\ar[r]^{\mathrm{global}}}(-)\ar[r]\ar[r] &\mathrm{deRham}_{\mathrm{sComm}_\mathrm{simplicial}\varphi_q\mathrm{Ind}\mathrm{Banach}(\{\mathrm{Robba}^\mathrm{extended}_{{R_0,-,I,q}}\}_I)}(-), 
}
\]
\[\displayindent=-0.4in
\xymatrix@R+6pc@C+0pc{
\mathrm{deRham}_{\mathrm{sComm}_\mathrm{simplicial}\mathrm{Ind}^m\mathrm{Banach}(\mathcal{O}_{X_{R_k,-,q}})\ar[r]^{\mathrm{global}}}(-)\ar[d]\ar[d]\ar[d]\ar[d]\ar[r]\ar[r] &\mathrm{deRham}_{\mathrm{sComm}_\mathrm{simplicial}\varphi_q\mathrm{Ind}^m\mathrm{Banach}(\{\mathrm{Robba}^\mathrm{extended}_{{R_k,-,I,q}}\}_I)}(-)\ar[d]\ar[d]\ar[d]\ar[d]\\
\mathrm{deRham}_{\mathrm{sComm}_\mathrm{simplicial}\mathrm{Ind}^m\mathrm{Banach}(\mathcal{O}_{X_{\mathbb{Q}_p(p^{1/p^\infty})^{\wedge\flat},-,q}})\ar[r]^{\mathrm{global}}}(-)\ar[r]\ar[r] &\mathrm{deRham}_{\mathrm{sComm}_\mathrm{simplicial}\varphi_q\mathrm{Ind}^m\mathrm{Banach}(\{\mathrm{Robba}^\mathrm{extended}_{{R_0,-,I,q}}\}_I)}(-). 
}
\]

\item Then we have the following a functor (global section) of $K$-group $(\infty,1)$-spectrum from \cite{BGT}:
\[
\xymatrix@R+6pc@C+0pc{
\mathrm{K}^\mathrm{BGT}_{\mathrm{sComm}_\mathrm{simplicial}\mathrm{Ind}\mathrm{Banach}(\mathcal{O}_{X_{R_k,-,q}})\ar[r]^{\mathrm{global}}}(-)\ar[d]\ar[d]\ar[d]\ar[d]\ar[r]\ar[r] &\mathrm{K}^\mathrm{BGT}_{\mathrm{sComm}_\mathrm{simplicial}\varphi_q\mathrm{Ind}\mathrm{Banach}(\{\mathrm{Robba}^\mathrm{extended}_{{R_k,-,I,q}}\}_I)}(-)\ar[d]\ar[d]\ar[d]\ar[d]\\
\mathrm{K}^\mathrm{BGT}_{\mathrm{sComm}_\mathrm{simplicial}\mathrm{Ind}\mathrm{Banach}(\mathcal{O}_{X_{\mathbb{Q}_p(p^{1/p^\infty})^{\wedge\flat},-,q}})\ar[r]^{\mathrm{global}}}(-)\ar[r]\ar[r] &\mathrm{K}^\mathrm{BGT}_{\mathrm{sComm}_\mathrm{simplicial}\varphi_q\mathrm{Ind}\mathrm{Banach}(\{\mathrm{Robba}^\mathrm{extended}_{{R_0,-,I,q}}\}_I)}(-), 
}
\]
\[
\xymatrix@R+6pc@C+0pc{
\mathrm{K}^\mathrm{BGT}_{\mathrm{sComm}_\mathrm{simplicial}\mathrm{Ind}^m\mathrm{Banach}(\mathcal{O}_{X_{R_k,-,q}})\ar[r]^{\mathrm{global}}}(-)\ar[d]\ar[d]\ar[d]\ar[d]\ar[r]\ar[r] &\mathrm{K}^\mathrm{BGT}_{\mathrm{sComm}_\mathrm{simplicial}\varphi_q\mathrm{Ind}^m\mathrm{Banach}(\{\mathrm{Robba}^\mathrm{extended}_{{R_k,-,I,q}}\}_I)}(-)\ar[d]\ar[d]\ar[d]\ar[d]\\
\mathrm{K}^\mathrm{BGT}_{\mathrm{sComm}_\mathrm{simplicial}\mathrm{Ind}^m\mathrm{Banach}(\mathcal{O}_{X_{\mathbb{Q}_p(p^{1/p^\infty})^{\wedge\flat},-,q}})\ar[r]^{\mathrm{global}}}(-)\ar[r]\ar[r] &\mathrm{K}^\mathrm{BGT}_{\mathrm{sComm}_\mathrm{simplicial}\varphi_q\mathrm{Ind}^m\mathrm{Banach}(\{\mathrm{Robba}^\mathrm{extended}_{{R_0,-,I,q}}\}_I)}(-). 
}
\]
The definition is given by the following:
\[
\xymatrix@R+0pc@C+0pc{
\mathrm{homotopycolimit}_i(\square)(\mathcal{O}_i),  
}
\]
each $\mathcal{O}_i$ is just as $\mathbb{Q}_p\left<C_1,...,C_\ell\right>,\ell=1,q,...$ and $\square$ is the relative diagram of $\infty$-functors.

\end{itemize}

\
\indent Then we consider further equivariance by considering the arithmetic profinite fundamental groups and actually its $q$-th power $\mathrm{Gal}(\overline{{Q}_p\left<T_1^{\pm 1},...,T_k^{\pm 1}\right>}/R_k)^{\times q}$ through the following diagram:\\

\[
\xymatrix@R+6pc@C+0pc{
\mathbb{Z}_p^k=\mathrm{Gal}(R_k/{\mathbb{Q}_p(p^{1/p^\infty})^\wedge\left<T_1^{\pm 1},...,T_k^{\pm 1}\right>})\ar[d]\ar[d]\ar[d]\ar[d] \ar[r]\ar[r] \ar[r]\ar[r] &\mathrm{Gal}(\overline{{Q}_p\left<T_1^{\pm 1},...,T_k^{\pm 1}\right>}/R_k) \ar[d]\ar[d]\ar[d] \ar[r]\ar[r] &\Gamma_{\mathbb{Q}_p} \ar[d]\ar[d]\ar[d]\ar[d]\\
(\mathbb{Z}_p^k=\mathrm{Gal}(R_k/{\mathbb{Q}_p(p^{1/p^\infty})^\wedge\left<T_1^{\pm 1},...,T_k^{\pm 1}\right>}))^{\times q} \ar[r]\ar[r] \ar[r]\ar[r] &\Gamma_k^{\times q}:=\mathrm{Gal}(R_k/{\mathbb{Q}_p\left<T_1^{\pm 1},...,T_k^{\pm 1}\right>})^{\times q} \ar[r] \ar[r]\ar[r] &\Gamma_{\mathbb{Q}_p}^{\times q}.
}
\]

\

We then have the correspond arithmetic profinite fundamental groups equivariant versions:

\begin{itemize}
\item (\text{Proposition}) There is a functor (global section) between the $\infty$-prestacks of inductive Banach quasicoherent presheaves:
\[
\xymatrix@R+6pc@C+0pc{
\mathrm{Ind}\mathrm{Banach}_{\Gamma_{k}^{\times q}}(\mathcal{O}_{X_{\mathbb{Q}_p(p^{1/p^\infty})^{\wedge}\left<T_1^{\pm 1/p^\infty},...,T_k^{\pm 1/p^\infty}\right>^\flat,-,q}})\ar[d]\ar[d]\ar[d]\ar[d] \ar[r]^{\mathrm{global}}\ar[r]\ar[r] &\varphi_q\mathrm{Ind}\mathrm{Banach}_{\Gamma_{k}^{\times q}}(\{\mathrm{Robba}^\mathrm{extended}_{{R_k,-,I,q}}\}_I) \ar[d]\ar[d]\ar[d]\ar[d].\\
\mathrm{Ind}\mathrm{Banach}_{\Gamma_{0}^{\times q}}(\mathcal{O}_{X_{\mathbb{Q}_p(p^{1/p^\infty})^{\wedge\flat},-,q}})\ar[r]^{\mathrm{global}}\ar[r]\ar[r] &\varphi_q\mathrm{Ind}\mathrm{Banach}_{\Gamma_{0}^{\times q}}(\{\mathrm{Robba}^\mathrm{extended}_{{R_0,-,I,q}}\}_I).\\ 
}
\]
The definition is given by the following:
\[
\xymatrix@R+0pc@C+0pc{
\mathrm{homotopycolimit}_i(\square)(\mathcal{O}_i),  
}
\]
each $\mathcal{O}_i$ is just as $\mathbb{Q}_p\left<C_1,...,C_\ell\right>,\ell=1,q,...$ and $\square$ is the relative diagram of $\infty$-functors.

\item (\text{Proposition}) There is a functor (global section) between the $\infty$-prestacks of monomorphic inductive Banach quasicoherent presheaves:
\[
\xymatrix@R+6pc@C+0pc{
\mathrm{Ind}^m\mathrm{Banach}_{\Gamma_{k}^{\times q}}(\mathcal{O}_{X_{R_k,-,q}})\ar[r]^{\mathrm{global}}\ar[d]\ar[d]\ar[d]\ar[d]\ar[r]\ar[r] &\varphi_q\mathrm{Ind}^m\mathrm{Banach}_{\Gamma_{k}^{\times q}}(\{\mathrm{Robba}^\mathrm{extended}_{{R_k,-,I,q}}\}_I)\ar[d]\ar[d]\ar[d]\ar[d]\\
\mathrm{Ind}^m\mathrm{Banach}_{\Gamma_{0}^{\times q}}(\mathcal{O}_{X_{\mathbb{Q}_p(p^{1/p^\infty})^{\wedge\flat},-,q}})\ar[r]^{\mathrm{global}}\ar[r]\ar[r] &\varphi_q\mathrm{Ind}^m\mathrm{Banach}_{\Gamma_{0}^{\times q}}(\{\mathrm{Robba}^\mathrm{extended}_{{R_0,-,I,q}}\}_I).\\  
}
\]
The definition is given by the following:
\[
\xymatrix@R+0pc@C+0pc{
\mathrm{homotopycolimit}_i(\square)(\mathcal{O}_i),  
}
\]
each $\mathcal{O}_i$ is just as $\mathbb{Q}_p\left<C_1,...,C_\ell\right>,\ell=1,q,...$ and $\square$ is the relative diagram of $\infty$-functors.

\item (\text{Proposition}) There is a functor (global section) between the $\infty$-stacks of inductive Banach quasicoherent commutative algebra $E_\infty$ objects:
\[\displayindent=-0.4in
\xymatrix@R+6pc@C+0pc{
\mathrm{sComm}_\mathrm{simplicial}\mathrm{Ind}\mathrm{Banach}_{\Gamma_{k}^{\times q}}(\mathcal{O}_{X_{R_k,-,q}})\ar[d]\ar[d]\ar[d]\ar[d]\ar[r]^{\mathrm{global}}\ar[r]\ar[r] &\mathrm{sComm}_\mathrm{simplicial}\varphi_q\mathrm{Ind}\mathrm{Banach}_{\Gamma_{k}^{\times q}}(\{\mathrm{Robba}^\mathrm{extended}_{{R_k,-,I,q}}\}_I)\ar[d]\ar[d]\ar[d]\ar[d]\\
\mathrm{sComm}_\mathrm{simplicial}\mathrm{Ind}\mathrm{Banach}_{\Gamma_{0}^{\times q}}(\mathcal{O}_{X_{\mathbb{Q}_p(p^{1/p^\infty})^{\wedge\flat},-,q}})\ar[r]^{\mathrm{global}}\ar[r]\ar[r] &\mathrm{sComm}_\mathrm{simplicial}\varphi_q\mathrm{Ind}\mathrm{Banach}_{\Gamma_{0}^{\times q}}(\{\mathrm{Robba}^\mathrm{extended}_{{R_0,-,I,q}}\}_I).  
}
\]
The definition is given by the following:
\[
\xymatrix@R+0pc@C+0pc{
\mathrm{homotopycolimit}_i(\square)(\mathcal{O}_i),  
}
\]
each $\mathcal{O}_i$ is just as $\mathbb{Q}_p\left<C_1,...,C_\ell\right>,\ell=1,q,...$ and $\square$ is the relative diagram of $\infty$-functors.

\item (\text{Proposition}) There is a functor (global section) between the $\infty$-prestacks of monomorphic inductive Banach quasicoherent commutative algebra $E_\infty$ objects:
\[\displayindent=-0.4in
\xymatrix@R+6pc@C+0pc{
\mathrm{sComm}_\mathrm{simplicial}\mathrm{Ind}^m\mathrm{Banach}_{\Gamma_{k}^{\times q}}(\mathcal{O}_{X_{R_k,-,q}})\ar[d]\ar[d]\ar[d]\ar[d]\ar[r]^{\mathrm{global}}\ar[r]\ar[r] &\mathrm{sComm}_\mathrm{simplicial}\varphi_q\mathrm{Ind}^m\mathrm{Banach}_{\Gamma_{k}^{\times q}}(\{\mathrm{Robba}^\mathrm{extended}_{{R_k,-,I,q}}\}_I)\ar[d]\ar[d]\ar[d]\ar[d]\\
 \mathrm{sComm}_\mathrm{simplicial}\mathrm{Ind}^m\mathrm{Banach}_{\Gamma_{0}^{\times q}}(\mathcal{O}_{X_{\mathbb{Q}_p(p^{1/p^\infty})^{\wedge\flat},-,q}})\ar[r]^{\mathrm{global}}\ar[r]\ar[r] &\mathrm{sComm}_\mathrm{simplicial}\varphi_q\mathrm{Ind}^m\mathrm{Banach}_{\Gamma_{0}^{\times q}}(\{\mathrm{Robba}^\mathrm{extended}_{{R_0,-,I,q}}\}_I). 
}
\]
The definition is given by the following:
\[
\xymatrix@R+0pc@C+0pc{
\mathrm{homotopycolimit}_i(\square)(\mathcal{O}_i),  
}
\]
each $\mathcal{O}_i$ is just as $\mathbb{Q}_p\left<C_1,...,C_\ell\right>,\ell=1,q,...$ and $\square$ is the relative diagram of $\infty$-functors.

\item Then parallel as in \cite{LBV} we have a functor (global section) of the de Rham complex after \cite[Definition 5.9, Section 5.2.1]{KKM}:
\[\displayindent=-0.6in
\xymatrix@R+6pc@C+0pc{
\mathrm{deRham}_{\mathrm{sComm}_\mathrm{simplicial}\mathrm{Ind}\mathrm{Banach}_{\Gamma_{k}^{\times q}}(\mathcal{O}_{X_{R_k,-,q}})\ar[r]^{\mathrm{global}}}(-)\ar[d]\ar[d]\ar[d]\ar[d]\ar[r]\ar[r] &\mathrm{deRham}_{\mathrm{sComm}_\mathrm{simplicial}\varphi_q\mathrm{Ind}\mathrm{Banach}_{\Gamma_{k}^{\times q}}(\{\mathrm{Robba}^\mathrm{extended}_{{R_k,-,I,q}}\}_I)}(-)\ar[d]\ar[d]\ar[d]\ar[d]\\
\mathrm{deRham}_{\mathrm{sComm}_\mathrm{simplicial}\mathrm{Ind}\mathrm{Banach}_{\Gamma_{0}^{\times q}}(\mathcal{O}_{X_{\mathbb{Q}_p(p^{1/p^\infty})^{\wedge\flat},-,q}})\ar[r]^{\mathrm{global}}}(-)\ar[r]\ar[r] &\mathrm{deRham}_{\mathrm{sComm}_\mathrm{simplicial}\varphi_q\mathrm{Ind}\mathrm{Banach}_{\Gamma_{0}^{\times q}}(\{\mathrm{Robba}^\mathrm{extended}_{{R_0,-,I,q}}\}_I)}(-), 
}
\]
\[\displayindent=-0.7in
\xymatrix@R+6pc@C+0pc{
\mathrm{deRham}_{\mathrm{sComm}_\mathrm{simplicial}\mathrm{Ind}^m\mathrm{Banach}_{\Gamma_{k}^{\times q}}(\mathcal{O}_{X_{R_k,-,q}})\ar[r]^{\mathrm{global}}}(-)\ar[d]\ar[d]\ar[d]\ar[d]\ar[r]\ar[r] &\mathrm{deRham}_{\mathrm{sComm}_\mathrm{simplicial}\varphi_q\mathrm{Ind}^m\mathrm{Banach}_{\Gamma_{k}^{\times q}}(\{\mathrm{Robba}^\mathrm{extended}_{{R_k,-,I,q}}\}_I)}(-)\ar[d]\ar[d]\ar[d]\ar[d]\\
\mathrm{deRham}_{\mathrm{sComm}_\mathrm{simplicial}\mathrm{Ind}^m\mathrm{Banach}_{\Gamma_{0}^{\times q}}(\mathcal{O}_{X_{\mathbb{Q}_p(p^{1/p^\infty})^{\wedge\flat},-,q}})\ar[r]^{\mathrm{global}}}(-)\ar[r]\ar[r] &\mathrm{deRham}_{\mathrm{sComm}_\mathrm{simplicial}\varphi_q\mathrm{Ind}^m\mathrm{Banach}_{\Gamma_{0}^{\times q}}(\{\mathrm{Robba}^\mathrm{extended}_{{R_0,-,I,q}}\}_I)}(-). 
}
\]
The definition is given by the following:
\[
\xymatrix@R+0pc@C+0pc{
\mathrm{homotopycolimit}_i(\square)(\mathcal{O}_i),  
}
\]
each $\mathcal{O}_i$ is just as $\mathbb{Q}_p\left<C_1,...,C_\ell\right>,\ell=1,q,...$ and $\square$ is the relative diagram of $\infty$-functors.

\item Then we have the following a functor (global section) of $K$-group $(\infty,1)$-spectrum from \cite{BGT}:
\[
\xymatrix@R+6pc@C+0pc{
\mathrm{K}^\mathrm{BGT}_{\mathrm{sComm}_\mathrm{simplicial}\mathrm{Ind}\mathrm{Banach}_{\Gamma_{k}^{\times q}}(\mathcal{O}_{X_{R_k,-,q}})\ar[r]^{\mathrm{global}}}(-)\ar[d]\ar[d]\ar[d]\ar[d]\ar[r]\ar[r] &\mathrm{K}^\mathrm{BGT}_{\mathrm{sComm}_\mathrm{simplicial}\varphi_q\mathrm{Ind}\mathrm{Banach}_{\Gamma_{k}^{\times q}}(\{\mathrm{Robba}^\mathrm{extended}_{{R_k,-,I,q}}\}_I)}(-)\ar[d]\ar[d]\ar[d]\ar[d]\\
\mathrm{K}^\mathrm{BGT}_{\mathrm{sComm}_\mathrm{simplicial}\mathrm{Ind}\mathrm{Banach}_{\Gamma_{0}^{\times q}}(\mathcal{O}_{X_{\mathbb{Q}_p(p^{1/p^\infty})^{\wedge\flat},-,q}})\ar[r]^{\mathrm{global}}}(-)\ar[r]\ar[r] &\mathrm{K}^\mathrm{BGT}_{\mathrm{sComm}_\mathrm{simplicial}\varphi_q\mathrm{Ind}\mathrm{Banach}_{\Gamma_{0}^{\times q}}(\{\mathrm{Robba}^\mathrm{extended}_{{R_0,-,I,q}}\}_I)}(-), 
}
\]
\[
\xymatrix@R+6pc@C+0pc{
\mathrm{K}^\mathrm{BGT}_{\mathrm{sComm}_\mathrm{simplicial}\mathrm{Ind}^m\mathrm{Banach}_{\Gamma_{k}^{\times q}}(\mathcal{O}_{X_{R_k,-,q}})\ar[r]^{\mathrm{global}}}(-)\ar[d]\ar[d]\ar[d]\ar[d]\ar[r]\ar[r] &\mathrm{K}^\mathrm{BGT}_{\mathrm{sComm}_\mathrm{simplicial}\varphi_q\mathrm{Ind}^m\mathrm{Banach}_{\Gamma_{k}^{\times q}}(\{\mathrm{Robba}^\mathrm{extended}_{{R_k,-,I,q}}\}_I)}(-)\ar[d]\ar[d]\ar[d]\ar[d]\\
\mathrm{K}^\mathrm{BGT}_{\mathrm{sComm}_\mathrm{simplicial}\mathrm{Ind}^m\mathrm{Banach}_{\Gamma_{0}^{\times q}}(\mathcal{O}_{X_{\mathbb{Q}_p(p^{1/p^\infty})^{\wedge\flat},-,q}})\ar[r]^{\mathrm{global}}}(-)\ar[r]\ar[r] &\mathrm{K}^\mathrm{BGT}_{\mathrm{sComm}_\mathrm{simplicial}\varphi_q\mathrm{Ind}^m\mathrm{Banach}_{\Gamma_{0}^{\times q}}(\{\mathrm{Robba}^\mathrm{extended}_{{R_0,-,I,q}}\}_I)}(-). 
}
\]
The definition is given by the following:
\[
\xymatrix@R+0pc@C+0pc{
\mathrm{homotopycolimit}_i(\square)(\mathcal{O}_i),  
}
\]
each $\mathcal{O}_i$ is just as $\mathbb{Q}_p\left<C_1,...,C_\ell\right>,\ell=1,q,...$ and $\square$ is the relative diagram of $\infty$-functors.

\end{itemize}

\

\begin{remark}
\noindent We can certainly consider the quasicoherent sheaves in \cite[Lemma 7.11, Remark 7.12]{1BBK}, therefore all the quasicoherent presheaves and modules will be those satisfying \cite[Lemma 7.11, Remark 7.12]{1BBK} if one would like to consider the the quasicoherent sheaves. That being all as this said, we would believe that the big quasicoherent presheaves are automatically quasicoherent sheaves (namely satisfying the corresponding \v{C}ech $\infty$-descent as in \cite[Section 9.3]{KKM} and \cite[Lemma 7.11, Remark 7.12]{1BBK}) and the corresponding global section functors are automatically equivalence of $\infty$-categories. \\
\end{remark}

\

\indent In Clausen-Scholze formalism we have the following\footnote{Certainly the homotopy colimit in the rings side will be within the condensed solid animated analytic rings from \cite{1CS2}.}:

\begin{itemize}
\item (\text{Proposition}) There is a functor (global section) between the $\infty$-prestacks of inductive Banach quasicoherent sheaves:
\[
\xymatrix@R+0pc@C+0pc{
{\mathrm{Modules}_\circledcirc}(\mathcal{O}_{X_{R,-,q}})\ar[r]^{\mathrm{global}}\ar[r]\ar[r] &\varphi_q{\mathrm{Modules}_\circledcirc}(\{\mathrm{Robba}^{\mathrm{extended},q}_{{R,-,I}}\}_I).  
}
\]
The definition is given by the following:
\[
\xymatrix@R+0pc@C+0pc{
\mathrm{homotopycolimit}_i({\mathrm{Modules}_\circledcirc}(\mathcal{O}_{X_{R,-,q}})\ar[r]^{\mathrm{global}}\ar[r]\ar[r] &\varphi_q{\mathrm{Modules}_\circledcirc}(\{\mathrm{Robba}^{\mathrm{extended},q}_{{R,-,I}}\}_I))(\mathcal{O}_i),  
}
\]
each $\mathcal{O}_i$ is just as $\mathbb{Q}_p\left<C_1,...,C_\ell\right>,\ell=1,2,...$.

\item (\text{Proposition}) There is a functor (global section) between the $\infty$-prestacks of inductive Banach quasicoherent sheaves:
\[
\xymatrix@R+0pc@C+0pc{
{\mathrm{Modules}_\circledcirc}(\mathcal{O}_{X_{R,-,q}})\ar[r]^{\mathrm{global}}\ar[r]\ar[r] &\varphi_q{\mathrm{Modules}_\circledcirc}(\{\mathrm{Robba}^{\mathrm{extended},q}_{{R,-,I}}\}_I).  
}
\]
The definition is given by the following:
\[
\xymatrix@R+0pc@C+0pc{
\mathrm{homotopycolimit}_i({\mathrm{Modules}_\circledcirc}(\mathcal{O}_{X_{R,-,q}})\ar[r]^{\mathrm{global}}\ar[r]\ar[r] &\varphi_q{\mathrm{Modules}_\circledcirc}(\{\mathrm{Robba}^{\mathrm{extended},q}_{{R,-,I}}\}_I))(\mathcal{O}_i),  
}
\]
each $\mathcal{O}_i$ is just as $\mathbb{Q}_p\left<C_1,...,C_\ell\right>,\ell=1,2,...$.

\item (\text{Proposition}) There is a functor (global section) between the $\infty$-prestacks of inductive Banach quasicoherent commutative algebra $E_\infty$ objects\footnote{Here $\circledcirc=\text{solidquasicoherentsheaves}$.}:
\[
\xymatrix@R+0pc@C+0pc{
\mathrm{sComm}_\mathrm{simplicial}{\mathrm{Modules}_\circledcirc}(\mathcal{O}_{X_{R,-,q}})\ar[r]^{\mathrm{global}}\ar[r]\ar[r] &\mathrm{sComm}_\mathrm{simplicial}\varphi_q{\mathrm{Modules}_\circledcirc}(\{\mathrm{Robba}^{\mathrm{extended},q}_{{R,-,I}}\}_I).  
}
\]
The definition is given by the following:
\[
\xymatrix@R+0pc@C+0pc{
\mathrm{homotopycolimit}_i(\mathrm{sComm}_\mathrm{simplicial}{\mathrm{Modules}_\circledcirc}(\mathcal{O}_{X_{R,-,q}})\ar[r]^{\mathrm{global}}\ar[r]\ar[r] &\mathrm{sComm}_\mathrm{simplicial}\varphi_q{\mathrm{Modules}_\circledcirc}(\{\mathrm{Robba}^{\mathrm{extended},q}_{{R,-,I}}\}_I))(\mathcal{O}_i),  
}
\]
each $\mathcal{O}_i$ is just as $\mathbb{Q}_p\left<C_1,...,C_\ell\right>,\ell=1,2,...$.
\item Then as in \cite{LBV} we have a functor (global section ) of the de Rham complex after \cite[Definition 5.9, Section 5.2.1]{KKM}\footnote{Here $\circledcirc=\text{solidquasicoherentsheaves}$.}:
\[
\xymatrix@R+0pc@C+0pc{
\mathrm{deRham}_{\mathrm{sComm}_\mathrm{simplicial}{\mathrm{Modules}_\circledcirc}(\mathcal{O}_{X_{R,-,q}})\ar[r]^{\mathrm{global}}}(-)\ar[r]\ar[r] &\mathrm{deRham}_{\mathrm{sComm}_\mathrm{simplicial}\varphi_q{\mathrm{Modules}_\circledcirc}(\{\mathrm{Robba}^{\mathrm{extended},q}_{{R,-,I}}\}_I)}(-), 
}
\]
The definition is given by the following:
\[
\xymatrix@R+0pc@C+0pc{
\mathrm{homotopycolimit}_i\\
(\mathrm{deRham}_{\mathrm{sComm}_\mathrm{simplicial}{\mathrm{Modules}_\circledcirc}(\mathcal{O}_{X_{R,-,q}})\ar[r]^{\mathrm{global}}}(-)\ar[r]\ar[r] &\mathrm{deRham}_{\mathrm{sComm}_\mathrm{simplicial}\varphi_q{\mathrm{Modules}_\circledcirc}(\{\mathrm{Robba}^{\mathrm{extended},q}_{{R,-,I}}\}_I)}(-))(\mathcal{O}_i),  
}
\]
each $\mathcal{O}_i$ is just as $\mathbb{Q}_p\left<C_1,...,C_\ell\right>,\ell=1,2,...$.\item Then we have the following a functor (global section) of $K$-group $(\infty,1)$-spectrum from \cite{BGT}\footnote{Here $\circledcirc=\text{solidquasicoherentsheaves}$.}:
\[
\xymatrix@R+0pc@C+0pc{
\mathrm{K}^\mathrm{BGT}_{\mathrm{sComm}_\mathrm{simplicial}{\mathrm{Modules}_\circledcirc}(\mathcal{O}_{X_{R,-,q}})\ar[r]^{\mathrm{global}}}(-)\ar[r]\ar[r] &\mathrm{K}^\mathrm{BGT}_{\mathrm{sComm}_\mathrm{simplicial}\varphi_q{\mathrm{Modules}_\circledcirc}(\{\mathrm{Robba}^{\mathrm{extended},q}_{{R,-,I}}\}_I)}(-). 
}
\]
The definition is given by the following:
\[
\xymatrix@R+0pc@C+0pc{
\mathrm{homotopycolimit}_i(\mathrm{K}^\mathrm{BGT}_{\mathrm{sComm}_\mathrm{simplicial}{\mathrm{Modules}_\circledcirc}(\mathcal{O}_{X_{R,-,q}})\ar[r]^{\mathrm{global}}}(-)\ar[r]\ar[r] &\mathrm{K}^\mathrm{BGT}_{\mathrm{sComm}_\mathrm{simplicial}\varphi_q{\mathrm{Modules}_\circledcirc}(\{\mathrm{Robba}^{\mathrm{extended},q}_{{R,-,I}}\}_I)}(-))(\mathcal{O}_i),  
}
\]
each $\mathcal{O}_i$ is just as $\mathbb{Q}_p\left<C_1,...,C_\ell\right>,\ell=1,2,...$.
\end{itemize}

\noindent Now let $R=\mathbb{Q}_p(p^{1/p^\infty})^{\wedge\flat}$ and $R_k=\mathbb{Q}_p(p^{1/p^\infty})^{\wedge}\left<T_1^{\pm 1/p^{\infty}},...,T_k^{\pm 1/p^{\infty}}\right>^\flat$ we have the following Galois theoretic results with naturality along $f:\mathrm{Spa}(\mathbb{Q}_p(p^{1/p^\infty})^{\wedge}\left<T_1^{\pm 1/p^\infty},...,T_k^{\pm 1/p^\infty}\right>^\flat)\rightarrow \mathrm{Spa}(\mathbb{Q}_p(p^{1/p^\infty})^{\wedge\flat})$:

\begin{itemize}
\item (\text{Proposition}) There is a functor (global section) between the $\infty$-prestacks of inductive Banach quasicoherent sheaves\footnote{Here $\circledcirc=\text{solidquasicoherentsheaves}$.}:
\[
\xymatrix@R+6pc@C+0pc{
{\mathrm{Modules}_\circledcirc}(\mathcal{O}_{X_{\mathbb{Q}_p(p^{1/p^\infty})^{\wedge}\left<T_1^{\pm 1/p^\infty},...,T_k^{\pm 1/p^\infty}\right>^\flat,-,q}})\ar[d]\ar[d]\ar[d]\ar[d] \ar[r]^{\mathrm{global}}\ar[r]\ar[r] &\varphi_q{\mathrm{Modules}_\circledcirc}(\{\mathrm{Robba}^{\mathrm{extended},q}_{{R_k,-,I}}\}_I) \ar[d]\ar[d]\ar[d]\ar[d].\\
{\mathrm{Modules}_\circledcirc}(\mathcal{O}_{X_{\mathbb{Q}_p(p^{1/p^\infty})^{\wedge\flat},-,q}})\ar[r]^{\mathrm{global}}\ar[r]\ar[r] &\varphi_q{\mathrm{Modules}_\circledcirc}(\{\mathrm{Robba}^{\mathrm{extended},q}_{{R_0,-,I}}\}_I).\\ 
}
\]
The definition is given by the following:
\[
\xymatrix@R+0pc@C+0pc{
\mathrm{homotopycolimit}_i(\square)(\mathcal{O}_i),  
}
\]
each $\mathcal{O}_i$ is just as $\mathbb{Q}_p\left<C_1,...,C_\ell\right>,\ell=1,2,...$ and $\square$ is the relative diagram of $\infty$-functors.

\item (\text{Proposition}) There is a functor (global section) between the $\infty$-prestacks of inductive Banach quasicoherent commutative algebra $E_\infty$ objects\footnote{Here $\circledcirc=\text{solidquasicoherentsheaves}$.}:
\[
\xymatrix@R+6pc@C+0pc{
\mathrm{sComm}_\mathrm{simplicial}{\mathrm{Modules}_\circledcirc}(\mathcal{O}_{X_{R_k,-,q}})\ar[d]\ar[d]\ar[d]\ar[d]\ar[r]^{\mathrm{global}}\ar[r]\ar[r] &\mathrm{sComm}_\mathrm{simplicial}\varphi_q{\mathrm{Modules}_\circledcirc}(\{\mathrm{Robba}^{\mathrm{extended},q}_{{R_k,-,I}}\}_I)\ar[d]\ar[d]\ar[d]\ar[d]\\
\mathrm{sComm}_\mathrm{simplicial}{\mathrm{Modules}_\circledcirc}(\mathcal{O}_{X_{\mathbb{Q}_p(p^{1/p^\infty})^{\wedge\flat},-,q}})\ar[r]^{\mathrm{global}}\ar[r]\ar[r] &\mathrm{sComm}_\mathrm{simplicial}\varphi_q{\mathrm{Modules}_\circledcirc}(\{\mathrm{Robba}^{\mathrm{extended},q}_{{R_0,-,I}}\}_I).  
}
\]
The definition is given by the following:
\[
\xymatrix@R+0pc@C+0pc{
\mathrm{homotopycolimit}_i(\square)(\mathcal{O}_i),  
}
\]
each $\mathcal{O}_i$ is just as $\mathbb{Q}_p\left<C_1,...,C_\ell\right>,\ell=1,2,...$ and $\square$ is the relative diagram of $\infty$-functors.

\item Then as in \cite{LBV} we have a functor (global section) of the de Rham complex after \cite[Definition 5.9, Section 5.2.1]{KKM}\footnote{Here $\circledcirc=\text{solidquasicoherentsheaves}$.}:
\[\displayindent=-0.4in
\xymatrix@R+6pc@C+0pc{
\mathrm{deRham}_{\mathrm{sComm}_\mathrm{simplicial}{\mathrm{Modules}_\circledcirc}(\mathcal{O}_{X_{R_k,-,q}})\ar[r]^{\mathrm{global}}}(-)\ar[d]\ar[d]\ar[d]\ar[d]\ar[r]\ar[r] &\mathrm{deRham}_{\mathrm{sComm}_\mathrm{simplicial}\varphi_q{\mathrm{Modules}_\circledcirc}(\{\mathrm{Robba}^{\mathrm{extended},q}_{{R_k,-,I}}\}_I)}(-)\ar[d]\ar[d]\ar[d]\ar[d]\\
\mathrm{deRham}_{\mathrm{sComm}_\mathrm{simplicial}{\mathrm{Modules}_\circledcirc}(\mathcal{O}_{X_{\mathbb{Q}_p(p^{1/p^\infty})^{\wedge\flat},-,q}})\ar[r]^{\mathrm{global}}}(-)\ar[r]\ar[r] &\mathrm{deRham}_{\mathrm{sComm}_\mathrm{simplicial}\varphi_q{\mathrm{Modules}_\circledcirc}(\{\mathrm{Robba}^{\mathrm{extended},q}_{{R_0,-,I}}\}_I)}(-), 
}
\]

\item Then we have the following a functor (global section) of $K$-group $(\infty,1)$-spectrum from \cite{BGT}\footnote{Here $\circledcirc=\text{solidquasicoherentsheaves}$.}:
\[
\xymatrix@R+6pc@C+0pc{
\mathrm{K}^\mathrm{BGT}_{\mathrm{sComm}_\mathrm{simplicial}{\mathrm{Modules}_\circledcirc}(\mathcal{O}_{X_{R_k,-,q}})\ar[r]^{\mathrm{global}}}(-)\ar[d]\ar[d]\ar[d]\ar[d]\ar[r]\ar[r] &\mathrm{K}^\mathrm{BGT}_{\mathrm{sComm}_\mathrm{simplicial}\varphi_q{\mathrm{Modules}_\circledcirc}(\{\mathrm{Robba}^{\mathrm{extended},q}_{{R_k,-,I}}\}_I)}(-)\ar[d]\ar[d]\ar[d]\ar[d]\\
\mathrm{K}^\mathrm{BGT}_{\mathrm{sComm}_\mathrm{simplicial}{\mathrm{Modules}_\circledcirc}(\mathcal{O}_{X_{\mathbb{Q}_p(p^{1/p^\infty})^{\wedge\flat},-,q}})\ar[r]^{\mathrm{global}}}(-)\ar[r]\ar[r] &\mathrm{K}^\mathrm{BGT}_{\mathrm{sComm}_\mathrm{simplicial}\varphi_q{\mathrm{Modules}_\circledcirc}(\{\mathrm{Robba}^{\mathrm{extended},q}_{{R_0,-,I}}\}_I)}(-), 
}
\]

The definition is given by the following:
\[
\xymatrix@R+0pc@C+0pc{
\mathrm{homotopycolimit}_i(\square)(\mathcal{O}_i),  
}
\]
each $\mathcal{O}_i$ is just as $\mathbb{Q}_p\left<C_1,...,C_\ell\right>,\ell=1,2,...$ and $\square$ is the relative diagram of $\infty$-functors.

\end{itemize}

\
\indent Then we consider further equivariance by considering the arithmetic profinite fundamental groups $\Gamma_{\mathbb{Q}_p}$ and $\mathrm{Gal}(\overline{{Q}_p\left<T_1^{\pm 1},...,T_k^{\pm 1}\right>}/R_k)$ through the following diagram:

\[
\xymatrix@R+0pc@C+0pc{
\mathbb{Z}_p^k=\mathrm{Gal}(R_k/{\mathbb{Q}_p(p^{1/p^\infty})^\wedge\left<T_1^{\pm 1},...,T_k^{\pm 1}\right>}) \ar[r]\ar[r] \ar[r]\ar[r] &\Gamma_k^{\times q}:=\mathrm{Gal}(R_k/{\mathbb{Q}_p\left<T_1^{\pm 1},...,T_k^{\pm 1}\right>}) \ar[r] \ar[r]\ar[r] &\Gamma_{\mathbb{Q}_p}.
}
\]

\begin{itemize}
\item (\text{Proposition}) There is a functor (global section) between the $\infty$-prestacks of inductive Banach quasicoherent sheaves\footnote{Here $\circledcirc=\text{solidquasicoherentsheaves}$.}:
\[
\xymatrix@R+6pc@C+0pc{
{\mathrm{Modules}_\circledcirc}_{\Gamma_k^{\times q}}(\mathcal{O}_{X_{\mathbb{Q}_p(p^{1/p^\infty})^{\wedge}\left<T_1^{\pm 1/p^\infty},...,T_k^{\pm 1/p^\infty}\right>^\flat,-,q}})\ar[d]\ar[d]\ar[d]\ar[d] \ar[r]^{\mathrm{global}}\ar[r]\ar[r] &\varphi_q{\mathrm{Modules}_\circledcirc}_{\Gamma_k^{\times q}}(\{\mathrm{Robba}^{\mathrm{extended},q}_{{R_k,-,I}}\}_I) \ar[d]\ar[d]\ar[d]\ar[d].\\
{\mathrm{Modules}_\circledcirc}(\mathcal{O}_{X_{\mathbb{Q}_p(p^{1/p^\infty})^{\wedge\flat},-,q}})\ar[r]^{\mathrm{global}}\ar[r]\ar[r] &\varphi_q{\mathrm{Modules}_\circledcirc}(\{\mathrm{Robba}^{\mathrm{extended},q}_{{R_0,-,I}}\}_I).\\ 
}
\]
The definition is given by the following:
\[
\xymatrix@R+0pc@C+0pc{
\mathrm{homotopycolimit}_i(\square)(\mathcal{O}_i),  
}
\]
each $\mathcal{O}_i$ is just as $\mathbb{Q}_p\left<C_1,...,C_\ell\right>,\ell=1,2,...$ and $\square$ is the relative diagram of $\infty$-functors.

\item (\text{Proposition}) There is a functor (global section) between the $\infty$-prestacks of inductive Banach quasicoherent commutative algebra $E_\infty$ objects\footnote{Here $\circledcirc=\text{solidquasicoherentsheaves}$.}:
\[\displayindent=-0.4in
\xymatrix@R+6pc@C+0pc{
\mathrm{sComm}_\mathrm{simplicial}{\mathrm{Modules}_\circledcirc}_{\Gamma_k^{\times q}}(\mathcal{O}_{X_{R_k,-,q}})\ar[d]\ar[d]\ar[d]\ar[d]\ar[r]^{\mathrm{global}}\ar[r]\ar[r] &\mathrm{sComm}_\mathrm{simplicial}\varphi_q{\mathrm{Modules}_\circledcirc}_{\Gamma_k^{\times q}}(\{\mathrm{Robba}^{\mathrm{extended},q}_{{R_k,-,I}}\}_I)\ar[d]\ar[d]\ar[d]\ar[d]\\
\mathrm{sComm}_\mathrm{simplicial}{\mathrm{Modules}_\circledcirc}_{\Gamma_0^{\times q}}(\mathcal{O}_{X_{\mathbb{Q}_p(p^{1/p^\infty})^{\wedge\flat},-,q}})\ar[r]^{\mathrm{global}}\ar[r]\ar[r] &\mathrm{sComm}_\mathrm{simplicial}\varphi_q{\mathrm{Modules}_\circledcirc}_{\Gamma_0^{\times q}}(\{\mathrm{Robba}^{\mathrm{extended},q}_{{R_0,-,I}}\}_I).  
}
\]
The definition is given by the following:
\[
\xymatrix@R+0pc@C+0pc{
\mathrm{homotopycolimit}_i(\square)(\mathcal{O}_i),  
}
\]
each $\mathcal{O}_i$ is just as $\mathbb{Q}_p\left<C_1,...,C_\ell\right>,\ell=1,2,...$ and $\square$ is the relative diagram of $\infty$-functors.

\item Then as in \cite{LBV} we have a functor (global section) of the de Rham complex after \cite[Definition 5.9, Section 5.2.1]{KKM}\footnote{Here $\circledcirc=\text{solidquasicoherentsheaves}$.}:
\[\displayindent=-0.6in
\xymatrix@R+6pc@C+0pc{
\mathrm{deRham}_{\mathrm{sComm}_\mathrm{simplicial}{\mathrm{Modules}_\circledcirc}_{\Gamma_k^{\times q}}(\mathcal{O}_{X_{R_k,-,q}})\ar[r]^{\mathrm{global}}}(-)\ar[d]\ar[d]\ar[d]\ar[d]\ar[r]\ar[r] &\mathrm{deRham}_{\mathrm{sComm}_\mathrm{simplicial}\varphi_q{\mathrm{Modules}_\circledcirc}_{\Gamma_k^{\times q}}(\{\mathrm{Robba}^{\mathrm{extended},q}_{{R_k,-,I}}\}_I)}(-)\ar[d]\ar[d]\ar[d]\ar[d]\\
\mathrm{deRham}_{\mathrm{sComm}_\mathrm{simplicial}{\mathrm{Modules}_\circledcirc}_{\Gamma_0^{\times q}}(\mathcal{O}_{X_{\mathbb{Q}_p(p^{1/p^\infty})^{\wedge\flat},-,q}})\ar[r]^{\mathrm{global}}}(-)\ar[r]\ar[r] &\mathrm{deRham}_{\mathrm{sComm}_\mathrm{simplicial}\varphi_q{\mathrm{Modules}_\circledcirc}_{\Gamma_0^{\times q}}(\{\mathrm{Robba}^{\mathrm{extended},q}_{{R_0,-,I}}\}_I)}(-), 
}
\]

The definition is given by the following:
\[
\xymatrix@R+0pc@C+0pc{
\mathrm{homotopycolimit}_i(\square)(\mathcal{O}_i),  
}
\]
each $\mathcal{O}_i$ is just as $\mathbb{Q}_p\left<C_1,...,C_\ell\right>,\ell=1,2,...$ and $\square$ is the relative diagram of $\infty$-functors.

\item Then we have the following a functor (global section) of $K$-group $(\infty,1)$-spectrum from \cite{BGT}\footnote{Here $\circledcirc=\text{solidquasicoherentsheaves}$.}:
\[
\xymatrix@R+6pc@C+0pc{
\mathrm{K}^\mathrm{BGT}_{\mathrm{sComm}_\mathrm{simplicial}{\mathrm{Modules}_\circledcirc}_{\Gamma_k^{\times q}}(\mathcal{O}_{X_{R_k,-,q}})\ar[r]^{\mathrm{global}}}(-)\ar[d]\ar[d]\ar[d]\ar[d]\ar[r]\ar[r] &\mathrm{K}^\mathrm{BGT}_{\mathrm{sComm}_\mathrm{simplicial}\varphi_q{\mathrm{Modules}_\circledcirc}_{\Gamma_k^{\times q}}(\{\mathrm{Robba}^{\mathrm{extended},q}_{{R_k,-,I}}\}_I)}(-)\ar[d]\ar[d]\ar[d]\ar[d]\\
\mathrm{K}^\mathrm{BGT}_{\mathrm{sComm}_\mathrm{simplicial}{\mathrm{Modules}_\circledcirc}_{\Gamma_0^{\times q}}(\mathcal{O}_{X_{\mathbb{Q}_p(p^{1/p^\infty})^{\wedge\flat},-,q}})\ar[r]^{\mathrm{global}}}(-)\ar[r]\ar[r] &\mathrm{K}^\mathrm{BGT}_{\mathrm{sComm}_\mathrm{simplicial}\varphi_q{\mathrm{Modules}_\circledcirc}_{\Gamma_0^{\times q}}(\{\mathrm{Robba}^{\mathrm{extended},q}_{{R_0,-,I}}\}_I)}(-), 
}
\]

The definition is given by the following:
\[
\xymatrix@R+0pc@C+0pc{
\mathrm{homotopycolimit}_i(\square)(\mathcal{O}_i),  
}
\]
each $\mathcal{O}_i$ is just as $\mathbb{Q}_p\left<C_1,...,C_\ell\right>,\ell=1,2,...$ and $\square$ is the relative diagram of $\infty$-functors.\\

\end{itemize}

\

\begin{proposition}
All the global functors from \cite[Proposition 13.8, Theorem 14.9, Remark 14.10]{1CS2} achieve the equivalences on both sides.	
\end{proposition}

\newpage
\subsection{$\infty$-Categorical Analytic Stacks and Descents V}

Here we consider the corresponding archimedean picture, after \cite[Problem A.4, Kedlaya's Lecture]{1CBCKSW}. Recall for any algebraic variety $R$ over $\mathbb{R}$ this $X_R(\mathbb{C})$ is defined to be the corresponding quotient:
\begin{align}
X_{R}(\mathbb{C}):=R(\mathbb{C})\times P^1(\mathbb{C})/\varphi,\\
Y_R(\mathbb{C}):=R(\mathbb{C})\times P^1(\mathbb{C}).	
\end{align}
The Hodge structure is given by $\varphi$. We define the relative version by considering a further algebraic variety over $\mathbb{C}$, say $A$ as in the following:
\begin{align}
X_{R,A}(\mathbb{C}):=R(\mathbb{C})\times P^1(\mathbb{C})\times A(\mathbb{C})/\varphi,\\
Y_{R,A}(\mathbb{C}):=R(\mathbb{C})\times P^1(\mathbb{C})\times A(\mathbb{C}).	
\end{align}

\indent We then take $q$-th self power to achieve $X_{R,q}(\mathbb{C})$ as 
\begin{align}
X_{R,q}(\mathbb{C}):=(R(\mathbb{C})\times P^1(\mathbb{C}))^q/\varphi_q,\\
Y_{R,q}(\mathbb{C}):=(R(\mathbb{C})\times P^1(\mathbb{C}))^q.	
\end{align}
The multi hyperk\"ahler Hodge structure is given by $\varphi_q$. We define the relative version by considering a further algebraic variety over $\mathbb{C}$, say $A$ as in the following:
\begin{align}
X_{R,A}(\mathbb{C}):=(R(\mathbb{C})\times P^1(\mathbb{C}))^q\times A(\mathbb{C})/\varphi_q,\\
Y_{R,A}(\mathbb{C}):=(R(\mathbb{C})\times P^1(\mathbb{C}))^q\times A(\mathbb{C}).	
\end{align}

Then by \cite{1BBK} and \cite{1CS2} we have the corresponding $\infty$-category of $\infty$-sheaves of simplicial ind-Banach quasicoherent modules which in our situation will be assumed to the modules in \cite{1BBK}, as well as the corresponding associated Clausen-Scholze spaces:
\begin{align}
X_{R}(\mathbb{C}):=R(\mathbb{C})\times P^1(\mathbb{C})^\blacksquare/\varphi,\\
Y_R(\mathbb{C}):=R(\mathbb{C})\times P^1(\mathbb{C})^\blacksquare.	
\end{align}
\begin{align}
X_{R,A}(\mathbb{C}):=R(\mathbb{C})\times P^1(\mathbb{C})\times A(\mathbb{C})^\blacksquare/\varphi,\\
Y_{R,A}(\mathbb{C}):=R(\mathbb{C})\times P^1(\mathbb{C})\times A(\mathbb{C})^\blacksquare,	
\end{align}
with the $\infty$-category of $\infty$-sheaves of simplicial liquid quasicoherent modules, liquid vector bundles and liquid perfect complexes, with further descent \cite[Proposition 13.8, Theorem 14.9, Remark 14.10]{1CS2}.\\
 
\indent We then take $q$-th self power to achieve $X_{R,q}(\mathbb{C})$ as 
\begin{align}
X_{R,q}(\mathbb{C}):=(R(\mathbb{C})\times P^1(\mathbb{C}))^{q,\blacksquare}/\varphi_q,\\
Y_{R,q}(\mathbb{C}):=(R(\mathbb{C})\times P^1(\mathbb{C}))^{q,\blacksquare}.	
\end{align}
The multi hyperk\"ahler Hodge structure is given by $\varphi_q$. We define the relative version by considering a further algebraic variety over $\mathbb{C}$, say $A$ as in the following:
\begin{align}
X_{R,A}(\mathbb{C}):=(R(\mathbb{C})\times P^1(\mathbb{C}))^{q,\blacksquare}\times A(\mathbb{C})/\varphi_q,\\
Y_{R,A}(\mathbb{C}):=(R(\mathbb{C})\times P^1(\mathbb{C}))^{q,\blacksquare}\times A(\mathbb{C}).	
\end{align}

We call the resulting global sections are the corresponding $c$-equivariant Hodge Modules. Then we have the following direct analogy:

\begin{itemize}
\item (\text{Proposition}) There is an equivalence between the $\infty$-categories of inductive Banach quasicoherent presheaves:
\[
\xymatrix@R+0pc@C+0pc{
\mathrm{Ind}\mathrm{Banach}(\mathcal{O}_{X_{R,A,q}})\ar[r]^{\mathrm{equi}}\ar[r]\ar[r] &\varphi_q\mathrm{Ind}\mathrm{Banach}(\mathcal{O}_{Y_{R,A,q}}).  
}
\]
\item (\text{Proposition}) There is an equivalence between the $\infty$-categories of monomorphic inductive Banach quasicoherent presheaves:
\[
\xymatrix@R+0pc@C+0pc{
\mathrm{Ind}^m\mathrm{Banach}(\mathcal{O}_{X_{R,A,q}})\ar[r]^{\mathrm{equi}}\ar[r]\ar[r] &\varphi_q\mathrm{Ind}^m\mathrm{Banach}(\mathcal{O}_{Y_{R,A,q}}).  
}
\]
\end{itemize}

\begin{itemize}

\item (\text{Proposition}) There is an equivalence between the $\infty$-categories of inductive Banach quasicoherent presheaves:
\[
\xymatrix@R+0pc@C+0pc{
\mathrm{Ind}\mathrm{Banach}(\mathcal{O}_{X_{R,A,q}})\ar[r]^{\mathrm{equi}}\ar[r]\ar[r] &\varphi_q\mathrm{Ind}\mathrm{Banach}(\mathcal{O}_{Y_{R,A,q}}).  
}
\]
\item (\text{Proposition}) There is an equivalence between the $\infty$-categories of monomorphic inductive Banach quasicoherent presheaves:
\[
\xymatrix@R+0pc@C+0pc{
\mathrm{Ind}^m\mathrm{Banach}(\mathcal{O}_{X_{R,A,q}})\ar[r]^{\mathrm{equi}}\ar[r]\ar[r] &\varphi_q\mathrm{Ind}^m\mathrm{Banach}(\mathcal{O}_{Y_{R,A,q}}).  
}
\]
\item (\text{Proposition}) There is an equivalence between the $\infty$-categories of inductive Banach quasicoherent commutative algebra $E_\infty$ objects:
\[
\xymatrix@R+0pc@C+0pc{
\mathrm{sComm}_\mathrm{simplicial}\mathrm{Ind}\mathrm{Banach}(\mathcal{O}_{X_{R,A,q}})\ar[r]^{\mathrm{equi}}\ar[r]\ar[r] &\mathrm{sComm}_\mathrm{simplicial}\varphi_q\mathrm{Ind}\mathrm{Banach}(\mathcal{O}_{Y_{R,A,q}}).  
}
\]
\item (\text{Proposition}) There is an equivalence between the $\infty$-categories of monomorphic inductive Banach quasicoherent commutative algebra $E_\infty$ objects:
\[
\xymatrix@R+0pc@C+0pc{
\mathrm{sComm}_\mathrm{simplicial}\mathrm{Ind}^m\mathrm{Banach}(\mathcal{O}_{X_{R,A,q}})\ar[r]^{\mathrm{equi}}\ar[r]\ar[r] &\mathrm{sComm}_\mathrm{simplicial}\varphi_q\mathrm{Ind}^m\mathrm{Banach}(\mathcal{O}_{Y_{R,A,q}}).  
}
\]

\item Then parallel as in \cite{LBV} we have the equivalence of the de Rham complex after \cite[Definition 5.9, Section 5.2.1]{KKM}:
\[
\xymatrix@R+0pc@C+0pc{
\mathrm{deRham}_{\mathrm{sComm}_\mathrm{simplicial}\mathrm{Ind}\mathrm{Banach}(\mathcal{O}_{X_{R,A,q}})\ar[r]^{\mathrm{equi}}}(-)\ar[r]\ar[r] &\mathrm{deRham}_{\mathrm{sComm}_\mathrm{simplicial}\varphi_q\mathrm{Ind}\mathrm{Banach}(\mathcal{O}_{Y_{R,A,q}})}(-), 
}
\]
\[
\xymatrix@R+0pc@C+0pc{
\mathrm{deRham}_{\mathrm{sComm}_\mathrm{simplicial}\mathrm{Ind}^m\mathrm{Banach}(\mathcal{O}_{X_{R,A,q}})\ar[r]^{\mathrm{equi}}}(-)\ar[r]\ar[r] &\mathrm{deRham}_{\mathrm{sComm}_\mathrm{simplicial}\varphi_q\mathrm{Ind}^m\mathrm{Banach}(\mathcal{O}_{Y_{R,A,q}})}(-). 
}
\]

\item Then we have the following equivalence of $K$-group $(\infty,1)$-spectrum from \cite{BGT}:
\[
\xymatrix@R+0pc@C+0pc{
\mathrm{K}^\mathrm{BGT}_{\mathrm{sComm}_\mathrm{simplicial}\mathrm{Ind}\mathrm{Banach}(\mathcal{O}_{X_{R,A,q}})\ar[r]^{\mathrm{equi}}}(-)\ar[r]\ar[r] &\mathrm{K}^\mathrm{BGT}_{\mathrm{sComm}_\mathrm{simplicial}\varphi_q\mathrm{Ind}\mathrm{Banach}(\mathcal{O}_{Y_{R,A,q}})}(-), 
}
\]
\[
\xymatrix@R+0pc@C+0pc{
\mathrm{K}^\mathrm{BGT}_{\mathrm{sComm}_\mathrm{simplicial}\mathrm{Ind}^m\mathrm{Banach}(\mathcal{O}_{X_{R,A,q}})\ar[r]^{\mathrm{equi}}}(-)\ar[r]\ar[r] &\mathrm{K}^\mathrm{BGT}_{\mathrm{sComm}_\mathrm{simplicial}\varphi_q\mathrm{Ind}^m\mathrm{Banach}(\mathcal{O}_{Y_{R,A,q}})}(-). 
}
\]
\end{itemize}

\begin{assumption}\label{assumtionpresheaves}
All the functors of modules or algebras below are Clausen-Scholze sheaves \cite[Proposition 13.8, Theorem 14.9, Remark 14.10]{1CS2}. 	
\end{assumption}

\begin{itemize}
\item (\text{Proposition}) There is an equivalence between the $\infty$-categories of inductive liquid sheaves:
\[
\xymatrix@R+0pc@C+0pc{
\mathrm{Module}_\circledcirc(\mathcal{O}_{X_{R,A,q}})\ar[r]^{\mathrm{equi}}\ar[r]\ar[r] &\varphi_q\mathrm{Module}_\circledcirc(\mathcal{O}_{Y_{R,A,q}}).  
}
\]
\end{itemize}

\begin{itemize}

\item (\text{Proposition}) There is an equivalence between the $\infty$-categories of inductive Banach quasicoherent commutative algebra $E_\infty$ objects:
\[
\xymatrix@R+0pc@C+0pc{
\mathrm{sComm}_\mathrm{simplicial}\mathrm{Module}_{\text{liquidquasicoherentsheaves}}(\mathcal{O}_{X_{R,A,q}})\ar[r]^{\mathrm{equi}}\ar[r]\ar[r] &\mathrm{sComm}_\mathrm{simplicial}\varphi_q\mathrm{Module}_{\text{liquidquasicoherentsheaves}}(\mathcal{O}_{Y_{R,A,q}}).  
}
\]

\item Then as in \cite{LBV} we have the equivalence of the de Rham complex after \cite[Definition 5.9, Section 5.2.1]{KKM}\footnote{Here $\circledcirc=\text{liquidquasicoherentsheaves}$.}:
\[
\xymatrix@R+0pc@C+0pc{
\mathrm{deRham}_{\mathrm{sComm}_\mathrm{simplicial}\mathrm{Module}_\circledcirc(\mathcal{O}_{X_{R,A,q}})\ar[r]^{\mathrm{equi}}}(-)\ar[r]\ar[r] &\mathrm{deRham}_{\mathrm{sComm}_\mathrm{simplicial}\varphi_q\mathrm{Module}_\circledcirc(\mathcal{O}_{Y_{R,A,q}})}(-). 
}
\]

\item Then we have the following equivalence of $K$-group $(\infty,1)$-spectrum from \cite{BGT}\footnote{Here $\circledcirc=\text{liquidquasicoherentsheaves}$.}:
\[
\xymatrix@R+0pc@C+0pc{
\mathrm{K}^\mathrm{BGT}_{\mathrm{sComm}_\mathrm{simplicial}\mathrm{Module}_\circledcirc(\mathcal{O}_{X_{R,A,q}})\ar[r]^{\mathrm{equi}}}(-)\ar[r]\ar[r] &\mathrm{K}^\mathrm{BGT}_{\mathrm{sComm}_\mathrm{simplicial}\varphi_q\mathrm{Module}_\circledcirc(\mathcal{O}_{Y_{R,A,q}})}(-). 
}
\]
\end{itemize}

\bibliographystyle{splncs}

\newpage\chapter{Hodge-Iwasawa Theory I}

\newpage\section{Introduction}

\subsection{Hodge-Iwasawa Theory and Higher Dimensional Modeling of the Weil Conjectures}

This paper deals with some simultaneous generalization of the noncommutative Iwasawa theory after Kedlaya-Pottharst and the relative $p$-adic Hodge theory in the style of Kedlaya-Liu, with further sophisticated philosophy in mind dated back to Fukaya-Kato and Kato. To make the motivations and the applications more serious and clear to the readers, let us start from recalling what happens along the intersection of the two pictures above.

\begin{example} \mbox{\textbf{(After Kedlaya-Pottharst)}}
Let $K$	be a finite extension of $\mathbb{Q}_p$. In \cite{2KP1}, the authors proposed that one could consider a family $\{L\}_L$ parametrized by perfectoid Galois fields of intrinsic descriptions of the category $\widetilde{\Phi\Gamma}_{K,A}$ (with the notations in \cite{2KP1}) in terms of $\varphi$-modules over $\widetilde{\Pi}_{L,A}$ in the family with the further action of the groups $G(L/K)$ for each $L$ in the family. These actually for a fixed $L/K$  will give rise to some perfectoid $p$-adic Hodge theory immediately. Then they considered some geometrization which is just some application of the theory in \cite[Chapter 9]{2KL15}, namely one first forgets the action of the group $G_{L/K}$ but replaces this by pro-\'etale topology over $\mathrm{Spa}(K,\mathcal{O}_K)$. Then one defines an object $\mathcal{M}$ in the category $\widetilde{\Phi\Gamma}_{K,A}$ which has no a priori relation to any tower $L/K$. One recovers the Iwasawa theory from the sheaf attached to the Iwasawa algebra attached to the group $G(L/K)$. This makes more sense if we know clearly the structure of $G(L/K)$. The memorization of the Iwasawa structure back from the Iwasawa deformation (instead of taking the limit throughout some rigidized towers) will be some key observation in our development. 
\end{example}

\indent On the other hand, we have the following example:

\begin{example}
Let $A$ be now an Artin algebra over $\mathbb{Q}_p$, in this situation one could actually consider the local systems $\underline{A}$ over some specific smooth proper rigid analytic space $X/\mathbb{Q}_p$. Therefore one could consider the situation where this deforms from some $\mathbb{Q}_p$-local systems. Then one should be able to relate constructible $\underline{A}$-local systems to the corresponding deformation of the corresponding $\mathbb{Q}_p$ representation of the corresponding \'etale fundamental groups. One could then as in \cite{2KL15} and \cite{2KL16} consider the corresponding Galois groups of affinoids in both rigid geometry, adic geometry or perfectoid geometry. Then in these situations, one could consider the corresponding deformation problem. Furthermore one could then consider more general rigid affinoids $A/\mathbb{Q}_p$ and consider the corresponding families of representations of the Galois groups mentioned above over $A$. In order to study these one could then construct the corresponding $A$-relative relative period rings perfect generalizing \cite{2KL15} and \cite{2KL16}. 	
\end{example}

\indent Therefore now we consider the combination of the two pictures and the corresponding simultaneous generalization of the pictures in the relative $p$-adic Hodge theory and the Iwasawa theory in the sense mentioned above. Recall that the toric chart in \cite{2KL16} could admit some cyclotomic deformation with the corresponding $\Gamma$-action where $\Gamma$ will be just $\mathbb{Z}_p^{d}$ for some rank $d$. Then one could consider the corresponding picture in the Iwasawa theory to consider the corresponding $(\varphi,\Gamma)$-cohomology of the corresponding multivariate cyclotomic deformations. We will call this generalized Iwasawa cohomology and the generalized Iwasawa theory, since over the perfectoid covering of the base space over $\widehat{\mathbb{Q}_p(p^{1/p^{\infty}})}$, one has the corresponding relative de Rham period rings. Therefore one could consider the generalized Bloch-Kato dual exponential maps and so on.\\

\indent One should somehow regard that Iwasawa theory is mimicking the celebrated Weil conjectures, namely the Riemann Hypothesis for the arithmetic of the algebraic curves. The conjectures themselves actually do not have too much sort of motivic motivation a priori, for instance one could just ask the questions on the counting problems within the arithmetic of algebraic curves. However the conjectures to some extent are definitely admitting deep and fundamental motivic point of views. Since the issue is that one could define the corresponding $L$-functions in terms of some Hodge structure for algebraic varieties over finite fields, namely the Frobenius actions on the $\ell$-adic \'etale cohomologies. Also one could actually directly reinterpret the corresponding number theoretic properties in term of the corresponding cohomological properties, for instance the zeros of zeta functions and the functional equations, which is further adopted by Perrin-Riou in her celebrated work for instance. Iwasawa established the corresponding point of view successfully in the context of number fields, which is actually to some extent very fascinating and deep. Although the original pictures from Iwasawa admit more class group formulation and the corresponding motivation (namely more Galois theoretic or more \'etale), which is similar to the original motivation of Weil conjectures, but following later philosophy like from Kato everything could be formulated quite cohomologically, namely due to some deep relationship and similarization between the \'etale cohomology and Galois cohomology.\\

\[\small
\xymatrix@R+4.5pc@C+0pc{
\text{Weil Conjectures for Curves over $k/\mathbb{F}_p$}
\ar[r]\ar[r]\ar[r]\ar[d]\ar[d]\ar[d] &\text{Deligne-Kedlaya-Caro-Abe Weil II for Varieties} \ar[d]\ar[d]\ar[d]\\
\text{\'Etale Iwasawa Theory for $\mathbb{Q}$}
\ar[u]\ar[u]\ar[u]\ar[r]\ar[r]\ar[r]\ar[d]\ar[d]\ar[d] &\text{??? for Varieties over $\mathbb{Q}$} \ar[u]\ar[u]\ar[u]\ar[d]\ar[d]\ar[d].\\
\text{Non-\'etale Iwasawa Theory}
\ar[r]\ar[r]\ar[r] &\text{??? for Varieties over $\mathbb{Q}_p$}.
}
\]
\\

\indent As shown in the diagram above, one could see that actually Weil conjectures are quite general, in the sense that first they are actually above varieties over finite fields, which was originally proved by Deligne in terms of \'etale cohomology. By $p$-adic method, Kedlaya reproved the conjectures with actually more general point of view, since there one considers general smooth coefficients which is to say the isocrystals. The key issue is to establish some sort of finiteness of the rigid cohomology with such smooth non-\'etale coefficients. The method adopted by Kedlaya is through deep $p$-adic local monodromy theorem established in \cite{2Ked1} which proves the affine curve cases directly, then through clever geometric method to reduce everything to affine spaces, and then to affine curves with induction. Actually the picture can be made more general in the sense that one could consider more general $p$-adic constructible coefficients, which was considered by Abe and etc, for instance the arithmetic $D$-modules.\\

\indent In the mixed-characteristic situation, Nakamura established some general Iwasawa theory for non-\'etale coefficients (by using Kedlaya-Liu's language these are pro-\'etale non-\'etale Hodge-Frobenius sheaves) over $X=\mathrm{Spa}(\mathbb{Q}_p,\mathcal{O}_{\mathbb{Q}_p})$. Through Kedlaya-Pottharst's point of view, namely the analytic cohomology of the cyclotomic deformation of pro-\'etale Hodge structures. We will give our understanding on the possible $???$ in the entries of the graph above, which is already represented at the beginning of this subsection. Replacing the cyclotomic towers or toric towers by pro-\'etale sites allows us to define the generalized Iwasawa cohomology. Therefore deformation of the Hodge sheaves over pro-\'etale sites will remember the Iwasawa theoretic consideration, which is observed as mentioned above by Kedlaya-Pottharst. \\

%%%%%%%%%%%%%%%%%%%%%%main result one:

\indent  The first main topic exhibited in our main body of work is on the deformation over a reduced affinoid algebra $A$ the corresponding categories of vector bundles and the corresponding Frobenius modules (bundles) over the ring $\widetilde{\Pi}^\infty_{R,A}$ and $\widetilde{\Pi}_{R,A}$ (with the notations in \cref{2corollary7.6} as below):

\begin{theorem}
The categories of $A$-relative quasicoherent finite locally free sheaves over the schematic relative Fargues-Fontaine curve $\mathrm{Proj}P_{R,A}$, the Frobenius modules over $\widetilde{\Pi}^\infty_{R,A}$, the Frobenius modules over $\widetilde{\Pi}_{R,A}$ and the Frobenius bundles over $\widetilde{\Pi}_{R,A}$ are equivalent. Moreover for those Fr\'echet-Stein distribution algebra $A_\infty(G)$ attached to a $p$-adic Lie group $G$ in the style of the Hodge-Iwasawa consideration assumed in \cref{setting3.8}, the corresponding statement holds as well (see \cref{2section3.2}).\\
\end{theorem}

\indent This is actually considered in \cite{2KP1} in the situation where $R$ is a perfectoid field. The main goal in our mind at this moment is essentially some belief (encoded in some possible application) that the deformations (both from some representation reason or Iwasawa-Tamagawa reason) share some rigid equivalence among the vector bundles over relative Fargues-Fontaine curves, $B$-pairs in families and finally the Frobenius-Hodge modules in relative sense. Of course if one focuses on the corresponding Galois representation context, then this belief will be somehow complicated while we believe the categories mentioned here are big enough to eliminate the difference from the deformations and the original absolute situation, without touching the \'etaleness. There are many contexts where such deformation on the level of $\varphi$-modules has been already proved to be important and convenient, such as \cite{2Nak4} and \cite{2KPX}. \cite{2KPX} considered essentially the deformation theory of $\varphi$-modules by using the tools from Fr\'echet algebras. Note that $B$-pairs are generalizations of Galois representations from the point of view of the coefficients, Frobenius modules over ind-Fr\'echet sheaves are generalizations of Galois representations from the point of view of the slope theory and the corresponding purity. The relative deformation of the Hodge structures in such generality looks challenging, for instance the parallel story to Deligne-Laumon, \cite{2AM}, \cite{2G1} and \cite{2G2} looks quite mysterious.\\

\indent We actually discuss in this paper the corresponding pseudocoherent objects as well in the sense of \cite{2SGA6} and the corresponding pseudocoherent Frobenius modules studied in \cite{2KL16}:

\begin{theorem}
The category of the pseudocoherent sheaves over the schematic Fargues-Fontaine curve $\mathrm{Proj}P_{R,A}$ and the category of the pseudocoherent modules over $\widetilde{\Pi}_{R,A}$ endowed with isomorphisms coming from the Frobenius pullbacks are equivalent.	
\end{theorem}

\indent One of our current \textit{main conjectures} (in the \'etale setting) in this article is quite general over a rigid analytic space $X$ (over $\mathbb{Q}_p$) separated and of finite type (in the category of adic spaces) with respect to an adic ring $T$ as in \cite{2Wit1}, \cite{2Wit2} and \cite[1.4]{2Wit3}, and we consider the corresponding categories $\mathbb{D}_\mathrm{perf}(X_\sharp,T)~ (\sharp=\text{\'et},\text{pro\'et})$ of compatible families $(\mathcal{F}^\bullet_I)_{I\subset T}$ complexes perfect and dg-flat of $T/I$ constructible flat \'etale or pro-\'etale sheaves and we conjecture that we have the corresponding Waldhausen enrichment. To be more precise:

\begin{conjecture}
Assume now that $p$ is a unit in the ring $T$. Assume that we have $\mathbb{D}_\mathrm{perf}(X_\sharp,T)~(\sharp=\text{\'et},\text{pro\'et})$ are Waldhausen categories and we have the corresponding continuous maps from the attached $K$-theory topological spaces to the one attached to $\mathbb{D}_\mathrm{perf}(T)$ induced by the direct image functor (which will make the space $X$ more restrictive). Then we have that the maps are null-homotopic in some canonical way. (See \cref{conjecture9.13})
\end{conjecture}

\indent Using relatively homotopical language, this could be directly related to the formulation adopted by Fukaya-Kato in Deligne's category when we focus on the corresponding \textit{spaces of 1 type} where the corresponding homotopy groups will vanish above the degree one. This has actually been considered extensively already in \cite{2Wit1} in the context in terms of SGAIV and V. The corresponding categories we mentioned above looks very complicated to study in the general situation. Since our consideration is highly globalized due to the fact that we are considering the whole categories of all the corresponding interesting objects, instead of each single object. We have discussed some further development beyond the \'etale setting (in the $p$-adic setting) at the end of the paper. \\

\subsection{Why Relative $p$-adic Hodge theory?}

\indent Since our motivation comes from the deformation of representations of \'etale fundamental groups of higher dimensional spaces and the geometrization of Iwasawa theory, it is important to look inside the intrinsic Hodge structures through some coherent tools, which is why we need to consider the corresponding relative $p$-adic Hodge theory in some deformed setting. The main goal of the study around the fundamental groups as illustrated above actually is usually of the following aspects:

\begin{setting}\mbox{}\\
I. Try to understand the intrinsic geometric structures and its reconstruction effect on the varieties. \\
II. Try to understand the representations.	\\
III. Try to understand the deformation of representations.\\
IV. Try to understand the Iwasawa sheaves with coefficients in some equivariant setting (where the equivariance sometimes factors through the \'etale fundamental groups).\\
\end{setting}

\indent Essentially and especially in our setting, the generalization is actually more complicated than the usual situation, namely the representations are relative in the sense that it is varying over a variety which is somewhat as general as possible with some generalized coefficients. In the usual setting, one could regard everything as local systems in most of desired cases. Then using Scholze's pro-\'etale sites and pushing down one could define the corresponding period functors. In our setting one can definitely rely on the advantage of the original ideas of Kedlaya-Liu \cite{2KL15} and \cite{2KL16}, namely transforming the problem of representations of fundamental groups not to the local systems over pro-\'etale sites but to the modules over relative Robba rings over quasi pro-\'etale sites. Under general philosophy from Berger, everything on the $p$-adic Hodge module structure will be encoded in the module structure over relative Robba rings.

\subsection{Further Study}
In our current establishment, we have touched some foundational results along our consideration towards general Hodge-Iwasawa theory. Deformations of the corresponding Hodge structures in our consideration (namely the analytic module structures over Robba sheaves and rings) have beed shown by us to have some well-behaved properties. The aspect of deformation of such general Hodge structures will be definitely one of our long-term goal which should be parallel to those current existing comparison theorems and the existing various types of equivariant cohomological projects (for instance in the setting of \'etale cohomology and rigid cohomology). Over a rigid analytic space, the investigation of the corresponding families of the representations of the \'etale fundamental group seems to be a very hard problem for instance, especially when one would like to consider the equivariant setting, namely the corresponding equivariance coming from the quotient of the \'etale fundamental groups. We actually expect our insight should work in more general setting where one considers some general enough quotients here.\\

On the Iwasawa side, we expect that we can follow and generalize the original ideas in \cite{2KP1} to geometrize the corresponding interesting towers in Tamagawa-Iwasawa theory. This is actually different from some standard construction in the original Iwasawa theory and $\Lambda$-adic Hodge theory. We expect this will give us a more uniform approach to recover many resulting Iwasawa theory of the general $p$-adic Hodge structures, for instance the corresponding exponential maps and the dual exponential maps, and the corresponding Iwasawa interpolation throughout some interesting character varieties, carrying some relativization. One another thing we would like to mention is that actually we have already motivated the corresponding study of the Iwasawa theory in the pseudocoherent setting beyond the vector bundle non-\'etale context.\\

Maybe a tour through the topological method in $p$-adic analysis should be possible in a long term consideration. Carrying some reasonable type of relativization and equivariance, one can apply higher categorical approach to study the corresponding relative $p$-adic Hodge theory in our context. We only here touched upon a piece of the story since our main goal here is Iwasawa deformation of general $p$-adic Hodge structure. But there is no reason to only restrict the consideration to this context. In other word the higher categorical investigation we limited here could motivate some further deeper understanding. Especially when we introduce the corresponding non-commutativity and higher-dimensional equivariance, where the context of the study of the rational $p$-adic Hodge theory looks quite challenging. On the other hand, actually the corresponding discussion above may also shed some light on the corresponding equivariant program in other contexts, such as the arithmetic $D$-modules. We would like to mention that actually the higher categorical aspects of $p$-adic Hodge theory are already encoded in some contexts in the literature, for instance the work \cite{2BMS}.\\

%And over an Artin stack one may also wonder if the corresponding $\infty$-categorical tools could be used to study the local systems (after Kedlaya-Liu) along the story we are considering here (as Scholze applied $\infty$-categorical tools from Lurie to have considered some problems around the cohomology of diamonds and even more general stacks, see \cite{2Scho1}). \\

\subsection{Notations}

\begin{center}
\begin{longtable}{p{7.0cm}p{8cm}}
Notation & Description \\
\hline
$p$ & A prime number.\\
%$Y:=\{y\}$ & A finite set.\\
$(R,R^+)$ & A perfect adic Banach algebra uniform over $\mathbb{F}_{p^h}$. \\
$A$ & A reduced affinoid algebra over $\mathbb{Q}_p$ in the rigid geometry after Tate. In general we consider the spectral norms on $A$ unless specified otherwise. \\
%Or we will consider the complete Iwasawa algebra $\mathbb{Z}_p[[G]]\otimes_{\mathbb{Z}_p}\mathbb{Q}_p$. In general we consider the spectral norms on $A$ unless specified otherwise. \\
$T$ & An adic ring. Note that this is considered in the work for instance \cite{2Wit1}, \cite{2Wit2} and \cite[Section 1.4]{2Wit3}.\\
$A_\infty(G)$ & Fr\'echet-Stein algebra attached to some nice $p$-adic group $G$ of type Lie. \\
$\Omega^\mathrm{int}_{R},	\Omega_{R}$ & Kedlaya-Liu's perfect period rings, over $R$. See \cite{2KL15} and \cite{2KL16}.\\
$\widetilde{\Pi}^{\mathrm{int},r}_{R},\widetilde{\Pi}^{\mathrm{bd},r}_{R},
\widetilde{\Pi}^{I}_{R},\widetilde{\Pi}^{+}_{R},\widetilde{\Pi}^{\infty}_{R},\widetilde{\Pi}_{R}$ & Kedlaya-Liu's perfect period rings, over $R$. See \cite{2KL15} and \cite{2KL16}.\\

$\Omega^\mathrm{int}_{R,A}$ & $A$-relative version of Kedlaya-Liu's perfect period rings, $\Omega^\mathrm{int}_{R,A}$.\\
$\Omega_{R,A}$ & $A$-relative version of Kedlaya-Liu's perfect period rings, $\Omega_{R,A}$.\\

$\widetilde{\Pi}^{I}_{R,A}$  & $A$-relative version of Kedlaya-Liu's relative perfect Robba rings.\\
$\widetilde{\Pi}^{r}_{R,A}$  & $A$-relative version of Kedlaya-Liu's relative perfect Robba rings.\\
$\widetilde{\Pi}^{\infty}_{R,A}$  & $A$-relative version of Kedlaya-Liu's relative perfect Robba rings.\\
$\widetilde{\Pi}^{+}_{R,A}$  & $A$-relative version of Kedlaya-Liu's relative perfect Robba rings.\\

$\widetilde{\Pi}^{\mathrm{int},r}_{R,A},\widetilde{\Pi}^{\mathrm{int}}_{R,A},
\widetilde{\Pi}^{\mathrm{int},+}_{R,A}$ & $A$-relative version of Kedlaya-Liu's relative perfect integral Robba rings.\\
$\widetilde{\Pi}^{\mathrm{bd},r}_{R,A},\widetilde{\Pi}^{\mathrm{bd}}_{R,A},
\widetilde{\Pi}^{\mathrm{bd},+}_{R,A}$ & $A$-relative version of Kedlaya-Liu's relative perfect bounded Robba rings.\\

$\underline{\underline{\Omega}}^\mathrm{int}_A,\underline{\underline{\Omega}}_A$ & $A$-relative version of Kedlaya-Liu's relative perfect sheaves over $\mathrm{Spa}(R,R^+)$.\\ 

$\underline{\underline{\widetilde{\Pi}}}_A^{\mathrm{int},r}, \underline{\underline{\widetilde{\Pi}}}_A^{\mathrm{int}}$ & $A$-relative version of Kedlaya-Liu's relative perfect sheaves over $\mathrm{Spa}(R,R^+)$. \\
$\underline{\underline{\widetilde{\Pi}}}_A^{\mathrm{bd},r},\underline{\underline{\widetilde{\Pi}}}_A^{\mathrm{bd}}$ & $A$-relative version of Kedlaya-Liu's relative perfect sheaves over $\mathrm{Spa}(R,R^+)$. \\

$\underline{\underline{\widetilde{\Pi}}}_A^{I},\underline{\underline{\widetilde{\Pi}}}_A^{\infty},\underline{\underline{\widetilde{\Pi}}}_A,\underline{\underline{\widetilde{\Pi}}}_A^{r}$ & $A$-relative version of Kedlaya-Liu's relative perfect sheaves over $\mathrm{Spa}(R,R^+)$.\\

%%%%%%%%%%%%%%%%%%%%%%%%%%%%%%%%%%%%%%%

$\Omega^\mathrm{int}_{R,A_\infty(G)}$ & $A_\infty(G)$-relative version of Kedlaya-Liu's perfect period rings, $\Omega^\mathrm{int}_{R,A_\infty(G)}$.\\
$\Omega_{R,A_\infty(G)}$ & $A_\infty(G)$-relative version of Kedlaya-Liu's perfect period rings, $\Omega_{R,A_\infty(G)}$.\\

$\widetilde{\Pi}^{I}_{R,A_\infty(G)}$  & $A_\infty(G)$-relative version of Kedlaya-Liu's relative perfect Robba rings.\\
$\widetilde{\Pi}^{r}_{R,A_\infty(G)}$  & $A_\infty(G)$-relative version of Kedlaya-Liu's relative perfect Robba rings.\\
$\widetilde{\Pi}^{\infty}_{R,A_\infty(G)}$  & $A_\infty(G)$-relative version of Kedlaya-Liu's relative perfect Robba rings.\\
$\widetilde{\Pi}^{+}_{R,A_\infty(G)}$  & $A_\infty(G)$-relative version of Kedlaya-Liu's relative perfect Robba rings.\\

$\widetilde{\Pi}^{\mathrm{int},r}_{R,A_\infty(G)},\widetilde{\Pi}^{\mathrm{int}}_{R,A_\infty(G)},
\widetilde{\Pi}^{\mathrm{int},+}_{R,A_\infty(G)}$ & $A_\infty(G)$-relative version of Kedlaya-Liu's relative perfect integral Robba rings.\\
$\widetilde{\Pi}^{\mathrm{bd},r}_{R,A_\infty(G)},\widetilde{\Pi}^{\mathrm{bd}}_{R,A_\infty(G)},
\widetilde{\Pi}^{\mathrm{bd},+}_{R,A_\infty(G)}$ & $A_\infty(G)$-relative version of Kedlaya-Liu's relative perfect bounded Robba rings.\\

$\underline{\underline{\Omega}}^\mathrm{int}_{A_\infty(G)},\underline{\underline{\Omega}}_{A_\infty(G)}$ & $A_\infty(G)$-relative version of Kedlaya-Liu's relative perfect sheaves over $\mathrm{Spa}(R,R^+)$.\\ 

$\underline{\underline{\widetilde{\Pi}}}_{A_\infty(G)}^{\mathrm{int},r}, \underline{\underline{\widetilde{\Pi}}}_{A_\infty(G)}^{\mathrm{int}}$ & $A_\infty(G)$-relative version of Kedlaya-Liu's relative perfect sheaves over $\mathrm{Spa}(R,R^+)$. \\
$\underline{\underline{\widetilde{\Pi}}}_{A_\infty(G)}^{\mathrm{bd},r},\underline{\underline{\widetilde{\Pi}}}_{A_\infty(G)}^{\mathrm{bd}}$ & $A_\infty(G)$-relative version of Kedlaya-Liu's relative perfect sheaves over $\mathrm{Spa}(R,R^+)$. \\

$\underline{\underline{\widetilde{\Pi}}}_{A_\infty(G)}^{I},\underline{\underline{\widetilde{\Pi}}}_{A_\infty(G)}^{\infty},\underline{\underline{\widetilde{\Pi}}}_{A_\infty(G)},\underline{\underline{\widetilde{\Pi}}}_{A_\infty(G)}^{r}$ & $A_\infty(G)$-relative version of Kedlaya-Liu's relative perfect sheaves over $\mathrm{Spa}(R,R^+)$.\\

%%%%%%%%%%%%%%%%%%%%%%%%

$X$ & Locally $v$-ringed spaces or preadic spaces.\\

$X,X_{\text{\'et}},X_{\text{pro\'et}}$ & The corresponding analytic sites, \'etale sites, pro-\'etale sites of $X$.\\

\end{longtable}
\end{center}

%%%%%%\newpage

\newpage\section{Frobenius Modules over Ind-Fr\'echet Algebras}

\indent Now we define the Frobenius modules in our setting, which are in some sense the generalized Hodge structures. Following \cite{2KL15} and \cite{2KL16} we know that actually the main tool in the study will be to find the links between the local systems or representations of fundamental groups with the Frobenius modules. Note that this is actually a common feature of $p$-adic cohomology theories or if you wish the $p$-adic weights after Deligne et al, which is to say in more detail: initially one has not the chance to see the intrinsic Hodge structures from the coefficient systems or the representations which is to say for instance the Frobenius structure, instead one can only have the chance to observe this after the application of Fontaine's functors. This is actually not the same as the situation in the initial Hodge theory and the $\ell$-adic weights theory in the sense again Deligne et al.\\

\indent So let us consider the following notations of Frobenius modules and Frobenius bundles, where the latter will be some collection of kind of sections with respect to each interval.

\begin{setting}
Following \cite[Chapter 5]{2KL15}, in the notation of both \cite{2KL16} and \cite{2KL15} we will in our context consider the situation where $h=1$ in our initial consideration on the finite field $\mathbb{F}_{p^h}$. And note that we could have another parameter $a$ in our initial setting of the definitions, so we are going to consider $a$-th power Frobenius $\varphi^a$ in the following. This means that in the characteristic $0$ situation we will consider $\mathbb{Q}_{p^a}$ and in the characteristic $p>0$ we consider working over the finite field $\mathbb{F}_{p^a}$ for some chosen positive integer $a\geq 1$.
%which will be first of all generalization of the situation as in  \cite{2KL15} to our situation.
\end{setting}

\begin{definition}
For general $A$ which is a reduced affinoid algebra over $\mathbb{Q}_p$ with integral subring $\mathfrak{o}_A$, we define the following various $A$-relative version of the corresponding Period rings from \cite[Chapter 5]{2KL15}:
\begin{align}
\widetilde{\Pi}_{R,A}:=\bigcup_{r}	\widetilde{\Pi}^{r}_{R,A},\widetilde{\Pi}^{r}_{R,A}:=	\widetilde{\Pi}^{r}_{R}\widehat{\otimes}_{\mathbb{Q}_p}A,
\end{align}
\begin{align}
\widetilde{\Pi}^{I}_{R,A}:=	\widetilde{\Pi}^{I}_{R}\widehat{\otimes}_{\mathbb{Q}_p}A,\widetilde{\Pi}^\infty_{R,A}:=	\widetilde{\Pi}^\infty_{R}\widehat{\otimes}_{\mathbb{Q}_p}A.
\end{align}
\begin{align}
\widetilde{\Pi}^{\mathrm{bd}}_{R,A}:=\bigcup_{r}	\widetilde{\Pi}^{\mathrm{bd},r}_{R}\widehat{\otimes}_{\mathbb{Q}_p}A,\widetilde{\Pi}^{\mathrm{bd},r}_{R,A}:=	\widetilde{\Pi}^{\mathrm{bd},r}_{R}\widehat{\otimes}_{\mathbb{Q}_p}A,
\end{align}
\begin{align}
\widetilde{\Pi}^{\mathrm{int}}_{R,A}:=\bigcup_{r}	\widetilde{\Pi}^{\mathrm{int},r}_{R}\widehat{\otimes}_{\mathbb{Z}_p}\mathfrak{o}_A,\widetilde{\Pi}^{\mathrm{int},r}_{R,A}:=	\widetilde{\Pi}^{\mathrm{int},r}_{R}\widehat{\otimes}_{\mathbb{Z}_p}\mathfrak{o}_A.
\end{align}

\end{definition}

\begin{remark}
One can also consider more deformed version of the period rings as in \cite{2KL15} in some natural way which we do not explicitly write.
\end{remark}

\begin{example}
To be more concrete, in the situation where $A$ is $\mathbb{Q}_p\left<T_1,...,T_d\right>$ for some $d$, we can consider the following more explicit definitions. First form the ring of Witt vectors $W(R)$, then we consider the algebra $W(R)_A$ which is defined by taking the corresponding completed tensor product whose elements have the form of:
\begin{displaymath}
\sum_{n\geq 0,i_1\geq 0,...,i_n \geq 0}p^n[\overline{x}_{n,i_1,...,i_d}]T_1^{i_1}T_2^{i_2}...T_d^{i_d}.
\end{displaymath}
Over $W(R)_A$ we have the Gauss norm $\left\|.\right\|_{\alpha^r,A}$ for each $r>0$ which is defined by 
\begin{displaymath}
\left\|.\right\|_{\alpha^r,A}(\sum_{n\geq 0,i_1\geq 0,...,i_n \geq 0}p^n[\overline{x}_{n,i_1,...,i_d}]T_1^{i_1}T_2^{i_2}...T_d^{i_d}):=\sup_{n\geq 0,i_1\geq 0,...,i_n \geq 0}\{p^{-n}\alpha(\overline{x}_{n,i_1,...,i_d})^r\}.
\end{displaymath}
Then we define $\widetilde{\Pi}_{R,A}^{\mathrm{int},r}$ to be the subring of $W(R)_A$ consisting of all the elements satisfying:
\begin{displaymath}
\lim_{n,i_1,...,i_d}\left\|.\right\|_{\alpha^r,A}(p^{-n}\alpha(\overline{x}_{n,i_1,...,i_d})^r)=0.	
\end{displaymath}
And we put $\widetilde{\Pi}_{R,A}^{\mathrm{int}}$ to be the union of $\widetilde{\Pi}_{R,A}^{\mathrm{int},r}$ for all $r>0$. We then define $\widetilde{\Pi}_{R,A}^{\mathrm{bd},r}$ to be $\widetilde{\Pi}_{R,A}^{\mathrm{int},r}[1/p]$, and furthermore define the ring $\widetilde{\Pi}_{R,A}^{\mathrm{bd}}$ to be the corresponding union throughout all $r>0$. Then we define the Robba ring $\widetilde{\Pi}^I_{R,A}$ for some interval $I$ as the Fr\'echet completion of the ring $W(R)_A[1/p]$ by using the family of norms $\left\|.\right\|_{\alpha^t,A}$ for all $t\in I$. Then one could take specific intervals $(0,r]$ or $(0,\infty)$ to define the corresponding Robba rings $\widetilde{\Pi}^r_{R,A}$ and $\widetilde{\Pi}^\infty_{R,A}$. Finally we can define the whole Robba ring $\widetilde{\Pi}_{R,A}$ as the union of all the $\widetilde{\Pi}^r_{R,A}$ throughout all $r>0$.

\end{example}

\begin{setting}
We need to specify the corresponding Frobenius in our setting, where we will consider the Frobenius action which is $A$-linear but acts trivially over the $A$-part in the products appeared as above in the definitions. 	
\end{setting}

\begin{proposition}
For $A$ a reduced affinoid algebra over $\mathbb{Q}_p$ as defined above, we have the corresponding equality:
\begin{displaymath}
\widetilde{\Pi}^{[s_{1},r_{1}]}_{R,A}\bigcap\widetilde{\Pi}^{[s_{2},r_{2}]}_{R,A}=\widetilde{\Pi}^{[s_{1},r_{2}]}_{R,A}.	
\end{displaymath}
Here we assume that $0<s_1\leq s_2\leq r_1\leq r_2$.	
\end{proposition}

\begin{proof}
This could be proved by using the strategy which mimicks \cite[Section 5.2]{2KL15}, or one could consider any representation taking the form of $\sum_{i}r_i\otimes a_i$ with $r_i$ in the non-relative Robba rings involved and $a_i\in A$, then the result follows from the situation without considering the algebra $A$.	
\end{proof}

\begin{definition} \label{Def2.6}
Consider the $A$-relative ring of periods 
\begin{center}
$\triangle:\Omega_{R,A},\widetilde{\Pi}^\mathrm{bd}_{R,A},\widetilde{\Pi}^\mathrm{int}_{R,A}$, $\widetilde{\Pi}_{R,A},\widetilde{\Pi}^+_{R,A},\widetilde{\Pi}^\infty_{R,A}$. 
\end{center}
Then we are going to define the $A$-relative $\varphi^a$-module over $\triangle$ to be a finite locally free $\triangle$-module carrying semilinear Frobenius action of $\varphi^a$. Moreover now for $\triangle:\widetilde{\Pi}^\mathrm{bd}_{R,A},\widetilde{\Pi}^\mathrm{int}_{R,A},\widetilde{\Pi}_{R,A}$ we have for each radius $r>0$ a module $M_{r}$ over $\triangle^{r}$ descending the module $M$ and the isomorphism 
\begin{center}
$\varphi^*M_r  \overset{\sim}{\rightarrow}M_r\otimes_{\triangle^{r}}\triangle^{r/q}$.
\end{center}
We will call the module $M_{r}$ over $\triangle^{r}$ a model of $M$. One can also define any $A$-relative Frobenius module to be one as above coming from some base change of some model $M_r$.
\end{definition}

\indent Then one has the corresponding bundles in our setting as in the following generalizing the situation in \cite[Section 6.1]{2KL15}.

\begin{definition}
Consider the period ring $\triangle:=\widetilde{\Pi}_{R,A}$. We define a $\varphi^a$-module $M_{I}$ over $\triangle_{I}$ for some closed interval $I\subset (0,\infty)$ having the form of $[s,r]$ with $0<s\leq r/q$ to be a finite locally free module over $\triangle^{I}$ carrying the semilinear Frobenius and the isomorphism 
\begin{displaymath}
\varphi^*M_{I}\otimes_{\triangle^{[s/q,r/q]}}\triangle^{[s,r/q]}\overset{\sim}{\rightarrow}M_{I}\otimes_{\triangle^{[s,r]}}\triangle^{[s,r/q]}.
\end{displaymath}

Then under this restriction on the interval we define a $\varphi^a$-bundle to be a collection $\{M_{I}\}_{I}$ of finite locally free modules $M_{I}$ over each $\triangle^{I}$ such that for each pair of intervals such that $I\subset I'$ we have an isomorphism $\psi_{I',I}:M_{I'}\otimes \triangle^{I}\overset{\sim}{\rightarrow} M_{I}$ and for a triple of intervals such that $I\subset I'\subset I''$ we have the obvious cocycle condition $\psi_{I',I}\circ\psi_{I'',I'}=\psi_{I'',I}$. Note that we have only finished the definitions under the restriction on the intervals. Now as in \cite[Section 6.1]{2KL15} we extend this to any more general interval $I$ by taking any pair $I\subset I'$ of  intervals and extracting some module $M_{I'}$ where $I'$ satisfies the restriction in our previous definition. Then also we extend the definition of $\varphi^a$-bundles to all the general intervals in our setting. Finally we define global section of a $\varphi^a$-bundle to be a collection of elements having the general form $\{v_{I}\}_{I}$ such that for each pair $I\subset I'$ we have that the image of $v_{I'}$ under the image of $\psi_{I',I}$ is $v_{I}$.
\end{definition}

\indent Now as in the usual situation it is natural to now consider the finiteness issue as in the following:

\begin{proposition}\mbox{\bf (After Kedlaya-Liu \cite[Lemma 6.1.4]{2KL15})} \label{proposition2.8}
Consider now an arbitrary $\varphi^a$-bundle $\{M_{I}\}$ over $\widetilde{\Pi}_{R,A}$. Suppose now that for any interval $I\subset (0,\infty)$  we have that there exists a finite number generating elements $\mathbf{e}_1,...,\mathbf{e}_n$ for some uniform $n\in \mathbb{Z}$ for the finitely locally free module $M_{I}$. Then this set of elements could be promoted to be a generating set of the global section as a module but now over $\widetilde{\Pi}^\infty_{R,A}$.
\end{proposition}

\begin{proof}
This is generalization of the corresponding results in \cite[Lemma 6.1.4]{2KL15} in the usual Hodge structures. Basically the main issue in the proof is to exhibit the desired finitely generatedness by some approximating process as in the usual situation. As in the usual situation we consider first the morphism $M_{[p^{-\ell},p^\ell]}\rightarrow (\widetilde{\Pi}^{[p^{-\ell},p^\ell]}_{R,A})^n$ by using the corresponding basis as in the assumption whose composite with the inverse $(\widetilde{\Pi}^{[p^{-\ell},p^\ell]}_{R,A})^n\rightarrow M_{[p^{-\ell},p^\ell]}$ is the identity map. And one could then bound the gap between the subspace norm and the quotient norm of $M_{[p^{-\ell},p^\ell]^n}$ by some constant $e_\ell$ as in \cite[Lemma 6.1.4]{2KL15}. Now let $\mathbf{v}\in M$. Then the point is to extract for each $i=1,...,n$ and $\ell=0,1,...$ the following element $B_{i,\ell}\in \widetilde{\Pi}^{[p^{-\ell},p^\ell]}_{R,A}$ and $A_{i,\ell}\in \widetilde{\Pi}^\infty_{R,A}$ by the following fashion following \cite[Lemma 6.1.4]{2KL15}, first for any $j<\ell$ we choose the corresponding $B_{i,j},j<\ell$ such that $\mathbf{v}-\sum_{i=1}^n\sum_{j\leq \ell}B_{i,j}\mathbf{e}_{i}=\sum_{i=1}A_{i,\ell}\mathbf{e}_i$ after which we choose $B_{i,\ell}$ to ensure that $\left\|.\right\|_{\alpha^{t},A}(B_{i,\ell}-A_{i,\ell})\leq p^{-1}e_\ell^{-1}\left\|.\right\|_{\alpha^{t},A}(A_{i,\ell})$ for each $i=1,...,n$ with $t \in [p^{-\ell},p^\ell]$. Then one could then finish as in \cite[Lemma 6.1.4]{2KL15} since by the suitable choices we will have desired convergence which gives rise to a desired expression of $\mathbf{v}$ in terms of the basis given.

\end{proof}

\begin{lemma}
Consider the corresponding projection map from a $\varphi^a$-modules $M$ over $\widetilde{\Pi}_{R,A}$ to the corresponding module $M_{I}$ over $\widetilde{\Pi}_{R,A}^{I}$ for some interval $I$. Then we have that this operation is a tensor equivalence. (Here we assume that the ring $\widetilde{\Pi}_{R,A}^{I}$ is sheafy.) 	
\end{lemma}

\begin{proof}
See \cite[Lemma 6.1.5]{2KL15}. 
\end{proof}

\indent Then we study the Frobenius invariance in our setting which generalizes the absolute situation where we do not have any deformation at all. Here we are going to use the notation $M(n)$ to denote the corresponding twist of $M$ where the twist is defined to be the one where $\varphi$ acts by $p^{-n}$.

\begin{proposition}\label{2prop2.9} \mbox{\bf (After Kedlaya-Liu \cite[Proposition 6.2.2]{2KL15})} 
Now we consider a $\varphi^a$-bundle $M$ over $\widetilde{\Pi}_{R,A}$. Then we have that there exists some integer $N\geq 1$ such that for $n\geq N$ we have that $\varphi^a-1:M_{[s,rq]}(n)\rightarrow M_{[s,r]}(n)$ is surjective. One could also have the chance to take the integer to be $1$ if the module could be derived from some module defined over $\widetilde{\Pi}_{R,A}^\mathrm{int}$. 	
\end{proposition}

\begin{proof}
Following \cite[Proposition 6.2.2]{2KL15}, we consider the quotient $r/s\leq q^{1/2}$. Then by considering the Frobenius actions one could further assume that $r\in [1,q]$. Then we will use the notation $A_{i,j}$ and $B_{i,j}$ to denote the corresponding matrices for the actions of $\varphi^{-a}$ and $\varphi^a$ respectively. Here $A_{i,j}$ has entries in $\widetilde{\Pi}^{[s,rq]}_{R,A}$ while the latter $B_{i,j}$ has entries in $\widetilde{\Pi}^{[s/q,r]}_{R,A}$. Then we put:
\begin{align}
c_1=\sup\{\left\|.\right\|_{\alpha^{t},A}(A_{i,j})|t\in [s,rq]\},\\
c_2=\sup\{\left\|.\right\|_{\alpha^{t},A}(B_{i,j})|t\in [s/q,r]\}.	
\end{align}
The goal is to extract from any element $\mathbf{w}$ in the target module a preimage $\mathbf{v}$ which amounts to solving the following equation:
\begin{displaymath}
\varphi^a\mathbf{v}-\mathbf{v}=\mathbf{w}.	
\end{displaymath}
Strictly speaking here we need some twisted version of this equation. By using our new $c_1$ and $c_2$ one chooses a suitable integer $N$ which shares the same property as the one in \cite[Proposition 6.2.2]{2KL15}. Then we adapt the argument therein to our $A$-relative setting. We start by setting the coordinate of $\mathbf{w}$ to be $(x_1,...,x_r)$ where $r$ is the rank. First we take suitable constant $c>0$ such that we have for all $n\geq N$:
\begin{displaymath}
\varepsilon:=\sup\{p^{-n}c_1c^{-(q-1)r/q},p^{n}c_2c^{(q-1)s}\}<1.	
\end{displaymath}
Then as in \cite[Proposition 6.2.2]{2KL15} one seeks suitable decomposition of $x_i$ taking the form of $x_i:=y_i+z_i$ where $y_i\in \widetilde{\Pi}_{R,A}^{[s/q,r]}$ and $z_i\in \widetilde{\Pi}_{R,A}^{[s,rq]}$ such that for the well-located $t$ as \cite[Lemma 5.2.9]{2KL15} with the estimates therein through the relative version of the equalities in \cite[Lemma 5.2.9]{2KL15}:
\begin{align}
\left\|.\right\|_{\alpha^t,A}(\varphi^{-a}(y_i))\leq c^{-(q-1)t/q}\left\|.\right\|_{\alpha^t,A}(x_i),\\
\left\|.\right\|_{\alpha^t,A}(\varphi^{a}(z_i))\leq c^{(q-1)t}\left\|.\right\|_{\alpha^t,A}(x_i),\\	
\end{align}
with the corresponding:
\begin{displaymath}
\left\|.\right\|_{\alpha^t,A}(y_i),	\left\|.\right\|_{\alpha^t,A}(z_i) \leq \left\|.\right\|_{\alpha^t,A}(x_i).
\end{displaymath}

(Note that this is actually $A$-relative version of the \cite[Lemma 5.2.9]{2KL15}, where one considers the product norm $\left\|.\right\|_{\alpha^{t},A}$.) 

Then we set as in \cite[Proposition 6.2.2]{2KL15} and compute:
\begin{align}
\sum_{i=1}^rx'_i\mathbf{f}_i &:=p^n\varphi^{-a}(\sum_{i=1}^ry_i\mathbf{f}_i)+p^{-n}\varphi^{a}(\sum_{i=1}^rz_i\mathbf{f}_i),	
\end{align}

which could also be deduced from the kind of definition by using the matrices $A_{i,j}$ and $B_{i,j}$. As in \cite[Proposition 6.2.2]{2KL15} we further consider the corresponding coordinate vectors $y,z,x'$ as functions of $x$. And set the following iteration $x_{(\ell+1)}:=x'(x_{(\ell)})$ for each $\ell\geq 0$. To extract the desired $\mathbf{v}$ one just sets now:
\begin{displaymath}
\mathbf{v}:=\sum_{\ell\geq 0}-p^{n}\varphi^{-a}(\sum_{i=1}^ry_i(x_{(\ell)})\mathbf{f}_i)+\sum_{i=1}^rz_i(x_{(\ell)})\mathbf{f}_i.	
\end{displaymath}
As in the usual situation and by our construction above we have that this element will converge to an desired element in the domain module since we have 
\begin{center}
$\mathbf{v}-p^{-n}\varphi^{a}\mathbf{v}=\sum_{i=1}^rx_{(0),i}\mathbf{f}_i$.
\end{center} 
Finally one could then follow the strategy in \cite[Proposition 6.2.2]{2KL15} to tackle the situation in more general setting by taking some suitable
\begin{center}
 $t=\max\{r/q^{1/2},s\}$ 
\end{center}
for those radii which do not satisfy the condition as above and extract suitable preimage $\mathbf{v}$ and then show that one can put this to be the desired element. To be more precise one considers the corresponding radius $t=\max\{r/q^{1/2},s\}$ where the corresponding radii $s,r$ may not satisfy the corresponding conditions as above. Then we have $t\geq r/q^{1/2}$ which implies that $r/t\leq q^{1/2}$ in this situation, which further implies (by the previous situation) one can extract a suitable element $\mathbf{v}_0$. To finish we still need to show that this element $\mathbf{v}_0$ is a desired preimage in our situation. Actually what we know in our situation is that the corresponding element $\mathbf{v}_0$ lives in the corresponding section $M_{[t,rq]}$, where this gives rise to the solution of the equation:
\begin{center}
$\mathbf{v}-p^{-n}\varphi^{a}\mathbf{v}=\sum_{i=1}^rx_{(0),i}\mathbf{f}_i$,
\end{center}
which further gives rise to the corresponding element in the corresponding section:
\begin{displaymath}
M_{[t/q,r]}.	
\end{displaymath}
Suppose now we a priori have the corresponding situation that:
\begin{displaymath}
t/q>s,	
\end{displaymath}
then we take the corresponding intersection in this situation with respect to the pair:
\begin{displaymath}
[t/q,r],[t,rq]	
\end{displaymath}
which gives rise to the fact that:
\begin{displaymath}
v_0\in M_{[t/q,rq]}	
\end{displaymath}
where if we have $t/q=s$ then we are done, otherwise we repeat the construction above again to produce the corresponding element gives rise to the situation:
\begin{align}
v_0\in M_{[t/q^2,rq]},\\
v_0\in M_{[t/q^3,rq]},\\
...,\\
v_0\in M_{[t/q^n,rq]},...,	
\end{align}
until there is some integer $N\geq 1$ such that we have $t/q^N\leq s$ then we are done.

\end{proof}

\begin{lemma}
For $0<s\leq r$, we consider the Robba ring $\Pi^{[s,r]}_{R,A}$. Any multiplicative seminorm of this ring is bounded by $\max\{\|.\|_{\alpha^s,A},\|.\|_{\alpha^s,A}\}$	if and only if this seminorm is bounded by the norm $\|.\|_{\alpha^t,A}$ for some unique $t\in [s,r]$ coming from \cite[Lemma 5.3.4]{2KL15}, by restricting the seminorm to the ring $\Pi^{[s,r]}_{R}$.
\end{lemma}

\begin{proof}
As in \cite[Lemma 5.3.4]{2KL15}, one direction is trivial by the property of the norm. So suppose now that we have that the seminorm is bounded by $\max\{\|.\|_{\alpha^s,A},\|.\|_{\alpha^s,A}\}$. Suppose we denote the seminorm by $|.|$, then we have that $|.|([\overline{r}]\otimes1)$ for any $\overline{r}\in R$ is bounded by $\|.\|_{\alpha^t,A}([\overline{r}])$ which could proved in the same way as in \cite[Lemma 5.3.4]{2KL15} (where the number $t$ is determined uniquely by \cite[Lemma 5.3.4]{2KL15}), while for any $a\in A$ we have $|.|(1\otimes a)$ is bounded by $\max\{\|.\|_{\alpha^s,A},\|.\|_{\alpha^s,A}\}(1\otimes a)$ which is just $\|.\|_A(a)$. Therefore we have $|.|([\overline{r}]\otimes 1)|.|(1\otimes a)=|.|([\overline{r}]\otimes a)$ is then bounded by the product of $\|.\|_{\alpha^t,A}([\overline{r}])$ and $\|.\|_A(a)$, which gives the result.
\end{proof}

\begin{proposition} \label{proposition1} \mbox{\bf (After Kedlaya-Liu \cite[Proposition 6.2.4]{2KL15})} 
Now again\\
 working over the ring $\widetilde{\Pi}_{R,A}$ and consider a $\varphi^a$-bundle which is denoted by $M$. Then we have that there exists some positive integer $N$ such that for all $n\geq N$ one could find finitely many $\varphi^a$-invariant global sections of $M(n)$ that generate $M$.	
\end{proposition}

\begin{proof}
We follow \cite[Proposition 6.2.4]{2KL15}. In our generalized setting, first we take $s$ to be $rq^{-1/2}$. Then as in the situation of \cite[Proposition 6.2.4]{2KL15} one chooses an element $\overline{\pi}$ whose norm is inside $(0,1)$. Then we consider the corresponding product which we will denote it by $[\overline{\pi}]$, and choose suitable rational number $z\in \mathbb{Z}[1/p]$ such that $[\overline{\pi}^{z}]$ gives rise to the value $c$ in the proof of the previous proposition. Then also as in the previous proposition one constructs (using $\varphi^a$ instead) the corresponding element $\mathbf{v}$ but in our situation we will denote this by $\mathbf{v}'_i$ for each $i$ suitable. To be more precise from the data above we put $\mathbf{v}'_i$ as the definition of $\mathbf{v}$ by using $x_{(0)}=0$ and $(y_{(0),j},z_{(0),j})=(-[\overline{\pi}^{z}],[\overline{\pi}^{z}])$ for $j=i$ and zero identically otherwise.\\
\indent Then we have the corresponding set of elements $\{\mathbf{v}'_i\}$ in $M_{[rq^{-1/2},rq]}$ consisting those elements which are invariant under the action of $\varphi^a$. We repeat the argument in the proof of the previous proposition as in the following. First one uses the notation $X_{i,j}$ to denote the corresponding matrix of the operator $\varphi^{-a}$ (certainly under some chosen set of generators with the specific rank $r$ in our situation) while one uses the notation $Y_{i,j}$ to denote the corresponding matrix of the operator $\varphi^a$, where from the corresponding matrix elements one could compute the following two numbers:
\begin{align}
b_1:=\mathrm{sup}_t\{\left\|.\right\|_{\alpha^t,A}(X_{i,j}),i,j\in \{1,2,...,r\},t\in [s,rq]\}\\
b_2:=\mathrm{sup}_t\{\left\|.\right\|_{\alpha^t,A}(Y_{i,j}),i,j\in \{1,2,...,r\},t\in [s/q,r]\}.	
\end{align}
Now by using our new scalars $b_1,b_2$ one could have the chance to extract some suitable $N\geq 1$ as in our previous proposition and in the corresponding result of \cite[Proposition 6.2.2]{2KL15}. Then we consider the following construction through the corresponding iterated induction as in the following. The corresponding initial coordinates have been constructed in our situation as mentioned in the previous paragraph. Then we consider some general coordinate $(x_1,...,x_r)$ for $r$ the rank. Then as in \cite[Lemma 5.2.9]{2KL15} for each $i=1,...,r$ one can make the following decomposition, namely putting $x_i=y_i+z_i$ with the following estimate with some chosen constant $c>0$:
\begin{align}
\left\|.\right\|_{\alpha^t,A}(\varphi^{-a}(y_i))\leq c^{-(q-1)t/q}\left\|.\right\|_{\alpha^t,A}(x_i),\\
\left\|.\right\|_{\alpha^t,A}(\varphi^{a}(z_i))\leq c^{(q-1)t}\left\|.\right\|_{\alpha^t,A}(x_i),\\	
\end{align}
and:
\begin{displaymath}
\left\|.\right\|_{\alpha^t,A}(y_i),	\left\|.\right\|_{\alpha^t,A}(z_i) \leq \left\|.\right\|_{\alpha^t,A}(x_i).
\end{displaymath}
This will be the key step in the iterating process, while the scalar $c$ will be chosen to guarantee that we have the following situation:
\begin{displaymath}
\delta:=\mathrm{max}\{p^{-n}b_1c^{-(q-1)r/q},p^{n}b_2c^{(q-1)s}\}<1.
\end{displaymath}
Then the corresponding iterating process will be further conducted in the following way, namely we consider the corresponding (well-defined due to our construction above):
\begin{displaymath}
\sum_{i=1}^rx_i'\mathbf{f}_i:=p^{n}\varphi^{-a}	\sum_{i=1}^ry_i\mathbf{f}_i+p^{-n}\varphi^{a}	\sum_{i=1}^rz_i\mathbf{f}_i.
\end{displaymath}
Note here that from each $x_i,i=1,...,r$ one can produce the corresponding elements $x_i',i=1,...,r$ which further produces the corresponding elements $y_i,z_i,i=1,...,r$, which implies that one can regard the corresponding elements $y_i,z_i,i=1,...,r$ as functions $y_i(x_1,...,x_r)$, $z_i(x_1,...,x_r),i=1,...,r$ of the original elements $x_i,i=1,...,r$. The whole iterating process is conducted through putting $x_{(\ell+1)}:=x'(x_{(\ell)})$ for each $\ell=0,1,...$ (note that in this situation $x_{(1)}$ is constructed directly from $x_{(0)}$), which produces the desired new elements by considering the following series:
\begin{displaymath}
\mathbf{v}'_i:=\sum_{\ell\geq 0}-p^{n}\varphi^{-a}	(\sum_{i=1}^ry_i(x_{(\ell)})\mathbf{f}_i) +\sum_{i=1}^rz_i(x_{(\ell)})\mathbf{f}_i.	
\end{displaymath}
As in our previous situation this process produces the desired elements which are invariant under the Frobenius operators in the sense that:
\begin{displaymath}
\mathbf{v}'_i-p^{-n}\varphi^a\mathbf{v}'_i=\sum_{i=0}^rx_{(0),i}\mathbf{f}_i.	
\end{displaymath}
Then we write $\mathbf{v}'_i$ as $[\overline{\pi}^{s}]\mathbf{f}_i+M\mathbf{f}$ with some matrix $\{M_{i,j}\}$ satisfying the condition that $\left\|.\right\|_{\alpha^{t},A}(M_{i,j})\leq \varepsilon \alpha(\overline{\pi})^{st}$ for each $t\in [rq^{-1/2},r]$. Therefore up to here we conclude that the matrix $[\overline{\pi}^s]+M$ will then consequently be invertible over the period ring $\widetilde{\Pi}_{R,A}^{[rq^{-1/2},r]}$ which implies the corresponding finite generating property for the interval $[rq^{-1},r]$ (by further repeating the corresponding construction above for the interval $[rq^{-1},rq^{-1/2}]$). Then to finish one could then iterate as in \cite[Proposition 6.2.4]{2KL15} and apply the Frobenius action from $\varphi^a$. Note that this will rely on the analog in our situation of \cite[Lemma 5.3.4]{2KL15} in the $A$-relative setting (in our setting this could be resolved by considering the corresponding seminorms on the period rings in the usual sense and $A$, as in the previous lemma).
\end{proof}

\indent Then we consider the Frobenius invariants acting on an exact sequence, generalizing the corresponding result from \cite[Corollary 6.2.3]{2KL15}:

\begin{corollary} \label{Corol2.12}
Suppose now we have an exact sequence of $\varphi^a$-modules over the corresponding period ring $\widetilde{\Pi}_{R,A}$ taking the form of $0\rightarrow M_\alpha\rightarrow M \rightarrow M_\beta\rightarrow 0$. Then there exists a positive integer $N$ such that for all $n\geq N$ we have the following exact sequence:
\[
\xymatrix@R+0pc@C+0pc{
0\ar[r]\ar[r]\ar[r] &M_\alpha(n)^{\varphi^a=1} \ar[r]\ar[r]\ar[r] &M(n)^{\varphi^a=1}
\ar[r]\ar[r]\ar[r] &M_\beta(n)^{\varphi^a=1}
\ar[r]\ar[r]\ar[r] &0
}.
\]
	
\end{corollary}

\begin{proof}
This is a direct consequence of the previous \cref{2prop2.9}.	
\end{proof}

\indent The finite projective objects are considered in the results above, actually one can naturally extend the discussion above to the setting of pseudocoherent objects. Following \cite[4.6.9]{2KL16} we consider the following results around the finite generated objects:

\begin{proposition} \label{propo2.13}
Let $M$ be a Frobenius module defined over $\widetilde{\Pi}_{R,A}$ (namely endowed with a semilinear action of the Frobenius operator) which is assumed to be finitely generated. Then we have that one can find then an integer $N\geq 0$ such that for all $n\geq N$, $H^0_{\varphi^a}(M(n))$ generates the module $M$ itself and we have that the $H^1_{\varphi^a}(M(n))$ vanishes. Also let $M$ be a Frobenius module defined over $\widetilde{\Pi}_{R,A}$ (namely endowed with a semilinear action of the Frobenius operator) which is assumed to be finite projective. Then we have that one can find then an integer $N\geq 0$ such that for all $n\geq N$, $H^0_{\varphi^a}(M(n))$ generates the module $M$ itself and we have that the $H^1_{\varphi^a}(M(n))$ vanishes. 
\end{proposition}

\begin{proof}
In this situation one can follow the proof of \cite[4.6.9]{2KL16} and the strategy in the proof of the previous two propositions.	
\end{proof}

%%%%%%\newpage

\newpage\section{Contact with Schematic Relative Fargues-Fontaine Curves}

\subsection{Rigid Analytic Deformation of Schematic Relative Fargues-Fontaine curves}

\indent In this section, we study some relationship between the vector bundles over Fargues-Fontaine curves and the generalized Frobenius modules defined in the previous section. This is relative version of the established results essentially in \cite{2KL15}. The original picture could be dated back to the work of Fargues-Fontaine, where they established the classification of the Galois equivariant vector bundles over the so-called Fargues-Fontaine curves. For the convenience of the reader we first recall the definition of the scheme $\mathrm{Proj}P$ in the current picture.

\begin{setting}
Recall some algebraic geometry from \cite[Definition 6.3.1]{2KL15} that first the ring $P$ is a graded commutative ring taking the form of $\bigoplus_{n\geq 0}P_n$ where each $P_{n}$ is defined to be a subring of $\widetilde{\Pi}^+_{R}$ 
%(here $R$ is compatible with our notation in the sense that here $Y$ is just a singleton) 
which consists of all the elements which is invariant under the Frobenius $\varphi^a$ (up to scalars $p^n$ for each $n$). Then for each element $f$ of degree $d>0$ in $P_d$ we consider the local affine scheme $\mathrm{Spec}(P[1/f]_0)$, then we glue these to get the projective spectrum of $P$ as in \cite[Definition 6.3.1]{2KL15} which is denoted by $\mathrm{Proj}P$. For more algebraic geometric discussion see \cite[Definition 6.3.1]{2KL15}. 
\end{setting}

\indent Then we could generalize the situation in \cite[Definition 6.3.1]{2KL15} to our $A$-relative situation.

\begin{setting}
In our situation, we consider the ring $\widetilde{\Pi}_{R,A}^+$ or the ring $\widetilde{\Pi}_{R,A}$ which is used to defined the corresponding graded ring as in the previous setting (namely, taking the corresponding $\varphi^a=p^n$ invariants for each $n\geq 0$), which will be denoted by $P_{R,A}$ or even just $P_A$.	To be slightly more precise $P_{R,A}$ is defined to be the corresponding direct sum $\bigoplus_{n\geq 0}P_{n,A}$.
\end{setting}

\begin{setting}
First we are going to use the notation $f$ to denote both the corresponding elements with specific degrees $d$ in $P_{d}$ and the corresponding elements in the tensor product of the local affine rings. And then we consider now the following invariance for any $M$, which is any $\varphi^a$-bundle over $\widetilde{\Pi}_{R,A}$. The invariance $M_{f}$ is now defined in our setting as $M[1/f]^{\varphi^a}$. To be more explicit we have $M_{f}=\bigcup_{n\geq 0}f^{-n}M(dn)^{\varphi^a}$. As in the usual situation one could regard this as a module over the tensor product $P_A[1/f]_0$, which then for any $I$ an interval induces the following map:

\begin{displaymath}
M_{f}\otimes_{P_A[1/f]_0}\widetilde{\Pi}^{I}_{R,A}[1/f]\rightarrow M_{f}\otimes_{\widetilde{\Pi}^{I}_{R,A}}\widetilde{\Pi}^{I}_{R,A}[1/f].	
\end{displaymath}
\end{setting}

\indent To study more, we first consider the following corollary after Kedlaya-Liu:

\begin{corollary}
Consider an arbitrary exact sequence of $\varphi^a$-modules taking the form of
\[
\xymatrix@R+0pc@C+0pc{
0\ar[r]\ar[r]\ar[r] &M_A \ar[r]\ar[r]\ar[r] &M
\ar[r]\ar[r]\ar[r] &M_B
\ar[r]\ar[r]\ar[r] &0.
}
\]
Then we have the derived corresponding exact sequence:
\[
\xymatrix@R+0pc@C+0pc{
0\ar[r]\ar[r]\ar[r] &M_{A,f} \ar[r]\ar[r]\ar[r] &M_{f}
\ar[r]\ar[r]\ar[r] &M_{B,f}
\ar[r]\ar[r]\ar[r] &0.
}
\]
			
\end{corollary}

\begin{proof}
Just consider the structure of $M_{f}$, then one could derive the results from our previous discussion.	
\end{proof}

\indent Then we have the following relative version of \cite[Corollary 6.3.4]{2KL15}:

\begin{lemma}
Let $M$ be a Frobenius bundle as above, and take $M_1$ and $M_2$ two Frobenius subbundles whose summation is also a subobject. Then we actually have an equality $(M_1+M_2)_f\overset{\sim}{\rightarrow} M_{1,f}+M_{2,f}$.
\end{lemma}

\begin{proof}
This is basically an $A$-relative version of the corresponding result in \cite[Corollary 6.3.4]{2KL15}.	
\end{proof}

\indent Then we consider the following key proposition for our further study:

\begin{proposition} \mbox{\bf (After Kedlaya-Liu \cite[Theorem 6.3.9]{2KL15})} 
Now let $M$ again be a $\varphi^a$-bundle over $\widetilde{\Pi}_{R,A}$. Then we have that
\begin{displaymath}
M_{f}\otimes_{P_A[1/f]_0}\widetilde{\Pi}^{I}_{R,A}[1/f]\rightarrow M_{f}\otimes_{\widetilde{\Pi}^{I}_{R,A}}\widetilde{\Pi}^{I}_{R,A}[1/f]	
\end{displaymath}  
is a bijection, and $M_{f}$ is finite projective over $P_A[1/f]_0$.
\end{proposition}

\begin{proof}
We follow \cite[Theorem 6.3.9]{2KL15} which gives basically from \cref{proposition1} that there exists $\mathbf{e}_1,...,\mathbf{e}_k$ which are the global sections of $M(dn)$ for some $n$ which generate $M$ itself. Then we have that considering $f^{-n}\mathbf{e}_1,...,f^{-n}\mathbf{e}_k$ could give rise to the corresponding surjectivity as in \cite[Theorem 6.3.9]{2KL15}. Then we consider the corresponding presentation by using these generating elements from $\widetilde{\Pi}_{R,A}(-dn)^k(dn)$ mapping to $M(dn)$ which gives rise to a mapping which is surjective mapping from $M_1\rightarrow M$ (certainly here we have to consider the corresponding bundle in our context corresponding to the module $\widetilde{\Pi}_{R,A}(-dn)^k$ and then twist to establish the corresponding map as above), which furthermore gives rise to (by the corresponding observation parallel to \cite[Theorem 6.3.9]{2KL15}) an exact sequence taking the form of
\[
\xymatrix@R+0pc@C+0pc{
0\ar[r]\ar[r]\ar[r] &M_2 \ar[r]\ar[r]\ar[r] &M_1
\ar[r]\ar[r]\ar[r] &M
\ar[r]\ar[r]\ar[r] &0
}.
\]
Then we consider the corresponding induced exact sequence taking the form of:
\[
\xymatrix@R+0pc@C+0pc{
0\ar[r]\ar[r]\ar[r] &M_{2,f} \ar[r]\ar[r]\ar[r] &M_{1,f}
\ar[r]\ar[r]\ar[r] &M_{f}
\ar[r]\ar[r]\ar[r] &0
},
\]
from the discussion above. Then to finish we look at the corresponding commutative diagram as in \cite[Theorem 6.3.9]{2KL15} in the following:
\[\tiny
\xymatrix@R+5pc@C+0.7pc{
 &M_{2,f}\otimes_{P_A[1/f]_0}\widetilde{\Pi}^{I}_{R,A}[1/f] \ar[d]\ar[d]\ar[d] \ar[r]\ar[r]\ar[r] &M_{1,f}\otimes_{P_A[1/f]_0}\widetilde{\Pi}^{I}_{R,A}[1/f] \ar[d]\ar[d]\ar[d]
\ar[r]\ar[r]\ar[r] &M_{f}\otimes_{P_A[1/f]_0}\widetilde{\Pi}^{I}_{R,A}[1/f] \ar[d]\ar[d]\ar[d]
\ar[r]\ar[r]\ar[r] &0\\
0\ar[r]\ar[r]\ar[r] &M_{2,f}\otimes_{\widetilde{\Pi}^{I}_{R,A}}\widetilde{\Pi}^{I}_{R,A}[1/f]	 \ar[r]\ar[r]\ar[r] &M_{1,f}\otimes_{\widetilde{\Pi}^{I}_{R,A}}\widetilde{\Pi}^{I}_{R,A}[1/f]	
\ar[r]\ar[r]\ar[r] &M_{f}\otimes_{\widetilde{\Pi}^{I}_{R,A}}\widetilde{\Pi}^{I}_{R,A}[1/f]	
\ar[r]\ar[r]\ar[r] &0
},
\]
which finishes the proof on the bijectivity as in the absolute situation in \cite[Theorem 6.3.9]{2KL15}. The finite projectiveness could be prove exactly by the same fashion in \cite[Theorem 6.3.9]{2KL15} (also one can similarly derive the corresponding finitely-presentedness). To be more precise one can look at the exact sequence constructed above:
\[
\xymatrix@R+0pc@C+0pc{
0\ar[r]\ar[r]\ar[r] &M_2 \ar[r]\ar[r]\ar[r] &M_1
\ar[r]\ar[r]\ar[r] &M
\ar[r]\ar[r]\ar[r] &0
}.
\]
Then by taking suitable higher degree twist in our situation (by considering the corresponding vanishing of the cohomology $\mathrm{H}^1(M^\vee\otimes M_2(dN))$) one can get the following push-out diagram just as in \cite[Theorem 6.3.9]{2KL15}:
\[
\xymatrix@R+2pc@C+2pc{
 M_2 \ar[d]\ar[d]\ar[d] \ar[r]\ar[r]\ar[r] &M_1 \ar[d]\ar[d]\ar[d]
\ar[r]\ar[r]\ar[r] &M \ar[d]\ar[d]\ar[d]\\
f^{-N}M_2	 \ar[r]\ar[r]\ar[r] &M'_{1}
\ar[r]\ar[r]\ar[r] &M	
},
\]
which gives rise to that $M_{1,f}\overset{\sim}{\rightarrow} M'_{1,f}$. Then by our construction, since $M_{1,f}$ is free we are done.

\end{proof}

\indent We then have the following results after suitable localization.

\begin{corollary} \label{corollary7.6}
Pick the element $f$ as in the previous proposition, we have that the following four categories are equivalent: \\
I. The category of finite locally free $P_A[1/f]_0$-modules over the relative coordinate ring of the localization of the Fargues-Fontaine curve $\mathrm{Proj}P_A$;\\
II.The category of $A$-relative $\varphi^a$-modules over $\widetilde{\Pi}^\infty_{R,A}[1/f]$;\\
III. The category of $A$-relative $\varphi^a$-modules over $\widetilde{\Pi}_{R,A}[1/f]$;\\
IV. The category of $A$-relative $\varphi^a$-bundles over $\widetilde{\Pi}_{R,A}[1/f]$.	
\end{corollary}

\begin{proof}
This is essentially by iteratedly applying the base change functor and the previous proposition.	
\end{proof}

%\begin{remark}
%\indent Then one could consider more global situation by glueing the construction in some local charts together, the globalization here is actually a bit complicated since we do not geometrize the ring $A$ as above (instead we treat this part as some external structures). The first issue is the definition of the sheaf $\mathcal{O}_{\mathrm{Spec}(P[1/f]_0)}\widehat{\otimes}A$ over again localized Fargues-Fontaine curves. We define first the corresponding presheaf in the obvious way. Then we believe that this might be directly automatically a sheaf, therefore we make the assumption that $A$ is external sheafy in the sense that the presheaves (defined in natural and obvious way):
%\begin{align}
%\mathcal{O}_{\mathrm{Spec}(P[1/f]_0)}\widehat{\otimes}A,\\
%\mathcal{O}_{\mathrm{Proj}(P)}\widehat{\otimes}A,\\
%\mathcal{O}_{\mathrm{Spec}(\widetilde{\Pi}^\infty_{R}[1/f])}\widehat{\otimes}A,\\
%\mathcal{O}_{\mathrm{Spec}(\widetilde{\Pi}^\infty_{R})}\widehat{\otimes}A	
%\end{align}
%are automatically sheaves. 	
%\end{remark}

%%%%%%%%%%%%%%%%%%%%%%%%%%%%%%%%%%revised.

\begin{proposition} \mbox{\bf (After Kedlaya-Liu \cite[Theorem 6.3.12]{2KL15})}     \label{2corollary7.6}
We have that the following four categories are equivalent: \\
I. The category of quasicoherent finite locally free sheaves of $\mathcal{O}_{\mathrm{Proj}P_A}$-modules (namely the $A$-relative vector bundles) over the Fargues-Fontaine curve $\mathrm{Proj}P_A$;\\
II.The category of $A$-relative $\varphi^a$-modules over $\widetilde{\Pi}^\infty_{R,A}$;\\
III. The category of $A$-relative $\varphi^a$-modules over $\widetilde{\Pi}_{R,A}$;\\
IV. The category of $A$-relative $\varphi^a$-bundles over $\widetilde{\Pi}_{R,A}$.	
\end{proposition}

\begin{proof}

\indent Let us mention briefly the corresponding functors involved since essentially the process is as in the same fashion of the construction of \cite[Theorem 6.3.12]{2KL15}. First from II to III, this is just the base change, and the III to IV is the functor which maps any $\varphi^a$-module to the associated $\varphi^a$-bundle. Then for the functor from I to II one considers the process which associates any vector bundle $V$ in the first category some finite locally free sheaf $\mathcal{V}$ over the corresponding space associated to $\widetilde{\Pi}^\infty_{R,A}$ by considering the localization and glueing through the consideration by $f$ mentioned above, for the glueing we apply the direct analog of \cite[Lemma 6.3.7]{2KL15} by considering \cref{proposition2.8} and \cref{proposition1}. Then one takes the global section to get desired module in the second category. Then in our situation the previous proposition gives us the final equivalence after the well-established glueing process under the consideration of the further relativization coming from the algebra $A$.
\end{proof}

\indent One can then further discuss the corresponding results on the comparison of pseudocoherent objects (just as in \cite[Definition 4.4.4]{2KL16}), as in \cite[Section 4.6]{2KL16}. The results for the pseudocoherent objects are bit more complicated than the results established above in the context of just vector bundles. We first have the following analog of \cite[Corollary 4.6.10]{2KL16}:

\begin{proposition}
Let $M_\alpha,M,M_\beta$ be three finitely generated Frobenius modules over the ring $\widetilde{\Pi}_{R,A}$ (namely endowed with a semilinear action from the Frobenius operator). And now we put the modules then in an exact sequence:
\[
\xymatrix@R+0pc@C+0pc{
0\ar[r]\ar[r]\ar[r] &M_\alpha \ar[r]\ar[r]\ar[r] &M
\ar[r]\ar[r]\ar[r] &M_\beta
\ar[r]\ar[r]\ar[r] &0.
}
\]
Then we have first the following exact sequence for sufficiently large integer $n\geq 0$:
\[
\xymatrix@R+0pc@C+0pc{
0\ar[r]\ar[r]\ar[r] &M_\alpha(n)^{\varphi^a} \ar[r]\ar[r]\ar[r] &M(n)^{\varphi^a}
\ar[r]\ar[r]\ar[r] &M_\beta(n)^{\varphi^a}
\ar[r]\ar[r]\ar[r] &0.
}
\]
Then we have for any element $f$ which is a $\varphi^a=p^\ell$ invariance, after forming the following module for each $M_*:=M_\alpha,M,M_\beta$:
\begin{displaymath}
M_{*,f}:=\bigcup_{n\in \mathbb{Z}}f^{-n}M_{*}(\ell n)^{\varphi^a},	
\end{displaymath}
the following exact sequence of Frobenius modules over $(\widetilde{\Pi}_{R,A}[1/f])^{\varphi^a}$:
\[
\xymatrix@R+0pc@C+0pc{
0\ar[r]\ar[r]\ar[r] &M_{\alpha,f} \ar[r]\ar[r]\ar[r] &M_f
\ar[r]\ar[r]\ar[r] &M_{\beta,f}
\ar[r]\ar[r]\ar[r] &0.
}
\]
Eventually we have the fact that if $M'$ is then a pseudocoherent Frobenius module over $\widetilde{\Pi}_{R,A}$ then the module $M_f$ is also a pseudocoherent module over $(\widetilde{\Pi}_{R,A}[1/f])^{\varphi^a}$.
\end{proposition}

\begin{proof}
Apply \cref{propo2.13} as in \cref{Corol2.12} to derive the first two consequences. As in \cite[Corollary 4.6.10]{2KL16} one can further prove the last statement by taking the corresponding projective resolutions and by using the previous statement repeatedly. 
\end{proof}

\begin{proposition}\mbox{\bf (After Kedlaya-Liu \cite[Theorem 4.6.12]{2KL16})} \label{2proposition2.3.10} Taking\\ pullbacks along the map from the scheme associated to $\widetilde{\Pi}^\infty_{R,A}$ (by the analog of \cite[Lemma 6.3.7]{2KL15} by using \cref{proposition1} and \cref{proposition2.8}) to the (schematic) Fargues-Fontaine curve (as in \cite[Definition 4.6.11]{2KL16}) gives rise to an exact functor which establishes the equivalence between the category of pseudocoherent sheaves over the Fargues-Fontaine curve and the category of pseudocoherent Frobenius modules over $\widetilde{\Pi}_{R,A}$ (with isomorphism by the Frobenius pullbacks). This equivalence respects the corresponding comparison on the sheaf cohomologies and the Frobenius cohomologies.
	
\end{proposition}

\begin{proof}
We follow the strategy of the proof of \cite[Theorem 4.6.12]{2KL16} to prove this. First in our context we have that the above functor is actually just initially known to be right exact under the base change. Now suppose $V$ is a pseudocoherent sheaf over the schematic Fargues-Fontaine curve in our context. As in the beginning of the proof of \cite[Theorem 4.6.12]{2KL16} one has the following exact sequence:
\[
\xymatrix@R+0pc@C+0pc{
0\ar[r]\ar[r]\ar[r] &V_1 \ar[r]\ar[r]\ar[r] &V_2
\ar[r]\ar[r]\ar[r] &V
\ar[r]\ar[r]\ar[r] &0,
}
\]
where $V_1$ is pseudocoherent and $V_2$ is assumed and set to be a vector bundle. Then we apply the pullback construction (as mentioned in the statement of the proposition) we have the corresponding (after taking the further base change) exact sequence of the pseudocoherent Frobenius modules over the ring $\widetilde{\Pi}_{R,A}$:
\[
\xymatrix@R+0pc@C+0pc{
W_1 \ar[r]\ar[r]\ar[r] &W_2
\ar[r]\ar[r]\ar[r] &W
\ar[r]\ar[r]\ar[r] &0,
}
\] 
where $W_2$ is finite projective and $W_1$ is finitely generated. We then have the following exact sequence:
\[
\xymatrix@R+0pc@C+0pc{
0 \ar[r]\ar[r]\ar[r] &\mathrm{*} \ar[r]\ar[r]\ar[r] &W_1 \ar[r]\ar[r]\ar[r] &W_2
\ar[r]\ar[r]\ar[r] &W
\ar[r]\ar[r]\ar[r] &0,
}
\]
where $*$ is the kernel of the map $W_1\rightarrow \mathrm{Kernel}(W_2\rightarrow W)$. Then choose some $f$ with some specific degree as what we did before, by the previous proposition we have the following exact sequence:
\[
\xymatrix@R+0pc@C+0pc{
0 \ar[r]\ar[r]\ar[r] &\mathrm{*}_f \ar[r]\ar[r]\ar[r] &W_{1,f} \ar[r]\ar[r]\ar[r] &W_{2,f}
\ar[r]\ar[r]\ar[r] &W_f
\ar[r]\ar[r]\ar[r] &0,
}
\]
then by taking the corresponding section we have the following commutative diagram:
\[
\xymatrix@R+5pc@C+1.2pc{
 &0 \ar[r]\ar[r]\ar[r] &V_1|_{\widetilde{\Pi}_{R,A}[1/f]^{\varphi^a}} \ar[r]\ar[r]\ar[r] \ar[d]\ar[d]\ar[d] &V_2|_{\widetilde{\Pi}_{R,A}[1/f]^{\varphi^a}} \ar[d]\ar[d]\ar[d]
\ar[r]\ar[r]\ar[r] &V|_{\widetilde{\Pi}_{R,A}[1/f]^{\varphi^a}} \ar[d]\ar[d]\ar[d]
\ar[r]\ar[r]\ar[r] &0 \\
0 \ar[r]\ar[r]\ar[r] &\mathrm{*}_f \ar[r]\ar[r]\ar[r] &W_{1,f} \ar[r]\ar[r]\ar[r] &W_{2,f}
\ar[r]\ar[r]\ar[r] &W_f
\ar[r]\ar[r]\ar[r] &0
},
\]
where by our above results on the comparison on the finite projective objects we have that the second vertical morphism is isomorphism which implies that the third vertical map is surjective. Then we repeat the corresponding argument and the construction by using $V_1$ as our $V$ we can also derive the fact that the first vertical map is also at this situation surjective. Then as in \cite[Theorem 4.6.12]{2KL16} the five lemma implies that the third vertical map is also injective. Then we apply the same construction and argument to the situation where we use $V_1$ to be our $V$ we can derive the fact that the first vertical map is also injective. Then we have that all the vertical maps are then in this situation isomorphism, which further implies that $*_f$ is trivial, as in \cite[Theorem 4.6.12]{2KL16} by using the previous proposition we have that $*$ is trivial. The functor send finite projective objects to the corresponding finite projective objects so then we have the pseudocoherent objects in the corresponding essential image. Then we have a well-defined exact functor in our situation where the quasi-inverse is just taking the corresponding Frobenius invariance. The functor from the modules over the Robba ring to the sheaves over the Fargues-Fontaine curve will have the corresponding composition with the quasi-inverse being equivalence. On the other hand, to show that the functor from the sheaves over the Fargues-Fontaine curve will have the corresponding composition with the quasi-inverse being equivalence will eventually reduce to the corresponding statement in the finite projective setting as what we did before.
\end{proof}

\subsection{Fr\'echet-Stein Deformation of Schematic Relative Fargues-Fontaine curves} \label{2section3.2}

\indent For our purpose and motivation, we would like to consider the corresponding situation where $A$ is replaced by the distribution algebras attached to some specific $p$-adic Lie groups or some mixed-type algebra by taking the product of $A$ with these interesting distribution algebras. Let us recall the following construction:

\begin{setting} \label{setting3.8}
Let $G$ be a $p$-adic Lie group in the fashion picked in the next setting, then we will use the notation $\mathcal{O}_K[[G]]$ to denote the corresponding integral completed group algebra over some finite extension $K/\mathbb{Q}_p$. Then we use the notation $K[[G]]$ to denote the base change of the integral ring to $K$. From here we consider the following integral version of the distribution we are considering:
\begin{displaymath}
A_n:=\mathcal{O}_K[[G]][\mathfrak{m}^n/p]^\wedge_{(p)}[1/p]	
\end{displaymath}
then take the inverse limit we get the ring $A_\infty$, where $\mathfrak{m}$ is the Jacobson radical. However this is not Noetherian, which makes life complicated. It is expected that one should consider this ring in order to do the equivariant or even more general noncommutative Iwasawa theory and the corresponding Tamagawa conjectures, in stead of just the usual Iwasawa algebra. For instance in the situation where $G$ is just $\mathbb{Z}_p^\times$ and $K=\mathbb{Q}_p$ this is just the ring $\Pi^\infty(\Gamma_K)$ which is just the Robba ring corresponding to the radius $\infty$. Also see \cite{2Zah1} where this is discussed in more detail.
\end{setting}

\begin{setting}
In our current context we are going to focus on those groups in the following form:
\begin{displaymath}
G=\mathbb{Z}_p^d,\mathbb{Z}_p^\times,\mathbb{Z}_p^\times \ltimes \mathbb{Z}_p^n.	
\end{displaymath}
In this situation we just define the corresponding period rings deformed over the distribution algebra $A_\infty(G)$:
\begin{align}
&\widetilde{\Omega}^\mathrm{int}_{R,A_\infty(G)}:=\varprojlim_{n\rightarrow \infty}\widetilde{\Omega}^\mathrm{int}_{R,A^\mathrm{int}_n(G)},\widetilde{\Omega}_{R,A_\infty(G)}:=\varprojlim_{n\rightarrow \infty}\widetilde{\Omega}_{R,A_n(G)},\\
&\widetilde{\Pi}^{\mathrm{int},r}_{R,A_\infty(G)}:=\varprojlim_{n\rightarrow \infty}\widetilde{\Pi}^{\mathrm{int},r}_{R,A^\mathrm{int}_n(G)},\widetilde{\Pi}^\mathrm{int}_{R,A_\infty(G)}:=\varprojlim_{n\rightarrow \infty}\widetilde{\Pi}^\mathrm{int}_{R,A^\mathrm{int}_n(G)},\\	
&\widetilde{\Pi}^{\mathrm{bd},r}_{R,A_\infty(G)}:=\varprojlim_{n\rightarrow \infty}\widetilde{\Pi}^{\mathrm{bd},r}_{R,A_n(G)},\widetilde{\Pi}^\mathrm{bd}_{R,A_\infty(G)}:=\varprojlim_{n\rightarrow \infty}\widetilde{\Pi}^\mathrm{bd}_{R,A_n(G)},\\
&\widetilde{\Pi}^I_{R,A_\infty(G)}:=\varprojlim_{n\rightarrow \infty}\widetilde{\Pi}^I_{R,A_n(G)},	
\end{align}
from which one can define furthermore the period rings 
\begin{align}
\widetilde{\Pi}^r_{R,A_\infty(G)},\widetilde{\Pi}^\infty_{R,A_\infty(G)},\widetilde{\Pi}_{R,A_\infty(G)}.	
\end{align}
Then we work over the corresponding finite locally free modules or bundles over the corresponding period rings defined above with the same definitions we used in \cref{Def2.6} to define those modules over the period rings associated to $A_\infty(G)$. We also endow the rings defined above with partial Frobenius where the action is induced form the Witt vector part, not the coefficient algebra part. We then define the corresponding families of $A_\infty(G)$-relative $\varphi^a$-module $M$ over the ring $\widetilde{\Pi}_{R,A_\infty(G)}^\infty$ or $\widetilde{\Pi}_{R,A_\infty(G)}$ to be a projective system $(M_n)_{n}$ where each $M_n$ for each $n$ is a corresponding $A_{n}(G)$-relative $\varphi^a$-modules over $\widetilde{\Pi}_{R,A_n(G)}^\infty$ or $\widetilde{\Pi}_{R,A_n(G)}$ such that we have the base change requirement for each $n$:
\begin{displaymath}
A_{n}(G)\otimes_{A_{n+1}(G)}M_{n+1}\overset{\sim}{\rightarrow}M_{n}.	
\end{displaymath}
Similarly we define the Frobenius bundles in this compatible way. And finally we define the corresponding global section of any families of the Frobenius modules over the corresponding period rings to be the corresponding projective limit taking the form of $\varprojlim_{n\rightarrow\infty}M_n$.
%Therefore by considering all the above proof but replacing $A$ with more general Banach algebras such as the ring $A_\infty(G)$ (with the restriction on the group $G$), we have the following.
\end{setting}

\begin{remark}
One can safely extend the assumption on the group $G$ to more general setting such that the corresponding algebra $A_\infty(G)$ are inverse limit of reduced (commutative) affinoid algebras.	
\end{remark}

\begin{definition}
We will in this situation use the notation $P_{A_\infty(G)}$ to mean the ring $\varprojlim_n P_{A_n(G)}$, therefore in our situation we have the corresponding ind-scheme $\mathrm{Proj}P_{A_\infty(G)}$ at the infinite level. Therefore we define the corresponding vector bundles over $\mathrm{Proj}P_{A_\infty(G)}$ to be the corresponding quasi-coherent finite locally free sheaves over the infinite level scheme $\mathrm{Proj}P_{A_\infty(G)}$. We then define the families of vector bundles 
to be those families taking the form of $(M_n)_n$ where each $M_n$ is a quasicoherent finite locally free sheaves $M_n$ over each $\mathrm{Proj}P_{A_n(G)}$ in the compatible way as in the previous setting. In this situation we define the corresponding global section of a family of vector bundles to be $\varprojlim_{n\rightarrow \infty} M_n$.	
\end{definition}

\begin{remark}
Here we actually have different notions of the corresponding modules over the corresponding Stein style. Since as in the usual story of the Iwasawa theory we have to study the corresponding derived categories, $K$-group spectra, the corresponding determinant category and so on, so we need to be careful when we would like to use the corresponding notions defined above. Here the main subtle point in the consideration is that actually it is nontrivial if the global sections of the corresponding sheaves are finitely generated. 	
\end{remark}

\indent For the corresponding families of the modules and the bundles we have the following comparison theorem:

\begin{proposition}
We have that the following four categories are equivalent: \\
I. The category of families of quasicoherent finite locally free sheaves of $\mathcal{O}_{\mathrm{Proj}P_{A_\infty(G)}}$-modules (namely the $A_\infty(G)$-relative vector bundles) over the Fargues-Fontaine curve $\mathrm{Proj}P_{A_\infty(G)}$;\\
II.The category of families of $A_\infty(G)$-relative $\varphi^a$-modules over $\widetilde{\Pi}^\infty_{R,A_\infty(G)}$;\\
III. The category of families of $A_\infty(G)$-relative $\varphi^a$-modules over $\widetilde{\Pi}_{R,A_\infty(G)}$;\\
IV. The category of families of $A_\infty(G)$-relative $\varphi^a$-bundles over $\widetilde{\Pi}_{R,A_\infty(G)}$.	
\end{proposition}

\begin{proof}
This will be further application of our established comparison above for $A_n(G)$ for each $n$.
\end{proof}

\indent And furthermore one could then further consider the following generality:

\begin{corollary}
We have that the following four categories are equivalent: \\
I. The category of families of quasicoherent finite locally free sheaves of $\mathcal{O}_{\mathrm{Proj}P_{A_\infty(G)\widehat{\otimes}_{\mathbb{Q}_p} A}}$-modules (namely the $A_\infty(G)\widehat{\otimes}_{\mathbb{Q}_p} A$-relative vector bundles) over the Fargues-Fontaine curve $\mathrm{Proj}P_{A_\infty(G)\widehat{\otimes}_{\mathbb{Q}_p} A}$;\\
II.The category of families of $A_\infty(G)\widehat{\otimes}_{\mathbb{Q}_p} A$-relative $\varphi^a$-modules over $\widetilde{\Pi}^\infty_{R,A_\infty(G)\widehat{\otimes}_{\mathbb{Q}_p} A}$;\\
III. The category of families of $A_\infty(G)\widehat{\otimes}_{\mathbb{Q}_p} A$-relative $\varphi^a$-modules over $\widetilde{\Pi}_{R,A_\infty(G)\widehat{\otimes}_{\mathbb{Q}_p} A}$;\\
IV. The category of families of $A_\infty(G)\widehat{\otimes}_{\mathbb{Q}_p} A$-relative $\varphi^a$-bundles over $\widetilde{\Pi}_{R,A_\infty(G)\widehat{\otimes}_{\mathbb{Q}_p} A}$.
%Here the group $G$ is of some equivariant type which means that $G$ is just a finite abelian group and we take $A$ to be $\mathbb{Q}_p$.
\end{corollary}

\indent Then we will study in some generality established below the cohomology which will generalize the situation in \cite[Proposition 6.3.19]{2KL15}. In our setting actually one could consider more general $A$. The algebras we are considering in the relative setting will be again the general Iwasawa theoretic level, namely the Fr\'echet algebra $A\widehat{\otimes}A_\infty(G)$. One can study the corresponding complexes of the Frobenius modules, or one can study the corresponding complexes of the corresponding families of the Frobenius modules and bundles.

\begin{proposition}
Consider now a $\varphi^a$-module $M$ over $\widetilde{\Pi}^\infty_{R,A}$ (which is always assumed to be finite projective). Then we have the following statements in our situation which generalize the corresponding results in \cite[Proposition 6.3.19]{2KL15}:\\
I. Now consider the following commutative diagram:
\[
\xymatrix@R+2pc@C+2pc{
0 \ar[r]\ar[r]\ar[r] &M \ar[d]\ar[d]\ar[d]
\ar[r]^{\varphi^a-1}\ar[r]\ar[r] &M \ar[d]\ar[d]\ar[d]
\ar[r]\ar[r]\ar[r] &0\\
0 \ar[r]\ar[r]\ar[r] &M_r
\ar[r]^{\varphi^a-1}\ar[r]\ar[r] &M_{r/q}
\ar[r]\ar[r]\ar[r] &0.
}
\]
Then we have that the vertical maps give rise to the quasi-isomorphism of the two complexes involved;\\
II. In our situation, the natural base change $M\rightarrow M\otimes_{\widetilde{\Pi}^\infty_{R,A}}\widetilde{\Pi}_{R,A}$ gives rise to a quasi-isomorphism;\\
III. Consider the following commutative diagram for some interval $[s,r]$:
\[
\xymatrix@R+2pc@C+2pc{
0 \ar[r]\ar[r]\ar[r] &M \ar[d]\ar[d]\ar[d]
\ar[r]^{\varphi^a-1}\ar[r]\ar[r] &M \ar[d]\ar[d]\ar[d]
\ar[r]\ar[r]\ar[r] &0\\
0 \ar[r]\ar[r]\ar[r] &M_{[s,r]}
\ar[r]^{\varphi^a-1}\ar[r]\ar[r] &M_{[s,r/q]}
\ar[r]\ar[r]\ar[r] &0.
}
\]
Then again we have that the vertical map is a quasi-isomorphism of the two complexes involved. Here the corresponding radii satisfy the condition that $0<s\leq r/q$.
%IV. When we only consider the algebra $A$, the map $\varphi^a-1$ is strict for $M$, therefore we have that the cohomology have the structure of a sheaf of Banach space over the Banach algebra $A$.	
\end{proposition}

\begin{proof}
Largely this follows the same idea as in \cite[Proposition 6.3.19]{2KL15} as below. The resulting complexes in our situation are those in the derived category of modules over ring $A$. 
%Since here we are working over some Fr\'echet-Stein space, we need to consider the situation over some affinoid $\mathcal{O}_{X_n}=A_n(G)$. 
We only need to prove I and III, then taking the limit will give us the second statement. Now for the first statement, as in the situation of \cite[Proposition 6.3.19]{2KL15} the injectivity of both the maps
\begin{displaymath}
H_{\varphi^a_{}}^0(M)\rightarrow H_{\varphi^a_{}}^0(M_{r}),
\end{displaymath}
\begin{displaymath}
H_{\varphi^a_{}}^0(M)\rightarrow H^0_{\varphi^a}(M_{I})	
\end{displaymath}
follows due to the property of the finite projectiveness. Then we consider the corresponding surjectivity. In this case, we just need to consider as in \cite[Proposition 6.3.19]{2KL15} the equation:
\begin{displaymath}
v=\varphi^{-na}v	
\end{displaymath}
for each $n\geq 1$, for instance considering those $v\in H_{\varphi^a_{}}^0(M_{r})$ and considering this will give rise to some element in $M_{q^{n}r}$, then taking the bootstrap like this to derived the result, which is the same for the situation in III. Then for the first order cohomology, one could just follow the idea in \cite[Proposition 6.3.19]{2KL15} with the help of the previous comparison to lift the corresponding extension. 
%Then one could finish as \cite{2KL15} to prove the last statement over the corresponding affinoid spaces. 
We do not repeat the argument again.
\end{proof}

\begin{remark}
Here we only consider the corresponding comparison carrying the deformation to $A$ instead of more general situation where the deformation is over simultaneously $A$ and $A_\infty(G)$. But one could easily have the idea on what should be established in more general setting.	
\end{remark}

%%%%%%%%%%%%%%%%%%%%Here needs some modification. 

%\newpage\section{$A$-relative $B$-pairs}
%
%
%\indent In what we have considered, we achieved some results under the $A$-deformation but remained in only the algebraic situation. That made us have added some localization with respect to some element $f$ within some ambient grading, for instance one can choose the corresponding element $t$ in our situation. Then as in the usual situation, one can ask the corresponding relation to the corresponding $B$-pairs but with some auxiliary structure involving the $A$-relativization. First let us recall the corresponding constructions in the usual setting in the language of \cite{2KL15} and \cite{2KL16}.
%
%
%
%\begin{setting}
%dd	
%\end{setting}

%%%%%%\newpage

\newpage\section{Hodge-Iwasawa Sheaves}

\subsection{Constructible $p$-adic Iwasawa Sheaves}

\indent It is actually quite natural to consider the corresponding constructible $p$-adic sheaves in Iwasawa theory as those considered in \cite{2Wit1} on the interpolation of $L$-functions after Grothendieck from \cite{2SGA5}. On the other hand, in the situation where one considers the Weil conjectures, it is actually more natural to consider both constructible and non-\'etale objects as in all types of Weil-II cnsidered by Deligne \cite{2DWe1} and \cite{2DWe2}, Kedlaya \cite{2KWe2}, Caro, Abe \cite{2CAWe2} and etc. Therefore we actually would like to consider the commutative equivariant version of the Hodge sheaves with Hodge structures considered in \cite{2KL15}. We will not consider the generality just as in the previous section here, but we will consider the generality instead more general than \cite{2KP1}.

\begin{setting}
All the flat constructible $p$-adic sheaves (both the $\mathbb{Z}_p$ ones or the $\mathbb{Z}_p$-isogeny ones), and all the Frobenius modules over the period Hodge sheaves will be considered in terms of some Iwasawa deformations, with the deformation in the corresponding reduced affinoids $A$ over $\mathbb{Q}_p$, or some adic ring $T$ (we will consider the situation considered by Witte in his noncommutative projects \cite{2Wit1} and \cite{2Wit3}), or a general Fr\'echet-Stein algebra $A_\infty(G)$ attached to some nice group $G$. In this section we assume that each adic ring is defined over $\mathbb{Z}_p$ and each quotient $T/I$ by some open two-sided ideal $I$ is of order a power $p^n$ of $p$.
%By the basic framework considered in \cite{2Wit1} the corresponding algebra $\mathbb{Z}_p[[G]]$ is an adic ring in the sense considered in \cite{2Wit1}.
\end{setting}

\begin{remark}
The Fr\'echet-Stein deformation over a point is essentially the one proposed by Kedlaya-Pottharst \cite{2KP1}, which means that actually it is very important to consider the corresponding $\mathrm{DfmLie}$ namely the $p$-adic Lie deformation of the corresponding relative $A$ Hodge-Frobenius sheaves over the corresponding pro-\'etale site of a pre-adic space $X$, which will be a natural generalization of the corresponding constructions proposed in \cite{2KP1} and established in \cite{2KL15}.
\end{remark}

\indent We now start from the spaces in our setting, namely the pre-adic spaces $X$ over $\mathbb{Q}_p$ in the sense of those considered in \cite{2KL15}. Then we will consider the corresponding homotopical categories of $p$-adic sheaves over the sites $X_{\text{\'et}}$ and $X_{\text{pro\'et}}$. As the setting considered by Witte in \cite{2Wit1} and \cite{2Wit3}, we will also consider the families of complexes with finite coefficients parametrized by the corresponding set of radicals $\mathfrak{J}$ of an adic ring $T$. \\

\indent In the style of the setting in Grothendieck's SGA IV and V, Witte studies the corresponding Iwasawa theory in the Fukaya-Kato setting, which is actually noncommutative. This is really general since the corresponding coefficients of the sheaves are actually in the adic ring mentioned above. The corresponding interpolation happens in the corresponding algebraic $K$-groups associated to some Waldhausen category. Waldhausen categories are natural places where one could do algebraic $K$-theory through Waldhausen $S$-constructions. We work out some generalization to the setting in our situations. We will consider from the following generality:

\begin{definition} \mbox{\bf (After Witte, \cite{2Wit1},\cite[Chapter 5]{2Wit3})}
Let $\mathfrak{J}$ be the set of all the two-sided ideals open of the adic ring $T$. Let $\mathbb{D}_\mathrm{perf}(X,T)$ denote the following category. Each object $(\mathcal{M}^\bullet_I)_{I\in \mathfrak{J}}$ of the category is now a family (to be more precise an inverse system) of the corresponding perfect (which are quasi-isomorphic to those strictly-perfect ones) complexes of constructible flat \'etale sheaves over $X_{\text{\'et}}$ parametrized by the corresponding open two-sided ideals of the adic ring $T$, such that: I. For each member of such family $(\mathcal{M}^\bullet_I)_{I\in \mathfrak{J}}$, we have that the complex $\mathcal{M}^\bullet_I$ is now supposed to be dg-flat over the ring $T/I$ (degreewise flat and the tensor product with any acyclic complex will be again acyclic, note here that the corresponding tensor product is taken over the ring $T/I$ when we are talking about each individual complex in a single family); II. We have the corresponding transition map $\psi_{I,J}:\mathcal{M}^\bullet_I\rightarrow \mathcal{M}^\bullet_J$ for any two open two-sided ideals $I\subset J$ with the basic base change requirement:
\begin{displaymath}
\mathcal{M}^\bullet_I\otimes_{T/I}T/J\overset{\sim}{\longrightarrow} \mathcal{M}^\bullet_J.	
\end{displaymath}
%We then will use the corresponding notation $\mathbb{D}^\infty_\mathrm{perf}(X,\mathbb{Z}_p[[G]])$ to denote the corresponding Waldhausen $(\infty,1)$-enrichement, for instance see \cite{2Bar1}.
All the modules over any quotient $T/I$ is assumed to be left $T/I$-modules, while note that then the tensor product will happen when we have another right $T/I$-module.
\end{definition}

\begin{remark}
The information of the homotopy type comes from the weak-equivalence which is defined to be the quasi-isomorphism and the cofibration which is defined to be the injection with kernel in the original category.	
\end{remark}

\indent Note that Fukaya-Kato originally considered Deligne's virtual categories (strictly speaking through more transparent approach) where they formulated their celebrated $\zeta$-isomorphism conjecture and the local $\varepsilon$-conjectures. Witte used Waldhausen construction to give more $K$-theoretic pictures, which could be regarded as some reinterpretation of the Fukaya-Kato's picture in the corresponding homotopical categories. We are somehow deeply inspired by this point of view. Concentrating everything in the zero-th degree we have the following definition:

\begin{definition}
We consider the subcategory $\mathbb{D}^\mathrm{const}_{\mathrm{perf}}(X,T)$ of $\mathbb{D}_{\mathrm{perf}}(X,T)$ which consists of all the families of complexes degenerated to zeroth degree of flat construcible sheaves over $X_\text{\'et}$, again parametrized by elements in $\mathfrak{J}$, which to be more precise means that the corresponding coefficients will live in the coefficient ring taking the form of $T/I$. For each element we will then use the corresponding notation $(\mathcal{F}^\bullet_I)_{I\in \mathfrak{J}}$ to denote each family of flat constructible \'etale sheaves.
\end{definition}

This subcategory plays the role exactly the same as those flat constructible $\ell$-adic sheaves considered by Witte's picture following Deligne's construction and equivariant Abe-Caro's arithmetic $D$-module point of views. If one is more Hodge theoretic, then one should degenerate again to define:

\begin{definition}
We will use the notation $\mathbb{D}^\mathrm{const,smooth}_{\mathrm{perf}}(X,T)$ to denote the subcategory of $\mathbb{D}^\mathrm{const}_{\mathrm{perf}}(X,T)$ consists of all the families of complexes degenerated to zeroth degree of \'etale local systems $(\mathbb{L}_I)_{I\in \mathfrak{J}}$ parametrized by the elements in $\mathfrak{J}$.
\end{definition}

\indent These will be generalization of our consideration in the previous section, but restricting to the integral $p$-adic Hodge theory.

\begin{definition}
We call each object $(\mathbb{L}_I)_{I\in \mathfrak{J}}$ in the corresponding category
\begin{displaymath}
\mathbb{D}^\mathrm{const,smooth}_{\mathrm{perf}}(X,T)
\end{displaymath}
a noncommutative \'etale local system.
\end{definition}

\indent Then we could use this point of view to generalize the corresponding Fontaine's style equivalence in the previous section over perfectoid spaces immediately after \cite[Section 8.5]{2KL15}:

\begin{proposition}
We have the following categories are equivalent for $X=\mathrm{Spa}(R,R^+)$ where $R$ is a perfect adic Banach uniform algebra over $\mathbb{F}_p$ with associated perfectoid space $X'=\mathrm{Spa}(A,A^+)$ under the perfectoid correspondence:\\
I. The category $\mathbb{D}^\mathrm{const,smooth}_{\mathrm{perf}}(X,T)$ of the noncommutative \'etale local systems over $X$;\\
II. The category $\mathbb{D}^\mathrm{const,smooth}_{\mathrm{perf}}(X',T)$ of the noncommutative \'etale local systems over $X'$.\\
%III. The category of all the families $(M_I)_I$ of $\varphi$-modules over $\Omega_R^\mathrm{int}$ with deformation to $T$ where each $\mathcal{M}_{I}$ is a $\varphi$-module over $\Omega_R^\mathrm{int}\otimes_{\mathbb{Z}_p}(T/I)$ in the usual sense such that we have the transition map $\psi_{I,J}$ as above with the additional base change properties. \\
\indent Moreover, we have that there is a fully faithful embedding from the category\\
 $\mathbb{D}^\mathrm{const,smooth}_{\mathrm{perf}}(X,T)$ (or $\mathbb{D}^\mathrm{const,smooth}_{\mathrm{perf}}(X',T)$ respectively) of all families of noncommmutative \'etale local systems over $X$ or $X'$ into the corresponding category of the families $(M_I)_I$ of $\varphi$-modules over $\widetilde{\Pi}_R^\mathrm{int}$ with deformation to $T$ where each $\mathcal{M}_{I}$ is a $\varphi$-module over $\widetilde{\Pi}_R^\mathrm{int}\otimes_{\mathbb{Z}_p}(T/I)$ in the usual sense such that we have the transition map $\psi_{I,J}$ as above with the additional base change properties. \\
%IV. The category of all the families $(\mathcal{M}_{I})_{I\in \mathfrak{J}}$ of \'etale $\varphi$-modules over $\widetilde{\Pi}_R^\mathrm{int}\widehat{\otimes}_{\mathbb{Z}_p}T$ where each $\mathcal{M}_{I}$ is an \'etale $\varphi$-module over $\widetilde{\Pi}_R^\mathrm{int}\otimes_{\mathbb{Z}_p}(T/I)$ in the usual sense such that we have the transition map $\psi_{I,J}$ as above with the additional base change properties.\\ 
\end{proposition}

\begin{proof}
In this situation this is derived from the corresponding result for the categories without Iwasawa deformation as proved in \cite[Section 8.5]{2KL15}.	
\end{proof}

\begin{definition}
We define the following key categories over the site $X_{\text{pro\'et}}$ just as in the way as above by replacing $X$ with $X_{\text{pro\'et}}$:
\begin{displaymath}
\mathbb{D}^\mathrm{const,smooth}_{\mathrm{perf}}(X_{\text{pro\'et}},T)\subset \mathbb{D}^\mathrm{const}_{\mathrm{perf}}(X_{\text{pro\'et}},T)\subset \mathbb{D}_{\mathrm{perf}}(X_{\text{pro\'et}},T).\\
\end{displaymath}	
%We also denote the $\infty$-enhancement by $\mathbb{D}^\infty_{\mathrm{perf}}(X_{\text{pro\'et}},T)$ (when this is well-established).
\end{definition}

\indent Then we glue everything in the previous two propositions over the corresponding pro-\'etale sites. Let $Y\in X_{\text{pro\'et}}$ be a perfectoid subdomain in the basis of the neighbourhood of this site, as in \cite[Section 9.2]{2KL15} taking the form of $U=\varprojlim_{n\rightarrow \infty} U_n$. Then we use the corresponding notations $\mathcal{O}_X$, $\mathcal{O}^+_X$, $\mathcal{O}^\circ_X$ to denote the corresponding sheaves evaluating over $U$ as in \cite[Section 9.2]{2KL15} as below:
\begin{align}
\mathcal{O}_X(U):=\varinjlim_{n} \mathcal{O}_{U_n}(U_n)\\
\mathcal{O}^+_X(U):=\varinjlim_{n} \mathcal{O}^+_{U_n}(U_n)\\
\mathcal{O}^\circ_X(U):=\varinjlim_{n} \mathcal{O}^\circ_{U_n}(U_n).
\end{align}

Recall that from \cite[Section 9.2]{2KL15} we have the corresponding completed sheaves which we use the notations $\widehat{\mathcal{O}}_X,\widehat{\mathcal{O}}_X^+$ to denote them, then under the perfectoid correspondence we have the corresponding version of the sheaves over the pro-\'etale sites. Recall that these are denoted in \cite[Section 9.2]{2KL15} by $\overline{\mathcal{O}}_X$ and $\overline{\mathcal{O}}_X^+$ following Andreatta-Iovita, by taking suitable tilting under the evaluation on any perfectoid subdomain $U$ as above. Then we are now at the position to glue the previous proposition over $X_{\text{pro\'et}}$.

\begin{definition}
Following \cite[Section 9.3]{2KL15}, we are going to use the notation $\widetilde{\mathbf{A}}_X\widehat{\otimes}_{\mathbb{Z}_p}T$ and the notation $\widetilde{\mathbf{A}}_X\otimes_{\mathbb{Z}_p}(T/I)$ to be the sheafifications of the corresponding presheaves over $X_{\text{pro\'et}}$ sending each $U\in X_{\text{pro\'et}}$ to
\begin{align} 
W(\overline{\mathcal{O}}_X(U))\widehat{\otimes}_{\mathbb{Z}_p}T,
W(\overline{\mathcal{O}}_X(U))\otimes_{\mathbb{Z}_p}(T/I),
\end{align}
for each element $I\in \mathfrak{J}$. And we use the notations:
\begin{displaymath}
\widetilde{\mathbf{A}}^\dagger_X\widehat{\otimes}_{\mathbb{Z}_p}T,\widetilde{\mathbf{A}}^\dagger_X\otimes_{\mathbb{Z}_p}(T/I)
\end{displaymath}
to denote the corresponding sheafification of the presheaves sending each element $U$ in $X_{\text{pro\'et}}$ to:
\begin{align} 
\widetilde{\Pi}^\mathrm{int}_{\overline{\mathcal{O}}_X(U)}\widehat{\otimes}_{\mathbb{Z}_p}T,
\widetilde{\Pi}^\mathrm{int}_{\overline{\mathcal{O}}_X(U)}\otimes_{\mathbb{Z}_p}(T/I).
\end{align}
\end{definition}

\indent Then we have the following generalization of our previous proposition to more general setting over the pro-\'etale sites, again immediately after \cite[Theorem 9.3.7]{2KL15} and its proof.

\begin{proposition}
We have the following categories are equivalent for $X$ defined above (as in the previous proposition):\\
I. The category $\mathbb{D}^\mathrm{const,smooth}_{\mathrm{perf}}(X_{\text{pro\'et}},T)$ of the noncommutative \'etale local systems over $X_{\text{pro\'et}}$;\\
II. The category $\mathbb{D}^\mathrm{const,smooth}_{\mathrm{perf}}(X'_{\text{pro\'et}},T)$ of the noncommutative \'etale local systems over $X_{\text{pro\'et}}'$.\\
 %III. The category of all the families $(M_I)_I$ of $\varphi$-modules (finite locally free) over $\widetilde{\mathbf{A}}_X$ with deformation to $T$ where each $\mathcal{M}_{I}$ is a $\varphi$-module over $\widetilde{\mathbf{A}}_X\otimes_{\mathbb{Z}_p}(T/I)$ in the usual sense such that we have the transition map $\psi_{I,J}$ as above with the additional base change properties.\\
 \indent Moreover, we have that there is a fully faithful embedding from the category \\
 $\mathbb{D}^\mathrm{const,smooth}_{\mathrm{perf}}(X_{\text{pro\'et}},T)$ (or $\mathbb{D}^\mathrm{const,smooth}_{\mathrm{perf}}(X'_{\text{pro\'et}},T)$ respectively) of all families of noncommmutative pro-\'etale local systems over $X$ or $X'$ into the corresponding category of the families $(M_I)_I$ of $\varphi$-modules over $\widetilde{\mathbf{A}}^\dagger_X$ with deformation to $T$ where each $\mathcal{M}_{I}$ is a $\varphi$-module over $\widetilde{\mathbf{A}}^\dagger_X\otimes_{\mathbb{Z}_p}(T/I)$ in the usual sense such that we have the transition map $\psi_{I,J}$ as above with the additional base change properties. \\
 %IV. The category of all the families $(\mathcal{M}_{I})_{I\in \mathfrak{J}}$ of \'etale $\varphi$-modules (finite locally free) over $\widetilde{\mathbf{A}}^\dagger_X\widehat{\otimes}_{\mathbb{Z}_p}T$ where each $\mathcal{M}_{I}$ is an \'etale $\varphi$-module over $\widetilde{\mathbf{A}}^\dagger_X\otimes_{\mathbb{Z}_p}(T/I)$ in the usual sense such that we have the transition map $\psi_{I,J}$ as above with the additional base change properties.\\
\end{proposition}

\subsection{Noncommutative-Equivariant $K$-Theory in Waldhausen C
ategories}

\indent Witte used $K$-theory to have had formulated some conjectures related to the corresponding Fukaya-Kato's key conjectures in Deligne's virtual category. Note that we defined key categories in the previous subsection for our further study on the Hodge-Iwasawa theory in the integral setting:
\begin{displaymath}
\mathbb{D}^\mathrm{const,smooth}_{\mathrm{perf}}(X,T)\subset \mathbb{D}^\mathrm{const}_{\mathrm{perf}}(X,T)\subset \mathbb{D}_{\mathrm{perf}}(X,T).
\end{displaymath}

\indent Also we have the following pro-\'etale version of the corresponding categories:

\begin{displaymath}
\mathbb{D}^\mathrm{const,smooth}_{\mathrm{perf}}(X_{\text{pro\'et}},T)\subset \mathbb{D}^\mathrm{const}_{\mathrm{perf}}(X_{\text{pro\'et}},T)\subset \mathbb{D}_{\mathrm{perf}}(X_{\text{pro\'et}},T).
\end{displaymath}

In the corresponding $K$-theory it is actually natural to ask for instance whether we have the corresponding special isomorphisms or homotopies as considered by Fukaya-Kato and Witte in the category of rigid analytic spaces or more general Huber's adic spaces. For instance one can consider the corresponding category $\mathbb{D}_{\mathrm{perf}}(\mathrm{Spa}(\mathbb{Q}_p,\mathfrak{o}_{\mathbb{Q}_p}),T)$ over a point. We actually put our observation into the following conjecture which actually predicts that over rigid analytic spaces we have the corresponding well-defined Waldhausen theory where in the algebraic geometry this is conjectured and studied by Witte in \cite{2Wit1}, \cite{2Wit2} and \cite[Section 1.4]{2Wit3}. First let us consider the following definition:

\begin{definition} \mbox{\bf (After Witte, \cite{2Wit1},\cite[Chapter 5]{2Wit3})}
Let $\mathfrak{J}$ be the set of all the two-sided ideals open of the adic ring $T$ where we assume that $p$ is a unit in $T$. Let $\mathbb{D}_\mathrm{perf}(X,T)$ denote the following category. Each object $(\mathcal{M}^\bullet_I)_{I\in \mathfrak{J}}$ of the category is now a family (to be more precise an inverse system) of the corresponding perfect (which are quasi-isomorphic to those strictly-perfect ones) complexes of constructible flat \'etale sheaves over $X_{\text{\'et}}$ parametrized by the corresponding open two-sided ideals of the adic ring $T$, such that: I. For each member of such family $(\mathcal{M}^\bullet_I)_{I\in \mathfrak{J}}$, we have that the complex $\mathcal{M}^\bullet_I$ is now supposed to be dg-flat over the ring $T/I$ (degreewise flat and the tensor product with any acyclic complex will be again acyclic, note here that the corresponding tensor product is taken over the ring $T/I$ when we are talking about each individual complex in a single family); II. We have the corresponding transition map $\psi_{I,J}:\mathcal{M}^\bullet_I\rightarrow \mathcal{M}^\bullet_J$ for any two open two-sided ideals $I\subset J$ with the basic base change requirement:
\begin{displaymath}
\mathcal{M}^\bullet_I\otimes_{T/I}T/J\overset{\sim}{\longrightarrow} \mathcal{M}^\bullet_J.	
\end{displaymath}
%We then will use the corresponding notation $\mathbb{D}^\infty_\mathrm{perf}(X,\mathbb{Z}_p[[G]])$ to denote the corresponding Waldhausen $(\infty,1)$-enrichement, for instance see \cite{2Bar1}.
All the modules over any quotient $T/I$ is assumed to be left $T/I$-modules, while note that then the tensor product will happen when we have another right $T/I$-module. And we also have the corresponding categories in the pro-\'etale setting.
\end{definition}

\begin{remark}
Note that here we consider the corresponding situation where $p$ is a unit in the adic ring $T$ which is different from the previous section.	
\end{remark}

\begin{conjecture} \label{conjecture9.13} \mbox{\bf (After Witte)}
With the assumption in the previous definition, for $X$ a rigid analytic space over the $p$-adic field $\mathbb{Q}_p$ which is separated and of finite type (when considered as some adic space locally of finite type, see \cite{2Huber1}) such that for $\sharp=\text{\'et},\text{pro\'et}$ the categories $\mathbb{D}_{\mathrm{perf}}(X_\sharp,T)$ could be endowed with the structure of Waldhausen categories, and we have there is a Waldhausen exact functor $R\Gamma(X_\sharp,.)$ (induced by the corresponding direct image functor as in \cite[Chapter 4-5]{2Wit3}) as below for $\sharp=\text{\'et},\text{pro\'et}$:
\[
\xymatrix@R+0pc@C+3pc{
\mathbb{D}_{\mathrm{perf}}(X_\sharp,T)\ar[r]^{R\Gamma(X_\sharp,.)}\ar[r]\ar[r] &\mathbb{D}_{\mathrm{perf}}(T)
}
\]
which induces the corresponding the morphism $\mathbb{K}R\Gamma(X_\sharp,.)$ as:
\[
\xymatrix@R+0pc@C+3pc{
\mathbb{K}\mathbb{D}_{\mathrm{perf}}(X_\sharp,T)\ar[r]^{\mathbb{K}R\Gamma(X_\sharp,.)}\ar[r]\ar[r] & \mathbb{K}\mathbb{D}_{\mathrm{perf}}(T),
}
\]	
we conjecture that there is a homotopy between this morphism and the zero map:
\[
\xymatrix@R+0pc@C+3pc{
\mathbb{K}\mathbb{D}_{\mathrm{perf}}(X_\sharp,T)\ar[r]^{0}\ar[r]\ar[r] & \mathbb{K}\mathbb{D}_{\mathrm{perf}}(T).
}
\]
%This actually could be enriched to the $\infty$-level in the sense that we can assume the functor $R\Gamma_c(X_\sharp,.)$ as below for $\sharp=\text{\'et},\text{pro\'et}$:
%\[
%\xymatrix@R+0pc@C+3pc{
%\mathbb{D}_{\mathrm{perf}}(X_\sharp,T)\ar[r]^{R\Gamma_c(X_\sharp,.)}\ar[r]\ar[r] &\mathbb{D}_{\mathrm{perf}}(T)
%}
%\]
%which induces the corresponding morphism $\mathbb{K}R\Gamma_c(X_\sharp,.)$ on the $K$-theory spectrum:
%\[
%\xymatrix@R+0pc@C+3pc{
%K\mathbb{D}_{\mathrm{perf}}(X_\sharp,T)\ar[r]^{KR\Gamma_c(X_\sharp,.)}\ar[r]\ar[r] & K\mathbb{D}_{\mathrm{perf}}(T).
%}
%\]		
%We conjecture that there is a homotopy between this morphism and the zero map:
%\[
%\xymatrix@R+0pc@C+3pc{
%\mathbb{K}\mathbb{D}^\infty_{\mathrm{perf}}(X_\sharp,T)\ar[r]^{0}\ar[r]\ar[r] & \mathbb{K}\mathbb{D}^\infty_{\mathrm{perf}}(T).
%}
%\]
And we conjecture this will be compatible with the higher direct image functor $Rf_!$ between different spaces (see Huber's book \cite[Chapter 5]{2Huber1}) and the changing coefficient morphism for a pair $T,T'$ of adic rings.
\end{conjecture}

\begin{remark}
We are actually not quite for sure how much higher categorical information one could see through the potential enrichment. But one might want to ask this even in the scheme situation. 
\end{remark}

\indent We would like to start from here (in the future consideration) to study the corresponding $K$-theoretic construction after the work of \cite{2Wit1} and \cite{2Wit3}. We would also like to start from in the future consideration the corresponding schematic version of the construction in Witte's thesis for instance in \cite{2Wit3}, where \cite{2Wit3} has already constructed the corresponding categories of \'etale sheaves and the corresponding functors and properties.

\begin{example}
Suppose now everything is just over the point associated to the point $\mathbb{Q}_p$, namely the space $\mathrm{Spa}(\mathbb{Q}_p,\mathfrak{o}_{\mathbb{Q}_p})$. Then the story collapses to the corresponding Galois cohomology and the corresponding Waldhausen category of the derived Galois representations, which is the Waldhausen categorical generalization of the usual Galois theoretic Tamagawa-Iwasawa theory. To be more precise we have first of all that the categories 
\begin{displaymath}
\mathbb{D}_{\mathrm{perf}}(\mathrm{Spa}(\mathbb{Q}_p,\mathfrak{o}_{\mathbb{Q}_p})_\sharp,T)	
\end{displaymath}
could be endowed with structure of Waldhausen categories as shown in \cite[Chapter 5, 6]{2Wit3}. 
%One can then apply the construction for instance in \cite{2Bar1} to enrich the categories to the level of Waldhausen $\infty$-categories:
%\begin{displaymath}
%\mathbb{D}^\infty_{\mathrm{perf}}(\mathrm{Spa}(\mathbb{Q}_p,\mathfrak{o}_{\mathbb{Q}_p})_\sharp,T).	
%\end{displaymath}
Moreover we have the following Waldhausen exact functor (see \cite[Chapter 5]{2Wit3}):
\[
\xymatrix@R+0pc@C+0pc{
R\Gamma(\mathrm{Spa}(\mathbb{Q}_p,\mathfrak{o}_{\mathbb{Q}_p})_\sharp,.):\mathbb{D}_{\mathrm{perf}}(\mathrm{Spa}(\mathbb{Q}_p,\mathfrak{o}_{\mathbb{Q}_p})_\sharp,T)
\ar[r]\ar[r]\ar[r] &\mathbb{D}_{\mathrm{perf}}(T)
}
\]
which gives rise to the following map on the spectra:
\[
\xymatrix@R+0pc@C+0pc{
KR\Gamma(\mathrm{Spa}(\mathbb{Q}_p,\mathfrak{o}_{\mathbb{Q}_p})_\sharp,.):K\mathbb{D}_{\mathrm{perf}}(\mathrm{Spa}(\mathbb{Q}_p,\mathfrak{o}_{\mathbb{Q}_p})_\sharp,T)\ar[r]\ar[r]\ar[r] &K\mathbb{D}_{\mathrm{perf}}(T).
}
\]	

\end{example}

\begin{remark}
Recall that actually in the categories of the virtual objects, Nakamura obtains partial results on the existence of the corresponding $\varepsilon$-isomorphism for \'etale $(\varphi,\Gamma)$-modules, for more detail see \cite{2Nak1} and \cite{2Nak2}. We only expect that one can construct the corresponding homotopy for the Waldhausen categorical consideration. Therefore in the situation where $p$ is not a unit in $T$ we still conjecture that there should be such Waldhausen category and the corresponding induced maps on the $K$-theoretic spaces, and even the null-homotopy. Note that this is some generalization of the corresponding picture proposed in \cite{2KP1} in the \'etale situation. Here we are considering derived \'etale $\varphi$-sheaves over the pro-\'etale sites in the corresponding Waldhausen categories. The situation of a point as above corresponding to $\mathbb{Q}_p$ is already interesting to some extent.	We will continue our discussion in this situation  below.
\end{remark}

\subsection{Hodge-Iwasawa Period Sheaves and Hodge-Iwasawa Vector Bundles}

\indent We consider here vector bundles over the corresponding adic version of the Fargues-Fontaine curves and consider the relationship between these and the ones we encountered in the last section. We are going to first recall the basics about the Fargues-Fontaine curves in the adic setting. We start from the following picture in the local situation, namely over some perfectoid ring $R$ which is some member in the basis of topology of an adic spaces in some well-established category.

\begin{setting}
We recall some adic geometry from \cite[Definition 8.7.4]{2KL15} around the corresponding geometry of the adic version of the Fargues-Fontaine curve with respect to the algebra $R$ as in our discussion above, which is to be denoted as in \cite[Definition 8.7.4]{2KL15} $\mathrm{FF}_R$. For the convenience of the readers we recall the basic construction of the curves on the level of the adic spaces. In this local setting the curve $\mathrm{FF}_R$ is defined to be the quotient of the corresponding adic space:
\begin{displaymath}
O_R:=\bigcup_{0<s<r}\mathrm{Spa}(\widetilde{\Pi}^{[s,r]}_R,\widetilde{\Pi}^{[s,r],+}_R).	
\end{displaymath}
See \cite[Definition 8.7.4]{2KL15}. Then by the basic construction we could have the chance to find a link between the adic Fargues-Fontaine curve and our originally considered schematic Fargues-Fontaine curve. To be more precise, we consider the following map:
\begin{displaymath}
P_R\rightarrow \widetilde{\Pi}^I_R	
\end{displaymath}
which naturally induces a morphism in the category of locally ringed spaces:
\begin{displaymath}
\mathrm{FF}_R\rightarrow \mathrm{Proj}P_R	
\end{displaymath}
coming from the one attached to the space $O_R$. Then we could consider the Iwasawa deformation of the vector bundles over these different locally-ringed spaces in the large categories of locally-ringed spaces.
\end{setting}

%\indent We are going to consider the corresponding deformed version of the construction above, we need to first discuss the sheafiness of the rings involved. Thanks Professor Kedlaya for the suggestion.
%
%\begin{proposition}\mbox{\bf (\cite[Kedlaya Theorem 1.4.20]{2AWS2017})}
%The sufficient and necessary condition for the quotient $L/I$ of a sheafy ring $L$ by some pseudocoherent closed ideal $I$	to be sheafy is that $L/I$ as an $L$-module is pseudocoherent.
%\end{proposition}
%
%\begin{theorem}
%For any closed interval $I$ we have that the ring $\widetilde{\Pi}^{[s,r]}_{R,A}$ is sheafy.	
%\end{theorem}
%
%\begin{proof}
%We are going to use the corresponding previous proposition to prove this. First in our situation we let $T$ be the Tate algebra, and let the ideal defining $A$ be $I$. First note that $\widetilde{\Pi}^{[s,r]}_{R,T}$ is sheafy since $T$ is pseudocoherent and we can take the product of the pseudocoherent resolution of $T$ with $\widetilde{\Pi}^{[s,r]}_{R}$ (note that here the corresponding ring $\widetilde{\Pi}^{[s,r]}_{R}$ preserves the corresponding pseudocoherence). $I'$ as an ideal generated by $I$ in the tensor product $\widetilde{\Pi}^{[s,r]}_{R,T}$ is also pseudocoherent since $I$ itself is pseudocoherent (which admits some finite free resolutation) and again the ring $\widetilde{\Pi}^{[s,r]}_{R}$ preserves the pseudocoherence. Then it suffices to check that the quotient $\widetilde{\Pi}^{[s,r]}_{R,A}$ is pseudocoherent this is again by acting the ring $\widetilde{\Pi}^{[s,r]}_{R,A}$ onto the pseudocoherent resolution of $T/I$. Then we are done.
%\end{proof}

\indent The corresponding deformed version of the adic Fargues-Fontaine curves could be defined in the same way. But note that here we need to require the corresponding deformed period rings are sheafy in the category of adic rings. This is in general not a trivial statement. For instance if we are in the situation where the affinoid algebra $A$ is a Tate algebra then we know that this is indeed the case. However in the full generality this is not that easy to see to be true. We now choose to work around this by just putting some assumption on the algebra $A$ which is already applicable enough in some interesting enough situations. Recall from \cite{2KH} that when a ring is sousperfectoid then we should have the following nice result. Recall that (from \cite[Introduction]{2KH}) a Huber ring $S$ is called sousperfectoid if it admits a split continuous morphism to a perfectoid Huber ring $P$. Note that this morphism is considered in the category of $S$-modules which are topological.

\begin{proposition} \mbox{\bf (Kedlaya-Hansen)}
The sousperfectoid rings are sheafy.	
\end{proposition}

\begin{proof}
See \cite[Corollary 7.4]{2KH}.	
\end{proof}

\begin{example}
Sousperfectoid rings are not rare, for instance one can just take the corresponding perfectoid field $\mathbb{Q}_p(p^{p^{-\infty}})^\wedge$ then $\mathbb{Q}_p$ will be then sousperfectoid in this larger perfectoid ring which admits a splitting. Also the ring $\widetilde{\Pi}^I_{R}$ is then sousperfectoid with respect to $\widetilde{\Pi}^I_{R}\widehat{\otimes}_{\mathbb{Q}_p}\mathbb{Q}_p(p^{p^{-\infty}})^\wedge$. More general, smooth affinoid algebras are sousperfectoid.	
\end{example}

\begin{assumption}
From now on, we assume that $A$ is sousperfectoid.	
\end{assumption}

\begin{corollary}
For any closed interval $I$ we have that the ring $\widetilde{\Pi}^{[s,r]}_{R,A}$ is sheafy.	
\end{corollary}

\begin{remark} Base on the observation above, we actually make the following conjecture:
\begin{conjecture}
The ring $\widetilde{\Pi}^{[s,r]}_{R,B}$ is stably uniform for any reduced affinoid algebra $B$, where the corresponding observation is as in the following. 
\end{conjecture}
\indent Note that a Tate affinoid algebra over a uniform Banach ring is uniform, so we have that the ring $\widetilde{\Pi}^{[s,r]}_{R,T_d}$ is uniform, but it is not the case that the rational localization of $\widetilde{\Pi}^{[s,r]}_{R,T_d}$ is the product of rational localizations of $\widetilde{\Pi}^{[s,r]}_{R}$ and $T_d$, which is to say that one does not have the chance to prove this directly along this. However since the ring $\widetilde{\Pi}^{[s,r]}_{R,T_d}$ is sousperfectoid,  
%This rational localization is actually uniform since we have that it is a product of preperfectoid ring with a Banach algebra where we can apply \cite[Proposition 3.7.3]{2KL15} to get this. 
so we have that the corresponding ring $\widetilde{\Pi}^{[s,r]}_{R,T_d}$ in the situation of Tate algebra is stably uniform. Now for general $A$ we suspect this is also true.
\end{remark}

\begin{setting}
In the same fashion, in our situation we consider the corresponding space 
\begin{displaymath}
O_{R,A}:=\bigcup_{0<s<r}\mathrm{Spa}(\widetilde{\Pi}^{[s,r]}_{R,A},\widetilde{\Pi}^{[s,r],+}_{R,A})	
\end{displaymath}	
which naturally gives rise to the corresponding quotient $\mathrm{FF}_{R,A}$ as mentioned above. This will allow us to consider the corresponding comparison in our generalized context as in the following.
\end{setting}

\begin{setting}
In our current situation, we have our adic version of the corresponding Fargues-Fontaine curves in some rigid analytic deformation context, which includes the situation where the rigid analytic spaces is some Stein space. In our situation the corresponding space will be:
\begin{displaymath}
O_{R,A_\infty(G)}:=\varinjlim_{n\rightarrow \infty}\bigcup_{0<s<r}\mathrm{Spa}(\widetilde{\Pi}^{[s,r]}_{R,A_n(G)},\widetilde{\Pi}^{[s,r],+}_{R,A_n(G)}).	
\end{displaymath}
Then we have the obvious notion of vector bundles over this space, which is Stein in some generalized sense, this is because the corresponding building rings are not noetherian. Then through taking the corresponding quotient we have the corresponding quotient space $\mathrm{FF}_{R,A_\infty(G)}$. We can also define the corresponding vector bundles over the $\infty$-level space while moreover we have the corresponding notions of families of vector bundles $(M_n)_n$ where each $M_n$ is a quasicoherent finite locally free sheaf over $\mathrm{FF}_{R,A_n(G)}$.	
\end{setting}

\indent Then we have the following relative version of \cite[Theorem 8.7.7]{2KL15}:

\begin{proposition}
The categories of $A$-relative vector bundles over the two spaces $\mathrm{FF}_{R,A}$ and $\mathrm{Proj}P_{R,A}$ defined above are equivalent. The categories of families of $A$-relative vector bundles over the two spaces $\mathrm{FF}_{R,A_\infty(G)}$ and $\mathrm{Proj}P_{R,A_\infty(G)}$ defined above are equivalent. 
%We have that the corresponding statements on the families of vector bundles also hold as well. 	
\end{proposition}

\begin{proof}
One could derive this by the same way as in \cite[Theorem 8.7.7]{2KL15}.	
\end{proof}

\begin{remark}
One interesting question here is that we can study the corresponding global sections taking the form of $\varprojlim_{n\rightarrow \infty}M_n$ of the families of vector bundles and the corresponding families of Frobenius modules through the equivalences we established so far. 
%Thank Professor Kedlaya for the discussion on the issue of the finiteness of the global sections.
\end{remark}

\indent First we have the following corollary of the previous proposition:

\begin{corollary}
Categories of the families of Frobenius modules over $\widetilde{\Pi}_{R,A_\infty(G)}$, the families of Frobenius bundles over $\widetilde{\Pi}_{R,A_\infty(G)}$ and the families of vector bundles over the adic Fargues-Fontaine curve $\mathrm{FF}_{R,A_\infty(G)}$ are equivalent.	
\end{corollary}

\begin{remark}
\indent The main issue here is that the corresponding families of Frobenius modules and vector bundles might not have finite generated global sections since in the whole process of taking inverse limit the numbers of the generators might be blowing up.
At this moment we do not know whether the global sections of a family of vector bundles over $\mathrm{FF}_{R,A_\infty(G)}$ are finitely generated. In the situation where one can get the chance to control the largest number of generators during the whole process of taking the inverse limit, one might be able to have the chance to control the finiteness of the global sections. However, this is not trivially correct. The Fargues-Fontaine curves involved are defined over some Fr\'echet-Stein algebras which behave well due to the well-established theory of Fr\'echet-Stein algebras in the noetherian setting. However the Robba rings in the relative setting are highly not noetherian. The corresponding cohomologies of these big modules or families of big modules will eventually then be the corresponding modules or families of modules over Fr\'echet-Stein algebras, which will be definitely easier to control and study since the theory of Fr\'echet-Stein algebras should be able to help us to construct the corresponding Iwasawa theoretic or $K$-theoretic objects. 
%Here we consider the corresponding parallel story of the corresponding quasi-Stein spaces from \cite[Section 2.6]{2KL16}. 
%We can use (after choosing some suitable radii $s,r$) the corresponding covering $\{\mathrm{Spa}(\widetilde{\Pi}^{[sq^n,rq^n]}_{\mathrm{R},A_{\alpha(G)+1}},\widetilde{\Pi}^{[sq^n,rq^n],\mathrm{Gr}}_{\mathrm{R},A_{\alpha+1}(G)})-\overline{\mathrm{Spa}(\widetilde{\Pi}^{[sq^n,rq^n]}_{\mathrm{R},A_{\alpha(G)}},\widetilde{\Pi}^{[sq^n,rq^n],\mathrm{Gr}}_{\mathrm{R},A_{\alpha}(G)})}\}_{n,\alpha}$ to cover the corresponding adic space $O_R$, which finishes the proof by \cite[Proposition 2.6.17]{2KL16}. 	

\end{remark}

\indent Furthermore we consider the following Hodge-Iwasawa construction under the idea from Kedlaya-Pottharst through the perfectoid subdomain covering mentioned above over any preadic space $X$ defined over $\mathbb{Q}_p$.

\begin{definition}
We define the following $A$-relative sheaves of period rings after \cite[Definition 9.3.3]{2KL15}:
\begin{align}
\underline{\underline{\Omega}}_{X,A},\underline{\underline{\Omega}}_{X,A}^\mathrm{int},	\underline{\underline{\widetilde{\Pi}}}_{X,A}^{\mathrm{int},r},\underline{\underline{\widetilde{\Pi}}}_{X,A}^\mathrm{int},\underline{\underline{\widetilde{\Pi}}}_{X,A}^\mathrm{bd,r},\underline{\underline{\widetilde{\Pi}}}_{X,A}^\mathrm{bd},
\underline{\underline{\widetilde{\Pi}}}_{X,A}^r,\underline{\underline{\widetilde{\Pi}}}_{X,A}^I,
\underline{\underline{\widetilde{\Pi}}}_{X,A}^\infty,
\underline{\underline{\widetilde{\Pi}}}_{X,A},
\end{align}
by locally taking the suitable completed product and then glueing. Similarly we have the following Iwasawa sheaves as well:
\begin{align}
\underline{\underline{\Omega}}_{X,A\widehat{\otimes}A_\infty(G)},\underline{\underline{\Omega}}_{X,A\widehat{\otimes}A_\infty(G)}^\mathrm{int},	\underline{\underline{\widetilde{\Pi}}}_{X,A\widehat{\otimes}A_\infty(G)}^{\mathrm{int},r},\underline{\underline{\widetilde{\Pi}}}_{X,A\widehat{\otimes}A_\infty(G)}^\mathrm{int},\underline{\underline{\widetilde{\Pi}}}_{X,A\widehat{\otimes}A_\infty(G)}^{\mathrm{bd},r},\underline{\underline{\widetilde{\Pi}}}_{X,A\widehat{\otimes}A_\infty(G)}^\mathrm{bd},
\end{align}
\begin{align}
\underline{\underline{\widetilde{\Pi}}}_{X,A\widehat{\otimes}A_\infty(G)}^r,\underline{\underline{\widetilde{\Pi}}}_{X,A\widehat{\otimes}A_\infty(G)}^I,
\underline{\underline{\widetilde{\Pi}}}_{X,A\widehat{\otimes}A_\infty(G)}^\infty,
\underline{\underline{\widetilde{\Pi}}}_{X,A\widehat{\otimes}A_\infty(G)}.
\end{align}
\end{definition}

Then we have the following $A$-relative version of \cite[Theorem 9.3.12]{2KL15}:

\begin{proposition}\mbox{}\\
The categories of $A$-relative vector bundles over $\mathrm{FF}_{X,A}$ and Frobenius-equivariant $A$-relative Hodge-Iwasawa finite locally free sheaves over the period sheaves $\underline{\underline{\widetilde{\Pi}}}_{X,A}^\infty,
\underline{\underline{\widetilde{\Pi}}}_{X,A}$ are all equivalent.
\end{proposition}

\begin{proof}
This is by considering the perfectoid covering and the local equivalence we established before.	
\end{proof}

\subsection{Higher Homotopical Geometrized Tamagawa-Iwasawa Theory of Hodge-Iwasawa Modules}
\indent The structure of the corresponding Hodge-Iwasawa modules in our mind has the potential to give us the chance to use them to study the corresponding rational Iwasawa theory. The corresponding relationship we established above could show us the corresponding Iwasawa theory of the relative Frobenius modules over the pro-\'etale site is equivalent to that of the corresponding $B$-pairs (we will also discuss this more in our further study). For instance the latter gives us directly some kind of family version of exponential maps and dual-exponential maps in the higher dimensional situation. We now give some further discussion around our main goals in mind. We first start from the following context:

\begin{setting}
Over a separated rigid analytic space over $\mathbb{Q}_p$ (we assume it to be of finite type as an adic space) $X$ we have the category of the pseudocoherent $\varphi$-sheaves over the sheaf of ring $\underline{\underline{\widetilde{\Pi}}}_{X}$ over the pro-\'etale site $X_\text{pro-\'etale}$. We then use the notation $ChQCoh\Pi\Phi_{X_\text{pro-\'etale}}$ to denote the category of all the chain complexes of objects in the category $QCoh\Pi\Phi_{X_\text{pro-\'etale}}$ which we define it to be the category of all the quasicoherent $\varphi$-sheaves (as in the pseudo-coherent situations) over the ring $\underline{\underline{\widetilde{\Pi}}}_{X}$. 
\end{setting}

\begin{remark}
One should actually regard the category $QCoh\Pi\Phi_{X_\text{pro-\'etale}}$ as the category of all the direct limits of the objects in the category of all the pseudocoherent $\varphi$-sheaves over the sheaf of ring $\underline{\underline{\widetilde{\Pi}}}_{X}$ over the pro-\'etale site $X_\text{pro-\'etale}$. During taking the direct limit we assume that we are actually always carrying the corresponding Frobenius action.  
\end{remark}

\begin{definition}
A complex of objects of $ChQCoh\Pi\Phi_{X_\text{pro-\'etale}}$  taking the form of
\[
\xymatrix@R+0pc@C+0pc{
\ar[r]\ar[r]\ar[r] &... \ar[r]\ar[r]\ar[r] &M^{n-1} \ar[r]\ar[r]\ar[r] &M^n
\ar[r]\ar[r]\ar[r] &M^{n+1}
\ar[r]\ar[r]\ar[r] &...
}
\]	
is called pseudocoherent if it is quasi-isomorphic to a bounded above complex of finite projective objects in $\Pi\Phi_{X_\text{pro-\'etale}}$ (the category of all the pseudo-coherent $\varphi$-sheaves over $\underline{\underline{\widetilde{\Pi}}}_{X}$)
\[
\xymatrix@R+0pc@C+0pc{
\ar[r]\ar[r]\ar[r] &... \ar[r]\ar[r]\ar[r] &F^{m-2} \ar[r]\ar[r]\ar[r] &F^{m-1}
\ar[r]\ar[r]\ar[r] &F^m
\ar[r]\ar[r]\ar[r] &0.
}
\] 
A complex of objects of $ChQCoh\Pi\Phi_{X_\text{pro-\'etale}}$  taking the form of
\[
\xymatrix@R+0pc@C+0pc{
\ar[r]\ar[r]\ar[r] &... \ar[r]\ar[r]\ar[r] &M^{n-1} \ar[r]\ar[r]\ar[r] &M^n
\ar[r]\ar[r]\ar[r] &M^{n+1}
\ar[r]\ar[r]\ar[r] &...
}
\]	
is called perfect if it is quasi-isomorphic to a bounded complex of finite projective objects in $\Pi\Phi_{X_\text{pro-\'etale}}$ (the category of all the pseudo-coherent $\varphi$-sheaves over $\underline{\underline{\widetilde{\Pi}}}_{X}$ as those in \cite[Chapter 8]{2KL16})
\[
\xymatrix@R+0pc@C+0pc{
\ar[r]\ar[r]\ar[r] &... \ar[r]\ar[r]\ar[r] &F^{m-2} \ar[r]\ar[r]\ar[r] &F^{m-1}
\ar[r]\ar[r]\ar[r] &F^m
\ar[r]\ar[r]\ar[r] &0.
}
\] 
We then use the notation $D_{\mathrm{pseudo}}\Pi\Phi_{X_\text{pro-\'etale}}$ to denote the category (not the derived one) of all the pseudocoherent complexes of objects in $QCoh\Pi\Phi_{X_\text{pro-\'etale}}$ in the above sense. And we use the notation $D_{\mathrm{perf}}\Pi\Phi_{X_\text{pro-\'etale}}$ for the perfect complexes. We use the notation $D^{\mathrm{cb}}_{\mathrm{pseudo}}\Pi\Phi_{X_\text{pro-\'etale}}$ to denote the category of cohomology bounded pseudocoherent complexes. And we use the notations $D^{\text{dg-flat}}_{\mathrm{perf}}\Pi\Phi_{X_\text{pro-\'etale}}$ to denote the the subcategories consisting of dg-flat complexes. And we use the notations $D^{\text{str}}_{\mathrm{perf}}\Pi\Phi_{X_\text{pro-\'etale}}$ to denote the subcategories consisting of strictly perfect complexes.
\end{definition}

\begin{proposition}
There is a Waldhausen structure over each of
\begin{align}
&D_{\mathrm{pseudo}}\Pi\Phi_{X_\text{pro-\'etale}},D^{\mathrm{cb}}_{\mathrm{pseudo}}\Pi\Phi_{X_\text{pro-\'etale}},\\
&D_{\mathrm{perf}}\Pi\Phi_{X_\text{pro-\'etale}},D^{\mathrm{dg-flat}}_{\mathrm{perf}}\Pi\Phi_{X_\text{pro-\'etale}},	D^{\mathrm{str}}_{\mathrm{perf}}\Pi\Phi_{X_\text{pro-\'etale}}.
\end{align} 	
\end{proposition}

\begin{proof}
We only show this for the first category. We are talking about the Waldhausen category in the sense considered by \cite[Proposition 3.1.1]{2Wit3} and \cite{2TT1} (namely the corresponding complicial biWaldhausen categories). Recalling from the work \cite[Proposition 3.1.1]{2Wit3}, the corresponding criterion is to check that the corresponding smaller full additive subcategory of complexes extracted from an abelian category has the stability under the corresponding shifts and the extension through the exact sequences. But we are talking about pseudocoherent complexes, two out of three property follows from \cite[Tag, 064R]{2SP} where one can show this directly by applying \cite[Derived Category, Lemma 16.11]{2SP}. Of course note that we are not talking about derived category. In this situation we define the corresponding cofibration to be the corresponding degreewise monomorphism with kernel in the original category in which we are considering, while we define the corresponding weak-equivalences to be quasi-isomorphisms. Then one can check the corresponding conditions in the criterion of Waldhausen category in the sense of \cite[Proposition 3.1.1]{2Wit3} and \cite{2TT1} satisfy.
\end{proof}

%\indent We then have the following conjectures which we believe should be some kind of generalization of the pictures in the \'etale setting and in the Deligne's category of virtual objects.

%\begin{conjecture}
%Let $X$ be $\mathrm{Spa}(\mathbb{Q}_p,\mathfrak{o}_{\mathbb{Q}_p})$. Then we have that the corresponding direct image functor induces a map from the $K$-theoretic space $\mathbb{K}D^{\mathrm{str}}_{\mathrm{perf}}\Pi\Phi_{X_\text{pro-\'etale}}$ or $\mathbb{K}D^{\mathrm{dg-flat}}_{\mathrm{perf}}\Pi\Phi_{X_\text{pro-\'etale}}$ to the corresponding space $\mathbb{K}D^{\mathrm{str}}_{\mathrm{perf}}(\mathbb{Q}_p)$ or $\mathbb{K}D^{\mathrm{dg-flat}}_{\mathrm{perf}}(\mathbb{Q}_p)$ respectively  (the $K$-theoretic space associated to the category of strictly perfect or dg-flat perfect complexes of all the modules over $\mathbb{Q}_p$). And we conjecture that this map is homotopic to zero.	
%\end{conjecture}

\indent We now consider the corresponding objects over the schematic Fargues-Fontaine curves. By our previous comparison theorems (and the corresponding theorems in \cite{2KL16}) on the deformations of the objects we can only consider the scheme theory to extract the corresponding information on the Tamagawa-Iwasawa deformations for the relative $p$-adic Hodge theory over the Robba rings which to some extent very complicated. We can now make the following discussion:

\begin{setting}
Consider any uniform perfect adic Banach algebra $R$ as in our previous consideration. We have the schematic version of the Fargues-Fontaine curve which is denoted by $\mathrm{Proj}P_{R}$. We now use the notation $Mod\mathcal{O}_{\mathrm{Proj}P_{R}}$ to denote the corresponding category of all the sheaves of $\mathcal{O}_{\mathrm{Proj}P_{R}}$-modules over the scheme $\mathrm{Proj}P_{R}$. Then we will use the notation $ChMod\mathcal{O}_{\mathrm{Proj}P_{R}}$ to denote the category of all the complexes of sheaves of $\mathcal{O}_{\mathrm{Proj}P_{R}}$-modules over the scheme $\mathrm{Proj}P_{R}$, with the corresponding derived category $DMod\mathcal{O}_{\mathrm{Proj}P_{R}}$.
\end{setting}

\begin{definition}
A complex of objects of $Mod\mathcal{O}_{\mathrm{Proj}P_{R}}$ taking the form of
\[
\xymatrix@R+0pc@C+0pc{
\ar[r]\ar[r]\ar[r] &... \ar[r]\ar[r]\ar[r] &M^{n-1} \ar[r]\ar[r]\ar[r] &M^n
\ar[r]\ar[r]\ar[r] &M^{n+1}
\ar[r]\ar[r]\ar[r] &...
}
\]	
is called pseudocoherent if it is quasi-isomorphic to a bounded above complex of finite projective objects in the category of all the sheaves of $\mathcal{O}_{\mathrm{Proj}P_{R}}$-modules over the schematic Fargues-Fontaine curve
\[
\xymatrix@R+0pc@C+0pc{
\ar[r]\ar[r]\ar[r] &... \ar[r]\ar[r]\ar[r] &F^{m-2} \ar[r]\ar[r]\ar[r] &F^{m-1}
\ar[r]\ar[r]\ar[r] &F^m
\ar[r]\ar[r]\ar[r] &0.
}
\]
We then use the notation $D_{\mathrm{pseudo}}\mathrm{Proj}P_{R}$ to denote the category (not the derived one) of all the pseudocoherent complexes of objects in $Mod\mathcal{O}_{\mathrm{Proj}P_{R}}$ in the above sense. A complex of objects of $Mod\mathcal{O}_{\mathrm{Proj}P_{R}}$ taking the form of
\[
\xymatrix@R+0pc@C+0pc{
\ar[r]\ar[r]\ar[r] &... \ar[r]\ar[r]\ar[r] &M^{n-1} \ar[r]\ar[r]\ar[r] &M^n
\ar[r]\ar[r]\ar[r] &M^{n+1}
\ar[r]\ar[r]\ar[r] &...
}
\]	
is called perfect if it is quasi-isomorphic to a bounded complex of finite projective objects in the category of all the sheaves of $\mathcal{O}_{\mathrm{Proj}P_{R}}$-modules over the schematic Fargues-Fontaine curve
\[
\xymatrix@R+0pc@C+0pc{
\ar[r]\ar[r]\ar[r] &... \ar[r]\ar[r]\ar[r] &F^{m-2} \ar[r]\ar[r]\ar[r] &F^{m-1}
\ar[r]\ar[r]\ar[r] &F^m
\ar[r]\ar[r]\ar[r] &0.
}
\]
We then use the notation $D_{\mathrm{perf}}\mathrm{Proj}P_{R}$ to denote the category (not the derived one) of all the perfect complexes of objects in $Mod\mathcal{O}_{\mathrm{Proj}P_{R}}$ in the above sense.	
\end{definition}

\begin{proposition}
The category $D_{\mathrm{pseudo}}\mathrm{Proj}P_{R}$ is a category admitting Waldhausen structure.	
\end{proposition}

\begin{proof}
This is standard such as in \cite[Section 1-3]{2TT1}. We present the proof as above for the convenience of the readers. We are now talking about again the context of \cite[Proposition 3.1.1]{2Wit3} and \cite{2TT1}. The corresponding criterion in \cite{2Wit3} requires the corresponding stability of the corresponding shifts and extension in the category of the corresponding complexes in our context. Again we can consider the corresponding results of \cite[Tag, 064R]{2SP} to prove this. Then the corresponding cofibration will be the corresponding monomorphism with kernel in the original category of the complexes while the corresponding weak-equivalences are taken to be the corresponding quasi-isomorphisms. 	
\end{proof}

Then we can define the corresponding category of pseudocoherent sheaves over the adic Fargues-Fontaine curves:

\begin{setting}
After \cite{2KL16}, we consider the cateogory $Mod\mathcal{O}_{\mathrm{FF}_{\widetilde{R}_\psi}}$ of all sheaves of $\mathcal{O}_{\mathrm{FF}_{\widetilde{R}_\psi}}$-modules. Here the ring $\widetilde{R}_\psi$ is the perfect adic uniform Banach ring attached to the toric tower considered in \cite{2KL16}. We then use the notation $ChMod\mathcal{O}_{\mathrm{FF}_{\widetilde{R}_\psi}}$ to denote the category of all the chain complexes of objects in the previous category $Mod\mathcal{O}_{\mathrm{FF}_{\widetilde{R}_\psi}}$. 
\end{setting}

\begin{definition}
A complex of objects of $ChMod\mathcal{O}_{\mathrm{FF}_{\widetilde{R}_\psi}}$ taking the form of
\[
\xymatrix@R+0pc@C+0pc{
\ar[r]\ar[r]\ar[r] &... \ar[r]\ar[r]\ar[r] &M^{n-1} \ar[r]\ar[r]\ar[r] &M^n
\ar[r]\ar[r]\ar[r] &M^{n+1}
\ar[r]\ar[r]\ar[r] &...
}
\]	
is called pseudocoherent if it is quasi-isomorphic to a bounded above complex of finite projective objects in the category of all the pseudocoherent sheaves over the Fargues-Fontaine curve
\[
\xymatrix@R+0pc@C+0pc{
\ar[r]\ar[r]\ar[r] &... \ar[r]\ar[r]\ar[r] &F^{m-2} \ar[r]\ar[r]\ar[r] &F^{m-1}
\ar[r]\ar[r]\ar[r] &F^m
\ar[r]\ar[r]\ar[r] &0.
}
\] 
We then use the notation $D_{\mathrm{pseudo,alg}}\mathrm{FF}_{\widetilde{R}_\psi}$ to denote the category (not the derived one) of all the pseudocoherent complexes of objects in $Mod\mathcal{O}_{\mathrm{FF}_{\widetilde{R}_\psi}}$ in the above sense. We also have the category $D^\mathrm{cb}_{\mathrm{pseudo,alg}}\mathrm{FF}_{\widetilde{R}_\psi}$ consisting of all the cohomology bounded pseudocoherent complexes.\\
A complex of objects of $ChMod\mathcal{O}_{\mathrm{FF}_{\widetilde{R}_\psi}}$ taking the form of
\[
\xymatrix@R+0pc@C+0pc{
\ar[r]\ar[r]\ar[r] &... \ar[r]\ar[r]\ar[r] &M^{n-1} \ar[r]\ar[r]\ar[r] &M^n
\ar[r]\ar[r]\ar[r] &M^{n+1}
\ar[r]\ar[r]\ar[r] &...
}
\]	
is called perfect if it is quasi-isomorphic to a bounded complex of finite projective objects in the category of all the pseudocoherent sheaves over the Fargues-Fontaine curve
\[
\xymatrix@R+0pc@C+0pc{
\ar[r]\ar[r]\ar[r] &... \ar[r]\ar[r]\ar[r] &F^{m-2} \ar[r]\ar[r]\ar[r] &F^{m-1}
\ar[r]\ar[r]\ar[r] &F^m
\ar[r]\ar[r]\ar[r] &0.
}
\] 
We then use the notation $D_{\mathrm{perf,alg}}\mathrm{FF}_{\widetilde{R}_\psi}$ to denote the category (not the derived one) of all the perfect complexes of objects in $Mod\mathcal{O}_{\mathrm{FF}_{\widetilde{R}_\psi}}$ in the above sense. We also have the subcategories $D^{\mathrm{dg-flat}}_{\mathrm{perf,alg}}\mathrm{FF}_{\widetilde{R}_\psi}$ and $D^{\mathrm{str}}_{\mathrm{perf,alg}}\mathrm{FF}_{\widetilde{R}_\psi}$ of dg-flat objects and strictly perfect objects.
\end{definition}

\begin{proposition}
The categories
\begin{align}
&D_{\mathrm{pseudo,alg}}\mathrm{FF}_{\widetilde{R}_\psi},D^\mathrm{cb}_{\mathrm{pseudo,alg}}\mathrm{FF}_{\widetilde{R}_\psi}\\
&D_{\mathrm{perf,alg}}\mathrm{FF}_{\widetilde{R}_\psi},D^\mathrm{dg-flat}_{\mathrm{perf,alg}}\mathrm{FF}_{\widetilde{R}_\psi},D^\mathrm{str}_{\mathrm{perf,alg}}\mathrm{FF}_{\widetilde{R}_\psi}	
\end{align}
admit structure of Waldhausen categories.	
\end{proposition}

%\begin{proposition}
%If the tower is usual cyclotomic one, then the conjecture holds true.	
%\end{proposition}

\begin{proof}
The corresponding category involved in this situation is not that far from being equivalent to the situation where we considered in the Frobenius sheaves situation. On the other hand, one can then following the proof to give the chance to prove the corresponding statement here.	
\end{proof}

\indent Now we consider the more general analytic objects:

\begin{setting}
Following \cite{2KL16}, we have the category $QCoh\mathcal{O}_{\mathrm{FF}_{\widetilde{R}_\psi}}$ of the quasi-coherent sheaves (direct limits of the pseudocoherent objects) over the Fargues-Fontaine curve $\mathrm{FF}_{\widetilde{R}_\psi}$. We then use the notation $ChQCoh\mathcal{O}_{\mathrm{FF}_{\widetilde{R}_\psi}}$ to denote the category of all the chain complexes of objects in the previous category $QCoh\mathcal{O}_{\mathrm{FF}_{\widetilde{R}_\psi}}$. 
\end{setting}

\begin{definition}
A complex of objects of $ChQCoh\mathcal{O}_{\mathrm{FF}_{\widetilde{R}_\psi}}$ taking the form of
\[
\xymatrix@R+0pc@C+0pc{
\ar[r]\ar[r]\ar[r] &... \ar[r]\ar[r]\ar[r] &M^{n-1} \ar[r]\ar[r]\ar[r] &M^n
\ar[r]\ar[r]\ar[r] &M^{n+1}
\ar[r]\ar[r]\ar[r] &...
}
\]	
is called pseudocoherent if it is quasi-isomorphic to a bounded above complex of finite projective objects in the category of all the pseudocoherent sheaves over the Fargues-Fontaine curve:
\[
\xymatrix@R+0pc@C+0pc{
\ar[r]\ar[r]\ar[r] &... \ar[r]\ar[r]\ar[r] &F^{m-2} \ar[r]\ar[r]\ar[r] &F^{m-1}
\ar[r]\ar[r]\ar[r] &F^m
\ar[r]\ar[r]\ar[r] &0.
}
\]
We then use the notation $D_{\mathrm{pseudo}}\mathrm{FF}_{\widetilde{R}_\psi}$ to denote the category (not the derived one) of all the pseudocoherent complexes of objects in $QCoh\mathcal{O}_{\mathrm{FF}_{\widetilde{R}_\psi}}$ in the above sense. And we use the notations $D^\mathrm{cb}_{\mathrm{pseudo}}\mathrm{FF}_{\widetilde{R}_\psi}$ to denote the subcategories of cohomology bounded complexes.
A complex of objects of $ChQCoh\mathcal{O}_{\mathrm{FF}_{\widetilde{R}_\psi}}$ taking the form of
\[
\xymatrix@R+0pc@C+0pc{
\ar[r]\ar[r]\ar[r] &... \ar[r]\ar[r]\ar[r] &M^{n-1} \ar[r]\ar[r]\ar[r] &M^n
\ar[r]\ar[r]\ar[r] &M^{n+1}
\ar[r]\ar[r]\ar[r] &...
}
\]	
is called perfect if it is quasi-isomorphic to a bounded complex of finite projective objects in the category of all the pseudocoherent sheaves over the Fargues-Fontaine curve:
\[
\xymatrix@R+0pc@C+0pc{
\ar[r]\ar[r]\ar[r] &... \ar[r]\ar[r]\ar[r] &F^{m-2} \ar[r]\ar[r]\ar[r] &F^{m-1}
\ar[r]\ar[r]\ar[r] &F^m
\ar[r]\ar[r]\ar[r] &0.
}
\]
We then use the notation $D_{\mathrm{perf}}\mathrm{FF}_{\widetilde{R}_\psi}$ to denote the category (not the derived one) of all the perfect complexes of objects in $QCoh\mathcal{O}_{\mathrm{FF}_{\widetilde{R}_\psi}}$ in the above sense. And we use the notations $D^\mathrm{dg-flat}_{\mathrm{perf}}\mathrm{FF}_{\widetilde{R}_\psi}$ and $D^\mathrm{str}_{\mathrm{perf}}\mathrm{FF}_{\widetilde{R}_\psi}$ to denote the subcategories of dg-flat complexes and strictly perfect complexes.
\end{definition}

\begin{proposition}
The categories
\begin{align}
&D_{\mathrm{pseudo}}\mathrm{FF}_{\widetilde{R}_\psi},D^{\mathrm{cb}}_{\mathrm{pseudo}}\mathrm{FF}_{\widetilde{R}_\psi},\\
&D_{\mathrm{perf}}\mathrm{FF}_{\widetilde{R}_\psi},D^{\mathrm{dg-flat}}_{\mathrm{perf}}\mathrm{FF}_{\widetilde{R}_\psi},D^{\mathrm{str}}_{\mathrm{perf}}\mathrm{FF}_{\widetilde{R}_\psi}
\end{align}
admit the structure of Waldhausen categories.	
\end{proposition}

\begin{proof}
See the proof for the previous proposition.	
\end{proof}

\indent Now one can further discuss the corresponding deformed version of the picture above, we mainly then focus on the schematic version of the Fargues-Fontaine curve:

\begin{setting}
Let $A$ be a reduced affinoid algebra as before. Consider any uniform perfect adic Banach algebra $R$ as in our previous consideration. We have the schematic version of the Fargues-Fontaine curve which is denoted by $\mathrm{Proj}P_{R,A}$. We now use the notation $Mod\mathcal{O}_{\mathrm{Proj}P_{R,A}}$ to denote the corresponding category of all the sheaves of $\mathcal{O}_{\mathrm{Proj}P_{R,A}}$-modules over the scheme $\mathrm{Proj}P_{R}$. Then we will use the notation $ChMod\mathcal{O}_{\mathrm{Proj}P_{R,A}}$ to denote the category of all the complexes of sheaves of $\mathcal{O}_{\mathrm{Proj}P_{R,A}}$-modules over the scheme $\mathrm{Proj}P_{R,A}$, with the corresponding derived category $DMod\mathcal{O}_{\mathrm{Proj}P_{R,A}}$.
\end{setting}

\begin{definition}
A complex of objects of $Mod\mathcal{O}_{\mathrm{Proj}P_{R,A}}$ taking the form of
\[
\xymatrix@R+0pc@C+0pc{
\ar[r]\ar[r]\ar[r] &... \ar[r]\ar[r]\ar[r] &M^{n-1} \ar[r]\ar[r]\ar[r] &M^n
\ar[r]\ar[r]\ar[r] &M^{n+1}
\ar[r]\ar[r]\ar[r] &...
}
\]	
is called pseudocoherent if it is quasi-isomorphic to a bounded above complex of finite projective objects in the category of all the quasicoherent sheaves over the schematic Fargues-Fontaine curve
\[
\xymatrix@R+0pc@C+0pc{
\ar[r]\ar[r]\ar[r] &... \ar[r]\ar[r]\ar[r] &F^{m-2} \ar[r]\ar[r]\ar[r] &F^{m-1}
\ar[r]\ar[r]\ar[r] &F^m
\ar[r]\ar[r]\ar[r] &0.
}
\]
We then use the notation $D_{\mathrm{pseudo}}\mathrm{Proj}P_{R,A}$ to denote the category (not the derived one) of all the pseudocoherent complexes of objects in $Mod\mathcal{O}_{\mathrm{Proj}P_{R,A}}$ in the above sense. We can also define the corresponding category $D^{\mathrm{cb}}_{\mathrm{pseudo}}\mathrm{Proj}P_{R,A}$ as above. \\	
\indent A complex of objects of $Mod\mathcal{O}_{\mathrm{Proj}P_{R,A}}$ taking the form of
\[
\xymatrix@R+0pc@C+0pc{
\ar[r]\ar[r]\ar[r] &... \ar[r]\ar[r]\ar[r] &M^{n-1} \ar[r]\ar[r]\ar[r] &M^n
\ar[r]\ar[r]\ar[r] &M^{n+1}
\ar[r]\ar[r]\ar[r] &...
}
\]	
is called perfect if it is quasi-isomorphic to a bounded complex of finite projective objects in the category of all the quasicoherent sheaves over the schematic Fargues-Fontaine curve
\[
\xymatrix@R+0pc@C+0pc{
\ar[r]\ar[r]\ar[r] &... \ar[r]\ar[r]\ar[r] &F^{m-2} \ar[r]\ar[r]\ar[r] &F^{m-1}
\ar[r]\ar[r]\ar[r] &F^m
\ar[r]\ar[r]\ar[r] &0.
}
\]
We then use the notation $D_{\mathrm{perf}}\mathrm{Proj}P_{R,A}$ to denote the category (not the derived one) of all the perfect complexes of objects in $Mod\mathcal{O}_{\mathrm{Proj}P_{R,A}}$ in the above sense. We can also define the corresponding category $D^{\mathrm{dg-flat}}_{\mathrm{perf}}\mathrm{Proj}P_{R,A}$ of dg-flat perfect complexes and the corresponding category $D^{\mathrm{str}}_{\mathrm{pseudo}}\mathrm{Proj}P_{R,A}$ of strictly perfect complexes.
\end{definition}

\begin{proposition}
The categories
\begin{align}
&D_{\mathrm{pseudo}}\mathrm{Proj}P_{R,A},D^{\mathrm{cb}}_{\mathrm{pseudo}}\mathrm{Proj}P_{R,A},\\
&D_{\mathrm{perf}}\mathrm{Proj}P_{R,A},D^{\mathrm{dg-flat}}_{\mathrm{perf}}\mathrm{Proj}P_{R,A},D^{\mathrm{str}}_{\mathrm{perf}}\mathrm{Proj}P_{R,A},	
\end{align}
are categories admitting Waldhausen structure.	
\end{proposition}

\begin{proof}
We only prove the statement for the first category. See the proof above for $D_{\mathrm{pseudo}}\mathrm{Proj}P_{R}$ without the deformation.	
\end{proof}

\begin{conjecture}
Assume that $\psi$ is the cyclotomic tower. Then the corresponding total derived section functor (by using the Godement resolution) induces a map on the $K$-theory space $\mathbb{K}D^\mathrm{dg-flat}_{\mathrm{perf}}\mathrm{Proj}P_{\widetilde{R}_\psi,A}$ to that of the category $\mathbb{K}D^{\mathrm{dg-flat}}_{\mathrm{perf}}(A)$ of all the dg-flat perfect complexes of $A$-modules. And we conjecture that this is homotopic to zero. \\
\indent Assume that $\psi$ is the cyclotomic tower. Then the corresponding total derived section functor (by using the Godement resolution) induces a map on the $K$-theory space $\mathbb{K}D^\mathrm{str}_{\mathrm{perf}}\mathrm{Proj}P_{\widetilde{R}_\psi,A}$ to that of the category $\mathbb{K}D^\mathrm{str}_{\mathrm{perf}}(A)$ of all the strictly perfect complexes of $A$-modules. 
And we conjecture that in this situation this is homotopic to zero.
\end{conjecture}

\indent Although not equivalent to the previous picture in the deformed setting we should still consider the following picture over the adic Fargues-Fontaine curve:

\begin{setting}
Let $A$ be as above such that we have the well-defined adic space $\mathrm{FF}_{\widetilde{R}_\psi,A}$	where the corresponding tower $\psi$ is the same as the previous one. Then we have the corresponding category $\mathrm{Mod}\mathcal{O}_{\mathrm{FF}_{\widetilde{R}_\psi,A}}$ as above by taking the corresponding sheaves of $\mathcal{O}_{\mathrm{FF}_{\widetilde{R}_\psi,A}}$-modules over the deformed version of the Fargues-Fontaine curve. Then we have the corresponding category of all the chain complexes namely the category $Ch\mathrm{Mod}\mathcal{O}_{\mathrm{FF}_{\widetilde{R}_\psi,A}}$ and the corresponding derived one $D\mathrm{Mod}\mathcal{O}_{\mathrm{FF}_{\widetilde{R}_\psi,A}}$.
\end{setting}

\begin{definition}
A complex in the category $Ch\mathrm{Mod}\mathcal{O}_{\mathrm{FF}_{\widetilde{R}_\psi,A}}$ taking the form of 
\[
\xymatrix@R+0pc@C+0pc{
\ar[r]\ar[r]\ar[r] &... \ar[r]\ar[r]\ar[r] &M^{n-1} \ar[r]\ar[r]\ar[r] &M^n
\ar[r]\ar[r]\ar[r] &M^{n+1}
\ar[r]\ar[r]\ar[r] &...
}
\]		
is now called pseudocoherent if it is quasi-isomorphic to a bounded above complex of finite projective objects in the category of all the pseudocoherent objects defined over the deformed version of adic version of the Fargues-Fontaine curve in the current situation.
Correspondingly we use the notation $D_\mathrm{pseudo,alg}\mathrm{FF}_{\widetilde{R}_\psi,A}$ to denote the corresponding category of all the pseudocoherent complexes defined above. We also have the corresponding category as above $D^\mathrm{cb}_\mathrm{pseudo,alg}\mathrm{FF}_{\widetilde{R}_\psi,A}$. 
A complex in the category $Ch\mathrm{Mod}\mathcal{O}_{\mathrm{FF}_{\widetilde{R}_\psi,A}}$ taking the form of 
\[
\xymatrix@R+0pc@C+0pc{
\ar[r]\ar[r]\ar[r] &... \ar[r]\ar[r]\ar[r] &M^{n-1} \ar[r]\ar[r]\ar[r] &M^n
\ar[r]\ar[r]\ar[r] &M^{n+1}
\ar[r]\ar[r]\ar[r] &...
}
\]		
is now called perfect if it is quasi-isomorphic to a bounded complex of finite projective objects in the category of all the pseudocoherent objects defined over the deformed version of adic version of the Fargues-Fontaine curve in the current situation.
Correspondingly we use the notation $D_\mathrm{perf,alg}\mathrm{FF}_{\widetilde{R}_\psi,A}$ to denote the corresponding category of all the pseudocoherent complexes defined above. We also have the corresponding categories of dg-flat complexes and the strictly perfect complexes which will be denoted by $D^\mathrm{dg-flat}_\mathrm{perf,alg}\mathrm{FF}_{\widetilde{R}_\psi,A}$ and $D^\mathrm{str}_\mathrm{perf,alg}\mathrm{FF}_{\widetilde{R}_\psi,A}$. 
\end{definition}

\begin{proposition}
The categories 
\begin{align}
&D_\mathrm{pseudo,alg}\mathrm{FF}_{\widetilde{R}_\psi,A},D^\mathrm{cb}_\mathrm{pseudo,alg}\mathrm{FF}_{\widetilde{R}_\psi,A},\\
&D_\mathrm{perf,alg}\mathrm{FF}_{\widetilde{R}_\psi,A},D^\mathrm{dg-flat}_\mathrm{perf,alg}\mathrm{FF}_{\widetilde{R}_\psi,A},D^\mathrm{str}_\mathrm{perf,alg}\mathrm{FF}_{\widetilde{R}_\psi,A}	
\end{align}
admit Waldhausen structure.	
\end{proposition}

\begin{proof}
See the proof before where we do not have the corresponding deformation with respect to the algebra $A$.
\end{proof}

\begin{setting}
We use the notation $QCoh\mathrm{FF}_{\widetilde{R}_\psi,A}$ to denote the corresponding quasi-coherent sheaves of $\mathcal{O}_{\mathrm{FF}_{\widetilde{R}_\psi,A}}$-modules (direct limits of the pseudocoherent ones which are assumed to form an abelian category) over the adic space $\mathrm{FF}_{\widetilde{R}_\psi,A}$. Then we use the notation $ChQCoh\mathrm{FF}_{\widetilde{R}_\psi,A}$ to denote the category of chain complexes consisting of all the objects in the category $QCoh\mathrm{FF}_{\widetilde{R}_\psi,A}$. And we use the notation $DQCoh\mathrm{FF}_{\widetilde{R}_\psi,A}$ to denote the corresponding derived category. 
\end{setting}

\begin{definition}
A complex in the category $ChQCoh\mathrm{FF}_{\widetilde{R}_\psi,A}$ taking the form of 
\[
\xymatrix@R+0pc@C+0pc{
\ar[r]\ar[r]\ar[r] &... \ar[r]\ar[r]\ar[r] &M^{n-1} \ar[r]\ar[r]\ar[r] &M^n
\ar[r]\ar[r]\ar[r] &M^{n+1}
\ar[r]\ar[r]\ar[r] &...
}
\]
is called then pseudocoherent if it is quasi-isomorphic to a bounded above complex of finite projective objects in the category of all the pseudocoherent objects over the deformed version of the adic version of the Fargues-Fontaine curve. The corresponding whole category of all such objects is denoted now by the notation $D_{\mathrm{pseudo}}\mathrm{FF}_{\widetilde{R}_\psi,A}$. We also have the category as above $D^{\mathrm{cb}}_{\mathrm{pseudo}}\mathrm{FF}_{\widetilde{R}_\psi,A}$.	
\end{definition}

\begin{definition}
A complex in the category $ChQCoh\mathrm{FF}_{\widetilde{R}_\psi,A}$ taking the form of 
\[
\xymatrix@R+0pc@C+0pc{
\ar[r]\ar[r]\ar[r] &... \ar[r]\ar[r]\ar[r] &M^{n-1} \ar[r]\ar[r]\ar[r] &M^n
\ar[r]\ar[r]\ar[r] &M^{n+1}
\ar[r]\ar[r]\ar[r] &...
}
\]
is called then perfect if it is quasi-isomorphic to a bounded complex of finite projective objects in the category of all the pseudocoherent objects over the deformed version of the adic version of the Fargues-Fontaine curve. The corresponding whole category of all such objects is denoted now by the notation $D_{\mathrm{perf}}\mathrm{FF}_{\widetilde{R}_\psi,A}$. We also have the categories of dg-flat complexes and strictly perfect complexes, namely $D^{\mathrm{dg-flat}}_{\mathrm{perf}}\mathrm{FF}_{\widetilde{R}_\psi,A}$ and $D^{\mathrm{str}}_{\mathrm{perf}}\mathrm{FF}_{\widetilde{R}_\psi,A}$.	
\end{definition}

\begin{remark}
The consideration here is in some sense different from the corresponding consideration in the schematic Fargues-Fontaine curve situation where we use the machinery of scheme theory and it is sometimes not stable under the analytic operations. But the latter should have its own interests from algebraic geometric point of view, in the general situation. For instance this gives us the chance to study the corresponding mysterious algebraic modules over the Robba rings with Frobenius pullbacks as isomorphisms. We remind the readers of the fact that over the cyclotomic tower many difficulties disappears. For instance more general algebraic approaches could be available.
\end{remark}

\begin{proposition}
The categories 
\begin{align}
&D_\mathrm{pseudo}\mathrm{FF}_{\widetilde{R}_\psi,A},D^\mathrm{cb}_\mathrm{pseudo}\mathrm{FF}_{\widetilde{R}_\psi,A},\\
&D_\mathrm{perf}\mathrm{FF}_{\widetilde{R}_\psi,A},D^\mathrm{dg-flat}_\mathrm{perf}\mathrm{FF}_{\widetilde{R}_\psi,A},D^\mathrm{str}_\mathrm{perf}\mathrm{FF}_{\widetilde{R}_\psi,A}	
\end{align}
admit Waldhausen structure.		
\end{proposition}

\begin{proof}
See the proof of the previous propositions around the same issue.	
\end{proof}

\begin{conjecture}
Assume that $\psi$ is the cyclotomic tower. Then the corresponding derived section functor induces a map on the $K$-theory space $\mathbb{K}D^\mathrm{dg-flat}_{\mathrm{perf}}\mathrm{FF}_{\widetilde{R}_\psi,A}$ to that of the category $D^\mathrm{dg-flat}_{\mathrm{perf}}(A)$ of all the dg-flat perfect complexes of $A$-modules. We conjecture that this is homotopic to zero.
 
Assume that $\psi$ is the cyclotomic tower. Then the corresponding derived section functor induces a map on the $K$-theory space $\mathbb{K}D^\mathrm{str}_{\mathrm{perf}}\mathrm{FF}_{\widetilde{R}_\psi,A}$ to that of the category $D^\mathrm{str}_{\mathrm{perf}}(A)$ of all the strictly perfect complexes of $A$-modules. We conjecture that this is homotopic to zero. 
\end{conjecture}

\begin{remark}
These are conjectured in the situation where $\psi$ is just the cyclotomic tower, in general one could easily guess what the conjectures look like.	
\end{remark}

\indent We make some discussion here on the picture in the \'etale setting. The conjecture is much easier in the $\ell$-adic situation. This represents the difference between the $\ell$-adic situation and the $p$-adic setting in some situation. For instance let us try to consider the following proposition after Witte (see \cite[Proposition 6.1.5]{2Wit3}):

\begin{proposition}
Let $T$ be $\mathbb{Z}_\ell$ where $\ell\neq p$. For $X$ a rigid analytic space over the $p$-adic field $\mathbb{Q}_p$ which is separated and of finite type (when considered as some adic space locally of finite type, see \cite{2Huber1}), suppose we have that for $\sharp=\text{\'et}$ the category $\mathbb{D}_{\mathrm{perf}}(X_\sharp,T)$ could be endowed with the structure of Waldhausen categories and there is a Waldhausen exact functor $R\Gamma_?(X_\sharp,.)$ (induced by the direct image functor in this context as in \cite[Chapter 4-5]{2Wit3}) as below for $\sharp=\text{\'et}$ and $?=\emptyset$:
\[
\xymatrix@R+0pc@C+3pc{
\mathbb{D}_{\mathrm{perf}}(X_\sharp,T)\ar[r]^{R\Gamma_?(X_\sharp,.)}\ar[r]\ar[r] &\mathbb{D}_{\mathrm{perf}}(T)
}
\]
which induces the corresponding the morphism $\mathbb{K}R\Gamma_?(X_\sharp,.)$ as:
\[
\xymatrix@R+0pc@C+3pc{
\mathbb{K}\mathbb{D}_{\mathrm{perf}}(X_\sharp,T)\ar[r]^{\mathbb{K}R\Gamma_?(X_\sharp,.)}\ar[r]\ar[r] & \mathbb{K}\mathbb{D}_{\mathrm{perf}}(T).
}
\]	
Then there is a homotopy between this morphism and the zero map:
\[
\xymatrix@R+0pc@C+3pc{
\mathbb{K}\mathbb{D}_{\mathrm{perf}}(X_\sharp,T)\ar[r]^{0}\ar[r]\ar[r] & \mathbb{K}\mathbb{D}_{\mathrm{perf}}(T).
}
\]
\end{proposition}

\begin{proof}
We follow the strategy of \cite[Proposition 6.1.5]{2Wit3} to prove this. The idea is to pullback all the things back along the following morphism:
\begin{displaymath}
f:\mathrm{Spa}(\mathbb{Q}^\mathrm{ur}_p,\mathfrak{o}_{\mathbb{Q}^\mathrm{ur}_p})\rightarrow \mathrm{Spa}(\mathbb{Q}_p,\mathfrak{o}_{\mathbb{Q}_p}).	
\end{displaymath}
Now for each complex $M^\bullet$ involved we put:
\begin{align}
R\Gamma(X_{\mathbb{Q}^\mathrm{ur}_p,\sharp},M^\bullet):=\Gamma(\mathrm{Spa}(\mathbb{Q}^\mathrm{ur}_p,\mathfrak{o}_{\mathbb{Q}^\mathrm{ur}_p}),f^*R\pi_* M^\bullet),
%R\Gamma_c(X_{\mathbb{Q}^\mathrm{ur}_p,\sharp},M^\bullet):=\Gamma(\mathrm{Spa}(\mathbb{Q}^\mathrm{ur}_p,\mathfrak{o}_{\mathbb{Q}^\mathrm{ur}_p}),f^*R\pi_! M^\bullet),
\end{align}
where the morphism $\pi:X\rightarrow \mathrm{Spa}(\mathbb{Q}_p,\mathfrak{o}_{\mathbb{Q}_p})$ is the corresponding structure morphism of $X$. We have the corresponding finiteness due to our assumption here. One then looks at the corresponding Frobenius $\mathrm{Fr}$ over the corresponding field $\mathbb{Q}^\mathrm{ur}_p$. Then as in \cite[Proposition 6.1.5]{2Wit3} we have that then the following cofibration sequence for each complex $M^\bullet$:
\[
\xymatrix@R+0pc@C+3pc{
R\Gamma_?(X_{\mathbb{Q}_p,\sharp},M^\bullet)\ar[r]\ar[r]\ar[r] &R\Gamma_?(X_{\mathbb{Q}^\mathrm{ur}_p,\sharp},M^\bullet)\ar[r]^{1-\mathrm{Fr}}\ar[r] &R\Gamma_?(X_{\mathbb{Q}^\mathrm{ur}_p,\sharp},M^\bullet),
}
\]
as in \cite[Proposition 6.1.5]{2Wit3} (note that here we need the corresponding flabbiness of the corresponding derived direct images), by the corresponding Waldhausen's additive theorem we have the corresponding homotopy:
\begin{align}
\mathbb{K}R\Gamma_?(X_{\mathbb{Q}^\mathrm{ur}_p,\sharp},M^\bullet)\oplus \mathbb{K}R\Gamma_?(X_{\mathbb{Q}_p,\sharp},M^\bullet)\leadsto \mathbb{K}R\Gamma_?(X_{\mathbb{Q}^\mathrm{ur}_p,\sharp},M^\bullet)	
\end{align}
which implies that we have then as in \cite[Proposition 6.1.5]{2Wit3}
\begin{align}
\pi_i\mathbb{K}R\Gamma_?(X_\sharp,.)=0.	
\end{align}
\end{proof}

\indent We now show how this recovers the results of Witte at a point:

\begin{corollary}
Suppose now the space $X$ is just the point $\mathrm{Spa}(\mathbb{Q}_p,\mathfrak{o}_{\mathbb{Q}_p})$, and keep the assumption that $T$ is the ring $\mathbb{Z}_\ell$. Then we have for $\sharp=\text{\'et}$ the Waldhausen category $\mathbb{D}_{\mathrm{perf}}(X_\sharp,T)$ which induces the Waldhausen exact functor (induced from the total derived direct image functor) $R\Gamma_?(X_\sharp,.)$ where $?=\emptyset$ and $\sharp$ as above:
\[
\xymatrix@R+0pc@C+3pc{
\mathbb{D}_{\mathrm{perf}}(X_\sharp,T)\ar[r]^{R\Gamma_?(X_\sharp,.)}\ar[r]\ar[r] &\mathbb{D}_{\mathrm{perf}}(T)
}
\]
which induces the corresponding the morphism $\mathbb{K}R\Gamma_?(X_\sharp,.)$ as:
\[
\xymatrix@R+0pc@C+3pc{
\mathbb{K}\mathbb{D}_{\mathrm{perf}}(X_\sharp,T)\ar[r]^{\mathbb{K}R\Gamma_?(X_\sharp,.)}\ar[r]\ar[r] & \mathbb{K}\mathbb{D}_{\mathrm{perf}}(T).
}
\]	
And moreover we have that this map is homotopic to zero:
\[
\xymatrix@R+0pc@C+3pc{
\mathbb{K}\mathbb{D}_{\mathrm{perf}}(X_\sharp,T)\ar[r]^{0}\ar[r]\ar[r] & \mathbb{K}\mathbb{D}_{\mathrm{perf}}(T).
}
\]
	
\end{corollary}

\begin{proof}
One can apply results of \cite{2Wit3} to recover the corresponding structures here or construct the corresponding structures here directly. Then we can apply the previous proposition. 
\end{proof}

\indent Let us analyze the corresponding $p$-adic setting in the following sense. Since in the proof used above the corresponding results on the finiteness the cohomological dimension over the space attached to the maximal unramified extension of $\mathbb{Q}_p$ are crucial. We first consider the integral version of the picture we considered above.

\begin{setting}
We are going to use the notation $QCoh\Omega^{\mathrm{int}}\Phi^a_{X_{\text{pro-\'etale}}}$ (where $X$ is assumed to be separated and of finite type as an adic space) to denote the corresponding quasi-coherent (direct limits of pseudocoherent ones) $\varphi^a$-sheaves over the sheaf of ring $\underline{\underline{{\Omega}}}^\mathrm{int}_{X_{\text{pro-\'etale}}}$, and we use the corresponding notation $ChQCoh\Omega^{\mathrm{int}}\Phi^a_{X_{\text{pro-\'etale}}}$ to denote the corresponding category of the chain complexes and use the corresponding notation $DQCoh\Omega^{\mathrm{int}}\Phi^a_{X_{\text{pro-\'etale}}}$ to denote the corresponding derived category.	
\end{setting}

\begin{definition}
A complex considered in the previous setting taking the form of 
\[
\xymatrix@R+0pc@C+0pc{
\ar[r]\ar[r]\ar[r] &... \ar[r]\ar[r]\ar[r] &M^{n-1} \ar[r]\ar[r]\ar[r] &M^n
\ar[r]\ar[r]\ar[r] &M^{n+1}
\ar[r]\ar[r]\ar[r] &...
}
\]	
is called then pseudocoherent if it is quasi-isomorphic to a bounded above complex of finite projective objects in the category of pseudocoherent $\varphi^a$-modules over $\underline{\underline{{\Omega}}}^\mathrm{int}_{X}$.	
\end{definition}

%\begin{conjecture}
%If we denote the category of all the pseudocoherent complexes in the previous definition by $D_{\mathrm{pseudo}}\Omega^{\mathrm{int}}\Phi^a_{X_{\text{pro-\'etale}}}$. Then we have that the category admits Waldhausen structure.	
%\end{conjecture}

%\begin{proof}
%See the proof above in the rational situation.	
%\end{proof}

%\begin{conjecture}
%We conjecture that the corresponding derived section functor in this context induces a map from the $K$-theory space $\mathbb{K}D_{\mathrm{pseudo}}\Omega^{\mathrm{int}}\Phi^a_{X_{\text{pro-\'etale}}}$	to the $K$-theory space $\mathbb{K}D_{\mathrm{pseudo}}(\mathbb{Q}_p)$. And we conjecture that this map is homotopic to zero in the situation where $X$ is just the point attached to $\mathbb{Q}_p$.
%\end{conjecture}

The corresponding setting of local systems gives us the following consideration after Witte:

\begin{conjecture}
Assume that $T$ is an adic ring over $\mathbb{Z}_p$ such that we can find a two-sided ideal $I$ such that each $T/I^n$ for $n\geq 0$ is finite of order a power of $p$. Suppose that $X$ is a rigid analytic space which is assumed to be separated and of finite type over the $p$-adic number field $\mathbb{Q}_p$. Assume that the corresponding category $\mathbb{D}_{\mathrm{perf}}(X_\sharp,T)$ for $\sharp=\text{\'et}$ admits Waldhausen structure and we have the corresponding Waldhausen exact functor (induced from the corresponding direct image functor as in \cite[Chapter 4-5]{2Wit3}) $R\Gamma_?(X_\sharp,.)$ where $?=\emptyset$:
 	\[
\xymatrix@R+0pc@C+3pc{
\mathbb{D}_{\mathrm{perf}}(X_\sharp,T)\ar[r]^{R\Gamma_?(X_\sharp,.)}\ar[r]\ar[r] &\mathbb{D}_{\mathrm{perf}}(T)
}
\]
which induces the corresponding the morphism $\mathbb{K}R\Gamma_?(X_\sharp,.)$ as:
\[
\xymatrix@R+0pc@C+3pc{
\mathbb{K}\mathbb{D}_{\mathrm{perf}}(X_\sharp,T)\ar[r]^{\mathbb{K}R\Gamma_?(X_\sharp,.)}\ar[r]\ar[r] & \mathbb{K}\mathbb{D}_{\mathrm{perf}}(T).
}
\]	
Then we have that this map is homotopic to the zero map:
\[
\xymatrix@R+0pc@C+3pc{
\mathbb{K}\mathbb{D}_{\mathrm{perf}}(X_\sharp,T)\ar[r]^{0}\ar[r]\ar[r] & \mathbb{K}\mathbb{D}_{\mathrm{perf}}(T).
}
\]

\end{conjecture}

\begin{remark}	

The proof should be actually the same as in the previous result in the corresponding $\ell$-adic situation, the difference here is that one uses the corresponding finiteness and cohomological dimension results in the $p$-adic setting by applying \cite[8.1,10.1]{2KL} for instance, which is to say in our setting the corresponding finiteness and cohomological dimension results over the point attached to $\mathbb{Q}_p^\mathrm{ur}$ would produce us the corresponding cofibration sequence involved. Also we need to consider a general ring $T$ here. Following the idea in \cite[Proposition 2.1.3]{2FK}, we quotient out the corresponding Jacobson radical of each finite quotient $T/I$ which is a finite ring. We can then reduce to semi-simple case, and then reduce to finite field case. 
\end{remark}

\begin{conjecture}
When $T$ is as in the previous proposition the statement in the previous corollary holds as well in our current situation, which is to say over $\mathrm{Spa}(\mathbb{Q}_p,\mathfrak{o}_{\mathbb{Q}_p})$ we have the Waldhausen category $\mathbb{D}_{\mathrm{perf}}(X_\sharp,T)$ for $\sharp=\text{\'et}$ with the corresponding Waldhausen exact functor $R\Gamma_?(X_\sharp,.)$ for $?=\emptyset$:
\[
\xymatrix@R+0pc@C+3pc{
\mathbb{D}_{\mathrm{perf}}(X_\sharp,T)\ar[r]^{R\Gamma_?(X_\sharp,.)}\ar[r]\ar[r] &\mathbb{D}_{\mathrm{perf}}(T)
}
\]
which induces the corresponding the morphism $\mathbb{K}R\Gamma_?(X_\sharp,.)$ as:
\[
\xymatrix@R+0pc@C+3pc{
\mathbb{K}\mathbb{D}_{\mathrm{perf}}(X_\sharp,T)\ar[r]^{\mathbb{K}R\Gamma_?(X_\sharp,.)}\ar[r]\ar[r] & \mathbb{K}\mathbb{D}_{\mathrm{perf}}(T).
}
\]	
And we have that this map is homotopic to the zero map:
\[
\xymatrix@R+0pc@C+3pc{
\mathbb{K}\mathbb{D}_{\mathrm{perf}}(X_\sharp,T)\ar[r]^{0}\ar[r]\ar[r] & \mathbb{K}\mathbb{D}_{\mathrm{perf}}(T).
}
\]
\end{conjecture}

%\begin{proof}
%See the proof of the previous corollary. 	
%\end{proof}

\begin{remark}
We used the corresponding context of Waldhausen categories after Witte, but certain cases of the corresponding $\epsilon$-isomorphism conjectures in the \'etale setting are due to many people who studied the corresponding $\epsilon$-isomorphism conjectures, mainly in \cite{2FK}, \cite{2F},\cite{2BB}, \cite{2LVZ} and \cite{2Nak3}.	
\end{remark}

\

This chapter is based on the following paper, where the author of this dissertation is the main author:
\begin{itemize}
\item Tong, Xin. "Hodge-Iwasawa Theory I." arXiv preprint arXiv:2006.03692 (2020). 
\end{itemize}

\newpage

\bibliographystyle{ams}

%%%%%%\newpage

%%%%%%\newpage

\newpage\chapter{Hodge-Iwasawa Theory II}

\newpage\section{Introduction}

\subsection{Introduction}

\indent In our previous paper \cite{3T1} on the Hodge-Iwasawa theory we studied some deformation of $p$-adic Hodge structures over some general perfectoid spaces and adic spaces. We mainly focused on the corresponding perfect setting which to some extent gives us some suitable and desired functoriality with respect to some general adic spaces. This reflects the corresponding motivation in some geometric and categorical way since our initial goal on the corresponding Hodge-Iwasawa consideration is to establish some possibility to deform the functorial geometric structures.\\

\indent In our current situation we still emphasize higher dimensional geometric structures. We will generalize the corresponding context of \cite{3KL16} to our context which on one hand produces more general tools to study the deformation of Hodge structures and on the other hand we will be more flexible over some specific towers both coming from geometry and Iwasawa theory. The context of \cite{3KL16} is already over some deep generalization. First off the towers considered in \cite{3KL16} are axiomatic before some more concrete study. Also the corresponding context is within the consideration of generalized Witt vectors. Also the context considers the corresponding pseudocoherent modules and sheaves instead of the objects initially considered in the context of \cite{3KL15}. We consider the corresponding deformations of these generalized structures over some affinoid algebras.\\

\indent In our current consideration we will do some kind of globalization to some extent which will generalize the corresponding story in some more rigid framework. To be more precise we will consider the situation where the deformation happens over some affinoid rigid analytic space. The sheafiness will be non-trivial here but will be true in this specific context.\\

\indent We will work in more general setting than the settings of our previous work and \cite{3KL16}. Although the setting of \cite{3KL16} is already more general than \cite{3KL15}, we choose to work in the setting where the coefficient field $E$ is ultrametric with perfect residue field of characteristic $p>0$ (even in the nondiscrete valued situation). When $E$ is not discrete valued this is actually quite hard to deal with, where we only consider the situation where we do have some stable unique expression in the variable $\pi$ (for elements in the ring of Witt vectors), which we hope go back to this later to work within more general setting as long as one can define the corresponding power-multiplicative norms to do various completion under the hypothesis that $R$ is at least uniform Banach. Therefore our work will provide some additional reference to the study on this level of generality. Please be careful enough on our notations.

\subsection{Main Results}

\indent We have managed to do some generalization of our previous work \cite{3T1}, and the corresponding context in \cite{3KL15}, \cite{3KL16}, \cite{3KP}. As mentioned in the previous subsection the corresponding imperfect setting of the Hodge-Iwasawa theory is one main goal of our development. Also we would like to discuss to some extent the corresponding globalization of the deformed information which the rings in our consideration are carrying. To summarize we have:\\

I. We have now considered the corresponding equal characteristic analog of the corresponding comparison between the vector bundles, pseudocoherent sheaves over the deformed schematic Fargues-Fontaine curves and the corresponding Frobenius modules over the deformed Robba rings, generalizing \cite{3T1}, \cite{3KL16} and \cite{3KP}. The corresponding motivation comes from certainly not only \cite{3KP}, but also the corresponding $t$-motivic Hodge theory in the relative setting such as in Hartl-Kim \cite{3HK}.\\

II. We have now considered the corresponding equal characteristic analog of the corresponding comparison between the vector bundles, pseudocoherent sheaves over the deformed adic Fargues-Fontaine curves and the corresponding Frobenius modules over the deformed Robba rings \cite{3T1}, \cite{3KL15}, \cite{3KL16}. Again the corresponding motivation comes from certainly not only \cite{3KP}, but also the corresponding $t$-motivic Hodge theory in the relative setting such as in Hartl-Kim \cite{3HK}. This will lead to some sort of globalization in the category of adic spaces.\\ 

%III. We include some \'etaleness spreading properties after \cite{3KL}, over the corresponding rigid families of Frobenius modules mainly in the equal characteristic setting. This is definitely motivated from the work from Hartl-Kim \cite{3HK}. Hartl-Kim proposed that one should be able to study the analog of Hellmann's stacks in \cite{3Hel1}, so it is interesting and should be very significant to study the corresponding \'etale locus living in the corresponding moduli stack of the corresponding isocrystals, where the corresponding analog of Rapoport-Zink conjecture should be formulated. \\

III. We have discussed the corresponding relationship between Frobenius modules over the perfect period rings and the Frobenius modules over the corresponding imperfection of the period rings in our deformed setting. This is to some extent important since the corresponding imperfection might lead to the corresponding noetherian objects which might be easy to control. This will have further application to the study of our deformation of the \'etale local systems.\\

IV. We have discussed the corresponding relationship between Frobenius modules over the perfect period rings and the Frobenius modules over the corresponding imperfection of the period rings in our deformed setting, but over the corresponding noncommutative Banach rings. This is to some extent important since the corresponding imperfection might lead to the corresponding noncommutative noetherian objects which might be easy to control. This will have further application to the study of our noncommutative deformation of the \'etale local systems. It is very natural to consider noncommutative setting in our development since we are considering the corresponding deformation into the noncommutative Fr\'echet-Stein algebras.\\

\subsection{Further Consideration}

We should mention that actually the corresponding ideas we mentioned within the globalization are not quite new since this bears some similarity to Pappas-Rapoport \cite{3PR1} and Hellmann's work \cite{3Hel1} on the arithmetic stacks of Frobenius modules. The corresponding comparison in our context could be morally conveyed in the following diagram:
\[
\xymatrix@R+6pc@C+0pc{
\text{Pappas-Rapoport-Hellmann}
\ar[r]\ar[r]\ar[r]\ar[d]\ar[d]\ar[d] &\text{Rapoport-Zink Spaces}\ar[l]\ar[l]\ar[l] \ar[d]\ar[d]\ar[d]\\
\text{???} \ar[r]\ar[r]\ar[r] &\text{Moduli Spaces of Local Shtukas after Scholze} \ar[l]\ar[l]\ar[l] .\\
}
\]

For instance in our consideration we will consider some arithmetic deformation of vector bundles over the arithmetic deformation of Fargues-Fontaine curves. It is very natural to consider the corresponding moduli problem around the corresponding vector bundles over these deformed spaces. Certainly these are not directly diamonds, but at least they should be stacks over the corresponding categories of certain adic spaces.\\

%\begin{conjecture}
%The moduli space $M_{\mathrm{FF}}$ of the corresponding $\mathrm{GL}_n$-bundle over the adic Fargues-Fontaine curve (in the setting where the base is defined from $\mathbb{Q}_p$ or $\mathbb{F}_p((t))$) in equal characteristic case or mixed characteristic case is an Artin stack fibered over the category of rigid analytic spaces in the fpqc topology. And one should have the corresponding openness result on the suitable \'etale locus.	
%\end{conjecture}

In this paper, we have built some well-posed equivariant versions of many results in the relative $p$-adic Hodge theory in \cite{3KL15} and \cite{3KL16} which have already reached beyond our previous consideration. We expect everything could be further applied to the corresponding study on the cohomologies of all the corresponding objects we defined and studied here. For instance if one applies these to the $t$-motivic setting, then immediately one has the chance to construct the corresponding $t$-adic local Tamagawa number conjecture after Nakamura \cite{3Nakamura1} (also inspired by the recent work \cite{3FGHP}), which certainly is a relative version in rigid family (again after Nakamura relying on the results from Kedlaya-Pottharst-Xiao on the finiteness of the corresponding cohomology of $(\varphi,\Gamma)$-modules over relative Robba rings).\\

We have reached some noncommutative Iwasawa deformation of the relative $p$-adic Hodge structures. At least the corresponding hope is directly targeted at the corresponding noncommutative Tamagawa number conjectures after Burns-Flach-Fukaya-Kato \cite{3BF1}, \cite{3BF2} and \cite{3FK1},  Nakamura and Z\"ahringer \cite{3Zah1}. Our belief is that although we have not shown that coherent sheaves and modules over noncommutative deformation of full Robba rings are equivalent, but at least for coherent sheaves one could have a coherent Iwasawa theory in the noncommutative setting. \\

%%\newpage

\subsection{Lists of Notations}

\indent The arrangement of the corresponding notations getting involved in this paper is bit complicated, so we have made some indication to indicate where they mainly emerge into the corresponding discussion.\\

\begin{center}
\begin{longtable}{p{7.0cm}p{8cm}}
Notation & Description (mainly in section 2, section 3, setion 4, section 6.1, section 6.2, setion 6.3) \\
\hline
$E$ & A complete discrete valued field which is nonarchimedean with residue field $k$ perfect in characteristic $p>0$.\\
$A$ & A reduced affinoid algebra in rigid analytic geometry in the sense of Tate, or a noncommutative Banach affinoid algebra.\\
$W_\pi$ & The ring of generalized Witt vectors with respect to some finite extension $E$ of $\mathbb{Q}_p$ with pseudo-uniformizer $\pi$. \\
$(R,R^+)$   & Adic perfect uniform algebra over $\mathbb{F}_{p^h}$ for some $h\geq 1$. \\
$\widetilde{\Omega}_{R,A}^\mathrm{int}$ & Deformed version of the period rings in the style of \cite{3KL15} and \cite{3KL16}.\\
$\widetilde{\Omega}_{R,A}$ & Deformed version of the period rings in the style of \cite{3KL15} and \cite{3KL16}.\\
$\widetilde{\Pi}_{R,A}^{\mathrm{int},r}$ & Deformed version of the period rings in the style of \cite{3KL15} and \cite{3KL16}.\\
$\widetilde{\Pi}_{R,A}^\mathrm{int}$ & Deformed version of the period rings in the style of \cite{3KL15} and \cite{3KL16}.\\
$\widetilde{\Pi}_{R,A}^{\mathrm{bd},r}$ & Deformed version of the period rings in the style of \cite{3KL15} and \cite{3KL16}.\\
$\widetilde{\Pi}_{R,A}^\mathrm{bd}$ & Deformed version of the period rings in the style of \cite{3KL15} and \cite{3KL16}.\\
$\widetilde{\Pi}_{R,A}^{r}$ & Deformed version of the period rings in the style of \cite{3KL15} and \cite{3KL16}.\\
$\widetilde{\Pi}_{R,A}$ & Deformed version of the period rings in the style of \cite{3KL15} and \cite{3KL16}.\\
$\widetilde{\Pi}_{R,A}^{\infty}$ & Deformed version of the period rings in the style of \cite{3KL15} and \cite{3KL16}.\\
$\widetilde{\Pi}_{R,A}^I$ & Deformed version of the period rings in the style of \cite{3KL15} and \cite{3KL16}.\\

$\widetilde{\Omega}_{*,A}^\mathrm{int}$ & Deformed version of the period sheaves in the style of \cite{3KL15} and \cite{3KL16}.\\
$\widetilde{\Omega}_{*,A}$ & Deformed version of the period sheaves in the style of \cite{3KL15} and \cite{3KL16}.\\
$\widetilde{\Pi}_{*,A}^{\mathrm{int},r}$ & Deformed version of the period sheaves in the style of \cite{3KL15} and \cite{3KL16}.\\
$\widetilde{\Pi}_{*,A}^\mathrm{int}$ & Deformed version of the period sheaves in the style of \cite{3KL15} and \cite{3KL16}.\\
$\widetilde{\Pi}_{*,A}^{\mathrm{bd},r}$ & Deformed version of the period sheaves in the style of \cite{3KL15} and \cite{3KL16}.\\
$\widetilde{\Pi}_{*,A}^\mathrm{bd}$ & Deformed version of the period sheaves in the style of \cite{3KL15} and \cite{3KL16}.\\
$\widetilde{\Pi}_{*,A}^{r}$ & Deformed version of the period sheaves in the style of \cite{3KL15} and \cite{3KL16}.\\
$\widetilde{\Pi}_{*,A}$ & Deformed version of the period sheaves in the style of \cite{3KL15} and \cite{3KL16}.\\
$\widetilde{\Pi}_{*,A}^{\infty}$ & Deformed version of the period sheaves in the style of \cite{3KL15} and \cite{3KL16}.\\
$\widetilde{\Pi}_{*,A}^I$ & Deformed version of the period sheaves in the style of \cite{3KL15} and \cite{3KL16}.\\

$W_{\pi,\infty}$ & The ring of generalized Witt vectors with respect to some finite extension $E$ of $\mathbb{Q}_p$ with pseudo-uniformizer $\pi$, but with base change to $E_\infty$. In this case each element admits unique expression $\sum_{n\in \mathbb{Z}[1/p]_{\geq 0}}\pi^n[\overline{x}_n]$. \\

$\widetilde{\Omega}_{R,\infty}^\mathrm{int},\widetilde{\Omega}_{*,\infty}^\mathrm{int}$ & Period rings or sheaves in the style of \cite{3KL15} and \cite{3KL16}, with base change to $E_\infty$.\\
$\widetilde{\Omega}_{R,\infty},\widetilde{\Omega}_{*,\infty}$ & Period rings or sheaves in the style of \cite{3KL15} and \cite{3KL16}, with base change to $E_\infty$.\\
$\widetilde{\Pi}_{R,\infty}^{\mathrm{int},r},\widetilde{\Pi}_{*,\infty}^{\mathrm{int},r}$ & Period rings or sheaves in the style of \cite{3KL15} and \cite{3KL16}, with base change to $E_\infty$.\\
$\widetilde{\Pi}_{R,\infty}^\mathrm{int},\widetilde{\Pi}_{*,\infty}^\mathrm{int}$ & Period rings or sheaves in the style of \cite{3KL15} and \cite{3KL16}, with base change to $E_\infty$.\\
$\widetilde{\Pi}_{R,\infty}^{\mathrm{bd},r},\widetilde{\Pi}_{*,\infty}^{\mathrm{bd},r}$ & Period rings or sheaves in the style of \cite{3KL15} and \cite{3KL16}, with base change to $E_\infty$.\\
$\widetilde{\Pi}_{R,\infty}^\mathrm{bd},\widetilde{\Pi}_{*,\infty}^\mathrm{bd}$ & Period rings or sheaves in the style of \cite{3KL15} and \cite{3KL16}, with base change to $E_\infty$.\\
$\widetilde{\Pi}_{R,\infty}^{r},\widetilde{\Pi}_{*,\infty}^{r}$ & Period rings or sheaves in the style of \cite{3KL15} and \cite{3KL16}, with base change to $E_\infty$.\\
$\widetilde{\Pi}_{R,\infty},\widetilde{\Pi}_{*,\infty}$ & Period rings or sheaves in the style of \cite{3KL15} and \cite{3KL16}, with base change to $E_\infty$.\\
$\widetilde{\Pi}_{R,\infty}^{\infty},\widetilde{\Pi}_{*,\infty}^{\infty}$ & Period rings or sheaves in the style of \cite{3KL15} and \cite{3KL16}, with base change to $E_\infty$.\\
$\widetilde{\Pi}_{R,\infty}^I,\widetilde{\Pi}_{*,\infty}^I$ & Period rings or sheaves in the style of \cite{3KL15} and \cite{3KL16}, with base change to $E_\infty$.\\

$\widetilde{\Omega}_{R,\infty,A}^\mathrm{int},\widetilde{\Omega}_{*,\infty,A}^\mathrm{int}$ & Deformed period rings or sheaves in the style of \cite{3KL15} and \cite{3KL16}, with base change to $E_\infty$.\\
$\widetilde{\Omega}_{R,\infty,A},\widetilde{\Omega}_{*,\infty,A}$ & Deformed period rings or sheaves in the style of \cite{3KL15} and \cite{3KL16}, with base change to $E_\infty$.\\
$\widetilde{\Pi}_{R,\infty,A}^{\mathrm{int},r},\widetilde{\Pi}_{*,\infty,A}^{\mathrm{int},r}$ & Deformed period rings or sheaves in the style of \cite{3KL15} and \cite{3KL16}, with base change to $E_\infty$.\\
$\widetilde{\Pi}_{R,\infty,A}^\mathrm{int},\widetilde{\Pi}_{*,\infty,A}^\mathrm{int}$ & Deformed period rings or sheaves in the style of \cite{3KL15} and \cite{3KL16}, with base change to $E_\infty$.\\
$\widetilde{\Pi}_{R,\infty,A}^{\mathrm{bd},r},\widetilde{\Pi}_{*,\infty,A}^{\mathrm{bd},r}$ & Deformed period rings or sheaves in the style of \cite{3KL15} and \cite{3KL16}, with base change to $E_\infty$.\\
$\widetilde{\Pi}_{R,\infty,A}^\mathrm{bd},\widetilde{\Pi}_{*,\infty,A}^\mathrm{bd}$ & Deformed period rings or sheaves in the style of \cite{3KL15} and \cite{3KL16}, with base change to $E_\infty$.\\
$\widetilde{\Pi}_{R,\infty,A}^{r},\widetilde{\Pi}_{*,\infty,A}^{r}$ & Deformed period rings or sheaves in the style of \cite{3KL15} and \cite{3KL16}, with base change to $E_\infty$.\\
$\widetilde{\Pi}_{R,\infty,A},\widetilde{\Pi}_{*,\infty,A}$ & Deformed period rings or sheaves in the style of \cite{3KL15} and \cite{3KL16}, with base change to $E_\infty$.\\
$\widetilde{\Pi}_{R,\infty,A}^{\infty},\widetilde{\Pi}_{*,\infty,A}^{\infty}$ & Deformed period rings or sheaves in the style of \cite{3KL15} and \cite{3KL16}, with base change to $E_\infty$.\\
$\widetilde{\Pi}_{R,\infty,A}^I,\widetilde{\Pi}_{*,\infty,A}^I$ & Deformed period rings or sheaves in the style of \cite{3KL15} and \cite{3KL16}, with base change to $E_\infty$.\\

%$\widetilde{\Omega}_{R,\mathbb{Q}_p\{T_1,...,T_d\}}^\mathrm{int}$ & Deformed version of the period rings in the style of \cite{3KL15} and \cite{3KL16}.\\
%$\widetilde{\Omega}_{R,\mathbb{Q}_p\{T_1,...,T_d\}}$ & Deformed version of the period rings in the style of \cite{3KL15} and \cite{3KL16}.\\
%$\widetilde{\Pi}_{R,\mathbb{Q}_p\{T_1,...,T_d\}}^{\mathrm{int},r}$ & Deformed version of the period rings in the style of \cite{3KL15} and \cite{3KL16}.\\
%$\widetilde{\Pi}_{R,\mathbb{Q}_p\{T_1,...,T_d\}}^\mathrm{int}$ & Deformed version of the period rings in the style of \cite{3KL15} and \cite{3KL16}.\\
%$\widetilde{\Pi}_{R,\mathbb{Q}_p\{T_1,...,T_d\}}^{\mathrm{bd},r}$ & Deformed version of the period rings in the style of \cite{3KL15} and \cite{3KL16}.\\
%$\widetilde{\Pi}_{R,\mathbb{Q}_p\{T_1,...,T_d\}}^\mathrm{bd}$ & Deformed version of the period rings in the style of \cite{3KL15} and \cite{3KL16}.\\
%$\widetilde{\Pi}_{R,\mathbb{Q}_p\{T_1,...,T_d\}}^{r}$ & Deformed version of the period rings in the style of \cite{3KL15} and \cite{3KL16}.\\
%$\widetilde{\Pi}_{R,\mathbb{Q}_p\{T_1,...,T_d\}}$ & Deformed version of the period rings in the style of \cite{3KL15} and \cite{3KL16}.\\
%$\widetilde{\Pi}_{R,\mathbb{Q}_p\{T_1,...,T_d\}}^{\infty}$ & Deformed version of the period rings in the style of \cite{3KL15} and \cite{3KL16}.\\
%$\widetilde{\Pi}_{R,\mathbb{Q}_p\{T_1,...,T_d\}}^I$ & Deformed version of the period rings in the style of \cite{3KL15} and \cite{3KL16}.\\
%
%
\end{longtable}
\end{center}

\begin{center}
\begin{longtable}{p{6.6cm}p{8cm}}
Notation & Description (mainly in section 5, section 6.4, section 6.5, section 6.6, section 6.7) \\
\hline
$(H_\bullet,H^+_\bullet)$ & A tower as in \cite[Chapter 5]{3KL16}.\\
$\Pi^{\mathrm{int},r}_{H,A}$ & Imperfect relative period ring corresponding to $(H_\bullet,H^+_\bullet)$ which is deformed analog of Kedlaya-Liu's imperfect ring.\\
$\Pi^{\mathrm{int},\dagger}_{H,A}$ & Imperfect relative period ring corresponding to $(H_\bullet,H^+_\bullet)$ which is deformed analog of Kedlaya-Liu's imperfect ring.\\
$\Omega^\mathrm{int}_{H,A}$ & Imperfect relative period ring corresponding to $(H_\bullet,H^+_\bullet)$ which is deformed analog of Kedlaya-Liu's imperfect ring.\\
$\Omega_{H,A}$ & Imperfect relative period ring corresponding to $(H_\bullet,H^+_\bullet)$ which is deformed analog of Kedlaya-Liu's imperfect ring.\\
$\Pi^{\mathrm{bd},r}_{H,A}$ & Imperfect relative period ring corresponding to $(H_\bullet,H^+_\bullet)$ which is deformed analog of Kedlaya-Liu's imperfect ring.\\
$\Pi^{\mathrm{bd},\dagger}_{H,A}$ & Imperfect relative period ring corresponding to $(H_\bullet,H^+_\bullet)$ which is deformed analog of Kedlaya-Liu's imperfect ring.\\
$\Pi^{[s,r]}_{H,A}$ & Imperfect relative period ring corresponding to $(H_\bullet,H^+_\bullet)$ which is deformed analog of Kedlaya-Liu's imperfect ring.\\
$\Pi^r_{H,A}$ & Imperfect relative period ring corresponding to $(H_\bullet,H^+_\bullet)$ which is deformed analog of Kedlaya-Liu's imperfect ring.\\
$\Pi_{H,A}$ & Imperfect relative period ring corresponding to $(H_\bullet,H^+_\bullet)$ which is deformed analog of Kedlaya-Liu's imperfect ring.\\

$\breve{\Pi}^{\mathrm{int},r}_{H,A}$ & Imperfect relative period ring corresponding to $(H_\bullet,H^+_\bullet)$ which is deformed analog of Kedlaya-Liu's imperfect ring.\\
$\breve{\Pi}^{\mathrm{int},\dagger}_{H,A}$ & Imperfect relative period ring corresponding to $(H_\bullet,H^+_\bullet)$ which is deformed analog of Kedlaya-Liu's imperfect ring.\\
$\breve{\Omega}^\mathrm{int}_{H,A}$ & Imperfect relative period ring corresponding to $(H_\bullet,H^+_\bullet)$ which is deformed analog of Kedlaya-Liu's imperfect ring.\\
$\breve{\Omega}_{H,A}$ & Imperfect relative period ring corresponding to $(H_\bullet,H^+_\bullet)$ which is deformed analog of Kedlaya-Liu's imperfect ring.\\
$\breve{\Pi}^{\mathrm{bd},r}_{H,A}$ & Imperfect relative period ring corresponding to $(H_\bullet,H^+_\bullet)$ which is deformed analog of Kedlaya-Liu's imperfect ring.\\
$\breve{\Pi}^{\mathrm{bd},\dagger}_{H,A}$ & Imperfect relative period ring corresponding to $(H_\bullet,H^+_\bullet)$ which is deformed analog of Kedlaya-Liu's imperfect ring.\\
$\breve{\Pi}^{[s,r]}_{H,A}$ & Imperfect relative period ring corresponding to $(H_\bullet,H^+_\bullet)$ which is deformed analog of Kedlaya-Liu's imperfect ring.\\
$\breve{\Pi}^r_{H,A}$ & Imperfect relative period ring corresponding to $(H_\bullet,H^+_\bullet)$ which is deformed analog of Kedlaya-Liu's imperfect ring.\\
$\breve{\Pi}_{H,A}$ & Imperfect relative period ring corresponding to $(H_\bullet,H^+_\bullet)$ which is deformed analog of Kedlaya-Liu's imperfect ring.\\

$\widehat{\Pi}^{\mathrm{int},r}_{H,A}$ & Imperfect relative period ring corresponding to $(H_\bullet,H^+_\bullet)$ which is deformed analog of Kedlaya-Liu's imperfect ring.\\
$\widehat{\Pi}^{\mathrm{int},\dagger}_{H,A}$ & Imperfect relative period ring corresponding to $(H_\bullet,H^+_\bullet)$ which is deformed analog of Kedlaya-Liu's imperfect ring.\\
$\widetilde{\Omega}^\mathrm{int}_{H,A}$ & Imperfect relative period ring corresponding to $(H_\bullet,H^+_\bullet)$ which is deformed analog of Kedlaya-Liu's imperfect ring.\\
$\widetilde{\Omega}_{H,A}$ & Imperfect relative period ring corresponding to $(H_\bullet,H^+_\bullet)$ which is deformed analog of Kedlaya-Liu's imperfect ring.\\
$\widehat{\Pi}^{\mathrm{bd},r}_{H,A}$ & Imperfect relative period ring corresponding to $(H_\bullet,H^+_\bullet)$ which is deformed analog of Kedlaya-Liu's imperfect ring.\\
$\widehat{\Pi}^{\mathrm{bd},\dagger}_{H,A}$ & Imperfect relative period ring corresponding to $(H_\bullet,H^+_\bullet)$ which is deformed analog of Kedlaya-Liu's imperfect ring.\\

\end{longtable}
\end{center}

%\newpage

\newpage\section{Arithmetic Families of Period Rings}

\subsection{Basic Settings and Basic Definitions of Period Rings}

\indent In this section we first define the corresponding period rings and period sheaves involved in our study which is generalized version of the study in our previous paper, where we consider the corresponding ring of usual $p$-typical Witt vectors. For the convenience of the readers we present the corresponding construction in detail to some extent. First following \cite[Hypothesis 3.1.1]{3KL16} we consider the following setting:

\begin{setting}
In our context we consider the corresponding generalized Witt vectors in the style of \cite[3.1]{3KL16} where we first fix some complete nonarchimedean (normalized as in \cite[Hypothesis 3.1.1]{3KL16}) discrete valued field $E$ where we will use the notation $\pi$ to denote a chosen pseudo-uniformizer of $E$, and we use the notation $\mathfrak{o}_E$ to denote the corresponding ring of integers of $E$. Slightly generalizing the context \cite[Hypothesis 3.1.1]{3KL16} of \cite{3KL16} we consider the corresponding situation where $k$ is a perfect field of characteristic $p>0$, and we assume the field $k$ to contain $\mathbb{F}_{p^h}$. We now assume that the ultrametric field $E$ has residue field $k$. Again as in our previous paper we will use the notation $(R,R^+)$ to denote a perfect adic Banach uniform algebra over $k$. Following \cite{3KL16} we will use the notation $\alpha$ to denote the spectral seminorms on perfect adic Banach uniform algebras. When $E$ is of equal-characteristic we assume that $E$ contains the field $\mathbb{F}_p((\eta))$. 
\end{setting}

\begin{remark}
The field $E$ is assumed to be discrete valued. One can actually consider more general $E$. The most general situation (namely just assume $E$ is ultrametric field) is actually a little bit hard to manage due to the fact that the element in the genralized Witt vector will have no unique expression as a power series, but in some nice situations this is not that hard for instance one can consider the base change to the perfectoid field $E_\infty$.	
\end{remark}

\indent Then we could consider the following definitions of the period rings in our situation which is just the corresponding deformation over $A$ of the corresponding rings defined in \cite[Definition 4.1.1]{3KL16} in some generalized fashion:

\begin{setting}
Recall that in our situation we have the corresponding rings of generalized Witt vectors associated to some perfect ring $R$ or $R^+$. Generalizing the situation in \cite[Definition 4.1.1]{3KL16} and \cite[Setting 3.1.1]{3KL16} we denote them by $W_\pi(R)$ or $W_\pi(R^+)$ where $\pi$ is the chosen pseudo-uniformizer as defined. To be more precise here $W_\pi(R)$ (same to $W_\pi(R^+)$) is defined to be
\begin{displaymath}
W_\pi(R):=W(R)\otimes_{W(k)}\mathcal{O}_E.	
\end{displaymath}
Then we have the following ring of period rings from this notation:
\begin{displaymath}
\widetilde{\Omega}^\mathrm{int}_{R}=W_\pi(R),	\widetilde{\Omega}_{R}=W_\pi(R)[1/\pi].
\end{displaymath}
Each element of $W_\pi(R)$ could be expressed as some power series:
\begin{displaymath}
\sum_{k\geq 0}\pi^k[\overline{x}_k]	
\end{displaymath}
where one has the corresponding Gauss norm coming from the norm $\alpha$ in the style of for each $r>0$ 
\begin{displaymath}
\left\|.\right\|_{\alpha^r}(\sum_{k\geq 0}\pi^k[\overline{x}_k]):=\sup_{k\geq 0}\{p^{-k}\alpha(\overline{x}_k)^r\}.	
\end{displaymath}
Then one can define the corresponding integral Robba ring $\widetilde{\Pi}^{\mathrm{int},r}_{R}$ by taking the completion of the ring $W_\pi(R^+)[[R]]$ by using the norm defined above.
%Then recall that by using this norm one can define the corresponding subring $\widetilde{\Pi}_R^{\mathrm{int},r}$ of $\widetilde{\Omega}_R^\mathrm{int}$ which is called the integral Robba ring consisting of all the elements $\sum_{k\geq 0}\pi^k[\overline{x}_k]$ satisfying the following condition:
%\begin{displaymath}
%\lim_{k\rightarrow \infty}\left\|.\right\|_{\alpha^r}(\pi^k[\overline{x}_k])=0.	
%\end{displaymath}
From this one has the corresponding bounded Robba ring $\widetilde{\Pi}_R^{\mathrm{bd},r}$ which is defined just to be $\widetilde{\Pi}_R^{\mathrm{int},r}[1/\pi]$. Then we consider the corresponding Robba $\widetilde{\Pi}^I_R$ with respect to some interval $I\subset (0,\infty)$ which is defined to be the corresponding ind-Fr\'echet completion of $W_\pi(R^+)[[R]][1/\pi]$ (namely $W_\pi(R^+)[[r],r\in R][1/\pi]$) with respect to the family of norms $\left\|.\right\|_{\alpha^r}$, $\forall r\in I$. By taking specific intervals like $(0,r]$ or $(0,\infty)$ we have the corresponding Robba ring $\widetilde{\Pi}_R^\mathrm{r}$ and $\widetilde{\Pi}_R^\mathrm{\infty}$. Then by taking the corresponding union throughout all the radius $r>0$ we have the corresponding Robba rings $\widetilde{\Pi}^\mathrm{int}_R,\widetilde{\Pi}^\mathrm{bd}_R,\widetilde{\Pi}_R$. As in \cite[4.1]{3KL16} we have also the corresponding integral versions of some of the rings defined above.
\end{setting}

\indent Our consideration is to extend the scope of the consideration mentioned above by extending the corresponding dimension of the ring involved to some extent. To be more precise one considers the corresponding power or Laurent series over the corresponding functional analytic algebras, and then take suitable quotient.

\begin{setting}
In characteristic zero situation, we are going to use the notation\\ $\mathbb{Q}_p\{T_1,...T_d\}$ to denote a Tate algebra in rigid analytic geometry after Tate, and use in general $A$ to represent the corresponding quotients of the Tate algebras over $\mathbb{Q}_p$ as above. Throughout we assume $A$ to be reduced, carrying the spectral seminorm. While in positive characteristic situation we are going to use the notation $A$ to denote a general affinoid algebra which is assumed to be reduced.
\end{setting}

We first deform the corresponding constructions in \cite[Definition 4.1.1]{3KL16}:

\begin{definition}
We first consider the corresponding deformation of the above rings over $\mathbb{Q}_p\{T_1,...T_d\}$. We are going to use the notation $W_\pi(R)_{\mathbb{Q}_p\{T_1,...T_d\}}$ to denote the complete tensor product of $W_\pi(R)$ with the Tate algebra $\mathbb{Q}_p\{T_1,...T_d\}$ consisting of all the element taking the form as:
\begin{displaymath}
\sum_{k\geq 0,i_1\geq 0,...,i_d\geq 0}\pi^k[\overline{x}_{k,i_1,...,i_d}]T_1^{i_1}T_2^{i_2}...T_d^{i_d}	
\end{displaymath}
over which we have the Gauss norm $\left\|.\right\|_{\alpha^r,\mathbb{Q}_p\{T_1,...T_d\}}$ for any $r>0$ which is defined by:
\begin{displaymath}
\left\|.\right\|_{\alpha^r,\mathbb{Q}_p\{T_1,...T_d\}}(\sum_{k\geq 0,i_1\geq 0,...,i_d\geq 0}\pi^k[\overline{x}_{k,i_1,...,i_d}]T_1^{i_1}T_2^{i_2}...T_d^{i_d}):=\sup_{k\geq 0,i_1\geq 0,...,i_d\geq 0}p^{-k}\alpha(\overline{x}_{k,i_1,...,i_d})^r.	
\end{displaymath}
Then we define the corresponding convergent rings:
\begin{displaymath}
\widetilde{\Omega}^\mathrm{int}_{R,\mathbb{Q}_p\{T_1,...T_d\}}=W_\pi(R)_{\mathbb{Q}_p\{T_1,...T_d\}},	\widetilde{\Omega}_{R,\mathbb{Q}_p\{T_1,...T_d\}}=W_\pi(R)_{\mathbb{Q}_p\{T_1,...T_d\}}[1/\pi].	
\end{displaymath}
Then as in the previous setting we define the corresponding $\mathbb{Q}_p\{T_1,...,T_d\}$-relative integral Robba ring $\widetilde{\Pi}_{R,\mathbb{Q}_p\{T_1,...T_d\}}^{\mathrm{int},r}$ as the completion of the ring $W_{\pi}(R^+)_{\mathbb{Q}_p\{T_1,...,T_d\}}[[R]]$ by using the Gauss norm defined above.
%Then as in the previous settings we can define the corresponding integral Robba ring by considering the subring $\widetilde{\Pi}_{R,\mathbb{Q}_p\{T_1,...T_d\}}^{\mathrm{int},r}$ of $\widetilde{\Omega}_{R,\mathbb{Q}_p\{T_1,...T_d\}}^{\mathrm{int}}$	consisting of all the elements satisfying:
%\begin{displaymath}
%\lim_{k\geq 0,i_1\geq 0,...,i_d\geq 0}\left\|.\right\|_{\alpha^r,\mathbb{Q}_p\{T_1,...T_d\}}(\pi^{k}[\overline{x}_{k,i_1,...,i_d}])=0.	
%\end{displaymath}
Then one can take the corresponding union of all $\widetilde{\Pi}_{R,\mathbb{Q}_p\{T_1,...T_d\}}^{\mathrm{int},r}$ throughout all $r>0$ to define the integral Robba ring $\widetilde{\Pi}_{R,\mathbb{Q}_p\{T_1,...T_d\}}^{\mathrm{int}}$. For the bounded Robba rings we just set $\widetilde{\Pi}_{R,\mathbb{Q}_p\{T_1,...T_d\}}^{\mathrm{bd},r}=\widetilde{\Pi}_{R,\mathbb{Q}_p\{T_1,...T_d\}}^{\mathrm{int},r}[1/\pi]$ and take the union throughout all $r>0$ to define the corresponding ring $\widetilde{\Pi}_{R,\mathbb{Q}_p\{T_1,...T_d\}}^{\mathrm{bd}}$. Then we define the corresponding Robba ring $\widetilde{\Pi}_{R,\mathbb{Q}_p\{T_1,...T_d\}}^{I}$ with respect to some interval $I\subset (0,\infty)$ by taking the Fr\'echet completion of ${W_\pi(R^+)}_{\mathbb{Q}_p\{T_1,...T_d\}}[[R]][1/\pi]$ with respect to all the norms $\left\|.\right\|_{\alpha^r,\mathbb{Q}_p\{T_1,...T_d\}}$ for all $r\in I$ which means that the corresponding equivalence classes in the completion procedure will be simultaneously Cauchy with respect to all the norms $\left\|.\right\|_{\alpha^r,\mathbb{Q}_p\{T_1,...T_d\}}$ for all $r\in I$. Then we take suitable intervals such as $(0,r]$ and $(0,\infty)$ to define the corresponding Robba rings $\widetilde{\Pi}_{R,\mathbb{Q}_p\{T_1,...T_d\}}^{r}$ and $\widetilde{\Pi}_{R,\mathbb{Q}_p\{T_1,...T_d\}}^{\infty}$, respectively. Then taking the union throughout all $r>0$ one can define the corresponding Robba ring $\widetilde{\Pi}_{R,\mathbb{Q}_p\{T_1,...T_d\}}$. Again as in \cite{3KL16} one can define the corresponding integral rings of the similar types.
\end{definition}

\indent Then one can define the corresponding period rings in our context deformed over some affinoid $A$ which is isomorphic to some quotient of $\mathbb{Q}_p\{T_1,...T_d\}$, again deforming from the context of \cite[Definition 4.1.1]{3KL16}.

\begin{definition}
In the characteristic $0$, we define the following period rings:
\begin{displaymath}
\widetilde{\Omega}^\mathrm{int}_{R,A},\widetilde{\Omega}_{R,A},\widetilde{\Pi}^\mathrm{int,r}_{R,A},\widetilde{\Pi}^\mathrm{bd,r}_{R,A}, \widetilde{\Pi}^I_{R,A},\widetilde{\Pi}^r_{R,A},\widetilde{\Pi}^\infty_{R,A}	
\end{displaymath}
by taking the suitable quotient of the following period rings defined above: 
\begin{align}
\widetilde{\Omega}^\mathrm{int}_{R,\mathbb{Q}_p\{T_1,...T_d\}},&\widetilde{\Omega}_{R,\mathbb{Q}_p\{T_1,...T_d\}},\widetilde{\Pi}^\mathrm{int,r}_{R,\mathbb{Q}_p\{T_1,...T_d\}},\widetilde{\Pi}^\mathrm{bd,r}_{R,\mathbb{Q}_p\{T_1,...T_d\}},\widetilde{\Pi}^I_{R,\mathbb{Q}_p\{T_1,...T_d\}},\\
 &\widetilde{\Pi}^r_{R,\mathbb{Q}_p\{T_1,...T_d\}},\widetilde{\Pi}^\infty_{R,\mathbb{Q}_p\{T_1,...T_d\}}	
\end{align}
with respect to the structure of the affinoid algebra $A$ in the sense of Tate. The corresponding constructions do not depend on the corresponding choices of the presentations. Therefore they carry the corresponding quotient seminorms of the above Gauss norms defined in the previous definition, note that these are not something induced from the corresponding spectral seminorms from $A$ chosen at the very beginning of our study. We use the notation $\overline{\left\|.\right\|}_{\alpha^r,\mathbb{Q}_p\{T_1,...T_d\}}$ to denote the corresponding quotient Gauss norm which induces then the corresponding spectral seminorm $\left\|.\right\|_{\alpha^r,A}$ for each $r>0$. Then we define the corresponding period rings:
\begin{align}
\widetilde{\Pi}^\mathrm{int}_{R,A},\widetilde{\Pi}^\mathrm{bd}_{R,A},\widetilde{\Pi}_{R,A}	
\end{align}
by taking suitable union throughout all $r>0$.
\end{definition}

\begin{definition}
In positive characteristic situation, when we are working over general affinoid algebra $A$, we use the same notations as in the previous definition, but by using $W(R)_{\mathbb{F}_p[[\eta]]\{T_1,...,T_d\}}$ and $W(R)_{A}$ as the starting rings, namely here $A$ is isomorphic to a quotient of $\mathbb{F}_p((\eta))\{T_1,...,T_d\}$.	
\end{definition}

\indent We can then define the corresponding rings over the ring $E_\infty$ which is defined to be the completion of the ring $E(\pi^{p^{{-\infty}}})$. We first consider the corresponding undeformed versions (after \cite[Definition 4.1.1]{3KL16}):

\begin{definition}
Consider now the base change $W_{\pi,\infty}(R)$ which is defined now to be the completed tensor product $W_{\pi}(R)\widehat{\otimes}_{\mathcal{O}_E}\mathcal{O}_{E_\infty}$. Then the point is that each element in this ring admits a unique expression taking the form of $\sum_{n\in \mathbb{Z}[1/p]_{\geq 0}}\pi^n[\overline{x}_n]$ which allows us to perform the construction mentioned above. First we can define for some radius $r>0$ the corresponding period ring $\widetilde{\Pi}^{\mathrm{int},r}_{R,\infty}$ by taking the completion of the ring
\begin{displaymath}
W_{\pi,\infty}(R^+)[[R]]	
\end{displaymath}
with respect to the following Gauss type norm:
\begin{displaymath}
\left\|.\right\|_{\alpha^r}(\sum_{n\in \mathbb{Z}[1/p]_{\geq 0}}\pi^n[\overline{x}_n]):=\sup_{n\in \mathbb{Z}[1/p]_{\geq 0}}\{p^{-n}\alpha(\overline{x}_n)^r\}.	
\end{displaymath}
Then we can define the union $\widetilde{\Pi}^{\mathrm{int}}_{R,\infty}$ throughout all the radius $r>0$. Then we just define the bounded Robba ring $\widetilde{\Pi}^{\mathrm{bd},r}_{R,\infty}$ by $\widetilde{\Pi}^{\mathrm{int},r}_{R,\infty}[1/\pi]$ and also we could define the union $\widetilde{\Pi}^{\mathrm{bd}}_{R,\infty}$ throughout all the radius $r>0$. Then for any interval in $(0,\infty)$ which is denoted by $I$ we can define the corresponding Robba rings $\widetilde{\Pi}^{I}_{R,\infty}$ by taking the Fr\'echet completion of 
\begin{displaymath}
W_{\pi,\infty}(R^+)[[R]][1/\pi]	
\end{displaymath}
with respect to all the norms $\left\|.\right\|_{\alpha^t}$ for all $t\in I$. Then by taking suitable specified intervals one can define the rings $\widetilde{\Pi}^{r}_{R,\infty}$ and $\widetilde{\Pi}^{\infty}_{R,\infty}$ as before, and finally one can define the corresponding union $\widetilde{\Pi}_{R,\infty}$ throughout all the radius $r>0$. Again we have the corresponding integral version of the rings defined over $E_\infty$ as \cite{3KL16}.

\end{definition}

\indent Then we can define the following affinoid deformations (after \cite[Definition 4.1.1]{3KL16} in the flavor as above):

\begin{definition}
Consider now the base change $W_{\pi,\infty,\mathbb{Q}_p\{T_1,...,T_d\}}(R)$ which is defined now to be the completed tensor product $(W_{\pi}(R)\widehat{\otimes}_{\mathcal{O}_E}\mathcal{O}_{E_\infty})\widehat{\otimes}_{\mathbb{Q}_p}\mathbb{Q}_p\{T_1,...,T_d\}$. Then the point is that each element in this ring admits a unique expression taking the form of 
\begin{center}
$\sum_{n\in \mathbb{Z}[1/p]_{\geq 0},i_1\geq 0,...,i_d\geq 0}\pi^n[\overline{x}_{n,i_1,...,i_d}]T_1^{i_1}...T_d^{i_d}$ 	
\end{center}
which allows us to perform the construction mentioned above. First we can define for some radius $r>0$ the corresponding period ring $\widetilde{\Pi}^{\mathrm{int},r}_{R,\infty,\mathbb{Q}_p\{T_1,...,T_d\}}$ by taking the completion of the ring
\begin{displaymath}
W_{\pi,\infty,\mathbb{Q}_p\{T_1,...,T_d\}}(R^+)[[R]]	
\end{displaymath}
with respect to the following Gauss type norm:
\begin{displaymath}
\left\|.\right\|_{\alpha^r}(\sum_{n\in \mathbb{Z}[1/p]_{\geq 0},i_1\geq 0,...,i_d\geq 0}\pi^n[\overline{x}_{n,i_1,...,i_d}]T_1^{i_1}...T_d^{i_d}):=\sup_{n\in \mathbb{Z}[1/p]_{\geq 0},i_1\geq 0,...,i_d\geq 0}\{p^{-n}\alpha(\overline{x}_{n,i_1,...,i_d})^r\}.	
\end{displaymath}
Then we can define the union $\widetilde{\Pi}^{\mathrm{int}}_{R,\infty,\mathbb{Q}_p\{T_1,...,T_d\}}$ throughout all the radius $r>0$. Then we just define the bounded Robba ring $\widetilde{\Pi}^{\mathrm{bd},r}_{R,\infty,\mathbb{Q}_p\{T_1,...,T_d\}}$ by $\widetilde{\Pi}^{\mathrm{int},r}_{R,\infty,\mathbb{Q}_p\{T_1,...,T_d\}}[1/\pi]$ and also we could define the union $\widetilde{\Pi}^{\mathrm{bd}}_{R,\infty,\mathbb{Q}_p\{T_1,...,T_d\}}$ throughout all the radius $r>0$. Then for any interval in $(0,\infty)$ which is denoted by $I$ we can define the corresponding Robba rings $\widetilde{\Pi}^{I}_{R,\infty,\mathbb{Q}_p\{T_1,...,T_d\}}$ by taking the Fr\'echet completion of 
\begin{displaymath}
W_{\pi,\infty,\mathbb{Q}_p\{T_1,...,T_d\}}(R^+)[[R]][1/\pi]	
\end{displaymath}
with respect to all the norms $\left\|.\right\|_{\alpha^t}$ for all $t\in I$. Then by taking suitable specified intervals one can define the rings $\widetilde{\Pi}^{r}_{R,\infty,\mathbb{Q}_p\{T_1,...,T_d\}}$ and $\widetilde{\Pi}^{\infty}_{R,\infty,\mathbb{Q}_p\{T_1,...,T_d\}}$ as before, and finally one can define the corresponding union $\widetilde{\Pi}_{R,\infty,\mathbb{Q}_p\{T_1,...,T_d\}}$ throughout all the radius $r>0$. Again we have the corresponding integral version of the rings defined over $E_\infty$ as \cite{3KL16}. Finally over $A$ we can define in the same way as above to deform all the rings over $E_\infty$, and we do not repeat the construction again.

\end{definition}

\subsection{Basic Properties of Period Rings}	
	
\indent Then we do some reality checks over the investigation of the properties of the above period rings in the style taken in \cite[Section 5.2]{3KL15} and \cite{3KL16}. 
%We will focus on the situation where $A$ is just the ring $\mathbb{Q}_p\{T_1,...,T_d\}$ for some $d\geq 1$.

\begin{proposition} \mbox{\bf{(After Kedlaya-Liu \cite[Lemma 5.2.1]{3KL15})}}
The function $t\mapsto \left\|x\right\|_{\alpha^t,\mathbb{Q}_p\{T_1,...,T_d\}}$ for $x\in \widetilde{\Pi}^r_{R,\mathbb{Q}_p\{T_1,...,T_d\}}$ is continuous log convex for the corresponding variable $t\in (0,r]$ for any $r>0$.
\end{proposition}

\begin{proof}
Adapt the corresponding argument in the proof of 5.2.1 of \cite[Lemma 5.2.1]{3KL15} to our situation we then first look at the situation where the element is just of the form of $\pi^k[\overline{x}_{k,i_1,...,i_d}]T_1^{i_1}...T_d^{i_d}$ where the corresponding norm function in terms of $t>0$ is just affine. Then one focuses on the finite sums of these kind of elements which gives rise to to the log convex directly. Finally by taking the approximation we get the desired result.	
\end{proof}

\begin{proposition}\mbox{\bf{(After Kedlaya-Liu \cite[Lemma 5.2.2]{3KL15})}}
For any element $x\in \widetilde{\Pi}_{R,\mathbb{Q}_p\{T_1,...,T_d\}}$ we have that $x\in \widetilde{\Pi}^\mathrm{bd}_{R,\mathbb{Q}_p\{T_1,...,T_d\}}$ if and only if we have the situation where $x$ actually lives in $\widetilde{\Pi}^r_{R,\mathbb{Q}_p\{T_1,...,T_d\}}$ (for some specific $r>0$) such that $x$ itself is bounded under the norm $\left\|.\right\|_{\alpha^t,\mathbb{Q}_p\{T_1,...,T_d\}}$ for each $t\in (0,r]$.	
\end{proposition}

\begin{proof}
One direction of the proof is easy, so we only choose to present the proof of the implication in the other direction as in the original proof of 5.2.2 of \cite[Lemma 5.2.2]{3KL15} as in the following. First choose some radius $r>0$ such that the element could be assumed to be living in the ring $\widetilde{\Pi}^r_{R,\mathbb{Q}_p\{T_1,...,T_d\}}$. The idea is to transfer the original question to the question about showing the integrality of $x$ when we add some hypothesis on the norm by taking suitable powers of $p$ (since the norm is bounded for each $t\in (0,r]$ so we are reduced to the situation where the norm is bounded by $1$). Then we argue as in \cite[Lemma 5.2.2]{3KL15} to choose some approximating sequence $\{x_i\}$ living in $\widetilde{\Pi}^\mathrm{bd}_{R,\mathbb{Q}_p\{T_1,...,T_d\}}$ of $x$. Therefore we have for any $j\geq 1$ one can find then some integer $N_j\geq 1$ such that for any $i\geq N_j$ we have the estimate:
\begin{displaymath}
\left\|.\right\|_{\alpha^t,\mathbb{Q}_p\{T_1,...,T_d\}}(x_i-x)\leq p^{-j}, \forall t\in [p^{-j}r,r].	
\end{displaymath}
Then the idea is to consider the integral decomposition of the element $x_i$ which has the form of $\sum_{k=n(x_i),i_1\geq 0,i_2\geq 0,...,i_d\geq 0}\pi^k[\overline{x}_{i,k,i_1,...,i_d}]$ into the following two parts:
\begin{align}
x_i &:= y_i+z_i\\
	&:=\sum_{k=0,i_1\geq 0,i_2\geq 0,...,i_d\geq 0}\pi^k[\overline{x}_{i,k,i_1,...,i_d}]T_1^{i_1}...T_d^{i_d}+z_i
\end{align}
from which we actually have the corresponding estimate over the residual part of the decomposition above:
\begin{displaymath}
\left\|.\right\|_{\alpha^{p^{-j}r},\mathbb{Q}_p\{T_1,...,T_d\}}(\pi^k[\overline{x}_{i,k,i_1,...,i_d}]T_1^{i_1}...T_d^{i_d})\leq 1, \forall k<0	
\end{displaymath}
which implies by direct computation:
\begin{align}
\alpha^{p^{-j}r}(\overline{x}_{i,k,i_1,...,i_d})&\leq p^k\\
\alpha(\overline{x}_{i,k,i_1,...,i_d})&\leq p^{kp^{j}/r}	
\end{align}
which implies that we have the following estimate: 
\begin{align}
\left\|.\right\|_{\alpha^{r},\mathbb{Q}_p\{T_1,...,T_d\}}(x_i-y_i)&\leq p^{-k}p^{kp^{j}r/r}\\
&\leq p^{1-p^j}
\end{align}
which implies that $y_i\rightarrow x$ under the norm $\left\|.\right\|_{\alpha^{r},\mathbb{Q}_p\{T_1,...,T_d\}}$ which furthermore under the norm $\left\|.\right\|_{\alpha^{r},\mathbb{Q}_p\{T_1,...,T_d\}}$ due the property of the norm, which finishes the proof of the desired result.
\end{proof}

\begin{proposition}\mbox{\bf{(After Kedlaya-Liu \cite[Corollary 5.2.3]{3KL15})}} We have the following identity:
\begin{displaymath}
(\widetilde{\Pi}_{R,\mathbb{Q}_p\{T_1,...,T_d\}})^\times=(\widetilde{\Pi}^\mathrm{bd}_{R,\mathbb{Q}_p\{T_1,...,T_d\}})^\times.
\end{displaymath}
\end{proposition}

\begin{proof}
See 5.2.3 of \cite[Corollary 5.2.3]{3KL15}.	
\end{proof}

\begin{proposition}\mbox{\bf{(After Kedlaya-Liu \cite[Lemma 5.2.6]{3KL15})}}
For any $0< r_1\leq r_2$ we have the following equality on the corresponding period rings:
\begin{displaymath}
\widetilde{\Pi}^{\mathrm{int},r_1}_{R,\mathbb{Q}_p\{T_1,...,T_d\}}\bigcap	\widetilde{\Pi}^{[r_1,r_2]}_{R,\mathbb{Q}_p\{T_1,...,T_d\}}=\widetilde{\Pi}^{\mathrm{int},r_2}_{R,\mathbb{Q}_p\{T_1,...,T_d\}}.
\end{displaymath}	
\end{proposition}

\begin{proof}
We adapt the argument in \cite[Lemma 5.2.6]{3KL15} 5.2 to prove this in the situation where $r_1<r_2$ (otherwise this is trivial), again one direction is easy where we only present the implication in the other direction. We take any element
\begin{center}
 $x\in \widetilde{\Pi}^{\mathrm{int},r_1}_{R,\mathbb{Q}_p\{T_1,...,T_d\}}\bigcap	\widetilde{\Pi}^{[r_1,r_2]}_{R,\mathbb{Q}_p\{T_1,...,T_d\}}$
 \end{center}
  and take suitable approximating elements $\{x_i\}$ living in the bounded Robba ring such that for any $j\geq 1$ one can find some integer $N_j\geq 1$ we have whenever $i\geq N_j$ we have the following estimate:
\begin{displaymath}
\left\|.\right\|_{\alpha^{t},\mathbb{Q}_p\{T_1,...,T_d\}}(x_i-x) \leq p^{-j}, \forall t\in [r_1,r_2].	
\end{displaymath}
Then we consider the corresponding decomposition of $x_i$ for each $i=1,2,...$ into a form having integral part and the rational part $x_i=y_i+z_i$ by setting
\begin{center}
 $y_i=\sum_{k=0,i_1,...,i_d}\pi^k[\overline{x}_{i,k,i_1,...,i_d}]T_1^{i_1}...T_d^{i_d}$ 
\end{center} 
out of
\begin{center} 
$x_i=\sum_{k=n(x_i),i_1,...,i_d}\pi^k[\overline{x}_{i,k,i_1,...,i_d}]T_1^{i_1}...T_d^{i_d}$.
\end{center}
Note that by our initial hypothesis we have that the element $x$ lives in the ring $\widetilde{\Pi}^{\mathrm{int},r_1}_{R,\mathbb{Q}_p\{T_1,...,T_d\}}$ which further implies that 
\begin{displaymath}
\left\|.\right\|_{\alpha^{r_1},\mathbb{Q}_p\{T_1,...,T_d\}}(\pi^k[\overline{x}_{i,k,i_1,...,i_d}]T_1^{i_1}...T^{i_d}_d)	\leq p^{-j}.
\end{displaymath}
Therefore we have $\alpha(\overline{x}_{i,k,i_1,...,i_d})\leq p^{(k-j)/r_1}$ directly from this through computation, which implies that then:
\begin{align}
\left\|.\right\|_{\alpha^{r_2},\mathbb{Q}_p\{T_1,...,T_d\}}(\pi^k[\overline{x}_{i,k,i_1,...,i_d}]T_1^{i_1}...T^{i_d}_d)	&\leq p^{-k}p^{(k-j)r_2/r_1}\\
	&\leq p^{1+(1-j)r_1/r_1}.
\end{align}
Then one can read off the result directly from this estimate since under this estimate we can have the chance to modify the original approximating sequence $\{x_i\}$ by $\{y_i\}$ which are initially chosen to be in the integral Robba ring, which implies that actually the element $x$ lives in the right-hand side of the identity in the statement of the proposition.
\end{proof}

\begin{proposition} \mbox{\bf{(After Kedlaya-Liu \cite[Lemma 5.2.6]{3KL15})}}
For any $0< r_1\leq r_2$ we have the following equality on the corresponding period rings:
\begin{displaymath}
\widetilde{\Pi}^{\mathrm{int},r_1}_{R,A}\bigcap	\widetilde{\Pi}^{[r_1,r_2]}_{R,A}=\widetilde{\Pi}^{\mathrm{int},r_2}_{R,A}.
\end{displaymath}	
Here $A$ is some reduced affinoid algebra over $\mathbb{Q}_p$.	
\end{proposition}

\begin{proof}
See the proof of \cref{proposition5.7}.	
\end{proof}

%%%%%%%%%%%%%%%%%%%%%%%%%%%%%%%%%%%%%%%%%%%%%%%%%%%%%%

\indent Then we have the following analog of the corresponding result of \cite[Lemma 5.2.8]{3KL15}:

\begin{proposition} \mbox{\bf{(After Kedlaya-Liu \cite[Lemma 5.2.8]{3KL15})}}
Consider now in our situation the radii $0< r_1\leq r_2$, and consider any element $x\in \widetilde{\Pi}^{[r_1,r_2]}_{R,\mathbb{Q}_p\{T_1,...,T_d\}}$. Then we have that for each $n\geq 1$ one can decompose $x$ into the form of $x=y+z$ such that $y\in \pi^n\widetilde{\Pi}^{\mathrm{int},r_2}_{R,\mathbb{Q}_p\{T_1,...,T_d\}}$ with $z\in \bigcap_{r\geq r_2}\widetilde{\Pi}^{[r_1,r]}_{R,\mathbb{Q}_p\{T_1,...,T_d\}}$ with the following estimate for each $r\geq r_2$:
\begin{displaymath}
\left\|.\right\|_{\alpha^r,\mathbb{Q}_p\{T_1,...,T_d\}}(z)\leq p^{(1-n)(1-r/r_2)}\left\|.\right\|_{\alpha^{r_2},\mathbb{Q}_p\{T_1,...,T_d\}}(z)^{r/r_2}.	
\end{displaymath}

\end{proposition}

\begin{proof}
As in \cite[Lemma 5.2.8]{3KL15} and in the proof of our previous proposition we first consider those elements $x$ living in the bounded Robba rings which could be expressed in general as
\begin{center}
 $\sum_{k=n(x),i_1,...,i_d}\pi^k[\overline{x}_{k,i_1,...,i_d}]T_1^{i_1}...T_d^{i_d}$.
 \end{center}	
In this situation the corresponding decomposition is very easy to come up with, namely we consider the corresponding $y_i$ as the corresponding series:
\begin{displaymath}
\sum_{k\geq n,i_1,...,i_d}\pi^k[\overline{x}_{k,i_1,...,i_d}]T_1^{i_1}...T_d^{i_d}	
\end{displaymath}
which give us the desired result since we have in this situation when focusing on each single term:
\begin{align}
\left\|.\right\|_{\alpha^r,\mathbb{Q}_p\{T_1,...,T_d\}}&(\pi^k[\overline{x}_{k,i_1,...,i_d}]T_1^{i_1}...T_d^{i_d})=p^{-k}\alpha(\overline{x}_{k,i_1,...,i_d})^r\\
&=p^{-k(1-r/r_2)}\left\|.\right\|_{\alpha^{r_2},\mathbb{Q}_p\{T_1,...,T_d\}}(\pi^k[\overline{x}_{k,i_1,...,i_d}]T_1^{i_1}...T_d^{i_d})^{r/r_2}\\
&\leq p^{(1-n)(1-r/r_2)}\left\|.\right\|_{\alpha^{r_2},\mathbb{Q}_p\{T_1,...,T_d\}}(\pi^k[\overline{x}_{k,i_1,...,i_d}]T_1^{i_1}...T_d^{i_d})^{r/r_2}
\end{align}
for all those suitable $k$. Then to tackle the more general situation we consider the approximating sequence consisting of all the elements in the bounded Robba ring as in the usual situation considered in \cite[Lemma 5.2.8]{3KL15}, namely we inductively construct the following approximating sequence just as:
\begin{align}
\left\|.\right\|_{\alpha^r,\mathbb{Q}_p\{T_1,...,T_d\}}(x-x_0-...-x_i)\leq p^{-i-1}	\left\|.\right\|_{\alpha^r,\mathbb{Q}_p\{T_1,...,T_d\}}(x), i=0,1,..., r\in [r_1,r_2].
\end{align}
Here all the elements $x_i$ for each $i=0,1,...$ are living in the bounded Robba ring, which immediately gives rise to the suitable decomposition as proved in the previous case namely we have for each $i$ the decomposition $x_i=y_i+z_i$ with the desired conditions as mentioned in the statement of the proposition. We first take the series summing all the elements $y_i$ up for all $i=0,1,...$, which first of all converges under the norm $\left\|.\right\|_{\alpha^r,\mathbb{Q}_p\{T_1,...,T_d\}}$ for all the radius $r\in [r_1,r_2]$, and also note that all the elements $y_i$ within the infinite sum live inside the corresponding integral Robba ring $\widetilde{\Pi}^{\mathrm{int},r_2}_{R,\mathbb{Q}_p\{T_1,...,T_d\}}$, which further implies the corresponding convergence ends up in $\widetilde{\Pi}^{\mathrm{int},r_2}_{R,\mathbb{Q}_p\{T_1,...,T_d\}}$. For the elements $z_i$ where $i=0,1,...$ also sum up to a converging series in the desired ring since combining all the estimates above we have:
\begin{displaymath}
\left\|.\right\|_{\alpha^r,\mathbb{Q}_p\{T_1,...,T_d\}}(z_i)\leq p^{(1-n)(1-r/r_2)}\left\|.\right\|_{\alpha^{r_2},\mathbb{Q}_p\{T_1,...,T_d\}}(x)^{r/r_2}.	
\end{displaymath}
\end{proof}

%\begin{corollary}
%Consider now in our situation the radii $0< r_1\leq r_2$, and consider any element $x\in \widetilde{\Pi}^{[r_1,r_2]}_{R,A}$. Then we have that for each $n\geq 1$ one can decompose $x$ into the form of $x=y+z$ such that $y\in \pi^n\widetilde{\Pi}^{\mathrm{int},r_2}_{R,A}$ with $z\in \bigcap_{r\geq r_2}\widetilde{\Pi}^{[r_1,r]}_{R,A}$ with the following estimate for each $r\geq r_2$:
%\begin{displaymath}
%\left\|.\right\|_{\alpha^r,A}(z)\leq p^{(1-n)(1-r/r_2)}\left\|.\right\|_{\alpha^{r_2},A}(z)^{r/r_2}.	
%\end{displaymath}	
%\end{corollary}
%
%\begin{proof}
%This could be derived directly from the previous proposition, by considering the corresponding residual norms and then the spectral norms.	
%\end{proof}

\begin{proposition} \mbox{\bf{(After Kedlaya-Liu \cite[Lemma 5.2.10]{3KL15})}}
We have the following identity:
\begin{displaymath}
\widetilde{\Pi}^{[s_1,r_1]}_{R,\mathbb{Q}_p\{T_1,...,T_d\}}\bigcap\widetilde{\Pi}^{[s_2,r_2]}_{R,\mathbb{Q}_p\{T_1,...,T_d\}}=\widetilde{\Pi}^{[s_1,r_2]}_{R,\mathbb{Q}_p\{T_1,...,T_d\}},
\end{displaymath}
here the radii satisfy $<s_1\leq s_2 \leq r_1 \leq r_2$.
\end{proposition}

\begin{proof}
In our situation one direction is obvious while on the other hand we consider any element $x$ in the intersection on the left, then by the previous proposition we	have the decomposition $x=y+z$ where $y\in \widetilde{\Pi}^{\mathrm{int},r_1}_{R,\mathbb{Q}_p\{T_1,...,T_d\}}$ and $z\in \widetilde{\Pi}^{[s_1,r_2]}_{R,\mathbb{Q}_p\{T_1,...,T_d\}}$. Then as in \cite[Lemma 5.2.10]{3KL15} section 5.2 we look at $y=x-z$ which lives in the intersection:
\begin{displaymath}
\widetilde{\Pi}^{\mathrm{int},r_1}_{R,\mathbb{Q}_p\{T_1,...,T_d\}}\bigcap	\widetilde{\Pi}^{[s_2,r_2]}_{R,\mathbb{Q}_p\{T_1,...,T_d\}}=\widetilde{\Pi}^{\mathrm{int},r_2}_{R,\mathbb{Q}_p\{T_1,...,T_d\}}
\end{displaymath}
which finishes the proof.
\end{proof}

\begin{proposition} \mbox{\bf{(After Kedlaya-Liu \cite[Lemma 5.2.10]{3KL15})}}
We have the following identity:
\begin{displaymath}
\widetilde{\Pi}^{[s_1,r_1]}_{R,A}\bigcap\widetilde{\Pi}^{[s_2,r_2]}_{R,A}=\widetilde{\Pi}^{[s_1,r_2]}_{R,A},
\end{displaymath}
here the radii satisfy $<s_1\leq s_2 \leq r_1 \leq r_2$.
	
\end{proposition}

\begin{proof}
See the proof of \cref{proposition5.7}.	
\end{proof}

\begin{remark}
This is subsection is finished so far only for the situation where $E$ is of mixed characteristic. But everything uniformly carries over for our original assumption on the field $E$ and $A$. We will not repeat the proof again.	
\end{remark}

%\newpage

\newpage\section{Period Modules and Period Sheaves}

\indent We now consider the corresponding Frobenius modules over the corresponding period rings defined in the previous section. Also one can consider the corresponding period sheaves in the style of \cite{3KL16}. We would like to point out that actually the corresponding sheaves in our context could mean the following two different objects. First the corresponding period rings defined in the previous section themselves are sheafy, which means that one can study the corresponding analytic geometry over the relative affinoid spaces over for instance the ring $\widetilde{\Pi}_{R}^{\mathrm{int},r}$ or the ring $\widetilde{\Pi}_{R}^{\mathrm{int},r}$ which has its own interests. On the other hand we have the situation where one can glue along the algebra $R$ over corresponding \'etale, corresponding pro-\'etale and corresponding v-sites but leaving the algebra $A$ unglued. We point out that both could have the chance to be compared in a coherent way, to the corresponding representation theoretic consideration in the pseudocoherent setting.
%%%%%%%%%%%%%%%%%%%%%%%%%%%%%%%%%%%%%%%%%%%%%%%%%%%%

\begin{setting}
We will work in the categories of the pseudocoherent, fpd and finite projective modules over the period rings defined above. First we specify the Frobenius in our setting. The rings involved are:
\begin{align}
\widetilde{\Omega}^\mathrm{int}_{R,A},\widetilde{\Omega}_{R,A}, \widetilde{\Pi}^\mathrm{int}_{R,A}, \widetilde{\Pi}^{\mathrm{int},r}_{R,A},\widetilde{\Pi}^\mathrm{bd}_{R,A},\widetilde{\Pi}^{\mathrm{bd},r}_{R,A}, \widetilde{\Pi}_{R,A}, \widetilde{\Pi}^r_{R,A},\widetilde{\Pi}^+_{R,A}, \widetilde{\Pi}^\infty_{R,A},\widetilde{\Pi}^I_{R,A}.
\end{align}
We are going to endow these rings with the Frobenius induced by continuation from the Witt vector part only, which is to say the corresponding Frobenius induced by the $p^h$-power absolute Frobenius over $R$. Note all the rings above are defined by taking the product of $\triangle$ where each $\triangle$ representing one of the following rings (over $E$):
\begin{align}
\widetilde{\Omega}^\mathrm{int}_{R},\widetilde{\Omega}_{R}, \widetilde{\Pi}^\mathrm{int}_{R}, \widetilde{\Pi}^{\mathrm{int},r}_{R},\widetilde{\Pi}^\mathrm{bd}_{R},\widetilde{\Pi}^{\mathrm{bd},r}_{R}, \widetilde{\Pi}_{R}, \widetilde{\Pi}^r_{R},\widetilde{\Pi}^+_{R}, \widetilde{\Pi}^\infty_{R},\widetilde{\Pi}^I_{R}
\end{align}
with the affinoid ring $A$. The Frobenius acts on $A$ trivially and we assume that the action is $A$-linear.  
\end{setting}

\indent First we consider the following sheafification as in \cite[Definition 4.4.2]{3KL16}:

\begin{setting} 
Consider the space $X=\mathrm{Spa}(R,R^+)$, over this perfectoid space there were sheaves:
\begin{align}
\widetilde{\Omega}^\mathrm{int}_{},\widetilde{\Omega}_{}, \widetilde{\Pi}^\mathrm{int}_{}, \widetilde{\Pi}^{\mathrm{int},r}_{},\widetilde{\Pi}^\mathrm{bd}_{},\widetilde{\Pi}^{\mathrm{bd},r}_{}, \widetilde{\Pi}_{}, \widetilde{\Pi}^r_{},\widetilde{\Pi}^+_{}, \widetilde{\Pi}^\infty_{},\widetilde{\Pi}^I_{}.
\end{align}
defined over this space through the corresponding adic, \'etale, pro-\'etale or $v$-topology, we consider the corresponding sheaves defined over the same Grothendieck sites but with further deformed consideration:
\begin{align}
\widetilde{\Omega}^\mathrm{int}_{*,A},\widetilde{\Omega}_{*,A}, \widetilde{\Pi}^\mathrm{int}_{*,A}, \widetilde{\Pi}^{\mathrm{int},r}_{*,A},\widetilde{\Pi}^\mathrm{bd}_{*,A},\widetilde{\Pi}^{\mathrm{bd},r}_{*,A}, \widetilde{\Pi}_{*,A}, \widetilde{\Pi}^r_{*,A},\widetilde{\Pi}^+_{*,A}, \widetilde{\Pi}^\infty_{*,A},\widetilde{\Pi}^I_{*,A}.
\end{align}
\end{setting}

\begin{remark}
The consideration in the previous setting in some sense reflects some kind of rigidity. Since the corresponding construction is just that the ring $A$ acts through the completed tensor product directly on the sheaves after \cite{3KL16}, in most cases when one would like to compare those sheaves of modules over the sheaves of rings above and the corresponding modules over the period rings defined above, one needs to work through the corresponding vanishing results and so on in this deformed context. We work around this by considering the corresponding globalized version of this. 	
\end{remark}

\begin{setting}
In our context we can make the following parallel discussion to that in \cite{3KL16}. First we can use the corresponding notation $E_\infty$ to denote the corresponding completion of $E(\pi^{1/{p^\infty}})$ and we denote the corresponding integral ring by $\mathcal{O}_{E_\infty}$. Then we have could have that the corresponding splitting of the ring $\mathcal{O}_E$ in this bigger perfectoid ring. Then we can consider the corresponding extended period rings $\widetilde{\Omega}^\mathrm{int}_{R}\widehat{\otimes}_{\mathcal{O}_E}\mathcal{O}_{E_\infty}$ and $\widetilde{\Pi}^{\mathrm{int},r}_{R}\widehat{\otimes}_{\mathcal{O}_E}\mathcal{O}_{E_\infty}$	in the same fashion as in \cite[Definition 4.1.11]{3KL16}.
\end{setting}

\begin{proposition}\mbox{}\\ \mbox{\bf{(After Kedlaya-Liu \cite[Lemma 4.1.12, Corollary 4.1.14]{3KL16})}}\\
 The ring $\widetilde{\Pi}^{\mathrm{int},r}_{R}\widehat{\otimes}_{\mathcal{O}_E}\mathcal{O}_{E_\infty}$ is perfectoid, so stably uniform and sheafy. By using the same notation in \cite[Lemma 4.1.12, Corollary 4.1.14]{3KL16} we have that the corresponding tilting perfectoid ring is $R\{(\overline{\pi}/p^{-1/r})^{1/p^\infty}\}$ which is the completion of $R[(\overline{\pi}/p^{-1/r})^{1/p^\infty}]$.	
\end{proposition}

\begin{proof}
This is the same as the proof for \cite[Lemma 4.1.12, Corollary 4.1.14]{3KL16} after we choose a suitable topologically nilpotent element.
\end{proof}

\begin{proposition}\mbox{}\\
\mbox{\bf{(After Kedlaya-Liu \cite[Proposition 4.1.13, Corollary 4.1.14]{3KL16})}}\\
The ring $\widetilde{\Pi}_R^I\widehat{\otimes}_{\mathcal{O}_E}\mathcal{O}_{E_\infty}$ for some closed interval $I=[s,r]$ is perfectoid so stably uniform and sheafy. Then in our situation this perfectoid admits tilting ring
\begin{center}
 $R\{(\overline{\pi}/p^{-1/r})^{1/p^\infty},(p^{-1/s}/\overline{\pi})^{1/p^\infty}\}$
\end{center} 
which is the completion of $R[(\overline{\pi}/p^{-1/r})^{1/p^\infty},(p^{-1/s}/\overline{\pi})^{1/p^\infty}]$.	
\end{proposition}

\begin{proof}
This is the same as the proof of \cite[Proposition 4.1.13, Corollary 4.1.14]{3KL16} by considering the previous proposition.	
\end{proof}

\begin{setting}
We are going to use the notation $\mathfrak{X}$ to denote a rigid analytic space in rigid geometry. Then we can consider the corresponding relative period rings or sheaves over $\mathcal{O}_{\mathfrak{X}}$. The period rings and the corresponding period sheaves over the sheaf $\mathcal{O}_{\mathfrak{X}}$ are defined to be the following sheaves (one takes the complete tensor product of the undeformed period rings with the corresponding exact sequence for the sheaf $?$, which again gives the corresponding exact sequence due to the fact that the undeformed period rings admit Schauder bases as in \cite[Definition 6.1]{3KP}):
%\begin{align}
%\widetilde{\Omega}^\mathrm{int}_{R,?},\widetilde{\Omega}_{R,?}, \widetilde{\Pi}^\mathrm{int}_{R,?}, \widetilde{\Pi}^{\mathrm{int},r}_{R,?},\widetilde{\Pi}^\mathrm{bd}_{R,?},\widetilde{\Pi}^{\mathrm{bd},r}_{R,?}, \widetilde{\Pi}_{R,?}, \widetilde{\Pi}^r_{R,?},\widetilde{\Pi}^+_{R,?}, \widetilde{\Pi}^\infty_{R,?},\widetilde{\Pi}^I_{R,?},
%\end{align}
\begin{align}
\widetilde{\Pi}_{R,?}, \widetilde{\Pi}^r_{R,?},\widetilde{\Pi}^\infty_{R,?},\widetilde{\Pi}^I_{R,?},
\end{align}
with
%\begin{align}
%\widetilde{\Omega}^\mathrm{int}_{*,?},\widetilde{\Omega}_{*,?}, \widetilde{\Pi}^\mathrm{int}_{*,?}, \widetilde{\Pi}^{\mathrm{int},r}_{*,?},\widetilde{\Pi}^\mathrm{bd}_{*,?},\widetilde{\Pi}^{\mathrm{bd},r}_{*,?}, \widetilde{\Pi}_{*,?}, \widetilde{\Pi}^r_{*,?},\widetilde{\Pi}^+_{*,?}, \widetilde{\Pi}^\infty_{*,?},\widetilde{\Pi}^I_{*,?},
%\end{align}
\begin{align}
\widetilde{\Pi}_{*,?}, \widetilde{\Pi}^r_{*,?},\widetilde{\Pi}^\infty_{*,?},\widetilde{\Pi}^I_{*,?},
\end{align}
where $*$ is some \'etale site, pro\'et site or $v$-site and etc, and $?$ represents $\mathcal{O}_{\mathfrak{X}}$ which implies that we treat these sheaves of rings as sheaves over $*$ in some relative sense. One should understand this as the corresponding sheaves essentially with respect $*$.
\end{setting}

%%%%%%%%%%%%%%%%%%%%%%%%%%do we really need the sheafification here?

\begin{remark}
It is hard to consider the corresponding comparison between the representation theories over these sheaves (deformed) with the rings (without any sheafified consideration), since the corresponding global section functor is a little bit hard to study. But on the other hand comparing this among the sheaves themselves and with the Fargues-Fontaine curves are very natural and interesting.
\end{remark}

\begin{assumption}
From now on, we are going to assume that the affinoid ring $A$ is splitting in some perfectoid covering $A^\mathrm{perf}$, but only when the corresponding sheafiness of the deformed period rings is essentially relevant. For instance if $A$ is just some analytic field as those in \cite{3KL16} then this is satisfied. This will guarantee the corresponding period rings over $A$ (we now assume them to be sousperfectoid) will then be sheafy by Kedlaya-Hansen's sheafiness criterion \cite{3KH}.	
\end{assumption}

\indent We first consider the following result which is slight generalization of the corresponding result from \cite[4.3.1]{3KL16}. Note that our $E$ is general than the basic setting of \cite[Setting 3.1.1]{3KL16}.

\begin{remark}
In what follows, there will be some situation we will consider the perfectoid deformed version of the period rings and sheaves. All the definitions are parallel to the base level definitions (including the corresponding Frobenius modules and sheaves).	
\end{remark}

\begin{lemma} \mbox{\bf{(After Kedlaya-Liu \cite[4.3.1]{3KL16})}} \label{lemma3.11}
Let $M$ be any \'etale-stably pseudocoherent or fpd module defined over $\widetilde{\Pi}^{\mathrm{int},r}_{R}$ and $\widetilde{\Pi}^{\mathrm{int},r}_{R,\infty}$ for some $r>0$ or $\widetilde{\Pi}^I_{R}$ and $\widetilde{\Pi}^I_{R,\infty}$ for some $I$ which is assumed to be closed, and let $\widetilde{M}$ be the corresponding sheaf attached to the module $M$ over the sheaf $\widetilde{\Pi}^{\mathrm{int},r}_{*}$ for some $r>0$ or $\widetilde{\Pi}^I_{*}$ for some $I$ which is assumed to be closed, over the \'etale and the pro-\'etale sites of $\mathrm{Spa}(R,R^+)$. Then we have the following statement:
\begin{align}
H^i(X,\widetilde{M})=H^i(X_\text{\'et},\widetilde{M})=H^i(X_\text{pro\'et},\widetilde{M})	
\end{align}
vanish beyond the $0$-th degree and give the module $M$ at the degree $0$. Also furthermore we have that in the projective situation we have that the same holds under the $v$-topology. 
\end{lemma}

\begin{proof}
The corresponding results could be transformed to perfectoid spaces as in \cite[4.3.1]{3KL16} by considering the splitting after passing to perfectoid field $E_\infty$ and perfectoid ring $\mathcal{O}_{\mathfrak{X}^\mathrm{perf}}$, then by applying the results in \cite[2.5.1,2.5.11,3.4.6,3.5.6]{3KL16} we can get the desired results.	
\end{proof}

\begin{definition}
Now working over some perfect uniform and adic space over $k$, as in \cite[Definition 4.3.2]{3KL16} we define the corresponding pseudocoherent sheaves or fpd sheaves over $\widetilde{\Pi}^{\mathrm{int},r}_*$ or $\widetilde{\Pi}^{I}_*$ (where $I$ is closed) where $*$ represents one of $X$, $X_\text{\'et}$ and $X_{\text{pro\'et}}$ to be the sheaves associated to \'etale-stably pseudocoherent or fpd modules in some local manner. 	
\end{definition}

\begin{proposition} \mbox{\bf{(After Kedlaya-Liu \cite[Theorem 4.3.3]{3KL16})}} \label{proposition3.13}
The global section functor establishes an equivalence between the  categories of the pseudocoherent sheaves over $\widetilde{\Pi}^{\mathrm{int},r}_*$ or $\widetilde{\Pi}^{I}_*$ (and $\widetilde{\Pi}^{\mathrm{int},r}_{*,\infty}$ or $\widetilde{\Pi}^{I}_{*,\infty}$) over the pro-\'etale site of $X=\mathrm{Spa}(R,R^+)$ and the \'etale-stably pseudocoherent modules over $\widetilde{\Pi}^{\mathrm{int},r}_R$ or $\widetilde{\Pi}^{I}_R$ (and $\widetilde{\Pi}^{\mathrm{int},r}_{R,\infty}$ or $\widetilde{\Pi}^{I}_{R,\infty}$). The global section functor establishes an equivalence between the  categories of the fpd sheaves over $\widetilde{\Pi}^{\mathrm{int},r}_*$ or $\widetilde{\Pi}^{I}_*$ over the pro-\'etale site of $X=\mathrm{Spa}(R,R^+)$ and the \'etale-stably fpd modules over $\widetilde{\Pi}^{\mathrm{int},r}_R$ or $\widetilde{\Pi}^{I}_R$. The following morphisms are effective descent morphisms for pseudocoherent Banach modules: I. $\widetilde{\Pi}_{R,A}^r\rightarrow \widetilde{\Pi}_{R,A}^s\oplus \widetilde{\Pi}_{R,A}^{[s,r]}$; II. $\widetilde{\Pi}_{R,A}^I\rightarrow \oplus_{i=1}^k \widetilde{\Pi}_{R,A}^{I_i}$. here the corresponding set of interval $\{I_i\}_{i=1,...,k}$ consists of finitely many closed intervals covering the interval $I$.
 
\end{proposition}

\begin{proof}
As in \cite[Theorem 4.3.3]{3KL16}, by taking the splitting base change to $\mathcal{O}_{E_\infty}$	and considering the corresponding properties of being perfectoid we can transform the proofs around the global section functor to the corresponding statement for juts perfectoid spaces, which could be finished by \cite[2.5.5,2.5.14,3.4.9]{3KL16}. The rest statements are further consequences of the fact that the period rings $\widetilde{\Pi}_{R,A}^r$ and $\widetilde{\Pi}_{R,A}^I$ are sheafy, see \cite[2.5.5,2.5.14]{3KL16}.
\end{proof}

\indent Then we generalize the corresponding definitions of Frobenius action and Frobenius modules in \cite[Definition 4.4.2-4.4.4]{3KL16} to our situation.

\begin{definition}
In our situation we consider the corresponding Frobenius action on the following period rings and sheaves:
\begin{align}
\widetilde{\Omega}^\mathrm{int}_{R,A},\widetilde{\Omega}^\mathrm{int}_{R,A},\widetilde{\Pi}^\mathrm{int}_{R,A},\widetilde{\Pi}^\mathrm{bd}_{R,A},\widetilde{\Pi}_{R,A},\widetilde{\Pi}^\infty_{R,A}, \widetilde{\Omega}^\mathrm{int}_{*,A}, \widetilde{\Pi}^\mathrm{int}_{*,A}, \widetilde{\Pi}^+_{*,A}, \widetilde{\Omega}_{*,A}, \widetilde{\Pi}^\mathrm{bd}_{*,A}, \widetilde{\Pi}^\infty_{*,A}, \widetilde{\Pi}^+_{*,A}, \widetilde{\Pi}_{*,A},   
\end{align}
which is defined by considering the corresponding lift of the absolute Frobenius of $p^h$-power in characteristic $p>0$ induced from $R$, which will be denoted by $\varphi$. We then introduce more general consideration by taking $E_a$ to be some unramified extension of $E$ of degree $a$ divisible by $h$, the corresponding Frobenius will be denoted by $\varphi^a$.
\end{definition}

%%%%%%%%%%%%%%%%%%%%%%%%%%%%%%%%%%%%%%%%%%%%%%%%T
	
\indent Then we generalize the corresponding Frobenius modules in \cite[Definition 4.4.4]{3KL16} to our situation as in the following.

\begin{definition}
Over the period rings and sheaves (each is denoted by $\triangle$ in this definition) defined in the previous definition we define as in \cite[Definition 4.4.4]{3KL16} the corresponding $\varphi^a$-modules over $\triangle$ which are respectively projective, pseudocoherent or fpd to be the corresponding finite projetive, pseudocoherent or fpd modules over $\triangle$ with further assigned semilinear action of the operator $\varphi^a$. Here we define in our situation the corresponding $\varphi^a$-cohomology to be the (hyper)-cohomology of the following complex:
\[
\xymatrix@R+0pc@C+0pc{
0\ar[r]\ar[r]\ar[r] &M
\ar[r]^{\varphi-1}\ar[r]\ar[r] &M
\ar[r]\ar[r]\ar[r] &0.
}
\]
We also require that the modules are complete for the natural topology involved in our situation and for any module over $\widetilde{\Pi}_{*,A}$ to be some base change of some module over $\widetilde{\Pi}^r_{*,A}$ (which will be defined in the following) over each perfectoid subdomain $Y$ (in this situation we are considering the pro-\'etale site). 
\end{definition}

\indent Now we define the corresponding modules over the rings which are the domains in the following morphisms induced from the Frobenius map $\varphi^a$:
\begin{align} 
\widetilde{\Pi}^{\mathrm{int},r}_{R,A}\rightarrow \widetilde{\Pi}^{\mathrm{int},rp^{-ha}}_{R,A},\widetilde{\Pi}^{\mathrm{bd},r}_{R,A}\rightarrow \widetilde{\Pi}^{\mathrm{bd},rp^{-ha}}_{R,A},\widetilde{\Pi}^{r}_{R,A}\rightarrow \widetilde{\Pi}^{rp^{-ha}}_{R,A}\\
\widetilde{\Pi}^{\mathrm{int},r}_{*,A}\rightarrow \widetilde{\Pi}^{\mathrm{int},rp^{-ha}}_{*,A},\widetilde{\Pi}^{\mathrm{bd},r}_{*,A}\rightarrow \widetilde{\Pi}^{\mathrm{bd},rp^{-ha}}_{*,A},\widetilde{\Pi}^{r}_{*,A}\rightarrow \widetilde{\Pi}^{rp^{-ha}}_{*,A}.	
\end{align}

\begin{definition}
Over each rings $\triangle$ which are the domains in the morphisms as mentioned just before this definition, we define the corresponding projective, pseudocoherent or fpd $\varphi^a$-module over any $\triangle$ listed above to be the corresponding finite projective, pseudocoherent or fpd module $M$ over $\triangle$ with additionally endowed semilinear Frobenius action from $\varphi^a$ such that we have the isomorphism $\varphi^{a*}M\overset{\sim}{\rightarrow}M\otimes \square$ where the ring $\square$ is one of the targets listed above. Also as in \cite[Definition 4.4.4]{3KL16} we assume that the module over $\widetilde{\Pi}^I_{*,A}$ is then complete for the natural topology and the corresponding base change to $\widetilde{\Pi}^I_{*,A}$ for any interval which is assumed to be closed $I\subset [0,r)$ gives rise to a module over $\widetilde{\Pi}^I_{*,A}$ with specified conditions which will be specified below. Also the cohomology of any module under this definition will be defined to be the (hyper)cohomology of the complex in the following form:
\[
\xymatrix@R+0pc@C+0pc{
0\ar[r]\ar[r]\ar[r] &M
\ar[r]^{\varphi-1}\ar[r]\ar[r] &M\otimes_\triangle \square
\ar[r]\ar[r]\ar[r] &0.
}
\]
\end{definition}

\indent Then we consider the following morphisms of specific period rings induced by the Frobenius. 
\begin{align}
\widetilde{\Pi}_{R,A}^{[s,r]}	\rightarrow \widetilde{\Pi}_{R,A}^{[sp^{-ah},rp^{-ah}]}\\
\widetilde{\Pi}_{*,A}^{[s,r]}	\rightarrow \widetilde{\Pi}_{*,A}^{[sp^{-ah},rp^{-ah}]}
\end{align}

with the corresponding morphisms in the following:

\begin{align}
\widetilde{\Pi}_{R,A}^{[s,r]}	\rightarrow \widetilde{\Pi}_{R,A}^{[s,rp^{-ah}]}\\
\widetilde{\Pi}_{*,A}^{[s,r]}	\rightarrow \widetilde{\Pi}_{*,A}^{[s,rp^{-ah}]}
\end{align}

\begin{definition}
Again as in \cite[Definition 4.4.4]{3KL16}, we define the corresponding projective, pseudocoherent and fpd $\varphi^a$-modules over the domain rings or sheaves of rings in the morphisms just before this definition to be the finite projective, pseudocoherent and fpd modules (which will be denoted by $M$) over the domain rings in the morphism just before this definition additionally endowed with semilinear Frobenius action from $\varphi^a$ with the following isomorphisms:
\begin{align}
\varphi^{a*}M\otimes_{\widetilde{\Pi}_{R,A}^{[sp^{-ah},rp^{-ah}]}}\widetilde{\Pi}_{R,A}^{[s,rp^{-ah}]}\overset{\sim}{\rightarrow}M\otimes_{\widetilde{\Pi}_{R,A}^{[s,r]}}\widetilde{\Pi}_{R,A}^{[s,rp^{-ah}]},\\
\varphi^{a*}M\otimes_{\widetilde{\Pi}_{*,A}^{[sp^{-ah},rp^{-ah}]}}\widetilde{\Pi}_{*,A}^{[s,rp^{-ah}]}\overset{\sim}{\rightarrow}M\otimes_{\widetilde{\Pi}_{*,A}^{[s,r]}}\widetilde{\Pi}_{*,A}^{[s,rp^{-ah}]}.
\end{align}
We now assume that the modules are complete with respect to the natural topology. And we assume that for any perfectoid subdomain $Y$ in the corresponding topology (defining the sheaves in this situation) the corresponding global sections over $Y$ give rise to \'etale-stably pseudocoherent modules.
\end{definition}

\indent Also one can further define the corresponding bundles carrying semilinear Frobenius in our context as in the situation of \cite[Definition 4.4.4]{3KL16}:

\begin{definition}
Over the ring $\widetilde{\Pi}_{R,A}$ we define a corresponding projective, pseudocoherent and fpd Frobenius bundle to be a family $(M_I)_I$ of finite projective, \'etale stably pseudocoherent and \'etale stably fpd modules over each $\widetilde{\Pi}^I_{R,A}$ carrying the natural Frobenius action coming from the operator $\varphi^a$ such that for any two involved intervals having the relation $I\subset J$ we have:
\begin{displaymath}
M_J\otimes_{\widetilde{\Pi}^J_{R,A}}\widetilde{\Pi}^I_{R,A}\overset{\sim}{\rightarrow}	M_I
\end{displaymath}
with the obvious cocycle condition. Here we have to propose condition on the intervals that for each $I=[s,r]$ involved we have $s\leq rp^{ah}$.
\end{definition}

\indent We then have the following analog of \cite[Lemma 4.4.8]{3KL16}:

\begin{proposition} \mbox{\bf{(After Kedlaya \cite[Lemma 4.4.8]{3KL15})}} \label{prop3.18}
Consider the correspo-
nding finite generated Frobenius modules over the following period rings or sheaves defined above:
\begin{displaymath}
\widetilde{\Omega}^\mathrm{int}_{R,A},\widetilde{\Omega}^\mathrm{int}_{R,A},\widetilde{\Pi}^\mathrm{int}_{R,A},\widetilde{\Pi}^\mathrm{bd}_{R,A},\widetilde{\Pi}_{R,A},\widetilde{\Pi}^\infty_{R,A}, \widetilde{\Omega}^\mathrm{int}_{*,A}, \widetilde{\Pi}^\mathrm{int}_{*,A}, \widetilde{\Pi}^+_{*,A}, \widetilde{\Omega}_{*,A}, \widetilde{\Pi}^\mathrm{bd}_{*,A}, \widetilde{\Pi}^\infty_{*,A}, \widetilde{\Pi}^+_{*,A}, \widetilde{\Pi}_{*,A},   
\end{displaymath}
and
\begin{align} 
\widetilde{\Pi}^{\mathrm{int},r}_{R,A}\rightarrow \widetilde{\Pi}^{\mathrm{int},rp^{-ha}}_{R,A},\widetilde{\Pi}^{\mathrm{bd},r}_{R,A}\rightarrow \widetilde{\Pi}^{\mathrm{bd},rp^{-ha}}_{R,A},\widetilde{\Pi}^{r}_{R,A}\rightarrow \widetilde{\Pi}^{rp^{-ha}}_{R,A}\\
\widetilde{\Pi}^{\mathrm{int},r}_{*,A}\rightarrow \widetilde{\Pi}^{\mathrm{int},rp^{-ha}}_{*,A},\widetilde{\Pi}^{\mathrm{bd},r}_{*,A}\rightarrow \widetilde{\Pi}^{\mathrm{bd},rp^{-ha}}_{*,A},\widetilde{\Pi}^{r}_{*,A}\rightarrow \widetilde{\Pi}^{rp^{-ha}}_{*,A}.	
\end{align} 
Then we have these are the quotients of finite projective ones again endowed with the corresponding Frobenius actions.
\end{proposition}

\begin{proof}
See 1.5.2 of \cite[Lemma 1.5.2]{3KL15}.	
\end{proof}

%\newpage
%%%%%%%%%%%%%%%%%%%%%%%%%%%%%%%%%%%%%%T

\newpage\section{Comparison Theorems}

\subsection{The Comparison between the Local Systems and the Period Modules}

\indent We now in our context study the corresponding local systems in the generalized setting. The objects to compare will be definitely the Frobenius modules defined in the previous section. However this is not that far from the situation studied in \cite{3KL16} since the undeformed period rings are actually finite projective over the rings defined in \cite{3KL16} in the most simplified situation. We will also use the notation $\underline{L}$ to denote the local system associated to some topological ring $L$. We will consider the pro-\'etale site $X_\text{pro\'et}$. As in \cite[Definition 4.5.1]{3KL16} we have the notion of projective, pseudocoherent and fpd $L$-local systems in our context over the site $X_\text{pro\'et}$.

\begin{proposition} \mbox{\bf{(After Kedlaya-Liu \cite[Lemma 4.5.4]{3KL16})}} \label{proposition4.1}
The Frobenius invariances give rise to the following exact sequences:
\[
\xymatrix@R+0pc@C+0pc{
0\ar[r]\ar[r]\ar[r] &\underline{\mathcal{O}_{E_a^{\varphi^a}}}
\ar[r]\ar[r]\ar[r] &\widetilde{\Omega}_{*}^\mathrm{int}
\ar[r]\ar[r]\ar[r] &\widetilde{\Omega}_{*}^\mathrm{int} \ar[r]\ar[r]\ar[r] &0,
}
\]
\[
\xymatrix@R+0pc@C+0pc{
0\ar[r]\ar[r]\ar[r] &\underline{\mathcal{O}_{E_a^{\varphi^a}}}
\ar[r]\ar[r]\ar[r] &\widetilde{\Pi}_{*}^\mathrm{int}
\ar[r]\ar[r]\ar[r] &\widetilde{\Pi}_{*}^\mathrm{int} \ar[r]\ar[r]\ar[r] &0,
}
\]
\[
\xymatrix@R+0pc@C+0pc{
0\ar[r]\ar[r]\ar[r] &\underline{E_a}^{\varphi^a}
\ar[r]\ar[r]\ar[r] &\widetilde{\Pi}_{*}^\mathrm{bd}
\ar[r]\ar[r]\ar[r] &\widetilde{\Pi}_{*}^\mathrm{bd} \ar[r]\ar[r]\ar[r] &0,
}
\]
\[
\xymatrix@R+0pc@C+0pc{
0\ar[r]\ar[r]\ar[r] &\underline{E_a}^{\varphi^a}
\ar[r]\ar[r]\ar[r] &\widetilde{\Pi}_{*}
\ar[r]\ar[r]\ar[r] &\widetilde{\Pi}_{*} \ar[r]\ar[r]\ar[r] &0,
}
\]
where in each exact sequence the third arrow represents the morphism $\varphi^a-1$.
	
\end{proposition}

\begin{proof}
It is actually \cite[Lemma 4.5.4]{3KL16}, although we consider some bigger base $k$, we reproduce the argument for the convenience of the reader. So it is obviously we have the exactness at the first positions which are not zero. The exactness in the middle could actually be derived from \cite[Lemma 5.2.4]{3KL15}. Then for the exactness in between the third and the last arrows one may follow the proof in \cite[Lemma 4.5.3]{3KL16} to prove this. For the first sequence see \cite[Lemma 4.5.3]{3KL16}. For the second sequence, relying on the first sequence one only needs to prove that suppose image of some $y\in \widetilde{\Omega}^\mathrm{int}_R$ under the map $\varphi^a-1$ is $x$ living in $\widetilde{\Pi}^\mathrm{int}_R\subset \widetilde{\Omega}^\mathrm{int}_R$ then we have that the preimage $y$ is not only in $\widetilde{\Omega}^\mathrm{int}_R$ but also actually in $\widetilde{\Pi}^\mathrm{int}_R$. This is because we have for each $n$ and the corresponding expressions $x=\sum_{n\geq 0}\pi^n[\overline{x}_n]$ and $y=\sum_{n\geq 0}\pi^n[\overline{y}_n]$, the coefficients for $y$ could be controlled by those of $x$. The corresponding statement for the last sequence is proved by reducing to the second one, namely put:
\begin{displaymath}
x=y+(\varphi^a-1)(-\sum_{k=0}\sum_{\ell=0}\pi^{i_k}[\overline{x}_{i_k}^{p_\ell}]),	
\end{displaymath}
where one may beforehand put $x$ into a summation of some element in the bounded Robba ring and some element of the form of $\sum_{n_k}\pi^{n_k}[\overline{x}_{i_k}]$ with $\overline{x}_{i_k}$ living in $R^{\circ\circ}$.
\end{proof}

We then make the following discussion as well, as in \cite[Lemma 4.5.4]{3KL16} in the following:

\begin{proposition} \mbox{\bf{(After Kedlaya-Liu \cite[Lemma 4.5.4]{3KL16})}}\\
We consider the finitely presented Frobenius modules defined over the ring $\widetilde{\Omega}^\mathrm{int}_{R}$ or $\widetilde{\Pi}^\mathrm{int}_R$. Then in our context, we have: (1) First the Fitting ideal of $M$ is generated by the elements which are in general form of $\pi^n[\overline{x}_n]$ where $\overline{x}_n$ is some idempotent. (2) And we have that the $\pi$-torsion submodule $T$ of $M$ is finitely generated, therefore could be annihilated by some power $\pi^n$ with some integer $n\geq 0$. (3) And the quotient $M/T$ is then projective which implies that the module $T$ is finitely presented. (4) And the submodule $T$ is then finite union of Frobenius modules taking the form of $?/\pi^k?$ for integer $k\geq 0$, where $?$ represents the ring $\widetilde{\Omega}^\mathrm{int}_{R}$ or $\widetilde{\Pi}^\mathrm{int}_R$. (5) Finally we have that the corresponding module $M$ is strictly pseudocoherent. 	
\end{proposition}

\begin{proof}
Again this is basically \cite[Lemma 4.5.4]{3KL16}, although we consider a larger base $k$. We reproduce the argument for the convenience of the reader. Then corresponding property of the Fitting ideal follows from the applying \cite[Proposition 3.2.13]{3KL15}. Then we apply \cref{prop3.18} to choose a projective covering of $M$ taking the form of $P\rightarrow M$ which further induces a map from $P\rightarrow M/T$, then we use the notation $N$ to denote the kernel. Taking quotient by $\pi$ we can form the exact sequence in our context due to the $\pi$-torsion freeness:
\[
\xymatrix@R+0pc@C+0pc{
0\ar[r]\ar[r]\ar[r] &N/\pi N
\ar[r]\ar[r]\ar[r] &P/\pi P \ar[r]\ar[r]\ar[r] &(M/T)/\pi (M/T) \ar[r]\ar[r]\ar[r] &0.
}
\]
Then as in \cite[Lemma 4.5.4]{3KL16} one can then in the situation where the base is $\widetilde{\Omega}_R^\mathrm{int}$ under the $\pi$-adic topology lift the generators of $N/\pi N$ by using the converging series of finite sums of the generators from the quotient. Then in the situation where the base is $\widetilde{\Pi}_R^\mathrm{int}$ we use the Fr\'echet topology to run the same argument. $T$ is a union of each $T_n$ for each $n$, which are the corresponding $\pi^n$-torsion components of $T$. Then by using \cite[Lemma 1.1.5]{3KL16} one derives that the module $M/T$ is then finitely presented, which further implies that the corresponding module $T$ is finite, which is just equal to some $T_n$. Then we apply \cite[Lemma 1.1.5]{3KL16} to finish the proof.
\end{proof}

\begin{corollary}
Let $M$ be a pseudocoherent Frobenius module over the period ring $\widetilde{\Omega}_R^\mathrm{int}$ or the period ring $\widetilde{\Pi}_R^\mathrm{int}$. Then we have that $M$ admits a projective resolution of length at most $1$.  
\end{corollary}

\indent Then we generalize the comparison of local systems and Frobenius modules in \cite[Theorem 4.5.7]{3KL16} to our context as in the following:

\begin{theorem} \mbox{\bf{(After Kedlaya-Liu \cite[Theorem 4.5.7]{3KL16})}} \label{theorem4.4}
The following categories are equivalent:\\
I. The category of projective, pseudocoherent or fpd $\mathcal{O}_{E_a^{\varphi^a}}$-local systems over $X_\text{pro\'et}$;\\
II. The category of projective, pseudocoherent or fpd Frobenius modules over the period ring $\widetilde{\Omega}_R^\mathrm{int}$;\\
III. The category of projective, pseudocoherent or fpd Frobenius modules over the period ring $\widetilde{\Pi}_R^\mathrm{int}$;\\
IV.  The category of projective, pseudocoherent or fpd Frobenius modules over the period sheaf $\widetilde{\Omega}_*^\mathrm{int}$ over $X$, $X_\text{\'et}$, $X_\text{pro\'et}$;\\
V.  The category of projective, pseudocoherent or fpd Frobenius modules over the period sheaf $\widetilde{\Pi}_*^\mathrm{int}$ over $X$, $X_\text{\'et}$, $X_\text{pro\'et}$.\\ 
Under further deformation we do not have the corresponding equivalence as above but we have the following statement:\\
1. There is a fully faithful embedding functor from the category of the projective, pseudocoherent or fpd $\mathcal{O}_{E_a^{\varphi^a}}\widehat{\otimes}\mathcal{O}_A$-local systems over $X_\text{pro\'et}$ to the following two categories:\\
1(a).  The category of projective, pseudocoherent or fpd Frobenius modules over the period sheaf $\widetilde{\Omega}_{*,A}^\mathrm{int}$ over $X$, $X_\text{\'et}$, $X_\text{pro\'et}$;\\
1(b).  The category of projective, pseudocoherent or fpd Frobenius modules over the period sheaf $\widetilde{\Pi}_{*,A}^\mathrm{int}$ over $X$, $X_\text{\'et}$, $X_\text{pro\'et}$.\\ 
2. We have that the following two categories are equivalent:\\
2(a). The category of projective, pseudocoherent or fpd Frobenius modules over the period ring $\widetilde{\Omega}_{R,\infty}^\mathrm{int}$;\\ 
2(b). The category of projective, pseudocoherent or fpd Frobenius modules over the period sheaf $\widetilde{\Omega}_{*,\infty}^\mathrm{int}$ over $X$, $X_\text{\'et}$, $X_\text{pro\'et}$.\\
3. Furthermore we have the following two categories are equivalent:\\
3(a). The category of projective, pseudocoherent or fpd Frobenius modules over the period ring $\widetilde{\Pi}_{R,\infty}^\mathrm{int}$;\\
3(b). The category of projective, pseudocoherent or fpd Frobenius modules over the period sheaf $\widetilde{\Pi}_{*,\infty}^\mathrm{int}$ over $X$, $X_\text{\'et}$, $X_\text{pro\'et}$.

\end{theorem}

\begin{proof}
The first main statement, see \cite[Theorem 4.5.7]{3KL16}. Then we proceed to consider the rest three main statements along the line of the proof of \cite[Theorem 4.5.7]{3KL16}. For the first statement the corresponding functor is just the base change functor:
\begin{displaymath}
\mathbb{L}\mapsto \mathbb{L}\otimes_{\underline{\mathcal{O}_{E_a^{\varphi^a}}\widehat{\otimes}\mathcal{O}_A}}\widetilde{\Omega}_{*,A}^\mathrm{int}	
\end{displaymath}
and 
\begin{displaymath}
\mathbb{L}\mapsto \mathbb{L}\otimes_{\underline{\mathcal{O}_{E_a^{\varphi^a}}\widehat{\otimes}\mathcal{O}_A}}\widetilde{\Pi}_{*,A}^\mathrm{int},	
\end{displaymath}
these are basically fully faithful by applying the corresponding deformed versions (which is just directly taking the complete tensor product with $A$) of the corresponding exact sequences in \cref{proposition4.1}. For the second and the last statements we follow the strategy of the proof of \cite[Theorem 4.5.7]{3KL16} to proceed as in the following. First for the second one, this is direct consequence of \cref{lemma3.11} and \cref{proposition3.13}. Then the third statement could be reduced to the statement by applying the second statement to the corresponding completion of $R[T^{\pm p^{-\infty}}]$ or $R[T^{p^{-\infty}}]$. 
\end{proof}

\indent Then as in \cite[Corollary 4.5.8]{3KL16} we consider the setting of more general adic spaces:

\begin{proposition} \mbox{\bf{(After Kedlaya-Liu \cite[Corollary 4.5.8]{3KL16})}}
1. Let $X$	be a preadic space over $E_a$. Then we have that there is fully faithful embedding functor from the category of all the $\mathcal{O}_{E_{a}^{\varphi^a}}\widehat{\otimes}\mathcal{O}_A$-local systems (projective, pseudocoherent or fpd) over $X$, $X_{\text{\'et}}$ and $X_\text{pro\'et}$ to the category of all the projective, pseudocoherent or fpd Frobenis modules over the sheaves over $\widetilde{\Pi}_{R,A}^\mathrm{int}$ over the corresponding sites;\\
2. Let $X$	be a preadic space over $E_{\infty,a}$. Then we have that there is fully faithful embedding functor from the category of all the $\mathcal{O}_{E_{\infty,a}^{\varphi^a}}\widehat{\otimes}\mathcal{O}_A$-local systems (projective, pseudocoherent or fpd) over $X$, $X_{\text{\'et}}$ and $X_\text{pro\'et}$ to the category of all the projective, pseudocoherent or fpd Frobenius modules over the sheaves over $\widetilde{\Pi}_{R,\infty,A}^\mathrm{int}$ over the corresponding sites.\\
\end{proposition}

\begin{proof}
For one, apply the previous theorem. For two, repeat the argument in the proof of the previous theorem.	
\end{proof}

\indent The following definition is kind of generalization of the corresponding one in \cite[Definition 4.5.9]{3KL16}:

\begin{definition}
Now over the ring $\widetilde{\Pi}_{R,A}$ or $\widetilde{\Pi}^\mathrm{bd}_{R,A}$	we call the corresponding Frobenius modules globally \'etale if they arise from the Frobenius modules over the ring $\widetilde{\Pi}^\mathrm{int}_{R,A}$. Now over the ring $\widetilde{\Pi}_{R,\infty,A}$ or $\widetilde{\Pi}^\mathrm{bd}_{R,\infty,A}$	we call the corresponding Frobenius modules globally \'etale if they arise from the Frobenius modules over the ring $\widetilde{\Pi}^\mathrm{int}_{R,\infty,A}$. \\
Now over the sheaf $\widetilde{\Pi}_{*,A}$ or $\widetilde{\Pi}^\mathrm{bd}_{*,A}$	we call the corresponding Frobenius modules globally \'etale if they arise from the Frobenius modules over the sheaf $\widetilde{\Pi}^\mathrm{int}_{*,A}$. Now over the ring $\widetilde{\Pi}_{*,\infty,A}$ or $\widetilde{\Pi}^\mathrm{bd}_{*,\infty,A}$	we call the corresponding Frobenius modules globally \'etale if they arise from the Frobenius modules over the sheaf $\widetilde{\Pi}^\mathrm{int}_{*,\infty,A}$.

\end{definition}

\begin{proposition} \mbox{\bf{(After Kedlaya-Liu \cite[Theorem 4.5.11]{3KL16})}}
1. Let $X$ be a preadic space over $E_{,a}$. Then we have that there is a fully faithful embedding of the category of the corresponding $E^{\varphi^a}_{,a}\widehat{\otimes}A$-local systems (in the projective setting) into the category of the corresponding projective Frobenius $\varphi^a$-modules over $\widetilde{\Pi}_{*,A}$;\\
2. Let $X$ be a preadic space over $E_{\infty,a}$. Then we have that there is a fully faithful embedding of the category of the corresponding $E^{\varphi^a}_{\infty,a}\widehat{\otimes}A$-local systems (in the projective setting) into the category of the corresponding projective Frobenius $\varphi^a$-modules over $\widetilde{\Pi}_{*,A}$;\\
\end{proposition}

\begin{proof}
As in \cite[Theorem 4.5.11]{3KL16}, we consider the corresponding base change of the corresponding exact sequence in \cref{proposition4.1} which reflects an exact sequence on the sheaves.
\end{proof}

\begin{remark}
We actually did not get the corresponding $v$-topology involved here, which is due to the fact that the pseudocoherent modules might possibly behavior annoyingly in such topology. But on the other hand, at least one will have the corresponding parallel statements as above over $v$-topology for finite projective objects.\\	
\end{remark}

\subsection{The Comparison beyond the \'Etale Setting}

\indent In this subsection we are going to study the corresponding relationship between the corresponding sheaves over the certain Fargues-Fontaine curves and the corresponding Frobenius modules, which is beyond the corresponding consideration of just \'etale objects as what we considered in the previous subsection. This is to some extend in our situation important since we would like to use different point of views to study essentially equivalent objects. The following result is actually directly reduced to \cite[Theorem 4.6.1]{3KL16}:

\begin{theorem} \mbox{\bf (\cite[Theorem 4.6.1]{3KL16})}\\
In the setting where the space $X$ is $\mathrm{Spa}(R,R^+)$ we have the following categories are equivalent:\\
1.  The category of all the pseudocohrent Frobenius $\varphi^a$-modules over the sheaf $\widetilde{\Pi}^\infty_{*}$;\\
2.  The category of all the pseudocohrent Frobenius $\varphi^a$-modules over the sheaf $\widetilde{\Pi}_{*}$;\\
3.  The category of all the pseudocohrent Frobenius $\varphi^a$-modules over the sheaf $\widetilde{\Pi}^I_{*}$, for some closed interval $I$;\\
4.  The category of all the strictly-pseudocoherent Frobenius $\varphi^a$-modules over $\widetilde{\Pi}^\infty_{R}$ such that for any closed interval $I$ the corresponding base change of any such module $M$ to the module over $\widetilde{\Pi}^I_{R}$ is \'etale-stably pseudocoherent;\\
5.  The category of all the pseudocoherent Frobenius $\varphi^a$-modules over $\widetilde{\Pi}_{R}$ such that each $M$ of such modules comes from a base change from some strictly-pseudocoherent Frobenius $\varphi^a$-module $M'$ over $\widetilde{\Pi}^r_{R}$ for some radius $r>0$ where the corresponding base change to (for any closed interval $I\subset (0,r]$) the module over $\widetilde{\Pi}^I_{R}$ is assumed to be \'etale-stably pseudocoherent;\\
6.  The category of all the pseudocoherent $\varphi^a$-bundles over $\widetilde{\Pi}_{R}$;\\
7.  The category of all the Frobenius $\varphi^a$-modules over $\widetilde{\Pi}^{[s,r]}_{R}$ where we have $0<s\leq r/p^{ha}$;\\
8.  The category of all the pseudocoherent sheaves over the adic version Fargues-Fontaine curve $\mathrm{FF}_{X,\text{\'et}}$ in the \'etale setting;\\
9.  The category of all the pseudocoherent sheaves over the adic version Fargues-Fontaine curve $\mathrm{FF}_{X,\text{pro\'et}}$ in the pro-\'etale setting.\\
\end{theorem}

\begin{proof}
This is \cite[Theorem 4.6.1]{3KL16}. We remind the readers of the corresponding functors involved. First the functors from 4 to 5 and to 7, from 4 to 1 and to 3 are base changes, and the functor from 9 to 8 is the pullback along the corresponding morphism of the sites. The functors from 9 to 7 and 6 are restrictions of the corresponding objects involved. 	
\end{proof}

\indent We do the first generalization to the general ultrametric field:

\begin{theorem} \mbox{\bf{(After Kedlaya-Liu \cite[Theorem 4.6.1]{3KL16})}}
In the setting where the space $X$ is $\mathrm{Spa}(R,R^+)$ we have the following categories are equivalent:\\
1.  The category of all the pseudocohrent Frobenius $\varphi^a$-modules over the sheaf $\widetilde{\Pi}^\infty_{*,\infty}$;\\
2.  The category of all the pseudocohrent Frobenius $\varphi^a$-modules over the sheaf $\widetilde{\Pi}_{*,\infty}$;\\
3.  The category of all the pseudocohrent Frobenius $\varphi^a$-modules over the sheaf $\widetilde{\Pi}^I_{*,\infty}$, for some closed interval $I$;\\
4.  The category of all the strictly-pseudocoherent Frobenius $\varphi^a$-modules over $\widetilde{\Pi}^\infty_{R,\infty}$ such that for any closed interval $I$ the corresponding base change of any such module $M$ to the module over $\widetilde{\Pi}^I_{R,\infty}$ is \'etale-stably pseudocoherent;\\
5.  The category of all the pseudocoherent Frobenius $\varphi^a$-modules over $\widetilde{\Pi}_{R,\infty}$ such that each $M$ of such modules comes from a base change from some strictly-pseudocoherent Frobenius $\varphi^a$-module $M'$ over $\widetilde{\Pi}^r_{R,\infty}$ for some radius $r>0$ where the corresponding base change to (for any closed interval $I\subset (0,r]$) the module over $\widetilde{\Pi}^I_{R,\infty}$ is assumed to be \'etale-stably pseudocoherent;\\
6.  The category of all the pseudocoherent $\varphi^a$-bundles over $\widetilde{\Pi}_{R,\infty}$;\\
7.  The category of all the Frobenius $\varphi^a$-modules over $\widetilde{\Pi}^{[s,r]}_{R,\infty}$ where we have $0<s\leq r/p^{ha}$;\\
8.  The category of all the pseudocoherent sheaves over the adic version Fargues-Fontaine curve $\mathrm{FF}_{X,\infty,\text{\'et}}$ in the \'etale setting;\\
9.  The category of all the pseudocoherent sheaves over the adic version Fargues-Fontaine curve $\mathrm{FF}_{X,\infty,\text{pro\'et}}$ in the pro-\'etale setting.
\end{theorem}

\begin{proof}
The functors are the same as in the previous theorem. The proof is actually encoded into the corresponding proof of the following theorem.	
\end{proof}

\indent We are then going to focus on the deformed version of the corresponding correspondences established above:

\begin{theorem}  \mbox{\bf{(After Kedlaya-Liu \cite[Theorem 4.6.1]{3KL16})}}
In the setting where the space $X$ is $\mathrm{Spa}(R,R^+)$ we have the following first group of categories are equivalent:\\
1.  The category of all the pseudocohrent Frobenius $\varphi^a$-modules over the sheaf $\widetilde{\Pi}^\infty_{*,A}$;\\
2.  The category of all the pseudocohrent Frobenius $\varphi^a$-modules over the sheaf $\widetilde{\Pi}_{*,A}$;\\
3.  The category of all the pseudocohrent Frobenius $\varphi^a$-modules over the sheaf $\widetilde{\Pi}^I_{*,A}$, for some closed interval $I$.\\
\indent Then we have the second group of categories which are equvalent:\\
4.  The category of all the strictly-pseudocoherent Frobenius $\varphi^a$-modules over $\widetilde{\Pi}^\infty_{R,A}$ such that for any closed interval $I$ the corresponding base change of any such module $M$ to the module over $\widetilde{\Pi}^I_{R,A}$ is \'etale-stably pseudocoherent;\\
5.  The category of all the pseudocoherent Frobenius $\varphi^a$-modules over $\widetilde{\Pi}_{R,A}$ such that each $M$ of such modules comes from a base change from some strictly-pseudocoherent Frobenius $\varphi^a$-module $M'$ over $\widetilde{\Pi}^r_{R,A}$ for some radius $r>0$ where the corresponding base change to (for any closed interval $I\subset (0,r]$) the module over $\widetilde{\Pi}^I_{R,A}$ is assumed to be \'etale-stably pseudocoherent;\\
6.  The category of all the pseudocoherent $\varphi^a$-bundles over $\widetilde{\Pi}_{R,A}$;\\
7.  The category of all the Frobenius $\varphi^a$-modules over $\widetilde{\Pi}^{[s,r]}_{R,A}$ where we have $0<s\leq r/p^{ha}$;\\
8.  The category of all the pseudocoherent sheaves over the adic version Fargues-Fontaine curve $\mathrm{FF}_{X,A,\text{\'et}}$ in the \'etale setting;\\
9.  The category of all the pseudocoherent sheaves over the adic version Fargues-Fontaine curve $\mathrm{FF}_{X,A,\text{pro\'et}}$ in the pro-\'etale setting.\\
\end{theorem}

\begin{proof}
First the functors are the same as in the proof of the previous theorem. For the proof, we can follow the idea of the proof of \cite[Theorem 4.6.1]{3KL16} to do so. First the equivalences between 9,8 and 6 could be proved by considering \cref{proposition3.13}. By considering the corresponding Frobenius action, one can establish the corresponding equivalences among them and further 7. Now as in the proof of \cite[Theorem 4.6.1]{3KL16}, we can now prove the equivalence between 4 and 6 as in the following. First as in the proof of \cite[Theorem 4.6.1]{3KL16} one considers the exact base change from the key period ring in 4 to the corresponding key period ring in 6, which implies that the corresponding base change functor is then fully faithful. To show it is also essentially surjective, we pick now an arbitrary Frobenius bundle in the category 6, and consider the corresponding reified quasi-Stein covering of the total space by spaces taking the general form as $\mathrm{Spa}(\Pi^{[sp^{nah},rp^{nah}]}_{R,A},\Pi^{[sp^{nah},rp^{nah}],\mathrm{Gr}}_{R,A})$ with well-located radii $r$ and $s$. This will imply that the corresponding global section of the Frobenius bundle is finitely generated after applying the \cite[Proposition 2.7.16]{3KL16} since the application of Frobenius action will imply the uniform bound of the numbers of generators over each such covering subspace as in \cite[Theorem 4.6.1]{3KL16}. Then having shown this we can as in \cite[Thereom 4.6.1]{3KL16} extract the corresponding desired property of being pseudocoherent for the corresponding desired object in the category 4 by considering the kernel of some finite free covering of the Frobenius bundle. Then as in \cite[Theorem 4.6.1]{3KL16} one can show that we have now the categories 4,5 and 6 are equivalent, and furthermore in the same way one may prove that the categories 1,2 and 3 are equivalent.
\end{proof}

\begin{remark}
We did not manage to compare among then the two groups of categories within this theorem, but we believe this could be achieved by further careful analysis.	Please note that the corresponding sheaves in the categories 1,2,3 are actually required to be \'etale-stably pseudocoherent over the ring $A$ as well. 
\end{remark}

\indent Now over smooth Fr\'echet-Stein algebra $A_\infty(G)$ (as in our previous paper) in the commutative setting where each $A_n(G)$ in the family for each $n\geq 0$ is smooth reduced affinoid algebra in rigid analytic space, we have the following comparison:

\begin{theorem} \mbox{\bf{(After Kedlaya-Liu \cite[Theorem 4.6.1]{3KL16})}} 
In the setting where the space $X$ is $\mathrm{Spa}(R,R^+)$ and the field $E$ is of mixed-characteristic we have the following categories are equivalent:\\
1.  The category of all the families of all the pseudocohrent Frobenius $\varphi^a$-modules over the sheaf $\widetilde{\Pi}^\infty_{*,A_\infty(G)}$;\\
2.  The category of all the families of all the pseudocohrent Frobenius $\varphi^a$-modules over the sheaf $\widetilde{\Pi}_{*,A_\infty(G)}$;\\
3.  The category of all the families of all the pseudocohrent Frobenius $\varphi^a$-modules over the sheaf $\widetilde{\Pi}^I_{*,A_\infty(G)}$, for some closed interval $I$.\\
\indent Then we have the second group of categories which are equvalent:\\
4.  The category of all the families of all the strictly-pseudocoherent Frobenius $\varphi^a$-modules over $\widetilde{\Pi}^\infty_{R,A_\infty(G)}$ such that for any closed interval $I$ the corresponding base change of any such module $M$ to the module over $\widetilde{\Pi}^I_{R,A_\infty(G)}$ is \'etale-stably pseudocoherent;\\
5.  The category of all the families of all the pseudocoherent Frobenius $\varphi^a$-modules over $\widetilde{\Pi}_{R,A_\infty(G)}$ such that each $M$ of such modules comes from a base change from some strictly-pseudocoherent Frobenius $\varphi^a$-module $M'$ over $\widetilde{\Pi}^r_{R,A_\infty(G)}$ for some radius $r>0$ where the corresponding base change to (for any closed interval $I\subset (0,r]$) the module over $\widetilde{\Pi}^I_{R,A_\infty(G)}$ is assumed to be \'etale-stably pseudocoherent;\\
6.  The category of all the families of all the pseudocoherent $\varphi^a$-bundles over $\widetilde{\Pi}_{R,A_\infty(G)}$;\\
7.  The category of all the families of all the Frobenius $\varphi^a$-modules over $\widetilde{\Pi}^{[s,r]}_{R,A_\infty(G)}$ where we have $0<s\leq r/p^{ha}$;\\
8.  The category of all the families of all the pseudocoherent sheaves over the adic version Fargues-Fontaine curve $\mathrm{FF}_{X,A_\infty(G),\text{\'et}}$ in the \'etale setting;\\
9.  The category of all the families of all the pseudocoherent sheaves over the adic version Fargues-Fontaine curve $\mathrm{FF}_{X,A_\infty(G),\text{pro\'et}}$ in the pro-\'etale setting.
\end{theorem}

\begin{proof}
This is the direct consequence of the previous theorem if one considers the corresponding systems of all the objects defined over each $A_n(G)$.	
\end{proof}

\indent We now work over a perfectoid field $R$ and drop the assumption on the affinoid algebras we imposed before:

\begin{theorem} \mbox{\bf{(After Kedlaya-Liu \cite[Theorem 4.6.1]{3KL16})}} 
In the setting where the space $X$ is $\mathrm{Spa}(R,R^+)$ where $R$ is further assumed to be an analytic field, and we let $A$ be a general reduced affinoid algebra, then we have the following categories are equivalent:\\
1.  The category of all the pseudocohrent Frobenius $\varphi^a$-modules over the sheaf $\widetilde{\Pi}^\infty_{*,A}$;\\
2.  The category of all the pseudocohrent Frobenius $\varphi^a$-modules over the sheaf $\widetilde{\Pi}_{*,A}$;\\
3.  The category of all the pseudocohrent Frobenius $\varphi^a$-modules over the sheaf $\widetilde{\Pi}^I_{*,A}$, for some closed interval $I$.\\
\indent Then we have the second group of categories which are equvalent:\\
4.  The category of all the strictly-pseudocoherent Frobenius $\varphi^a$-modules over $\widetilde{\Pi}^\infty_{R,A}$ such that for any closed interval $I$ the corresponding base change of any such module $M$ to the module over $\widetilde{\Pi}^I_{R,A}$ is \'etale-stably pseudocoherent;\\
5.  The category of all the pseudocoherent Frobenius $\varphi^a$-modules over $\widetilde{\Pi}_{R,A}$ such that each $M$ of such modules comes from a base change from some strictly-pseudocoherent Frobenius $\varphi^a$-module $M'$ over $\widetilde{\Pi}^r_{R,A}$ for some radius $r>0$ where the corresponding base change to (for any closed interval $I\subset (0,r]$) the module over $\widetilde{\Pi}^I_{R,A}$ is assumed to be \'etale-stably pseudocoherent;\\
6.  The category of all the pseudocoherent $\varphi^a$-bundles over $\widetilde{\Pi}_{R,A}$;\\
7.  The category of all the Frobenius $\varphi^a$-modules over $\widetilde{\Pi}^{[s,r]}_{R,A}$ where we have $0<s\leq r/p^{ha}$;\\
8.  The category of all the pseudocoherent sheaves over the adic version Fargues-Fontaine curve $\mathrm{FF}_{X,A,\text{\'et}}$ in the \'etale setting;\\
9.  The category of all the pseudocoherent sheaves over the adic version Fargues-Fontaine curve $\mathrm{FF}_{X,A,\text{pro\'et}}$ in the pro-\'etale setting.\\
\end{theorem}

\begin{proof}
See the proof of the previous theorem.	\\
\end{proof}

%%%%%%%%%%%%%%%%%%%%%%%%%%%%%%%%%%%%%%%%%%%%%revise the citation in the following:

\subsection{The Comparison with the Schematic Fargues-Fontaine Cur
ve}

\indent In this section we are going to study the corresponding relationship between the Frobenius modules over period rings and the certain sheaves over the deformed version of schematic Fargues-Fontaine curves. We actually discussed in our previous paper the corresponding results in the setting where $E$ is just $\mathbb{Q}_p$. Now we consider the more general setting and the corresponding story in the equal characteristic. First we consider the following key lemmas which are analogs in our context of \cite[6.2.2-6.2.4]{3KL15}, \cite[Lemma 4.6.9]{3KL16} and \cite[Proposition 2.11]{3T1}:

\begin{lemma} \mbox{\bf{(After Kedlaya-Liu \cite[6.2.2-6.2.4]{3KL15}, \cite[Lemma 4.6.9]{3KL16}})}\\ \mbox{\bf{(And also see \cite[Proposition 2.11]{3T1})}}
For any Frobenius $\varphi^a$-bundle $M$ over $\widetilde{\Pi}_{R,A}$, then we have that for any interval $I=[s,r]$ where $0<s\leq r$ the map $\varphi^a-1: M_{[s,rq]}\rightarrow M_{[s,r]}$ is surjective after taking some Frobenius twist as in \cite[6.2.2]{3KL15} and our previous work namely the new morphism $\varphi^a-1: M(n)_{[s,rq]}\rightarrow M(n)_{[s,r]}$ for sufficiently large $n\geq 1$. As in the previous established version one may have the chance to take the number to be $1$ if the bundle initially comes from the corresponding base change from the integral Robba ring.	
\end{lemma}

\begin{proof}
See the proof of \cite[Proposition 2.11]{3T1}.	
\end{proof}

\begin{lemma} \mbox{\bf{(After Kedlaya-Liu \cite[6.2.2-6.2.4]{3KL15}, \cite[Lemma 4.6.9]{3KL16}})}\\ \mbox{\bf{(And also see \cite[Proposition 2.11]{3T1})}}
For any finitely generated bundle $M$ carrying the Frobenius action over $\widetilde{\Pi}_{R,A}$, then we have that for any interval $I=[s,r]$ where $0<s\leq r$ the map $\varphi^a-1: M_{[s,rq]}\rightarrow M_{[s,r]}$ is surjective after taking some Frobenius twist as in \cite[6.2.2]{3KL15} and our previous work namely the new morphism $\varphi^a-1: M(n)_{[s,rq]}\rightarrow M(n)_{[s,r]}$ for sufficiently large $n\geq 1$. As in the previous established version one may have the chance to take the number to be $1$ if the bundle initially comes from the corresponding base change from the integral Robba ring.	
\end{lemma}

\begin{proof}
See the proof of \cite[Proposition 2.11]{3T1}.	
\end{proof}

\begin{lemma} \mbox{\bf{(After Kedlaya-Liu \cite[6.2.2-6.2.4]{3KL15}, \cite[Lemma 4.6.9]{3KL16}})}\\ \mbox{\bf{(And also see \cite[Proposition 2.12]{3T1})}}
For any Frobenius $\varphi^a$-bundle $M$ over $\widetilde{\Pi}_{R,A}$, then we have that for sufficiently large number $n\geq 0$ the module $M$ could be generated by finitely many Frobenius $\varphi^a$-invariant elements of the global sections of $M(n)$.	
\end{lemma}

\begin{proof}
See the proof of \cite[Proposition 2.12]{3T1}.	
\end{proof}

\begin{lemma} \mbox{\bf{(After Kedlaya-Liu \cite[6.2.2-6.2.4]{3KL15}, \cite[Lemma 4.6.9]{3KL16}})}\\ \mbox{\bf{(And also see \cite[Proposition 2.12, Proposition 2.11]{3T1})}}
For any finitely generated bundle $M$ carrying Frobenius action over $\widetilde{\Pi}_{R,A}$, then we have that for sufficiently large number $n\geq 0$ the module $M$ could be generated by finitely many Frobenius $\varphi^a$-invariant elements of the global sections of $M(n)$.	
\end{lemma}

\begin{proof}
See the proof of \cite[Proposition 2.12, Proposition 2.11]{3T1}.	
\end{proof}

\begin{lemma} \mbox{\bf{(After Kedlaya-Liu \cite[6.2.2-6.2.4]{3KL15}, \cite[Lemma 4.6.9]{3KL16}})}\\ \mbox{\bf{(And also see \cite[Proposition 2.12, Proposition 2.11, Proposition 2.14]{3T1})}}\\
For any Frobenius $\varphi^a$-module $M$ over $\widetilde{\Pi}_{R,A}$, then we have that for sufficiently large number $n\geq 1$ the space $H^0_{\varphi^a}(M(n))$ generates $M$ and the space $H^1_{\varphi^a}(M(n))$ vanishes.
\end{lemma}

\begin{proof}
See the proof of \cite[Proposition 2.12, Proposition 2.11, Proposition 2.14]{3T1}.	
\end{proof}

\begin{lemma} \mbox{\bf{(After Kedlaya-Liu \cite[6.2.2-6.2.4]{3KL15}, \cite[Lemma 4.6.9]{3KL16}})}\\ \mbox{\bf{(And also see \cite[Proposition 2.12, Proposition 2.11, Proposition 2.14]{3T1})}}\\
For any finitely generated module $M$ carrying Frobenius action over $\widetilde{\Pi}_{R,A}$, then we have that for sufficiently large number $n\geq 1$ the space $H^0_{\varphi^a}(M(n))$ generates $M$ and the space $H^1_{\varphi^a}(M(n))$ vanishes.
\end{lemma}

\begin{proof}
See the proof of \cite[Proposition 2.12, Proposition 2.11, Proposition 2.14]{3T1}.	
\end{proof}

\indent Then we have the following key corollary which is analog of \cite[Corollary 6.2.3, Lemma 6.3.3]{3KL15}, \cite[Corollary 4.6.10]{3KL16} and \cite[Corollary 2.13, Corollary 3.4, Proposition 3.9]{3T1}:

\begin{corollary} \label{coro4.16}
I. When $M_\alpha,M,M_\beta$ are three Frobenius $\varphi^a$-bundles over the ring $\widetilde{\Pi}_{R,A}$ then we have that for sufficiently large $n\geq 0$ we have the following exact sequence:
\[
\xymatrix@R+0pc@C+0pc{
0\ar[r]\ar[r]\ar[r] &M_{\alpha,I}(n)^{\varphi^a}
\ar[r]\ar[r]\ar[r] &M_I(n)^{\varphi^a}
\ar[r]\ar[r]\ar[r] &M_{\beta,I}(n)^{\varphi^a} \ar[r]\ar[r]\ar[r] &0,
}
\]	
for each interval $I$ (same holds true for finite objects). \\
II. When $M_\alpha,M,M_\beta$ are three Frobenius $\varphi^a$-modules over the ring $\widetilde{\Pi}_{R,A}$ then we have that for sufficiently large $n\geq 0$ we have the following exact sequence:
\[
\xymatrix@R+0pc@C+0pc{
0\ar[r]\ar[r]\ar[r] &M_{\alpha,I}(n)^{\varphi^a}
\ar[r]\ar[r]\ar[r] &M_I(n)^{\varphi^a}
\ar[r]\ar[r]\ar[r] &M_{\beta,I}(n)^{\varphi^a} \ar[r]\ar[r]\ar[r] &0,
}
\]	
for each interval $I$ (same holds true for finite objects).\\
III. For a Frobenius $\varphi^a$-bundle $M$ over $\widetilde{\Pi}_{R,A}$, and for each module $M_I$ over some $\widetilde{\Pi}_{R,A}^{I}$, we put for any element $f$ such that $\varphi^af=p^df$:
\begin{displaymath}
M_{I,f}:=\bigcup_{n\in \mathbb{Z}}f^{-n}M_I(dn)^{\varphi^a}.	
\end{displaymath}
Then with this convention suppose we have three Frobenius $\varphi^a$-bundles taking the form of $M_\alpha,M,M_\beta$ over
$\widetilde{\Pi}_{R,A}$, then for each closed interval $I$ we have the following is an exact sequence:
\[
\xymatrix@R+0pc@C+0pc{
0\ar[r]\ar[r]\ar[r] &M_{\alpha,I,f}\ar[r]\ar[r]\ar[r] &M_{I,f}
\ar[r]\ar[r]\ar[r] &M_{\beta,I,f} \ar[r]\ar[r]\ar[r] &0,
}
\]
where each module in the exact sequence is now a module over $\widetilde{\Pi}_{R,A}[1/f]^{\varphi^a}$ (same holds true for finite objects). \\
IV. For a Frobenius $\varphi^a$-module $M$ over $\widetilde{\Pi}_{R,A}$, we put for any element $f$ such that $\varphi^af=p^df$:
\begin{displaymath}
M_{f}:=\bigcup_{n\in \mathbb{Z}}f^{-n}M(dn)^{\varphi^a}.	
\end{displaymath}
Then with this convention suppose we have three Frobenius $\varphi^a$-modules taking the form of $M_\alpha,M,M_\beta$ over
$\widetilde{\Pi}_{R,A}$, then we have the following is an exact sequence:
\[
\xymatrix@R+0pc@C+0pc{
0\ar[r]\ar[r]\ar[r] &M_{\alpha,f}\ar[r]\ar[r]\ar[r] &M_{f}
\ar[r]\ar[r]\ar[r] &M_{\beta,f} \ar[r]\ar[r]\ar[r] &0,
}
\]
where each module in the exact sequence is now a module over $\widetilde{\Pi}_{R,A}[1/f]^{\varphi^a}$ (same holds true for finite objects). \\
V. Suppose $M$ is now a pseudocoherent Frobenius $\varphi^a$-module over $\widetilde{\Pi}_{R,A}$, then we have that the corresponding module $M_f$ then is also pseudocoherent over the ring $\widetilde{\Pi}_{R,A}[1/f]^{\varphi^a}$. For a Frobenius $\varphi^a$-bundle $M$ over $\widetilde{\Pi}_{R,A}$, and for each module $M_I$ over some $\widetilde{\Pi}_{R,A}^{I}$, we put for any element $f$ such that $\varphi^af=p^df$:
\begin{displaymath}
M_{I,f}:=\bigcup_{n\in \mathbb{Z}}f^{-n}M_I(dn)^{\varphi^a}.	
\end{displaymath}
Then with this convention suppose we have three Frobenius $\varphi^a$-bundles taking the form of $M_\alpha,M,M_\beta$ over
$\widetilde{\Pi}_{R,A}$, then for each closed interval $I$ we have the following is an exact sequence:
\[
\xymatrix@R+0pc@C+0pc{
0\ar[r]\ar[r]\ar[r] &M_{\alpha,I,f}\ar[r]\ar[r]\ar[r] &M_{I,f}
\ar[r]\ar[r]\ar[r] &M_{\beta,I,f} \ar[r]\ar[r]\ar[r] &0,
}
\]
where each module in the exact sequence is now in our situation a pseudocoherent module over $\widetilde{\Pi}_{R,A}[1/f]^{\varphi^a}$.\\

 \end{corollary}

\begin{proof}
See the proof of \cite[Corollary 6.2.3, Lemma 6.3.3]{3KL15}, \cite[Corollary 4.6.10]{3KL16} and \cite[Corollary 2.13, Corollary 3.4, Proposition 3.9]{3T1}.	
\end{proof}

\indent Then we recall the following construction both from our previous consideration and the corresponding consideration in \cite[Definition 6.3.1]{3KL15} and \cite[Definition 4.6.11]{3KL16}:

\begin{setting}
Now we consider the graded commutative ring $P_{R,A}$ which is constructed from the ring $\widetilde{\Pi}_{R,A}^+$, or $\widetilde{\Pi}_{R,A}$, or $\widetilde{\Pi}_{R,A}^\infty$ by taking each $\pi^d$-eigenfunction of the operator $\varphi^a$, where $d$ corresponds to the degree of each element in $P_{d,R,A}$. Then for each element $f$ in this graded commutative ring, one considers the affine charts taking the form of $\mathrm{Spec}P_{R,A}[1/f]_0$ which further glues to a projective scheme which is denoted by $\mathrm{Proj}P_{R,A}$ and called the schematic deformed version of the schematic Fargues-Fontaine curve.
\end{setting}

And we have the following analog of \cite{3KL15} and the corresponding construction in \cite[Setting 3.3]{3T1}:

\begin{setting}
Starting from any $\varphi^a$-bundle $M$ over the ring $\widetilde{\Pi}_{R,A}$, one has by the construction in the previous corollary the corresponding module $M_{I,f}$ over the ring $P_{R,A}[1/f]_0$, which then by considering the natural base change construction gives rise to the map:
\begin{displaymath}
M_{I,f}\otimes_{P_{R,A}[1/f]_0}	\widetilde{\Pi}^I_{R,A}[1/f]\rightarrow 
M_{I,f}\otimes_{\widetilde{\Pi}^I_{R,A}}\widetilde{\Pi}^I_{R,A}[1/f]
\end{displaymath}
for each closed interval $I$.
\end{setting}

\begin{proposition} \mbox{\bf{(After Kedlaya-Liu \cite[Theorem 6.3.9]{3KL15})}}
The map defined in the previous setting is a bijection and the corresponding module $M_{I,f}$ for each specific interval $I$ is then in our situation projective of finite type.
\end{proposition}

\begin{proof}
See \cite[Proposition 3.6]{3T1}.	
\end{proof}

\begin{proposition} \mbox{\bf{(After Kedlaya-Liu \cite[Theorem 6.3.12]{3KL15})}} \label{3propsition3.4.25}
We have the following categories are equivalent:\\
I. The category of all the $\varphi^a$-bundles over $\widetilde{\Pi}_{R,A}$;\\
II. The category of all the $\varphi^a$-modules over $\widetilde{\Pi}_{R,A}$;\\
III. The category of all the $\varphi^a$-modules over $\widetilde{\Pi}^\infty_{R,A}$;\\
IV.  The category of all the quasicoherent finite locally free sheaves over the scheme\\ $\mathrm{Proj}P_{R,A}$.	
\end{proposition}

\begin{proof}
Based on the previous proposition, we only need to recall the corresponding functor involved. Starting from an object in the category $IV$ one locally considers the corresponding pullbacks along the natural map from the localized scheme attached to the key ring in the category $III$, which globally glues to a well-defined object in the corresponding module with Frobenius structure in the category $III$, by applying the corresponding analog of \cite[Lemma 6.3.7]{3KL15}. Then the natural base change functor sends the resulting object to some well-defined object in the category $II$. Finally we associate the corresponding $\varphi^a$-bundle to derive the corresponding object in the first category $I$. 
\end{proof}

\begin{proposition} \mbox{\bf{(After Kedlaya-Liu \cite[Theorem 4.6.12]{3KL16})}} \label{3propsition3.4.26}
The natural pullback functor from the Fargues-Fontaine curve (which is schematic) to the scheme associated to the Robba ring $\widetilde{\Pi}^\infty_{R,A}$ establishes an equivalence between the category of all the pseudocoherent sheaves over the deformed version of the schematic Fargues-Fontaine curve (which is schematic) and the category of all the pseudocoherent modules over $\widetilde{\Pi}_{R,A}$ with isomorphisms established by the Frobenius pullbacks. 	
\end{proposition}

\begin{proof}
As in \cite[Theorem 4.6.12]{3KL16} and our previous paper \cite[Proposition 3.10]{3T1}, we start from an object $V$ in the category of the pseudocoherent sheaves over the schematic Fargues-Fontaine curve, which gives rise to the following exact sequence:
\[
\xymatrix@R+0pc@C+0pc{
0\ar[r]\ar[r]\ar[r] &V_2\ar[r]\ar[r]\ar[r] &V_1
\ar[r]\ar[r]\ar[r] &V \ar[r]\ar[r]\ar[r] &0,
}
\]	
where the sheaf $V_1$ is projective of finite type. Then the corresponding functor mentioned in the statement of the proposition establishes then the following exact sequence:
\[
\xymatrix@R+0pc@C+0pc{
 W_2\ar[r]\ar[r]\ar[r] &W_1
\ar[r]\ar[r]\ar[r] &W \ar[r]\ar[r]\ar[r] &0,
}
\]	
in the category of all the pseudocoherent modules over the Robba ring with isomorphisms established by the Frobenius pullbacks, where then $W_1$ is finite projective and $W_2$ is finite generated. Then as in \cite[Theorem 4.6.12]{3KL16} and \cite[Proposition 3.10]{3T1} we consider the kernel $K$ of the map $W_2\rightarrow \mathrm{Ker}(W_1\rightarrow W)$, which gives rise to the following exact sequence:
\[
\xymatrix@R+0pc@C+0pc{
0\ar[r]\ar[r]\ar[r] &K\ar[r]\ar[r]\ar[r] &W_2\ar[r]\ar[r]\ar[r] &W_1
\ar[r]\ar[r]\ar[r] &W \ar[r]\ar[r]\ar[r] &0,
}
\]
which then gives rise to the following exact sequence by \cref{coro4.16}:
\[
\xymatrix@R+0pc@C+0pc{
0\ar[r]\ar[r]\ar[r] &K_f\ar[r]\ar[r]\ar[r] &W_{2,f}\ar[r]\ar[r]\ar[r] &W_{1,f}
\ar[r]\ar[r]\ar[r] &W_f \ar[r]\ar[r]\ar[r] &0,
}
\]
over the ring $\widetilde{\Pi}_{R,A}[1/f]^{\varphi^a}$. Then by taking the corresponding section over the ring $P_{R,A}[1/f]_0$ we have the following commutative diagram:
\[
\xymatrix@R+2pc@C+1.5pc{
&0\ar[r]\ar[r]\ar[r] &V_{2}\Big|_{P_{R,A}[1/f]_0}\ar[d]\ar[d]\ar[d]\ar[r]\ar[r]\ar[r] &V_{1}\Big|_{P_{R,A}[1/f]_0}
\ar[d]\ar[d]\ar[d]\ar[r]\ar[r]\ar[r] &V \Big|_{P_{R,A}[1/f]_0} \ar[d]\ar[d]\ar[d]\ar[r]\ar[r]\ar[r] &0\\
0\ar[r]\ar[r]\ar[r] &K_f\ar[r]\ar[r]\ar[r] &W_{2,f}\ar[r]\ar[r]\ar[r] &W_{1,f}
\ar[r]\ar[r]\ar[r] &W_f \ar[r]\ar[r]\ar[r] &0.\\
}
\]
Now the middle vertical arrow is isomorphism which shows that the right vertical arrow is surjective. Then apply the corresponding construction to the situation where we replace $V$ by $V_2$ we will have the parallel result which imply that the first vertical arrow is also surjective. Now by five lemma we have that the rightmost vertical arrow is injective. Then by applying the corresponding consideration to the situation where we replace $V$ by $V_2$ we could derive the fact that actually the first vertical arrow is also injective. Up to here, we have that the corresponding module $K_f$ is zero which implies by considering \cref{coro4.16} $K$ itself is zero. Then by the previous results on the corresponding finite projective objects we have that functor in our situation sends the pseudocoherent objects to pseudocoherent objects. And from the argument above we have that the functor is exact. Then to finish the composition of the functors from category of Frobenius modules to the sheaves and back from sheaves to the category of Frobenius modules is an equivalence. And the other direction could be directly deduced from the results on the vector bundles.
\end{proof}

%\newpage

\newpage\section{Deformation of Imperfect Period Rings}

\subsection{Key Rings}

\indent We now consider the corresponding imperfectization of the corresponding constructions we considered above, after \cite{3KL16}. The construction in \cite[Chapter 5]{3KL16} is actually quite general, which is sort of direct imperfection of the corresponding perfect period rings and modules off the tower. We will consider the corresponding towers in \cite[Chapter 5]{3KL16}, so we recall the setting as in the following:

\begin{assumption}
In this section, we are going to assume that $k$ is just $\mathbb{F}_{p^h}$.
\end{assumption}

\begin{setting}
Recall from \cite[Chapter 5]{3KL16}, we have the following setting up. The base will be a Banach adic ring $(H,H^+)$ over $\mathfrak{o}_E$ which is assumed to be uniform and carrying the corresponding spectral norm $\alpha$. Then we consider the tower:
\[
\xymatrix@R+0pc@C+0pc{
...\ar[r]\ar[r]\ar[r] &H_{n-1}\ar[r]\ar[r]\ar[r] &H_n
\ar[r]\ar[r]\ar[r] &H_{n+1} \ar[r]\ar[r]\ar[r] &...,
}
\]	
where each map linking two adjacent rings is the corresponding morphism of the corresponding Banach uniform adic rings whose corresponding induced map on the adic spaces is surjective as in the corresponding construction in \cite[Definition 5.1.1]{3KL16} with the corresponding spectral norm $\alpha_n$. The infinite level of the tower will be denoted by $H_\infty$. This is not actually necessarily complete under the multiplicative extension of all the finite level spectral norms $\alpha_n$ namely $\alpha_\infty$, so we need to consider the completed ring $\overline{H}_\infty$. In the situation where the tower is Fontaine perfectoid, we have that following \cite[Definition 5.1.1]{3KL16} that this is called the corresponding perfectoid tower, which gives rise the equal characteristic counterpart $\overline{H'}_\infty$ wit the spectral norm $\overline{\alpha}_\infty$ correspondingly under the perfectoid correspondence in the sense of \cite[Theorem 3.3.8]{3KL16}. Recall that from \cite[Definition 5.1.1]{3KL16} the tower is called finite \'etale if each transition map is finite \'etale.
\end{setting}

\indent Now we describe the corresponding imperfect period rings which we will deform. These rings are those introduced in \cite[Definition 5.2.1]{3KL16} by using series of imperfection processes. Recall in more detail we have:

\begin{setting} \label{setting6.3}
Fix a perfectoid tower $(H_\bullet,H^+_\bullet)$ Recall from \cite[Definition 5.2.1]{3KL16} we have the following different imperfect constructions:\\
A. First we have the ring $\overline{H'}_\infty$, which could give us the corresponding ring $\widetilde{\Omega}^{\mathrm{int}}_{\overline{H'}_\infty}$; \\
B. We then have the ring $\widetilde{\Pi}^{\mathrm{int},r}_{\overline{H'}_\infty}$ coming from $\overline{H'}_\infty$;\\
C. We then have the ring $\widetilde{\Pi}^{\mathrm{int}}_{\overline{H'}_\infty}$ coming from $\overline{H'}_\infty$;\\
D. We then have the ring $\widetilde{\Pi}^{\mathrm{bd},r}_{\overline{H'}_\infty}$ coming from $\overline{H'}_\infty$;\\
E. We then have the ring $\widetilde{\Pi}^{\mathrm{bd}}_{\overline{H'}_\infty}$ coming from $\overline{H'}_\infty$;\\
F. We then have the ring $\widetilde{\Pi}^{r}_{\overline{H'}_\infty}$ coming from $\overline{H'}_\infty$;\\
G. We then have the ring $\widetilde{\Pi}^{[s,r]}_{\overline{H'}_\infty}$ coming from $\overline{H'}_\infty$;\\
H. We then have the ring $\widetilde{\Pi}_{\overline{H'}_\infty}$ coming from $\overline{H'}_\infty$;\\
I. $\Pi^{\mathrm{int},r}_{H}$ comes from the ring $\widetilde{\Pi}^{\mathrm{int},r}_{\overline{H'}_\infty}$ consisting of those elements of $\widetilde{\Pi}^{\mathrm{int},r}_{\overline{H'}_\infty}$ with the requirement that whenever we have $n$ an integer such that $nh>-\log_pr$ we have then $\theta(\varphi^{-n}(x))\in H_n$;\\
J. $\Pi^{\mathrm{int},\dagger}_{H}$ is defined to be the corresponding union of the rings in $I$;\\
K. $\Omega^\mathrm{int}_H$ is defined to be the corresponding period ring coming from the corresponding $\pi$-adic completion of the ring $\Pi^{\mathrm{int},\dagger}_{H}$ in $J$;\\
L. $\breve{\Omega}^{\mathrm{int}}_{H}$ is the ring which is defined to be the union of all the $\varphi^{-n}\Omega^\mathrm{int}_H$;\\
M. $\breve{\Pi}^{\mathrm{int},r}_H$ is then the ring which is defined to be the union of all the $\varphi^{-n}\Pi^{\mathrm{int},p^{hn}r}_H$;\\
N. $\breve{\Pi}^{\mathrm{int},\dagger}_H$ is define to be union of all the $\varphi^{-n}\breve{\Pi}^{\mathrm{int},\dagger}_H$;\\
O. $\widehat{\Omega}^{\mathrm{int}}_{H}$ is defined to be the $\pi$-completion of $\breve{\Omega}^{\mathrm{int}}_{H}$;\\
P. $\widehat{\Pi}^{\mathrm{int},r}_H$ is defined to be the $\mathrm{max}\{\|.\|_{\overline{\alpha}_\infty^r},\|.\|_{\pi-\text{adic}}\}$-completion of $\breve{\Pi}^{\mathrm{int},r}_H$;\\
Q. We then have $\widehat{\Pi}^{\mathrm{int},\dagger}_H$ by taking the union over $r>0$;\\
R. Correspondingly we have $\breve{\Omega}_{H}$, $\breve{\Pi}^{\mathrm{bd},r}_H$, $\breve{\Pi}^{\mathrm{bd},\dagger}_H$ by inverting the element $\pi$;\\
S. Correspondingly we also have $\widehat{\Omega}_{H}$, $\widehat{\Pi}^{\mathrm{bd},r}_H$, $\widehat{\Pi}^{\mathrm{bd},\dagger}_H$ by inverting the corresponding element $\pi$;\\
T. We also have $\Omega_{H}$, $\Pi^{\mathrm{bd},r}_H$, $\Pi^{\mathrm{bd},\dagger}_H$ again by inverting the element $\pi$;\\
U. Taking the $\max\{\|.\|_{\overline{\alpha}_\infty^s},\|.\|_{\overline{\alpha}_\infty^r}\}$ (for $0<s\leq r$) completion of the ring $\Pi^{\mathrm{bd},r}_H$ we have the ring $\Pi^{[s,r]}_{H}$, while taking the Fr\'echet completion with respect to the norm $\|.\|_{\overline{\alpha}^t_\infty}$ for each $0<t\leq r$ we have the ring $\Pi^r_{H}$;\\
V. Taking the union we have the ring $\Pi_H$;\\
W. We use the notation $\breve{\Pi}_H$ to denote the corresponding union of all the $\varphi^{-n}\Pi_H$;\\
X. We use the notation $\breve{\Pi}^{[s,r]}_{H}$ to be the corresponding union of all the $\varphi^{-n}\Pi^{[p^{-hn}s,p^{-hn}r]}_H$;\\
Y. We use the notation $\breve{\Pi}_H^r$ to be the union of all the $\varphi^{-n}\breve{\Pi}_H^{p^{-hn}r}$.
\end{setting}

\indent Then we have the following direct analog of the relative version of the ring defined above (here as before the ring $A$ denotes a reduced affinoid algebra):

\begin{setting} \label{setting6.4}
Now we consider the deformation of the rings above:\\
I. We have the first group of the period rings in the deformed setting:
\begin{align}
\Pi^{\mathrm{int},r}_{H,A},\Pi^{\mathrm{int},\dagger}_{H,A},\Omega^\mathrm{int}_{H,A}, \Omega_{H,A}, \Pi^{\mathrm{bd},r}_{H,A}, \Pi^{\mathrm{bd},\dagger}_{H,A},\Pi^{[s,r]}_{H,A}, \Pi^r_{H,A}, \Pi_{H,A}.
\end{align}
II. We also have the second group of desired rings in the desired setting:
\begin{align}
\breve{\Pi}^{\mathrm{int},r}_{H,A},\breve{\Pi}^{\mathrm{int},\dagger}_{H,A},\breve{\Omega}^\mathrm{int}_{H,A}, \breve{\Omega}_{H,A}, \breve{\Pi}^{\mathrm{bd},r}_{H,A}, \breve{\Pi}^{\mathrm{bd},\dagger}_{H,A},\breve{\Pi}^{[s,r]}_{H,A}, \breve{\Pi}^r_{H,A}, \breve{\Pi}_{H,A}.	
\end{align}
III. We also have the third group:
\begin{align}
\widehat{\Pi}^{\mathrm{int},r}_{H,A},\widehat{\Pi}^{\mathrm{int},\dagger}_{H,A},\widehat{\Omega}^\mathrm{int}_{H,A}, \widehat{\Omega}_{H,A}, \widehat{\Pi}^{\mathrm{bd},r}_{H,A}, \widehat{\Pi}^{\mathrm{bd},\dagger}_{H,A}.	
\end{align}
\end{setting}

\indent Then we have the following globalization as what we did before in the perfect setting:

\begin{definition}
Now consider the sheaf $\mathcal{O}_{\mathfrak{X}}$ of a rigid analytic space $\mathfrak{X}$ in rigid analytic geometry. We have the following three groups of sheaves of rings over $\mathfrak{X}$:
%\begin{align}
%&\Pi^{\mathrm{int},r}_{H,\mathfrak{X}},\Pi^{\mathrm{int},\dagger}_{H,\mathfrak{X}},\Omega^\mathrm{int}_{H,\mathfrak{X}}, \Omega_{H,\mathfrak{X}}, \Pi^{\mathrm{bd},r}_{H,\mathfrak{X}}, \Pi^{\mathrm{bd},\dagger}_{H,\mathfrak{X}},\Pi^{[s,r]}_{H,\mathfrak{X}}, \Pi^r_{H,\mathfrak{X}}, \Pi_{H,\mathfrak{X}},\\
%&\breve{\Pi}^{\mathrm{int},r}_{H,\mathfrak{X}},\breve{\Pi}^{\mathrm{int},\dagger}_{H,\mathfrak{X}},\breve{\Omega}^\mathrm{int}_{H,\mathfrak{X}}, \breve{\Omega}_{H,\mathfrak{X}}, \breve{\Pi}^{\mathrm{bd},r}_{H,\mathfrak{X}}, \breve{\Pi}^{\mathrm{bd},\dagger}_{H,\mathfrak{X}},\breve{\Pi}^{[s,r]}_{H,\mathfrak{X}}, \breve{\Pi}^r_{H,\mathfrak{X}}, \breve{\Pi}_{H,\mathfrak{X}},\\
%&\widehat{\Pi}^{\mathrm{int},r}_{H,\mathfrak{X}},\widehat{\Pi}^{\mathrm{int},\dagger}_{H,\mathfrak{X}},\widehat{\Omega}^\mathrm{int}_{H,\mathfrak{X}}, \widehat{\Omega}_{H,\mathfrak{X}}, \widehat{\Pi}^{\mathrm{bd},r}_{H,\mathfrak{X}}, \widehat{\Pi}^{\mathrm{bd},\dagger}_{H,\mathfrak{X}},	
%\end{align}
\begin{align}
&\Pi^{[s,r]}_{H,\mathfrak{X}}, \Pi^r_{H,\mathfrak{X}}, \Pi_{H,\mathfrak{X}},\\
&\breve{\Pi}^{[s,r]}_{H,\mathfrak{X}}, \breve{\Pi}^r_{H,\mathfrak{X}}, \breve{\Pi}_{H,\mathfrak{X}}
\end{align}
to be:
\begin{align}
&\Pi^{[s,r]}_{H,\mathcal{O}_\mathfrak{X}}, \Pi^r_{H,\mathcal{O}_\mathfrak{X}}, \Pi_{H,\mathcal{O}_\mathfrak{X}},\\
&\breve{\Pi}^{[s,r]}_{H,\mathcal{O}_\mathfrak{X}}, \breve{\Pi}^r_{H,\mathcal{O}_\mathfrak{X}}, \breve{\Pi}_{H,\mathcal{O}_\mathfrak{X}}.
\end{align}	
%\begin{align}
%&\Pi^{\mathrm{int},r}_{H,\mathcal{O}_\mathfrak{X}},\Pi^{\mathrm{int},\dagger}_{H,\mathcal{O}_\mathfrak{X}},\Omega^\mathrm{int}_{H,\mathcal{O}_\mathfrak{X}}, \Omega_{H,\mathcal{O}_\mathfrak{X}}, \Pi^{\mathrm{bd},r}_{H,\mathcal{O}_\mathfrak{X}}, \Pi^{\mathrm{bd},\dagger}_{H,\mathcal{O}_\mathfrak{X}},\Pi^{[s,r]}_{H,\mathcal{O}_\mathfrak{X}}, \Pi^r_{H,\mathcal{O}_\mathfrak{X}}, \Pi_{H,\mathcal{O}_\mathfrak{X}},\\
%&\breve{\Pi}^{\mathrm{int},r}_{H,\mathcal{O}_\mathfrak{X}},\breve{\Pi}^{\mathrm{int},\dagger}_{H,\mathcal{O}_\mathfrak{X}},\breve{\Omega}^\mathrm{int}_{H,\mathcal{O}_\mathfrak{X}}, \breve{\Omega}_{H,\mathcal{O}_\mathfrak{X}}, \breve{\Pi}^{\mathrm{bd},r}_{H,\mathcal{O}_\mathfrak{X}}, \breve{\Pi}^{\mathrm{bd},\dagger}_{H,\mathcal{O}_\mathfrak{X}},\breve{\Pi}^{[s,r]}_{H,\mathcal{O}_\mathfrak{X}}, \breve{\Pi}^r_{H,\mathcal{O}_\mathfrak{X}}, \breve{\Pi}_{H,\mathcal{O}_\mathfrak{X}},\\
%&\widehat{\Pi}^{\mathrm{int},r}_{H,\mathcal{O}_\mathfrak{X}},\widehat{\Pi}^{\mathrm{int},\dagger}_{H,\mathcal{O}_\mathfrak{X}},\widehat{\Omega}^\mathrm{int}_{H,\mathcal{O}_\mathfrak{X}}, \widehat{\Omega}_{H,\mathcal{O}_\mathfrak{X}}, \widehat{\Pi}^{\mathrm{bd},r}_{H,\mathcal{O}_\mathfrak{X}}, \widehat{\Pi}^{\mathrm{bd},\dagger}_{H,\mathcal{O}_\mathfrak{X}}.
%\end{align}	
\end{definition}

\indent Now we discuss some properties of the corresponding deformed version of the imperfect rings in our context, which is parallel to the corresponding discussion we made in the perfect setting and generalizing the corresponding discussion in \cite{3KL16}:

\begin{proposition}\mbox{\bf{(After Kedlaya-Liu \cite[Lemma 5.2.10]{3KL16})}}
For any $0< r_1\leq r_2$ we have the following equality on the corresponding period rings:
\begin{displaymath}
\Pi^{\mathrm{int},r_1}_{H,\mathbb{Q}_p\{T_1,...,T_d\}}\bigcap	\Pi^{[r_1,r_2]}_{H,\mathbb{Q}_p\{T_1,...,T_d\}}=\Pi^{\mathrm{int},r_2}_{H,\mathbb{Q}_p\{T_1,...,T_d\}}.
\end{displaymath}	
\end{proposition}

\begin{proof}
We adapt the argument in \cite[Lemma 5.2.10]{3KL16} to prove this in the situation where $r_1<r_2$ (otherwise this is trivial), again one direction is easy where we only present the implication in the other direction. We take any element
\begin{center}
 $x\in \Pi^{\mathrm{int},r_1}_{H,\mathbb{Q}_p\{T_1,...,T_d\}}\bigcap	\Pi^{[r_1,r_2]}_{H,\mathbb{Q}_p\{T_1,...,T_d\}}$
\end{center} 
and take suitable approximating elements $\{x_i\}$ living in the bounded Robba ring such that for any $j\geq 1$ one can find some integer $N_j\geq 1$ we have whenever $i\geq N_j$ we have the following estimate:
\begin{displaymath}
\left\|.\right\|_{\overline{\alpha}_\infty^{t},\mathbb{Q}_p\{T_1,...,T_d\}}(x_i-x) \leq p^{-j}, \forall t\in [r_1,r_2].	
\end{displaymath}
Then we consider the corresponding decomposition of $x_i$ for each $i=1,2,...$ into a form having integral part and the rational part $x_i=y_i+z_i$ by setting
\begin{center}
 $y_i=\sum_{k=0,i_1,...,i_d}\pi^kx_{i,k,i_1,...,i_d}T_1^{i_1}...T_d^{i_d}$ 
\end{center} 
out of
\begin{center} 
$x_i=\sum_{k=n(x_i),i_1,...,i_d}\pi^kx_{i,k,i_1,...,i_d}T_1^{i_1}...T_d^{i_d}$.
\end{center}
Note that by our initial hypothesis we have that the element $x$ lives in the ring
\begin{center} 
 $\Pi^{\mathrm{int},r_1}_{H,\mathbb{Q}_p\{T_1,...,T_d\}}$ 
\end{center} 
which further implies that 
\begin{displaymath}
\left\|.\right\|_{\overline{\alpha}_\infty^{r_1},\mathbb{Q}_p\{T_1,...,T_d\}}(\pi^kx_{i,k,i_1,...,i_d}T_1^{i_1}...T^{i_d}_d)	\leq p^{-j}.
\end{displaymath}
Therefore we have ${\overline{\alpha}_\infty}(\overline{x}_{i,k,i_1,...,i_d})\leq p^{(k-j)/r_1},\forall k< 0$ directly from this through computation, which implies that then:
\begin{align}
\left\|.\right\|_{\overline{\alpha}_\infty^{r_2},\mathbb{Q}_p\{T_1,...,T_d\}}(\pi^kx_{i,k,i_1,...,i_d}T_1^{i_1}...T^{i_d}_d)	&\leq p^{-k}p^{(k-j)r_2/r_1}\\
	&\leq p^{1+(1-j)r_1/r_1}.
\end{align}
Then one can read off the result directly from this estimate since under this estimate we can have the chance to modify the original approximating sequence $\{x_i\}$ by $\{y_i\}$ which are initially chosen to be in the integral Robba ring, which implies that actually the element $x$ lives in the right-hand side of the identity in the statement of the proposition.
\end{proof}

\begin{proposition} \mbox{\bf{(After Kedlaya-Liu \cite[Lemma 5.2.10]{3KL16})}}\label{proposition5.7}
For any $0< r_1\leq r_2$ we have the following equality on the corresponding period rings:
\begin{displaymath}
\Pi^{\mathrm{int},r_1}_{H,A}\bigcap	\Pi^{[r_1,r_2]}_{H,A}=\Pi^{\mathrm{int},r_2}_{H,A}.
\end{displaymath}	
Here $A$ is some reduced affinoid algebra over $\mathbb{Q}_p$.	
\end{proposition}

\begin{proof}
This is actually not a direct corollary from the corresponding result as in the previous proposition. But since the map from the corresponding Tate algebra to $A$ is strict, which will remain to be so when tensor (completely) with the corresponding Robba rings in the previous proposition (see \cite[2.1.8, Proposition 6]{3BGR}). Then for any element $\overline{x}$ in the intersection on the left hand side, one can consider the corresponding construction to any lift $x$ of $\overline{x}$ with the corresponding approximating estimate by some sequence:
\begin{displaymath}
\left\|.\right\|_{\overline{\alpha}_\infty^{t},\mathbb{Q}_p\{T_1,...,T_d\}}(x_i-x) \leq p^{-j}, \forall t\in [r_1,r_2],	
\end{displaymath}
for any $i\geq N_j$ with some $N_j$ when $j$ is arbitrarily chosen. Then we to each $x_i$ we associate $y_i$ and $z_i$ as in the proof of the previous proposition. Also for $x$ we have that it is living in the integral Robba ring which implies that we have:
\begin{displaymath}
\left\|.\right\|_{\overline{\alpha}_\infty^{r_1},\mathbb{Q}_p\{T_1,...,T_d\}}(\pi^kx_{i,k,i_1,...,i_d}T_1^{i_1}...T^{i_d}_d)	\leq p^{-j}.
\end{displaymath}
Therefore we have ${\overline{\alpha}_\infty}(\overline{x}_{i,k,i_1,...,i_d})\leq p^{(k-j)/r_1},\forall k< 0$ directly from this through computation, which implies that then:
\begin{align}
\left\|.\right\|_{\overline{\alpha}_\infty^{r_2},\mathbb{Q}_p\{T_1,...,T_d\}}(\pi^kx_{i,k,i_1,...,i_d}T_1^{i_1}...T^{i_d}_d)	&\leq p^{-k}p^{(k-j)r_2/r_1}\\
	&\leq p^{1+(1-j)r_1/r_1}.
\end{align}
This shows that actually one can rearrange the corresponding lifting to be some lifting with respect to the ring on the right hand side, which will finish the proof then by the strictness again after the corresponding projection.	
\end{proof}

%%%%%%%%%%%%%%%%%%%%%%%%%%%%%%%%%%%%%%%%%%%%%%%%%

\begin{proposition} \mbox{\bf{(After Kedlaya-Liu \cite[Lemma 5.2.8]{3KL15} and \cite{3KL16})}}
Consider now in our situation the radii $0< r_1\leq r_2$, and consider any element $x\in \Pi^{[r_1,r_2]}_{H,\mathbb{Q}_p\{T_1,...,T_d\}}$. Then we have that for each $n\geq 1$ one can decompose $x$ into the form of $x=y+z$ such that $y\in \pi^n\Pi^{\mathrm{int},r_2}_{H,\mathbb{Q}_p\{T_1,...,T_d\}}$ with $z\in \bigcap_{r\geq r_2}\Pi^{[r_1,r]}_{H,\mathbb{Q}_p\{T_1,...,T_d\}}$ with the following estimate for each $r\geq r_2$:
\begin{displaymath}
\left\|.\right\|_{\overline{\alpha}_\infty^r,\mathbb{Q}_p\{T_1,...,T_d\}}(z)\leq p^{(1-n)(1-r/r_2)}\left\|.\right\|_{\overline{\alpha}_\infty^{r_2},\mathbb{Q}_p\{T_1,...,T_d\}}(z)^{r/r_2}.	
\end{displaymath}

\end{proposition}

\begin{proof}
As in \cite[Lemma 5.2.8]{3KL15} and \cite{3KL16} and in the proof of our previous proposition we first consider those elements $x$ lives in the bounded Robba rings which could be expressed in general as
\begin{center}
 $\sum_{k=n(x),i_1,...,i_d}\pi^kx_{k,i_1,...,i_d}T_1^{i_1}...T_d^{i_d}$.
 \end{center}	
In this situation the corresponding decomposition is very easy to come up with, namely we consider the corresponding $y_i$ as the corresponding series:
\begin{displaymath}
\sum_{k\geq n,i_1,...,i_d}\pi^kx_{k,i_1,...,i_d}T_1^{i_1}...T_d^{i_d}	
\end{displaymath}
which give us the desired result since we have in this situation when focusing on each single term:
\begin{align}
\left\|.\right\|_{\overline{\alpha}_\infty^r,\mathbb{Q}_p\{T_1,...,T_d\}}&(\pi^kx_{k,i_1,...,i_d}T_1^{i_1}...T_d^{i_d})=p^{-k}\overline{\alpha}_\infty(\overline{x}_{k,i_1,...,i_d})^r\\
&=p^{-k(1-r/r_2)}\left\|.\right\|_{\overline{\alpha}_\infty^{r_2},\mathbb{Q}_p\{T_1,...,T_d\}}(\pi^kx_{k,i_1,...,i_d}T_1^{i_1}...T_d^{i_d})^{r/r_2}\\
&\leq p^{(1-n)(1-r/r_2)}\left\|.\right\|_{\overline{\alpha}_\infty^{r_2},\mathbb{Q}_p\{T_1,...,T_d\}}(\pi^kx_{k,i_1,...,i_d}T_1^{i_1}...T_d^{i_d})^{r/r_2}
\end{align}
for all those suitable $k$. Then to tackle the more general situation we consider the approximating sequence consisting of all the elements in the bounded Robba ring as in the usual situation considered in \cite[Lemma 5.2.8]{3KL15} and \cite{3KL16}, namely we inductively construct the following approximating sequence just as:
\begin{align}
\left\|.\right\|_{\overline{\alpha}_\infty^r,\mathbb{Q}_p\{T_1,...,T_d\}}(x-x_0-...-x_i)\leq p^{-i-1}	\left\|.\right\|_{\overline{\alpha}_\infty^r,\mathbb{Q}_p\{T_1,...,T_d\}}(x), i=0,1,..., r\in [r_1,r_2].
\end{align}
Here all the elements $x_i$ for each $i=0,1,...$ are living in the bounded Robba ring, which immediately gives rise to the suitable decomposition as proved in the previous case namely we have for each $i$ the decomposition $x_i=y_i+z_i$ with the desired conditions as mentioned in the statement of the proposition. We first take the series summing all the elements $y_i$ up for all $i=0,1,...$, which first of all converges under the norm $\left\|.\right\|_{\overline{\alpha}_\infty^r,\mathbb{Q}_p\{T_1,...,T_d\}}$ for all the radius $r\in [r_1,r_2]$, and also note that all the elements $y_i$ within the infinite sum live inside the corresponding integral Robba ring $\Pi^{\mathrm{int},r_2}_{H,\mathbb{Q}_p\{T_1,...,T_d\}}$, which further implies the corresponding convergence ends up in $\Pi^{\mathrm{int},r_2}_{H,\mathbb{Q}_p\{T_1,...,T_d\}}$. For the elements $z_i$ where $i=0,1,...$ also sum up to a converging series in the desired ring since combining all the estimates above we have:
\begin{displaymath}
\left\|.\right\|_{\overline{\alpha}_\infty^r,\mathbb{Q}_p\{T_1,...,T_d\}}(z_i)\leq p^{(1-n)(1-r/r_2)}\left\|.\right\|_{\overline{\alpha}_\infty^{r_2},\mathbb{Q}_p\{T_1,...,T_d\}}(x)^{r/r_2}.	
\end{displaymath}
\end{proof}

%\begin{corollary}
%Consider now in our situation the radii $0< r_1\leq r_2$, and consider any element $x\in \widetilde{\Pi}^{[r_1,r_2]}_{R,A}$. Then we have that for each $n\geq 1$ one can decompose $x$ into the form of $x=y+z$ such that $y\in \pi^n\widetilde{\Pi}^{\mathrm{int},r_2}_{R,A}$ with $z\in \bigcap_{r\geq r_2}\widetilde{\Pi}^{[r_1,r]}_{R,A}$ with the following estimate for each $r\geq r_2$:
%\begin{displaymath}
%\left\|.\right\|_{\alpha^r,A}(z)\leq p^{(1-n)(1-r/r_2)}\left\|.\right\|_{\alpha^{r_2},A}(z)^{r/r_2}.	
%\end{displaymath}	
%\end{corollary}
%
%\begin{proof}
%This could be derived directly from the previous proposition, by considering the corresponding residual norms and then the spectral norms.	
%\end{proof}

\begin{proposition} \mbox{\bf{(After Kedlaya-Liu \cite[Lemma 5.2.10]{3KL15})}}
We have the following identity:
\begin{displaymath}
\Pi^{[s_1,r_1]}_{H,\mathbb{Q}_p\{T_1,...,T_d\}}\bigcap\Pi^{[s_2,r_2]}_{H,\mathbb{Q}_p\{T_1,...,T_d\}}=\Pi^{[s_1,r_2]}_{H,\mathbb{Q}_p\{T_1,...,T_d\}},
\end{displaymath}
here the radii satisfy $<s_1\leq s_2 \leq r_1 \leq r_2$.
\end{proposition}

\begin{proof}
In our situation one direction is obvious while on the other hand we consider any element $x$ in the intersection on the left, then by the previous proposition we	have the decomposition $x=y+z$ where $y\in \Pi^{\mathrm{int},r_1}_{H,\mathbb{Q}_p\{T_1,...,T_d\}}$ and $z\in \Pi^{[s_1,r_2]}_{H,\mathbb{Q}_p\{T_1,...,T_d\}}$. Then as in \cite[Lemma 5.2.10]{3KL15} section 5.2 we look at $y=x-z$ which lives in the intersection:
\begin{displaymath}
\Pi^{\mathrm{int},r_1}_{H,\mathbb{Q}_p\{T_1,...,T_d\}}\bigcap	\Pi^{[s_2,r_2]}_{H,\mathbb{Q}_p\{T_1,...,T_d\}}=\Pi^{\mathrm{int},r_2}_{H,\mathbb{Q}_p\{T_1,...,T_d\}}
\end{displaymath}
which finishes the proof.
\end{proof}

\begin{proposition} \mbox{\bf{(After Kedlaya-Liu \cite[Lemma 5.2.10]{3KL15})}}
We have the following identity:
\begin{displaymath}
\Pi^{[s_1,r_1]}_{H,A}\bigcap \Pi^{[s_2,r_2]}_{H,A}=\Pi^{[s_1,r_2]}_{H,A},
\end{displaymath}
here the radii satisfy $<s_1\leq s_2 \leq r_1 \leq r_2$.
	
\end{proposition}

\begin{proof}
See the proof of \cref{proposition5.7}.	
\end{proof}

\begin{remark}
Again one can follow the same strategy to deal with the corresponding equal-characteristic situation.	
\end{remark}

\subsection{Modules and Bundles}

Now we consider the modules and bundles over the rings introduced in the previous subsection. First we make the following assumption:

\begin{setting}
Recall that from \cite[Definition 5.2.3]{3KL16} any tower $(H_\bullet,H_\bullet^+)$	is called weakly decompleting if we have that first the density of the perfection of $H_{\infty}$ in $\overline{H}_\infty$. Here the ring $H_\infty$ is the ring coming from the mod-$\pi$ construction of the ring $\Omega^\mathrm{int}_{H}$, also at the same time one can find some $r>0$ such that the corresponding modulo $\pi$ operation from the ring $\Omega^\mathrm{int}_{H}$ to the ring $H_\infty$ is actually surjective strictly. 
\end{setting}

\begin{assumption}
We now assume that we are basically in the situation where $(H_\bullet,H_\bullet^+)$ is actually weakly decompleting. Also as in \cite[Lemma 5.2.7]{3KL16} we assume we fix some radius $r_0>0$, for instance this will correspond to the corresponding index in the situation we have the corresponding noetherian tower. Recall that a tower is called noetherian if we have some specific radius as above such that we have the strongly noetherian property on the ring $\Pi^{[s,r]}_{H}$ with $[s,r]\subset (0,r_0]$. Under this condition we have that actually we have the corresponding strongly noetherian property on the ring $\Pi^{[s,r]}_{H,A}$ with $[s,r]\subset (0,r_0]$, therefore consequently we have that the ring $(\Pi^{[s,r]}_{H,A},\Pi^{[s,r],+}_{H,A})$ is sheafy. Here the ring $\Pi^{[s,r],+}_{H,A}$ is defined by taking the product construction between $A$ and the ring $\Pi^{[s,r],+}_{H}$ which is defined in \cite[Definition 5.3.2]{3KL16}. We now assume that the tower is then noetherian. For more examples, see \cite[5.3.3]{3KL16}.
\end{assumption}

\indent Then we can start to discuss the corresponding modules over the rings in our deformed setting, first as in \cite[Lemma 5.3.3]{3KL16} the following result should be derived from our construction:

\begin{lemma}\mbox{\bf{(After Kedlaya-Liu \cite[Lemma 5.3.3]{3KL16})}}
We have the following strict isomorphism in our setting, where the corresponding notations are as in \cite[Lemma 5.3.3]{3KL16}:
\begin{align}
\Pi^{[s,r]}_{H,A}\{X/p^{-t}\}/(\pi X-1)\rightarrow \Pi^{[t,r]}_{H,A}\\
\Pi^{[s,r]}_{H,A}\{X/p^{-t}\}/(X-\pi)\rightarrow \Pi^{[s,t]}_{H,A}\\
\Pi^{r}_{H,A}\{X/p^{-s}\}/(X-\pi)\rightarrow \Pi^{s}_{H,A}\\
\Pi^{r}_{H,A}\{X/p^{-s}\}/(\pi X-1)\rightarrow \Pi^{[s,r]}_{H,A},
\end{align}
where the corresponding radii satisfy $0<s\leq t\leq r\leq r_{0}$.
	
\end{lemma}

\begin{proof}
See \cite[Lemma 5.3.3]{3KL16}.	
\end{proof}

We then as in \cite[Lemma 5.3.4]{3KL16} have the following:

\begin{lemma} \mbox{\bf{(After Kedlaya-Liu \cite[Lemma 5.3.4]{3KL16})}}\\
We have the 2-pseudoflatness of the following maps:
\begin{align}
\Pi^{[r_1,r_2]}_{H,A}\rightarrow\Pi^{[r_1,t]}_{H,A}, \Pi^{[r_1,r_2]}_{H,A}\rightarrow\Pi^{[t,r_2]}_{H,A} 	
\end{align}
where we have $0<r_1\leq t\leq r_2$.
	
\end{lemma}

\begin{proof}
See \cite[Lemma 5.3.4]{3KL16}.	
\end{proof}

\begin{remark}
In effect, one can make more strong statement here, due to the fact that we are working in the noetherian setting. To be more precise see the corresponding result (and the proof) of \cite[Corollary 2.6.9]{3KL16}, one can actually have even the flatness of the corresponding maps in the previous lemma as long as we are working over noetherian rings.	
\end{remark}

\indent Then we have the following implication:

\begin{proposition}\label{proposition6.15} \mbox{\bf{(After Kedlaya-Liu \cite[Proposition 5.3.5]{3KL16})}} 
Suppose that we have a module $M$ over the ring $\Pi_{H,A}^{I}$ (which is assumed to be stably pseudocoherent) where the closed interval $I$ admits a covering by finite many closed intervals $I=\cup_{i}I_i$. Then we have in our situation the following exact augmented \v{C}ech complex:
\begin{displaymath}
0\rightarrow M \rightarrow \bigoplus_{i} M_i \rightarrow...,	
\end{displaymath}
where the module $M_i$ is defined to be the base change of the module $M$ to the ring $\Pi_{H,A}^{I_i}$ with respect to each $i$. Also we have that the morphism $\Pi_{H,A}^{I}\rightarrow \bigoplus_i \Pi_{H,A}^{I_i}$ is in our situation effective descent with respect to the corresponding categories of \'etale-stably pseudocoherent modules in the Banach setting and with respect to the corresponding categories of finite projective modules in the Banach setting. 
\end{proposition}

\begin{proof}
As in \cite[Proposition 5.3.5]{3KL16} we only have to to look at the situation of two intervals. Then as in \cite[Proposition 5.3.5]{3KL16} by using the previous lemma one can argue as in \cite[Proposition 5.3.5]{3KL16}.	
\end{proof}

\indent Then we deform the basic notation of bundles in \cite[Definition 5.3.6]{3KL16}:

\begin{definition} \mbox{\bf{(After Kedlaya-Liu \cite[Definition 5.3.6]{3KL16})}}
We define the bundle over the ring $\Pi^{r_0}_{H,A}$ to be a collection $(M_I)_I$ of finite projective modules over each $\Pi_{H,A}^{I}$ with $I\subset (0,r_0]$ closed subintervals of $(0,r_0]$ such that we have the following requirement in the glueing fashion. First for any $I_1\subset I_2$ two closed intervals we have $M_{I_2}\otimes_{\Pi_{H,A}^{I_2}}\Pi_{H,A}^{I_1}\overset{\sim}{\rightarrow} M_{I_1}$ with the obvious cocycle condition with respect to three closed subintervals of $(0,r_0]$ namely taking the form of $I_1\subset I_2\subset I_3$.\\
\indent We define the pseudocoherent sheaf over the ring $\Pi^{r_0}_{H,A}$ to be a collection $(M_I)_I$ of \'etale-stably pseudocoherent modules over each $\Pi_{H,A}^{I}$ with $I\subset (0,r_0]$ closed subintervals of $(0,r_0]$ such that we have the following requirement in the glueing fashion. First for any $I_1\subset I_2$ two closed intervals we have $M_{I_2}\otimes_{\Pi_{H,A}^{I_2}}\Pi_{H,A}^{I_1}\overset{\sim}{\rightarrow} M_{I_1}$ with the obvious cocycle condition with respect to three closed subintervals of $(0,r_0]$ namely taking the form of $I_1\subset I_2\subset I_3$. 	
\end{definition}

\indent We make the following discussion around the corresponding module and sheaf structures defined above.

\begin{lemma} \mbox{\bf{(After Kedlaya-Liu \cite[Lemma 5.3.8]{3KL16})}}
We have the isomorphism between the ring $\Pi^r_{H,A}$ and the inverse limit of the ring $\Pi^{[s,r]}_{H,A}$ with respect to the radius $s$ by the map $\Pi^r_{H,A}\rightarrow \Pi^{[s,r]}_{H,A}$.	
\end{lemma}
 
\begin{proof}
As in \cite[Lemma 5.3.8]{3KL16} it is injective by the isometry, and then use the corresponding elements $x_n,n=0,1,...$ in the dense ring $\Pi^{r}_{H,A}$ to approximate any element $x$ in the ring $\Pi^{[s,r]}_{H,A}$ in the same way as in \cite[Lemma 5.3.8]{3KL16}:
\begin{displaymath}
\|.\|_{\overline{\alpha}_\infty^t,A}(x-x_n)\leq p^{n}	
\end{displaymath}
for any radius $t$ now living in the corresponding interval $[r2^{-n},r]$. This will establish Cauchy sequence which finishes the proof as in \cite[Lemma 5.3.8]{3KL16}.
\end{proof}

\begin{proposition} \mbox{\bf{(After Kedlaya-Liu \cite[Lemma 5.3.9]{3KL16})}}
For some radius $r\in (0,r_0]$. Suppose we have that $M$ is a vector bundle in the general setting or $M$ is a pseudocoherent sheaf in the setting where the tower is noetherian. Then we have that the corresponding global section is actually dense in each section with respect to some closed interval. And then we have the corresponding vanishing result of the first derived inverse limit functor.	
\end{proposition}

\begin{proof}
See \cite[Lemma 5.3.9]{3KL16}.	
\end{proof}

\indent The interesting issue here as in \cite{3KL16} is the corresponding finitely generatedness of the global section of a pseudocoherent sheaf which is actually not guaranteed in general. Therefore as in \cite{3KL16} we have to distinguish the corresponding well-behaved sheaves out from the corresponding category of all the corresponding pseudocoherent sheaves.

\begin{proposition} \mbox{\bf{(After Kedlaya-Liu \cite[Lemma 5.3.10]{3KL16})}}
As in the previous proposition we choose some $r\in (0,r_0]$. Now assume that the corresponding tower $(H_\bullet,H_\bullet^+)$ is noetherian. Now for any pseudocoherent sheaf $M$ defined above we have the following three statements are equivalent. A. The first statement is that one can find a sequence of positive integers $x_1,x_2,...$ such that for any closed interval living inside $(0,r]$ the section of the sheaf with respect this closed interval admits a projective resolution of modules with corresponding ranks bounded by the sequence of integer $x_1,x_2,...$. B. The second statement is that for any locally finite covering of the corresponding interval $(0,r]$ which takes the corresponding form of $\{I_i\}$ one can find a sequence of positive integers $x_1,x_2,...$ such that for any closed interval living inside $\{I_i\}$ the section of the sheaf with respect this closed interval admits a projective resolution of modules with corresponding ranks bounded by the sequence of integer $x_1,x_2,...$. C. Lastly the third statement is that the corresponding global section is a pseudocoherent module over the ring $\Pi^r_{H,A}$. 	
\end{proposition}

\begin{proof}
See \cite[Lemma 5.3.10]{3KL16}.	
\end{proof}

\indent As in \cite[Definition 5.3.11]{3KL16} we call the sheaf satisfying the corresponding equivalent conditions in the proposition above uniform pseudocoherent sheaf. Then we have the following analog of \cite[Lemma 5.3.12]{3KL16}:

\begin{proposition} \mbox{\bf{(After Kedlaya-Liu \cite[Lemma 5.3.12]{3KL16})}}
The global section functor defines the corresponding equivalence between the categories of the following two sorts of objects. The first ones are the corresponding uniform pseudocoherent sheaves over $\Pi^r_{H,A}$. The second ones are those pseudocoherent modules over the ring $\Pi^r_{H,A}$. 	
\end{proposition}

\begin{proof}
See \cite[Lemma 5.3.12]{3KL16}.	
\end{proof}

\subsection{Frobenius Structure and $\Gamma$-Structure on Hodge-Iwasawa M
odules}

\indent Now we consider the corresponding Frobenius actions over the corresponding imperfect rings we defined before, note that the corresponding Frobenius actions are induced from the corresponding imperfect rings in the undeformed situation from \cite{3KL16} which is to say that the Frobenius action on the ring $A$ is actually trivial.

\indent First we consider the corresponding Frobenius modules:

\begin{definition} \mbox{\bf{(After Kedlaya-Liu \cite[Definition 5.4.2]{3KL16})}}
Over the period rings $\Pi_{H,A}$ or $\breve{\Pi}_{H,A}$  (which is denoted by $\triangle$ in this definition) we define the corresponding $\varphi^a$-modules over $\triangle$ which are respectively projective, pseudocoherent or fpd to be the corresponding finite projective, pseudocoherent or fpd modules over $\triangle$ with further assigned semilinear action of the operator $\varphi^a$. We also require that the modules are complete for the natural topology involved in our situation and for any module over $\Pi_{H,A}$ to be some base change of some module over $\Pi^r_{H,A}$ (which will be defined in the following). We also require that the modules are complete for the natural topology involved in our situation and any module over $\widetilde{\Pi}_{H,A}$ to be some base change of some module over $\widetilde{\Pi}^r_{H,A}$ (which will be defined in the following).
\end{definition}

\begin{definition} \mbox{\bf{(After Kedlaya-Liu \cite[Definition 5.4.2]{3KL16})}}
Over each rings $\triangle=\Pi^r_{H,A},\breve{\Pi}^r_{H,A}$ we define the corresponding projective, pseudocoherent or fpd $\varphi^a$-module over any $\triangle$ to be the corresponding finite projective, pseudocoherent or fpd module $M$ over $\triangle$ with additionally endowed semilinear Frobenius action from $\varphi^a$ such that we have the isomorphism $\varphi^{a*}M\overset{\sim}{\rightarrow}M\otimes \square$ where the ring $\square$ is one $\triangle=\Pi^{r/p}_{H,A},\breve{\Pi}^{r/p}_{H,A}$. Also as in \cite[Definition 5.4.2]{3KL16} we assume that the module over $\Pi^r_{H,A}$ is then complete for the natural topology and the corresponding base change to $\Pi^I_{H,A}$ for any interval which is assumed to be closed $I\subset [0,r)$ gives rise to a module over $\Pi^I_{H,A}$ with specified conditions which will be specified below. Also as in \cite[Definition 5.4.2]{3KL16} we assume that the module over $\widetilde{\Pi}^r_{H,A}$ is then complete for the natural topology and the corresponding base change to $\widetilde{\Pi}^I_{H,A}$ for any interval which is assumed to be closed $I\subset [0,r)$ gives rise to a module over $\widetilde{\Pi}^I_{H,A}$ with specified conditions which will be specified below.

\end{definition}

\begin{definition} \mbox{\bf{(After Kedlaya-Liu \cite[Definition 5.4.2]{3KL16})}}
Again as in\\ \cite[Definition 5.4.2]{3KL16}, we define the corresponding projective, pseudocoherent and fpd $\varphi^a$-modules over ring $\Pi^{[s,r]}_{H,A}$ or $\breve{\Pi}^{[s,r]}_{H,A}$  to be the finite projective, pseudocoherent and fpd modules (which will be denoted by $M$) over $\Pi^{[s,r]}_{H,A}$ or $\breve{\Pi}^{[s,r]}_{H,A}$ respectively additionally endowed with semilinear Frobenius action from $\varphi^a$ with the following isomorphisms:
\begin{align}
\varphi^{a*}M\otimes_{\Pi_{H,A}^{[sp^{-ah},rp^{-ah}]}}\Pi_{H,A}^{[s,rp^{-ah}]}\overset{\sim}{\rightarrow}M\otimes_{\Pi_{H,A}^{[s,r]}}\Pi_{H,A}^{[s,rp^{-ah}]}
\end{align}
and
\begin{align}
\varphi^{a*}M\otimes_{\breve{\Pi}_{R,A}^{[sp^{-ah},rp^{-ah}]}}\breve{\Pi}_{R,A}^{[s,rp^{-ah}]}\overset{\sim}{\rightarrow}M\otimes_{\breve{\Pi}_{R,A}^{[s,r]}}\breve{\Pi}_{R,A}^{[s,rp^{-ah}]}.\\
\end{align}
We now assume that the modules are complete with respect to the natural topology and \'etale stably pseudocoherent.
\end{definition}

\indent Also one can further define the corresponding bundles carrying semilinear Frobenius in our context as in the situation of \cite[Definition 5.4.10]{3KL16}:

\begin{definition} \mbox{\bf{(After Kedlaya-Liu \cite[Definition 5.4.10]{3KL16})}}
Over the ring $\Pi^r_{H,A}$ or $\breve{\Pi}^r_{H,A}$ we define a corresponding projective, pseudocoherent and fpd Frobenius bundle to be a family $(M_I)_I$ of finite projective, \'etale stably pseudocoherent and \'etale stably fpd modules over each $\Pi^I_{H,A}$ or $\breve{\Pi}^I_{H,A}$ respectively carrying the natural Frobenius action coming from the operator $\varphi^a$ such that for any two involved intervals having the relation $I\subset J$ we have:
\begin{displaymath}
M_J\otimes_{\Pi^J_{H,A}}\Pi^I_{H,A}\overset{\sim}{\rightarrow}	M_I
\end{displaymath}
and 
\begin{displaymath}
M_J\otimes_{\breve{\Pi}^J_{H,A}}\breve{\Pi}^I_{H,A}\overset{\sim}{\rightarrow}	M_I
\end{displaymath}
with the obvious cocycle condition. Here we have to propose condition on the intervals that for each $I=[s,u]$ involved we have $s\leq u/p^{ah}$. Then as in \cite[Definition 5.4.10]{3KL16}, we can consider the corresponding direct 2-limit to achieve the corresponding objects in the categories over full Robba rings.
\end{definition}

\indent We can then compare the corresponding objects defined above:

\begin{proposition} \mbox{\bf{(After Kedlaya-Liu \cite[Lemma 5.4.11]{3KL16})}}\\
I. Consider the following objects for some radius $r_0$ in our situation. The first group of objects are those finite projective $\varphi^a$-modules over the Robba ring $\Pi^{r_0}_{H,A}$. The second group of objects are those finite projective $\varphi^a$-bundles over the Robba ring $\Pi^{r_0}_{H,A}$. The third group of objects are those finite projective $\varphi^a$-modules over the Robba ring $\Pi^{[s,r]}_{H,A}$ for some $[s,r]\in (0,r_0)$. Then we have that the corresponding categories of the three groups of objects are equivalent. \\
II. Consider the following objects for some radius $r_0$ in our situation. The first group of objects are those pseudocoherent $\varphi^a$-modules over the Robba ring $\Pi^{r_0}_{H,A}$. The second group of objects are those pseudocoherent $\varphi^a$-bundles over the Robba ring $\Pi^{r_0}_{H,A}$. The third group of objects are those pseudocoherent $\varphi^a$-modules over the Robba ring $\Pi^{[s,r]}_{H,A}$ for some $[s,r]\in (0,r_0)$. Then we have that the corresponding categories of the three groups of objects are equivalent. \\
III. Consider the following objects for some radius $r_0$ in our situation. The first group of objects are those finite projective dimensional $\varphi^a$-modules over the Robba ring $\Pi^{r_0}_{H,A}$. The second group of objects are those finite projective dimensional $\varphi^a$-bundles over the Robba ring $\Pi^{r_0}_{H,A}$. The third group of objects are those finite projective dimensional $\varphi^a$-modules over the Robba ring $\Pi^{[s,r]}_{H,A}$ for some $[s,r]\in (0,r_0)$. Then we have that the corresponding categories of the three groups of objects are equivalent.  	
\end{proposition}

\begin{proof}
See the proof of \cite[Lemma 5.4.11]{3KL16}. In our situation we need to use the corresponding \cref{proposition6.15} to prove the corresponding equivalences between the corresponding categories of bundles in the second groups and modules in the third groups. The corresponding equivalences between the corresponding categories of bundles in the second groups and modules in the first groups are following \cite[Lemma 5.4.11]{3KL16} after applying the Frobenius pullbacks to compare those sections over different intervals, and note that we have the corresponding uniform pseudocoherent objects in our development above which actually restricts further to the finite projective objects as in \cite[Proposition 2.7.16]{3KL16}.
\end{proof}

\indent Now we define the corresponding $\Gamma$-modules over the period rings attached to the tower $(H_\bullet,H_\bullet^+)$. The corresponding structures are actually abstractly defined in the same way as in \cite{3KL16}. First we consider the deformation of the corresponding complex $*_{H^\bullet}$ for any ring $*$ in \cref{setting6.4}.

\begin{assumption}
Recall that the corresponding tower is called decompleting if it is weakly decompleting, finite \'etale on each finite level and having the exact sequence $\overline{\varphi}^{-1}H'_{H^\bullet}/H'_{H^\bullet}$ is exact. We now assume that the tower $(H_\bullet,H^+_\bullet)$ is then decompleting.	
\end{assumption}

%\begin{definition}
%Recall from \cite{3KL16} we have that the corresponding notation $*_{H^k}$ for some corresponding nonnegative integer $k>0$ to denote the corresponding completed tensor product of the corresponding periods rings in $\label{setting6.4}$ of order $k+1$ over the given tower with itself. Then one extend the corresponding definition to our situation denoted by $*_{(H^k,A^k}$. We now denote the map $\sharp_{p,k} *_{H^k,A} \rightarrow *_{H^k+1,A}$ by the map induced from the corresponding map 
%\end{definition}

\begin{setting} 
Assume now $\Gamma$ is a topological group as in \cite[Definition 5.5.5]{3KL16} acting on the corresponding period rings in the \cref{setting6.3} in the original context of \cref{setting6.3}. Then we consider the corresponding induced continous action over the corresponding deformed version in our context namely in \cref{setting6.4}. Assume now that the tower is Galois with the corresponding Galois group $\Gamma$.	
\end{setting}

\begin{definition}
We now consider the corresponding inhomogeneous continuous cocycles of the group $\Gamma$, as in \cite[Definition 5.5.5]{3KL16} we use the following notation to denote the corresponding complex extracted from a single tower for a given period ring $*_{H,A}$ in \cref{setting6.4} for each $k>0$:
\begin{displaymath}
*_{H^k,A}:=C_\mathrm{con}(\Gamma^k,*_{H})\widehat{\otimes}_{\sharp} ?
\end{displaymath}
where $\sharp=\mathbb{Q}_p,\mathbb{Z}_p$ and $?=A,\mathfrak{o}_A$ respectively, which forms the corresponding complex $(*_{H^\bullet,A},d^\bullet)$ with the corresponding differential as in \cite[Definition 5.5.5]{3KL16} in the sense of continuous group cohomology.	
\end{definition}

\begin{definition}
Having established the corresponding meaning of the $\Gamma$-structure we now consider the corresponding definition of $\Gamma$-modules. Such modules called the corresponding $\Gamma$-modules are defined over the corresponding rings in \cref{setting6.4}. Again we allow the corresponding modules to be finite projective, or pseudocoherent or fpd over the rings in \cref{setting6.4}. And the modules are defined to be carrying the corresponding continuous semilinear action from the group $\Gamma$.	
\end{definition}

\begin{proposition} \mbox{\bf{(After Kedlaya-Liu \cite[Corollary 5.6.5]{3KL16})}}\\
The complex $\varphi^{-1}\Pi^{[sp^{-h},rp^{-h}]}_{H^\bullet_{\geq n},A}/\Pi^{[s,r]}_{H^\bullet_{\geq n},A}$ and the complex $\widetilde{\Pi}^{[s,r]}_{H^\bullet_{\geq n},A}/\Pi^{[s,r]}_{H^\bullet_{\geq n},A}$ are strict exact for any truncation index $n$. The corresponding radii satisfy the corresponding relation $0<s\leq r\leq r_0$.
\end{proposition}

\begin{proof}
See \cite[Corollary 5.6.5]{3KL16}, and consider the corresponding Schauder basis of $A$.	
\end{proof}

\begin{proposition} \mbox{\bf{(After Kedlaya-Liu \cite[Lemma 5.6.6]{3KL16})}}
The complex 
\begin{displaymath}
M\otimes_{\Pi^{[s,r]}_{H,A}} \varphi^{-(\ell+1)}\Pi^{[sp^{-h(\ell+1)},rp^{-h(\ell+1)}]}_{H^\bullet,A}/ \varphi^{-\ell}\Pi^{[sp^{-h\ell},rp^{-h\ell}]}_{H^\bullet,A}	
\end{displaymath}
and the complex 
\begin{displaymath}
M\otimes_{\Pi^{[s,r]}_{H,A}} \widetilde{\Pi}^{[s,r]}_{H^\bullet,A}/\varphi^{-\ell}\Pi^{[sp^{-h(\ell)},rp^{-h(\ell)}]}_{H^\bullet,A}
\end{displaymath}
are strict exact for any truncation index $n$. The corresponding radii satisfy the corresponding relation $0<s\leq r\leq r_0$, and $\ell$ is bigger than some existing truncated integer $\ell_0\geq  0$. Here $M$ is any $\Gamma$-module in our context.
\end{proposition}

\begin{proof}
See \cite[Lemma 5.6.6]{3KL16}.	
\end{proof}

\begin{proposition} \mbox{\bf{(After Kedlaya-Liu \cite[Lemma 5.6.9]{3KL16})}}
With the corresponding notations as above we have that the corresponding base change from the ring taking the form of $\breve{\Pi}^{[s,r]}_{H,A}$ to the ring taking the form of $\widetilde{\Pi}^{[s,r]}_{H,A}$ establishes the corresponding equivalence on the corresponding categories of $\Gamma$-modules.	
\end{proposition}

\begin{proof}
This is a relative version of the corresponding result in \cite[Lemma 5.6.9]{3KL16} we adapt the corresponding argument here. Indeed the corresponding fully faithfulness comes from the previous proposition the idea to prove the corresponding essential surjectivity comes from writing the module $M$ (the corresponding differential) over the ring taking form of $\widetilde{\Pi}^{[s,r]}_{H,A}$ as the base change from $\varphi^{-k}(\Pi^{[s,r]}_{H,A})$ of a module $M_0$ (the corresponding differential) after the corresponding \cite[Lemma 5.6.8]{3KL16}.	Then as in \cite{3KL16} we consider the corresponding norms on the corresponding $M$ and the base change of $M_0$ which could be controlled up to some constant which could be further modified to be zero by reducing each time positive amount of constant from the constant represented by the difference of the norms.
\end{proof}

\begin{definition} \mbox{\bf{(After Kedlaya-Liu \cite[Definition 5.7.2]{3KL16})}}
Over the period rings $\Pi_{H,A}$ or $\breve{\Pi}_{H,A}$  (which is denoted by $\triangle$ in this definition) we define the corresponding $(\varphi^a,\Gamma)$-modules over $\triangle$ which are respectively projective, pseudocoherent or fpd to be the corresponding finite projective, pseudocoherent or fpd $\Gamma$-modules over $\triangle$ with further assigned semilinear action of the operator $\varphi^a$ with the isomorphism defined by using the Frobenius. We also require that the modules are complete for the natural topology involved in our situation and for any module over $\Pi_{H,A}$ to be some base change of some module over $\Pi^r_{H,A}$ (which will be defined in the following). 
\end{definition}

\begin{definition} \mbox{\bf{(After Kedlaya-Liu \cite[Definition 5.7.2]{3KL16})}}
Over each rings $\triangle=\Pi^r_{H,A},\breve{\Pi}^r_{H,A}$ we define the corresponding projective, pseudocoherent or fpd $(\varphi^a,\Gamma)$-module over any $\triangle$ to be the corresponding finite projective, pseudocoherent or fpd $\Gamma$-module $M$ over $\triangle$ with additionally endowed semilinear Frobenius action from $\varphi^a$ such that we have the isomorphism $\varphi^{a*}M\overset{\sim}{\rightarrow}M\otimes \square$ where the ring $\square$ is one $\triangle=\Pi^{r/p}_{H,A},\breve{\Pi}^{r/p}_{H,A}$. Also as in \cite[Definition 5.7.2]{3KL16} we assume that the module over $\Pi^r_{H,A}$ is then complete for the natural topology and the corresponding base change to $\Pi^I_{H,A}$ for any interval which is assumed to be closed $I\subset [0,r)$ gives rise to a module over $\Pi^I_{H,A}$ with specified conditions which will be specified below. 

\end{definition}

\begin{definition} \mbox{\bf{(After Kedlaya-Liu \cite[Definition 5.7.2]{3KL16})}}
Again as in\\ \cite[Definition 5.7.2]{3KL16}, we define the corresponding projective, pseudocoherent and fpd $(\varphi^a,\Gamma)$-modules over ring $\Pi^{[s,r]}_{H,A}$ or $\breve{\Pi}^{[s,r]}_{H,A}$ to be the finite projective, pseudocoherent and fpd $\Gamma$-modules (which will be denoted by $M$) over $\Pi^{[s,r]}_{H,A}$ or $\breve{\Pi}^{[s,r]}_{H,A}$ respectively additionally endowed with semilinear Frobenius action from $\varphi^a$ with the following isomorphisms:
\begin{align}
\varphi^{a*}M\otimes_{\Pi_{H,A}^{[sp^{-ah},rp^{-ah}]}}\Pi_{H,A}^{[s,rp^{-ah}]}\overset{\sim}{\rightarrow}M\otimes_{\Pi_{H,A}^{[s,r]}}\Pi_{H,A}^{[s,rp^{-ah}]}
\end{align}
and
\begin{align}
\varphi^{a*}M\otimes_{\breve{\Pi}_{R,A}^{[sp^{-ah},rp^{-ah}]}}\breve{\Pi}_{R,A}^{[s,rp^{-ah}]}\overset{\sim}{\rightarrow}M\otimes_{\breve{\Pi}_{R,A}^{[s,r]}}\breve{\Pi}_{R,A}^{[s,rp^{-ah}]}.\\
\end{align}
We also require the corresponding topological conditions as we considered in the Frobenius module situation. 
\end{definition}

\begin{definition} \mbox{\bf{(After Kedlaya-Liu \cite[Definition 5.7.2]{3KL16})}}
Over the ring\\ $\Pi^r_{H,A}$ or $\breve{\Pi}^r_{H,A}$ we define a corresponding projective, pseudocoherent and fpd $(\varphi^a,\Gamma)$ bundle to be a family $(M_I)_I$ of finite projective, pseudocoherent and fpd $\Gamma$-modules over each $\Pi^I_{H,A}$ or $\breve{\Pi}^I_{H,A}$ respectively carrying the natural Frobenius action coming from the operator $\varphi^a$ such that for any two involved intervals having the relation $I\subset J$ we have:
\begin{displaymath}
M_J\otimes_{\Pi^J_{H,A}}\Pi^I_{H,A}\overset{\sim}{\rightarrow}	M_I
\end{displaymath}
and 
\begin{displaymath}
M_J\otimes_{\breve{\Pi}^J_{H,A}}\breve{\Pi}^I_{H,A}\overset{\sim}{\rightarrow}	M_I
\end{displaymath}
with the obvious cocycle condition. Here we have to propose condition on the intervals that for each $I=[v,u]$ involved we have $v\leq u/p^{ah}$. We put the corresponding topological conditions as before when we consider the corresponding Frobenius bundles. And one can take the corresponding 2-limit in the direct sense to define the corresponding objects over the full Robba rings.
\end{definition}

\begin{definition} \mbox{\bf{(After Kedlaya-Liu \cite[Definition 5.7.2]{3KL16})}}
Over the period rings $\widetilde{\Pi}_{H,A}$ (which is denoted by $\triangle$ in this definition) we define the corresponding $(\varphi^a,\Gamma)$-modules over $\triangle$ which are respectively projective, pseudocoherent or fpd to be the corresponding finite projective, pseudocoherent or fpd $\Gamma$-modules over $\triangle$ with further assigned semilinear action of the operator $\varphi^a$ with the isomorphism defined by using the Frobenius. We also require that the modules are complete for the natural topology involved in our situation and for any module over $\widetilde{\Pi}_{H,A}$ to be some base change of some module over $\widetilde{\Pi}^r_{H,A}$ (which will be defined in the following).
\end{definition}

\begin{definition} \mbox{\bf{(After Kedlaya-Liu \cite[Definition 5.7.2]{3KL16})}}
Over each ring $\triangle=\widetilde{\Pi}^r_{H,A}$ we define the corresponding projective, pseudocoherent or fpd $(\varphi^a,\Gamma)$-module over any $\triangle$ to be the corresponding finite projective, pseudocoherent or fpd $\Gamma$-module $M$ over $\triangle$ with additionally endowed semilinear Frobenius action from $\varphi^a$ such that we have the isomorphism $\varphi^{a*}M\overset{\sim}{\rightarrow}M\otimes \square$ where the ring $\square$ is one $\triangle=\widetilde{\Pi}^{r/p}_{H,A}$. Also as in \cite[Definition 5.7.2]{3KL16} we assume that the module over $\widetilde{\Pi}^r_{H,A}$ is then complete for the natural topology and the corresponding base change to $\widetilde{\Pi}^I_{H,A}$ for any interval which is assumed to be closed $I\subset [0,r)$ gives rise to a module over $\widetilde{\Pi}^I_{H,A}$ with specified conditions which will be specified below.

\end{definition}

\begin{definition} \mbox{\bf{(After Kedlaya-Liu \cite[Definition 5.7.2]{3KL16})}}
Again as in \cite[Definition 5.7.2]{3KL16}, we define the corresponding projective, pseudocoherent and fpd $(\varphi^a,\Gamma)$-modules over ring $\widetilde{\Pi}^{[s,r]}_{H,A}$ to be the finite projective, pseudocoherent and fpd $\Gamma$-modules (which will be denoted by $M$) over $\widetilde{\Pi}^{[s,r]}_{H,A}$ additionally endowed with semilinear Frobenius action from $\varphi^a$ with the following isomorphisms:
\begin{align}
\varphi^{a*}M\otimes_{\widetilde{\Pi}_{H,A}^{[sp^{-ah},rp^{-ah}]}}\widetilde{\Pi}_{H,A}^{[s,rp^{-ah}]}\overset{\sim}{\rightarrow}M\otimes_{\widetilde{\Pi}_{H,A}^{[s,r]}}\widetilde{\Pi}_{H,A}^{[s,rp^{-ah}]}.
\end{align}
We also require the corresponding topological conditions as we considered in the Frobenius module situation. 
\end{definition}

\begin{definition} \mbox{\bf{(After Kedlaya-Liu \cite[Definition 5.7.2]{3KL16})}}
Over the ring\\ $\widetilde{\Pi}^r_{H,A}$ we define a corresponding projective, pseudocoherent and fpd $(\varphi^a,\Gamma)$ bundle to be a family $(M_I)_I$ of finite projective, pseudocoherent and fpd $\Gamma$-modules over each $\widetilde{\Pi}^I_{H,A}$ carrying the natural Frobenius action coming from the operator $\varphi^a$ such that for any two involved intervals having the relation $I\subset J$ we have:
\begin{displaymath}
M_J\otimes_{\widetilde{\Pi}^J_{H,A}}\widetilde{\Pi}^I_{H,A}\overset{\sim}{\rightarrow}	M_I
\end{displaymath}
%and 
%\begin{displaymath}
%M_J\otimes_{\widetilde{\Pi}^J_{H,A}}\widetilde{\Pi}^I_{H,A}\overset{\sim}{\rightarrow}	M_I
%\end{displaymath}
with the obvious cocycle condition. Here we have to propose condition on the intervals that for each $I=[v,u]$ involved we have $v\leq u/p^{ah}$. We put the corresponding topological conditions as before when we consider the corresponding Frobenius bundles. And one can take the corresponding 2-limit in the direct sense to define the corresponding objects over the full Robba rings.
\end{definition}

\begin{remark}
In the following proposition, we assume that the corresponding ring $\widetilde{\Pi}^{[s,r]}_{\overline{H}'_\infty,A}$ is sheafy.	
\end{remark}

\begin{proposition} \mbox{\bf{(After Kedlaya-Liu \cite[Theorem 5.7.5]{3KL16})}}
We have now the following categories are equivalence for the corresponding radii $0< s\leq r\leq r_0$ (with the further requirement as in \cite[Theorem 5.7.5]{3KL16} that $s\in (0,r/q]$):\\
1. The category of all the finite projective sheaves over the ring $\widetilde{\Pi}_{\mathrm{Spa}(H_0,H_0^+),A}$, carrying the $\varphi^a$-action;\\
2. The category of all the finite projective sheaves over the ring $\widetilde{\Pi}^r_{\mathrm{Spa}(H_0,H_0^+),A}$, carrying the $\varphi^a$-action;\\
3. The category of all the finite projective sheaves over the ring $\widetilde{\Pi}^{[s,r]}_{\mathrm{Spa}(H_0,H_0^+),A}$, carrying the $\varphi^a$-action.\\
%1b. The category of all the finite projective sheaves over the ring $\widetilde{\Pi}_{\mathrm{Spa}(\overline{H}'_\infty,\overline{H}^+_\infty),A}$, carrying the $\varphi^a$-action;\\
%2b. The category of all the finite projective sheaves over the ring $\widetilde{\Pi}^r_{\mathrm{Spa}(\overline{H}'_\infty,\overline{H}^+_\infty),A}$, carrying the $\varphi^a$-action;\\
%3b. The category of all the finite projective sheaves over the ring $\widetilde{\Pi}^{[s,r]}_{\mathrm{Spa}(\overline{H}'_\infty,\overline{H}^+_\infty),A}$, carrying the $\varphi^a$-action.\\
\indent Then we have the second group of categories which are equvalent:\\
4. The category of all the finite projective quasi-coherent sheaves over corresponding adic Fargues-Fontaine curve in the deformed setting $\mathrm{FF}_{\overline{H}'_\infty,A}$, carrying the corresponding action from the group $\Gamma$ which is assumed to be semilinear and continuous over each section over any affinoid subspace of the whole space which is assumed to be $\Gamma$-invariant;\\
5. The category of all the finite projective modules over the ring $\Pi_{H,A}$, carrying the $(\varphi^a,\Gamma)$-action;\\
6. The category of all the finite projective bundles over the ring $\Pi_{H,A}$, carrying the $(\varphi^a,\Gamma)$-action;\\
7. The category of all the finite projective modules over the ring $\breve{\Pi}_{H,A}$, carrying the $(\varphi^a,\Gamma)$-action;\\
8. The category of all the finite projective bundles over the ring $\breve{\Pi}_{H,A}$, carrying the $(\varphi^a,\Gamma)$-action;\\
9. The category of all the finite projective modules over the ring $\breve{\Pi}^{[s,r]}_{H,A}$, carrying the $(\varphi^a,\Gamma)$-action;\\
10. The category of all the finite projective modules over the ring $\widetilde{\Pi}_{H,A}$, carrying the $(\varphi^a,\Gamma)$-action;\\
11. The category of all the finite projective bundles over the ring $\widetilde{\Pi}_{H,A}$, carrying the $(\varphi^a,\Gamma)$-action;\\
12. The category of all the finite projective modules over the ring $\widetilde{\Pi}^{[s,r]}_{H,A}$, carrying the $(\varphi^a,\Gamma)$-action.\\	
\end{proposition}

\begin{proof}
The corresponding comparisons on the sheaves and bundles in 1-4 and 10-12 are derived in the corresponding context in the perfect setting as what we did in the previous sections. The rest ones are proved exactly the same as \cite[Theorem 5.7.5]{3KL16} by using our development.	
\end{proof}

\begin{remark}
This proposition generalizes the corresponding results in \cite{3KP} including the situation therein considered by Chojecki-Gaisin, while note that we have also included the situation in the equal characteristic situation.	
\end{remark}

\indent Furthermore if one considers the corresponding context where the ring $H'_{H}$ is further assumed to be $F$-(finite projective) then we can discuss the level of pseudocoherent objects.

\begin{lemma}\mbox{\bf{(After Kedlaya-Liu \cite[Lemma 5.8.7]{3KL16})}}
For any radii in our situation namely $0<s\leq r\leq r_0$ we have the following isomorphism:
\begin{displaymath}
\widetilde{\Pi}^{[s,r]}_{H,A}\overset{\sim}{\rightarrow} \Pi_{H,A}^{[s,r]}\oplus (\oplus_{\ell=0}^\infty \varphi^{-(\ell+1)
}\Pi_{H,A}^{[sp^{-h(\ell+1)},rp^{-h(\ell+1)}]}/\varphi^{-\ell
}\Pi_{H,A}^{[sp^{-h\ell},rp^{-h\ell}]})^\wedge.	
\end{displaymath}	
\end{lemma}

\begin{proof}
This is by considering the Schauder basis of $A$. See \cite[Lemma 5.8.7]{3KL16}.	
\end{proof}

\begin{corollary}\mbox{}\\
\mbox{\bf{(After Kedlaya-Liu \cite[Corollary 5.8.8, Corollary 5.8.11]{3KL16})}}
The map for the radii as in the previous lemma is 2-pseudoflat:
\begin{displaymath}
\Pi_{H,A}^{[s,r]}\rightarrow \widetilde{\Pi}^{[s,r]}_{H,A}.	
\end{displaymath}
Also we have that the corresponding base change along this map will preserve the corresponding \'etale stably pseudocoherence.
	
\end{corollary}

\begin{lemma} \mbox{\bf{(After Kedlaya-Liu \cite[Lemma 5.8.14]{3KL16})}}
The complex 
\begin{displaymath}
M\otimes_{\Pi^{[s,r]}_{H,A}} \varphi^{-(\ell+1)}\Pi^{[sp^{-h(\ell+1)},rp^{-h(\ell+1)}]}_{H^\bullet,A}/ \varphi^{-\ell}\Pi^{[sp^{-h\ell},rp^{-h\ell}]}_{H^\bullet,A}	
\end{displaymath}
and the complex 
\begin{displaymath}
M\otimes_{\Pi^{[s,r]}_{H,A}} \widetilde{\Pi}^{[s,r]}_{H^\bullet,A}/\varphi^{-\ell}\Pi^{[sp^{-h(\ell)},rp^{-h(\ell)}]}_{H^\bullet,A}
\end{displaymath}
are strict exact for any truncation index $n$. The corresponding radii satisfy the corresponding relation $0<s\leq r\leq r_0$, and $\ell$ is bigger than some existing truncated integer $\ell_0\geq  0$. Here $M$ is any pseudocoherent $\Gamma$-module in our context.	
\end{lemma}

\begin{proposition}\mbox{\bf{(After Kedlaya-Liu \cite[Lemma 5.8.17]{3KL16})}}
The corresponding base change along the following map is fully faithful for the corresponding pseudocoherent modules carrying the corresponding structure of $\Gamma$-action:
\begin{displaymath}
\breve{\Pi}_{H,A}^{[s,r]}\rightarrow \widetilde{\Pi}^{[s,r]}_{H,A}.	
\end{displaymath}
The image in the essential sense consists of those modules descending to the corresponding ring in the domain of this map when forgetting the corresponding $\Gamma$-action.	
\end{proposition}

\begin{proof}
See \cite[Lemma 5.8.17]{3KL16}.	
\end{proof}

\begin{proposition} \mbox{\bf{(After Kedlaya-Liu \cite[Lemma 5.9.2]{3KL16})}}
Keep the notations as above. For any finitely generated in general module over the ring $\widetilde{\Pi}^{[s,r]}_{H,A}$ carrying the action of $\Gamma$, one can find another module over $\breve{\Pi}^{[s,r]}_{H,A}$ which is now pseudocoherent carrying the action of $\Gamma$ which covers the previous module through a surjective map after taking the corresponding base change to the ring $\widetilde{\Pi}^{[s,r]}_{H,A}$.	
\end{proposition}

\begin{proof}
As in the proof of \cite[Lemma 5.9.2]{3KL16}, one can prove this in the similar fashion by considering the corresponding isomorphism respecting the cocycle requirement coming from the free module $\widetilde{\Lambda}$ in the presentation of a given module $\widetilde{\Delta}$ over the perfect Robba ring:
\begin{displaymath}
\widetilde{\Lambda}\otimes_{i_{0,0}} \widetilde{\Pi}^{[s,r]}_{H^0,A}\rightarrow \widetilde{\Lambda} \otimes_{i_{0,1}} \widetilde{\Pi}^{[s,r]}_{H^1,A}	
\end{displaymath}
coming from our prescribed $\Gamma$-action on the module $\widetilde{\Delta}$ over the perfect Robba ring in the way that one considers a converging sequence of different desired lifts (as in \cite[Lemma 5.6.9]{3KL16}) of 
\begin{displaymath}
\Lambda\otimes_{i_{0,0}} \breve{\Pi}^{[s,r]}_{H^0,A}\rightarrow \Lambda \otimes_{i_{0,1}} \breve{\Pi}^{[s,r]}_{H^1,A}.	
\end{displaymath}	
To get the corresponding desired covering of $\widetilde{\Delta}$ one considers the corresponding finitely generated submodule of the corresponding kernel of the map from $\Lambda$ to $\widetilde{\Delta}$ then take the quotient of $\Lambda$ through this corresponding submodule. Note that in our situation the ring $\breve{\Pi}^{[s,r]}_{H,A}$ is also coherent which finishes the proof as in \cite[Lemma 5.9.2]{3KL16}. 
\end{proof}

\begin{proposition} \mbox{\bf{(After Kedlaya-Liu \cite[Theorem 5.9.4]{3KL16})}} \label{3proposition3.5.51}
We have now the following categories are equivalence for the corresponding radii $0< s\leq r\leq r_0$ (with the further requirement as in \cite[Theorem 5.7.5]{3KL16} that $s\in (0,r/q]$):\\
1. The category of all the pseudocoherent sheaves over the ring $\widetilde{\Pi}_{\mathrm{Spa}(H_0,H_0^+),A}$, carrying the $\varphi^a$-action;\\
2. The category of all the pseudocoherent sheaves over the ring $\widetilde{\Pi}^r_{\mathrm{Spa}(H_0,H_0^+),A}$, carrying the $\varphi^a$-action;\\
3. The category of all the pseudocoherent sheaves over the ring $\widetilde{\Pi}^{[s,r]}_{\mathrm{Spa}(\overline{H}_\infty,\overline{H}^+_\infty),A}$, carrying the $\varphi^a$-action.\\
%1b. The category of all the pseudocoherent sheaves over the ring $\widetilde{\Pi}_{\mathrm{Spa}(\overline{H}'_\infty,\overline{H}^+_\infty),A}$, carrying the $\varphi^a$-action;\\
%2b. The category of all the pseudocoherent sheaves over the ring $\widetilde{\Pi}^r_{\mathrm{Spa}(\overline{H}'_\infty,\overline{H}^+_\infty),A}$, carrying the $\varphi^a$-action;\\
%3b. The category of all the pseudocoherent sheaves over the ring $\widetilde{\Pi}^{[s,r]}_{\mathrm{Spa}(\overline{H}'_\infty,\overline{H}^+_\infty),A}$, carrying the $\varphi^a$-action.\\
\indent Then we have the second group of categories which are equvalent:\\
4. The category of all the pseudocoherent quasi-coherent sheaves over corresponding adic Fargues-Fontaine curve in the deformed setting $\mathrm{FF}_{\overline{H}'_\infty,A}$, carrying the corresponding action from the group $\Gamma$ which is assumed to be semilinear and continuous over each section over any affinoid subspace of the whole space which is assumed to be $\Gamma$-invariant;\\
5. The category of all the pseudocoherent modules over the ring $\Pi_{H,A}$, carrying the $(\varphi^a,\Gamma)$-action;\\
6. The category of all the pseudocoherent bundles over the ring $\Pi_{H,A}$, carrying the $(\varphi^a,\Gamma)$-action;\\
7. The category of all the pseudocoherent modules over the ring $\breve{\Pi}_{H,A}$, carrying the $(\varphi^a,\Gamma)$-action;\\
8. The category of all the pseudocoherent bundles over the ring $\breve{\Pi}_{H,A}$, carrying the $(\varphi^a,\Gamma)$-action;\\
9. The category of all the pseudocoherent modules over the ring $\breve{\Pi}^{[s,r]}_{H,A}$, carrying the $(\varphi^a,\Gamma)$-action;\\
10. The category of all the pseudocoherent modules over the ring $\widetilde{\Pi}_{H,A}$, carrying the $(\varphi^a,\Gamma)$-action;\\
11. The category of all the pseudocoherent bundles over the ring $\widetilde{\Pi}_{H,A}$, carrying the $(\varphi^a,\Gamma)$-action;\\
12. The category of all the pseudocoherent modules over the ring $\widetilde{\Pi}^{[s,r]}_{H,A}$, carrying the $(\varphi^a,\Gamma)$-action.\\	
\end{proposition}

\begin{proof}
This is by using the $A$-relative version of \cite[Lemma 5.9.3]{3KL16} which states that actually the base change functor along
\begin{displaymath}
\breve{\Pi}_{H,A}^{[s,r]}\rightarrow \widetilde{\Pi}^{[s,r]}_{H,A}	
\end{displaymath}	
is then not only fully faithful on the category of pseudocoherent objects as above but also essential surjective. Indeed by the previous proposition in our situation we have that a surjective covering of any module $\widetilde{\Lambda}$ over $\widetilde{\Pi}^{[s,r]}_{H,A}$:
\begin{displaymath}
\widetilde{\Lambda}'\rightarrow \widetilde{\Lambda}\rightarrow 0	
\end{displaymath}
where $\widetilde{\Lambda}'$ descend to the ring $\breve{\Pi}_{H,A}^{[s,r]}$ carrying $\Gamma$-action. Then by applying the same process above to the corresponding kernel of the covering above we have another exact sequence in the following form:
\begin{displaymath}
\widetilde{\Lambda}''\rightarrow \widetilde{\Lambda}'\rightarrow \widetilde{\Lambda}\rightarrow 0	
\end{displaymath}
Here $\widetilde{\Lambda}'$ and $ \widetilde{\Lambda}$ are finitely presented as in \cite[Lemma 5.9.3]{3KL16}. Then we can now consider the corresponding cokernel of the map $\widetilde{\Lambda}''\rightarrow \widetilde{\Lambda}'$ which will present a corresponding desired object carrying $\Gamma$-action over the ring $\breve{\Pi}_{H,A}^{[s,r]}$ which is again pseudocoherent in our context.
\end{proof}

%\newpage

\newpage\section{Organizaion of Noncommutative Setting}

\subsection{Noncommutative Period Rings in Perfect Setting}

\indent Now we assume $A$ to be some noncommutative Banach affinoid algebra over the local fields we consider above. Examples of such rings could be coming from the corresponding context of \cite{3Soi1}, or the corresponding noncommutative Fr\'echet-Stein algebras as in \cite{3ST1}. We will use the notation $E\{Z_1,...,Z_n\}$ to denote the corresponding noncommutative non-noetherian Tate algebra as in the commutative setting.

\begin{remark}
The noncommutative consideration is not quite new, for instance the aspects rooted in \cite{3Zah1}, \cite{3Wit1} and \cite{3Wit2} (and even \cite{3KL16}), but obviously the corresponding noncommutative setting is more complicated than the corresponding commutative setting which we discussed extensively in the previous context, therefore we do not have the chance to reach all the corresponding noncommutative version of the results above. 	
\end{remark}

\begin{remark}
We choose to closely in some parallel way present the corresponding construction in the noncommutative setting, which is parallel in some aspects to the commutative setting we presented above. We definitely won't have the chance to see the complete picture as in the commutative setting and even more consideration and effort has to be made during the corresponding development.	
\end{remark}

%Also we will consider the following setting:
%
%\begin{setting}
%As above $A$ is a noncommutative Banach affinoid algebra in the sense of \cite{3Soi1}, and we use the notation $A_\infty$ F\'echet-Stein algebra which could be written as limit of the rings as $A$.
%\end{setting}
We first deform the corresponding constructions in \cite[Definition 4.1.1]{3KL16}:

\begin{definition}
We first consider the corresponding deformation of the above rings over $\mathbb{Q}_p\{Z_1,...,Z_d\}$. We are going to use the notation $W_\pi(R)_{\mathbb{Q}_p\{Z_1,...,Z_d\}}$ to denote the complete tensor product of $W_\pi(R)$ with the Tate algebra $\mathbb{Q}_p\{Z_1,...,Z_d\}$ consisting of all the element taking the form as:
\begin{displaymath}
\sum_{k\geq 0,i_1\geq 0,...,i_d\geq 0}\pi^k[\overline{x}_{k,i_1,...,i_d}]Z_1^{i_1}Z_2^{i_2}...Z_d^{i_d}	
\end{displaymath}
over which we have the Gauss norm $\left\|.\right\|_{\alpha^r,\mathbb{Q}_p\{Z_1,...Z_d\}}$ for any $r>0$ which is defined by:
\begin{displaymath}
\left\|.\right\|_{\alpha^r,\mathbb{Q}_p\{Z_1,...Z_d\}}(\sum_{k\geq 0,i_1\geq 0,...,i_d\geq 0}\pi^k[\overline{x}_{k,i_1,...,i_d}]Z_1^{i_1}Z_2^{i_2}...Z_d^{i_d}):=\sup_{k\geq 0,i_1\geq 0,...,i_d\geq 0}p^{-k}\alpha(\overline{x}_{k,i_1,...,i_d})^r.	
\end{displaymath}
Then we define the corresponding convergent rings:
\begin{displaymath}
\widetilde{\Omega}^\mathrm{int}_{R,\mathbb{Q}_p\{Z_1,...,Z_d\}}=W_\pi(R)_{\mathbb{Q}_p\{Z_1,...,Z_d\}},	\widetilde{\Omega}_{R,\mathbb{Q}_p\{Z_1,...,Z_d\}}=W_\pi(R)_{\mathbb{Q}_p\{Z_1,...,Z_d\}}[1/\pi].	
\end{displaymath}
Then as in the previous setting we define the corresponding $\mathbb{Q}_p\{Z_1,...,Z_d\}$-relative integral Robba ring $\widetilde{\Pi}_{R,\mathbb{Q}_p\{Z_1,...Z_d\}}^{\mathrm{int},r}$ as the completion of the ring $W_{\pi}(R^+)_{\mathbb{Q}_p\{Z_1,...,Z_d\}}[[R]]$ by using the Gauss norm defined above.
%Then as in the previous settings we can define the corresponding integral Robba ring by considering the subring $\widetilde{\Pi}_{R,\mathbb{Q}_p\{T_1,...T_d\}}^{\mathrm{int},r}$ of $\widetilde{\Omega}_{R,\mathbb{Q}_p\{T_1,...T_d\}}^{\mathrm{int}}$	consisting of all the elements satisfying:
%\begin{displaymath}
%\lim_{k\geq 0,i_1\geq 0,...,i_d\geq 0}\left\|.\right\|_{\alpha^r,\mathbb{Q}_p\{T_1,...T_d\}}(\pi^{k}[\overline{x}_{k,i_1,...,i_d}])=0.	
%\end{displaymath}
Then one can take the corresponding union of all $\widetilde{\Pi}_{R,\mathbb{Q}_p\{Z_1,...,Z_d\}}^{\mathrm{int},r}$ throughout all $r>0$ to define the integral Robba ring $\widetilde{\Pi}_{R,\mathbb{Q}_p\{Z_1,...,Z_d\}}^{\mathrm{int}}$. For the bounded Robba rings we just set $\widetilde{\Pi}_{R,\mathbb{Q}_p\{Z_1,...,Z_d\}}^{\mathrm{bd},r}=\widetilde{\Pi}_{R,\mathbb{Q}_p\{Z_1,...,Z_d\}}^{\mathrm{int},r}[1/\pi]$ and taking the union throughout all $r>0$ to define the corresponding ring $\widetilde{\Pi}_{R,\mathbb{Q}_p\{Z_1,...,Z_d\}}^{\mathrm{bd}}$. Then we define the corresponding Robba ring $\widetilde{\Pi}_{R,\mathbb{Q}_p\{Z_1,...,Z_d\}}^{I}$ with respect to some interval $I\subset (0,\infty)$ by taking the Fr\'echet completion of ${W_\pi(R^+)}_{\mathbb{Q}_p\{Z_1,...,Z_d\}}[[R]][1/\pi]$ with respect to all the norms $\left\|.\right\|_{\alpha^r,\mathbb{Q}_p\{Z_1,...,Z_d\}}$ for all $r\in I$ which means that the corresponding equivalence classes in the completion procedure will be simultaneously Cauchy with respect to all the norms $\left\|.\right\|_{\alpha^r,\mathbb{Q}_p\{Z_1,...,Z_d\}}$ for all $r\in I$. Then we take suitable intervals such as $(0,r]$ and $(0,\infty)$ to define the corresponding Robba rings $\widetilde{\Pi}_{R,\mathbb{Q}_p\{Z_1,...,Z_d\}}^{r}$ and $\widetilde{\Pi}_{R,\mathbb{Q}_p\{Z_1,...,Z_d\}}^{\infty}$, respectively. Then taking the union throughout all $r>0$ one can define the corresponding Robba ring $\widetilde{\Pi}_{R,\mathbb{Q}_p\{Z_1,...,Z_d\}}$. Again as in \cite{3KL16} one can define the corresponding integral rings of the similar types.
\end{definition}

\indent Then one can define the corresponding period rings in our context deformed over some affinoid $A$ which is isomorphic to some quotient of $\mathbb{Q}_p\{Z_1,...Z_d\}$, again deforming from the context of \cite[Definition 4.1.1]{3KL16}.

\begin{definition}
In the characteristic $0$, we define the following period rings:
\begin{displaymath}
\widetilde{\Omega}^\mathrm{int}_{R,A},\widetilde{\Omega}_{R,A},\widetilde{\Pi}^\mathrm{int,r}_{R,A},\widetilde{\Pi}^\mathrm{bd,r}_{R,A}, \widetilde{\Pi}^I_{R,A},\widetilde{\Pi}^r_{R,A},\widetilde{\Pi}^\infty_{R,A}	
\end{displaymath}
by taking the suitable quotient of the following period rings defined above: 
\begin{align}
\widetilde{\Omega}^\mathrm{int}_{R,\mathbb{Q}_p\{Z_1,...,Z_d\}},&\widetilde{\Omega}_{R,\mathbb{Q}_p\{Z_1,...,Z_d\}},\widetilde{\Pi}^\mathrm{int,r}_{R,\mathbb{Q}_p\{Z_1,...,Z_d\}},\widetilde{\Pi}^\mathrm{bd,r}_{R,\mathbb{Q}_p\{Z_1,...,Z_d\}},\widetilde{\Pi}^I_{R,\mathbb{Q}_p\{Z_1,...,Z_d\}},\\
 &\widetilde{\Pi}^r_{R,\mathbb{Q}_p\{Z_1,...,Z_d\}},\widetilde{\Pi}^\infty_{R,\mathbb{Q}_p\{Z_1,...,Z_d\}}	
\end{align}
with respect to the structure of the affinoid algebra $A$ in the sense of Tate. Therefore they carry the corresponding quotient seminorms of the above Gauss norms defined in the previous definition, note that these are not something induced from the corresponding spectral seminorms from $A$ chosen at the very beginning of our study. We use the notation $\overline{\left\|.\right\|}_{\alpha^r,\mathbb{Q}_p\{Z_1,...,Z_d\}}$ to denote the corresponding quotient Gauss norm which induces then the corresponding spectral seminorm $\left\|.\right\|_{\alpha^r,A}$ for each $r>0$. Then we define the corresponding period rings:
\begin{align}
\widetilde{\Pi}^\mathrm{int}_{R,A},\widetilde{\Pi}^\mathrm{bd}_{R,A},\widetilde{\Pi}_{R,A}	
\end{align}
by taking suitable union throughout all $r>0$.
\end{definition}

\begin{definition}
In positive characteristic situation, when we are working over general affinoid algebra $A$, we use the same notations as in the previous definition, but by using $W(R)_{\mathbb{F}_p[[\eta]]\{Z_1,...,Z_d\}}$ and $W(R)_{A}$ as the starting rings, namely here $A$ is isomorphic to a quotient of $\mathbb{F}_p((\eta))\{Z_1,...,Z_d\}$. Again the corresponding Tate algebra is defined in term of free variables.	
\end{definition}

\indent Then we can define the following affinoid deformations (after \cite[Definition 4.1.1]{3KL16} in the flavor as above):

\begin{definition}
Consider now the base change $W_{\pi,\infty,\mathbb{Q}_p\{Z_1,...,Z_d\}}(R)$ which is defined now to be the completed tensor product $(W_{\pi}(R)\widehat{\otimes}_{\mathcal{O}_E}\mathcal{O}_{E_\infty})\widehat{\otimes}_{\mathbb{Q}_p}\mathbb{Q}_p\{Z_1,...,Z_d\}$. Then the point is that each element in this ring admits a unique expression taking the form of 
\begin{center}
$\sum_{n\in \mathbb{Z}[1/p]_{\geq 0},i_1\geq 0,...,i_d\geq 0}\pi^n[\overline{x}_{n,i_1,...,i_d}]Z_1^{i_1}...Z_d^{i_d}$ 	
\end{center}
which allows us to perform the construction mentioned above. First we can define for some radius $r>0$ the corresponding period ring $\widetilde{\Pi}^{\mathrm{int},r}_{R,\infty,\mathbb{Q}_p\{Z_1,...,Z_d\}}$ by taking the completion of the ring
\begin{displaymath}
W_{\pi,\infty,\mathbb{Q}_p\{Z_1,...,Z_d\}}(R^+)[[R]]	
\end{displaymath}
with respect to the following Gauss type norm:
\begin{displaymath}
\left\|.\right\|_{\alpha^r}(\sum_{n\in \mathbb{Z}[1/p]_{\geq 0},i_1\geq 0,...,i_d\geq 0}\pi^n[\overline{x}_{n,i_1,...,i_d}]Z_1^{i_1}...Z_d^{i_d}):=\sup_{n\in \mathbb{Z}[1/p]_{\geq 0},i_1\geq 0,...,i_d\geq 0}\{p^{-n}\alpha(\overline{x}_{n,i_1,...,i_d})^r\}.	
\end{displaymath}
Then we can define the union $\widetilde{\Pi}^{\mathrm{int}}_{R,\infty,\mathbb{Q}_p\{Z_1,...,Z_d\}}$ throughout all the radius $r>0$. Then we just define the bounded Robba ring $\widetilde{\Pi}^{\mathrm{bd},r}_{R,\infty,\mathbb{Q}_p\{Z_1,...,Z_d\}}$ by $\widetilde{\Pi}^{\mathrm{int},r}_{R,\infty,\mathbb{Q}_p\{Z_1,...,Z_d\}}[1/\pi]$ and also we could define the union $\widetilde{\Pi}^{\mathrm{bd}}_{R,\infty,\mathbb{Q}_p\{Z_1,...,Z_d\}}$ throughout all the radius $r>0$. Then for any interval in $(0,\infty)$ which is denoted by $I$ we can define the corresponding Robba rings $\widetilde{\Pi}^{I}_{R,\infty,\mathbb{Q}_p\{Z_1,...,Z_d\}}$ by taking the Fr\'echet completion of 
\begin{displaymath}
W_{\pi,\infty,\mathbb{Q}_p\{Z_1,...,Z_d\}}(R^+)[[R]][1/\pi]	
\end{displaymath}
with respect to all the norms $\left\|.\right\|_{\alpha^t}$ for all $t\in I$. Then by taking suitable specified intervals one can define the rings $\widetilde{\Pi}^{r}_{R,\infty,\mathbb{Q}_p\{Z_1,...,Z_d\}}$ and $\widetilde{\Pi}^{\infty}_{R,\infty,\mathbb{Q}_p\{Z_1,...,Z_d\}}$ as before, and finally one can define the corresponding union $\widetilde{\Pi}_{R,\infty,\mathbb{Q}_p\{Z_1,...,Z_d\}}$ throughout all the radius $r>0$. Again we have the corresponding integral version of the rings defined over $E_\infty$ as \cite{3KL16}. Finally over $A$ we can define in the same way as above to deform all the rings over $E_\infty$, and we do not repeat the construction again.

\end{definition}

\subsection{Basic Properties of Period Rings}	
	
\indent Then we do some reality checks over the investigation of the properties of the above period rings in the style taken in \cite[Section 5.2]{3KL15} and \cite{3KL16}. 
%We will focus on the situation where $A$ is just the ring $\mathbb{Q}_p\{T_1,...,T_d\}$ for some $d\geq 1$.

\begin{proposition} \mbox{\bf{(After Kedlaya-Liu \cite[Lemma 5.2.1]{3KL15})}}
The function $t\mapsto \left\|x\right\|_{\alpha^t,\mathbb{Q}_p\{Z_1,...,Z_d\}}$ for $x\in \widetilde{\Pi}^r_{R,\mathbb{Q}_p\{Z_1,...,Z_d\}}$ is continuous log convex for the corresponding variable $t\in (0,r]$ for any $r>0$.
\end{proposition}

\begin{proof}
Adapt the corresponding argument in the proof of 5.2.1 of \cite[Lemma 5.2.1]{3KL15} to our situation we then first look at the situation where the element is just of the form of $\pi^k[\overline{x}_{k,i_1,...,i_d}]Z_1^{i_1}...Z_d^{i_d}$ where the corresponding norm function in terms of $t>0$ is just affine. Then one focuses on the finite sums of these kind of elements which gives rise to to the log convex directly. Finally by taking the approximation we get the desired result.	
\end{proof}

\begin{proposition}\mbox{\bf{(After Kedlaya-Liu \cite[Lemma 5.2.2]{3KL15})}}
For any element $x\in \widetilde{\Pi}_{R,\mathbb{Q}_p\{Z_1,...,Z_d\}}$ we have that $x\in \widetilde{\Pi}^\mathrm{bd}_{R,\mathbb{Q}_p\{Z_1,...,Z_d\}}$ if and only if we have the situation where $x$ actually lives in $\widetilde{\Pi}^r_{R,\mathbb{Q}_p\{Z_1,...,Z_d\}}$ (for some specific $r>0$) such that $x$ itself is bounded under the norm $\left\|.\right\|_{\alpha^t,\mathbb{Q}_p\{Z_1,...,Z_d\}}$ for each $t\in (0,r]$.	
\end{proposition}

\begin{proof}
One direction of the proof is easy, so we only choose to present the proof of the implication in the other direction as in the original proof of 5.2.2 of \cite[Lemma 5.2.2]{3KL15} as in the following. First choose some radius $r>0$ such that the element could be assumed to be living in the ring $\widetilde{\Pi}^r_{R,\mathbb{Q}_p\{Z_1,...,Z_d\}}$. The idea is to transfer the original question to the question about showing the integrality of $x$ when we add some hypothesis on the norm by taking suitable powers of $p$ (since the norm is bounded for each $t\in (0,r]$ so we are reduced to the situation where the norm is bounded by $1$). Then we argue as in \cite[Lemma 5.2.2]{3KL15} to choose some approximating sequence $\{x_i\}$ living in $\widetilde{\Pi}^\mathrm{bd}_{R,\mathbb{Q}_p\{Z_1,...,Z_d\}}$ of $x$. Therefore we have for any $j\geq 1$ one can find then some integer $N_j\geq 1$ such that for any $i\geq N_j$ we have the estimate:
\begin{displaymath}
\left\|.\right\|_{\alpha^t,\mathbb{Q}_p\{Z_1,...,Z_d\}}(x_i-x)\leq p^{-j}, \forall t\in [p^{-j}r,r].	
\end{displaymath}
Then the idea is to consider the integral decomposition of the element $x_i$ which has the form of $\sum_{k=n(x_i),i_1\geq 0,i_2\geq 0,...,i_d\geq 0}\pi^k[\overline{x}_{i,k,i_1,...,i_d}]$ into the following two parts:
\begin{align}
x_i &:= y_i+z_i\\
	&:=\sum_{k=0,i_1\geq 0,i_2\geq 0,...,i_d\geq 0}\pi^k[\overline{x}_{i,k,i_1,...,i_d}]Z_1^{i_1}...Z_d^{i_d}+z_i
\end{align}
from which we actually have the corresponding estimate over the residual part of the decomposition above:
\begin{displaymath}
\left\|.\right\|_{\alpha^{p^{-j}r},\mathbb{Q}_p\{Z_1,...,Z_d\}}(\pi^k[\overline{x}_{i,k,i_1,...,i_d}]Z_1^{i_1}...Z_d^{i_d})\leq 1, \forall k<0	
\end{displaymath}
which implies by direct computation:
\begin{align}
\alpha^{p^{-j}r}(\overline{x}_{i,k,i_1,...,i_d})&\leq p^k\\
\alpha(\overline{x}_{i,k,i_1,...,i_d})&\leq p^{kp^{j}/r}	
\end{align}
which implies that we have the following estimate: 
\begin{align}
\left\|.\right\|_{\alpha^{r},\mathbb{Q}_p\{Z_1,...,Z_d\}}(x_i-y_i)&\leq p^{-k}p^{kp^{j}r/r}\\
&\leq p^{1-p^j}
\end{align}
which implies that $y_i\rightarrow x$ under the norm $\left\|.\right\|_{\alpha^{r},\mathbb{Q}_p\{Z_1,...,Z_d\}}$ which furthermore under the norm $\left\|.\right\|_{\alpha^{r},\mathbb{Q}_p\{Z_1,...,Z_d\}}$ due the property of the norm, which finishes the proof of the desired result.
\end{proof}

\begin{proposition}\mbox{\bf{(After Kedlaya-Liu \cite[Corollary 5.2.3]{3KL15})}} We have the following identity:
\begin{displaymath}
(\widetilde{\Pi}_{R,\mathbb{Q}_p\{Z_1,...,Z_d\}})^\times=(\widetilde{\Pi}^\mathrm{bd}_{R,\mathbb{Q}_p\{Z_1,...,Z_d\}})^\times.
\end{displaymath}
\end{proposition}

\begin{proof}
See 5.2.3 of \cite[Corollary 5.2.3]{3KL15}.	
\end{proof}

\begin{proposition}\mbox{\bf{(After Kedlaya-Liu \cite[Lemma 5.2.6]{3KL15})}}
For any $0< r_1\leq r_2$ we have the following equality on the corresponding period rings:
\begin{displaymath}
\widetilde{\Pi}^{\mathrm{int},r_1}_{R,\mathbb{Q}_p\{Z_1,...,Z_d\}}\bigcap	\widetilde{\Pi}^{[r_1,r_2]}_{R,\mathbb{Q}_p\{Z_1,...,Z_d\}}=\widetilde{\Pi}^{\mathrm{int},r_2}_{R,\mathbb{Q}_p\{Z_1,...,Z_d\}}.
\end{displaymath}	
\end{proposition}

\begin{proof}
We adapt the argument in \cite[Lemma 5.2.6]{3KL15} 5.2 to prove this in the situation where $r_1<r_2$ (otherwise this is trivial), again one direction is easy where we only present the implication in the other direction. We take any element
\begin{center}
 $x\in \widetilde{\Pi}^{\mathrm{int},r_1}_{R,\mathbb{Q}_p\{Z_1,...,Z_d\}}\bigcap	\widetilde{\Pi}^{[r_1,r_2]}_{R,\mathbb{Q}_p\{Z_1,...,Z_d\}}$
\end{center} 
and take suitable approximating elements $\{x_i\}$ living in the bounded Robba ring such that for any $j\geq 1$ one can find some integer $N_j\geq 1$ we have whenever $i\geq N_j$ we have the following estimate:
\begin{displaymath}
\left\|.\right\|_{\alpha^{t},\mathbb{Q}_p\{Z_1,...,Z_d\}}(x_i-x) \leq p^{-j}, \forall t\in [r_1,r_2].	
\end{displaymath}
Then we consider the corresponding decomposition of $x_i$ for each $i=1,2,...$ into a form having integral part and the rational part $x_i=y_i+z_i$ by setting
\begin{center}
 $y_i=\sum_{k=0,i_1,...,i_d}\pi^k[\overline{x}_{i,k,i_1,...,i_d}]Z_1^{i_1}...Z_d^{i_d}$ 
\end{center} 
out of
\begin{center} 
$x_i=\sum_{k=n(x_i),i_1,...,i_d}\pi^k[\overline{x}_{i,k,i_1,...,i_d}]Z_1^{i_1}...Z_d^{i_d}$.
\end{center}
Note that by our initial hypothesis we have that the element $x$ lives in the ring $\widetilde{\Pi}^{\mathrm{int},r_1}_{R,\mathbb{Q}_p\{Z_1,...,Z_d\}}$ which further implies that 
\begin{displaymath}
\left\|.\right\|_{\alpha^{r_1},\mathbb{Q}_p\{Z_1,...,Z_d\}}(\pi^k[\overline{x}_{i,k,i_1,...,i_d}]Z_1^{i_1}...Z^{i_d}_d)	\leq p^{-j}.
\end{displaymath}
Therefore we have $\alpha(\overline{x}_{i,k,i_1,...,i_d})\leq p^{(k-j)/r_1}$ directly from this through computation, which implies that then:
\begin{align}
\left\|.\right\|_{\alpha^{r_2},\mathbb{Q}_p\{Z_1,...,Z_d\}}(\pi^k[\overline{x}_{i,k,i_1,...,i_d}]Z_1^{i_1}...Z^{i_d}_d)	&\leq p^{-k}p^{(k-j)r_2/r_1}\\
	&\leq p^{1+(1-j)r_1/r_1}.
\end{align}
Then one can read off the result directly from this estimate since under this estimate we can have the chance to modify the original approximating sequence $\{x_i\}$ by $\{y_i\}$ which are initially chosen to be in the integral Robba ring, which implies that actually the element $x$ lives in the right-hand side of the identity in the statement of the proposition.
\end{proof}

\begin{proposition} \mbox{\bf{(After Kedlaya-Liu \cite[Lemma 5.2.6]{3KL15})}}
For any $0< r_1\leq r_2$ we have the following equality on the corresponding period rings:
\begin{displaymath}
\widetilde{\Pi}^{\mathrm{int},r_1}_{R,A}\bigcap	\widetilde{\Pi}^{[r_1,r_2]}_{R,A}=\widetilde{\Pi}^{\mathrm{int},r_2}_{R,A}.
\end{displaymath}	
Here $A$ is some noncommutative Banach affinoid algebra over $\mathbb{Q}_p$.	
\end{proposition}

\begin{proof}
See the proof of \cref{proposition5.7}.	
\end{proof}

\indent Then we have the following analog of the corresponding result of \cite[Lemma 5.2.8]{3KL15}:

\begin{proposition} \mbox{\bf{(After Kedlaya-Liu \cite[Lemma 5.2.8]{3KL15})}}
Consider now in our situation the radii $0< r_1\leq r_2$, and consider any element $x\in \widetilde{\Pi}^{[r_1,r_2]}_{R,\mathbb{Q}_p\{Z_1,...,Z_d\}}$. Then we have that for each $n\geq 1$ one can decompose $x$ into the form of $x=y+z$ such that $y\in \pi^n\widetilde{\Pi}^{\mathrm{int},r_2}_{R,\mathbb{Q}_p\{Z_1,...,Z_d\}}$ with $z\in \bigcap_{r\geq r_2}\widetilde{\Pi}^{[r_1,r]}_{R,\mathbb{Q}_p\{Z_1,...,Z_d\}}$ with the following estimate for each $r\geq r_2$:
\begin{displaymath}
\left\|.\right\|_{\alpha^r,\mathbb{Q}_p\{Z_1,...,Z_d\}}(z)\leq p^{(1-n)(1-r/r_2)}\left\|.\right\|_{\alpha^{r_2},\mathbb{Q}_p\{Z_1,...,Z_d\}}(z)^{r/r_2}.	
\end{displaymath}

\end{proposition}

\begin{proof}
As in \cite[Lemma 5.2.8]{3KL15} and in the proof of our previous proposition we first consider those elements $x$ lives in the bounded Robba rings which could be expressed in general as
\begin{center}
 $\sum_{k=n(x),i_1,...,i_d}\pi^k[\overline{x}_{k,i_1,...,i_d}]Z_1^{i_1}...Z_d^{i_d}$.
 \end{center}	
In this situation the corresponding decomposition is very easy to come up with, namely we consider the corresponding $y_i$ as the corresponding series:
\begin{displaymath}
\sum_{k\geq n,i_1,...,i_d}\pi^k[\overline{x}_{k,i_1,...,i_d}]Z_1^{i_1}...Z_d^{i_d}	
\end{displaymath}
which give us the desired result since we have in this situation when focusing on each single term:
\begin{align}
\left\|.\right\|_{\alpha^r,\mathbb{Q}_p\{Z_1,...,Z_d\}}&(\pi^k[\overline{x}_{k,i_1,...,i_d}]Z_1^{i_1}...Z_d^{i_d})=p^{-k}\alpha(\overline{x}_{k,i_1,...,i_d})^r\\
&=p^{-k(1-r/r_2)}\left\|.\right\|_{\alpha^{r_2},\mathbb{Q}_p\{Z_1,...,Z_d\}}(\pi^k[\overline{x}_{k,i_1,...,i_d}]Z_1^{i_1}...Z_d^{i_d})^{r/r_2}\\
&\leq p^{(1-n)(1-r/r_2)}\left\|.\right\|_{\alpha^{r_2},\mathbb{Q}_p\{Z_1,...,Z_d\}}(\pi^k[\overline{x}_{k,i_1,...,i_d}]Z_1^{i_1}...Z_d^{i_d})^{r/r_2}
\end{align}
for all those suitable $k$. Then to tackle the more general situation we consider the approximating sequence consisting of all the elements in the bounded Robba ring as in the usual situation considered in \cite[Lemma 5.2.8]{3KL15}, namely we inductively construct the following approximating sequence just as:
\begin{align}
\left\|.\right\|_{\alpha^r,\mathbb{Q}_p\{Z_1,...,Z_d\}}(x-x_0-...-x_i)\leq p^{-i-1}	\left\|.\right\|_{\alpha^r,\mathbb{Q}_p\{Z_1,...,Z_d\}}(x), i=0,1,..., r\in [r_1,r_2].
\end{align}
Here all the elements $x_i$ for each $i=0,1,...$ are living in the bounded Robba ring, which immediately gives rise to the suitable decomposition as proved in the previous case namely we have for each $i$ the decomposition $x_i=y_i+z_i$ with the desired conditions as mentioned in the statement of the proposition. We first take the series summing all the elements $y_i$ up for all $i=0,1,...$, which first of all converges under the norm $\left\|.\right\|_{\alpha^r,\mathbb{Q}_p\{Z_1,...,Z_d\}}$ for all the radius $r\in [r_1,r_2]$, and also note that all the elements $y_i$ within the infinite sum live inside the corresponding integral Robba ring $\widetilde{\Pi}^{\mathrm{int},r_2}_{R,\mathbb{Q}_p\{Z_1,...,Z_d\}}$, which further implies the corresponding convergence ends up in $\widetilde{\Pi}^{\mathrm{int},r_2}_{R,\mathbb{Q}_p\{Z_1,...,Z_d\}}$. For the elements $z_i$ where $i=0,1,...$ also sum up to a converging series in the desired ring since combining all the estimates above we have:
\begin{displaymath}
\left\|.\right\|_{\alpha^r,\mathbb{Q}_p\{Z_1,...,Z_d\}}(z_i)\leq p^{(1-n)(1-r/r_2)}\left\|.\right\|_{\alpha^{r_2},\mathbb{Q}_p\{Z_1,...,Z_d\}}(x)^{r/r_2}.	
\end{displaymath}
\end{proof}

%\begin{corollary}
%Consider now in our situation the radii $0< r_1\leq r_2$, and consider any element $x\in \widetilde{\Pi}^{[r_1,r_2]}_{R,A}$. Then we have that for each $n\geq 1$ one can decompose $x$ into the form of $x=y+z$ such that $y\in \pi^n\widetilde{\Pi}^{\mathrm{int},r_2}_{R,A}$ with $z\in \bigcap_{r\geq r_2}\widetilde{\Pi}^{[r_1,r]}_{R,A}$ with the following estimate for each $r\geq r_2$:
%\begin{displaymath}
%\left\|.\right\|_{\alpha^r,A}(z)\leq p^{(1-n)(1-r/r_2)}\left\|.\right\|_{\alpha^{r_2},A}(z)^{r/r_2}.	
%\end{displaymath}	
%\end{corollary}
%
%\begin{proof}
%This could be derived directly from the previous proposition, by considering the corresponding residual norms and then the spectral norms.	
%\end{proof}

\begin{proposition} \mbox{\bf{(After Kedlaya-Liu \cite[Lemma 5.2.10]{3KL15})}}
We have the following identity:
\begin{displaymath}
\widetilde{\Pi}^{[s_1,r_1]}_{R,\mathbb{Q}_p\{Z_1,...,Z_d\}}\bigcap\widetilde{\Pi}^{[s_2,r_2]}_{R,\mathbb{Q}_p\{Z_1,...,Z_d\}}=\widetilde{\Pi}^{[s_1,r_2]}_{R,\mathbb{Q}_p\{Z_1,...,Z_d\}},
\end{displaymath}
here the radii satisfy $<s_1\leq s_2 \leq r_1 \leq r_2$.
\end{proposition}

\begin{proof}
In our situation one direction is obvious while on the other hand we consider any element $x$ in the intersection on the left, then by the previous proposition we	have the decomposition $x=y+z$ where $y\in \widetilde{\Pi}^{\mathrm{int},r_1}_{R,\mathbb{Q}_p\{Z_1,...,Z_d\}}$ and $z\in \widetilde{\Pi}^{[s_1,r_2]}_{R,\mathbb{Q}_p\{Z_1,...,Z_d\}}$. Then as in \cite[Lemma 5.2.10]{3KL15} section 5.2 we look at $y=x-z$ which lives in the intersection:
\begin{displaymath}
\widetilde{\Pi}^{\mathrm{int},r_1}_{R,\mathbb{Q}_p\{Z_1,...,Z_d\}}\bigcap	\widetilde{\Pi}^{[s_2,r_2]}_{R,\mathbb{Q}_p\{Z_1,...,Z_d\}}=\widetilde{\Pi}^{\mathrm{int},r_2}_{R,\mathbb{Q}_p\{Z_1,...,Z_d\}}
\end{displaymath}
which finishes the proof.
\end{proof}

\begin{proposition}\mbox{\bf{(After Kedlaya-Liu \cite[Lemma 5.2.10]{3KL15})}}
We have the following identity:
\begin{displaymath}
\widetilde{\Pi}^{[s_1,r_1]}_{R,A}\bigcap\widetilde{\Pi}^{[s_2,r_2]}_{R,A}=\widetilde{\Pi}^{[s_1,r_2]}_{R,A},
\end{displaymath}
here the radii satisfy $<s_1\leq s_2 \leq r_1 \leq r_2$.
	
\end{proposition}

\begin{proof}
See the proof of \cref{proposition5.7}.		
\end{proof}

\begin{remark}
Again this subsection is finished so far only for the situation where $E$ is of mixed characteristic. But everything uniformly carries over for our original assumption on the field $E$ and $A$. We will not repeat the proof again.	
\end{remark}

\subsection{Noncommutative Period Rings and Period Sheaves}

\begin{setting}
We will work in the categories of the pseudocoherent, fpd and finite projective modules over the period rings defined above. First we specify the Frobenius in our setting. The rings involved are:
\begin{align}
\widetilde{\Omega}^\mathrm{int}_{R,A},\widetilde{\Omega}_{R,A}, \widetilde{\Pi}^\mathrm{int}_{R,A}, \widetilde{\Pi}^{\mathrm{int},r}_{R,A},\widetilde{\Pi}^\mathrm{bd}_{R,A},\widetilde{\Pi}^{\mathrm{bd},r}_{R,A}, \widetilde{\Pi}_{R,A}, \widetilde{\Pi}^r_{R,A},\widetilde{\Pi}^+_{R,A}, \widetilde{\Pi}^\infty_{R,A},\widetilde{\Pi}^I_{R,A}.
\end{align}
We are going to endow these rings with the Frobenius induced by continuation from the Witt vector part only, which is to say the corresponding Frobenius induced by the $p^h$-power absolute Frobenius over $R$. Note all the rings above are defined by taking the product of $\triangle$ where each $\triangle$ representing one of the following rings (over $E$):
\begin{align}
\widetilde{\Omega}^\mathrm{int}_{R},\widetilde{\Omega}_{R}, \widetilde{\Pi}^\mathrm{int}_{R}, \widetilde{\Pi}^{\mathrm{int},r}_{R},\widetilde{\Pi}^\mathrm{bd}_{R},\widetilde{\Pi}^{\mathrm{bd},r}_{R}, \widetilde{\Pi}_{R}, \widetilde{\Pi}^r_{R},\widetilde{\Pi}^+_{R}, \widetilde{\Pi}^\infty_{R},\widetilde{\Pi}^I_{R}
\end{align}
with the affinoid ring $A$. The Frobenius acts on $A$ trivially and we assume that the action is $A$-linear.  
\end{setting}

\indent First we consider the following sheafification as in \cite[Definition 4.4.2]{3KL16}:

\begin{setting} 
Consider the space $X=\mathrm{Spa}(R,R^+)$, over this perfectoid space there were sheaves:
\begin{align}
\widetilde{\Omega}^\mathrm{int}_{},\widetilde{\Omega}_{}, \widetilde{\Pi}^\mathrm{int}_{}, \widetilde{\Pi}^{\mathrm{int},r}_{},\widetilde{\Pi}^\mathrm{bd}_{},\widetilde{\Pi}^{\mathrm{bd},r}_{}, \widetilde{\Pi}_{}, \widetilde{\Pi}^r_{},\widetilde{\Pi}^+_{}, \widetilde{\Pi}^\infty_{},\widetilde{\Pi}^I_{}.
\end{align}
defined over this space through the corresponding adic, \'etale, pro-\'etale or $v$-topology, we consider the corresponding sheaves defined over the same Grothendieck sites but with further deformed consideration:
\begin{align}
\widetilde{\Omega}^\mathrm{int}_{*,A},\widetilde{\Omega}_{*,A}, \widetilde{\Pi}^\mathrm{int}_{*,A}, \widetilde{\Pi}^{\mathrm{int},r}_{*,A},\widetilde{\Pi}^\mathrm{bd}_{*,A},\widetilde{\Pi}^{\mathrm{bd},r}_{*,A}, \widetilde{\Pi}_{*,A}, \widetilde{\Pi}^r_{*,A},\widetilde{\Pi}^+_{*,A}, \widetilde{\Pi}^\infty_{*,A},\widetilde{\Pi}^I_{*,A}.
\end{align}
\end{setting}

\begin{remark}
This is well defined, since we have orthogonal basis for the ring $A$.	
\end{remark}

\begin{definition}
In our situation we consider the corresponding Frobenius action on the following period rings and sheaves:
\begin{align}
\widetilde{\Omega}^\mathrm{int}_{R,A},\widetilde{\Omega}^\mathrm{int}_{R,A},\widetilde{\Pi}^\mathrm{int}_{R,A},\widetilde{\Pi}^\mathrm{bd}_{R,A},\widetilde{\Pi}_{R,A},\widetilde{\Pi}^\infty_{R,A}, \widetilde{\Omega}^\mathrm{int}_{*,A}, \widetilde{\Pi}^\mathrm{int}_{*,A}, \widetilde{\Pi}^+_{*,A}, \widetilde{\Omega}_{*,A}, \widetilde{\Pi}^\mathrm{bd}_{*,A}, \widetilde{\Pi}^\infty_{*,A}, \widetilde{\Pi}^+_{*,A}, \widetilde{\Pi}_{*,A},   
\end{align}
which is defined by considering the corresponding lift of the absolute Frobenius of $p^h$-power in characteristic $p>0$ induced from $R$, which will be denoted by $\varphi$. We then introduce more general consideration by taking $E_a$ to be some unramified extension of $E$ of degree $a$ divisible by $h$, the corresponding Frobenius will be denoted by $\varphi^a$.
\end{definition}

\indent Then we generalize the corresponding Frobenius modules in \cite[Definition 4.4.4]{3KL16} to our situation as in the following.

\begin{definition}
Over the period rings and sheaves (each is denoted by $\triangle$ in this definition) defined in the previous definition we define as in \cite[Definition 4.4.4]{3KL16} the corresponding $\varphi^a$-modules over $\triangle$ which are respectively projective to be the corresponding finite projetive modules over $\triangle$ with further assigned semilinear action of the operator $\varphi^a$. Here we define in our situation the corresponding $\varphi^a$-cohomology to be the (hyper)-cohomology of the following complex:
\[
\xymatrix@R+0pc@C+0pc{
0\ar[r]\ar[r]\ar[r] &M
\ar[r]^{\varphi-1}\ar[r]\ar[r] &M
\ar[r]\ar[r]\ar[r] &0.
}
\]
All modules are right over the noncommutative rings. 
\end{definition}

\indent Now we define the corresponding modules over the rings which are the domains in the following morphisms induced from the Frobenius map $\varphi^a$:
\begin{align} 
\widetilde{\Pi}^{\mathrm{int},r}_{R,A}\rightarrow \widetilde{\Pi}^{\mathrm{int},rp^{-ha}}_{R,A},\widetilde{\Pi}^{\mathrm{bd},r}_{R,A}\rightarrow \widetilde{\Pi}^{\mathrm{bd},rp^{-ha}}_{R,A},\widetilde{\Pi}^{r}_{R,A}\rightarrow \widetilde{\Pi}^{rp^{-ha}}_{R,A}\\
\widetilde{\Pi}^{\mathrm{int},r}_{*,A}\rightarrow \widetilde{\Pi}^{\mathrm{int},rp^{-ha}}_{*,A},\widetilde{\Pi}^{\mathrm{bd},r}_{*,A}\rightarrow \widetilde{\Pi}^{\mathrm{bd},rp^{-ha}}_{*,A},\widetilde{\Pi}^{r}_{*,A}\rightarrow \widetilde{\Pi}^{rp^{-ha}}_{*,A}.	
\end{align}

\begin{definition}
Over each rings $\triangle$ which are the domains in the morphisms as mentioned just before this definition, we define the corresponding projective $\varphi^a$-module over any $\triangle$ listed above to be the corresponding finite projective module $M$ over $\triangle$ with additionally endowed semilinear Frobenius action from $\varphi^a$ such that we have the isomorphism $\varphi^{a*}M\overset{\sim}{\rightarrow}M\otimes \square$ where the ring $\square$ is one of the targets listed above. Also the cohomology of any module under this definition will be defined to be the (hyper)cohomology of the complex in the following form:
\[
\xymatrix@R+0pc@C+0pc{
0\ar[r]\ar[r]\ar[r] &M
\ar[r]^{\varphi-1}\ar[r]\ar[r] &M\otimes_\triangle \square
\ar[r]\ar[r]\ar[r] &0.
}
\]
All modules are right over the noncommutative rings.
\end{definition}

\indent Then we consider the following morphisms of specific period rings induced by the Frobenius. 
\begin{align}
\widetilde{\Pi}_{R,A}^{[s,r]}	\rightarrow \widetilde{\Pi}_{R,A}^{[sp^{-ah},rp^{-ah}]}\\
\widetilde{\Pi}_{*,A}^{[s,r]}	\rightarrow \widetilde{\Pi}_{*,A}^{[sp^{-ah},rp^{-ah}]}
\end{align}

with the corresponding morphisms in the following:

\begin{align}
\widetilde{\Pi}_{R,A}^{[s,r]}	\rightarrow \widetilde{\Pi}_{R,A}^{[s,rp^{-ah}]}\\
\widetilde{\Pi}_{*,A}^{[s,r]}	\rightarrow \widetilde{\Pi}_{*,A}^{[s,rp^{-ah}]}
\end{align}

\begin{definition}
Again as in \cite[Definition 4.4.4]{3KL16}, we define the corresponding projective $\varphi^a$-modules over the domain rings or sheaves of rings in the morphisms just before this definition to be the finite projective modules (which will be denoted by $M$) over the domain rings in the morphism just before this definition additionally endowed with semilinear Frobenius action from $\varphi^a$ with the following isomorphisms:
\begin{align}
\varphi^{a*}M\otimes_{\widetilde{\Pi}_{R,A}^{[sp^{-ah},rp^{-ah}]}}\widetilde{\Pi}_{R,A}^{[s,rp^{-ah}]}\overset{\sim}{\rightarrow}M\otimes_{\widetilde{\Pi}_{R,A}^{[s,r]}}\widetilde{\Pi}_{R,A}^{[s,rp^{-ah}]},\\
\varphi^{a*}M\otimes_{\widetilde{\Pi}_{*,A}^{[sp^{-ah},rp^{-ah}]}}\widetilde{\Pi}_{*,A}^{[s,rp^{-ah}]}\overset{\sim}{\rightarrow}M\otimes_{\widetilde{\Pi}_{*,A}^{[s,r]}}\widetilde{\Pi}_{*,A}^{[s,rp^{-ah}]}.
\end{align}
All modules are right over the noncommutative rings.
\end{definition}

\indent Also one can further define the corresponding bundles carrying semilinear Frobenius in our context as in the situation of \cite[Definition 4.4.4]{3KL16}:

\begin{definition}
Over the ring $\widetilde{\Pi}_{R,A}$ we define a corresponding projective Frobenius bundle to be a family $(M_I)_I$ of finite projective modules over each $\widetilde{\Pi}^I_{R,A}$ carrying the natural Frobenius action coming from the operator $\varphi^a$ such that for any two involved intervals having the relation $I\subset J$ we have:
\begin{displaymath}
M_J\otimes_{\widetilde{\Pi}^J_{R,A}}\widetilde{\Pi}^I_{R,A}\overset{\sim}{\rightarrow}	M_I
\end{displaymath}
with the obvious cocycle condition. Here we have to propose condition on the intervals that for each $I=[s,r]$ involved we have $s\leq rp^{ah}$. All modules are right over the noncommutative rings.
\end{definition}

\indent Then we generalize the comparison of local systems and Frobenius modules in \cite[Theorem 4.5.7]{3KL16} to our noncommutative context as in the following:

\begin{theorem} \mbox{\bf{(After Kedlaya-Liu \cite[Theorem 4.5.7]{3KL16})}}
1. There is a fully faithful embedding functor from the category of the projective right $\mathcal{O}_{E_a^{\varphi^a}}\widehat{\otimes}\mathcal{O}_A$-local systems over $X_\text{pro\'et}$ to the following two categories:\\
1(a).  The category of projective right Frobenius modules over the period sheaf $\widetilde{\Omega}_{*,A}^\mathrm{int}$ over $X$, $X_\text{\'et}$, $X_\text{pro\'et}$;\\
1(b).  The category of projective right Frobenius modules over the period sheaf $\widetilde{\Pi}_{*,A}^\mathrm{int}$ over $X$, $X_\text{\'et}$, $X_\text{pro\'et}$.
\end{theorem}

\begin{proof}
See \cref{theorem4.4}.	
\end{proof}

\indent Then as in \cite[Corollary 4.5.8]{3KL16} we consider the setting of more general adic spaces:

\begin{proposition} \mbox{\bf{(After Kedlaya-Liu \cite[Corollary 4.5.8]{3KL16})}}
1. Let $X$	be a preadic space over $E_a$. Then we have that there is fully faithful embedding functor from the category of all the right $\mathcal{O}_{E_{a}^{\varphi^a}}\widehat{\otimes}\mathcal{O}_A$-local systems (projective) over $X$, $X_{\text{\'et}}$ and $X_\text{pro\'et}$ to the category of all the projective right Frobenius modules over the sheaves over $\widetilde{\Pi}_{R,A}^\mathrm{int}$ over the corresponding sites;\\
2. Let $X$	be a preadic space over $E_{\infty,a}$. Then we have that there is fully faithful embedding functor from the category of all the right $\mathcal{O}_{E_{\infty,a}^{\varphi^a}}\widehat{\otimes}\mathcal{O}_A$-local systems (projective) over $X$, $X_{\text{\'et}}$ and $X_\text{pro\'et}$ to the category of all the projective right Frobenius modules over the sheaves over $\widetilde{\Pi}_{R,\infty,A}^\mathrm{int}$ over the corresponding sites.
\end{proposition}

\begin{proof}
For one, apply the previous theorem. For two, repeat the argument in the proof of the previous theorem.	
\end{proof}

\indent The following definition is kind of generalization of the corresponding one in \cite[Definition 4.5.9]{3KL16}:

\begin{definition}
Now over the ring $\widetilde{\Pi}_{R,A}$ or $\widetilde{\Pi}^\mathrm{bd}_{R,A}$	we call the corresponding Frobenius modules globally \'etale if they arise from the Frobenius modules over the ring $\widetilde{\Pi}^\mathrm{int}_{R,A}$. Now over the ring $\widetilde{\Pi}_{R,\infty,A}$ or $\widetilde{\Pi}^\mathrm{bd}_{R,\infty,A}$	we call the corresponding Frobenius modules globally \'etale if they arise from the Frobenius modules over the ring $\widetilde{\Pi}^\mathrm{int}_{R,\infty,A}$. \\
Now over the sheaf $\widetilde{\Pi}_{*,A}$ or $\widetilde{\Pi}^\mathrm{bd}_{*,A}$	we call the corresponding Frobenius modules globally \'etale if they arise from the Frobenius modules over the sheaf $\widetilde{\Pi}^\mathrm{int}_{*,A}$. Now over the ring $\widetilde{\Pi}_{*,\infty,A}$ or $\widetilde{\Pi}^\mathrm{bd}_{*,\infty,A}$	we call the corresponding Frobenius modules globally \'etale if they arise from the Frobenius modules over the sheaf $\widetilde{\Pi}^\mathrm{int}_{*,\infty,A}$.

\end{definition}

\begin{proposition} \mbox{\bf{(After Kedlaya-Liu \cite[Theorem 4.5.11]{3KL16})}}
1. Let $X$ be a preadic space over $E_{a}$. Then we have that there is a fully faithful embedding of the category of the corresponding right $E^{\varphi^a}_{a}\widehat{\otimes}A$-local systems (in the projective setting) into the category of the corresponding projective right Frobenius $\varphi^a$-modules over $\widetilde{\Pi}_{*,A}$;\\
2. Let $X$ be a preadic space over $E_{\infty,a}$. Then we have that there is a fully faithful embedding of the category of the corresponding right $E^{\varphi^a}_{\infty,a}\widehat{\otimes}A$-local systems (in the projective setting) into the category of the corresponding projective right  Frobenius $\varphi^a$-modules over $\widetilde{\Pi}_{*,A}$;\\
\end{proposition}

\begin{proof}
As in \cite[Theorem 4.5.11]{3KL16}, we consider the corresponding base change of the corresponding exact sequence in \cref{proposition4.1} which reflects an exact sequence on the sheaves.
\end{proof}

\begin{lemma} \mbox{\bf{(After Kedlaya-Liu \cite[6.2.2-6.2.4]{3KL15}, \cite[Lemma 4.6.9]{3KL16}})}\\ \mbox{\bf{(And also see \cite[Proposition 2.11]{3T1})}}
For any Frobenius $\varphi^a$-bundle $M$ over $\widetilde{\Pi}_{R,A}$, then we have that for any interval $I=[s,r]$ where $0<s\leq r$ the map $\varphi^a-1: M_{[s,rq]}\rightarrow M_{[s,r]}$ is surjective after taking some Frobenius twist as in \cite[6.2.2]{3KL15} and our previous work namely the new morphism $\varphi^a-1: M(n)_{[s,rq]}\rightarrow M(n)_{[s,r]}$ for sufficiently large $n\geq 1$. As in the previous established version one may have the chance to take the number to be $1$ is the bundle initially comes from the corresponding base change from the integral Robba ring. 
\end{lemma}

\begin{proof}
See the proof of \cite[Proposition 2.11]{3T1}.

\end{proof}

%\begin{lemma} \mbox{\bf{(After Kedlaya-Liu \cite[6.2.2-6.2.4]{3KL15}, \cite[Lemma 4.6.9]{3KL16}, and also see \cite{3T1})}}
%For any Frobenius $\varphi^a$-bundle $M$ over $\widetilde{\Pi}_{R,A}$, then we have that for sufficiently large number $n\geq 0$ the module $M$ could be generated by finitely many Frobenius $\varphi^a$-invariant elements of the global sections of $M(n)$.	
%\end{lemma}
%
%\begin{proof}
%See the proof of \cite{3T1}.	
%\end{proof}

\begin{lemma} \mbox{\bf{(After Kedlaya-Liu \cite[6.2.2-6.2.4]{3KL15}, \cite[Lemma 4.6.9]{3KL16}})}\\ \mbox{\bf{(And also see \cite[Proposition 2.11, Proposition 2.14]{3T1})}}\\
For any Frobenius $\varphi^a$-module $M$ over $\widetilde{\Pi}_{R,A}$, we have that for suifficiently large number $n\geq 1$ the space $H^1_{\varphi^a}(M(n))$ vanishes.	
\end{lemma}

\begin{proof}
See the proof of \cite[Proposition 2.11, Proposition 2.14]{3T1}.	
\end{proof}

\indent Then we have the following key corollary which is analog of \cite[Corollary 6.2.3, Lemma 6.3.3]{3KL15}, \cite[Corollary 4.6.10]{3KL16} and \cite[Corollary 2.13, Corollary 3.4, Proposition 3.9]{3T1}:

\begin{corollary} \label{2coro4.16}
I. When $M_\alpha,M,M_\beta$ are three Frobenius $\varphi^a$-bundles over the ring $\widetilde{\Pi}_{R,A}$ then we have that for sufficiently large $n\geq 0$ we have the following exact sequence:
\[
\xymatrix@R+0pc@C+0pc{
0\ar[r]\ar[r]\ar[r] &M_{\alpha,I}(n)^{\varphi^a}
\ar[r]\ar[r]\ar[r] &M_I(n)^{\varphi^a}
\ar[r]\ar[r]\ar[r] &M_{\beta,I}(n)^{\varphi^a} \ar[r]\ar[r]\ar[r] &0,
}
\]	
for each interval $I$.	\\
II. When $M_\alpha,M,M_\beta$ are three Frobenius $\varphi^a$-modules over the ring $\widetilde{\Pi}_{R,A}$ then we have that for sufficiently large $n\geq 0$ we have the following exact sequence:
\[
\xymatrix@R+0pc@C+0pc{
0\ar[r]\ar[r]\ar[r] &M_{\alpha,I}(n)^{\varphi^a}
\ar[r]\ar[r]\ar[r] &M_I(n)^{\varphi^a}
\ar[r]\ar[r]\ar[r] &M_{\beta,I}(n)^{\varphi^a} \ar[r]\ar[r]\ar[r] &0,
}
\]	
for each interval $I$.\\
III. For a Frobenius $\varphi^a$-bundle $M$ over $\widetilde{\Pi}_{R,A}$, and for each module $M_I$ over some $\widetilde{\Pi}_{R,A}^{I}$, we put for any element $f$ such that $\varphi^af=p^df$:
\begin{displaymath}
M_{I,f}:=\bigcup_{n\in \mathbb{Z}}f^{-n}M_I(dn)^{\varphi^a}.	
\end{displaymath}
Then with this convention suppose we have three Frobenius $\varphi^a$-bundles taking the form of $M_\alpha,M,M_\beta$ over
$\widetilde{\Pi}_{R,A}$, then for each closed interval $I$ we have the following is an exact sequence:
\[
\xymatrix@R+0pc@C+0pc{
0\ar[r]\ar[r]\ar[r] &M_{\alpha,I,f}\ar[r]\ar[r]\ar[r] &M_{I,f}
\ar[r]\ar[r]\ar[r] &M_{\beta,I,f} \ar[r]\ar[r]\ar[r] &0,
}
\]
where each module in the exact sequence is now a module over $\widetilde{\Pi}_{R,A}[1/f]^{\varphi^a}$;\\
IV. For a Frobenius $\varphi^a$-module $M$ over $\widetilde{\Pi}_{R,A}$, we put for any element $f$ such that $\varphi^af=p^df$:
\begin{displaymath}
M_{f}:=\bigcup_{n\in \mathbb{Z}}f^{-n}M(dn)^{\varphi^a}.	
\end{displaymath}
Then with this convention suppose we have three Frobenius $\varphi^a$-modules taking the form of $M_\alpha,M,M_\beta$ over
$\widetilde{\Pi}_{R,A}$, then we have the following is an exact sequence:
\[
\xymatrix@R+0pc@C+0pc{
0\ar[r]\ar[r]\ar[r] &M_{\alpha,f}\ar[r]\ar[r]\ar[r] &M_{f}
\ar[r]\ar[r]\ar[r] &M_{\beta,f} \ar[r]\ar[r]\ar[r] &0,
}
\]
where each module in the exact sequence is now a module over $\widetilde{\Pi}_{R,A}[1/f]^{\varphi^a}$.\\
%V. Suppose $M$ is now a pseudocoherent Frobenius $\varphi^a$-bundle over $\widetilde{\Pi}_{R,A}$, then we have that the corresponding module  module $M_I$ over some $\widetilde{\Pi}_{R,A}^{I}$, we put for any element $f$ such that $\varphi^af=p^df$:
%\begin{displaymath}
%M_{I,f}:=\bigcup_{n\in \mathbb{Z}}f^{-n}M_I(dn)^{\varphi^a}.	
%\end{displaymath}
%Then with this convention suppose we have three Frobenius $\varphi^a$-bundles taking the form of $M_\alpha,M,M_\beta$ over
%$\widetilde{\Pi}_{R,A}$, then for each closed interval $I$ we have the following is an exact sequence:
%\[
%\xymatrix@R+0pc@C+0pc{
%0\ar[r]\ar[r]\ar[r] &M_{\alpha,I,f}\ar[r]\ar[r]\ar[r] &M_{I,f}
%\ar[r]\ar[r]\ar[r] &M_{\beta,I,f} \ar[r]\ar[r]\ar[r] &0,
%}
%\]
%where each module in the exact sequence is now in our situation a pseudocoherent module over $\widetilde{\Pi}_{R,A}[1/f]^{\varphi^a}$.\\

\end{corollary}

\begin{proof}
See the proof of \cite[Corollary 6.2.3, Lemma 6.3.3]{3KL15}, \cite[Corollary 4.6.10]{3KL16} and \cite[Corollary 2.13, Corollary 3.4, Proposition 3.9]{3T1}.	
\end{proof}

\subsection{Noncommutative Imperfect Setting}

\indent We now consider the corresponding imperfectization of the corresponding constructions we considered above, after \cite{3KL16}. We will consider the corresponding towers in \cite[Chapter 5]{3KL16}, so we keep the assumption on the towers as in the commutative setting, namely:

\begin{assumption}
In this section, we are going to assume that $k$ is just $\mathbb{F}_{p^h}$.
\end{assumption}

%\begin{setting}
%Recall from \cite[Chapter 5]{3KL16}, we have the following setting up. The base will be a Banach adic ring $(H,H^+)$ over $\mathfrak{o}_E$ which is assumed to be uniform and carrying the corresponding spectral norm $\alpha$. Then we consider the tower:
%\[
%\xymatrix@R+0pc@C+0pc{
%...\ar[r]\ar[r]\ar[r] &H_{n-1}\ar[r]\ar[r]\ar[r] &H_n
%\ar[r]\ar[r]\ar[r] &H_{n+1} \ar[r]\ar[r]\ar[r] &...,
%}
%\]	
%where each map linking two adjacent rings are the corresponding morphisms of the corresponding Banach uniform adic rings whose corresponding induced maps on the adic spaces are surjective as in the corresponding construction in \cite[Definition 5.1.1]{3KL16} with the corresponding spectral norm $\alpha_n$. The infinite level of the tower will be denoted by $H_\infty$. This is not actually necessarily complete under the multiplicative extension of all the finite level spectral norms $\alpha_n$ namely $\alpha_\infty$, so we need to consider the completed ring $\overline{H}_\infty$. In the situation where the tower is Fontaine perfectoid, we have that following \cite[Definition 5.1.1]{3KL16} that this is called the corresponding perfectoid tower, which gives rise the equal characteristic counterpart $\overline{H'}_\infty$ wit the spectral norm $\overline{\alpha}_\infty$ correspondingly under the perfectoid correspondence in the sense of \cite[Theorem 3.3.8]{3KL16}. Recall that from \cite[Definition 5.1.1]{3KL16} the tower is called finite \'etale if each transition map is finite \'etale.
%\end{setting}

\indent Now we describe the corresponding imperfect period rings which we will deform. These rings are those introduced in \cite[Definition 5.2.1]{3KL16} by using series of imperfection processes. Recall in more detail in our noncommutative setting we have:

\begin{setting} \label{2setting6.3}
Fix a perfectoid tower $(H_\bullet,H^+_\bullet)$ Recall from \cite[Definition 5.2.1]{3KL16} we have the following different imperfect constructions:\\
A. First we have the ring $\overline{H'}_\infty$, which could give us the corresponding ring $\widetilde{\Omega}^{\mathrm{int}}_{\overline{H'}_\infty}$. Taking the corresponding product with $A$ we have the corresponding deformed period ring;\\
B. We then have the ring $\widetilde{\Pi}^{\mathrm{int},r}_{\overline{H'}_\infty}$ coming from $\overline{H'}_\infty$. Taking the corresponding product with $A$ we have the corresponding deformed period ring;\\
C. We then have the ring $\widetilde{\Pi}^{\mathrm{int}}_{\overline{H'}_\infty}$ coming from $\overline{H'}_\infty$. Taking the corresponding product with $A$ we have the corresponding deformed period ring;\\
D. We then have the ring $\widetilde{\Pi}^{\mathrm{bd},r}_{\overline{H'}_\infty}$ coming from $\overline{H'}_\infty$. Taking the corresponding product with $A$ we have the corresponding deformed period ring;\\
E. We then have the ring $\widetilde{\Pi}^{\mathrm{bd}}_{\overline{H'}_\infty}$ coming from $\overline{H'}_\infty$. Taking the corresponding product with $A$ we have the corresponding deformed period ring;\\
F. We then have the ring $\widetilde{\Pi}^{r}_{\overline{H'}_\infty}$ coming from $\overline{H'}_\infty$. Taking the corresponding product with $A$ we have the corresponding deformed period ring;\\
G. We then have the ring $\widetilde{\Pi}^{[s,r]}_{\overline{H'}_\infty}$ coming from $\overline{H'}_\infty$. Taking the corresponding product with $A$ we have the corresponding deformed period ring;\\
H. We then have the ring $\widetilde{\Pi}_{\overline{H'}_\infty}$ coming from $\overline{H'}_\infty$. Taking the corresponding product with $A$ we have the corresponding deformed period ring;\\
I. $\Pi^{\mathrm{int},r}_{H}$ comes from the ring $\widetilde{\Pi}^{\mathrm{int},r}_{\overline{H'}_\infty}$ consisting of those elements of $\widetilde{\Pi}^{\mathrm{int},r}_{\overline{H'}_\infty}$ with the requirement that whenever we have $n$ an integer such that $nh>-\log_pr$ we have then $\theta(\varphi^{-n}(x))\in H_n$. Taking the corresponding product with $A$ we have the corresponding deformed period ring;\\
J. $\Pi^{\mathrm{int},\dagger}_{H}$ is defined to be the corresponding union of the rings in $I$. Taking the corresponding product with $A$ we have the corresponding deformed period ring;\\
K. $\Omega^\mathrm{int}_H$ is defined to be the corresponding period ring coming from the corresponding $\pi$-adic completion of the ring $\Pi^{\mathrm{int},\dagger}_{H}$ in $J$. Taking the corresponding product with $A$ we have the corresponding deformed period ring;\\
L. $\breve{\Omega}^{\mathrm{int}}_{H}$ is the ring which is defined to be the union of all the $\varphi^{-n}\Omega^\mathrm{int}_H$. Taking the corresponding product with $A$ we have the corresponding deformed period ring;\\
M. $\breve{\Pi}^{\mathrm{int},r}_H$ is then the ring which is defined to be the union of all the $\varphi^{-n}\Pi^{\mathrm{int},p^{hn}r}_H$. Taking the corresponding product with $A$ we have the corresponding deformed period ring;\\
N. $\breve{\Pi}^{\mathrm{int},\dagger}_H$ is define to be union of all the $\varphi^{-n}\breve{\Pi}^{\mathrm{int},\dagger}_H$. Taking the corresponding product with $A$ we have the corresponding deformed period ring;\\
O. $\widehat{\Omega}^{\mathrm{int}}_{H}$ is defined to be the $\pi$-completion of $\breve{\Omega}^{\mathrm{int}}_{H}$. Taking the corresponding product with $A$ we have the corresponding deformed period ring;\\
P. $\widehat{\Pi}^{\mathrm{int},r}_H$ is defined to be the $\mathrm{max}\{\|.\|_{\overline{\alpha}_\infty^r},\|.\|_{\pi-\text{adic}}\}$-completion of $\breve{\Pi}^{\mathrm{int},r}_H$. Taking the corresponding product with $A$ we have the corresponding deformed period ring;\\
Q. We then have $\widehat{\Pi}^{\mathrm{int},\dagger}_H$ by taking the union over $r>0$. Taking the corresponding product with $A$ we have the corresponding deformed period ring;\\
R. Correspondingly we have $\breve{\Omega}_{H}$, $\breve{\Pi}^{\mathrm{bd},r}_H$, $\breve{\Pi}^{\mathrm{bd},\dagger}_H$ by inverting the element $\pi$. Taking the corresponding product with $A$ we have the corresponding deformed period ring;\\
S. Correspondingly we also have $\widehat{\Omega}_{H}$, $\widehat{\Pi}^{\mathrm{bd},r}_H$, $\widehat{\Pi}^{\mathrm{bd},\dagger}_H$ by inverting the corresponding element $\pi$. Taking the corresponding product with $A$ we have the corresponding deformed period ring;\\
T. We also have $\Omega_{H}$, $\Pi^{\mathrm{bd},r}_H$, $\Pi^{\mathrm{bd},\dagger}_H$ again by inverting the element $\pi$. Taking the corresponding product with $A$ we have the corresponding deformed period ring;\\
U. Taking the $\max\{\|.\|_{\overline{\alpha}_\infty^s},\|.\|_{\overline{\alpha}_\infty^r}\}$ (for $0<s\leq r$) completion of the ring $\Pi^{\mathrm{bd},r}_H$ we have the ring $\Pi^{[s,r]}_{H}$, while taking the Fr\'echet completion with respect to the norm $\|.\|_{\overline{\alpha}^t_\infty}$ for each $0<t\leq r$ we have the ring $\Pi^r_{H}$. Taking the corresponding product with $A$ we have the corresponding deformed period ring;\\
V. Taking the union we have the ring $\Pi_H$. Taking the corresponding product with $A$ we have the corresponding deformed period ring;\\
W. We use the notation $\breve{\Pi}_H$ to denote the corresponding union of all the $\varphi^{-n}\Pi_H$. Taking the corresponding product with $A$ we have the corresponding deformed period ring;\\
X. We use the notation $\breve{\Pi}^{[s,r]}_{H}$ to be the corresponding union of all the $\varphi^{-n}\Pi^{[p^{-hn}s,p^{-hn}r]}_H$. Taking the corresponding product with $A$ we have the corresponding deformed period ring;\\
Y. We use the notation $\breve{\Pi}_H^r$ to be the union of all the $\varphi^{-n}\breve{\Pi}_H^{p^{-hn}r}$. Taking the corresponding product with $A$ we have the corresponding deformed period ring.\\
\end{setting}

\indent Then we have the following direct analog of the relative version of the ring defined above (here as before the ring $A$ denotes a Banach affinoid algebra in the noncommutative setting):

\begin{setting} \label{2setting6.4}
Now we consider the deformation of the rings above:\\
I. We have the first group of the period rings in the deformed setting:
\begin{align}
\Pi^{\mathrm{int},r}_{H,A},\Pi^{\mathrm{int},\dagger}_{H,A},\Omega^\mathrm{int}_{H,A}, \Omega_{H,A}, \Pi^{\mathrm{bd},r}_{H,A}, \Pi^{\mathrm{bd},\dagger}_{H,A},\Pi^{[s,r]}_{H,A}, \Pi^r_{H,A}, \Pi_{H,A}.
\end{align}
II. We also have the second group of desired rings in the desired setting:
\begin{align}
\breve{\Pi}^{\mathrm{int},r}_{H,A},\breve{\Pi}^{\mathrm{int},\dagger}_{H,A},\breve{\Omega}^\mathrm{int}_{H,A}, \breve{\Omega}_{H,A}, \breve{\Pi}^{\mathrm{bd},r}_{H,A}, \breve{\Pi}^{\mathrm{bd},\dagger}_{H,A},\breve{\Pi}^{[s,r]}_{H,A}, \breve{\Pi}^r_{H,A}, \breve{\Pi}_{H,A}.	
\end{align}
III. We also have the third group:
\begin{align}
\widehat{\Pi}^{\mathrm{int},r}_{H,A},\widehat{\Pi}^{\mathrm{int},\dagger}_{H,A},\widehat{\Omega}^\mathrm{int}_{H,A}, \widehat{\Omega}_{H,A}, \widehat{\Pi}^{\mathrm{bd},r}_{H,A}, \widehat{\Pi}^{\mathrm{bd},\dagger}_{H,A}.	
\end{align}
\end{setting}

\indent Now we discuss some properties of the corresponding deformed version of the imperfect rings in our context, which is parallel to the corresponding discussion we made in the perfect setting and generalizing the corresponding discussion in \cite{3KL16}, again in our noncommutative setting:

\begin{proposition}\mbox{\bf{(After Kedlaya-Liu \cite[Lemma 5.2.10]{3KL16})}}
For any $0< r_1\leq r_2$ we have the following equality on the corresponding period rings:
\begin{displaymath}
\Pi^{\mathrm{int},r_1}_{H,\mathbb{Q}_p\{Z_1,...,Z_d\}}\bigcap	\Pi^{[r_1,r_2]}_{H,\mathbb{Q}_p\{Z_1,...,Z_d\}}=\Pi^{\mathrm{int},r_2}_{H,\mathbb{Q}_p\{Z_1,...,Z_d\}}.
\end{displaymath}	
\end{proposition}

\begin{proof}
We adapt the argument in \cite[Lemma 5.2.10]{3KL16} to prove this in the situation where $r_1<r_2$ (otherwise this is trivial), again one direction is easy where we only present the implication in the other direction. We take any element
\begin{center}
 $x\in \Pi^{\mathrm{int},r_1}_{H,\mathbb{Q}_p\{Z_1,...,Z_d\}}\bigcap	\Pi^{[r_1,r_2]}_{H,\mathbb{Q}_p\{Z_1,...,Z_d\}}$ 
\end{center} 
and take suitable approximating elements $\{x_i\}$ living in the bounded Robba ring such that for any $j\geq 1$ one can find some integer $N_j\geq 1$ we have whenever $i\geq N_j$ we have the following estimate:
\begin{displaymath}
\left\|.\right\|_{\overline{\alpha}_\infty^{t},\mathbb{Q}_p\{Z_1,...,Z_d\}}(x_i-x) \leq p^{-j}, \forall t\in [r_1,r_2].	
\end{displaymath}
Then we consider the corresponding decomposition of $x_i$ for each $i=1,2,...$ into a form having integral part and the rational part $x_i=y_i+z_i$ by setting
\begin{center}
 $y_i=\sum_{k=0,i_1,...,i_d}\pi^kx_{i,k,i_1,...,i_d}Z_1^{i_1}...Z_d^{i_d}$ 
\end{center} 
out of
\begin{center} 
$x_i=\sum_{k=n(x_i),i_1,...,i_d}\pi^kx_{i,k,i_1,...,i_d}Z_1^{i_1}...Z_d^{i_d}$.
\end{center}
Note that by our initial hypothesis we have that the element $x$ lives in the ring $\Pi^{\mathrm{int},r_1}_{H,\mathbb{Q}_p\{Z_1,...,Z_d\}}$ which further implies that 
\begin{displaymath}
\left\|.\right\|_{\overline{\alpha}_\infty^{r_1},\mathbb{Q}_p\{Z_1,...,Z_d\}}(\pi^kx_{i,k,i_1,...,i_d}Z_1^{i_1}...Z^{i_d}_d)	\leq p^{-j}.
\end{displaymath}
Therefore we have ${\overline{\alpha}_\infty}(\overline{x}_{i,k,i_1,...,i_d})\leq p^{(k-j)/r_1},\forall k< 0$ directly from this through computation, which implies that then:
\begin{align}
\left\|.\right\|_{\overline{\alpha}_\infty^{r_2},\mathbb{Q}_p\{Z_1,...,Z_d\}}(\pi^kx_{i,k,i_1,...,i_d}Z_1^{i_1}...Z^{i_d}_d)	&\leq p^{-k}p^{(k-j)r_2/r_1}\\
	&\leq p^{1+(1-j)r_1/r_1}.
\end{align}
Then one can read off the result directly from this estimate since under this estimate we can have the chance to modify the original approximating sequence $\{x_i\}$ by $\{y_i\}$ which are initially chosen to be in the integral Robba ring, which implies that actually the element $x$ lives in the right-hand side of the identity in the statement of the proposition.
\end{proof}

\begin{remark}
The following result cannot be derived from the previous proposition but the proof could be made essentially in the fashion, which is the same as the corresponding situation we encountered in the commutative situation, for instance see  \cref{proposition5.7}.	
\end{remark}

\begin{proposition}\mbox{\bf{(After Kedlaya-Liu \cite[Lemma 5.2.10]{3KL16})}}
For any $0< r_1\leq r_2$ we have the following equality on the corresponding period rings:
\begin{displaymath}
\Pi^{\mathrm{int},r_1}_{H,A}\bigcap	\Pi^{[r_1,r_2]}_{H,A}=\Pi^{\mathrm{int},r_2}_{H,A}.
\end{displaymath}	
Here $A$ is some Banach affinoid algebra over $\mathbb{Q}_p$ in the noncommutative setting.	
\end{proposition}

\begin{proof}
See the proof of \cref{proposition5.7}.		
\end{proof}

%%%%%%%%%%%%%%%%%%%%%%%%%%%%%%%%%%%%%%%%%%%%%%%%%

\begin{proposition} \mbox{\bf{(After Kedlaya-Liu \cite[Lemma 5.2.8]{3KL15} and \cite{3KL16})}}
Consider now in our situation the radii $0< r_1\leq r_2$, and consider any element 
\begin{center}
$x\in \Pi^{[r_1,r_2]}_{H,\mathbb{Q}_p\{Z_1,...,Z_d\}}$. 
\end{center}
Then we have that for each $n\geq 1$ one can decompose $x$ into the form of $x=y+z$ such that $y\in \pi^n\Pi^{\mathrm{int},r_2}_{H,\mathbb{Q}_p\{Z_1,...,Z_d\}}$ with $z\in \bigcap_{r\geq r_2}\Pi^{[r_1,r]}_{H,\mathbb{Q}_p\{Z_1,...,Z_d\}}$ with the following estimate for each $r\geq r_2$:
\begin{displaymath}
\left\|.\right\|_{\overline{\alpha}_\infty^r,\mathbb{Q}_p\{Z_1,...,Z_d\}}(z)\leq p^{(1-n)(1-r/r_2)}\left\|.\right\|_{\overline{\alpha}_\infty^{r_2},\mathbb{Q}_p\{Z_1,...,Z_d\}}(z)^{r/r_2}.	
\end{displaymath}

\end{proposition}

\begin{proof}
As in \cite[Lemma 5.2.8]{3KL15} and \cite{3KL16} and in the proof of our previous proposition we first consider those elements $x$ lives in the bounded Robba rings which could be expressed in general as
\begin{center}
 $\sum_{k=n(x),i_1,...,i_d}\pi^kx_{k,i_1,...,i_d}Z_1^{i_1}...Z_d^{i_d}$.
 \end{center}	
In this situation the corresponding decomposition is very easy to come up with, namely we consider the corresponding $y_i$ as the corresponding series:
\begin{displaymath}
\sum_{k\geq n,i_1,...,i_d}\pi^kx_{k,i_1,...,i_d}Z_1^{i_1}...Z_d^{i_d}	
\end{displaymath}
which give us the desired result since we have in this situation when focusing on each single term:
\begin{align}
\left\|.\right\|_{\overline{\alpha}_\infty^r,\mathbb{Q}_p\{Z_1,...,Z_d\}}&(\pi^kx_{k,i_1,...,i_d}Z_1^{i_1}...Z_d^{i_d})=p^{-k}\overline{\alpha}_\infty(\overline{x}_{k,i_1,...,i_d})^r\\
&=p^{-k(1-r/r_2)}\left\|.\right\|_{\overline{\alpha}_\infty^{r_2},\mathbb{Q}_p\{Z_1,...,Z_d\}}(\pi^kx_{k,i_1,...,i_d}Z_1^{i_1}...Z_d^{i_d})^{r/r_2}\\
&\leq p^{(1-n)(1-r/r_2)}\left\|.\right\|_{\overline{\alpha}_\infty^{r_2},\mathbb{Q}_p\{Z_1,...,Z_d\}}(\pi^kx_{k,i_1,...,i_d}Z_1^{i_1}...Z_d^{i_d})^{r/r_2}
\end{align}
for all those suitable $k$. Then to tackle the more general situation we consider the approximating sequence consisting of all the elements in the bounded Robba ring as in the usual situation considered in \cite[Lemma 5.2.8]{3KL15} and \cite{3KL16}, namely we inductively construct the following approximating sequence just as:
\begin{align}
\left\|.\right\|_{\overline{\alpha}_\infty^r,\mathbb{Q}_p\{Z_1,...,Z_d\}}(x-x_0-...-x_i)\leq p^{-i-1}	\left\|.\right\|_{\overline{\alpha}_\infty^r,\mathbb{Q}_p\{Z_1,...,Z_d\}}(x), i=0,1,..., r\in [r_1,r_2].
\end{align}
Here all the elements $x_i$ for each $i=0,1,...$ are living in the bounded Robba ring, which immediately gives rise to the suitable decomposition as proved in the previous case namely we have for each $i$ the decomposition $x_i=y_i+z_i$ with the desired conditions as mentioned in the statement of the proposition. We first take the series summing all the elements $y_i$ up for all $i=0,1,...$, which first of all converges under the norm $\left\|.\right\|_{\overline{\alpha}_\infty^r,\mathbb{Q}_p\{Z_1,...,Z_d\}}$ for all the radius $r\in [r_1,r_2]$, and also note that all the elements $y_i$ within the infinite sum live inside the corresponding integral Robba ring $\Pi^{\mathrm{int},r_2}_{H,\mathbb{Q}_p\{Z_1,...,Z_d\}}$, which further implies the corresponding convergence ends up in $\Pi^{\mathrm{int},r_2}_{H,\mathbb{Q}_p\{Z_1,...,Z_d\}}$. For the elements $z_i$ where $i=0,1,...$ also sum up to a converging series in the desired ring since combining all the estimates above we have:
\begin{displaymath}
\left\|.\right\|_{\overline{\alpha}_\infty^r,\mathbb{Q}_p\{Z_1,...,Z_d\}}(z_i)\leq p^{(1-n)(1-r/r_2)}\left\|.\right\|_{\overline{\alpha}_\infty^{r_2},\mathbb{Q}_p\{Z_1,...,Z_d\}}(x)^{r/r_2}.	
\end{displaymath}
\end{proof}

%\begin{corollary}
%Consider now in our situation the radii $0< r_1\leq r_2$, and consider any element $x\in \widetilde{\Pi}^{[r_1,r_2]}_{R,A}$. Then we have that for each $n\geq 1$ one can decompose $x$ into the form of $x=y+z$ such that $y\in \pi^n\widetilde{\Pi}^{\mathrm{int},r_2}_{R,A}$ with $z\in \bigcap_{r\geq r_2}\widetilde{\Pi}^{[r_1,r]}_{R,A}$ with the following estimate for each $r\geq r_2$:
%\begin{displaymath}
%\left\|.\right\|_{\alpha^r,A}(z)\leq p^{(1-n)(1-r/r_2)}\left\|.\right\|_{\alpha^{r_2},A}(z)^{r/r_2}.	
%\end{displaymath}	
%\end{corollary}
%
%\begin{proof}
%This could be derived directly from the previous proposition, by considering the corresponding residual norms and then the spectral norms.	
%\end{proof}

\begin{proposition} \mbox{\bf{(After Kedlaya-Liu \cite[Lemma 5.2.10]{3KL15})}}
We have the following identity:
\begin{displaymath}
\Pi^{[s_1,r_1]}_{H,\mathbb{Q}_p\{Z_1,...,Z_d\}}\bigcap\Pi^{[s_2,r_2]}_{H,\mathbb{Q}_p\{Z_1,...,Z_d\}}=\Pi^{[s_1,r_2]}_{H,\mathbb{Q}_p\{Z_1,...,Z_d\}},
\end{displaymath}
here the radii satisfy $<s_1\leq s_2 \leq r_1 \leq r_2$.
\end{proposition}

\begin{proof}
In our situation one direction is obvious while on the other hand we consider any element $x$ in the intersection on the left, then by the previous proposition we	have the decomposition $x=y+z$ where $y\in \Pi^{\mathrm{int},r_1}_{H,\mathbb{Q}_p\{Z_1,...,Z_d\}}$ and $z\in \Pi^{[s_1,r_2]}_{H,\mathbb{Q}_p\{Z_1,...,Z_d\}}$. Then as in \cite[Lemma 5.2.10]{3KL15} section 5.2 we look at $y=x-z$ which lives in the intersection:
\begin{displaymath}
\Pi^{\mathrm{int},r_1}_{H,\mathbb{Q}_p\{Z_1,...,Z_d\}}\bigcap	\Pi^{[s_2,r_2]}_{H,\mathbb{Q}_p\{Z_1,...,Z_d\}}=\Pi^{\mathrm{int},r_2}_{H,\mathbb{Q}_p\{Z_1,...,Z_d\}}
\end{displaymath}
which finishes the proof.
\end{proof}

\begin{proposition} \mbox{\bf{(After Kedlaya-Liu \cite[Lemma 5.2.10]{3KL15})}}
We have the following identity:
\begin{displaymath}
\Pi^{[s_1,r_1]}_{H,A}\bigcap \Pi^{[s_2,r_2]}_{H,A}=\Pi^{[s_1,r_2]}_{H,A},
\end{displaymath}
here the radii satisfy $<s_1\leq s_2 \leq r_1 \leq r_2$.
	
\end{proposition}

\begin{proof}
See the proof of \cref{proposition5.7}.		
\end{proof}

\begin{remark}
Again one can follow the same strategy to deal with the corresponding equal-characteristic situation.	
\end{remark}

\subsection{Modules and Bundles over Noncommutative Rings}

Now we consider the modules and bundles over the rings introduced in the previous subsection. First we make the following assumption:

\begin{setting}
Recall that from \cite[Definition 5.2.3]{3KL16} any tower $(H_\bullet,H_\bullet^+)$ is called weakly decompleting if we have that first the density of the perfection of $H_{\infty}$ in $\overline{H}_\infty$. Here the ring $H_\infty$ is the ring coming from the mod-$\pi$ construction of the ring $\Omega^\mathrm{int}_{H}$, also at the same time one can find some $r>0$ such that the corresponding modulo $\pi$ operation from the ring $\Omega^\mathrm{int}_{H}$ to the ring $H_\infty$ is actually surjective strictly. 
\end{setting}

\begin{assumption} \label{assumption642}
We now assume that we are basically in the situation where $(H_\bullet,H_\bullet^+)$ is actually weakly decompleting. Also as in \cite[Lemma 5.2.7]{3KL16} we assume we fix some radius $r_0>0$, for instance this will correspond to the corresponding index in the situation we have the corresponding noetherian tower. Recall that a tower is called noetherian if we have some specific radius as above such that we have the strongly noetherian property on the ring $\Pi^{[s,r]}_{H}$ with $[s,r]\subset (0,r_0]$. We now assume that the tower is then noetherian and any ring $\Pi^{[s,r]}_{H,A}$ with $[s,r]\subset (0,r_0]$ is noetherian in the noncommutative setting. 
\end{assumption}

\begin{example}
This assumption could be achieved as in the following. First naively we can consider a noncommutative affinoid algebra which is finite over some commutative affinoid algebra over $\mathbb{Q}_p$, then tensor with the corresponding the Robba ring $\Pi^{[s,r]}_{H}$ which is assumed to be strongly noetherian. Then for some deeper case, we can take the Robba ring $D_{[a,b]}(\mathbb{Z}_p,\mathbb{Q}_p)$ attached to the group $\mathbb{Z}_p$ in the sense of \cite[Proposition 4.1]{3Z1} (namely the usual Robba ring over $\mathbb{Q}_p$, here as in \cite[Proposition 4.1]{3Z1} we assume that $a,b\in p^\mathbb{Q}$), and we take the corresponding complete tensor product with the local chart ring of the distribution algebra attached to a uniform $p$-adic Lie group $G$ as in \cite[Proposition 4.1]{3Z1}:
\begin{displaymath}
D_\rho(G,\mathbb{Q}_p)	
\end{displaymath}
with some radius $\rho>0$ living in $p^\mathbb{Q}$. The whole product:
\begin{displaymath}
D_{[a,b]}(\mathbb{Z}_p,\mathbb{Q}_p)\widehat{\otimes}D_\rho(G,\mathbb{Q}_p)	
\end{displaymath}
is noetherian. Indeed, as in \cite[Proposition 4.1]{3Z1} we apply the criterion in \cite[Proposition I.7.1.2]{3LVO} to check that actually the graded ring induced by the product norm:
\begin{displaymath}
gr^\bullet_{\|.\|_\otimes}(D_{[a,b]}(\mathbb{Z}_p,\mathbb{Q}_p)\widehat{\otimes}D_\rho(G,\mathbb{Q}_p)	)\overset{\sim}{\rightarrow} gr^\bullet_{\|.\|_\otimes}(D_{[a,b]}(\mathbb{Z}_p,\mathbb{Q}_p)\otimes D_\rho(G,\mathbb{Q}_p)	)
\end{displaymath}
admits a surjection map from:
\begin{align}
gr^\bullet_{\|.\|_{[a,b]}}(D_{[a,b]}(\mathbb{Z}_p,\mathbb{Q}_p))\otimes_{gr^\bullet\mathbb{Q}_p} gr^\bullet_{\|.\|_\rho}(D_\rho(G,\mathbb{Q}_p)	)
\end{align}
As in the proof of \cite[Proposition 4.1]{3Z1} one can show that this is a tensor product of a Laurent polynomial ring and a polynomial ring, which is noetherian.
\end{example}

\indent Then we can start to discuss the corresponding modules over the rings in our deformed setting, first as in \cite[Lemma 5.3.3]{3KL16} the following result should be derived from our construction:

\indent Then we deform the basic notation of bundles in \cite[Definition 5.3.6]{3KL16}:

\begin{definition} \mbox{\bf{(After Kedlaya-Liu \cite[Definition 5.3.6]{3KL16})}}
We define the bundle over the ring $\Pi^{r_0}_{H,A}$ to be a collection $(M_I)_I$ of finite projective right modules over each $\Pi_{H,A}^{I}$ with $I\subset (0,r_0]$ closed subintervals of $(0,r_0]$ such that we have the following requirement in the glueing fashion. First for any $I_1\subset I_2$ two closed intervals we have $M_{I_2}\otimes_{\Pi_{H,A}^{I_2}}\Pi_{H,A}^{I_1}\overset{\sim}{\rightarrow} M_{I_1}$ with the obvious cocycle condition with respect to three closed subintervals of $(0,r_0]$ namely taking the form of $I_1\subset I_2\subset I_3$.\\
\indent We define the pseudocoherent sheaf over the ring $\Pi^{r_0}_{H,A}$ to be a collection $(M_I)_I$ of pseudocoherent right modules over each $\Pi_{H,A}^{I}$ with $I\subset (0,r_0]$ closed subintervals of $(0,r_0]$ such that we have the following requirement in the glueing fashion. First for any $I_1\subset I_2$ two closed intervals we have $M_{I_2}\otimes_{\Pi_{H,A}^{I_2}}\Pi_{H,A}^{I_1}\overset{\sim}{\rightarrow} M_{I_1}$ with the obvious cocycle condition with respect to three closed subintervals of $(0,r_0]$ namely taking the form of $I_1\subset I_2\subset I_3$. 	
\end{definition}

\indent We make the following discussion around the corresponding module and sheaf structures defined above.

\begin{lemma} \mbox{\bf{(After Kedlaya-Liu \cite[Lemma 5.3.8]{3KL16})}}
We have the isomorphism between the ring $\Pi^r_{H,A}$ and the inverse limit of the ring $\Pi^{[s,r]}_{H,A}$ with respect to the radius $s$ by the map $\Pi^r_{H,A}\rightarrow \Pi^{[s,r]}_{H,A}$.	
\end{lemma}
 
\begin{proof}
As in \cite[Lemma 5.3.8]{3KL16} it is injective by the isometry, and then use the corresponding elements $x_n,n=0,1,...$ in the dense ring $\Pi^{r}_{H,A}$ to approximate any element $x$ in the ring $\Pi^{[s,r]}_{H,A}$ in the same way as in \cite[Lemma 5.3.8]{3KL16}:
\begin{displaymath}
\|.\|_{\overline{\alpha}_\infty^t,A}(x-x_n)\leq p^{n}	
\end{displaymath}
for any radius $t$ now living in the corresponding interval $[r2^{-n},r]$. This will establish Cauchy sequence which finishes the proof as in \cite[Lemma 5.3.8]{3KL16}.
\end{proof}

\begin{proposition} \mbox{\bf{(After Kedlaya-Liu \cite[Lemma 5.3.9]{3KL16})}}
For some radius $r\in (0,r_0]$. Suppose we have that $M$ is a vector bundle in the general setting or $M$ is a pseudocoherent sheaf in the setting where the tower is noetherian. Then we have that the corresponding global section is actually dense in each section with respect to some closed interval. And then we have the corresponding vanishing result of the first derived inverse limit functor.	
\end{proposition}

\begin{proof}
See \cite[Lemma 5.3.9]{3KL16}. In our setting, note that we are actually in the formalism of the corresponding (noncommutative) Fr\'echet-Stein algebras as in \cite{3ST1}.	
\end{proof}

\indent The interesting issue here as in \cite{3KL16} is the corresponding finitely generated of the global section of a pseudocoherent sheaf which is actually not guaranteed in general. Therefore as in \cite{3KL16} we have to distinguish the corresponding well-behaved sheaves out from the corresponding category of all the corresponding pseudocoherent sheaves. But we are not quite for sure how more complicated the noncommutative situation is than the commutative setting. For the completeness of the presentation we present some conjectures in our mind.

\begin{conjecture} \mbox{\bf{(After Kedlaya-Liu \cite[Lemma 5.3.10]{3KL16})}}
As in the previous proposition we choose some $r\in (0,r_0]$. Now assume that the corresponding tower $(H_\bullet,H_\bullet^+)$ is noetherian. Now for any pseudocoherent sheaf $M$ defined above we have the following three statements are equivalent. A. The first statement is that one can find a sequence of positive integers $x_1,x_2,...$ such that for any closed interval living inside $(0,r]$ the section of the sheaf with respect this closed interval admits a projective resolution of modules with corresponding ranks bounded by the sequence of integer $x_1,x_2,...$. B. The second statement is that for any locally finite covering of the corresponding interval $(0,r]$ which takes the corresponding form of $\{I_i\}$ one can find a sequence of positive integers $x_1,x_2,...$ such that for any closed interval living inside $\{I_i\}$ the section of the sheaf with respect this closed interval admits a projective resolution of modules with corresponding ranks bounded by the sequence of integer $x_1,x_2,...$. C. Lastly the third statement is that the corresponding global section is a pseudocoherent module over the ring $\Pi^r_{H,A}$. 	
\end{conjecture}

%\begin{proposition} \mbox{\bf{(After Kedlaya-Liu \cite[Lemma 5.3.10]{3KL16})}}
%As in the previous proposition we choose some $r\in (0,r_0]$. Now assume that the corresponding tower $(H_\bullet,H_\bullet^+)$ is noetherian. And we assume that Now for any pseudocoherent sheaf $M$ defined above we have the following second statement implies the corresponding first statement. The first statement is that one can find a sequence of positive integers $x_1,x_2,...$ such that for any closed interval living inside $(0,r]$ the section of the sheaf with respect this closed interval admits a projective resolution of modules with corresponding ranks bounded by the sequence of integer $x_1,x_2,...$. The second statement is that the corresponding global section is a pseudocoherent module over the ring $\Pi^r_{H,A}$. 	
%\end{proposition}
%
%
%\begin{proof}
%See \cite[Lemma 5.3.10]{3KL16}.	
%\end{proof}

\indent As in \cite[Definition 5.3.11]{3KL16} we call the sheaf satisfies the corresponding equivalent conditions in the proposition above uniform pseudocoherent sheaf. Then we have the following analog of \cite[Lemma 5.3.12]{3KL16}:

\begin{conjecture} \mbox{\bf{(After Kedlaya-Liu \cite[Lemma 5.3.12]{3KL16})}}
The global section functor defines the corresponding equivalence between the categories of the following two sorts of objects. The first ones are the corresponding uniform pseudocoherent sheaves over $\Pi^r_{H,A}$. The second ones are those pseudocoherent modules over the ring $\Pi^r_{H,A}$. 	
\end{conjecture}

%\begin{proof}
%See \cite[Lemma 5.3.12]{3KL16}.	
%\end{proof}

\begin{conjecture} \mbox{\bf{(After Kedlaya-Liu \cite[Lemma 5.3.13]{3KL16})}}
The global section functor defines the corresponding equivalence between the categories of the following two sorts of objects. The first ones are the corresponding finite projective uniform pseudocoherent sheaves over $\Pi^r_{H,A}$. The second ones are those finite projective modules over the ring $\Pi^r_{H,A}$. 	
\end{conjecture}

\subsection{$\varphi$-Module Structures and $\Gamma$-Module Structures over Noncommutative Rings}

\indent Now we consider the corresponding Frobenius actions over the corresponding imperfect rings we defined before, note that the corresponding Frobenius actions are induced from the corresponding imperfect rings in the undeformed situation from \cite{3KL16} which is to say that the Frobenius action on the ring $A$ is actually trivial.

\begin{assumption}
We now drop the corresponding requirement on $A$ which makes the noncommutative deformation of the Robba ring over some interval both left and right noetherian, in \cref{assumption642}. 	
\end{assumption}

\indent First we consider the corresponding Frobenius modules:

\begin{definition} \mbox{\bf{(After Kedlaya-Liu \cite[Definition 5.4.2]{3KL16})}}
Over the period rings $\Pi_{H,A}$ or $\breve{\Pi}_{H,A}$  (which is denoted by $\triangle$ in this definition) we define the corresponding $\varphi^a$-modules over $\triangle$ which are respectively projective to be the corresponding finite projective right modules over $\triangle$ with further assigned semilinear action of the operator $\varphi^a$. 
\end{definition}

\begin{definition} \mbox{\bf{(After Kedlaya-Liu \cite[Definition 5.4.2]{3KL16})}}
Over each rings $\triangle=\Pi^r_{H,A},\breve{\Pi}^r_{H,A}$ we define the corresponding projective $\varphi^a$-module over any $\triangle$ to be the corresponding finite projective right module $M$ over $\triangle$ with additionally endowed semilinear Frobenius action from $\varphi^a$ such that we have the isomorphism $\varphi^{a*}M\overset{\sim}{\rightarrow}M\otimes \square$ where the ring $\square$ is one $\triangle=\Pi^{r/p}_{H,A},\breve{\Pi}^{r/p}_{H,A}$.   
\end{definition}

\begin{definition} \mbox{\bf{(After Kedlaya-Liu \cite[Definition 5.4.2]{3KL16})}}
Again as in\\ \cite[Definition 5.4.2]{3KL16}, we define the corresponding projective $\varphi^a$-modules over ring $\Pi^{[s,r]}_{H,A}$ or $\breve{\Pi}^{[s,r]}_{H,A}$  to be the finite projective right modules (which will be denoted by $M$) over $\Pi^{[s,r]}_{H,A}$ or $\breve{\Pi}^{[s,r]}_{H,A}$ respectively additionally endowed with semilinear Frobenius action from $\varphi^a$ with the following isomorphisms:
\begin{align}
\varphi^{a*}M\otimes_{\Pi_{H,A}^{[sp^{-ah},rp^{-ah}]}}\Pi_{H,A}^{[s,rp^{-ah}]}\overset{\sim}{\rightarrow}M\otimes_{\Pi_{H,A}^{[s,r]}}\Pi_{H,A}^{[s,rp^{-ah}]}
\end{align}
and
\begin{align}
\varphi^{a*}M\otimes_{\breve{\Pi}_{R,A}^{[sp^{-ah},rp^{-ah}]}}\breve{\Pi}_{R,A}^{[s,rp^{-ah}]}\overset{\sim}{\rightarrow}M\otimes_{\breve{\Pi}_{R,A}^{[s,r]}}\breve{\Pi}_{R,A}^{[s,rp^{-ah}]}.\\
\end{align} 
\end{definition}

\indent Also one can further define the corresponding bundles carrying semilinear Frobenius in our context as in the situation of \cite[Definition 5.4.10]{3KL16}:

\begin{definition} \mbox{\bf{(After Kedlaya-Liu \cite[Definition 5.4.10]{3KL16})}}
Over the ring $\Pi^t_{H,A}$ or $\breve{\Pi}^t_{H,A}$ we define a corresponding projective Frobenius bundle to be a family $(M_I)_I$ of finite projective right modules over each $\Pi^I_{H,A}$ or $\breve{\Pi}^I_{H,A}$ respectively carrying the natural Frobenius action coming from the operator $\varphi^a$ such that for any two involved intervals having the relation $I\subset J$ we have:
\begin{displaymath}
M_J\otimes_{\Pi^J_{H,A}}\Pi^I_{H,A}\overset{\sim}{\rightarrow}	M_I
\end{displaymath}
and 
\begin{displaymath}
M_J\otimes_{\breve{\Pi}^J_{H,A}}\breve{\Pi}^I_{H,A}\overset{\sim}{\rightarrow}	M_I
\end{displaymath}
with the obvious cocycle condition. Here we have to propose condition on the intervals that for each $I=[s,r]$ involved we have $s\leq r/p^{ah}$. As in \cite[Definition 5.4.10]{3KL16} one can take the corresponding 2-limit in the direct sense to define the corresponding objects over the full Robba rings.
\end{definition}

\indent We can then compare the corresponding objects defined above:

\begin{conjecture} \mbox{\bf{(After Kedlaya-Liu \cite[Lemma 5.4.11]{3KL16})}}\\
I. Consider the following objects for some radius $r_0$ in our situation. The first group of objects are those finite projective $\varphi^a$-modules over the Robba ring $\Pi^{r_0}_{H,A}$. The second group of objects are those finite projective $\varphi^a$-bundles over the Robba ring $\Pi^{r_0}_{H,A}$. The third group of objects are those finite projective $\varphi^a$-modules over the Robba ring $\Pi^{[s,r]}_{H,A}$ for some radii $0<s\leq r\leq r_0$. Then we have that the corresponding categories of the two groups of objects are equivalent.
%II. Consider the following objects for some radius $r_0$ in our situation. The first group of objects are those pseudocoherent $\varphi^a$-modules over the Robba ring $\Pi^{r_0}_{H,A}$. The second group of objects are those pseudocoherent $\varphi^a$-bundles over the Robba ring $\Pi^{r_0}_{H,A}$. Then we have that the corresponding categories of the two groups of objects are equivalent. \\
%III. Consider the following objects for some radius $r_0$ in our situation. The first group of objects are those finite projective dimensional $\varphi^a$-modules over the Robba ring $\Pi^{r_0}_{H,A}$. The second group of objects are those finite projective dimensional $\varphi^a$-bundles over the Robba ring $\Pi^{r_0}_{H,A}$. Then we have that the corresponding categories of the two groups of objects are equivalent.  	
\end{conjecture}

%
%\begin{proof}
%See the proof of \cite[Lemma 5.4.11]{3KL16}. This needs the noncommutative version of \cite[Proposition 2.7.16]{3KL16}.
%\end{proof}

\indent Now we define the corresponding $\Gamma$-modules over the period rings attached to the tower $(H_\bullet,H_\bullet^+)$. The corresponding structures are actually abstractly defined in the same way as in \cite{3KL16} and parallel to our commutative setting. First we consider the deformation of the corresponding complex $*_{H^\bullet}$ for any ring $*$ in \cref{2setting6.4}.

\begin{assumption}
Recall that the corresponding tower is called decompleting if it is weakly decompleting, finite \'etale on each finite level and having the exact sequence $\overline{\varphi}^{-1}H'_{H^\bullet}/H'_{H^\bullet}$ is exact. We now assume that the tower $(H_\bullet,H^+_\bullet)$ is then decompleting.	
\end{assumption}

%\begin{definition}
%Recall from \cite{3KL16} we have that the corresponding notation $*_{H^k}$ for some corresponding nonnegative integer $k>0$ to denote the corresponding completed tensor product of the corresponding periods rings in $\label{setting6.4}$ of order $k+1$ over the given tower with itself. Then one extend the corresponding definition to our situation denoted by $*_{(H^k,A^k}$. We now denote the map $\sharp_{p,k} *_{H^k,A} \rightarrow *_{H^k+1,A}$ by the map induced from the corresponding map 
%\end{definition}

\begin{setting}
Assume now $\Gamma$ is a topological group as in \cite[Definition 5.5.5]{3KL16} acting on the corresponding period rings in the \cref{2setting6.3} in the original context of \cref{2setting6.3}. Then we consider the corresponding induced continous action over the corresponding deformed version in our context namely in \cref{2setting6.4}. Assume now that the tower is Galois with the corresponding Galois group $\Gamma$.	
\end{setting}

\begin{definition}
We now consider the corresponding inhomogeneous continuous cocycles of the group $\Gamma$, as in \cite[Definition 5.5.5]{3KL16} we use the following notation to denote the corresponding complex extracted from a single tower for a given period ring $*_{H,A}$ in \cref{2setting6.4} for each $k>0$:
\begin{displaymath}
*_{H^k,A}:=C_\mathrm{con}(\Gamma^k,*_{H})\widehat{\otimes}_{\sharp} ?
\end{displaymath}
where $\sharp=\mathbb{Q}_p,\mathbb{Z}_p$ and $?=A,\mathfrak{o}_A$ respectively, which forms the corresponding complex $(*_{H^\bullet,A},d^\bullet)$ with the corresponding differential as in \cite[Definition 5.5.5]{3KL16} in the sense of continuous group cohomology.	
\end{definition}

\begin{definition}
Having established the corresponding meaning of the $\Gamma$-structure we now consider the corresponding definition of $\Gamma$-modules. Such modules called the corresponding right $\Gamma$-modules are defined over the corresponding rings in \cref{2setting6.4}. Again we allow that the corresponding right modules to be finite projective, or pseudocoherent or fpd over the rings in \cref{2setting6.4}. And the modules are defined to be carrying the corresponding continuous semilinear action from the group $\Gamma$.	
\end{definition}

\begin{proposition} \mbox{\bf{(After Kedlaya-Liu \cite[Corollary 5.6.5]{3KL16})}}\\
The complex $\varphi^{-1}\Pi^{[sp^{-h},rp^{-h}]}_{H^\bullet_{\geq n},A}/\Pi^{[s,r]}_{H^\bullet_{\geq n},A}$ and the complex $\widetilde{\Pi}^{[s,r]}_{H^\bullet_{\geq n},A}/\Pi^{[s,r]}_{H^\bullet_{\geq n},A}$ are strictly exact for any truncation index $n$. The corresponding radii satisfy the corresponding relation $0<s\leq r\leq r_0$.
\end{proposition}

\begin{proof}
See \cite[Corollary 5.6.5]{3KL16}, and consider the corresponding Schauder basis of $A$.	
\end{proof}

\begin{proposition} \mbox{\bf{(After Kedlaya-Liu \cite[Lemma 5.6.6]{3KL16})}}
The complex 
\begin{displaymath}
M\otimes_{\Pi^{[s,r]}_{H,A}} \varphi^{-(\ell+1)}\Pi^{[sp^{-h(\ell+1)},rp^{-h(\ell+1)}]}_{H^\bullet,A}/ \varphi^{-\ell}\Pi^{[sp^{-h\ell},rp^{-h\ell}]}_{H^\bullet,A}	
\end{displaymath}
and the complex 
\begin{displaymath}
M\otimes_{\Pi^{[s,r]}_{H,A}} \widetilde{\Pi}^{[s,r]}_{H^\bullet,A}/\varphi^{-\ell}\Pi^{[sp^{-h(\ell)},rp^{-h(\ell)}]}_{H^\bullet,A}
\end{displaymath}
are strictly exact for any truncation index $n$. The corresponding radii satisfy the corresponding relation $0<s\leq r\leq r_0$, and $\ell$ is bigger than some existing truncated integer $\ell_0\geq  0$.
\end{proposition}

\begin{proof}
See \cite[Lemma 5.6.6]{3KL16}.	
\end{proof}

\begin{lemma} \mbox{\bf{(After Kedlaya-Liu \cite[Lemma 5.6.8]{3KL16})}} Suppose we have a direct system of Banach rings $(B_i,\iota_i)$ where the corresponding map $\iota_i$ as in \cite[Lemma 5.6.8]{3KL16} is submetric for each $i$. As in \cite[Lemma 5.6.8]{3KL16} we assume the direct limit is endowed with the corresponding infimum seminorm. Now assume that $B\rightarrow C$ is a isometry and of the corresponding image which is assumed to be dense into a Banach ring $C$. Then we have that one can have the chance to have that any finitely generated right projective module over the ring $C$ could come from the corresponding one finitely generated right projective module over some ring $B_i$ after the corresponding base change.

\end{lemma}

\begin{proof}
We need to work in the corresponding noncommutative setting. We adapt the argument in the proof of \cite[Lemma 5.6.8]{3KL16} with some possible modification due to the fact that the rings are noncommutative. Let $O$ be the matrix attached to the corresponding projector associated to the module over the ring $C$. Now choose $P$ such that $\|O-P\|< \|O\|^{-3}$ (and we need to guarantee that this $P$ could live over some $B_i$). Now we consider the corresponding iteration where $Q_0:=P$ and set $Q_{k+1}=3Q_k^2-2Q_k^3$, so we have:
\begin{align}
Q_{k+1}-Q_k&=(1-2Q_k)(Q_k^2-Q_k),\\
Q^2_{k+1}-Q_{k+1}&=(4Q_k^2-4Q_k-3)(Q_k^2-Q_k)^2,	
\end{align}
which implies that by induction (see the estimate on $P^2-P$ below):
\begin{align}
\|Q_{k+1}-O\|&\leq \|O\|\|P^2-P\|,\\
\|Q^2_{k}-Q_{k}\|&\leq (\|P^2-P\|\|O\|^2)^{2^k}\|O\|^{-2}.	
\end{align}
%For instance for $Q_k-O$ we could choose suitable $P$ such that we have:
%\begin{align}
%\|Q_0-O\|&=\|P-O\|\\
%         &\leq \|O\|\|P^2-P\|	
%\end{align}
%and if the statement is true for $k$ then we have:
%\begin{align}
%\|Q_{k+1}-O\|&=\|Q_{k+1}-Q_k+Q_k-O\|\\
%             &\leq \sup\{\|Q_{k+1}-Q_k\|,\|Q_k-O\|\}\\
%             &\leq \sup\{\|Q^2_{k}-Q_k\|,\|O\|\|P^2-P\|\}\\	
%             &\leq \sup\{\}
%\end{align}
%
We estimate $P^2-P$ in the following way:
\begin{align}
P^2-P &=P^2-P-O^2+O	\\
      &=(P-O)(P+O-1)+OP-PO\\
      &=(P-O)(P+O-1)+OP-O^2+O^2-PO\\
      &=(P-O)(P+O-1)+O(P-O)+(O-P)O
\end{align}
which implies that actually:
\begin{align}
\|P^2-P\| &\leq \sup\{\|(P-O)(P+O-1)\|,\|O(P-O)\|,\|(O-P)O\|\}\\
	      &\leq \sup\{\|(P-O)(P+O-1)\|,\|O\|\|(P-O)\|,\|(O-P)\|\|O\|\}\\
	      &< \|O\|^{-3}\|O\|\\
	      &< \|O\|^{-2}
\end{align}
which by taking the limit gives rise to the desired matrix $Q$ which proves the result. Indeed one can consider the corresponding projection induced by $Q$ from some free module over $B_i$, which gives rise to the desired projective module over $B_i$ whose base change to $C$ gives the initially chosen finite projective module 
over $C$ since one can compute:
\begin{align}
OQ+(1-O)(1-Q)-1&=OQ-O^2-Q^2+OQ\\
	           &=O(Q-O)-(Q-O)Q
\end{align}
which implies that $\|OQ+(1-O)(1-Q)-1\|\leq 1$.

\end{proof}

\begin{remark}
The bound for $Q_k-O$ is not the same as in \cite[Lemma 5.6.8]{3KL16}, we would like to thank Professor Kedlaya for letting us know this correction.	
\end{remark}

\begin{proposition} \mbox{\bf{(After Kedlaya-Liu \cite[Lemma 5.6.9]{3KL16})}}
With the corresponding notations as above we have that the corresponding base change from the ring taking the form of $\breve{\Pi}^{[s,r]}_{H,A}$ to the ring taking the form of $\widetilde{\Pi}^{[s,r]}_{H,A}$ establishes the corresponding equivalence on the corresponding categories of $\Gamma$-modules.	
\end{proposition}

\begin{proof}
This is a relative version of the corresponding result in \cite[Lemma 5.6.9]{3KL16} we adapt the corresponding argument here. Indeed the corresponding fully faithfulness comes from the previous proposition the idea to prove the corresponding essential surjectivity comes from writing the module $M$ over the ring taking form of $\widetilde{\Pi}^{[s,r]}_{H,A}$ as the base change from $\varphi^{-k}(\Pi^{[s,r]}_{H,A})$ of a module $M_0$ after the corresponding analog of \cite[Lemma 5.6.8]{3KL16} as above.	Then as in \cite[Lemma 5.6.9]{3KL16} we consider the corresponding norms on the corresponding $M$ (the corresponding differentials) and the base change of $M_0$ (the corresponding differentials) which could be controlled up to some constant which could be further modified to be zero by reducing each time positive amount of constant from the constant represented by the difference of the norms.
\end{proof}

\begin{definition} \mbox{\bf{(After Kedlaya-Liu \cite[Definition 5.7.2]{3KL16})}}
Over the period rings $\Pi_{H,A}$ or $\breve{\Pi}_{H,A}$  (which is denoted by $\triangle$ in this definition) we define the corresponding $(\varphi^a,\Gamma)$-modules over $\triangle$ which are respectively projective to be the corresponding finite projective right $\Gamma$-modules over $\triangle$ with further assigned semilinear action of the operator $\varphi^a$ with the isomorphism defined by using the Frobenius.  
\end{definition}

\begin{definition} \mbox{\bf{(After Kedlaya-Liu \cite[Definition 5.7.2]{3KL16})}}
Over each rings $\triangle=\Pi^r_{H,A},\breve{\Pi}^r_{H,A}$ we define the corresponding projective $(\varphi^a,\Gamma)$-module over any $\triangle$ to be the corresponding finite projective right $\Gamma$-module $M$ over $\triangle$ with additionally endowed semilinear Frobenius action from $\varphi^a$ such that we have the isomorphism $\varphi^{a*}M\overset{\sim}{\rightarrow}M\otimes \square$ where the ring $\square$ is one $\triangle=\Pi^{r/p}_{H,A},\breve{\Pi}^{r/p}_{H,A}$.  
\end{definition}

\begin{definition} \mbox{\bf{(After Kedlaya-Liu \cite[Definition 5.7.2]{3KL16})}}
Again as in \cite[Definition 5.7.2]{3KL16}, we define the corresponding projective $(\varphi^a,\Gamma)$-modules over ring $\Pi^{[s,r]}_{H,A}$ or $\breve{\Pi}^{[s,r]}_{H,A}$ to be the finite projective right $\Gamma$-modules (which will be denoted by $M$) over $\Pi^{[s,r]}_{H,A}$ or $\breve{\Pi}^{[s,r]}_{H,A}$ respectively additionally endowed with semilinear Frobenius action from $\varphi^a$ with the following isomorphisms:
\begin{align}
\varphi^{a*}M\otimes_{\Pi_{H,A}^{[sp^{-ah},rp^{-ah}]}}\Pi_{H,A}^{[s,rp^{-ah}]}\overset{\sim}{\rightarrow}M\otimes_{\Pi_{H,A}^{[s,r]}}\Pi_{H,A}^{[s,rp^{-ah}]}
\end{align}
and
\begin{align}
\varphi^{a*}M\otimes_{\breve{\Pi}_{R,A}^{[sp^{-ah},rp^{-ah}]}}\breve{\Pi}_{R,A}^{[s,rp^{-ah}]}\overset{\sim}{\rightarrow}M\otimes_{\breve{\Pi}_{R,A}^{[s,r]}}\breve{\Pi}_{R,A}^{[s,rp^{-ah}]}.
\end{align}
\end{definition}

\begin{definition} \mbox{\bf{(After Kedlaya-Liu \cite[Definition 5.7.2]{3KL16})}}\\
Over the ring $\Pi^t_{H,A}$ or $\breve{\Pi}^t_{H,A}$ we define a corresponding projective $(\varphi^a,\Gamma)$ bundle to be a family $(M_I)_I$ of finite projective right $\Gamma$-modules over each $\widetilde{\Pi}^I_{H,A}$ carrying the natural Frobenius action coming from the operator $\varphi^a$ such that for any two involved intervals having the relation $I\subset J$ we have:
\begin{displaymath}
M_J\otimes_{\Pi^J_{H,A}}\Pi^I_{H,A}\overset{\sim}{\rightarrow}	M_I
\end{displaymath}
and 
\begin{displaymath}
M_J\otimes_{\breve{\Pi}^J_{H,A}}\breve{\Pi}^I_{H,A}\overset{\sim}{\rightarrow}	M_I
\end{displaymath}
with the obvious cocycle condition. Here we have to propose condition on the intervals that for each $I=[s,r]$ involved we have $s\leq r/p^{ah}$. We require the corresponding topological conditions as we did for the corresponding Frobenius bundles. one can take the corresponding 2-limit in the direct sense to define the corresponding objects over the full Robba rings.
\end{definition}

\begin{definition} \mbox{\bf{(After Kedlaya-Liu \cite[Definition 5.7.2]{3KL16})}}
Over the period rings $\widetilde{\Pi}_{H,A}$ (which is denoted by $\triangle$ in this definition) we define the corresponding $(\varphi^a,\Gamma)$-modules over $\triangle$ which are respectively projective to be the corresponding finite projective right $\Gamma$-modules over $\triangle$ with further assigned semilinear action of the operator $\varphi^a$ with the isomorphism defined by using the Frobenius. 
\end{definition}

\begin{definition} \mbox{\bf{(After Kedlaya-Liu \cite[Definition 5.7.2]{3KL16})}}
Over each ring $\triangle=\widetilde{\Pi}^r_{H,A}$ we define the corresponding projective $(\varphi^a,\Gamma)$-module over any $\triangle$ to be the corresponding finite projective right $\Gamma$-module $M$ over $\triangle$ with additionally endowed semilinear Frobenius action from $\varphi^a$ such that we have the isomorphism $\varphi^{a*}M\overset{\sim}{\rightarrow}M\otimes \square$ where the ring $\square$ is one $\triangle=\widetilde{\Pi}^{r/p}_{H,A}$.

\end{definition}

\begin{definition} \mbox{\bf{(After Kedlaya-Liu \cite[Definition 5.7.2]{3KL16})}}
Again as in \cite[Definition 5.7.2]{3KL16}, we define the corresponding projective $(\varphi^a,\Gamma)$-modules over ring $\widetilde{\Pi}^{[s,r]}_{H,A}$ to be the finite projective right $\Gamma$-modules (which will be denoted by $M$) over $\widetilde{\Pi}^{[s,r]}_{H,A}$ additionally endowed with semilinear Frobenius action from $\varphi^a$ with the following isomorphisms:
\begin{align}
\varphi^{a*}M\otimes_{\widetilde{\Pi}_{H,A}^{[sp^{-ah},rp^{-ah}]}}\widetilde{\Pi}_{H,A}^{[s,rp^{-ah}]}\overset{\sim}{\rightarrow}M\otimes_{\widetilde{\Pi}_{H,A}^{[s,r]}}\widetilde{\Pi}_{H,A}^{[s,rp^{-ah}]}.
\end{align}
\end{definition}

\begin{definition} \mbox{\bf{(After Kedlaya-Liu \cite[Definition 5.7.2]{3KL16})}}\\
Over the ring $\widetilde{\Pi}^t_{H,A}$ we define a corresponding projective $(\varphi^a,\Gamma)$ bundle to be a family $(M_I)_I$ of finite projective right $\Gamma$-modules over each $\widetilde{\Pi}^I_{H,A}$ carrying the natural Frobenius action coming from the operator $\varphi^a$ such that for any two involved intervals having the relation $I\subset J$ we have:
\begin{displaymath}
M_J\otimes_{\widetilde{\Pi}^J_{H,A}}\widetilde{\Pi}^I_{H,A}\overset{\sim}{\rightarrow}	M_I
\end{displaymath}
%and 
%\begin{displaymath}
%M_J\otimes_{\widetilde{\Pi}^J_{H,A}}\widetilde{\Pi}^I_{H,A}\overset{\sim}{\rightarrow}	M_I
%\end{displaymath}
with the obvious cocycle condition. Here we have to propose condition on the intervals that for each $I=[s,r]$ involved we have $s\leq r/p^{ah}$. We require the corresponding topological conditions as we did for the corresponding Frobenius bundles. one can take the corresponding 2-limit in the direct sense to define the corresponding objects over the full Robba rings.
\end{definition}

%\begin{remark}
%In the following, we assume that the corresponding ring $\Pi^{[s,r]}_{R,A}$ is sheafy.	
%\end{remark}

\begin{conjecture} \mbox{\bf{(After Kedlaya-Liu \cite[Theorem 5.7.5]{3KL16})}} \label{conjecture672}
We have now the following categories are equivalence for the corresponding radii $0< s\leq r\leq r_0$ (with the further requirement as in \cite[Theorem 5.7.5]{3KL16} that $s\in (0,r/q]$):\\
1. The category of all the finite projective sheaves over the ring $\widetilde{\Pi}_{\mathrm{Spa}(H_0,H_0^+),A}$, carrying the $\varphi^a$-action;\\
2. The category of all the finite projective sheaves over the ring $\widetilde{\Pi}^r_{\mathrm{Spa}(H_0,H_0^+),A}$, carrying the $\varphi^a$-action;\\
3. The category of all the finite projective sheaves over the ring $\widetilde{\Pi}^{[s,r]}_{\mathrm{Spa}(H_0,H_0^+),A}$, carrying the $\varphi^a$-action;\\
%1b. The category of all the finite projective sheaves over the ring $\widetilde{\Pi}_{\mathrm{Spa}(\overline{H}'_\infty,\overline{H}^+_\infty),A}$, carrying the $\varphi^a$-action;\\
%2b. The category of all the finite projective sheaves over the ring $\widetilde{\Pi}^r_{\mathrm{Spa}(\overline{H}'_\infty,\overline{H}^+_\infty),A}$, carrying the $\varphi^a$-action;\\
%3b. The category of all the finite projective sheaves over the ring $\widetilde{\Pi}^{[s,r]}_{\mathrm{Spa}(\overline{H}'_\infty,\overline{H}^+_\infty),A}$, carrying the $\varphi^a$-action;\\
%4. The categories of all the finite projective quasi-coherent sheaves over corresponding adic Fargues-Fontaine curve in the deformed setting $\mathrm{FF}_{\overline{H}'_\infty,A}$, carrying the corresponding action from the group $\Gamma$ which is assumed to be semilinear and continuous over each section over any affinoid subspace of the whole space which is assumed to be $\Gamma$-invariant;\\
4. The category of all the finite projective modules over the ring $\Pi_{H,A}$, carrying the $(\varphi^a,\Gamma)$-action;\\
5. The category of all the finite projective bundles over the ring $\Pi_{H,A}$, carrying the $(\varphi^a,\Gamma)$-action;\\
6. The category of all the finite projective modules over the ring $\breve{\Pi}_{H,A}$, carrying the $(\varphi^a,\Gamma)$-action;\\
7. The category of all the finite projective bundles over the ring $\breve{\Pi}_{H,A}$, carrying the $(\varphi^a,\Gamma)$-action;\\
8. The category of all the finite projective modules over the ring $\breve{\Pi}^{[s,r]}_{H,A}$, carrying the $(\varphi^a,\Gamma)$-action;\\
9. The category of all the finite projective modules over the ring $\widetilde{\Pi}_{H,A}$, carrying the $(\varphi^a,\Gamma)$-action;\\
10. The category of all the finite projective bundles over the ring $\widetilde{\Pi}_{H,A}$, carrying the $(\varphi^a,\Gamma)$-action;\\
11. The category of all the finite projective modules over the ring $\widetilde{\Pi}^{[s,r]}_{H,A}$, carrying the $(\varphi^a,\Gamma)$-action.
\end{conjecture}

%\begin{proof}
%The corresponding comparisons on the sheaves and bundles in 1-4 and 10-12 are derived in the corresponding context in the perfect setting as what we did in the previous sections. The rest ones are proved exactly the same as \cite[Theorem 5.7.5]{3KL16} by using our development.	
%\end{proof}

\begin{proposition} \mbox{\bf{(After Kedlaya-Liu \cite[Theorem 5.7.5]{3KL16})}}
We have now the corresponding equivalence among categories described as below:\\
1. The category of all the finite projective modules over the ring $\Pi_{H,A}$, carrying the $(\varphi^a,\Gamma)$-action;\\
%2. The categories of all the finite projective bundles over the ring $\Pi_{H,A}$, carrying the $(\varphi^a,\Gamma)$-action;\\
2. The category of all the finite projective modules over the ring $\breve{\Pi}_{H,A}$, carrying the $(\varphi^a,\Gamma)$-action;\\
%4. The categories of all the finite projective bundles over the ring $\breve{\Pi}_{H,A}$, carrying the $(\varphi^a,\Gamma)$-action;\\
3. The category of all the finite projective modules over the ring $\widetilde{\Pi}_{H,A}$, carrying the $(\varphi^a,\Gamma)$-action.

\end{proposition}

\begin{proof}
These are proved exactly the same as \cite[Theorem 5.7.5]{3KL16} by using our development.
\end{proof}

\subsection{Noncommutative Descent}

\indent Motivated by the corresponding comparison between the bundles and modules carrying the corresponding action coming from the Frobenius and the group $\Gamma$ we study some noncommutative descent, which takes the form which is more representation theoretic and topos theoretic. Note that the corresponding localization and descent are very complicated issues within the corresponding noncommutative geometry.

\indent Here we study the corresponding noncommutative descent in the style coming from \cite[Section 1.3]{3KL15}. The corresponding story in our setting will be a noncommutative analog of \cite[Section 1.3]{3KL15}. Also we mention that \cite{3DLLZ1} has already considered some interesting generalization of \cite[Section 1.3]{3KL15} along the other direction within the study of the corresponding logarithmic spaces.

\begin{setting}\mbox{\bf{(Noncommutaive Glueing)}}
Consider the following noncommutative analog of \cite[Section 1.3]{3KL15}. First we have a square diagram taking the following form:
\[
\xymatrix@R+2pc@C+2pc{
\Pi \ar[r]\ar[r]\ar[r]\ar[d]\ar[d]\ar[d] &\Pi_2 \ar[d]\ar[d]\ar[d]\\
\Pi_1 \ar[r]\ar[r]\ar[r] &\Pi_{12},\\
}
\]
which expands to the corresponding short exact sequence:
\[
\xymatrix@R+1pc@C+1pc{
0 \ar[r]\ar[r]\ar[r] &\Pi \ar[r]\ar[r]\ar[r] &\Pi_1\oplus \Pi_2 \ar[r]^{*-*}\ar[r]\ar[r] &\Pi_{12}  \ar[r]\ar[r]\ar[r] &0
}
\]	
in the corresponding category of $\Pi$-modules. The corresponding datum consists also of three right modules $M_1,M_2,M_{12}$ over the corresponding rings $\Pi_1,\Pi_2,\Pi_{12}$ respectively, with the corresponding base change isomorphisms from the module $M_1,M_2$ to the module $M_{12}$. Take the corresponding kernel we have the corresponding exact sequence:
\[
\xymatrix@R+1pc@C+1pc{
0 \ar[r]\ar[r]\ar[r] &M \ar[r]\ar[r]\ar[r] &M_1\oplus M_2 \ar[r]\ar[r]\ar[r] &M_{12}.  
}
\]
We are going to call the corresponding datum in our situation to be \textit{coherent}, \textit{pseudocoherent}, $ finite$, $finitely~presented$, $finite~projective$ if the corresponding modules involved are so over the corresponding rings $\Pi_1,\Pi_2,\Pi_{12}$, which is to say $coherent$, $pseudocoherent$, $ finite$, $finitely~presented$, $finite~projective$.
\end{setting}

%\begin{assumption}
%Assume all the corresponding rings in this section to be unital noetherian rings.	
%\end{assumption}

\begin{lemma}\mbox{\bf{(After Kedlaya-Liu \cite[Lemma 1.3.8]{3KL15})}}\label{lemma677} 
For finite datum defined in the corresponding sense, suppose we further assume that the corresponding map $M\otimes\Pi_1\rightarrow M_1$ is surjective. Then we have that the corresponding maps $M\otimes\Pi_2\rightarrow M_2$ and $f_1-f_2:M_1\oplus M_2\rightarrow M_{12}$ are also being surjective, and one can find a finite submodule of $M$ such that the corresponding base change to the $R_1$ or $R_2$ maps through surjective maps to the corresponding module $M_1$ or $M_2$ respectively.
\end{lemma}

\begin{proof}
The proof in our noncommutative setting is actually parallel to \cite[Lemma 1.3.8]{3KL15}. To be more precise the map $f_1-f_2:M_1\oplus M_2\rightarrow M_{12}$ is then surjective is just because we can derive this from the corresponding surjectivity of $M\otimes\Pi_1\rightarrow M_1$. The surjectivity of $M\otimes\Pi_1\rightarrow M_1$ implies that the corresponding surjectivity of $M\otimes (\Pi_1\oplus \Pi_2)\rightarrow M_{12}$ which factor through the corresponding desired map in our situation. Then after proving this one can show the corresponding surjectivity of the map $M\otimes\Pi_2\rightarrow M_2$. Indeed suppose we start from some element $m\in M_2$ then by taking the corresponding base change to $\Pi_{12}$ we can have corresponding elements $m_1$ and $m_2$ living in the image of the maps:
\begin{align}
g_1:M\otimes \Pi_1\rightarrow M_1,\\
g_2:M\otimes \Pi_2\rightarrow M_2,	
\end{align}
which gives the corresponding element $f_1(m_1)-f_2(m_2)$ which is just the corresponding image of the element $m$ in $\Pi_{12}$. Then look at the element $(m_2,m_2+m)$ in the direct sum, then one can see that the image of the map in our mind will contain actually the element $m_2$ and $m+m_2$ at the same times, which implies the result. Finally the last statement then follows from this.
\end{proof}

\begin{lemma} \mbox{\bf{(After Kedlaya-Liu \cite[Lemma 1.3.9]{3KL15})}} \label{lemma678}
For finitely projective datum, suppose that we have the surjectivity of the map $M\otimes \Pi_1\rightarrow M_1$. Then we have the kernel $M$ is finitely presented and we have the isomorphisms $M\otimes \Pi_1\rightarrow M_1$ and $M\otimes \Pi_2 \rightarrow M_2$.
\end{lemma}

\begin{proof}
By applying the previous lemma one can have the chance to find a submodule $M'$ of $M$ which is finite which admits a covering from some finite free module $G$. Then set $G_1,G_2,G_{12}$ to be the corresponding base change of $G$ to the corresponding rings $\Pi_1,\Pi_2,\Pi_{12}$. Then we have the following commutative diagram:
\[
\xymatrix@R+2pc@C+2pc{
     &0\ar[d]\ar[d]\ar[d] &0\ar[d]\ar[d]\ar[d] &0\ar[d]\ar[d]\ar[d] \\ 
0 \ar[r]\ar[r]\ar[r] &F \ar[r]\ar[r]\ar[r]\ar[d]\ar[d]\ar[d] &F_1\oplus F_2 \ar[d]\ar[d]\ar[d] \ar[r]\ar[r]\ar[r] &F_{12} \ar[d]\ar[d]\ar[d]\\
0 \ar[r]\ar[r]\ar[r] &G \ar[r]\ar[r]\ar[r]\ar[d]\ar[d]\ar[d] &G_1\oplus G_2 \ar[d]\ar[d]\ar[d] \ar[r]\ar[r]\ar[r] &G_{12} \ar[d]\ar[d]\ar[d] \ar[r]\ar[r]\ar[r] \ar[d]\ar[d]\ar[d] &0\\
0 \ar[r]\ar[r]\ar[r] &M \ar[r]\ar[r]\ar[r]  &M_1\oplus M_2  \ar[r]\ar[r]\ar[r] \ar[d]\ar[d]\ar[d] &M_{12} \ar[d]\ar[d]\ar[d] \ar[r]\ar[r]\ar[r] \ar[d]\ar[d]\ar[d] &0  \\
&  &0  &0.
}
\]
Here the module $F,F_1,F_2,F_{12}$ are the corresponding kernel of the corresponding maps in the diagram, namely $G\rightarrow M$, $G_1\rightarrow M_1$, $G_2\rightarrow M_2$, $G_{12}\rightarrow M_{12}$ respectively. Now for each $j=\{1,2\}$ we have the corresponding exact sequence:
\[
\xymatrix@R+1pc@C+1pc{
0 \ar[r]\ar[r]\ar[r] &F_j\otimes_{\Pi_j}\Pi_{12} \ar[r]\ar[r]\ar[r] &G_{12} \ar[r]\ar[r]\ar[r] &M_{12}  \ar[r]\ar[r]\ar[r] &0
}
\]
since the corresponding module $M_j$ ($j\in\{1,2\}$) is finite projective. We then have the situation that $F_j\otimes_{\Pi_j}\Pi_{12}\overset{\sim}{\rightarrow} F_{12}$ for each $j\in \{1,2\}$ and the corresponding modules $F_1,F_2$ are finitely generated and projective. Now apply the previous lemma we have that $F_1\oplus F_2\rightarrow F_{12}$ is then surjective which gives the following commutative diagram:
\[
\xymatrix@R+2pc@C+2pc{
     &0\ar[d]\ar[d]\ar[d] &0\ar[d]\ar[d]\ar[d] &0\ar[d]\ar[d]\ar[d] \\ 
0 \ar[r]\ar[r]\ar[r] &F \ar[r]\ar[r]\ar[r]\ar[d]\ar[d]\ar[d] &F_1\oplus F_2 \ar[d]\ar[d]\ar[d] \ar[r]^{\gamma}\ar[r]\ar[r] &F_{12} \ar[d]\ar[d]\ar[d] \ar[r]\ar[r]\ar[r] \ar[d]\ar[d]\ar[d] &0\\
0 \ar[r]\ar[r]\ar[r] &G \ar[r]\ar[r]\ar[r]\ar[d]\ar[d]\ar[d] &G_1\oplus G_2 \ar[d]\ar[d]\ar[d] \ar[r]^{\beta}\ar[r]\ar[r] &G_{12} \ar[d]\ar[d]\ar[d] \ar[r]\ar[r]\ar[r] \ar[d]\ar[d]\ar[d] &0\\
0 \ar[r]\ar[r]\ar[r] &M \ar[r]\ar[r]\ar[r]  &M_1\oplus M_2  \ar[r]^{\alpha}\ar[r]\ar[r] \ar[d]\ar[d]\ar[d] &M_{12} \ar[d]\ar[d]\ar[d] \ar[r]\ar[r]\ar[r] \ar[d]\ar[d]\ar[d] &0  \\
&  &0  &0.
}
\]
Then we claim that then the map $G\rightarrow M$ is surjective. Indeed this is by direct diagram chasing. First take any element in $M_1\oplus M_2$, which is denoted by $(m_1,m_2)$ for which we assume that $(m_1,m_2)\in \mathrm{Ker}\alpha$, then take any element in $(f_1,f_2)\in G_1\oplus G_2$ lifting this element. Let $\overline{f_1-f_2}$ be the image of $(f_1,f_2)$ under $\beta$. By the commutativity of the diagram we have $\overline{f_1-f_2}$ dies in $M_{12}$, so we can find $g_{12}\in F_{12}$ whose image in $G_{12}$ is $\overline{f_1-f_2}$. Then take any element element $(g_1,g_2)$ in $F_1\oplus F_2$ which lifts $f_{12}$. Then since the diagram is commutative we then have that $(g_1,g_2)$ takes image in $G_1\oplus G_2$ which must be living in same equivalence class with $(f_1,f_2)$ with respect to $G$, so we have $(f_1,f_2)=(g_1,g_2)\oplus f$ where $f\in G$. But then we have $f$ which is an element in the kernel of $\beta$ maps to $(m_1,m_2)$ which finishes the proof of the claim. Then repeat this we can have the chance to derive the finiteness of the module $F$. Then to finish we look at the corresponding commutative diagram:
\[
\xymatrix@R+2pc@C+2pc{
 &\mathrm{Ker}(G\rightarrow M)\otimes \Pi_i \ar[r]\ar[r]\ar[r]\ar[d]\ar[d]\ar[d] &G_i \ar[d]\ar[d]\ar[d] \ar[r]\ar[r]\ar[r] &M\otimes \Pi_i \ar[d]\ar[d]\ar[d] \ar[r]\ar[r]\ar[r] \ar[d]\ar[d]\ar[d] &0\\
0 \ar[r]\ar[r]\ar[r] &\mathrm{Ker}(G_i\rightarrow M_i) \ar[r]\ar[r]\ar[r]&G_i  \ar[r]\ar[r]\ar[r] &M_i  \ar[r]\ar[r]\ar[r]  &0,
}
\]
which implies that the corresponding map from $M\otimes \Pi_i$ to $M_i$ is injective as well by five lemma (note that the corresponding map $\mathrm{Ker}(G\rightarrow M)\otimes \Pi_i$ to $\mathrm{Ker}(G_i\rightarrow M_i)$ is surjective). So we have the corresponding desired isomorphisms.	
\end{proof}

\indent The following is a similar noncommutative version of \cite[Lemma A.3]{3DLLZ1}:

\begin{lemma}\mbox{\bf{(After Kedlaya-Liu \cite[Lemma 1.3.9]{3KL15})}} \label{lemma679}
For finitely presented datum, suppose that we have the surjectivity of the map $M\otimes \Pi_1\rightarrow M_1$. And suppose we have the corresponding additional requirement that $\Pi_1\rightarrow \Pi_{12}$ and $\Pi_2\rightarrow \Pi_{12}$ are now assumed to be flat. Then we have the kernel $M$ is finitely presented and we have the isomorphisms $M\otimes \Pi_1\rightarrow M_1$ and $M\otimes \Pi_2 \rightarrow M_2$.
\end{lemma}

\begin{proof}
By applying the previous lemma one can have the chance to find a submodule $M'$ of $M$ which is finite which admits a covering from some finite free module $G$. Then set $G_1,G_2,G_{12}$ to be the corresponding base change of $G$ to the corresponding rings $\Pi_1,\Pi_2,\Pi_{12}$. Then we have the following commutative diagram:
\[
\xymatrix@R+2pc@C+2pc{
     &0\ar[d]\ar[d]\ar[d] &0\ar[d]\ar[d]\ar[d] &0\ar[d]\ar[d]\ar[d] \\ 
0 \ar[r]\ar[r]\ar[r] &F \ar[r]\ar[r]\ar[r]\ar[d]\ar[d]\ar[d] &F_1\oplus F_2 \ar[d]\ar[d]\ar[d] \ar[r]\ar[r]\ar[r] &F_{12} \ar[d]\ar[d]\ar[d]\\
0 \ar[r]\ar[r]\ar[r] &G \ar[r]\ar[r]\ar[r]\ar[d]\ar[d]\ar[d] &G_1\oplus G_2 \ar[d]\ar[d]\ar[d] \ar[r]\ar[r]\ar[r] &G_{12} \ar[d]\ar[d]\ar[d] \ar[r]\ar[r]\ar[r] \ar[d]\ar[d]\ar[d] &0\\
0 \ar[r]\ar[r]\ar[r] &M \ar[r]\ar[r]\ar[r]  &M_1\oplus M_2  \ar[r]\ar[r]\ar[r] \ar[d]\ar[d]\ar[d] &M_{12} \ar[d]\ar[d]\ar[d] \ar[r]\ar[r]\ar[r] \ar[d]\ar[d]\ar[d] &0  \\
&  &0  &0.
}
\]
Here the module $F,F_1,F_2,F_{12}$ are the corresponding kernel of the corresponding maps in the diagram, namely $G\rightarrow M$, $G_1\rightarrow M_1$, $G_2\rightarrow M_2$, $G_{12}\rightarrow M_{12}$ respectively. Now for each $j=\{1,2\}$ we have the corresponding exact sequence:
\[
\xymatrix@R+1pc@C+1pc{
0 \ar[r]\ar[r]\ar[r] &F_j\otimes_{\Pi_j}\Pi_{12} \ar[r]\ar[r]\ar[r] &G_{12} \ar[r]\ar[r]\ar[r] &M_{12}  \ar[r]\ar[r]\ar[r] &0
}
\]
since we have the corresponding additional requirement that $\Pi_1\rightarrow \Pi_{12}$ and $\Pi_2\rightarrow \Pi_{12}$ are now assumed to be flat. We then have the situation that $F_j\otimes_{\Pi_j}\Pi_{12}\overset{\sim}{\rightarrow} F_{12}$ for each $j\in \{1,2\}$. The corresponding modules $F_1,F_2$ are finitely generated by assumption. Now applying the previous lemma we have that $F_1\oplus F_2\rightarrow F_{12}$ is then surjective which gives the following commutative diagram:
\[
\xymatrix@R+2pc@C+2pc{
     &0\ar[d]\ar[d]\ar[d] &0\ar[d]\ar[d]\ar[d] &0\ar[d]\ar[d]\ar[d] \\ 
0 \ar[r]\ar[r]\ar[r] &F \ar[r]\ar[r]\ar[r]\ar[d]\ar[d]\ar[d] &F_1\oplus F_2 \ar[d]\ar[d]\ar[d] \ar[r]^{\gamma}\ar[r]\ar[r] &F_{12} \ar[d]\ar[d]\ar[d] \ar[r]\ar[r]\ar[r] \ar[d]\ar[d]\ar[d] &0\\
0 \ar[r]\ar[r]\ar[r] &G \ar[r]\ar[r]\ar[r]\ar[d]\ar[d]\ar[d] &G_1\oplus G_2 \ar[d]\ar[d]\ar[d] \ar[r]^{\beta}\ar[r]\ar[r] &G_{12} \ar[d]\ar[d]\ar[d] \ar[r]\ar[r]\ar[r] \ar[d]\ar[d]\ar[d] &0\\
0 \ar[r]\ar[r]\ar[r] &M \ar[r]\ar[r]\ar[r]  &M_1\oplus M_2  \ar[r]^{\alpha}\ar[r]\ar[r] \ar[d]\ar[d]\ar[d] &M_{12} \ar[d]\ar[d]\ar[d] \ar[r]\ar[r]\ar[r] \ar[d]\ar[d]\ar[d] &0  \\
&  &0  &0.
}
\]
Then we claim that then the map $G\rightarrow M$ is surjective. Indeed this is by direct diagram chasing. First take any element in $M_1\oplus M_2$, which is denoted by $(m_1,m_2)$ for which we assume that $(m_1,m_2)\in \mathrm{Ker}\alpha$, then take any element in $(f_1,f_2)\in G_1\oplus G_2$ lifting this element. Let $\overline{f_1-f_2}$ be the image of $(f_1,f_2)$ under $\beta$. By the commutativity of the diagram we have $\overline{f_1-f_2}$ dies in $M_{12}$, so we can find $g_{12}\in F_{12}$ whose image in $G_{12}$ is $\overline{f_1-f_2}$. Then take any element element $(g_1,g_2)$ in $F_1\oplus F_2$ which lifts $f_{12}$. Then since the diagram is commutative we then have that $(g_1,g_2)$ takes image in $G_1\oplus G_2$ which must be living in same equivalence class with $(f_1,f_2)$ with respect to $G$, so we have $(f_1,f_2)=(g_1,g_2)\oplus f$ where $f\in G$. But then we have $f$ which is an element in the kernel of $\beta$ maps to $(m_1,m_2)$ which finishes the proof of the claim. Then repeat this we can have the chance to derive the finiteness of the module $F$. Then to finish we look at the corresponding commutative diagram:
\[
\xymatrix@R+2pc@C+2pc{
 &\mathrm{Ker}(G\rightarrow M)\otimes \Pi_i \ar[r]\ar[r]\ar[r]\ar[d]\ar[d]\ar[d] &G_i \ar[d]\ar[d]\ar[d] \ar[r]\ar[r]\ar[r] &M\otimes \Pi_i \ar[d]\ar[d]\ar[d] \ar[r]\ar[r]\ar[r] \ar[d]\ar[d]\ar[d] &0\\
0 \ar[r]\ar[r]\ar[r] &\mathrm{Ker}(G_i\rightarrow M_i) \ar[r]\ar[r]\ar[r]&G_i  \ar[r]\ar[r]\ar[r] &M_i  \ar[r]\ar[r]\ar[r]  &0,
}
\]
which implies that the corresponding map from $M\otimes \Pi_i$ to $M_i$ is injective as well by five lemma (note that the corresponding map $\mathrm{Ker}(G\rightarrow M)\otimes \Pi_i$ to $\mathrm{Ker}(G_i\rightarrow M_i)$ is surjective). So we have the corresponding desired isomorphisms.	
	
\end{proof}

\indent Now we consider the corresponding Banach structures as well. In this case we consider the following glueing process. First we have the following analog of \cite[Lemma 2.7.1]{3KL15}:

\begin{lemma}\mbox{\bf{(After Kedlaya-Liu \cite[Lemma 2.7.1]{3KL15})}}
Suppose we have the following two homomorphism of Banach rings (which are not necessarily commutative in our situation):
\begin{align}
f_1: \Pi_1\rightarrow \Gamma\\
f_2: \Pi_2\rightarrow \Gamma	
\end{align}
such that they are bounded and we have the following sum map induced which is surjective in the strict sense:
\begin{displaymath}
\Pi_1\oplus \Pi_2\rightarrow \Gamma.	
\end{displaymath}
Then one can find some constant $c$ with the fact that for any positive integer $n\geq 1$, we have the corresponding decomposition of any invertible matrix $I\in \Gamma^{n\times n}$ such that $\|I-1\|\leq c$ in the following form:
\begin{align}
f_1(I_1)f_2(I_2)	
\end{align}
with invertible matrices $I_1,I_2$ living in $\Pi_1^{n\times n}$ and $\Pi_2^{n\times n}$	respectively.
\end{lemma}

\begin{proof}
Following 	\cite[Lemma 2.7.1]{3KL15}, we consider the corresponding lifting process. First by the conditions of the lemma we have that one can lift $I-1$ to a pair matrices $A,B$ living in $\mathrm{GL}(\Pi_1)$ and $\mathrm{GL}(\Pi_2)$ with the corresponding estimate on the entries:
\begin{displaymath}
\|A_{ij}\|\leq c_0 \|(I-1)_{ij}\|\\
\|B_{ij}\|\leq c_0 \|(I-1)_{ij}\|\\	
\end{displaymath}
for some constant $c_0$. Then consider the corresponding iteration process, namely from what we have set $I'$ to be $f_1(1-A)If_2(1-B)$ such that we have the corresponding estimate in the form of $\|I'-1\|\leq d\|I-1\|^2$. Then to finish run the construction process above step by step, we will be done.
\end{proof}

\indent Now in the following setting, we set all rings involved to be Banach:

\begin{setting}\mbox{\bf{(Noncommutaive Glueing)}}
Consider the following noncommutative analog of \cite[Section 1.3]{3KL15}. First we have a square diagram taking the following form which is commutative:
\[
\xymatrix@R+2pc@C+2pc{
\Pi \ar[r]\ar[r]\ar[r]\ar[d]\ar[d]\ar[d] &\Pi_2 \ar[d]\ar[d]\ar[d]\\
\Pi_1 \ar[r]\ar[r]\ar[r] &\Pi_{12},\\
}
\]
which expands to the corresponding short exact sequence:
\[
\xymatrix@R+1pc@C+1pc{
0 \ar[r]\ar[r]\ar[r] &\Pi \ar[r]\ar[r]\ar[r] &\Pi_1\oplus \Pi_2 \ar[r]^{*-*}\ar[r]\ar[r] &\Pi_{12}  \ar[r]\ar[r]\ar[r] &0
}
\]	
in the corresponding category of $\Pi$-modules. The corresponding datum consists also of three right modules $M_1,M_2,M_{12}$ over the corresponding rings $\Pi_1,\Pi_2,\Pi_{12}$ respectively, with the corresponding base change isomorphisms from the module $M_1,M_2$ to the module $M_{12}$. Take the corresponding kernel we have the corresponding exact sequence:
\[
\xymatrix@R+1pc@C+1pc{
0 \ar[r]\ar[r]\ar[r] &M \ar[r]\ar[r]\ar[r] &M_1\oplus M_2 \ar[r]\ar[r]\ar[r] &M_{12}.  
}
\]
We are going to call the corresponding datum in our situation to be \textit{coherent}, \textit{pseudocoherent}, $ finite$, $finitely~presented$, $finite~projective$ if the corresponding modules involved are so over the corresponding rings $\Pi_1,\Pi_2,\Pi_{12}$, which is to say $coherent$, $pseudocoherent$, $finite$, $finitely~presented$, $finite~projective$, and we require that the corresponding exact sequence:
\[
\xymatrix@R+1pc@C+1pc{
0 \ar[r]\ar[r]\ar[r] &\Pi \ar[r]\ar[r]\ar[r] &\Pi_1\oplus \Pi_2 \ar[r]^{*-*}\ar[r]\ar[r] &\Pi_{12}  \ar[r]\ar[r]\ar[r] &0
}
\]
is surjective in the strict sense, and we require that the map:
\begin{displaymath}
\Pi_2\rightarrow \Pi_{12}	
\end{displaymath}
has the image which is assumed in our situation to be dense.
\end{setting}

\begin{lemma} \mbox{\bf{(After Kedlaya-Liu \cite[Lemma 2.7.3]{3KL15})}}
Consider a finite glueing datum as above, we then have the fact that 
$M\otimes \Pi_1 \rightarrow M_{1}$ and $M\otimes \Pi_2 \rightarrow M_{2}$ are surjective and the map $M_1\oplus M_2\rightarrow M_{12}$ is also surjective.
\end{lemma}

\begin{proof}
This is a noncommutative version of \cite[Lemma 2.7.3]{3KL15}. We just need to check that in our situation the condition of \cref{lemma677} holds. Now consider the bases of $M_1$ and $M_2$ respectively:
\begin{align}
e_1,e_2,...,e_m,\\
f_1,f_2,...,f_m.	
\end{align}
And we consider the corresponding expansion with respect to these bases after taking the image under the maps $s_1:M_1\otimes \Pi_{12}\rightarrow M_{12}$ and $s_2:M_2\otimes \Pi_{12}\rightarrow M_{12}$:
\begin{displaymath}
s_2(f_j)=\sum_{i} I_{ij} s_1(e_i), s_1(f_j)=\sum_{i} J_{ij} s_2(f_i).	
\end{displaymath}
Then we have from the assumption on the dense image some matrix $J'$ such that $I(J'-J)$ could be satisfying the condition $\|I(J'-J)\|\leq c$ where $c$ is the constant in the previous lemma. Then by using the previous lemma we have the following decomposition:
\begin{displaymath}
1+I(J'-J)=K_1K_2^{-1}	
\end{displaymath}
with $K_1\in \mathrm{GL}(\Pi_1),K_2\in \mathrm{GL}(\Pi_2)$. Then we consider the element:
\begin{displaymath}
(a_{j},b_j):=(\sum_{i} K_{1,ij} e_i,\sum_{i} (J'K_2)_{ij} f_i).	
\end{displaymath}
We compute as in the following:
\begin{align}
s_1(a_{j})-s_2(b_j)&=\sum_{i} K_{1,ij} s_1(e_i)-	\sum_{i} J'K_{2,ij} s_2(f_i))\\
                   &=\sum_i K_{1,ij} s_1(e_i)-(IJ'K_{2})_{ij}s_1(e_i))\\
                   &=\sum_i ((1-IJ)K_2)_{ij} s_1(e_i))\\
                   &=0.
\end{align}
Therefore the element we defined is now living in the kernel $M$. Consider just the first component we can see that the map $M\otimes \Pi_1\rightarrow M_1$ is now surjective.	
\end{proof}

%\begin{lemma}\mbox{\bf{(After Kedlaya-Liu \cite[Lemma 2.7.4]{3KL15})}}
%Assume we are woring over noetherian rings. For finitely presented datum, we then in this situation have the surjectivity of the map $M\otimes \Pi_1\rightarrow M_1$. Therefore we have the kernel $M$ is finitely presented and we have the isomorphisms $M\otimes \Pi_1\rightarrow M_1$ and $M\otimes \Pi_2 \rightarrow M_2$.
%\end{lemma}
%
%
%\begin{proof}
%This is by the previous lemma and \cref{lemma677}.	
%\end{proof}

\begin{lemma}\mbox{\bf{(After Kedlaya-Liu \cite[Lemma 2.7.4]{3KL15})}}
For finitely projective datum, we then in our situation have the surjectivity of the map $M\otimes \Pi_1\rightarrow M_1$. Then we have the kernel $M$ is finitely presented and we have the isomorphisms $M\otimes \Pi_1\rightarrow M_1$ and $M\otimes \Pi_2 \rightarrow M_2$.
\end{lemma}

\begin{proof}
This is by the previous lemma and \cref{lemma678}.	
\end{proof}

\indent The following is a similar noncommutative version of \cite[Lemma A.7]{3DLLZ1}:

\begin{lemma}\mbox{\bf{(After Kedlaya-Liu \cite[Lemma 2.7.4]{3KL15})}}
For finitely presented datum, we then in our situation have the surjectivity of the map $M\otimes \Pi_1\rightarrow M_1$. Furthermore suppose we are in the situation where $\Pi_1\rightarrow \Pi_{12}$ and $\Pi_2\rightarrow \Pi_{12}$ are simultaneously assumed to be flat. Then we have the kernel $M$ is finitely presented and we have the isomorphisms $M\otimes \Pi_1\rightarrow M_1$ and $M\otimes \Pi_2 \rightarrow M_2$.
\end{lemma}

\begin{proof}
See lemma above and \cref{lemma679}.	
\end{proof}

%\indent Our hope is merely to extend the corresponding comparison of pseudocoherent hodge structures to the noncommutative situation, definitely in the non-noetherian situation this will be complicated.

\begin{definition} \mbox{\bf{(After Kedlaya-Liu \cite[Definition 5.7.2]{3KL16})}}
Over the period rings $\Pi_{H,A}$ (which is denoted by $\triangle$ in this definition) we define the corresponding $(\varphi^a,\Gamma)$-modules over $\triangle$ which are pseudocoherent or fpd to be the corresponding pseudocoherent or fpd right $\Gamma$-modules over $\triangle$ with further assigned semilinear action of the operator $\varphi^a$ with the isomorphism defined by using the Frobenius. We also require that the modules are complete for the natural topology involved in our situation and for any module over $\Pi_{H,A}$ to be some base change of some module over $\Pi^r_{H,A}$ (which will be defined in the following).
\end{definition}

\begin{definition} \mbox{\bf{(After Kedlaya-Liu \cite[Definition 5.7.2]{3KL16})}}
Over each ring $\triangle=\Pi^r_{H,A}$ we define the corresponding pseudocoherent or fpd $(\varphi^a,\Gamma)$-module over any $\triangle$ to be the corresponding pseudocoherent or fpd right $\Gamma$-module $M$ over $\triangle$ with additionally endowed semilinear Frobenius action from $\varphi^a$ such that we have the isomorphism $\varphi^{a*}M\overset{\sim}{\rightarrow}M\otimes \square$ where the ring $\square$ is one $\triangle=\Pi^{r/p}_{H,A}$. Also as in \cite[Definition 5.7.2]{3KL16} we assume that the module over $\Pi^r_{H,A}$ is then complete for the natural topology and the corresponding base change to $\Pi^I_{H,A}$ for any interval which is assumed to be closed $I\subset [0,r)$ gives rise to a module over $\Pi^I_{H,A}$ with specified conditions which will be specified below.

\end{definition}

\begin{definition} \mbox{\bf{(After Kedlaya-Liu \cite[Definition 5.7.2]{3KL16})}}
Again as in \cite[Definition 5.7.2]{3KL16}, we define the corresponding pseudocoherent and fpd $(\varphi^a,\Gamma)$-modules over ring $\Pi^{[s,r]}_{H,A}$ to be the pseudocoherent and fpd right $\Gamma$-modules (which will be denoted by $M$) over $\Pi^{[s,r]}_{H,A}$ additionally endowed with semilinear Frobenius action from $\varphi^a$ with the following isomorphisms:
\begin{align}
\varphi^{a*}M\otimes_{\Pi_{H,A}^{[sp^{-ah},rp^{-ah}]}}\Pi_{H,A}^{[s,rp^{-ah}]}\overset{\sim}{\rightarrow}M\otimes_{\Pi_{H,A}^{[s,r]}}\Pi_{H,A}^{[s,rp^{-ah}]}.
\end{align}
We now assume that the modules are complete with respect to the natural topology. 
\end{definition}

\begin{definition} \mbox{\bf{(After Kedlaya-Liu \cite[Definition 5.7.2]{3KL16})}}\\
Over the ring $\Pi^t_{H,A}$ we define a corresponding pseudocoherent and fpd $(\varphi^a,\Gamma)$ bundle to be a family $(M_I)_I$ of pseudocoherent and fpd right $\Gamma$-modules over each $\Pi^I_{H,A}$ carrying the natural Frobenius action coming from the operator $\varphi^a$ such that for any two involved intervals having the relation $I\subset J$ we have:
\begin{displaymath}
M_J\otimes_{\Pi^J_{H,A}}\Pi^I_{H,A}\overset{\sim}{\rightarrow}	M_I
\end{displaymath}
%and 
%\begin{displaymath}
%M_J\otimes_{\Pi^J_{H,A}}\Pi^I_{H,A}\overset{\sim}{\rightarrow}	M_I
%\end{displaymath}
with the obvious cocycle condition. Here we have to propose condition on the intervals that for each $I=[s,r]$ involved we have $s\leq r/p^{ah}$. We require the corresponding topological conditions as above. one can take the corresponding 2-limit in the direct sense to define the corresponding objects over the full Robba rings.
\end{definition}

\begin{remark}
We have the obvious notion of coherent objects in the noetherian setting. Also note that in the noetherian setting we do not have to worry about the corresponding completion issue, which to say the finiteness is all enough. In our current situation and furthermore also in the noetherian setting, by using the corresponding open mapping, one can also mimick the proof in \cite[3.7.3/2,3]{3BGR} to prove this by considering the corresponding presentation and take the corresponding kernel.
\end{remark}

\begin{proposition}
Keeping \cref{assumption642} and under the assumption of \cref{conjecture672}, consider the following categories:\\
1. The category of all the pseudocoherent bundles over the ring $\Pi_{H,A}$, carrying the $(\varphi^a,\Gamma)$-action;\\
%%6. The category of all the pseudocoherent modules over the ring $\breve{\Pi}_{H,A}$, carrying the $(\varphi^a,\Gamma)$-action;\\
%%2. The category of all the pseudocoherent bundles over the ring $\breve{\Pi}_{H,A}$, carrying the $(\varphi^a,\Gamma)$-action;\\
%%3. The category of all the pseudocoherent modules over the ring $\breve{\Pi}^{[s,r]}_{H,A}$, carrying the $(\varphi^a,\Gamma)$-action;\\
%%9. The category of all the pseudocoherent modules over the ring $\widetilde{\Pi}_{H,A}$, carrying the $(\varphi^a,\Gamma)$-action;\\
%%4. The category of all the pseudocoherent bundles over the ring $\widetilde{\Pi}_{H,A}$, carrying the $(\varphi^a,\Gamma)$-action;\\
%%5. The category of all the pseudocoherent modules over the ring $\widetilde{\Pi}^{[s,r]}_{H,A}$, carrying the $(\varphi^a,\Gamma)$-action.\\	
2. The category of all the pseudocoherent modules over the ring $\Pi^{[s,r]}_{H,A}$, carrying the $(\varphi^a,\Gamma)$-action.\\
%%3. The category of all the pseudocoherent modules over the ring $\breve{\Pi}^{[s,r]}_{H,A}$, carrying the $(\varphi^a,\Gamma)$-action;\\
%%4. The category of all the pseudocoherent modules over the ring $\widetilde{\Pi}^{[s,r]}_{H,A}$, carrying the $(\varphi^a,\Gamma)$-action.\\
\indent Then in our situation they are equivalent, if we have for any intervals permitted $I\subset I'$ we have $\Pi^{I'}_{H,A}\rightarrow \Pi^I_{H,A}$ is flat.	
\end{proposition}

\begin{proof}
First note that we are working over noetherian rings. The functor from the category of the Frobenius bundles to that of the Frobenius modules is just projection of the family of Frobenius modules to the corresponding Frobenius module with respect to the specific interval. In our our situation we have 
\[
\xymatrix@R+1pc@C+1pc{
0 \ar[r]\ar[r]\ar[r] &\Pi \ar[r]\ar[r]\ar[r] &\Pi_1\oplus \Pi_2 \ar[r]^{*-*}\ar[r]\ar[r] &\Pi_{12}  \ar[r]\ar[r]\ar[r] &0
}
\]
is exact in the strict sense due to the existence of Schauder basis, and we have the image of $\Pi_2\rightarrow \Pi_{12}$ is dense as well in the same way. Therefore we can see that the corresponding categories are equivalent since we can show the corresponding essential surjectivity by using the Frobenius action and the glueing process established above to reach any closed interval involved.	
\end{proof}

%\begin{remark}
%For instance if $A$ is commutative, then this will give a new proof of the corresponding equivalence involved which we proved before.	
%\end{remark}

\

This chapter is based on the following paper, where the author of this dissertation is the main author:
\begin{itemize}
\item Tong, Xin. "Hodge-Iwasawa Theory II." arXiv preprint arXiv:2010.06093 (2020).
\end{itemize}

\newpage\chapter{Analytic Geometry and Hodge-Frobenius Structure}

\newpage\section{Introduction}

\subsection{Higher Dimensional Analytic Geometry}

Higher dimensional analytic geometry in the nonarchimedean setting and corresponding $p$-adic analysis received very rapid development in the past decades, dated back to the last century to some point. Its complex analog, which is to say the complex analysis in higher dimension of several variables or the study of complex analytic varieties could be dated back even earlier. For instance searching for the automorphism groups of higher dimensional domains is still attracting extensive attention. On the other hand in the $p$-adic setting, considering higher dimensional analytic domains will be a very natural generalization (especially in the rigid analytic geometry, or more general analytic geometry). However, the approaches do not generalize very smoothly to the higher dimension, even, one will immediately need to find new approaches. This is actually already clear in the complex setting, for instance one may believe that complex analysis of several variables over complex polydiscs is natural generalization of one variable situation, for instance the Cauchy integration formula and maximal value issues. However in the $p$-adic setting, this might be more complicated as discussed in the following.\\

\indent The famous $p$-adic local monodromy theorem, proved by Tsuzuki in the unit-root setting, and then proved by Andr\'e-Kedlaya-Mebkhout \cite{4Ked1}, \cite{4And1} and \cite{4Me} in the full general setting, is already a very difficult theorem over one-dimensional $p$-adic annulus, which roughly speaks that over such annulus, each finite locally free sheaf with connection and Frobenius (not necessarily absolute) structure could be made simpler with respect to the connection structure. For instance in the unit-root setting, the picture is even simpler, since then we only need to consider Galois representations with finite monodromy which is to say potentially unramified somewhat: we have an equivalence between the category of unit-root overconvergent isocrystals and the one of all the potentially unramified \'etale $\mathbb{Q}_p$-local systems, over some space. Actually Kedlaya's approach established in \cite{4Ked1} used very deep theorem on the slope filtration over one-dimensional $p$-adic annulus to reduce to Tsuzuki's approach.\\

\indent To generalize to the higher dimensional situation, we have many ways to do so which are encoded in the following projects. First one could consider the picture in the relative $p$-adic Hodge theory, for instance established in \cite{4KL1} and \cite{4KL2}, or the multidimensional program of \cite{4CKZ18} and \cite{4PZ19}. Also we would like to mention some relativization programs in \cite{4BC1}, \cite{4Ked2}, \cite{4Liu1} where some attempt to generalize the slope filtration theorem and the local monodromy theorem to the relative Frobenius situation had been successful. The applications from \cite{4Ked1} to the rigid cohomology and isocrystals have already considered some sort of multidimensional construction as in \cite{4Ked3},\cite{4Ked4},\cite{4Ked5},\cite{4Ked6},\cite{4Ked7},\cite{4Ked8}. However, one could still only get information deduced from the monodromy theorem or the slope filtration theorem essentially from the one mentioned above. The theory of partial differential equations have been regarded as non-trivial and significant generalization of the story relevant above. To be more precise, these are reflected in the following work: \cite{4Xiao1}, \cite{4Xiao2} and \cite{4KX1}, \cite{4Ked10}, \cite{4Ked11}, \cite{4Ked12}.\\

\indent That being all said, we will only focus on mixed-characteristic situation in this paper. Here let us do some summarization in a uniform way on the known pictures which one might want to follow being local or global, being \'etale or non-\'etale:

\begin{setting}
Mixed-characteristic local \'etale: Berger's theorem for potentially semi-stable Galois representations \cite{4Ber1}.	
\end{setting}

\begin{setting}
Mixed-characteristic local non-\'etale: Again Berger's theorem for $(\varphi,\Gamma)$-modules.
\end{setting}

\begin{setting}
Mixed-characteristic global \'etale: Relative $p$-adic Hodge theory after Faltings, Andreatta, Iovita, Brinon, Scholze, Kedlaya, Liu and so on, in \cite{4Fal}, \cite{4A1}, \cite{4AI1}, \cite{4AI2}, \cite{4AB1}, \cite{4AB2}, \cite{4AB3}, \cite{4Sch1}, \cite{4KL1} and \cite{4KL2}.
\end{setting}

\begin{setting}
Mixed-characteristic global non-\'etale: Relative $p$-adic Hodge theory in the style of Kedlaya-Liu \cite{4KL1} and \cite{4KL2}.	
\end{setting}

\indent In this paper, we first consider the program of \cite{4CKZ18} and \cite{4PZ19}, and we will work in the generality of \cite{4KPX}.  The generalization itself is already interesting enough and actually we need to address some fundamental issues around. The multidimensionalization we studied here represents some similarity to the globalization and multidimensionalization mentioned above in the relative $p$-adic Hodge theory and some globalization of the study of isocrystals in the style of Kedlaya's theorem on Shiho's conjecture. Actually we consider slightly different actions within the module structures, which is to say that we will have multi-Frobenius actions, multi-Lie group actions, multi-differentials. 

\subsection{Results and Applications}

\indent Here is our first main result, while the corresponding proof will be given in \cref{section4.1}:

\begin{theorem} \label{4theorem4.1.5}
With the same assumptions as in \cref{assumption1} as below, we set $M$ to be the global section of a locally free coherent sheaf of module over the sheaf of the Robba ring $\Pi_{\mathrm{an},\mathrm{con},I,X}(\pi_{K_I})$ over an affinoid domain $X$ attached to an affinoid algebra $A$ in the rigid analytic geometry, carrying mutually commuting actions of $(\varphi_I,\Gamma_{K_I})$. Then we have the corresponding Herr complex $C^\bullet_{\varphi_I,\Gamma_{K_I}}(M)$ lives then in the derived category $\mathbb{D}^\flat_{\mathrm{perf}}(A)$, and we have that the corresponding complex $C^\bullet_{\psi_I}(M)$ has all the cohomology groups which are coadmissible over $\Pi_{\mathrm{an},\infty,I,A}(\Gamma_{K_I})$. Moreover when $A=\mathbb{Q}_p$, $C^\bullet_{\psi_I}(M)$ lives in $\mathbb{D}^\flat_{\mathrm{perf}}(\Pi_{\mathrm{an},\infty,I,A}(\Gamma_{K_I}))$.
\end{theorem}

\indent Another motivation in our situation comes from the equivariant Iwasawa theory proposed by Kedlaya-Pottharst in \cite{4KP1}. Here the idea will be working over more general distribution algebra after Berthelot, Schneider and Teitelbaum in the following sense. The idea is from the following observation. Now suppose $G$ is a nice $p$-adic Lie group for instance as above we considered the product of $\mathbb{Z}_p$ or $\mathbb{Z}_p^\times$. Then we consider the corresponding distribution algebra which we will use the notation $\Pi_{\mathrm{an},\infty,A}(G)$ to denote this ring (note that here we have already considered the even more general relative setting). Now consider $L/K$ a $p$-adic Lie extension. Then we have that one could consider the projection $G_{K}\rightarrow G(L/K)$ to define the corresponding Galois representation with coefficients in $\Pi_{\mathrm{an},\infty,A}(G)$, which then gives rise to a $(\varphi,\Gamma_K)$-module over $\Pi_{\mathrm{an},\infty,A}(G)$. Then we take the external tensor product of any $(\varphi,\Gamma)$-module $M$ with this $p$-adic Lie deformation $\mathbf{DfmLie}$ we will get a module $M\boxtimes \mathbf{DfmLie}$ which is actually $(\varphi,\Gamma)$-module over $\Pi_{\mathrm{an},\infty,A}(G)$. Then we are going to define the Herr complex of this module, which will be called the $G$-equivariant Iwasawa cohomology. It is natural to investigate the fundamental properties of this cohomology which is over $\mathrm{Max}A\times \mathrm{Max}\Pi_{\mathrm{an},\infty}(G)$. For instance over $\mathbb{Z}_p^{\times,d}$ as above, we will have to study the cohomology of the big $(\varphi,\Gamma)$-module over the product of the ring $\Pi_{\mathrm{an},\mathrm{con},A}(\pi_K)$ and the ring $\Pi_{\mathrm{an},\infty}(G)$. Therefore through the previous theorem we have the following coherence over more general Fr\'echet-Stein algebras, while the corresponding proof will be given in \cref{section4.1}:

\begin{theorem}   \label{4theorem4.1.6}
Let $X_\infty(G)$ be the quasi-Stein rigid analytic space attached to the $p$-adic Lie group $G:=\Gamma_{K_J}$ for some set $J$ possibly different from $I$ as in \cref{4setting3.1}. Let $A$ be an affinoid algebra over $\mathbb{Q}_p$. Suppose $\mathcal{F}$ is a locally free coherent sheaf over the sheaf of Robba ring $\Pi_{\mathrm{an},\mathrm{con},I,X_\infty(G)\widehat{\otimes} A}(\pi_{K_I})$, where we use the same notation to denote the global section and we assume that the global section is finite projective, carrying mutually commuting actions of $(\varphi_I,\Gamma_{K_I})$. Then we have that the cohomology groups of the Herr complex $C^\bullet_{\varphi_I,\Gamma_I}(\mathcal{F})$ are coadmissible over $X_\infty(G)\widehat{\otimes}A$, and the corresponding cohomology groups of Iwasawa complex $C^\bullet_{\psi_I} (\mathcal{F})$ are coadmissible over $\Pi_{\mathrm{an},\infty,I,X_\infty(G)\widehat{\otimes} A}(\Gamma_{K_I})$. Moreover when $A=\mathbb{Q}_p$, $C^\bullet_{\psi_I}(\mathcal{F})$ lives in $\mathbb{D}^\flat_{\mathrm{perf}}(\Pi_{\mathrm{an},\infty,I,X_\infty(G)\widehat{\otimes} A}(\Gamma_{K_I}))$ and $C^\bullet_{\varphi_I,\Gamma_I}(\mathcal{F})$ lives in $\mathbb{D}^\flat_{\mathrm{perf}}(X_\infty(G)\widehat{\otimes} A)$.
\end{theorem}

\begin{remark}
The two theorems are around commutative period rings. Actually we conjectured that the corresponding analogs in the noncommutative setting should hold. In other words not only we would expect that one can follow the ideas in this work to establish some finiteness theorems when we have the coefficients that are noetherian noncommutative Banach rings as those discussed in \cref{4section5.1}, but also we should as well have some analogs of the results on the (noncommutative) equivariant Iwasawa cohomology. For more detail, see \cref{conjecture5.22}, \cref{conjecture5.23}, \cref{conjecture5.24}. 
\end{remark}

\indent It is actually not expected (due to some difficulty in high dimensional $p$-adic analysis) that one could derive the finiteness or perfectness in some ambient derived categories by using analogue of Kedlaya's key ingredient to prove the theorem mentioned in setting 1.2 above, which is to say there is no d\'evissage possibility transparent for us to use in order to reduce all the proof in the non-\'etale setting to the main results of \cite{4CKZ18} and \cite{4PZ19}, in the \'etale setting. The coherence of the cohomology is significant in our theory, for instance if one would like to study the global triangulation of $(\varphi,\Gamma)$-modules over any rigid analytic spaces, for more application see the discussion below.\\

\indent Now we discuss the applications of our results. The multidimensionalization represented in the program of \cite{4PZ19} and \cite{4CKZ18} is aimed at some motivation from $p$-adic Local Langlands program and relative Fukaya-Kato-Nakamura's local Tamagawa number conjecture or the $\varepsilon$-isomorphism conjecture. One could somehow regard the two conjectures as closely related to each other to some extent. Roughly speaking we have the following conjecture:

\begin{conjecture} \mbox{\bf (Nakamura \cite{4Nak2})} In the determinant category attached to any $(\varphi,\Gamma)$-modules over relative Robba rings, for each such $M$ there exists an isomorphism between the object $\mathbf{1}$ with the fundamental line attached to $M$. Moreover this isomorphism is compatible with algebraic functional equation, duality and the de Rham isomorphism.
\end{conjecture}

The program initiated in \cite{4PZ19} and \cite{4CKZ18} gives rise to some possible approach to this conjecture, by first generalizing this to higher dimensional situation, and then consider some sort of diagonal embedding in the group theoretic consideration. Again there should be a corresponding relative version of the conjectures of Pal-Z\'abr\'adi \cite{4PZ19}. \\

%\begin{center}
%\begin{tabularx}{\linewidth}{lX}
%Setting & Known or not\\
%\hline
%1.1    & Not known to us.\\
%1.2    & Not known to us.\\
%
%1.3    & Not known to us.\\
%
%1.4    & Not known to us.\\
%
%1.5    & Fukaya-Kato, meanwhile we have variants of Iwasawa main conjectures.\\
%
%1.6    & Fukaya-Kato-Nakamura, meanwhile we have variants of Iwasawa main conjectures.\\
%
%1.7    & We conjecture that for nice rigid analytic spaces (for instance admitting toric charts) there should be a variant.\\
%
%1.8    & We conjecture that for nice spaces this could be having a variant.\\
%
%1.9    & Not quite known to us.\\
%
%1.10    & Not quite known to us.\\
%
%1.11    & Not quite known to us.\\
%
%1.12    & Not quite known to us.\\
%
%1.13    & Z\'abr\'adi-Pal's conjecture.\\
%
%1.14  &Z\'abr\'adi-Pal's conjecture.\\
%
%1.15   & We conjecture there could be some variant.\\
%1.16   & We conjecture there could be some variant.
%\end{tabularx}
%\end{center}

\indent Another application will be deformation theory and the variation of the triangulation property represented already in the one-dimensional situation, see \cite{4KPX}. As mentioned above, local triangulation or global triangulation are significant properties in the study of Galois theoretic eigenvarieties. We should mention that actually our approach which is parallel to \cite{4KPX} is different from the one in \cite{4Liu1}. The latter used essentially some local spreading properties of the slope filtrations, which is tested by considering the sensibility of the neighbourhood around a good point on the slope filtrations. All of these are not directly expected in the context of \cite{4KPX} and the one we are considering in this paper. For the extended and detailed exposition, see the corresponding \cref{triangulation}.\\

\subsection{What's Next}

The framework we studied here is just some partial aspects of the whole story in the higher dimensional analytic geometry. One can definitely study more general analytic geometry even in the archimedean setting where one can search for more very typical multi-Hodge structures such as those coming from the hyper-K\"ahler geometry. We have only worked over polydisks and polyannuli in the strictly affinoid setting in regular rigid analytic geometry which we believe of course should be the most similar setting to the classical complex analysis of several variables. But that being said, these restrictions are not necessary at all. One can obviously ask how general the spaces which we could handle are. We certainly believe in some very general sense that the spaces we could potentially handle are those allowing us to do the suitable $p$-adic functional analysis and $p$-adic complex analysis of several variables.\\

The corresponding construction we established here targets directly at the noncommutative Tamagawa number conjectures. We should then follow \cite{4Nak1}, \cite{4Nak2}, \cite{4Zah1}, \cite{4PZ19} to amplify our programs in our mind from the philosophy in \cite{4BF1}, \cite{4BF2} and \cite{4FK}. And we should believe that our approach could be also  applied in many other contexts where \cite{4Kie1} and \cite{4KL3} will provide insights and applications, which certainly might not have to be Iwasawa theoretic. Certainly in the context of \cite{4PZ19} one still has to investigate in more detail the corresponding $p$-adic Hodge theoretic properties of objects studied in this paper, such as the corresponding possible analogs of the main results in \cite{4Ked1}, \cite{4And1}, \cite{4Me} and \cite{4Ber1}.\\

%\newpage 

\begin{center}
\begin{tabularx}{\linewidth}{lX}
Notation & Description (These are mainly starting participating in the discussion in Chapter 2 in the commutative setting, while in Chapter 5 in the noncommutative setting)\\
\hline
$p$    & A prime number which is $>2$.\\
$I$    & A finite set. \\
$K_\alpha$   & Finite extension of $\mathbb{Q}_p$ for each $\alpha\in I$, with chosen $\pi_{K,\alpha}$.\\
$G_{K_{\alpha}}$ & The corresponding Galois group of $\mathrm{Gal}(\overline{\mathbb{Q}}_p/K_\alpha)$ for each $\alpha\in I$.\\
$\Gamma_{K_{\alpha}}$ & The corresponding $\Gamma$ group of $K_\alpha$ for each $\alpha\in I$.\\
$s_I,r_I$ & Multidimensional radii $0<s_I\leq r_I$ with coordinates which are rational numbers.\\
$A$  & An affinoid algebra (which could be noncommutative).\\
$\Pi_{r_I,I}$ & Multidimensional convergent ring with respect to the radius $r_I$.\\
$\Pi_{[s_I,r_I],I}$ & Multidimensional analytic ring with respect to the radii in $[s_I,r_I]$.\\
$\Pi_{\mathrm{an},r_I,I}$ & Multidimensional analytic ring with respect to the radius $r_{I}$.\\
$\Pi_{\mathrm{an},\mathrm{con},I}$ & Multidimensional Robba rings.\\
$\Pi_{r_I,I,A}$ & Multidimensional convergent ring with respect to the radius $r_I$ with coefficients in an affinoid algebra $A$.\\
$\Pi_{[s_I,r_I],I,A}$ & Multidimensional analytic ring with respect to the radii in $[s_I,r_I]$ with coefficients in an affinoid algebra $A$.\\
$\Pi_{\mathrm{an},r_I,I,A}$ & Multidimensional analytic ring with respect to the radius $r_{I}$ with coefficients in an affinoid algebra $A$.\\
$\Pi_{\mathrm{an},\mathrm{con},I,A}$ & Multidimensional Robba rings with coefficients in an affinoid algebra $A$.\\

$(C_\alpha,l_\alpha)_{\alpha\in I}$  & Collections of pair $(C_\alpha,l_\alpha)$ each of which consists of a multiplicative topological group which is abelian with finite torsion part and free part isomorphic to $\mathbb{Z}_p$, and a continuous homomorphism from $C_\alpha$ to $\mathbb{Q}_p$ for each $\alpha\in I$.  \\
$\Pi_{\mathrm{an},r_I,I}(C_I)$ & Multidimensional analytic ring with respect to the radius $r_{I}$, with new variables from $C_I$.\\
$\Pi_{\mathrm{an},\mathrm{con},I}(C_I)$ & Multidimensional Robba rings, with new variables from $C_I$.\\ 	
$\Pi_{\mathrm{an},r_I,I,A}(C_I)$ & Multidimensional analytic ring with respect to the radius $r_{I}$ with coefficients in an affinoid algebra $A$, with new variables from $C_I$.\\
$\Pi_{\mathrm{an},\mathrm{con},I,A}(C_I)$ & Multidimensional Robba rings with coefficients in an affinoid algebra $A$, with new variables from $C_I$.\\ 	
$\Pi_{\mathrm{an},\infty,I}(C_I)$ & Multidimensional analytic ring with respect to the radius $r_{I}$, with new variables from $C_I$.\\
$\Pi_{\mathrm{an},\infty,I,A}(C_I)$ & Multidimensional analytic ring with respect to the radius $r_{I}$ with coefficients in an affinoid algebra $A$, with new variables from $C_I$.\\

\end{tabularx}
\end{center}

\begin{center}
\begin{tabularx}{\linewidth}{lX}
Notation & Description (These are mainly starting participating in the discussion in Chapter 3 and 4 in the commutative setting, while in Chapter 5 in the noncommutative setting)\\
\hline

$\Pi_{\mathrm{an},r_I,I}(\pi_{K,I})$ & Multidimensional analytic ring with respect to the radius $r_{I}$, with new variables in $\pi_{K,I}$.\\
$\Pi_{[s_I,r_I],I}(\pi_{K,I})$ & Multidimensional analytic ring with respect to the radii in $[s_I,r_{I}]$, with new variables in $\pi_{K,I}$.\\
$\Pi_{\mathrm{an},\mathrm{con},I}(\pi_{K,I})$ & Multidimensional Robba rings, with new variables in $\pi_{K,I}$.\\ 	
$\Pi_{\mathrm{an},r_I,I,A}(\pi_{K,I})$ & Multidimensional analytic ring with respect to the radius $r_{I}$ with coefficients in an affinoid algebra $A$, with new variables in $\pi_{K,I}$.\\
$\Pi_{[s_I,r_I],I,A}(\pi_{K,I})$ & Multidimensional analytic ring with respect to the radii in $[s_I,r_{I}]$ with coefficients in an affinoid algebra $A$, with new variables in $\pi_{K,I}$.\\
$\Pi_{\mathrm{an},\mathrm{con},I,A}(\pi_{K,I})$ & Multidimensional Robba rings with coefficients in an affinoid algebra $A$, with new variables in $\pi_{K,I}$.\\

$\Pi_{\mathrm{an},\infty,I}(\Gamma_{K,I})$ & Multidimensional analytic ring with respect to the radius $r_{I}$, with new variables from $\Gamma_K$.\\

$\Pi_{\mathrm{an},\infty,I,A}(\Gamma_{K,I})$ & Multidimensional analytic ring with respect to the radius $r_{I}$ with coefficients in an affinoid algebra $A$, with new variables from $\Gamma_K$.\\

\end{tabularx}
\end{center}

%\newpage

\newpage\section{Fr\'echet Rings and Sheaves in the Commutative Ca
se}

\subsection{Fr\'echet Objects}

\indent In this subsection we are going to define the main objects we are going to study ring theoretically. Everything will be modeling the one-dimensional situation. As indicated in the title of the section, we will consider Fr\'echet-Stein objects in the sense of \cite{4ST1}. We first proceed by defining the following rings:

\begin{definition} \mbox{\bf{(After KPX \cite[Notation 2.1.1]{4KPX})}}
We consider the multi\\ radii as mentioned in the notation list, where $0<s_I\leq r_I$ are all rational numbers or $\infty$. Under this assumption first let $I$ is a finite set. Then we consider the multidimensional polynomial rings generalizing the one-dimensional situation with variables $T_1^{\pm 1},..,T_I^{\pm 1}$ which is to say the ring:
\begin{displaymath}
\mathbb{Q}_p[T_1^{\pm 1},..,T_I^{\pm 1}].	
\end{displaymath}
Now we endow this rings with some specific Gauss valuations and norms as in the following way for any multi radii $t_I>0$:
\begin{displaymath}
v_{t_I}(f)=v_{t_I}(\sum_{n_I}a_{n_I}T_I^{i_I}):= \inf_{n_I}\left\{v_p(a_{n_I})+\sum_{\alpha\in I}t_\alpha i_\alpha\right\},
\end{displaymath}
\begin{displaymath}
\left\|f\right\|_{t_I}=\left\|\sum_{n_I}a_{n_I}T_I^{i_I}\right\|_{t_I}:= \sup_{n_I}\left\{|a_{n_I}|p^{-\sum_{\alpha\in I}t_\alpha i_\alpha/(p-1)}\right\}.
\end{displaymath}
We are going to use the Gauss valuations and norms with respect to the multi-interval $[s_I,r_I]$ in the following sense:
\begin{displaymath}
v_{[s_I,r_I]}(f)=v_{[s_I,r_I]}(\sum_{n_I}a_{n_I}T_I^{i_I}):= \min_{t_\alpha\in\{s_I,r_I\}}\left\{v_p(a_{n_I})+\sum_{\alpha\in I}t_\alpha i_\alpha\right\},
\end{displaymath}
\begin{displaymath}
\left\|f\right\|_{[s_I,r_I]}=\left\|\sum_{n_I}a_{n_I}T_I^{i_I}\right\|_{[s_I,r_I]}:= \max_{t_\alpha\in\{s_I,r_I\}}\left\|\sum_{n_I}a_{n_I}T_I^{i_I}\right\|_{t_\alpha}.	
\end{displaymath}
The latter will give rise to the following definition:
\begin{displaymath}
\Pi_{[s_I,r_I],I}	
\end{displaymath}
which is defined to the completion of the ring $\mathbb{Q}_p[T_1^{\pm 1},..,T_I^{\pm 1}]$ with respect to the norms defined above with respect to the multi-interval $[s_I,r_I]$. For the polynomial rings:
\begin{displaymath}
\mathbb{Q}_p[T_1,..,T_I]	
\end{displaymath}
we then instead define for some $t_I>0$:

\begin{displaymath}
v_{t_I}(f)=v_{t_I}(\sum_{n_I}a_{n_I}T_I^{i_I}):= \inf_{n_I}\left\{v_p(a_{n_I})+\sum_{\alpha\in I}t_\alpha i_\alpha\right\},
\end{displaymath}
\begin{displaymath}
\left\|f\right\|_{t_I}=\left\|\sum_{n_I}a_{n_I}T_I^{i_I}\right\|_{t_I}:= \sup_{n_I}\left\{|a_{n_I}|p^{-\sum_{\alpha\in I}t_\alpha i_\alpha/(p-1)}\right\}.
\end{displaymath}
This will then give rise to the rings corresponding to the closed polydiscs and open polydiscs:
\begin{displaymath}
\Pi_{[t_I,\infty_I],I}, \Pi_{\mathrm{an},\infty_I,I}:=\bigcap_{t_I>0}\Pi_{[t_I,\infty_I],I}.	
\end{displaymath}
Then correspondingly we define the following analytic rings corresponding to the full polyannuli:
\begin{displaymath}
\Pi_{\mathrm{an},r_I,I}:=\bigcap_{0<s_I\leq r_I} \Pi_{[s_I,r_I],I},\Pi_{\mathrm{an},\mathrm{con},I}:=\bigcup_{r_I>0}\bigcap_{0<s_I\leq r_I} \Pi_{[s_I,r_I],I}. 
\end{displaymath}

\end{definition}

\begin{remark}
One could actually define the rings by considering all the (rigid) analytic functions over the corresponding open polydiscs, closed polydiscs, polyannuli and full polyannuli, which will be more transparent.
\end{remark}

\begin{remark}
Actually one could see that in the definition, one could find that especially in the multi-interval situation, the rightmost radii could be mixed-type with respect to the infinity which is to say that one could have a genuine proper subset of indexes where the rightmost radii are exactly infinite. We will mainly focus on the situation where such proper subset does not emerge in our discussion later, if not specified otherwise.	
\end{remark}

\indent Now we consider the relative version of the definition as in the one dimensional situation in \cite[Notation 2.1.1]{4KPX}.

\begin{definition} \mbox{\bf{(After KPX \cite[Notation 2.1.1]{4KPX})}}
Let $A$ be an affinoid algebra over $\mathbb{Q}_p$ in the sense of rigid analytic spaces. Then we are going to define the relative version of the rings defined as above in the absolute situation:
\begin{displaymath}
\Pi_{[s_I,r_I],I,A}, 
\end{displaymath}
as the analytic functions relative to $A$ over the products:
\begin{displaymath}
\mathrm{Max}A\times A^I[s_I,r_I],0<s_I\leq r_I.
\end{displaymath}
Then we define:
\begin{displaymath}
\Pi_{\mathrm{an},r_I,I,A}:=\bigcap_{0<s_I\leq r_I} \Pi_{[s_I,r_I],I,A},\Pi_{\mathrm{an},\mathrm{con},I,A}:=\bigcup_{r_I>0}\bigcap_{0<s_I\leq r_I} \Pi_{[s_I,r_I],I,A}.	
\end{displaymath}
Correspondingly almost in the parallel and similar way one could define the following ring of analytic functions relative to the space $A$:
\begin{displaymath}
\Pi_{[s_I,\infty_I],I,A},\Pi_{\mathrm{an},\infty_I,I,A}:=\bigcap_{s_I>0}\Pi_{[s_I,\infty_I],I,A}.	
\end{displaymath}
\end{definition}

\indent Eventually we are going to work within the framework of Fr\'echet-Stein algebras over rigid analytic spaces. Now we are going to consider the construction over the affinoid algebras. First let $r_{I,0}$ be a fixed radii (where we do not consider the situation that there is a finite proper set of $I$ where the radii approach the infinity). Then we consider a decreasing sequence of multi radii $\{r_{I,n}\}_{n\geq 0}$ which is to say that for each $\alpha\in I$, $r_{\alpha,n+1}\leq r_{\alpha,n}$, and when $n=0$ the radii is the one fixed above. And we require that for each $\alpha\in I$, $r_{\alpha,1}< r_{\alpha,0}$.

\begin{remark}
Actually this looks a little bit strange, but one could eventually choose in the following way that for each $n\geq 1$ one just set $r_{\alpha,n}=r_n$ for each $\alpha$.	
\end{remark}

\indent Then we consider the following result which in some non-trivial way generalizes the one dimensional situation:

\begin{proposition}
For all $n\geq 1$, the rings $\Pi_{[r_{I,n},r_{I,0}],I,A}$ and $\Pi_{[r_{I,n},\infty_I],I,A}$ are noetherian as Banach $A$-algebras.	
\end{proposition}

\begin{proof}
We are going to follow ideas in \cite{4PZ19} and \cite{4Ked9} to prove this. This will be some application of a theorem due to Fulton in \cite{4Ful1}. We are going to show the statement in the situation where $A$ is just the field $\mathbb{Q}_p$, and leave the general situation to the reader. We are going to only choose to prove the statement for the ring $\Pi_{[r_{I,n},r_{I,0}],I,A}$, and we leave the check for the other ring to the readers. We are going to show that the ring $\Pi_{[r_{I,n},r_{I,0}],I,A}$ is weakly complete algebras over the corresponding adic completion of the ring $\mathbb{Z_p}[T_\alpha,\alpha\in I]$ with the ideal $(p)$ generated by the element $p$ , and a theorem due to Fulton will imply that the ring $\Pi_{[r_{I,n},r_{I,0}],I,A}$ is actually noetherian.\\
\indent Following Monsky-Washnitzer, we recall as in the following the basic definition of a ring $R$ being relatively weakly complete algebras over $A$. This means that for any element $f\in R$ and any degree $d\geq 0$, one could find $k$ elements $t_1,...,t_k$ and a polynomial $p_d$ with $k$-variables with the coefficients in $(p)^d$, such that there exists a constant $C>0$ with $\mathrm{deg}p_d\leq C(d+1)$ and with that:
\begin{displaymath}
f=\sum_{d\geq 0}p_d(t_1,...,t_k).	
\end{displaymath}
Therefore we are going to find some essential estimate over the degree of the polynomial with respect to some given element $f$. To get start for any such given $f$ and any given $d$ in the condition recalled as above, we consider a degree $\partial_d$ positive which is maximal in the sense that $f$ is inside $T_I^{-\partial_d}\mathbb{Z}_p/(p)^d$ and the index is the maximal throughout all such indexes. Then we are going to use the notation $\partial_I$ to denote the corresponding multi-index (which is possibly not the same as $-\partial_d^I$). Now we take a constant $c_1>1$ such that $\left\|f\right\|_{[r_{I,n},r_{I,0}]}\leq c_1$. Then we compute:
\begin{displaymath}
c_1\geq \left\|f\right\|_{[r_{I,n},r_{I,0}]} \geq \left\|f\right\|_{r_{I,0}} \geq |\mathrm{coeff}_{\partial_I}(f)|p^{\partial_d r_{\alpha,0}/(p-1)}\prod_{\beta\in I\backslash \{\alpha\}}p^{-r_{\beta,0} i_\beta/(p-1)}
\end{displaymath}
which gives rise to that for some $c_2>0$:
\begin{displaymath}
|\mathrm{coeff}_{\partial_I}(f)|p^{\partial_d r_{\alpha,0}/(p-1)}\leq c_2.	
\end{displaymath}
Then we compute as in the following:
\begin{align}
|\mathrm{coeff}_{\partial_I}(f)|p^{\partial_d r_{\alpha,0}/(p-1)}&\leq c_2\\	
\mathrm{log}_p|\mathrm{coeff}_{\partial_I}(f)|+\partial_d r_{\alpha,0}/(p-1)&\leq \mathrm{log}_pc_2\\
(1-d)+\partial_d r_{\alpha,0}/(p-1)&\leq \mathrm{log}_pc_2\\
\partial_d r_{\alpha,0}&\leq (p-1)\mathrm{log}_pc_2+(d-1)(p-1)\\
\partial_d &\leq (p-1)\mathrm{log}_pc_2/r_{\alpha,0}+(d-1)(p-1)/r_{\alpha,0}
\end{align}
which gives rise to the desired estimate, which then further finishes the proof by considering Fulton's theorem in \cite{4Ful1} which says that the desired relative weakly complete algebra is now noetherian.
\end{proof}

 \indent Then we consider the following statement which sends us to the framework of \cite{4ST1}, which generalizes the one dimensional situation in \cite[Section 2.1]{4KPX}:
 
\begin{proposition}
For all $n\geq 1$, the rings $\Pi_{\mathrm{an},r_{I,0},I,A}$ and $\Pi_{\mathrm{an},\infty_I,I,A}$ are Fr\'echet-Stein in the sense of \cite{4ST1}.
\end{proposition}

\begin{proof}
We keep the notations simple so we just consider the proof for the first ring and leave the task for checking the second case to the reader. Now indeed, we first note that $\Pi_{[r_{I,n},r_{I,0}],I,A}$ is Fr\'echet algebra in the sense that it is the ring of the analytic functions over a union of the relative annuli as in the following:
\begin{displaymath}
\bigcup_{0<r_{I,n}\leq r_{I,0}} \mathrm{Max}A\times A^I[r_{I,n},r_{I,0}]. 
\end{displaymath}
Furthermore we note that this is also the increasing union which is of the form of:
\begin{displaymath}
\bigcup_{0<r_{I,n}\leq r_{I,0}} \mathrm{Max}A\times A^I[r_{I,n},r_{I,0}],
\end{displaymath}
where $ \mathrm{Max}A\times A^I[r_{I,n},r_{I,0}]\subset  \mathrm{Max}A\times A^I[r_{I,n+1},r_{I,0}]$. The corresponding transition map over the algebras then with respect to this increasing union is actually flat since it is \'etale by the production construction and the base change properties. Finally the transition map is of topologically dense image again by the product construction.	
\end{proof}

\indent Before we proceed to study more in more detail on the representation and structure of the algebras defined above, we first define the vector bundles over polyannuli generalizing the one dimensional result in \cite[Definition 2.1.3]{4KPX}. We mimick the definition, which gives rise to the following:

\begin{definition} \mbox{\bf{(After KPX \cite[Definition 2.1.3]{4KPX})}}
We define coherent\\
 sheaves of modules and the global sections over the ring $\Pi_{\mathrm{an},r_{I_0},I,A}$ in the following way. First a coherent sheaf $\mathcal{M}$ over the relative multidimensional Robba ring $\Pi_{\mathrm{an},r_{I_0},I,A}$ is defined to be a family $(M_{[s_I,r_I]})_{[s_I,r_I]\subset (0,r_{I,0}]}$	of modules of finite type over each relative $\Pi_{[s_I,r_I],I,A}$ satisfying the following two conditions as in the more classical situation:
\begin{displaymath}
M_{[s_I',r_I']}\otimes_{\Pi_{[s'_I,r'_I],I,A}}\Pi_{[s_I'',r_I''],I,A}\overset{\sim}{\rightarrow} M_{[s''_I,r''_I]}	
\end{displaymath}
for any multi radii satisfying $0<s_I'\leq s_I''\leq r_I''\leq r_I'\leq r_{I,0}$, with furthermore 
\begin{displaymath}
(M_{[s'_I,r'_I]}\otimes_{\Pi_{[s'_I,r'_I],I,A}}\Pi_{[s''_I,r''_I],I,A}	)\otimes_{\Pi_{[s''_I,r''_I],I,A}}\Pi_{[s'''_I,r'''_I],I,A}\overset{\sim}{\rightarrow}\Pi_{[s'''_I,r'''_I],I,A},
\end{displaymath}
for any multi radii satisfying the following condition:
\begin{displaymath}
0<s_I'\leq s_I''\leq s_I'''\leq r_I'''\leq r_I''\leq r_I'\leq r_{I,0}.	
\end{displaymath}
Then we define the corresponding global section of the coherent sheaf $\mathcal{M}$, which is usually denoted by $M$ which is defined to be the following inverse limit:
\begin{displaymath}
M:=\varprojlim_{s_I\rightarrow 0_I^+} M_{[s_I,r_{I,0}]}.	
\end{displaymath}
We are going to follow \cite[Definition 2.1.3]{4KPX} to call any module defined over $\Pi_{\mathrm{an},r_{I_0},I,A}$ coadmissible it comes from a global section of a coherent sheaf $\mathcal{M}$ in the sense defined above. 
\end{definition}

\indent Then we could collect the following properties form more general theory from \cite{4ST1} as in \cite[Lemma 2.1.4]{4KPX} in the one dimensional situation.

\begin{proposition} \mbox{\bf{(After KPX \cite[Lemma 2.1.4]{4KPX})}} \label{4prop2.9}
We have the following properties for the coherent sheaves defined above over the relative multi-dimensional Robba rings $\Pi_{\mathrm{an},r_{I_0},I,A}$:\\
I.	For each multi radii $0< s_I\leq r_I\leq r_{I,0}$ as above, we have that the global section $M$ is dense in the corresponding section $M_{[s_I,r_I]}$ where $r_{I,0}$ is finite, where the statement is also true for $M[1/T_I]$ for $r_{I,0}$ infinite but with $r_I$ not;\\
II. For each multi radii $0< s_I\leq r_I\leq r_{I,0}$ as above, we have the following comparison:
\begin{displaymath}
M\otimes_{\Pi_{\mathrm{an},r_{I_0},I,A}}\Pi_{[s_I,r_I],I,A}\overset{\sim}{\rightarrow}	M_{[s_I,r_I]};\\
\end{displaymath}
III. For each multi radii $0< s_I \leq r_{I,0}$ we have the following higher vanishing result:
\begin{displaymath}
R^1\varprojlim_{s_I\rightarrow 0^+}M_{[s_I,r_{I,0}]}=0;\\	
\end{displaymath}
IV. The kernel or cokernel of a morphism between coadmissible modules over $\Pi_{\mathrm{an},r_{I_0},I,A}$ is again coadmissible;\\
V. Any finite presented module over $\Pi_{\mathrm{an},r_{I_0},I,A}$ is coadmissible;\\
VI. Any submodule of finite type of any coadmissible  $\Pi_{\mathrm{an},r_{I_0},I,A}$-module is again coadmissible;\\ 
VII. For each multi radii $0< s_I\leq r_I\leq r_{I,0} $, and for the global section defined above, we have that $M\otimes_{\Pi_{\mathrm{an},r_{I_0},I,A}}\Pi_{[s_I,r_I],I,A}$ is flat over $\Pi_{[s_I,r_I],I,A}$. 

\end{proposition}

\begin{proof}
See \cite[Lemma 2.1.4]{4KPX} and \cite{4ST1}. Again let us mention that actually in this setting initially we only have the corresponding results for $r_I$ fixed which could for instance be $r_{I,0}$, however this could be overcome by the same method as in the one dimensional situation in \cite[Lemma 2.1.4]{4KPX}.
\end{proof}

\indent We then follow the ideas and approaches from \cite{4KPX} to 
study further the structure theory of the coherent sheaves defined above. The following is then a generalization of the corresponding result in the one-dimensional situation in \cite[Corollary 2.1.5]{4KPX}.

\begin{lemma} \mbox{\bf{(After KPX \cite[Corollary 2.1.5]{4KPX})}}
The rings $\Pi_{\mathrm{an},r_{I_0},I,A}$, and $\Pi_{\mathrm{an},\mathrm{con},I,A}$ are flat over $A$ as in \cite[Corollary 2.1.5]{4KPX}.	
\end{lemma}

\begin{proof}
We adapt the corresponding method of proof in \cite[Corollary 2.1.5]{4KPX} to prove this generalization. Indeed one could then prove this by the same method of \cite[Corollary 2.1.5]{4KPX} since then by taking a suitable basis we will have $\Pi_{\mathrm{an},r_{I_0},I,\mathbb{Q}_p}\overset{\simeq}{\rightarrow}\widehat{\oplus}_{\eta}\mathbb{Q}_pe_\eta$. Then one could just finish as in \cite[Corollary 2.1.5]{4KPX}. 	
\end{proof}

\begin{definition} \mbox{\bf{(After KPX \cite[Definition 2.1.3]{4KPX})}}\\
Consider a coherent sheaf $M_{r_{I,0}}$ over $\Pi_{\mathrm{an},r_{I_0},I,A}$ defined as above, then we are going to call this sheaf a vector bundle if for each $[s_I,r_I]$ suitably located multi interval we have $M_{[s_I,r_I]}$ is flat over $\Pi_{\mathrm{an},[s_I,r_I],I,A}$.	
\end{definition}

\indent This will give us a chance to consider the following result:

\begin{proposition} \mbox{\bf{(After KPX \cite[Lemma 2.1.6]{4KPX})}} \label{prop2.12}
For any coherent sheaf $M_{r_{I,0}}$ over $\Pi_{\mathrm{an},r_{I_0},I,A}$ which admits a vector bundle structure in our situation, we then have that the global section $M$ of it is finitely projective if and only if it is finitely presented. (Actually one can even derive the corresponding sufficient and necessary condition to be just being finitely generated by applying \cite[Corollary 2.6.8]{4KL2}). 	
\end{proposition}

\begin{proof}
Basic representation theory tells us that the module involved is finitely projective if and only if it is finitely presented and flat. Therefore back to our proposition one direction of the implications is trivial, which is to say that we could now assume that the corresponding module is finitely presented, and it suffices to show that it is flat. Therefore as in \cite[Lemma 2.1.6]{4KPX} we consider the map $I\otimes_{\Pi_{\mathrm{an},r_{I_0},I,A}} M\rightarrow M$ for some finitely generated ideal $I$. Then we are going to show that this is injective. First choose some presentation of $M$ which is of type $(m,n)$ over the base ring $\Pi_{\mathrm{an},r_{I_0},I,A}$ then translate this to the corresponding statement for the ideal $I$ which is to say $I^m\rightarrow I^n$. This actually gives us a chance to derive the statement that $I\otimes_{\Pi_{\mathrm{an},r_{I_0},I,A}} M$ is also coadmissible by the \cref{4prop2.9}. Now again by \cref{4prop2.9} we have that it suffices to show that:
\begin{displaymath}
I\otimes_{\Pi_{\mathrm{an},r_{I_0},I,A}} M\rightarrow M	
\end{displaymath}
is injective over some base change to some $\Pi_{\mathrm{an},[s_I,r_I],I,A} $. Then as in \cite[Lemma 2.1.6]{4KPX} one could finish the proof since this is nothing but the following map:
\begin{displaymath}
(I\otimes_{\Pi_{\mathrm{an},r_{I_0},I,A}} \Pi_{\mathrm{an},[s_I,r_I],I,A})\otimes  M_{[s_I,r_I]} \rightarrow M_{[s_I,r_I]}	
\end{displaymath}
which is injective by the flatness from \cite[Lemma 2.1.6]{4KPX}.
\end{proof}

\indent We then have the following corollary which is a direct analog of the corresponding results in \cite[Corollary 2.1.7]{4KPX}:

\begin{corollary}
Let $A$ be an affinoid algebra which is reduced. Now consider a module $M$ over the Robba ring with respect to a specific radius $r_I$, defined over $\mathrm{Max}(A)\times \mathbb{A}^n(0,r_I]$ and we assume that for any maximal ideal $\mathfrak{m}_\eta$ of $A$ the fiber $M_{\mathfrak{m}_\eta}$ has the same rank. Then we have that $M$ is finite projective if we further assume that $M$ is finitely presented (as mentioned in the previous proposition one can even assume that it is finitely generated).	
\end{corollary}

\begin{proof}
This is a corollary of the previous proposition and the corresponding noetherian property as in \cite[Corollary 2.1.7]{4KPX}.	
\end{proof}

\indent Before we study the relationship between the finiteness of the global sections and the coherent sheaves, we generalize the uniform finiteness in the one dimensional situation which is established in \cite[Definition 2.1.9]{4KPX}.

\begin{definition} \mbox{\bf{(After KPX \cite[Definition 2.1.9]{4KPX})}} 
We first generalize the notion of admissible covering in our setting, which is higher dimensional generalization of the situation established in \cite[Definition 2.1.9]{4KPX}. For any covering $\{[s_I,r_I]\}$ of $(0,r_{I,0}]$, we are going to specify those admissible coverings if the given covering admits refinement finite locally (and each corresponding member in the covering has the corresponding interior part which is assumed to be not empty with respect to each $\alpha\in I$). Then it is very natural in our setting that we have the corresponding notations of $(m,n)$-finitely presentedness and $n$-finitely generatedness for any $m,n$ positive integers. Then based on these definitions we have the following generalization of the corresponding uniform finiteness. First we are going to call a coherent sheaf $\{M_{s_I,r_I}\}_{[s_I,r_I]}$ over $\Pi_{\mathrm{an},r_{I,0},I,A}$ uniformly $(m,n)$-finitely presented if there exists an admissible covering $\{[s_I,r_I]\}$ and a pair of positive integers $(m,n)$ such that each module $M_{[s_I,r_I]}$ defined over $\Pi_{\mathrm{an},[s_I,r_I],I,A}$ is $(m,n)$-finitely presented. Also we have the notion of uniformly $n$-finitely generated for any positive integer $n$ by defining that to mean under the existence of an admissible covering we have that each member $M_{[s_I,r_I]}$ in the family defined over $\Pi_{\mathrm{an},[s_I,r_I],I,A}$ is $n$-finitely generated. Sometimes we are also going to use the notions of being uniformly finitely generated and being uniformly finitely presented to mean the corresponding objects when the corresponding uniform numbers of the generators are well-understood.
\end{definition}

\indent After defining these notions one could consider the following generalization of the corresponding results in \cite[Lemma 2.1.10]{4KPX} specifically in the one-dimensional situation. Before establishing this kind of generalization we first consider the following:

\begin{lemma} \mbox{\bf{(After KPX \cite[Lemma 2.1.10]{4KPX})}} \label{lemma2.15}
Consider now the following short exact sequence of coherent sheaves over $\Pi_{\mathrm{an},r_{I,0},I,A}$:
\[
\xymatrix@R+2pc@C+0pc{
0\ar[r] \ar[r] &(M^1_{[s_I,r_I]})_{\{[s_I,r_I]\}} \ar[r] \ar[r] &(M_{[s_I,r_I]})_{\{[s_I,r_I]\}} \ar[r] \ar[r] &(M^2_{[s_I,r_I]})_{\{[s_I,r_I]\}} \ar[r] \ar[r] &0.} 
\]
Then we have that if $(M^1_{[s_I,r_I]})_{\{[s_I,r_I]\}}$ and one of the rest two sheaves are uniformly finitely presented so is the third, moreover we have that if $(M_{[s_I,r_I]})_{\{[s_I,r_I]\}},(M^2_{[s_I,r_I]})_{\{[s_I,r_I]\}}$ are uniformly finitely presented then the third one will be uniformly finitely generated.
\end{lemma}

\begin{proof}
See \cite[Lemma 2.1.10]{4KPX}. 	
\end{proof}

%%%%%%%%%%%%%%%%%%%%%%%%%%%%%%%%%%%%%%revise this if necessary.

\begin{lemma}  \mbox{\bf{(After KPX \cite[Lemma 2.1.11]{4KPX})}}
Let $\{M_{[s_I,r_I]}\}$	be a coherent sheaf over $\Pi_{\mathrm{an},r_{I,0},I,A}$ with global section which we will denote it by $M$ (with respect to an admissible covering in our context). Suppose that we have a set of generators $\{\mathbf{e}_1,...,\mathbf{e}_n\}$ which generates $M_{[s_I,r_I]}$ for each $[s_I,r_I]$ involved, then we have that this set of generators will generate $M$ as well.

\end{lemma}

\begin{proof}
As in the one dimensional situation (see \cite[Lemma 2.1.11]{4KPX}) one could prove this by first exhibiting a map by using the given basis $(\Pi_{\mathrm{an},r_{I,0},I,A})^n\rightarrow M$, then arguing that this will give rise to the desired generating set by considering the corresponding projection which is again a well-defined surjection and comparison under the assumption mentioned above.
\end{proof}

\begin{proposition} \mbox{\bf{(After KPX \cite[Proposition 2.1.13]{4KPX})}} \label{4prop1}
Consider any arbitrary coherent sheaf $\{M_{[s_I,r_I]}\}$ over $\Pi_{\mathrm{an},r_{I,0},I,A}$ with the global section $M$. Then we have the following corresponding statements:\\
I.	The coherent sheaf $\{M_{[s_I,r_I]}\}$ is uniformly finitely generated iff the global section $M$ is finitely generated;\\
II. The coherent sheaf $\{M_{[s_I,r_I]}\}$ is uniformly finitely generated iff the global section $M$ is finitely presented;\\
III. The coherent sheaf $\{M_{[s_I,r_I]}\}$ is uniformly finitely presented vector bundle iff the global section $M$ is finite projective.
\end{proposition}

\begin{proof}
One could derive the the second and the third statements by using the \cref{lemma2.15} and \cref{prop2.12}. Therefore it suffices for us to prove the first statement. Indeed one could see that one direction in the statement is straightforward, therefore it suffices now to consider the other direction. So now we assume that the coherent sheaf $\{M_{[s_I,r_I]}\}$ is uniformly finitely generated for some $n$ for instance. Then we are going to show that one could find corresponding generators for the global section $M$ from the given uniform finitely generatedness. 

Now by definition we consider an admissible covering taking the form of $\{[s_I,r_I]\}$. In \cite[Proposition 2.1.13]{4KPX} the corresponding argument goes by using special decomposition of the covering in some refined way, mainly by reorganizing the initial admissible covering into two single ones, where each of them will be made into the one consisting of those well-established intervals with no overlap. In our setting, we need to use instead the corresponding upgraded nice decomposition of the given admissible covering in our setting, which is extensively discussed in \cite[2.6.14-2.6.17]{4KL2}. Namely as in \cite[2.6.14-2.6.17]{4KL2} there is a chance for us to extract $2^{|I|}$ families $\{\{[s_{I,\delta_i},r_{I,\delta_i}]\}_{i=0,1,...}\}_{\{1_i\},\{2_i\},...,\{N_i\},}$ ($N=2^{|I|}$) of intervals, where there is a chance to have the situation where for each such family the corresponding intervals in it could be made to be disjoint in pairs. Then we are in the situation of \cite[Proposition 2.6.17]{4KL2}, namely whenever we have a $2^{|I|}$-uniform covering, then the uniformly finiteness throughout all the coverings involved will guarantee the corresponding finiteness of the corresponding global section which proves the proposition. To be more precise we make the argument in \cite[Proposition 2.6.17]{4KL2} more precise in our setting. First we write the corresponding space $\Lambda^I(0,r_{I,0}]_A$ into the following form:
\begin{displaymath}
\Lambda^I(0,r_{I,0}]_A:=\bigcup_{\alpha\geq 0} \Lambda^I[s_{I,\alpha},r_{I,0}]_A.	
\end{displaymath}
For $\Lambda^I[s_{I,-1},r_{I,0}]_A$ we choose that to be empty. Then for a chosen family $\{[s_{I,\delta_i},r_{I,\delta_i}]\}_{i=0,1,...}$, we use the notation $\{[s_{I,\gamma},r_{I,\gamma}]\}_{\Gamma}$ to denote this specific family where $\Gamma$ is endowed with the corresponding well ordering as in \cite[proposition 2.6.17]{4KL2}. Now we consider for each $\gamma$ the corresponding index $\alpha(\gamma)$ such that the intersection between $\Lambda^I[s_{I,\alpha(\gamma)},r_{I,0}]_A$ and $\Lambda^I[s_{I,\gamma},r_{I,\gamma}]_A$ is empty and meanwhile this index is the largest among all the candidates such that the intersection like this is empty. Then we consider the corresponding subspace in the following way:
\begin{displaymath}
B_\gamma:=	\Lambda^I[s_{I,\alpha(\gamma)},r_{I,0}]_A\cup \bigcup_{\gamma'<\gamma,\gamma'\in \Gamma}\Lambda^I[s_{I,\gamma},r_{I,\gamma}]_A.
\end{displaymath}
Then we can see that one can find an element in the section $\mathcal{O}(B_\gamma)$ which by considering the corresponding restriction gives rise to some element which represents an invertible of an element which is topologically nilpotent over the section over $\Lambda^I[s_{I,\gamma},r_{I,\gamma}]_A$ but meanwhile vanishes over the section over $B_\gamma$. Then by the density of the global section one can see that one can upgrade this element to be some element in $\mathcal{O}(\Lambda^I(0,r_{I,0}]_A)$ which which represents an invertible of an element which is topologically nilpotent over the section over $\Lambda^I[s_{I,\gamma},r_{I,\gamma}]_A$ but meanwhile represents an element which is topologically nilpotent over the section over $B_\gamma$. We denote this element by $x_\gamma$. Now we can perform the corresponding iteration process in \cite[proposition 2.6.17]{4KL2} to generate the following families of elements in $M$:
\begin{align}
&...,\\
&e_{\gamma,1},...,e_{\gamma,k},\\
&...	
\end{align}
whose limit will be the corresponding finite generating set of the global section. We first start from a generating set of $M_{[s_{I,\gamma},r_{I,\gamma}]}$ such $e_1,...,e_k$, then we are considering will be the sequence defined in the following way:
\begin{displaymath}
e_{\gamma,i}:=e_{\gamma,i}+x_\gamma^\eta e_i,  	
\end{displaymath}
where $\eta$ is sufficiently large integer. This process guarantee that we can have converging limit set of generators. The generators actually generate the corresponding section $M_{[s_{I,\gamma},r_{I,\gamma}]}$ since in our situation we can control the corresponding difference between the restriction of our produced set of elements and the actually generator over such section as in \cite[proposition 2.6.17]{4KL2} and \cite[Lemma 2.1.12]{4KPX}.

\end{proof}

\subsection{Further Discussion and Residue Pairings}

\indent We collect some further continuation of the discussion in the previous section before we consider the residue pairings in our setting, these are higher dimensional generalization of the corresponding results in \cite{4KPX}.

\begin{lemma} \mbox{\bf{(After KPX \cite[Lemma 2.1.16]{4KPX})}}
Consider a morphism $f:M\rightarrow N$ of modules over $\Pi_{\mathrm{an},r_{I,0},I,A}$ with kernel and cokernel of finite types. Suppose $f:M\rightarrow N$ is injective over $\Pi_{\mathrm{an},\mathrm{con},I,A}$. Then $f$ is injective over $\Pi_{\mathrm{an},r_{I},I,A}$ for some radii $r_I\leq r_{I,0}$. Furthermore if $f:M\rightarrow N$ is surjective over $\Pi_{\mathrm{an},con,I,A}$, then $f$ is surjective over $\Pi_{\mathrm{an},r_{I},I,A}$ for some radii $r_I\leq r_{I,0}$.  	
\end{lemma}

\begin{proof}
See the proof of \cite[Lemma 2.1.16]{4KPX} which is a reinterpretation of the properties of being of finite type of the corresponding kernel and cokernel involved.	
\end{proof}

\begin{proposition} \mbox{\bf{(After KPX \cite[Lemma 2.1.18]{4KPX})}}
Suppose now that $M$ is a module defined over $\Pi_{\mathrm{an},r_{I,0},I,A}$ and is assumed to be finite over $A$. Then we have that there exists some multi radii $r_I\leq r_{I,0}$ such that the base change of $M$ to $\Pi_{\mathrm{an},r_{I},I,A}$ vanishes. 	
\end{proposition}

\begin{proof}
This is a higher dimensional generalization of the corresponding results in \cite[Lemma 2.1.18]{4KPX}. The argument is parallel, we adapt the argument to our setting as in the following. First we consider the situation where $A$ is Artin, then it is easy to see that the proposition is true. Then more generally we consider the copy $A\left<T_I\right>\subset \Pi_{\mathrm{an},r_{I,0},I,A}$, which allows us to go back to the previous situation by taking the reduction. Note that we still satisfy the condition since then the module is also finite over $A\left<T_I\right>$ since it is so over $A$. Then one could finish as in the one dimensional situation, see \cite[Lemma 2.1.18]{4KPX}.	
\end{proof}

\indent Then we define the following higher dimensional analog of the differential $\Omega$ over the one dimensional Robba ring relative to the ring $A$. This amounts to considering a higher dimensional residue formula for a function with $I$ variables.

\begin{definition} \mbox{\bf{(After KPX)}}
We first generalize the definition of the differentials of the Robba rings to higher dimensional situation. First consider our full Robba rings $\Pi_{\mathrm{an},\mathrm{con},I,A}$, then we are going to define the corresponding $n$-th differential of this ring relative to $A$ as the following:
\begin{displaymath}
\Omega_{\Pi_{\mathrm{an},\mathrm{con},I,A}/A}:=\wedge^{|I|} \Omega^1_{\Pi_{\mathrm{an},\mathrm{con},I,A}/A}	
\end{displaymath}
which is nothing but the $|I|$-th differential of the ring $\Pi_{\mathrm{an},\mathrm{con},I,A}$ with respect to the affinoid $A$. Here the differentials are assumed to be the continuous ones.	
\end{definition}

\indent Then for any element in the $|I|$-th differential we could define the corresponding residue of this chosen element as in the following:

\begin{definition} \mbox{\bf{(After KPX)}}
For any $f\in \Omega_{\Pi_{\mathrm{an},\mathrm{con},I,A}/A}$ of the general form of $f=\sum_{n_I}a_{n_1,...,n_I}T_1^{n_1}...T_I^{n_I}$, we define the residue of this differential of $|I|$-th order by the following formula:
\begin{displaymath}
\mathrm{res}(f):=\mathrm{res}(\sum_{n_I}a_{n_1,...,n_I}T_1^{n_1}...T_I^{n_I}):= a_{-\mathbf{1}}.	
\end{displaymath}
Here the symbol $\mathbf{1}$ is defined to be the index of $(1,...,1)$ in a uniform way.	
\end{definition}

\begin{proposition} \mbox{\bf{(After KPX \cite[Lemma 2.1.19]{4KPX})}}
The residue map defined above in our setting is well-defined, and the differential module $\Omega_{\Pi_{\mathrm{an},\mathrm{con},I,A}/A}$ is free of rank one. We define the following differential pairing over the relative differential module $\Omega_{\Pi_{\mathrm{an},\mathrm{con},I,A}/A}$:
\begin{align}
h:{\Pi_{\mathrm{an},\mathrm{con},I,A}}\times {\Pi_{\mathrm{an},\mathrm{con},I,A}} &\rightarrow A\\
(f,g)	&\mapsto \mathrm{res}(f(T_I)g(T_I)dT_1\wedge...\wedge dT_I).
\end{align}
Then the pairing now induces the following canonical isomorphisms (where the topology depends on the choices of the variables):
\begin{align}
\mathrm{Hom}_{\mathrm{con}}(\Pi_{\mathrm{an},\mathrm{con},I,A},A)\overset{\sim}{\rightarrow}	\Pi_{\mathrm{an},\mathrm{con},I,A}\\
\mathrm{Hom}_{\mathrm{con}}(\Pi_{\mathrm{an},\mathrm{con},I,A}/\sum_{k}\Pi_{\mathrm{an},k,I,A} ,A)\overset{\sim}{\rightarrow}	\Pi_{\mathrm{an},\infty,I,A}.\\
\end{align}
Here the ring $\Pi_{\mathrm{an},k,I,A}$ is defined as
\begin{align}
\Pi_{\mathrm{an},k,I,A}:=\bigcup_{r_i>0,i\neq k, r_k=\infty,k=1,...,|I|}\Pi_{\mathrm{an},r_I,I,A}.	
\end{align}

\end{proposition}

\begin{proof}
As in the one dimensional situation, from the very definitions, one could derive that the residue map is well-defined and the fact that the differential module is generated by $dT_I$. Consider our product construction and the descent process in the one dimensional situation we have that one could then check the issue of well-definedness over the corresponding multivariate convergent rings $\mathcal{E}_{\Delta}$ in \cite{4CKZ18} and \cite{4PZ19}. Then for the first statement in the isomorphisms, first one could see that one direction (from right to left) which is induced directly from the residue pairing. On the other hand, we consider the following map which is a generalization of the one dimensional situation:
\begin{align}
\mathrm{Hom}_{\mathrm{con}}(\Pi_{\mathrm{an},\mathrm{con},I,A},A) &\overset{}{\rightarrow}	\Pi_{\mathrm{an},\mathrm{con},I,A}\\\mu_I	&\mapsto \sum_{n_I}\mu_I(T_I^{-\mathbf{1}-n_I})T_1^{n_1}...T_I^{n_I}, 
\end{align}
which is easy to see that this gives an inverse mapping by the definitions of the residue pairings. For the second isomorphism, under our construction, see \cite[Proposition 5.5]{4Cre1} and the corresponding \cite[Lemma 8.5.1]{4Ked3}.

\end{proof}

\indent Then as below we are going to contact with the multi Lie groups structures which is a higher dimensional generalization of the usual pictures in \cite[Definition 2.1.20]{4KPX}. To achieve so we need to generalize the constructions in \cite[Definition 2.1.20]{4KPX} in the following sense:

\begin{definition} \mbox{\bf{(After KPX \cite[Definition 2.1.20]{4KPX})}} \label{definition1}
Uniformizing the notations as in our list of notations at the beginning, we use the notation $(C_\alpha,\ell_\alpha)_{\alpha\in I}$ to denote $I$ pairs consisting of $I$ commutative topological groups with finite torsion parts and with the free part uniformly isomorphic to $\mathbb{Z}_p$, and $I$ morphism from $\ell_\alpha:C_\alpha\rightarrow \mathbb{Q}_p$ for each $\alpha\in I$. Then based on these we consider the following generalization of the constructions in \cite[Definition 2.1.20]{4KPX}. First we are going to choose now for each $\alpha\in I$ a generator taking the form of $c_\alpha$ of the free part of the $\alpha$-th topological group for all $\alpha\in I$. Then we generalize the corresponding definitions in the usual situation as in the following. First assume that we only have torsion free topological groups in our setting for each $\alpha\in I$. Then we will use the notations $\Pi_{\mathrm{an},r_{I},I,A}(C_I),\Pi_{[s_I,r_I],I,A}(C_I)$ respectively to denote the corresponding higher dimensional period rings defined above namely $\Pi_{\mathrm{an},r_{I},I,A},\Pi_{[s_I,r_I],I,A}$ formally modified by considering the replacement of the variables with $c_\alpha-1$ for each dimension $\alpha\in I$. As in the usual situations in this situation we do not have the issue of the dependence of the choices of the generators. Also we know since by the corresponding one dimensional result we will have that
\begin{displaymath}
\left\|f(...,(1+T_\alpha)^p-1,...)\right\|_{r_I}=\left\|f(T_1,...,T_I)\right\|_{...,pr_\alpha,...}
\end{displaymath} 
for each $\alpha\in I$ which implies that certainly for all the $\alpha\in I$ together. Therefore from this fact we know that for each $\alpha\in I$:
\begin{displaymath}
\Pi_{\mathrm{an},(...,pr_{\alpha},...),I,A}(...,pC_\alpha,...)\otimes_{\mathbb{Z}_p[(...,pC_\alpha,...)]}	\mathbb{Z}_p[C_I]\overset{\sim}{\rightarrow} \Pi_{\mathrm{an},r_{I},I,A}(C_I),
\end{displaymath}
with
\begin{displaymath}
\Pi_{\mathrm{an},(...,p[s_{\alpha},r_{\alpha}],...),I,A}(...,pC_\alpha,...)\otimes_{\mathbb{Z}_p[(...,pC_\alpha,...)]}	\mathbb{Z}_p[C_I]\overset{\sim}{\rightarrow} \Pi_{[s_I,r_{I}],I,A}(C_I).
\end{displaymath}
As in \cite[Definition 2.1.20]{4KPX} one could then consider more general situation where $C_I$ may not be torsion free. In this case one considers simultaneously all small $r_I<\mathbf{1}$ or all small $s_I$ when $r_I$ are all infinity, and any $D_\alpha$ a torsion free open subgroup of $C_\alpha$ for each $\alpha\in I$, one has:
\begin{displaymath}
\Pi_{\mathrm{an},([C_1:D_1]r_1,...,[C_I:D_I]r_I),I,A}(D_I)\otimes_{\mathbb{Z}_p[D_I]}	\mathbb{Z}_p[C_I]\overset{\sim}{\rightarrow} \Pi_{\mathrm{an},r_{I},I,A}(C_I),
\end{displaymath}
with
\begin{displaymath}
\Pi_{\mathrm{an},(...,[C_\alpha:D_\alpha][s_{\alpha},r_{\alpha}],...),I,A}(...,D_\alpha,...)\otimes_{\mathbb{Z}_p[(...,D_\alpha,...)]}	\mathbb{Z}_p[C_I]\overset{\sim}{\rightarrow} \Pi_{[s_I,r_{I}],I,A}(C_I).
\end{displaymath}
Then we consider the situation for some general open subgroups $D_I'$ of the groups $D_I$ which are torsion free then for each $\alpha\in I$, and for each $r_\alpha\leq 1/\sharp D_\mathrm{tor}$, or for each $s_\alpha\leq 1/\sharp D_\mathrm{tor}$ for $r_I$ infinite, we define:
\begin{displaymath}
\Pi_{\mathrm{an},([D_1:D'_1]r_1,...,[D_I:D'_I]r_I),I,A}(D'_I)\otimes_{\mathbb{Z}_p[D'_I]}	\mathbb{Z}_p[D_I]\overset{\sim}{\rightarrow} \Pi_{\mathrm{an},r_{I},I,A}(D_I),
\end{displaymath}
with
\begin{displaymath}
\Pi_{\mathrm{an},(...,[D_\alpha:D'_\alpha][s_{\alpha},r_{\alpha}],...),I,A}(...,D'_\alpha,...)\otimes_{\mathbb{Z}_p[(...,D'_\alpha,...)]}	\mathbb{Z}_p[D_I]\overset{\sim}{\rightarrow} \Pi_{[s_I,r_{I}],I,A}(D_I).
\end{displaymath}
As in the one dimensional situation everything is well-defined with respect to the choices involved. Then based on this observation we start to consider the generalization of the differentials defined above as in the following. In our situation we consider the differential $\omega_{\ell_I}$ which is defined to be the wedge product $\omega_{\ell_1}\wedge...\wedge\omega_{\ell_I}$, where each $\omega_{\ell_\alpha}$ is defined as in the situation in \cite[Definition 2.1.20]{4KPX}. Then we could then consider any element $f\in \Omega_{\Pi_{\mathrm{an},\mathrm{con},I,A}(C_I)/A}$, and define:
\begin{displaymath}
\mathrm{res}(f):=a_{-\mathbf{1}}.	
\end{displaymath}
Then we define pairing as above in more general setting as:
\begin{align}
h:{\Pi_{\mathrm{an},\mathrm{con},I,A}(C_I)}\times {\Pi_{\mathrm{an},\mathrm{con},I,A}(C_I)} &\rightarrow A\\
(f,g) &\mapsto \mathrm{res}(fg\omega_{\ell_I}).	
\end{align}
And we also have the following analog:
\begin{align}
\mathrm{Hom}_{\mathrm{con}}({\Pi_{\mathrm{an},\mathrm{con},I,A}(C_I)},A)\overset{\sim}{\longrightarrow}	{\Pi_{\mathrm{an},\mathrm{con},I,A}(C_I)}.
%\mathrm{Hom}_{\mathrm{con}}(\Pi_{\mathrm{an},\mathrm{con},I,A}(C_I)/\Pi_{\mathrm{an},\infty,I,A}(C_I),A)\overset{\sim}{\rightarrow}	\Pi_{\mathrm{an},\infty,I,A}(C_I).
\end{align}
\end{definition}

%\newpage

\newpage

\section{Relative Frobenius Sheaves over Affinoids and Fr\'e
chet Rings}

\indent Now we introduce main objects in our study and our consideration. We first generalize the setting in \cite{4KPX} to the situation we have multi Hodge-Frobenius structure which is quite natural since we consider already some higher dimensional spaces. Then we generalize everything to the setting in \cite{4KP1}, which is also very natural in the consideration of noncommutative Iwasawa theory in the style of \cite{4KP1}. The latter will be around the Iwasawa cohomology of some $p$-adic Lie deformation of the usual $(\varphi,\Gamma)$-modules. Therefore in the second situation we will consider the usual Hodge structure but we consider the $p$-adic Lie extension, for instance one could consider the product of $\mathbb{Z}_p^\times$ with some power of $\mathbb{Z}_p^\times$ or $\mathbb{Z}_p$. This will be then highly relative since we will work over some relative Robba rings or sheaves taking the form of $\Pi_{\mathrm{an},\mathrm{con},I,\Pi_{\mathrm{an},\infty,A}(G)}$, where $\Pi_{\mathrm{an},\infty}(G)$ for instance could be Berthelot's $A_\infty$-algebra, namely the distribution over $G$.\\

\begin{setting} \label{4setting3.1}
We will work over the ring $\Pi_{\mathrm{an},\infty,A}(G)$ which is just the ring defined in the previous section which is abelian and which is the corresponding completed tensor product of the corresponding rings. To be more precise we will consider the Frobenius modules over $\Pi_{\mathrm{an},\mathrm{con},I,\Pi_{\mathrm{an},\infty,J,A}(\Gamma_{K_J})}$. Note that here $J$ could be different from $I$, but we assume that we have the relation $I\subset J$, and for those $K_\alpha$ ($\alpha\in J\backslash I$) we have that $K_\alpha=\mathbb{Q}_p$. Now set $G$ be $\Gamma_{K_J}$ (see the discussion below for the definitions). By considering the philosophy of Burns-Flach-Fukaya-Kato, this is actually already some noncommutative type Hodge structures.	
\end{setting}

\begin{remark}
The noncommutative setting is considered in the 2017 thesis \cite{4Zah1}. 	
\end{remark}

\indent Let $K_I$ be $I$ finite extension of the $p$-adic number field $\mathbb{Q}_p$. Now for each $\alpha\in I$, we consider the cyclotomic field $K_{\alpha}(\mu_{p^\infty})$, which gives rise to the the extension $K_{\alpha}(\mu_{p^\infty})/K_\alpha$ with the Galois group $G(K_{\alpha}(\mu_{p^\infty})/K_\alpha)$ which is just the corresponding $\Gamma_{K_\alpha}$ in our situation. Note that for each $K_\alpha$ we have the associated cyclotomic character $\chi_{\alpha}$ whose kernel will be denoted as usual by $H_{K_\alpha}$, for each $\alpha\in I$. And for each $\alpha$ we will use the notation $\widetilde{e}_{K_{\alpha,\infty}}$ to denote the ramification indexes along the towers, namely for each $\alpha$ this is the quotient $[K_{\alpha}(\mu_{p^\infty}):\mathbb{Q}_p(\mu_{p^\infty})]/[k_{K_{\alpha}(\mu_{p^\infty})}:k_{K_\alpha}]$. Then first step now is to add some Hodge-Frobenius structures over the period rings we established in the previous section. To do so we consider the following situation:

\begin{definition} \mbox{\bf{(After KPX \cite[Definition 2.2.2]{4KPX})}}
For each $\alpha\in I$, we\\ choose suitable uniformizer $\pi_{K_\alpha}$ in our consideration. Then we are going to use the notation $\Pi_{\mathrm{an},\mathrm{con},I,A}(\pi_{K_I})$ to denote the corresponding period ring constructed from $\Pi_{\mathrm{an},\mathrm{con},I,A}$ just by directly replacing the variables by the corresponding uniformizers as above. And similarly for other period rings we use the corresponding notations taking the same form namely $\Pi_{\mathrm{an},\mathrm{con},I,A}(\pi_{K_I})$, $\Pi_{\mathrm{an},r_{I,0},I,A}(\pi_{K_I})$, $\Pi_{\mathrm{an},[s_I,r_I],I,A}(\pi_{K_I})$. As in the one dimensional situation we consider the situation where the radii are all sufficiently small. Then one could define the multiple Frobenius actions from the multi Frobenius $\varphi_I$ over each the ring mentioned above. For suitable radii $r_I$ we define a $\varphi_I$-module to be a finite projective $\Pi_{\mathrm{an},r_{I},I,A}(\pi_{K_I})$-module with the requirement that for each $\alpha\in I$ we have $\varphi_\alpha^*M_{r_I}\overset{\sim}{\rightarrow}M_{...,r_\alpha/p,...}$ (after suitable base changes). Then we define $M:=M_{r_I}\otimes_{\Pi_{\mathrm{an},r_{I},I,A}(\pi_{K_I})}\Pi_{\mathrm{an},\mathrm{con},I,A}(\pi_{K_I})$ to define a $\varphi_I$-module over the full relative Robba ring in our situation. Furthermore we have the notion of $\varphi_I$-bundles in our situation which consists of a family of $\varphi_I$-modules $\{M_{[s_I,r_I]}\}$ where each module $M_{[s_I,r_I]}$ is defined to be finite projective over $\Pi_{\mathrm{an},[s_I,r_I],I,A}(\pi_{K_I})$ satisfying the action formula taking the form of $\varphi_\alpha^*M_{[s_I,r_I]}\overset{\sim}{\rightarrow}M_{...,[s_\alpha/p,r_\alpha/p],...}$ (after suitable base changes). Furthermore for each $\alpha$ we have the corresponding operator $\varphi_\alpha:M_{[s_I,r_I]}\rightarrow M_{...,[s_\alpha/p,r_\alpha/p],...}$ and we have the corresponding operator $\psi_\alpha$ which is defined to be $p^{-1}\varphi_\alpha^{-1}\circ\mathrm{Trace}_{M_{...,[s_\alpha/p,r_\alpha/p],...}/\varphi_\alpha(M_{[s_I,r_I]})}$. Certainly we have the corresponding operator $\psi_\alpha$ over the global section $M_{r_I}$. Note that here we require that the Hodge-Frobenius structures are commutative in the sense that all the Frobenius are commuting with each other, and they are semilinear.
\end{definition}

\indent So now based on the definitions above we can consider the corresponding objects over the ring $\Pi_{\mathrm{an},\infty,A}(G)$, which are actually complicated. The basic reason is that over the infinite level we will lose some control of the corresponding finiteness, also very significantly we lose the control of the behavior of the radius throughout the whole variation over the Stein space attached to $\Pi_{\mathrm{an},\infty,A}(G)$. Let $X(\Pi_{\mathrm{an},\infty,A}(G))$ be the corresponding Fr\'echet-Stein space attached to the group $G$, which gives rise to the following represetation:
\begin{align}
X(\Pi_{\mathrm{an},\infty,A}(G))=\varinjlim_k X(\Pi_{\mathrm{an},\infty,A}(G))_k,\\
\mathcal{O}_{X(\Pi_{\mathrm{an},\infty,A}(G))}= \varprojlim_k \mathcal{O}_{X(\Pi_{\mathrm{an},\infty,A}(G))_k	},\\
\Pi_{\mathrm{an},\infty,A}(G)=: \varprojlim_k \Pi_{\mathrm{an},\infty,A}(G)_k.
\end{align}

\begin{definition} 
Over $\square:=\Pi_{\mathrm{an},\mathrm{con},I,\Pi_{\mathrm{an},\infty,A}(G)}(\pi_{K_I})$, or\\ $\Pi_{\mathrm{an},r_{I,0},I,\Pi_{\mathrm{an},\infty,A}(G)}(\pi_{K_I})$, or $\Pi_{\mathrm{an},[s_I,r_I],I,\Pi_{\mathrm{an},\infty,A}(G)}(\pi_{K_I})$ as in the above (with the notations for the radii as considered above) we have the notion of $\varphi_I$-modules over $\square$. Be careful that we will define a $\varphi_I$-module $M$ over $\square$ to be the following projective limit:
\begin{align}
M:=\varprojlim_k M_k	
\end{align}
where $M_k$ is a corresponding object over $\square_k$ as defined above as a corresponding $\varphi_I$-module $M$ (finite projective) over $\square_k$. Here $\square_k:=\Pi_{\mathrm{an},\mathrm{con},I,\Pi_{\mathrm{an},\infty,A}(G)_k}(\pi_{K_I})$, or\\ $\Pi_{\mathrm{an},r_{I,0},I,\Pi_{\mathrm{an},\infty,A}(G)_k}(\pi_{K_I})$, or $\Pi_{\mathrm{an},[s_I,r_I],I,\Pi_{\mathrm{an},\infty,A}(G)_k}(\pi_{K_I})$ respectively. One defines the $\varphi_I$-bundles over
\begin{center}
$\Pi_{\mathrm{an},\mathrm{con},I,\Pi_{\mathrm{an},\infty,A}(G)}(\pi_{K_I})$ 
\end{center} 
in the same fashion.
\end{definition}

\indent We also have the corresponding notions of families of the corresponding objects considered above in the Fr\'echet-Stein situation.

\begin{definition} 
Over $\square:=\Pi_{\mathrm{an},\mathrm{con},I,X(\Pi_{\mathrm{an},\infty,A}(G))}(\pi_{K_I})$, or $\Pi_{\mathrm{an},r_{I,0},I,X(\Pi_{\mathrm{an},\infty,A}(G))}(\pi_{K_I})$, or\\
 $\Pi_{\mathrm{an},[s_I,r_I],I,X(\Pi_{\mathrm{an},\infty,A}(G))}(\pi_{K_I})$ as in the above (with the notations for the radii as considered above) we have the notion of $\varphi_I$-modules over $\square$. Be careful that we will define a $\varphi_I$-module $M$ over $\square$ to be the following the sheaf taking each quasi-compact $Y\subset X(\Pi_{\mathrm{an},\infty,A}(G))$ to $M(Y)$ such that this is a corresponding object over $\square_Y$ as defined above as a corresponding $\varphi_I$-module $M$ (finite projective) over $\square_Y$. Here $\square_Y:=\Pi_{\mathrm{an},\mathrm{con},I,\mathcal{O}_{X(\Pi_{\mathrm{an},\infty,A}(G))}(Y)}(\pi_{K_I})$, or $\Pi_{\mathrm{an},r_{I,0},I,\mathcal{O}_{X(\Pi_{\mathrm{an},\infty,A}(G))}(Y)}(\pi_{K_I})$, or\\
 $\Pi_{\mathrm{an},[s_I,r_I],I,\mathcal{O}_{X(\Pi_{\mathrm{an},\infty,A}(G))}(Y)}(\pi_{K_I})$ respectively. One defines the families of $\varphi_I$-bundles over $\Pi_{\mathrm{an},\mathrm{con},I,X(\Pi_{\mathrm{an},\infty,A}(G))}(\pi_{K_I})$ in the same fashion. Here one can take the corresponding $Y$ as some member in the filtration:
\begin{displaymath}
X(\Pi_{\mathrm{an},\infty,A}(G))=\varinjlim_k X(\Pi_{\mathrm{an},\infty,A}(G))_k.	
\end{displaymath}
 
\end{definition}

\begin{proposition} \mbox{\bf{(After KPX \cite[Proposition 2.2.7]{4KPX})}} \label{proposition3.5}
We have that natural functor from families of $\varphi_I$-modules (in the sense of the previous definition without the finite condition on the global section) over $\Pi_{\mathrm{an},r_{I,0},I,\Pi_{\mathrm{an},\infty,A}(G)}(\pi_{K_I})$ to the corresponding families of $\varphi_I$-bundles (in the sense of the previous definition without the finite condition on the global section) is an equivalence.	Here by a family of $\varphi_I$-module we mean a $\varphi_I$-module over $\Pi_{\mathrm{an},r_{I,0},I,X(\Pi_{\mathrm{an},\infty,A}(G))}(\pi_{K_I})$, which is the same to a family of $\varphi_I$-bundles over $\Pi_{\mathrm{an},r_{I,0},I,\Pi_{\mathrm{an},\infty,A}(G)}(\pi_{K_I})$.
\end{proposition}

\begin{proof}
Here we need to consider the corresponding Fr\'echet sheaves as in the previous definition since in our situation the relative Robba rings are defined over some non-noetherian Stein space. The statement is obviously true for each $M_n$ as in the usual situation. First the faithfulness is straightforward, then for the surjectivity part, we only need to show that the global section is finite projective, which is a directly consequence of \cref{4prop1} as long as in our situation one considers the simultaneous application of our multiple Frobenius to shift around to vary through all the intervals.	
\end{proof}

\begin{definition} \mbox{\bf{(After KPX \cite[Definition 2.2.12]{4KPX})}}
Then we add the corresponding more Lie group action in our setting, namely the multi semiliear $\Gamma_{K_I}$-action. Again in this situation we require that all the actions of the $\Gamma_{K_I}$ are commuting with each other and with all the semilinear Frobenius actions defined above. We define the corresponding $(\varphi_I,\Gamma_I)$-modules over $\Pi_{\mathrm{an},r_{I,0},I,\Pi_{\mathrm{an},\infty,A}(G)}(\pi_{K_I})$ to be finite projective module with mutually commuting semilinear actions of $\varphi_I$ and $\Gamma_{K_I}$. For the latter action we require that to be continuous.
\end{definition}

\indent The continuity condition here could be made more explicit, as in \cite[Proposition 2.2.14]{4KPX}. In our setting, we just need to generalize the result in \cite[Proposition 2.2.14]{4KPX} to multiple actions from $\Gamma_{K_I}$.

\begin{lemma} \mbox{\bf{(After KPX \cite[Proposition 2.2.14]{4KPX})}}\\
For suitable $r_{I,0}$, consider $M_{r_{I,0}}$ a module (carrying a semilinear action of $\Gamma_{K_I}$) of finite type over $\Pi_{\mathrm{an},r_{I,0},I,\Pi_{\mathrm{an},\infty,A}(G)}(\pi_{K_I})$ with the generating set $\mathbf{e}_1,...,\mathbf{e}_m$ (namely we look at a sheaf over $\Pi_{\mathrm{an},r_{I,0},I,X(\Pi_{\mathrm{an},\infty,A}(G))}(\pi_{K_I})$ with global section finitely generated). We will follow \cite[Proposition 2.2.14]{4KPX} to differentiate the corresponding Banach module norm $|.|_{[s_I,r_{I,0}]}$ and the operator norm	$\left\|.\right\|_{[s_I,r_{I,0}]}$ when we consider the module over the polydisc with respect to the corresponding multi-interval (and the corresponding member in the Fr\'echet family). Then we have: I. The action of the group $\Gamma_{K_I}$ is continuous if and only if that for each $\gamma\in \Gamma_{K_I}$ and each element in the generating set we have $\lim_{\gamma\rightarrow 1} \gamma \mathbf{e}_j=\mathbf{e}_j$; The continuity of the action implies that: 
II. $\lim_{\gamma\rightarrow 1}\left\|\gamma-1\right\|_{[s_I,r_{I,0}]}=0$ for any $s_I$; 
III. Whenever we have the subgroups (with the notations in \cite[Proposition 2.2.14]{4KPX}) $\Gamma_{n_I}\subset \Gamma_{K_I}$ ($\forall n_I\geq 0$), one can find separately $g_{\alpha,s_I}$ with the fact that $\left\|(\gamma'_{n_\alpha}-1)^{g_{\alpha,s_I}}\right\|_{[s_I,r_{I,0}]}< p^{-1/(p-1)}$ for each $\alpha\in I$.
\end{lemma}

\begin{proof}
First we need to work over some affinoid domain with respect to the Fr\'echet-Stein algebra in the coefficient, but this will be reduced to the corresponding affinoid situation. We then follow the strategy of \cite[Proposition 2.2.14]{4KPX}, where we just need to consider each group in the product $\Gamma_{K_I}$. The first statement could be derived in the same way as in \cite[Proposition 2.2.14]{4KPX}. For the second statement therefore we assume the limit identity in I, then prove the corresponding results on the identity in II. As in the proof in \cite[Proposition 2.2.14]{4KPX}, for I and II, we reduce ourselves to the situation where the base field is $\mathbb{Q}_p$ for each $K_\alpha$ from the induction properties, and correspondingly we reduce us to the situation where all the uniformizers $\pi_I$ are just same as a product taking the form of $\pi^I$ as considered in the proof in \cite[Proposition 2.2.14]{4KPX}. First is the following observation on the Newton-Leibniz style formula for each decomposition of elements with respect to given basis $\sum r_i e_i$ and for any $\gamma\in \Gamma_{K_I}$:
\begin{displaymath}
(\gamma-1)(r_i e_i)=\gamma(r_i)	(\gamma-1)(e_i)+(\gamma-1)r_i e_i.
\end{displaymath}
It suffices to look at the corresponding last term. We then further reduce ourselves to the situation where the module is just free of rank one over the base ring, which directly implies it suffices to consider the element in the ring taking the general form of $\pi_1^{n_1}...\pi_I^{n_I}$. Since we consider further the group element $\gamma_\alpha$ separately for each $\alpha\in I$, the result follows from then directly observation that $\gamma_\alpha(\pi_\beta)-\pi_\beta=f(\pi_\beta)((1+\pi_\beta)^{ap^k}-1)$ as well as $\gamma_\alpha(\pi_\beta^{-1})-\pi_\beta^{-1}=f'(\pi_\beta)((1+\pi_\beta)^{ap^k}-1)$ ($\alpha=\beta$ otherwise we have trivial action) where $f'$ actually depends on $k$ but could be controlled by some $g'$ which does not depend. Then from this for each $\beta,\alpha$ one gets the result. For the third statement we again argue by separately considering the corresponding group $\Gamma_\alpha$ where $\alpha\in I$. By first stage of reduction, we reduce ourselves to the situation where the base fields are all $\mathbb{Q}_p$ by induction and the uniformizer are all $\pi$ for each $\alpha \in I$. Since for each $\alpha$, as in \cite[Proposition 2.2.14]{4KPX} the Banach module norm $|.|_{[s_I,r_{I,0}]}$ is invariant under the action of $\Gamma_{n_\alpha}$ for some $n_\alpha$. Note that this is based on the corresponding observation that:
\begin{align}
\gamma_\alpha(\pi^{i_1}_1...\pi_\alpha...\pi^{i_I}_I)-\pi^{i_1}_1...\pi_\alpha...\pi^{i_I}_I=\pi^{i_1}_1...f(\pi_\alpha)((1+\pi_\alpha)^{ap^k}-1)...\pi^{i_I}_I\\
\gamma_\alpha(\pi^{i_1}_1...\pi^{-1}_\alpha...\pi^{i_I}_I)-\pi^{i_1}_1...\pi_\alpha^{-1}...\pi^{i_I}_I=\pi^{i_1}_1...f'(\pi_\alpha)((1+\pi_\alpha)^{ap^k}-1)...\pi^{i_I}_I,	
\end{align}
and note that $|.|_{[s_I,r_{I,0}]}(\pi^{i_1}_1...((1+\pi_\alpha)^{ap^k}-1)...\pi^{i_I}_I)\leq |.|_{[s_I,r_{I,0}]}(\pi^{i_1}_1...\pi_\alpha...\pi^{i_I}_I)$. Then we consider the corresponding average which is invariant as in \cite[Proposition 2.2.14]{4KPX} for some $n_\alpha'\geq n_\alpha$:
\begin{displaymath}
|.|^{-}_{[s_I,r_{I,0}]}:=\frac{1}{[\Gamma_{n_\alpha}:\Gamma_{n_\alpha}']}\sum_{\widetilde{\gamma}_\alpha\in \Gamma_{n_\alpha}/\Gamma_{n_\alpha}'}|.|_{[s_I,r_{I,0}]}\circ \widetilde{\gamma}_\alpha	
\end{displaymath}
which is equivalent to the original one:
\begin{displaymath}
C^{-1}|.|^{-}_{[s_I,r_{I,0}]}\leq |.|_{[s_I,r_{I,0}]} \leq C|.|^{-}_{[s_I,r_{I,0}]}	
\end{displaymath}
for some positive number $C>0$. First for the average norm we have that by taking suitable $n_\alpha'$ one could get $\left\|\gamma'_{n_\alpha'}-1\right\|^{-}_{[s_I,r_{I,0}]}$ is small enough, then take the corresponding $p^{n_\alpha'-n_\alpha}$-power of the element $\gamma'_{n_\alpha}-1$ one could further make this term smaller than $p^{-1/(p-1)}$ under the corresponding operator norm, namely the term $(\gamma'_{n_\alpha'}-1)^{p^{n_\alpha'-n_\alpha}}$ as in \cite[Proposition 2.2.14]{4KPX}. Then by the inequality above between the original norms and the average ones, we have that by taking suitable power $g_{s_\alpha}$ one could make $\left\|(\gamma'_{n_\alpha}-1)^{g_{s_\alpha}}\right\|_{[s_I,r_{I,0}]}$ smaller than $p^{-1/(p-1)}$. To be more precise by raising the power $g_{n_\alpha}$ one will have $\left\|(\gamma'_{n_\alpha}-1)^{g_{s_\alpha}}\right\|_{[s_I,r_{I,0}]}\leq C\left\|(\gamma'_{n_\alpha}-1)^{g_{s_\alpha}}\right\|_{[s_I,r_{I,0}]}^{-}\leq CC_{g_{s_\alpha}}$ where $C_{g_{s_\alpha}}$ could be made sufficiently small along this process.
\end{proof}

%\newpage

\newpage\section{Der Endlichkeitssatz}

\subsection{Finiteness through $p$-adic Functional Analaysis} \label{section4.1}

\indent In this section we are going to study the cohomology of our relative Frobenius modules defined in the previous sections. We start from a $(\varphi_I,\Gamma_I)$-module over $\Pi_{\mathrm{an},\mathrm{con},I,A}(\pi_{K_I})$ as well a $(\psi_I,\Gamma_I)$-module over $\Pi_{\mathrm{an},\mathrm{con},I,A}(\pi_{K_I})$. We are going to use $M$ to denote such a module. First we are going to use the notation $\triangle_I$ to denote the $p$-torsion part of $\Gamma_{I}$ in the product form, which will be used in the definitions of the corresponding cohomologies as below.

\begin{assumption} \label{assumption1}
We choose to consider the story in which $I$ consists of two elements. We note that actually everything could be carried over to more general setting. We will then use the notation $\{1,2\}$ to denote our set which indicates the corresponding Frobenius actions and Lie group actions with respect to each index.	
\end{assumption}

\begin{definition} \mbox{\bf{(After KPX \cite[2.3]{4KPX})}}
For any $(\varphi_I,\Gamma_I)$-module $M$ over\\ $\Pi_{\mathrm{an},\mathrm{con},I,A}(\pi_{K_I})$ we define the complex $C^\bullet_{\Gamma_I}(M)$ of $M$ to be the total complex of the following complex through the induction:
\[
\xymatrix@R+0pc@C+0pc{
[ C^\bullet_{\Gamma_{I\backslash\{I\}}}(M) \ar[r]^{\Gamma_I} \ar[r] \ar[r]  & C^\bullet_{\Gamma_I\backslash\{I\}}(M)]
.}
\]
Then we define the corresponding double complex $C^\bullet_{\varphi_I}C^\bullet_{\Gamma_I}(M)$ again by taking the corresponding totalization of the following complex through induction:
\[
\xymatrix@R+0pc@C+0pc{
[C^\bullet_{\varphi_{I\backslash\{I\}}}C^\bullet_{\Gamma_{I}}(M) \ar[r]^{\varphi_I} \ar[r] \ar[r]  & C^\bullet_{\varphi_{I\backslash\{I\}}}C^\bullet_{\Gamma_I}(M)]
.}
\]
Then we define the complex $C^\bullet_{\varphi_I,\Gamma_I}(M)$ to be the totalization of the double complex defined above, which is called the complex of a $(\varphi_I,\Gamma_I)$-module $M$. Similarly we define the corresponding complex $C^\bullet_{\psi_I,\Gamma_I}(M)$ to be the totalization of the double complex define  d in the same way as above by replacing $\varphi_I$ by then the operators $\psi_I$. For later comparison we also need the corresponding $\psi_I$-cohomology, which is to say the complex $C^\bullet_{\psi_I}(M)$ which could be defined by the totalization of the following complex through induction:
\[
\xymatrix@R+0pc@C+0pc{
[ C^\bullet_{\psi_{I\backslash\{I\}}}(M) \ar[r]^{\psi_I} \ar[r] \ar[r]  & C^\bullet_{\psi_I\backslash\{I\}}(M)]
.}
\]
\end{definition}

\indent One can actually compare the two total complexes defined above in the following way, as in the situation where $I$ is a singleton. First for the $(\varphi_I,\Gamma_I)$-complex we have the following initial double complex:
\[
\xymatrix@R+6pc@C+0pc{
 &0\ar[d] \ar[d] \ar[d]  &0 \ar[d] \ar[d] \ar[d]  & 0 \ar[d] \ar[d] \ar[d] \\
0\ar[r] \ar[r] \ar[r] &M^{\triangle_I}\ar[d]^{(\gamma_1-1,\gamma_2-1)} \ar[d] \ar[d] \ar[r]^{(\varphi_1-1,\varphi_2-1)} \ar[r] \ar[r] &M^{\triangle_I}\oplus M^{\triangle_I}\ar[d] \ar[d] \ar[d]\ar[r]^{(\varphi_2-1)+(1-\varphi_1)} \ar[r] \ar[r]  & M^{\triangle_I}\ar[d] \ar[d] \ar[d] \ar[r] \ar[r] \ar[r] &0\\
0\ar[r] \ar[r] \ar[r] &M^{\triangle_I}\oplus M^{\triangle_I}\ar[d]^{(\gamma_2-1)+(1-\gamma_1)} \ar[d] \ar[d] \ar[r] \ar[r] \ar[r] &M^{\triangle_I}\oplus M^{\triangle_I}\oplus M^{\triangle_I}\oplus M^{\triangle_I}\ar[d] \ar[d] \ar[d]\ar[r] \ar[r] \ar[r]  & M^{\triangle_I}\oplus M^{\triangle_I}\ar[d] \ar[d] \ar[d] \ar[r] \ar[r] \ar[r] &0\\
0\ar[r] \ar[r] \ar[r] &M^{\triangle_I}\ar[d] \ar[d] \ar[d] \ar[r] \ar[r] \ar[r] &M^{\triangle_I}\oplus M^{\triangle_I}\ar[d] \ar[d] \ar[d]\ar[r] \ar[r] \ar[r]  & M^{\triangle_I}\ar[d] \ar[d] \ar[d] \ar[r] \ar[r] \ar[r] &0\\
&0 &0  & 0
.}
\]
Here we followed the corresponding convention as in \cite{4KPX}. And then correspondingly we have the following double complex for the corresponding $(\psi_I,\Gamma_I)$-module:
\begin{center}
\[
\xymatrix@R+6pc@C+0pc{
 &0\ar[d] \ar[d] \ar[d]  &0 \ar[d] \ar[d] \ar[d]  & 0 \ar[d] \ar[d] \ar[d] \\
0\ar[r] \ar[r] \ar[r] &M^{\triangle_I}\ar[d]^{(\gamma_1-1,\gamma_2-1)} \ar[d] \ar[d] \ar[r]^{(\psi_1-1,\psi_2-1)} \ar[r] \ar[r] &M^{\triangle_I}\oplus M^{\triangle_I}\ar[d] \ar[d] \ar[d]\ar[r]^{(\psi_2-1)+(1-\psi_1)} \ar[r] \ar[r]  & M^{\triangle_I}\ar[d] \ar[d] \ar[d] \ar[r] \ar[r] \ar[r] &0\\
0\ar[r] \ar[r] \ar[r] &M^{\triangle_I}\oplus M^{\triangle_I}\ar[d]^{(\gamma_2-1)+(1-\gamma_1)} \ar[d] \ar[d] \ar[r] \ar[r] \ar[r] &M^{\triangle_I}\oplus M^{\triangle_I}\oplus M^{\triangle_I}\oplus M^{\triangle_I}\ar[d] \ar[d] \ar[d]\ar[r] \ar[r] \ar[r]  & M^{\triangle_I}\oplus M^{\triangle_I}\ar[d] \ar[d] \ar[d] \ar[r] \ar[r] \ar[r] &0\\
0\ar[r] \ar[r] \ar[r] &M^{\triangle_I}\ar[d] \ar[d] \ar[d] \ar[r] \ar[r] \ar[r] &M^{\triangle_I}\oplus M^{\triangle_I}\ar[d] \ar[d] \ar[d]\ar[r] \ar[r] \ar[r]  & M^{\triangle_I}\ar[d] \ar[d] \ar[d] \ar[r] \ar[r] \ar[r] &0\\
&0 &0  & 0
.}
\]
\end{center}

\indent In our situation in order to relate the two double complexes we consider the corresponding map induced from the following maps for the two first lines of two big complexes where one complex is given as above, while the other is the following one:

\begin{center}
\[
\xymatrix@R+6pc@C+0pc{
 &0\ar[d] \ar[d] \ar[d]  &0 \ar[d] \ar[d] \ar[d]  & 0 \ar[d] \ar[d] \ar[d] \\
0\ar[r] \ar[r] \ar[r] &M^{\triangle_I}\ar[d]^{(\gamma_1-1,\gamma_2-1)} \ar[d] \ar[d] \ar[r]^{(\varphi_1-1,\psi_2-1)} \ar[r] \ar[r] &M^{\triangle_I}\oplus M^{\triangle_I}\ar[d] \ar[d] \ar[d]\ar[r]^{(\psi_2-1)+(1-\varphi_1)} \ar[r] \ar[r]  & M^{\triangle_I}\ar[d] \ar[d] \ar[d] \ar[r] \ar[r] \ar[r] &0\\
0\ar[r] \ar[r] \ar[r] &M^{\triangle_I}\oplus M^{\triangle_I}\ar[d]^{(\gamma_2-1)+(1-\gamma_1)} \ar[d] \ar[d] \ar[r] \ar[r] \ar[r] &M^{\triangle_I}\oplus M^{\triangle_I}\oplus M^{\triangle_I}\oplus M^{\triangle_I}\ar[d] \ar[d] \ar[d]\ar[r] \ar[r] \ar[r]  & M^{\triangle_I}\oplus M^{\triangle_I}\ar[d] \ar[d] \ar[d] \ar[r] \ar[r] \ar[r] &0\\
0\ar[r] \ar[r] \ar[r] &M^{\triangle_I}\ar[d] \ar[d] \ar[d] \ar[r] \ar[r] \ar[r] &M^{\triangle_I}\oplus M^{\triangle_I}\ar[d] \ar[d] \ar[d]\ar[r] \ar[r] \ar[r]  & M^{\triangle_I}\ar[d] \ar[d] \ar[d] \ar[r] \ar[r] \ar[r] &0\\
&0 &0  & 0
}
\]
\end{center}

 (which is a direct generalization of the corresponding morphisms in one dimensional situation in \cite[Definition 2.3.3]{4KPX}):

\[
\xymatrix@R+5pc@C+4pc{
0\ar[r] \ar[r] \ar[r] &M^{\triangle_I}\ar[d]^{\mathrm{Id}} \ar[d] \ar[d] \ar[r]^{(\varphi_1-1,\gamma_1-1)} \ar[r] \ar[r] &M^{\triangle_I}\oplus M^{\triangle_I}\ar[d]^{(-\psi_1,\mathrm{Id})} \ar[d] \ar[d]\ar[r]^{(\gamma_1-1)+(1-\varphi_1)} \ar[r] \ar[r]  & M^{\triangle_I}\ar[d]^{-\psi_1} \ar[d] \ar[d] \ar[r] \ar[r] \ar[r] &0\\
0\ar[r] \ar[r] \ar[r] &M^{\triangle_I} \ar[r]^{(\psi_1-1,\gamma_1-1)} \ar[r] \ar[r] &M^{\triangle_I}\oplus M^{\triangle_I}\ar[r]^{(\gamma_1-1)+(1-\psi_1)} \ar[r] \ar[r]  & M^{\triangle_I} \ar[r] \ar[r] \ar[r] &0
.}
\]

We will denote the corresponding morphism in between the two total complex by:

\[
\xymatrix@R+0pc@C+0pc{
\Psi_1: C^\bullet_{\varphi_1,\psi_2,\Gamma_1,\Gamma_2}(M) \ar[r] \ar[r] \ar[r]  & C^\bullet_{\psi_1,\psi_2,\Gamma_1,\Gamma_2}(M)
.}
\]

\begin{proposition} \mbox{\bf{(After KPX \cite[Proposition 2.3.6,Lemma 3.1.2]{4KPX})}}
The morphism  
\[
\xymatrix@R+0pc@C+0pc{
\Psi_1': C^\bullet_{\varphi_1,\Gamma_1}(M) \ar[r] \ar[r] \ar[r]  & C^\bullet_{\psi_1,\Gamma_1}(M)
}
\]
is a quasi-isomorphism. This will imply that $\Psi_1$ is a quasi-isomorphism.	
\end{proposition}

\begin{proof}
First the problem is now reduced to the checking of the injectivity of the morphism $\Psi_1'$ defined above. Note now that the kernel of the above morphism is now just the corresponding complex which consists of the invariance taking the form of $M^{\triangle_I,\psi_1=0}$, so it suffices to prove that the operator $\gamma_1-1$ in our situation bas the corresponding property of being invertible. Therefore we will follow the argument in \cite[Proposition 2.3.6]{4KPX}, which will be just a direct generalization of the corresponding one dimensional result in \cite[Proposition 2.3.6]{4KPX}. First we consider the reduction of the proof to the situation where all the finite extensions $K_I$ over $\mathbb{Q}_p$ are all just the fields $\mathbb{Q}_p$, which could be checked through the induction. Then we assume that we have a corresponding $\Pi_{\mathrm{an},r_{I,0},I,A}(\pi^I)$-module $M_{r_{I,0}}^{\triangle_{\mathbb{Q}_p}}$ and the corresponding invariance $M_{r_{I,0}}^{\triangle_{\mathbb{Q}_p},\psi_1=0}$.
Then we consider the suitable multi integers $n_I\geq 1$ which gives rise to (for suitable multi radii $s_I$ and $r_I$ smaller than the corresponding original multi radii $r_{I,0}$) the corresponding decomposition of the section 
\begin{displaymath}
M_{r_{I,0}}^{\triangle_{\mathbb{Q}_p},\psi_1=0}\otimes_{\Pi_{\mathrm{an},r_{I,0},I,A}(\pi^I)} \Pi_{\mathrm{an},[s_I/p^{n_I},r_I/p^{n_I}],I,A}(\pi^I)	
\end{displaymath}
into the corresponding summation of the corresponding components taking the form of:
\begin{displaymath}
(1+\pi_1)^{i_1}\varphi_1^{n_1}\varphi_2^{n_2}M_{[s_I,r_I]}^{\triangle_{\mathbb{Q}_p}},	
\end{displaymath}
for suitable integers $i_1$ in our generalized situation. Then it now suffices to show that for each $n_I$, the action of the corresponding operators $\gamma_{n_I}-1$ are invertible over each member in the decomposition as above, which is to say the action of each $\gamma_\alpha-1$ is invertible over each $(1+\pi_1)\varphi_1^{n_1}M_{[s_I,r_I]}^{\triangle_{\mathbb{Q}_p}}$, which implies that the $A[\gamma_{n_I}]$-module structure could be promoted by continuation to the corresponding $\Pi_{\mathrm{an},r_{I,0},I,A}(\Gamma_{n_I})$-module structure. Now we look at each member in the decomposition as above, and do the standard computation as in \cite[Lemma 3.1.2]{4KPX} as below (which is generalization since we are working with higher dimensional situation) for each $m$ in each member participating in the decomposition as above:
\begin{align}
(\gamma_{n_1}-1)&(1+\pi_1)\varphi_1^{n_1}\varphi_2^{n_2}(m)\\
&=(\gamma_{n_1}-1)(1+\pi_1)\varphi_1^{n_1}\varphi_2^{n_2}(m)\\
&=(1+\pi_1)^{p^{n_1+1}}\varphi_1^{n_1}\varphi_2^{n_2}\gamma_{n_1}(m)-(1+\pi_1)\varphi_1^{n_1}\varphi_2^{n_2}(m)\\
&=(1+\pi_1)((1+\pi_1)^{p^{n_1}}\varphi_1^{n_1}\varphi_2^{n_2}\gamma_{n_1}-\varphi_1^{n_1}\varphi_2^{n_2})(m)\\
&=(1+\pi_1)\varphi_1^{n_1}((1+\pi_1)\varphi_2^{n_2}\gamma_{n_1}-\varphi_2^{n_2})(m)\\
&=(1+\pi_1)\varphi_1^{n_1}\varphi_2^{n_2}\pi_1(1+(1+\pi_1)/\pi_1(\gamma_{n_1}-1))(m)
\end{align}
which implies that the inverse could be taken to be:
\begin{align}
[(1+\pi_1)\varphi_1^{n_1}&\varphi_2^{n_2}\pi_1(1+(1+\pi_1)/\pi_1(\gamma_{n_1}-1))]^{-1}(m)	\\
&=\frac{1}{(1+\pi_1)\varphi_1^{n_1}\varphi_2^{n_2}\pi_1}\cdot \frac{1}{1+(1+\pi_1)/\pi_1(\gamma_{n_1}-1)}m\\
&=\frac{1}{(1+\pi_1)\varphi_1^{n_1}\varphi_2^{n_2}\pi_1}\cdot \sum_{k\geq 0}(-1)^k((1+\pi_1)/\pi_1(\gamma_{n_1}-1))^k m.
\end{align}
\end{proof}

\indent Now in the same fashion one can define the corresponding morphism:
\[
\xymatrix@R+0pc@C+0pc{
\Psi_2: C^\bullet_{\varphi_1,\varphi_2,\Gamma_I}(M) \ar[r] \ar[r] \ar[r]  & C^\bullet_{\varphi_1,\psi_2,\Gamma_I}(M)
}
\]

\begin{lemma}
This morphism $\Psi_2$ is also a quasi-isomorphism. Therefore by composition we have the following quasi-isomorphism:
\[
\xymatrix@R+0pc@C+0pc{
\Psi_I: C^\bullet_{\varphi_1,\varphi_2,\Gamma_I}(M) \ar[r] \ar[r] \ar[r]  & C^\bullet_{\psi_1,\psi_2,\Gamma_I}(M).
}
\]
%induced by the following diagram:
%\[
%\xymatrix@R+6pc@C+4pc{
%0\ar[r] \ar[r] \ar[r] &M^{\triangle_I}\ar[d]^{\mathrm{Id}} \ar[d] \ar[d] \ar[r]^{(\varphi_1-1,\varphi_2-1)} \ar[r] \ar[r] &M^{\triangle_I}\oplus M^{\triangle_I}\ar[d]^{(-\psi_2,\mathrm{Id})} \ar[d] \ar[d]\ar[r]^{(\varphi_2-1)+(1-\varphi_1)} \ar[r] \ar[r]  & M^{\triangle_I}\ar[d]^{-\psi_2} \ar[d] \ar[d] \ar[r] \ar[r] \ar[r] &0\\
%0\ar[r] \ar[r] \ar[r] &M^{\triangle_I}\ar[d]^{\mathrm{Id}} \ar[d] \ar[d] \ar[r]^{(\varphi_1-1,\psi_2-1)} \ar[r] \ar[r] &M^{\triangle_I}\oplus M^{\triangle_I}\ar[d]^{(-\psi_1,\mathrm{Id})} \ar[d] \ar[d]\ar[r]^{(\psi_2-1)+(1-\varphi_1)} \ar[r] \ar[r]  & M^{\triangle_I}\ar[d]^{-\psi_1} \ar[d] \ar[d] \ar[r] \ar[r] \ar[r] &0\\
%0\ar[r] \ar[r] \ar[r] &M^{\triangle_I} \ar[r]^{(\psi_1-1,\psi_2-1)} \ar[r] \ar[r] &M^{\triangle_I}\oplus M^{\triangle_I}\ar[r]^{(\psi_2-1)+(1-\psi_1)} \ar[r] \ar[r]  & M^{\triangle_I} \ar[r] \ar[r] \ar[r] &0
%.}
%\]

\end{lemma}

\indent Then now we are at the position to consider the finiteness of the corresponding involved cohomologies. In our situation we consider the following derived categories:

\begin{definition} \mbox{\bf{(After KPX \cite[Notation 4.1.2]{4KPX})}}
Now we consider the following derived categories. The first one is the sub-derived category consisting of all the objects in $\mathbb{D}(A)$ which are quasi-isomorphic to those bounded above complexes of finite projective modules over the ring $A$. We denote this category (which is defined in the same way as in \cite[Notation 4.1.2]{4KPX}) by $\mathbb{D}^-_\mathrm{perf}(A)$. Over the larger ring $\Pi_{\mathrm{an},\infty,I,A}(\Gamma_{K_I})$ we also have the sub-derived category of $\mathbb{D}(\Pi_{\mathrm{an},\infty,I,A}(\Gamma_{K_I}))$ consisting of all those objects in $\mathbb{D}(\Pi_{\mathrm{an},\infty,I,A}(\Gamma_{K_I}))$ which are quasi-isomorphic to those bounded above complexes of finite projective modules now over the ring $\Pi_{\mathrm{an},\infty,I,A}(\Gamma_{K_I})$. We denote this by $\mathbb{D}^-_\mathrm{perf}(\Pi_{\mathrm{an},\infty,I,A}(\Gamma_{K_I}))$.
\end{definition}

\indent Then we will investigate the complexes attached to any finite projective $(\varphi_I,\Gamma_I)$-module $M$ over $\Pi_{\mathrm{an},\mathrm{con},I,A}(\pi_{K_I})$:
\begin{displaymath}
C^\bullet_{\varphi_I,\Gamma_I}(M), C^\bullet_{\psi_I,\Gamma_I}(M),C^\bullet_{\psi_I}(M),	
\end{displaymath}
which are living inside the corresponding derived categories:
\begin{displaymath}
\mathbb{D}(A),\mathbb{D}(A),\mathbb{D}(\Pi_{\mathrm{an},\infty,I,A}(\Gamma_{K_I})).	
\end{displaymath}

\begin{proposition} \mbox{\bf{(After KPX \cite[Proposition 4.2.3]{4KPX})}} \label{proposition4.4}
Now for any multiplicative character $\eta_I$ of the group $\Gamma_{K_I}$, we have the following isomorphism in the corresponding derived categories defined above:
\begin{align}
C^\bullet_{\psi_I}(M)\otimes^{\mathbb{L}}\Pi_{\mathrm{an},\infty,I,A}(\Gamma_{K_I})/\mathfrak{m}_\eta \Pi_{\mathrm{an},\infty,I,A}(\Gamma_{K_I})&\overset{\sim}{\longrightarrow} C^\bullet_{\psi_I,\Gamma_I}(M\otimes M_{\eta^{-1}})\\
	&\overset{\sim}{\longrightarrow}C^\bullet_{\varphi_I,\Gamma_I}(M\otimes M_{\eta^{-1}}).
\end{align}	
\end{proposition}

\indent Then we are going to study the finiteness of the cohomologies defined above in more details. The methods in \cite{4KPX} requires some deep results on the slope filtrations for modules over the one-dimensional Robba rings, where allows \cite{4KPX} to reduce the proof in some sense to the \'etale situation. It looks like we have no chance to apply the ideas in \cite{4KPX}. Here we adopt the ideas in some development from Kedlaya-Liu in \cite{4KL3}. To be more precise we are going to apply some ideas in \cite{4KL3} around the application of some $p$-adic functional analysis.

%\begin{proposition}
%With the notation and setting up as above we have $C^\bullet_{\psi_I}(M)$ lives in $ \mathbb{D}_\mathrm{perf}^\flat(\Pi_{\mathrm{an},\infty,I,A}(\Gamma_{K_I}))$.	
%\end{proposition}

\begin{proposition} \label{proposition4.7}
With the notation and setting up as above we have $h^\bullet(C_{\psi_I}(M))$ are coadmissible over $\Pi_{\mathrm{an},\infty,I,A}(\Gamma_{K_I})$.	
\end{proposition}

\begin{proof}
The original situation for $A=\mathbb{Q}_p$ and $M$ is \'etale is proved by \cite{4CKZ18} and \cite{4PZ19}, which is a direct consequence of Galois cohomologies. Here we adopt the Cartan-Serre method involving the completely continuous maps within the $p$-adic functional analysis. By simultaneous induction to the base cases for the field $\mathbb{Q}_p$ from all $K_I$, we can now assume that we are considering the situation where all $K_I$ are $\mathbb{Q}_p$. Then for the complex $C^\bullet_{\psi_I}(M)$, by definition one could then regard this complex as the double complex taking the following form (after considering the corresponding totalization):
\[
\xymatrix@R+5pc@C+5pc{
M \ar[r]^{\psi_1-1} \ar[r] \ar[r] \ar[d]^{\psi_2-1} \ar[d] \ar[d] & M \ar[d] \ar[d] \ar[d]\\
M \ar[r] \ar[r] \ar[r]  & M. 
}
\]

%\[
%\xymatrix@R+0pc@C+0pc{
%0\ar[r] \ar[r] \ar[r]  &C^\bullet_{\psi_{I\backslash \{I\}}}(M) \ar[r]^{\psi_I-1} \ar[r] \ar[r]  & C^\bullet_{\psi_{I\backslash \{I\}}}(M)\ar[r] \ar[r] \ar[r] &0
%.}
%\]

\indent The boundedness of the complex could be observed from the definition and the corresponding ideas in \cite[Theorem 3.3, Theorem 3.5]{4KL3}. For the finiteness we construct suitable completely continuous maps comparing the corresponding quasi-isomorphic complexes for $M$ as below. First we consider some model $M_{r_I}$ over some $\Pi_{\mathrm{an},r_{I,0},I,A}(\pi_{K_I})$. Then as in \cite[Theorem 3.3, Theorem 3.5]{4KL3} we could regard the corresponding complex for $M_{r_I}$ as equivalently the corresponding complex over $\Pi_{\mathrm{an},r_{I,0},I,\Pi_{\mathrm{an},\infty,I,A}(\Gamma_{K_I})}(\pi_{K_I})$. We are then reduced to considering the corresponding complex over $\Pi_{\mathrm{an},r_{I_0},I,\Pi_{\mathrm{an},[t_I,\infty],I,A}(\Gamma_{K_I})}(\pi_{K_I})$ for some $t_I$.\\
%\[
%\xymatrix@R+0pc@C+0pc{
%0\ar[r] \ar[r] \ar[r]  &C^\bullet_{\psi_{I\backslash \{I\}}}(M_{...,r_I}) \ar[r]^{\psi_I-1} \ar[r] \ar[r]  & C^\bullet_{\psi_{I\backslash \{I\}}}(M_{...,r_I})\ar[r] \ar[r] \ar[r] &0
%.}
%\]
\indent Then we consider the complex for some section over some multi interval taking the form of $[s_I,r_I]$ which gives rise to the corresponding quasi-isomorphic complex now over 
\begin{center}
$\Pi_{\mathrm{an},[s_I,r_{I_0}],I,\Pi_{\mathrm{an},[t_I,\infty],I,A}(\Gamma_{K_I})}(\pi_{K_I})$ 
\end{center}
which gives rise to the complex:
\[
\xymatrix@R+5pc@C+2pc{
M_{[s_1,r_1]\times[s_2,r_2]} \ar[r]^{\psi_1-1} \ar[r] \ar[r] \ar[d]^{\psi_2-1} \ar[d] \ar[d] & M_{[ps_1,r_1]\times[s_2,r_2]} \ar[d] \ar[d] \ar[d]\\
M_{[s_1,r_1]\times[ps_2,r_2]}\ar[r] \ar[r] \ar[r]  & M_{[ps_1,r_1]\times[ps_2,r_2]}. 
}
\]
Here the corresponding quasi-isomorphism among different boxes (e.g. multi intervals) is not actually trivial. For the serious argument after the knowledge on the equivalence between the categories of $\varphi_I$-modules over the full Robba ring, $\varphi_I$-Frobenius bundles and $\varphi_I$-modules over the Robba ring with respect to some single multi interval, which could be realized in our context by using \cref{proposition3.5} as in \cite[Lemma 2.10]{4KP1}, we refer to \cite[Lemma 5.2]{4KP1} and \cite[Proposition 5.12, Theorem 5.7.11]{4KL2} for the corresponding comparison on the cohomology groups. \\
\indent To show that the cohomology groups do not change when we switch from a multi-interval $[s_1',r_1']\times [s_2',r_2']$ to a multi-interval $[s_1'',r_1'']\times [s_2'',r_2'']$, we choose to look at an intermediate multi-interval $[s_1'',r_1'']\times [s_2',r_2']$ namely we consider the following in order rings over multi-intervals:
\begin{align}
\Pi_{\mathrm{an},[s'_1,r'_1]\times [s'_2,r'_2],I,\Pi_{\mathrm{an},[t_I,\infty],I,A}(\Gamma_{K_I})}(\pi_{K_I})\\
\Pi_{\mathrm{an},[s''_1,r''_1]\times [s'_2,r'_2],I,\Pi_{\mathrm{an},[t_I,\infty],I,A}(\Gamma_{K_I})}(\pi_{K_I})\\
\Pi_{\mathrm{an},[s''_1,r''_1]\times [s''_2,r''_2],I,\Pi_{\mathrm{an},[t_I,\infty],I,A}(\Gamma_{K_I})}(\pi_{K_I})	
\end{align}
and compare the first two and the last two. To compare with respect to the first two multi-intervals we first consider the $\psi_1$-complex, the corresponding cohomology groups could be regarded as the corresponding Yoneda extension groups over the corresponding twisted polynomial ring in variable of $\psi_1$ with coefficients in the Robba rings, which shows that the $\psi_I$-cohomology groups do not change when we change the boxes. Then for the remaining comparison we argue in the parallel way. \\
\indent To show that the cohomology groups do not change when we switch from Robba ring over some fixed multi-radii $r_{1}\times r_{2}$ to some interval $[s_1,r_1]\times [s_2,r_2]$, we consider the corresponding intermediate Robba ring:
\begin{align}
\bigcap_{s_2>0} \Pi_{\mathrm{an},[s_1,r_1]\times [s_2,r_{2}],I,\Pi_{\mathrm{an},[t_I,\infty],I,A}(\Gamma_{K_I})}(\pi_{K_I})	
\end{align}
and then in order the following rings:
\begin{align}
\bigcap_{s_1>0,s_2>0} \Pi_{\mathrm{an},[s_1,r_{1}]\times [s_2,r_{2}],I,\Pi_{\mathrm{an},[t_I,\infty],I,A}(\Gamma_{K_I})}(\pi_{K_I}),\\
\bigcap_{s_2>0} \Pi_{\mathrm{an},[s_1,r_{1}]\times [s_2,r_{2}],I,\Pi_{\mathrm{an},[t_I,\infty],I,A}(\Gamma_{K_I})}(\pi_{K_I}),\\
 \Pi_{\mathrm{an},[s_1,r_{1}]\times [s_2,r_{2}],I,\Pi_{\mathrm{an},[t_I,\infty],I,A}(\Gamma_{K_I})}(\pi_{K_I}).
\end{align}
Note that finite projective modules over these are equivalent when carrying the corresponding $\varphi_I$-structure. To compare with respect to the first two situations, we first consider the $\psi_1$-cohomology complex, then by using the corresponding Yoneda extension we have the corresponding $\psi_1$-cohomology groups are the same, which implies further that the corresponding $\psi_I$-cohomology groups are the same. Then for the remaining comparison we argue in the parallel way.

%\[
%\xymatrix@R+0pc@C+0pc{
% 0\ar[r] \ar[r] \ar[r]  &C^\bullet_{\psi_{I\backslash \{I\}}}(M_{...,[s_I,r_I]}) \ar[r]^{\psi_I-1} \ar[r] \ar[r]  & C^\bullet_{\psi_{I\backslash \{I\}}}(M_{...,[ps_I,r_I]})\ar[r] \ar[r] \ar[r] &0
%.}
%\]
Then to finish, we choose two well-located multi intervals such that one embeds into another (which exists trivially in our situation which is a little bit different from the situation in \cite{4KL3}):
\[
\xymatrix@R+6pc@C+3pc{
 0\ar[r] \ar[r] \ar[r]  &C^\bullet_{\psi_{I\backslash \{I\}}}(M_{...,[s'_I,r'_I]}) \ar[r]^{\psi_I-1} \ar[r] \ar[r] \ar[d]\ar[d]\ar[d]  & C^\bullet_{\psi_{I\backslash \{I\}}}(M_{...,[ps'_I,r'_I]})\ar[r] \ar[r] \ar[r] \ar[d]\ar[d]\ar[d]  &0\\
% 0\ar[r] \ar[r] \ar[r]  &C^\bullet_{\psi_{I\backslash \{I\}}}(M_{...,[s_I,r_I]}) \ar[r]^{\psi_I-1} \ar[r] \ar[r] \ar[d]\ar[d]\ar[d]  & C^\bullet_{\psi_{I\backslash \{I\}}}(M_{...,[ps_I,r_I]})\ar[r] \ar[r] \ar[r] \ar[d]\ar[d]\ar[d]  &0\\
 0\ar[r] \ar[r] \ar[r]  &C^\bullet_{\psi_{I\backslash \{I\}}}(M_{...,[s''_I,r''_I]}) \ar[r]^{\psi_I-1} \ar[r] \ar[r]  & C^\bullet_{\psi_{I\backslash \{I\}}}(M_{...,[ps''_I,r''_I]})\ar[r] \ar[r] \ar[r] &0
,}
\]

or more explicit:

\[
\xymatrix@R+7pc@C+3pc{
M_{[s'_1,r'_1]\times[s'_2,r'_2]}\ar[r] \ar[r] \ar[r] \ar[d] \ar[d] \ar[d] & M_{[ps'_1,r'_1]\times[s'_2,r'_2]}\oplus M_{[s'_1,r'_1]\times[ps'_2,r'_2]}\ar[r] \ar[r] \ar[r]\ar[d] \ar[d] \ar[d] & M_{[ps'_1,r'_1]\times[ps'_2,r'_2]} \ar[d]\ar[d]\ar[d] \\
%M_{[s_1,r_1]\times[s_2,r_2]}\ar[r] \ar[r] \ar[r] \ar[d] \ar[d] \ar[d] & M_{[ps_1,r_1]\times[s_2,r_2]}\oplus M_{[s_1,r_1]\times[ps_2,r_2]} \ar[r] \ar[r] \ar[r] \ar[d] \ar[d] \ar[d] &M_{[ps_1,r_1]\times[ps_2,r_2]} \ar[d]\ar[d]\ar[d]\\
M_{[s''_1,r''_1]\times[s''_2,r''_2]} \ar[r] \ar[r] \ar[r]  & M_{[ps''_1,r''_1]\times[s''_2,r''_2]}\oplus M_{[s''_1,r''_1]\times[ps''_2,r''_2]}  \ar[r] \ar[r] \ar[r] & M_{[ps''_1,r''_1]\times[ps''_2,r''_2]}.\\
}
\]
Here the corresponding maps between objects over different multi intervals are induced exactly from the corresponding inclusion of multi intervals. Then by the fundamental lemma for in \cite[Satz 5.2]{4Kie1} and  \cite[Lemma 1.10]{4KL3} on the continuous completely maps we get the corresponding finiteness when we restrict to $[t_I,\infty]$, since the map according to the embeddings of the multi intervals are quasi-isomorphisms. Then the corresponding cohomology groups could be arranged to be a coherent sheaf over each section with respect to each such $[t_I,\infty]$, then the coadmissibility follows.
\end{proof}

\begin{proposition}
Based on the previous proposition we have:
\begin{displaymath}
C_{\psi_I}(M)\in  \mathbb{D}^-_\mathrm{perf}(\Pi_{\mathrm{an},\infty,I,A}(\Gamma_{K_I})),	
\end{displaymath}
for $A=\mathbb{Q}_p$.	
\end{proposition}

\begin{proof}
We have the corresponding coadmissibility from the previous proposition, then we work in the category of all the coadmissible modules over the ring $\Pi_{\mathrm{an},\infty,I,A}(\Gamma_{K_I})$. Then by applying the corresponding result \cite[Theorem 3.10]{4Zab1} we can show that the complex is quasi-isomorphic to a corresponding complex being finite projective at each degree in the following. By base change to some $\Pi_{\mathrm{an},[t_I,\infty],I,A}(\Gamma_{K_I})$ for some $t_I$ we have that any such complex over $\Pi_{\mathrm{an},[t_I,\infty],I,A}(\Gamma_{K_I})$ is quasi-isomorphic to a complex being finitely generated at each degree. Therefore the corresponding original complex is quasi-isomorphic to a corresponding complex being coadmissible at each degree. Then in fact at each degree this latter complex is also projective (in the category of all the coadmissible modules over $\Pi_{\mathrm{an},\infty,I,A}(\Gamma_{K_I})$), since this will be further equivalent to that each section over $\Pi_{\mathrm{an},[t_I,\infty],I,A}(\Gamma_{K_I})$ is projective for some $t_I$, which is the case in our situation. Then we are done by \cite[Theorem 3.10]{4Zab1} which implies that coadmissible and projective modules are finitely generated.

%we have that it suffices to check that the corresponding cohomology modules (global sections over $\Pi_{\mathrm{an},\infty,I,A}(\Gamma_{K_I})$) are actually projective in the category of all the coadmissible modules over $\Pi_{\mathrm{an},\infty,I,A}(\Gamma_{K_I})$. 

%But this will be further equivalent to checking that each section over $\Pi_{\mathrm{an},[t_I,\infty],I,A}(\Gamma_{K_I})$ is projective, which is the case in our situation.
\end{proof}

\indent We then give a direct proof on the finiteness of the cohomology of $(\varphi_I,\Gamma_I)$-modules:

\begin{proposition} \label{4proposition4.9}
With the notations above, we have that 
\begin{displaymath}
C^\bullet_{\varphi_I,\Gamma_{I}}(M)\in \mathbb{D}_\mathrm{perf}^-(A), C^\bullet_{\psi_I,\Gamma_{I}}(M)\in \mathbb{D}_\mathrm{perf}^-(A).
\end{displaymath}
\end{proposition}

\begin{proof}
This is proved in the \'etale setting in \cite{4PZ19} and \cite{4CKZ18} in the absolute situation. We use the corresponding Cartan-Serre method to prove this in the following fashion, modeled on the corresponding context of \cite{4KL3}. By simultaneous induction to the base cases for the field $\mathbb{Q}_p$ from all $K_I$, we can now assume that we are considering the situation where all $K_I$ are $\mathbb{Q}_p$. We choose to prove the corresponding result for $(\psi_I,\Gamma_I)$-modules. In our situation we look at the corresponding continuous complexes coming from the corresponding $\Gamma_I$-action, which gives us the corresponding complex $C^\bullet_{\mathrm{con}}(\Gamma_{K_I},M)$. Then the idea is to look at the following double complex of the complex $C^\bullet_{\mathrm{con}}(\Gamma_{K_I},M)$ (to be more precise we need to look at the corresponding totalization):
\[
\xymatrix@R+6pc@C+2pc{
C^\bullet_{\mathrm{con}}(\Gamma_{K_I},M_{[s_1,r_1]\times[s_2,r_2]}) \ar[r]^{\psi_1-1} \ar[r] \ar[r] \ar[d]^{\psi_2-1} \ar[d] \ar[d] & C^\bullet_{\mathrm{con}}(\Gamma_{K_I},M_{[ps_1,r_1]\times[s_2,r_2]}) \ar[d] \ar[d] \ar[d]\\
C^\bullet_{\mathrm{con}}(\Gamma_{K_I},M_{[s_1,r_1]\times[ps_2,r_2]})\ar[r] \ar[r] \ar[r]  &C^\bullet_{\mathrm{con}}(\Gamma_{K_I},M_{[ps_1,r_1]\times[ps_2,r_2]}). 
}
\]
Here the corresponding quasi-isomorphism among different boxes (e.g. multi intervals) is not actually trivial. For the serious argument after the knowledge on the equivalence between the categories of $\varphi_I$-modules over the full Robba ring, $\varphi_I$-Frobenius bundles and $\varphi_I$-modules over the Robba ring with respect to some single multi interval, which could be realized in our context by using \cref{proposition3.5} as in \cite[Lemma 2.10]{4KP1}, we refer to \cite[Lemma 5.2]{4KP1} and \cite[Proposition 5.12, Theorem 5.7.11]{4KL2} for the corresponding comparison on the cohomology groups. \\
\indent To show that the cohomology groups do not change when we switch from a multi-interval $[s_1',r_1']\times [s_2',r_2']$ to a multi-interval $[s_1'',r_1'']\times [s_2'',r_2'']$, we choose to look at an intermediate multi-interval $[s_1'',r_1'']\times [s_2',r_2']$ namely we consider the following in order rings over multi-intervals:
\begin{align}
\Pi_{\mathrm{an},[s'_1,r'_1]\times [s'_2,r'_2],I,A}(\pi_{K_I})\\
\Pi_{\mathrm{an},[s''_1,r''_1]\times [s'_2,r'_2],I,A}(\pi_{K_I})\\
\Pi_{\mathrm{an},[s''_1,r''_1]\times [s''_2,r''_2],I,A}(\pi_{K_I})	
\end{align}
and compare the first two and the last two. To compare with respect to the first two multi-intervals we first consider the $\psi_1$-complex, the corresponding cohomology groups could be regarded as the corresponding Yoneda extension groups over the corresponding twisted polynomial ring in variable of $\psi_1$ with coefficients in the Robba rings, which shows that the $\psi_I$-cohomology groups do not change when we change the boxes. Then for the remaining comparison we argue in the parallel way. \\
\indent To show that the cohomology groups do not change when we switch from Robba ring over some fixed multi-radii $r_{1}\times r_{2}$ to some interval $[s_1,r_1]\times [s_2,r_2]$, we consider the corresponding intermediate Robba ring:
\begin{align}
\bigcap_{s_2>0} \Pi_{\mathrm{an},[s_1,r_1]\times [s_2,r_{2}],I,A}(\pi_{K_I})	
\end{align}
and then in order the following rings:
\begin{align}
\bigcap_{s_1>0,s_2>0} \Pi_{\mathrm{an},[s_1,r_{1}]\times [s_2,r_{2}],I,A}(\pi_{K_I}),\\
\bigcap_{s_2>0} \Pi_{\mathrm{an},[s_1,r_{1}]\times [s_2,r_{2}],I,A}(\pi_{K_I}),\\
 \Pi_{\mathrm{an},[s_1,r_{1}]\times [s_2,r_{2}],I,A}(\pi_{K_I}).
\end{align}
Note that finite projective modules over these are equivalent when carrying the corresponding $\varphi_I$-structure. To compare with respect to the first two situations, we first consider the $\psi_1$-cohomology complex, then by using the corresponding Yoneda extension we have the corresponding $\psi_1$-cohomology groups are the same, which implies further that the corresponding $\psi_I$-cohomology groups are the same. Then for the remaining comparison we argue in the parallel way.\\
\indent Then in order to apply the corresponding fundamental lemma namely \cite[Lemma 1.10]{4KL3}, we need to choose suitable internals in our multi setting where one embeds into the other which gives rise to the following diagram:
\[
\xymatrix@R+6pc@C+0pc{
 0\ar[r] \ar[r] \ar[r]  &C^\bullet_{\psi_{I\backslash \{I\}}}(C^\bullet_{\mathrm{con}}(\Gamma_{K_I},M_{...,[s'_I,r'_I]})) \ar[r]^{\psi_I-1} \ar[r] \ar[r] \ar[d]\ar[d]\ar[d]  & C^\bullet_{\psi_{I\backslash \{I\}}}(C^\bullet_{\mathrm{con}}(\Gamma_{K_I},M_{...,[ps'_I,r'_I]}))\ar[r] \ar[r] \ar[r] \ar[d]\ar[d]\ar[d]  &0\\
 %0\ar[r] \ar[r] \ar[r]  &C^\bullet_{\psi_{I\backslash \{I\}}}(C^\bullet_{\mathrm{con}}(\Gamma_{K_I},M_{...,[s_I,r_I]})) \ar[r]^{\psi_I-1} \ar[r] \ar[r] \ar[d]\ar[d]\ar[d]  & C^\bullet_{\psi_{I\backslash \{I\}}}(C^\bullet_{\mathrm{con}}(\Gamma_{K_I},M_{...,[ps_I,r_I]}))\ar[r] \ar[r] \ar[r] \ar[d]\ar[d]\ar[d]  &0\\
 0\ar[r] \ar[r] \ar[r]  &C^\bullet_{\psi_{I\backslash \{I\}}}(C^\bullet_{\mathrm{con}}(\Gamma_{K_I},M_{...,[s''_I,r''_I]})) \ar[r]^{\psi_I-1} \ar[r] \ar[r]  & C^\bullet_{\psi_{I\backslash \{I\}}}(C^\bullet_{\mathrm{con}}(\Gamma_{K_I},M_{...,[ps''_I,r''_I]}))\ar[r] \ar[r] \ar[r] &0
.}
\]
To be more explicit this gives rise to the following diagram:
\[\tiny
\xymatrix@R+7pc@C+0pc{
C^\bullet_{\Gamma_I}(M_{[s'_1,r'_1]\times[s'_2,r'_2]}) \ar[r] \ar[r] \ar[r] \ar[d] \ar[d] \ar[d] & C^\bullet_{\Gamma_I}(M_{[ps'_1,r'_1]\times[s'_2,r'_2]})\oplus C^\bullet_{\Gamma_I}(M_{[s'_1,r'_1]\times[ps'_2,r'_2]})\ar[r] \ar[r] \ar[r]\ar[d] \ar[d] \ar[d] & C^\bullet_{\Gamma_I}(M_{[ps'_1,r'_1]\times[ps'_2,r'_2]})\ar[d]\ar[d]\ar[d] \\
%C^\bullet_{\Gamma_I}(M_{[s_1,r_1]\times[s_2,r_2]}) \ar[r] \ar[r] \ar[r] \ar[d] \ar[d] \ar[d] & C^\bullet_{\Gamma_I}(M_{[ps_1,r_1]\times[s_2,r_2]})\oplus C^\bullet_{\Gamma_I}(M_{[s_1,r_1]\times[ps_2,r_2]})\ar[r] \ar[r] \ar[r]\ar[d] \ar[d] \ar[d] & C^\bullet_{\Gamma_I}(M_{[ps_1,r_1]\times[ps_2,r_2]})\ar[d]\ar[d]\ar[d]\\
C^\bullet_{\Gamma_I}(M_{[s''_1,r''_1]\times[s''_2,r''_2]}) \ar[r] \ar[r] \ar[r] & C^\bullet_{\Gamma_I}(M_{[ps''_1,r''_1]\times[s''_2,r''_2]})\oplus C^\bullet_{\Gamma_I}(M_{[s''_1,r''_1]\times[ps''_2,r''_2]})\ar[r] \ar[r] \ar[r] & C^\bullet_{\Gamma_I}(M_{[ps''_1,r''_1]\times[ps''_2,r''_2]}).\\
}
\]
Here the corresponding maps between objects over different multi intervals are induced exactly from the corresponding inclusion of multi intervals. Then we consider the corresponding fundamental lemma \cite[Satz 5.2]{4Kie1} and \cite[Lemma 1.10]{4KL3} which implies the corresponding finiteness. And by our basic construction the result follows.
\end{proof}

\indent Now we consider a finite projective $(\varphi_I,\Gamma_I)$-module $M$ over the ring\\ $\Pi_{\mathrm{an},\mathrm{con},I,\Pi_{\mathrm{an},\infty,A}(G)}$, namely we have:
\begin{displaymath}
M=\varprojlim_k M_k,	
\end{displaymath}
where each object $M_k$ satisfies the corresponding framework in the above discussion.

\begin{definition} \mbox{\bf{(After KPX \cite[2.3]{4KPX})}}
For any $(\varphi_I,\Gamma_I)$-module $M$ over\\ $\Pi_{\mathrm{an},\mathrm{con},I,\Pi_{\mathrm{an},\infty,A}(G)}$ we define the complex $C^\bullet_{\Gamma_I}(M)$ of $M$ to be the total complex of the following complex through the induction:
\[
\xymatrix@R+0pc@C+0pc{
[ C^\bullet_{\Gamma_{I\backslash\{I\}}}(M) \ar[r]^{\Gamma_I} \ar[r] \ar[r]  & C^\bullet_{\Gamma_I\backslash\{I\}}(M)]
.}
\]
Then we define the corresponding double complex $C^\bullet_{\varphi_I}C^\bullet_{\Gamma_I}(M)$ again by taking the corresponding totalization of the following complex through induction:
\[
\xymatrix@R+0pc@C+0pc{
[C^\bullet_{\varphi_{I\backslash\{I\}}}C^\bullet_{\Gamma_{I}}(M) \ar[r]^{\varphi_I} \ar[r] \ar[r]  & C^\bullet_{\varphi_{I\backslash\{I\}}}C^\bullet_{\Gamma_I}(M)]
.}
\]
Then we define the complex $C^\bullet_{\varphi_I,\Gamma_I}(M)$ to be the totalization of the double complex defined above, which is called the complex of a $(\varphi_I,\Gamma_I)$-module $M$. Similarly we define the corresponding complex $C^\bullet_{\psi_I,\Gamma_I}(M)$ to be the totalization of the double complex define  d in the same way as above by replacing $\varphi_I$ by then the operators $\psi_I$. For later comparison we also need the corresponding $\psi_I$-cohomology, which is to say the complex $C^\bullet_{\psi_I}(M)$ which could be defined by the totalization of the following complex through induction:
\[
\xymatrix@R+0pc@C+0pc{
[ C^\bullet_{\psi_{I\backslash\{I\}}}(M) \ar[r]^{\psi_I} \ar[r] \ar[r]  & C^\bullet_{\psi_I\backslash\{I\}}(M)]
.}
\]
\end{definition}

\begin{proposition} 
With the notation and setting up as above we have $h^\bullet(C_{\psi_I}(M))$ are coadmissible over $\Pi_{\mathrm{an},\infty,A \widehat{\otimes} \Pi_{\mathrm{an},\infty}(\Gamma_{K_I})}(G)$. 	
\end{proposition}

\begin{proof}
By simultaneous induction to the base cases for the field $\mathbb{Q}_p$ from all $K_I$, we can now assume that we are considering the situation where all $K_I$ are $\mathbb{Q}_p$. Then for the complex $C^\bullet_{\psi_I}(M)$, by definition one could then regard this complex as the double complex taking the following form (after considering the corresponding totalization):
\[
\xymatrix@R+5pc@C+5pc{
M \ar[r]^{\psi_1-1} \ar[r] \ar[r] \ar[d]^{\psi_2-1} \ar[d] \ar[d] & M \ar[d] \ar[d] \ar[d]\\
M \ar[r] \ar[r] \ar[r]  & M. 
}
\]

%\[
%\xymatrix@R+0pc@C+0pc{
%0\ar[r] \ar[r] \ar[r]  &C^\bullet_{\psi_{I\backslash \{I\}}}(M) \ar[r]^{\psi_I-1} \ar[r] \ar[r]  & C^\bullet_{\psi_{I\backslash \{I\}}}(M)\ar[r] \ar[r] \ar[r] &0
%.}
%\]

\indent The boundedness of the complex could be observed from the definition and the corresponding ideas in \cite[Theorem 3.3, Theorem 3.5]{4KL3}. For the finiteness we construct suitable completely continuous maps comparing the corresponding quasi-isomorphic complexes for $M$ as below. Write $\Pi_{\mathrm{an},\infty,A}(G)$ as the following form:
\begin{displaymath}
\Pi_{\mathrm{an},\infty,A }(G):=\varprojlim_k \Pi_{\mathrm{an},\infty,A}(G)_k,	
\end{displaymath}
as a Fr\'echet-Stein algebra. Then take some $\Pi_{\mathrm{an},\infty,A}(G)_k$ in the family with the corresponding $M_k$. By our result above for $M_k$ \cref{proposition4.7} we have that $h^\bullet(C_{\psi_I}(M_k))$ are coadmissible over $\Pi_{\mathrm{an},\infty,\Pi_{\mathrm{an},\infty,A}(G)_k}(\Gamma_{K_I})$. Therefore we have the the corresponding cohomology groups $h^\bullet(C_{\psi_I}(M))$ are actually coming from (the global sections of) some coherent sheaves over $\Pi_{\mathrm{an},\infty,\Pi_{\mathrm{an},\infty,A}(G)}(\Gamma_{K_I})$.

\end{proof}

\begin{proposition}
Based on the previous proposition we have:
\begin{displaymath}
C_{\psi_I}(M)\in  \mathbb{D}^-_\mathrm{perf}(\Pi_{\mathrm{an},\infty,A \widehat{\otimes} \Pi_{\mathrm{an},\infty}(\Gamma_{K_I})}(G)),	
\end{displaymath}
for $A=\mathbb{Q}_p$.	
\end{proposition}

\begin{proof}
We have the corresponding coadmissibility from the previous proposition, then we work in the category of all the coadmissible modules over the ring $\Pi_{\mathrm{an},\infty,I,A}(\Gamma_{K_I})$. Then by applying the corresponding result \cite[Theorem 3.10]{4Zab1} we can show that the complex is quasi-isomorphic to a corresponding complex being finite projective at each degree in the following. Set:
\begin{align}
\Pi_{\mathrm{an},\infty,A \widehat{\otimes} \Pi_{\mathrm{an},\infty}(\Gamma_{K_I})}(G):=\varprojlim_k \Pi_{\mathrm{an},\infty,A \widehat{\otimes} \Pi_{\mathrm{an},\infty}(\Gamma_{K_I})}(G)_k.	
\end{align}

By base change to some $\Pi_{\mathrm{an},\infty,A \widehat{\otimes} \Pi_{\mathrm{an},\infty}(\Gamma_{K_I})}(G)_k$ for some $k$ we have that any such complex over $ \Pi_{\mathrm{an},\infty,A \widehat{\otimes} \Pi_{\mathrm{an},\infty}(\Gamma_{K_I})}(G)_k$ is quasi-isomorphic to a complex being finitely generated at each degree. Therefore the corresponding original complex is quasi-isomorphic to a corresponding complex being coadmissible at each degree. Then in fact at each degree this latter complex is also projective (in the category of all the coadmissible modules over $\varprojlim_k \Pi_{\mathrm{an},\infty,A \widehat{\otimes} \Pi_{\mathrm{an},\infty}(\Gamma_{K_I})}(G)_k$), since this will be further equivalent to that each section over $\varprojlim_k \Pi_{\mathrm{an},\infty,A \widehat{\otimes} \Pi_{\mathrm{an},\infty}(\Gamma_{K_I})}(G)_k$ is projective for some $k$, which is the case in our situation. Then we are done by \cite[Theorem 3.10]{4Zab1} which implies that coadmissible and projective modules are finitely generated.

%we have that it suffices to check that the corresponding cohomology modules (global sections over $\Pi_{\mathrm{an},\infty,I,A}(\Gamma_{K_I})$) are actually projective in the category of all the coadmissible modules over $\Pi_{\mathrm{an},\infty,I,A}(\Gamma_{K_I})$. 

%But this will be further equivalent to checking that each section over $\Pi_{\mathrm{an},[t_I,\infty],I,A}(\Gamma_{K_I})$ is projective, which is the case in our situation.
\end{proof}

\indent We then give a direct proof on the finiteness of the cohomology of $(\varphi_I,\Gamma_I)$-modules:

\begin{proposition}
With the notations above, we have that 
\begin{displaymath}
C^\bullet_{\varphi_I,\Gamma_{I}}(M)\in \mathbb{D}_\mathrm{perf}^-(\Pi_{\mathrm{an},\infty,A }(G)), C^\bullet_{\psi_I,\Gamma_{I}}(M)\in \mathbb{D}_\mathrm{perf}^-(\Pi_{\mathrm{an},\infty,A }(G)).
\end{displaymath}
Here again we assume that the corresponding ring $A$ is $\mathbb{Q}_p$.
\end{proposition}

\begin{proof}
We choose to prove the corresponding result for $(\psi_I,\Gamma_I)$-modules. As above we only have to consider the corresponding situation where $K_I$ are all $\mathbb{Q}_p$ by induction in our current context. Write $\Pi_{\mathrm{an},\infty,A }(G)$ as the limit:
\begin{displaymath}
\Pi_{\mathrm{an},\infty,A }(G)=\varprojlim_k \Pi_{\mathrm{an},\infty,A }(G)_k.	
\end{displaymath}
Working over the ring $\Pi_{\mathrm{an},\infty,A }(G)_k$ we have the corresponding finiteness by the result from \cref{4proposition4.9}, which will show that $C^\bullet_{\psi_I,\Gamma_{I}}(M)$ are all coadmissible modules over $\varprojlim_k \Pi_{\mathrm{an},\infty,A }(G)_k$. By base change to some $\Pi_{\mathrm{an},\infty,A }(G)_k$ for some $k$ we have that any such complex over $ \Pi_{\mathrm{an},\infty,A }(G)_k$ is quasi-isomorphic to a complex being finitely generated at each degree. Therefore the corresponding original complex is quasi-isomorphic to a corresponding complex being coadmissible at each degree. Then in fact at each degree this latter complex is also projective (in the category of all the coadmissible modules over $\varprojlim_k \Pi_{\mathrm{an},\infty,A }(G)_k$), since this will be further equivalent to that each section over $\varprojlim_k \Pi_{\mathrm{an},\infty,A }(G)_k$ is projective for some $k$, which is the case in our situation. Then we are done by \cite[Theorem 3.10]{4Zab1} which implies that coadmissible and projective modules are finitely generated.\\ 
\end{proof}

\subsection{Application to Triangulation over Rigid Analytic Spaces} \label{triangulation}

\indent We are now going to establish some application of the finiteness of the corresponding complexes with respect to the action $(\varphi_I,\Gamma_I)$. We will generalize the interesting globalization of the triangulation established in \cite{4KPX} to our generalized setting. This has its own interests since it will be very natural to ask about the behavior of some triangulation in our context over the corresponding interesting rigid families of the generalized Hodge structures.\\

\indent Since we only proved the corresponding finiteness in the situation where the set $I$ consists of two elements, therefore in this section we are going to make the same assumption. Then we consider the corresponding finite free rank one $(\varphi_I,\Gamma_I)$-modules coming from the corresponding algebraic characters from the product of $\Gamma_{K_I}$, and moreover from the product of the groups $K_I^\times$. Actually in the one dimensional situation this consideration finishes the corresponding construction of all the finite free rank one $(\varphi_I,\Gamma_I)$-modules in a reasonably well-defined way. We follow largely the corresponding notations in the one dimensional situation. We will assume that over some strong (since we do not have a well-established classifications as mentioned above) pointwisely we have the corresponding strong triangulation (meaning that we have the triangulations by using the parameters strictly coming from the corresponding continuous characters of the group $\Gamma_{K_I}^\times$), then we would like to study the corresponding behavior of the corresponding global triangulations. This means that in order to have some coherent triangulations and spreading effect, one only needs to check things pointwisely.\\

\begin{setting}
In this section we are going to work over a general rigid analytic space $X$ defined over suitable finite extension of $L$, where we assume that $L$ is large enough with respect to the fields $K_I$ generalizing the corresponding situation in \cite[Chapter 6]{4KPX}.	
\end{setting}

\indent We start from the following geometric setting:

\begin{definition} \mbox{\bf{(After KPX \cite[Definition 6.2.1]{4KPX})}}
First we are going to consider the corresponding sheaves $\Pi_{\mathrm{an},\mathrm{con},I,X}(\pi_{K_I})$ with respect to the analytic space $X$ over $L$, which is defined to be glueing the rings of analytic functions over each affinoid $\mathrm{Max}A$ for the space $X$, namely glueing the rings having the form of $\Pi_{\mathrm{an},\mathrm{con},I,A}(\pi_{K_I})$ with respect to this affinoid. Also as in \cite[Definition 6.2.1]{4KPX} one can define the corresponding sheaf with respect to some specific radii which is to say the sheaves taking the form of $\Pi_{\mathrm{an},r_{I},I,X}(\pi_{K_I})$. Then we define the corresponding $(\varphi_I,\Gamma_I)$-modules over these two kinds of sheaves by considering the corresponding families of finite locally free sheaves of modules over $\Pi_{\mathrm{an},\mathrm{con},I,X}(\pi_{K_I})$ or $\Pi_{\mathrm{an},r_{I},I,X}(\pi_{K_I})$ in the well-defined and compatible sense as in the one-dimensional situation in \cite[Definition 6.2.1]{4KPX}.

\end{definition}

\begin{setting}
We choose to now consider the following situation where all the fields $K_I$ are assumed to be $\mathbb{Q}_p$. Therefore now the $p$-adic Lie groups we are considering take the following form namely the product of $\mathbb{Q}_p^\times$.	
\end{setting}

\indent Under this assumption namely where all the fields $K_I$ are $\mathbb{Q}_p$ one could directly consider the corresponding construction as generalized below from the one dimensional situation considered in \cite[Notation 6.2.2]{4KPX}:

\begin{definition} \mbox{\bf{(After KPX \cite[Notation 6.2.2]{4KPX})}}
Consider the product of the $p$-adic Lie group $\mathbb{Q}_p^\times$ which is to say $\prod_{\alpha\in I}\mathbb{Q}_p^\times$. We are going to define the corresponding free rank one $(\varphi_I,\Gamma_I)$-modules of the character types attached to any character taking the form of $\delta_I:\prod_{\alpha\in I}\mathbb{Q}_p^\times\rightarrow \Gamma(X,\mathcal{O}_X)^\times$, to be the free of rank one module $\Pi_{\mathrm{an},\mathrm{con},I,X}(\pi_{K_I})\mathbf{e}$ such that we have for each $\alpha\in I$, the Frobenius $\varphi_\alpha$ acts via $\delta_I(p_\alpha)$ while the group $\Gamma_\alpha$ acts via $\delta_I(\chi_\alpha(\gamma_\alpha))$. Note that when we have that the set $I$ is a singleton the corresponding definition recovers the one defined in \cite[Notation 6.2.2]{4KPX}. We will use the notation $\Pi_{\mathrm{an},\mathrm{con},I,X}(\pi_{K_I})(\delta_I)$ to denote such free rank one object in the general setting.	
\end{definition}

\indent On the other hand actually we have no direct classification literally as in the one dimensional situation on the free rank one objects at least up to this definition. Therefore we have to at this moment to fix our objects considered, which is to say that we will focus on all the free rank one objects of character types as defined above, then to study the corresponding geometrization and the variation of the corresponding triangulation. We start from the following basic definition which is a direct generalization of the corresponding one in the one dimensional situation.

\begin{definition} \mbox{\bf{(After KPX \cite[Definition 6.3.1]{4KPX})}}
We define the \textit{triangulated} $(\varphi_I,\Gamma_I)$-modules over the higher dimensional Robba rings over a rigid analytic space $X$ taking the form of $\Pi_{\mathrm{an},\mathrm{con},I,X}(\pi_{K_I})$. Such a module $M$ is defined to be a usual $(\varphi_I,\Gamma_I)$-module over $\Pi_{\mathrm{an},\mathrm{con},I,X}(\pi_{K_I})$ carrying a triangulation with $n$-parameters $\delta_1,...,\delta_n$ which are $n$ characters (we will assume that these are continuous ones) of the group $\prod_{I}\mathbb{Q}_p^\times$ with the values in the corresponding group $\Gamma(X,\mathcal{O}_X)^\times$, which means that we have the following filtration on $M$ by the corresponding sub $(\varphi_I,\Gamma_I)$-modules over the sheaf $\Pi_{\mathrm{an},\mathrm{con},I,X}(\pi_{K_I})$:
\begin{displaymath}
0=M_0\subset M_1\subset...\subset M_n=M	
\end{displaymath}
where $n$ is assumed in our situation to be the rank of $M$, moreover we have that the corresponding each graded piece $M_{i}/M_{i-1}$ for $i=1,...,n$ is a rank one $(\varphi_I,\Gamma_I)$-module with the corresponding parameter $\delta_i$ defined above with some twisted line bundle $\mathbb{L}_i,i=1,...,n$:
\begin{displaymath}
M_{i}/M_{i-1}\overset{\sim}{\longrightarrow}\Pi_{\mathrm{an},\mathrm{con},I,X}(\pi_{K_I})(\delta_i)\otimes \mathbb{L}_i,i=1,...,n.	
\end{displaymath}
\end{definition}

\begin{definition} \mbox{\bf{(After KPX \cite[Definition 6.3.1]{4KPX})}}
As in \cite[Definition 6.3.1]{4KPX} we generalize the corresponding $strictly~triangulated~modules$ to the higher dimensional situation by defining them to be over $X=\mathrm{Max}(L)$ the triangulated modules as above in a manner of building from each $i$-th filtration to the $i+1$-th uniquely through the corresponding parameters.
\end{definition}

\indent We are going to consider the corresponding integrated behavior of a dense subset of pointwise triangulations. First recall that a subspace $U$ is $Zariski~dense$ if for each member in an admissible covering of the space $X$ the space $U$ is dense inside (see \cite[Chapter 6.3]{4KPX}). Then we can now consider the corresponding densely pointwise strictly triangulated families:

\begin{definition} \mbox{\bf{(After KPX \cite[Definition 6.3.2]{4KPX})}}
We are going to call a $(\varphi_I,\Gamma_I)$-module over the sheaf $\Pi_{\mathrm{an},\mathrm{con},I,X}(\pi_{K_I})$ a \textit{densely pointwise strictly triangulated family} if there exists a Zariski-dense subspace $X_{dpt}$ inside $X$ such that the corresponding the restriction to $X_{dpt}$ of the sheaf $M$ is pointwise strictly triangulated in the sense defined above.	
\end{definition}

\begin{lemma}\mbox{\bf{(After KPX \cite[6.3.3]{4KPX})}}
In our situation, we have the cohomology groups $H^\bullet_{\varphi_I,\Gamma_I}(M)$ are coherent sheaves over the rigid analytic space $X$ for $M$ any $(\varphi_I,\Gamma_I)$-modules over the Frobenius sheaves $\Pi_{\mathrm{an},\mathrm{con},I,X}(\pi_{K_I})$ with the corresponding non-vanishing throughout the degrees in $[0,2I]$.	
\end{lemma}

\begin{proof}
The proof of this result is parallel to the corresponding result in the one dimensional situation, see \cite[6.3.3]{4KPX}. Essentially one has this result as a direct consequence of the results in the local setting in \cref{section4.1}.
%One derives the corresponding coherence of $H^\bullet_{\varphi_I,\Gamma_I}(M/t_\alpha)$ from that of the former cohomology groups.
\end{proof}

%\begin{lemma} \mbox{\bf{(After KPX \cite[Corollary 6.2.9, Lemma 6.2.10]{4KPX})}}
%For any character $\chi_I:\prod \mathbb{Q}_p^\times\rightarrow E^\times$, the corresponding cohomology group $H^0(\Pi_{\mathrm{an},\mathrm{con},I,X}(\pi_{K_I})(\chi_I))$ vanishes unless we have the situation where the character could be written as $x_1^{-k_1}...x_I^{-k_I}$ where the integer $k_\alpha >0 $ for each $\alpha\in I$, and in the latter situation we have that the corresponding cohomology group is of one dimension generated by $t_1^{k_1}...t_I^{k_I}$. Consequently the corresponding submodules of $\Pi_{\mathrm{an},\mathrm{con},I,X}(\pi_{K_I})(\chi_I)$ will take the form of $t_1^{\ell_1}...t_I^{\ell_I}\Pi_{\mathrm{an},\mathrm{con},I,X}(\pi_{K_I})(\chi_I)$ for $\ell_\alpha>0$ for each $\alpha\in I$. Finally take any element in the space $H^0(M)$ for any $(\varphi_I,\Gamma_{I})$-module $M$ and regard this element as a map from the Robba ring to the underlying module of $M$, then we have that this map realizes a saturation only when we have that the corresponding image of this element does not vanish in the corresponding zeroth cohomology group of $M/t_\alpha$ for each $\alpha\in I$.	
%\end{lemma}
%
%
%\begin{proof}
%The calculation on the cohomology group is straightforward, see \cite[Lemma I.5]{4Col1}. The statement on the submodules follows from this. For the last statement one can follow the strategy of \cite[Lemma 6.2.10]{4KPX} to prove this, the only difference is just one has to use different $t_\alpha$ for each $\alpha\in I$. 
%\end{proof}

\indent We then consider our main results in this section, which is a direct generalization of the results of \cite{4KPX} (chapter 6) to the situation where we have higher dimensional Hodge-Frobenius structure. 

%\begin{theorem}\mbox{\bf{(After KPX \cite[Theorem 6.3.9, Corollary 6.3.10]{4KPX})}}
%Now suppose that we have a $(\varphi_I,\Gamma_I)$-sheaf over $\Pi_{\mathrm{an},\mathrm{con},I,X}(\pi_{K_I})$ which we denote by $M$, where the space $X$ is a reduced rigid analytic space with irreducible connected components. 
%and let $\delta$ be a character from the group $\prod_{I}\mathbb{Q}_p^\times$ to the group $\Gamma(X,\mathcal{O}_X)^\times$. 
%Assume now that $M$ is a pointwise densely strictly triangulated $(\varphi_I,\Gamma_I)$-sheaf over $\Pi_{\mathrm{an},\mathrm{con},I,X}(\pi_{K_I})$. Then there exists a blow-up (birational, proper, reduced) $b:Y\rightarrow X$ such that the pullback of $M$ along this admits global triangulation (namely, a global triangulated filtration by coherent subsheaves as in \cite[Corollary 6.3.10]{4KPX}).	
%\end{theorem}

\begin{proposition} \mbox{\bf{(After KPX \cite[Theorem 6.3.9]{4KPX})}}
Suppose that we have a $(\varphi_I,\Gamma_I)$-module $M$ over the ring $\Pi_{\mathrm{an},\mathrm{con},I,X}(\pi_{K_I})$, where we assume that this module is of rank $k$. Suppose that we have a continuous character from the group $\prod_I \mathbb{Q}_p^\times$ into the group $\Gamma(X,\mathcal{O}_X)^\times$. Suppose that there exists a subset $X_{pdt}\subset X$ over which we have that the cohomology group $H^0_{\varphi_I,\Gamma_I}(M_x^\vee(\delta_x))$ is of dimension one (which is to say for each $x\in X_{pdt}$), with the further assumption that the image of the corresponding ring $\Pi_{\mathrm{an},\mathrm{con},I,X}(\pi_{K_I})$ is saturated under this assumption (namely under the corresponding consideration of the basis). Assume this subset is Zariski-dense. Assume now the space $X$ is reduced, then one could be able to find a morphism from another blow-up $b:Y\rightarrow X$ which is birational and proper and a corresponding map $f:b^*M\rightarrow \Pi_{\mathrm{an},\mathrm{con},I,X}(\pi_{K_I})(\delta)\otimes_{\mathcal{O}_X}\mathbb{L}$ (with some well-defined line bundle $\mathbb{L}$ as in the situation of \cite[Theorem 6.3.9]{4KPX}). 

%with the following further fact that the corresponding kernel $\mathbf{Ker}f$ is of rank just $k-1$.:\\

%\indent A.  We have the corresponding subset $Y_1^c\subset Y$ of all the points $y$ where the corresponding map $f_y:M_y\rightarrow \Pi_{\mathrm{an},\mathrm{con},I,\kappa_y}(\pi_{K_I})(\delta_y)\otimes_{\mathcal{O}_X}\mathbb{L}_y$ is surjective (as the original context of \cite{4KPX} we need to assume that the corresponding element generates the zero-th cohomology of $M_y^\vee(\delta_y)$ as well), is open and dense in the Zariski topology; B.  The corresponding kernel $\mathbf{Ker}f$ is of rank just $k-1$. And the corresponding cokernel $\mathbf{CoKer}f$ is killed by some power taking the form of $t_1^{k_1}...t_I^{k_I}$ for some integers $k_\alpha\geq 1$ ($\forall \alpha\in I$) with support on $Y_1$. 	
\end{proposition}

\begin{proof}

Indeed, we can adapt the corresponding strategy of \cite[Theorem 6.3.9]{4KPX} to prove this. First by taking the well-established (see \cite[Theorem 6.3.9]{4KPX} the geometry actually) normalization we are reduced to assume that the space $X$ is normal and connected, and consequently the ranks of the corresponding coherent sheaves preserve. Then we consider the corresponding local coherence of the corresponding cohomology complex $C^\bullet_{\varphi_I,\Gamma_I}(M)$ (for all $\alpha\in I$), and the corresponding local base change results (namely the direct analog of the theorem 4.4.3 of \cite{4KPX}), together with the 6.3.6 of \cite{4KPX} we have the corresponding groups $H^\bullet_{\varphi_I,\Gamma_I}(M)$ are flat for $\bullet=0$ and having Tor dimension less than or equal to 1 above, after we replace $M$ by some $M_0$ which is the pullback of $M^\vee(\delta)$ along some $f_0:Y_0\rightarrow X$ which is birational and proper. We then glue along the corresponding local construction above as in \cite[Theorem 6.3.9]{4KPX}. The first task for us then will be the corresponding construction of the morphism over this $Y_0$ which will be later be promoted to be the desired $Y$ by taking the corresponding blow-up. The construction goes as in the following way which is modeled on \cite[Theorem 6.3.9]{4KPX}. One first considers the corresponding subset $Y_1'^c\subset Y_0$ that consists of all the points at which the corresponding forth term in the following exact sequence (by considering the corresponding base change spectral sequence)
\[
\xymatrix@R+0pc@C+0pc{
0\ar[r] \ar[r] \ar[r]  &H^0_{\varphi_I,\Gamma_I}(M_0)\otimes\kappa_y \ar[r] \ar[r] \ar[r]  &H^0_{\varphi_I,\Gamma_I}(M_{0,y}) \ar[r] \ar[r] \ar[r] &\mathrm{Tor}_1(H^1_{\varphi_I,\Gamma_I}(M_0),\kappa_y) \ar[r] \ar[r] \ar[r] &0
.}
\]
vanishes in our situation. Note that this is dense and open in the Zariski topology. Over the subset $Y_1'$ the resulting exact sequence degenerates just to be:
\[
\xymatrix@R+0pc@C+0pc{
0\ar[r] \ar[r] \ar[r]  &H^0_{\varphi_I,\Gamma_I}(M_0)\otimes\kappa_y \ar[r] \ar[r] \ar[r]  &H^0_{\varphi_I,\Gamma_I}(M_{0,y}) \ar[r] \ar[r] \ar[r] &\mathrm{Tor}_1(H^1_{\varphi_I,\Gamma_I}(M_0),\kappa_y) \ar[r] \ar[r] \ar[r] &0
,}
\]
where the forth term does not vanish. Then we note that over the subset $Y_1'^c$ we have that the corresponding resulting cohomology group $H^0_{\varphi_I,\Gamma_I}(M_0)$ is then of dimensional one, which implies that by the topological property of the subspace it is coherent of rank one throughout the whole space $Y_0$. We then have the chance to extract the desired morphism in our situation as below: the corresponding line bundle $\mathbb{L}$ is chosen to be the corresponding dual line bundle. Then we have the result morphism $\Pi_{\mathrm{an},\mathrm{con},I,X}(\pi_{K_I})\otimes \mathbb{L}^\vee\rightarrow M_0$ and consequently we have the corresponding desired morphism which is to say $f_0^*M\rightarrow \Pi_{\mathrm{an},\mathrm{con},I,X}(\pi_{K_I})(\delta)\otimes \mathbb{L}$. \\
\indent Then up to this point, one can then consider 6.3.6 of \cite{4KPX} to extract the desired blow-up $Y\rightarrow Y_0$ and define the desired data through the pullback along blow-up, which is to say the corresponding $\mathbb{L}$, $f$ (see \cite[Theorem 6.3.9]{4KPX}). 

%Then to finish we consider the corresponding properties of the kernel. We also in our situation consider the corresponding base change spectral sequence which gives rise to the following exact sequence:
%\[
%\xymatrix@R+0pc@C+0pc{
%0\ar[r] \ar[r] \ar[r]  &\mathrm{Tor}^Y_{2}(\mathbf{Ker},\kappa_y) \ar[r] \ar[r] \ar[r]  &\mathbf{Ker}_y \ar[r] \ar[r] \ar[r] &\mathbf{Ker} \ar[r] \ar[r] \ar[r] &\mathrm{Tor}^Y_{1}(\mathbf{Ker},\kappa_y) \ar[r] \ar[r] \ar[r]  &0 
%,}
%\]
%which further implies that the corresponding property that $\mathbf{Ker}_y$ is being of rank $k-1$, since the second $\mathrm{Tor}$ is zero in our situation as well. However this is not the desired eventual result for the kernel in our situation. To get the desired property of our kernel we need to integrate the results at each point $y$ to get the desired result for the kernel in our situation, but essentially this is not different from the corresponding proof in \cite[Theorem 6.3.9]{4KPX} by considering the 2.1.8 of \cite{4KPX} and the previous lemmas in this section, which is to say we just need to work locally and then consider the corresponding models of our modules with respect to some radii $r_I$, then integrate the results at all the points $y$ as above. Then we are done.
\end{proof}

\begin{remark}
The corresponding morphism established by the previous proposition is actually expected to satisfy the corresponding properties parallel to the corresponding situation in \cite[Theorem 6.3.9, Corollary 6.3.10]{4KPX}. Namely in the context mentioned in the proof of the previous proposition we are now at the position to (thanks to the previous proposition) conjecture: \\
\indent A. We have the corresponding subset $Y_1^c\subset Y$ of all the points $y$ where the corresponding map $f_y:M_y\rightarrow \Pi_{\mathrm{an},\mathrm{con},I,\kappa_y}(\pi_{K_I})(\delta_y)\otimes_{\mathcal{O}_X}\mathbb{L}_y$ is surjective (as the original context of \cite{4KPX} we need to assume that the corresponding element generates the zero-th cohomology of $M_y^\vee(\delta_y)$ as well), is open and dense in the Zariski topology; B. The corresponding kernel $\mathbf{Ker}f$ is of rank just $k-1$.

\end{remark}

%\newpage

\newpage\section{Noncommutative Coefficients}

\subsection{Frobenius Modules over Noncommutative Rings} \label{4section5.1}

\indent We now consider the coefficient where $A$ is a noetherian noncommutative Banach affinoid algebra over $\mathbb{Q}_p$ namely quotient of noncommutative Tate algebra with free variables. The finiteness even in the univariate situation with noncommutative coefficients is actually not known in \cite{4Zah1} which is key ingredient in the noncommutative local $p$-adic Tamagawa number conjecture in \cite{4Zah1}. We make the corresponding discussion in our situation.

%We give a proof by using noncommutative $p$-adic functional analysis in the following development.

\begin{example} \label{example5.1}
A nontrivial example of such noetherian noncommutative Banach affinoid algebra could be the corresponding noncommutative analog of the commutative Tate algebra $\mathbb{T}_2$ considered in \cite[Section 3]{4So1}. To be more precise one considers first the corresponding twisted polynomial ring $\mathbb{Q}_p[X_1,X_2]'$ such that $X_1X_2=aX_2X_1$ for some nonzero element $a\in \mathbb{Q}_p$, then takes the corresponding completion with respect to the Gauss norm $\|.\|$ which is defined by:
\begin{displaymath}
\|\sum_{i_1,i_2}a_{i_1,i_2}X_1^{i_1}X_2^{i_2}\|:=\sup_{i_1,i_2} \{|a_{i_1,i_2}|_p\}.	
\end{displaymath}
Note that as in \cite[Proposition 3.1]{4So1} the corresponding categories of finite left $A$-modules and the correspnding Banach ones are actually equivalent, under the forgetful functor. 	
\end{example}

\indent We then consider the families of the corresponding Frobenius modules with coefficients in noncommutative affinoid algebra $A$ considered as above.

\begin{setting}
For the Robba rings with respect to some closed interval $[s_I,r_I]\subset (0,\infty]^{|I|}$ we defined above namely $\Pi_{\mathrm{an},[s_I,r_I],I,\mathbb{Q}_p}$, we take the completed tensor product with the corresponding ring $A$ in this section to define the corresponding noncommutative version of the Robba rings. We keep all the notations compatible with the ones in the commutative setting in the previous discussion in the previous sections.	
\end{setting}

\begin{definition} \mbox{\bf{(After KPX \cite[Definition 2.2.2]{4KPX})}}
We keep all the notations (except the coefficient $A$ which is now noncommutative) compatible with the ones in the commutative setting in the previous discussion in the previous sections. For each $\alpha\in I$, we choose suitable uniformizer $\pi_{K_\alpha}$ in our consideration. Then we are going to use the notation $\Pi_{\mathrm{an},\mathrm{con},I,A}(\pi_{K_I})$ to denote the corresponding period ring constructed from $\Pi_{\mathrm{an},\mathrm{con},I,A}$ just by directly replacing the variables by the corresponding uniformizers as above. And similarly for other period rings we use the corresponding notations taking the same form namely $\Pi_{\mathrm{an},\mathrm{con},I,A}(\pi_{K_I})$, $\Pi_{\mathrm{an},r_{I,0},I,A}(\pi_{K_I})$, $\Pi_{\mathrm{an},[s_I,r_I],I,A}(\pi_{K_I})$. As in the one dimensional situation we consider the situation where the radii are all sufficiently small. Then one could define the multiple Frobenius actions from the multi Frobenius $\varphi_I$ over each the ring mentioned above. For suitable radii $r_I$ we define a $\varphi_I$-module to be a finite projective left $\Pi_{\mathrm{an},r_{I},I,A}(\pi_{K_I})$-module with the requirement that for each $\alpha\in I$ we have $\varphi_\alpha^*M_{r_I}\overset{\sim}{\rightarrow}M_{...,r_\alpha/p,...}$ (after suitable base changes). Then we define $M:=\Pi_{\mathrm{an},\mathrm{con},I,A}(\pi_{K_I})\otimes_{\Pi_{\mathrm{an},r_{I},I,A}(\pi_{K_I})}M_{r_I}$ to define a $\varphi_I$-module over the full relative Robba ring in our situation. Furthermore we have the notion of $\varphi_I$-bundles in our situation which consists of a family of $\varphi_I$-modules $\{M_{[s_I,r_I]}\}$ where each module $M_{[s_I,r_I]}$ is defined to be finite projective over $\Pi_{\mathrm{an},[s_I,r_I],I,A}(\pi_{K_I})$ satisfying the action formula taking the form of $\varphi_\alpha^*M_{[s_I,r_I]}\overset{\sim}{\rightarrow}M_{...,[s_\alpha/p,r_\alpha/p],...}$ (after suitable base changes). Furthermore for each $\alpha$ we have the corresponding operator $\varphi_\alpha:M_{[s_I,r_I]}\rightarrow M_{...,[s_\alpha/p,r_\alpha/p],...}$ and we have the corresponding operator $\psi_\alpha$ which is defined to be $p^{-1}\varphi_\alpha^{-1}\circ\mathrm{Trace}_{M_{...,[s_\alpha/p,r_\alpha/p],...}/\varphi_\alpha(M_{[s_I,r_I]})}$. Certainly we have the corresponding operator $\psi_\alpha$ over the global section $M_{r_I}$. Note that here we require that the Hodge-Frobenius structures are commutative in the sense that all the Frobenius are commuting with each other, and they are semilinear.
\end{definition}

\begin{definition} \mbox{\bf{(After KPX \cite[Definition 2.2.2]{4KPX})}} \label{definition6.4}
We keep all the notations (except the coefficient $A$ which is now noncommutative) compatible with the ones in the commutative setting in the previous discussion in the previous sections. For each $\alpha\in I$, we choose suitable uniformizer $\pi_{K_\alpha}$ in our consideration. Then we are going to use the notation $\Pi_{\mathrm{an},\mathrm{con},I,A}(\pi_{K_I})$ to denote the corresponding period ring constructed from $\Pi_{\mathrm{an},\mathrm{con},I,A}$ just by directly replacing the variables by the corresponding uniformizers as above. And similarly for other period rings we use the corresponding notations taking the same form namely $\Pi_{\mathrm{an},\mathrm{con},I,A}(\pi_{K_I})$, $\Pi_{\mathrm{an},r_{I,0},I,A}(\pi_{K_I})$, $\Pi_{\mathrm{an},[s_I,r_I],I,A}(\pi_{K_I})$. As in the one dimensional situation we consider the situation where the radii are all sufficiently small. Then one could define the multiple Frobenius actions from the multi Frobenius $\varphi_I$ over each the ring mentioned above. For suitable radii $r_I$ we define a coherent $\varphi_I$-module to be a finitely presented left $\Pi_{\mathrm{an},r_{I},I,A}(\pi_{K_I})$-module with the requirement that for each $\alpha\in I$ we have $\varphi_\alpha^*M_{r_I}\overset{\sim}{\rightarrow}M_{...,r_\alpha/p,...}$ (after suitable base changes). Then we define $M:=\Pi_{\mathrm{an},\mathrm{con},I,A}(\pi_{K_I})\otimes_{\Pi_{\mathrm{an},r_{I},I,A}(\pi_{K_I})}M_{r_I}$ to define a coherent $\varphi_I$-module over the full relative Robba ring in our situation. Furthermore we have the notion of coherent $\varphi_I$-bundles in our situation which consists of a family of coherent $\varphi_I$-modules $\{M_{[s_I,r_I]}\}$ where each module $M_{[s_I,r_I]}$ is defined to be finitely presented over $\Pi_{\mathrm{an},[s_I,r_I],I,A}(\pi_{K_I})$ satisfying the action formula taking the form of $\varphi_\alpha^*M_{[s_I,r_I]}\overset{\sim}{\rightarrow}M_{...,[s_\alpha/p,r_\alpha/p],...}$ (after suitable base changes). Furthermore for each $\alpha$ we have the corresponding operator $\varphi_\alpha:M_{[s_I,r_I]}\rightarrow M_{...,[s_\alpha/p,r_\alpha/p],...}$ and we have the corresponding operator $\psi_\alpha$ which is defined to be $p^{-1}\varphi_\alpha^{-1}\circ\mathrm{Trace}_{M_{...,[s_\alpha/p,r_\alpha/p],...}/\varphi_\alpha(M_{[s_I,r_I]})}$. Certainly we have the corresponding operator $\psi_\alpha$ over the global section $M_{r_I}$. Note that here we require that the Hodge-Frobenius structures are commutative in the sense that all the Frobenius are commuting with each other, and they are semilinear.
\end{definition}

\begin{definition} \mbox{\bf{(After KPX \cite[Definition 2.2.12]{4KPX})}}
Then we add the corresponding more Lie group action in our setting, namely the multi semiliear $\Gamma_{K_I}$-action. Again in this situation we require that all the actions of the $\Gamma_{K_I}$ are commuting with each other and with all the semilinear Frobenius actions defined above. We define the corresponding $(\varphi_I,\Gamma_I)$-modules over $\Pi_{\mathrm{an},r_{I,0},I,A}(\pi_{K_I})$ to be finite projective module with mutually commuting semilinear actions of $\varphi_I$ and $\Gamma_{K_I}$. For the latter action we require that to be continuous. We can also define the corresponding actions for coherent Frobenius modules.
\end{definition}

%\begin{assumption} \label{assumption1}
%We choose to consider the story in which $I$ consists of two element. We note that actually everything could be carried over to more general setting. We will then use the notation $\{1,2\}$ to denote our set which indicates the corresponding Frobenius actions and Lie group actions with respect to each index.	
%\end{assumption}

\begin{definition} \mbox{\bf{(After KPX \cite[2.3]{4KPX})}}
For any $(\varphi_I,\Gamma_I)$-module $M$ over\\
 $\Pi_{\mathrm{an},\mathrm{con},I,A}(\pi_{K_I})$ we define the complex $C^\bullet_{\Gamma_I}(M)$ of $M$ to be the total complex of the following complex through the induction:
\[
\xymatrix@R+0pc@C+0pc{
[ C^\bullet_{\Gamma_{I\backslash\{I\}}}(M) \ar[r]^{\Gamma_I} \ar[r] \ar[r]  & C^\bullet_{\Gamma_I\backslash\{I\}}(M)]
.}
\]
Then we define the corresponding double complex $C^\bullet_{\varphi_I}C^\bullet_{\Gamma_I}(M)$ again by taking the corresponding totalization of the following complex through induction:
\[
\xymatrix@R+0pc@C+0pc{
[C^\bullet_{\varphi_{I\backslash\{I\}}}C^\bullet_{\Gamma_{I}}(M) \ar[r]^{\varphi_I} \ar[r] \ar[r]  & C^\bullet_{\varphi_{I\backslash\{I\}}}C^\bullet_{\Gamma_I}(M)]
.}
\]
Then we define the complex $C^\bullet_{\varphi_I,\Gamma_I}(M)$ to be the totalization of the double complex defined above, which is called the complex of a $(\varphi_I,\Gamma_I)$-module $M$. Similarly we define the corresponding complex $C^\bullet_{\psi_I,\Gamma_I}(M)$ to be the totalization of the double complex define  d in the same way as above by replacing $\varphi_I$ by then the operators $\psi_I$. For later comparison we also need the corresponding $\psi_I$-cohomology, which is to say the complex $C^\bullet_{\psi_I}(M)$ which could be defined by the totalization of the following complex through induction:
\[
\xymatrix@R+0pc@C+0pc{
[ C^\bullet_{\psi_{I\backslash\{I\}}}(M) \ar[r]^{\psi_I} \ar[r] \ar[r]  & C^\bullet_{\psi_I\backslash\{I\}}(M)]
.}
\]
\end{definition}

\indent In order to proceed we now make the following assumption:

\begin{assumption} \label{assumption6.12}
We assume that the ring $A$ in our situation satisfies the following assumption. For any closed interval $[s_I,r_I]$ we assume that the Robba ring $\Pi_{\mathrm{an},[s_I,r_I],I,A}$ is left and right noetherian, namely noetherian in the noncommutative situation. In the following we will discuss some samples.
\end{assumption}

\begin{example} 
First note we point out that the corresponding Tate algebra with several free variables is actually not noetherian, even the corresponding polynomial rings with several free variables are not noetherian at all. For instance in the situation where we have two variables, as in \cite[Introduction]{4K1} and \cite[Exercise 1E]{4GW1} we know that the ring $k[x,y]^\mathrm{nc}$ noncommutative polynomial ring contains specific corresponding ideal generated by $xy^k$ for each nonnegative integer power $k\geq 0$. Even in the commutative setting, suppose we have a commutative noetherian Banach algebra $B$, it is also definitely not trivial that the affinoid algebra:
\begin{displaymath}
B\{T_1/a_1,...,T_m/a_m,b_1/T_1,...b_m/T_m\}	
\end{displaymath}
is noetherian or not. Actually there does exist some counterexample which makes this not true in general as one might not expect, as in \cite[8.3]{4FGK}. 
\end{example}

\begin{example}
One can actually produce another examplification of the\\ \cref{assumption6.12}, namely we consider the following situation. We look at the corresponding (generalized version of) twisted polynomial rings in the sense of \cref{example5.1} (see \cite[Section 3]{4So1}). For instance we take a such ring with three variables:
\begin{align}
\mathbb{Q}_p[X_1,X_2,X_3]'	
\end{align}
with:
\begin{align}
X_iX_j=a_{i,j}X_jX_i,i,j\in\{1,2,3\}	
\end{align}
such that $a_{3,1}=a_{1,3}=1$, $a_{3,2}=a_{2,3}=1$ and $a_{2,1}=a_{1,2}\neq 1$. Then after taking the corresponding completion under the Gauss norm we will have the corresponding noncommutative affinoid which is noetherian and satisfies the corresponding \cref{assumption6.12} by suitable rational localization.
\end{example}

Now consider the corresponding rings in the style considered in \cite[Proposition 4.1]{4Zab1}:
\begin{displaymath}
D_{[s,r]}(\mathbb{Z}_p,\mathbb{Q}_p), D_{[\rho_1,\rho_2]}(G,K)	
\end{displaymath}
with $G$ a pro-$p$ group as in \cite[Chapter 4]{4Zab1} (which is to say we will also assume $G$ to be uniform). The corresponding ring $D_{[\rho_1,\rho_2]}(G,K)$ could be defined through suitable microlocalization with respect to some noncommutative polyannulus with respect to the radii $\rho_1,\rho_2$ which is direct noncommutative generalization of the usual Robba ring we considered in this paper including obviously the corresponding multivariate ones. With the notation in \cite[Proposition 4.1]{4Zab1} namely for sufficiently large $s,r,\rho_1,\rho_2\in p^\mathbb{Q}$ and $s,r,\rho_1,\rho_2<1$ we have:

\begin{theorem} \mbox{\bf{(Z\'abr\'adi \cite[Proposition 4.1]{4Zab1})}}\\
The noncommutative Robba rings $D_{[\rho_1,\rho_2]}(G,K)$ is noetherian.	
\end{theorem}

We now take the corresponding completed tensor product we have then:
\begin{displaymath}
D_{[s,r]}(\mathbb{Z}_p,\mathbb{Q}_p)\widehat{\otimes}_{\mathbb{Q}_p}D_{[\rho_1,\rho_2]}(G,K)	
\end{displaymath}
which corresponds directly to our interested situation, namely the noncommutative deformation of the usual Robba ring in one variable.

\begin{proposition}
For suitable radii $s,r,\rho_1,\rho_2\in p^\mathbb{Q}$ such that $s,r,\rho_1,\rho_2<1$ and:
\begin{displaymath}
D_{[s,r]}(\mathbb{Z}_p,\mathbb{Q}_p), D_{[\rho_1,\rho_2]}(G,K)	
\end{displaymath}
are noetherian as in \cite[Proposition 4.1]{4Zab1}. Then we have that the corresponding completed tensor product 
\begin{displaymath}
D_{[s,r]}(\mathbb{Z}_p,\mathbb{Q}_p)\widehat{\otimes}_{\mathbb{Q}_p}D_{[\rho_1,\rho_2]}(G,K)	
\end{displaymath}
is also noetherian.

\end{proposition}
	
\begin{remark}
Before the proof, let us mention that this actually provides a very typical example for situation required in \cref{assumption6.12}.	
\end{remark}

\begin{proof}
One can follow the proof of \cite[Proposition 4.1]{4Zab1} to do so by applying the criterion in \cite[Proposition I.7.1.2]{4LVO}, namely we look at the the graded ring of 
\begin{displaymath}
gr^\bullet_{\|.\|_{\otimes}} \left(D_{[s,r]}(\mathbb{Z}_p,\mathbb{Q}_p)\widehat{\otimes}_{\mathbb{Q}_p}D_{[\rho_1,\rho_2]}(G,K)\right),	
\end{displaymath}
defined by using the bounds of the norm $\|.\|_{\otimes}$, and if this is noetherian then we will have by \cite[Proposition I.7.1.2]{4LVO} that the ring $D_{[s,r]}(\mathbb{Z}_p,\mathbb{Q}_p)\widehat{\otimes}D_{[\rho_1,\rho_2]}(G,K)$ is also noetherian. 
It suffices to look at the corresponding graded ring of the dense subring of tensor product (see \cite[Lemma 4.3]{4ST1}) 
\begin{displaymath}
gr^\bullet_{\|.\|_{\otimes}} \left(D_{[s,r]}(\mathbb{Z}_p,\mathbb{Q}_p)\otimes_{\mathbb{Q}_p} D_{[\rho_1,\rho_2]}(G,K)\right).	
\end{displaymath}
We do not have the corresponding isomorphism in mind to split the corresponding components to get a product of two single graded rings but we have the corresponding surjective map (see \cite[I.6.13]{4LVO}):
\begin{displaymath}
gr^\bullet \left(D_{[s,r]}(\mathbb{Z}_p,\mathbb{Q}_p))\otimes_{gr^\bullet \mathbb{Q}_p} gr^\bullet( D_{[\rho_1,\rho_2]}(G,K)\right)\rightarrow gr_{\|.\|_{\otimes}}^\bullet \left(D_{[s,r]}(\mathbb{Z}_p,\mathbb{Q}_p)\otimes_{\mathbb{Q}_p} D_{[\rho_1,\rho_2]}(G,K)\right).	
\end{displaymath}
This is defined by:
\begin{displaymath}
\sigma_k(x)\otimes \sigma_l(y)\mapsto \sigma_{k+l}(x\otimes y),	
\end{displaymath}
where the principal symbol $\sigma_k(x)$ is defined to be $x+F^{k+} D_{[s,r]}(\mathbb{Z}_p,\mathbb{Q}_p)$ if we we have $x\in F^k D_{[s,r]}(\mathbb{Z}_p,\mathbb{Q}_p) \backslash F^{k+} D_{[s,r]}(\mathbb{Z}_p,\mathbb{Q}_p)$, in the corresponding algebraic microlocalization. The construction for $\sigma_l(y)$ and $ \sigma_{k+l}(x\otimes y)$ is the same. Therefore we look at then the corresponding product:
\begin{displaymath}
gr^\bullet \left(D_{[s,r]}(\mathbb{Z}_p,\mathbb{Q}_p)\right)\otimes gr^\bullet \left( D_{[\rho_1,\rho_2]}(G,K)\right)	
\end{displaymath}
The corresponding elements of $gr^\bullet \left( D_{[\rho_1,\rho_2]}(G,K)\right)$ are those Laurent polynomials over $gr^\bullet K=k[\pi_2,\pi_2^{-1}]$ with variables:
\begin{displaymath}
t_1,...,t_d	
\end{displaymath}
where $t_i=\sigma(b_i)/\pi_2^l$. The corresponding elements of $gr^\bullet \left(D_{[s,r]}(\mathbb{Z}_p,\mathbb{Q}_p) \right)$ are those Laurent polynomials over $gr^\bullet \mathbb{Q}_p=\mathbb{F}_p[\pi_1,\pi_1^{-1}]$ with variable:
\begin{displaymath}
t=\sigma(b)/\pi_1^{k}.
\end{displaymath}
We followed the notations in \cite[Proposition 4.1]{4Zab1}, where $b$ is the corresponding local coordinate for $D_{[s,r]}(\mathbb{Z}_p,\mathbb{Q}_p)$ while $\sigma$ represents the corresponding principal symbol in the corresponding nonarchimedean microlocalization, while we have that $b_1,...,b_d$ are local coordinates for $D_{[\rho_1,\rho_2]}(G,K)$. The idea is that this ring is actually noetherian then we have the ring
\begin{displaymath}
gr^\bullet \left(D_{[s,r]}(\mathbb{Z}_p,\mathbb{Q}_p)\otimes D_{[\rho_1,\rho_2]}(G,K)\right)	
\end{displaymath}
is also noetherian. Then one can endow 
\begin{displaymath}
D_{[s,r]}(\mathbb{Z}_p,\mathbb{Q}_p)\widehat{\otimes} D_{[\rho_1,\rho_2]}(G,K)
\end{displaymath}
with the corresponding filtration which will imply that the ring 
\begin{displaymath}
D_{[s,r]}(\mathbb{Z}_p,\mathbb{Q}_p)\widehat{\otimes}  D_{[\rho_1,\rho_2]}(G,K)
\end{displaymath}
is noetherian as well by \cite[Proposition I.7.1.2]{4LVO}. Indeed as in \cite[Proposition 4.1]{4Zab1} what we could do is to consider the corresponding projection to the corresponding Laurent polynomial ring over $k$ in the variables $\pi_2,t,t_1,...,t_d$ as in the proof of \cite[Proposition 4.1]{4Zab1}:
\begin{displaymath}
\mathbb{F}_p[\pi_1,\pi_1^{-1},t,t^{-1}]\otimes_{\mathbb{F}_p[\pi_1,\pi_1^{-1}]} k[\pi_2,\pi_2^{-1},t_1,t_1t_2^{-1},...,t_1t_d^{-1},t_2t_1^{-1},...,t_dt_1^{-1}]	
\end{displaymath}
which is noetherian, then the corresponding ring
\begin{displaymath}
gr^\bullet \left(D_{[s,r]}(\mathbb{Z}_p,\mathbb{Q}_p))\otimes_{gr^\bullet\mathbb{Q}_p} gr^\bullet(D_{[\rho_1,\rho_2]}(G,K)\right)	
\end{displaymath}
is noetherian since one can then lift the generators for the corresponding projection of any ideal through this construction as in the same way of \cite[Proposition 4.1]{4Zab1}.
\end{proof}

%While we conjecture that any affinoid algebra taking the form of 
%\begin{displaymath}
%A\{T_1/a_1,...,T_m/a_m,b_1/T_1,...b_m/T_m\}	
%\end{displaymath}
%with mutually commuting variables $T_1,...,T_m$ which are also assumed to be commuting with any element in $A$ is noetherian (note that the Hilbert basis theorem holds for $A[[T_1,...,T_m]]$), one can actually look at the corresponding ring of noncommutative Laurent series as in \cite{4ST1} and \cite{4Zab1} over the noncommutative polyannuli with respect to the $p$-adic Lie group:
%\begin{displaymath}
%\mathbb{Z}_p^{|I|}\times G	
%\end{displaymath}
%for suitable compact $p$-adic Lie group $G$. One can vary the corresponding annuli with respect to the factor $\mathbb{Z}_p^{|I|}$ and fix the corresponding noncommutative annuli with respect to the factor $G$. Then this will give us the corresponding examples in our situation.

\begin{setting}
Now under the corresponding noetherian assumption in \cref{assumption6.12}, we actually have that the corresponding multivariate Robba ring $\Pi_{\mathrm{an},r_{I,0},I,A}$	is Fr\'echet-Stein in the sense of \cite[Chapter 3]{4ST1}. 
\end{setting}

\begin{definition} \mbox{\bf{(After KPX \cite[Definition 2.1.3]{4KPX})}}\\
We define coherent sheaves of left modules and the global sections over the ring $\Pi_{\mathrm{an},r_{I_0},I,A}$ in the following way. First a coherent sheaf $\mathcal{M}$ over the relative multidimensional Robba ring $\Pi_{\mathrm{an},r_{I_0},I,A}$ is defined to be a family $(M_{[s_I,r_I]})_{[s_I,r_I]\subset (0,r_{I,0}]}$of left modules of finite type over each relative $\Pi_{[s_I,r_I],I,A}$ satisfying the following two conditions as in the more classical situation:
\begin{displaymath}
\Pi_{[s_I'',r_I''],I,A}\otimes_{\Pi_{[s'_I,r'_I],I,A}}M_{[s_I',r_I']}\overset{\sim}{\rightarrow} M_{[s''_I,r''_I]}	
\end{displaymath}
for any multi radii satisfying $0<s_I'\leq s_I''\leq r_I''\leq r_I'\leq r_{I,0}$, with furthermore 
\begin{displaymath}
\Pi_{[s'''_I,r'''_I],I,A}\otimes_{\Pi_{[s''_I,r''_I],I,A}}(\Pi_{[s''_I,r''_I],I,A}	\otimes_{\Pi_{[s'_I,r'_I],I,A}}M_{[s'_I,r'_I]}) \overset{\sim}{\rightarrow}\Pi_{[s'''_I,r'''_I],I,A},
\end{displaymath}
for any multi radii satisfying the following condition:
\begin{displaymath}
0<s_I'\leq s_I''\leq s_I'''\leq r_I'''\leq r_I''\leq r_I'\leq r_{I,0}.	
\end{displaymath}
Then we define the corresponding global section of the coherent sheaf $\mathcal{M}$, which is usually denoted by $M$ which is defined to be the following inverse limit:
\begin{displaymath}
M:=\varprojlim_{s_I\rightarrow 0_I^+} M_{[s_I,r_{I,0}]}.	
\end{displaymath}
We are going to follow \cite[Definition 2.1.3]{4KPX} to call any module defined over $\Pi_{\mathrm{an},r_{I_0},I,A}$ coadmissible it comes from a global section of a coherent sheaf $\mathcal{M}$ in the sense defined above. 
\end{definition}

\begin{remark}
Since we are in the framework of \cite[Chapter 3]{4ST1}, so now by \cite[Chapter 3, Theorem]{4ST1} for $r_{I,0}$ finite (namely for all $\alpha\in I$ the radius $r_{\alpha,0}$ is finite) we have that the corresponding global section in the previous definition is dense in each member participating in the coherent sheaf over some multi interval. We also have the same result for the situation where $r_{I,0}=\infty^{|I|}$ after inverting the corresponding variables in the multivariate Robba rings as in the commutative situation we considered in \cref{4prop2.9}. Also we have the higher vanishing of the corresponding derived inverse limit on coherent sheaves in the sense defined above.	
\end{remark}

\begin{definition} \mbox{\bf{(After KPX \cite[Definition 2.1.9]{4KPX})}} 
We first generalize the notion of admissible covering in our setting, which is higher dimensional generalization of the situation established in \cite[Definition 2.1.9]{4KPX}. For any covering $\{[s_I,r_I]\}$ of $(0,r_{I,0}]$, we are going to specify those admissible coverings if the given covering admits refinement finite locally (and each corresponding member in the covering has the corresponding interior part which is assumed to be not empty with respect to each $\alpha\in I$). Then it is very natural in our setting that we have the corresponding notations of $(m,n)$-finitely presentedness and $n$-finitely generatedness for any $m,n$ positive integers. Then based on these definitions we have the following generalization of the corresponding uniform finiteness. First we are going to call a coherent sheaf $\{M_{s_I,r_I}\}_{[s_I,r_I]}$ over $\Pi_{\mathrm{an},r_{I,0},I,A}$ uniformly $(m,n)$-finitely presented if there exists an admissible covering $\{[s_I,r_I]\}$ and a pair of positive integers $(m,n)$ such that each module $M_{[s_I,r_I]}$ defined over $\Pi_{\mathrm{an},[s_I,r_I],I,A}$ is $(m,n)$-finitely presented. Also we have the notion of uniformly $n$-finitely generated for any positive integer $n$ by defining that to mean under the existence of an admissible covering we have that each member $M_{[s_I,r_I]}$ in the family defined over $\Pi_{\mathrm{an},[s_I,r_I],I,A}$ is $n$-finitely generated. 
\end{definition}

\begin{lemma}\mbox{\bf{(After KPX \cite[Lemma 2.1.11]{4KPX})}} Suppose that for a coherent sheaf $\{M_{s_I,r_I}\}_{[s_I,r_I]}$ defined above and for an admissible covering $\{[s_I,r_I]\}$ defined above in our noncommutative setting we have that any module $M_{s_I,r_I}$ is finitely generated uniformly by finite many elements of global section. Then the global section will be also finitely generated by these set of elements as well. 
	
\end{lemma}

\begin{proof}
As in \cite[Lemma 2.1.11]{4KPX}, we use the set of elements to present a map from the free module generated by them to the global section, then after base change to each $\Pi_{\mathrm{an},[s_I,r_I],I,A}$ we do have that the map is corresponding presentation over $\Pi_{\mathrm{an},[s_I,r_I],I,A}$, then this will give rise to the corresponding presentation for the global section since we have the vanishing of the derived inverse limit functor.	
\end{proof}

\begin{lemma} \mbox{\bf{(After KPX \cite[Lemma 2.1.12]{4KPX})}} \label{lemma6.16}
Let $H$ be a noncommutative Banach ring, and consider some finite module $M$ over $H$ generated by $\{f_1,...,f_n\}$ and any another set $\{f'_1,...,f'_n\}$ of elements such that $f_i-f_i'$ could be made sufficiently small uniformly for each $i=1,...,n$. Then we have that the set $\{f'_1,...,f'_n\}$ could also finitely generate the module $M$.	
\end{lemma}

\begin{proof}
As in \cite[Lemma 2.1.12]{4KPX}, one looks at the corresponding matrix $G_{i,j}$ defined by $f_i-f_i'=\sum_{k}G_{i,k}f_k$. Then by using the open mapping theorem, one can control the norm of each entry $G_{i,j}$ to make them small enough below a uniform small bound, which makes the corresponding matrix $1+G$ invertible.	
\end{proof}

\begin{proposition} \mbox{\bf{(After KPX \cite[Proposition 2.1.13]{4KPX})}} \label{4prop1}
Assume now $I$ consists of exactly one index, consider any arbitrary coherent sheaf $\{M_{[s_I,r_I]}\}$ over $\Pi_{\mathrm{an},r_{I,0},I,A}$ with the global section $M$. Then we have the following corresponding statements:\\
I.	The coherent sheaf $\{M_{[s_I,r_I]}\}$ is uniformly finitely generated iff the global section $M$ is finitely generated;\\
II. The coherent sheaf $\{M_{[s_I,r_I]}\}$ is uniformly finitely presented iff the global section $M$ is finitely presented.
%III. The coherent sheaf $\{M_{[s_I,r_I]}\}$ is uniformly finitely presented vector bundle iff the global section $M$ is finite projective. (Here the vector bundles will be defined to be those corresponding families of flat modules as in the commutative case.) 
\end{proposition}

\begin{proof}
One could derive the the second statement by using the analog of \cref{lemma2.15} (see \cite[Lemma 4.54]{4L}) and the first statement. Therefore it suffices for us to prove the first statement. Indeed one could see that one direction in the statement is straightforward, therefore it suffices now to consider the other direction. So now we assume that the coherent sheaf $\{M_{[s_I,r_I]}\}$ is uniformly finitely generated for some $n$ for instance. Then we are going to show that one could find corresponding generators for the global section $M$ from the given uniform finitely generatedness. 

Now by definition we consider an admissible covering taking the form of $\{[s_I,r_I]\}$. In \cite[Proposition 2.1.13]{4KPX} the corresponding argument goes by using special decomposition of the covering in some refined way, mainly by reorganizing the initial admissible covering into two single ones, where each of them will be made into the one consisting of those well-established intervals with no overlap. We can also use instead the corresponding upgraded nice decomposition of the given admissible covering in our setting, which is extensively discussed in \cite[2.6.14-2.6.17]{4KL2}. Namely as in \cite[2.6.14-2.6.17]{4KL2} there is a chance for us to extract $2^{|I|}$ families $\{\{[s_{I,\delta_i},r_{I,\delta_i}]\}_{i=0,1,...}\}_{\{1_i\},\{2_i\},...,\{N_i\},}$ ($N=2^{|I|}$) of intervals, where there is a chance to have the situation where for each such family the corresponding intervals in it could be made to be disjoint in pairs.

Then for a chosen family $\{[s_{I,\delta_i},r_{I,\delta_i}]\}_{i=0,1,...}$, we now for each $i$ involved we consider the corresponding generating set $\mathbf{e}_{\delta_i,1},...,\mathbf{e}_{\delta_i,n}$ where $i=0,1,...$ and $n\in \{1,2,...,n\}$. Note that these should be generated over the integral rings, which is not the case when we work with the infinite radii. In the latter situation one should first invert $T_I$ then consider some large powers which could eliminate the different denominators. Then by the previous results we know that one could actually arrange these generators to be global sections by the density of the global sections in the sections over any chosen interval. Now we consider the following step where we for each $j=1,...,n$ form the following sum with coefficients to be chosen:
\begin{displaymath}
\mathbf{e}_j:=\sum_{i=0}^\infty a_{\delta_i,j}T_I^{\lambda_{I,\delta_i,j}}\mathbf{e}_{\delta_i,j}.
\end{displaymath}
Note that in our situation actually $I$ is a singleton, so we can actually regard $\delta_i$ as $2i$.

To check the suitable convergence in this situation we consider the following estimate with bounds given as below to be chosen from the conditions:
\begin{align}
\forall i''<\delta_i-1, \left\|a_{\delta_i,j}T^{\lambda_{I,\delta_i,j}}\mathbf{e}_{\delta_i,j}\right\|_{[s_{I,i''},r_{I,i''}]}\leq p^{-i},\\
\forall i'<i, \left\|a_{\delta_i,j}T^{\lambda_{I,\delta_i,j}}\mathbf{e}_{\delta_i,j}/(a_{\delta_{i'},j}T^{\lambda_{I,\delta_{i'},j}})\right\|_{[s_{I,\delta_{i'}},r_{I,\delta_{i'}}]} \leq \varepsilon_{\delta_{i'},j}\\
\forall i'<i, \left\|a_{\delta_{i'},j}T^{\lambda_{I,\delta_{i'},j}}\mathbf{e}_{\delta_{i'},j}/(a_{\delta_i,j}T^{\lambda_{I,\delta_i,j}})\right\|_{[s_{I,\delta_{i}},r_{I,\delta_{i}}]} \leq \varepsilon_{\delta_i,j}.	
\end{align}
As in the commutative situation in \cite[Proposition 2.1.13]{4KPX} in order to get the desired quantity expressed in the estimate above we seek the following situation from what we presented above:
\begin{align}
\forall i''<\delta_i-1, \left\|a_{\delta_i,j}T^{\lambda_{I,\delta_i,j}}\mathbf{e}_{\delta_i,j}\right\|_{[s_{I,i''},r_{I,i''}]}\leq p^{-i},\\
\forall i'<i, \left\|a_{\delta_i,j}T^{\lambda_{I,\delta_i,j}}\mathbf{e}_{\delta_i,j}\right\|_{[s_{I,\delta_{i'}},r_{I,\delta_{i'}}]} \leq \varepsilon_{\delta_{i'},j} \left\|a_{\delta_{i'},j}T^{\lambda_{I,\delta_{i'},j}}\right\|_{[s_{I,\delta_{i'}},r_{I,\delta_{i'}}]}\\
\forall i'<i, \left\|a_{\delta_{i'},j}T^{\lambda_{I,\delta_{i'},j}}\mathbf{e}_{\delta_{i'},j}\right\|_{[s_{I,\delta_i},r_{I,\delta_i}]} \leq \varepsilon_{\delta_i,j} \left\|a_{\delta_i,j}T^{\lambda_{I,\delta_i,j}}\right\|_{[s_{I,\delta_i},r_{I,\delta_i}]}.	
\end{align} 	
Then what we can achieve up to now is to consider sufficiently large powers on the both sides in order to make the inequalities hold, which is possible due to our initial choice on the intervals, where the first statement could also be proved in the similar way (also see \cite[Proposition 2.1.13]{4KPX}). To be more precise, one chooses the corresponding coefficients and the powers of the monomials by induction, where the initial coefficient for $i=0$ is set to be 1 with order 0. Then one enlarges the corresponding difference between the power of the monomials on the both sides of the inequalities above while evaluating the first and second inequalities at the left most radii and evaluating the last one at right most radii. Then actually this will prove the statement since up to some powers the difference between $\mathbf{e}_{\delta_i,j}$ and $\mathbf{e}_j$ will have sufficiently small norm, which could give rise to suitable finitely generatedness through the \cref{lemma6.16}. Finally one could then finish the whole proof in the odd situation in the same fashion.

\end{proof}

\begin{proposition} \mbox{\bf{(After KPX \cite[Proposition 2.2.7]{4KPX})}}\\
In the situation where $I$ consists of exactly one element. We have that natural functor from coherent $\varphi_I$-modules (in the sense of the previous \cref{definition6.4}) over $\Pi_{\mathrm{an},r_{I,0},I,A}(\pi_{K_I})$ to the corresponding coherent $\varphi_I$-bundles (in the sense of the previous \cref{definition6.4}) is an equivalence.	
\end{proposition}

\begin{proof}
See \cref{proposition3.5}, by using the previous proposition.	
\end{proof}

\begin{definition} \mbox{\bf{(After KPX \cite[Notation 4.1.2]{4KPX})}}
Now we consider the following derived categories. All modules will be assumed to be left. The first one is the sub-derived category consisting of all the objects in $\mathbb{D}_\mathrm{left}(A)$ which are quasi-isomorphic to those bounded above complexes of finite projective modules over the ring $A$. We denote this category (which is defined in the same way as in \cite[Notation 4.1.2]{4KPX}) by $\mathbb{D}^-_\mathrm{perf,left}(A)$. Over the larger ring $\Pi_{\mathrm{an},\infty,I,A}(\Gamma_{K_I})$ we also have the sub-derived category of $\mathbb{D}_\mathrm{left}(\Pi_{\mathrm{an},\infty,I,A}(\Gamma_{K_I}))$ consisting of all those objects in $\mathbb{D}_\mathrm{left}(\Pi_{\mathrm{an},\infty,I,A}(\Gamma_{K_I}))$ which are quasi-isomorphic to those bounded above complexes of finite projective modules now over the ring $\Pi_{\mathrm{an},\infty,I,A}(\Gamma_{K_I})$. We denote this by $\mathbb{D}^-_\mathrm{perf,left}(\Pi_{\mathrm{an},\infty,I,A}(\Gamma_{K_I}))$.
\end{definition}

\indent Then we will investigate the complexes attached to any finite projective or coherent $(\varphi_I,\Gamma_I)$-module $M$ over $\Pi_{\mathrm{an},\mathrm{con},I,A}(\pi_{K_I})$:
\begin{displaymath}
C^\bullet_{\varphi_I,\Gamma_I}(M), C^\bullet_{\psi_I,\Gamma_I}(M),C^\bullet_{\psi_I}(M),	
\end{displaymath}
which are living inside the corresponding derived categories:
\begin{displaymath}
\mathbb{D}_\mathrm{left}(A),\mathbb{D}_\mathrm{left}(A),\mathbb{D}_\mathrm{left}(\Pi_{\mathrm{an},\infty,I,A}(\Gamma_{K_I})).	
\end{displaymath}

\begin{conjecture} \label{conjecture5.22}
With the notation and setting up as above we have $h^\bullet(C_{\psi_I}(M))$ are coadmissible over $\Pi_{\mathrm{an},\infty,I,A}(\Gamma_{K_I})$.	
\end{conjecture}

%\begin{proof}
%By our development, see \cref{proposition4.7}.	
%\end{proof}

\begin{conjecture} \label{conjecture5.23}
With the notations above, we have that 
\begin{displaymath}
C^\bullet_{\varphi_I,\Gamma_{I}}(M)\in \mathbb{D}_\mathrm{perf,left}^-(A), C^\bullet_{\psi_I,\Gamma_{I}}(M)\in \mathbb{D}_\mathrm{perf,left}^-(A).
\end{displaymath}
\end{conjecture}

\begin{conjecture} \label{conjecture5.24}
Let $A_\infty(H)$ be the Fr\'echet-Stein algebra attached to any compact $p$-adic Lie group $H$ which could be written as inverse limit of noncommutative noetherian Banach affinoid algebras over $\mathbb{Q}_p$. Suppose $\mathcal{F}$ is a family of finite projective left modules over the Robba ring $\Pi_{\mathrm{an},\mathrm{con},I,A_\infty(H)}(\pi_{K_I})$, where we use the same notation to denote the global section and we assume that the global section is finite projective, carrying mutually commuting actions of $(\varphi_I,\Gamma_{K_I})$. Then we have that $C^\bullet_{\psi_I}(\mathcal{F})$ lives in $\mathbb{D}^\flat_{\mathrm{perf,left}}(\Pi_{\mathrm{an},\infty,I,A_\infty(H)}(\Gamma_{K_I}))$ and $C^\bullet_{\varphi_I,\Gamma_I}(\mathcal{F})$ lives in $\mathbb{D}^\flat_{\mathrm{perf,left}}(A_\infty(H))$. 
\end{conjecture}

%\begin{proof}
%By our development, see \cref{proposition4.9}.	
%\end{proof}

\begin{remark}
We conjecture that the corresponding Cartan-Serre approach holds in the noncommutative setting as in \cite[Lemma1.10]{4KL3}, which may possibly provide proof of the conjectures above.	
\end{remark}

\

This chapter is based on the following paper, where the author of this dissertation is the main author:
\begin{itemize}
\item Tong, Xin. "Analytic Geometry and Hodge-Frobenius Structure." arXiv preprint arXiv:2011.08358 (2020). 
\end{itemize}

\newpage

\bibliographystyle{ams}

\newpage\chapter{Analytic Geometry and Hodge-Frobenius Structure Continued}

\newpage\section{Introduction}

\subsection{Introduction and Summary}

\indent Generally speaking analytic geometry studies very general manifolds or varieties which are locally related to analytic functions of several variables. In complex analytic geometry, the special domains (like those classical domains) form a very important subject in the corresponding understanding of higher dimensional analytic geometry and their application to other subjects such as number theory. One definitely would like to understand the corresponding story in nonarchimedean geometry, namely what are the corresponding significant special domains. Another issue is that if we could have the chance to find the links with other geometries, such as formal ones and pure algebraic geometry. In the study of analytic geometry, certainly not only the analytic methods will be definitely applied, but also we have algebraic ones. Algebraic ones certainly work well in some very significant cases, for instance those in the Stein cases in classical complex geometry, quasi-Stein cases within the generality of Kedlaya-Liu \cite{5KL2} and Fr\'echet-Stein cases after Schneider-Teitelbaum \cite{5ST}.\\

\indent While the corresponding nonarchimedean analytic geometry has its own very significant geometric rights, many useful applications also exist extensively as well. For instance in \cite{5T1} we followed \cite{5CKZ} and \cite{5PZ} to have investigated some relative cohomologies corresponding to very specific multi Frobenius structures. The picture in \cite{5T1} was a very practical first step towards more general pictures and more interesting pictures.\\

\indent The tools in \cite{5T1} extend to our current paper since in more general framework of geometry as in \cite{5KL1} and \cite{5KL2} the methods and ideas could be applied directly to our situation with possible suitable more detailization. For instance in our situation we could consider nonnoetherian base rings. Note that this does not mean they are necessarily preperfectoid or perfect, since in our situation we can consider some very interesting mixed-type Robba rings in very general sense. For instance if we have two variables, we could have one living in some preperfectoid part while the other one living in the some rigid analytic affinoid component. Note that we do not have to maintain in the corresponding strictly analytic situation, especially if one would like to work with Kedlaya's reified adic spaces \cite{5Ked1}.\\

\indent We will also consider the corresponding generalized $B$-pairs after Berger and Nakamura. However for some deep study we could introduce some mixed-type Hodge structures. Namely we could introduce $(\varphi,\Gamma)$-$B$-objects. For instance in the situation where we have two factors we can defined the corresponding period rings:
\begin{displaymath}
\mathrm{B}_{\mathrm{dR},\mathbb{Q}_p}\widehat{\otimes}_{\mathbb{Q}_p}\Pi_{\mathrm{an},\mathrm{con},A},\mathrm{B}^+_{\mathrm{dR},\mathbb{Q}_p}\widehat{\otimes}_{\mathbb{Q}_p}\Pi_{\mathrm{an},\mathrm{con},A}	
\end{displaymath}
with finite projective objects defined over them, carrying one partial action from the Galois groups and carrying one partial action from the corresponding $(\varphi,\Gamma)$ operators.\\

\indent To make summary, we have:\\

\noindent A. Over really general Banach rings, we have defined many useful mixed type Robba rings by perfection along some partial variables. In fact we also have defined some $\infty$-period rings following the recent work from Bambozzi-Kremnizer \cite{5BK}. The corresponding issue we encountered is certainly some common issue around sheafiness in adic geometry in some very general classical sense. These directly will be expected to give the chance for one to compare Hodge structures over multivariate imperfect Robba rings with the corresponding Hodge structures over multivariate perfect Robba rings as those with accents $\breve{.}$ and $\widetilde{.}$ in \cite{5KL2}, see sections from 2-3. \\

\noindent B. Over really general Banach rings, we defined many useful  $(\varphi_I,\Gamma_I)$-modules over the mixed type Robba rings by perfection along some partial variables. This directly gives the chance for one to compare multivariate $(\varphi_I,\Gamma_I)$-modules over partially imperfect Robba rings with the corresponding multivariate $(\varphi_I,\Gamma_I)$-modules over partially perfect Robba rings as those with accents $\breve{.}$ and $\widetilde{.}$ in \cite{5KL2}, see section 4.\\

\noindent C. Over really general Banach rings over $\mathbb{Q}_p$, we defined many useful mixed type big period rings by taking product with Fontaine's de Rham period ring along some partial variables. First of all, we then immediately have the definition for some $B_I'$-$(\varphi_I,\Gamma_I)$-modules. Therefore relying on these rings (although they are very interesting in their rights) proved the equivalence between multivariate $B_I'$-pairs and multivariate $(\varphi_I,\Gamma_I)$-modules generalizing Berger's work and \cite{5KP}, see section 5.\\

\noindent D. In sections 6-8, we defined the corresponding $B_I'$-$(\varphi_I,\Gamma_I)$-cohomologies for $B_I'$-$(\varphi_I,\Gamma_I)$-modules, and we promote the equivalence to some quasi-isomorphisms within the derived category $\mathrm{D}^\flat(A)$, where $A$ is the base Banach relative ring. \\

\subsection{Comments on the Notation}

Our notations on the corresponding multivariate mixed type big Hodge structures are inspired by essentially some Langlands programs in $\ell$-adic situation such as in \cite{5VL}, namely the Drinfeld's Lemma. However the work \cite{5PZ} ad \cite{5CKZ} use $\Delta$ which is inspired by essentially some Langlands programs in $p$-adic situation such as that rooted in the work of Z\'abr\'adi and is related to reductive datum. We remind the readers that these are actually not the corresponding intervals for the Robba rings in order to eliminate the corresponding possible confusion. \\

\subsection{Future Study}

\indent The current geometric discussion covered in this paper is basically around the commutative analytic geometry. We have not made it to add the corresponding discussion on the noncommutative setting, but this will be pushed to our further study. Certainly the corresponding noncommutative deformation will require some further well-established understanding on the foundational issues, such as the corresponding noncommutative descent as in the commutative situation from \cite{5KL1} and \cite{5KL2}. Definitely any good understanding on the noncommutative settings of these sorts will be essential to the corresponding good understanding on noncommutative analytic geometry and noncommutative Tamagawa number conjectures after \cite{5FK1}, \cite{5BF1} and \cite{5BF2}.  \\

\indent Since our dreams will be really those where we can handle very general analytic spaces. Certainly adic spaces need extensively restrictive requirement on the sheafiness of the Banach rings. However this might be resolved completely by considering \cite{5BK} (or possibly equivalently the work of Clausen-Scholze \cite{5CS}). We have already defined many interesting $\infty$-analytic stacks and the $\infty$-Robba rings. We will study Kedlaya-Liu glueing on this level after \cite{5KL1} and \cite{5KL2} in future work.\\

\indent We have discussed the corresponding higher dimensional $B$-pairs by introducing the corresponding higher dimensional de Rham period rings. We certainly hope to amplify the corresponding discussion in more $p$-adic Hodge theoretic sense. However one could definitely maintain in the world of $(\varphi_I,\Gamma_I)$-modules in some very flexible way, literally after Berger \cite{5Ber2}. We believe that we have more rigidity in the current context. In fact, we should mention that this higher dimensionalization of \cite{5Ber2} is literally motivated by the work \cite{5KL1} and \cite{5KL2}, as well as the work \cite{5CKZ} and \cite{5PZ}.

%%\newpage

\newpage\section{Big Robba Rings over Rigid Analytic Affinoids and Fr\'echet Objects in Mixed-characteristic Case}

\subsection{Big Robba Rings over Rigid Analytic Affinoids}

\indent We follow \cite{5T1} to give the thorough definition and discussion of those very big period rings we will need in the further discussions in the following body of the paper, where we consider as in \cite{5T1} a finite set $I$ with some subset $J$. Our Robba rings in this paper will be depending on the $I$ and $J$ simultaneously.

\begin{definition}
Let $A$ be any affinoid algebra over $\mathbb{Q}_p$ in rigid analytic geometry. We consider the corresponding multi intervals $[\omega^{r_I},\omega^{s_I}]$. Recall we have the corresponding Robba rings defined in \cite[Definition 2.4]{5T1}:
\begin{displaymath}
\Pi_{[s_I,r_I],I,A}	
\end{displaymath}
which is defined to be the corresponding affinoid:
\begin{displaymath}
A\widehat{\otimes}_{\mathbb{Q}_p}\mathbb{Q}_p\{\omega^{r_1}/T_1,...,\omega^{r_I}/T_I,T_1/\omega^{s_1},...,T_I/\omega^{s_I}\}.	
\end{displaymath}
Then we have the corresponding rings:
\begin{displaymath}
\Pi_{\mathrm{an},r_I,I,A}:= \varprojlim_{s_I} \Pi_{[s_I,r_I],I,A}.			
\end{displaymath}
with 
\begin{displaymath}
\Pi_{\mathrm{an,con},I,A}:= \bigcup_{r_I}\varprojlim_{s_I} \Pi_{[s_I,r_I],I,A}.			
\end{displaymath}
However we will in this paper to consider some more complicated version of the rings. We will use some partial Frobenius to perfectize partially the rings defined above. Therefore we will in some more uniform way to denote the rings in the following different way:
\begin{align}
\Pi_{[s_I,r_I],I,I,\emptyset,A}:=\Pi_{[s_I,r_I],I,A}\\
\Pi_{\mathrm{an},r_I,I,I,\emptyset,A}:=\Pi_{\mathrm{an},r_I,I,A}\\
\Pi_{\mathrm{an,con},I,I,\emptyset,A}:=\Pi_{\mathrm{an,con},I,\emptyset,A}.	
\end{align}
	
\end{definition}

\indent Now we follow the idea in \cite[Definition 5.2.1]{5Ked2} to define some extended version of the rings. We will have the following rings to be:

\begin{align}
%\Pi_{[s_I,r_I],I,J,I\backslash J,A},\\	
\Pi_{[s_I,r_I],I,\breve{J},I\backslash J,A},\\	
\Pi_{[s_I,r_I],I,\widetilde{J},I\backslash J,A},\\
\Pi_{[s_I,r_I],I,J,\breve{I\backslash J},A},\\	
\Pi_{[s_I,r_I],I,\breve{J},\breve{I\backslash J},A},\\	
\Pi_{[s_I,r_I],I,\widetilde{J},\breve{I\backslash J},A},\\
\Pi_{[s_I,r_I],I,J,\widetilde{I\backslash J},A},\\	
\Pi_{[s_I,r_I],I,\breve{J},\widetilde{I\backslash J},A},\\	
\Pi_{[s_I,r_I],I,\widetilde{J},\widetilde{I\backslash J},A}.	
\end{align}

and 

\begin{align}
%\Pi_{\mathrm{an},r_I,I,J,I\backslash J,A},\\	
\Pi_{\mathrm{an},r_I,I,\breve{J},I\backslash J,A},\\	
\Pi_{\mathrm{an},r_I,I,\widetilde{J},I\backslash J,A},\\
\Pi_{\mathrm{an},r_I,I,J,\breve{I\backslash J},A},\\	
\Pi_{\mathrm{an},r_I,I,\breve{J},\breve{I\backslash J},A},\\	
\Pi_{\mathrm{an},r_I,I,\widetilde{J},\breve{I\backslash J},A},\\
\Pi_{\mathrm{an},r_I,I,J,\widetilde{I\backslash J},A},\\	
\Pi_{\mathrm{an},r_I,I,\breve{J},\widetilde{I\backslash J},A},\\	
\Pi_{\mathrm{an},r_I,I,\widetilde{J},\widetilde{I\backslash J},A}.	
\end{align}

and 

\begin{align}
%\Pi_{\mathrm{an},\mathrm{con},I,J,I\backslash J,A},\\	
\Pi_{\mathrm{an},\mathrm{con},I,\breve{J},I\backslash J,A},\\	
\Pi_{\mathrm{an},\mathrm{con},I,\widetilde{J},I\backslash J,A},\\
\Pi_{\mathrm{an},\mathrm{con},I,J,\breve{I\backslash J},A},\\	
\Pi_{\mathrm{an},\mathrm{con},I,\breve{J},\breve{I\backslash J},A},\\
\Pi_{\mathrm{an},\mathrm{con},I,\widetilde{J},\breve{I\backslash J},A},\\
\Pi_{\mathrm{an},\mathrm{con},I,J,\widetilde{I\backslash J},A},\\	
\Pi_{\mathrm{an},\mathrm{con},I,\breve{J},\widetilde{I\backslash J},A},\\	
\Pi_{\mathrm{an},\mathrm{con},I,\widetilde{J},\widetilde{I\backslash J},A}.	
\end{align}

\begin{definition}\mbox{\bf{(After Kedlaya-Liu, \cite[Definition 5.2.1]{5KL2})}}
We first define the corresponding first group of the rings. The corresponding rings in groups as mentioned above are defined by using the corresponding partial Frobenius $\varphi_1,...,\varphi_I$ and the corresponding Fr\'echet completion.	For the ring $\Pi_{[s_I,r_I],I,\breve{J},I\backslash J,A}$, this is defined by:
\begin{align}
\Pi_{[s_I,r_I],I,\breve{J},I\backslash J,A}:=\varinjlim_{n_\alpha\geq 0,\alpha  \in J}\prod_{\alpha\in J}\varphi_\alpha^{n_\alpha}\Pi_{[s_I,r_I],I,{J},I\backslash J,A}.
\end{align}
Note that for the corresponding rings getting involved in the corresponding definition above we consider the corresponding various Fr\'echet norms for each $t_I>0$:
\begin{displaymath}
\|.\|_{\prod_{\alpha\in J}\varphi_\alpha^{n_\alpha}\Pi_{[s_I,r_I],I,{J},I\backslash J,A},t_I}	
\end{displaymath}
Then we define the corresponding ring $\Pi_{[s_I,r_I],I,\widetilde{J},I\backslash J,A}$, this is defined by the following Fr\'echet completion process:
\begin{align}
\Pi_{[s_I,r_I],I,\widetilde{J},I\backslash J,A}:=\left(\varinjlim_{n_\alpha\geq 0,\alpha  \in J}\prod_{\alpha\in J}\varphi_\alpha^{n_\alpha}\Pi_{[s_I,r_I],I,{J},I\backslash J,A}\right)^\wedge_{\|.\|_{\prod_{\alpha\in J}\varphi_\alpha^{n_\alpha}\Pi_{[s_I,r_I],I,{J},I\backslash J,A},t_I},t_I\in [s_I,r_I]}.
\end{align}

\end{definition}

%%%%%%%%%%%%%%%%%%%%%%%%%%%%%%%%%%!!!!!!!!!!!!!!!!!!

\indent Then in the corresponding symmetric way we have the following definition:

\begin{definition}\mbox{\bf{(After Kedlaya-Liu, \cite[Definition 5.2.1]{5KL2})}}\\
For the ring $\Pi_{[s_I,r_I],I,{J},\breve{I\backslash J},A}$, this is defined by:
\begin{align}
\Pi_{[s_I,r_I],I,{J},\breve{I\backslash J},A}:=\varinjlim_{n_\alpha\geq 0,\alpha  \in I\backslash J}\prod_{\alpha\in I\backslash J}\varphi_\alpha^{n_\alpha}\Pi_{[s_I,r_I],I,{J},I\backslash J,A}.
\end{align}
Note that for the corresponding rings getting involved in the corresponding definition above we consider the corresponding various Fr\'echet norms for each $t_I>0$:
\begin{displaymath}
\|.\|_{\prod_{\alpha\in I\backslash J}\varphi_\alpha^{n_\alpha}\Pi_{[s_I,r_I],I,{J},I\backslash J,A},t_I}	
\end{displaymath}
Then we define the corresponding ring $\Pi_{[s_I,r_I],I,{J},\widetilde{I\backslash J},A}$, this is defined by the following Fr\'echet completion process:
\begin{align}
\Pi_{[s_I,r_I],I,{J},\widetilde{I\backslash J},A}:=\left(\varinjlim_{n_\alpha\geq 0,\alpha  \in I\backslash J}\prod_{\alpha\in I\backslash J}\varphi_\alpha^{n_\alpha}\Pi_{[s_I,r_I],I,\breve{J},\breve{I\backslash J},A}\right)^\wedge_{\|.\|_{\prod_{\alpha\in I}\varphi_\alpha^{n_\alpha}\Pi_{[s_I,r_I],I,{J},I\backslash J,A},t_I},t_I\in [s_I,r_I]}.
\end{align}

\end{definition}

\indent Then we do the following one:

\begin{definition}\mbox{\bf{(After Kedlaya-Liu, \cite[Definition 5.2.1]{5KL2})}}\\
For the ring $\Pi_{[s_I,r_I],I,\breve{J},\breve{I\backslash J},A}$, this is defined by:
\begin{align}
\Pi_{[s_I,r_I],I,\breve{J},\breve{I\backslash J},A}:=\varinjlim_{n_\alpha\geq 0,\alpha  \in I}\prod_{\alpha\in I}\varphi_\alpha^{n_\alpha}\Pi_{[s_I,r_I],I,{J},I\backslash J,A}.
\end{align}
Note that for the corresponding rings getting involved in the corresponding definition above we consider the corresponding various Fr\'echet norms for each $t_I>0$:
\begin{displaymath}
\|.\|_{\prod_{\alpha\in I}\varphi_\alpha^{n_\alpha}\Pi_{[s_I,r_I],I,{J},I\backslash J,A},t_I}	
\end{displaymath}
Then we define the corresponding ring $\Pi_{[s_I,r_I],I,\widetilde{J},\widetilde{I\backslash J},A}$, this is defined by the following Fr\'echet completion process:
\begin{align}
\Pi_{[s_I,r_I],I,\widetilde{J},\widetilde{I\backslash J},A}:=\left(\varinjlim_{n_\alpha\geq 0,\alpha  \in I}\prod_{\alpha\in I}\varphi_\alpha^{n_\alpha}\Pi_{[s_I,r_I],I,{J},I\backslash J,A}\right)^\wedge_{\|.\|_{\prod_{\alpha\in I}\varphi_\alpha^{n_\alpha}\Pi_{[s_I,r_I],I,{J},I\backslash J,A},t_I},t_I\in [s_I,r_I]}.
\end{align}

\end{definition}

\indent Then we consider the following definition building on the definitions above:

\begin{definition}\mbox{\bf{(After Kedlaya-Liu, \cite[Definition 5.2.1]{5KL2})}}\\
For the ring $\Pi_{[s_I,r_I],I,\widetilde{J},\breve{I\backslash J},A}$, this is defined by:
\begin{align}
\Pi_{[s_I,r_I],I,\widetilde{J},\breve{I\backslash J},A}:=\varinjlim_{n_\alpha\geq 0,\alpha  \in I\backslash J}\prod_{\alpha\in I\backslash J}\varphi_\alpha^{n_\alpha}\Pi_{[s_I,r_I],I,\widetilde{J},I\backslash J,A}.
\end{align}
For the ring $\Pi_{[s_I,r_I],I,\breve{J},\widetilde{I\backslash J},A}$, this is defined by:
\begin{align}
\Pi_{[s_I,r_I],I,\breve{J},\widetilde{I\backslash J},A}:=\varinjlim_{n_\alpha\geq 0,\alpha  \in J}\prod_{\alpha\in J}\varphi_\alpha^{n_\alpha}\Pi_{[s_I,r_I],I,{J},\widetilde{I\backslash J},A}.
\end{align}

\end{definition}

\indent Then we have the following definitions:

\begin{definition} \mbox{\bf{(After Kedlaya-Liu, \cite[Definition 5.2.1]{5KL2})}}
\begin{align}
%\Pi_{\mathrm{an},r_I,I,J,I\backslash J,A},\\	
\Pi_{\mathrm{an},r_I,I,\breve{J},I\backslash J,A}:=\varprojlim_{s_I}\Pi_{[s_I,r_I],I,\breve{J},I\backslash J,A},\\	
\Pi_{\mathrm{an},r_I,I,\widetilde{J},I\backslash J,A}:=\varprojlim_{s_I} \Pi_{[s_I,r_I],I,\widetilde{J},I\backslash J,A},\\
\Pi_{\mathrm{an},r_I,I,J,\breve{I\backslash J},A}:=\varprojlim_{s_I}\Pi_{[s_I,r_I],I,J,\breve{I\backslash J},A},\\	
\Pi_{\mathrm{an},r_I,I,\breve{J},\breve{I\backslash J},A}:=\varprojlim_{s_I} \Pi_{[s_I,r_I],I,\breve{J},\breve{I\backslash J},A},\\	
\Pi_{\mathrm{an},r_I,I,\widetilde{J},\breve{I\backslash J},A}:=\varprojlim_{s_I} \Pi_{[s_I,r_I],I,\widetilde{J},\breve{I\backslash J},A},\\
\Pi_{\mathrm{an},r_I,I,J,\widetilde{I\backslash J},A}:=\varprojlim_{s_I} \Pi_{[s_I,r_I],I,J,\widetilde{I\backslash J},A},\\	
\Pi_{\mathrm{an},r_I,I,\breve{J},\widetilde{I\backslash J},A}:=\varprojlim_{s_I} \Pi_{[s_I,r_I],I,\breve{J},\widetilde{I\backslash J},A},\\	
\Pi_{\mathrm{an},r_I,I,\widetilde{J},\widetilde{I\backslash J},A}:=\varprojlim_{s_I} \Pi_{[s_I,r_I],I,\widetilde{J},\widetilde{I\backslash J},A}.	
\end{align}

and 

\begin{align}
%\Pi_{\mathrm{an},\mathrm{con},I,J,I\backslash J,A},\\	
\Pi_{\mathrm{an},\mathrm{con},I,\breve{J},I\backslash J,A}:=\varinjlim_{r_I}\varprojlim_{s_I}\Pi_{[s_I,r_I],I,\breve{J},I\backslash J,A},\\	
\Pi_{\mathrm{an},\mathrm{con},I,\widetilde{J},I\backslash J,A}:=\varinjlim_{r_I}\varprojlim_{s_I} \Pi_{[s_I,r_I],I,\widetilde{J},I\backslash J,A},\\
\Pi_{\mathrm{an},\mathrm{con},I,J,\breve{I\backslash J},A}:=\varinjlim_{r_I}\varprojlim_{s_I}\Pi_{[s_I,r_I],I,J,\breve{I\backslash J},A},\\	
\Pi_{\mathrm{an},\mathrm{con},I,\breve{J},\breve{I\backslash J},A}:=\varinjlim_{r_I}\varprojlim_{s_I} \Pi_{[s_I,r_I],I,\breve{J},\breve{I\backslash J},A},\\
\Pi_{\mathrm{an},\mathrm{con},I,\widetilde{J},\breve{I\backslash J},A}:=\varinjlim_{r_I}\varprojlim_{s_I} \Pi_{[s_I,r_I],I,\widetilde{J},\breve{I\backslash J},A},\\
\Pi_{\mathrm{an},\mathrm{con},I,J,\widetilde{I\backslash J},A}:=\varinjlim_{r_I}\varprojlim_{s_I} \Pi_{[s_I,r_I],I,J,\widetilde{I\backslash J},A},\\	
\Pi_{\mathrm{an},\mathrm{con},I,\breve{J},\widetilde{I\backslash J},A}:=\varinjlim_{r_I}\varprojlim_{s_I} \Pi_{[s_I,r_I],I,\breve{J},\widetilde{I\backslash J},A},\\	
\Pi_{\mathrm{an},\mathrm{con},I,\widetilde{J},\widetilde{I\backslash J},A}:=\varinjlim_{r_I}\varprojlim_{s_I} \Pi_{[s_I,r_I],I,\widetilde{J},\widetilde{I\backslash J},A}.\\	
\end{align}	
\end{definition}

%%%%%%%%%%%%%%%%%%%%%%%%%%%%%%%%!!!!!!!!!!!!!!

\begin{setting}
The corresponding construction we established above is also motivated from the corresponding construction of multivariate Robba rings in \cite{5Ked2}.\\	
\end{setting}

\subsection{Fr\'echet Objects}

\indent Now we perform some of the corresponding construction parallel to the corresponding Fr\'echet-Stein construction used originally in \cite{5KPX}, which is essentially developed in \cite[Section 2.6]{5KL2}. In effect, the corresponding construction was already used in \cite{5T1}, although the paper behaves as if we are not fully after \cite{5KL2}. The reason of such impression left for the readers is essentially due to the fact that we are in the noetherian situation. Now we are not definitely in the noetherian situation, but \cite[Section 2.6]{5KL2} has already tackled this issue.

\begin{setting} \mbox{\bf{(After Kedlaya-Liu, \cite[Definition 2.6.2]{5KL2})}}
Following \cite[Definition 2.6.2]{5KL2} we call a Banach uniform adic ring $(R,R^+)$ quasi-Stein if it could be written as the following inverse limit with the corresponding transition map of dense image for some inverse system $\{\alpha\}$ of indexes:
\begin{displaymath}
(R,R^+):=\varprojlim_{\alpha}  (R_\alpha,R_\alpha^+).	
\end{displaymath}
And we call the ring ind-Fr\'echet-Stein if it is further could be written the following injective-projective limit:
\begin{displaymath}
(R,R^+):=\varinjlim_{\alpha'}\varprojlim_{\alpha}(R_{\alpha,\alpha'} ,R_{\alpha,\alpha'}^+).	
\end{displaymath}
	
\end{setting}

\indent After this axiomization we could now study the corresponding sheaves over these rings. In what follow, we consider the corresponding radii living in the set of all the rational numbers.

\begin{definition} \mbox{\bf{(After KPX, \cite[Definition 2.1.3]{5KPX})}}
Over $\Pi_{\mathrm{an},r_{I,0},I,\breve{J},I\backslash J,A}$ we define the corresponding stably pseudocoherent sheaves to mean a collection of stably pseudocoherent modules $(M_{[s_I,r_I]})$	over each $\Pi_{[s_I,r_I],I,\breve{J},I\backslash J,A}$ satisfying the corresponding compatibility condition and the obvious cocycle condition with respect to the family of the corresponding multi-intervals $\{[s_I,r_I]\}$.
\end{definition}

\begin{definition} \mbox{\bf{(After KPX, \cite[Definition 2.1.3]{5KPX})}}
Over $\Pi_{\mathrm{an},r_{I,0},I,\widetilde{J},I\backslash J,A}$ we define the corresponding stably pseudocoherent sheaves to mean a collection of stably pseudocoherent modules $(M_{[s_I,r_I]})$	over each $\Pi_{[s_I,r_I],I,\widetilde{J},I\backslash J,A}$ satisfying the corresponding compatibility condition and the obvious cocycle condition with respect to the family of the corresponding multi-intervals $\{[s_I,r_I]\}$.
\end{definition}

\begin{definition} \mbox{\bf{(After KPX, \cite[Definition 2.1.3]{5KPX})}}
Over $\Pi_{\mathrm{an},r_{I,0},I,{J},\breve{I\backslash J},A}$ we define the corresponding stably pseudocoherent sheaves to mean a collection of stably pseudocoherent modules $(M_{[s_I,r_I]})$	over each $\Pi_{[s_I,r_I],I,{J},\breve{I\backslash J},A}$ satisfying the corresponding compatibility condition and the obvious cocycle condition with respect to the family of the corresponding multi-intervals $\{[s_I,r_I]\}$.
\end{definition}

\begin{definition} \mbox{\bf{(After KPX, \cite[Definition 2.1.3]{5KPX})}}
Over $\Pi_{\mathrm{an},r_{I,0},I,\breve{J},\breve{I\backslash J},A}$ we define the corresponding stably pseudocoherent sheaves to mean a collection of stably pseudocoherent modules $(M_{[s_I,r_I]})$	over each $\Pi_{[s_I,r_I],I,\breve{J},\breve{I\backslash J},A}$ satisfying the corresponding compatibility condition and the obvious cocycle condition with respect to the family of the corresponding multi-intervals $\{[s_I,r_I]\}$.
\end{definition}

\begin{definition} \mbox{\bf{(After KPX, \cite[Definition 2.1.3]{5KPX})}}
Over $\Pi_{\mathrm{an},r_{I,0},I,\widetilde{J},\breve{I\backslash J},A}$ we define the corresponding stably pseudocoherent sheaves to mean a collection of stably pseudocoherent modules $(M_{[s_I,r_I]})$ over each $\Pi_{[s_I,r_I],I,\widetilde{J},\breve{I\backslash J},A}$ satisfying the corresponding compatibility condition and the obvious cocycle condition with respect to the family of the corresponding multi-intervals $\{[s_I,r_I]\}$.
\end{definition}

\begin{definition} \mbox{\bf{(After KPX, \cite[Definition 2.1.3]{5KPX})}}
Over $\Pi_{\mathrm{an},r_{I,0},I,{J},\widetilde{I\backslash J},A}$ we define the corresponding stably pseudocoherent sheaves to mean a collection of stably pseudocoherent modules $(M_{[s_I,r_I]})$	over each $\Pi_{[s_I,r_I],I,{J},\widetilde{I\backslash J},A}$ satisfying the corresponding compatibility condition and the obvious cocycle condition with respect to the family of the corresponding multi-intervals $\{[s_I,r_I]\}$.
\end{definition}

\begin{definition} \mbox{\bf{(After KPX, \cite[Definition 2.1.3]{5KPX})}}
Over $\Pi_{\mathrm{an},r_{I,0},I,\breve{J},\widetilde{I\backslash J},A}$ we define the corresponding stably pseudocoherent sheaves to mean a collection of stably pseudocoherent modules $(M_{[s_I,r_I]})$ over each $\Pi_{[s_I,r_I],I,\breve{J},\widetilde{I\backslash J},A}$ satisfying the corresponding compatibility condition and the obvious cocycle condition with respect to the family of the corresponding multi-intervals $\{[s_I,r_I]\}$.
\end{definition}

\begin{definition} \mbox{\bf{(After KPX, \cite[Definition 2.1.3]{5KPX})}}
Over $\Pi_{\mathrm{an},r_{I,0},I,\widetilde{J},\widetilde{I\backslash J},A}$ we define the corresponding stably pseudocoherent sheaves to mean a collection of stably pseudocoherent modules $(M_{[s_I,r_I]})$ over each $\Pi_{[s_I,r_I],I,\widetilde{J},\widetilde{I\backslash J},A}$ satisfying the corresponding compatibility condition and the obvious cocycle condition with respect to the family of the corresponding multi-intervals $\{[s_I,r_I]\}$.
\end{definition}

\begin{remark}
There is some overlap and repeating on some objects within the eight categories. Therefore in the following we are going to then work with only 1st, 2nd, 4th, 5th, 8th categories. 	
\end{remark}

\begin{proposition} \mbox{\bf{(After Kedlaya-Liu, \cite[Corollary 2.6.8]{5KL2})}}
The corresponding pseudocoherent finitely projective bundles over
\begin{align}
\Pi_{\mathrm{an},r_{I,0},I,\breve{J},I\backslash J,A},\\
\Pi_{\mathrm{an},r_{I,0},I,\widetilde{J},I\backslash J,A},\\
\Pi_{\mathrm{an},r_{I,0},I,\breve{J},\breve{I\backslash J},A},\\
\Pi_{\mathrm{an},r_{I,0},I,\widetilde{J},\breve{I\backslash J},A},\\
\Pi_{\mathrm{an},r_{I,0},I,\widetilde{J},\widetilde{I\backslash J},A}
\end{align}
defined as above have global sections which are finite projective if and only if the global sections are finitely generated. 	
\end{proposition}

\begin{proof}
Just apply \cite[Corollary 2.6.8]{5KL2}.	
\end{proof}

\begin{proposition} \mbox{\bf{(After Kedlaya-Liu, \cite[Proposition 2.6.17]{5KL2})}} \label{5proposition2.18}
The corresponding pseudocoherent sheaves over
\begin{align}
\Pi_{\mathrm{an},r_{I,0},I,\breve{J},I\backslash J,A},\\
\Pi_{\mathrm{an},r_{I,0},I,\widetilde{J},I\backslash J,A},\\
\Pi_{\mathrm{an},r_{I,0},I,\breve{J},\breve{I\backslash J},A},\\
\Pi_{\mathrm{an},r_{I,0},I,\widetilde{J},\breve{I\backslash J},A},\\
\Pi_{\mathrm{an},r_{I,0},I,\widetilde{J},\widetilde{I\backslash J},A}
\end{align}
defined as above have global sections which are finitely generated as long as we have the uniform bound on the rank of the bundles over quasi-compacts with respect to closed multi-intervals taking the general form of $[s_I,r_I]$. 	
\end{proposition}

\begin{proof}
This is actually a direct consequence of \cite[Proposition 2.6.17]{5KL2}, where the space $(0,r_{I,0}]$ admits $2^{|I|}$-uniform covering.	
\end{proof}

\begin{proposition}\mbox{}\\
 \mbox{\bf{(After Kedlaya-Liu, \cite[Corollary 2.6.8, Proposition 2.6.17]{5KL2})}} \label{proposition2.19}
The corresponding pseudocoherent sheaves over
\begin{align}
\Pi_{\mathrm{an},r_{I,0},I,\breve{J},I\backslash J,A},\\
\Pi_{\mathrm{an},r_{I,0},I,\widetilde{J},I\backslash J,A},\\
\Pi_{\mathrm{an},r_{I,0},I,\breve{J},\breve{I\backslash J},A},\\
\Pi_{\mathrm{an},r_{I,0},I,\widetilde{J},\breve{I\backslash J},A},\\
\Pi_{\mathrm{an},r_{I,0},I,\widetilde{J},\widetilde{I\backslash J},A}
\end{align}
defined as above have global sections which are finite projective as long as we have the uniform bound on the rank of the bundle over each quasi-compact with respect to each closed multi-interval $[s_I,r_I]$, and we have that the sheaves admits section actually finite projective over each quasi-compact with respect to each closed multi-interval $[s_I,r_I]$. 	
\end{proposition}

\begin{proof}
This is actually a direct corollary of the previous two propositions. 	
\end{proof}

%%\newpage

\newpage\section{Big Robba Rings over General Banach Affinoids and Fr\'echet Objects in Mixed-characteristic Case}

\subsection{Big Robba Rings over General Banach Affinoids}

\indent Since we have already considered the corresponding foundation from \cite{5KL2} on the quasi-Stein nonnoetherian adic Banach uniform algebra over $\mathbb{Q}_p$, we hope then now study more general $p$-adic analysis of several variables. Certainly as mentioned in \cite{5KPX} one could carry some strongly noetherian coefficients, where everything is sheafy, but one might be very curious about the situation where we do not have so strong condition on the noetherianness. Actually we could then apply the derived analytic geometry in \cite{5BK}.

\begin{definition}
Let $A$ be any commutative Banach algebra over $\mathbb{Q}_p$. We consider the corresponding multi intervals $[\omega^{r_I},\omega^{s_I}]$. We have the corresponding Robba rings defined as in \cite[Definition 2.4]{5T1}:
\begin{displaymath}
\Pi_{[s_I,r_I],I,A}	
\end{displaymath}
which is defined to be the corresponding affinoid:
\begin{displaymath}
A\widehat{\otimes}_{\mathbb{Q}_p}\mathbb{Q}_p\{\omega^{r_1}/T_1,...,\omega^{r_I}/T_I,T_1/\omega^{s_1},...,T_I/\omega^{s_I}\}.	
\end{displaymath}
Then we have the corresponding rings:
\begin{displaymath}
\Pi_{\mathrm{an},r_I,I,A}:= \varprojlim_{s_I} \Pi_{[s_I,r_I],I,A}.			
\end{displaymath}
with 
\begin{displaymath}
\Pi_{\mathrm{an,con},I,A}:= \bigcup_{r_I}\varprojlim_{s_I} \Pi_{[s_I,r_I],I,A}.			
\end{displaymath}
However we will in this paper to consider some more complicated version of the rings. We will use some partial Frobenius to perfectize partially the rings defined above. Therefore we will in some more uniform way to denote the rings in the following different way:
\begin{align}
\Pi_{[s_I,r_I],I,I,\emptyset,A}:=\Pi_{[s_I,r_I],I,A}\\
\Pi_{\mathrm{an},r_I,I,I,\emptyset,A}:=\Pi_{\mathrm{an},r_I,I,A}\\
\Pi_{\mathrm{an,con},I,I,\emptyset,A}:=\Pi_{\mathrm{an,con},I,\emptyset,A}.	
\end{align}
	
\end{definition}

\indent Now we follow the idea in \cite[Definition 5.2.1]{5Ked2} to define some extended version of the rings. We will have the following rings to be:

\begin{align}
%\Pi_{[s_I,r_I],I,J,I\backslash J,A},\\	
\Pi_{[s_I,r_I],I,?,?',A},?=J,\widetilde{J},\breve{J},?'=I\backslash J,\widetilde{I\backslash J},\breve{I\backslash J},\\		
\end{align}

and 

\begin{align}
%\Pi_{\mathrm{an},r_I,I,J,I\backslash J,A},\\	
\Pi_{\mathrm{an},r_I,I,?,?',A},?=J,\widetilde{J},\breve{J},?'=I\backslash J,\widetilde{I\backslash J},\breve{I\backslash J},\\ 	
\end{align}

and 

\begin{align}
%\Pi_{\mathrm{an},\mathrm{con},I,J,I\backslash J,A},\\	
\Pi_{\mathrm{an},\mathrm{con},I,?,?',A},?=J,\widetilde{J},\breve{J},?'=I\backslash J,\widetilde{I\backslash J},\breve{I\backslash J}.\\	
\end{align}

\begin{definition}\mbox{\bf{(After Kedlaya-Liu, \cite[Definition 5.2.1]{5KL2})}}
We first define the corresponding first group of the rings. The corresponding rings in groups as mentioned above are defined by using the corresponding partial Frobenius $\varphi_1,...,\varphi_I$ and the corresponding Fr\'echet completion.	For the ring $\Pi_{[s_I,r_I],I,\breve{J},I\backslash J,A}$, this is defined by:
\begin{align}
\Pi_{[s_I,r_I],I,\breve{J},I\backslash J,A}:=\varinjlim_{n_\alpha\geq 0,\alpha  \in J}\prod_{\alpha\in J}\varphi_\alpha^{n_\alpha}\Pi_{[s_I,r_I],I,{J},I\backslash J,A}.
\end{align}
Note that for the corresponding rings getting involved in the corresponding definition above we consider the corresponding various Fr\'echet norms for each $t_I>0$:
\begin{displaymath}
\|.\|_{\prod_{\alpha\in J}\varphi_\alpha^{n_\alpha}\Pi_{[s_I,r_I],I,{J},I\backslash J,A},t_I}	
\end{displaymath}
Then we define the corresponding ring $\Pi_{[s_I,r_I],I,\widetilde{J},I\backslash J,A}$, this is defined by the following Fr\'echet completion process:
\begin{align}
\Pi_{[s_I,r_I],I,\widetilde{J},I\backslash J,A}:=\left(\varinjlim_{n_\alpha\geq 0,\alpha  \in J}\prod_{\alpha\in J}\varphi_\alpha^{n_\alpha}\Pi_{[s_I,r_I],I,{J},I\backslash J,A}\right)^\wedge_{\|.\|_{\prod_{\alpha\in J}\varphi_\alpha^{n_\alpha}\Pi_{[s_I,r_I],I,{J},I\backslash J,A},t_I},t_I\in [s_I,r_I]}.
\end{align}

\end{definition}

\indent Then in the corresponding symmetric way we have the following definition:

\begin{definition}\mbox{\bf{(After Kedlaya-Liu, \cite[Definition 5.2.1]{5KL2})}}\\
For the ring $\Pi_{[s_I,r_I],I,{J},\breve{I\backslash J},A}$, this is defined by:
\begin{align}
\Pi_{[s_I,r_I],I,{J},\breve{I\backslash J},A}:=\varinjlim_{n_\alpha\geq 0,\alpha  \in I\backslash J}\prod_{\alpha\in I\backslash J}\varphi_\alpha^{n_\alpha}\Pi_{[s_I,r_I],I,{J},I\backslash J,A}.
\end{align}
Note that for the corresponding rings getting involved in the corresponding definition above we consider the corresponding various Fr\'echet norms for each $t_I>0$:
\begin{displaymath}
\|.\|_{\prod_{\alpha\in I\backslash J}\varphi_\alpha^{n_\alpha}\Pi_{[s_I,r_I],I,{J},I\backslash J,A},t_I}	
\end{displaymath}
Then we define the corresponding ring $\Pi_{[s_I,r_I],I,{J},\widetilde{I\backslash J},A}$, this is defined by the following Fr\'echet completion process:
\begin{align}
\Pi_{[s_I,r_I],I,{J},\widetilde{I\backslash J},A}:=\left(\varinjlim_{n_\alpha\geq 0,\alpha  \in I\backslash J}\prod_{\alpha\in I\backslash J}\varphi_\alpha^{n_\alpha}\Pi_{[s_I,r_I],I,\breve{J},\breve{I\backslash J},A}\right)^\wedge_{\|.\|_{\prod_{\alpha\in I}\varphi_\alpha^{n_\alpha}\Pi_{[s_I,r_I],I,{J},I\backslash J,A},t_I},t_I\in [s_I,r_I]}.
\end{align}

\end{definition}

\indent Then we do the following one:

\begin{definition}\mbox{\bf{(After Kedlaya-Liu, \cite[Definition 5.2.1]{5KL2})}}\\
For the ring $\Pi_{[s_I,r_I],I,\breve{J},\breve{I\backslash J},A}$, this is defined by:
\begin{align}
\Pi_{[s_I,r_I],I,\breve{J},\breve{I\backslash J},A}:=\varinjlim_{n_\alpha\geq 0,\alpha  \in I}\prod_{\alpha\in I}\varphi_\alpha^{n_\alpha}\Pi_{[s_I,r_I],I,{J},I\backslash J,A}.
\end{align}
Note that for the corresponding rings getting involved in the corresponding definition above we consider the corresponding various Fr\'echet norms for each $t_I>0$:
\begin{displaymath}
\|.\|_{\prod_{\alpha\in I}\varphi_\alpha^{n_\alpha}\Pi_{[s_I,r_I],I,{J},I\backslash J,A},t_I}	
\end{displaymath}
Then we define the corresponding ring $\Pi_{[s_I,r_I],I,\widetilde{J},\widetilde{I\backslash J},A}$, this is defined by the following Fr\'echet completion process:
\begin{align}
\Pi_{[s_I,r_I],I,\widetilde{J},\widetilde{I\backslash J},A}:=\left(\varinjlim_{n_\alpha\geq 0,\alpha  \in I}\prod_{\alpha\in I}\varphi_\alpha^{n_\alpha}\Pi_{[s_I,r_I],I,{J},\breve{I\backslash J},A}\right)^\wedge_{\|.\|_{\prod_{\alpha\in I}\varphi_\alpha^{n_\alpha}\Pi_{[s_I,r_I],I,{J},I\backslash J,A},t_I},t_I\in [s_I,r_I]}.
\end{align}

\end{definition}

\indent Then we consider the following definition building on the definitions above:

\begin{definition}\mbox{\bf{(After Kedlaya-Liu, \cite[Definition 5.2.1]{5KL2})}}\\
For the ring $\Pi_{[s_I,r_I],I,\widetilde{J},\breve{I\backslash J},A}$, this is defined by:
\begin{align}
\Pi_{[s_I,r_I],I,\widetilde{J},\breve{I\backslash J},A}:=\varinjlim_{n_\alpha\geq 0,\alpha  \in I\backslash J}\prod_{\alpha\in I\backslash J}\varphi_\alpha^{n_\alpha}\Pi_{[s_I,r_I],I,\widetilde{J},I\backslash J,A}.
\end{align}
For the ring $\Pi_{[s_I,r_I],I,\breve{J},\widetilde{I\backslash J},A}$, this is defined by:
\begin{align}
\Pi_{[s_I,r_I],I,\breve{J},\widetilde{I\backslash J},A}:=\varinjlim_{n_\alpha\geq 0,\alpha  \in J}\prod_{\alpha\in J}\varphi_\alpha^{n_\alpha}\Pi_{[s_I,r_I],I,{J},\widetilde{I\backslash J},A}.
\end{align}

\end{definition}

\indent Then we have the following definitions:

\begin{definition} \mbox{\bf{(After Kedlaya-Liu, \cite[Definition 5.2.1]{5KL2})}}
\begin{align}
%\Pi_{\mathrm{an},r_I,I,J,I\backslash J,A},\\	
\Pi_{\mathrm{an},r_I,I,?,?',A}:=\varprojlim_{s_I}\Pi_{[s_I,r_I],I,?,?',A},?=J,\widetilde{J},\breve{J},?'=I\backslash J,\widetilde{I\backslash J},\breve{I\backslash J},\\	
\end{align}

and 

\begin{align}
%\Pi_{\mathrm{an},\mathrm{con},I,J,I\backslash J,A},\\	
\Pi_{\mathrm{an},\mathrm{con},I,?,?',A}:=\varinjlim_{r_I}\varprojlim_{s_I}\Pi_{[s_I,r_I],I,?,?',A},?=J,\widetilde{J},\breve{J},?'=I\backslash J,\widetilde{I\backslash J},\breve{I\backslash J}.\\	
\end{align}	
\end{definition}

\subsection{$\infty$-Robba Rings over General Banach Affinoids}

\indent We now apply the construction of \cite{5BK} to the rings defined in the previous section. Recall from \cite{5BK}, for any Banach adic algebra $R$ over $\mathbb{Q}_p$ we have the derived spectrum $\mathrm{Spa}^h(R):=\mathrm{Spa}^h_{\mathrm{Rat}}(R)$.

\begin{definition}
Consider the following rings we defined in the previous section:
\begin{align}
%\Pi_{[s_I,r_I],I,J,I\backslash J,A},\\	
\Pi_{[s_I,r_I],I,\breve{J},I\backslash J,A},\\	
\Pi_{[s_I,r_I],I,\widetilde{J},I\backslash J,A},\\
\Pi_{[s_I,r_I],I,J,\breve{I\backslash J},A},\\	
\Pi_{[s_I,r_I],I,\breve{J},\breve{I\backslash J},A},\\	
\Pi_{[s_I,r_I],I,\widetilde{J},\breve{I\backslash J},A},\\
\Pi_{[s_I,r_I],I,J,\widetilde{I\backslash J},A},\\	
\Pi_{[s_I,r_I],I,\breve{J},\widetilde{I\backslash J},A},\\	
\Pi_{[s_I,r_I],I,\widetilde{J},\widetilde{I\backslash J},A}.	
\end{align}
We then take the corresponding derived spectrum from Bambozzi-Kremnizer to defined the following $\infty$-analytic stacks:
\begin{align}
%\Pi_{[s_I,r_I],I,J,I\backslash J,A},\\	
\mathrm{Spa}^h\Pi_{[s_I,r_I],I,\breve{J},I\backslash J,A},\\	
\mathrm{Spa}^h\Pi_{[s_I,r_I],I,\widetilde{J},I\backslash J,A},\\
\mathrm{Spa}^h\Pi_{[s_I,r_I],I,J,\breve{I\backslash J},A},\\	
\mathrm{Spa}^h\Pi_{[s_I,r_I],I,\breve{J},\breve{I\backslash J},A},\\	
\mathrm{Spa}^h\Pi_{[s_I,r_I],I,\widetilde{J},\breve{I\backslash J},A},\\
\mathrm{Spa}^h\Pi_{[s_I,r_I],I,J,\widetilde{I\backslash J},A},\\	
\mathrm{Spa}^h\Pi_{[s_I,r_I],I,\breve{J},\widetilde{I\backslash J},A},\\	
\mathrm{Spa}^h\Pi_{[s_I,r_I],I,\widetilde{J},\widetilde{I\backslash J},A}.	
\end{align}
Taking the global section we have the following ring spectra:
\begin{align}
%\Pi_{[s_I,r_I],I,J,I\backslash J,A},\\	
\Pi^h_{[s_I,r_I],I,\breve{J},I\backslash J,A},\\	
\Pi^h_{[s_I,r_I],I,\widetilde{J},I\backslash J,A},\\
\Pi^h_{[s_I,r_I],I,J,\breve{I\backslash J},A},\\	
\Pi^h_{[s_I,r_I],I,\breve{J},\breve{I\backslash J},A},\\	
\Pi^h_{[s_I,r_I],I,\widetilde{J},\breve{I\backslash J},A},\\
\Pi^h_{[s_I,r_I],I,J,\widetilde{I\backslash J},A},\\	
\Pi^h_{[s_I,r_I],I,\breve{J},\widetilde{I\backslash J},A},\\	
\Pi^h_{[s_I,r_I],I,\widetilde{J},\widetilde{I\backslash J},A}.
\end{align}
\end{definition}

\indent Then we have the following definitions:

\begin{definition} \mbox{\bf{(After Kedlaya-Liu, \cite[Definition 5.2.1]{5KL2})}}
\begin{align}
%\Pi_{\mathrm{an},r_I,I,J,I\backslash J,A},\\	
\Pi^h_{\mathrm{an},r_I,I,\breve{J},I\backslash J,A}:=\varprojlim_{s_I}\Pi^h_{[s_I,r_I],I,\breve{J},I\backslash J,A},\\	
\Pi^h_{\mathrm{an},r_I,I,\widetilde{J},I\backslash J,A}:=\varprojlim_{s_I} \Pi^h_{[s_I,r_I],I,\widetilde{J},I\backslash J,A},\\
\Pi^h_{\mathrm{an},r_I,I,J,\breve{I\backslash J},A}:=\varprojlim_{s_I}\Pi^h_{[s_I,r_I],I,J,\breve{I\backslash J},A},\\	
\Pi^h_{\mathrm{an},r_I,I,\breve{J},\breve{I\backslash J},A}:=\varprojlim_{s_I} \Pi^h_{[s_I,r_I],I,\breve{J},\breve{I\backslash J},A},\\	
\Pi^h_{\mathrm{an},r_I,I,\widetilde{J},\breve{I\backslash J},A}:=\varprojlim_{s_I} \Pi^h_{[s_I,r_I],I,\widetilde{J},\breve{I\backslash J},A},\\
\Pi^h_{\mathrm{an},r_I,I,J,\widetilde{I\backslash J},A}:=\varprojlim_{s_I} \Pi^h_{[s_I,r_I],I,J,\widetilde{I\backslash J},A},\\	
\Pi^h_{\mathrm{an},r_I,I,\breve{J},\widetilde{I\backslash J},A}:=\varprojlim_{s_I} \Pi^h_{[s_I,r_I],I,\breve{J},\widetilde{I\backslash J},A},\\	
\Pi^h_{\mathrm{an},r_I,I,\widetilde{J},\widetilde{I\backslash J},A}:=\varprojlim_{s_I} \Pi^h_{[s_I,r_I],I,\widetilde{J},\widetilde{I\backslash J},A}.	
\end{align}

and 

\begin{align}
%\Pi_{\mathrm{an},\mathrm{con},I,J,I\backslash J,A},\\	
\Pi^h_{\mathrm{an},\mathrm{con},I,\breve{J},I\backslash J,A}:=\varinjlim_{r_I}\varprojlim_{s_I}\Pi^h_{[s_I,r_I],I,\breve{J},I\backslash J,A},\\	
\Pi^h_{\mathrm{an},\mathrm{con},I,\widetilde{J},I\backslash J,A}:=\varinjlim_{r_I}\varprojlim_{s_I} \Pi^h_{[s_I,r_I],I,\widetilde{J},I\backslash J,A},\\
\Pi^h_{\mathrm{an},\mathrm{con},I,J,\breve{I\backslash J},A}:=\varinjlim_{r_I}\varprojlim_{s_I}\Pi^h_{[s_I,r_I],I,J,\breve{I\backslash J},A},\\	
\Pi^h_{\mathrm{an},\mathrm{con},I,\breve{J},\breve{I\backslash J},A}:=\varinjlim_{r_I}\varprojlim_{s_I} \Pi^h_{[s_I,r_I],I,\breve{J},\breve{I\backslash J},A},\\
\Pi^h_{\mathrm{an},\mathrm{con},I,\widetilde{J},\breve{I\backslash J},A}:=\varinjlim_{r_I}\varprojlim_{s_I} \Pi^h_{[s_I,r_I],I,\widetilde{J},\breve{I\backslash J},A},\\
\Pi^h_{\mathrm{an},\mathrm{con},I,J,\widetilde{I\backslash J},A}:=\varinjlim_{r_I}\varprojlim_{s_I} \Pi^h_{[s_I,r_I],I,J,\widetilde{I\backslash J},A},\\	
\Pi^h_{\mathrm{an},\mathrm{con},I,\breve{J},\widetilde{I\backslash J},A}:=\varinjlim_{r_I}\varprojlim_{s_I} \Pi^h_{[s_I,r_I],I,\breve{J},\widetilde{I\backslash J},A},\\	
\Pi^h_{\mathrm{an},\mathrm{con},I,\widetilde{J},\widetilde{I\backslash J},A}:=\varinjlim_{r_I}\varprojlim_{s_I} \Pi^h_{[s_I,r_I],I,\widetilde{J},\widetilde{I\backslash J},A}.\\	
\end{align}	
\end{definition}

\subsection{Fr\'echet Objects in the Sheafy Situation}

\begin{assumption}
Assume the followings are sheafy adic Banach uniform algebra over $\mathbb{Q}_p$:
\begin{align}
%\Pi_{[s_I,r_I],I,J,I\backslash J,A},\\	
\Pi_{[s_I,r_I],I,?,?',A},?=J,\widetilde{J},\breve{J},?'=I\backslash J,\widetilde{I\backslash J},\breve{I\backslash J}.	
\end{align}	
\end{assumption}

\indent In what follow, we consider the corresponding radii living in the set of all the rational numbers.

\begin{definition} \mbox{\bf{(After KPX, \cite[Definition 2.1.3]{5KPX})}}
Over $\Pi_{\mathrm{an},r_{I,0},I,\breve{J},I\backslash J,A}$ we define the corresponding stably pseudocoherent sheaves to mean a collection of stably pseudocoherent modules $(M_{[s_I,r_I]})$	over each $\Pi_{[s_I,r_I],I,\breve{J},I\backslash J,A}$ satisfying the corresponding compatibility condition and the obvious cocycle condition with respect to the family of the corresponding multi-intervals $\{[s_I,r_I]\}$.
\end{definition}

\begin{definition} \mbox{\bf{(After KPX, \cite[Definition 2.1.3]{5KPX})}}
Over $\Pi_{\mathrm{an},r_{I,0},I,\widetilde{J},I\backslash J,A}$ we define the corresponding stably pseudocoherent sheaves to mean a collection of stably pseudocoherent modules $(M_{[s_I,r_I]})$	over each $\Pi_{[s_I,r_I],I,\widetilde{J},I\backslash J,A}$ satisfying the corresponding compatibility condition and the obvious cocycle condition with respect to the family of the corresponding multi-intervals $\{[s_I,r_I]\}$.
\end{definition}

\begin{definition} \mbox{\bf{(After KPX, \cite[Definition 2.1.3]{5KPX})}}
Over $\Pi_{\mathrm{an},r_{I,0},I,{J},\breve{I\backslash J},A}$ we define the corresponding stably pseudocoherent sheaves to mean a collection of stably pseudocoherent modules $(M_{[s_I,r_I]})$	over each $\Pi_{[s_I,r_I],I,{J},\breve{I\backslash J},A}$ satisfying the corresponding compatibility condition and the obvious cocycle condition with respect to the family of the corresponding multi-intervals $\{[s_I,r_I]\}$.
\end{definition}

\begin{definition} \mbox{\bf{(After KPX, \cite[Definition 2.1.3]{5KPX})}}
Over $\Pi_{\mathrm{an},r_{I,0},I,\breve{J},\breve{I\backslash J},A}$ we define the corresponding stably pseudocoherent sheaves to mean a collection of stably pseudocoherent modules $(M_{[s_I,r_I]})$	over each $\Pi_{[s_I,r_I],I,\breve{J},\breve{I\backslash J},A}$ satisfying the corresponding compatibility condition and the obvious cocycle condition with respect to the family of the corresponding multi-intervals $\{[s_I,r_I]\}$.
\end{definition}

\begin{definition} \mbox{\bf{(After KPX, \cite[Definition 2.1.3]{5KPX})}}
Over $\Pi_{\mathrm{an},r_{I,0},I,\widetilde{J},\breve{I\backslash J},A}$ we define the corresponding stably pseudocoherent sheaves to mean a collection of stably pseudocoherent modules $(M_{[s_I,r_I]})$ over each $\Pi_{[s_I,r_I],I,\widetilde{J},\breve{I\backslash J},A}$ satisfying the corresponding compatibility condition and the obvious cocycle condition with respect to the family of the corresponding multi-intervals $\{[s_I,r_I]\}$.
\end{definition}

\begin{definition} \mbox{\bf{(After KPX, \cite[Definition 2.1.3]{5KPX})}}
Over $\Pi_{\mathrm{an},r_{I,0},I,{J},\widetilde{I\backslash J},A}$ we define the corresponding stably pseudocoherent sheaves to mean a collection of stably pseudocoherent modules $(M_{[s_I,r_I]})$	over each $\Pi_{[s_I,r_I],I,{J},\widetilde{I\backslash J},A}$ satisfying the corresponding compatibility condition and the obvious cocycle condition with respect to the family of the corresponding multi-intervals $\{[s_I,r_I]\}$.
\end{definition}

\begin{definition} \mbox{\bf{(After KPX, \cite[Definition 2.1.3]{5KPX})}}
Over $\Pi_{\mathrm{an},r_{I,0},I,\breve{J},\widetilde{I\backslash J},A}$ we define the corresponding stably pseudocoherent sheaves to mean a collection of stably pseudocoherent modules $(M_{[s_I,r_I]})$ over each $\Pi_{[s_I,r_I],I,\breve{J},\widetilde{I\backslash J},A}$ satisfying the corresponding compatibility condition and the obvious cocycle condition with respect to the family of the corresponding multi-intervals $\{[s_I,r_I]\}$.
\end{definition}

\begin{definition} \mbox{\bf{(After KPX, \cite[Definition 2.1.3]{5KPX})}}
Over $\Pi_{\mathrm{an},r_{I,0},I,\widetilde{J},\widetilde{I\backslash J},A}$ we define the corresponding stably pseudocoherent sheaves to mean a collection of stably pseudocoherent modules $(M_{[s_I,r_I]})$ over each $\Pi_{[s_I,r_I],I,\widetilde{J},\widetilde{I\backslash J},A}$ satisfying the corresponding compatibility condition and the obvious cocycle condition with respect to the family of the corresponding multi-intervals $\{[s_I,r_I]\}$.
\end{definition}

\begin{remark}
There is some overlap and repeating on some objects within the eight categories. Therefore in the following we are going to then work with only 1st, 2nd, 4th, 5th, 8th categories. 	
\end{remark}

\begin{proposition} \mbox{\bf{(After Kedlaya-Liu, \cite[Corollary 2.6.8]{5KL2})}}
The corresponding pseudocoherent finitely projective bundles over
\begin{align}
\Pi_{\mathrm{an},r_{I,0},I,\breve{J},I\backslash J,A},\\
\Pi_{\mathrm{an},r_{I,0},I,\widetilde{J},I\backslash J,A},\\
\Pi_{\mathrm{an},r_{I,0},I,\breve{J},\breve{I\backslash J},A},\\
\Pi_{\mathrm{an},r_{I,0},I,\widetilde{J},\breve{I\backslash J},A},\\
\Pi_{\mathrm{an},r_{I,0},I,\widetilde{J},\widetilde{I\backslash J},A}
\end{align}
defined as above have global sections which are finite projective if and only if the global sections are finitely generated. 	
\end{proposition}

\begin{proof}
Just apply \cite[Corollary 2.6.8]{5KL2}.	
\end{proof}

\begin{proposition} \mbox{\bf{(After Kedlaya-Liu, \cite[Proposition 2.6.17]{5KL2})}}
The corresponding pseudocoherent sheaves over
\begin{align}
\Pi_{\mathrm{an},r_{I,0},I,\breve{J},I\backslash J,A},\\
\Pi_{\mathrm{an},r_{I,0},I,\widetilde{J},I\backslash J,A},\\
\Pi_{\mathrm{an},r_{I,0},I,\breve{J},\breve{I\backslash J},A},\\
\Pi_{\mathrm{an},r_{I,0},I,\widetilde{J},\breve{I\backslash J},A},\\
\Pi_{\mathrm{an},r_{I,0},I,\widetilde{J},\widetilde{I\backslash J},A}
\end{align}
defined as above have global sections which are finitely generated as long as we have the uniform bound on the rank of the bundles over quasi-compacts with respect to closed multi-intervals taking the general form of $[s_I,r_I]$. 	
\end{proposition}

\begin{proof}
This is actually a direct consequence of \cite[Proposition 2.6.17]{5KL2}, where the space $(0,r_{I,0}]$ admits $2^{|I|}$-uniform covering.	
\end{proof}

\begin{proposition}\mbox{}\\ \mbox{\bf{(After Kedlaya-Liu, \cite[Corollary 2.6.8, Proposition 2.6.17]{5KL2})}}
The corresponding pseudocoherent sheaves over
\begin{align}
\Pi_{\mathrm{an},r_{I,0},I,\breve{J},I\backslash J,A},\\
\Pi_{\mathrm{an},r_{I,0},I,\widetilde{J},I\backslash J,A},\\
\Pi_{\mathrm{an},r_{I,0},I,\breve{J},\breve{I\backslash J},A},\\
\Pi_{\mathrm{an},r_{I,0},I,\widetilde{J},\breve{I\backslash J},A},\\
\Pi_{\mathrm{an},r_{I,0},I,\widetilde{J},\widetilde{I\backslash J},A}
\end{align}
defined as above have global sections which are finite projective as long as we have the uniform bound on the rank of the bundle over each quasi-compact with respect to each closed multi-interval $[s_I,r_I]$, and we have that the sheaves admits section actually finite projective over each quasi-compact with respect to each closed multi-interval $[s_I,r_I]$. 	
\end{proposition}

\begin{proof}
This is actually a direct corollary of the previous two propositions. 	
\end{proof}

%%\newpage

\newpage\section{Cyclotomic Multivariate $(\varphi_I,\Gamma_I)$-Modules over Rigid Analytic Affinoids in Mixed-characteristic Case}

\subsection{Fundamental Definitions}

\indent In the situation where $A$ is a rigid analytic affinoid over $\mathbb{Q}_p$. Recall from \cite{5T1}, suppose we have $|I|$ finite extensions of $\mathbb{Q}_p$, where we denote them as $K_1,...,K_I$. Then we have the corresponding uniformizers $\pi_{K_1},...,\pi_{K_I}$, the corresponding Frobenius operators $\varphi_1,...,\varphi_I$ and the corresponding groups $\Gamma_{K_1},...,\Gamma_{K_I}$. Recall from \cite{5T1}, by adding the corresponding variables from $\pi_{K_1},...,\pi_{K_I}$ and $\Gamma_{K_1},...,\Gamma_{K_I}$ we have the following rings:

\begin{align}
\Pi_{[s_I,r_I],I,I,\emptyset,A}(\pi_{K_I}):=\Pi_{[s_I,r_I],I,A}(\pi_{K_I})\\
\Pi_{\mathrm{an},r_I,I,I,\emptyset,A}(\pi_{K_I}):=\Pi_{\mathrm{an},r_I,I,A}(\pi_{K_I})\\
\Pi_{\mathrm{an,con},I,I,\emptyset,A}(\pi_{K_I}):=\Pi_{\mathrm{an,con},I,\emptyset,A}(\pi_{K_I})	
\end{align}
and
\begin{align}
\Pi_{[s_I,r_I],I,I,\emptyset,A}(\Gamma_{K_I}):=\Pi_{[s_I,r_I],I,A}(\Gamma_{K_I})\\
\Pi_{\mathrm{an},r_I,I,I,\emptyset,A}(\Gamma_{K_I}):=\Pi_{\mathrm{an},r_I,I,A}(\Gamma_{K_I})\\
\Pi_{\mathrm{an,con},I,I,\emptyset,A}(\Gamma_{K_I}):=\Pi_{\mathrm{an,con},I,\emptyset,A}(\Gamma_{K_I}).	
\end{align}

\begin{definition}
By taking the direct base change we have the following rings in mixed-characteristic situation:
\begin{align}
%\Pi_{[s_I,r_I],I,J,I\backslash J,A},\\	
\Pi_{[s_I,r_I],I,\breve{J},I\backslash J,A}(\pi_{K_I}),\\	
\Pi_{[s_I,r_I],I,\widetilde{J},I\backslash J,A}(\pi_{K_I}),\\
\Pi_{[s_I,r_I],I,J,\breve{I\backslash J},A}(\pi_{K_I}),\\	
\Pi_{[s_I,r_I],I,\breve{J},\breve{I\backslash J},A}(\pi_{K_I}),\\	
\Pi_{[s_I,r_I],I,\widetilde{J},\breve{I\backslash J},A}(\pi_{K_I}),\\
\Pi_{[s_I,r_I],I,J,\widetilde{I\backslash J},A}(\pi_{K_I}),\\	
\Pi_{[s_I,r_I],I,\breve{J},\widetilde{I\backslash J},A}(\pi_{K_I}),\\	
\Pi_{[s_I,r_I],I,\widetilde{J},\widetilde{I\backslash J},A}(\pi_{K_I}).	
\end{align}
By taking the direct base change we have the following rings in mixed-characteristic situation:
\begin{align}
%\Pi_{[s_I,r_I],I,J,I\backslash J,A},\\	
\Pi_{[s_I,r_I],I,\breve{J},I\backslash J,A}(\Gamma_{K_I}),\\	
\Pi_{[s_I,r_I],I,\widetilde{J},I\backslash J,A}(\Gamma_{K_I}),\\
\Pi_{[s_I,r_I],I,J,\breve{I\backslash J},A}(\Gamma_{K_I}),\\	
\Pi_{[s_I,r_I],I,\breve{J},\breve{I\backslash J},A}(\Gamma_{K_I}),\\
\Pi_{[s_I,r_I],I,\widetilde{J},\breve{I\backslash J},A}(\Gamma_{K_I}),\\
\Pi_{[s_I,r_I],I,J,\widetilde{I\backslash J},A}(\Gamma_{K_I}),\\	
\Pi_{[s_I,r_I],I,\breve{J},\widetilde{I\backslash J},A}(\Gamma_{K_I}),\\	
\Pi_{[s_I,r_I],I,\widetilde{J},\widetilde{I\backslash J},A}(\Gamma_{K_I}).	
\end{align}	
\end{definition}

\begin{definition} \mbox{\bf{(After Kedlaya-Liu, \cite[Definition 5.2.1]{5KL2})}}
\begin{align}
%\Pi_{\mathrm{an},r_I,I,J,I\backslash J,A},\\	
\Pi_{\mathrm{an},r_I,I,\breve{J},I\backslash J,A}(\pi_{K_I}):=\varprojlim_{s_I}\Pi_{[s_I,r_I],I,\breve{J},I\backslash J,A}(\pi_{K_I}),\\	
\Pi_{\mathrm{an},r_I,I,\widetilde{J},I\backslash J,A}(\pi_{K_I}):=\varprojlim_{s_I} \Pi_{[s_I,r_I],I,\widetilde{J},I\backslash J,A}(\pi_{K_I}),\\
\Pi_{\mathrm{an},r_I,I,J,\breve{I\backslash J},A}(\pi_{K_I}):=\varprojlim_{s_I}\Pi_{[s_I,r_I],I,J,\breve{I\backslash J},A}(\pi_{K_I}),\\	
\Pi_{\mathrm{an},r_I,I,\breve{J},\breve{I\backslash J},A}(\pi_{K_I}):=\varprojlim_{s_I} \Pi_{[s_I,r_I],I,\breve{J},\breve{I\backslash J},A}(\pi_{K_I}),\\	
\Pi_{\mathrm{an},r_I,I,\widetilde{J},\breve{I\backslash J},A}(\pi_{K_I}):=\varprojlim_{s_I} \Pi_{[s_I,r_I],I,\widetilde{J},\breve{I\backslash J},A}(\pi_{K_I}),\\
\Pi_{\mathrm{an},r_I,I,J,\widetilde{I\backslash J},A}(\pi_{K_I}):=\varprojlim_{s_I} \Pi_{[s_I,r_I],I,J,\widetilde{I\backslash J},A}(\pi_{K_I}),\\	
\Pi_{\mathrm{an},r_I,I,\breve{J},\widetilde{I\backslash J},A}(\pi_{K_I}):=\varprojlim_{s_I} \Pi_{[s_I,r_I],I,\breve{J},\widetilde{I\backslash J},A}(\pi_{K_I}),\\	
\Pi_{\mathrm{an},r_I,I,\widetilde{J},\widetilde{I\backslash J},A}(\pi_{K_I}):=\varprojlim_{s_I} \Pi_{[s_I,r_I],I,\widetilde{J},\widetilde{I\backslash J},A}(\pi_{K_I}).	
\end{align}

and 

\begin{align}
%\Pi_{\mathrm{an},\mathrm{con},I,J,I\backslash J,A},\\	
\Pi_{\mathrm{an},\mathrm{con},I,\breve{J},I\backslash J,A}(\pi_{K_I}):=\varinjlim_{r_I}\varprojlim_{s_I}\Pi_{[s_I,r_I],I,\breve{J},I\backslash J,A}(\pi_{K_I}),\\	
\Pi_{\mathrm{an},\mathrm{con},I,\widetilde{J},I\backslash J,A}(\pi_{K_I}):=\varinjlim_{r_I}\varprojlim_{s_I} \Pi_{[s_I,r_I],I,\widetilde{J},I\backslash J,A}(\pi_{K_I}),\\
\Pi_{\mathrm{an},\mathrm{con},I,J,\breve{I\backslash J},A}(\pi_{K_I}):=\varinjlim_{r_I}\varprojlim_{s_I}\Pi_{[s_I,r_I],I,J,\breve{I\backslash J},A}(\pi_{K_I}),\\	
\Pi_{\mathrm{an},\mathrm{con},I,\breve{J},\breve{I\backslash J},A}(\pi_{K_I}):=\varinjlim_{r_I}\varprojlim_{s_I} \Pi_{[s_I,r_I],I,\breve{J},\breve{I\backslash J},A}(\pi_{K_I}),\\
\Pi_{\mathrm{an},\mathrm{con},I,\widetilde{J},\breve{I\backslash J},A}(\pi_{K_I}):=\varinjlim_{r_I}\varprojlim_{s_I} \Pi_{[s_I,r_I],I,\widetilde{J},\breve{I\backslash J},A}(\pi_{K_I}),\\
\Pi_{\mathrm{an},\mathrm{con},I,J,\widetilde{I\backslash J},A}(\pi_{K_I}):=\varinjlim_{r_I}\varprojlim_{s_I} \Pi_{[s_I,r_I],I,J,\widetilde{I\backslash J},A}(\pi_{K_I}),\\	
\Pi_{\mathrm{an},\mathrm{con},I,\breve{J},\widetilde{I\backslash J},A}(\pi_{K_I}):=\varinjlim_{r_I}\varprojlim_{s_I} \Pi_{[s_I,r_I],I,\breve{J},\widetilde{I\backslash J},A}(\pi_{K_I}),\\	
\Pi_{\mathrm{an},\mathrm{con},I,\widetilde{J},\widetilde{I\backslash J},A}(\pi_{K_I}):=\varinjlim_{r_I}\varprojlim_{s_I} \Pi_{[s_I,r_I],I,\widetilde{J},\widetilde{I\backslash J},A}(\pi_{K_I}).\\	
\end{align}	
\end{definition}

\begin{definition} \mbox{\bf{(After Kedlaya-Liu, \cite[Definition 5.2.1]{5KL2})}}
\begin{align}
%\Pi_{\mathrm{an},r_I,I,J,I\backslash J,A},\\	
\Pi_{\mathrm{an},r_I,I,\breve{J},I\backslash J,A}(\Gamma_{K_I}):=\varprojlim_{s_I}\Pi_{[s_I,r_I],I,\breve{J},I\backslash J,A}(\Gamma_{K_I}),\\	
\Pi_{\mathrm{an},r_I,I,\widetilde{J},I\backslash J,A}(\Gamma_{K_I}):=\varprojlim_{s_I} \Pi_{[s_I,r_I],I,\widetilde{J},I\backslash J,A}(\Gamma_{K_I}),\\
\Pi_{\mathrm{an},r_I,I,J,\breve{I\backslash J},A}(\Gamma_{K_I}):=\varprojlim_{s_I}\Pi_{[s_I,r_I],I,J,\breve{I\backslash J},A}(\Gamma_{K_I}),\\	
\Pi_{\mathrm{an},r_I,I,\breve{J},\breve{I\backslash J},A}(\Gamma_{K_I}):=\varprojlim_{s_I} \Pi_{[s_I,r_I],I,\breve{J},\breve{I\backslash J},A}(\Gamma_{K_I}),\\	
\Pi_{\mathrm{an},r_I,I,\widetilde{J},\breve{I\backslash J},A}(\Gamma_{K_I}):=\varprojlim_{s_I} \Pi_{[s_I,r_I],I,\widetilde{J},\breve{I\backslash J},A}(\Gamma_{K_I}),\\
\Pi_{\mathrm{an},r_I,I,J,\widetilde{I\backslash J},A}(\Gamma_{K_I}):=\varprojlim_{s_I} \Pi_{[s_I,r_I],I,J,\widetilde{I\backslash J},A}(\Gamma_{K_I}),\\	
\Pi_{\mathrm{an},r_I,I,\breve{J},\widetilde{I\backslash J},A}(\Gamma_{K_I}):=\varprojlim_{s_I} \Pi_{[s_I,r_I],I,\breve{J},\widetilde{I\backslash J},A}(\Gamma_{K_I}),\\	
\Pi_{\mathrm{an},r_I,I,\widetilde{J},\widetilde{I\backslash J},A}(\Gamma_{K_I}):=\varprojlim_{s_I} \Pi_{[s_I,r_I],I,\widetilde{J},\widetilde{I\backslash J},A}(\Gamma_{K_I}).	
\end{align}

and 

\begin{align}
%\Pi_{\mathrm{an},\mathrm{con},I,J,I\backslash J,A},\\	
\Pi_{\mathrm{an},\mathrm{con},I,\breve{J},I\backslash J,A}(\Gamma_{K_I}):=\varinjlim_{r_I}\varprojlim_{s_I}\Pi_{[s_I,r_I],I,\breve{J},I\backslash J,A}(\Gamma_{K_I}),\\	
\Pi_{\mathrm{an},\mathrm{con},I,\widetilde{J},I\backslash J,A}(\Gamma_{K_I}):=\varinjlim_{r_I}\varprojlim_{s_I} \Pi_{[s_I,r_I],I,\widetilde{J},I\backslash J,A}(\Gamma_{K_I}),\\
\Pi_{\mathrm{an},\mathrm{con},I,J,\breve{I\backslash J},A}(\Gamma_{K_I}):=\varinjlim_{r_I}\varprojlim_{s_I}\Pi_{[s_I,r_I],I,J,\breve{I\backslash J},A}(\Gamma_{K_I}),\\	
\Pi_{\mathrm{an},\mathrm{con},I,\breve{J},\breve{I\backslash J},A}(\Gamma_{K_I}):=\varinjlim_{r_I}\varprojlim_{s_I} \Pi_{[s_I,r_I],I,\breve{J},\breve{I\backslash J},A}(\Gamma_{K_I}),\\
\Pi_{\mathrm{an},\mathrm{con},I,\widetilde{J},\breve{I\backslash J},A}(\Gamma_{K_I}):=\varinjlim_{r_I}\varprojlim_{s_I} \Pi_{[s_I,r_I],I,\widetilde{J},\breve{I\backslash J},A}(\Gamma_{K_I}),\\
\Pi_{\mathrm{an},\mathrm{con},I,J,\widetilde{I\backslash J},A}(\Gamma_{K_I}):=\varinjlim_{r_I}\varprojlim_{s_I} \Pi_{[s_I,r_I],I,J,\widetilde{I\backslash J},A}(\Gamma_{K_I}),\\	
\Pi_{\mathrm{an},\mathrm{con},I,\breve{J},\widetilde{I\backslash J},A}(\Gamma_{K_I}):=\varinjlim_{r_I}\varprojlim_{s_I} \Pi_{[s_I,r_I],I,\breve{J},\widetilde{I\backslash J},A}(\Gamma_{K_I}),\\	
\Pi_{\mathrm{an},\mathrm{con},I,\widetilde{J},\widetilde{I\backslash J},A}(\Gamma_{K_I}):=\varinjlim_{r_I}\varprojlim_{s_I} \Pi_{[s_I,r_I],I,\widetilde{J},\widetilde{I\backslash J},A}(\Gamma_{K_I}).	
\end{align}	
\end{definition}

\indent Now we consider the corresponding definition of the $(\varphi_I,\Gamma_I)$-modules over the corresponding period rings  defined above.

\begin{setting}
For the corresponding $\varphi_I$-modules over the rings:
\begin{align}
%\Pi_{\mathrm{an},r_I,I,J,I\backslash J,A},\\	
\Pi_{\mathrm{an},r_I,I,\breve{J},I\backslash J,A}(\pi_{K_I}):=\varprojlim_{s_I}\Pi_{[s_I,r_I],I,\breve{J},I\backslash J,A}(\pi_{K_I}),\\	
\Pi_{\mathrm{an},r_I,I,\widetilde{J},I\backslash J,A}(\pi_{K_I}):=\varprojlim_{s_I} \Pi_{[s_I,r_I],I,\widetilde{J},I\backslash J,A}(\pi_{K_I}),\\
\Pi_{\mathrm{an},r_I,I,J,\breve{I\backslash J},A}(\pi_{K_I}):=\varprojlim_{s_I}\Pi_{[s_I,r_I],I,J,\breve{I\backslash J},A}(\pi_{K_I}),\\	
\Pi_{\mathrm{an},r_I,I,\breve{J},\breve{I\backslash J},A}(\pi_{K_I}):=\varprojlim_{s_I} \Pi_{[s_I,r_I],I,\breve{J},\breve{I\backslash J},A}(\pi_{K_I}),\\	
\Pi_{\mathrm{an},r_I,I,\widetilde{J},\breve{I\backslash J},A}(\pi_{K_I}):=\varprojlim_{s_I} \Pi_{[s_I,r_I],I,\widetilde{J},\breve{I\backslash J},A}(\pi_{K_I}),\\
\Pi_{\mathrm{an},r_I,I,J,\widetilde{I\backslash J},A}(\pi_{K_I}):=\varprojlim_{s_I} \Pi_{[s_I,r_I],I,J,\widetilde{I\backslash J},A}(\pi_{K_I}),\\	
\Pi_{\mathrm{an},r_I,I,\breve{J},\widetilde{I\backslash J},A}(\pi_{K_I}):=\varprojlim_{s_I} \Pi_{[s_I,r_I],I,\breve{J},\widetilde{I\backslash J},A}(\pi_{K_I}),\\	
\Pi_{\mathrm{an},r_I,I,\widetilde{J},\widetilde{I\backslash J},A}(\pi_{K_I}):=\varprojlim_{s_I} \Pi_{[s_I,r_I],I,\widetilde{J},\widetilde{I\backslash J},A}(\pi_{K_I})	
\end{align}	
we assume we have the sufficiently small radius as in \cite[Definition 2.2.6]{5KPX}. We will keep this assumption in all similar situation involving the corresponding $\varphi_I$-modules.
\end{setting}

\begin{definition} \mbox{\bf{(After KPX \cite[Definition 2.2.6]{5KPX})}}
We define in the following way the corresponding $\varphi_I$-modules over the following period rings:
\begin{align}
%\Pi_{\mathrm{an},\mathrm{con},I,J,I\backslash J,A},\\	
\Pi_{\mathrm{an},\mathrm{con},I,\breve{J},I\backslash J,A}(\pi_{K_I}):=\varinjlim_{r_I}\varprojlim_{s_I}\Pi_{[s_I,r_I],I,\breve{J},I\backslash J,A}(\pi_{K_I}),\\	
\Pi_{\mathrm{an},\mathrm{con},I,\widetilde{J},I\backslash J,A}(\pi_{K_I}):=\varinjlim_{r_I}\varprojlim_{s_I} \Pi_{[s_I,r_I],I,\widetilde{J},I\backslash J,A}(\pi_{K_I}),\\
\Pi_{\mathrm{an},\mathrm{con},I,J,\breve{I\backslash J},A}(\pi_{K_I}):=\varinjlim_{r_I}\varprojlim_{s_I}\Pi_{[s_I,r_I],I,J,\breve{I\backslash J},A}(\pi_{K_I}),\\	
\Pi_{\mathrm{an},\mathrm{con},I,\breve{J},\breve{I\backslash J},A}(\pi_{K_I}):=\varinjlim_{r_I}\varprojlim_{s_I} \Pi_{[s_I,r_I],I,\breve{J},\breve{I\backslash J},A}(\pi_{K_I}),\\
\Pi_{\mathrm{an},\mathrm{con},I,\widetilde{J},\breve{I\backslash J},A}(\pi_{K_I}):=\varinjlim_{r_I}\varprojlim_{s_I} \Pi_{[s_I,r_I],I,\widetilde{J},\breve{I\backslash J},A}(\pi_{K_I}),\\
\Pi_{\mathrm{an},\mathrm{con},I,J,\widetilde{I\backslash J},A}(\pi_{K_I}):=\varinjlim_{r_I}\varprojlim_{s_I} \Pi_{[s_I,r_I],I,J,\widetilde{I\backslash J},A}(\pi_{K_I}),\\	
\Pi_{\mathrm{an},\mathrm{con},I,\breve{J},\widetilde{I\backslash J},A}(\pi_{K_I}):=\varinjlim_{r_I}\varprojlim_{s_I} \Pi_{[s_I,r_I],I,\breve{J},\widetilde{I\backslash J},A}(\pi_{K_I}),\\	
\Pi_{\mathrm{an},\mathrm{con},I,\widetilde{J},\widetilde{I\backslash J},A}(\pi_{K_I}):=\varinjlim_{r_I}\varprojlim_{s_I} \Pi_{[s_I,r_I],I,\widetilde{J},\widetilde{I\backslash J},A}(\pi_{K_I}).	
\end{align}	
These are the corresponding base change of the corresponding $\varphi_I$-modules coming from the the ones over the following rings:
\begin{align}
%\Pi_{\mathrm{an},r_I,I,J,I\backslash J,A},\\	
\Pi_{\mathrm{an},r_I,I,\breve{J},I\backslash J,A}(\pi_{K_I}):=\varprojlim_{s_I}\Pi_{[s_I,r_I],I,\breve{J},I\backslash J,A}(\pi_{K_I}),\\	
\Pi_{\mathrm{an},r_I,I,\widetilde{J},I\backslash J,A}(\pi_{K_I}):=\varprojlim_{s_I} \Pi_{[s_I,r_I],I,\widetilde{J},I\backslash J,A}(\pi_{K_I}),\\
\Pi_{\mathrm{an},r_I,I,J,\breve{I\backslash J},A}(\pi_{K_I}):=\varprojlim_{s_I}\Pi_{[s_I,r_I],I,J,\breve{I\backslash J},A}(\pi_{K_I}),\\	
\Pi_{\mathrm{an},r_I,I,\breve{J},\breve{I\backslash J},A}(\pi_{K_I}):=\varprojlim_{s_I} \Pi_{[s_I,r_I],I,\breve{J},\breve{I\backslash J},A}(\pi_{K_I}),\\	
\Pi_{\mathrm{an},r_I,I,\widetilde{J},\breve{I\backslash J},A}(\pi_{K_I}):=\varprojlim_{s_I} \Pi_{[s_I,r_I],I,\widetilde{J},\breve{I\backslash J},A}(\pi_{K_I}),\\
\Pi_{\mathrm{an},r_I,I,J,\widetilde{I\backslash J},A}(\pi_{K_I}):=\varprojlim_{s_I} \Pi_{[s_I,r_I],I,J,\widetilde{I\backslash J},A}(\pi_{K_I}),\\	
\Pi_{\mathrm{an},r_I,I,\breve{J},\widetilde{I\backslash J},A}(\pi_{K_I}):=\varprojlim_{s_I} \Pi_{[s_I,r_I],I,\breve{J},\widetilde{I\backslash J},A}(\pi_{K_I}),\\	
\Pi_{\mathrm{an},r_I,I,\widetilde{J},\widetilde{I\backslash J},A}(\pi_{K_I}):=\varprojlim_{s_I} \Pi_{[s_I,r_I],I,\widetilde{J},\widetilde{I\backslash J},A}(\pi_{K_I}).	
\end{align}
For this latter group of period rings, we define a corresponding pseudocoherent or finite projective $\varphi_I$-module to be a corresponding pseudocoherent or finite projective module $M$ over this latter group of period rings carrying the corresponding semilinear partial Frobenius action coming from each Frobenius operator $\varphi_\alpha,\alpha\in I$ such that for each $\alpha$ we have 
\begin{align}
\varphi_\alpha^*M\otimes_{\Pi_{\mathrm{an},\{...,r_\alpha/p,...\},I,*,*,A}(\pi_{K_I})} &\Pi_{\mathrm{an},\{...,r_\alpha/p,...\},I,*,*,A}(\pi_{K_I})\\
&\overset{\sim}{\longrightarrow} M\otimes_{\Pi_{\mathrm{an},\{...,r_\alpha,...\},I,*,*,A}(\pi_{K_I})}\Pi_{\mathrm{an},\{...,r_\alpha/p,...\},I,*,*,A}(\pi_{K_I}).	
\end{align}
And we assume that altogether the partial Frobenius operators are commuting with each other. We assume all the modules involved are complete for the natural topology (mainly in the pseudocoherent situation) whose base changes to the following rings:
\begin{align}
%\Pi_{[s_I,r_I],I,J,I\backslash J,A},\\	
\Pi_{[s_I,r_I],I,\breve{J},I\backslash J,A}(\pi_{K_I}),\\	
\Pi_{[s_I,r_I],I,\widetilde{J},I\backslash J,A}(\pi_{K_I}),\\
\Pi_{[s_I,r_I],I,J,\breve{I\backslash J},A}(\pi_{K_I}),\\	
\Pi_{[s_I,r_I],I,\breve{J},\breve{I\backslash J},A}(\pi_{K_I}),\\
\Pi_{[s_I,r_I],I,\widetilde{J},\breve{I\backslash J},A}(\pi_{K_I}),\\
\Pi_{[s_I,r_I],I,J,\widetilde{I\backslash J},A}(\pi_{K_I}),\\	
\Pi_{[s_I,r_I],I,\breve{J},\widetilde{I\backslash J},A}(\pi_{K_I}),\\
\Pi_{[s_I,r_I],I,\widetilde{J},\widetilde{I\backslash J},A}(\pi_{K_I})
\end{align}	
give rise to the corresponding $\varphi_I$-modules defined over these rings which defined in the following way.	
\indent We define a corresponding pseudocoherent or finite projective $\varphi_I$-module over:
\begin{align}
%\Pi_{[s_I,r_I],I,J,I\backslash J,A},\\	
\Pi_{[s_I,r_I],I,\breve{J},I\backslash J,A}(\pi_{K_I}),\\	
\Pi_{[s_I,r_I],I,\widetilde{J},I\backslash J,A}(\pi_{K_I}),\\
\Pi_{[s_I,r_I],I,J,\breve{I\backslash J},A}(\pi_{K_I}),\\	
\Pi_{[s_I,r_I],I,\breve{J},\breve{I\backslash J},A}(\pi_{K_I}),\\
\Pi_{[s_I,r_I],I,\widetilde{J},\breve{I\backslash J},A}(\pi_{K_I}),\\
\Pi_{[s_I,r_I],I,J,\widetilde{I\backslash J},A}(\pi_{K_I}),\\	
\Pi_{[s_I,r_I],I,\breve{J},\widetilde{I\backslash J},A}(\pi_{K_I}),\\
\Pi_{[s_I,r_I],I,\widetilde{J},\widetilde{I\backslash J},A}(\pi_{K_I})	
\end{align}
to be a corresponding stably-pseudocoherent or finite projective module $M$ over this latter group of period rings carrying the corresponding semilinear partial Frobenius action coming from each Frobenius operator $\varphi_\alpha,\alpha\in I$ such that for each $\alpha$ we have 
\begin{align}
\varphi_\alpha^*M\otimes_{\Pi_{\mathrm{an},...,[s_\alpha/p,r_\alpha/p],...,I,*,*,A}(\pi_{K_I})} &\Pi_{\mathrm{an},...,[s_\alpha,r_\alpha/p],...,I,*,*,A}(\pi_{K_I})\\
&\overset{\sim}{\longrightarrow} M\otimes_{\Pi_{\mathrm{an},...,[s_\alpha,r_\alpha],...,I,*,*,A}(\pi_{K_I})}\Pi_{\mathrm{an},...,[s_\alpha,r_\alpha/p],...,I,*,*,A}(\pi_{K_I}).	
\end{align}
And we assume that altogether the partial Frobenius operators are commuting with each other. We assume all the modules involved are complete for the natural topology (mainly in the pseudocoherent situation).\\
\indent Finally we define the pseudocoherent or finite projective $\varphi_I$-modules over the corresponding period rings:
\begin{align}
%\Pi_{\mathrm{an},\mathrm{con},I,J,I\backslash J,A},\\	
\Pi_{\mathrm{an},\mathrm{con},I,\breve{J},I\backslash J,A}(\pi_{K_I}):=\varinjlim_{r_I}\varprojlim_{s_I}\Pi_{[s_I,r_I],I,\breve{J},I\backslash J,A}(\pi_{K_I}),\\	
\Pi_{\mathrm{an},\mathrm{con},I,\widetilde{J},I\backslash J,A}(\pi_{K_I}):=\varinjlim_{r_I}\varprojlim_{s_I} \Pi_{[s_I,r_I],I,\widetilde{J},I\backslash J,A}(\pi_{K_I}),\\
\Pi_{\mathrm{an},\mathrm{con},I,J,\breve{I\backslash J},A}(\pi_{K_I}):=\varinjlim_{r_I}\varprojlim_{s_I}\Pi_{[s_I,r_I],I,J,\breve{I\backslash J},A}(\pi_{K_I}),\\	
\Pi_{\mathrm{an},\mathrm{con},I,\breve{J},\breve{I\backslash J},A}(\pi_{K_I}):=\varinjlim_{r_I}\varprojlim_{s_I} \Pi_{[s_I,r_I],I,\breve{J},\breve{I\backslash J},A}(\pi_{K_I}),\\
\Pi_{\mathrm{an},\mathrm{con},I,\widetilde{J},\breve{I\backslash J},A}(\pi_{K_I}):=\varinjlim_{r_I}\varprojlim_{s_I} \Pi_{[s_I,r_I],I,\widetilde{J},\breve{I\backslash J},A}(\pi_{K_I}),\\
\Pi_{\mathrm{an},\mathrm{con},I,J,\widetilde{I\backslash J},A}(\pi_{K_I}):=\varinjlim_{r_I}\varprojlim_{s_I} \Pi_{[s_I,r_I],I,J,\widetilde{I\backslash J},A}(\pi_{K_I}),\\	
\Pi_{\mathrm{an},\mathrm{con},I,\breve{J},\widetilde{I\backslash J},A}(\pi_{K_I}):=\varinjlim_{r_I}\varprojlim_{s_I} \Pi_{[s_I,r_I],I,\breve{J},\widetilde{I\backslash J},A}(\pi_{K_I}),\\	
\Pi_{\mathrm{an},\mathrm{con},I,\widetilde{J},\widetilde{I\backslash J},A}(\pi_{K_I}):=\varinjlim_{r_I}\varprojlim_{s_I} \Pi_{[s_I,r_I],I,\widetilde{J},\widetilde{I\backslash J},A}(\pi_{K_I})
\end{align}
to be the corresponding base changes of some $\varphi_I$-modules over the rings:
\begin{align}
%\Pi_{\mathrm{an},r_I,I,J,I\backslash J,A},\\	
\Pi_{\mathrm{an},r_I,I,\breve{J},I\backslash J,A}(\pi_{K_I}):=\varprojlim_{s_I}\Pi_{[s_I,r_I],I,\breve{J},I\backslash J,A}(\pi_{K_I}),\\	
\Pi_{\mathrm{an},r_I,I,\widetilde{J},I\backslash J,A}(\pi_{K_I}):=\varprojlim_{s_I} \Pi_{[s_I,r_I],I,\widetilde{J},I\backslash J,A}(\pi_{K_I}),\\
\Pi_{\mathrm{an},r_I,I,J,\breve{I\backslash J},A}(\pi_{K_I}):=\varprojlim_{s_I}\Pi_{[s_I,r_I],I,J,\breve{I\backslash J},A}(\pi_{K_I}),\\	
\Pi_{\mathrm{an},r_I,I,\breve{J},\breve{I\backslash J},A}(\pi_{K_I}):=\varprojlim_{s_I} \Pi_{[s_I,r_I],I,\breve{J},\breve{I\backslash J},A}(\pi_{K_I}),\\	
\Pi_{\mathrm{an},r_I,I,\widetilde{J},\breve{I\backslash J},A}(\pi_{K_I}):=\varprojlim_{s_I} \Pi_{[s_I,r_I],I,\widetilde{J},\breve{I\backslash J},A}(\pi_{K_I}),\\
\Pi_{\mathrm{an},r_I,I,J,\widetilde{I\backslash J},A}(\pi_{K_I}):=\varprojlim_{s_I} \Pi_{[s_I,r_I],I,J,\widetilde{I\backslash J},A}(\pi_{K_I}),\\	
\Pi_{\mathrm{an},r_I,I,\breve{J},\widetilde{I\backslash J},A}(\pi_{K_I}):=\varprojlim_{s_I} \Pi_{[s_I,r_I],I,\breve{J},\widetilde{I\backslash J},A}(\pi_{K_I}),\\	
\Pi_{\mathrm{an},r_I,I,\widetilde{J},\widetilde{I\backslash J},A}(\pi_{K_I}):=\varprojlim_{s_I} \Pi_{[s_I,r_I],I,\widetilde{J},\widetilde{I\backslash J},A}(\pi_{K_I}),	
\end{align}
with the corresponding requirement that they are basically complete with respect to the natural topology and the partial Frobenius operators are commuting with each other.

\end{definition}

\begin{definition} \mbox{\bf{(After KPX \cite[Definition 2.2.6]{5KPX})}}
Then we define a corresponding pseudocoherent or finite projective $\varphi_I$-sheaf $F$ over one of the following period rings:
\begin{align}
%\Pi_{\mathrm{an},r_I,I,J,I\backslash J,A},\\	
\Pi_{\mathrm{an},r_{I,0},I,\breve{J},I\backslash J,A}(\pi_{K_I}):=\varprojlim_{s_I}\Pi_{[s_I,r_I],I,\breve{J},I\backslash J,A}(\pi_{K_I}),\\	
\Pi_{\mathrm{an},r_{I,0},I,\widetilde{J},I\backslash J,A}(\pi_{K_I}):=\varprojlim_{s_I} \Pi_{[s_I,r_I],I,\widetilde{J},I\backslash J,A}(\pi_{K_I}),\\
\Pi_{\mathrm{an},r_{I,0},I,J,\breve{I\backslash J},A}(\pi_{K_I}):=\varprojlim_{s_I}\Pi_{[s_I,r_I],I,J,\breve{I\backslash J},A}(\pi_{K_I}),\\	
\Pi_{\mathrm{an},r_{I,0},I,\breve{J},\breve{I\backslash J},A}(\pi_{K_I}):=\varprojlim_{s_I} \Pi_{[s_I,r_I],I,\breve{J},\breve{I\backslash J},A}(\pi_{K_I}),\\	
\Pi_{\mathrm{an},r_{I,0},I,\widetilde{J},\breve{I\backslash J},A}(\pi_{K_I}):=\varprojlim_{s_I} \Pi_{[s_I,r_I],I,\widetilde{J},\breve{I\backslash J},A}(\pi_{K_I}),\\
\Pi_{\mathrm{an},r_{I,0},I,J,\widetilde{I\backslash J},A}(\pi_{K_I}):=\varprojlim_{s_I} \Pi_{[s_I,r_I],I,J,\widetilde{I\backslash J},A}(\pi_{K_I}),\\	
\Pi_{\mathrm{an},r_{I,0},I,\breve{J},\widetilde{I\backslash J},A}(\pi_{K_I}):=\varprojlim_{s_I} \Pi_{[s_I,r_I],I,\breve{J},\widetilde{I\backslash J},A}(\pi_{K_I}),\\	
\Pi_{\mathrm{an},r_{I,0},I,\widetilde{J},\widetilde{I\backslash J},A}(\pi_{K_I}):=\varprojlim_{s_I} \Pi_{[s_I,r_I],I,\widetilde{J},\widetilde{I\backslash J},A}(\pi_{K_I})	
\end{align}	
to be the corresponding compatible family of the $\varphi_I$-modules over any one $\Pi_{[s_I,r_I],I,*,*,A}(\pi_{K_I})$ of the following rings:
\begin{align}
%\Pi_{[s_I,r_I],I,J,I\backslash J,A},\\	
\Pi_{[s_I,r_I],I,\breve{J},I\backslash J,A}(\pi_{K_I}),\\	
\Pi_{[s_I,r_I],I,\widetilde{J},I\backslash J,A}(\pi_{K_I}),\\
\Pi_{[s_I,r_I],I,J,\breve{I\backslash J},A}(\pi_{K_I}),\\	
\Pi_{[s_I,r_I],I,\breve{J},\breve{I\backslash J},A}(\pi_{K_I}),\\
\Pi_{[s_I,r_I],I,\widetilde{J},\breve{I\backslash J},A}(\pi_{K_I}),\\
\Pi_{[s_I,r_I],I,J,\widetilde{I\backslash J},A}(\pi_{K_I}),\\	
\Pi_{[s_I,r_I],I,\breve{J},\widetilde{I\backslash J},A}(\pi_{K_I}),\\
\Pi_{[s_I,r_I],I,\widetilde{J},\widetilde{I\backslash J},A}(\pi_{K_I})	
\end{align}
satisfying the corresponding restriction requirement and the corresponding cocycle condition as in \cite[Definition 2.2.6]{5KPX}, such that $[s_I,r_I]\subset (0,r_{I,0}]$.
\end{definition}

\begin{definition}  \mbox{\bf{(After KPX \cite[Definition 2.2.12]{5KPX})}}
As in \cite[Definition 2.2.12]{5KPX} we impose the corresponding $\Gamma_I$-structure by adding the corresponding semilinear continuous action of $\Gamma_I$ on the modules induced from that on the period rings, which are assumed to be commuting with the action from $\varphi_I$.	
\end{definition}

\subsection{The Comparison Theorems}

\indent Now we establish some results on the comparison on the corresponding $\varphi_I$-modules we defined above.

\begin{theorem}\mbox{\bf{(After KPX \cite[Proposition 2.2.7]{5KPX})}} 
Consider the following categories:\\
1. The category of all the finite projective $\varphi_I$-modules over the ring $\Pi_{\mathrm{an},r_{I,0},I,?,?,A}(\pi_{K_I})$;\\
2. The category of all the finite projective $\varphi_I$-sheaves over the ring $\Pi_{\mathrm{an},r_{I,0},I,?,?,A}(\pi_{K_I})$.\\
Then we have that the two categories are equivalent.	
\end{theorem}

\begin{proof}
The base change gives rise to the corresponding fully faithful functor from the first category to the second one, while to show the corresponding essential surjectivity, consider the corresponding multi-interval $[r_{1,0}/p,r_{1,0}]\times...\times [r_{I,0}/p,r_{I,0}]$ and use the corresponding Frobenius to reach all the corresponding intervals taking the general form of:
\begin{align}
[r_{1,0}/p^{k_1},r_{1,0}/p^{k_1-1}]\times...\times [r_{I,0}/p^{k_I},r_{I,0}/p^{k_I-1}],k_\alpha=1,2,...,\forall\alpha\in I.	
\end{align}
This forms a $2^{|I|}$-uniform covering of the whole space. And the corresponding uniform finiteness of the modules over each 
\begin{align}
[r_{1,0}/p^{k_1},r_{1,0}/p^{k_1-1}]\times...\times [r_{I,0}/p^{k_I},r_{I,0}/p^{k_I-1}],k_\alpha=1,2,...,\forall\alpha\in I.	
\end{align}	
could be achieved by using the corresponding partial Frobenius actions. Then we are done by applying \cref{proposition2.19}.
\end{proof}

\begin{theorem}\mbox{\bf{(After KPX \cite[Proposition 2.2.7]{5KPX})}} 
Consider the following categories:\\
1. The category of all the pseudocoherent $\varphi_I$-modules over the ring $\Pi_{\mathrm{an},r_{I,0},I,?,?,A}(\pi_{K_I})$;\\
2. The category of all the pseudocoherent $\varphi_I$-sheaves over the ring $\Pi_{\mathrm{an},r_{I,0},I,?,?,A}(\pi_{K_I})$.\\
Then we have that the two categories are equivalent.	
\end{theorem}

\begin{proof}
The base change gives rise to the corresponding fully faithful functor from the first category to the second one, while to show the corresponding essential surjectivity, consider the corresponding multi-interval $[r_{1,0}/p,r_{1,0}]\times...\times [r_{I,0}/p,r_{I,0}]$ and use the corresponding Frobenius to reach all the corresponding intervals taking the general form of:
\begin{align}
[r_{1,0}/p^{k_1},r_{1,0}/p^{k_1-1}]\times...\times [r_{I,0}/p^{k_I},r_{I,0}/p^{k_I-1}],k_\alpha=1,2,...,\forall\alpha\in I.	
\end{align}
This forms a $2^{|I|}$-uniform covering of the whole space. And the corresponding uniform finiteness of the modules over each 
\begin{align}
[r_{1,0}/p^{k_1},r_{1,0}/p^{k_1-1}]\times...\times [r_{I,0}/p^{k_I},r_{I,0}/p^{k_I-1}],k_\alpha=1,2,...,\forall\alpha\in I.	
\end{align}	
could be achieved by using the corresponding partial Frobenius actions. Then we are done by applying \cref{5proposition2.18}. Then one choose finite free covering to promote the finiteness to pseudocoherence as in \cite[Theorem 4.6.1, Lemma 5.4.11]{5KL2}.
\end{proof}

%%%%%%%%%%%%%%%!!!!!!!!!!!!!!!!!!

\indent We now consider the corresponding vertical comparison for $I=\{1,2\}$:

\begin{theorem}\mbox{\bf{(After Kedlaya-Liu \cite[Theorem 5.7.5]{5KL2})}} 
Consider the following categories:\\
1.The category of all the finite projective $(\varphi_I,\Gamma_I)$-modules over the ring $\Pi_{[s_I,r_{I}],I,J,I\backslash J,A}(\pi_{K_I})$;\\
2.The category of all the finite projective $(\varphi_I,\Gamma_I)$-modules over the ring $\Pi_{[s_I,r_{I}],I,J,\breve{I\backslash J},A}(\pi_{K_I})$;\\
3.The category of all the finite projective $(\varphi_I,\Gamma_I)$-modules over the ring $\Pi_{[s_I,r_{I}],I,J,\widetilde{I\backslash J},A}(\pi_{K_I})$.\\
Then we have that these categories are equivalent. Here $0< s_\alpha<r_\alpha<\infty$ for any $\alpha\in I$.
\end{theorem}

\indent This is the consequence of the following theorem:

\begin{theorem} \mbox{\bf{(After Kedlaya-Liu \cite[Theorem 5.7.5]{5KL2})}}
Consider the following categories:\\
1.The category of all the finite projective $(\varphi_2,\Gamma_2)$-modules over the ring $\Pi_{[s_I,r_{I}],I,{1},I\backslash J,A}(\pi_{K_I})$;\\
2.The category of all the finite projective $(\varphi_2,\Gamma_2)$-modules over the ring $\Pi_{[s_I,r_{I}],I,{1},\breve{I\backslash J},A}(\pi_{K_I})$;\\
3.The category of all the finite projective $(\varphi_2,\Gamma_2)$-modules over the ring $\Pi_{[s_I,r_{I}],I,{1},\widetilde{I\backslash J},A}(\pi_{K_I})$.\\
Then we have that these categories are equivalent. Here $0< s_\alpha<r_\alpha<\infty$ for any $\alpha\in I$.
\end{theorem}

\begin{proof}
In this situation this is just the corresponding relative comparison for finite projective $(\varphi,\Gamma)$ modules, which is \cite[Theorem 4.4]{5KP}.	
\end{proof}

\indent We now consider the corresponding vertical comparison in the following context for $I=\{1,2\}$:

\begin{theorem}\mbox{\bf{(After Kedlaya-Liu \cite[Theorem 5.7.5]{5KL2})}}
Consider the following categories:\\
1.The category of all the pseudocoherent $(\varphi_I,\Gamma_I)$-modules over the ring $\Pi_{[s_I,r_{I}],I,J,I\backslash J,A}(\pi_{K_I})$;\\
2.The category of all the pseudocoherent $(\varphi_I,\Gamma_I)$-modules over the ring $\Pi_{[s_I,r_{I}],I,J,\breve{I\backslash J},A}(\pi_{K_I})$;\\
3.The category of all the pseudocoherent $(\varphi_I,\Gamma_I)$-modules over the ring $\Pi_{[s_I,r_{I}],I,J,\widetilde{I\backslash J},A}(\pi_{K_I})$.\\
Then we have that these categories are equivalent. Here $0< s_\alpha<r_\alpha<\infty$ for any $\alpha\in I$.
\end{theorem}

\indent This is the consequence of the following theorem:

\begin{theorem} \mbox{\bf{(After Kedlaya-Liu \cite[Theorem 5.7.5]{5KL2})}}
Consider the following categories:\\
1.The category of all the pseudocoherent $(\varphi_2,\Gamma_2)$-modules over the ring $\Pi_{[s_I,r_{I}],I,{1},I\backslash J,A}(\pi_{K_I})$;\\
2.The category of all the pseudocoherent $(\varphi_2,\Gamma_2)$-modules over the ring $\Pi_{[s_I,r_{I}],I,{1},\breve{I\backslash J},A}(\pi_{K_I})$;\\
3.The category of all the pseudocoherent $(\varphi_2,\Gamma_2)$-modules over the ring $\Pi_{[s_I,r_{I}],I,{1},\widetilde{I\backslash J},A}(\pi_{K_I})$.\\
Then we have that these categories are equivalent. Here $0< s_\alpha<r_\alpha<\infty$ for any $\alpha\in I$.
\end{theorem}

\begin{proof}
In this situation this is just the corresponding relative comparison for finite projective $(\varphi,\Gamma)$ modules, which is \cite[Proposition 5.44, Proposition 5.51]{5T2}.	
\end{proof}

\begin{corollary}
Assume $I=\{1,2\}$. Consider the following categories:\\
1.The category of all the pseudocoherent $(\varphi_I,\Gamma_I)$-modules over the ring $\Pi_{[s_I,r_{I}],I,J,I\backslash J,A}(\pi_{K_I})$;\\
2.The category of all the pseudocoherent $(\varphi_I,\Gamma_I)$-modules over the ring $\Pi_{[s_I,r_{I}],I,J,\breve{I\backslash J},A}(\pi_{K_I})$;\\
3.The category of all the pseudocoherent $(\varphi_I,\Gamma_I)$-modules over the ring $\Pi_{[s_I,r_{I}],I,J,\widetilde{I\backslash J},A}(\pi_{K_I})$.\\
Then we have that these categories are equivalent. Here $0< s_\alpha\leq r_\alpha/p<\infty$ for any $\alpha\in I$. \\
Consider the following categories:\\
1.The category of all the finite projective $(\varphi_I,\Gamma_I)$-modules over the ring $\Pi_{[s_I,r_{I}],I,J,I\backslash J,A}(\pi_{K_I})$;\\
2.The category of all the finite projective $(\varphi_I,\Gamma_I)$-modules over the ring $\Pi_{[s_I,r_{I}],I,J,\breve{I\backslash J},A}(\pi_{K_I})$;\\
3.The category of all the finite projective $(\varphi_I,\Gamma_I)$-modules over the ring $\Pi_{[s_I,r_{I}],I,J,\widetilde{I\backslash J},A}(\pi_{K_I})$.\\
Then we have that these categories are equivalent. Here $0< s_\alpha\leq r_\alpha/p<\infty$ for any $\alpha\in I$.	
\end{corollary}

%%\newpage

\newpage\section{$B_{I}$-Pairs and Intermediate Objects in Rigid Family}

\subsection{Mixed-type Hodge Structures}

\indent Now we work with the corresponding $B$-pairs and some mixed-type objects as in \cite{5Ber1} and \cite{5Nak1}.

\begin{definition}
Define $B^+_{\mathrm{dR},I}:=\mathbb{C}_p[[t_1,...,t_I]]$, and define\\  $B_{\mathrm{dR},I}:=\mathbb{C}_p[[t_1,...,t_I]][t_1^{-1},...,t_I^{-1}]$, and similarly we have the obvious higher dimensional analog $B_{e,I}$ of the corresponding period ring $B_{e}$.

.	
\end{definition}

\indent Then after Bloch-Kato \cite{5BK1} we have the following fundamental sequence:

\begin{proposition}
We have the higher dimensional generalization of the corresponding Bloch-Kato fundamental sequence:
\[
\xymatrix@C+0pc@R+0pc{
0 \ar[r] \ar[r] \ar[r] &\mathbb{Q}_p \ar[r] \ar[r] \ar[r]  &B_{e,I}\bigoplus B^+_{\mathrm{dR},I} \ar[r] \ar[r] \ar[r] &B_{\mathrm{dR},I} \ar[r] \ar[r] \ar[r] &0
}
\]
induced by the corresponding short exact sequence:
\[
\xymatrix@C+0pc@R+0pc{
0 \ar[r] \ar[r] \ar[r] &\mathbb{Q}_p \ar[r] \ar[r] \ar[r]  &B_{e,I}\bigoplus \mathbb{C}_p[[t_1,...,t_I]] \ar[r] \ar[r] \ar[r] &\mathbb{C}_p[[t_1,...,t_I]][t_1^{-1},...,t_I^{-1}] \ar[r] \ar[r] \ar[r] &0.
}
\]	
\end{proposition}

\begin{setting}
In what follows, we assume that the corresponding $A$ to be a rigid affinoid in rigid analytic geometry over $\mathbb{Q}_p$.	
\end{setting}

\begin{definition}
We now consider the following rings:
\begin{align}
B^+_{\mathrm{dR},I'}	\widehat{\otimes}\Pi_{[s_I,r_I],I,J,I\backslash J,A},\\	
B^+_{\mathrm{dR},I'}	\widehat{\otimes}\Pi_{[s_I,r_I],I,\breve{J},I\backslash J,A},\\	
B^+_{\mathrm{dR},I'}	\widehat{\otimes}\Pi_{[s_I,r_I],I,\widetilde{J},I\backslash J,A},\\
B^+_{\mathrm{dR},I'}	\widehat{\otimes}\Pi_{[s_I,r_I],I,J,\breve{I\backslash J},A},\\	
B^+_{\mathrm{dR},I'}	\widehat{\otimes}\Pi_{[s_I,r_I],I,\breve{J},\breve{I\backslash J},A},\\
B^+_{\mathrm{dR},I'}	\widehat{\otimes}\Pi_{[s_I,r_I],I,\widetilde{J},\breve{I\backslash J},A},\\
B^+_{\mathrm{dR},I'}	\widehat{\otimes}\Pi_{[s_I,r_I],I,J,\widetilde{I\backslash J},A},\\	
B^+_{\mathrm{dR},I'}	\widehat{\otimes}\Pi_{[s_I,r_I],I,\breve{J},\widetilde{I\backslash J},A},\\
B^+_{\mathrm{dR},I'}	\widehat{\otimes}\Pi_{[s_I,r_I],I,\widetilde{J},\widetilde{I\backslash J},A}
\end{align}

with

\begin{align}
B^+_{\mathrm{dR},I'}	\widehat{\otimes}\Pi_{[s_I,r_I],I,J,I\backslash J,A}[t_1^{-1},...,t_{I'}^{-1}],\\	
B^+_{\mathrm{dR},I'}	\widehat{\otimes}\Pi_{[s_I,r_I],I,\breve{J},I\backslash J,A}[t_1^{-1},...,t_{I'}^{-1}],\\	
B^+_{\mathrm{dR},I'}	\widehat{\otimes}\Pi_{[s_I,r_I],I,\widetilde{J},I\backslash J,A}[t_1^{-1},...,t_{I'}^{-1}],\\
B^+_{\mathrm{dR},I'}	\widehat{\otimes}\Pi_{[s_I,r_I],I,J,\breve{I\backslash J},A}[t_1^{-1},...,t_{I'}^{-1}],\\	
B^+_{\mathrm{dR},I'}	\widehat{\otimes}\Pi_{[s_I,r_I],I,\breve{J},\breve{I\backslash J},A}[t_1^{-1},...,t_{I'}^{-1}],\\
B^+_{\mathrm{dR},I'}	\widehat{\otimes}\Pi_{[s_I,r_I],I,\widetilde{J},\breve{I\backslash J},A}[t_1^{-1},...,t_{I'}^{-1}],\\
B^+_{\mathrm{dR},I'}	\widehat{\otimes}\Pi_{[s_I,r_I],I,J,\widetilde{I\backslash J},A}[t_1^{-1},...,t_{I'}^{-1}],\\	
B^+_{\mathrm{dR},I'}	\widehat{\otimes}\Pi_{[s_I,r_I],I,\breve{J},\widetilde{I\backslash J},A}[t_1^{-1},...,t_{I'}^{-1}],\\
B^+_{\mathrm{dR},I'}	\widehat{\otimes}\Pi_{[s_I,r_I],I,\widetilde{J},\widetilde{I\backslash J},A}[t_1^{-1},...,t_{I'}^{-1}],
\end{align}

with 

\begin{align}
B_{e,I'}	\widehat{\otimes}\Pi_{[s_I,r_I],I,J,I\backslash J,A},\\	
B_{e,I'}	\widehat{\otimes}\Pi_{[s_I,r_I],I,\breve{J},I\backslash J,A},\\	
B_{e,I'}	\widehat{\otimes}\Pi_{[s_I,r_I],I,\widetilde{J},I\backslash J,A},\\
B_{e,I'}	\widehat{\otimes}\Pi_{[s_I,r_I],I,J,\breve{I\backslash J},A},\\	
B_{e,I'}	\widehat{\otimes}\Pi_{[s_I,r_I],I,\breve{J},\breve{I\backslash J},A},\\
B_{e,I'}	\widehat{\otimes}\Pi_{[s_I,r_I],I,\widetilde{J},\breve{I\backslash J},A},\\
B_{e,I'}	\widehat{\otimes}\Pi_{[s_I,r_I],I,J,\widetilde{I\backslash J},A},\\	
B_{e,I'}	\widehat{\otimes}\Pi_{[s_I,r_I],I,\breve{J},\widetilde{I\backslash J},A},\\
B_{e,I'} \widehat{\otimes}\Pi_{[s_I,r_I],I,\widetilde{J},\widetilde{I\backslash J},A}
\end{align}

with

\begin{align}
\varprojlim_{s_I}  B^+_{\mathrm{dR},I'}	\widehat{\otimes}\Pi_{[s_I,r_I],I,J,I\backslash J,A},\\	
\varprojlim_{s_I}  B^+_{\mathrm{dR},I'}	\widehat{\otimes}\Pi_{[s_I,r_I],I,\breve{J},I\backslash J,A},\\	
\varprojlim_{s_I} B^+_{\mathrm{dR},I'}	\widehat{\otimes}\Pi_{[s_I,r_I],I,\widetilde{J},I\backslash J,A},\\
\varprojlim_{s_I} B^+_{\mathrm{dR},I'}	\widehat{\otimes}\Pi_{[s_I,r_I],I,J,\breve{I\backslash J},A},\\	
\varprojlim_{s_I} B^+_{\mathrm{dR},I'}	\widehat{\otimes}\Pi_{[s_I,r_I],I,\breve{J},\breve{I\backslash J},A},\\
\varprojlim_{s_I} B^+_{\mathrm{dR},I'}	\widehat{\otimes}\Pi_{[s_I,r_I],I,\widetilde{J},\breve{I\backslash J},A},\\
\varprojlim_{s_I} B^+_{\mathrm{dR},I'}	\widehat{\otimes}\Pi_{[s_I,r_I],I,J,\widetilde{I\backslash J},A},\\	
\varprojlim_{s_I} B^+_{\mathrm{dR},I'}	\widehat{\otimes}\Pi_{[s_I,r_I],I,\breve{J},\widetilde{I\backslash J},A},\\
\varprojlim_{s_I} B^+_{\mathrm{dR},I'}	\widehat{\otimes}\Pi_{[s_I,r_I],I,\widetilde{J},\widetilde{I\backslash J},A}
\end{align}

with

\begin{align}
\varprojlim_{s_I} B^+_{\mathrm{dR},I'}	\widehat{\otimes}\Pi_{[s_I,r_I],I,J,I\backslash J,A}[t_1^{-1},...,t_{I'}^{-1}],\\	
\varprojlim_{s_I} B^+_{\mathrm{dR},I'}	\widehat{\otimes}\Pi_{[s_I,r_I],I,\breve{J},I\backslash J,A}[t_1^{-1},...,t_{I'}^{-1}],\\	
\varprojlim_{s_I} B^+_{\mathrm{dR},I'}	\widehat{\otimes}\Pi_{[s_I,r_I],I,\widetilde{J},I\backslash J,A}[t_1^{-1},...,t_{I'}^{-1}],\\
\varprojlim_{s_I} B^+_{\mathrm{dR},I'}	\widehat{\otimes}\Pi_{[s_I,r_I],I,J,\breve{I\backslash J},A}[t_1^{-1},...,t_{I'}^{-1}],\\	
\varprojlim_{s_I} B^+_{\mathrm{dR},I'}	\widehat{\otimes}\Pi_{[s_I,r_I],I,\breve{J},\breve{I\backslash J},A}[t_1^{-1},...,t_{I'}^{-1}],\\
\varprojlim_{s_I} B^+_{\mathrm{dR},I'}	\widehat{\otimes}\Pi_{[s_I,r_I],I,\widetilde{J},\breve{I\backslash J},A}[t_1^{-1},...,t_{I'}^{-1}],\\
\varprojlim_{s_I} B^+_{\mathrm{dR},I'}	\widehat{\otimes}\Pi_{[s_I,r_I],I,J,\widetilde{I\backslash J},A}[t_1^{-1},...,t_{I'}^{-1}],\\	
\varprojlim_{s_I} B^+_{\mathrm{dR},I'}	\widehat{\otimes}\Pi_{[s_I,r_I],I,\breve{J},\widetilde{I\backslash J},A}[t_1^{-1},...,t_{I'}^{-1}],\\
\varprojlim_{s_I} B^+_{\mathrm{dR},I'}	\widehat{\otimes}\Pi_{[s_I,r_I],I,\widetilde{J},\widetilde{I\backslash J},A}[t_1^{-1},...,t_{I'}^{-1}],
\end{align}

with 

\begin{align}
\varprojlim_{s_I} B_{e,I'}	\widehat{\otimes}\Pi_{[s_I,r_I],I,J,I\backslash J,A},\\	
\varprojlim_{s_I} B_{e,I'}	\widehat{\otimes}\Pi_{[s_I,r_I],I,\breve{J},I\backslash J,A},\\	
\varprojlim_{s_I} B_{e,I'}	\widehat{\otimes}\Pi_{[s_I,r_I],I,\widetilde{J},I\backslash J,A},\\
\varprojlim_{s_I} B_{e,I'}	\widehat{\otimes}\Pi_{[s_I,r_I],I,J,\breve{I\backslash J},A},\\	
\varprojlim_{s_I} B_{e,I'}	\widehat{\otimes}\Pi_{[s_I,r_I],I,\breve{J},\breve{I\backslash J},A},\\
\varprojlim_{s_I} B_{e,I'}	\widehat{\otimes}\Pi_{[s_I,r_I],I,\widetilde{J},\breve{I\backslash J},A},\\
\varprojlim_{s_I} B_{e,I'}	\widehat{\otimes}\Pi_{[s_I,r_I],I,J,\widetilde{I\backslash J},A},\\	
\varprojlim_{s_I} B_{e,I'}	\widehat{\otimes}\Pi_{[s_I,r_I],I,\breve{J},\widetilde{I\backslash J},A},\\
\varprojlim_{s_I} B_{e,I'} \widehat{\otimes}\Pi_{[s_I,r_I],I,\widetilde{J},\widetilde{I\backslash J},A}
\end{align}

with

\begin{align}
\varinjlim_{r_I}\varprojlim_{s_I}  B^+_{\mathrm{dR},I'}	\widehat{\otimes}\Pi_{[s_I,r_I],I,J,I\backslash J,A},\\	
\varinjlim_{r_I}\varprojlim_{s_I}  B^+_{\mathrm{dR},I'}	\widehat{\otimes}\Pi_{[s_I,r_I],I,\breve{J},I\backslash J,A},\\	
\varinjlim_{r_I}\varprojlim_{s_I} B^+_{\mathrm{dR},I'}	\widehat{\otimes}\Pi_{[s_I,r_I],I,\widetilde{J},I\backslash J,A},\\
\varinjlim_{r_I}\varprojlim_{s_I} B^+_{\mathrm{dR},I'}	\widehat{\otimes}\Pi_{[s_I,r_I],I,J,\breve{I\backslash J},A},\\	
\varinjlim_{r_I}\varprojlim_{s_I} B^+_{\mathrm{dR},I'}	\widehat{\otimes}\Pi_{[s_I,r_I],I,\breve{J},\breve{I\backslash J},A},\\
\varinjlim_{r_I}\varprojlim_{s_I} B^+_{\mathrm{dR},I'}	\widehat{\otimes}\Pi_{[s_I,r_I],I,\widetilde{J},\breve{I\backslash J},A},\\
\varinjlim_{r_I}\varprojlim_{s_I} B^+_{\mathrm{dR},I'}	\widehat{\otimes}\Pi_{[s_I,r_I],I,J,\widetilde{I\backslash J},A},\\	
\varinjlim_{r_I}\varprojlim_{s_I} B^+_{\mathrm{dR},I'}	\widehat{\otimes}\Pi_{[s_I,r_I],I,\breve{J},\widetilde{I\backslash J},A},\\
\varinjlim_{r_I}\varprojlim_{s_I} B^+_{\mathrm{dR},I'}	\widehat{\otimes}\Pi_{[s_I,r_I],I,\widetilde{J},\widetilde{I\backslash J},A}
\end{align}

with

\begin{align}
\varinjlim_{r_I}\varprojlim_{s_I} B^+_{\mathrm{dR},I'}	\widehat{\otimes}\Pi_{[s_I,r_I],I,J,I\backslash J,A}[t_1^{-1},...,t_{I'}^{-1}],\\	
\varinjlim_{r_I}\varprojlim_{s_I} B^+_{\mathrm{dR},I'}	\widehat{\otimes}\Pi_{[s_I,r_I],I,\breve{J},I\backslash J,A}[t_1^{-1},...,t_{I'}^{-1}],\\	
\varinjlim_{r_I}\varprojlim_{s_I} B^+_{\mathrm{dR},I'}	\widehat{\otimes}\Pi_{[s_I,r_I],I,\widetilde{J},I\backslash J,A}[t_1^{-1},...,t_{I'}^{-1}],\\
\varinjlim_{r_I}\varprojlim_{s_I} B^+_{\mathrm{dR},I'}	\widehat{\otimes}\Pi_{[s_I,r_I],I,J,\breve{I\backslash J},A}[t_1^{-1},...,t_{I'}^{-1}],\\	
\varinjlim_{r_I}\varprojlim_{s_I} B^+_{\mathrm{dR},I'}	\widehat{\otimes}\Pi_{[s_I,r_I],I,\breve{J},\breve{I\backslash J},A}[t_1^{-1},...,t_{I'}^{-1}],\\
\varinjlim_{r_I}\varprojlim_{s_I} B^+_{\mathrm{dR},I'}	\widehat{\otimes}\Pi_{[s_I,r_I],I,\widetilde{J},\breve{I\backslash J},A}[t_1^{-1},...,t_{I'}^{-1}],\\
\varinjlim_{r_I}\varprojlim_{s_I} B^+_{\mathrm{dR},I'}	\widehat{\otimes}\Pi_{[s_I,r_I],I,J,\widetilde{I\backslash J},A}[t_1^{-1},...,t_{I'}^{-1}],\\	
\varinjlim_{r_I}\varprojlim_{s_I} B^+_{\mathrm{dR},I'}	\widehat{\otimes}\Pi_{[s_I,r_I],I,\breve{J},\widetilde{I\backslash J},A}[t_1^{-1},...,t_{I'}^{-1}],\\
\varinjlim_{r_I}\varprojlim_{s_I} B^+_{\mathrm{dR},I'}	\widehat{\otimes}\Pi_{[s_I,r_I],I,\widetilde{J},\widetilde{I\backslash J},A}[t_1^{-1},...,t_{I'}^{-1}],
\end{align}

with 

\begin{align}
\varinjlim_{r_I}\varprojlim_{s_I} B_{e,I'}	\widehat{\otimes}\Pi_{[s_I,r_I],I,J,I\backslash J,A},\\	
\varinjlim_{r_I}\varprojlim_{s_I} B_{e,I'}	\widehat{\otimes}\Pi_{[s_I,r_I],I,\breve{J},I\backslash J,A},\\	
\varinjlim_{r_I}\varprojlim_{s_I} B_{e,I'}	\widehat{\otimes}\Pi_{[s_I,r_I],I,\widetilde{J},I\backslash J,A},\\
\varinjlim_{r_I}\varprojlim_{s_I} B_{e,I'}	\widehat{\otimes}\Pi_{[s_I,r_I],I,J,\breve{I\backslash J},A},\\	
\varinjlim_{r_I}\varprojlim_{s_I} B_{e,I'}	\widehat{\otimes}\Pi_{[s_I,r_I],I,\breve{J},\breve{I\backslash J},A},\\
\varinjlim_{r_I}\varprojlim_{s_I} B_{e,I'}	\widehat{\otimes}\Pi_{[s_I,r_I],I,\widetilde{J},\breve{I\backslash J},A},\\
\varinjlim_{r_I}\varprojlim_{s_I} B_{e,I'}	\widehat{\otimes}\Pi_{[s_I,r_I],I,J,\widetilde{I\backslash J},A},\\	
\varinjlim_{r_I}\varprojlim_{s_I} B_{e,I'}	\widehat{\otimes}\Pi_{[s_I,r_I],I,\breve{J},\widetilde{I\backslash J},A},\\
\varinjlim_{r_I}\varprojlim_{s_I} B_{e,I'} \widehat{\otimes}\Pi_{[s_I,r_I],I,\widetilde{J},\widetilde{I\backslash J},A}.
\end{align}\\
	
\end{definition}

\indent Now we combine the construction in the following coherent way following \cite[Section 2]{5Ber1}, \cite[Definition 2.2]{5Nak1} and \cite[Definition 2.2.6]{5KPX}.

\begin{definition} 
We define a $B_{I'}$-$(\varphi_I,\Gamma_I)$-module over
\begin{center}
$(B^+_{\mathrm{dR},I'}	\widehat{\otimes}\Pi_{[s_I,r_I],I,J,I\backslash J,A},B^+_{\mathrm{dR},I'}	\widehat{\otimes}\Pi_{[s_I,r_I],I,*,*,A}[t_1^{-1},...,t_I'^{-1}],B_{e,I'}	\widehat{\otimes}\Pi_{[s_I,r_I],I,*,*,A})$	
\end{center}
to be a triplet of finite projective modules:
\begin{align}
(M_e,M_{\mathrm{dR}},M^+_{\mathrm{dR}})	
\end{align}
over 
\begin{center}
$(B^+_{\mathrm{dR},I'}	\widehat{\otimes}\Pi_{[s_I,r_I],I,*,*,A},B^+_{\mathrm{dR},I'}	\widehat{\otimes}\Pi_{[s_I,r_I],I,*,*,A}[t_1^{-1},...,t_I'^{-1}],B_{e,I'}	\widehat{\otimes}\Pi_{[s_I,r_I],I,*,*,A})$	
\end{center}
such that we have glueing datum along :
\begin{align}
B^+_{\mathrm{dR},I'}	\widehat{\otimes}\Pi_{[s_I,r_I],I,*,*,A}\rightarrow B^+_{\mathrm{dR},I'}	\widehat{\otimes}\Pi_{[s_I,r_I],I,*,*,A}[t_1^{-1},...,t_I'^{-1}] \leftarrow B_{e,I'}	\widehat{\otimes}\Pi_{[s_I,r_I],I,*,*,A}.	
\end{align}
And this carry the corresponding relative Galois action of:
\begin{align}
\mathrm{Gal}_{\mathbb{Q}_p,1}\times...\times \mathrm{Gal}_{\mathbb{Q}_p,I'}	
\end{align}
on the multi de Rham period rings which is semilinear. And we have that the three modules involved are relative $(\varphi_I,\Gamma_I)$-modules relative to
\begin{align}
B^+_{\mathrm{dR},I'},B_{\mathrm{dR},I'},B_{e,I'}.	
\end{align}

\end{definition}

\begin{definition}
We define a pseudocoherent $B_{I'}$-$(\varphi_I,\Gamma_I)$-module over
\begin{center}
$(B^+_{\mathrm{dR},I'}	\widehat{\otimes}\Pi_{[s_I,r_I],I,J,I\backslash J,A},B^+_{\mathrm{dR},I'}	\widehat{\otimes}\Pi_{[s_I,r_I],I,*,*,A}[t_1^{-1},...,t_I'^{-1}],B_{e,I'}	\widehat{\otimes}\Pi_{[s_I,r_I],I,*,*,A})$	
\end{center}
to be a triplet of stably-pseudocoherent modules:
\begin{align}
(M_e,M_{\mathrm{dR}},M^+_{\mathrm{dR}})	
\end{align}
over 
\begin{center}
$(B^+_{\mathrm{dR},I'}	\widehat{\otimes}\Pi_{[s_I,r_I],I,*,*,A},B^+_{\mathrm{dR},I'}	\widehat{\otimes}\Pi_{[s_I,r_I],I,*,*,A}[t_1^{-1},...,t_I'^{-1}],B_{e,I'}	\widehat{\otimes}\Pi_{[s_I,r_I],I,*,*,A})$	
\end{center}
such that this is glueing datum along :
\begin{align}
B^+_{\mathrm{dR},I'}	\widehat{\otimes}\Pi_{[s_I,r_I],I,*,*,A}\rightarrow B^+_{\mathrm{dR},I'}	\widehat{\otimes}\Pi_{[s_I,r_I],I,*,*,A}[t_1^{-1},...,t_I'^{-1}] \leftarrow B_{e,I'}	\widehat{\otimes}\Pi_{[s_I,r_I],I,*,*,A}.	
\end{align}
And this carry the corresponding relative Galois action of:
\begin{align}
\mathrm{Gal}_{\mathbb{Q}_p,1}\times...\times \mathrm{Gal}_{\mathbb{Q}_p,I'}		
\end{align}
on the multi de Rham period rings which is semilinear. And we have that the three modules involved are relative pseudocoherent $(\varphi_I,\Gamma_I)$-modules relative to
\begin{align}
B^+_{\mathrm{dR},I'},B_{\mathrm{dR},I'},B_{e,I'}.	
\end{align}

\end{definition}

\begin{definition}
We define a $B_{I'}$-$(\varphi_I,\Gamma_I)$-module over
\begin{center}
$(B^+_{\mathrm{dR},I'}	\widehat{\otimes}\Pi_{\mathrm{an},r_{I},I,J,I\backslash J,A},B^+_{\mathrm{dR},I'}	\widehat{\otimes}\Pi_{[s_I,r_I],I,*,*,A}[t_1^{-1},...,t_I'^{-1}],B_{e,I'}	\widehat{\otimes}\Pi_{\mathrm{an},r_{I},I,*,*,A})$	
\end{center}
to be a triplet of finite projective modules:
\begin{align}
(M_e,M_{\mathrm{dR}},M^+_{\mathrm{dR}})	
\end{align}
over 
\begin{center}
$(B^+_{\mathrm{dR},I'}	\widehat{\otimes}\Pi_{\mathrm{an},r_{I},I,*,*,A},B^+_{\mathrm{dR},I'}	\widehat{\otimes}\Pi_{\mathrm{an},r_{I},I,*,*,A}[t_1^{-1},...,t_I'^{-1}],B_{e,I'}	\widehat{\otimes}\Pi_{\mathrm{an},r_{I},I,*,*,A})$	
\end{center}
such that this is glueing datum along :
\begin{align}
B^+_{\mathrm{dR},I'}	\widehat{\otimes}\Pi_{\mathrm{an},r_{I},I,*,*,A}\rightarrow B^+_{\mathrm{dR},I'}	\widehat{\otimes}\Pi_{\mathrm{an},r_{I},I,*,*,A}[t_1^{-1},...,t_I'^{-1}] \leftarrow B_{e,I'}	\widehat{\otimes}\Pi_{\mathrm{an},r_{I},I,*,*,A}.	
\end{align}
And this carry the corresponding relative Galois action of:
\begin{align}
\mathrm{Gal}_{\mathbb{Q}_p,1}\times...\times \mathrm{Gal}_{\mathbb{Q}_p,I'}		
\end{align}
on the multi de Rham period rings which is semilinear. And we have that the three modules involved are relative $(\varphi_I,\Gamma_I)$-modules relative to
\begin{align}
B^+_{\mathrm{dR},I'},B_{\mathrm{dR},I'},B_{e,I'}.	
\end{align}

\end{definition}

\begin{definition}
We define a pseudocoherent $B_{I'}$-$(\varphi_I,\Gamma_I)$-module over
\begin{center}
$(B^+_{\mathrm{dR},I'}	\widehat{\otimes}\Pi_{\mathrm{an},r_{I},I,J,I\backslash J,A},B^+_{\mathrm{dR},I'}	\widehat{\otimes}\Pi_{[s_I,r_I],I,*,*,A}[t_1^{-1},...,t_I'^{-1}],B_{e,I'}	\widehat{\otimes}\Pi_{\mathrm{an},r_{I},I,*,*,A})$	
\end{center}
to be a triplet of stably pseudocoherent modules:
\begin{align}
(M_e,M_{\mathrm{dR}},M^+_{\mathrm{dR}})	
\end{align}
over 
\begin{center}
$(B^+_{\mathrm{dR},I'}	\widehat{\otimes}\Pi_{\mathrm{an},r_{I},I,*,*,A},B^+_{\mathrm{dR},I'}	\widehat{\otimes}\Pi_{\mathrm{an},r_{I},I,*,*,A}[t_1^{-1},...,t_I'^{-1}],B_{e,I'}	\widehat{\otimes}\Pi_{\mathrm{an},r_{I},I,*,*,A})$	
\end{center}
such that this is glueing datum along :
\begin{align}
B^+_{\mathrm{dR},I'}	\widehat{\otimes}\Pi_{\mathrm{an},r_{I},I,*,*,A}\rightarrow B^+_{\mathrm{dR},I'}	\widehat{\otimes}\Pi_{\mathrm{an},r_{I},I,*,*,A}[t_1^{-1},...,t_I'^{-1}] \leftarrow B_{e,I'}	\widehat{\otimes}\Pi_{\mathrm{an},r_{I},I,*,*,A}.	
\end{align}
And this carry the corresponding relative Galois action of:
\begin{align}
\mathrm{Gal}_{\mathbb{Q}_p,1}\times...\times \mathrm{Gal}_{\mathbb{Q}_p,I'}		
\end{align}
on the multi de Rham period rings which is semilinear. And we have that the three modules involved are relative pseudocoherent $(\varphi_I,\Gamma_I)$-modules relative to
\begin{align}
B^+_{\mathrm{dR},I'},B_{\mathrm{dR},I'},B_{e,I'}.	
\end{align}\\

\end{definition}

\subsection{Fundamental Comparison on the Mixed-Type Objects}

\begin{proposition} \mbox{\bf{(After Berger \cite[Th\'eor\`eme A]{5Ber1})}}
Let $I'$ be a set consisting of two elements and $I$ is empty, then we have that the category of all the $B_{I'}$-$(\varphi_I,\Gamma_I)$ modules over 
\begin{center}
$(B^+_{\mathrm{dR},I'}	\widehat{\otimes}\Pi_{[s_I,r_I],I,J,I\backslash J,A},B^+_{\mathrm{dR},I'}	\widehat{\otimes}\Pi_{[s_I,r_I],I,J,I\backslash J,A}[t_1^{-1},...,t_I'^{-1}],B_{e,I'}	\widehat{\otimes}\Pi_{[s_I,r_I],I,J,I\backslash J,A})$	
\end{center}
is equivalent to the category of all the $(\varphi_{I'},\Gamma_{I'})$-modules in the finite projective setting.
\end{proposition}
\begin{proof}
One has the result after the following two propositions.	
\end{proof}

\indent We first consider the following comparison:
\begin{proposition}
Let $I'=\{1,2\}$ be a set consisting of two elements and $I$ is empty, then we have that the category of all the $B_{I'}$-$(\varphi_I,\Gamma_I)$ modules over 
\begin{center}
$(B^+_{\mathrm{dR},I'}	\widehat{\otimes}\Pi_{[s_I,r_I],I,J,I\backslash J,A},B^+_{\mathrm{dR},I'}	\widehat{\otimes}\Pi_{[s_I,r_I],I,J,I\backslash J,A}[t_1^{-1},...,t_I'^{-1}],B_{e,I'}	\widehat{\otimes}\Pi_{[s_I,r_I],I,J,I\backslash J,A})$	
\end{center}
is equivalent to the category of all the $B_{\{1\}}$-$(\varphi_{\{2\}},\Gamma_{\{2\}})$ modules.
\end{proposition}
\begin{proof}
This will be the corresponding consequence of the following. Let $I'=\{1,2\}$ be a set consisting of two elements and $I$ is empty, then we have that the category of all the $B_{I'}$-$(\varphi_I,\Gamma_I)$ modules over 
\begin{center}
$(B^+_{\mathrm{dR},I'}	\widehat{\otimes}\Pi_{[s_I,r_I],I,\widetilde{J},\widetilde{I\backslash J},A},B^+_{\mathrm{dR},I'}	\widehat{\otimes}\Pi_{[s_I,r_I],I,\widetilde{J},\widetilde{I\backslash J},A}[t_1^{-1},...,t_I'^{-1}],B_{e,I'}	\widehat{\otimes}\Pi_{[s_I,r_I],I,\widetilde{J},\widetilde{I\backslash J},A})$	
\end{center}
is equivalent to the category of all the $B_{\{1\}}$-$(\varphi_{\{2\}},\Gamma_{\{2\}})$ modules over the corresponding perfected rings with $\widetilde{.}$ accent. However this could be proved as in \cite[Theorem 2.18]{5KP} as long as one works with mod $t^k,k\in \mathbb{Z}$ coefficients (also see the corresponding proof of \cite[Proposition 3.8]{5T3}). To be more precise first we consider the corresponding base change of any $B_{I'}$-$(\varphi_I,\Gamma_I)$ module over 
\begin{center}
$(B^+_{\mathrm{dR},I'}	\widehat{\otimes}\Pi_{[s_I,r_I],I,\widetilde{J},\widetilde{I\backslash J},A},B^+_{\mathrm{dR},I'}	\widehat{\otimes}\Pi_{[s_I,r_I],I,\widetilde{J},\widetilde{I\backslash J},A}[t_1^{-1},...,t_I'^{-1}],B_{e,I'}	\widehat{\otimes}\Pi_{[s_I,r_I],I,\widetilde{J},\widetilde{I\backslash J},A})$	
\end{center}
to $B_\mathrm{dR,\{1\}}$, which then by the strategy above could be associated a $B_{\{1\}}$-$(\varphi_{\{2\}},\Gamma_{\{2\}})$ module over the corresponding perfected rings with $\widetilde{.}$ accent (see \cite[Theorem 2.18]{5KP}, \cite[Proposition 3.8]{5T3}). As in \cite[Theorem 2.18]{5KP} we will have the situation where the category of all the $B_{I'}$-$(\varphi_I,\Gamma_I)$ modules over 
\begin{center}
$(B^+_{\mathrm{dR},I'}	\widehat{\otimes}\Pi_{[s_I,r_I],I,\widetilde{J},\widetilde{I\backslash J},A},B^+_{\mathrm{dR},I'}	\widehat{\otimes}\Pi_{[s_I,r_I],I,\widetilde{J},\widetilde{I\backslash J},A}[t_1^{-1},...,t_I'^{-1}],B_{e,I'}	\widehat{\otimes}\Pi_{[s_I,r_I],I,\widetilde{J},\widetilde{I\backslash J},A})$	
\end{center}
is equivalent to the category of all the $B_{\{1\}}$-$(\varphi_{\{2\}},\Gamma_{\{2\}})$ modules over the corresponding perfected rings with $\widetilde{.}$ accent. However by the proof of \cite[Theorem 4.4]{5KP} we further have the situation where the category of all the $B_{I'}$-$(\varphi_I,\Gamma_I)$ modules over 
\begin{center}
$(B^+_{\mathrm{dR},I'}	\widehat{\otimes}\Pi_{[s_I,r_I],I,\widetilde{J},\widetilde{I\backslash J},A},B^+_{\mathrm{dR},I'}	\widehat{\otimes}\Pi_{[s_I,r_I],I,\widetilde{J},\widetilde{I\backslash J},A}[t_1^{-1},...,t_I'^{-1}],B_{e,I'}	\widehat{\otimes}\Pi_{[s_I,r_I],I,\widetilde{J},\widetilde{I\backslash J},A})$	
\end{center}
is equivalent to the category of all the $B_{\{1\}}$-$(\varphi_{\{2\}},\Gamma_{\{2\}})$ modules over the corresponding perfected rings with no accent. 
\end{proof}

\begin{proposition}
Let $I=\{1,2\}$ be a set consisting of two elements and $I'$ is empty, then we have that the category of all the $B_{I'}$-$(\varphi_I,\Gamma_I)$ modules over 
\begin{center}
$(B^+_{\mathrm{dR},I'}	\widehat{\otimes}\Pi_{[s_I,r_I],I,J,I\backslash J,A},B^+_{\mathrm{dR},I'}	\widehat{\otimes}\Pi_{[s_I,r_I],I,J,I\backslash J,A}[t_1^{-1},...,t_I'^{-1}],B_{e,I'}	\widehat{\otimes}\Pi_{[s_I,r_I],I,J,I\backslash J,A})$	
\end{center}
is equivalent to the category of all the $B_{\{1\}}$-$(\varphi_{\{2\}},\Gamma_{\{2\}})$ modules.

\end{proposition}

\begin{proof}
This will be the corresponding consequence of the following. Let $I=\{1,2\}$ be a set consisting of two elements and $I'$ is empty, then we have that the category of all the $B_{I'}$-$(\varphi_I,\Gamma_I)$ modules over 
\begin{center}
$(B^+_{\mathrm{dR},I'}	\widehat{\otimes}\Pi_{[s_I,r_I],I,\widetilde{J},\widetilde{I\backslash J},A},B^+_{\mathrm{dR},I'}	\widehat{\otimes}\Pi_{[s_I,r_I],I,\widetilde{J},\widetilde{I\backslash J},A}[t_1^{-1},...,t_I'^{-1}],B_{e,I'}	\widehat{\otimes}\Pi_{[s_I,r_I],I,\widetilde{J},\widetilde{I\backslash J},A})$	
\end{center}
is equivalent to the category of all the $B_{\{1\}}$-$(\varphi_{\{2\}},\Gamma_{\{2\}})$ modules over the corresponding perfected rings with $\widetilde{.}$ accent. However this is proved in \cite[Theorem 2.18]{5KP} as in the proof of the previous proposition.

\end{proof}

%\begin{proposition} \mbox{\bf{(After Berger \cite[Th\'eor\`eme A]{5Ber1})}}
%Let $I$ be empty set, then we have that the category of all the $B_{I'}$-$(\varphi_I,\Gamma_I)$ modules over 
%\begin{center}
%$(B^+_{\mathrm{dR},I'}	\widehat{\otimes}\Pi_{[s_I,r_I],I,J,I\backslash J,A},B^+_{\mathrm{dR},I'}	\widehat{\otimes}\Pi_{[s_I,r_I],I,J,I\backslash J,A}[t_1^{-1},...,t_I'^{-1}],B_{e,I'}	\widehat{\otimes}\Pi_{[s_I,r_I],I,J,I\backslash J,A})$	
%\end{center}
%is equivalent to the category of all the $(\varphi_{I'},\Gamma_{I'})$-modules in the finite projective setting.
%
%
%\end{proposition}
%
%\begin{proof}
%By induction one can prove this by using the strategy of the above.
%	
%\end{proof}
%

\indent Now we combine the construction in the following coherent way following \cite[Section 2]{5Ber1}, \cite[Definition 2.2]{5Nak1} and \cite[Definition 2.2.6]{5KPX}.

\begin{definition}
We define a $B_{I'}$-$(\varphi_I,\Gamma_I)$-bundle over
\begin{center}
$(B^+_{\mathrm{dR},I'}	\widehat{\otimes}\Pi_{\mathrm{an},r_{I},I,J,I\backslash J,A},B^+_{\mathrm{dR},I'}	\widehat{\otimes}\Pi_{[s_I,r_I],I,*,*,A}[t_1^{-1},...,t_I'^{-1}],B_{e,I'}	\widehat{\otimes}\Pi_{\mathrm{an},r_{I},I,*,*,A})$	
\end{center}
to be a compatible family (with respect to the Robba rings) of triplets of finite projective modules:
\begin{align}
(M_e,M_{\mathrm{dR}},M^+_{\mathrm{dR}})	
\end{align}
over 
\begin{center}
$(B^+_{\mathrm{dR},I'}	\widehat{\otimes}\Pi_{\mathrm{an},r_{I},I,*,*,A},B^+_{\mathrm{dR},I'}	\widehat{\otimes}\Pi_{\mathrm{an},r_{I},I,*,*,A}[t_1^{-1},...,t_I'^{-1}],B_{e,I'}	\widehat{\otimes}\Pi_{\mathrm{an},r_{I},I,*,*,A})$	
\end{center}
such that this is glueing datum along :
\begin{align}
B^+_{\mathrm{dR},I'}	\widehat{\otimes}\Pi_{\mathrm{an},r_{I},I,*,*,A}\rightarrow B^+_{\mathrm{dR},I'}	\widehat{\otimes}\Pi_{\mathrm{an},r_{I},I,*,*,A}[t_1^{-1},...,t_I'^{-1}] \leftarrow B_{e,I'}	\widehat{\otimes}\Pi_{\mathrm{an},r_{I},I,*,*,A}.	
\end{align}
And this carry the corresponding relative Galois action of:
\begin{align}
\mathrm{Gal}_{\mathbb{Q}_p,1}\times...\times \mathrm{Gal}_{\mathbb{Q}_p,I'}		
\end{align}
on the multi de Rham period rings which is semilinear. And we have that the three modules involved are relative $(\varphi_I,\Gamma_I)$-bundles relative to
\begin{align}
B^+_{\mathrm{dR},I'},B_{\mathrm{dR},I'},B_{e,I'}.	
\end{align}

\end{definition}

\begin{definition}
We define a pseudocoherent $B_{I'}$-$(\varphi_I,\Gamma_I)$-bundle over
\begin{center}
$(B^+_{\mathrm{dR},I'}	\widehat{\otimes}\Pi_{\mathrm{an},r_{I},I,J,I\backslash J,A},B^+_{\mathrm{dR},I'}	\widehat{\otimes}\Pi_{[s_I,r_I],I,*,*,A}[t_1^{-1},...,t_I'^{-1}],B_{e,I'}	\widehat{\otimes}\Pi_{\mathrm{an},r_{I},I,*,*,A})$	
\end{center}
to be a compatible family of triplets of stably pseudocoherent modules:
\begin{align}
(M_e,M_{\mathrm{dR}},M^+_{\mathrm{dR}})	
\end{align}
over 
\begin{center}
$(B^+_{\mathrm{dR},I'}	\widehat{\otimes}\Pi_{\mathrm{an},r_{I},I,*,*,A},B^+_{\mathrm{dR},I'}	\widehat{\otimes}\Pi_{\mathrm{an},r_{I},I,*,*,A}[t_1^{-1},...,t_I'^{-1}],B_{e,I'}	\widehat{\otimes}\Pi_{\mathrm{an},r_{I},I,*,*,A})$	
\end{center}
such that this is glueing datum along :
\begin{align}
B^+_{\mathrm{dR},I'}	\widehat{\otimes}\Pi_{\mathrm{an},r_{I},I,*,*,A}\rightarrow B^+_{\mathrm{dR},I'}	\widehat{\otimes}\Pi_{\mathrm{an},r_{I},I,*,*,A}[t_1^{-1},...,t_I'^{-1}] \leftarrow B_{e,I'}	\widehat{\otimes}\Pi_{\mathrm{an},r_{I},I,*,*,A}.	
\end{align}
And this carry the corresponding relative Galois action of:
\begin{align}
\mathrm{Gal}_{\mathbb{Q}_p,1}\times...\times \mathrm{Gal}_{\mathbb{Q}_p,I'}		
\end{align}
on the multi de Rham period rings which is semilinear. And we have that the three modules involved are relative pseudocoherent $(\varphi_I,\Gamma_I)$-bundles relative to
\begin{align}
B^+_{\mathrm{dR},I'},B_{\mathrm{dR},I'},B_{e,I'}.	
\end{align}

\end{definition}

\begin{proposition} \mbox{\bf{(After KPX \cite[Proposition 2.2.7]{5KPX})}}
The category of all the finite projective $B_{I'}$-$(\varphi_I,\Gamma_I)$-bundle over 
\begin{center}
$(B^+_{\mathrm{dR},I'}	\widehat{\otimes}\Pi_{\mathrm{an},r_{I},I,*,*,A},B^+_{\mathrm{dR},I'}	\widehat{\otimes}\Pi_{\mathrm{an},r_{I},I,*,*,A}[t_1^{-1},...,t_I'^{-1}],B_{e,I'}	\widehat{\otimes}\Pi_{\mathrm{an},r_{I},I,*,*,A})$	
\end{center}
is equivalent to the category of all the finite projective $B_{I'}$-$(\varphi_I,\Gamma_I)$-modules over 
\begin{center}
$(B^+_{\mathrm{dR},I'}	\widehat{\otimes}\Pi_{\mathrm{an},r_{I},I,*,*,A},B^+_{\mathrm{dR},I'}	\widehat{\otimes}\Pi_{\mathrm{an},r_{I},I,*,*,A}[t_1^{-1},...,t_I'^{-1}],B_{e,I'}	\widehat{\otimes}\Pi_{\mathrm{an},r_{I},I,*,*,A})$.	
\end{center}
\end{proposition}

\begin{proof}
Without considering the corresponding Galois actions for the $B_{I'}$-pair components we could prove this as in the relative situation carrying just $A$-coefficient. To be more precise, the base change gives rise to the corresponding fully faithful functor from the first category to the second one, while to show the corresponding essential surjectivity, consider the corresponding multi-interval $[r_{1,0}/p,r_{1,0}]\times...\times [r_{I,0}/p,r_{I,0}]$ and use the corresponding Frobenius to reach all the corresponding intervals taking the general form of:
\begin{align}
[r_{1,0}/p^{k_1},r_{1,0}/p^{k_1-1}]\times...\times [r_{I,0}/p^{k_I},r_{I,0}/p^{k_I-1}],k_\alpha=1,2,...,\forall\alpha\in I.	
\end{align}
This forms a $2^{|I|}$-uniform covering of the whole space. And the corresponding uniform finiteness of the modules over each 
\begin{align}
[r_{1,0}/p^{k_1},r_{1,0}/p^{k_1-1}]\times...\times [r_{I,0}/p^{k_I},r_{I,0}/p^{k_I-1}],k_\alpha=1,2,...,\forall\alpha\in I.	
\end{align}	
could be achieved by using the corresponding partial Frobenius actions. Then we are done by applying \cref{proposition2.19}. 	
\end{proof}

\begin{proposition} \mbox{\bf{(After KPX \cite[Proposition 2.2.7]{5KPX})}}
The category of all the pseudocoherent $B_{I'}$-$(\varphi_I,\Gamma_I)$-bundle over 
\begin{center}
$(B^+_{\mathrm{dR},I'}	\widehat{\otimes}\Pi_{\mathrm{an},r_{I},I,*,*,A},B^+_{\mathrm{dR},I'}	\widehat{\otimes}\Pi_{\mathrm{an},r_{I},I,*,*,A}[t_1^{-1},...,t_I'^{-1}],B_{e,I'}	\widehat{\otimes}\Pi_{\mathrm{an},r_{I},I,*,*,A})$	
\end{center}
is equivalent to the category of all the pseudocoherent $B_{I'}$-$(\varphi_I,\Gamma_I)$-modules over 
\begin{center}
$(B^+_{\mathrm{dR},I'}	\widehat{\otimes}\Pi_{\mathrm{an},r_{I},I,*,*,A},B^+_{\mathrm{dR},I'}	\widehat{\otimes}\Pi_{\mathrm{an},r_{I},I,*,*,A}[t_1^{-1},...,t_I'^{-1}],B_{e,I'}	\widehat{\otimes}\Pi_{\mathrm{an},r_{I},I,*,*,A})$.	
\end{center}
\end{proposition}

\begin{proof}
Without considering the corresponding Galois actions for the $B_{I'}$-pair components we could prove this as in the relative situation carrying just $A$-coefficient. To be more precise, the base change gives rise to the corresponding fully faithful functor from the first category to the second one, while to show the corresponding essential surjectivity, consider the corresponding multi-interval $[r_{1,0}/p,r_{1,0}]\times...\times [r_{I,0}/p,r_{I,0}]$ and use the corresponding Frobenius to reach all the corresponding intervals taking the general form of:
\begin{align}
[r_{1,0}/p^{k_1},r_{1,0}/p^{k_1-1}]\times...\times [r_{I,0}/p^{k_I},r_{I,0}/p^{k_I-1}],k_\alpha=1,2,...,\forall\alpha\in I.	
\end{align}
This forms a $2^{|I|}$-uniform covering of the whole space. And the corresponding uniform finiteness of the modules over each 
\begin{align}
[r_{1,0}/p^{k_1},r_{1,0}/p^{k_1-1}]\times...\times [r_{I,0}/p^{k_I},r_{I,0}/p^{k_I-1}],k_\alpha=1,2,...,\forall\alpha\in I.	
\end{align}	
could be achieved by using the corresponding partial Frobenius actions. Then we are done by applying \cref{5proposition2.18}. 	
\end{proof}

%%\newpage

\newpage\section{Cohomologies of Cyclotomic Multivariate $(\varphi_I,\Gamma_I)$-Modules over Rigid Analytic Affinoids in Mixed-characteristic Case}

\indent Now we define the corresponding cohomologies of the multivariate $(\varphi_I,\Gamma_I)$-modules over the following groups of rings:

\begin{align}
%\Pi_{\mathrm{an},r_I,I,J,I\backslash J,A},\\	
\Pi_{\mathrm{an},r_I,I,\breve{J},I\backslash J,A}(\pi_{K_I}):=\varprojlim_{s_I}\Pi_{[s_I,r_I],I,\breve{J},I\backslash J,A}(\pi_{K_I}),\\	
\Pi_{\mathrm{an},r_I,I,\widetilde{J},I\backslash J,A}(\pi_{K_I}):=\varprojlim_{s_I} \Pi_{[s_I,r_I],I,\widetilde{J},I\backslash J,A}(\pi_{K_I}),\\
\Pi_{\mathrm{an},r_I,I,J,\breve{I\backslash J},A}(\pi_{K_I}):=\varprojlim_{s_I}\Pi_{[s_I,r_I],I,J,\breve{I\backslash J},A}(\pi_{K_I}),\\	
\Pi_{\mathrm{an},r_I,I,\breve{J},\breve{I\backslash J},A}(\pi_{K_I}):=\varprojlim_{s_I} \Pi_{[s_I,r_I],I,\breve{J},\breve{I\backslash J},A}(\pi_{K_I}),\\	
\Pi_{\mathrm{an},r_I,I,\widetilde{J},\breve{I\backslash J},A}(\pi_{K_I}):=\varprojlim_{s_I} \Pi_{[s_I,r_I],I,\widetilde{J},\breve{I\backslash J},A}(\pi_{K_I}),\\
\Pi_{\mathrm{an},r_I,I,J,\widetilde{I\backslash J},A}(\pi_{K_I}):=\varprojlim_{s_I} \Pi_{[s_I,r_I],I,J,\widetilde{I\backslash J},A}(\pi_{K_I}),\\	
\Pi_{\mathrm{an},r_I,I,\breve{J},\widetilde{I\backslash J},A}(\pi_{K_I}):=\varprojlim_{s_I} \Pi_{[s_I,r_I],I,\breve{J},\widetilde{I\backslash J},A}(\pi_{K_I}),\\	
\Pi_{\mathrm{an},r_I,I,\widetilde{J},\widetilde{I\backslash J},A}(\pi_{K_I}):=\varprojlim_{s_I} \Pi_{[s_I,r_I],I,\widetilde{J},\widetilde{I\backslash J},A}(\pi_{K_I}).	
\end{align}

and

\begin{align}
%\Pi_{[s_I,r_I],I,J,I\backslash J,A},\\	
\Pi_{[s_I,r_I],I,\breve{J},I\backslash J,A}(\Gamma_{K_I}),\\	
\Pi_{[s_I,r_I],I,\widetilde{J},I\backslash J,A}(\Gamma_{K_I}),\\
\Pi_{[s_I,r_I],I,J,\breve{I\backslash J},A}(\Gamma_{K_I}),\\	
\Pi_{[s_I,r_I],I,\breve{J},\breve{I\backslash J},A}(\Gamma_{K_I}),\\
\Pi_{[s_I,r_I],I,\widetilde{J},\breve{I\backslash J},A}(\Gamma_{K_I}),\\
\Pi_{[s_I,r_I],I,J,\widetilde{I\backslash J},A}(\Gamma_{K_I}),\\	
\Pi_{[s_I,r_I],I,\breve{J},\widetilde{I\backslash J},A}(\Gamma_{K_I}),\\	
\Pi_{[s_I,r_I],I,\widetilde{J},\widetilde{I\backslash J},A}(\Gamma_{K_I}).	
\end{align}

\begin{definition} \mbox{\bf{(After KPX, \cite[Definition 2.3.3]{5KPX})}}
We define by induction the corresponding $\varphi_I$-complex $C^\bullet_{\varphi_I}$ of a corresponding $(\varphi_I,\Gamma_I)$-module $M$ over 
\begin{align}
%\Pi_{\mathrm{an},r_I,I,J,I\backslash J,A},\\	
\Pi_{\mathrm{an},r_I,I,\breve{J},I\backslash J,A}(\pi_{K_I}):=\varprojlim_{s_I}\Pi_{[s_I,r_I],I,\breve{J},I\backslash J,A}(\pi_{K_I}),\\	
\Pi_{\mathrm{an},r_I,I,\widetilde{J},I\backslash J,A}(\pi_{K_I}):=\varprojlim_{s_I} \Pi_{[s_I,r_I],I,\widetilde{J},I\backslash J,A}(\pi_{K_I}),\\
\Pi_{\mathrm{an},r_I,I,J,\breve{I\backslash J},A}(\pi_{K_I}):=\varprojlim_{s_I}\Pi_{[s_I,r_I],I,J,\breve{I\backslash J},A}(\pi_{K_I}),\\	
\Pi_{\mathrm{an},r_I,I,\breve{J},\breve{I\backslash J},A}(\pi_{K_I}):=\varprojlim_{s_I} \Pi_{[s_I,r_I],I,\breve{J},\breve{I\backslash J},A}(\pi_{K_I}),\\	
\Pi_{\mathrm{an},r_I,I,\widetilde{J},\breve{I\backslash J},A}(\pi_{K_I}):=\varprojlim_{s_I} \Pi_{[s_I,r_I],I,\widetilde{J},\breve{I\backslash J},A}(\pi_{K_I}),\\
\Pi_{\mathrm{an},r_I,I,J,\widetilde{I\backslash J},A}(\pi_{K_I}):=\varprojlim_{s_I} \Pi_{[s_I,r_I],I,J,\widetilde{I\backslash J},A}(\pi_{K_I}),\\	
\Pi_{\mathrm{an},r_I,I,\breve{J},\widetilde{I\backslash J},A}(\pi_{K_I}):=\varprojlim_{s_I} \Pi_{[s_I,r_I],I,\breve{J},\widetilde{I\backslash J},A}(\pi_{K_I}),\\	
\Pi_{\mathrm{an},r_I,I,\widetilde{J},\widetilde{I\backslash J},A}(\pi_{K_I}):=\varprojlim_{s_I} \Pi_{[s_I,r_I],I,\widetilde{J},\widetilde{I\backslash J},A}(\pi_{K_I})	
\end{align}
to be the corresponding totalization of the following complex:
\[
\xymatrix@C+0pc@R+0pc{
0 \ar[r] \ar[r] \ar[r] &C^\bullet_{\varphi_{I\backslash |I|}}(M_{...,r_{|I|}}) \ar[r]^{\varphi_{|I|}-1}\ar[r]\ar[r] &C^\bullet_{\varphi_{I\backslash |I|}}(M_{...,r_{|I|}/p}) \ar[r] \ar[r] \ar[r] &0
}
\]
as long as $C^\bullet_{\varphi_{I\backslash |I|}}(M)$ is constructed. We define by induction the corresponding $\varphi_I$-complex $C^\bullet_{\varphi_I}(M)$ of a corresponding $(\varphi_I,\Gamma_I)$-module $M$ over
\begin{align}
%\Pi_{[s_I,r_I],I,J,I\backslash J,A},\\	
\Pi_{[s_I,r_I],I,\breve{J},I\backslash J,A}(\pi_{K_I}),\\	
\Pi_{[s_I,r_I],I,\widetilde{J},I\backslash J,A}(\pi_{K_I}),\\
\Pi_{[s_I,r_I],I,J,\breve{I\backslash J},A}(\pi_{K_I}),\\	
\Pi_{[s_I,r_I],I,\breve{J},\breve{I\backslash J},A}(\pi_{K_I}),\\
\Pi_{[s_I,r_I],I,\widetilde{J},\breve{I\backslash J},A}(\pi_{K_I}),\\
\Pi_{[s_I,r_I],I,J,\widetilde{I\backslash J},A}(\pi_{K_I}),\\	
\Pi_{[s_I,r_I],I,\breve{J},\widetilde{I\backslash J},A}(\pi_{K_I}),\\
\Pi_{[s_I,r_I],I,\widetilde{J},\widetilde{I\backslash J},A}(\pi_{K_I})	
\end{align}	
to be the corresponding totalization of the following complex:
\[
\xymatrix@C+0pc@R+0pc{
0 \ar[r] \ar[r] \ar[r] &C^\bullet_{\varphi_{I\backslash |I|}}(M_{...,[s_{|I|},r_{|I|}]}) \ar[r]^{\varphi_{|I|}-1}\ar[r]\ar[r] &C^\bullet_{\varphi_{I\backslash |I|}}(M_{...,[s_{|I|},r_{|I|}/p]}) \ar[r] \ar[r] \ar[r] &0
}
\]
as long as $C^\bullet_{\varphi_{I\backslash |I|}}(M)$ is constructed. Here we assume $0< s_\alpha \leq r_\alpha/p$ for each $\alpha\in I$.	
\end{definition}

\begin{definition} \mbox{\bf{(After KPX, \cite[Definition 2.3.3]{5KPX})}}
We define by induction the corresponding $\psi_I$-complex $C^\bullet_{\psi_I}$ of a corresponding $(\varphi_I,\Gamma_I)$-module $M$ over 
\begin{align}
%\Pi_{\mathrm{an},r_I,I,J,I\backslash J,A},\\	
\Pi_{\mathrm{an},r_I,I,\breve{J},I\backslash J,A}(\pi_{K_I}):=\varprojlim_{s_I}\Pi_{[s_I,r_I],I,\breve{J},I\backslash J,A}(\pi_{K_I}),\\	
\Pi_{\mathrm{an},r_I,I,\widetilde{J},I\backslash J,A}(\pi_{K_I}):=\varprojlim_{s_I} \Pi_{[s_I,r_I],I,\widetilde{J},I\backslash J,A}(\pi_{K_I}),\\
\Pi_{\mathrm{an},r_I,I,J,\breve{I\backslash J},A}(\pi_{K_I}):=\varprojlim_{s_I}\Pi_{[s_I,r_I],I,J,\breve{I\backslash J},A}(\pi_{K_I}),\\	
\Pi_{\mathrm{an},r_I,I,\breve{J},\breve{I\backslash J},A}(\pi_{K_I}):=\varprojlim_{s_I} \Pi_{[s_I,r_I],I,\breve{J},\breve{I\backslash J},A}(\pi_{K_I}),\\	
\Pi_{\mathrm{an},r_I,I,\widetilde{J},\breve{I\backslash J},A}(\pi_{K_I}):=\varprojlim_{s_I} \Pi_{[s_I,r_I],I,\widetilde{J},\breve{I\backslash J},A}(\pi_{K_I}),\\
\Pi_{\mathrm{an},r_I,I,J,\widetilde{I\backslash J},A}(\pi_{K_I}):=\varprojlim_{s_I} \Pi_{[s_I,r_I],I,J,\widetilde{I\backslash J},A}(\pi_{K_I}),\\	
\Pi_{\mathrm{an},r_I,I,\breve{J},\widetilde{I\backslash J},A}(\pi_{K_I}):=\varprojlim_{s_I} \Pi_{[s_I,r_I],I,\breve{J},\widetilde{I\backslash J},A}(\pi_{K_I}),\\	
\Pi_{\mathrm{an},r_I,I,\widetilde{J},\widetilde{I\backslash J},A}(\pi_{K_I}):=\varprojlim_{s_I} \Pi_{[s_I,r_I],I,\widetilde{J},\widetilde{I\backslash J},A}(\pi_{K_I})	
\end{align}
to be the corresponding totalization of the following complex:
\[
\xymatrix@C+0pc@R+0pc{
0 \ar[r] \ar[r] \ar[r] &C^\bullet_{\psi_{I\backslash |I|}}(M_{...,r_{|I|}}) \ar[r]^{\psi_{|I|}-1}\ar[r]\ar[r] &C^\bullet_{\psi_{I\backslash |I|}}(M_{...,pr_{|I|}}) \ar[r] \ar[r] \ar[r] &0
}
\]
as long as $C^\bullet_{\psi_{I\backslash |I|}}(M)$ is constructed. We define by induction the corresponding $\psi_I$-complex $C^\bullet_{\psi_I}(M)$ of a corresponding $(\varphi_I,\Gamma_I)$-module $M$ over
\begin{align}
%\Pi_{[s_I,r_I],I,J,I\backslash J,A},\\	
\Pi_{[s_I,r_I],I,\breve{J},I\backslash J,A}(\pi_{K_I}),\\	
\Pi_{[s_I,r_I],I,\widetilde{J},I\backslash J,A}(\pi_{K_I}),\\
\Pi_{[s_I,r_I],I,J,\breve{I\backslash J},A}(\pi_{K_I}),\\	
\Pi_{[s_I,r_I],I,\breve{J},\breve{I\backslash J},A}(\pi_{K_I}),\\
\Pi_{[s_I,r_I],I,\widetilde{J},\breve{I\backslash J},A}(\pi_{K_I}),\\
\Pi_{[s_I,r_I],I,J,\widetilde{I\backslash J},A}(\pi_{K_I}),\\	
\Pi_{[s_I,r_I],I,\breve{J},\widetilde{I\backslash J},A}(\pi_{K_I}),\\
\Pi_{[s_I,r_I],I,\widetilde{J},\widetilde{I\backslash J},A}(\pi_{K_I})	
\end{align}	
to be the corresponding totalization of the following complex:
\[
\xymatrix@C+0pc@R+0pc{
0 \ar[r] \ar[r] \ar[r] &C^\bullet_{\psi_{I\backslash |I|}}(M_{...,[s_{|I|},r_{|I|}]}) \ar[r]^{\psi_{|I|}-1}\ar[r]\ar[r] &C^\bullet_{\psi_{I\backslash |I|}}(M_{...,[ps_{|I|},r_{|I|}]}) \ar[r] \ar[r] \ar[r] &0
}
\]
as long as $C^\bullet_{\psi_{I\backslash |I|}}(M)$ is constructed. Here we assume $0< s_\alpha \leq r_\alpha/p$ for each $\alpha\in I$.	
\end{definition}

\begin{definition} \mbox{\bf{(After KPX, \cite[Definition 2.3.3]{5KPX})}}
We define by induction the corresponding $\Gamma_I$-complex $C^\bullet_{\Gamma_I}$ of a corresponding $(\varphi_I,\Gamma_I)$-module $M$ over 
\begin{align}
%\Pi_{\mathrm{an},r_I,I,J,I\backslash J,A},\\	
\Pi_{\mathrm{an},r_I,I,\breve{J},I\backslash J,A}(\pi_{K_I}):=\varprojlim_{s_I}\Pi_{[s_I,r_I],I,\breve{J},I\backslash J,A}(\pi_{K_I}),\\	
\Pi_{\mathrm{an},r_I,I,\widetilde{J},I\backslash J,A}(\pi_{K_I}):=\varprojlim_{s_I} \Pi_{[s_I,r_I],I,\widetilde{J},I\backslash J,A}(\pi_{K_I}),\\
\Pi_{\mathrm{an},r_I,I,J,\breve{I\backslash J},A}(\pi_{K_I}):=\varprojlim_{s_I}\Pi_{[s_I,r_I],I,J,\breve{I\backslash J},A}(\pi_{K_I}),\\	
\Pi_{\mathrm{an},r_I,I,\breve{J},\breve{I\backslash J},A}(\pi_{K_I}):=\varprojlim_{s_I} \Pi_{[s_I,r_I],I,\breve{J},\breve{I\backslash J},A}(\pi_{K_I}),\\	
\Pi_{\mathrm{an},r_I,I,\widetilde{J},\breve{I\backslash J},A}(\pi_{K_I}):=\varprojlim_{s_I} \Pi_{[s_I,r_I],I,\widetilde{J},\breve{I\backslash J},A}(\pi_{K_I}),\\
\Pi_{\mathrm{an},r_I,I,J,\widetilde{I\backslash J},A}(\pi_{K_I}):=\varprojlim_{s_I} \Pi_{[s_I,r_I],I,J,\widetilde{I\backslash J},A}(\pi_{K_I}),\\	
\Pi_{\mathrm{an},r_I,I,\breve{J},\widetilde{I\backslash J},A}(\pi_{K_I}):=\varprojlim_{s_I} \Pi_{[s_I,r_I],I,\breve{J},\widetilde{I\backslash J},A}(\pi_{K_I}),\\	
\Pi_{\mathrm{an},r_I,I,\widetilde{J},\widetilde{I\backslash J},A}(\pi_{K_I}):=\varprojlim_{s_I} \Pi_{[s_I,r_I],I,\widetilde{J},\widetilde{I\backslash J},A}(\pi_{K_I})	
\end{align}
to be the corresponding totalization of the following complex:
\[
\xymatrix@C+0pc@R+0pc{
0 \ar[r] \ar[r] \ar[r] &C^\bullet_{\Gamma_{I\backslash |I|}}(M) \ar[r]^{\gamma_{|I|}-1}\ar[r]\ar[r] &C^\bullet_{\Gamma_{I\backslash |I|}}(M) \ar[r] \ar[r] \ar[r] &0
}
\]
as long as $C^\bullet_{\Gamma_{I\backslash |I|}}(M)$ is constructed. We define by induction the corresponding $\Gamma_I$-complex $C^\bullet_{\Gamma_I}(M)$ of a corresponding $(\varphi_I,\Gamma_I)$-module $M$ over
\begin{align}
%\Pi_{[s_I,r_I],I,J,I\backslash J,A},\\	
\Pi_{[s_I,r_I],I,\breve{J},I\backslash J,A}(\pi_{K_I}),\\	
\Pi_{[s_I,r_I],I,\widetilde{J},I\backslash J,A}(\pi_{K_I}),\\
\Pi_{[s_I,r_I],I,J,\breve{I\backslash J},A}(\pi_{K_I}),\\	
\Pi_{[s_I,r_I],I,\breve{J},\breve{I\backslash J},A}(\pi_{K_I}),\\
\Pi_{[s_I,r_I],I,\widetilde{J},\breve{I\backslash J},A}(\pi_{K_I}),\\
\Pi_{[s_I,r_I],I,J,\widetilde{I\backslash J},A}(\pi_{K_I}),\\	
\Pi_{[s_I,r_I],I,\breve{J},\widetilde{I\backslash J},A}(\pi_{K_I}),\\
\Pi_{[s_I,r_I],I,\widetilde{J},\widetilde{I\backslash J},A}(\pi_{K_I})	
\end{align}	
to be the corresponding totalization of the following complex:
\[
\xymatrix@C+0pc@R+0pc{
0 \ar[r] \ar[r] \ar[r] &C^\bullet_{\Gamma_{I\backslash |I|}}(M) \ar[r]^{\gamma_{|I|}-1}\ar[r]\ar[r] &C^\bullet_{\Gamma_{I\backslash |I|}}(M) \ar[r] \ar[r] \ar[r] &0
}
\]
as long as $C^\bullet_{\Gamma_{I\backslash |I|}}(M)$ is constructed. Here we assume $0< s_\alpha \leq r_\alpha/p$ for each $\alpha\in I$.	
\end{definition}

\begin{definition} \mbox{\bf{(After KPX, \cite[Definition 2.3.3]{5KPX})}}
For any $(\varphi_I,\Gamma_I)$-module $M$ over 
\begin{align}
%\Pi_{\mathrm{an},r_I,I,J,I\backslash J,A},\\	
\Pi_{\mathrm{an},r_I,I,\breve{J},I\backslash J,A}(\pi_{K_I}):=\varprojlim_{s_I}\Pi_{[s_I,r_I],I,\breve{J},I\backslash J,A}(\pi_{K_I}),\\	
\Pi_{\mathrm{an},r_I,I,\widetilde{J},I\backslash J,A}(\pi_{K_I}):=\varprojlim_{s_I} \Pi_{[s_I,r_I],I,\widetilde{J},I\backslash J,A}(\pi_{K_I}),\\
\Pi_{\mathrm{an},r_I,I,J,\breve{I\backslash J},A}(\pi_{K_I}):=\varprojlim_{s_I}\Pi_{[s_I,r_I],I,J,\breve{I\backslash J},A}(\pi_{K_I}),\\	
\Pi_{\mathrm{an},r_I,I,\breve{J},\breve{I\backslash J},A}(\pi_{K_I}):=\varprojlim_{s_I} \Pi_{[s_I,r_I],I,\breve{J},\breve{I\backslash J},A}(\pi_{K_I}),\\	
\Pi_{\mathrm{an},r_I,I,\widetilde{J},\breve{I\backslash J},A}(\pi_{K_I}):=\varprojlim_{s_I} \Pi_{[s_I,r_I],I,\widetilde{J},\breve{I\backslash J},A}(\pi_{K_I}),\\
\Pi_{\mathrm{an},r_I,I,J,\widetilde{I\backslash J},A}(\pi_{K_I}):=\varprojlim_{s_I} \Pi_{[s_I,r_I],I,J,\widetilde{I\backslash J},A}(\pi_{K_I}),\\	
\Pi_{\mathrm{an},r_I,I,\breve{J},\widetilde{I\backslash J},A}(\pi_{K_I}):=\varprojlim_{s_I} \Pi_{[s_I,r_I],I,\breve{J},\widetilde{I\backslash J},A}(\pi_{K_I}),\\	
\Pi_{\mathrm{an},r_I,I,\widetilde{J},\widetilde{I\backslash J},A}(\pi_{K_I}):=\varprojlim_{s_I} \Pi_{[s_I,r_I],I,\widetilde{J},\widetilde{I\backslash J},A}(\pi_{K_I})	
\end{align} 
we define the corresponding complex $C^\bullet_{\varphi_I,\Gamma_I}(M)$ to be the corresponding totalization of $C^\bullet_{\varphi_I}C^\bullet_{\Gamma_I}(M)$.
For any $(\varphi_I,\Gamma_I)$-module $M$ over 
\begin{align}
%\Pi_{[s_I,r_I],I,J,I\backslash J,A},\\	
\Pi_{[s_I,r_I],I,\breve{J},I\backslash J,A}(\pi_{K_I}),\\	
\Pi_{[s_I,r_I],I,\widetilde{J},I\backslash J,A}(\pi_{K_I}),\\
\Pi_{[s_I,r_I],I,J,\breve{I\backslash J},A}(\pi_{K_I}),\\	
\Pi_{[s_I,r_I],I,\breve{J},\breve{I\backslash J},A}(\pi_{K_I}),\\
\Pi_{[s_I,r_I],I,\widetilde{J},\breve{I\backslash J},A}(\pi_{K_I}),\\
\Pi_{[s_I,r_I],I,J,\widetilde{I\backslash J},A}(\pi_{K_I}),\\	
\Pi_{[s_I,r_I],I,\breve{J},\widetilde{I\backslash J},A}(\pi_{K_I}),\\
\Pi_{[s_I,r_I],I,\widetilde{J},\widetilde{I\backslash J},A}(\pi_{K_I})	
\end{align}	
we define the corresponding complex $C^\bullet_{\varphi_I,\Gamma_I}(M)$ to be the corresponding totalization of $C^\bullet_{\varphi_I}C^\bullet_{\Gamma_I}(M)$. Here we assume $0< s_\alpha \leq r_\alpha/p$ for each $\alpha\in I$.	
	
\end{definition}

%%%%%%%%%%%%%%%%%%%%%%%%%%%%%%%%%!!!!!!!!!!!!!!!

\begin{definition} \mbox{\bf{(After KPX, \cite[Definition 2.3.3]{5KPX})}}
For any $(\psi_I,\Gamma_I)$-module $M$ over 
\begin{align}
%\Pi_{\mathrm{an},r_I,I,J,I\backslash J,A},\\	
\Pi_{\mathrm{an},r_I,I,\breve{J},I\backslash J,A}(\pi_{K_I}):=\varprojlim_{s_I}\Pi_{[s_I,r_I],I,\breve{J},I\backslash J,A}(\pi_{K_I}),\\	
\Pi_{\mathrm{an},r_I,I,\widetilde{J},I\backslash J,A}(\pi_{K_I}):=\varprojlim_{s_I} \Pi_{[s_I,r_I],I,\widetilde{J},I\backslash J,A}(\pi_{K_I}),\\
\Pi_{\mathrm{an},r_I,I,J,\breve{I\backslash J},A}(\pi_{K_I}):=\varprojlim_{s_I}\Pi_{[s_I,r_I],I,J,\breve{I\backslash J},A}(\pi_{K_I}),\\	
\Pi_{\mathrm{an},r_I,I,\breve{J},\breve{I\backslash J},A}(\pi_{K_I}):=\varprojlim_{s_I} \Pi_{[s_I,r_I],I,\breve{J},\breve{I\backslash J},A}(\pi_{K_I}),\\	
\Pi_{\mathrm{an},r_I,I,\widetilde{J},\breve{I\backslash J},A}(\pi_{K_I}):=\varprojlim_{s_I} \Pi_{[s_I,r_I],I,\widetilde{J},\breve{I\backslash J},A}(\pi_{K_I}),\\
\Pi_{\mathrm{an},r_I,I,J,\widetilde{I\backslash J},A}(\pi_{K_I}):=\varprojlim_{s_I} \Pi_{[s_I,r_I],I,J,\widetilde{I\backslash J},A}(\pi_{K_I}),\\	
\Pi_{\mathrm{an},r_I,I,\breve{J},\widetilde{I\backslash J},A}(\pi_{K_I}):=\varprojlim_{s_I} \Pi_{[s_I,r_I],I,\breve{J},\widetilde{I\backslash J},A}(\pi_{K_I}),\\	
\Pi_{\mathrm{an},r_I,I,\widetilde{J},\widetilde{I\backslash J},A}(\pi_{K_I}):=\varprojlim_{s_I} \Pi_{[s_I,r_I],I,\widetilde{J},\widetilde{I\backslash J},A}(\pi_{K_I})	
\end{align} 
we define the corresponding complex $C^\bullet_{\psi_I,\Gamma_I}(M)$ to be the corresponding totalization of $C^\bullet_{\psi_I}C^\bullet_{\Gamma_I}(M)$.
For any $(\varphi_I,\Gamma_I)$-module $M$ over 
\begin{align}
%\Pi_{[s_I,r_I],I,J,I\backslash J,A},\\	
\Pi_{[s_I,r_I],I,\breve{J},I\backslash J,A}(\pi_{K_I}),\\	
\Pi_{[s_I,r_I],I,\widetilde{J},I\backslash J,A}(\pi_{K_I}),\\
\Pi_{[s_I,r_I],I,J,\breve{I\backslash J},A}(\pi_{K_I}),\\	
\Pi_{[s_I,r_I],I,\breve{J},\breve{I\backslash J},A}(\pi_{K_I}),\\
\Pi_{[s_I,r_I],I,\widetilde{J},\breve{I\backslash J},A}(\pi_{K_I}),\\
\Pi_{[s_I,r_I],I,J,\widetilde{I\backslash J},A}(\pi_{K_I}),\\	
\Pi_{[s_I,r_I],I,\breve{J},\widetilde{I\backslash J},A}(\pi_{K_I}),\\
\Pi_{[s_I,r_I],I,\widetilde{J},\widetilde{I\backslash J},A}(\pi_{K_I})	
\end{align}	
we define the corresponding complex $C^\bullet_{\psi_I,\Gamma_I}(M)$ to be the corresponding totalization of $C^\bullet_{\psi_I}C^\bullet_{\Gamma_I}(M)$. Here we assume $0< s_\alpha \leq r_\alpha/p$ for each $\alpha\in I$.	
	
\end{definition}

%%\newpage

\newpage\section{Cohomologies of $B_I$-pairs and Mixed-Type Objects over Rigid Analytic Affinoids in Mixed-characterist
ic Case}

\subsection{Partial $(\varphi_I,\Gamma_I)$-Cohomology and Partial 
$(\psi_I,\Gamma_I)$
-Cohomol
ogy}

\indent Now we define the corresponding cohomologies of the multivariate $(\varphi_I,\Gamma_I)$-modules over the following two groups of rings:

\begin{align}
\varprojlim_{s_I}  B^+_{\mathrm{dR},I'}	\widehat{\otimes}\Pi_{[s_I,r_I],I,?,?',A},?=J,\widetilde{J},\breve{J},?'=I\backslash J,\widetilde{I\backslash J},\breve{I\backslash J}\\
\end{align}

with

\begin{align}
\varprojlim_{s_I} B^+_{\mathrm{dR},I'}	\widehat{\otimes}\Pi_{[s_I,r_I],I,?,?',A}[t_1^{-1},...,t_{I'}^{-1}],?=J,\widetilde{J},\breve{J},?'=I\backslash J,\widetilde{I\backslash J},\breve{I\backslash J},\\
\end{align}

with 

\begin{align}
\varprojlim_{s_I} B_{e,I'}	\widehat{\otimes}\Pi_{[s_I,r_I],I,?,?',A}, ?=J,\widetilde{J},\breve{J},?'=I\backslash J,\widetilde{I\backslash J},\breve{I\backslash J},\\
\end{align}

and

\begin{align}
B^+_{\mathrm{dR},I'}	\widehat{\otimes}\Pi_{[s_I,r_I],I,?,?',A},?=J,\widetilde{J},\breve{J},?'=I\backslash J,\widetilde{I\backslash J},\breve{I\backslash J}
\end{align}

with

\begin{align}
B^+_{\mathrm{dR},I'}	\widehat{\otimes}\Pi_{[s_I,r_I],I,?,?',A}[t_1^{-1},...,t_{I'}^{-1}],?=J,\widetilde{J},\breve{J},?'=I\backslash J,\widetilde{I\backslash J},\breve{I\backslash J}
\end{align}

with 

\begin{align}
B_{e,I'}	\widehat{\otimes}\Pi_{[s_I,r_I],I,?,?',A},?=J,\widetilde{J},\breve{J},?'=I\backslash J,\widetilde{I\backslash J},\breve{I\backslash J}.
\end{align}

\begin{definition} \mbox{\bf{(After KPX, \cite[Definition 2.3.3]{5KPX})}}
We define by induction the corresponding $\varphi_I$-complex $C^\bullet_{\varphi_I}$ of a corresponding $B_{I'}$-$(\varphi_I,\Gamma_I)$-module $M$ over the first three groups of rings
to be the corresponding totalization of the following complex:
\[
\xymatrix@C+0pc@R+0pc{
0 \ar[r] \ar[r] \ar[r] &C^\bullet_{\varphi_{I\backslash |I|}}(M_{...,r_{|I|}}) \ar[r]^{\varphi_{|I|}-1}\ar[r]\ar[r] &C^\bullet_{\varphi_{I\backslash |I|}}(M_{...,r_{|I|}/p}) \ar[r] \ar[r] \ar[r] &0
}
\]
as long as $C^\bullet_{\varphi_{I\backslash |I|}}(M)$ is constructed. We define by induction the corresponding $\varphi_I$-complex $C^\bullet_{\varphi_I}(M)$ of a corresponding $B_{I'}$-$(\varphi_I,\Gamma_I)$-module $M$ over the second three groups of rings to be the corresponding totalization of the following complex:
\[
\xymatrix@C+0pc@R+0pc{
0 \ar[r] \ar[r] \ar[r] &C^\bullet_{\varphi_{I\backslash |I|}}(M_{...,[s_{|I|},r_{|I|}]}) \ar[r]^{\varphi_{|I|}-1}\ar[r]\ar[r] &C^\bullet_{\varphi_{I\backslash |I|}}(M_{...,[s_{|I|},r_{|I|}/p]}) \ar[r] \ar[r] \ar[r] &0
}
\]
as long as $C^\bullet_{\varphi_{I\backslash |I|}}(M)$ is constructed. Here we assume $0< s_\alpha \leq r_\alpha/p$ for each $\alpha\in I$.	
\end{definition}

\begin{definition} \mbox{\bf{(After KPX, \cite[Definition 2.3.3]{5KPX})}}
We define by induction the corresponding $\psi_I$-complex $C^\bullet_{\psi_I}$ of a corresponding $B_{I'}$-$(\varphi_I,\Gamma_I)$-module $M$ over the first three groups of rings 
to be the corresponding totalization of the following complex:
\[
\xymatrix@C+0pc@R+0pc{
0 \ar[r] \ar[r] \ar[r] &C^\bullet_{\psi_{I\backslash |I|}}(M_{...,r_{|I|}}) \ar[r]^{\psi_{|I|}-1}\ar[r]\ar[r] &C^\bullet_{\psi_{I\backslash |I|}}(M_{...,pr_{|I|}}) \ar[r] \ar[r] \ar[r] &0
}
\]
as long as $C^\bullet_{\psi_{I\backslash |I|}}(M)$ is constructed. We define by induction the corresponding $\psi_I$-complex $C^\bullet_{\psi_I}(M)$ of a corresponding $B_{I'}$-$(\varphi_I,\Gamma_I)$-module $M$ over the second three groups of rings
to be the corresponding totalization of the following complex:
\[
\xymatrix@C+0pc@R+0pc{
0 \ar[r] \ar[r] \ar[r] &C^\bullet_{\psi_{I\backslash |I|}}(M_{...,[s_{|I|},r_{|I|}]}) \ar[r]^{\psi_{|I|}-1}\ar[r]\ar[r] &C^\bullet_{\psi_{I\backslash |I|}}(M_{...,[ps_{|I|},r_{|I|}]}) \ar[r] \ar[r] \ar[r] &0
}
\]
as long as $C^\bullet_{\psi_{I\backslash |I|}}(M)$ is constructed. Here we assume $0< s_\alpha \leq r_\alpha/p$ for each $\alpha\in I$.	
\end{definition}

\begin{definition} \mbox{\bf{(After KPX, \cite[Definition 2.3.3]{5KPX})}}
We define by induction the corresponding $\Gamma_I$-complex $C^\bullet_{\Gamma_I}$ of a corresponding $B_{I'}$-$(\varphi_I,\Gamma_I)$-module $M$ over 
the first three groups of rings
to be the corresponding totalization of the following complex:
\[
\xymatrix@C+0pc@R+0pc{
0 \ar[r] \ar[r] \ar[r] &C^\bullet_{\Gamma_{I\backslash |I|}}(M) \ar[r]^{\gamma_{|I|}-1}\ar[r]\ar[r] &C^\bullet_{\Gamma_{I\backslash |I|}}(M) \ar[r] \ar[r] \ar[r] &0
}
\]
as long as $C^\bullet_{\Gamma_{I\backslash |I|}}(M)$ is constructed. We define by induction the corresponding $\Gamma_I$-complex $C^\bullet_{\Gamma_I}(M)$ of a corresponding $B_{I'}$-$(\varphi_I,\Gamma_I)$-module $M$ over
the second three groups of rings 	
to be the corresponding totalization of the following complex:
\[
\xymatrix@C+0pc@R+0pc{
0 \ar[r] \ar[r] \ar[r] &C^\bullet_{\Gamma_{I\backslash |I|}}(M) \ar[r]^{\gamma_{|I|}-1}\ar[r]\ar[r] &C^\bullet_{\Gamma_{I\backslash |I|}}(M) \ar[r] \ar[r] \ar[r] &0
}
\]
as long as $C^\bullet_{\Gamma_{I\backslash |I|}}(M)$ is constructed. Here we assume $0< s_\alpha \leq r_\alpha/p$ for each $\alpha\in I$.	
\end{definition}

\begin{definition} \mbox{\bf{(After KPX, \cite[Definition 2.3.3]{5KPX})}}\\
For any $B_{I'}$-$(\varphi_I,\Gamma_I)$-module $M$ over 
the first three groups of rings 
we define the corresponding complex $C^\bullet_{\varphi_I,\Gamma_I}(M)$ to be the corresponding totalization of $C^\bullet_{\varphi_I}C^\bullet_{\Gamma_I}(M)$.
For any $(\varphi_I,\Gamma_I)$-module $M$ over 
the second three groups of rings
we define the corresponding complex $C^\bullet_{\varphi_I,\Gamma_I}(M)$ to be the corresponding totalization of $C^\bullet_{\varphi_I}C^\bullet_{\Gamma_I}(M)$. Here we assume $0< s_\alpha \leq r_\alpha/p$ for each $\alpha\in I$.	
	
\end{definition}

\begin{definition} \mbox{\bf{(After KPX, \cite[Definition 2.3.3]{5KPX})}}\\
For any $B_{I'}$-$(\varphi_I,\Gamma_I)$-module $M$ over 
the first three groups of rings
we define the corresponding complex $C^\bullet_{\psi_I,\Gamma_I}(M)$ to be the corresponding totalization of $C^\bullet_{\psi_I}C^\bullet_{\Gamma_I}(M)$.
For any $B_{I'}$-$(\varphi_I,\Gamma_I)$-module $M$ over 
the second three groups of rings
we define the corresponding complex $C^\bullet_{\psi_I,\Gamma_I}(M)$ to be the corresponding totalization of $C^\bullet_{\psi_I}C^\bullet_{\Gamma_I}(M)$. Here we assume $0< s_\alpha \leq r_\alpha/p$ for each $\alpha\in I$.	
	
\end{definition}

\subsection{Partial $B_{I'}$-Cohomology}

\indent We now define the corresponding partial $B_{I'}$-cohomology  in the situation

\begin{definition}  \mbox{\bf{(After Nakamura, \cite[Appendix 5]{5Nak1})}} 
For any\\ $B_{I'}$-$(\varphi_I,\Gamma_I)$-module $M$ over 
\begin{align}
\varprojlim_{s_I}  B^+_{\mathrm{dR},I'}	\widehat{\otimes}\Pi_{[s_I,r_I],I,?,?',A},?=J,\widetilde{J},\breve{J},?'=I\backslash J,\widetilde{I\backslash J},\breve{I\backslash J},\\
\end{align}

with

\begin{align}
\varprojlim_{s_I} B^+_{\mathrm{dR},I'}	\widehat{\otimes}\Pi_{[s_I,r_I],I,?,?',A}[t_1^{-1},...,t_{I'}^{-1}],?=J,\widetilde{J},\breve{J},?'=I\backslash J,\widetilde{I\backslash J},\breve{I\backslash J},\\
\end{align}

with 

\begin{align}
\varprojlim_{s_I} B_{e,I'}	\widehat{\otimes}\Pi_{[s_I,r_I],I,?,?',A},?=J,\widetilde{J},\breve{J},?'=I\backslash J,\widetilde{I\backslash J},\breve{I\backslash J},\\
\end{align}

and

\begin{align}
B^+_{\mathrm{dR},I'}	\widehat{\otimes}\Pi_{[s_I,r_I],I,J,I\backslash J,A},\\	
B^+_{\mathrm{dR},I'}	\widehat{\otimes}\Pi_{[s_I,r_I],I,\breve{J},I\backslash J,A},\\	
B^+_{\mathrm{dR},I'}	\widehat{\otimes}\Pi_{[s_I,r_I],I,\widetilde{J},I\backslash J,A},\\
B^+_{\mathrm{dR},I'}	\widehat{\otimes}\Pi_{[s_I,r_I],I,J,\breve{I\backslash J},A},\\	
B^+_{\mathrm{dR},I'}	\widehat{\otimes}\Pi_{[s_I,r_I],I,\breve{J},\breve{I\backslash J},A},\\
B^+_{\mathrm{dR},I'}	\widehat{\otimes}\Pi_{[s_I,r_I],I,\widetilde{J},\breve{I\backslash J},A},\\
B^+_{\mathrm{dR},I'}	\widehat{\otimes}\Pi_{[s_I,r_I],I,J,\widetilde{I\backslash J},A},\\	
B^+_{\mathrm{dR},I'}	\widehat{\otimes}\Pi_{[s_I,r_I],I,\breve{J},\widetilde{I\backslash J},A},\\
B^+_{\mathrm{dR},I'}	\widehat{\otimes}\Pi_{[s_I,r_I],I,\widetilde{J},\widetilde{I\backslash J},A}
\end{align}

with

\begin{align}
B^+_{\mathrm{dR},I'}	\widehat{\otimes}\Pi_{[s_I,r_I],I,J,I\backslash J,A}[t_1^{-1},...,t_{I'}^{-1}],\\	
B^+_{\mathrm{dR},I'}	\widehat{\otimes}\Pi_{[s_I,r_I],I,\breve{J},I\backslash J,A}[t_1^{-1},...,t_{I'}^{-1}],\\	
B^+_{\mathrm{dR},I'}	\widehat{\otimes}\Pi_{[s_I,r_I],I,\widetilde{J},I\backslash J,A}[t_1^{-1},...,t_{I'}^{-1}],\\
B^+_{\mathrm{dR},I'}	\widehat{\otimes}\Pi_{[s_I,r_I],I,J,\breve{I\backslash J},A}[t_1^{-1},...,t_{I'}^{-1}],\\	
B^+_{\mathrm{dR},I'}	\widehat{\otimes}\Pi_{[s_I,r_I],I,\breve{J},\breve{I\backslash J},A}[t_1^{-1},...,t_{I'}^{-1}],\\
B^+_{\mathrm{dR},I'}	\widehat{\otimes}\Pi_{[s_I,r_I],I,\widetilde{J},\breve{I\backslash J},A}[t_1^{-1},...,t_{I'}^{-1}],\\
B^+_{\mathrm{dR},I'}	\widehat{\otimes}\Pi_{[s_I,r_I],I,J,\widetilde{I\backslash J},A}[t_1^{-1},...,t_{I'}^{-1}],\\	
B^+_{\mathrm{dR},I'}	\widehat{\otimes}\Pi_{[s_I,r_I],I,\breve{J},\widetilde{I\backslash J},A}[t_1^{-1},...,t_{I'}^{-1}],\\
B^+_{\mathrm{dR},I'}	\widehat{\otimes}\Pi_{[s_I,r_I],I,\widetilde{J},\widetilde{I\backslash J},A}[t_1^{-1},...,t_{I'}^{-1}],
\end{align}

with 

\begin{align}
B_{e,I'}	\widehat{\otimes}\Pi_{[s_I,r_I],I,J,I\backslash J,A},\\	
B_{e,I'}	\widehat{\otimes}\Pi_{[s_I,r_I],I,\breve{J},I\backslash J,A},\\	
B_{e,I'}	\widehat{\otimes}\Pi_{[s_I,r_I],I,\widetilde{J},I\backslash J,A},\\
B_{e,I'}	\widehat{\otimes}\Pi_{[s_I,r_I],I,J,\breve{I\backslash J},A},\\	
B_{e,I'}	\widehat{\otimes}\Pi_{[s_I,r_I],I,\breve{J},\breve{I\backslash J},A},\\
B_{e,I'}	\widehat{\otimes}\Pi_{[s_I,r_I],I,\widetilde{J},\breve{I\backslash J},A},\\
B_{e,I'}	\widehat{\otimes}\Pi_{[s_I,r_I],I,J,\widetilde{I\backslash J},A},\\	
B_{e,I'}	\widehat{\otimes}\Pi_{[s_I,r_I],I,\breve{J},\widetilde{I\backslash J},A},\\
B_{e,I'} \widehat{\otimes}\Pi_{[s_I,r_I],I,\widetilde{J},\widetilde{I\backslash J},A}
\end{align}
we define the corresponding $B_{I'}$-complex $C^\bullet_{B_{I'}}(M)$ to be 
\[
\xymatrix@C+0pc@R+0pc{
0 \ar[r] \ar[r] \ar[r] &{*}_1 \times {*}_2  \ar[r]^{?-?'}\ar[r]\ar[r] &{*}_3\ar[r] \ar[r] \ar[r] &0,
}
\]
where
\begin{align}
*_1:=C^\bullet(\mathrm{Gal}_{\mathbb{Q}_p,1}\times...\times \mathrm{Gal}_{\mathbb{Q}_p,I'}, M^+_\mathrm{dR}),\\
*_2:=C^\bullet(\mathrm{Gal}_{\mathbb{Q}_p,1}\times...\times \mathrm{Gal}_{\mathbb{Q}_p,I'}, M_e),\\
*_3:=C^\bullet(\mathrm{Gal}_{\mathbb{Q}_p,1}\times...\times \mathrm{Gal}_{\mathbb{Q}_p,I'}, M_\mathrm{dR}).	
\end{align}

\end{definition}

\begin{definition}
Now we consider any $B_{I'}$-$(\varphi_I,\Gamma_I)$-module $M$ over 

\begin{align}
\varprojlim_{s_I}  B^+_{\mathrm{dR},I'}	\widehat{\otimes}\Pi_{[s_I,r_I],I,?,?',A},?=J,\widetilde{J},\breve{J},?'=I\backslash J,\widetilde{I\backslash J},\breve{I\backslash J}\\
\end{align}

with

\begin{align}
\varprojlim_{s_I} B^+_{\mathrm{dR},I'}	\widehat{\otimes}\Pi_{[s_I,r_I],I,?,?',A}[t_1^{-1},...,t_{I'}^{-1}],?=J,\widetilde{J},\breve{J},?'=I\backslash J,\widetilde{I\backslash J},\breve{I\backslash J},\\
\end{align}

with 

\begin{align}
\varprojlim_{s_I} B_{e,I'}	\widehat{\otimes}\Pi_{[s_I,r_I],I,?,?',A}, ?=J,\widetilde{J},\breve{J},?'=I\backslash J,\widetilde{I\backslash J},\breve{I\backslash J},\\
\end{align}

and

\begin{align}
B^+_{\mathrm{dR},I'}	\widehat{\otimes}\Pi_{[s_I,r_I],I,?,?',A},?=J,\widetilde{J},\breve{J},?'=I\backslash J,\widetilde{I\backslash J},\breve{I\backslash J}
\end{align}

with

\begin{align}
B^+_{\mathrm{dR},I'}	\widehat{\otimes}\Pi_{[s_I,r_I],I,?,?',A}[t_1^{-1},...,t_{I'}^{-1}],?=J,\widetilde{J},\breve{J},?'=I\backslash J,\widetilde{I\backslash J},\breve{I\backslash J}
\end{align}

with 

\begin{align}
B_{e,I'}	\widehat{\otimes}\Pi_{[s_I,r_I],I,?,?',A},?=J,\widetilde{J},\breve{J},?'=I\backslash J,\widetilde{I\backslash J},\breve{I\backslash J}
\end{align}

Then we can define the corresponding $C^\bullet_{B_{I'},\varphi_I,\Gamma_I}$-cohomology complex by taking the corresponding totalization of the corresponding double complex $C^\bullet_{B_{I'}}C^\bullet_{\varphi_I,\Gamma_I}(M)$.
\end{definition}

%%\newpage

\newpage\section{The Results on the Cohomologies}

\subsection{Comparisons for $C^\bullet_{\varphi_I,\Gamma_I},C^\bullet_{\psi_I,\Gamma_I},C^\bullet_{\psi_I}$}

\indent We now consider the following categories for $A$ a rigid affinoid over $\mathbb{Q}_p$:\\

\noindent A. The corresponding category of all the $(\varphi_I,\Gamma_I)$-modules over the corresponding rings carrying the corresponding cohomologies $C^\bullet_{(\varphi_I,\Gamma_I)},C^\bullet_{(\psi_I,\Gamma_I)},C^\bullet_{\psi_I}$ (for sufficiently small $r_{I,0}$): \\

\begin{align}
%\Pi_{\mathrm{an},r_I,I,J,I\backslash J,A},\\	
\Pi_{\mathrm{an},r_{I,0},I,\breve{J},I\backslash J,A}(\pi_{K_I}):=\varprojlim_{s_I}\Pi_{[s_I,r_{I,0}],I,\breve{J},I\backslash J,A}(\pi_{K_I}),\\	
\Pi_{\mathrm{an},r_{I,0},I,\widetilde{J},I\backslash J,A}(\pi_{K_I}):=\varprojlim_{s_I} \Pi_{[s_I,r_{I,0}],I,\widetilde{J},I\backslash J,A}(\pi_{K_I}),\\
\Pi_{\mathrm{an},r_{I,0},I,J,\breve{I\backslash J},A}(\pi_{K_I}):=\varprojlim_{s_I}\Pi_{[s_I,r_{I,0}],I,J,\breve{I\backslash J},A}(\pi_{K_I}),\\	
\Pi_{\mathrm{an},r_{I,0},I,\breve{J},\breve{I\backslash J},A}(\pi_{K_I}):=\varprojlim_{s_I} \Pi_{[s_I,r_{I,0}],I,\breve{J},\breve{I\backslash J},A}(\pi_{K_I}),\\	
\Pi_{\mathrm{an},r_{I,0},I,\widetilde{J},\breve{I\backslash J},A}(\pi_{K_I}):=\varprojlim_{s_I} \Pi_{[s_I,r_{I,0}],I,\widetilde{J},\breve{I\backslash J},A}(\pi_{K_I}),\\
\Pi_{\mathrm{an},r_{I,0},I,J,\widetilde{I\backslash J},A}(\pi_{K_I}):=\varprojlim_{s_I} \Pi_{[s_I,r_{I,0}],I,J,\widetilde{I\backslash J},A}(\pi_{K_I}),\\	
\Pi_{\mathrm{an},r_{I,0},I,\breve{J},\widetilde{I\backslash J},A}(\pi_{K_I}):=\varprojlim_{s_I} \Pi_{[s_I,r_{I,0}],I,\breve{J},\widetilde{I\backslash J},A}(\pi_{K_I}),\\	
\Pi_{\mathrm{an},r_{I,0},I,\widetilde{J},\widetilde{I\backslash J},A}(\pi_{K_I}):=\varprojlim_{s_I} \Pi_{[s_I,r_{I,0}],I,\widetilde{J},\widetilde{I\backslash J},A}(\pi_{K_I}).	
\end{align}

\noindent B. The corresponding category of all the $(\varphi_I,\Gamma_I)$-modules over the corresponding rings carrying the corresponding cohomologies $C^\bullet_{(\varphi_I,\Gamma_I)},C^\bullet_{(\psi_I,\Gamma_I)},C^\bullet_{\psi_I}$ (where $0<s_\alpha\leq r_\alpha/p\leq  r_{\alpha,0}$ for each $\alpha\in I$): \\

\begin{align}
%\Pi_{[s_I,r_I],I,J,I\backslash J,A},\\	
\Pi_{[s_I,r_I],I,\breve{J},I\backslash J,A}(\pi_{K_I}),\\	
\Pi_{[s_I,r_I],I,\widetilde{J},I\backslash J,A}(\pi_{K_I}),\\
\Pi_{[s_I,r_I],I,J,\breve{I\backslash J},A}(\pi_{K_I}),\\	
\Pi_{[s_I,r_I],I,\breve{J},\breve{I\backslash J},A}(\pi_{K_I}),\\	
\Pi_{[s_I,r_I],I,\widetilde{J},\breve{I\backslash J},A}(\pi_{K_I}),\\
\Pi_{[s_I,r_I],I,J,\widetilde{I\backslash J},A}(\pi_{K_I}),\\	
\Pi_{[s_I,r_I],I,\breve{J},\widetilde{I\backslash J},A}(\pi_{K_I}),\\	
\Pi_{[s_I,r_I],I,\widetilde{J},\widetilde{I\backslash J},A}(\pi_{K_I}).	
\end{align}

\begin{theorem}
Let $I$ be a set of two elements. Let $M$ be some object over\\
 $\Pi_{\mathrm{an},r_{I,0},I,*,*,A}(\pi_{K_I})$ in the corresponding category $A$ and let $M_{[s_I,r_I]}$ be the corrresponding (under the horizontal equivalence of the categories for the rings of the same type) object over $\Pi_{[s_I,r_{I}],I,*,*,A}(\pi_{K_I})$ in the category $B$. Then we have the following quasi-isomorphisms:
\begin{align}
C^\bullet_{\varphi_I,\Gamma_I}(M)\overset{\sim}{\rightarrow}C^\bullet_{\varphi_I,\Gamma_I}(M_{[s_I,r_I]}), \\
C^\bullet_{\psi_I,\Gamma_I}(M)\overset{\sim}{\rightarrow}C^\bullet_{\psi_I,\Gamma_I}(M_{[s_I,r_I]}), \\
C^\bullet_{\psi_I}(M)\overset{\sim}{\rightarrow}C^\bullet_{\psi_I}(M_{[s_I,r_I]}), \\	
\end{align}
in the bounded derived category of $A$-modules $D^\flat(A)$.	
\end{theorem}

\begin{proof}
The proof relies on the following intermediate $(\varphi_I,\Gamma_I)$-modules. Let $I=\{1,2\}$, we now defined the following rings:
\begin{align}
%\Pi_{\mathrm{an},r_I,I,J,I\backslash J,A},\\	
\Pi_{\mathrm{an},r_{1,0},[s_2,r_2],I,\breve{J},I\backslash J,A}(\pi_{K_I}):=\varprojlim_{s_1}\Pi_{[s_1,r_{1,0}]\times[s_2,r_{2}],I,\breve{J},I\backslash J,A}(\pi_{K_I}),\\	
\Pi_{\mathrm{an},r_{1,0},[s_2,r_2],I,\widetilde{J},I\backslash J,A}(\pi_{K_I}):=\varprojlim_{s_1} \Pi_{[s_1,r_{1,0}]\times[s_2,r_{2}],I,\widetilde{J},I\backslash J,A}(\pi_{K_I}),\\
\Pi_{\mathrm{an},r_{1,0},[s_2,r_2],I,J,\breve{I\backslash J},A}(\pi_{K_I}):=\varprojlim_{s_1}\Pi_{[s_1,r_{1,0}]\times[s_2,r_{2}],I,J,\breve{I\backslash J},A}(\pi_{K_I}),\\	
\Pi_{\mathrm{an},r_{1,0},[s_2,r_2],I,\breve{J},\breve{I\backslash J},A}(\pi_{K_I}):=\varprojlim_{s_1} \Pi_{[s_1,r_{1,0}]\times[s_2,r_{2}],I,\breve{J},\breve{I\backslash J},A}(\pi_{K_I}),\\	
\Pi_{\mathrm{an},r_{1,0},[s_2,r_2],I,\widetilde{J},\breve{I\backslash J},A}(\pi_{K_I}):=\varprojlim_{s_1} \Pi_{[s_1,r_{1,0}]\times[s_2,r_{2}],I,\widetilde{J},\breve{I\backslash J},A}(\pi_{K_I}),\\
\Pi_{\mathrm{an},r_{1,0},[s_2,r_2],I,J,\widetilde{I\backslash J},A}(\pi_{K_I}):=\varprojlim_{s_1} \Pi_{[s_1,r_{1,0}]\times[s_2,r_{2}],I,J,\widetilde{I\backslash J},A}(\pi_{K_I}),\\	
\Pi_{\mathrm{an},r_{1,0},[s_2,r_2],I,\breve{J},\widetilde{I\backslash J},A}(\pi_{K_I}):=\varprojlim_{s_1} \Pi_{[s_1,r_{1,0}]\times[s_2,r_{2}],I,\breve{J},\widetilde{I\backslash J},A}(\pi_{K_I}),\\	
\Pi_{\mathrm{an},r_{1,0},[s_2,r_2],I,\widetilde{J},\widetilde{I\backslash J},A}(\pi_{K_I}):=\varprojlim_{s_1} \Pi_{[s_1,r_{1,0}]\times[s_2,r_{2}],I,\widetilde{J},\widetilde{I\backslash J},A}(\pi_{K_I}).	
\end{align}
This category serves as a corresponding intermediate category which factors through the corresponding original equivalence on the categories involved. From $A$ to this category we have the equivalence on the cohomology groups by considering the $(\varphi_2,\Gamma_2)$-module structure, then by regarding the corresponding $(\varphi_2,\Gamma_2)$-cohomology groups as Yoneda extension groups we have the isomorphism on $(\varphi_2,\Gamma_2)$-cohomology groups which further shows the equivalence on the full cohomology groups. Similarly one compares this category with the category $B$ to finish.

\end{proof}

\

This chapter is based on the following paper, where the author of this dissertation is the main author:
\begin{itemize}
\item Tong, Xin. "Analytic Geometry and Hodge-Frobenius Structure Continued." arXiv preprint arXiv:2012.07336 (2020).
\end{itemize}

\newpage\chapter{Category and Cohomology of Hodge-Iwasawa Modules}

\newpage\section{Introduction}

\subsection{The Main Scope of the Discussion}

\indent The corresponding Hodge-Iwasawa theory was studied in our papers \cite{6XT1} and \cite{6XT2}. We in the corresponding deformed setting discussed the corresponding relative $p$-adic Hodge theory after \cite{6KP}, \cite{6KL15} and \cite{6KL16}. We call the theory partially Iwasawa in the corresponding sense that the corresponding deformation will reflect the corresponding Iwasawa theoretic information which was observed as in \cite{6KP}. On the Hodge-theoretic side, the corresponding cohomology in the deformed setting will also reflect the corresponding towers related to quotients (even nonabelian) of fundamental groups.\\

\indent Therefore in our current consideration in this paper, we apply all we have developed to the corresponding categories and the corresponding cohomologies. For instance, very interesting application would be the corresponding pro-\'etale cohomology of general local systems over some smooth proper rigid analytic spaces as in \cite{6KL3}. Definitely we expect more applications in the corresponding Iwasawa theoretic consideration. \\

\indent The other interesting things we would like to pursue in the application of our previously developed theory is the corresponding categorical and $K$-theoretic study on the corresponding $(\varphi,\Gamma)$-modules over very relative Robba rings. The corresponding categories should form some well-established abelian ones which could further have more derived enrichment (carrying some six functors).\\

\indent The second goal of this paper is also to consider the Hodge-Structure on more general spaces. Although we will not systematically consider some general framework, since it will impossible to handle very general adic spaces. To be more general than the corresponding rigid analytic spaces, we are going to consider the corresponding spaces admitting atlas of the so-called $k_\Delta$ affinoid spaces as in \cite{6TC} and \cite{6DFN}, where $\Delta$ is some valuation group larger than $|k^\times|$.\\

\indent We briefly mention now the corresponding structures of the paper. In the second section we are going to study the corresponding Hodge-Iwasawa modules over logarithmic towers which generalizes for instance in the mixed-characteristic case the corresponding consideration in \cite[Chapter 7]{6KL16}. Then in the third section we study after that the corresponding categories and cohomologies of the corresponding Hodge-Iwasawa modules involved. Then we are going to study the corresponding cohomology of more general analytic spaces. In fourth section we considere the Hodge-Iwasawa module over rigid analytic spaces, $k_\Delta$-analytic spaces and we contact the corresponding arithmetic Riemann-Hilbert correspondence after \cite{6LZ}. The fifth section is for the corresponding equivariant Iwasawa theory of de Rham $(\varphi,\Gamma)$-modules which is a equivariant version of the corresponding Iwasawa theory in \cite{6Nakamura1}. In such context, we still emphasize the corresponding geometrization aspect along \cite{6KP} on the reconstruction effect of the corresponding Iwasawa deformation in different contexts. Note that the story is equivariant and relative inspired by Nakamura \cite{6Nakamura3}, Kedlaya-Pottharst \cite{6KP}, as well as some further motivation from \cite{6BF1}, \cite{6BF2} and \cite{6FK}.\\

\subsection{Results Involved}

\indent We searched many interesting directions in this paper to apply the Hodge-Iwasawa consideration in \cite{6XT1} and \cite{6XT2}. We present in some form in this current introduction part of this paper some results along the main body of this paper below. First motivation comes from directly following \cite{6KL16} where the relative $p$-adic Hodge theory is essentially applied in the corresponding context of rigid analytic spaces in the strictly situation. \cite{6KL16} proved many very deep results around the corresponding categories of pseudocoherent relative $(\varphi,\Gamma)$-modules, namely \cite{6KL16} proved that these are abelian categories although working in Banach categories does cause essential difficulties. In our situation we seek the same goal when we deform the corresponding Frobenius modules in arithmetic family. In the context and notation in \cref{proposition6.3} we have:

\begin{proposition} \mbox{\bf{(After Kedlaya-Liu, \cite[Theorem 8.10.6]{6KL16})}} 
Working over period ring $\widetilde{\Pi}_{X,A}$, we have the corresponding category of $(\varphi^a,\Gamma)$-modules is abelian, which is basically compatible with the corresponding category of all the sheaves of modules over $\widetilde{\Pi}_{X,A}$ when one would like to form the corresponding kernels and cokernels. $A$ is assumed to be sousperfectoid as those considered in \cite{6KH}.
\end{proposition}

\indent Then following Temkin's notes \cite[Part I Chapter I]{6DFN} we consider the corresponding analytic spaces which could admit an atlas consisting of all the $k_\Delta$-affinoids namely those quotients of $k_\Delta$ strictly rational localizations of Tate algebras. Over such space $X$: 

\begin{proposition} \mbox{\bf{(After Kedlaya-Liu, \cite[Theorem 8.10.6]{6KL16})}} 
Working over period ring $\widetilde{\Pi}_{X,A}$, we have the corresponding category of $(\varphi^a,\Gamma)$-modules is abelian, which is basically compatible with the corresponding category of all the sheaves of modules over $\widetilde{\Pi}_{X,A}$ when one would like to form the corresponding kernel and cokernel. Here $A$ is $\mathbb{Q}_p$.
\end{proposition}

\begin{proposition} \mbox{\bf{(After Kedlaya-Liu, \cite[Theorem 8.10.6]{6KL16})}} 
Working over period ring $\widetilde{\Pi}_{X,A}$, we have the corresponding category of $(\varphi^a,\Gamma)$-modules is abelian, which is basically compatible with the corresponding category of all the sheaves of modules over $\widetilde{\Pi}_{X,A}$ when one would like to form the corresponding kernel and cokernel. Here $X$ is assumed to be smooth. $A$ is assumed to be sousperfectoid as those considered in \cite{6KH}.\\
\end{proposition}

\indent Then we consider arithmetic family version of Higgs bundles after \cite{6LZ} where a version of Riemann-Hilbert correspondence in the arithmetic setting is achieved. The idea in our situation will be definitely the corresponding version of such correspondence for general $A$-relative $(\varphi,\Gamma)$-modules $M$ in the 'geometric' setting. We also consider the corresponding equivalent version of $B$-pairs in the $A$-relative setting after \cite{6KL16}. With the notation in later discussion (see \cref{proposition4.45} and \cref{proposition4.42}) we have for a de Rham $A$-relative $(\varphi,\Gamma)$-module $M$ or equivalently a de Rham $A$-relative $B$-pair $M$:

\begin{proposition} \mbox{\bf{(After Kedlaya-Liu, \cite[below Definition 10.10]{6KL3})}}
Over a general rigid analyic space $X$, consider a de Rham $A$-relative $(\varphi,\Gamma)$-modules $M$ or equivalently a de Rham $A$-relative $B$-pair $M$. We then have that the corresponding each higher de Rham derived cohomology group $D^i_\mathrm{dR}(M)$ for each $i\geq 0$ is coherent sheaf over $X$. The projectivity could achieved if we assume that $X$ is smooth. $A$ is assumed to be sousperfectoid as those considered in \cite{6KH}.

\end{proposition}

\begin{proposition} \mbox{\bf{(After Kedlaya-Liu, \cite[below Definition 10.10]{6KL3})}}
Over a general rigid analyic space $X$, consider a de Rham $A$-relative $(\varphi,\Gamma)$-modules $M$ or equivalently a de Rham $A$-relative $B$-pair $M$. We then have that the corresponding each higher de Rham derived cohomology group $E^i_\mathrm{dR}(M)$ for each $i\geq 0$ is coherent sheaf over $X$. The projectivity could achieved if we assume that $X$ is smooth. $A$ is assumed to be sousperfectoid as those considered in \cite{6KH}.\\

\end{proposition}

\indent The other corresponding application we want to search for while again along some idea proposed in \cite{6KP} applications in equivariant Iwasawa theory. We choose to consider equivariant version of the Iwasawa theory of de Rham $(\varphi,\Gamma)$-module $M$. 

\begin{remark}
We believe (while motivated by \cite{6KP}) our idea should eventually produce more than cyclotomic consideration by deforming Berger's $(\varphi,\nabla,\Gamma)$-modules over pro-\'etale site which could absorb the $\Gamma$-action (which could be reconstructed through the Iwasawa deformation).
\end{remark}

\begin{definition}\mbox{}\\ \mbox{\bf{(After Perrin-Riou and Nakamura, \cite[Definition 3.7]{6Nakamura1})}} 
With the notation in \cref{6section5.1} especially the corresponding abelian Fr\'echet-Stein space we are working on, we have the following derived big Perrin-Riou-Nakamura Exponential map:
\begin{displaymath}
\mathrm{Exp}^{\Pi^\infty(\Gamma'),\bullet}_{\mathcal{F}(M),m}: R^\bullet\Gamma_\mathrm{sheaf}(\varprojlim\mathcal{F}(\Theta_{\mathrm{Ber,dif}}(M))_p) \rightarrow R^\bullet\Gamma_\mathrm{sheaf}(\varprojlim\mathcal{F}(M)_p),
\end{displaymath}
which is corresponding to the one for $(\varphi,\Gamma)$-modules:
\begin{displaymath}
\mathrm{Exp}^{\Pi^\infty(\Gamma'),\bullet}_{M,m}: R^\bullet\Gamma_{\varphi,\Gamma}(\varprojlim \Theta_{\mathrm{Ber,dif}}(M)_p) \rightarrow R^\bullet\Gamma_{\varphi,\Gamma}(\varprojlim M_p).
\end{displaymath}
\end{definition}

\indent We actually conjectured that the corresponding this equivariant\\ Perrin-Riou-Nakamura map in algebraic way produces the equivariant $p$-adic $L$-functions (in some not necessarily explicit way)  which relate directly to the corresponding characteristic ideals. See \cref{section5.2}.\\

\subsection{Future Study}

There are many interesting possible extensions of the current paper. We would like to study more in the future, although here we mainly emphasize the corresponding Iwasawa theoretic aspects. The first scope of topics we would like to study further is the relative $p$-adic Hodge theory in rigid family and the arithmetic Riemann-Hilbert correspondence in rigid family for $B$-pairs after \cite{6LZ}, \cite{6TT} and \cite{6Shi}, along the corresponding foundation we established here.\\

The second scope of the topics we would like to study in the future is to extend the corresponding discussion to more general analytic spaces. In fact here we do not go beyond \cite[Chapter 8]{6KL16}, although we touched the corresponding $k_\Delta$-analytic spaces as those considered in Temkin's lecture notes \cite[Part I Chapter I]{6DFN}. Namely we still restrict ourselves to the corresponding varieties covered by strictly or non-strictly affinoids. \\

The Iwasawa consideration we considered here is actually also worthwhile to be amplified and to be noncommutativized. Certainly we want to study noncommutative deformation of the sheaves in the future not only in all kinds of Iwasawa theories, but also we would like to study the noncommutative deformation of the corresponding the relative $p$-adic Hodge theory in rigid family and the arithmetic Riemann-Hilbert correspondence in rigid family for $B$-pairs after \cite{6LZ}, \cite{6TT} and \cite{6Shi} in the sense of many existing noncommutative geometric contexts.

%%\newpage

\newpage\section{Logarithmic Hodge-Iwasawa Structures in Mixed-Characteristic Case}\label{section3}

\subsection{$\Gamma$-Modules over Ramified Towers}

\begin{setting}
Let $E$ be $\mathbb{Q}_p$. We use the corresponding $\pi$ to denote the corresponding uniformizer in a uniform notation. Let $A$ be an affinoid algebra over the corresponding field $E$ in our current situation defined through some strict quotient from the corresponding Tate algebras.
\end{setting}

\indent We first consider the corresponding standard ones for $(P,P^+)$ a perfectoid adic Banach uniform pair over $\mathcal{O}_E$, containing the perfectoid field $\mathbb{Q}_p(\zeta_{p^\infty})^\wedge$.

\begin{setting}
In the mixed characteristic situation, we consider the setting in the following fashion:
\begin{align}
H_0&=P\{x_1/t_1,...,x_k/t_k,t_1/y_1,...,t_\ell/y_\ell\},\\
H_0^+&=P^+\{x_1/t_1,...,x_k/t_k,t_1/y_1,...,t_\ell/y_\ell\},\\
H_n&=P\{(x_1/t_1)^{1/p^n},...,(x_k/t_k)^{1/p^n},(t_1/y_1)^{1/p^n},...,(t_\ell/y_\ell)^{1/p^n}\},\\
n&=0,1,..., 0<k\leq\ell.	
\end{align}		
And here the corresponding quantities $x_i,y_j,i=1,...,k,j=1,...,\ell$ could be allowed to be real numbers.
\end{setting}

\begin{lemma} \mbox{\bf{(Kedlaya-Liu \cite[Lemma 7.3.3]{6KL16})}}
The corresponding standard ramified toric towers defined as above are weakly decompleting.	
\end{lemma}

\begin{proof}
See \cite[Lemma 7.3.3]{6KL16}.	
\end{proof}

\indent The corresponding tower is not finite \'etale, so the corresponding machinery in \cite[Chapter 5]{6KL16} needs to be modified, which is already done in the corresponding context in \cite[Chapter 7]{6KL16}. As in \cite[Chapter 7, Section 3, Section 4]{6KL16} one can also consider direct semilinear $\Gamma$-action and define the corresponding $\Gamma$-modules to be certain ones carrying such semilinear action.

\indent  When we in some situations have the corresponding topological group $\Gamma$ (just as in the corresponding logarithmic setting we are considering here) we can simply consider the corresponding definition as in \cite[Definition 7.3.5]{6KL16}:

\begin{definition} \mbox{\bf{(After Kedlaya-Liu \cite[Definition 7.3.5]{6KL16})}}
In our current situation, we consider the corresponding action of the corresponding group $\Gamma$ in our situation, namely the $\ell$-fold product of the corresponding additive group $\mathbb{Z}_p$ on the corresponding period rings taking the form of $*_{H,A}$. We now define the corresponding $\Gamma$-modules in the finite projective, pseudocoherent or fpd setting to the corresponding modules over the corresponding period rings as above carrying the corresponding action coming from the group $\Gamma$ which is assumed to as in \cite[Definition 7.3.5]{6KL16} particularly semilinear. 	
\end{definition}

\begin{remark}
Again the deformation only happens over the rings in mixed-characteristic situation.	
\end{remark}

\begin{proposition} \mbox{\bf{(Kedlaya-Liu \cite[Lemma 7.3.6]{6KL16})}} \label{proposition4.8}
Let $\Gamma_n:=p^n\Gamma$ for $n\geq 0$, then we have in our situation the corresponding result that the corresponding complex $C^\bullet(\Gamma_n,\overline{\varphi}^{-1}R_H/R_H)$ is strict exact.
\end{proposition}

\begin{proof}
See \cite[Lemma 7.3.6]{6KL16}.
\end{proof}

\begin{proposition} \mbox{\bf{(Kedlaya-Liu \cite[Lemma 5.6.4]{6KL16})}} \label{6proposition4.9}
Let $\Gamma_n:=p^n\Gamma$ for $n\geq 0$, then we have in our situation the corresponding result that the corresponding complex $C^\bullet(\Gamma_n,\overline{R}_H/R_H)$ is strict exact.
\end{proposition}

\begin{proof}
See \cite[Lemma 5.6.4]{6KL16}. 
\end{proof}

\begin{lemma} \mbox{\bf{(After Kedlaya-Liu \cite[Corollary 5.6.5]{6KL16})}} \label{lemma4.10}
We need to fix some suitable $r_0>0$ as in \cite[Corollary 5.6.5]{6KL16} and for the corresponding radii $s,r$ in $(0,r_0]$ we will have the corresponding complexes $C^\bullet(\Gamma_n,\varphi^{-1}\Pi^{\mathrm{int},r/p}_{H}/\Pi^{\mathrm{int},r}_{H})_A$ and $C^\bullet(\Gamma_n,\widetilde{\Pi}^{\mathrm{int},r}_{H}/\Pi^{\mathrm{int},r}_{H})_A$ are strict exact. Moreover we have that the corresponding complexes\\
 $C^\bullet(\Gamma_n,\varphi^{-1}\Pi^{[s/p,r/p]}_{H}/\Pi^{[s,r]}_{H})_A$ and $C^\bullet(\Gamma_n,\widetilde{\Pi}^{[s/p,r/p]}_{H}/\Pi^{[s,r]}_{H})_A$ are then in our situation strict exact as well.  	
\end{lemma}

\begin{proof}
See \cite[Corollary 5.6.5]{6KL16}.	
\end{proof}

\begin{lemma} \mbox{\bf{(After Kedlaya-Liu \cite[Lemma 5.6.6]{6KL16})}} \label{lemma4.11}
Let $M$ be a $\Gamma$-module defined over the period rings involved below. Then one can find a sufficiently large integer $\ell^*\geq 0$ such that we have that for any $\ell\geq \ell^*$, the following complexes are correspondingly strict exact: 
\begin{align}
&M\otimes_{R_H} C^\bullet(\Gamma_n,\overline{\varphi}^{-\ell-1}R_H/\overline{\varphi}^{-\ell}R_H), M\otimes_{R_H} C^\bullet(\Gamma_n,\overline{R}_H/\overline{\varphi}^{-\ell}R_H),\\
&M\otimes_{\Pi^{\mathrm{int},r}_{H,A}} C^\bullet(\Gamma_n,\varphi^{-(\ell+1)}\Pi^{\mathrm{int},r/p^{\ell+1}}_{H}/\varphi^{-\ell}\Pi^{\mathrm{int},r/p^\ell}_{H})_A, M\otimes_{\Pi^{\mathrm{int},r}_{H}} C^\bullet(\Gamma_n,\widetilde{\Pi}^{\mathrm{int},r}_{H}/\varphi^{-\ell}\Pi^{\mathrm{int},r/p^\ell}_{H})_A\\ 
&M\otimes_{\Pi_{H,A}} C^\bullet(\Gamma_n,\varphi^{-(\ell+1)}\Pi^{[s/p^{\ell+1},r/p^{\ell+1}]}_{H}/\varphi^{-\ell}\Pi^{[s/p^{\ell},r/p^{\ell}]}_{H})_A, \\
&M\otimes_{\Pi^{\mathrm{int},r}_{H}} C^\bullet(\Gamma_n,\widetilde{\Pi}_{H}/\varphi^{-\ell}\Pi^{[s/p^{\ell},r/p^{\ell}]}_{H})_A.
\end{align}

\end{lemma}

\begin{proof}
This is the corresponding corollary of \cref{proposition4.8}, \cref{6proposition4.9} and \cref{lemma4.10}, along the corresponding argument in \cite[Lemma 5.6.6]{6KL16}. 	
\end{proof}

\begin{theorem}\mbox{\bf{(After Kedlaya-Liu \cite[Lemma 5.6.9]{6KL16})}} \label{6theorem4.12}
The corresponding base change functor from the ring $\breve{R}_H$ to $\widetilde{R}_H$ is in our situation an equivalence between the corresponding $\Gamma$-modules.	 The corresponding base change functor from the ring $\breve{\Pi}^{\mathrm{int}}_{H,A}$ to $\widetilde{\Pi}^\mathrm{int}_{H,A}$ is in our situation an equivalence between the corresponding $\Gamma$-modules. The corresponding base change functor from the ring $\breve{\Pi}_{H,A}$ to $\widetilde{\Pi}_{H,A}$ is in our situation an equivalence between the corresponding $\Gamma$-modules.
\end{theorem}

\begin{proof}
We briefly mention the corresponding argument since this is just analog of the corresponding argument for the proof of the corresponding result \cite[Lemma 5.6.9]{6KL16}. The corresponding fully faithfulness comes from the \cref{lemma3.10}. Then we apply the corresponding \cite[Lemma 5.6.8]{6KL16} to descend the corresponding modules over the bigger rings with accent $\widetilde{*}$ to the base rings after some Frobenius action. Then the corresponding differentials on the original modules over the bigger rings could be made have norms just a small quantity twist away from the corresponding differentials on the corresponding base changes from the corresponding base rings after Frobenius twists. Then to finish we apply the corresponding \label{lemma3.10} to modify this quantity to finish the proof as in \cite[Lemma 5.6.9]{6KL16}.
\end{proof}

\indent For the relative logarithmic toric towers, we consider the corresponding following mixed characteristic situation. We consider the following tower $(H,H^+)$:
\begin{align}
H_0&=\mathbb{Q}_p\{x_1/t_1,...,x_k/t_k,t_1/y_1,...,t_\ell/y_\ell\},\\
H_0^+&=H_0^\circ,\\
H_n&=\mathbb{Q}_p(\zeta_{p^n})\{(x_1/t_1)^{1/p^n},...,(x_k/t_k)^{1/p^n},(t_1/y_1)^{1/p^n},...,(t_k/y_k)^{1/p^n}\}, n=0,1,....	
\end{align}

The corresponding parameters $x_i,y_i$ are in $p^\mathbb{Q}$, and $\ell\geq k\geq 0$. The corresponding tower is again not Galois. Again when we in some situations have the corresponding topological group $\Gamma$ (just as in the corresponding logarithmic setting we are considering here) we can simply consider the corresponding definition as in \cite[Definition 7.3.5]{6KL16}:

\begin{definition} \mbox{\bf{(After Kedlaya-Liu \cite[Definition 7.3.5]{6KL16})}}
In our current situation, we consider the corresponding action of the corresponding group $\Gamma$ in our situation, namely the semidirect product taking the form of $\mathbb{Z}^\times_p\ltimes \mathbb{Z}_p^\ell$ on the corresponding period rings taking the form of $*_{H,A}$. We now define the corresponding $\Gamma$-modules in the finite projective, pseudocoherent or fpd setting to the corresponding modules over the corresponding period rings as above carrying the corresponding action coming from the group $\Gamma$ which is assumed to as in \cite[Definition 7.3.5]{6KL16} particularly semilinear. 	
\end{definition}

\begin{proposition} \mbox{\bf{(After Kedlaya-Liu \cite[Lemma 7.3.6]{6KL16})}} 
Let $\Gamma_n:=p^n\Gamma$ for $n\geq 0$, then we have in our situation the corresponding result that the corresponding complex $C^\bullet(\Gamma_n,\overline{\varphi}^{-1}R_H/R_H)$ is strict exact.
\end{proposition}

\begin{proof}
See \cref{proposition4.8}.
\end{proof}

\begin{proposition} \mbox{\bf{(After Kedlaya-Liu \cite[Lemma 5.6.4]{6KL16})}} 
Let $\Gamma_n:=p^n\Gamma$ for $n\geq 0$, then we have in our situation the corresponding result that the corresponding complex $C^\bullet(\Gamma_n,\overline{R}_H/R_H)$ is strict exact.
\end{proposition}

\begin{proof}
See \cref{6proposition4.9}.
\end{proof}

\begin{lemma} \mbox{\bf{(After Kedlaya-Liu \cite[Corollary 5.6.5]{6KL16})}} 
We need to fix some suitable $r_0>0$ as in \cite[Corollary 5.6.5]{6KL16} and for the corresponding radii $s,r$ in $(0,r_0]$ we will have the corresponding complexes $C^\bullet(\Gamma_n,\varphi^{-1}\Pi^{\mathrm{int},r/p}_{H}/\Pi^{\mathrm{int},r}_{H})_A$ and $C^\bullet(\Gamma_n,\widetilde{\Pi}^{\mathrm{int},r}_{H}/\Pi^{\mathrm{int},r}_{H})_A$ are strict exact. Moreover we have that the corresponding complexes\\ $C^\bullet(\Gamma_n,\varphi^{-1}\Pi^{[s/p,r/p]}_{H}/\Pi^{[s,r]}_{H})_A$ and $C^\bullet(\Gamma_n,\widetilde{\Pi}^{[s/p,r/p]}_{H}/\Pi^{[s,r]}_{H})_A$ are then in our situation strict exact as well.  	
\end{lemma}

\begin{proof}
See \cref{lemma4.10}.	
\end{proof}

\begin{lemma} \mbox{\bf{(After Kedlaya-Liu \cite[Lemma 5.6.6]{6KL16})}} 
Let $M$ be a $\Gamma$-module defined over the period rings involved below. Then one can find a sufficiently large integer $\ell^*\geq 0$ such that we have that for any $\ell\geq \ell^*$, the following complexes are correspondingly strict exact: 
\begin{align}
&M\otimes_{R_H} C^\bullet(\Gamma_n,\overline{\varphi}^{-\ell-1}R_H/\overline{\varphi}^{-\ell}R_H), M\otimes_{R_H} C^\bullet(\Gamma_n,\overline{R}_H/\overline{\varphi}^{-\ell}R_H),\\
&M\otimes_{\Pi^{\mathrm{int},r}_{H,A}} C^\bullet(\Gamma_n,\varphi^{-(\ell+1)}\Pi^{\mathrm{int},r/p^{\ell+1}}_{H}/\varphi^{-\ell}\Pi^{\mathrm{int},r/p^\ell}_{H})_A, M\otimes_{\Pi^{\mathrm{int},r}_{H}} C^\bullet(\Gamma_n,\widetilde{\Pi}^{\mathrm{int},r}_{H}/\varphi^{-\ell}\Pi^{\mathrm{int},r/p^\ell}_{H})_A\\ 
&M\otimes_{\Pi_{H,A}} C^\bullet(\Gamma_n,\varphi^{-(\ell+1)}\Pi^{[s/p^{\ell+1},r/p^{\ell+1}]}_{H}/\varphi^{-\ell}\Pi^{[s/p^{\ell},r/p^{\ell}]}_{H})_A,\\
&M\otimes_{\Pi^{\mathrm{int},r}_{H}} C^\bullet(\Gamma_n,\widetilde{\Pi}_{H}/\varphi^{-\ell}\Pi^{[s/p^{\ell},r/p^{\ell}]}_{H})_A.
\end{align}

\end{lemma}

\begin{proof}
See \cref{lemma4.11}.	
\end{proof}

\begin{theorem}\mbox{\bf{(After Kedlaya-Liu \cite[Lemma 5.6.9]{6KL16})}} 
The corresponding base change functor from the ring $\breve{R}_H$ to $\widetilde{R}_H$ is in our situation an equivalence between the corresponding $\Gamma$-modules.	 The corresponding base change functor from the ring $\breve{\Pi}^{\mathrm{int}}_{H,A}$ to $\widetilde{\Pi}^\mathrm{int}_{H,A}$ is in our situation an equivalence between the corresponding $\Gamma$-modules. The corresponding base change functor from the ring $\breve{\Pi}_{H,A}$ to $\widetilde{\Pi}_{H,A}$ is in our situation an equivalence between the corresponding $\Gamma$-modules.
\end{theorem}

\begin{proof}
See \cref{6theorem4.12}.
\end{proof}

\begin{remark} \label{remark2.17}
The corresponding statements for the mixed characteristic version of the relative log towers are established in \cite[Chapter 7.3, Chapter 7.4]{6KL16} without touching the corresponding large coefficient $A$ here, therefore for more detailed discussion see \cite[Chapter 7.3, Chapter 7.4]{6KL16}. We should mention that the corresponding \cite[Chapter 7.3, Chapter 7.4]{6KL16} target at the corresponding finite projective objects, essentially in the generality above we should come across essentially similar difficulties to possible ones in \cite{6KL16} when one would like to consider the pseudocoherent sheaves. Things in our current section are basically the complementary discussion to the discussion we made in \cite{6XT2} for the unramified towers.\\  	
\end{remark}

%%%%%%%%%%%%%%%%%%%%%%%%%%%%%%%%%%!!!!!!!!!!!!!!!!!!!!!!

\subsection{$(\varphi,\Gamma)$-Modules}

\indent Now we study the corresponding $\Gamma$-modules carrying further Frobenius action, but we will again mainly focus on the corresponding finite projective objects (see \cref{remark2.17}), which is parallel to the corresponding consideration in \cite{6KL16}. First recall the following definition:

\begin{setting} \mbox{\bf{(After Kedlaya-Liu \cite[Definition 7.5.3]{6KL16})}} \label{6setting2.18}
We will consider the corresponding Kummer \'etale topology and the corresponding pro-Kummer \'etale topology considered in \cite[Definition 7.5.3]{6KL16}. Let $(H,H^+)$ be either the standard ramified toric tower or the relative ramified toric tower we considered as obave. Here we recall that a covering of the adic space $\mathrm{Spa}(H_0,H_0^+)$ is called \'etale if it is so after considering the corresponding pullback along for some integer $v\geq 0$ the map $T_1\mapsto T_1^{p^{-v}},...,T_\ell\mapsto T_\ell^{p^{-v}}$. And we recall that a covering of the adic space $\mathrm{Spa}(H_0,H_0^+)$ is called finite \'etale if it is so after considering the corresponding pullback along for some integer $v\geq 0$ the map $T_1\mapsto T_1^{p^{-v}},...,T_\ell\mapsto T_\ell^{p^{-v}}$. And we recall that a covering of the adic space $\mathrm{Spa}(H_0,H_0^+)$ is called faithfully finite \'etale if it is so after considering the corresponding pullback along for some integer $v\geq 0$ the map $T_1\mapsto T_1^{p^{-v}},...,T_\ell\mapsto T_\ell^{p^{-v}}$.	
\end{setting}

\begin{definition}\mbox{\bf{(After Kedlaya-Liu \cite[Definition 5.7.2]{6KL16})}}
Consider the corresponding framework in the previous subsection. Consider the corresponding period rings in our situation (with the corresponding notations in \cite{6XT1} and \cite{6XT2}):
\begin{align}
\Omega^{\mathrm{int}}_{H,A},\breve{\Omega}_{H,A}^{\mathrm{int}},\widehat{\Omega}_{H,A}^{\mathrm{int}},\widetilde{\Omega}_{H,A}^{\mathrm{int}},\\
\Pi^{\mathrm{int}}_{H,A},\breve{\Pi}_{H,A}^{\mathrm{int}},\widehat{\Pi}_{H,A}^{\mathrm{int}},\widetilde{\Pi}_{H,A}^{\mathrm{int}},\\
\Pi^{\mathrm{bd}}_{H,A},\breve{\Pi}_{H,A}^{\mathrm{bd}},\widehat{\Pi}_{H,A}^{\mathrm{bd}},\widetilde{\Pi}_{H,A}^{\mathrm{bd}},\\
\Pi_{H,A},\breve{\Pi}_{H,A},\widetilde{\Pi}_{H,A},
\end{align}
we have the corresponding notions of finite projective or pseudocoherent or finite projective dimension $\Gamma$-modules, where we can then add the corresponding semilinear actions as in \cite[Definition 5.7.2]{6KL16} to define the corresponding $(\varphi,\Gamma)$-modules satisfying the corresponding condition on the corresponding Frobenius pullbacks. Then furthermore we can define the same objects for those rings having some radius $r>0$ of interval $[s,r]\in (0,\infty)$ as in \cite[Definition 5.7.2]{6KL16}, again satisfying the corresponding Frobenius pullback condition taking the general form:
\begin{align}
&\varphi^*\Delta \otimes_{*_{H,A}^{r/p}} *_{{H,A}}^{r/p} \overset{\sim}{\rightarrow} \Delta \otimes *_{{H,A}}^{r/p},\\
&\varphi^*\Delta \otimes_{*_{H,A}^{[s/p,r/p]}} *_{H,A}^{[s,r/p]} \overset{\sim}{\rightarrow} \Delta \otimes_{*_{H,A}^{[s,r]}} *_{H,A}^{[s,r/p]}.	
\end{align}
%We need to consider the corresponding topological conditions on the modules over the corresponding Robba rings $\Pi^r_{H,A},\breve{\Pi}^r_{H,A},\widetilde{\Pi}^r_{H,A},\Pi^{[s,r]}_{H,A},\breve{\Pi}^{[s,r]}_{H,A},\widetilde{\Pi}^{[s,r]}_{H,A}$. Namely all the modules at least will be complete for the natural topology inherited from the corresponding topological rings. Moreover modules over $\Pi^r_{H,A},\breve{\Pi}^r_{H,A},\widetilde{\Pi}^r_{H,A}$ will be assumed to base change to modules which are Kummer \'etale-stably pseudocoherent over the rings $\Pi^{[s,r]}_{H,A},\breve{\Pi}^{[s,r]}_{H,A},\widetilde{\Pi}^{[s,r]}_{H,A}$.
\end{definition}

%\begin{proposition} \mbox{\bf{(After Kedlaya-Liu \cite[Theorem 5.7.3]{6KL16})}}
%The categories of $(\varphi,\Gamma)$-modules (namely finite projective ones) over the rings:
%\begin{align}
%\Omega^{\mathrm{int}}_H,\breve{\Omega}_H^{\mathrm{int}},\widehat{\Omega}_H^{\mathrm{int}},\widetilde{\Omega}_H^{\mathrm{int}},\\
%\Pi^{\mathrm{int}}_H,\breve{\Pi}_H^{\mathrm{int}},\widehat{\Pi}_H^{\mathrm{int}},\widetilde{\Pi}_H^{\mathrm{int}}
%\end{align}	
%are equivalent to each other. And the categories of $(\varphi,\Gamma)$-modules (namely finite projective ones) over the rings:
%\begin{align}
%\Pi^{\mathrm{int}}_{H,A},\breve{\Pi}_{H,A}^{\mathrm{int}},\widehat{\Pi}_{H,A}^{\mathrm{int}},\widetilde{\Pi}_{H,A}^{\mathrm{int}}
%\end{align}	
%are equivalent to each other.
%\end{proposition}
%
%
%\begin{proof}
%See \cite[Theorem 5.7.3]{6KL16}.	
%\end{proof}

\begin{proposition} \mbox{\bf{(After Kedlaya-Liu \cite[Theorem 5.7.4]{6KL16})}} \label{proposition2.20}
Consider the following period rings in the rational setting:
\begin{align}
\Omega_H,\breve{\Omega}_H,\widehat{\Omega}_H,\widetilde{\Omega}_H,\\
\Pi^{\mathrm{bd}}_H,\breve{\Pi}_H^{\mathrm{bd}},\widehat{\Pi}_H^{\mathrm{bd}},\widetilde{\Pi}_H^{\mathrm{bd}},\\
\Pi_H,\breve{\Pi}_H,\widetilde{\Pi}_H.	
\end{align}
Then the corresponding rational finite projective modules over these rings carrying the corresponding $(\varphi,\Gamma)$-structures. For the deformed we look at the corresponding rings in the following:
\begin{align}
\Pi_{H,A},\breve{\Pi}_{H,A}
,\widetilde{\Pi}_{H,A}.	
\end{align}	
Then we have the corresponding parallel results as well.
\end{proposition}

\begin{proof}
See \cite[Theorem 5.7.4]{6KL16}.	
\end{proof}

\begin{remark}
One can certainly compare the corresponding $(\varphi,\Gamma)$-modules above with the corresponding Frobenius sheaves over the Kummer pro-\'etale topology in \cref{6setting2.18} as in \cite{6KL16} by adding additional effective conditions on the objects.
\end{remark}

%%%%%%%%%%%%%%%%%%%%%%%%%%%%%%%%%%!!!!!!!!!!!!!!!!!!!!

%%\newpage

\newpage\section{Hodge-Iwasawa Modules}

\subsection{The Categories of Hodge-Iwasawa Modules}

\indent We recall the corresponding Hodge-Iwasawa modules we defined in the previous papers \cite{6XT1} and \cite{6XT2}. We use our notations in \cite{6XT2} for the corresponding tower and the corresponding period rings.

\begin{setting}
We use the notation $(H,H^+)$ to denote the corresponding tower in our situation. Here we assume that $(H,H^+)$ is finite \'etale, noetherian, weakly-decompleting, decompleting, Galois with the Galois topological group $\Gamma$. We also then follow the corresponding setting of \cite[Hypothesis 3.1.1]{6KL16}. And fix some $a>0$ positive integer.	
\end{setting}

\begin{setting}
Recall that we have the corresponding period rings in the following, defined in \cite{6KL15} and \cite{6KL16}:
\begin{displaymath}
\widetilde{\Omega}_H^\mathrm{int},\widetilde{\Omega}_H,\widetilde{\Pi}_H^{\mathrm{int},r},	\widetilde{\Pi}_H^{\mathrm{int}},\widetilde{\Pi}_H^{\mathrm{bd},r},	\widetilde{\Pi}_H^{\mathrm{bd}},\widetilde{\Pi}_H^{r},	\widetilde{\Pi}_H^{I},\widetilde{\Pi}_H^\infty,	\widetilde{\Pi}_H
\end{displaymath}
form the corresponding group of perfect period rings along the tower $(H,H^+)$, here $r>0$ is some radius while $I$ corresponds to some closed interval. Then as in \cite[Definition 5.2.1]{6KL16} we have various imperfection of the corresponding rings listed above. We mainly need the corresponding Robba rings in the following fashion:
\begin{displaymath}
\breve{\Pi}_H^{r},	\breve{\Pi}_H^{I},\breve{\Pi}_H^\infty,	\breve{\Pi}_H,	
\end{displaymath}
and 
\begin{displaymath}
{\Pi}_H^{r},	{\Pi}_H^{I}, \Pi_H^\infty,	{\Pi}_H.	
\end{displaymath}
\end{setting}

\begin{setting}
The Hodge-Iwasawa modules we considered and studied are basically over the corresponding deformed version of the corresponding period rings we recalled above over some affinoid algebra $A$ over the base field $E$ as in \cite{6XT2}:
\begin{displaymath}
\breve{\Pi}_{H,A}^{r},	\breve{\Pi}_{H,A}^{I},\breve{\Pi}_{H,A}^\infty,	\breve{\Pi}_{H,A},	
\end{displaymath}
and 
\begin{displaymath}
{\Pi}_{H,A}^{r},	{\Pi}_{H,A}^{I}, \Pi_{H,A}^\infty,	{\Pi}_{H,A},	
\end{displaymath}	
and
\begin{displaymath}
\widetilde{\Pi}_{H,A}^{r},	\widetilde{\Pi}_{H,A}^{I},\widetilde{\Pi}_{H,A}^\infty,	\widetilde{\Pi}_{H,A}.
\end{displaymath}

\end{setting}

\indent Recall that we have the following corresponding key categories involved in our study (\cite{6XT2}). First we consider the corresponding finite projective objects for any suitable $s,r$ such that $0\leq s\leq r/p^{ah}$ (for any $A$ such that the corresponding rings building the corresponding Fargues-Fontaine curve namely the perfect Robba rings with respect to some intervals are sheafy as adic rings):

\noindent A. Look at the corresponding space $\mathrm{Spa}(H_0,H_0^+)$ which is the corresponding base of the tower, then we have the corresponding category of all the Frobenius sheaves over the pro-\'etale site of this space over the sheaf of ring $\widetilde{\Pi}_{\mathrm{Spa}(H_0,H_0^+)_\text{pro\'et},A}$. Here the Frobenius will be chosen in our situation to be $\varphi^a$;\\
\noindent B. We look at the same space as in the previous item, and we look at the corresponding category of the Frobenius bundles over the corresponding sheaf of ring $\widetilde{\Pi}_{\mathrm{Spa}(H_0,H_0^+)_\text{pro\'et},A}$;\\
\noindent C. We look at the same base space and then the corresponding sheaf of ring $\widetilde{\Pi}^{[s,r]}_{\mathrm{Spa}(H_0,H_0^+)_\text{pro\'et},A}$ for $s\leq r/p^{ah}$, we then have the category of the corresponding finite projective Frobenius modules over the period rings in this setting;\\
\noindent D. We then have the corresponding finite projective modules over the period ring ${\widetilde{\Pi}}_{H,A}$, carrying the corresponding action coming from $(\varphi^a,\Gamma)$;\\
\noindent E. We then have the corresponding finite projective bundles over the period ring ${\widetilde{\Pi}}_{H,A}$, carrying the corresponding action coming from $(\varphi^a,\Gamma)$;\\
\noindent F. We then have the corresponding finite projective modules over the period ring ${\widetilde{\Pi}}^{[s,r]}_{H,A}$, carrying the corresponding action coming from $(\varphi^a,\Gamma)$, $0\leq s\leq r/p^{ah}$;\\
\noindent G. We then have the corresponding finite projective modules over the period ring ${{\Pi}}_{H,A}$, carrying the corresponding action coming from $(\varphi^a,\Gamma)$;\\
\noindent H. We then have the corresponding finite projective bundles over the period ring ${{\Pi}}_{H,A}$, carrying the corresponding action coming from $(\varphi^a,\Gamma)$;\\
\noindent I. We then have the corresponding finite projective modules over the period ring ${{\Pi}}^{[s,r]}_{H,A}$, carrying the corresponding action coming from $(\varphi^a,\Gamma)$, $0\leq s\leq r/p^{ah}$;\\
\noindent J. We then have the corresponding finite projective modules over the period ring ${\breve{\Pi}}_{H,A}$, carrying the corresponding action coming from $(\varphi^a,\Gamma)$;\\
\noindent K. We then have the corresponding finite projective bundles over the period ring ${\breve{\Pi}}_{H,A}$, carrying the corresponding action coming from $(\varphi^a,\Gamma)$;\\
\noindent L. We then have the corresponding finite projective modules over the period ring ${\breve{\Pi}}^{[s,r]}_{H,A}$, carrying the corresponding action coming from $(\varphi^a,\Gamma)$, $0\leq s\leq r/p^{ah}$;\\
\noindent M. We have the corresponding adic version of the corresponding Fargues-Fontaine curve $FF_{\overline{H}'_\infty,A}$. And we can consider the corresponding category of $\Gamma$-equivariant quasi-coherent locally free sheaves over this adic version of Fargues-Fontaine curve, and we assume the corresponding $\Gamma$-equivariant action to be semilinear and continuous over each $\Gamma$-equivariant affinoid subspace.

\begin{remark}
One has to be careful when one would like to define the objects in $A,B,C$ which could be compared to the other objects. One needs to consider the corresponding locally finite projective objects instead of locally finite free objects.	
\end{remark}

\indent Then we have:

\begin{proposition}\mbox{}\\
\mbox{\bf{(See \cite[Proposition 5.44]{6XT2}, after Kedlaya-Liu \cite[Theorem 5.7.4]{6KL16})}}
The categories described above are equivalent to each other.	
\end{proposition}

%%%\\\\\\\\\\\\\\\\

\indent Then we consider the corresponding pseudocoherent objects, where we further assume that the corresponding ring $R_H$ is now $F$-(finite projective) (for any $A$ such that the corresponding rings building the corresponding Fargues-Fontaine curve namely the perfect Robba rings with respect to some intervals are sheafy as adic rings).

\noindent A. Look at the corresponding space $\mathrm{Spa}(H_0,H_0^+)$ which is the corresponding base of the tower, then we have the corresponding category of all the Frobenius sheaves over the pro-\'etale site of this space over the sheaf of ring $\widetilde{\Pi}_{\mathrm{Spa}(H_0,H_0^+)_\text{pro\'et},A}$. Here the Frobenius will be chosen in our situation to be $\varphi^a$;\\
\noindent B. We look at the same space as in the previous item, and we look at the corresponding category of the Frobenius bundles over the corresponding sheaf of ring $\widetilde{\Pi}_{\mathrm{Spa}(H_0,H_0^+)_\text{pro\'et},A}$;\\
\noindent C. We look at the same base space and then the corresponding sheaf of ring $\widetilde{\Pi}^{[s,r]}_{\mathrm{Spa}(H_0,H_0^+)_\text{pro\'et},A}$ for $s\leq r/p^{ah}$, we then have the category of the corresponding pseudocoherent Frobenius modules over the period rings in this setting;\\
\noindent D. We then have the corresponding pseudocoherent modules over the period ring ${\widetilde{\Pi}}_{H,A}$, carrying the corresponding action coming from $(\varphi^a,\Gamma)$;\\
\noindent E. We then have the corresponding pseudocoherent bundles over the period ring ${\widetilde{\Pi}}_{H,A}$, carrying the corresponding action coming from $(\varphi^a,\Gamma)$;\\
\noindent F. We then have the corresponding pseudocoherent modules over the period ring ${\widetilde{\Pi}}^{[s,r]}_{H,A}$, carrying the corresponding action coming from $(\varphi^a,\Gamma)$, $0\leq s\leq r/p^{ah}$;\\
\noindent G. We then have the corresponding pseudocoherent modules over the period ring ${{\Pi}}_{H,A}$, carrying the corresponding action coming from $(\varphi^a,\Gamma)$;\\
\noindent H. We then have the corresponding pseudocoherent bundles over the period ring ${{\Pi}}_{H,A}$, carrying the corresponding action coming from $(\varphi^a,\Gamma)$;\\
\noindent I. We then have the corresponding pseudocoherent modules over the period ring ${{\Pi}}^{[s,r]}_{H,A}$, carrying the corresponding action coming from $(\varphi^a,\Gamma)$, $0\leq s\leq r/p^{ah}$;\\
\noindent J. We then have the corresponding pseudocoherent modules over the period ring ${\breve{\Pi}}_{H,A}$, carrying the corresponding action coming from $(\varphi^a,\Gamma)$;\\
\noindent K. We then have the corresponding pseudocoherent bundles over the period ring ${\breve{\Pi}}_{H,A}$, carrying the corresponding action coming from $(\varphi^a,\Gamma)$;\\
\noindent L. We then have the corresponding pseudocoherent modules over the period ring ${\breve{\Pi}}^{[s,r]}_{H,A}$, carrying the corresponding action coming from $(\varphi^a,\Gamma)$, $0\leq s\leq r/p^{ah}$;\\
\noindent M. We have the corresponding adic version of the corresponding Fargues-Fontaine curve $FF_{\overline{H}'_\infty,A}$. And we can consider the corresponding category of pseudocoherent sheaves over this adic version of Fargues-Fontaine curve.

\begin{remark}
One has to be careful when one would like to define the objects in $A,B,C$ which could be compared to the other objects. One needs to consider the corresponding locally pseudocoherent objects instead of locally pseudocoherent objects. Moreover in this current setting locally we have to impose $p$-adic functional analytic conditions on the pseudocoherent objects when one would like to do the localization in different topologies.
\end{remark}

\indent Then we have:

\begin{proposition}\mbox{}\\
\mbox{\bf{(See \cite[Proposition 5.51]{6XT2}, after Kedlaya-Liu \cite[Theorem 5.9.4]{6KL16})}}
The categories described above are equivalent to each other.	
\end{proposition}

\begin{proposition} \mbox{\bf{(After Kedlaya-Liu \cite[Theorem 5.9.4]{6KL16})}}
The categories above are abelian.	
\end{proposition}

\begin{proof}
By our assumption the corresponding ring $\breve{\Pi}^{[s,r]}_{H,A}$ is coherent.	
\end{proof}

\indent At least in the corresponding finite projective situation, we also have some comparison on the corresponding objects in the logarithmic setting:

\begin{proposition} \mbox{\bf{(After Kedlaya-Liu \cite[Lemma 5.4.11]{6KL16})}}
Let $(H,H^+)$ be the one of the corresponding logarithmic towers we encountered in the previous section, we then have the following categories are basically equivalent:\\
\noindent A. Look at the corresponding space $\mathrm{Spa}(H_0,H_0^+)$ which is the corresponding base of the tower, then we have the corresponding category of all the Frobenius sheaves over the pro-\'etale site of this space over the sheaf of ring $\widetilde{\Pi}_{\mathrm{Spa}(H_0,H_0^+)_\text{prok\'et},A}$. Here the Frobenius will be chosen in our situation to be $\varphi^a$;\\
\noindent B. We look at the same space as in the previous item, and we look at the corresponding category of the Frobenius bundles over the corresponding sheaf of ring $\widetilde{\Pi}_{\mathrm{Spa}(H_0,H_0^+)_\text{prok\'et},A}$;\\
\noindent C. We look at the same base space and then the corresponding sheaf of ring $\widetilde{\Pi}^{[s,r]}_{\mathrm{Spa}(H_0,H_0^+)_\text{prok\'et},A}$ for $s\leq r/p^{ah}$, we then have the category of the corresponding finite projective Frobenius modules over the period rings in this setting.\\
\end{proposition}

\begin{proof}
This is because the corresponding log perfectoid affinoids form a basis of the corresponding pro-Kummer \'etale site (see the corresponding development in \cite[Chapter 5.3]{6DLLZ2}).
\end{proof}

\begin{proposition} \mbox{\bf{(After Kedlaya-Liu \cite[Theorem 5.7.4]{6KL16})}} Let $(H,H^+)$ be the one of the corresponding logarithmic towers we encountered in the previous section. Then we have the following categories of Hodge-Iwasawa modules are equivalent:\\
\noindent D. We then have the corresponding finite projective modules over the period ring ${\widetilde{\Pi}}_{H,A}$, carrying the corresponding action coming from $(\varphi^a,\Gamma)$;\\
\noindent E. We then have the corresponding finite projective bundles over the period ring ${\widetilde{\Pi}}_{H,A}$, carrying the corresponding action coming from $(\varphi^a,\Gamma)$;\\
\noindent F. We then have the corresponding finite projective modules over the period ring ${\widetilde{\Pi}}^{[s,r]}_{H,A}$, carrying the corresponding action coming from $(\varphi^a,\Gamma)$, $0\leq s\leq r/p^{ah}$;\\
\noindent G. We then have the corresponding finite projective modules over the period ring ${{\Pi}}_{H,A}$, carrying the corresponding action coming from $(\varphi^a,\Gamma)$;\\
\noindent H. We then have the corresponding finite projective bundles over the period ring ${{\Pi}}_{H,A}$, carrying the corresponding action coming from $(\varphi^a,\Gamma)$;\\
\noindent I. We then have the corresponding finite projective modules over the period ring ${{\Pi}}^{[s,r]}_{H,A}$, carrying the corresponding action coming from $(\varphi^a,\Gamma)$, $0\leq s\leq r/p^{ah}$;\\
\noindent J. We then have the corresponding finite projective modules over the period ring ${\breve{\Pi}}_{H,A}$, carrying the corresponding action coming from $(\varphi^a,\Gamma)$;\\
\noindent K. We then have the corresponding finite projective bundles over the period ring ${\breve{\Pi}}_{H,A}$, carrying the corresponding action coming from $(\varphi^a,\Gamma)$;\\
\noindent L. We then have the corresponding finite projective modules over the period ring ${\breve{\Pi}}^{[s,r]}_{H,A}$, carrying the corresponding action coming from $(\varphi^a,\Gamma)$, $0\leq s\leq r/p^{ah}$.\\
%\noindent M. We have the corresponding adic version of the corresponding Fargues-Fontaine curve $FF_{\overline{H}'_\infty,A}$. And we can consider the corresponding category of $\Gamma$-equivariant quasi-coherent locally free sheaves over this adic version of Fargues-Fontaine curve, and we assume the corresponding $\Gamma$-equivariant action to be semilinear and continuous over each $\Gamma$-equivariant affinoid subspace.
\end{proposition}

\begin{proof}
The vertical comparison along changing the corresponding ring with respect to the accents was established in \cref{proposition2.20}. Then to compare the corresponding objects in the remaining cases see \cite[Theorem 5.7.4]{6KL16} and \cite[Proposition 5.44]{6XT2}.	
\end{proof}

\indent In Iwaswa theory, the corresponding aspects of Hodge-Iwasawa theory will usually happen over some Fr\'echet-Stein algebras. Here we will use some geometric language which is different from the situation we considered in our previous work. We now first consider the corresponding commutative setting. We use the notation $X$ to denote a general Stein space:
\begin{displaymath}
X:= \bigcup_{n} X_n	
\end{displaymath}
regarded as a general quasi-Stein adic space in the corresponding context of \cite[Chapter 2.6]{6KL16}. Then we have the corresponding sheaf will be organized in the following sense:
\begin{align}
\mathcal{O}_X:=\varprojlim_{n} \mathcal{O}_{X_n}.	
\end{align}

Therefore we will correspondingly consider the period rings over the corresponding spaces taking the form of $X$ as above. We assume each $\mathcal{O}_{X_n}$ satisfies the corresponding assumption as that for the ring $A$ as discussed above. 

\begin{definition}
For any ring $R$ (commutative) which could be written as the corresponding inverse system $\varprojlim_n R_n$. We call the projective system of modules over $R$ as a projective system $\{\mathcal{M}_n\}_n$ of all modules $M_n$ over $R_n$ for each $n\geq 0$.	
\end{definition}

\noindent A. Look at the corresponding space $\mathrm{Spa}(H_0,H_0^+)$ which is the corresponding base of the tower, then we have the corresponding category of all the projective systems of the Frobenius sheaves over the pro-\'etale site of this space over the sheaf of ring $\widetilde{\Pi}_{\mathrm{Spa}(H_0,H_0^+)_\text{pro\'et},\varprojlim_n\mathcal{O}_{X_n}}$. Here the Frobenius will be chosen in our situation to be $\varphi^a$;\\
\noindent B. We look at the same space as in the previous item, and we look at the corresponding category of all the projective systems of the Frobenius bundles over the corresponding sheaf of ring $\widetilde{\Pi}_{\mathrm{Spa}(H_0,H_0^+)_\text{pro\'et},\varprojlim_n\mathcal{O}_{X_n}}$;\\
\noindent C. We look at the same base space and then the corresponding sheaf of ring\\ $\widetilde{\Pi}^{[s,r]}_{\mathrm{Spa}(H_0,H_0^+)_\text{pro\'et},\varprojlim_n\mathcal{O}_{X_n}}$ for $s\leq r/p^{ah}$, we then have the category of the projective systems of  corresponding finite projective Frobenius modules over the period rings in this setting;\\
\noindent D. We then have the corresponding projective systems of  finite projective modules over the period ring ${\widetilde{\Pi}}_{H,\varprojlim_n\mathcal{O}_{X_n}}$, carrying the corresponding action coming from $(\varphi^a,\Gamma)$;\\
\noindent E. We then have the corresponding projective systems of finite projective bundles over the period ring ${\widetilde{\Pi}}_{H,\varprojlim_n\mathcal{O}_{X_n}}$, carrying the corresponding action coming from $(\varphi^a,\Gamma)$;\\
\noindent F. We then have the corresponding projective systems of finite projective modules over the period ring ${\widetilde{\Pi}}^{[s,r]}_{H,\varprojlim_n\mathcal{O}_{X_n}}$, carrying the corresponding action coming from $(\varphi^a,\Gamma)$;\\
\noindent G. We then have the corresponding projective systems of finite projective modules over the period ring ${{\Pi}}_{H,\varprojlim_n\mathcal{O}_{X_n}}$, carrying the corresponding action coming from $(\varphi^a,\Gamma)$;\\
\noindent H. We then have the corresponding projective systems of finite projective bundles over the period ring ${{\Pi}}_{H,\varprojlim_n\mathcal{O}_{X_n}}$, carrying the corresponding action coming from $(\varphi^a,\Gamma)$;\\
\noindent I. We then have the corresponding projective systems of finite projective modules over the period ring ${{\Pi}}^{[s,r]}_{H,\varprojlim_n\mathcal{O}_{X_n}}$, carrying the corresponding action coming from $(\varphi^a,\Gamma)$;\\
\noindent J. We then have the corresponding projective systems of finite projective modules over the period ring ${\breve{\Pi}}_{H,\varprojlim_n\mathcal{O}_{X_n}}$, carrying the corresponding action coming from $(\varphi^a,\Gamma)$;\\
\noindent K. We then have the corresponding projective systems of finite projective bundles over the period ring ${\breve{\Pi}}_{H,\varprojlim_n\mathcal{O}_{X_n}}$, carrying the corresponding action coming from $(\varphi^a,\Gamma)$;\\
\noindent L. We then have the corresponding projective systems of finite projective modules over the period ring ${\breve{\Pi}}^{[s,r]}_{H,\varprojlim_n\mathcal{O}_{X_n}}$, carrying the corresponding action coming from $(\varphi^a,\Gamma)$;\\
\noindent M. We have the corresponding adic version of the corresponding Fargues-Fontaine curve $FF_{\overline{H}'_\infty,\varprojlim_n\mathcal{O}_{X_n}}$. And we can consider the corresponding category of the projective systems of $\Gamma$-equivariant quasi-coherent locally free sheaves over this adic version of Fargues-Fontaine curve, and we assume the corresponding $\Gamma$-equivariant action to be semilinear and continuous over each $\Gamma$-equivariant affinoid subspace.\\

\indent Then we have:

\begin{proposition}\mbox{}\\\mbox{\bf{(See \cite[Theorem 4.11]{6XT2}, after Kedlaya-Liu \cite[Theorem 5.7.4]{6KL16})}}
The categories described above are equivalent to each other. 	
\end{proposition}

\indent For the pseudocoherent setting we have the following categories:

\noindent A. Look at the corresponding space $\mathrm{Spa}(H_0,H_0^+)$ which is the corresponding base of the tower, then we have the corresponding category of all the projective systems of the pseudocoherent Frobenius sheaves over the pro-\'etale site of this space over the sheaf of ring $\widetilde{\Pi}_{\mathrm{Spa}(H_0,H_0^+)_\text{pro\'et},\varprojlim_n\mathcal{O}_{X_n}}$. Here the Frobenius will be chosen in our situation to be $\varphi^a$;\\
\noindent B. We look at the same space as in the previous item, and we look at the corresponding category of all the projective systems of the pseudocoherent Frobenius bundles over the corresponding sheaf of ring $\widetilde{\Pi}_{\mathrm{Spa}(H_0,H_0^+)_\text{pro\'et},\varprojlim_n\mathcal{O}_{X_n}}$;\\
\noindent C. We look at the same base space and then the corresponding sheaf of ring\\ $\widetilde{\Pi}^{[s,r]}_{\mathrm{Spa}(H_0,H_0^+)_\text{pro\'et},\varprojlim_n\mathcal{O}_{X_n}}$ for $s\leq r/p^{ah}$, we then have the category of the projective systems of corresponding pseudocoherent Frobenius modules over the period rings in this setting;\\
\noindent D. We then have the corresponding projective systems of  pseudocoherent modules over the period ring ${\widetilde{\Pi}}_{H,\varprojlim_n\mathcal{O}_{X_n}}$, carrying the corresponding action coming from $(\varphi^a,\Gamma)$;\\
\noindent E. We then have the corresponding projective systems of pseudocoherent bundles over the period ring ${\widetilde{\Pi}}_{H,\varprojlim_n\mathcal{O}_{X_n}}$, carrying the corresponding action coming from $(\varphi^a,\Gamma)$;\\
\noindent F. We then have the corresponding projective systems of pseudocoherent modules over the period ring ${\widetilde{\Pi}}^{[s,r]}_{H,\varprojlim_n\mathcal{O}_{X_n}}$, carrying the corresponding action coming from $(\varphi^a,\Gamma)$;\\
\noindent G. We then have the corresponding projective systems of pseudocoherent modules over the period ring ${{\Pi}}_{H,\varprojlim_n\mathcal{O}_{X_n}}$, carrying the corresponding action coming from $(\varphi^a,\Gamma)$;\\
\noindent H. We then have the corresponding projective systems of pseudocoherent bundles over the period ring ${{\Pi}}_{H,\varprojlim_n\mathcal{O}_{X_n}}$, carrying the corresponding action coming from $(\varphi^a,\Gamma)$;\\
\noindent I. We then have the corresponding projective systems of pseudocoherent modules over the period ring ${{\Pi}}^{[s,r]}_{H,\varprojlim_n\mathcal{O}_{X_n}}$, carrying the corresponding action coming from $(\varphi^a,\Gamma)$;\\
\noindent J. We then have the corresponding projective systems of pseudocoherent modules over the period ring ${\breve{\Pi}}_{H,\varprojlim_n\mathcal{O}_{X_n}}$, carrying the corresponding action coming from $(\varphi^a,\Gamma)$;\\
\noindent K. We then have the corresponding projective systems of pseudocoherent bundles over the period ring ${\breve{\Pi}}_{H,\varprojlim_n\mathcal{O}_{X_n}}$, carrying the corresponding action coming from $(\varphi^a,\Gamma)$;\\
\noindent L. We then have the corresponding projective systems of pseudocoherent modules over the period ring ${\breve{\Pi}}^{[s,r]}_{H,\varprojlim_n\mathcal{O}_{X_n}}$, carrying the corresponding action coming from $(\varphi^a,\Gamma)$;\\
\noindent M. We have the corresponding adic version of the corresponding Fargues-Fontaine curve $FF_{\overline{H}'_\infty,\varprojlim_n\mathcal{O}_{X_n}}$. And we can consider the corresponding category of the projective systems of $\Gamma$-equivariant pseudocoherent sheaves over this adic version of Fargues-Fontaine curve, and we assume the corresponding $\Gamma$-equivariant action to be semilinear and continuous over each $\Gamma$-equivariant affinoid subspace.\\

\indent Then we have:

\begin{proposition}\mbox{}\\
\mbox{\bf{(See \cite[Theorem 4.11]{6XT2}, after Kedlaya-Liu \cite[Theorem 5.9.6]{6KL16})}}
The categories described above are equivalent to each other. And they are actually abelian.\\	
\end{proposition}

%%%%%%%%%%%%%%%%%%%%%%%%%%%%%%%%%!!!!!!!!!!!!!!!!!!!!!!!!

\subsection{The Cohomologies of Hodge-Iwasawa Modules}

\indent We then define the corresponding Herr style cohomology in our current deformed setting for the corresponding. Note that the corresponding cohomology over the corresponding Fr\'echet-Stein algebras are very important in the study of Iwasawa conjectures and the corresponding $p$-adic Tamagawa number conjectures.

\indent First keep the following setting:

\begin{setting}
We use the notation $(H,H^+)$ to denote the corresponding tower in our situation. Here we assume that $(H,H^+)$ is finite \'etale, noetherian, weakly-decompleting, decompleting, Galois with the Galois topological group $\Gamma$. We also then follow the corresponding setting of \cite[Hypothesis 3.1.1]{6KL16}. And fix some $a>0$ positive integer.	
\end{setting}

\noindent A. Look at the corresponding space $\mathrm{Spa}(H_0,H_0^+)$ which is the corresponding base of the tower, then we have the corresponding category of all the Frobenius sheaves over the pro-\'etale site of this space over the sheaf of ring $\widetilde{\Pi}_{\mathrm{Spa}(H_0,H_0^+)_\text{pro\'et},A}$. Here the Frobenius will be chosen in our situation to be $\varphi^a$;\\
\begin{definition} \mbox{\bf{(After Kedlaya-Liu, \cite[Definition 4.4.4]{6KL16})}} For the category in the above (and the corresponding setting in pseudocoherent situation) we define the corresponding complex $C^\bullet_{\varphi}(M)$ of any object $\mathcal{M}$ as the corresponding hypercomplex of the following complex for just the Frobenius operator:
\[
\xymatrix@R+0pc@C+0pc{
0\ar[r] \ar[r] \ar[r] &\mathcal{M} \ar[r]^{\varphi^a-1} \ar[r] \ar[r]  & \mathcal{M} \ar[r] \ar[r] \ar[r]&0.
}
\]	
\end{definition}

\noindent B. We look at the same space as in the previous item, and we look at the corresponding category of the Frobenius bundles over the corresponding sheaf of ring $\widetilde{\Pi}_{\mathrm{Spa}(H_0,H_0^+)_\text{pro\'et},A}$;\\

\noindent C. We look at the same base space and then the corresponding sheaf of ring $\widetilde{\Pi}^{[s,r]}_{\mathrm{Spa}(H_0,H_0^+)_\text{pro\'et},A}$ for $s\leq r/p^{ah}$, we then have the category of the corresponding finite projective Frobenius modules over the period rings in this setting;

\begin{definition} \mbox{\bf{(After Kedlaya-Liu, \cite[Definition 4.4.4]{6KL16})}} For the three categories in the above (and the corresponding setting in pseudocoherent situation) we define the corresponding complex $C^\bullet_{\varphi}(M)$ of any object $\mathcal{M}$ as the corresponding hypercomplex of the following complex for just the Frobenius operator:
\[
\xymatrix@R+0pc@C+0pc{
0\ar[r] \ar[r] \ar[r] &\mathcal{M} \ar[r] \ar[r] \ar[r]  & \mathcal{M}\otimes_{\widetilde{\Pi}^{[s,r]}_{\mathrm{Spa}(H_0,H_0^+)_\text{pro\'et},A}}\widetilde{\Pi}^{[s,r/p^{ah}]}_{\mathrm{Spa}(H_0,H_0^+)_\text{pro\'et},A} \ar[r] \ar[r] \ar[r]&0.
}
\]

\end{definition}

\noindent D. We then have the corresponding finite projective modules over the period ring ${\widetilde{\Pi}}_{H,A}$, carrying the corresponding action coming from $(\varphi^a,\Gamma)$;\\
\noindent E. We then have the corresponding finite projective bundles over the period ring ${\widetilde{\Pi}}_{H,A}$, carrying the corresponding action coming from $(\varphi^a,\Gamma)$;\\
\noindent F. We then have the corresponding finite projective modules over the period ring ${\widetilde{\Pi}}^{[s,r]}_{H,A}$, carrying the corresponding action coming from $(\varphi^a,\Gamma)$;\\
\noindent G. We then have the corresponding finite projective modules over the period ring ${{\Pi}}_{H,A}$, carrying the corresponding action coming from $(\varphi^a,\Gamma)$;\\
\noindent H. We then have the corresponding finite projective bundles over the period ring ${{\Pi}}_{H,A}$, carrying the corresponding action coming from $(\varphi^a,\Gamma)$;\\
\noindent I. We then have the corresponding finite projective modules over the period ring ${{\Pi}}^{[s,r]}_{H,A}$, carrying the corresponding action coming from $(\varphi^a,\Gamma)$;\\
\noindent J. We then have the corresponding finite projective modules over the period ring ${\breve{\Pi}}_{H,A}$, carrying the corresponding action coming from $(\varphi^a,\Gamma)$;\\
\noindent K. We then have the corresponding finite projective bundles over the period ring ${\breve{\Pi}}_{H,A}$, carrying the corresponding action coming from $(\varphi^a,\Gamma)$;\\
\noindent L. We then have the corresponding finite projective modules over the period ring ${\breve{\Pi}}^{[s,r]}_{H,A}$, carrying the corresponding action coming from $(\varphi^a,\Gamma)$;\\

\begin{definition} \mbox{\bf{(After Kedlaya-Liu, \cite[Definition 5.7.9]{6KL16})}} For the categories (except for the corresponding bundles) in the above (and the corresponding setting in pseudocoherent situation) we define the corresponding complex $C^\bullet_{\varphi,\Gamma}(\mathcal{M})$ of any object $\mathcal{M}$ as the corresponding totalization of the following complex:
\[
\xymatrix@R+0pc@C+0pc{
0\ar[r] \ar[r] \ar[r] &C^\bullet_{\Gamma}(\mathcal{M})\ar[r]^{\varphi^a-1} \ar[r] \ar[r]  &C^\bullet_{\Gamma}(\mathcal{M}) \ar[r] \ar[r] \ar[r]&0.
}
\]	
or 
\[
\xymatrix@R+0pc@C+0pc{
0\ar[r] \ar[r] \ar[r] &C^\bullet_{\Gamma}(\mathcal{M})\ar[r]^{\varphi^a-1} \ar[r] \ar[r]  &C^\bullet_{\Gamma}(\mathcal{M})\otimes \Delta \ar[r] \ar[r] \ar[r]&0,
}
\]
respectively for different situations in the above where $\Delta\in \{{\breve{\Pi}}^{[s,r/p^{ah}]}_{H,A},{{\Pi}}^{[s,r/p^{ah}]}_{H,A},{\widetilde{\Pi}}^{[s,r/p^{ah}]}_{H,A}\}$.
\end{definition}

\noindent M. We have the corresponding adic version of the corresponding Fargues-Fontaine curve $FF_{\overline{H}'_\infty,A}$. And we can consider the corresponding category of $\Gamma$-equivariant quasi-coherent locally free sheaves over this adic version of Fargues-Fontaine curve (with the deformation from the algebra $A$), and we assume the corresponding $\Gamma$-equivariant action to be semilinear and continuous over each $\Gamma$-equivariant affinoid subspace.

\begin{definition} \mbox{\bf{(After Kedlaya-Liu, \cite{6KL16})}} For the categories in the above we define the corresponding complex $C^\bullet_{\Gamma}(\mathcal{M})$ of any object $\mathcal{M}$ as the corresponding hypercomplex of the complex $C^\bullet_{\Gamma}(\mathcal{M})$.	 In pseudocoherent situation, we have the parallel definition.
\end{definition}

\begin{remark}
Although we will not basically repeat the corresponding construction of the cohomology as above, but all the definitions on the corresponding cohomologies could be translated to the corresponding construction for the corresponding logarithmic setting and the corresponding Fr\'echet-Stein context we considered seriously above.	
\end{remark}

\indent Based on the corresponding comparison theorem above we can actually consider the corresponding comparison between the corresponding cohomologies. The corresponding statements will take the following forms:

\begin{proposition} \mbox{\bf{(After Kedlaya-Liu, \cite[Theorem 5.7.10]{6KL16})}}
Suppose we have a corresponding finite projective module $\mathcal{M}$ over $\Pi_{H,A}$ carrying the corresponding action of $(\varphi^a,\Gamma)$, and we use the notation $\breve{\mathcal{M}}$ and $\widetilde{M}$ to be the corresponding base change of $\mathcal{M}$ from the original base ring to the ring being one of $\breve{\Pi}_{H,A}$ and $\widetilde{\Pi}_{H,A}$ respectively. Then we have the following three complexes are mutually quasi-isomorphic:
\begin{align}
C_{\varphi,\Gamma}^\bullet(\mathcal{M}),C_{\varphi,\Gamma}^\bullet(\breve{\mathcal{M}}),C_{\varphi,\Gamma}^\bullet(\widetilde{\mathcal{M}}).
\end{align}

\end{proposition}

\begin{proposition} \mbox{\bf{(After Kedlaya-Liu, \cite[Theorem 5.7.10]{6KL16})}}
Suppose we have a corresponding finite projective module $\mathcal{M}$ over $\Pi_{H,\varprojlim_n\mathcal{O}_{X_n}}$ carrying the corresponding action of $(\varphi^a,\Gamma)$, and we use the notation $\breve{\mathcal{M}}$ and $\widetilde{M}$ to be the corresponding base change of $\mathcal{M}$ from the original base ring to the ring being one of $\breve{\Pi}_{H,\varprojlim_n\mathcal{O}_{X_n}}$ and $\widetilde{\Pi}_{H,\varprojlim_n\mathcal{O}_{X_n}}$ respectively. Then we have the following three complexes are mutually quasi-isomorphic:
\begin{align}
C_{\varphi,\Gamma}^\bullet(\mathcal{M}),C_{\varphi,\Gamma}^\bullet(\breve{\mathcal{M}}),C_{\varphi,\Gamma}^\bullet(\widetilde{\mathcal{M}}).
\end{align}

\end{proposition}

\indent As in \cite{6KL16} starting from a corresponding finite projective module carrying the corresponding $(\varphi^a,\Gamma)$-structure, one can compare the corresponding cohomologies of modules over the full Robba rings with those of those modules over some Robba rings with respect to some intervals.

\begin{proposition}\mbox{\bf{(After Kedlaya-Liu, \cite[Theorem 5.7.11]{6KL16})}}
Descend $\mathcal{M}$\\
 which is a finite projective module over $\Pi_{H,A}$ carrying the action of $(\varphi^a,\Gamma)$ to the corresponding Robba ring $\Pi^{r_0}_{H,A}$ for some specific positive radius $r_0>0$. Then as in the previous proposition we consider the corresponding base change of the module $\mathcal{M}$ to the corresponding Robba rings $\breve{\Pi}^{[s,r]}_{H,A}$ for $0<s\leq r\leq r_0$ with that $s\leq r/p^{ah}$, and to the corresponding Robba rings $\widetilde{\Pi}^{[s,r]}_{H,A}$. We denote these two modules by $\breve{\mathcal{M}}$ and $\widetilde{\mathcal{M}}$. Then for sufficiently large integer $m\geq 0$ suppose we consider the following three complexes:
\[
\xymatrix@R+0pc@C+0pc{
0\ar[r] \ar[r] \ar[r] &C^\bullet_{\Gamma}(\breve{\mathcal{M}}_{[s,r]})\ar[r]^{\varphi^a-1} \ar[r] \ar[r]  &C^\bullet_{\Gamma}(\breve{\mathcal{M}}_{[s,r/p^{ah}]}) \ar[r] \ar[r] \ar[r]&0,
}
\]
\[
\xymatrix@R+0pc@C+0pc{
0\ar[r] \ar[r] \ar[r] &C^\bullet_{\Gamma}(\widetilde{\mathcal{M}}_{[s,r]})\ar[r]^{\varphi^a-1} \ar[r] \ar[r]  &C^\bullet_{\Gamma}(\widetilde{\mathcal{M}}_{[s,r/p^{ah}]}) \ar[r] \ar[r] \ar[r]&0,
}
\]
\[
\xymatrix@R+0pc@C+0pc{
0\ar[r] \ar[r] \ar[r] &\varphi^{-am}C^\bullet_{\Gamma}(\mathcal{M}_{[s/p^{ahm},r/p^{ahm}]})\ar[r]^{\varphi^a-1} \ar[r] \ar[r]  &\varphi^{-am}C^\bullet_{\Gamma}(\mathcal{M}_{[s/p^{ahm},r/p^{ah(m+1)}]}) \ar[r] \ar[r] \ar[r]&0.
}
\]
Then we have that these three complexes and the corresponding complex $C^\bullet_{\varphi,\Gamma}(\mathcal{M})$ are quasi-isomorphic to each other.	
\end{proposition}

\begin{proof}
See \cite[Theorem 5.7.11]{6KL16}.	
\end{proof}

\begin{remark}
We want to mention here that actually the corresponding context should work in general for the tower which is finite \'etale but not actually Galois. Although in our previous work \cite{6XT2} we only considered the corresponding towers which are required to be Galois with respect to some Galois topological group $\Gamma$.	
\end{remark}

%%\newpage

\newpage\section{Applications to General Analytic Spaces}

\subsection{Contact with Rigid Analytic Spaces}

\indent We now consider the corresponding context of more global geometric objects, over smooth proper rigid analytic spaces. We certainly will consider very general coefficients namely the corresponding affinoid algebras in the same framework. The two considerations happen to carry the same context coming from the corresponding rigid analytic spaces, but more strictly speaking they will interact in certain situations such as when we look at some nonabelian quotient of the profinite fundamental groups.

%%%%%%%%%%%%%%%%%%%%%%%%%%%%%%%%%%%%%%!!!!!!!!!!!!!!!!!!!!!!

\indent We basically consider the following context, where $X$ will be a corresponding smooth proper rigid analytic space over some analytic field $K$. Here we assume that the corresponding analytic field $K$ is containing the $p$-adic number field $\mathbb{Q}_p$. The affinoid algebra $A$ over $\mathbb{Q}_p$ in the following discussion is assumed to be sousperfectoid as those considered in \cite{6KH}. Then we can consider the following category.

\begin{definition}\mbox{\bf{(After Kedlaya-Liu, \cite[Definition 8.10.1]{6KL16})}}
We use the corresponding notation $D_{\mathrm{pseudo},\widetilde{\Pi}^{[s,r]}_{X,A}}$ ($0<s\leq r/p^{a}$) to denote the corresponding category of all the pseudocoherent $\widetilde{\Pi}^{[s,r]}_{X,A}$ sheaves, which is regarded as a subcategory of the corresponding category of all the corresponding sheaves over $\widetilde{\Pi}^{[s,r]}_{X,A}$, in the corresponding pro-\'etale topology.	
\end{definition}

\indent Also we have the following category of all the $(\varphi^a,\Gamma)$ over the space $X$ as in the above:

\begin{definition}\mbox{\bf{(After Kedlaya-Liu, \cite[Definition 8.10.5]{6KL16})}} 
We use the corresponding notation $D_{\mathrm{pseudo},\varphi^a,\widetilde{\Pi}_{X,A}}$ to denote the corresponding category of all the pseudocoherent $\widetilde{\Pi}_{X,A}$ sheaves carrying the corresponding Frobenius operator $\varphi^a$, which is regarded as a subcategory of the corresponding category of all the corresponding sheaves over $\widetilde{\Pi}_{X,A}$, in the corresponding pro-\'etale topology.	
\end{definition}

\begin{proposition} \mbox{\bf{(After Kedlaya-Liu, \cite[Theorem 8.10.6]{6KL16})}} \label{proposition6.3}
The corresponding category $D_{\mathrm{pseudo},\varphi^a,\widetilde{\Pi}_{X,A}}$ is abelian, where the corresponding kernels and cokernels are compatible with the corresponding category of all the sheaves of $\widetilde{\Pi}_{X,A}$-modules over the corresponding pro-\'etale site.	
\end{proposition}

\begin{proof}
Locally the corresponding object in this category will be $\varphi^a$-sheaf in the corresponding category $D_{\mathrm{pseudo},\varphi^a,\widetilde{\Pi}_{X_0,A}}$ where $X_0$ is a smooth affinoid subspace. We then deduce the corresponding $\Gamma$-action will be realized by the corresponding Galois group of the corresponding toric tower. Then by the comparison we established in the above, we can have the chance to directly compare this to the corresponding $(\varphi^a,\Gamma)$-module over the ring $\breve{\Pi}_{H,A}$ for some restricted tower $H$, which implies the corresponding result since this ring is coherent. 	
\end{proof}

\subsection{Contact with $k_\Delta$-Analytic Spaces}

\indent In this section we are going to consider the corresponding $k_\Delta$-analytic space in the sense of \cite[Part I Chapter I]{6DFN}. The corresponding $\Delta$ was denoted by $H$ in \cite[Part I Chapter I]{6DFN}, namely a commutative multiplication group satisfying the corresponding relation:
\begin{displaymath}
|k^\times|\subset H \subset \mathbb{R}^\times_+.	
\end{displaymath}
Here for $k$ we assume that this is an analytic field of characteristic $0$. Then we are going to use the notation $\Delta$ to denote a corresponding commutative multiplicative group such that we have the following:
\begin{displaymath}
|k^\times|\subset \Delta \subset \mathbb{R}^\times_+.	
\end{displaymath}

\indent Recall from \cite[Part I Chapter I.3]{6DFN} we have the corresponding notation of $k_\Delta$-affinoid algebras and $k_\Delta$-affinoid spaces. These are basically the corresponding non-strictly generalization of the corresponding strictly affinoid algebras and the corresponding spectrum. For instance the corresponding $k_\Delta$-affinoid algebras are defined to be those taking the corresponding form coming form the admissible quotient from some algebra taking the form of:
\begin{displaymath}
k\{T_1/x_1,...,T_d/x_d\}	
\end{displaymath}
for $x_i\in \Delta$ and $i=1,...,d$. Then we have the corresponding notion from this of the affinoid spaces from some suitable geometric point of view, we will use the language of adic spaces as in \cite{6KL16}. So we define:

\begin{definition}\mbox{\bf{(After Temkin, \cite[Part I Chapter I.3]{6DFN})}}
We define the corresponding $k_\Delta$-affinoid spaces as the corresponding adic space associated some $k_\Delta$-affinoid algebra $R$. Note that this is endowed with the structure sheaf $\mathcal{O}_{\mathrm{Spa}(R,R^+)}$. 
\end{definition}

\indent Then as in \cite[Part I Chapter I]{6DFN} we can define the corresponding general $k_\Delta$-analytic spaces as in the following:

\begin{definition}\mbox{\bf{(After Temkin, \cite[Part I Chapter I.3]{6DFN})}} 
We define a general $k_\Delta$-analytic space to be an adic space $X$ which could be covered by an atlas forming by the corresponding $k_\Delta$-analytic affinoid spaces. 
	
\end{definition}

\indent We would like to study the corresponding properties of the corresponding sheaves over a general $k_\Delta$-analytic space $X$ as in \cite{6KL16}. And certainly eventually we contact with the corresponding Hodge-Iwasawa theory. But we would like to further make some further more detailed discussion around the current generalized context.

\begin{remark}
We should mention that we reproduce many key arguments following \cite[Chapter 8]{6KL16} in the following presentation of the parallel results to \cite[Chapter 8]{6KL16}. However, if one really wants to work with some larger base field (by taking the valuation group even up to $\mathbb{R}_+^\times$), then one should be able to simplified many arguments below.	
\end{remark}

\begin{remark} \label{remark6.6}
The corresponding desingularization result from Temkin (as in the context of \cite[Remark 8.1.2]{6KL16}) in our situation namely we can desingularize a corresponding $k_\Delta$-affinoid algebra $A$ by the corresponding excellence of this algebra as in \cite[Part I Chapter I Fact 3.1.2.1]{6DFN}. Namely in the category of schemes we can choose a map with smooth domain $Y\rightarrow \mathrm{Spec}A$ which will have the corresponding analytification in the categories of all the $k_\Delta$-affinoid algebras, again with smooth domain.	
\end{remark}

\indent Then we consider the corresponding correspondence with the Ax-Sen-Tate style result in our context by using the corresponding pro-\'etale site. So we use the notation $f_\text{pro\'et}$ to denote the corresponding morphism $f_\text{pro\'et}:X_\text{pro\'et}\rightarrow X$.

\begin{proposition} \mbox{\bf{(After Kedlaya-Liu, \cite[Theorem 8.2.3]{6KL16})}}
Suppose we have that the $k_\Delta$-analytic space $X$ is seminormal. Then we have the corresponding isomorphism $\mathcal{O}_X\overset{\sim}{\rightarrow}f_{\text{pro\'et},*}\widehat{\mathcal{O}}_X$. On the other hand this isomorphism will also imply that the corresponding space $X$ is essentially seminormal.	
\end{proposition}

\begin{proof}
We follow the proof of \cite[Theorem 8.2.3]{6KL16}, which proves the smooth situation first and then proves the corresponding general cases by considering the corresponding desingularization. If $k$ is perfectoid and the space $X$ is smooth and affinoid, which as in \cite[Theorem 8.2.3]{6KL16} could be directly proved by relating to the corresponding restricted toric towers, then the corresponding result follows from \cite[Theorem 8.2.3]{6KL16}. Then to remove the corresponding current assumption on the prefectoidness we just consider the corresponding prefectoidization $k_1$ of $k$ and the corresponding uniform completion of $k_2:=k_1\otimes_k k_1$, and we consider the base changes of $X$ (affinoid) to the corresponding $k_1$ and $k_2$ to form the following short exact sequence:
\[
\xymatrix@R+0pc@C+0pc{
0\ar[r] \ar[r] \ar[r] &X    \ar[r] \ar[r] \ar[r]  &X\otimes k_1 \ar[r] \ar[r] \ar[r] & X\otimes k_2.
}
\]
Then as in \cite[Theorem 8.2.3]{6KL16} one further reduces to the corresponding situation around the fields. Finally to consider the general space $X$ we consider the corresponding desingularization as in \cref{remark6.6} to consider some smooth space $Y$ which maps under $h$ say to $X$, then by considering the corresponding pullbacks of perfectoid subdomain one can basically up to such desingularization compare to have the desired isomorphism by the smooth case as above. After this one considers the fact that seminormality implies that we have the isomorphism $\mathcal{O}_{X}\overset{\sim}{\rightarrow} h_* \mathcal{O}_{Y}$, which finishes the proof in this case. Then conversely if we have the corresponding isomorphism in the statement of the proposition, then we choose the corresponding seminormalization $X^{\text{seminor}}$ of $X$ which then finishes the corresponding proof as in \cite[Theorem 8.2.3]{6KL16}.

\end{proof}

\begin{proposition}\mbox{\bf{(After Kedlaya-Liu, \cite[Corollary 8.2.4]{6KL16})}}
When we consider the corresponding vector bundle $V$ over $X$ which is assumed to be seminormal, we will have then the corresponding isomorphism:
\begin{displaymath}
V\overset{\sim}{\longrightarrow}f_{\text{pro\'et},*}f^*_{\text{pro\'et}}V.	
\end{displaymath}

\end{proposition}

\begin{proof}
From the previous proposition.	
\end{proof}

\begin{proposition}\mbox{\bf{(After Kedlaya-Liu, \cite[Lemma 8.2.7]{6KL16})}}
Let $X$ be seminormal, and consider any coherent subsheaf of any locally finite free $\mathcal{O}_X$-module $V$, which is denoted by $W$. We then have the following identity:
\begin{displaymath}
W=f_{\text{pro\'et},*}\mathrm{Image}(f^*_{\text{pro\'et}}W\rightarrow f^*_{\text{pro\'et}}V).	
\end{displaymath}
\end{proposition}

\begin{proof}
As in \cite[Lemma 8.2.7]{6KL16}, we basically realize the corresponding coherent subsheaf as the kernel of a map taking the form of $V\rightarrow V'$. Then we have that:
\[
\xymatrix@R+0pc@C+0pc{
0\ar[r] \ar[r] \ar[r] &\mathrm{Image}(f^*_{\text{pro\'et}}W\rightarrow f^*_{\text{pro\'et}}V)    \ar[r] \ar[r] \ar[r]  &f^*_{\text{pro\'et}}V \ar[r] \ar[r] \ar[r] & f^*_{\text{pro\'et}}V'.
}
\]
is exact sequence. Then the previous proposition directly show the corresponding identity after taking the corresponding pushforward along the canonical morphism for the corresponding pro-\'etale site of $X$.
\end{proof}

\begin{proposition}\mbox{\bf{(After Kedlaya-Liu, \cite[Lemma 8.2.8]{6KL16})}} Now we consider the following correspondence. Consider a corresponding locally finite free $\mathcal{O}_X$-module $V$ over $X$ which is assumed to be seminormal, and also consider the corresponding locally finite free $\mathcal{O}_X$-sheaf $f^*_{\text{pro\'et}}V$. Then the corresponding functor $W\rightarrow \mathrm{Image}(f^*_{\text{pro\'et}}W\rightarrow f^*_{\text{pro\'et}}V)$ will establish an equivalence between the coherent subsheaves of $V$ and the coherent subsheaves of its pullbacks.
	
\end{proposition}

\begin{proof}
One can finish the proof as in the proof of \cite[Lemma 8.2.8]{6KL16}. To recall the argument, we briefly present the proof for the readers. The corresponding fully faithfulness comes from the previous proposition. For the essential surjectivity as in \cite[Lemma 8.2.8]{6KL16} we first consider the corresponding smooth case, which could be then reduced to the \cite[Lemma 8.2.6]{6KL16} by relating to the corresponding restricted toric tower. In general the argument is by considering any decreasing sequence of ideal of $\mathcal{O}_X$ whose pullback will basically contain some given ideal $B$ in $\widehat{\mathcal{O}}_X$. We denote this sequence by $\{A_i\}_i$. Then for each $A_i$ choose some finite free $C_i$ which maps to $A_i$, and one can basically consider the corresponding inverse image of $B$ in the pullback of $C_i$ which could be denoted by $B_i$, which corresponds to by the previous smooth case a coherent subsheaf of the $C_i$. Now considering further pushforwarding this into $\widehat{\mathcal{O}}_X$ we will have $B_{i+1}$. Now we have if $A_{i+1}=A_i$ for some $i\geq 1$ then we have that $B=B_i$. However this will eventually be the case by the corresponding noetherian property. 	
\end{proof}

\begin{proposition}\mbox{\bf{(After Kedlaya-Liu)}}\\
\mbox{\bf{(\cite[Proposition 8.2.9, Corollary 8.2.10, Corollary 8.2.11]{6KL16})}}
A. When we have that $X$ is seminormal, we will have that the corresponding functor $f^*_{\text{pro\'et}}$ establishes embedding from the category of coherent sheaves over $\mathcal{O}_X$ to the category of pseudocoherent sheaves over $\widehat{\mathcal{O}}_X$, which is fully faithful and exact. B. And for any coherent sheaf $F$ over any analytic $k_\Delta$-space, then the tensoring with $F$ will then establish a functor on the category of all the pseudocoherent sheaves over $\widehat{\mathcal{O}}_X$ which is an endofunctor, here $X$ could be allowed to be more general by considering the corresponding seminormalization. C. Over $X$ which is again assumed to be seminormal, then we have the injective map $f^*_{\text{pro\'et}}f_{\text{pro\'et},*}V\hookrightarrow V$ from any locally finite free sheaf $V$ over $\widehat{O}_X$.  	
\end{proposition}

\begin{proof}
See \cite[Proposition 8.2.9, Corollary 8.2.10, Corollary 8.2.11]{6KL16}.	
\end{proof}

\indent Then building on the basic facts from \cite[Lemma 8.2.12, Corollary 8.2.15]{6KL16} on the corresponding restricted toric towers, we can now study the corresponding categories of the corresponding $\mathcal{O}_X$-sheaves.

\begin{setting}
We now assume we are working over $X=\mathrm{Spa}(L,L^+)$ which is seminormal and affinoid. And consider $(H_\bullet,H_\bullet^+)$ which is assumed to be a perfectoid and finite \'etale tower over the base $X$. Take any ideal $l\subset L$ we can form the corresponding quotient $\overline{L}:=L/JL$ and $\overline{\Omega}:=\overline{H}_\infty/J\overline{H}_\infty$.
\end{setting}

\indent We first have the following:

\begin{proposition}\mbox{\bf{(After Kedlaya-Liu, \cite[Lemma 8.3.3]{6KL16})}}
The map $L\rightarrow \overline{H}_\infty$ is faithfully flat. 
\end{proposition}

\begin{proof}
See \cite[Lemma 8.3.3]{6KL16}.	
\end{proof}

\begin{proposition}\mbox{\bf{(After Kedlaya-Liu, \cite[Lemma 8.3.4]{6KL16})}}
$\Gamma$-modules over the quotient $\overline{\Omega}$ satisfy the following property. The finitely generatedness of any $M$ of such corresponding modules as such could be promoted to be finite projective one if one works over $\overline{\Omega}_e$ with $M_e$, where $e$ exists as an element in $L$ which is not mapped to zero divisors in the quotients $\overline{L}$ and $\overline{\Omega}$.   	
\end{proposition}

\begin{proof}
The corresponding proof is basically parallel to \cite{6KL16}. To achieve so we consider the corresponding construction as in the following. First we consider the corresponding minimal ideals of $L$ containing the prescribed ideal $l$. Then for each of such ideal $\mathfrak{p}$ (finitely many) we look at the Fitting ideal $\mathrm{Fitt}_k$ no contained in this prime (for minimized such $k$), choose any element in $\mathrm{Fitt}_k-\mathfrak{p}$, then this will do the job that the corresponding base change of $M$ to $(L/l)_{e_\mathfrak{p}}$ will be finite projective the reasonable quotient. Then to integrate all the information along all such minimal primes we consider for each pair $\mathfrak{p}\neq \mathfrak{q}$ and choose any element $e_{\mathfrak{p}\neq \mathfrak{q}}$. Finally just put $e$ to be $\sum_{\mathfrak{p}}e_\mathfrak{p}\prod_{\mathfrak{p}\neq \mathfrak{q}}e_{\mathfrak{p}\neq \mathfrak{q}}$.	
\end{proof}

\indent The following will be crucial in the understanding of the corresponding category of all the $\widehat{\mathcal{O}}_X$-sheaves over some $k_\Delta$-analytic space $X$.

\begin{proposition}\mbox{\bf{(After Kedlaya-Liu, \cite[Proposition 8.3.5]{6KL16})}} \label{proposition615}
Suppose over the ring $\overline{H}_\infty$, we have a finite module which is carrying additionally the structure of $\Gamma$-module. Denote this by $M$. Then this is automatically pseudocoherent and more importantly \'etale-stably pseudocoherent.  
\end{proposition}

\begin{proof}
This is a nontrivial result around the $\Gamma$-modules. We recall the argument in \cite[Proposition 8.3.5]{6KL16} in the following, which will prove our result in our current generalized context. First the ideal is using the corresponding noetherian reduction along some prescribed ideal $l\subset L$. As in \cite[Proposition 8.3.5]{6KL16} we reduce ourselves to the situation where $\overline{L}$ is reduced and connected. To perform the argument, one first chooses some element $e\in L$ as in the previous proposition such that the statement of the previous proposition hold. Then consider all the corresponding $e$-torsions over $L$ which is denoted by $C$. After applying the corresponding induction hypothesis we have for $M/C$ that the corresponding $M/C/eM/C$ is then pseudocoherent (not known to be strictly or so on). And consider the fact that $(M/C)_e=M_e$ and $M/C[e]$ is trivial, we have then as in \cite[Proposition 8.3.5, Lemma 1.5.8]{6KL16}	that $M_C$ is then pseudocoherent (again not known to be strictly or so on). But then this implies the torsion part is then finitely generated, which is killed by some power of $e$ up to some power $h$ for instance. Then repeat considering using the corresponding induction hypothesis we can show that the whole module $M$ is then pseudocoherent (again not known to be strictly or so on) by considering the quotient $e^{x-1}C/e^{x}C$ for $x=0,...,h$. Then to consider the corresponding strictly pseudocoherence, we play with then the kernel along some finite free cover of $M$, which is assumed to be linear over the ring $L$. Then consider kernel and complete the kernel, then apply the corresponding previous paragraph to show that the quotient of $M$ by this kernel is pseudocoherent, which by \cite[Corollary 1.2.11]{6KL16} shows that this kernel is actually complete initially, which shows that the corresponding strictly pseudocoherent (not known to be \'etale-stably pseudocoherent). To finish relate the \'etale map out of the period ring to some step in the tower, which shows by the previous argument above we are done.
\end{proof}

\indent We now choose to consider the corresponding parallel discussion as above but for the Robba ring taking the form of $\widetilde{\Pi}_{H}$ and the corresponding related Robba rings. The order of our discussion is different from \cite{6KL16}, where we will study the corresponding categories of pseudocoherent sheaves over $\widehat{\mathcal{O}}_X$ and $\widetilde{\Pi}_{X}$  together after we finish more concrete local study. Consider at this moment $X$ which is affinoid defined over a perfectoid field $k$, and consider $(H_\bullet,H_\bullet^+)$ a finite \'etale perfectoid tower as in \cite{6KL16}.

\begin{proposition} \mbox{\bf{(After Kedlaya-Liu, \cite[Lemma 8.7.3]{6KL16})}}
Consider the ring $\widetilde{\Pi}^{[s,r]}_{H}$ and consider any $\Gamma$-stable nonzero ideal of this ring. Then we have that this is basically containing some nonzero element in $\widetilde{\Pi}_{H}$.	
\end{proposition}

\begin{proof}
Since over a corresponding restricted toric tower this is just the same as in the proof of \cite[Lemma 8.7.3]{6KL16} we will not repeat it. And then one can follow the strategy to tackle more general situation by looking at the corresponding $\widetilde{\Pi}^{[s,r]}_X$-sheaf of module associated to the quotient of $\widetilde{\Pi}^{[s,r]}_H$ by the ideal given. Then by applying the previous case one can finish the proof.
\end{proof}

\begin{proposition} \mbox{\bf{(After Kedlaya-Liu, \cite[Corollary 8.7.4]{6KL16})}}
Consider the ring $\widetilde{\Pi}^{[s,r]}_{H}$ and consider any $\Gamma$-module which is now assumed to be finite generated over this ring. Then we can find an element such that inverting this element over the base ring the module will be finite projective.	
\end{proposition}

\begin{proof}
See \cite[Corollary 8.7.4]{6KL16}.
\end{proof}

\begin{proposition} \mbox{\bf{(After Kedlaya-Liu, \cite[Proposition 8.8.9]{6KL16})}}
Under the same condition as in \cite[Proposition 8.8.9]{6KL16}, consider the ring $\widetilde{\Pi}^{[s,r]}_{H}$ and consider an element $t$ which is coprime to all the images of $t_\theta$ (see \cite[Proposition 8.8.9]{6KL16} for the definition of this element) under all the powers of the Frobenius. By modulo this element any nontrivial ideal of this ring, suppose we have that this ideal will be same under the two corresponding base change to $\widetilde{\Pi}^{[s,r]}_{H^1}$. Then we have that this is just a trivial ideal.
\end{proposition}

\begin{proof}
See \cite[Proposition 8.8.9]{6KL16}.
\end{proof}

\begin{proposition} \mbox{\bf{(After Kedlaya-Liu, \cite[Corollary 8.8.10]{6KL16})}}
Under the same condition as in \cite[Corollary 8.8.10]{6KL16}, consider the ring $\widetilde{\Pi}^{[s,r]}_{H}$ and consider an element $t$ which is coprime to all the images of $t_\theta$ (see \cite[Corollary 8.8.10]{6KL16} for the definition of this element) under all the powers of the Frobenius. Then by considering taking the corresponding quotient by this element we will have the chance to promote finitely generatedness to finite projectivity. Now to be more precise suppose over the quotient over $\widetilde{\Pi}^{[s,r]}_{H}$ by this element we have a module which is contained in a finitely generated one and which carries an action from $\Gamma$, then we have that actually this is projective.
\end{proposition}

\begin{proof}
See \cite[Corollary 8.8.10]{6KL16}.
\end{proof}

\begin{proposition}\mbox{\bf{(After Kedlaya-Liu, \cite[Proposition 8.9.2]{6KL16})}} \label{proposition6.20}
Suppose over the ring $\widetilde{\Pi}_H^{[s,r]}$ for some $0<s\leq r$, we have a finite module which is carrying additionally the structure of $\Gamma$-module. Denote this by $M$. Then this is automatically pseudocoherent and more importantly \'etale-stably pseudocoherent.  
\end{proposition}

\begin{proof}
See \cite[Proposition 8.9.2]{6KL16}. This is actually just parallel to \cref{proposition615} by using decomposition of $t$ to supercede the corresponding role of the element $e$.	
\end{proof}

\indent After these technical discussion we can now consider the corresponding categories of the corresponding pseudocoherent sheaves over the ring $\widehat{\mathcal{O}}_X$ and $\widetilde{\Pi}_X^{[s,r]}$ for $0<s\leq r$. Here now $X$ will be then a general $k_\Delta$-analytic space. The affinoid algebra $A$ over $\mathbb{Q}_p$ in the following discussion is assumed to be sousperfectoid as those considered in \cite{6KH}.

\begin{definition}\mbox{\bf{(After Kedlaya-Liu, \cite[Definition 8.4.1]{6KL16})}}
We use the corresponding notation $D_{\mathrm{pseudo},\widehat{\mathcal{O}}_X}$ to denote the corresponding category of all the pseudocoherent $\widehat{\mathcal{O}}_X$ sheaves over $X$ (certainly under the corresponding pro-\'etale site), which is regarded as a subcategory of the corresponding category of all the corresponding sheaves over $\widehat{\mathcal{O}}_X$, in the corresponding pro-\'etale topology.	
\end{definition}

\begin{proposition} \mbox{\bf{(After Kedlaya-Liu, \cite[Theorem 8.4.3]{6KL16})}}
The corresponding category $D_{\mathrm{pseudo},\widehat{\mathcal{O}}_X}$ is abelian, where the corresponding kernels and cokernels are compatible with the corresponding category of all the sheaves of $\widehat{\mathcal{O}}_X$-modules over the corresponding pro-\'etale site.	
\end{proposition}

\begin{proof}
This is a direct consequence of the corresponding properties in \cref{proposition615}.	
\end{proof}

\begin{definition}\mbox{\bf{(After Kedlaya-Liu, \cite[Definition 8.10.1]{6KL16})}}
We use the corresponding notation $D_{\mathrm{pseudo},\widetilde{\Pi}^{[s,r]}_{X,A}}$ ($0<s\leq r/p^{a}$) to denote the corresponding category of all the pseudocoherent $\widetilde{\Pi}^{[s,r]}_{X,A}$ sheaves, which is regarded as a subcategory of the corresponding category of all the corresponding sheaves over $\widetilde{\Pi}^{[s,r]}_{X,A}$, in the corresponding pro-\'etale topology.	
\end{definition}

\indent Also we have the following category of all the $(\varphi^a,\Gamma)$ over the space $X$ as in the above:

\begin{definition}\mbox{\bf{(After Kedlaya-Liu, \cite[Definition 8.10.5]{6KL16})}}
We use the corresponding notation $D_{\mathrm{pseudo},\varphi^a,\widetilde{\Pi}_{X,A}}$ to denote the corresponding category of all the pseudocoherent $\widetilde{\Pi}_{X,A}$ sheaves carrying the corresponding Frobenius operator $\varphi^a$, which is regarded as a subcategory of the corresponding category of all the corresponding sheaves over $\widetilde{\Pi}_{X,A}$, in the corresponding pro-\'etale topology.	
\end{definition}

\begin{proposition} \mbox{\bf{(After Kedlaya-Liu, \cite[Theorem 8.10.2]{6KL16})}} \label{6proposition6.4.26}
The corresponding category $D_{\mathrm{pseudo},\widetilde{\Pi}^{[s,r]}_{X,A}}$ is abelian, where the corresponding kernels and cokernels are compatible with the corresponding category of all the sheaves of $\widetilde{\Pi}^{[s,r]}_{X,A}$-modules over the corresponding pro-\'etale site. Here $A$ is just $\mathbb{Q}_p$.	
\end{proposition}

\begin{proof}
This is a direct consequence of the corresponding properties in \cref{proposition6.20}.	
\end{proof}

\begin{proposition} \mbox{\bf{(After Kedlaya-Liu, \cite[Theorem 8.10.6]{6KL16})}} \label{6proposition6.4.27}
The corresponding category $D_{\mathrm{pseudo},\varphi^a,\widetilde{\Pi}_{X,A}}$ is abelian, where the corresponding kernels and cokernels are compatible with the corresponding category of all the sheaves of $\widetilde{\Pi}_{X,A}$-modules over the corresponding pro-\'etale site. Here $A$ is just $\mathbb{Q}_p$.	
\end{proposition}

\begin{proof}
See the previous proposition.	
\end{proof}

\begin{proposition} \mbox{\bf{(After Kedlaya-Liu, \cite[Theorem 8.10.6]{6KL16})}} \label{6proposition6.4.28}\\
When $X$ is smooth. The corresponding category $D_{\mathrm{pseudo},\varphi^a,\widetilde{\Pi}_{X,A}}$ is abelian, where the corresponding kernels and cokernels are compatible with the corresponding category of all the sheaves of $\widetilde{\Pi}_{X,A}$-modules over the corresponding pro-\'etale site. 	
\end{proposition}

\begin{proof}
See \cref{proposition6.3}.	
\end{proof}

\subsection{Contact with Arithmetic Riemann-Hilbert Correspondence}

We now consider the corresponding setting where the base field is finite extension of $\mathbb{F}_p((t))$ or $\mathbb{Q}_p$. The affinoid algebra $A$ over $\mathbb{F}_p((t))$ or $\mathbb{Q}_p$ in the following discussion is assumed to be sousperfectoid as those considered in \cite{6KH}. And we consider the corresponding setting above on the corresponding Robba rings. \\

We would like to take this chance to develop a little bit on the corresponding Hodge-Iwasawa $B$-pairs in our context. This was considered by \cite{6KP} with some restriction in the representation theoretic context. Since we would like to discuss the corresponding some tower where the internal structure is a little bit complicated to handle (although this is definitely doable). Note that one might want also to have the chance to consider more general tower contacting with the corresponding Iwasawa theory such as the Lubin-Tate as proposed in \cite{6KP}.

%%%%%%%%%%%%%%%%%%%%%%%%%%%%%%%%%%%%%%%!!!!!!!!!!!!!!!!!!!

\begin{definition} \mbox{\bf{(After Kedlaya-Pottharst \cite[Definition 2.17]{6KP})}}
Now we first consider more general rigid analytic spaces. Let $X$ be a rigid analytic space over a finite extension $E$ of $\mathbb{Q}_p$ or $\mathbb{F}_p((t))$. Now we consider the corresponding sheaves $\mathbb{B}_{\mathrm{dR},X,A}$, $\mathbb{B}_{e,X,A}$ and $\mathbb{B}^+_{\mathrm{dR},X,A}$ which are defined in the following way. For each of the corresponding perfectoid subdomain $Y$ of the space $X$ in our situation the corresponding rings $\mathbb{B}_{\mathrm{dR},X,A}(Y)$, $\mathbb{B}_{e,X,A}(Y)$ and $\mathbb{B}^+_{\mathrm{dR},X,A}(Y)$ are defined as in the following. So the corresponding map $Y_A \rightarrow \mathrm{Proj} P_{X,A}(Y)$ will by considering the corresponding zero locus define a closed subscheme in the corresponding deformed schematic Fargues-Fontaine curve by $A$. Here in more detail this is just the corresponding map $Y_A \rightarrow \mathrm{Proj} P_{Y^\flat,A}$.
Then we consider the corresponding open complement subspace by removing this subscheme, whose coordinate ring will be the ring $\mathbb{B}_{e,X,A}(Y)$, then we take the corresponding completion of the corresponding deformed schematic Fargues-Fontaine curve along the corresponding closed subscheme to define the corresponding space whose coordinate ring is $\mathbb{B}^+_{\mathrm{dR},X,A}(Y)$, finally taking the corresponding space corresponding to the fiber product of the open affine subscheme and this completion we can take the corresponding coordinate ring which is denoted by $\mathbb{B}_{\mathrm{dR},X,A}(Y)$. Then we take the corresponding graded rings we can define the sheaves $\mathbb{B}^+_{\mathrm{HT},X,A}$ and $\mathbb{B}_{\mathrm{HT},X,A}$ as in \cite[Definition 8.6.5]{6KL16}. 
\end{definition}

\begin{definition}\mbox{}\\ \mbox{\bf{(After Tan-Tong, Shimizu \cite[Section 2.1]{6TT}, \cite[Section 2]{6Shi})}}
In the very parallel way one can consider the corresponding cristalline period sheaves (by taking the corresponding completion in the corresponding de Rham periods above) $\mathbb{B}^+_{\mathrm{cris},X,A}$ and $\mathbb{B}_{\mathrm{cris},X,A}$ as in \cite[Section 2.1]{6TT}. In the very parallel way one can consider the corresponding semistable period sheaves $\mathbb{B}^+_{\mathrm{st},X,A}$ and $\mathbb{B}_{\mathrm{st},X,A}$ as in \cite[Section 2]{6Shi}.
\end{definition}

%%%%%%%%%%%%%%%%%%%%%%%%%%%%%%%%%%%%%%%%%%%%%%%%!!!!!!!!!!!!!!

Then we have the following generalization of Bloch-Kato \cite{6BK1} and Nakamura's \cite{6Nakamura1} fundamental sequences in rigid family:

\begin{lemma}  
We then have the corresponding Bloch-Kato fundamental sequence which is defined by the following short exact sequence:
\[
\xymatrix@R+0pc@C+0pc{
0\ar[r] \ar[r] \ar[r] &\mathbb{B}_{e,X,A}\bigcap \mathbb{B}^+_{\mathrm{dR},X,A}    \ar[r] \ar[r] \ar[r]  &\mathbb{B}_{e,X,A}\bigoplus \mathbb{B}^+_{\mathrm{dR},X,A} \ar[r] \ar[r] \ar[r] &\mathbb{B}_{e,X,A}\bigcup \mathbb{B}^+_{\mathrm{dR},X,A} \ar[r] \ar[r] \ar[r] &0 .
}
\]
\end{lemma}

\begin{proof}
This comes from the corresponding construction directly, for instance see \cite[Definition 4.8.2]{6KL16}.	
\end{proof}

\begin{definition}\mbox{}\\ \mbox{\bf{(After Scholze, Kedlaya-Liu \cite[Definition 6.8]{6Sch1}, \cite[Definition 8.6.5]{6KL16})}}\\
Now we first consider more general rigid analytic spaces.  Let $X$ be a rigid analytic space over a finite extension $E$ of $\mathbb{Q}_p$ or $\mathbb{F}_p((t))$. Now we consider the corresponding sheaves $\mathcal{O}\mathbb{B}_{\mathrm{dR},X,A}$, $\mathcal{O}\mathbb{B}_{e,X,A}$ and $\mathcal{O}\mathbb{B}^+_{\mathrm{dR},X,A}$ which are defined in the following way. $\mathcal{O}\mathbb{B}^+_{\mathrm{dR},X,A}$ is the sheafication of the following presheaf. For any perfectoid subdomain taking the form of $(P_i)_i$ as rings we consider the corresponding product $P^+_i\otimes_{W(\kappa_E)}W(P^{\flat+})$, then take the corresponding $\varpi$-adic completion and invert the corresponding $\varphi$. Then we consider the corresponding completion along the ring $P_i$ and the kernel of the map $\theta$. Then taking the limit throughout all $i$ we can define the corresponding ring $\mathcal{O}\mathbb{B}^+_{\mathrm{dR},X,A}((P_i)_i)$. Then by considering the base change to $\mathbb{B}^+_{\mathrm{dR},X,A}$, one can define the corresponding ring $\mathcal{O}\mathbb{B}^+_{\mathrm{dR},X,A}$. Then by considering the base change to $\mathbb{B}_{\mathrm{dR},X,A}$, one can define the corresponding ring $\mathcal{O}\mathbb{B}_{\mathrm{dR},X,A}$. Then similarly one can define $\mathcal{O}\mathbb{B}_{e,X,A}$. Then taking the corresponding graded pieces we have $\mathcal{O}\mathbb{B}_{\mathrm{HT},X,A}$ and $\mathcal{O}\mathbb{B}^+_{\mathrm{HT},X,A}$ as in \cite[Definition 8.6.5]{6KL16}. 
\end{definition}

\begin{definition}\mbox{}\\ \mbox{\bf{(After Tan-Tong, Shimizu \cite[Section 2.1]{6TT}, \cite[Section 2]{6Shi})}}
In the very parallel way one can consider the corresponding cristalline period sheaves (by taking the corresponding completion in the corresponding de Rham periods above) $\mathcal{O}\mathbb{B}^+_{\mathrm{cris},X,A}$ and $\mathcal{O}\mathbb{B}_{\mathrm{cris},X,A}$ as in \cite[Section 2.1]{6TT}. In the very parallel way one can consider the corresponding semistable period sheaves $\mathcal{O}\mathbb{B}^+_{\mathrm{st},X,A}$ and $\mathcal{O}\mathbb{B}_{\mathrm{st},X,A}$ as in \cite[Section 2]{6Shi}.
\end{definition}

\indent Then we have the following generalization of Bloch-Kato \cite{6BK1} and Nakamura's \cite{6Nakamura1} fundamental sequences in rigid family:

\begin{lemma}  
For $X$ smooth, we then have the corresponding Bloch-Kato fundamental sequence which is defined by the following short exact sequence:
\[
\xymatrix@R+0pc@C+0pc{
0\ar[r] \ar[r] \ar[r] &\mathcal{O}\mathbb{B}_{e,X,A}\bigcap \mathcal{O}\mathbb{B}^+_{\mathrm{dR},X,A}    \ar[r] \ar[r] \ar[r]  &\mathcal{O}\mathbb{B}_{e,X,A}\bigoplus \mathcal{O}\mathbb{B}^+_{\mathrm{dR},X,A} \ar[r] \ar[r] \ar[r] &\mathcal{O}\mathbb{B}_{e,X,A}\bigcup \mathcal{O}\mathbb{B}^+_{\mathrm{dR},X,A} \ar[r] \ar[r] \ar[r] &0 .
}
\]
\end{lemma}

\begin{proof}
In the situation where $X$ is smooth, we do have the corresponding explicit expression for the pro-\'etale sheaves, for instance for $\mathcal{O}\mathbb{B}^+_{\mathrm{dR},X,A}$ and $\mathcal{O}\mathbb{B}_{\mathrm{dR},X,A}$ one could follow \cite[Proposition 6.10]{6Sch1} to have the chance to find the corresponding expressions which represents these rings as formal power series over the rings $\mathbb{B}^+_{\mathrm{dR},X,A}$ and $\mathbb{B}_{\mathrm{dR},X,A}$. For the ring $\mathcal{O}\mathbb{B}_{e,X,A}$	one could make the parallel argument we will not repeat this again.
\end{proof}

\indent Then we can discuss the corresponding $B$-pairs:

\begin{definition} \mbox{\bf{(After Kedlaya-Liu \cite[Definition 9.3.11]{6KL15})}}\\
	We now define a $B$-pair over the corresponding triplet $(\mathbb{B}_{\mathrm{dR},X,A}, \mathbb{B}_{e,X,A}, \mathbb{B}^+_{\mathrm{dR},X,A})$ is a triplet  $(\mathbb{M}_{\mathrm{dR},X,A}, \mathbb{M}_{e,X,A}, \mathbb{M}^+_{\mathrm{dR},X,A})$ which are the corresponding finite projective objects over the corresponding rings involved above, such that that we have the following isomorphisms:
\begin{displaymath}
\mathbb{M}_{\mathrm{dR},X,A}\overset{\sim}{\longrightarrow}\mathbb{M}_{\mathrm{dR},X,A}^+\otimes_{\mathbb{B}^+_{\mathrm{dR},X,A}}\mathbb{B}_{\mathrm{dR},X,A}\overset{\sim}{\longrightarrow}\mathbb{M}_{e,X,A}^+\otimes_{\mathbb{B}_{e,X,A}} \mathbb{B}_{\mathrm{dR},X,A}.
\end{displaymath}
Similarly as in \cite[Definition 9.4.4]{6KL15} we can define the corresponding cohomology of the $B$-pairs by considering the hypercohomology of the following sequence:
\[
\xymatrix@R+0pc@C+0pc{
0\ar[r] \ar[r] \ar[r] &\mathbb{M}_{e,X,A}\bigoplus \mathbb{M}^+_{\mathrm{dR},X,A}    \ar[r] \ar[r] \ar[r]  &\mathbb{M}_{\mathrm{dR},X,A} \ar[r] \ar[r] \ar[r] &0,
}
\]
which is given by subtraction from the first coordinate by the second coordinate in the first non-zero term in this sequence.

\end{definition}

\begin{remark}
The finite projective objects mean those ones which are locally finite projective. Since we are carrying the corresponding coefficient $A$ which is essentially some nontrivial deformation, therefore locally on the space $X$ we have to further glue along $A$ in order to really consider the corresponding finite projective objects which are interesting enough in our development.\\	
\end{remark}

\indent Then we can discuss the corresponding $\mathcal{O}B$-pairs:

\begin{definition} \mbox{\bf{(After Kedlaya-Liu \cite[Definition 9.3.11]{6KL15})}}\\
	We now define a $\mathcal{O}B$-pair over the corresponding triplet $(\mathcal{O}\mathbb{B}_{\mathrm{dR},X,A}, \mathcal{O}\mathbb{B}_{e,X,A}, \mathcal{O}\mathbb{B}^+_{\mathrm{dR},X,A})$ is a triplet $(\mathbb{M}_{\mathrm{dR},X,A}, \mathbb{M}_{e,X,A}, \mathbb{M}^+_{\mathrm{dR},X,A})$ which are the corresponding finite projective objects over the corresponding rings involved above, such that that we have the following isomorphisms:
\begin{displaymath}
\mathbb{M}_{\mathrm{dR},X,A}\overset{\sim}{\longrightarrow}\mathbb{M}_{\mathrm{dR},X,A}^+\otimes_{\mathcal{O}\mathbb{B}^+_{\mathrm{dR},X,A}}\mathcal{O}\mathbb{B}_{\mathrm{dR},X,A}\overset{\sim}{\longrightarrow}\mathbb{M}_{e,X,A}^+\otimes_{\mathcal{O}\mathbb{B}_{e,X,A}} \mathcal{O}\mathbb{B}_{\mathrm{dR},X,A}.
\end{displaymath}
Similarly as in \cite[Definition 9.4.4]{6KL15} we can define the corresponding cohomology of the $B$-pairs by considering the hypercohomology of the following sequence:
\[
\xymatrix@R+0pc@C+0pc{
0\ar[r] \ar[r] \ar[r] &\mathbb{M}_{e,X,A}\bigoplus \mathbb{M}^+_{\mathrm{dR},X,A}    \ar[r] \ar[r] \ar[r]  &\mathbb{M}_{\mathrm{dR},X,A} \ar[r] \ar[r] \ar[r] &0,
}
\]
which is given by subtraction from the first coordinate by the second coordinate in the first non-zero term in this sequence.

\end{definition}

We now construct the corresponding deformed geometrization of the corresponding Bloch-Kato exponentials. To do that we consider the following derived functors which are inspired by \cite[Theorem 1.5]{6LZ} which could be also dated back to Nakamura's work \cite[2.1]{6Nakamura2} in the nonderived situation:

\begin{definition} 
For any $B$-pair we set the corresponding derived de Rham functor $D^\bullet_{\mathrm{dR}}$ and derived Hodge-Tate functor $D^\bullet_{\mathrm{HT}}$ on the $B$-pairs as in the following.
For any $B$-pair we set the corresponding derived cristalline functor $D^\bullet_{\mathrm{cris}}$ and derived semistable functor $D^\bullet_{\mathrm{st}}$ on the $B$-pairs as in the following. Consider a $B$-pair $(\mathbb{M}_{\mathrm{dR},X,A}, \mathbb{M}_{e,X,A}, \mathbb{M}^+_{\mathrm{dR},X,A})$ we have the following functors:
\begin{displaymath}
D_{\mathrm{dR}}^\bullet(\mathbb{M}_{\mathrm{dR},X,A}):=R^\bullet f_{\text{pro\'et},*} \mathbb{M}_{\mathrm{dR},X,A}\otimes_{\mathbb{B}_{\mathrm{dR},X,A}}	 \mathcal{O}\mathbb{B}_{\mathrm{dR},X,A},
\end{displaymath}
and
\begin{displaymath}
D_{\mathrm{HT}}^\bullet(\mathbb{M}_{\mathrm{dR},X,A}):=R^\bullet f_{\text{pro\'et},*} \mathbb{M}_{\mathrm{dR},X,A}\otimes_{\mathbb{B}_{\mathrm{dR},X,A}}	 \mathcal{O}\mathbb{B}_{\mathrm{HT},X,A}
\end{displaymath}
and
\begin{displaymath}
D_{\mathrm{cris}}^\bullet(\mathbb{M}_{\mathrm{dR},X,A}):=R^\bullet f_{\text{pro\'et},*} \mathbb{M}_{\mathrm{dR},X,A}\otimes_{\mathbb{B}_{\mathrm{dR},X,A}}	 \mathcal{O}\mathbb{B}_{\mathrm{cris},X,A}
\end{displaymath}
and
\begin{displaymath}
D_{\mathrm{st}}^\bullet(\mathbb{M}_{\mathrm{dR},X,A}):=R^\bullet f_{\text{pro\'et},*} \mathbb{M}_{\mathrm{dR},X,A}\otimes_{\mathbb{B}_{\mathrm{dR},X,A}}	 \mathcal{O}\mathbb{B}_{\mathrm{st},X,A}.
\end{displaymath}
We call that the $B$-pair de Rham if we have the following isomorphism:
\begin{displaymath}
f^*_\text{pro\'et}(f_{\text{pro\'et},*} \mathbb{M}_{\mathrm{dR},X,A}\otimes_{\mathbb{B}_{\mathrm{dR},X,A}}	 \mathcal{O}\mathbb{B}_{\mathrm{dR},X,A})\otimes\mathcal{O}\mathbb{B}_{\mathrm{dR},X,A} \overset{\sim}{\rightarrow} \mathbb{M}_{\mathrm{dR},X,A}\otimes_{\mathbb{B}_{\mathrm{dR},X,A}}	 \mathcal{O}\mathbb{B}_{\mathrm{dR},X,A},
\end{displaymath}
and we call this Hodge-Tate if we have the corresponding similar isomorphism. We will use the corresponding notation $\mathbb{D}$ to denote the corresponding functor on the derived category level. And we call that the $B$-pair cristalline if we have the corresponding isomorphism:
\begin{displaymath}
f^*_\text{pro\'et}(f_{\text{pro\'et},*} \mathbb{M}_{\mathrm{dR},X,A}\otimes	 \mathcal{O}\mathbb{B}_{\mathrm{cris},X,A})\otimes\mathcal{O}\mathbb{B}_{\mathrm{cris},X,A} \overset{\sim}{\rightarrow} \mathbb{M}_{\mathrm{dR},X,A}\otimes_{\mathbb{B}_{\mathrm{dR},X,A}}	 \mathcal{O}\mathbb{B}_{\mathrm{cris},X,A}.
\end{displaymath}
And we call that the $B$-pair semi-stable if we have the corresponding isomorphism:
\begin{displaymath}
f^*_\text{pro\'et}(f_{\text{pro\'et},*} \mathbb{M}_{\mathrm{dR},X,A}\otimes	 \mathcal{O}\mathbb{B}_{\mathrm{st},X,A})\otimes\mathcal{O}\mathbb{B}_{\mathrm{st},X,A} \overset{\sim}{\rightarrow} \mathbb{M}_{\mathrm{dR},X,A}\otimes_{\mathbb{B}_{\mathrm{dR},X,A}}	 \mathcal{O}\mathbb{B}_{\mathrm{st},X,A}.
\end{displaymath}	
\end{definition}

\indent Now we address some issues around the corresponding finiteness of the construction getting involved as in the above. First when the space $X$ is smooth and proper over the corresponding base fields, in the situation where $A$ is just $\mathbb{Q}_p$ we have in this case $f_\text{pro\'et,*} \mathbb{M}_{\mathrm{dR},X,A}\otimes_{\mathbb{B}_{\mathrm{dR},X,A}}	 \mathcal{O}\mathbb{B}_{\mathrm{dR},X,A}$ of any $B$-pair is coherent sheaf over $\mathcal{O}_X$, which is the consequence of \cite[Theorem 8.6.2 (a)]{6KL16}. And moreover in this case since the corresponding space is smooth as noted in \cite[Definition 10.10]{6KL3} we have that the corresponding existence of the connection $\nabla$ which will ensure the corresponding coherent sheaves are also projective modules. We now first consider the situation where we have the ring $A$ not just $\mathbb{Q}_p$ in more general situation. Certainly to study the corresponding $B$-pairs it will be basically convenient to relate to the corresponding $(\varphi,\Gamma)$-modules as in the earlier work from Kedlaya-Liu and Nakamura. To do so we can first consider the following higher dimensional generalization of \cite[Theorem 2.18]{6KP}:

\begin{remark}
Note that at this moment here we assume that the corresponding base fields are of 0 characteristic.	
\end{remark}

\begin{proposition} \mbox{\bf{(After Kedlaya-Pottharst \cite[Theorem 2.18]{6KP})}}
We look at the two groups of objects. In the first consideration we look at the corresponding finite locally free sheaves over the corresponding $A$-deformed version of the schematic Fargues-Fontaine curves attache to any local perfectoid chart for some ring $R$. In the other second consideration we look at the corresponding $A$-deformed version of the corresponding $B$-pairs defined above but after taking the corresponding section over some perfectoid domain associated to some ring $R$. Then we have the corresponding equivalence between the corresponding categories respecting the two groups of objects.	
\end{proposition}

\begin{proof}
As in \cite[Theorem 2.18]{6KP} this is by Beauville-Laszlo descent.	
\end{proof}

\begin{remark}
\indent Then we could basically translate the corresponding results on the cohomology of $B$-pairs to those on the corresponding cohomology of the vector bundles over the Fargues-Fontaine curves and further more to those on that of the Frobenius modules over the pro-\'etale site. Then we have that this will be further related to $(\varphi,\Gamma)$-modules over the imperfect Robba rings with respect to the corresponding base spaces. So we move on to consider the corresponding $(\varphi,\Gamma)$-modules over the corresponding pro-\'etale site.	
\end{remark}

\indent Following Nakamura one can further generalize the corresponding Bloch-Kato long exact sequence which provides the corresponding exponentials:

\begin{definition} \mbox{\bf{(After Bloch-Kato-Nakamura)}}
For $X$ smooth, consider any $(\varphi,\Gamma)$-module $M$ over the pro-\'etale sheaf $\widetilde{\Pi}^\infty_{X,A}$ we then take the corresponding pro-\'etale invariance, namely we consider the following two functor with the one $D_\mathrm{dR}$:
\begin{displaymath}
D^\bullet_\mathcal{O}(M):= R^\bullet f_{\text{pro\'et},*}(M\otimes_{\widetilde{\Pi}^\infty_{X,A}} (\mathcal{O}\mathbb{B}_{e,X,A}\bigcap \mathcal{O}\mathbb{B}^+_{\mathrm{dR},X,A})),	
\end{displaymath}
\begin{displaymath}
D^\bullet_\oplus(M):= R^\bullet f_{\text{pro\'et},*}(M\otimes_{\widetilde{\Pi}^\infty_{X,A}} (\mathcal{O}\mathbb{B}_{e,X,A}\bigoplus \mathcal{O}\mathbb{B}^+_{\mathrm{dR},X,A}),	
\end{displaymath}
as well as the corresponding definitions on the derived category level:
\begin{displaymath}
\mathbb{D}^\bullet_\mathcal{O}(M):= R^\bullet f_{\text{pro\'et},*}(M\otimes_{\widetilde{\Pi}^\infty_{X,A}} (\mathcal{O}\mathbb{B}_{e,X,A}\bigcap \mathcal{O}\mathbb{B}^+_{\mathrm{dR},X,A})),	
\end{displaymath}
\begin{displaymath}
\mathbb{D}^\bullet_\oplus(M):= R^\bullet f_{\text{pro\'et},*}(M\otimes_{\widetilde{\Pi}^\infty_{X,A}} (\mathcal{O}\mathbb{B}_{e,X,A}\bigoplus \mathcal{O}\mathbb{B}^+_{\mathrm{dR},X,A}),	
\end{displaymath}

Taking the corresponding complexes induced and considering the corresponding mapping cone we have the following distinguished triangle in $D(\mathcal{O}_X\widehat{\otimes}_{\mathbb{Q}_p}A)$ (the derived category of all the $\mathcal{O}_X\widehat{\otimes}_{\mathbb{Q}_p}A$-modules):
\[
\xymatrix@R+0pc@C+0pc{
C^\bullet_\mathcal{O}(M)\ar[r] \ar[r] \ar[r] &C^\bullet_\oplus(M)    \ar[r] \ar[r] \ar[r]  &C^\bullet_\mathrm{dR}(M) \ar[r] \ar[r] \ar[r] &C^\bullet_\mathcal{O}(M)[1].
}
\]
Consider the corresponding degree $0$ and $1$ we will have geometric and arithmetic family version of the corresponding Bloch-Kato-Nakamura exponential map:
\begin{displaymath}
h^0(C^\bullet_\mathrm{dR}(M))\rightarrow h^1(C^\bullet_\mathcal{O}(M)),	
\end{displaymath}
as the map of sheaves over $\mathcal{O}_X\widehat{\otimes}_{\mathbb{Q}_p}A$. Here we assume that the corresponding $B$-pairs getting involved are de Rham. 

\end{definition}

\begin{definition} \mbox{\bf{(After Kedlaya-Liu \cite[Definition 8.6.5]{6KL16})}} 
For any $\varphi$-module $M$ over the sheaf of ring $\widetilde{\Pi}^\infty_{X,A}$ we set the corresponding derived de Rham functor $D^\bullet_{\mathrm{dR}}$ and derived Hodge-Tate functor $D^\bullet_{\mathrm{HT}}$ on the $\varphi$-modules $M$ over the sheaf of ring $\widetilde{\Pi}^\infty_{X,A}$ as in the following.  For any $\varphi$-module $M$ over the sheaf of ring $\widetilde{\Pi}^\infty_{X,A}$ we set the corresponding derived cristalline functor $D^\bullet_{\mathrm{cris}}$ and derived semistable functor $D^\bullet_{\mathrm{st}}$ on the $\varphi$-modules $M$ over the sheaf of ring $\widetilde{\Pi}^\infty_{X,A}$ as in the following. Consider a $\varphi$-module $M$ over the sheaf of ring $\widetilde{\Pi}^\infty_{X,A}$ we have the following functors:
\begin{displaymath}
D_{\mathrm{dR}}^\bullet(M):=R^\bullet f_{\text{pro\'et},*} M\otimes_{\widetilde{\Pi}^\infty_{X,A}}	 \mathcal{O}\mathbb{B}_{\mathrm{dR},X,A},
\end{displaymath}
and
\begin{displaymath}
D_{\mathrm{HT}}^\bullet(M):=R^\bullet f_{\text{pro\'et},*} M \otimes_{\widetilde{\Pi}^\infty_{X,A}}	 \mathcal{O}\mathbb{B}_{\mathrm{HT},X,A}
\end{displaymath}
and
\begin{displaymath}
D_{\mathrm{cris}}^\bullet(M):=R^\bullet f_{\text{pro\'et},*} M\otimes_{\widetilde{\Pi}^\infty_{X,A}}	 \mathcal{O}\mathbb{B}_{\mathrm{cris},X,A},
\end{displaymath}
and
\begin{displaymath}
D_{\mathrm{st}}^\bullet(M):=R^\bullet f_{\text{pro\'et},*} M \otimes_{\widetilde{\Pi}^\infty_{X,A}}	 \mathcal{O}\mathbb{B}_{\mathrm{st},X,A}.
\end{displaymath}
We call that the $B$-pair de Rham if we have the following isomorphism:
\begin{displaymath}
f^*_\text{pro\'et}(f_{\text{pro\'et},*} M\otimes_{\widetilde{\Pi}^\infty_{X,A}}	 \mathcal{O}\mathbb{B}_{\mathrm{dR},X,A})\otimes\mathcal{O}\mathbb{B}_{\mathrm{dR},X,A} \overset{\sim}{\rightarrow} M \otimes_{\widetilde{\Pi}^\infty_{X,A}}	 \mathcal{O}\mathbb{B}_{\mathrm{dR},X,A},
\end{displaymath}
and we call this Hodge-Tate if we have the corresponding similar isomorphism. We will use the corresponding notation $\mathbb{D}$ to denote the corresponding functor on the derived category level. And we define the corresponding properties of being cristalline and semistable in the same parallel way.
	
\end{definition}

\begin{proposition} \label{proposition4.42}
For general rigid analytic space $X$ (note that at this moment we are in the 0 characteristic situation), the corresponding cohomology $D^i_\mathrm{dR}(M)$ is coherent sheaf over $X$ for each $i$. In the situation where the space is smooth we have the further projectivity. 
\end{proposition}

\begin{proof}
This is a relative version of the result of \cite[below Definition 10.10]{6KL3}. Namely we consider the corresponding resolution of singularities and reduce to the corresponding smooth case. In the corresponding smooth case we further reduce to the corresponding case for the base of a corresponding toric tower in the restricted sense. Then we relate to pseudocoherent $\Gamma$-modules over $\overline{H}_\infty\widehat{\otimes}A$ in our context, which further relates to the the corresponding pseudocoherent $\Gamma$-modules over $\widetilde{\Pi}_A$ and then over $H_\infty\widehat{\otimes}A$. Then we descend to some finite level ring $H_k\widehat{\otimes}A$ for some $k$. Then we need to consider as in \cite[Corollary 5.9.5]{6KL16} that the corresponding p.f. descent will allow us to realize the corresponding $\Gamma$-cohomology as the abstract one by using the higher period rings. Then this will realize the corresponding finiteness. The projectivity follows from the existence of the connection.
\end{proof}

\begin{proposition} \label{6proposition6.4.45}
For general rigid analytic space $X$ (note that at this moment we are in the 0 characteristic situation), the corresponding cohomology $D^i_\mathrm{cris}(M),D^i_\mathrm{st}(M)$ are coherent sheaves over $X$ for each $i$. In the situation where the space is smooth we have the further projectivity.	
\end{proposition}

\begin{proof}
See the proof in the previous proposition.
\end{proof}

\begin{corollary}
The corresponding sheaf
\begin{center}
 $f_{\text{pro\'et},*} \mathbb{M}_{\mathrm{dR},X,A}\otimes_{\mathbb{B}_{\mathrm{dR},X,A}}	 \mathcal{O}\mathbb{B}_{\mathrm{dR},X,A}$\end{center}
over $\mathcal{O}_X\widehat{\otimes}A$ is basically a coherent sheaf.	When we have that the corresponding space $X$ is smooth we could further have the corresponding projectivity in our context.  
\end{corollary}

\begin{remark}
The construction above is generalization of \cite[Definition 8.6.5]{6KL16}, which is different slightly from \cite[Section 3.2]{6LZ}. But we can also push things down to the \'etale site along $g_\text{pro\'et}:X_\text{pro\'et}\rightarrow X_\text{\'et}$.	
\end{remark}

\begin{definition} 
For any $B$-pair we set the corresponding derived de Rham functor $E^\bullet_{\mathrm{dR}}$ and derived Hodge-Tate functor $E^\bullet_{\mathrm{HT}}$ on the $B$-pairs as in the following.
For any $B$-pair we set the corresponding derived cristalline functor $E^\bullet_{\mathrm{cris}}$ and derived semistable functor $E^\bullet_{\mathrm{st}}$ on the $B$-pairs as in the following. Consider a $B$-pair $(\mathbb{M}_{\mathrm{dR},X,A}, \mathbb{M}_{e,X,A}, \mathbb{M}^+_{\mathrm{dR},X,A})$ we have the following functors:
\begin{displaymath}
E_{\mathrm{dR}}^\bullet(\mathbb{M}_{\mathrm{dR},X,A}):=R^\bullet g_{\text{pro\'et},*} \mathbb{M}_{\mathrm{dR},X,A}\otimes_{\mathbb{B}_{\mathrm{dR},X,A}}	 \mathcal{O}\mathbb{B}_{\mathrm{dR},X,A},
\end{displaymath}
and
\begin{displaymath}
E_{\mathrm{HT}}^\bullet(\mathbb{M}_{\mathrm{dR},X,A}):=R^\bullet g_{\text{pro\'et},*} \mathbb{M}_{\mathrm{dR},X,A}\otimes_{\mathbb{B}_{\mathrm{dR},X,A}}	 \mathcal{O}\mathbb{B}_{\mathrm{HT},X,A}
\end{displaymath}
and
\begin{displaymath}
E_{\mathrm{cris}}^\bullet(\mathbb{M}_{\mathrm{dR},X,A}):=R^\bullet g_{\text{pro\'et},*} \mathbb{M}_{\mathrm{dR},X,A}\otimes_{\mathbb{B}_{\mathrm{dR},X,A}}	 \mathcal{O}\mathbb{B}_{\mathrm{cris},X,A}
\end{displaymath}
and
\begin{displaymath}
E_{\mathrm{st}}^\bullet(\mathbb{M}_{\mathrm{dR},X,A}):=R^\bullet g_{\text{pro\'et},*} \mathbb{M}_{\mathrm{dR},X,A}\otimes_{\mathbb{B}_{\mathrm{dR},X,A}}	 \mathcal{O}\mathbb{B}_{\mathrm{st},X,A}.
\end{displaymath}
We call that the $B$-pair de Rham if we have the following isomorphism:
\begin{displaymath}
g^*_\text{pro\'et}(g_{\text{pro\'et},*} \mathbb{M}_{\mathrm{dR},X,A}\otimes_{\mathbb{B}_{\mathrm{dR},X,A}}	 \mathcal{O}\mathbb{B}_{\mathrm{dR},X,A})\otimes\mathcal{O}\mathbb{B}_{\mathrm{dR},X,A} \overset{\sim}{\rightarrow} \mathbb{M}_{\mathrm{dR},X,A}\otimes_{\mathbb{B}_{\mathrm{dR},X,A}}	 \mathcal{O}\mathbb{B}_{\mathrm{dR},X,A},
\end{displaymath}
and we call this Hodge-Tate if we have the corresponding similar isomorphism. We will use the corresponding notation $\mathbb{E}$ to denote the corresponding functor on the derived category level. And we call that the $B$-pair cristalline if we have the corresponding isomorphism:
\begin{displaymath}
g^*_\text{pro\'et}(g_{\text{pro\'et},*} \mathbb{M}_{\mathrm{dR},X,A}\otimes	 \mathcal{O}\mathbb{B}_{\mathrm{cris},X,A})\otimes\mathcal{O}\mathbb{B}_{\mathrm{cris},X,A} \overset{\sim}{\rightarrow} \mathbb{M}_{\mathrm{dR},X,A}\otimes_{\mathbb{B}_{\mathrm{dR},X,A}}	 \mathcal{O}\mathbb{B}_{\mathrm{cris},X,A}.
\end{displaymath}
And we call that the $B$-pair semi-stable if we have the corresponding isomorphism:
\begin{displaymath}
g^*_\text{pro\'et}(g_{\text{pro\'et},*} \mathbb{M}_{\mathrm{dR},X,A}\otimes	 \mathcal{O}\mathbb{B}_{\mathrm{st},X,A})\otimes\mathcal{O}\mathbb{B}_{\mathrm{st},X,A} \overset{\sim}{\rightarrow} \mathbb{M}_{\mathrm{dR},X,A}\otimes_{\mathbb{B}_{\mathrm{dR},X,A}}	 \mathcal{O}\mathbb{B}_{\mathrm{st},X,A}.
\end{displaymath}	
\end{definition}

\begin{definition} \mbox{\bf{(After Kedlaya-Liu \cite[Definition 8.6.5]{6KL16})}} 
For any $\varphi$-module $M$ over the sheaf of ring $\widetilde{\Pi}^\infty_{X,A}$ we set the corresponding derived de Rham functor $E^\bullet_{\mathrm{dR}}$ and derived Hodge-Tate functor $E^\bullet_{\mathrm{HT}}$ on the $\varphi$-modules $M$ over the sheaf of ring $\widetilde{\Pi}^\infty_{X,A}$ as in the following.  For any $\varphi$-module $M$ over the sheaf of ring $\widetilde{\Pi}^\infty_{X,A}$ we set the corresponding derived cristalline functor $E^\bullet_{\mathrm{cris}}$ and derived semistable functor $E^\bullet_{\mathrm{st}}$ on the $\varphi$-modules $M$ over the sheaf of ring $\widetilde{\Pi}^\infty_{X,A}$ as in the following. Consider a $\varphi$-module $M$ over the sheaf of ring $\widetilde{\Pi}^\infty_{X,A}$ we have the following functors:
\begin{displaymath}
E_{\mathrm{dR}}^\bullet(M):=R^\bullet g_{\text{pro\'et},*} M\otimes_{\widetilde{\Pi}^\infty_{X,A}}	 \mathcal{O}\mathbb{B}_{\mathrm{dR},X,A},
\end{displaymath}
and
\begin{displaymath}
E_{\mathrm{HT}}^\bullet(M):=R^\bullet g_{\text{pro\'et},*} M \otimes_{\widetilde{\Pi}^\infty_{X,A}}	 \mathcal{O}\mathbb{B}_{\mathrm{HT},X,A}
\end{displaymath}
and
\begin{displaymath}
E_{\mathrm{cris}}^\bullet(M):=R^\bullet g_{\text{pro\'et},*} M\otimes_{\widetilde{\Pi}^\infty_{X,A}}	 \mathcal{O}\mathbb{B}_{\mathrm{cris},X,A},
\end{displaymath}
and
\begin{displaymath}
E_{\mathrm{st}}^\bullet(M):=R^\bullet g_{\text{pro\'et},*} M \otimes_{\widetilde{\Pi}^\infty_{X,A}}	 \mathcal{O}\mathbb{B}_{\mathrm{st},X,A}.
\end{displaymath}
We call that the $B$-pair de Rham if we have the following isomorphism:
\begin{displaymath}
g^*_\text{pro\'et}(g_{\text{pro\'et},*} M\otimes_{\widetilde{\Pi}^\infty_{X,A}}	 \mathcal{O}\mathbb{B}_{\mathrm{dR},X,A})\otimes\mathcal{O}\mathbb{B}_{\mathrm{dR},X,A} \overset{\sim}{\rightarrow} M \otimes_{\widetilde{\Pi}^\infty_{X,A}}	 \mathcal{O}\mathbb{B}_{\mathrm{dR},X,A},
\end{displaymath}
and we call this Hodge-Tate if we have the corresponding similar isomorphism. We will use the corresponding notation $\mathbb{E}$ to denote the corresponding functor on the derived category level. And we define the corresponding properties of being cristalline and semistable in the same parallel way.
	
\end{definition}

\begin{proposition} \label{proposition4.45}
For general rigid analytic space $X$ (note that at this moment we are in the 0 characteristic situation), the corresponding cohomology $E^i_\mathrm{dR}(M)$ for each $i$ is a coherent sheaf over $X$. In the situation where the space is smooth we have the further projectivity.	
\end{proposition}

\begin{proof}
One looks at the following composition:
\begin{align}
X_{\text{pro\'et}}\rightarrow X_{\text{\'et}} \rightarrow X_\mathrm{an}.	
\end{align}
We call the second functor $f_{\text{\'et},\text{an}}$. In effect, the second one realizes the corresponding equivalence between the category of all the coherent $\mathcal{O}_{X_{\text{\'et}}}$-modules and the category of all the coherent $\mathcal{O}_{X_{\text{an}}}$-modules. We have that $f_{\text{pro\'et},*}=f_{\text{\'et},\text{an},*}g_{\text{pro\'et},*}$. Then the result follows from the corresponding \cref{proposition4.42}.
\end{proof}

\begin{proposition}  \label{6proposition6.4.51}
For general rigid analytic space $X$ (note that at this moment we are in the 0 characteristic situation), the corresponding cohomology $E^i_\mathrm{cris}(M),E^i_\mathrm{st}(M)$ for each $i$ are coherent sheaves over $X$. In the situation where the space is smooth we have the further projectivity.	
\end{proposition}

\begin{proof}
See the previous proposition.	
\end{proof}

\begin{corollary}
The corresponding sheaf
\begin{center}
 $g_{\text{pro\'et},*} \mathbb{M}_{\mathrm{dR},X,A}\otimes_{\mathbb{B}_{\mathrm{dR},X,A}}	 \mathcal{O}\mathbb{B}_{\mathrm{dR},X,A}$
\end{center}
over $\mathcal{O}_X\widehat{\otimes}A$ is basically a coherent sheaf.	When we have that the corresponding space $X$ is smooth we could further have the corresponding projectivity in our context. 
\end{corollary}

\indent Now we consider the corresponding reconstruction process as in the philosophy of \cite{6KP}. 

\begin{setting}
Namely we take the space $X$ to be smooth, and we take the corresponding $A$ to be some local chart of some abelian Fr\'echet-Stein algebra $A_\infty$ (but the local quasi-compacts are assumed to be within our sousperfectoid requirement). Take a $B$-pair $M$ now over $\mathcal{O}\mathbb{B}_{\mathrm{dR},X,A_\infty}$ namely let us require the corresponding finiteness and projectivity at the infinite level, as a family of $B$-pair over $\mathcal{O}\mathbb{B}_{\mathrm{dR},X,A_n}$ for all $n$.
\end{setting}

\begin{proposition}
Let $X$ be the space attached to $E\{T_1,...,T_n\}$, and define the functor $\mathbb{D}^\bullet_\mathrm{dR}$ as above for a family $M$ in the previous setting over $\mathcal{O}\mathbb{B}_{\mathrm{dR},X,A_\infty}$. Then we have that $\mathbb{D}^\bullet_\mathrm{dR}$ consists of coadmissible and finite projective modules over $\mathcal{O}_X\widehat{\otimes}A_\infty$ where $A_\infty$ is assumed to come from a $p$-adic Lie group as in \cite[Theorem 3.10]{6Z1} (but the local quasi-compacts are assumed to be within our sousperfectoid requirement).
\end{proposition}

\begin{proof}
Our result above for $A_n$ for each $n=0,1,...$ shows the coadmissibility. The finite projectivity follows from the existence of the connection $\nabla_X$ and a theorem due to Z\'ab\'radi \cite[Theorem 3.10]{6Z1}.	
\end{proof}

\begin{conjecture}
If the space $X$ is seminormal, then we have the same result.
\end{conjecture}

\begin{remark}
This would reconstruct, generalize and geometrize the corresponding Berger and Fourquaux's Lubin-Tate Iwasawa construction \cite{6Berger}. And in the corresponding cyclotomic situation this also provides the corresponding generalization and geometrization of the corresponding construction after Nakamura \cite{6Nakamura1}.	
\end{remark}

%%\newpage

\newpage\section{Applications to Equivariant Iwasawa Theory}

\subsection{Equivariant Big Perrin-Riou-Nakamura Exponential Maps} \label{6section5.1}

\indent Now we establish the corresponding equivariant version of Nakamura's de Rham Iwasawa theory \cite{6Nakamura1} of $(\varphi,\Gamma)$-modules over the Robba ring in the corresponding cyclotomic situation. The idea and basic philosophy come from the corresponding reconstruction process observed by the work of \cite{6KPX} and \cite{6KP}, namely the corresponding Fontaine's original construction using $\psi$-operator of the corresponding Iwasawa cohomology of $(\varphi,\Gamma)$-modules could be reconstructed by using the corresponding Iwasawa analytic deformation and even Frobenius sheaves in the style of Kedlaya-Liu \cite{6KL15} and \cite{6KL16} over the pro-\'etale sites after Scholze.\\

\indent In our previous work in \cite{6T1}, we partially looked at some sort of generalized version of the picture above (although \cite{6T1} is really and essentially aimed at the corresponding direction after \cite{6PZ1} and \cite{6CKZ}), namely instead of just deforming the corresponding $(\varphi,\Gamma)$-modules by some Iwasawa deformation over $\Pi^\infty(\Gamma)$ which is the corresponding Fr\'echet-Stein algebra attached to the group $\Gamma$, we look at more general Iwasawa deformation over $\Pi^\infty(\Gamma)\widehat{\otimes}\Pi^\infty(\Gamma')$ namely we could look at a suitable Fr\'echet-Stein algebra attached to the group $\Gamma\times \Gamma'$ such that we have that it could be written as family of affinoid algebras in rigid analytic geometry. The corresponding $(\varphi,\Gamma)$ cohomology of such modules could be regarded as the corresponding genuine $\Gamma'$-equivariant Iwasawa cohomology.\\

\indent Definitely we work in the philosophy of the corresponding Burns-Flach-Fukaya-Kato's grand picture on the so-called noncommutative Tamagawa number conjectures as in \cite{6BF1}, \cite{6BF2} and \cite{6FK}. Nakamura's original Iwasawa theory on the $(\varphi,\Gamma)$-modules happens over the tower of Iwasawa tower. Namely the corresponding big Perrin-Riou-Nakamura exponential maps should be related to the corresponding characteristic ideals of the $(\varphi,\Gamma)$-cohomology of the corresponding Iwasawa deformation of any de Rham $(\varphi,\Gamma)$-module in some very precise way. One should regard our generalization as climbing even higher in some equivariant way. Namely the corresponding big Perrin-Riou-Nakamura exponential maps should be related to the corresponding characteristic ideals of the $(\varphi,\Gamma)$-cohomology of the corresponding higher dimensional equivariant deformation of any de Rham $(\varphi,\Gamma)$-module in some very precise way.\\

\indent We now follow \cite{6Nakamura1} to establish the equivariant version of corresponding big (dual-)exponential maps after Perrin-Riou-Nakamura. This happens at the infinite level of the corresponding generalized admissible abelian tower in our setting, which generalizes the corresponding setting considered in \cite{6Nakamura1}.

%%%%%%%%%%%%%%%%%%%%%%%!!!!!!!!!!!!!!

\begin{setting}
We consider an abelian $p$-adic Lie group $\Gamma'$ such that the Fr\'echet-Stein algebra $\Pi^\infty(\Gamma')$ in the sense of \cite{6ST1} could be written as a family of affinoid algebras $\Pi^\infty(\Gamma')_k,k=0,1,...$ in the sense of rigid analytic geometry:
\begin{displaymath}
\Pi^\infty(\Gamma')\overset{\sim}{\longrightarrow} \varprojlim_k \Pi^\infty(\Gamma')_k.
\end{displaymath}
If we use the notation $X$ to denote the Fr\'echet-Stein space attached to the algebra $\Pi^\infty(\Gamma')$ which could be written as:
\begin{displaymath}
X\overset{\sim}{\longrightarrow} \varinjlim_k X_k:	
\end{displaymath}	
with:
\begin{displaymath}
\mathcal{O}_X\overset{\sim}{\longrightarrow} \varprojlim_k \mathcal{O}_{X_k}:	
\end{displaymath}
\end{setting}

\indent Now we start to consider the corresponding relative Robba rings where our Iwasawa theory happens:

\begin{definition}
Many definitions of the Robba rings are actually compatible with \cite{6Nakamura1}, the only difference is that we will deform over some general algebras. We use the notations $\Pi^r_{\mathrm{an,con},K}$ and $\Pi_{\mathrm{an,con},K}$ to denote the corresponding Robba rings as in \cite[Section 2.1]{6Nakamura1}. Now we consider the corresponding Robba rings\\ $\Pi^r_{\mathrm{an,con},K,\Pi^\infty(\Gamma_K)\widehat{\otimes}\Pi^\infty(\Gamma')}$ and $\Pi_{\mathrm{an,con},K,\Pi^\infty(\Gamma_K)\widehat{\otimes}\Pi^\infty(\Gamma')}$ in the relative sense:
\begin{align}
\Pi^r_{\mathrm{an,con},K,\Pi^\infty(\Gamma_K)\widehat{\otimes}\Pi^\infty(\Gamma')}:=\Pi^r_{\mathrm{an,con},K}\widehat{\otimes}\Pi^\infty(\Gamma_K)\widehat{\otimes}\Pi^\infty(\Gamma'),\\
\Pi_{\mathrm{an,con},K,\Pi^\infty(\Gamma_K)\widehat{\otimes}\Pi^\infty(\Gamma')}:=\Pi_{\mathrm{an,con},K}\widehat{\otimes}\Pi^\infty(\Gamma_K)\widehat{\otimes}\Pi^\infty(\Gamma').\\	
\end{align}
One also has the corresponding definitions for Robba rings with respect to some interval $[s,r]$. And be very careful that in our current development the corresponding Robba ring $\Pi^?_{\mathrm{an,con},K}$ is the ring $\mathrm{B}^?_{\mathrm{rig},K}$ where $?=\empty,r,[s,r]$, instead of the ring $\Pi^?_{\mathrm{an,con}}(\pi_K)$ considered in \cite[Definition 2.2.2]{6KPX}.	
\end{definition}

\begin{setting}
Our motivic objects are still usual $(\varphi,\Gamma)$-modules over the Robba ring $\Pi_{\mathrm{an,con},K}$. Namely we consider those modules finite locally free over this big non-noetherian ring carrying commuting semilinear action from the $\varphi$ and $\Gamma$. 	
\end{setting}

\begin{remark}
The ideas we presented here actually are initiated in \cite{6KP}. One replaces the corresponding $p$-adic Lie group $\Gamma$ by a pro-\'etale site, and considers the deformation of the sheaves (both the corresponding motivic objects and the corresponding Berger's $p$-adic differential equations attached in whatever reasonable context like Lubin-Tate setting or some else) by the corresponding Iwasawa deformations, and then considers the corresponding big exponential maps by using the connections and their duals.	
\end{remark}

\indent To study Iwasawa theory we look at those geometric $(\varphi,\Gamma)$-modules, namely those geometric ones in the sense of $p$-adic Hodge theory after Berger's celebrated thesis. We do not go beyond \cite{6Nakamura1} on this, but we would like to recall the basic construction on the corresponding $p$-adic differential equation attached to a de Rham $(\varphi,\Gamma)$-module $M$.

\begin{setting}
Recall the basic definition of the corresponding $p$-adic differential equation attache to a de Rham $(\varphi,\Gamma)$-module $M$ after Berger. As in \cite[Section 3.2]{6Nakamura1} (we will use some different notations from \cite[Section 3.2]{6Nakamura1}), we consider the following constructions. Let $M$ be a $(\varphi,\Gamma)$-module $M$ over $\Pi_{\mathrm{an,con},K}$. In our situation we work with at least de Rham ones in the sense of Berger's celebrated paper \cite{6Ber1}. Then we have the corresponding relation for any integer $k\geq l'$ for some $l'(K,M)\geq 0$ fixed which is assumed to be sufficiently large depending on $K,M$:
\begin{align}
F_{\mathrm{dif},k}(M)\overset{\sim}{\longrightarrow} F_{\mathrm{dR},K}(M)\otimes K_k((t))	
\end{align}
with the corresponding lattice $F_{\mathrm{dR},K}(M)\otimes K_k[[t]]$ inside with reasonable rank. Then we have the corresponding continuous differential:
\begin{displaymath}
\widehat{\Omega}_{\Pi_{\mathrm{an,con},K}/K_0^\mathrm{ur}}=	\widehat{\Omega}_{\Pi_{\mathrm{an,con},K}}dT.
\end{displaymath}
 
\end{setting}

The corresponding differential equation in the sense of Berger is the following construction which is denoted by $\mathrm{N}_{\mathrm{dif}}(M)$, we will denote this by $\Theta_{\mathrm{Ber,dif}}(M)$ which is actually defined in the following way. Recall as in \cite[Theorem 3.5]{6Nakamura1} by \cite[5.10]{6Ber1} and \cite[3.2.3]{6Ber2}, first we have $\Theta_{\mathrm{Ber,dif},k}(M)$ which is defined to be the subset of all the elements in $M^{(k)}[1/t]$ such that the corresponding image of them under each map $\iota_k'$ for all $k'\geq k$ lives in $F_{\mathrm{dR},K}(M)\otimes K_k'[[t]]$. Recall from \cite[Section 3.2]{6Nakamura1} we have then the corresponding set $\Theta_{\mathrm{Ber,dif}}(M)=\varinjlim_k \Theta_{\mathrm{Ber,dif},k}(M)$ which could be endowed with furthermore the structure of $(\varphi,\Gamma)$-module over the Robba ring $\Pi_{\mathrm{an,con},K}$. And furthermore we have the corresponding properties that $d_0(\Theta_{\mathrm{Ber,dif}}(M))\subset t\Theta_{\mathrm{Ber,dif}}(M)$ where the differential operator $d_0$ is defined to be $\frac{\mathrm{log}(\gamma)}{\mathrm{log}(\chi(\gamma))}$ living in in our notation $\Pi^\infty(\Gamma_K)$.\\

\indent Therefore recall from \cite[Section 3.2]{6Nakamura1} then we can define the connection $\partial:=\frac{1}{t}d_0$ mapping the corresponding differential $(\varphi,\Gamma)$-module considered above to itself. The differential structures and the corresponding Frobenius structures and the corresponding Lie group action structures are compatible in the following sense:
\begin{align}
\partial\varphi-p\varphi \partial=0,\\
\partial\gamma-\chi(\gamma)\gamma \partial=0,	
\end{align}
and:
\begin{align}
\varphi(dT)-pdT=0,\\
\gamma(dT)-\chi(\gamma)dT=0.	
\end{align}

\begin{lemma}
Now we set $d_i:=d_0-i$ for $i\geq 0$ integer, we have for any $m$-th ($m\geq 0$) Hodge filtration of $M$ (which is assumed to be de Rham) which gives the full module:
\begin{displaymath}
d_{m-1}d_{m-2}d_{m-3}...d_1d_0(\Theta_{\mathrm{Ber,dif}}(M))\subset M.	
\end{displaymath}
	
\end{lemma}

\begin{proof}
This is just \cite[Lemma 3.6]{6Nakamura1}.	
\end{proof}

\indent We now define the corresponding equivariant big exponential map after Perrin-Riou-Nakamura.

\begin{assumption}
Now we assume that there is a quotient $G_K\rightarrow \Gamma'$ in our situation among all the other requirements as above.	
\end{assumption}

\begin{definition}
Take any sufficiently large $m\geq 0$ as in the previous lemma. First we consider the differential equation $\Theta_{\mathrm{Ber,dif}}(M)$ attached to $M$ which is assumed to be de Rham. We consider the $\Gamma'$-deformation of $M$ over the algebra $\Pi^\infty(\Gamma')$. This happens by considering the corresponding quotient $G_K\rightarrow \Gamma'$	which defines a corresponding Galois representation with value in $\Pi^\infty(\Gamma')$, which by corresponding embedding fully faithful gives rise to a corresponding $(\varphi,\Gamma)$-module over the relative Robba ring $\Pi_{\mathrm{an,con},K}\widehat{\otimes}\Pi^\infty(\Gamma')$. Denote this by $\mathrm{Dfm}_{\Gamma'}$ which could be written as the projective limit of the families $\{\mathrm{Dfm}_{\Gamma',p}\}_{p\geq 0}$ over $\Pi_{\mathrm{an,con},K}\widehat{\otimes}\Pi^\infty(\Gamma')_p$ for all $p\geq 0$. Then we consider the product $\mathrm{Dfm}_{\Gamma',p}\widehat{\otimes} M$ then take the corresponding projective limit over all $p\geq 0$. Then we consider the following definition which is an equivariant version of \cite[Definition 3.7]{6Nakamura1}:
\begin{displaymath}
\mathrm{Exp}^{\Pi^\infty(\Gamma')}_{M,m}: (\varprojlim \mathrm{Dfm}_{\Gamma',p}\widehat{\otimes} \Theta_{\mathrm{Ber,dif}}(M)	)^{\psi=1}\rightarrow (\varprojlim \mathrm{Dfm}_{\Gamma',p}\widehat{\otimes} M	)^{\psi=1}
\end{displaymath}
which is just defined by extension of $d_{m-1}d_{m-2}d_{m-3}...d_1d_0$. This is now $\Pi^\infty(\Gamma_K)\widehat{\otimes}\Pi^\infty(\Gamma')$-linear instead of just $\Pi^\infty(\Gamma_K)$. And the corresponding $(\varphi,\Gamma)$-modules here could be basically geometrized by the corresponding Hodge-Iwasawa construction above. We consider the corresponding family of vector bundle $\varprojlim\mathcal{F}(M)_p$ and $\varprojlim\mathcal{F}(\Theta_{\mathrm{Ber,dif}}(M))_p$ over the corresponding adic version of Fargues-Fontaine curve deformed by $\Pi^\infty(\Gamma_K)\widehat{\otimes}\Pi^\infty(\Gamma')$. We have now that the geometrized big exponential map by taking the corresponding sheaf cohomology:
\begin{displaymath}
\mathrm{Exp}^{\Pi^\infty(\Gamma')}_{\mathcal{F},m}: R^1\Gamma_\mathrm{sheaf}(\varprojlim\mathcal{F}(\Theta_{\mathrm{Ber,dif}}(M))_p) \rightarrow R^1\Gamma_\mathrm{sheaf}(\varprojlim\mathcal{F}(M)_p),
\end{displaymath}
and the corresponding derived version:
\begin{displaymath}
\mathrm{Exp}^{\Pi^\infty(\Gamma'),\bullet}_{\mathcal{F},m}: R^\bullet\Gamma_\mathrm{sheaf}(\varprojlim\mathcal{F}(\Theta_{\mathrm{Ber,dif}}(M))_p) \rightarrow R^\bullet\Gamma_\mathrm{sheaf}(\varprojlim\mathcal{F}(M)_p).
\end{displaymath}\\
\end{definition}

\subsection{Characteristic Ideal Sheaves} \label{section5.2}

\indent Now we follow the corresponding idea in \cite{6Nakamura1} to consider the corresponding characteristic ideal and the corresponding determinant of a given derived maps with respect to the corresponding equivariant Iwasawa cohomology. For this we consider the corresponding map defined above, and we look at the corresponding equivariant Iwasawa cohomology, i.e. the corresponding more general commutative deformation in our current context.\\

\indent First in our situation we make the following assumption:

\begin{assumption}
We now assume that the corresponding deformation happens over Krull noetherian rings. This could basically allow us to define the corresponding characteristic ideal over the sheaf $\Pi^\infty(\Gamma_K)\widehat{\otimes} X_{\Pi^\infty(\Gamma')}$ of any torsion module $M$, which we will denote it by $\mathrm{char}_{\Pi^\infty(\Gamma_K)\widehat{\otimes}X_{\Pi^\infty(\Gamma')}}(M)$ which is again a sheaf of module over $\Pi^\infty(\Gamma_K)\widehat{\otimes} X_{\Pi^\infty(\Gamma')}$.	
\end{assumption}

\indent In the current generalized context, following Nakamura \cite[Section 3.4]{6Nakamura1} we define the corresponding determinant of a general map $f:M\rightarrow N$ of coadmissible modules over the sheaf $\Pi^\infty(\Gamma_K)\widehat{\otimes}X_{\Pi^\infty(\Gamma')}$. First as above we have the corresponding characteristic ideals $\mathrm{char}_{\Pi^\infty(\Gamma_K)\widehat{\otimes}X_{\Pi^\infty(\Gamma')}}(M_\mathrm{torsion})$ and $\mathrm{char}_{\Pi^\infty(\Gamma_K)\widehat{\otimes}X_{\Pi^\infty(\Gamma')}}(N_\mathrm{torsion})$. Then quotienting out the corresponding torsion part we have:
\begin{align}
\overline{f}:M/M_{\mathrm{torsion}}\rightarrow 	N/N_{\mathrm{torsion}}.
\end{align}

We will regard this as a corresponding map on the sheaves over the the space $\Pi^\infty(\Gamma_K)\widehat{\otimes}X_{\Pi^\infty(\Gamma')}$, then for any $Y=\mathrm{Sp}(\Pi^\infty(\Gamma')_p)$ quasicompact of $X_{\Pi^\infty(\Gamma')}$ we consider the corresponding ring $\Pi^\infty(\Gamma_K)\widehat{\otimes}\Pi^\infty(\Gamma')_p$, then we have the notion of determinant of $\overline{f}$ which we will denote it by $\mathrm{det}_{\Pi^\infty(\Gamma_K)\widehat{\otimes}\Pi^\infty(\Gamma')_p}(\overline{f})$. We then organize this construction to get the determinant sheaf $\mathrm{det}_{\Pi^\infty(\Gamma_K)\widehat{\otimes}X_{\Pi^\infty(\Gamma')}}(\overline{f})$. Then now we define the determinant of $f$ as:
\begin{align}
&\mathrm{det}_{\Pi^\infty(\Gamma_K)\widehat{\otimes}X_{\Pi^\infty(\Gamma')}}(f)\\
&:=\mathrm{det}_{\Pi^\infty(\Gamma_K)\widehat{\otimes}X_{\Pi^\infty(\Gamma')}}(\overline{f})\mathrm{char}_{\Pi^\infty(\Gamma_K)\widehat{\otimes}X_{\Pi^\infty(\Gamma')}}(M_\mathrm{torsion})\mathrm{char}_{\Pi^\infty(\Gamma_K)\widehat{\otimes}X_{\Pi^\infty(\Gamma')}}(N_\mathrm{torsion})^{-1}.	
\end{align}

\indent On the corresponding derived level we consider after \cite[Section 3.4]{6Nakamura1} the following construction for the big exponential map:
\begin{displaymath}
\mathrm{Exp}^{\Pi^\infty(\Gamma'),\bullet}_{\mathcal{F},m}: R^\bullet\Gamma_\mathrm{sheaf}(\varprojlim\mathcal{F}(\Theta_{\mathrm{Ber,dif}}(M))_p) \rightarrow R^\bullet\Gamma_\mathrm{sheaf}(\varprojlim\mathcal{F}(M)_p).
\end{displaymath}

For each $p$ we consider the localized map:

\begin{displaymath}
\mathrm{Exp}^{\Pi^\infty(\Gamma')_p,\bullet}_{\mathcal{F},m}: R^\bullet\Gamma_\mathrm{sheaf}(\mathcal{F}(\Theta_{\mathrm{Ber,dif}}(M))_p) \rightarrow R^\bullet\Gamma_\mathrm{sheaf}(\mathcal{F}(M)_p).
\end{displaymath}

Then what we can do is to define the corresponding determinant of the this by:
\begin{align}
\mathrm{det}_{\Pi^\infty(\Gamma_K)\widehat{\otimes}\Pi^\infty(\Gamma')_p}&(\mathrm{Exp}^{\Pi^\infty(\Gamma')_p,\bullet}_{\mathcal{F},m}):=	\\
& \mathrm{det}_{\Pi^\infty(\Gamma_K)\widehat{\otimes}\Pi^\infty(\Gamma')_p}(R^1\Gamma_\mathrm{sheaf}(\mathcal{F}(\Theta_{\mathrm{Ber,dif}}(M))_p) \rightarrow R^1\Gamma_\mathrm{sheaf}(\mathcal{F}(M)_p))\cdot\\
&\mathrm{det}_{\Pi^\infty(\Gamma_K)\widehat{\otimes}\Pi^\infty(\Gamma')_p}(R^2\Gamma_\mathrm{sheaf}(\mathcal{F}(\Theta_{\mathrm{Ber,dif}}(M))_p) \rightarrow R^2\Gamma_\mathrm{sheaf}(\mathcal{F}(M)_p))^{-1}.
\end{align}

Then we can organize this to be the corresponding determinant sheaf:

\begin{align}
&\mathrm{det}_{\Pi^\infty(\Gamma_K)\widehat{\otimes}\Pi^\infty(\Gamma')}(\mathrm{Exp}^{X_{\Pi^\infty(\Gamma')},\bullet}_{\mathcal{F},m})	\\
&= \mathrm{det}_{\Pi^\infty(\Gamma_K)\widehat{\otimes}\Pi^\infty(\Gamma')}(R^1\Gamma_\mathrm{sheaf}(\mathcal{F}(\{\Theta_{\mathrm{Ber,dif}}(M))_p\}_p) \rightarrow R^1\Gamma_\mathrm{sheaf}(\{\mathcal{F}(M)_p\}_p))\cdot\\
&\mathrm{det}_{\Pi^\infty(\Gamma_K)\widehat{\otimes}\Pi^\infty(\Gamma')}(R^2\Gamma_\mathrm{sheaf}(\mathcal{F}(\{\Theta_{\mathrm{Ber,dif}}(M))_p\}_p) \rightarrow R^2\Gamma_\mathrm{sheaf}(\{\mathcal{F}(M)_p\}_p))^{-1}\\
&= \mathrm{det}_{\Pi^\infty(\Gamma_K)\widehat{\otimes}\Pi^\infty(\Gamma')}(R^1\Gamma_\mathrm{sheaf}(\mathcal{F}(\{\Theta_{\mathrm{Ber,dif}}(M))_p\}_p) \rightarrow R^1\Gamma_\mathrm{sheaf}(\{\mathcal{F}(M)_p\}_p))\cdot\\
&\mathrm{char}_{\Pi^\infty(\Gamma_K)\widehat{\otimes}\Pi^\infty(\Gamma')}(R^2\Gamma_\mathrm{sheaf}(\mathcal{F}(\{\Theta_{\mathrm{Ber,dif}}(M))_p\}_p))^{-1}\cdot\\
&\mathrm{char}_{\Pi^\infty(\Gamma_K)\widehat{\otimes}\Pi^\infty(\Gamma')}(R^2\Gamma_\mathrm{sheaf}(\{\mathcal{F}(M)_p\}_p)). 
\end{align}\\

\indent The following conjecture reads that in our current situation we have the corresponding equivariant Iwasawa main conjecture. But we will only consider the situation where $K=\mathbb{Q}_p$ and the corresponding full admissible $p$-adic Lie extension after \cite{6LZ1} namely we consider $K_\infty=\mathbb{Q}^\mathrm{ur}_p(\mu_{p^\infty})$. We first recall Nakamura's result in \cite{6Nakamura1}:

\begin{theorem} \mbox{\bf{(Nakamura's Iwasawa Main Conjecture)}}
In our current notations, let $\Gamma'$ be trivial. Consider a de Rham $(\varphi,\Gamma)$ module $M$ over the usual Robba ring. We assume that this module has rank $g$ with $g$ Hodge-Tate weights	$a_1,...a_g$. Then we have the following formula linking the determinant of the Perrin-Riou-Nakamura big exponential map and the characteristic ideals in the equivariant way:
\begin{align}
\prod_{i}\prod^{g-a_i-1}_{b_i} &\nabla_{a_i+b_i}\\
&= \mathrm{det}_{\Pi^\infty(\Gamma_K)\widehat{\otimes}\Pi^\infty(\Gamma')}(R^1\Gamma_\mathrm{sheaf}(\mathcal{F}(\{\Theta_{\mathrm{Ber,dif}}(M))_p\}_p) \rightarrow R^1\Gamma_\mathrm{sheaf}(\{\mathcal{F}(M)_p\}_p))\cdot\\
&\mathrm{char}_{\Pi^\infty(\Gamma_K)\widehat{\otimes}\Pi^\infty(\Gamma')}(R^2\Gamma_\mathrm{sheaf}(\mathcal{F}(\{\Theta_{\mathrm{Ber,dif}}(M))_p\}_p))^{-1}\cdot\\
&\mathrm{char}_{\Pi^\infty(\Gamma_K)\widehat{\otimes}\Pi^\infty(\Gamma')}(R^2\Gamma_\mathrm{sheaf}(\{\mathcal{F}(M)_p\}_p)).	
\end{align}
\end{theorem}

Then in this situation:

\begin{conjecture} \mbox{\bf{(Equivariant Iwasawa Main Conjecture, after Nakamura)}}\\
Consider a de Rham $(\varphi,\Gamma)$ module $M$ over the usual Robba ring. We assume that this module has rank $g$ with $g$ Hodge-Tate weights	$a_1,...a_g$. Then we conjecture we have the corresponding analog of the formula linking the determinant of the Perrin-Riou-Nakamura big exponential map and the characteristic ideals in the equivariant way.
\end{conjecture}

\begin{remark}
This is also inspired by the work of \cite{6LVZ} and \cite{6LZ1}, while we work in the generality after \cite{6Nakamura1}. Although we have not mentioned explicitly but we definitely conjecture that one has the parallel statement for Lubin-Tate differential equations after \cite{6Berger} in the equivariant setting. 	
\end{remark}

\

This chapter is based on the following paper, where the author of this dissertation is the main author:
\begin{itemize}
\item Tong, Xin. "Category and Cohomology of Hodge-Iwasawa Modules." arXiv preprint arXiv:2012.07308 (2020). 
\end{itemize}

\newpage\chapter{Period Rings with Big Coefficients and Applications I}

\newpage\section{Introduction}

\subsection{Some Motivation of Big Coefficients}

\indent We initiate our development on the corresponding period rings with big coefficients and application, with the corresponding motivation partially coming from our Hodge-Iwasawa theory (see \cite{7T1}, \cite{7T2}), although we do not want to restrict ourselves to Hodge-Iwasawa theory. In our previous work \cite{7T1}, \cite{7T2} and actually also \cite{7T3}, we developed very carefully the corresponding deformation of various type of relative $p$-adic Hodge structures over rigid analytic spaces. It is not necessary to just consider things over rigid analytic spaces namely those analytic spaces which admit atlas made up of strictly affinoids. We first want to mention that one could have the chance to generalize this.\\

\indent One of our target here is the corresponding adic space coefficient consideration, which is natural if one looks at the corresponding context of \cite{7CKZ} or \cite{7PZ}, where one has some multivariate Robba rings which admit multi partial Frobenius actions. For instance when we have two variables taking quotient by one partial Frobenius will basically give us some adic space coefficients. The remaining Frobenius gives rise to some relative Frobenius Hodge structures. \\

\indent The corresponding \cite{7CKZ} actually proposed that one can also even consider the product version of the corresponding context of \cite{7KL1} and \cite{7KL2}, which is our second goal here. We want to systematically consider deforming the Frobenius modules structures of \cite{7KL1} and \cite{7KL2} after our work \cite{7T1} and \cite{7T2}. We focus on the translating the corresponding results in \cite{7T1} and \cite{7T2} to the corresponding context in the current adic space consideration, again after \cite{7KL1}, \cite{7KL2} as well. We hope to work in full generality, namely over perfectoid uniform adic Banach coefficients and those adic space coefficients forming from the quotient of the perfectoid uniform adic Banach ones. \\

\indent In the situation where we have the corresponding Kedlaya-Liu's perfect Robba rings with general coefficients in some Banach rings, one might have to impose some of the corresponding sheafiness condition. But we would like to use the $\infty$-Huber spectrum given by \cite{7BK} to tackle this from some other derived perspective where we do have the corresponding derived sheafiness. Definitely this is very complicated to manipulate, but we initiate some of the discussion though not complete at all. We hope we could come back to this systematically later at some point.\\

\indent Finally we revisit the corresponding noncommutative Hodge-Iwasawa theory we initiated in the corresponding papers \cite{7T1} and \cite{7T2}. We use some observation coming from Kedlaya to tackle the corresponding noncommutative descent for finite projective datum.\\

\subsection{Results}

\indent In this section we make some summary on our main results covered in the main body of our current paper:\\

\noindent 1. The first scope of the discussion is around some perfectoid or preperfectoid deformation of essentially the Frobenius modules and the corresponding quasi-coherent sheaves over the deformed version of the adic Fargues-Fontaine curves by perfectoids or preperfectoids. In \cref{theorem3.11}, we consider the comparison between the following objects (with the corresponding notations in \cref{theorem3.11}):\\

\noindent A. The pseudocoherent sheaves over the adic relative Fargues-Fontaine curve $Y_{\mathrm{FF},R,A}$ where $A$ is a perfectoid adic Banach uniform algebra over $E$, here we consider \'etale topology;\\
\noindent B. The pseudocoherent sheaves over the adic relative Fargues-Fontaine curve $Y_{\mathrm{FF},R,A}$ where $A$ is a perfectoid adic Banach uniform algebra over $E$, here we consider pro-\'etale topology;\\
\noindent C. The pseudocoherent modules over the relative $A$-coefficient Robba ring $\widetilde{\Pi}_{R,A}$, carring the Frobenius action, where $A$ is a perfectoid adic Banach uniform algebra over $E$;\\
\noindent D. The pseudocoherent modules over the relative $A$-coefficient Robba ring $\widetilde{\Pi}^\infty_{R,A}$, carring the Frobenius action, where $A$ is a perfectoid adic Banach uniform algebra over $E$;\\ 
\noindent E. The pseudocoherent bundles over the relative $A$-coefficient Robba ring $\widetilde{\Pi}_{R,A}$, carring the Frobenius action, where $A$ is a perfectoid adic Banach uniform algebra over $E$;\\
\noindent F. The pseudocoherent bundles over the relative $A$-coefficient Robba ring $\widetilde{\Pi}^{[s,r]}_{R,A}$, carring the Frobenius action, where $A$ is a perfectoid adic Banach uniform algebra over $E$, $0<s\leq r/p^{ah}$.

\begin{theorem}\mbox{\bf{(After Kedlaya-Liu \cite[Theorem 4.6.1]{7KL2})}}
The categories mentioned above are equivalent to each other. (See \cref{theorem3.11}.)
	
\end{theorem}

\indent Now for an adic space $D$ in the pro-\'etale topology which could be covered by perfectoid subdomains over $E$, we have the following categories (objects defined over $D$ are actually compatible families):\\

\noindent A. The pseudocoherent sheaves over the adic relative Fargues-Fontaine curve $Y_{\mathrm{FF},R,D}$, here we consider \'etale topology;\\
\noindent B. The pseudocoherent sheaves over the adic relative Fargues-Fontaine curve $Y_{\mathrm{FF},R,D}$, here we consider pro-\'etale topology;\\
\noindent C. The pseudocoherent modules over the relative $D$-coefficient Robba ring $\widetilde{\Pi}_{R,D}$, carring the Frobenius action;\\
\noindent D. The pseudocoherent modules over the relative $D$-coefficient Robba ring $\widetilde{\Pi}^\infty_{R,D}$, carring the Frobenius action;\\ 
\noindent E. The pseudocoherent bundles over the relative $D$-coefficient Robba ring $\widetilde{\Pi}_{R,D}$, carring the Frobenius action;\\
\noindent F. The pseudocoherent bundles over the relative $D$-coefficient Robba ring $\widetilde{\Pi}^{[s,r]}_{R,D}$, carring the Frobenius action, $0<s\leq r/p^{ah}$.

\begin{theorem}\mbox{\bf{(After Kedlaya-Liu \cite[Theorem 4.6.1]{7KL2})}}
The categories mentioned above are equivalent to each other. (See \cref{theorem3.10}.)
	
\end{theorem}

\indent We then apply to the corresponding construction to the corresponding Fargues-Fontaine curve $Y_{\mathrm{FF},R}$ we have the following categories:\\

\noindent A. The pseudocoherent sheaves over the adic relative Fargues-Fontaine curve $Y_{\mathrm{FF},R,Y_{\mathrm{FF},R}}$, here we consider \'etale topology;\\
\noindent B. The pseudocoherent sheaves over the adic relative Fargues-Fontaine curve $Y_{\mathrm{FF},R,Y_{\mathrm{FF},R}}$, here we consider pro-\'etale topology;\\
\noindent C. The pseudocoherent modules over the relative $Y_{\mathrm{FF},R}$-coefficient Robba ring $\widetilde{\Pi}_{R,Y_{\mathrm{FF},R}}$, carring the Frobenius action;\\
\noindent D. The pseudocoherent modules over the relative $Y_{\mathrm{FF},R}$-coefficient Robba ring $\widetilde{\Pi}^\infty_{R,Y_{\mathrm{FF},R}}$, carring the Frobenius action;\\ 
\noindent E. The pseudocoherent bundles over the relative $Y_{\mathrm{FF},R}$-coefficient Robba ring $\widetilde{\Pi}_{R,Y_{\mathrm{FF},R}}$, carring the Frobenius action;\\
\noindent F. The pseudocoherent bundles over the relative $Y_{\mathrm{FF},R}$-coefficient Robba ring $\widetilde{\Pi}^{[s,r]}_{R,Y_{\mathrm{FF},R}}$, carring the Frobenius action, $0<s\leq r/p^{ah}$.

\begin{theorem}\mbox{\bf{(After Kedlaya-Liu \cite[Theorem 4.6.1]{7KL2})}}
The categories mentioned above are equivalent to each other. (See \cref{corollary3.12}.)\\
	
\end{theorem}

\noindent 2. For general Banach objects, we believe the corresponding picture will still be very interesting and have more potential applications. At least the relative analytic geometry (over the period rings) will be very interesting. That being said, the corresponding study of such relative geometry is really not that easy. The first difficulty will be essentially the sheafiness. Bambozzi-Kremnizer established some derived sheafiness when we replace the Huber's spectrum with the derived ones in \cite{7BK}. For any Banach ring $S$ we have (with the notation in the \cref{section4}) the $\infty$-analytic space $\mathrm{Spa}^h(S)$. We then apply to our period rings with general Banach coefficient $A$ we have the corresponding $\infty$-period rings. And we then look at the following categories (with the notation in the \cref{section4}):\\

\noindent A. The corresponding $f$-projective Frobenius $\varphi^a$-bundles over the $\infty$-Robba ring $\widetilde{\Pi}^{h}_{R,A}$, here $A$ is arbitrary Banach adic uniform algebra over $\mathbb{Q}_p$;\\
%\noindent B. The corresponding $f$-projective Frobenius $\varphi^a$-modules over the $\infty$-Robba ring $\widetilde{\Pi}^{\infty,h}_{R,A}$;\\
\noindent B. The corresponding $f$-projective Frobenius $\varphi^a$-modules over the $\infty$-Robba ring $\widetilde{\Pi}^{[s,r],h}_{R,A}$, here $0<s\leq r/p^{ah}$.

\begin{theorem}\mbox{\bf{(After Kedlaya-Liu \cite[Theorem 4.6.1]{7KL2})}}
The categories mentioned above are equivalent to each other. (See \cref{7theorem4.12}.)\\
	
\end{theorem}

\noindent 3. We revisit the corresponding noncommutative descent which was initiated in \cite{7T2}. We will noncommutativize some argument due to Kedlaya on the corresponding Kiehl's glueing properties without the noetherian assumption. This means that we could actually descent the corresponding finite projective modules in the corresponding noncommutative setting.

\begin{proposition} 
The descent for finite projective bimodules over noncommutative Banach rings holds under the conditions in \cite[Definition 2.7.3 (a),(b)]{7KL1}.	(See \cref{proposition5.11}.)\\
\end{proposition}

\indent Then we apply this to $p$-adic Hodge theory, we look at the following categories:\\
\noindent A. The corresponding finite projective Frobenius $\varphi^a$-bundles over the Robba ring $\widetilde{\Pi}_{R,A}$, here $A$ is arbitrary Banach algebra over $\mathbb{Q}_p$;\\
%\noindent B. The corresponding $f$-projective Frobenius $\varphi^a$-modules over the $\infty$-Robba ring $\widetilde{\Pi}^{\infty,h}_{R,A}$;\\
\noindent B. The corresponding finite projective Frobenius $\varphi^a$-modules over the Robba ring $\widetilde{\Pi}^{[s,r]}_{R,A}$, here $0<s\leq r/p^{ah}$.

\begin{proposition}\mbox{\bf{(After Kedlaya-Liu \cite[Theorem 4.6.1]{7KL2})}}
The categories mentioned above are equivalent to each other. (See \cref{proposition5.13}.)\\
\end{proposition}

\subsection{Convention and Notation}

\indent We make some convention here. When we say an $\infty$-analytic stack, we will mean a $\infty$-sheaf fibered over the corresponding category of all the seminormed commutative monoids satisfying the condition and framework in \cite{7BBBK}. Please note the difference between a derived analytic stack, an $\infty$-analytic stack and a stack, fibered over category of seminormed commutative monoids, as in \cite{7BBBK}.\\

\subsection{Remarks on Future Work}

\indent We will consider more general $\infty$-glueing in the style of Kedlaya-Liu \cite{7KL2} after the foundations we have applied from \cite{7BK} (or possibly equivalently work of Clausen-Scholze \cite{7CS}) whenever one would study more general quasi-coherent sheaves over some sheaves of Banach $\mathcal{O}$-algebras over reasonable sites essentially encoded in our current work.\\

\indent We have not considered the corresponding imperfect setting of the corresponding period rings with really big coefficients along the corresponding style and fashions we established in \cite{7T1} and \cite{7T2} in the context of reasonable towers in \cite{7KL1} and \cite{7KL2}. We would like to study the corresponding commutative and noncommutative setting in future work along these towers going upside and down. There will be some consequences for some specific analytic spaces. Also we could basically contact the corresponding $\infty$-context as well by considering the corresponding sheaves of $\infty$-Banach algebras over some interesting pro-\'etale site over adic spaces.\\

%%\newpage 

\newpage\section{Period Rings with Coefficients in General Adic Spa
ces}

\subsection{Period Rings with Coefficients in General Banach Adic Rin
gs}

\begin{setting}
We now consider the corresponding Kedlaya-Liu's period sheaves with the big coefficients in general adic spaces. We will consider those adic spaces which could be written as the corresponding quotients of the corresponding affinoid (pre-)perfectoid spaces. To be more precise we look at the following rings, we first consider some finite multi-intervals $[s_I,r_I]$ (with real bounds) with $0\leq s_\alpha \leq r_\alpha< \infty$ for each $\alpha \in I$ and we consider the corresponding base field $E$ which is complete nonarchimedian discrete valued carrying some Banach norm, taking the form of $\mathbb{Q}_p$ or $\mathbb{F}_p((\overline{t}))$ with the uniform notation $\pi$ for chosen uniformizer. We then consider over $E$ the corresponding $E$-strictly affinoids:
\begin{align}
E\{s_1/T_1,s_2/T_2,...,s_n/T_n,T_1/r_1,...,T_n/r_n\}	
\end{align}
and we have the corresponding Fr\'echet rings:
\begin{align}
\varprojlim_{r_\alpha \rightarrow 1,\forall \alpha\in I} E\{s_1/T_1,s_2/T_2,...,s_n/T_n,T_1/r_1,...,T_n/r_n\}
\end{align}
with the corresponding ind-Fr\'echet rings:
\begin{align}
\varinjlim_{s_\alpha \rightarrow 0^+s,\forall \alpha\in I}\varprojlim_{r_\alpha \rightarrow 1,\forall \alpha\in I} E\{s_1/T_1,s_2/T_2,...,s_n/T_n,T_1/r_1,...,T_n/r_n\}.
\end{align}
We then consider the following:
\begin{align}
&A_{[s_I,r_I]}\\
&=E(\pi^{1/p^\infty})\{(s_1/T_1)^{1/p^\infty},(s_2/T_2)^{1/p^\infty},...,(s_n/T_n)^{1/p^\infty},(T_1/r_1)^{1/p^\infty},...,(T_n/r_n)^{1/p^\infty}\}^\wedge,\\
&A_{[s_I,r_I]}^+\\
&=\mathcal{O}[\pi^{1/p^\infty}]\{(s_1/T_1)^{1/p^\infty},(s_2/T_2)^{1/p^\infty},...,(s_n/T_n)^{1/p^\infty},(T_1/r_1)^{1/p^\infty},...,(T_n/r_n)^{1/p^\infty}\}^\wedge.	
\end{align}
And we also consider the corresponding Fr\'echet ones:
\begin{align}
&A_{s_I}=\\
&\varprojlim_{r_\alpha \rightarrow 1,\forall \alpha\in I}E(\pi^{1/p^\infty})\{(s_1/T_1)^{1/p^\infty},(s_2/T_2)^{1/p^\infty},...,(s_n/T_n)^{1/p^\infty},(T_1/r_1)^{1/p^\infty},...,(T_n/r_n)^{1/p^\infty}\}^\wedge,\\
&A_{r_I}^+=\\
&\varprojlim_{r_\alpha \rightarrow 1,\forall \alpha\in I}\mathcal{O}[\pi^{1/p^\infty}]\{(s_1/T_1)^{1/p^\infty},(s_2/T_2)^{1/p^\infty},...,(s_n/T_n)^{1/p^\infty},(T_1/r_1)^{1/p^\infty},...,(T_n/r_n)^{1/p^\infty}\}^\wedge.	
\end{align}
\end{setting}

\begin{setting}
We use the corresponding notations as in the following to denote the corresponding preperfectoid ones:
\begin{align}
A'_{[s_I,r_I]}=E\{(s_1/T_1)^{1/p^\infty},(s_2/T_2)^{1/p^\infty},...,(s_n/T_n)^{1/p^\infty},(T_1/r_1)^{1/p^\infty},...,(T_n/r_n)^{1/p^\infty}\}^\wedge,\\
A_{[s_I,r_I]}^{'+}=\mathcal{O}\{(s_1/T_1)^{1/p^\infty},(s_2/T_2)^{1/p^\infty},...,(s_n/T_n)^{1/p^\infty},(T_1/r_1)^{1/p^\infty},...,(T_n/r_n)^{1/p^\infty}\}^\wedge.	
\end{align}
And we also consider the corresponding Fr\'echet ones:
\begin{align}
A'_{s_I}= \varprojlim_{r_\alpha \rightarrow 1,\forall \alpha\in I}E\{(s_1/T_1)^{1/p^\infty},(s_2/T_2)^{1/p^\infty},...,(s_n/T_n)^{1/p^\infty},(T_1/r_1)^{1/p^\infty},...,(T_n/r_n)^{1/p^\infty}\}^\wedge,\\
A_{r_I}^{'+}= \varprojlim_{r_\alpha \rightarrow 1,\forall \alpha\in I}\mathcal{O}\{(s_1/T_1)^{1/p^\infty},(s_2/T_2)^{1/p^\infty},...,(s_n/T_n)^{1/p^\infty},(T_1/r_1)^{1/p^\infty},...,(T_n/r_n)^{1/p^\infty}\}^\wedge.	
\end{align}	
\end{setting}

\indent Then we consider the corresponding context of \cite[Chapter 1]{7T2}, where with trivial coefficients we just consider the same context after \cite{7KL2}. First we have the corresponding consideration, let $F$ be the corresponding field which is a complete discrete valued nnonarchimedean field with uniformizer $\pi$ with normalized Banach norm $\|.\|_F$ such the corresponding evaluation on $\pi$ is just $p^{-1}$, and we assume the residue field of $F$ is just $\mathbb{F}_{p^h}$ for some integer $h>0$. Now we define the following big Robba rings:

\begin{definition} \mbox{\bf{(After Kedlaya-Liu \cite[Definition 4.1.1]{7KL2})}}
We let $A$ be a perfectoid adic uniform Banach algebra over $E$ with integral subring $\mathcal{O}_A$ over $\mathcal{O}_E$. Recall from the corresponding context in \cite{7KL1} we have the corresponding period rings in the relative setting. We follow the corresponding notations we used in \cite[Section 2.1]{7T2} for those corresponding period rings. We first have for a pair $(R,R^+)$ where $R$ is a uniform perfect adic Banach ring over the assumed base $\mathcal{O}_F$. Then we take the corresponding generalized Witt vectors taking the form of $W_{\mathcal{O}_E}(R)$, which is just the ring $\widetilde{\Omega}_R^\mathrm{int}$, and by inverting the corresponding uniformizer we have the ring $\widetilde{\Omega}_R$, then by taking the completed product with $A$ we have $\widetilde{\Omega}_{R,A}$. Now for some $r>0$ we consider the ring $\widetilde{\Pi}^{\mathrm{int},r}_{R}$ which is the completion of $W_{\mathcal{O}_E}(R^+)[[R]]$ (note that this is actually just $W_{\mathcal{O}_E}(R^+)[[x],x\in R]$) by the norm $\|.\|_{\alpha^r}$ defined by:
\begin{align}
\|.\|_{\alpha^r}(\sum_{n\geq 0}\pi^n[\overline{r}_n])=\sup_{n\geq 0}\{p^{-n}\alpha(\overline{r}_n)^r\}.	
\end{align}
Then we have the product $\widetilde{\Pi}^{\mathrm{bd},r}_{R,A}$ defined as the completion under the product norm $\|.\|_{\alpha^r}\otimes \|.\|_A$ of the corresponding ring $\widetilde{\Pi}^{\mathrm{bd},r}_{R}\otimes_{E}A$, which could be also defined from the corresponding integral Robba rings defined above. Then we define the corresponding Robba ring for some interval $I\subset (0,\infty)$ with coefficient in the perfectoid ring $A$ denoted by $\widetilde{\Pi}^{I}_{R,A}$ as the following complete tensor product:
\begin{displaymath}
\widetilde{\Pi}^{I}_{R}\widehat{\otimes}_{E} A	
\end{displaymath}
under the the corresponding tensor product norm $\|.\|_{\alpha^r}\otimes \|.\|_A$. Then we set $\widetilde{\Pi}^r_{R,A}$ as $\varprojlim_{s\rightarrow 0}\widetilde{\Pi}^{[s,r]}_{R,A}$, and then we define $\widetilde{\Pi}_{R,A}$ as $\varinjlim_{r\rightarrow \infty}\widetilde{\Pi}^{[s,r]}_{R,A}$, and we define $\widetilde{\Pi}^\infty_{R,A}$ as $\varprojlim_{r\rightarrow \infty}\widetilde{\Pi}^{r}_{R,A}$. And we also have the corresponding full integral Robba ring and the corresponding full bounded Robba ring by taking the corresponding union through all $r>0$.
\end{definition}

\indent One can even consider the corresponding coefficients in the corresponding Banach adic uniform algebra over the base field.

\begin{definition} \mbox{\bf{(After Kedlaya-Liu \cite[Definition 4.1.1]{7KL2})}}
We let $B$ be an adic uniform Banach algebra over $E$ with integral subring $\mathcal{O}_B$ over $\mathcal{O}_E$. Recall from the corresponding context in \cite{7KL1} we have the corresponding period rings in the relative setting. We follow the corresponding notations we used in \cite[Section 2.1]{7T2} for those corresponding period rings. We first have for a pair $(R,R^+)$ where $R$ is a uniform perfect adic Banach ring over the assumed base $\mathcal{O}_F$. Then we take the corresponding generalized Witt vectors taking the form of $W_{\mathcal{O}_E}(R)$, which is just the ring $\widetilde{\Omega}_R^\mathrm{int}$, and by inverting the corresponding uniformizer we have the ring $\widetilde{\Omega}_R$, then by taking the completed product with $B$ we have $\widetilde{\Omega}_{R,B}$. Now for some $r>0$ we consider the ring $\widetilde{\Pi}^{\mathrm{int},r}_{R}$ which is the completion of $W_{\mathcal{O}_E}(R^+)[[r]:r\in R] $ by the norm $\|.\|_{\alpha^r}$ defined by:
\begin{align}
\|.\|_{\alpha^r}(\sum_{n\geq 0}\pi^n[\overline{r}_n])=\sup_{n\geq 0}\{p^{-n}\alpha(\overline{r}_n)^r\}.	
\end{align}
Then we have the product $\widetilde{\Pi}^{\mathrm{bd},r}_{R,B}$ defined as the completion under the product norm $\|.\|_{\alpha^r}\otimes \|.\|_B$ of the corresponding ring $\widetilde{\Pi}^{\mathrm{bd},r}_{R}\otimes_{E }B$, which could be also defined from the corresponding integral Robba rings defined above. Then we define the corresponding Robba ring for some interval $I\subset (0,\infty)$ with coefficient in the perfectoid ring $B$ denoted by $\widetilde{\Pi}^{I}_{R,B}$ as the following complete tensor product:
\begin{displaymath}
\widetilde{\Pi}^{I}_{R}\widehat{\otimes}_{E} B	
\end{displaymath}
under the the corresponding tensor product norm $\|.\|_{\alpha^r}\otimes \|.\|_B$. Then we set $\widetilde{\Pi}^r_{R,B}$ as $\varprojlim_{s\rightarrow 0}\widetilde{\Pi}^{[s,r]}_{R,B}$, and then we define $\widetilde{\Pi}_{R,B}$ as $\varinjlim_{r\rightarrow \infty}\widetilde{\Pi}^{[s,r]}_{R,B}$, and we define $\widetilde{\Pi}^\infty_{R,B}$ as $\varprojlim_{r\rightarrow \infty}\widetilde{\Pi}^{r}_{R,B}$. And we also have the corresponding full integral Robba ring and the corresponding full bounded Robba ring by taking the corresponding union through all $r>0$.
\end{definition}

\indent In some situation we also need the corresponding preperfectoid coefficients for instance the local chart coming from the adic Fargues-Fontaine curves:

\begin{definition}  \mbox{\bf{(After Kedlaya-Liu \cite[Definition 4.1.1]{7KL2})}}
We let $A$ be a preperfectoid adic uniform Banach algebra over $E$ with integral subring $\mathcal{O}_A$ over $\mathcal{O}_E$. Recall from the corresponding context in \cite{7KL1} we have the corresponding period rings in the relative setting. We follow the corresponding notations we used in \cite[Section 2.1]{7T2} for those corresponding period rings. We first have for a pair $(R,R^+)$ where $R$ is a uniform perfect adic Banach ring over the assumed base $\mathcal{O}_F$. Then we take the corresponding generalized Witt vectors taking the form of $W_{\mathcal{O}_E}(R)$, which is just the ring $\widetilde{\Omega}_R^\mathrm{int}$, and by inverting the corresponding uniformizer we have the ring $\widetilde{\Omega}_R$, then by taking the completed product with $A$ we have $\widetilde{\Omega}_{R,A}$. Now for some $r>0$ we consider the ring $\widetilde{\Pi}^{\mathrm{int},r}_{R}$ which is the completion of $W_{\mathcal{O}_E}(R^+)[[r]:r\in R] $ by the norm $\|.\|_{\alpha^r}$ defined by:
\begin{align}
\|.\|_{\alpha^r}(\sum_{n\geq 0}\pi^n[\overline{r}_n])=\sup_{n\geq 0}\{p^{-n}\alpha(\overline{r}_n)^r\}.	
\end{align}
Then we have the product $\widetilde{\Pi}^{\mathrm{bd},r}_{R,A}$ defined as the completion under the product norm $\|.\|_{\alpha^r}\otimes \|.\|_A$ of the corresponding ring $\widetilde{\Pi}^{\mathrm{bd},r}_{R}\otimes_{E}A$, which could be also defined from the corresponding integral Robba rings defined above. Then we define the corresponding Robba ring for some interval $I\subset (0,\infty)$ with coefficient in the perfectoid ring $A$ denoted by $\widetilde{\Pi}^{I}_{R,A}$ as the following complete tensor product:
\begin{displaymath}
\widetilde{\Pi}^{I}_{R}\widehat{\otimes}_{E} A	
\end{displaymath}
under the the corresponding tensor product norm $\|.\|_{\alpha^r}\otimes \|.\|_A$. Then we set $\widetilde{\Pi}^r_{R,A}$ as $\varprojlim_{s\rightarrow 0}\widetilde{\Pi}^{[s,r]}_{R,A}$, and then we define $\widetilde{\Pi}_{R,A}$ as $\varinjlim_{r\rightarrow \infty}\widetilde{\Pi}^{[s,r]}_{R,A}$, and we define $\widetilde{\Pi}^\infty_{R,A}$ as $\varprojlim_{r\rightarrow \infty}\widetilde{\Pi}^{r}_{R,A}$. And we also have the corresponding full integral Robba ring and the corresponding full bounded Robba ring by taking the corresponding union through all $r>0$.
\end{definition}

\indent In the previous situation where we have the more explicit representation of the corresponding perfectoid coefficients we can basically consider more explicit construction:

\begin{definition}\mbox{\bf{(After Kedlaya-Liu \cite[Definition 4.1.1]{7KL2})}}
We let $A$ be a perfectoid adic uniform Banach algebra over $E$ with integral subring $\mathcal{O}_A$ over $\mathcal{O}_E$ taking the form of one of the rings $A_{[s_I,r_I]}$. Recall from the corresponding context in \cite{7KL1} we have the corresponding period rings in the relative setting. We follow the corresponding notations we used in \cite[Section 2.1]{7T2} for those corresponding period rings. We first have for a pair $(R,R^+)$ where $R$ is a uniform perfect adic Banach ring over the assumed base $\mathcal{O}_F$. Then we take the corresponding generalized Witt vectors taking the form of $W_{\mathcal{O}_E}(R)$, which is just the ring $\widetilde{\Omega}_R^\mathrm{int}$, and by inverting the corresponding uniformizer we have the ring $\widetilde{\Omega}_R$, then by taking the completed product with $A$ we have $\widetilde{\Omega}_{R,A}$. Now for some $r>0$ we consider the ring $\widetilde{\Pi}^{\mathrm{int},r}_{R,A}$ which is the completion of $W_{\mathcal{O}_E}(R^+)[[r]:r\in R]\otimes \mathcal{O}_A $ by the norm $\|.\|_{\alpha^r,A}$ defined by:
\begin{align}
\|.\|_{\alpha^r}&(\sum_{n,i_1,...,i_k,j_1,...,j_k\in \mathbb{Z}[1/p]_{\geq 0}}\pi^n[\overline{r}_n](s_1/T_1)^{i_1}...(s_k/T_k)^{i_k}(T_1/r_1)^{j_1}...(T_k/r_k)^{j_k})\\
&=\sup_{n\geq 0}\{p^{-n}\alpha(\overline{r}_n)^r\}.	
\end{align}
Then we have the product $\widetilde{\Pi}^{\mathrm{bd},r}_{R,A}$ defined as the completion under the product norm $\|.\|_{\alpha^r}\otimes \|.\|_A$ of the corresponding ring $\widetilde{\Pi}^{\mathrm{bd},r}_{R}\otimes_{E}A$, which could be also defined from the corresponding integral Robba rings defined above. Then we define the corresponding Robba ring for some interval $I\subset (0,\infty)$ with coefficient in the perfectoid ring $A$ denoted by $\widetilde{\Pi}^{I}_{R,A}$ as the following complete tensor product:
\begin{displaymath}
\widetilde{\Pi}^{I}_{R}\widehat{\otimes}_{E} E(\pi^{1/p^\infty})\{(s_1/T_1)^{1/p^\infty},(s_2/T_2)^{1/p^\infty},...,(s_n/T_n)^{1/p^\infty},(T_1/r_1)^{1/p^\infty},...,(T_n/r_n)^{1/p^\infty}\}^\wedge	
\end{displaymath}
under the the corresponding tensor product norm $\|.\|_{\alpha^r}\otimes \|.\|_A$. Then we set $\widetilde{\Pi}^r_{R,A}$ as $\varprojlim_{s\rightarrow 0}\widetilde{\Pi}^{[s,r]}_{R,A}$, and then we define $\widetilde{\Pi}_{R,A}$ as $\varinjlim_{r\rightarrow \infty}\widetilde{\Pi}^{[s,r]}_{R,A}$, and we define $\widetilde{\Pi}^\infty_{R,A}$ as $\varprojlim_{r\rightarrow \infty}\widetilde{\Pi}^{r}_{R,A}$. And we also have the corresponding full integral Robba ring and the corresponding full bounded Robba ring by taking the corresponding union through all $r>0$.\\
\end{definition}

%%%%%%%%%%%%%%%%%%%%%%%%%%%%%%%%%%!!!!!!!!!!!!!!!!!!!!!!!

\subsection{Properties}

\indent As in our work \cite{7T2} and \cite{7KL2}, we can discuss some of the key properties of the corresponding period rings defined above.

\begin{setting}
We will mainly consider the corresponding coefficients in 
\begin{center}
$E(\pi^{1/p^\infty})\{(s_1/T_1)^{1/p^\infty},(s_2/T_2)^{1/p^\infty},...,(s_n/T_n)^{1/p^\infty},(T_1/r_1)^{1/p^\infty},...,(T_n/r_n)^{1/p^\infty}\}^\wedge$ 
\end{center}
and their Banach uniform adic strict quotient through the following strict morphism:
\begin{displaymath}
A_{[s_I,r_I]}\rightarrow \overline{A_{[s_I,r_I]}}\rightarrow 0.	
\end{displaymath}
Also we consider the corresponding preperfectoid setting:
\begin{center}
$E\{(s_1/T_1)^{1/p^\infty},(s_2/T_2)^{1/p^\infty},...,(s_n/T_n)^{1/p^\infty},(T_1/r_1)^{1/p^\infty},...,(T_n/r_n)^{1/p^\infty}\}^\wedge$ 
\end{center}
and their Banach uniform adic strict quotient through the following strict morphism:
\begin{displaymath}
A'_{[s_I,r_I]}\rightarrow \overline{A'_{[s_I,r_I]}}\rightarrow 0.	
\end{displaymath} 
 	
\end{setting}

%\begin{proposition}\mbox{\bf{(After Kedlaya-Liu \cite{7KL2})}} 
%Suppose that the corresponding quotient $\overline{A_{[s_I,r_I]}}$ is reduced, then we have that the following identifications:
%\begin{displaymath}
%(\widetilde{\Pi}_{R,\overline{A_{[s_I,r_I]}}})^\times=	(\widetilde{\Pi}^\mathrm{bd}_{R,\overline{A_{[s_I,r_I]}}})^\times
%\end{displaymath}
	
%\end{proposition}

%\begin{proof}
	
%\end{proof}

\begin{proposition} \mbox{\bf{(After Kedlaya-Liu \cite[Lemma 5.2.6]{7KL1})}} \\
Under the corresponding definitions and the corresponding notations in our context we have the following identification:
\begin{align}
\widetilde{\Pi}^{\mathrm{int},s}_{R,A_{[s_I,r_I]}}\bigcap \widetilde{\Pi}^{[s,r]}_{R,A_{[s_I,r_I]}}=	\widetilde{\Pi}^{\mathrm{int},r}_{R,A_{[s_I,r_I]}}\\
\widetilde{\Pi}^{\mathrm{int},s}_{R,\overline{A_{[s_I,r_I]}}}\bigcap \widetilde{\Pi}^{[s,r]}_{R,\overline{A_{[s_I,r_I]}}}=	\widetilde{\Pi}^{\mathrm{int},r}_{R,\overline{A_{[s_I,r_I]}}}.
\end{align}
	
\end{proposition}

\begin{proof}
For the first corresponding identity above we perform the corresponding argument as in \cite[Lemma 5.2.6]{7KL2} and \cite{7T2}. First consider now arbitrary element $h$ in the corresponding intersection on the left. Then we use a corresponding sequence of elements in the bounded Robba ring to approximate the corresponding element $h$ above, say $h_1,h_2,...$. This will mean that for any $k$ we could find some sufficiently large number $P_k$ such that the following holds for any $t\in [s,r]$ and for any $i\geq P_k$:
\begin{align}
\|.\|_{\alpha^t,A_{[s_I,r_I]}}(h-h_i)\leq p^{-k}.	
\end{align}
Then we are going to consider the corresponding decomposition of each element $h_i$ in the following way:
\begin{align}
h_i=\sum_{n\in \mathbb{Z}[1/p],i_1,...,i_k,j_1,...,j_k\in \mathbb{Z}[1/p]_{\geq 0}}\pi^n[\overline{h}_{i,n}](s_1/T_1)^{i_1}...(s_k/T_k)^{i_k}(T_1/r_1)^{j_1}...(T_k/r_k)^{j_k}.	
\end{align}
Then we consider the corresponding integral part of this decomposition namely we consider:
\begin{align}
g_i=\sum_{n\in \mathbb{Z}[1/p]_{\geq 0},i_1,...,i_k,j_1,...,j_k\in \mathbb{Z}[1/p]_{\geq 0}}\pi^n[\overline{h}_{i,n}](s_1/T_1)^{i_1}...(s_k/T_k)^{i_k}(T_1/r_1)^{j_1}...(T_k/r_k)^{j_k}.	
\end{align}
Now since we have that the elements involved are also living in the integral rings we have then that in the expression of $h_i$ we have:
\begin{align}
\|.\|_{\alpha^s,A_{[s_I,r_I]}}(\pi^n[\overline{h}_{i,n}](s_1/T_1)^{i_1}...(s_k/T_k)^{i_k}(T_1/r_1)^{j_1}...(T_k/r_k)^{j_k})\leq p^{-k},\forall n<0.	
\end{align}
Then we have:
\begin{align}
{\alpha}(\overline{h}_{i,n})\leq p^{-(k-n)/s}.	
\end{align}
Then we have:
\begin{displaymath}
\|.\|_{\alpha^r,A_{[s_I,r_I]}}(h-g_i)\leq p^{-n}p^{-(i-n)r/s}.	
\end{displaymath}
This will prove the corresponding results since we can now use the corresponding elements in the integral Robba rings to approximate.\\
For the corresponding second corresponding identity we consider the corresponding strictness of the quotient in the original formation of the rings. Now take any element $\overline{h}$ in the quotient, then lift this to $h$. For the first corresponding identity above we perform the corresponding argument as in \cite[Lemma 5.2.6]{7KL2} and \cite{7T2}. First consider now arbitrary element $h$ in the corresponding intersection on the left. Then we use a corresponding sequence of elements in the bounded Robba ring to approximate the corresponding element $h$ above, say $h_1,h_2,...$. This will mean that for any $k$ we could find some sufficiently large number $P_k$ such that the following holds for any $t\in [s,r]$ and for any $i\geq P_k$:
\begin{align}
\|.\|_{\alpha^t,A_{[s_I,r_I]}}(h-h_i)\leq p^{-k}.	
\end{align}
Then we are going to consider the corresponding decomposition of each element $h_i$ in the following way:
\begin{align}
h_i=\sum_{n\in \mathbb{Z}[1/p],i_1,...,i_k,j_1,...,j_k\in \mathbb{Z}[1/p]_{\geq 0}}\pi^n[\overline{h}_{i,n}](s_1/T_1)^{i_1}...(s_k/T_k)^{i_k}(T_1/r_1)^{j_1}...(T_k/r_k)^{j_k}.	
\end{align}
Then we consider the corresponding integral part of this decomposition namely we consider:
\begin{align}
g_i=\sum_{n\in \mathbb{Z}[1/p]_{\geq 0},i_1,...,i_k,j_1,...,j_k\in \mathbb{Z}[1/p]_{\geq 0}}\pi^n[\overline{h}_{i,n}](s_1/T_1)^{i_1}...(s_k/T_k)^{i_k}(T_1/r_1)^{j_1}...(T_k/r_k)^{j_k}.	
\end{align}
Now since we have that the elements involved are also living in the integral rings we have then that in the expression of $h_i$ we have:
\begin{align}
\|.\|_{\alpha^s,A_{[s_I,r_I]}}(\pi^n[\overline{h}_{i,n}](s_1/T_1)^{i_1}...(s_k/T_k)^{i_k}(T_1/r_1)^{j_1}...(T_k/r_k)^{j_k})\leq p^{-k},\forall n<0.	
\end{align}
Then we have:
\begin{align}
{\alpha}(\overline{h}_{i,n})\leq p^{-(k-n)/s}.	
\end{align}
Then we have:
\begin{displaymath}
\|.\|_{\alpha^r,A_{[s_I,r_I]}}(h-g_i)\leq p^{-n}p^{-(i-n)r/s}.	
\end{displaymath}
This will prove the corresponding results since we can now use the corresponding elements in the integral Robba rings to approximate. Then we have that:
\begin{align}
\|.\|_{\alpha^r,\overline{A_{[s_I,r_I]}}}(\overline{h}-\overline{g}_i)\leq p^{-n}p^{-(i-n)r/s},	
\end{align}
which shows the corresponding approximating process in the strict quotient.

\end{proof}

\begin{proposition} \mbox{\bf{(After Kedlaya-Liu \cite[Lemma 5.2.6]{7KL1})}} 
Under the corresponding definitions and the corresponding notations in our context we have the following identification:
\begin{align}
\widetilde{\Pi}^{\mathrm{int},s}_{R,A'_{[s_I,r_I]}}\bigcap \widetilde{\Pi}^{[s,r]}_{R,A'_{[s_I,r_I]}}=	\widetilde{\Pi}^{\mathrm{int},r}_{R,A'_{[s_I,r_I]}}\\
\widetilde{\Pi}^{\mathrm{int},s}_{R,\overline{A_{[s_I,r_I]}}}\bigcap \widetilde{\Pi}^{[s,r]}_{R,\overline{A'_{[s_I,r_I]}}}=	\widetilde{\Pi}^{\mathrm{int},r}_{R,\overline{A'_{[s_I,r_I]}}}.
\end{align}
	
\end{proposition}

\begin{proof}
For the first corresponding identity above we perform the corresponding argument as in \cite[Lemma 5.2.6]{7KL2} and \cite{7T2}. First consider now arbitrary element $h$ in the corresponding intersection on the left. Then we use a corresponding sequence of elements in the bounded Robba ring to approximate the corresponding element $h$ above, say $h_1,h_2,...$. This will mean that for any $k$ we could find some sufficiently large number $P_k$ such that the following holds for any $t\in [s,r]$ and for any $i\geq P_k$:
\begin{align}
\|.\|_{\alpha^t,A'_{[s_I,r_I]}}(h-h_i)\leq p^{-k}.	
\end{align}
Then we are going to consider the corresponding decomposition of each element $h_i$ in the following way:
\begin{align}
h_i=\sum_{n\in \mathbb{Z},i_1,...,i_k,j_1,...,j_k\in \mathbb{Z}[1/p]_{\geq 0}}\pi^n[\overline{h}_{i,n}](s_1/T_1)^{i_1}...(s_k/T_k)^{i_k}(T_1/r_1)^{j_1}...(T_k/r_k)^{j_k}.	
\end{align}
Then we consider the corresponding integral part of this decomposition namely we consider:
\begin{align}
g_i=\sum_{n\in \mathbb{Z}_{\geq 0},i_1,...,i_k,j_1,...,j_k\in \mathbb{Z}[1/p]_{\geq 0}}\pi^n[\overline{h}_{i,n}](s_1/T_1)^{i_1}...(s_k/T_k)^{i_k}(T_1/r_1)^{j_1}...(T_k/r_k)^{j_k}.	
\end{align}
Now since we have that the elements involved are also living in the integral rings we have then that in the expression of $h_i$ we have:
\begin{align}
\|.\|_{\alpha^s,A'_{[s_I,r_I]}}(\pi^n[\overline{h}_{i,n}](s_1/T_1)^{i_1}...(s_k/T_k)^{i_k}(T_1/r_1)^{j_1}...(T_k/r_k)^{j_k})\leq p^{-k},\forall n<0.	
\end{align}
Then we have:
\begin{align}
{\alpha}(\overline{h}_{i,n})\leq p^{-(k-n)/s}.	
\end{align}
Then we have:
\begin{displaymath}
\|.\|_{\alpha^r,A'_{[s_I,r_I]}}(h-g_i)\leq p^{-n}p^{-(i-n)r/s}.	
\end{displaymath}
This will prove the corresponding results since we can now use the corresponding elements in the integral Robba rings to approximate.\\
For the corresponding second corresponding identity we consider the corresponding strictness of the quotient in the original formation of the rings. Now take any element $\overline{h}$ in the quotient, then lift this to $h$. For the first corresponding identity above we perform the corresponding argument as in \cite[Lemma 5.2.6]{7KL2} and \cite{7T2}. First consider now arbitrary element $h$ in the corresponding intersection on the left. Then we use a corresponding sequence of elements in the bounded Robba ring to approximate the corresponding element $h$ above, say $h_1,h_2,...$. This will mean that for any $k$ we could find some sufficiently large number $P_k$ such that the following holds for any $t\in [s,r]$ and for any $i\geq P_k$:
\begin{align}
\|.\|_{\alpha^t,A'_{[s_I,r_I]}}(h-h_i)\leq p^{-k}.	
\end{align}
Then we are going to consider the corresponding decomposition of each element $h_i$ in the following way:
\begin{align}
h_i=\sum_{n\in \mathbb{Z},i_1,...,i_k,j_1,...,j_k\in \mathbb{Z}[1/p]_{\geq 0}}\pi^n[\overline{h}_{i,n}](s_1/T_1)^{i_1}...(s_k/T_k)^{i_k}(T_1/r_1)^{j_1}...(T_k/r_k)^{j_k}.	
\end{align}
Then we consider the corresponding integral part of this decomposition namely we consider:
\begin{align}
g_i=\sum_{n\in \mathbb{Z}_{\geq 0},i_1,...,i_k,j_1,...,j_k\in \mathbb{Z}[1/p]_{\geq 0}}\pi^n[\overline{h}_{i,n}](s_1/T_1)^{i_1}...(s_k/T_k)^{i_k}(T_1/r_1)^{j_1}...(T_k/r_k)^{j_k}.	
\end{align}
Now since we have that the elements involved are also living in the integral rings we have then that in the expression of $h_i$ we have:
\begin{align}
\|.\|_{\alpha^s,A'_{[s_I,r_I]}}(\pi^n[\overline{h}_{i,n}](s_1/T_1)^{i_1}...(s_k/T_k)^{i_k}(T_1/r_1)^{j_1}...(T_k/r_k)^{j_k})\leq p^{-k},\forall n<0.	
\end{align}
Then we have:
\begin{align}
{\alpha}(\overline{h}_{i,n})\leq p^{-(k-n)/s}.	
\end{align}
Then we have:
\begin{displaymath}
\|.\|_{\alpha^r,A'_{[s_I,r_I]}}(h-g_i)\leq p^{-n}p^{-(i-n)r/s}.	
\end{displaymath}
This will prove the corresponding results since we can now use the corresponding elements in the integral Robba rings to approximate. Then we have that:
\begin{align}
\|.\|_{\alpha^r,\overline{A'_{[s_I,r_I]}}}(\overline{h}-\overline{g}_i)\leq p^{-n}p^{-(i-n)r/s},	
\end{align}
which shows the corresponding approximating process in the strict quotient.

\end{proof}

\begin{proposition} \mbox{\bf{(After Kedlaya-Liu \cite[Lemma 5.2.7]{7KL1})}}\\ 
We consider the Robba ring over some closed interval $\widetilde{\Pi}^{[s,r]}_{R,A_{[s_I,r_I]}}$ we have the following decomposition for each element $x$ in this ring into the following way:
\begin{displaymath}
x=y+z	
\end{displaymath}
where $y\in \pi^n \widetilde{\Pi}^{\mathrm{int},r}_{R,A_{[s_I,r_I]}}$ and $z\in \widetilde{\Pi}^{[s,r']}_{R,A_{[s_I,r_I]}}$ for each $r'\geq r$, and here $n>0$ lives in $\mathbb{Z}[1/p]$. With this decomposition we have the following estimate:
\begin{align}
\|.\|_{\alpha^t,A_{[s_I,r_I]}}(z)\leq p^{(1-n)(1-t/r)}\|.\|_{\alpha^r,A_{[s_I,r_I]}}(x)^{t/r}.	
\end{align}
And we consider the Robba ring over some closed interval $\widetilde{\Pi}^{[s,r]}_{R,\overline{A_{[s_I,r_I]}}}$ we have the following decomposition for each element $x$ in this ring into the following way:
\begin{displaymath}
x=y+z	
\end{displaymath}
where $y\in p^n \widetilde{\Pi}^{\mathrm{int},r}_{R,\overline{A_{[s_I,r_I]}}}$ and $z\in \widetilde{\Pi}^{[s,r']}_{R,\overline{A_{[s_I,r_I]}}}$ for each $r'\geq r$, and here $n>0$ lives in $\mathbb{Z}[1/p]$. With this decomposition we have the following estimate:
\begin{align}
\|.\|_{\alpha^t,\overline{A_{[s_I,r_I]}}}(z)\leq p^{(1-n)(1-t/r)}\|.\|_{\alpha^r,\overline{A_{[s_I,r_I]}}}(x)^{t/r}.	
\end{align}
	
\end{proposition}

\begin{proof}
For the first part, we first consider any element $r$ in the bounded Robba ring. Each such element admits the corresponding expansion by the linear combination of the Teichm\"uller elements, each term takes the corresponding form as $\pi^l[r_{l}]$ summing up from some negative integer. Then the corresponding decomposition just reads:
\begin{align}
y=\sum_{l\geq n,i_1,...,i_k,j_1,...,j_k\in \mathbb{Z}[1/p]_{\geq 0}}\pi^l[\overline{r}_{l}](s_1/T_1)^{i_1}...(s_k/T_k)^{i_k}(T_1/r_1)^{j_1}...(T_k/r_k)^{j_k},\\
z=r-y.	
\end{align}
Then in this case we have the corresponding desired estimate:
\begin{align}
\|.\|_{\alpha^t,A_{[s_I,r_I]}}(\pi^l[\overline{r}_{l}](s_1/T_1)^{i_1}...(s_k/T_k)^{i_k}(T_1/r_1)^{j_1}...(T_k/r_k)^{j_k})\\
\leq p^{(1-n)(1-t/r)}\|.\|_{\alpha^r,A_{[s_I,r_I]}}(\pi^l[\overline{r}_{l}](s_1/T_1)^{i_1}...(s_k/T_k)^{i_k}(T_1/r_1)^{j_1}...(T_k/r_k)^{j_k})^{t/r}.
\end{align}
 Then we choose a corresponding sequence of the corresponding elements in the bounded Robba ring to approximate any chosen arbitrary element in $\widetilde{\Pi}^{[s,r]}_{R,A_{[s_I,r_I]}}$ such that we have:
 \begin{displaymath}
 \|.\|_{\alpha^t,A_{[s_I,r_I]}}(r-r_0-...-r_i)\leq p^{-1-i}  \|.\|_{\alpha^t,A_{[s_I,r_I]}}(r),\forall t\in [s,r].	
 \end{displaymath}
Then apply the corresponding construction in the previous situation to each $r_i$ we have that there is a corresponding decomposition:
\begin{displaymath}
r_i=y_i+z_i	
\end{displaymath}
with the corresponding desired extracted sum $\sum_{i}y_i$ which converges to the desired element $y$. While for $z_i$ we consider the corresponding estimate as in the following:
\begin{align}
 \|.\|_{\alpha^t,A_{[s_I,r_I]}}(z_i) &\leq  p^{(1-n)(1-t/r)}\|.\|_{\alpha^r,A_{[s_I,r_I]}}(r_i)^{t/r}\\
 &\leq p^{-i}p^{(1-n)(1-t/r)}\|.\|_{\alpha^r,A_{[s_I,r_I]}}(r_i)^{t/r},	
\end{align}
which shows the corresponding sum $\sum_i z_i$ converges to desired element $z$ in the corresponding union of all the Robba rings $\widetilde{\Pi}^{[s,r']}_{R,A_{[s_I,r_I]}}$ for all $r'\geq r$. Then for the second statement we just perform as before to look at the corresponding lifting of elements through the corresponding strict quotient morphism.
 
\end{proof}

\begin{proposition} \mbox{\bf{(After Kedlaya-Liu \cite[Lemma 5.2.7]{7KL1})}} \\
We consider the Robba ring over some closed interval $\widetilde{\Pi}^{[s,r]}_{R,A'_{[s_I,r_I]}}$ we have the following decomposition for each element $x$ in this ring into the following way:
\begin{displaymath}
x=y+z	
\end{displaymath}
where $y\in \pi^n \widetilde{\Pi}^{\mathrm{int},r}_{R,A'_{[s_I,r_I]}}$ and $z\in \widetilde{\Pi}^{[s,r']}_{R,A'_{[s_I,r_I]}}$ for each $r'\geq r$, and here $n>0$ is some arbitrary given integer. With this decomposition we have the following estimate:
\begin{align}
\|.\|_{\alpha^t,A'_{[s_I,r_I]}}(z)\leq p^{(1-n)(1-t/r)}\|.\|_{\alpha^r,A'_{[s_I,r_I]}}(x)^{t/r}.	
\end{align}
And we consider the Robba ring over some closed interval $\widetilde{\Pi}^{[s,r]}_{R,\overline{A'_{[s_I,r_I]}}}$ we have the following decomposition for each element $x$ in this ring into the following way:
\begin{displaymath}
x=y+z	
\end{displaymath}
where $y\in p^n \widetilde{\Pi}^{\mathrm{int},r}_{R,\overline{A'_{[s_I,r_I]}}}$ and $z\in \widetilde{\Pi}^{[s,r']}_{R,\overline{A'_{[s_I,r_I]}}}$ for each $r'\geq r$, and here $n>0$ lives in $\mathbb{Z}$. With this decomposition we have the following estimate:
\begin{align}
\|.\|_{\alpha^t,\overline{A'_{[s_I,r_I]}}}(z)\leq p^{(1-n)(1-t/r)}\|.\|_{\alpha^r,\overline{A'_{[s_I,r_I]}}}(x)^{t/r}.	
\end{align}

\end{proposition}

\begin{proof}
For the first part, we first consider any element $r$ in the bounded Robba ring. Each such element admits the corresponding expansion by the linear combination of the Teichm\"uller elements, each term takes the corresponding form as $\pi^l[r_{l}]$ summing up from some negative integer. Then the corresponding decomposition just reads:
\begin{align}
y=\sum_{l\geq n,i_1,...,i_k,j_1,...,j_k\in \mathbb{Z}[1/p]_{\geq 0}}\pi^l[\overline{r}_{l}](s_1/T_1)^{i_1}...(s_k/T_k)^{i_k}(T_1/r_1)^{j_1}...(T_k/r_k)^{j_k},\\
z=r-y.	
\end{align}
Then in this case we have the corresponding desired estimate:
\begin{align}
\|.\|_{\alpha^t,A'_{[s_I,r_I]}}(\pi^l[\overline{r}_{l}](s_1/T_1)^{i_1}...(s_k/T_k)^{i_k}(T_1/r_1)^{j_1}...(T_k/r_k)^{j_k})\\
\leq p^{(1-n)(1-t/r)}\|.\|_{\alpha^r,A'_{[s_I,r_I]}}(\pi^l[\overline{r}_{l}](s_1/T_1)^{i_1}...(s_k/T_k)^{i_k}(T_1/r_1)^{j_1}...(T_k/r_k)^{j_k})^{t/r}.
\end{align}
 Then we choose a corresponding sequence of the corresponding elements in the bounded Robba ring to approximate any chosen arbitrary element in $\widetilde{\Pi}^{[s,r]}_{R,A'_{[s_I,r_I]}}$ such that we have:
 \begin{displaymath}
 \|.\|_{\alpha^t,A'_{[s_I,r_I]}}(r-r_0-...-r_i)\leq p^{-1-i}  \|.\|_{\alpha^t,A'_{[s_I,r_I]}}(r),\forall t\in [s,r].	
 \end{displaymath}
Then apply the corresponding construction in the previous situation to each $r_i$ we have that there is a corresponding decomposition:
\begin{displaymath}
r_i=y_i+z_i	
\end{displaymath}
with the corresponding desired extracted sum $\sum_{i}y_i$ which converges to the desired element $y$. While for $z_i$ we consider the corresponding estimate as in the following:
\begin{align}
 \|.\|_{\alpha^t,A'_{[s_I,r_I]}}(z_i) &\leq  p^{(1-n)(1-t/r)}\|.\|_{\alpha^r,A'_{[s_I,r_I]}}(r_i)^{t/r}\\
 &\leq p^{-i}p^{(1-n)(1-t/r)}\|.\|_{\alpha^r,A'_{[s_I,r_I]}}(r_i)^{t/r},	
\end{align}
which shows the corresponding sum $\sum_i z_i$ converges to desired element $z$ in the corresponding union of all the Robba rings $\widetilde{\Pi}^{[s,r']}_{R,A'_{[s_I,r_I]}}$ for all $r'\geq r$. Then for the second statement we just perform as before to look at the corresponding lifting of elements through the corresponding strict quotient morphism.	
\end{proof}

\begin{proposition} \mbox{\bf{(After Kedlaya-Liu \cite[Lemma 5.2.10]{7KL1})}}
We could have the following result in our new context around the corresponding period rings with respect to two intervals $[s_1,r_1],[s_2,r_2]$ such that we have:
\begin{align}
0<s_1\leq s_2\leq r_1 \leq r_2.	
\end{align}
To be more precise in our current situation, we have:
\begin{align}
\widetilde{\Pi}^{[s_1,r_1]}_{R,A_{[s_I,r_I]}}\bigcap \widetilde{\Pi}^{[s_2,r_2]}_{R,A_{[s_I,r_I]}}=	\widetilde{\Pi}^{\mathrm{int},r}_{R,A_{[s_1,r_2]}}	
\end{align}
and 
\begin{align}
\widetilde{\Pi}^{[s_1,r_1]}_{R,\overline{A_{[s_I,r_I]}}}\bigcap \widetilde{\Pi}^{[s_2,r_2]}_{R,\overline{A_{[s_I,r_I]}}}=	\widetilde{\Pi}^{\mathrm{int},r}_{R,\overline{A_{[s_I,r_I]}}}.	
\end{align}

\end{proposition}

\begin{proof}
See the proof of \cite[Lemma 5.2.10]{7KL1}.	
\end{proof}

\begin{proposition} \mbox{\bf{(After Kedlaya-Liu \cite[Lemma 5.2.10]{7KL1})}}
We could have the following result in our new context around the corresponding period rings with respect to two intervals $[s_1,r_1],[s_2,r_2]$ such that we have:
\begin{align}
0<s_1\leq s_2\leq r_1 \leq r_2.	
\end{align}
To be more precise in our current situation, we have:
\begin{align}
\widetilde{\Pi}^{[s_1,r_1]}_{R,A'_{[s_I,r_I]}}\bigcap \widetilde{\Pi}^{[s_2,r_2]}_{R,A'_{[s_I,r_I]}}=	\widetilde{\Pi}^{\mathrm{int},r}_{R,A'_{[s_1,r_2]}}	
\end{align}
and 
\begin{align}
\widetilde{\Pi}^{[s_1,r_1]}_{R,\overline{A'_{[s_I,r_I]}}}\bigcap \widetilde{\Pi}^{[s_2,r_2]}_{R,\overline{A'_{[s_I,r_I]}}}=	\widetilde{\Pi}^{\mathrm{int},r}_{R,\overline{A'_{[s_I,r_I]}}}.	
\end{align}

\end{proposition}

\begin{proof}
See the proof of the previous proposition.	
\end{proof}

%%%%%%%%%%%%%%%%%%%%%%%%%%%%%%%%%%!!!!!!!!!!!!!!!!!!!

%%\newpage

\newpage\section{Period Rings and Sheaves with coefficient in Adic Spaces}

\subsection{Fundamental Settings}

\indent Scholze's diamonds \cite{7Sch1} are very general adic stacks where the category is big and convenient enough to include the corresponding perfectoid spaces and the corresponding seminormal rigid analytic spaces. The two spaces happen to be our main interests at the same time. However, we do not actually use the language of diamonds.

\indent Now we use the notation $D$ to denote a general adic space over $E$ in the sense of \cite{7KL1} and \cite{7KL2}. Then we define the following period rings with coefficient in $D$ by using the corresponding Banach algebra structures over $E$.

\begin{definition}

We now consider the following sheaves in adic topology, \'etale topology and pro-\'etale topology (which are indeed not just presheaves by using the corresponding Schauder basis):
\begin{align}
\widetilde{\Pi}_{R,D},\widetilde{\Pi}^\infty_{R,D},\widetilde{\Pi}^I_{R,D},\widetilde{\Pi}^r_{R,D}	
\end{align}
which assign each affinoid subdomain or perfectoid subdomain $\mathrm{Spa}(A,A^+)$ of $D$ to the following rings:
\begin{align}
\widetilde{\Pi}_{R,A},\widetilde{\Pi}^\infty_{R,A},\widetilde{\Pi}^I_{R,A},\widetilde{\Pi}^r_{R,A}.	
\end{align}

\end{definition}

\indent To define Frobenius modules over $D$ we first define the following local picture:

\begin{definition} \mbox{\bf{(After Kedlaya-Liu \cite[Definition 4.4.4]{7KL2})}}
We define the corresponding $\varphi^a$-modules over the corresponding period rings above, which we could use some uniform notation $\triangle^?,?=\emptyset,I,r,\infty,\triangle=\widetilde{\Pi}_{R,B}$ to denote these. We define any $\varphi^a$-module over the period rings with index $\emptyset$ above to be a finitely generated projective module over $\triangle^?$ carrying semilinear action from the Frobenius operator $\varphi^a$ such that $\varphi^{a*}M\overset{\sim}{\rightarrow}M$. We define any $\varphi^a$-module over the period rings with index $r$ above to be a finitely generated projective module over $\triangle^?$ carrying semilinear action from the Frobenius operator $\varphi^a$ such that $\varphi^{a*}M\overset{\sim}{\rightarrow}M\otimes \triangle^{rp^{-ah}}$. We define any $\varphi^a$-module over the period rings with index $[s,r]$ above to be a finitely generated projective module over $\triangle^?$ carrying semilinear action from the Frobenius operator $\varphi^a$ such that $\varphi^{a*}M\otimes_{\triangle^{[sp^{-ha},rp^{-ha}]}} \triangle^{[s,rp^{-ha}]}\overset{\sim}{\rightarrow}M\otimes \triangle^{[s,rp^{-ha}]}$. For a corresponding object over $\triangle^\emptyset$, we assume the module is the base change from some module over $\triangle^{r_0}$ for some $r_0>0$.
\end{definition}

\begin{definition} \mbox{\bf{(After Kedlaya-Liu \cite[Definition 4.4.4]{7KL2})}}
We define the corresponding pseudocoherent $\varphi^a$-modules over the corresponding period rings above, which we could use some uniform notation $\triangle^?,?=\emptyset,I,r,\infty,\triangle=\widetilde{\Pi}_{R,B}$ to denote these. We define any pseudocoherent $\varphi^a$-module over the period rings with index $\empty$ above to be a pseudocoherent module over $\triangle^?$ carrying semilinear action from the Frobenius operator $\varphi^a$ such that $\varphi^{a*}M\overset{\sim}{\rightarrow}M$. We define any pseudocoherent $\varphi^a$-module over the period rings with index $r$ above to be a pseudocoherent module over $\triangle^?$ carrying semilinear action from the Frobenius operator $\varphi^a$ such that $\varphi^{a*}M\overset{\sim}{\rightarrow}M\otimes \triangle^{rp^{-ah}}$. We define any pseudocoherent $\varphi^a$-module over the period rings with index $[s,r]$ above to be a pseudocoherent module over $\triangle^?$ carrying semilinear action from the Frobenius operator $\varphi^a$ such that $\varphi^{a*}M\otimes_{\triangle^{[sp^{-ha},rp^{-ha}]}} \triangle^{[s,rp^{-ha}]}\overset{\sim}{\rightarrow}M\otimes \triangle^{[s,rp^{-ha}]}$. For a corresponding object over $\triangle^\emptyset$, we assume the module is the base change from some module over $\triangle^{r_0}$ for some $r_0>0$. As in \cite[Definition 4.4.4]{7KL2} we impose the corresponding topological condition on the corresponding modules by imposing that all the modules are complete with respect to the corresponding natural topology, and over the Robba rings with respect to some specific interval we assume that the corresponding modules are \'etale-stably pseudocoherent. And we assume the corresponding stability with respect to the base change as in \cite[Theorem 4.6.1]{7KL2} for modules over $\triangle^{r_0},\triangle^\infty$.
\end{definition}

%\begin{remark}
%We believe that one can discuss more general objects such as the corresponding coherent sheaves after \cite{7Lu1} and \cite{7Lu2}, but at this moment let us just focus on the derived vector bundles.	
%\end{remark}

\begin{definition} \mbox{\bf{(After Kedlaya-Liu \cite[Definition 4.4.6]{7KL2})}}
We now define the corresponding $\varphi^a$-bundles over the corresponding period rings above, which we could use some uniform notation $\triangle_{*,B}^?,?=\emptyset,r,*=R$ to denote these. We define a $\varphi^a$-bundle over $\triangle_{*,B}^?,?=\emptyset,r$ to be a compatible family of finitely generated projective $\varphi^a$-modules over $\triangle_{*,B}^?,?=[s',r']$ for suitable $[s',r']$ (namely contained in $(0,r]$ if we are looking at the ring $\triangle_{*,B}^r$) satisfying the corresponding restriction compatibility and cocycle condition.
\end{definition}

\begin{definition} \mbox{\bf{(After Kedlaya-Liu \cite[Definition 4.4.6]{7KL2})}}
We now define the corresponding pseudocoherent $\varphi^a$-bundles over the corresponding period rings above, which we could use some uniform notation $\triangle_{*,B}^?,?=\emptyset,r,*=R$ to denote these. We define a pseudocoherent $\varphi^a$-bundle over $\triangle_{*,B}^?,?=\emptyset,r$ to be a compatible family of \'etale-stably pseudocoherent $\varphi^a$-modules over $\triangle_{*,B}^?,?=[s',r']$ for suitable $[s',r']$ (namely contained in $(0,r]$ if we are looking at the ring $\triangle_{*,B}^r$) satisfying the corresponding restriction compatibility and cocycle condition.
\end{definition}

\begin{remark}
The $p$-adic functional analytic property we imposed on the algebraic pseudocoherent modules could be also translated to adic topology, which is the corresponding stably pseudocoherent one. This is very important whenever one works with the corresponding analytic topology instead of the corresponding pro-\'etale topology.	
\end{remark}

\indent Then we define things over the adic spaces.

\begin{definition}
In adic topology, we define the $\varphi^a$-module over any\\ $\triangle^?,?=\emptyset,I,r,\infty,\triangle=\widetilde{\Pi}_{R,D}$ to be a compatible family of $\varphi^a$-modules over $\triangle^?,?=\emptyset,I,r,\infty,\triangle=\widetilde{\Pi}_{R,A}$ for each affinoid subdomain $\mathrm{Spa}(A,A^+)$ of $D$.	
\end{definition}

\begin{definition}
In adic topology, we define a pseudocoherent $\varphi^a$-module over any $\triangle^?,?=\emptyset,I,r,\infty,\triangle=\widetilde{\Pi}_{R,D}$ to be a compatible family of pseudocoherent $\varphi^a$-modules over $\triangle^?,?=\emptyset,I,r,\infty,\triangle=\widetilde{\Pi}_{R,A}$ for each affinoid subdomain $\mathrm{Spa}(A,A^+)$ of $D$.	
\end{definition}

\begin{definition}
In adic topology, we define the $\varphi^a$-bundle over any\\ $\triangle^?,?=\emptyset,I,r,\infty,\triangle=\widetilde{\Pi}_{R,D}$ to be a compatible family of $\varphi^a$-bundles over $\triangle^?,?=\emptyset,I,r,\infty,\triangle=\widetilde{\Pi}_{R,A}$ for each affinoid subdomain $\mathrm{Spa}(A,A^+)$ of $D$.	
\end{definition}

\begin{definition}
In adic topology, we define a pseudocoherent $\varphi^a$-bundle over any $\triangle^?,?=\emptyset,I,r,\infty,\triangle=\widetilde{\Pi}_{R,D}$ to be a compatible family of pseudocoherent $\varphi^a$-bundles over $\triangle^?,?=\emptyset,I,r,\infty,\triangle=\widetilde{\Pi}_{R,A}$ for each affinoid subdomain $\mathrm{Spa}(A,A^+)$ of $D$.	
\end{definition}

\begin{definition}
In pro-\'etale topology, we define the $\varphi^a$-module over any $\triangle^?,?=\emptyset,I,r,\infty,\triangle=\widetilde{\Pi}_{R,D}$ to be a compatible family of $\varphi^a$-modules over $\triangle^?,?=\emptyset,I,r,\infty,\triangle=\widetilde{\Pi}_{R,A}$ for each perfectoid subdomain $\mathrm{Spa}(A,A^+)$ of $D$.	
\end{definition}

\begin{definition}
In pro-\'etale topology, we define a pseudocoherent $\varphi^a$-module over any $\triangle^?,?=\emptyset,I,r,\infty,\triangle=\widetilde{\Pi}_{R,D}$ to be a compatible family of pseudocoherent $\varphi^a$-modules over $\triangle^?,?=\emptyset,I,r,\infty,\triangle=\widetilde{\Pi}_{R,A}$ for each perfectoid subdomain $\mathrm{Spa}(A,A^+)$ of $D$.	
\end{definition}

\begin{definition}
In pro-\'etale topology, we define the $\varphi^a$-bundle over any $\triangle^?,?=\emptyset,I,r,\infty,\triangle=\widetilde{\Pi}_{R,D}$ to be a compatible family of $\varphi^a$-bundles over $\triangle^?,?=\emptyset,I,r,\infty,\triangle=\widetilde{\Pi}_{R,A}$ for each perfectoid subdomain $\mathrm{Spa}(A,A^+)$ of $D$.	
\end{definition}

\begin{definition}
In pro-\'etale topology, we define a pseudocoherent $\varphi^a$-bundle over any $\triangle^?,?=\emptyset,I,r,\infty,\triangle=\widetilde{\Pi}_{R,D}$ to be a compatible family of pseudocoherent $\varphi^a$-bundles over $\triangle^?,?=\emptyset,I,r,\infty,\triangle=\widetilde{\Pi}_{R,A}$ for each perfectoid subdomain $\mathrm{Spa}(A,A^+)$ of $D$.	
\end{definition}

\subsection{Fundamental Comparisons} \label{7section3.2}

\begin{theorem} \mbox{\bf{(After Kedlaya-Liu \cite[Theorem 4.6.1]{7KL2})}} \label{theorem3.10}
Consider the following categories:\\
1. The corresponding category of all the sheaves of pseudocoherent $\mathcal{O}_{\text{\'etale}}$-modules over the adic Fargues-Fontaine curve $Y_{\mathrm{FF},R,D}$ in the corresponding \'etale topology;\\
2. The corresponding category of all the sheaves of pseudocoherent $\mathcal{O}_{\text{pro-\'etale}}$-modules over the adic Fargues-Fontaine curve $Y_{\mathrm{FF},R,D}$ in the corresponding pro-\'etale topology;\\
3. The corresponding category of all the pseudocoherent bundles over the period ring $\widetilde{\Pi}_{R,D}$, carrying the corresponding Frobenius action from the operator $\varphi^a$;\\
4. The corresponding category of all the pseudocoherent modules over the period ring $\widetilde{\Pi}^{\infty}_{R,D}$, carrying the corresponding Frobenius action from the operator $\varphi^a$;\\
5. The corresponding category of all the pseudocoherent modules over the period ring $\widetilde{\Pi}^{}_{R,D}$, carrying the corresponding Frobenius action from the operator $\varphi^a$;\\
6. The corresponding category of all the finitely generated projective modules over the period ring $\widetilde{\Pi}^{[s,r]}_{R,D}$, carrying the corresponding Frobenius action from the operator $\varphi^a$, where $0<s\leq r/p^{ah}$.\\
%5. The corresponding category of all the finitely generated projective modules over the homotopical period ring $\widetilde{\Pi}^{[s,r],h}_{R,A}$, carrying the corresponding Frobenius action from the operator $\varphi^a$, where $0<s\leq r/p^{ah}$, such that the corresponding underlying modules are plat.\\
Then we have that they are equivalent.
\end{theorem}

This could be proved by the following theorem which is also considered in \cite{7T2} in the corresponding Hodge-Iwasawa theory:

\begin{theorem} \mbox{\bf{(After Kedlaya-Liu \cite[Theorem 4.6.1]{7KL2})}}  \label{theorem3.11}
Consider the following categories ($A$ is a perfectoid algebra corresponding to some local chart of $D$ in the previous theorem):\\
1. The corresponding category of all the sheaves of pseudocoherent $\mathcal{O}_{\text{\'etale}}$-modules over the adic Fargues-Fontaine curve $Y_{\mathrm{FF},R,A}$ in the corresponding \'etale topology;\\
2. The corresponding category of all the sheaves of pseudocoherent $\mathcal{O}_{\text{pro-\'etale}}$-modules over the adic Fargues-Fontaine curve $Y_{\mathrm{FF},R,A}$ in the corresponding pro-\'etale topology;\\
3. The corresponding category of all the pseudocoherent bundles over the period ring $\widetilde{\Pi}_{R,A}$, carrying the corresponding Frobenius action from the operator $\varphi^a$;\\
4. The corresponding category of all the pseudocoherent modules over the period ring $\widetilde{\Pi}^{\infty}_{R,A}$, carrying the corresponding Frobenius action from the operator $\varphi^a$;\\
5. The corresponding category of all the pseudocoherent modules over the period ring $\widetilde{\Pi}^{}_{R,A}$, carrying the corresponding Frobenius action from the operator $\varphi^a$ (note that here we assume that the module descends to some module carrying Frobenius over $\widetilde{\Pi}^{r}_{R,A}$ for some radius $r_0>0$, which further base changes to some module carrying Frobenius over $\widetilde{\Pi}^{[s,r]}_{R,A}$ for some interval $[s,r]$ which is assumed to be \'etale-stably pseudocoherent);\\
6. The corresponding category of all the finitely generated projective modules over the period ring $\widetilde{\Pi}^{[s,r]}_{R,A}$, carrying the corresponding Frobenius action from the operator $\varphi^a$, where $0<s\leq r/p^{ah}$.\\
%5. The corresponding category of all the finitely generated projective modules over the homotopical period ring $\widetilde{\Pi}^{[s,r],h}_{R,A}$, carrying the corresponding Frobenius action from the operator $\varphi^a$, where $0<s\leq r/p^{ah}$, such that the corresponding underlying modules are plat.\\
Then we have that they are equivalent.
\end{theorem}

\begin{proof}
The proof is parallel to \cite[Theorem 4.6.1]{7KL2}. From 1,2 to 3, the corresponding results follow from \cite[Theorem 4.4.3]{7KL2} due to the natural sheafiness since we are working over perfectoid spaces. Then to glue a bundle to get a module over the Robba ring over the corresponding Robba ring $\widetilde{\Pi}_{R,A}^\infty$, we use the corresponding reified adic space $\mathrm{Spra}(\widetilde{\Pi}_{R,A}^{[sp^{-kah},rp^{-kah}]},\widetilde{\Pi}_{R,A}^{[sp^{-kah},rp^{-kah}],\mathrm{Gr}})$ to cover the corresponding whole spaces, and use the corresponding Frobenius pullback to basically control the corresponding finiteness of the corresponding sections over each member in this covering, then we can derive the result from \cite[Proposition 2.6.17]{7KL2}. The same strategy allows us to go from 3 to 4,5. The corresponding equivalence between 3 and 6 could be proved by the following argument. First from 3 to 6 we just take the corresponding projection. Then back we consider the corresponding Frobenius to control finiteness over each $\mathrm{Spra}(\widetilde{\Pi}_{R,A}^{[sp^{-kah},rp^{-kah}]},\widetilde{\Pi}_{R,A}^{[sp^{-kah},rp^{-kah}],\mathrm{Gr}})$ for any $k\geq 0$, then glueing over any arbitrary interval by using the sheafiness of the period rings. One also checks the corresponding pseudocoherence as in \cite[Theorem 4.6.1]{7KL2}.
\end{proof}

\begin{corollary} \label{corollary3.12}
We have the following categories are equivalent predicted in \cite{7CKZ} partially and implicitly (here $R'$ is a perfect algebra with the same type as that of $R$ but needs not to be the same as $R$):	\\
1. The corresponding category of all the sheaves of pseudocoherent $\mathcal{O}_{\text{\'etale}}$-modules over the adic Fargues-Fontaine curve $Y_{\mathrm{FF},R,Y_{\mathrm{FF},R'}}$ in the corresponding \'etale topology;\\
2. The corresponding category of all the sheaves of pseudocoherent $\mathcal{O}_{\text{pro-\'etale}}$-modules over the adic Fargues-Fontaine curve $Y_{\mathrm{FF},R,Y_{\mathrm{FF},R'}}$ in the corresponding pro-\'etale topology;\\
3. The corresponding category of all the pseudocoherent bundles over the period ring $\widetilde{\Pi}_{R,Y_{\mathrm{FF},R'}}$, carrying the corresponding Frobenius action from the operator $\varphi^a$;\\
4. The corresponding category of all the pseudocoherent modules over the period ring $\widetilde{\Pi}^{\infty}_{R,Y_{\mathrm{FF},R'}}$, carrying the corresponding Frobenius action from the operator $\varphi^a$;\\
5. The corresponding category of all the pseudocoherent modules over the period ring $\widetilde{\Pi}^{}_{R,Y_{\mathrm{FF},R'}}$, carrying the corresponding Frobenius action from the operator $\varphi^a$;\\
6. The corresponding category of all the finitely generated projective modules over the period ring $\widetilde{\Pi}^{[s,r]}_{R,Y_{\mathrm{FF},R'}}$, carrying the corresponding Frobenius action from the operator $\varphi^a$, where $0<s\leq r/p^{ah}$.
\end{corollary}

\indent For some application one could consider the corresponding preperfectoid coefficients:

\begin{theorem} \mbox{\bf{(After Kedlaya-Liu \cite[Theorem 4.6.1]{7KL2})}}  \label{7theorem7.3.18}
Consider the following categories ($A$ is a preperfectoid algebra over $E$):\\
1. The corresponding category of all the sheaves of pseudocoherent $\mathcal{O}_{\text{\'etale}}$-modules over the adic Fargues-Fontaine curve $Y_{\mathrm{FF},R,A}$ in the corresponding \'etale topology;\\
2. The corresponding category of all the sheaves of pseudocoherent $\mathcal{O}_{\text{pro-\'etale}}$-modules over the adic Fargues-Fontaine curve $Y_{\mathrm{FF},R,A}$ in the corresponding pro-\'etale topology;\\
3. The corresponding category of all the pseudocoherent bundles over the period ring $\widetilde{\Pi}_{R,A}$, carrying the corresponding Frobenius action from the operator $\varphi^a$;\\
4. The corresponding category of all the pseudocoherent modules over the period ring $\widetilde{\Pi}^{\infty}_{R,A}$, carrying the corresponding Frobenius action from the operator $\varphi^a$;\\
5. The corresponding category of all the pseudocoherent modules over the period ring $\widetilde{\Pi}^{}_{R,A}$, carrying the corresponding Frobenius action from the operator $\varphi^a$ (note that here we assume that the module descends to some module carrying Frobenius over $\widetilde{\Pi}^{r}_{R,A}$ for some radius $r>0$, which further base changes to some module carrying Frobenius over $\widetilde{\Pi}^{[s,r]}_{R,A}$ for some interval $[s,r]$ which is assumed to be \'etale-stably pseudocoherent);\\
6. The corresponding category of all the finitely generated projective modules over the period ring $\widetilde{\Pi}^{[s,r]}_{R,A}$, carrying the corresponding Frobenius action from the operator $\varphi^a$, where $0<s\leq r/p^{ah}$.\\
%5. The corresponding category of all the finitely generated projective modules over the homotopical period ring $\widetilde{\Pi}^{[s,r],h}_{R,A}$, carrying the corresponding Frobenius action from the operator $\varphi^a$, where $0<s\leq r/p^{ah}$, such that the corresponding underlying modules are plat.\\
Then we have that they are equivalent.
\end{theorem}

\begin{proof}
The proof is parallel to \cite[Theorem 4.6.1]{7KL2}. From 1,2 to 3, the corresponding results follow from \cite[Theorem 4.4.3]{7KL2} due to the natural sheafiness since we are working over preperfectoid spaces. Then to glue a bundle to get a module over the Robba ring over the corresponding Robba ring $\widetilde{\Pi}_{R,A}^\infty$, we use the corresponding reified adic space $\mathrm{Spra}(\widetilde{\Pi}_{R,A}^{[sp^{-kah},rp^{-kah}]},\widetilde{\Pi}_{R,A}^{[sp^{-kah},rp^{-kah}],\mathrm{Gr}})$ to cover the corresponding whole spaces, and use the corresponding Frobenius pullback to basically control the corresponding finiteness of the corresponding sections over each member in this covering, then we can derive the result from \cite[Proposition 2.6.17]{7KL2}. The same strategy allows us to go from 3 to 4,5. The corresponding equivalence between 3 and 6 could be proved by the following argument. First from 3 to 6 we just take the corresponding projection. Then back we consider the corresponding Frobenius to control finiteness over each $\mathrm{Spra}(\widetilde{\Pi}_{R,A}^{[sp^{-kah},rp^{-kah}]},\widetilde{\Pi}_{R,A}^{[sp^{-kah},rp^{-kah}],\mathrm{Gr}})$ for any $k\geq 0$, then glueing over any arbitrary interval by using the sheafiness of the period rings. One also checks the corresponding pseudocoherence as in \cite[Theorem 4.6.1]{7KL2}.\\
\end{proof}

\begin{definition} \mbox{\bf{(After Kedlaya-Liu \cite[Definition 4.4.4]{7KL2})}}
We define the corresponding $\varphi^a$-$\varphi^{'a}$-modules over the corresponding period rings, which we could use some uniform notation $\triangle^?\widehat{\otimes}\triangle^{?'},?,?'=\emptyset,I,r,\infty,\triangle=\widetilde{\Pi}_{R}$ to denote these. We define any $\varphi^a$-$\varphi^{'a}$-module over the period rings with index $\emptyset$ above to be a finitely generated projective module over $\triangle^?\widehat{\otimes}\triangle^{?'}$ carrying semilinear action from the Frobenius operator $\varphi^a$ such that $\varphi^{a*}M\overset{\sim}{\rightarrow}M$ where $\varphi^a$ comes from the first factor, and carrying semilinear action from the Frobenius operator $\varphi^{'a}$ such that $\varphi^{'a*}M\overset{\sim}{\rightarrow}M$ where $\varphi^{'a}$ comes from the second factor, and we assume that this module descends to some $\triangle^r\widehat{\otimes}\triangle^{r'}$ for $0<r,r'<\infty$. We define any $\varphi^a$-$\varphi^{'a}$-module over the period rings with index $r_1,r_2$ above to be a finitely generated projective module over $\triangle^?\widehat{\otimes}\triangle^{?'}$ carrying semilinear action from the Frobenius operator $\varphi^a$ such that $\varphi^{a*}M\overset{\sim}{\rightarrow}M\otimes \triangle^{r_1p^{-ah}}\widehat{\otimes}\triangle^{?'}$, carrying semilinear action from the Frobenius operator $\varphi^{'a}$ such that $\varphi^{'a*}M\overset{\sim}{\rightarrow}M\otimes \triangle^?\widehat{\otimes}\triangle^{r_2p^{-ah}}$. We define any $\varphi^a$-$\varphi^{'a}$-module over the period rings with index $[s_1,r_1],[s_2,r_2]$ above to be a finitely generated projective module over $\triangle^?\widehat{\otimes}\triangle^{?'}$ carrying semilinear action from the Frobenius operator $\varphi^a$ such that $\varphi^{a*}M\otimes_{\triangle^{[s_1p^{-ha},r_1p^{-ha}]}\widehat{\otimes}\triangle^{?'}} \triangle^{[s_1,r_1p^{-ha}]} \widehat{\otimes}\triangle^{?'}\overset{\sim}{\rightarrow}M\otimes \triangle^{[s_1,r_1p^{-ha}]}\widehat{\otimes}\triangle^{?'}$, carrying semilinear action from the Frobenius operator $\varphi^{'a}$ such that $\varphi^{'a*}M\otimes_{\triangle^?\widehat{\otimes}\triangle^{[s_2p^{-ha},r_2p^{-ha}]}} \triangle^?\widehat{\otimes}\triangle^{[s_2,r_2p^{-ha}]}\overset{\sim}{\rightarrow}M\otimes \triangle^?\widehat{\otimes}\triangle^{[s_2,r_2p^{-ha}]}$. For a corresponding object over $\triangle^\emptyset\widehat{\otimes}\triangle^\emptyset$, we assume the module is the base change from some module over $\triangle^{r_0}\widehat{\otimes}\triangle^{r'_0}$ for some $r_0,r_0'>0$, as mentioned above. 
\end{definition}

\begin{definition} \mbox{\bf{(After Kedlaya-Liu \cite[Definition 4.4.4]{7KL2})}}
We define the corresponding pseudocoherent $\varphi^a$-$\varphi^{'a}$-modules over the corresponding period rings above, which we could use some uniform notation $\triangle^?\widehat{\otimes}\triangle^{?'},?,?'=\emptyset,I,r,\infty,\triangle=\widetilde{\Pi}_{R}$ to denote these. We define any pseudocoherent $\varphi^a$-$\varphi^{'a}$-module over the period rings with index $\emptyset$ above to be a pseudocoherent module over $\triangle^?\widehat{\otimes}\triangle^{?'}$ carrying semilinear action from the Frobenius operator $\varphi^a$ such that $\varphi^{a*}M\overset{\sim}{\rightarrow}M$ where $\varphi^a$ comes from the first factor, and carrying semilinear action from the Frobenius operator $\varphi^{'a}$ such that $\varphi^{'a*}M\overset{\sim}{\rightarrow}M$ where $\varphi^{'a}$ comes from the second factor, and we assume that this module descends to some $\triangle^r\widehat{\otimes}\triangle^{r'}$ for $0<r,r'<\infty$. We define any pseudocoherent $\varphi^a$-$\varphi^{'a}$-module over the period rings with index $r_1,r_2$ above to be a pseudocoherent module over $\triangle^?\widehat{\otimes}\triangle^{?'}$ carrying semilinear action from the Frobenius operator $\varphi^a$ such that $\varphi^{a*}M\overset{\sim}{\rightarrow}M\otimes \triangle^{r_1p^{-ah}}\widehat{\otimes}\triangle^{?'}$, carrying semilinear action from the Frobenius operator $\varphi^{'a}$ such that $\varphi^{'a*}M\overset{\sim}{\rightarrow}M\otimes \triangle^?\widehat{\otimes}\triangle^{r_2p^{-ah}}$. We define any pseudocoherent $\varphi^a$-$\varphi^{'a}$-module over the period rings with index $[s_1,r_1],[s_2,r_2]$ above to be a pseudocoherent module over $\triangle^?\widehat{\otimes}\triangle^{?'}$ carrying semilinear action from the Frobenius operator $\varphi^a$ such that $\varphi^{a*}M\otimes_{\triangle^{[s_1p^{-ha},r_1p^{-ha}]}\widehat{\otimes}\triangle^{?'}} \triangle^{[s_1,r_1p^{-ha}]} \widehat{\otimes}\triangle^{?'}\overset{\sim}{\rightarrow}M\otimes \triangle^{[s_1,r_1p^{-ha}]}\widehat{\otimes}\triangle^{?'}$, carrying semilinear action from the Frobenius operator $\varphi^{'a}$ such that $\varphi^{'a*}M\otimes_{\triangle^?\widehat{\otimes}\triangle^{[s_2p^{-ha},r_2p^{-ha}]}} \triangle^?\widehat{\otimes}\triangle^{[s_2,r_2p^{-ha}]}\overset{\sim}{\rightarrow}M\otimes \triangle^?\widehat{\otimes}\triangle^{[s_2,r_2p^{-ha}]}$. For a corresponding object over $\triangle^\emptyset\widehat{\otimes}\triangle^\emptyset$, we assume the module is the base change from some module over $\triangle^{r_0}\widehat{\otimes}\triangle^{r'_0}$ for some $r_0,r_0'>0$. As in \cite[Definition 4.4.4]{7KL2} we impose the corresponding topological condition on the corresponding modules by imposing that all the modules are complete with respect to the corresponding natural topology, and over the Robba rings with respect to some specific multi-interval we assume that the corresponding modules are \'etale-stably pseudocoherent.
\end{definition}

\begin{theorem} \label{7theorem7.3.21}
We have the following categories are equivalent predicted in \cite{7CKZ}:	\\
1. The corresponding category of all the sheaves of pseudocoherent $\mathcal{O}_{\text{\'etale}}$-modules over the adic Fargues-Fontaine curve $Y_{\mathrm{FF},R,Y_{\mathrm{FF},R}}$ in the corresponding \'etale topology;\\
2. The corresponding category of all the sheaves of pseudocoherent $\mathcal{O}_{\text{pro-\'etale}}$-modules over the adic Fargues-Fontaine curve $Y_{\mathrm{FF},R,Y_{\mathrm{FF},R}}$ in the corresponding pro-\'etale topology;\\
3. The corresponding category of all the pseudocoherent modules over the period ring $\widetilde{\Pi}^{[s,r]}_{R}\widehat{\otimes}\widetilde{\Pi}^{[s',r']}_{R}$, carrying the corresponding partial Frobenius action from the operators $\varphi^a$ and $\varphi^{'a}$ from the first $R$ and the second $R$ respectively, here we assume that $0<s\leq r/p^{ah}<\infty$ and $0<s'\leq r'/p^{ah}<\infty$ (which are certainly assumed to be \'etale-stably pseudocoherent).\\
\end{theorem}

\begin{proof}
From 1,2 to 3, we consider the corresponding projection to the subspace with respect to the first factor namely:
\begin{align}
\mathrm{Spa}(\widetilde{\Pi}^{[s,r]}_{R,Y_{\mathrm{FF},R}},\widetilde{\Pi}^{[s,r],+}_{R,Y_{\mathrm{FF},R}})\slash \varphi^\mathbb{Z}.	
\end{align}
By the quasi-compactness of $Y_{\mathrm{FF},R}$ we have that for some fixed nice $[s,r]$ the pseudocoherent $\mathcal{O}$-modules could be regarded in some equivalent way as pseudocoherent $\widetilde{\Pi}^{[s,r]}_{R,Y_{\mathrm{FF},R}}$-modules carrying the partial Frobenius action from $\varphi^a$ (namely the corresponding finiteness is preserved during this functor). Then consider the preperfectoid coefficients in $\widetilde{\Pi}^{[s,r]}_{R}$ we have that this category could be further related to the one of all the pseudocoherent $\widetilde{\Pi}^{[s,r]}_{R}\widehat{\otimes}\widetilde{\Pi}^{[s',r']}_{R}$-modules carrying the partial Frobenius action from $\varphi^a$ and the partial Frobenius action from $\varphi^{'a}$. We then start from any object in 3, we then use the corresponding reified adic space:
\begin{align}
\mathrm{Spra}(\widetilde{\Pi}^{[sp^{-kha},rp^{-kha}]}_{R}\widehat{\otimes}\widetilde{\Pi}^{[s'p^{-k'ha},r'p^{-k'ha}]}_{R},(\widetilde{\Pi}^{[sp^{-kha},rp^{-kha}]}_{R}\widehat{\otimes}\widetilde{\Pi}^{[s'p^{-k'ha},r'p^{-k'ha}]}_{R})^\mathrm{Gr})	
\end{align}
to cover the corresponding spaces showing the corresponding essential surjectivity of the functor from 1 to 3. Again we need to apply \cite[Proposition 2.6.17]{7KL2} to obtain the finiteness around the global section.\\
\end{proof}

%%%%%%%%%%%%%%%%%%%%%%%%%%%%%%%%%%%%%%%%!!!!!!!!!!!!!!!!!!!!!

%%\newpage

\newpage\section{Sheafiness and $\infty$-Period Rings and Sheaves} \label{section4}

\subsection{The $\infty$-Period Rings and Sheaves}

\indent Here we would like to mention something related to the corresponding sheafiness of the Robba rings in our context. Namely when we discuss the corresponding modules over period rings with big coefficients we might come across the issues around the sheafiness. However recent work \cite{7BK} could let us get around the issue by just considering the $\infty$-adic space $(\mathrm{Spa}^h(S), \mathcal{O}^\mathrm{der}_{\mathrm{Spa}^h(S)})$ attached a Banach algebra $S$, which is just an $\infty$-analytic stack.

\indent Following \cite[Introduction]{7BK} we define the following $\infty$-analytic space which are the corresponding coverings of the $\infty$-version of Fargues-Fontaine curves as those considered in \cite{7KL1}, \cite{7KL2} and \cite{7T2}. Here we drop the corresponding sheafiness condition on the ring $\widetilde{\Pi}^{[s,r]}_{R,A}$ for any closed subinterval of $(0,\infty)$. Then we define:
\begin{displaymath}
X^\infty_{\mathrm{FF},R,A}:= \varinjlim_{0<s\leq r< \infty} \mathrm{Spa}^h(\widetilde{\Pi}^{[s,r]}_{R,A})	
\end{displaymath}
which carries a corresponding derived sheaves of $\infty$-rings, namely we have ringed $\infty$-adic space:
\begin{align}
(\varinjlim_{0<s\leq r< \infty} \mathrm{Spa}^h(\widetilde{\Pi}^{[s,r]}_{R,A}), \varprojlim_{0<s\leq r< \infty}\mathcal{O}^\mathrm{der}_{\mathrm{Spa}^h(\widetilde{\Pi}^{[s,r]}_{R,A})}).	
\end{align}

\begin{remark}
In \cite{7BK}, there are two different notions of the corresponding $\infty$-Huber spectra attached to Banach algebras. We choose to focus on $(\mathrm{Spa}^h_\mathrm{Rat}(S), \mathcal{O}^\mathrm{der}_{\mathrm{Spa}_\mathrm{Rat}^h(S)})$ as our $(\mathrm{Spa}^h(S), \mathcal{O}^\mathrm{der}_{\mathrm{Spa}^h(S)})$ in this section. Recall the corresponding Tate acyclicity in this case still holds for derived standard rational localization.	
\end{remark}

\begin{setting}
Note (!) that in our current section we are going to use the notation $A$ to denote any commutative Banach uniform adic algebra over $\mathbb{Q}_p$.
\end{setting}

\begin{definition}
We define the corresponding $\infty$ deformed Fargues-Fontaine curves in general notation $Y^\infty_{\mathrm{FF},R,A}$ to be the corresponding quotient of the space by the $\mathbb{Z}$-power of the Frobenius operator which acts through the corresponding affinoids on the algebraic level.	
\end{definition}

\indent The presheaf of $\infty$-rings $\varprojlim_{0<s\leq r< \infty}\mathcal{O}^\mathrm{der}_{\mathrm{Spa}^h(\widetilde{\Pi}^{[s,r]}_{R,A})}$ are general limit of the corresponding sheaves $\mathcal{O}^\mathrm{der}_{\mathrm{Spa}^h(\widetilde{\Pi}^{[s,r]}_{R,A})}$, which provides upgrading of the discussion in \cite{7T2}. So now we consider similarly as above the following $\infty$-adic spaces:

\begin{align}
\mathrm{Spa}^h\widetilde{\Pi}_{R,A},\mathrm{Spa}^h\widetilde{\Pi}^\infty_{R,A},\mathrm{Spa}^h\widetilde{\Pi}^I_{R,A},\mathrm{Spa}^h\widetilde{\Pi}^r_{R,A}.	
\end{align}
%%%%%%%%%

\indent Then take the corresponding global section we have the following $\infty$-period rings:

\begin{align}
\widetilde{\Pi}^h_{R,A},\widetilde{\Pi}^{\infty,h}_{R,A},\widetilde{\Pi}^{I,h}_{R,A},\widetilde{\Pi}^{r,h}_{R,A}.	
\end{align}

\begin{definition}
We define the corresponding Frobenius (on the $\infty$-rings) as in the corresponding non-derived situation through the corresponding affinoids on the algebraic level.	
\end{definition}

%\indent In the following we define the corresponding finitely generated projective objects in more general sense from \cite[Definition 7.2.2.1, Definition 7.2.2.4]{7Lu1}. Recall that the corresponding $finitely~generated~free~ modules$ over any connective $\infty$ algebra $S$ are defined to ones which could be written as coproducts of $S$ in terms of some finite copies, and recall that a $projective$ $S$-module $M$ is defined to be a projective object in the corresponding $\infty$-category of connective $S$-modules.

%\begin{definition} \mbox{\bf{(After Lurie \cite[Corollary 7.2.2.9,Definition 7.2.2.10]{7Lu1})}}
%For general projective with finiteness condition we look at the following definition following \cite[Corollary 7.2.2.9]{7Lu1}. We will use the notion a $finitely~generated~projective$ module over a connective $\mathbb{E}_\infty$-algebra $S$ to mean a projective object $M$ in the corresponding $\infty$-category of all the connective $S$-modules such that we have $\pi_0M$ is finitely generated over $\pi_0 S$. We will use the notion \cite[Definition 7.2.2.10]{7Lu1} a $plat$ module over a connective $\mathbb{E}_\infty$-algebra $S$ to mean a connective $S$-modules such that we have $\pi_0M$ is flat over $\pi_0 S$ and the base change of $\pi_0 M$ to each $\pi_nS$ will be then $\pi_n M$ for each $n\geq 0$. 	
%\end{definition}

\begin{definition} \mbox{\bf{(After Lurie \cite[Corollary 7.2.2.9,Definition 7.2.2.10]{7Lu1})}}
For general projective with finiteness condition we look at the following definition following \cite[Corollary 7.2.2.9]{7Lu1}. We will use the notion a $f$-$projective$ module over a $\infty$-algebra $S$ to mean a projective object $M$ in the corresponding $\infty$-category of all the $S$-modules such that we have $\pi_0M$ is finite projective over $\pi_0 S$. 	
\end{definition}

%\begin{remark}
%The modules here should have the same vanishing on the homotopy groups w	
%\end{remark}

%\begin{example}
%By \cite[Proposition 7.2.2.18]{7Lu1} a connective module over a connective $\mathbb{E}_\infty$-ring $S$ is plat finitely generated projective if and only if $\pi_0M$ is finite projective over $\pi_0 S$.	
%\end{example}

\begin{definition} \mbox{\bf{(After Kedlaya-Liu \cite[Definition 4.4.4]{7KL2})}}
We define the corresponding $\varphi^a$-modules over the corresponding $\infty$-period rings above, which we could use some uniform notation $\triangle^{?,h},?=\emptyset,I,r,\infty$ to denote these. We define any $\varphi^a$-module over the $\infty$-period rings with index $\emptyset$ above to be a $f$-projective module over $\triangle^{?,h}$ carrying semilinear action from the Frobenius operator $\varphi^a$ such that $\varphi^{a*}M\overset{\sim}{\rightarrow}M$. We define any $\varphi^a$-module over the derived period rings with index $r$ above to be a $f$-projective module over $\triangle^{?,h}$ carrying semilinear action from the Frobenius operator $\varphi^a$ such that $\varphi^{a*}M\overset{\sim}{\rightarrow}M\otimes \triangle^{rp^{-ah},h}$. We define any $\varphi^a$-module over the $\infty$-period rings with index $[s,r]$ above to be a $f$-projective module over $\triangle^{?,h}$ carrying semilinear action from the Frobenius operator $\varphi^a$ such that $\varphi^{a*}M\otimes \triangle^{[s,rp^{-ha}],h}\overset{\sim}{\rightarrow}M\otimes_{\triangle^{[sp^{-ha},rp^{-ha}],h}} \triangle^{[s,rp^{-ha}],h}$. Here $\triangle=\widetilde{\Pi}_{R,A}$. For a corresponding object over $\triangle^{\emptyset,h}$, we assume the module is the base changes from some module over $\triangle^{r_0,h}$ for some $r_0>0$.
\end{definition}

\begin{remark}
We believe that one can discuss more general objects such as the corresponding coherent sheaves after \cite{7Lu1} and \cite{7Lu2}, but at this moment let us just focus on the derived vector bundles.	
\end{remark}

\begin{definition} \mbox{\bf{(After Kedlaya-Liu \cite[Definition 4.4.6]{7KL2})}}
We now define the corresponding $\varphi^a$-bundles over the corresponding $\infty$-period rings above, which we could use some uniform notation $\triangle_{*,A}^{?,h},?=\emptyset,I,r,\infty,*=R$ to denote these. We define a $\varphi^a$-bundle over $\triangle_{*,A}^{?,h},?=\emptyset,r$ to be a compatible family of $f$-projective $\varphi^a$-modules over $\triangle_{*,A}^{?,h},?=[s',r']$ for suitable $[s',r']$ (namely contained in $(0,r]$) satisfying the corresponding restriction compatibility and cocycle condition.
\end{definition}

\subsection{Some Results on the $\infty$-Descent}

\begin{conjecture} \mbox{\bf{(After Kedlaya-Liu \cite[Theorem 4.6.1]{7KL2})}}
Consider the following categories:\\
%1. The corresponding category of all the sheaves of $f$-projective $\mathcal{O}^\mathrm{der}_{Y^\infty_{\mathrm{FF},R,A}}$-modules over the derived adic Fargues-Fontaine curve $Y^\infty_{\mathrm{FF},R,A}$;\\
1. The corresponding category of all the bundles over the homotopical period ring $\widetilde{\Pi}^h_{R,A}$, carrying the corresponding Frobenius action from the operator $\varphi^a$;\\
2. The corresponding category of all the $f$-projective modules over the homotopical period ring $\widetilde{\Pi}^{\infty,h}_{R,A}$, carrying the corresponding Frobenius action from the operator $\varphi^a$;\\
%3. The corresponding category of all the $f$-projective modules over the homotopical period ring $\widetilde{\Pi}^{h}_{R,A}$, carrying the corresponding Frobenius action from the operator $\varphi^a$;\\
3. The corresponding category of all the $f$-projective modules over the homotopical period ring $\widetilde{\Pi}^{[s,r],h}_{R,A}$, carrying the corresponding Frobenius action from the operator $\varphi^a$, where $0<s\leq r/p^{ah}$.\\
%5. The corresponding category of all the finitely generated projective modules over the homotopical period ring $\widetilde{\Pi}^{[s,r],h}_{R,A}$, carrying the corresponding Frobenius action from the operator $\varphi^a$, where $0<s\leq r/p^{ah}$, such that the corresponding underlying modules are plat.\\
Then we have that they are equivalent.
\end{conjecture}

\indent In the corresponding sheafy case, we do have some good comparison:

\begin{theorem} \mbox{\bf{(After Kedlaya-Liu \cite[Theorem 4.6.1]{7KL2})}}
In the situation where the Robba ring $\widetilde{\Pi}^{[s,r]}_{R,A}$ for any closed interval $[s,r]$ is stably uniform. Consider the following categories:\\
1. The corresponding category of all the bundles over the homotopical period ring $\widetilde{\Pi}^h_{R,A}$, carrying the corresponding Frobenius action from the operator $\varphi^a$;\\
2. The corresponding category of all the $f$-projective modules over the homotopical period ring $\widetilde{\Pi}^{\infty,h}_{R,A}$, carrying the corresponding Frobenius action from the operator $\varphi^a$;\\
%3. The corresponding category of all the $f$-projective modules over the homotopical period ring $\widetilde{\Pi}^{h}_{R,A}$, carrying the corresponding Frobenius action from the operator $\varphi^a$;\\
3. The corresponding category of all the $f$-projective modules over the homotopical period ring $\widetilde{\Pi}^{[s,r],h}_{R,A}$, carrying the corresponding Frobenius action from the operator $\varphi^a$, where $0<s\leq r/p^{ah}$.
%5. The corresponding category of all the finitely generated projective modules over the homotopical period ring $\widetilde{\Pi}^{[s,r],h}_{R,A}$, carrying the corresponding Frobenius action from the operator $\varphi^a$, where $0<s\leq r/p^{ah}$, such that the corresponding underlying modules are plat.\\
Then we have that they are equivalent.
\end{theorem}

\begin{proof}
This amounts to the corresponding statement for modules over classical rings. The proof is parallel to \cite[Theorem 4.6.1]{7KL2}. Then to glue a bundle to get a module over the Robba ring over the corresponding Robba ring $\widetilde{\Pi}_{R,A}^\infty$, we use the corresponding reified adic space $\mathrm{Spra}(\widetilde{\Pi}_{R,A}^{[sp^{-kah},rp^{-kah}]},\widetilde{\Pi}_{R,A}^{[sp^{-kah},rp^{-kah}],\mathrm{Gr}})$ to cover the corresponding whole spaces, and use the corresponding Frobenius pullback to basically control the corresponding finiteness of the corresponding sections over each member in this covering, then we can derive the result from \cite[Proposition 2.6.17, Corollary 2.6.10]{7KL2}. This will show the equivalence between 1 and 2. Then the functor from 1 to 3 is just the corresponding projection. On the other hand the functor from 3 to 1, we use the corresponding Frobenius to reach any interval taking the form of $[sp^{-nha},rp^{-nha}]$ for any $n\in \mathbb{Z}$, then we still have to consider the corresponding extraction of a single module from two over some overlapped such specific intervals coming from the Frobenius. However this is achievable since for $\pi_0$ of the corresponding two modules we are done through \cite[Theorem 1.3.9]{7KL1}. 	
\end{proof}

\begin{remark}
The stably-uniform condition here relates directly to the construction of \cite{7BK}. 	
\end{remark}

\begin{theorem} \mbox{\bf{(After Kedlaya-Liu \cite[Theorem 4.6.1]{7KL2})}} \label{7theorem4.12}
Consider the following categories:\\
%1. The corresponding category of all the sheaves of finitely generated projective $\mathcal{O}^\mathrm{der}_{Y^\infty_{\mathrm{FF},R,A}}$-modules over the derived adic Fargues-Fontaine curve $Y^\infty_{\mathrm{FF},R,A}$, such that the corresponding underlying modules are plat;\\
1. The corresponding category of all the bundles over the homotopical period ring $\widetilde{\Pi}^h_{R,A}$, carrying the corresponding Frobenius action from the operator $\varphi^a$, such that the corresponding underlying modules are plat (here we assume that the modules are $f$-projective);\\
%2. The corresponding category of all the finitely generated projective modules over the homotopical period ring $\widetilde{\Pi}^{\infty,h}_{R,A}$, carrying the corresponding Frobenius action from the operator $\varphi^a$, such that the corresponding underlying modules are plat;\\
2. The corresponding category of all the $f$-projective modules over the homotopical period ring $\widetilde{\Pi}^{[s,r],h}_{R,A}$, carrying the corresponding Frobenius action from the operator $\varphi^a$, where $0<s\leq r/p^{ah}$, such that the corresponding underlying modules are plat.\\
%5. The corresponding category of all the finitely generated projective modules over the homotopical period ring $\widetilde{\Pi}^{[s,r],h}_{R,A}$, carrying the corresponding Frobenius action from the operator $\varphi^a$, where $0<s\leq r/p^{ah}$, such that the corresponding underlying modules are plat.\\
Then we have that they are equivalent.
\end{theorem}

\begin{proof}
Due to the fact that the $\infty$-presheaf $\mathcal{O}^\mathrm{der}_{\mathrm{Spa}^h(\widetilde{\Pi}^{I,h}_{R,A})}$ is a $\infty$-sheaf, the proof reduces to \cite[Theorem 4.6.1]{7KL2} which we also performed in \cite[Section 4.2]{7T2}. To compare, the functor from 1 to 2 is just the corresponding projection. On the other hand the functor from 2 to 1, we use the corresponding Frobenius to reach any interval taking the form of $[sp^{-nha},rp^{-nha}]$ for any $n\in \mathbb{Z}$, then we still have to consider the corresponding extraction of a single module from two over some overlapped such specific intervals coming from the Frobenius. However this is achievable since for $\pi_0$ of the corresponding two modules we are done through \cite[Theorem 1.3.9]{7KL1}.
\end{proof}

\begin{remark}
Professor Kedlaya has informed us that the condition on the spectrum in \cite[Theorem 1.3.9 (b)]{7KL1} could be removed (see our amplification in \cref{section5} below namely \cref{proposition5.11con}). This made us believe above holds. We want to mention here that since as mentioned in \cite{7BK} the authors of \cite{7BK} believe there is deep relationship between \cite{7BK} and Clausen-Scholze's work \cite{7CS}. Therefore we believe there is deep relationship between ours and the possible ones after Clausen-Scholze \cite{7CS}. 
\end{remark}

%%\newpage

\newpage\section{Noncommutative Descent Revisit} \label{section5}

\subsection{Big Noncommutative Coefficients}

\indent We now also consider the corresponding more general noncommutative Banach coefficients as in our previous work \cite{7T2}, which is targeted at the corresponding $K$-theory of the general Robba rings.

\begin{setting}
Note (!) that in our current section we are going to use the notation $B$ to denote any Banach algebra over $E$.
\end{setting}

\begin{definition} \mbox{\bf{(After Kedlaya-Liu \cite[Definition 4.1.1]{7KL2})}}
We let $B$ be a Banach algebra over $E$ with integral subring $\mathcal{O}_B$ over $\mathcal{O}_E$. Recall from the corresponding context in \cite{7KL1} we have the corresponding period rings in the relative setting. We follow the corresponding notations we used in \cite[Section 2.1]{7T2} for those corresponding period rings. We first have for a pair $(R,R^+)$ where $R$ is a uniform perfect adic Banach ring over the assumed base $\mathcal{O}_F$. Then we take the corresponding generalized Witt vectors taking the form of $W_{\mathcal{O}_E}(R)$, which is just the ring $\widetilde{\Omega}_R^\mathrm{int}$, and by inverting the corresponding uniformizer we have the ring $\widetilde{\Omega}_R$, then by taking the completed product with $B$ we have $\widetilde{\Omega}_{R,B}$. Now for some $r>0$ we consider the ring $\widetilde{\Pi}^{\mathrm{int},r}_{R}$ which is the completion of $W_{\mathcal{O}_E}(R^+)[[r]:r\in R] $ by the norm $\|.\|_{\alpha^r}$ defined by:
\begin{align}
\|.\|_{\alpha^r}(\sum_{n\geq 0}\pi^n[\overline{r}_n])=\sup_{n\geq 0}\{p^{-n}\alpha(\overline{r}_n)^r\}.	
\end{align}
Then we have the product $\widetilde{\Pi}^{\mathrm{bd},r}_{R,B}$ defined as the completion under the product norm $\|.\|_{\alpha^r}\otimes \|.\|_B$ of the corresponding ring $\widetilde{\Pi}^{\mathrm{bd},r}_{R}\otimes_{\mathbb{Q}_p}B$, which could be also defined from the corresponding integral Robba rings defined above. Then we define the corresponding Robba ring for some interval $I\subset (0,\infty)$ with coefficient in the perfectoid ring $B$ denoted by $\widetilde{\Pi}^{I}_{R,B}$ as the following complete tensor product:
\begin{displaymath}
\widetilde{\Pi}^{I}_{R}\widehat{\otimes}_{\mathbb{Q}_p} B	
\end{displaymath}
under the the corresponding tensor product norm $\|.\|_{\alpha^r}\otimes \|.\|_B$. Then we set $\widetilde{\Pi}^r_{R,B}$ as $\varprojlim_{s\rightarrow 0}\widetilde{\Pi}^{[s,r]}_{R,B}$, and then we define $\widetilde{\Pi}_{R,B}$ as $\varinjlim_{r\rightarrow \infty}\widetilde{\Pi}^{[s,r]}_{R,B}$, and we define $\widetilde{\Pi}^\infty_{R,B}$ as $\varprojlim_{r\rightarrow \infty}\widetilde{\Pi}^{r}_{R,B}$. And we also have the corresponding full integral Robba ring and the corresponding full bounded Robba ring by taking the corresponding union through all $r>0$.
\end{definition}

\indent We first define the corresponding right Frobenius modules and bundles:

\begin{definition} \mbox{\bf{(After Kedlaya-Liu \cite[Definition 4.4.4]{7KL2})}}
We define the corresponding right $\varphi^a$-modules over the corresponding period rings above, which we could use some uniform notation $\triangle^?,?=\emptyset,I,r,\infty,\triangle=\widetilde{\Pi}_{R,B}$ to denote these. We define any right $\varphi^a$-module over the period rings with index $\empty$ above to be a finitely generated projective right module over $\triangle^?$ carrying semilinear action from the Frobenius operator $\varphi^a$ such that $\varphi^{a*}M\overset{\sim}{\rightarrow}M$. We define any right $\varphi^a$-module over the period rings with index $r$ above to be a finitely generated projective right module over $\triangle^?$ carrying semilinear action from the Frobenius operator $\varphi^a$ such that $\varphi^{a*}M\overset{\sim}{\rightarrow}M\otimes \triangle^{rp^{-ah}}$. We define any right $\varphi^a$-module over the period rings with index $[s,r]$ above to be a finitely generated projective right module over $\triangle^?$ carrying semilinear action from the Frobenius operator $\varphi^a$ such that $\varphi^{a*}M\otimes_{\triangle^{[sp^{-ha},rp^{-ha}]}} \triangle^{[s,rp^{-ha}]}\overset{\sim}{\rightarrow}M\otimes \triangle^{[s,rp^{-ha}]}$. For a corresponding object over $\triangle^\emptyset$, we assume the module is the base changes from some module over $\triangle^{r_0}$ for some $r_0$.
\end{definition}

\begin{definition} \mbox{\bf{(After Kedlaya-Liu \cite[Definition 4.4.6]{7KL2})}}
We now define the corresponding right $\varphi^a$-bundles over the corresponding period rings above, which we could use some uniform notation $\triangle_{*,B}^?,?=\emptyset,I,r,\infty,*=R$ to denote these. We define a right $\varphi^a$-bundle over $\triangle_{*,B}^?,?=\emptyset,r$ to be a compatible family of finitely generated projective right $\varphi^a$-modules over $\triangle_{*,B}^?,?=[s',r']$ for suitable $[s',r']$ (namely contained in $(0,r]$) satisfying the corresponding restriction compatibility and cocycle condition.
\end{definition}

\indent We first define the corresponding left Frobenius modules and bundles:

\begin{definition} \mbox{\bf{(After Kedlaya-Liu \cite[Definition 4.4.4]{7KL2})}}
We define the corresponding left $\varphi^a$-modules over the corresponding period rings above, which we could use some uniform notation $\triangle^?,?=\emptyset,I,r,\infty,\triangle=\widetilde{\Pi}_{R,B}$ to denote these. We define any left $\varphi^a$-module over the period rings with index $\empty$ above to be a finitely generated projective left module over $\triangle^?$ carrying semilinear action from the Frobenius operator $\varphi^a$ such that $\varphi^{a*}M\overset{\sim}{\rightarrow}M$. We define any left $\varphi^a$-module over the period rings with index $r$ above to be a finitely generated projective left module over $\triangle^?$ carrying semilinear action from the Frobenius operator $\varphi^a$ such that $\varphi^{a*}M\overset{\sim}{\rightarrow} \triangle^{rp^{-ah}}\otimes M$. We define any left $\varphi^a$-module over the period rings with index $[s,r]$ above to be a finitely generated projective left module over $\triangle^?$ carrying semilinear action from the Frobenius operator $\varphi^a$ such that $\triangle^{[s,rp^{-ha}]} \otimes_{\triangle^{[sp^{-ha},rp^{-ha}]}} \varphi^{a*}M \overset{\sim}{\rightarrow} \triangle^{[s,rp^{-ha}]} \otimes M$. For a corresponding object over $\triangle^{\emptyset}$, we assume the module is the base changes from some module over $\triangle^{r_0}$ for some $r_0>0$.

\end{definition}

\begin{definition} \mbox{\bf{(After Kedlaya-Liu \cite[Definition 4.4.6]{7KL2})}}
We now define the corresponding left $\varphi^a$-bundles over the corresponding period rings above, which we could use some uniform notation $\triangle_{*,B}^?,?=\emptyset,I,r,\infty,*=R$ to denote these. We define a left $\varphi^a$-bundle over $\triangle_{*,B}^?,?=\emptyset,r$ to be a compatible family of finitely generated projective left $\varphi^a$-modules over $\triangle_{*,B}^?,?=[s',r']$ for suitable $[s',r']$ (namely contained in $(0,r]$) satisfying the corresponding restriction compatibility and cocycle condition.
\end{definition}

\indent We first define the corresponding Frobenius bimodules and bibundles:

\begin{definition} \mbox{\bf{(After Kedlaya-Liu \cite[Definition 4.4.4]{7KL2})}}
We define the corresponding $\varphi^a$-bimodules over the corresponding period rings above, which we could use some uniform notation $\triangle^?,?=\emptyset,I,r,\infty,\triangle=\widetilde{\Pi}_{R,B}$ to denote these. We define any $\varphi^a$-bimodule over the period rings with index $\empty$ above to be a finitely generated projective bimodule over $\triangle^?$ carrying semilinear action from the Frobenius operator $\varphi^a$ such that $\varphi^{a*}M\overset{\sim}{\rightarrow}M$. We define any $\varphi^a$-bimodule over the period rings with index $r$ above to be a finitely generated projective bimodule over $\triangle^?$ carrying semilinear action from the Frobenius operator $\varphi^a$ such that $\varphi^{a*}M\overset{\sim}{\rightarrow} \triangle^{rp^{-ah}}\otimes M$ and $\varphi^{a*}M\overset{\sim}{\rightarrow}M\otimes \triangle^{rp^{-ah}}$. We define any $\varphi^a$-bimodule over the period rings with index $[s,r]$ above to be a finitely generated projective bimodule over $\triangle^?$ carrying semilinear action from the Frobenius operator $\varphi^a$ such that $\triangle^{[s,rp^{-ha}]} \otimes_{\triangle^{[sp^{-ha},rp^{-ha}]}} \varphi^{a*}M \overset{\sim}{\rightarrow} \triangle^{[s,rp^{-ha}]} \otimes M$ and $\varphi^{a*}M\otimes_{\triangle^{[sp^{-ha},rp^{-ha}]}} \triangle^{[s,rp^{-ha}]}\overset{\sim}{\rightarrow}M\otimes \triangle^{[s,rp^{-ha}]}$. For a corresponding object over $\triangle^{\emptyset}$, we assume the module is the base changes from some module over $\triangle^{r_0}$ for some $r_0>0$.

\end{definition}

\begin{definition} \mbox{\bf{(After Kedlaya-Liu \cite[Definition 4.4.6]{7KL2})}}
We now define the corresponding $\varphi^a$-bibundles over the corresponding period rings above, which we could use some uniform notation $\triangle_{*,B}^?,?=\emptyset,I,r,\infty,*=R$ to denote these. We define a $\varphi^a$-bibundle over $\triangle_{*,B}^?,?=\emptyset,r$ to be a compatible family of finitely generated projective $\varphi^a$-bimodules over $\triangle_{*,B}^?,?=[s',r']$ for suitable $[s',r']$ (namely contained in $(0,r]$) satisfying the corresponding restriction compatibility and cocycle condition.\\
\end{definition}

\subsection{Glueing Noncommutative Vector Bundles}

\indent The following is what we achieved in the corresponding paper \cite[Lemma 6.82]{7T2}:

\begin{proposition}
Consider the following exact sequence of Banach algebras satisfying the corresponding conditions in \cite[Definition 2.7.3 (a),(b)]{7KL1}:
\begin{align}
0\rightarrow {\Pi}\rightarrow \Pi_1\bigoplus \Pi_2\rightarrow \Pi_{12}\rightarrow 0.	
\end{align}
And consider the corresponding right glueing datum $(M_1,M_2,M_{12})$ over the corresponding rings above, and we assume that the glueing datum is right finite projective. Then we have that the corresponding kernel $M$ of
\begin{align}
M_1\bigoplus M_2\rightarrow M_{12}	
\end{align}
as right module over $\Pi$ is finitely presented, with the corresponding isomorphisms:
\begin{align}
M\otimes \Pi_1 \overset{\sim}{\rightarrow}	M_1,\\
M\otimes \Pi_2 \overset{\sim}{\rightarrow}	M_2.
\end{align}

\end{proposition}

\begin{proof}
This is essentially proved in \cite[Lemma 6.82]{7T2}.	
\end{proof}

\indent By considering the parallel argument one has the following:

\begin{proposition}
Consider the following exact sequence of Banach algebras satisfying the corresponding conditions in \cite[Definition 2.7.3 (a),(b)]{7KL1}:
\begin{align}
0\rightarrow {\Pi}\rightarrow \Pi_1\bigoplus \Pi_2\rightarrow \Pi_{12}\rightarrow 0.	
\end{align}
And consider the corresponding left glueing datum $(M_1,M_2,M_{12})$ over the corresponding rings above, and we assume that the glueing datum is left finite projective. Then we have that the corresponding kernel $M$ of
\begin{align}
M_1\bigoplus M_2\rightarrow M_{12}	
\end{align}
as left module over $\Pi$ is finitely presented, with the corresponding isomorphisms:
\begin{align}
M\otimes \Pi_1 \overset{\sim}{\rightarrow}	M_1,\\
M\otimes \Pi_2 \overset{\sim}{\rightarrow}	M_2.
\end{align}

\end{proposition}

\indent Then we consider the corresponding bimodules, and by relying on some argument essentially due to Kedlaya we have the following:

\begin{proposition} \mbox{\bf{(Kedlaya)}} \label{proposition5.11con}
Assume now the Banach algebras are commutative. Consider the following exact sequence of Banach algebras satisfying the corresponding conditions in \cite[Definition 2.7.3 (a),(b)]{7KL1}:
\begin{align}
0\rightarrow {\Pi}\rightarrow \Pi_1\bigoplus \Pi_2\rightarrow \Pi_{12}\rightarrow 0.	
\end{align}
And consider the corresponding glueing datum $(M_1,M_2,M_{12})$ over the corresponding rings above, and we assume that the glueing datum is finite projective. Then we have that the corresponding kernel $M$ of
\begin{align}
M_1\bigoplus M_2\rightarrow M_{12}	
\end{align}
as module over $\Pi$ is finite projective, with the corresponding isomorphisms:
\begin{align}
M\otimes \Pi_1 \overset{\sim}{\rightarrow}	M_1,\\
M\otimes \Pi_2 \overset{\sim}{\rightarrow}	M_2.
\end{align}

\end{proposition}

\begin{proof}
This proof is due to Kedlaya. The corresponding finitely presentedness is given in \cite[Theorem 1.3.9(a)]{7KL1}, note here that we are under the conditions (a), (b) in \cite[Definition 2.7.3 (a),(b)]{7KL1}. To promote the corresponding property to that being finite projective. We look at the corresponding diagram coming from the desired presentation:
\[
\xymatrix@R+3pc@C+3pc{
B \ar[r] \ar[r] \ar[r] \ar[d] \ar[d] \ar[d] &B_1\oplus B_2 \ar[r] \ar[r] \ar[r] \ar[d] \ar[d] \ar[d] &B_{12} \ar[d] \ar[d] \ar[d]\\
A \ar[r] \ar[r] \ar[r] \ar[d] \ar[d] \ar[d] &A_1\oplus A_2 \ar[r] \ar[r] \ar[r] \ar[d] \ar[d] \ar[d] &A_{12} \ar[d] \ar[d] \ar[d]\\
M \ar[r] \ar[r] \ar[r]  &M_1\oplus M_2 \ar[r] \ar[r] \ar[r]  &M_{12}. \\
}
\]
Now we look at the sets:
\begin{align}
\mathrm{Hom}(M_1,B_1),\mathrm{Hom}(M_2,B_2)	
\end{align}
which actually are also modules over the corresponding base rings. By \cite[Theorem 1.3.9(a)]{7KL1} we can basically apply the construction to 
\begin{align}
\mathrm{Hom}(M_1,B_1),\mathrm{Hom}(M_2,B_2),\mathrm{Hom}(M_{12},B_{12}) 	
\end{align}	
which relates directly to a glueing square, which gives rise to the surjective morphism:
\begin{align}
\mathrm{Hom}(M_1,B_1)\oplus\mathrm{Hom}(M_2,B_2)\rightarrow \mathrm{Hom}(M_{12},B_{12}).	
\end{align}	
Then we choose the corresponding splittings $s_1,s_2$ for $M_1$ and $M_2$ respectively. Then map them to $\mathrm{Hom}(M_{12},B_{12})$. Then we have there is an element $s\in \mathrm{Hom}(M_{12},B_{12})$ such that $f_{1,12}(s_1)-f_{2,12}(s_2)=s$. Then choose the corresponding preimage $s'_1+s'_2$ of $s$ in $\mathrm{Hom}(M_1,B_1)\oplus\mathrm{Hom}(M_2,B_2)$ through:
\begin{align}
\mathrm{Hom}(M_1,B_1)\oplus\mathrm{Hom}(M_2,B_2)\rightarrow \mathrm{Hom}(M_{12},B_{12}).	
\end{align}
This makes the corresponding correct modification to the corresponding splitting $s_1+s_2$ to make sure that $s_1''+s_2''$ lives in the kernel of 
\begin{align}
\mathrm{Hom}(M_1,B_1)\oplus\mathrm{Hom}(M_2,B_2)\rightarrow \mathrm{Hom}(M_{12},B_{12}).	
\end{align}
Then this gives rise to a corresponding splitting for $M$.
\end{proof}

\begin{proposition} \mbox{\bf{(After Kedlaya)}} \label{proposition5.11}
Consider the following exact sequence of Banach algebras satisfying the corresponding conditions in \cite[Definition 2.7.3 (a),(b)]{7KL1}:
\begin{align}
0\rightarrow {\Pi}\rightarrow \Pi_1\bigoplus \Pi_2\rightarrow \Pi_{12}\rightarrow 0.	
\end{align}
And consider the corresponding glueing datum $(M_1,M_2,M_{12})$ over the corresponding rings above, and we assume that the glueing datum is finite projective. Then we have that the corresponding kernel $M$ of
\begin{align}
M_1\bigoplus M_2\rightarrow M_{12}	
\end{align}
as bimodule over $\Pi$ is finite projective, with the corresponding isomorphisms:
\begin{align}
M\otimes \Pi_1 \overset{\sim}{\rightarrow}	M_1,\\
M\otimes \Pi_2 \overset{\sim}{\rightarrow}	M_2.
\end{align}

\end{proposition}

\begin{proof}
The left module and the right module structures have already given us the corresponding finitely presentedness. To promote the corresponding property to that being finite projective. We look at the corresponding diagram coming from the desired presentation:
\[
\xymatrix@R+3pc@C+3pc{
B \ar[r] \ar[r] \ar[r] \ar[d] \ar[d] \ar[d] &B_1\oplus B_2 \ar[r] \ar[r] \ar[r] \ar[d] \ar[d] \ar[d] &B_{12} \ar[d] \ar[d] \ar[d]\\
A \ar[r] \ar[r] \ar[r] \ar[d] \ar[d] \ar[d] &A_1\oplus A_2 \ar[r] \ar[r] \ar[r] \ar[d] \ar[d] \ar[d] &A_{12} \ar[d] \ar[d] \ar[d]\\
M \ar[r] \ar[r] \ar[r]  &M_1\oplus M_2 \ar[r] \ar[r] \ar[r]  &M_{12}. \\
}
\]
Now we look at the sets:
\begin{align}
\mathrm{Hom}(M_1,B_1),\mathrm{Hom}(M_2,B_2)	
\end{align}
which actually are also left modules. By the left glueing process we consider before we can basically apply the construction to 
\begin{align}
\mathrm{Hom}(M_1,B_1),\mathrm{Hom}(M_2,B_2),\mathrm{Hom}(M_{12},B_{12}) 	
\end{align}	
which relates directly to a glueing square, which gives rise to the surjective morphism:
\begin{align}
\mathrm{Hom}(M_1,B_1)\oplus\mathrm{Hom}(M_2,B_2)\rightarrow \mathrm{Hom}(M_{12},B_{12}).	
\end{align}	
This is again as left module consideration. Then we choose the corresponding splittings $s_1,s_2$ for $M_1$ and $M_2$ respectively. Then map them to $\mathrm{Hom}(M_{12},B_{12})$. Then we have there is an element $s\in \mathrm{Hom}(M_{12},B_{12})$ such that $f_{1,12}(s_1)-f_{2,12}(s_2)=s$. Then choose the corresponding preimage $s'_1+s'_2$ of $s$ in $\mathrm{Hom}(M_1,B_1)\oplus\mathrm{Hom}(M_2,B_2)$ through:
\begin{align}
\mathrm{Hom}(M_1,B_1)\oplus\mathrm{Hom}(M_2,B_2)\rightarrow \mathrm{Hom}(M_{12},B_{12}).	
\end{align}
This makes the corresponding correct modification to the corresponding splitting $s_1+s_2$ to make sure that $s_1''+s_2''$ lives in the kernel of 
\begin{align}
\mathrm{Hom}(M_1,B_1)\oplus\mathrm{Hom}(M_2,B_2)\rightarrow \mathrm{Hom}(M_{12},B_{12}).	
\end{align}
Then this gives rise to a corresponding splitting for $M$ as left module. Then on the other hand, we look at the corresponding diagram coming from the desired presentation:
\[
\xymatrix@R+3pc@C+3pc{
B \ar[r] \ar[r] \ar[r] \ar[d] \ar[d] \ar[d] &B_1\oplus B_2 \ar[r] \ar[r] \ar[r] \ar[d] \ar[d] \ar[d] &B_{12} \ar[d] \ar[d] \ar[d]\\
A \ar[r] \ar[r] \ar[r] \ar[d] \ar[d] \ar[d] &A_1\oplus A_2 \ar[r] \ar[r] \ar[r] \ar[d] \ar[d] \ar[d] &A_{12} \ar[d] \ar[d] \ar[d]\\
M \ar[r] \ar[r] \ar[r]  &M_1\oplus M_2 \ar[r] \ar[r] \ar[r]  &M_{12}. \\
}
\]
Now we look at the sets:
\begin{align}
\mathrm{Hom}(M_1,B_1),\mathrm{Hom}(M_2,B_2)	
\end{align}
which actually are also right modules. By the right glueing process we consider before we can basically apply the construction to 
\begin{align}
\mathrm{Hom}(M_1,B_1),\mathrm{Hom}(M_2,B_2),\mathrm{Hom}(M_{12},B_{12}) 	
\end{align}	
which relates directly to a glueing square, which gives rise to the surjective morphism:
\begin{align}
\mathrm{Hom}(M_1,B_1)\oplus\mathrm{Hom}(M_2,B_2)\rightarrow \mathrm{Hom}(M_{12},B_{12}).	
\end{align}	
This is again as right module consideration. Then we choose the corresponding splittings $s_1,s_2$ for $M_1$ and $M_2$ respectively. Then map them to $\mathrm{Hom}(M_{12},B_{12})$. Then we have there is an element $s\in \mathrm{Hom}(M_{12},B_{12})$ such that $f_{1,12}(s_1)-f_{2,12}(s_2)=s$. Then choose the corresponding preimage $s'_1+s'_2$ of $s$ in $\mathrm{Hom}(M_1,B_1)\oplus\mathrm{Hom}(M_2,B_2)$ through:
\begin{align}
\mathrm{Hom}(M_1,B_1)\oplus\mathrm{Hom}(M_2,B_2)\rightarrow \mathrm{Hom}(M_{12},B_{12}).	
\end{align}
This makes the corresponding correct modification to the corresponding splitting $s_1+s_2$ to make sure that $s_1''+s_2''$ lives in the kernel of 
\begin{align}
\mathrm{Hom}(M_1,B_1)\oplus\mathrm{Hom}(M_2,B_2)\rightarrow \mathrm{Hom}(M_{12},B_{12}).	
\end{align}
Then this gives rise to a corresponding splitting for $M$ as right module as well.

\end{proof}

\indent Apply the corresponding general results developed above we have the following result for glueing finite projective Frobenius modules.

\begin{proposition} \label{proposition5.13}
Consider the following two categories. The first category is the corresponding finite projective $\varphi^a$-bibundles over the corresponding period ring $\widetilde{\Pi}_{R,B}$. And the corresponding second finite projective $\varphi^a$-bimodules over the period ring $\widetilde{\Pi}^{[s,r]}_{R,B}$	such that we have $0<s\leq rp^{-ha}$. Then we have that the corresponding involved categories are actually equivalent to each other.
\end{proposition}

\begin{proof}
The functor from the first category to the second one realizing this equivalence is just the corresponding projection. To show that for any bimodule over $\widetilde{\Pi}^{[s,r]}_{R,B}$ we have the corresponding essential surjectivity we lift the corresponding bimodule up to achieve a family of bimodules with the same rank by applying the corresponding Frobenius with respect to different intervals taking the form of $[sp^{-akh},rp^{-akh}]$. Then for general interval we are free to glue finite projective modules by applying \cref{proposition5.11}.\\
\end{proof}

\

This chapter is based on the following paper, where the author of this dissertation is the main author:
\begin{itemize}
\item Tong, Xin. "Period Rings with Big Coefficients and Applications I." arXiv preprint arXiv:2012.07338 (2020). 
\end{itemize}

\newpage\chapter{Period Rings with Big Coefficients and Applications II}

\newpage\section{Introduction}

\subsection{Noncommutative Topological Pseudocoherence after Illusie-Kedlaya-Liu}

\indent In our previous work on noncommutative deformation of period rings and application \cite{8T1}, \cite{8T2} and \cite{8T3}, we considered many useful deformation of the sheaves from many deep constructions coming from Kedlaya-Liu. The corresponding local pictures in both commutative setting and noncommutative setting were carefully developed. \\

\indent The globalization of the picture will happen both in the corresponding coefficients and in the corresponding underlying spaces. Both of them require very deep descent results, although in the commutative setting one do not have to modify much from \cite{8KL1} and \cite{8KL2}. However in the noncommutative setting, things will be very complicated. Since if one would like to really consider the corresponding localization with respect to the noncommutative coefficient Banach algebra, one has to study the corresponding noncommutative toposes. This will be really complicated since the corresponding topological issues are really  hard to manipulate on the level of toposes. That being said, with respect to some of the corresponding application in mind, we could actually do some extension of Kedlaya-Liu glueing to noncommutative setting by keeping track of the noncommutative spatial information in the following sense.\\

\indent Over $\mathbb{Q}_p$ or $\mathbb{F}_p((t))$, in our application in mind we consider the corresponding deformation of period rings in various context and the corresponding period sheaves over certain sites. One could deal with such very big Banach sheaves by considering the corresponding localization concentrated at the corresponding commutative part but release the corresponding deformation components. This will basically be slightly different from Kedlaya-Liu's original extension to Banach context \cite{8KL2} from Illusie's original notion of pseudocoherence in \cite{8SGAVI}. What is happening is that we will have some stability notion which is only sensible to the commutative spatial part. To be more precise when we consider the product of a commutative Banach adic algebra $A$ with a noncommutative Banach algebra $B$, we will consider just the corresponding stability with respect to the corresponding rational localization from $A$ (not at all for $B$). Therefore when we sheafify the modules with such stability we will have a sheaf locally attached to modules with such stability in the same fashion.\\

\indent In some of the applications one might want to consider the corresponding deformations by more general coefficients, for instance in noncommutative Iwasawa theory one might want to consider the corresponding noncommutative Iwasawa algebras. Therefore in such situation one could deform the structure sheaves by some pro-Banach algebras, motivated by Burns-Flach-Fukaya-Kato \cite{8BF1}, \cite{8BF2} and \cite{8FK1}. On the other hand suppose one considers deformations by some limit of Fr\'echet rings such as the corresponding ring $B_e$ in the theory of Bergers $B$-pairs in \cite{8Ber1} and the corresponding various types of Robba rings with respect to some radius $r>0$ or nothing (namely the full Robba rings) in \cite{8KL1} and \cite{8KL2}. Carrying these LF algebras one can consider some types of mixed-type of Hodge-structures along the work of \cite{8CKZ}, \cite{8PZ} and \cite{8Ked1} in more general sense or in the original context.\\

\indent In this paper, we will also consider the following discussion. First we would like to tackle the corresponding integral aspects of the corresponding story, especially when we have the adic spaces which are not Tate but just analytic (since we are dealing with spaces over $\mathbb{Z}_p$) in the sense of \cite{8Ked2}. Certainly we will carry some large coefficients. We will expect that this will have the potential application to integral $p$-adic Hodge theory such as in \cite{8BMS} and \cite{8BS}. That being said, the corresponding descent in this paper will definitely have its own interests. We work in detail out the corresponding descent of the corresponding pseudocoherent sheaves with some stability carrying the corresponding large coefficients in the corresponding uniform analytic Huber pair situation, while we also translated several results to the corresponding uniform analytic adic Banach rings after \cite{8KL2}. \\

\subsection{Main Results}

\indent  The corresponding notions of stably-pseudocoherence is generalized to noncommutative deformed setting from \cite{8KL2}. Over affinoids we realized the following theorems around taking global section of pseudocoherent sheaves. Namely with the corresponding notations in \cref{theorem2.14}:

\begin{theorem}\mbox{\bf{(After Kedlaya-Liu \cite[Theorem 2.5.5]{8KL2})}} \label{8theorem8.1.1} In analytic topology, taking the global section will realize an equivalence between the category of pseudocoherent $\mathcal{O}\widehat{\otimes} B$-sheaves and the category of $B$-stably-pseudocoherent $\mathcal{O}(X)\widehat{\otimes}B$ modules. Here $X$ is an adic affinoid $\mathrm{Spa}(A,A^+)$. 
	
\end{theorem}

\indent In further application, one can also establish the desired parallel results in the \'etale topology. Namely with the corresponding notations in \cref{theorem2.22}:

\begin{theorem}\mbox{\bf{(After Kedlaya-Liu \cite[Theorem 2.5.14]{8KL2})}} \label{8theorem8.1.2} In \'etale topology, taking the global section will realize an equivalence between the category of pseudocoherent $\mathcal{O}_\text{\'et}\widehat{\otimes} B$-sheaves and the category of $B$-\'etale-stably-pseudocoherent $\mathcal{O}_\text{\'et}(X)\widehat{\otimes}B$ modules. Here $X$ is an adic affinoid $\mathrm{Spa}(A,A^+)$.  
	
\end{theorem}

\indent As we mentioned one could also consider more general deformation such as the corresponding. The first scope of the consideration could be made in order to include the context of pro-Banach families associated to some $p$-adic rational Iwasawa algebra with Banach coefficients $B[[G]]$ attached to some $p$-adic Lie group $G$. Namely with the corresponding notations in \cref{theorem3.15}:

\begin{theorem}\mbox{\bf{(After Kedlaya-Liu \cite[Theorem 2.5.5]{8KL2})}} \label{8theorem8.1.3} In analytic topology, taking the global section will realize an equivalence between the category of pseudocoherent $\mathcal{O}\widehat{\otimes} B[[G]]$-sheaves and the category of $B[[G]]$-stably-pseudocoherent $\mathcal{O}(X)\widehat{\otimes}B[[G]]$ modules. Here $X$ is an adic affinoid $\mathrm{Spa}(A,A^+)$, and here all the modules are pro-systems over $\mathbb{Q}_p[[G]]$. 
	
\end{theorem}

\indent In further application, one can also establish the desired parallel results in the \'etale topology. Namely with the corresponding notations in \cref{theorem3.23}:

\begin{theorem}\mbox{\bf{(After Kedlaya-Liu \cite[Theorem 2.5.14]{8KL2})}} \label{8theorem8.1.4} In \'etale topology, taking the global section will realize an equivalence between the category of pseudocoherent $\mathcal{O}_\text{\'et}\widehat{\otimes} B[[G]]$-sheaves and the category of $B[[G]]$-\'etale-stably-pseudocoherent $\mathcal{O}_\text{\'et}(X)\widehat{\otimes}B[[G]]$ modules. Here $X$ is an adic affinoid $\mathrm{Spa}(A,A^+)$, and here all the modules are pro-systems over $\mathbb{Q}_p[[G]]$.	
\end{theorem}

\indent One can also even consider the limit of Fr\'echet coefficients with some further application to some equivariant sheaves over adic Fargues-Fontaine curves. Namely with the corresponding notations in \cref{theorem4.14}:

\begin{theorem}\mbox{\bf{(After Kedlaya-Liu \cite[Theorem 2.5.5]{8KL2})}} \label{8theorem8.1.5} In analytic topology, taking the global section will realize an equivalence between the category of pseudocoherent $\mathcal{O}\widehat{\otimes} B$-sheaves and the category of $B$-stably-pseudocoherent $\mathcal{O}(X)\widehat{\otimes}B$ modules. Here $X$ is an adic affinoid $\mathrm{Spa}(A,A^+)$, and here $B$ is injective limit of Banach algebras $\varinjlim_h B_h$. 
	
\end{theorem}

\indent In further application, one can also establish the desired parallel results in the \'etale topology. Namely with the corresponding notations in \cref{theorem4.22}:

\begin{theorem}\mbox{\bf{(After Kedlaya-Liu \cite[Theorem 2.5.14]{8KL2})}} \label{8theorem8.1.6} In \'etale topology, taking the global section will realize an equivalence between the category of pseudocoherent $\mathcal{O}_\text{\'et}\widehat{\otimes} B$-sheaves and the category of $B$-\'etale-stably-pseudocoherent $\mathcal{O}_\text{\'et}(X)\widehat{\otimes}B$ modules. Here $X$ is an adic affinoid $\mathrm{Spa}(A,A^+)$, and here $B$ is injective limit of Banach algebras $\varinjlim_h B_h$.  
	
\end{theorem}

\indent One can also consider the corresponding families of affinoids in some coherent way for the commutative parts, namely this is the notion of quasi-Stein spaces from \cite[Chapter 2.6]{8KL2}. Then with the corresponding notations in \cref{theorem2.31}, \cref{theorem3.32} and \cref{theorem4.31}:

\begin{theorem}\mbox{\bf{(After Kedlaya-Liu \cite[Proposition 2.6.17]{8KL2})}} \label{8theorem8.1.7} Carrying Banach, pro-Banach or limit of Fr\'echet coefficient $B$ one could have the finiteness of the global sections of $B$-stably-pseudocoherent sheaves over some quasi-Stein space $X$ as long as one has $m$-uniform covering as in \cite[Proposition 2.6.17]{8KL2}.  
	
\end{theorem}

\indent The application will be encoded in our last chapter, where we focused on finding relationship between the corresponding Frobenius equivariant sheaves over Fargues-Fontaine space:

\begin{align}
\bigcup \mathrm{Spa}(\widetilde{\Pi}_R^{[r_1,r_2]},\widetilde{\Pi}_R^{[r_1,r_2],+})	
\end{align}
and the corresponding Frobenius equivariant modules over the global sections:
\begin{align}
\widetilde{\Pi}_R^{[r_1,r_2]},\forall 0<s<r.	
\end{align}
in glueing fashion. The sheaves and modules here carrying some big Banach coefficients (noncommutative). And we also considered some very interesting pro-Banach situation and LF situation, see \cref{theorem5.18} and \cref{theorem5.24}.\\

\subsection{Further Consideration}

\indent One thing we have observed (although have not written so) is that suppose we apply our results to some locally noetherian spaces, then we will have no stability issue, and we will have not just pseudoflatness but even more. Instead, this will immediately mean that we are going to be more algebraic throughout. The corresponding story is very interesting in the situation of rigid analytic spaces.\\ 

\indent The discussion in the previous paragraph in fact implies that in the noetherian situation, one should be able to extend the discussion around $\infty$-glueing along the style of Kedlaya-Liu by applying Bambozzi-Kremnizer spectrum \cite{8BK1} and Clausen-Scholze space \cite{8CS} to the corresponding coherent sheaves which are interesting to us where the corresponding highly nontrivial stability and nonflatness issues in \cite{8KL2} disappear. Again note that we are still carrying some large coefficients.\\

%%\newpage 

\newpage\section{Foundations on Noncommutative Descent for Adic Spectra in Banach Case}

\subsection{Noncommutative Pseudocoherence in Analytic Topology}

\indent We now establish some foundations in the noncommutative setting on glueing noncommutative pseudocoherent modules after \cite[Chapter 2]{8KL2}. But we remind the readers that we will fix the base space:

\begin{setting} \label{setting2.1}
Consider a corresponding sheafy Banach adic uniform algebra $(A,A^+)$ over $\mathbb{Q}_p$ or $\mathbb{F}_p((t))$, we consider the base space $\mathrm{Spa}(A,A^+)$ as in \cite[Chapter 2]{8KL2}. And now we will consider a further noncommutative Banach algebra $(B,B^+)$ over $\mathbb{Q}_p$ or $\mathbb{F}_p((t))$. 	
\end{setting}

\begin{definition} \mbox{\bf{(After Kedlaya-Liu \cite[Definition 2.4.1]{8KL2})}}
For any left $A\widehat{\otimes}B$-module $M$, we call it $m$-$B$-stably pseudocoherent if we have that it is $m$-$B$-pseudocoherent, complete for the natural topology and for any morphism $A\rightarrow A'$ which is the corresponding rational localization, the base change of $M$ to $A'\widehat{\otimes} B$ is complete for the natural topology as the corresponding left $A'\widehat{\otimes} B$-module.	As in \cite[Definition 2.4.1]{8KL2} we call that the corresponding left $A\widehat{\otimes}B$-module $M$ just $B$-stably pseudocoherent if we have that it is simply just $\infty$-$B$-stably pseudocoherent.
\end{definition}

\begin{remark}
Throughout the whole paper, it is obvious we need to use the corresponding noncommutative version of the corresponding pseudocoherent modules on the algebraic level due to Illusie along the development of Grothendieck's SGAVI \cite{8SGAVI}, but the main difference is just that we consider the corresponding modules with action happening on one side.
\end{remark}

\indent Since our spaces are commutative, so the rational localization is not flat which is exactly the same situation as in \cite{8KL1} and \cite{8KL2}, therefore in our situation we need to define something weaker than the corresponding flatness which will be sufficient for us to apply in the following development: 

\begin{definition}\mbox{\bf{(After Kedlaya-Liu \cite[Definition 2.4.4]{8KL2})}}
For any Banach left $A\widehat{\otimes}B$-module $M$, we call it is $m$-$B$-pseudoflat if for any right $A\widehat{\otimes}B$-module $M'$ $m$-$B$-stably pseudocoherent we have $\mathrm{Tor}_1^{A\widehat{\otimes}B}(M',M)=0$. For any Banach right $A\widehat{\otimes}B$-module $M$, we call it is $m$-$B$-pseudoflat if for any left $A\widehat{\otimes}B$-module $M'$ $m$-$B$-stably pseudocoherent we have $\mathrm{Tor}_1^{A\widehat{\otimes}B}(M,M')=0$.
\end{definition}

\begin{definition}\mbox{\bf{(After Kedlaya-Liu \cite[Definition 2.4.6]{8KL2})}} We consider the corresponding notion of the corresponding pro-projective module. We define over $A$ the corresponding $B$-pro-projective module $M$ to be a corresponding left $A\widehat{\otimes}B$-module such that one could find a filtered sequence of projectors such that in the sense of taking the prolimit we have the corresponding projector will converge to any chosen element in $M$. Here we assume the module is complete with respect to natural topology, and we assume the projectors are $A\widehat{\otimes}B$-linear and we assume that the corresponding image of the projectors are modules which are also finitely generated and projective.
	
\end{definition}

\begin{lemma}\mbox{\bf{(After Kedlaya-Liu \cite[Lemma 2.4.7]{8KL2})}}
Suppose we have over the space $\mathrm{Spa}(A,A^+)$ a corresponding $B$-pro-projective left module $M$. And suppose that we have a $2$-$B$-pseudocoherent right module $C$ over $A$. And we assume that $C$ is complete with respect to the natural topology. Then we have that the corresponding product $C{\otimes}_{A\widehat{\otimes}B}M$ is then complete under the corresponding natural topology. And moreover we have that in our situation:
\begin{displaymath}
\mathrm{Tor}_1^{A\widehat{\otimes}B}(C,M)=0.	
\end{displaymath}

\end{lemma}

\begin{proof}
This will be just a noncommutative version of the corresponding  \cite[Lemma 2.4.7]{8KL2}. We adapt the corresponding argument in \cite[Lemma 2.4.7]{8KL2} to our situation. What we are going to consider in this case is then first consider the following presentation:
\[
\xymatrix@C+0pc@R+0pc{
M'  \ar[r]\ar[r]\ar[r] &(A\widehat{\otimes}B)^k \ar[r]\ar[r]\ar[r] & M\ar[r]\ar[r]\ar[r] &0,
}
\]
where the left module $M'$ is finitely presented. Then we consider the following commutative diagram:
\[
\xymatrix@C+0pc@R+3pc{
&C\otimes_{A\widehat{\otimes}B}M'  \ar[r]\ar[r]\ar[r]  \ar[d]\ar[d]\ar[d] &C\otimes_{A\widehat{\otimes}B}(A\widehat{\otimes}B)^k \ar[r]\ar[r]\ar[r] \ar[d]\ar[d]\ar[d] & C\otimes_{A\widehat{\otimes}B}M\ar[r]\ar[r]\ar[r] \ar[d]\ar[d]\ar[d] &0,\\
0 \ar[r]\ar[r]\ar[r] &C\widehat{\otimes}_{A\widehat{\otimes}B}M'  \ar[r]\ar[r]\ar[r] &C\widehat{\otimes}_{A\widehat{\otimes}B}(A\widehat{\otimes}B)^k \ar[r]\ar[r]\ar[r] & C\widehat{\otimes}_{A\widehat{\otimes}B}M\ar[r]\ar[r]\ar[r] &0.
}
\]\\
The first row is exact as in \cite[Lemma 2.4.7]{8KL2}, and we have the corresponding second row is also exact since we have that by hypothesis the corresponding module $M$ is $B$-pro-projective.  
Then as in \cite[Lemma 2.4.7]{8KL2} the results follow.	
\end{proof}

\begin{proposition} \mbox{\bf{(After Kedlaya-Liu \cite[Corollary 2.4.8]{8KL2})}}
Over the space\\ $\mathrm{Spa}(A,A^+)$, we have any $B$-pro-projective left module is $2$-$B$-pseudoflat.
\end{proposition}

\begin{proof}
This is a direct consequence of the previous lemma.	
\end{proof}

\begin{lemma} \mbox{\bf{(After Kedlaya-Liu \cite[Corollary 2.4.9]{8KL2})}}
Keep the notation above, suppose we are working over the adic space $
\mathrm{Spa}(A,A^+)$. Suppose now the left $A\widehat{\otimes}B$-module is finitely generated, then we have that the following natural maps:
\begin{align}
A\{T\}\widehat{\otimes}B\otimes_{A\widehat{\otimes}B} M	\rightarrow M\{T\},\\\
A\{T,T^{-1}\}\widehat{\otimes}B \otimes_{A\widehat{\otimes}B} M 	\rightarrow M\{T,T^{-1}\}	
\end{align}
are surjective. Suppose now the left $A\widehat{\otimes}B$-module is finitely presented, then we have that the following natural maps:
\begin{align}
A\{T\}\widehat{\otimes}B\otimes_{A\widehat{\otimes}B} M	\rightarrow M\{T\},\\\
A\{T,T^{-1}\}\widehat{\otimes}B \otimes_{A\widehat{\otimes}B} M 	\rightarrow M\{T,T^{-1}\}	
\end{align}
are bijective. Here all the modules are assumed to be compete with respect to the natural topology in our context.

\end{lemma}

\begin{proof}
This is the corresponding consequence of the previous lemma.
\end{proof}

\

\begin{proposition} \mbox{\bf{(After Kedlaya-Liu \cite[Lemma 2.4.12]{8KL2})}}
We keep the corresponding notations in \cite[Lemma 2.4.10]{8KL2}. Then in our current sense we have that the corresponding morphism $A\rightarrow B_2$ is then $2$-$B$-pseudoflat. 
\end{proposition}

\begin{proof}
The corresponding proof could be made parallel to the corresponding proof of \cite[Lemma 2.4.12]{8KL2}. We need to keep track the corresponding coefficient $B$ in our current context. To be more precise we consider the following commutative diagram from the corresponding chosen exact sequence (which is just the analog of the corresponding one in \cite[Lemma 2.4.12]{8KL2}, namely as below we choose arbitrary short exact sequence with $P$ as a left $A\widehat{\otimes}B$-module which is $2$-$B$-stably pseudocoherent):
\[
\xymatrix@C+0pc@R+0pc{
0   \ar[r]\ar[r]\ar[r] &M \ar[r]\ar[r]\ar[r] &N \ar[r]\ar[r]\ar[r] & P \ar[r]\ar[r]\ar[r] &0,
}
\]
which then induces the following corresponding big commutative diagram:
\[
\xymatrix@C+0pc@R+4pc{
& &0 \ar[d]\ar[d]\ar[d] &0 \ar[d]\ar[d]\ar[d] &,\\
0   \ar[r]\ar[r]\ar[r]  &{A\{T\}\widehat{\otimes}B}\otimes_{A\widehat{\otimes}B} M \ar[r]\ar[r]\ar[r] \ar[d]^{1-fT}\ar[d]\ar[d] &{A\{T\}\widehat{\otimes}B}\otimes_{A\widehat{\otimes}B} N \ar[r]\ar[r]\ar[r] \ar[d]^{1-fT}\ar[d]\ar[d] &{A\{T\}\widehat{\otimes}B} \otimes_{A\widehat{\otimes}B}P \ar[r]\ar[r]\ar[r] \ar[d]^{1-fT}\ar[d]\ar[d] &0,\\
0   \ar[r]\ar[r]\ar[r] &{A\{T\}\widehat{\otimes}B}\otimes_{A\widehat{\otimes}B}M  \ar[r]\ar[r]\ar[r] \ar[d]\ar[d]\ar[d] &{A\{T\}\widehat{\otimes}B} \otimes_{A\widehat{\otimes}B}N\ar[r]\ar[r]\ar[r] \ar[d]\ar[d]\ar[d] & {A\{T\}\widehat{\otimes}B}\otimes_{A\widehat{\otimes}B}P \ar[r]\ar[r]\ar[r] \ar[d]\ar[d]\ar[d] &0,\\
0  \ar[r]^?\ar[r]\ar[r] &{B_2\widehat{\otimes}B}\otimes_{A\widehat{\otimes}B} M \ar[r]\ar[r]\ar[r] \ar[d]\ar[d]\ar[d] &{B_2\widehat{\otimes}B}\otimes_{A\widehat{\otimes}B}N \ar[r]\ar[r]\ar[r] \ar[d]\ar[d]\ar[d] & {B_2\widehat{\otimes}B}\otimes_{A\widehat{\otimes}B} P \ar[r]\ar[r]\ar[r] \ar[d]\ar[d]\ar[d] &0,\\
&0&0&0
}
\]
with the notations in \cite[Lemma 2.4.10]{8KL2} where we actually have the corresponding exactness along the horizontal direction at the corner around ${B_2\widehat{\otimes}B}\otimes_{A\widehat{\otimes}B} M $ marked with $?$, by diagram chasing.
\end{proof}

\begin{proposition} \mbox{\bf{(After Kedlaya-Liu \cite[Lemma 2.4.13]{8KL2})}} \label{proposition2.9}
We keep the corresponding notations in \cite[Lemma 2.4.10]{8KL2}. Then in our current sense we have that the corresponding morphism $A\rightarrow B_1$ is then $2$-$B$-pseudoflat. And in our current sense we have that the corresponding morphism $A\rightarrow B_{12}$ is then $2$-$B$-pseudoflat. 
\end{proposition}

\begin{proof}
The statement for $A\rightarrow B_{12}$ could be proved by consider the composition $A\rightarrow B_2\rightarrow B_{12}$ where the $2$-$B$-pseudoflatness could be proved as in \cite[Lemma 2.4.13]{8KL2} by inverting the corresponding variables. For $A\rightarrow B_{1}$,
the corresponding proof could be made parallel to the corresponding proof of \cite[Lemma 2.4.13]{8KL2}. To be more precise we consider the following commutative diagram from the corresponding chosen exact sequence (which is just the analog of the corresponding one in \cite[Remark 2.4.5]{8KL2}, namely as below we choose arbitrary short exact sequence with $P$ as a left $A\widehat{\otimes}B$-module which is $2$-$B$-stably pseudocoherent):
\[
\xymatrix@C+0pc@R+0pc{
0   \ar[r]\ar[r]\ar[r] &M \ar[r]\ar[r]\ar[r] &N \ar[r]\ar[r]\ar[r] & P \ar[r]\ar[r]\ar[r] &0,
}
\]
which then induces the following corresponding big commutative diagram:
\[\tiny
\xymatrix@C+0pc@R+5pc{
& &0 \ar[d]\ar[d]\ar[d] &0 \ar[d]\ar[d]\ar[d] &,\\
0   \ar[r]\ar[r]\ar[r]  &M \ar[r]\ar[r]\ar[r] \ar[d]\ar[d]\ar[d] &N \ar[r]\ar[r]\ar[r] \ar[d]\ar[d]\ar[d] & P \ar[r]\ar[r]\ar[r] \ar[d]\ar[d]\ar[d] &0,\\
0   \ar[r]^?\ar[r]\ar[r] &({B_1\widehat{\otimes}B}\bigoplus {B_2\widehat{\otimes}B}) \otimes_{A\widehat{\otimes}B}M  \ar[r]\ar[r]\ar[r] \ar[d]\ar[d]\ar[d] &({B_1\widehat{\otimes}B}\bigoplus {B_2\widehat{\otimes}B})\otimes_{A\widehat{\otimes}B}N  \ar[r]\ar[r]\ar[r] \ar[d]\ar[d]\ar[d] &({B_1\widehat{\otimes}B}\bigoplus {B_2\widehat{\otimes}B})\otimes_{A\widehat{\otimes}B}P\ar[r]\ar[r]\ar[r] \ar[d]\ar[d]\ar[d] &0,\\
0  \ar[r]\ar[r]\ar[r] &{B_{12}\widehat{\otimes}B}\otimes_{A\widehat{\otimes}B}M \ar[r]\ar[r]\ar[r] \ar[d]\ar[d]\ar[d] &{B_{12}\widehat{\otimes}B}\otimes_{A\widehat{\otimes}B}N \ar[r]\ar[r]\ar[r] \ar[d]\ar[d]\ar[d] &{B_{12}\widehat{\otimes}B} \otimes_{A\widehat{\otimes}B} P \ar[r]\ar[r]\ar[r] \ar[d]\ar[d]\ar[d] &0,\\
&0&0&0
}
\]
with the notations in \cite[Lemma 2.4.10]{8KL2} where we actually have the corresponding exactness along the horizontal direction at the corner around $M\otimes_{A\widehat{\otimes}B} ({B_1\widehat{\otimes}B}\bigoplus {B_2\widehat{\otimes}B})$ marked with $?$, by diagram chasing.
\end{proof}

\indent After these foundational results, as in \cite{8KL2} we have the following proposition which is the corresponding noncommutative generalization of the corresponding result established in \cite[Theorem 2.4.15]{8KL2}.

\begin{proposition} \mbox{\bf{(After Kedlaya-Liu \cite[Theorem 2.4.15]{8KL2})}} \label{8proposition2.10}
In our current context we have that for any rational localization $\mathrm{Spa}(A',A^{',+})\rightarrow \mathrm{Spa}(A,A^+)$ we have that along this base change the corresponding $B$-stably pseudocoherence is preserved.\\
	
\end{proposition}

\indent Then we have the corresponding Tate's acyclicity in the noncommutative deformed setting:

\begin{theorem}\mbox{\bf{(After Kedlaya-Liu \cite[Theorem 2.5.1]{8KL2})}} \label{theorem2.11} Now suppose we have in our corresponding \cref{setting2.1} a corresponding $B$-stably pseudocoherent module $M$. Then we consider the corresponding assignment such that for any $U\subset \mathrm{Spa}(A,A^+)$ we define $\widetilde{M}(U)$ as in the following:
\begin{align}
\widetilde{M}(U):=\varprojlim_{\mathrm{Spa}(S,S^+)\subset U,\mathrm{rational}} S\widehat{\otimes}B\otimes_{A\widehat{\otimes}B}M.	
\end{align}
Then we have that for any $\mathfrak{B}$ which is a rational covering of $U=\mathrm{Spa}(S,S^+)\subset \mathrm{Spa}(A,A^+)$ (certainly this $U$ is also assumed to be rational) we have that the vanishing of the following two cohomology groups:
\begin{align}
H^i(U,\widetilde{M}), \check{H}^i(U:\mathfrak{B},\widetilde{M})
\end{align}
for any $i>0$. When concentrating at the degree zero we have:
\begin{align}
H^0(U,\widetilde{M})=S\widehat{\otimes}B\otimes_{A\widehat{\otimes}B}M, \check{H}^0(U:\mathfrak{B},\widetilde{M})=S\widehat{\otimes}B\otimes_{A\widehat{\otimes}B}M.
\end{align}
	
\end{theorem}

\begin{proof}
By \cite[Propositions 2.4.20-2.4.21]{8KL1}, we could then finish proof as in \cite[Theorem 2.5.1]{8KL2} by using our previous \cref{8proposition2.10} and \cref{proposition2.9} as above.	
\end{proof}

\indent We now consider the corresponding noncommutative deformed version of Kiehl's glueing properties for stably pseudocoherent modules after \cite{8KL2}:

\begin{definition} \mbox{\bf{(After Kedlaya-Liu \cite[Definition 2.5.3]{8KL2})}}
Consider in our context (over $\mathrm{Spa}(A,A^+)$) the corresponding sheaves $\mathcal{O}\widehat{\otimes}B$, we will then define the corresponding pseudocoherent sheaves over $\mathrm{Spa}(A,A^+)$ to be those locally defined by attaching stably-pseudocoherent modules over the section. 
\end{definition}

\begin{lemma} \mbox{\bf{(After Kedlaya-Liu \cite[Lemma 2.5.4]{8KL2})}}\label{8lemma2.13}
	Consider the corresponding notations in \cite[Lemma 2.4.10]{8KL2}, we have the corresponding morphism $A\rightarrow B_1\bigoplus B_2$. Then we have that this morphism is an descent morphism effective for the corresponding $B$-stably pseodocoherent Banach modules. 
\end{lemma}

\begin{proof}
This is a $B$-relative version of the \cite[Lemma 2.5.4]{8KL2}. We adapt the corresponding proof to our situation. First consider the corresponding notations on the commutative adic space $\mathrm{Spa}(A,A^+)$ as in \cite[Lemma 2.4.10]{8KL2}. And we have a coving $\{U_1,U_2\}$. Then consider a descent datum for this cover of $B$-stably pseudocoherent sheaves. By \cref{theorem2.11} we could realize this as a corresponding pseudocoherent sheaf over $\mathcal{O}\widehat{\otimes}B$. We denote this by $\mathcal{M}$. Then by (the proof of) \cite[Lemma 6.82]{8T2}, we have that this is finitely generated and we have the corresponding vanishing of $\check{H}^i(\mathrm{Spa}(A,A^+),\mathfrak{B};\mathcal{M})$ for $i>0$. Then by \cref{theorem2.11} again we have that the vanishing of $H^i(U_j;\mathcal{M})$ for $j\in \{1,2\},i>0$. Then as in \cite[Lemma 2.5.4]{8KL2} we have the corresponding vanishing of $H^i(\mathrm{Spa}(A,A^+);\mathcal{M})$ for $i>0$. Then as in \cite[Lemma 2.5.4]{8KL2} we choose a corresponding covering of some finite free object $\mathcal{F}$ and take its kernel $\mathcal{K}$ which is actually pseudocoherent by the corresponding \cref{8proposition2.10}. Then forming the following diagram:
\[\tiny
\xymatrix@C+0pc@R+6pc{
0 \ar[r]\ar[r]\ar[r] & (B_p\widehat{\otimes}B)\otimes \mathcal{K}(\mathrm{Spa}(A,A^+))  \ar[r]\ar[r]\ar[r] \ar[d]\ar[d]\ar[d] &(B_p\widehat{\otimes}B)\otimes \mathcal{F}(\mathrm{Spa}(A,A^+))\ar[r]\ar[r]\ar[r] \ar[d]\ar[d]\ar[d] & (B_p\widehat{\otimes}B)\otimes \mathcal{M}(\mathrm{Spa}(A,A^+))\ar[r]\ar[r]\ar[r]  \ar[d]\ar[d]\ar[d]&0,\\
0 \ar[r]\ar[r]\ar[r] &\mathcal{K}(U_p)   \ar[r]\ar[r]\ar[r] &\mathcal{F}(U_p) \ar[r]\ar[r]\ar[r] & \mathcal{M}(U_p) \ar[r]\ar[r]\ar[r] &0,\\
}
\]
where we have $p=1,2$. The first row is exact by $B$-pseudoflatness as in \cref{8proposition2.10} and the second is row is exact as well by \cref{theorem2.11}. The middle vertical arrow is an isomorphism while the rightmost vertical arrow is surjective and the leftmost vertical arrow is surjective. The surjectivity just here comes from actually \cite[Lemma 6.82]{8T2}. This will imply that the rightmost vertical arrow is injective as well. Then we have the corresponding isomorphism $\mathcal{M}\overset{\sim}{\rightarrow} \widetilde{\mathcal{M}(\mathrm{Spa}(A,A^+))}$. Then we consider the following diagram:

\[ \tiny
\xymatrix@C+0pc@R+6pc{
 &(C\widehat{\otimes}B)\otimes\mathcal{K}(\mathrm{Spa}(A,A^+))   \ar[r]\ar[r]\ar[r] \ar[d]\ar[d]\ar[d] &(C\widehat{\otimes}B)\otimes\mathcal{F}(\mathrm{Spa}(A,A^+)) \ar[r]\ar[r]\ar[r] \ar[d]\ar[d]\ar[d] & (C\widehat{\otimes}B)\otimes\mathcal{M}(\mathrm{Spa}(A,A^+))\ar[r]\ar[r]\ar[r]  \ar[d]\ar[d]\ar[d]&0,\\
0 \ar[r]\ar[r]\ar[r] &\mathcal{K}(\mathrm{Spa}(C,C^+))   \ar[r]\ar[r]\ar[r] &\mathcal{F}(\mathrm{Spa}(C,C^+)) \ar[r]\ar[r]\ar[r] & \mathcal{M}(\mathrm{Spa}(C,C^+)) \ar[r]\ar[r]\ar[r] &0,\\
}
\]
$C$ is any rational localization of $A$. Here the first row is exact by $B$-pseudoflatness and the second is row is exact as well. The middle vertical arrow is an isomorphism while the rightmost vertical arrow is surjective as in \cite[Lemma 2.5.4]{8KL2}. Then we apply the same argument to the leftmost vertical arrow (with a new diagram), we will have the corresponding surjectivity of this leftmost vertical arrow. This will imply that the rightmost vertical arrow is injective as well. Then we have the corresponding isomorphism for the rightmost vertical arrow.
\end{proof}

\begin{theorem}\mbox{\bf{(After Kedlaya-Liu \cite[Theorem 2.5.5]{8KL2})}} \label{theorem2.14}
Taking global section will realize the equivalence between the following two categories: A. The category of all the pseudocoherent $\mathcal{O}\widehat{\otimes}B$-sheaves; B. The category of all the $B$-stably pseudocoherent modules over $A\widehat{\otimes}B$. 	
\end{theorem}

\begin{proof}
See \cite[Theorem 2.5.5]{8KL2}. We need to still apply \cite[Proposition 2.4.20]{8KL1}, as long as one considers instead in our situation \cref{8lemma2.13} and \cref{theorem2.11}.\\	
\end{proof}

\subsection{Noncommutative Pseudocoherence in \'Etale Topology}

\indent We consider the extension of the corresponding discussion in the previous subsection to the \'etale topology. At this moment we consider the following same context on the geometric level as in \cite{8KL2}:

\begin{setting} \mbox{\bf{(After Kedlaya-Liu \cite[Hypothesis 2.5.8]{8KL2})}} \label{setting2.15}
As in the previous subsection we consider now the corresponding geometric setting namely an adic space $\mathrm{Spa}(A,A^+)$ where $A$ is assumed to be sheafy. Then we consider $B$ as above. And we consider the corresponding \'etale site $\mathrm{Spa}(A,A^+)_\text{\'et}$. And by keeping the corresponding setting in the geometry as in \cite[Hypothesis 2.5.8]{8KL2}, we assume that there is a stable basis $\mathfrak{B}$ containing the space $\mathrm{Spa}(A,A^+)$ itself.	
\end{setting}

\begin{definition} \mbox{\bf{(After Kedlaya-Liu \cite[Definition 2.5.9]{8KL2})}}
For any left $A\widehat{\otimes}B$-module $M$, we call it $m$-$B$-\'etale stably pseudocoherent with respect to $\mathfrak{B}$ if we have that it is $m$-$B$-pseudocoherent, complete for the natural topology and for any rational localization $A\rightarrow A'$, the base change of $M$ to $A'\widehat{\otimes} B$ is complete for the natural topology as the corresponding left $A'\widehat{\otimes} B$-module.	As in \cite[Definition 2.4.1]{8KL2} we call that the corresponding left $A\widehat{\otimes}B$-module $M$ just $B$-stably pseudocoherent if we have that it is simply just $\infty$-$B$-stably pseudocoherent.
\end{definition}

\begin{definition}\mbox{\bf{(After Kedlaya-Liu \cite[Below Definition 2.5.9]{8KL2})}}
For any Banach left $A\widehat{\otimes}B$-module $M$, we call it $m$-$B$-\'etale-pseudoflat if for any right $A\widehat{\otimes}B$-module $M'$ $m$-$B$-\'etale-stably pseudocoherent we have $\mathrm{Tor}_1^{A\widehat{\otimes}B}(M',M)=0$. For any Banach right $A\widehat{\otimes}B$-module $M$, we call it is $m$-$B$-\'etale-pseudoflat if for any left $A\widehat{\otimes}B$-module $M'$ $m$-$B$-\'etale-stably pseudocoherent we have $\mathrm{Tor}_1^{A\widehat{\otimes}B}(M,M')=0$.
\end{definition}

\indent The following proposition holds in our current setting.

\begin{proposition} \mbox{\bf{(After Kedlaya-Liu \cite[Lemma 2.5.10]{8KL2})}} \label{8proposition2.18}
One can actually find another basis $\mathfrak{C}$ in $\mathfrak{B}$ such that any morphism in $\mathfrak{C}$ could be $2$-$B$-\'etale-pseudoflat with respect to either $\mathfrak{C}$ or $\mathfrak{B}$.	
\end{proposition}

\begin{proof}
We follow the proof of \cite[Lemma 2.5.10]{8KL2} in our current $B$-relative situation, however not that much needs to be adapted by relying on the proof of \cite[Lemma 2.5.10]{8KL2}. To be more precise the corresponding selection of the new basis $\mathfrak{C}$ comes from including all the morphisms made up of some composition of rational localization and the corresponding finite \'etale ones. For these we have shown above the corresponding 2-$B$-\'etale pseudoflatness. Then the rest will be pure geometric for the analytic spaces which are just the same as the situation of \cite[Lemma 2.5.10]{8KL2}, therefore we just omit the corresponding argument, see \cite[Lemma 2.5.10]{8KL2}.	
\end{proof}

\indent Then we have the corresponding Tate's acyclicity in the noncommutative deformed setting in \'etale topology (here fix $\mathfrak{C}$ as above):

\begin{theorem}\mbox{\bf{(After Kedlaya-Liu \cite[Theorem 2.5.11]{8KL2})}} \label{theorem2.19} Now suppose we have in our corresponding \cref{setting2.1} a corresponding $B$-\'etale-stably pseudocoherent module $M$. Then we consider the corresponding assignment such that for any $U\subset \mathrm{Spa}(A,A^+)$ we define $\widetilde{M}(U)$ as in the following:
\begin{align}
\widetilde{M}(U):=\varprojlim_{\mathrm{Spa}(S,S^+)\subset U,\in \mathfrak{C}} S\widehat{\otimes}B\otimes_{A\widehat{\otimes}B}M.	
\end{align}
Then we have that for any $\mathfrak{B}$ which is a covering of $U=\mathrm{Spa}(S,S^+)\subset \mathrm{Spa}(A,A^+)$ (certainly this $U$ is also assumed to be in $\mathfrak{C}$, and we assume this covering is formed by using the corresponding members in $\mathfrak{C}$) we have that the vanishing of the following two cohomology groups:
\begin{align}
H^i(U,\widetilde{M}), \check{H}^i(U,\mathfrak{B},\widetilde{M})
\end{align}
for any $i>0$. When concentrating at the degree zero we have:
\begin{align}
H^0(U,\widetilde{M})=S\widehat{\otimes}B\otimes_{A\widehat{\otimes}B}M, \check{H}^0(U,\mathfrak{B},\widetilde{M})=S\widehat{\otimes}B\otimes_{A\widehat{\otimes}B}M.
\end{align}
	
\end{theorem}

\begin{proof}
By \cite[Proposition 8.2.21]{8KL1}, we could then finish proof as in \cite[Theorem 2.5.11]{8KL2} by using our previous \cref{8proposition2.18} and the corresponding faithfully flat descent as in the situation of \cite[Theorem 2.5.11]{8KL2}.	
\end{proof}

\indent We now consider the corresponding noncommutative deformed version of Kiehl's glueing property for stably pseudocoherent modules after \cite[Definition 2.5.12]{8KL2}:

\begin{definition} \mbox{\bf{(After Kedlaya-Liu \cite[Definition 2.5.12]{8KL2})}}
Consider in our context (over $\mathrm{Spa}(A,A^+)_{\text{\'et}}$) the corresponding sheaves $\mathcal{O}_{\mathrm{Spa}(A,A^+)_{\text{\'et}}}\widehat{\otimes}B$, we will then define the corresponding pseudocoherent sheaves over $\mathrm{Spa}(A,A^+)_{\text{\'et}}$ to be those locally defined by attaching \'etale-stably-pseudocoherent modules over the section. 
\end{definition}

\begin{lemma} \mbox{\bf{(After Kedlaya-Liu \cite[Lemma 2.5.13]{8KL2})}} \label{lemma2.21}
	Consider the corresponding notations in \cite[Lemma 2.4.10]{8KL2}, we have the corresponding morphism $A\rightarrow B_1\bigoplus B_2$. Then we have that this morphism is an descent morphism effective for the corresponding $B$-\'etale-stably pseodocoherent Banach modules. 
\end{lemma}

\begin{proof}
This is a $B$-relative version of the \cite[Lemma 2.5.13]{8KL2}. We adapt the corresponding proof to our situation. First consider the corresponding notations on the commutative adic space $\mathrm{Spa}(A,A^+)$ as in \cite[Lemma 2.4.10]{8KL2}. And we have a coving $\{U_1,U_2\}$. Then consider a descent datum for this cover of $B$-\'etale-stably pseudocoherent sheaves. By \cref{theorem2.19} we could realize this as a corresponding pseudocoherent sheaf over $\mathcal{O}_{\text{\'et}}\widehat{\otimes}B$. We denote this by $\mathcal{M}$. Then by consider the corresponding results in the rational localization situation above we have that the vanishing of $H^i(U_j;\mathcal{M})$ for $j\in \{1,2\},i>0$. And note that we have the corresponding equality on the degree zero namely $H^i_\text{\'et}(U_j;\mathcal{M})=H^i(U_j;\mathcal{M})$ for just $i$ which is zero. Then as in \cite[Lemma 2.5.13]{8KL2} we have the corresponding vanishing of $H^i(\mathrm{Spa}(A,A^+);\mathcal{M})$ for $i>0$. Namely this is due to the fact that the covering is already one in the analytic topology. Then as in \cite[Lemma 2.5.13]{8KL2} we choose a corresponding covering of some finite free object $\mathcal{F}$ and take its kernel $\mathcal{K}$ which gives rise to the corresponding exact sequence:
\[\tiny
\xymatrix@C+0pc@R+6pc{
0 \ar[r]\ar[r]\ar[r] &\mathcal{K}(\mathrm{Spa}(A,A^+))   \ar[r]\ar[r]\ar[r] &\mathcal{F}(\mathrm{Spa}(A,A^+)) \ar[r]\ar[r]\ar[r] & \mathcal{M}(\mathrm{Spa}(A,A^+)) \ar[r]\ar[r]\ar[r] &0,
}
\]
where we have $p=1,2$. We have already the corresponding stably-pseudocoherence of the global section. As in \cite[Lemma 2.5.13]{8KL2} we only need to show that it is stable under \'etale morphism. Then we consider the following diagram:

\[ \tiny
\xymatrix@C+0pc@R+6pc{
 &(C\widehat{\otimes}B)\otimes\mathcal{K}(\mathrm{Spa}(A,A^+))   \ar[r]\ar[r]\ar[r] \ar[d]\ar[d]\ar[d] &(C\widehat{\otimes}B)\otimes\mathcal{F}(\mathrm{Spa}(A,A^+)) \ar[r]\ar[r]\ar[r] \ar[d]\ar[d]\ar[d] & (C\widehat{\otimes}B)\otimes\mathcal{M}(\mathrm{Spa}(A,A^+))\ar[r]\ar[r]\ar[r]  \ar[d]\ar[d]\ar[d]&0,\\
0 \ar[r]\ar[r]\ar[r] &\mathcal{K}(\mathrm{Spa}(C,C^+))   \ar[r]\ar[r]\ar[r] &\mathcal{F}(\mathrm{Spa}(C,C^+)) \ar[r]\ar[r]\ar[r] & \mathcal{M}(\mathrm{Spa}(C,C^+)) \ar[r]\ar[r]\ar[r] &0,\\
}
\]
$A\rightarrow C$ is some member in $\mathfrak{C}$. Here the first row is exact by $B$-\'etale-pseudoflatness and the second is row is exact as well. The middle vertical arrow is an isomorphism while the rightmost vertical arrow is surjective as in \cite[Lemma 2.5.13]{8KL2}. Then we apply the same argument to the leftmost vertical arrow (with a new diagram), we will have the corresponding surjectivity of this leftmost vertical arrow. This will imply that the rightmost vertical arrow is injective as well. Then we have the corresponding isomorphism for the rightmost vertical arrow.
\end{proof}

\begin{theorem}\mbox{\bf{(After Kedlaya-Liu \cite[Theorem 2.5.14]{8KL2})}} \label{theorem2.22}
Taking global section will realize the equivalence between the following two categories: A. The category of all the pseudocoherent $\mathcal{O}_{\mathrm{Spa}(A,A^+)_{\text{\'et}}}\widehat{\otimes}B$-sheaves; B. The category of all the $B$-\'etale-stably pseudocoherent modules over $A\widehat{\otimes}B$. 	
\end{theorem}

\begin{proof}
See \cite[Theorem 2.5.14]{8KL2}. We need to still apply \cite[Theorem 8.2.22]{8KL1}, as long as one considers instead in our situation, and \cite[Tag 03OD]{8SP}.\\	
\end{proof}

\subsection{Noncommutative Deformation over Quasi-Stein Spaces}

\indent In this subsection we consider the corresponding noncommutative deformation over the corresponding context in \cite[Chapter 2.6]{8KL2}.

\begin{setting}
Let $X$ be a corresponding quasi-Stein adic space over $\mathbb{Q}_p$ or $\mathbb{F}_p((t))$ in the sense of \cite[Definition 2.6.2]{8KL2}. Recall that what is happening is that $X$ could be written as the corresponding direct limit of affinoids $X:=\varinjlim_i X_i$.\\ 	
\end{setting}

\begin{lemma} \mbox{\bf{(After Kedlaya-Liu \cite[Lemma 2.6.3]{8KL2})}} \label{lemma2.24}
We now consider the corresponding rings $A_i:=\mathcal{O}_{X_i}$ for all $i=0,1,...$, and in our current situation we consider the corresponding rings $A_i\widehat{\otimes}B$ (over $\mathbb{Q}_p$ or $\mathbb{F}_p((t))$) for all $i=0,1,...$. And in our situation we consider the corresponding modules $M_i$ over $A_i\widehat{\otimes}B$ for all $i=0,1,...$ with the same requirement as in \cite[Lemma 2.6.3]{8KL2} (namely those complete with respect to the natural topology). Suppose that we have bounded surjective map from $f_i:A_{i}\widehat{\otimes}B\otimes_{A_{i+1}\widehat{\otimes}B} M_{i+1}\rightarrow M_i,i=0,1,...$. Then we have first the density of the corresponding image of $\varprojlim_i M_i$ in each $M_i$ for any $i=0,1,2,...$. And we have as well the corresponding vanishing of $R^1\varprojlim_i M_i$.
\end{lemma}

\begin{proof}
This is the $B$-relative version of the result in \cite[Lemma 2.6.3]{8KL2}. For the first statement we just choose sequence of Banach norms on all the corresponding modules for all $i=0,1,...$ such that we have $\|f_i(x_{i+1})\|_i\leq 1/2\|x_{i+1}\|_{i+1}$ for any $x_{i+1}\in M_{i+1}$. Then for any $x_i\in M_i$ and any $\delta>0$, we now consider for any $j\geq 1$ the corresponding $x_{i+j}$ such that we have $\|x_{i+j}-f_{i+j+1}(x_{i+j+1})\|_{i+j}\leq \delta$. Then the sequence $x_{i+j+k},k=0,1,...$ will converge to some well-defined $y_{i+j}$ with in our situation the corresponding $y_{i}=f_{i}(y_{i+1})$. We then have $\|x_i-y_i\|_i\leq \delta$. This will prove the first statement. For the second statement as in \cite[Lemma 2.6.3]{8KL2} we form the product $M_0\times M_1\times M_2\times...$ and the consider the induced map $F$ from $M_{i+1}\rightarrow M_i$, and consider the corresponding cokernel of the map $1-F$ since this is just the corresponding limit we are considering. Then to show that the cokernel is zero we just look at the corresponding cokernel of the corresponding map on the corresponding completed direct summand which will project to the original one. But then we will have $\|f_i(v)\|_i\leq 1/2 \|v\|_i$, which produces an inverse to $1-F$ which will basically finish the proof for the second statement. 	
\end{proof}

\begin{proposition} \mbox{\bf{(After Kedlaya-Liu \cite[Lemma 2.6.4]{8KL2})}} In the same situation as above, suppose we have that the corresponding modules $M_i$ are basically $B$-stably pseudocoherent over the rings $A_i\widehat{\otimes}B$ for all $i=0,1,...$. Now we consider the situation where $f_i:A_i\widehat{\otimes}B\otimes_{A_{i+1}\widehat{\otimes}B}M_{i+1}\rightarrow M_i$ is an isomorphism. Then the conclusion in our situation is then that the corresponding projection from $\varprojlim M_i$ to $M_i$ for each $i=0,1,2,...$ is an isomorphism.

\end{proposition}

\begin{proof}
This is a $B$-relative version of the \cite[Lemma 2.6.4]{8KL2}. We adapt the argument to our situation as in the following. First we choose some finite free covering $T$ of the limit $M$ such that we have for each $i$ the corresponding map $T_i\rightarrow M_i$ is surjective. Then we consider the index $j\geq i$ and set the kernel of the map from $T_j$ to $M_j$ to be $S_j$. By the direct analog of \cite[Lemma 2.5.6]{8KL2} we have that $A_i\widehat{\otimes}B\otimes S_j\overset{}{\rightarrow}S_i$ realizes the isomorphism, and we have that the corresponding surjectivity of the corresponding map from $\varprojlim_i S_i$ projecting to the $S_i$. Then one could finish the proof by 5-lemma to the following commutative diagram as in \cite[Lemma 2.6.4]{8KL2}:
\[ \tiny
\xymatrix@C+0pc@R+6pc{
 &(A_i\widehat{\otimes}B)\otimes \varprojlim_i S_i \ar[r]\ar[r]\ar[r] \ar[d]\ar[d]\ar[d] &(A_i\widehat{\otimes}B)\otimes \varprojlim_i F_i \ar[r]\ar[r]\ar[r] \ar[d]\ar[d]\ar[d] &(A_i\widehat{\otimes}B)\otimes \varprojlim_i M_i\ar[r]\ar[r]\ar[r]  \ar[d]\ar[d]\ar[d]&0,\\
0 \ar[r]\ar[r]\ar[r] &S_i   \ar[r]\ar[r]\ar[r] &F_i \ar[r]\ar[r]\ar[r] &M_i \ar[r]\ar[r]\ar[r] &0.\\
}
\]
 
\end{proof}

\begin{proposition}  \mbox{\bf{(After Kedlaya-Liu \cite[Theorem 2.6.5]{8KL2})}}
For any \\quasi-compact adic affinoid space of $X$ which is denoted by $Y$, we have that the map $\mathcal{M}(X)\rightarrow \mathcal{M}(Y)$ is surjective for any $B$-stably pseudocoherent sheaf $\mathcal{M}$ over the sheaf $\mathcal{O}_X\widehat{\otimes}B$.	
\end{proposition}

\begin{proof}
This is just the corresponding corollary of the previous proposition.	
\end{proof}

\begin{proposition}  \mbox{\bf{(After Kedlaya-Liu \cite[Theorem 2.6.5]{8KL2})}}
We have that the stalk $\mathcal{M}_x$ is generated over the stalk $\mathcal{O}_{X,x}$ for any $x\in X$ by $M(X)$, for any $B$-stably pseudocoherent sheaf $\mathcal{M}$ over the sheaf $\mathcal{O}_X\widehat{\otimes}B$.	
\end{proposition}

\begin{proof}
This is just the corresponding corollary of the proposition before the previous proposition.	
\end{proof}

\begin{proposition}  \mbox{\bf{(After Kedlaya-Liu \cite[Theorem 2.6.5]{8KL2})}} \label{proposition2.28}
For any\\ quasi-compact adic affinoid space of $X$ which is denoted by $Y$, we have that the corresponding vanishing of the corresponding sheaf cohomology groups $H^k(X,\mathcal{M})$ of $\mathcal{M}$ for higher $k>0$, for any $B$-stably pseudocoherent sheaf $\mathcal{M}$ over the sheaf $\mathcal{O}_X\widehat{\otimes}B$.	
\end{proposition}

\begin{proof}
We follow the idea of the proof of \cite[Theorem 2.6.5]{8KL2} by comparing this to the corresponding \v{C}ech cohomology with some covering $\mathfrak{X}=\{X_1,...,X_N,...\}$:
\begin{align}
\breve{H}^k(X,\mathfrak{X}=\{X_1,...,X_N,...\};\mathcal{M})=H^k(X,\mathcal{M}),\\
\breve{H}^k(X_i,\mathfrak{X}=\{X_1,...,X_i\};\mathcal{M})=H^k(X_i,\mathcal{M})=0,k\geq 1.
\end{align}
Here we have applied the corresponding \cite[Tag 01EW]{8SP}. Now we consider the situation where $k>1$:
\begin{align}
\breve{H}^{k-1}(X_{j+1},\mathfrak{X}=\{X_1,...,X_{j+1}\};\mathcal{M})	\rightarrow \breve{H}^{k-1}(X_j,\mathfrak{X}=\{X_1,...,X_{j}\};\mathcal{M})\rightarrow 0,
\end{align}
which induces the following isomorphism by \cite[2.6 Hilfssatz]{8Kie1}:
\begin{align}
\varprojlim_{j\rightarrow \infty} \breve{H}^{k-1}(X_j,\mathfrak{X}=\{X_1,...,X_{j}\};\mathcal{M})	\overset{\sim}{\rightarrow} \breve{H}^k(X,\mathfrak{X}=\{X_1,...,X_N,...\};\mathcal{M}). 
\end{align}
Then we have the corresponding results for the index $k>1$. For $k=1$, one can as in \cite[Theorem 2.6.5]{8KL2} relate this to the corresponding $R^1\varprojlim_i$, which will finishes the proof in the same fashion.\\

\end{proof}

\begin{corollary}  \mbox{\bf{(After Kedlaya-Liu \cite[Corollary 2.6.6]{8KL2})}}
The corresponding functor from the corresponding $B$-deformed pseudocoherent sheaves over $X$ to the corresponding $B$-stably by taking the corresponding global section is an exact functor.
\end{corollary}

\begin{corollary}  \mbox{\bf{(After Kedlaya-Liu \cite[Corollary 2.6.8]{8KL2})}}
Consider a particular $\mathcal{O}_X\widehat{\otimes}B$-pseudocoherent sheaf $\mathcal{M}$ which is finite locally free throughout the whole space $X$. Then we have that the global section $\mathcal{M}(X)$ as $\mathcal{O}_X(X)\widehat{\otimes}B$ left module admits the corresponding structures of finite projective structure if and only if we have the corresponding global section is finitely generated.
\end{corollary}

\begin{proof}
As in \cite[Corollary 2.6.8]{8KL2} one could find some global splitting through the local splittings.
\end{proof}

\indent We now consider the following $B$-relative analog of \cite[Proposition 2.6.17]{8KL2}:

\begin{theorem} \mbox{\bf{(After Kedlaya-Liu \cite[Proposition 2.6.17]{8KL2})}} \label{theorem2.31}
Consider the following two statements for a particular $\mathcal{O}_X\widehat{\otimes}B$-pseudocoherent sheaf $\mathcal{M}$. First is that one can find finite many generators (the number is up to some uniform fixed integer $n\geq 0$) for each section of $\mathcal{M}(X_i)$ for each $i=1,2,...$. The second statement is that the global section $\mathcal{M}(X)$ is just finitely generated. Then in our situation the two statement is equivalent if we have that the corresponding space $X$ admits an $m$-uniform covering in the exact same sense of \cite[Proposition 2.6.17]{8KL2}. 	
\end{theorem}

\begin{proof}
One direction is obvious, the other direct could be proved in the same way as in \cite[Proposition 2.6.17]{8KL2} where essentially the information on $X$ does not change at all. Namely as in \cite[Proposition 2.6.17]{8KL2} we could basically consider one single subcovering $\{Y_u\}_{u\in U}$ indexed by $U$ of the covering which is $m$-uniform. For any $u\in U$ we consider the smallest $i$ such that $X_i\bigcap Y_u$, and we denote this by $i(u)\geq 0$. Then we form:
\begin{align}
V_{u}:=X_{i(u)}\bigcup_{v\in U,v\leq u} Y_{v}.	
\end{align}
Then as in \cite[Proposition 2.6.17]{8KL2} we can find some $x_u,u\in U$ (which are denoted by the general form of $x_i$ in \cite[Proposition 2.6.17]{8KL2}) in the global section of $\mathcal{O}_X(X)$ such that it restricts to some topologically nilpotent onto the space $V_u$ and it restricts to some inverse of some topologically nilpotent to $Y_u$. Then we can build up a corresponding generating set for the global section by the following approximation process just as in \cite[Proposition 2.6.17]{8KL2}. To be more precise what we have consider is to modify the corresponding generator for each $Y_u$ (for instance let them be denoted by $y_1,...,y_u$) to be $y_{u,1},...,y_{u,n}$ ($n$ is the corresponding uniform integer in the statement of the theorem) for all $u\in U$. It is achieved through induction, when we have $y_{u,1},...,y_{u,n}$ then we will set that to be $0$ otherwise in the situation where there exists some predecessor $u'$ of $u$ in the corresponding set $U$ we then set this to be $y_{u',1},...,y_{u',n}$. And we set as in \cite[Proposition 2.6.17]{8KL2}:
\begin{displaymath}
y_{u,j}:=	y_{u-1,j}+x_u^cy_u
\end{displaymath}
by lifting the corresponding power $c$ to be as large as possible. This will guarantee the convergence of:
\begin{align}
\lim_{u\rightarrow \infty}\{y_{u,1},...,y_{u,n}\},	
\end{align}
which gives the set of the corresponding global generators desired.
\end{proof}

%%\newpage 

\newpage\section{Foundations on Noncommutative Descent for Adic Spectra in Pro-Banach Case}

\subsection{Noncommutative Pseudocoherence in Analytic Topology}

\indent We now establish some foundations in the noncommutative setting on glueing noncommutative pseudocoherent modules after \cite[Chapter 2]{8KL2}. However the corresponding difference we will make here from the previous section is that we will consider some further Iwasawa deformation similar to the sense of Fukaya-Kato's adic rings. We first consider the simpler situation where we will deform over $B[[G]]$ where $B$ as before comes from the following beginning setting up and $G$ is a given $p$-adic Lie group:

\begin{setting} \label{8setting3.1}
Consider a corresponding sheafy Banach adic uniform algebra $(A,A^+)$ over $\mathbb{Q}_p$, we consider the base space $\mathrm{Spa}(A,A^+)$. And now we will consider a further noncommutative Banach algebra $(B,B^+)$ over $\mathbb{Q}_p$. 	
\end{setting}

\begin{setting}
In what follow when we mean a module $M$ over $A\widehat{\otimes}B[[G]]$ we will mean a pro-system $\{M_E\}_{E\subset G}$ of left module over $A\widehat{\otimes}B[[G]]$. And the corresponding definition of pseudocoherence is also defined on the adic system $\{M_E\}_{E\subset G}$ namely for each member in the system. We abuse the notation a bit when we write the notation $[[G]]/E$ which will really mean that $\mathbb{Q}_p[[G/E]]$. Here we have considered a corresponding representation of the algebra $B[[G]]$ by $\varprojlim_{E} B[[G/E]]$.	
\end{setting}

\begin{definition} \mbox{\bf{(After Kedlaya-Liu \cite[Definition 2.4.1]{8KL2})}}
For any left\\ $A\widehat{\otimes}B[[G]]$-module $M$, we call it $m$-$B[[G]]$-stably pseudocoherent if we have that it is $m$-$B[[G]]$-pseudocoherent, complete for the natural topology and for any morphism $A\rightarrow A'$ which is the corresponding rational localization, the base change of $M$ to $A'\widehat{\otimes} B[[G]]$ is complete for the natural topology as the corresponding left $A'\widehat{\otimes} B[[G]]$-module (namely each member in the family is still remained to be complete after the corresponding rational localization). As in \cite[Definition 2.4.1]{8KL2} we call that the corresponding left $A\widehat{\otimes}B[[G]]$-module $M$ just $B[[G]]$-stably pseudocoherent if we have that it is simply just $\infty$-$B[[G]]$-stably pseudocoherent.
\end{definition}

\indent Since our spaces are commutative, so the rational localization is not flat, therefore in our situation we need to define something weaker than the corresponding flatness which will be sufficient for us to apply in the following development: 

\begin{definition}\mbox{\bf{(After Kedlaya-Liu \cite[Definition 2.4.4]{8KL2})}}
For any Banach left $A\widehat{\otimes}B[[G]]$-module $M$, we call it is $m$-$B[[G]]$-pseudoflat if for any right $A\widehat{\otimes}B[[G]]$-module $M'$ $m$-$B[[G]]$-stably pseudocoherent we have $\mathrm{Tor}_1^{A\widehat{\otimes}B[[G]]}(M',M)=0$ (namely\\ $\{\mathrm{Tor}_1^{A\widehat{\otimes}B[[G]]/E}(M'_E,M_E)\}_{E\subset G}=\{0\}_{E\subset G}$). For any Banach right $A\widehat{\otimes}B[[G]]$-module $M$, we call it is $m$-$B[[G]]$-pseudoflat if for any left $A\widehat{\otimes}B[[G]]$-module $M'$ $m$-$B[[G]]$-stably pseudocoherent we have
\begin{center}
 $\mathrm{Tor}_1^{A\widehat{\otimes}B[[G]]}(M,M')=0$.
\end{center}
\end{definition}

\begin{definition}\mbox{\bf{(After Kedlaya-Liu \cite[Definition 2.4.6]{8KL2})}} We consider the corresponding notion of the corresponding pro-projective modules. We define over $A$ the corresponding $B[[G]]$-pro-projective module $M$ to be a corresponding left $A\widehat{\otimes}B[[G]]$-module such that one could find a filtered sequence of projectors such that in the sense of taking the prolimit we have the corresponding projector will converge to any chosen element in $M$. Here we assume the module is complete with respect to natural topology, and we assume the projectors are $A\widehat{\otimes}B[[G]]$-linear and we assume that the corresponding image of the projectors are modules which are also finitely generated and projective.
	
\end{definition}

\begin{lemma}\mbox{\bf{(After Kedlaya-Liu \cite[Lemma 2.4.7]{8KL2})}}
Suppose we have over the space $\mathrm{Spa}(A,A^+)$ a corresponding $B[[G]]$-pro-projective left module $M$. And suppose that we have a $2$-$B[[G]]$-pseudoflat right module $C$ over $A$. And here we assume that $C$ is complete. Then we have that the corresponding product $C{\otimes}_{A\widehat{\otimes}B[[G]]}M$ is then complete under the corresponding natural topology. And moreover we have that in our situation:
\begin{displaymath}
\mathrm{Tor}_1^{A\widehat{\otimes}B[[G]]}(C,M)=0.	
\end{displaymath}

\end{lemma}

\begin{proof}
This will be just a noncommutative version of the corresponding  \cite[Lemma 2.4.7]{8KL2}. We adapt the corresponding argument in \cite[Lemma 2.4.7]{8KL2} to our situation. We just consider each member in the families for the adic modules. What we are going to consider in this case is then first to consider the following presentation:
\[
\xymatrix@C+0pc@R+0pc{
M'  \ar[r]\ar[r]\ar[r] &(A\widehat{\otimes}B[[G]]/E)^k \ar[r]\ar[r]\ar[r] & M\ar[r]\ar[r]\ar[r] &0,
}
\]
where the left module $M'$ is finitely presented. Then we consider the following commutative diagram:
\[
\xymatrix@C+0pc@R+3pc{
&C\otimes_{A\widehat{\otimes}B[[G]]/E}M'  \ar[r]\ar[r]\ar[r]  \ar[d]\ar[d]\ar[d] &C\otimes_{A\widehat{\otimes}B[[G]]/E}(A\widehat{\otimes}B[[G]]/E)^k \ar[r]\ar[r]\ar[r] \ar[d]\ar[d]\ar[d] & C\otimes_{A\widehat{\otimes}B[[G]]/E}M\ar[r]\ar[r]\ar[r] \ar[d]\ar[d]\ar[d] &0,\\
0 \ar[r]\ar[r]\ar[r] &C\widehat{\otimes}_{A\widehat{\otimes}B[[G]]/E}M'  \ar[r]\ar[r]\ar[r] &C\widehat{\otimes}_{A\widehat{\otimes}B[[G]]/E}(A\widehat{\otimes}B[[G]]/E)^k \ar[r]\ar[r]\ar[r] & C\widehat{\otimes}_{A\widehat{\otimes}B[[G]]/E}M\ar[r]\ar[r]\ar[r] &0.
}
\]\\
The first row is exact as in \cite[Lemma 2.4.7]{8KL2}, and we have the corresponding second row is also exact since we have that by hypothesis the corresponding module $M$ is $[[G]]/E$-pro-projective.  
Then as in \cite[Lemma 2.4.7]{8KL2} the results follow.	
\end{proof}

\begin{proposition} \mbox{\bf{(After Kedlaya-Liu \cite[Corollary 2.4.8]{8KL2})}}
Over the space\\ $\mathrm{Spa}(A,A^+)$, we have any $B[[G]]$-pro-projective left module is $2$-$B[[G]]$-pseudoflat.
\end{proposition}

\begin{proof}
This is a direct consequence of the previous lemma.	
\end{proof}

\begin{lemma} \mbox{\bf{(After Kedlaya-Liu \cite[Corollary 2.4.9]{8KL2})}}
Keep the notation above, suppose we are working over the adic space $
\mathrm{Spa}(A,A^+)$. Suppose now the left $A\widehat{\otimes}B[[G]]$-module is finitely generated (again note that we are talking about a family of objects over the corresponding quotients), then we have that the following natural maps:
\begin{align}
A\{T\}\widehat{\otimes}B[[G]]\otimes_{A\widehat{\otimes}B[[G]]}M 	\rightarrow M\{T\},\\\
A\{T,T^{-1}\}\widehat{\otimes}B[[G]] \otimes_{A\widehat{\otimes}B[[G]]} M	\rightarrow M\{T,T^{-1}\}	
\end{align}
are surjective. Suppose now the left $A\widehat{\otimes}B$-module is finitely presented, then we have that the following natural maps:
\begin{align}
A\{T\}\widehat{\otimes}B[[G]]\otimes_{A\widehat{\otimes}B[[G]]}M 	\rightarrow M\{T\},\\\
A\{T,T^{-1}\}\widehat{\otimes}B[[G]] \otimes_{A\widehat{\otimes}B[[G]]} M	\rightarrow M\{T,T^{-1}\}	
\end{align}
are bijective. All modules here are assumed to be complete.

\end{lemma}

\begin{proof}
This is the corresponding consequence of the previous lemma.
\end{proof}

\begin{proposition} \mbox{\bf{(After Kedlaya-Liu \cite[Lemma 2.4.12]{8KL2})}}
We keep the corresponding notations in \cite[Lemma 2.4.10]{8KL2}. Then in our current sense we have that the corresponding morphism $A\rightarrow B_2$ is then $2$-$B[[G]]$-pseudoflat. 
\end{proposition}

\begin{proof}
The corresponding proof could be made parallel to the corresponding proof of \cite[Lemma 2.4.12]{8KL2}. To be more precise we consider the following commutative diagram from the corresponding chosen exact sequence (which is just the analog of the corresponding one in \cite[Lemma 2.4.12]{8KL2}, namely as below we choose arbitrary short exact sequence with $P$ as a left $A\widehat{\otimes}B[[G]]/E$-module which is $2$-$B[[G]]/E$-stably pseudocoherent):
\[
\xymatrix@C+0pc@R+0pc{
0   \ar[r]\ar[r]\ar[r] &M \ar[r]\ar[r]\ar[r] &N \ar[r]\ar[r]\ar[r] & P \ar[r]\ar[r]\ar[r] &0,
}
\]
which then induces the following corresponding big commutative diagram:
\[\tiny
\xymatrix@C+0pc@R+6pc{
& &0 \ar[d]\ar[d]\ar[d] &0 \ar[d]\ar[d]\ar[d] &,\\
0   \ar[r]\ar[r]\ar[r]  &{A\{T\}\widehat{\otimes}B[[G]]/E}\otimes_{A\widehat{\otimes}B[[G]]/E} M \ar[r]\ar[r]\ar[r] \ar[d]^{1-fT}\ar[d]\ar[d] &{A\{T\}\widehat{\otimes}B[[G]]/E}\otimes_{A\widehat{\otimes}B[[G]]/E} N \ar[r]\ar[r]\ar[r] \ar[d]^{1-fT}\ar[d]\ar[d] &{A\{T\}\widehat{\otimes}B[[G]]/E} \otimes_{A\widehat{\otimes}B[[G]]/E}P \ar[r]\ar[r]\ar[r] \ar[d]^{1-fT}\ar[d]\ar[d] &0,\\
0   \ar[r]\ar[r]\ar[r] &{A\{T\}\widehat{\otimes}B[[G]]/E}\otimes_{A\widehat{\otimes}B[[G]]/E}M  \ar[r]\ar[r]\ar[r] \ar[d]\ar[d]\ar[d] &{A\{T\}\widehat{\otimes}B[[G]]/E} \otimes_{A\widehat{\otimes}B[[G]]/E}N\ar[r]\ar[r]\ar[r] \ar[d]\ar[d]\ar[d] & {A\{T\}\widehat{\otimes}B[[G]]/E}\otimes_{A\widehat{\otimes}B[[G]]/E}P \ar[r]\ar[r]\ar[r] \ar[d]\ar[d]\ar[d] &0,\\
0  \ar[r]^?\ar[r]\ar[r] &{B_2\widehat{\otimes}B[[G]]/E}\otimes_{A\widehat{\otimes}B[[G]]/E} M \ar[r]\ar[r]\ar[r] \ar[d]\ar[d]\ar[d] &{B_2\widehat{\otimes}B[[G]]/E}\otimes_{A\widehat{\otimes}B[[G]]/E}N \ar[r]\ar[r]\ar[r] \ar[d]\ar[d]\ar[d] & {B_2\widehat{\otimes}B[[G]]/E}\otimes_{A\widehat{\otimes}B[[G]]/E} P \ar[r]\ar[r]\ar[r] \ar[d]\ar[d]\ar[d] &0,\\
&0&0&0
}
\]
with the notations in \cite[Lemma 2.4.10]{8KL2} where we actually have the corresponding exactness along the horizontal direction at the corner around $M\otimes_{A\widehat{\otimes}B[[G]]/E} {B_2\widehat{\otimes}B[[G]]/E}$ marked with $?$, by diagram chasing.
\end{proof}

\begin{proposition} \mbox{\bf{(After Kedlaya-Liu \cite[Lemma 2.4.13]{8KL2})}} \label{proposition3.10}
We keep the corresponding notations in \cite[Lemma 2.4.10]{8KL2}. Then in our current sense we have that the corresponding morphism $A\rightarrow B_1$ is then $2$-$B[[G]]$-pseudoflat. And in our current sense we have that the corresponding morphism $A\rightarrow B_{1,2}$ is then $2$-$B[[G]]$-pseudoflat. 
\end{proposition}

\begin{proof}
The statement for $A\rightarrow B_{1,2}$ could be proved by consider the composition $A\rightarrow B_2\rightarrow B_{12}$ where the $2$-$B[[G]]$-pseudoflatness could be proved as in \cite[Lemma 2.4.13]{8KL2} by inverting the corresponding variable. For $A\rightarrow B_{1}$,
the corresponding proof could be made parallel to the corresponding proof of \cite[Lemma 2.4.13]{8KL2}. To be more precise we consider the following commutative diagram from the corresponding chosen exact sequence (which is just the analog of the corresponding one in \cite[Remark 2.4.5]{8KL2}, namely as below we choose arbitrary short exact sequence with $P$ as a left $A\widehat{\otimes}B[[G]]/E$-module which is $2$-$B[[G]]/E$-stably pseudocoherent):
\[
\xymatrix@C+0pc@R+0pc{
0   \ar[r]\ar[r]\ar[r] &M \ar[r]\ar[r]\ar[r] &N \ar[r]\ar[r]\ar[r] & P \ar[r]\ar[r]\ar[r] &0,
}
\]
which then induces the following corresponding big commutative diagram:
\[\tiny
\xymatrix@C+0pc@R+6pc{
& &0 \ar[d]\ar[d]\ar[d] &0 \ar[d]\ar[d]\ar[d] &,\\
0   \ar[r]\ar[r]\ar[r]  &M \ar[r]\ar[r]\ar[r] \ar[d]\ar[d]\ar[d] &N \ar[r]\ar[r]\ar[r] \ar[d]\ar[d]\ar[d] & P \ar[r]\ar[r]\ar[r] \ar[d]\ar[d]\ar[d] &0,\\
0   \ar[r]^?\ar[r]\ar[r] &({B_1\widehat{\otimes}\frac{B[[G]]}{E}}\bigoplus {B_2\widehat{\otimes}\frac{B[[G]]}{E}}) \otimes M  \ar[r]\ar[r]\ar[r] \ar[d]\ar[d]\ar[d] &({B_1\widehat{\otimes}\frac{B[[G]]}{E}}\bigoplus {B_2\widehat{\otimes}\frac{B[[G]]}{E}})\otimes N  \ar[r]\ar[r]\ar[r] \ar[d]\ar[d]\ar[d] &({B_1\widehat{\otimes}\frac{B[[G]]}{E}}\bigoplus {B_2\widehat{\otimes}\frac{B[[G]]}{E}})\otimes P\ar[r]\ar[r]\ar[r] \ar[d]\ar[d]\ar[d] &0,\\
0  \ar[r]\ar[r]\ar[r] &{B_{12}\widehat{\otimes}\frac{B[[G]]}{E}}\otimes M \ar[r]\ar[r]\ar[r] \ar[d]\ar[d]\ar[d] &{B_{12}\widehat{\otimes}\frac{B[[G]]}{E}}\otimes N \ar[r]\ar[r]\ar[r] \ar[d]\ar[d]\ar[d] &{B_{12}\widehat{\otimes}\frac{B[[G]]}{E}} \otimes P \ar[r]\ar[r]\ar[r] \ar[d]\ar[d]\ar[d] &0,\\
&0&0&0
}
\]
with the notations in \cite[Lemma 2.4.10]{8KL2} where we actually have the corresponding exactness along the horizontal direction at the corner around
\begin{center}
 $M\otimes_{A\widehat{\otimes}B[[G]]/E} ({B_1\widehat{\otimes}B[[G]]/E}\bigoplus {B_2\widehat{\otimes}B[[G]]/E})$
\end{center} 
marked with $?$, by diagram chasing.
\end{proof}

\indent After these foundational results as in \cite{8KL2} we have the following proposition which is the corresponding noncommutative generalization of the corresponding result established in \cite[Theorem 2.4.15]{8KL2}.

\begin{proposition} \mbox{\bf{(After Kedlaya-Liu \cite[Theorem 2.4.15]{8KL2})}} \label{proposition3.11}
In our current context we have that for any rational localization $\mathrm{Spa}(A',A^{',+})\rightarrow \mathrm{Spa}(A,A^+)$ we have that along this base change the corresponding $B[[G]]$-stably pseudocoherence is preserved.	
\end{proposition}

\indent Then we have the corresponding Tate's acyclicity in the noncommutative deformed setting:

\begin{theorem}\mbox{\bf{(After Kedlaya-Liu \cite[Theorem 2.5.1]{8KL2})}} \label{theorem3.12} Now suppose we have in our corresponding \cref{setting2.1} a corresponding $B[[G]]$-stably pseudocoherent module $M$. Then we consider the corresponding assignment such that for any $U\subset \mathrm{Spa}(A,A^+)$ we define $\widetilde{M}(U)$ as in the following:
\begin{align}
\widetilde{M}(U):=\varprojlim_{\mathrm{Spa}(S,S^+)\subset U,\mathrm{rational}} S\widehat{\otimes}B[[G]]\otimes_{A\widehat{\otimes}B[[G]]}M.	
\end{align}
Then we have that for any $\mathfrak{B}$ which is a rational covering of $U=\mathrm{Spa}(S,S^+)\subset \mathrm{Spa}(A,A^+)$ (certainly this $U$ is also assumed to be rational) we have that the vanishing of the following two cohomology groups:
\begin{align}
H^i(U,\widetilde{M}), \check{H}^i(U:\mathfrak{B},\widetilde{M})
\end{align}
for any $i>0$. When concentrating at the degree zero we have:
\begin{align}
H^0(U,\widetilde{M})=S\widehat{\otimes}B[[G]]\otimes_{A\widehat{\otimes}B[[G]]}M, \check{H}^0(U:\mathfrak{B},\widetilde{M})=S\widehat{\otimes}B[[G]]\otimes_{A\widehat{\otimes}B[[G]]}M.
\end{align}
	
\end{theorem}

\begin{proof}
By \cite[Propositions 2.4.20-2.4.21]{8KL1}, we could then finish proof as in \cite[Theorem 2.5.1]{8KL2} by using our previous \cref{proposition3.11} and \cref{proposition3.10} as above.	
\end{proof}

\indent We now consider the corresponding noncommutative deformed version of Kiehl's glueing properties for stably pseudocoherent modules after \cite{8KL2}:

\begin{definition} \mbox{\bf{(After Kedlaya-Liu \cite[Definition 2.5.3]{8KL2})}}
Consider in our context (over $\mathrm{Spa}(A,A^+)$) the corresponding sheaves $\mathcal{O}\widehat{\otimes}B[[G]]$, we will then define the corresponding pseudocoherent sheaves over $\mathrm{Spa}(A,A^+)$ to be those locally defined by attaching stably-pseudocoherent modules over the section. 
\end{definition}

\begin{lemma} \mbox{\bf{(After Kedlaya-Liu \cite[Lemma 2.5.4]{8KL2})}}\label{lemma3.14}
	Consider the corresponding notations in \cite[Lemma 2.4.10]{8KL2}, we have the corresponding morphism $A\rightarrow B_1\bigoplus B_2$. Then we have that this morphism is an descent morphism effective for the corresponding $B[[G]]$-stably pseodocoherent Banach modules. 
\end{lemma}

\begin{proof}
This is a $B[[G]]$-relative version of the \cite[Lemma 2.5.4]{8KL2}. We adapt the corresponding proof to our situation. First consider the corresponding notations on the commutative adic space $\mathrm{Spa}(A,A^+)$ as in \cite[Lemma 2.4.10]{8KL2}. And we have a coving $\{U_1,U_2\}$. Then consider a descent datum for this cover of $B[[G]]$-stably pseudocoherent sheaves. By \cref{theorem3.12} we could realize this as a corresponding pseudocoherent sheaf over $\mathcal{O}\widehat{\otimes}B[[G]]$. We denote this by $\mathcal{M}$. Then by (the proof of) \cite[Lemma 6.82]{8T2}, we have that this is finitely generated and we have the corresponding vanishing of $\check{H}^i(\mathrm{Spa}(A,A^+),\mathfrak{B};\mathcal{M})$ for $i>0$. Then by \cref{theorem3.12} again we have that the vanishing of $H^i(U_j;\mathcal{M})$ for $j\in \{1,2\},i>0$. Then as in \cite[Lemma 2.5.4]{8KL2} we have the corresponding vanishing of $H^i(\mathrm{Spa}(A,A^+);\mathcal{M})$ for $i>0$. Then as in \cite[Lemma 2.5.4]{8KL2} we choose a corresponding covering of some finite free object $\mathcal{F}$ and take its kernel $\mathcal{K}$ which is actually pseudocoherent by the corresponding \cref{proposition3.11}. Then forming the following diagram:
\[\tiny
\xymatrix@C-0.4pc@R+6pc{
0 \ar[r]\ar[r]\ar[r] & (B_p\widehat{\otimes}B[[G]]/E)\otimes \mathcal{K}(\mathrm{Spa}(A,A^+))  \ar[r]\ar[r]\ar[r] \ar[d]\ar[d]\ar[d] &(B_p\widehat{\otimes}B[[G]]/E)\otimes \mathcal{F}(\mathrm{Spa}(A,A^+))\ar[r]\ar[r]\ar[r] \ar[d]\ar[d]\ar[d] & (B_p\widehat{\otimes}B[[G]]/E)\otimes \mathcal{M}(\mathrm{Spa}(A,A^+))\ar[r]\ar[r]\ar[r]  \ar[d]\ar[d]\ar[d]&0,\\
0 \ar[r]\ar[r]\ar[r] &\mathcal{K}(U_p)   \ar[r]\ar[r]\ar[r] &\mathcal{F}(U_p) \ar[r]\ar[r]\ar[r] & \mathcal{M}(U_p) \ar[r]\ar[r]\ar[r] &0,\\
}
\]
where we have $p=1,2$. The first row is exact by $B[[G]]/E$-pseudoflatness as in \cref{8proposition2.10} and the second is row is exact as well by \cref{theorem2.11}. The middle vertical arrow is an isomorphism while the rightmost vertical arrow is surjective and the leftmost vertical arrow is surjective. The surjectivity just here comes from actually \cite[Lemma 6.82]{8T2}. This will imply that the rightmost vertical arrow is injective as well. Then we have the corresponding isomorphism $\mathcal{M}\overset{\sim}{\rightarrow} \widetilde{\mathcal{M}(\mathrm{Spa}(A,A^+))}$. Then we consider the following diagram:

\[ \tiny
\xymatrix@C+0pc@R+6pc{
 &(C\widehat{\otimes}B[[G]]/E)\otimes\mathcal{K}(\mathrm{Spa}(A,A^+))   \ar[r]\ar[r]\ar[r] \ar[d]\ar[d]\ar[d] &(C\widehat{\otimes}B[[G]]/E)\otimes\mathcal{F}(\mathrm{Spa}(A,A^+)) \ar[r]\ar[r]\ar[r] \ar[d]\ar[d]\ar[d] & (C\widehat{\otimes}B[[G]]/E)\otimes\mathcal{M}(\mathrm{Spa}(A,A^+))\ar[r]\ar[r]\ar[r]  \ar[d]\ar[d]\ar[d]&0,\\
0 \ar[r]\ar[r]\ar[r] &\mathcal{K}(\mathrm{Spa}(C,C^+))   \ar[r]\ar[r]\ar[r] &\mathcal{F}(\mathrm{Spa}(C,C^+)) \ar[r]\ar[r]\ar[r] & \mathcal{M}(\mathrm{Spa}(C,C^+)) \ar[r]\ar[r]\ar[r] &0,\\
}
\]
$C$ is any rational localization of $A$. Here the first row is exact by $B[[G]]/E$-pseudoflatness and the second is row is exact as well. The middle vertical arrow is an isomorphism while the rightmost vertical arrow is surjective as in \cite[Lemma 2.5.4]{8KL2}. Then we apply the same argument to the leftmost vertical arrow (with a new diagram), we will have the corresponding surjectivity of this leftmost vertical arrow. This will imply that the rightmost vertical arrow is injective as well. Then we have the corresponding isomorphism for the rightmost vertical arrow.
\end{proof}

\begin{theorem}\mbox{\bf{(After Kedlaya-Liu \cite[Theorem 2.5.5]{8KL2})}} \label{theorem3.15}
Taking global section will realize the equivalence between the following two categories: A. The category of all the pseudocoherent $\mathcal{O}\widehat{\otimes}B[[G]]$-sheaves; B. The category of all the $B[[G]]$-stably pseudocoherent modules over $A\widehat{\otimes}B[[G]]$. 	
\end{theorem}

\begin{proof}
See \cite[Theorem 2.5.5]{8KL2}. We need to still apply \cite[Proposition 2.4.20]{8KL1}, as long as one considers instead in our situation \cref{lemma3.14} and \cref{theorem3.12}.\\	
\end{proof}

\subsection{Noncommutative Pseudocoherence in \'Etale Topology}

\indent We consider the extension of the corresponding discussion in the previous subsection to the \'etale topology. At this moment we consider the following same context on the geometric level as in \cite{8KL2}:

\begin{setting} \mbox{\bf{(After Kedlaya-Liu \cite[Hypothesis 2.5.8]{8KL2})}} \label{setting3.16}
As in the previous subsection we consider now the corresponding geometric setting namely an adic space $\mathrm{Spa}(A,A^+)$ where $A$ is assumed to be sheafy. Then we consider $B[[G]]$ as above. And we consider the corresponding \'etale site $\mathrm{Spa}(A,A^+)_\text{\'et}$. And by keeping the corresponding setting in the geometry as in \cite[Hypothesis 2.5.8]{8KL2}, we assume that there is a stable basis $\mathfrak{B}$ containing the space $\mathrm{Spa}(A,A^+)$ itself.	
\end{setting}

\begin{definition} \mbox{\bf{(After Kedlaya-Liu \cite[Definition 2.5.9]{8KL2})}}
For any left\\ $A\widehat{\otimes}B[[G]]$-module $M$, we call it $m$-$B[[G]]$-\'etale stably pseudocoherent with respect to $\mathfrak{B}$ if we have that it is $m$-$B[[G]]$-pseudocoherent, complete for the natural topology and for any rational localization $A\rightarrow A'$, the base change of $M$ to $A'\widehat{\otimes} B[[G]]$ is complete for the natural topology as the corresponding left $A'\widehat{\otimes} B[[G]]$-module.	As in \cite[Definition 2.4.1]{8KL2} we call that the corresponding left $A\widehat{\otimes}B[[G]]$-module $M$ just $B[[G]]$-stably pseudocoherent if we have that it is simply just $\infty$-$B[[G]]$-stably pseudocoherent.
\end{definition}

\begin{definition}\mbox{\bf{(After Kedlaya-Liu \cite[Below Definition 2.5.9]{8KL2})}}
For any Banach left $A\widehat{\otimes}B[[G]]$-module $M$, we call it is $m$-$B[[G]]$-\'etale-pseudoflat if for any right $A\widehat{\otimes}B[[G]]$-module $M'$ $m$-$B[[G]]$-\'etale-stably pseudocoherent we have $\mathrm{Tor}_1^{A\widehat{\otimes}B[[G]]}(M',M)=0$. For any Banach right $A\widehat{\otimes}B[[G]]$-module $M$, we call it is $m$-$B[[G]]$-\'etale-pseudoflat if for any left $A\widehat{\otimes}B[[G]]$-module $M'$ $m$-$B[[G]]$-\'etale-stably pseudocoherent we have
\begin{center}
 $\mathrm{Tor}_1^{A\widehat{\otimes}B[[G]]}(M,M')=0$.
\end{center}
\end{definition}

\indent The following proposition holds in our current setting.

\begin{proposition} \mbox{\bf{(After Kedlaya-Liu \cite[Lemma 2.5.10]{8KL2})}} \label{proposition3.19}
One can actually find another basis $\mathfrak{C}$ in $\mathfrak{B}$ such that any morphism in $\mathfrak{C}$ could be $2$-$B[[G]]$-\'etale-pseudoflat with respect to either $\mathfrak{C}$ or $\mathfrak{B}$.	
\end{proposition}

\begin{proof}
We follow the proof of \cite[Lemma 2.5.10]{8KL2} in our current $B[[G]]$-relative situation, however not that much needs to proof by relying on the proof of \cite[Lemma 2.5.10]{8KL2}. To be more precise the corresponding selection of the new basis $\mathfrak{C}$ comes from including all the morphism made up of some composition of rational localization and the corresponding finite \'etale ones. For these we have shown above the corresponding 2-$B[[G]]$-\'etale pseudoflatness. Then the rest will be pure geometric for the analytic spaces which are just the same as the situation of \cite[Lemma 2.5.10]{8KL2}, therefore we just omit the corresponding argument, see \cite[Lemma 2.5.10]{8KL2}.	
\end{proof}

\indent Then we have the corresponding Tate's acyclicity in the noncommutative deformed setting in \'etale topology (here fix $\mathfrak{C}$ as above):

\begin{theorem}\mbox{\bf{(After Kedlaya-Liu \cite[Theorem 2.5.11]{8KL2})}} \label{theorem3.20} Now suppose we have in our corresponding \cref{8setting3.1} a corresponding $B[[G]]$-\'etale-stably pseudocoherent module $M$. Then we consider the corresponding assignment such that for any $U\subset \mathrm{Spa}(A,A^+)$ we define $\widetilde{M}(U)$ as in the following:
\begin{align}
\widetilde{M}(U):=\varprojlim_{\mathrm{Spa}(S,S^+)\subset U,\in \mathfrak{C}} S\widehat{\otimes}B[[G]]\otimes_{A\widehat{\otimes}B[[G]]}M.	
\end{align}
Then we have that for any $\mathfrak{B}$ which is a covering of $U=\mathrm{Spa}(S,S^+)\subset \mathrm{Spa}(A,A^+)$ (certainly this $U$ is also assumed to be in $\mathfrak{C}$, and we assume this covering is formed by using the corresponding members in $\mathfrak{C}$) we have that the vanishing of the following two cohomology groups:
\begin{align}
H^i(U,\widetilde{M}), \check{H}^i(U:\mathfrak{B},\widetilde{M})
\end{align}
for any $i>0$. When concentrating at the degree zero we have:
\begin{align}
H^0(U,\widetilde{M})=S\widehat{\otimes}B[[G]]\otimes_{A\widehat{\otimes}B[[G]]}M, \check{H}^0(U:\mathfrak{B},\widetilde{M})=S\widehat{\otimes}B[[G]]\otimes_{A\widehat{\otimes}B[[G]]}M.
\end{align}
	
\end{theorem}

\begin{proof}
By \cite[Proposition 8.2.21]{8KL1}, we could then finish proof as in \cite[Theorem 2.5.11]{8KL2} by using our previous \cref{proposition3.19} and the corresponding faithfully flat descent as in the situation of \cite[Theorem 2.5.11]{8KL2}.	
\end{proof}

\indent We now consider the corresponding noncommutative deformed version of Kiehl's glueing properties for stably pseudocoherent modules after \cite{8KL2}:

\begin{definition} \mbox{\bf{(After Kedlaya-Liu \cite[Definition 2.5.12]{8KL2})}}
Consider in our context (over $\mathrm{Spa}(A,A^+)_{\text{\'et}}$) the corresponding sheaves $\mathcal{O}_{\mathrm{Spa}(A,A^+)_{\text{\'et}}}\widehat{\otimes}B[[G]]$, we will then define the corresponding pseudocoherent sheaves over $\mathrm{Spa}(A,A^+)_{\text{\'et}}$ to be those locally defined by attaching \'etale-stably-pseudocoherent modules over the section. 
\end{definition}

\begin{lemma} \mbox{\bf{(After Kedlaya-Liu \cite[Lemma 2.5.13]{8KL2})}}
	Consider the corresponding notations in \cite[Lemma 2.4.10]{8KL2}, we have the corresponding morphism $A\rightarrow B_1\bigoplus B_2$. Then we have that this morphism is an descent morphism effective for the corresponding $B[[G]]$-\'etale-stably pseodocoherent Banach modules. 
\end{lemma}

\begin{proof}
This is a $B[[G]]$-relative version of the \cite[Lemma 2.5.13]{8KL2}. We adapt the corresponding proof to our situation. First consider the corresponding notations on the commutative adic space $\mathrm{Spa}(A,A^+)$ as in \cite[Lemma 2.4.10]{8KL2}. And we have a coving $\{U_1,U_2\}$. Then consider a descent datum for this cover of $B[[G]]$-\'etale-stably pseudocoherent sheaves. By \cref{theorem2.19} we could realize this as a corresponding pseudocoherent sheaf over $\mathcal{O}\widehat{\otimes}B[[G]]$. We denote this by $\mathcal{M}$. Then by consider the corresponding results in the rational localization situation above we have that the vanishing of $H^i(U_j;\mathcal{M})$ for $j\in \{1,2\},i>0$. Then as in \cite[Lemma 2.5.13]{8KL2} we have the corresponding vanishing of $H^i(\mathrm{Spa}(A,A^+);\mathcal{M})$ for $i>0$. Then as in \cite[Lemma 2.5.13]{8KL2} we choose a corresponding covering of some finite free object $\mathcal{F}$ and take its kernel $\mathcal{K}$ which gives rise to the corresponding exact sequence:
\[\tiny
\xymatrix@C+0pc@R+6pc{
0 \ar[r]\ar[r]\ar[r] &\mathcal{K}(\mathrm{Spa}(A,A^+))   \ar[r]\ar[r]\ar[r] &\mathcal{F}(\mathrm{Spa}(A,A^+)) \ar[r]\ar[r]\ar[r] & \mathcal{M}(\mathrm{Spa}(A,A^+)) \ar[r]\ar[r]\ar[r] &0,
}
\]
where we have $p=1,2$. We have already the corresponding stably-pseudocoherence of the global section. As in \cite[Lemma 2.5.13]{8KL2} we only need to show that it is stable under \'etale morphism. Then we consider the following diagram:

\[ \tiny
\xymatrix@C+0pc@R+6pc{
 &(C\widehat{\otimes}B[[G]]/E)\otimes\mathcal{K}(\mathrm{Spa}(A,A^+))   \ar[r]\ar[r]\ar[r] \ar[d]\ar[d]\ar[d] &(C\widehat{\otimes}B[[G]]/E)\otimes\mathcal{F}(\mathrm{Spa}(A,A^+)) \ar[r]\ar[r]\ar[r] \ar[d]\ar[d]\ar[d] & (C\widehat{\otimes}B[[G]]/E)\otimes\mathcal{M}(\mathrm{Spa}(A,A^+))\ar[r]\ar[r]\ar[r]  \ar[d]\ar[d]\ar[d]&0,\\
0 \ar[r]\ar[r]\ar[r] &\mathcal{K}(\mathrm{Spa}(C,C^+))   \ar[r]\ar[r]\ar[r] &\mathcal{F}(\mathrm{Spa}(C,C^+)) \ar[r]\ar[r]\ar[r] & \mathcal{M}(\mathrm{Spa}(C,C^+)) \ar[r]\ar[r]\ar[r] &0,\\
}
\]
$A\rightarrow C$ is some member in $\mathfrak{C}$. Here the first row is exact by $B[[G]]/E$-\'etale-pseudoflatness and the second is row is exact as well. The middle vertical arrow is an isomorphism while the rightmost vertical arrow is surjective as in \cite[Lemma 2.5.13]{8KL2}. Then we apply the same argument to the leftmost vertical arrow (with a new diagram), we will have the corresponding surjectivity of this leftmost vertical arrow. This will imply that the rightmost vertical arrow is injective as well. Then we have the corresponding isomorphism for the rightmost vertical arrow.
\end{proof}

\begin{theorem}\mbox{\bf{(After Kedlaya-Liu \cite[Theorem 2.5.14]{8KL2})}} \label{theorem3.23}
Taking global section will realize the equivalence between the following two categories: A. The category of all the pseudocoherent $\mathcal{O}_{\mathrm{Spa}(A,A^+)_{\text{\'et}}}\widehat{\otimes}B[[G]]$-sheaves; B. The category of all the $B[[G]]$-\'etale-stably pseudocoherent modules over $A\widehat{\otimes}B[[G]]$. 	
\end{theorem}

\begin{proof}
See \cite[Theorem 2.5.14]{8KL2}. We need to still apply \cite[Theorem 8.2.22]{8KL1}, as long as one considers instead in our situation, and \cite[Tag 03OD]{8SP}.\\	
\end{proof}

\subsection{Noncommutative Deformation over Quasi-Stein Spaces}

\indent In this subsection we consider the corresponding noncommutative deformation over the corresponding context in \cite[Chapter 2.6]{8KL2}.

\begin{setting}
Let $X$ be a corresponding quasi-Stein adic space over $\mathbb{Q}_p$  in the sense of \cite[Definition 2.6.2]{8KL2}. Recall that what is happening is that $X$ could be written as the corresponding direct limit of affinoids $X:=\varinjlim_i X_i$. 	
\end{setting}

\begin{lemma} \mbox{\bf{(After Kedlaya-Liu \cite[Lemma 2.6.3]{8KL2})}}
We now consider the corresponding rings $A_i:=\mathcal{O}_{X_i}$ for all $i=0,1,...$, and in our current situation we consider the corresponding rings $A_i\widehat{\otimes}B[[G]]$ (over $\mathbb{Q}_p$) for all $i=0,1,...$. And in our situation we consider the corresponding modules $M_i$ over $A_i\widehat{\otimes}B[[G]]$ for all $i=0,1,...$ with the same requirement as in \cite[Lemma 2.6.3]{8KL2} (namely those complete with respect to the natural topology). Suppose that we have bounded surjective map from $f_i:A_{i}\widehat{\otimes}B[[G]]\otimes_{A_{i+1}\widehat{\otimes}B[[G]]} M_{i+1}\rightarrow M_i,i=0,1,...$. Then we have first the density of the corresponding image of $\varprojlim_i M_i$ in each $M_i$ for any $i=0,1,2,...$. And we have as well the corresponding vanishing of $R^1\varprojlim_i M_i$.
\end{lemma}

\begin{proof}
This is the $B[[G]]$-relative version of the result in \cite[Lemma 2.6.3]{8KL2}. One needs to consider the corresponding truncation through some specific open subgroup $E$. Then see \cref{lemma2.24}.	
\end{proof}

\begin{proposition} \mbox{\bf{(After Kedlaya-Liu \cite[Lemma 2.6.4]{8KL2})}} In the same situation as above, suppose we have that the corresponding modules $M_i$ are basically $B[[G]]$-stably pseudocoherent over the rings $A_i$ for all $i=0,1,...$. Now we consider the situation where $f_i:A_i\widehat{\otimes}B[[G]]\otimes_{A_{i+1}\widehat{\otimes}B[[G]]}M_{i+1}\rightarrow M_i$ is an isomorphism. Then the conclusion in our situation is then that the corresponding projection from $\varprojlim M_i$ to $M_i$ for each $i=0,1,2,...$ is an isomorphism.

\end{proposition}

\begin{proof}
This is a $B[[G]]$-relative version of the \cite[Lemma 2.6.4]{8KL2}. We adapt the argument to our situation as in the following. First we choose some finite free covering $T$ of the limit $M$ such that we have for each $i$ the corresponding map $T_i\rightarrow M_i$ is surjective. Then we consider the index $j\geq i$ and set the kernel of the map from $T_j$ to $M_j$ to be $S_j$. By the direct analog of \cite[Lemma 2.5.6]{8KL2} we have that $A_i\widehat{\otimes}B[[G]]/E\otimes S_j\overset{}{\rightarrow}S_i$ realizes the isomorphism, and we have that the corresponding surjectivity of the corresponding map from $\varprojlim_i S_i$ projecting to the $S_i$. Then one could finish the proof by 5-lemma to the following commutative diagram as in \cite[Lemma 2.6.4]{8KL2}:
\[ \tiny
\xymatrix@C+0pc@R+6pc{
 &(A_i\widehat{\otimes}B[[G]]/E)\otimes \varprojlim_i S_i \ar[r]\ar[r]\ar[r] \ar[d]\ar[d]\ar[d] &(A_i\widehat{\otimes}B[[G]]/E)\otimes \varprojlim_i F_i \ar[r]\ar[r]\ar[r] \ar[d]\ar[d]\ar[d] &(A_i\widehat{\otimes}B[[G]]/E)\otimes \varprojlim_i M_i\ar[r]\ar[r]\ar[r]  \ar[d]\ar[d]\ar[d]&0,\\
0 \ar[r]\ar[r]\ar[r] &S_i   \ar[r]\ar[r]\ar[r] &F_i \ar[r]\ar[r]\ar[r] &M_i \ar[r]\ar[r]\ar[r] &0.\\
}
\]
 
\end{proof}

\begin{proposition}  \mbox{\bf{(After Kedlaya-Liu \cite[Theorem 2.6.5]{8KL2})}}
For any\\ quasi-compact adic affinoid space of $X$ which is denoted by $Y$, we have that the map $\mathcal{M}(X)\rightarrow \mathcal{M}(Y)$ is surjective for any $B[[G]]$-stably pseudocoherent sheaf $\mathcal{M}$ over the sheaf $\mathcal{O}_X\widehat{\otimes}B[[G]]$.	
\end{proposition}

\begin{proof}
This is just the corresponding corollary of the previous proposition.	
\end{proof}

\begin{proposition}  \mbox{\bf{(After Kedlaya-Liu \cite[Theorem 2.6.5]{8KL2})}}
We have that the stalk $\mathcal{M}_x$ is finitely generated over the stalk $\mathcal{O}_{X,x}$ for any $x\in X$ by $M(X)$, for any $B[[G]]$-stably pseudocoherent sheaf $\mathcal{M}$ over the sheaf $\mathcal{O}_X\widehat{\otimes}B[[G]]$.	
\end{proposition}

\begin{proof}
This is just the corresponding corollary of the proposition before the previous proposition.	
\end{proof}

\begin{proposition}  \mbox{\bf{(After Kedlaya-Liu \cite[Theorem 2.6.5]{8KL2})}}
For any\\ quasi-compact adic affinoid space of $X$ which is denoted by $Y$, we have that the corresponding vanishing of the corresponding sheaf cohomology groups $H^k(X,\mathcal{M})$ of $\mathcal{M}$ for higher $k>0$, for any $B[[G]]$-stably pseudocoherent sheaf $\mathcal{M}$ over the sheaf $\mathcal{O}_X\widehat{\otimes}B[[G]]$.	
\end{proposition}

\begin{proof}
See the proof of \cref{proposition2.28}. Again in our situation we only have to work over some quotient.

\end{proof}

\begin{corollary}  \mbox{\bf{(After Kedlaya-Liu \cite[Corollary 2.6.6]{8KL2})}}
The corresponding functor from the corresponding $B[[G]]$-deformed pseudocoherent sheaves over $X$ to the corresponding $B[[G]]$-stably by taking the corresponding global section is an exact functor.
\end{corollary}

\begin{corollary}  \mbox{\bf{(After Kedlaya-Liu \cite[Corollary 2.6.8]{8KL2})}}
Consider a particular $\mathcal{O}_X\widehat{\otimes}B[[G]]$-pseudocoherent sheaf $\mathcal{M}$ which is finite locally free throughout the whole space $X$. Then we have that the global section $\mathcal{M}(X)$ as $\mathcal{O}_X(X)\widehat{\otimes}B[[G]]$ left module admits the corresponding structures of finite projective structure if and only if we have the corresponding global section is finitely generated.
\end{corollary}

\begin{proof}
As in \cite[Corollary 2.6.8]{8KL2} one could find some global splitting through the local splittings.
\end{proof}

\indent We now consider the following $B[[G]]$-relative analog of \cite[Proposition 2.6.17]{8KL2}:

\begin{theorem} \mbox{\bf{(After Kedlaya-Liu \cite[Proposition 2.6.17]{8KL2})}} \label{theorem3.32}
Consider the following two statements for a particular $\mathcal{O}_X\widehat{\otimes}B[[G]]$-pseudocoherent sheaf $\mathcal{M}$. First is that one can find finite many generators (the number is up to some uniform fixed integer $n\geq 0$) for each section of $\mathcal{M}(X_i)$ for each $i=1,2,...$. The second statement is that the global section $\mathcal{M}(X)$ is just finitely generated. Then in our situation the two statement is equivalent if we have that the corresponding space $X$ admits a $m$-uniformed covering in the exact same sense of \cite[Proposition 2.6.17]{8KL2}. 	
\end{theorem}

\begin{proof}
See the proof of \cref{theorem2.31}.
\end{proof}

%%\newpage 

\newpage\section{Foundations on Noncommutative Descent for Adic Spectra in Limit of Fr\'echet Case}

\subsection{Noncommutative Pseudocoherence in Analytic Topology}

\indent We now establish some foundations in the noncommutative setting on glueing noncommutative pseudocoherent modules after \cite[Chapter 2]{8KL2}. But we remind the readers that we will fix the base space:

\begin{setting} \label{setting4.1}
Consider a corresponding sheafy Banach adic uniform algebra $(A,A^+)$ over $\mathbb{Q}_p$ or $\mathbb{F}_p((t))$, we consider the base space $\mathrm{Spa}(A,A^+)$. And now we will consider a further noncommutative Banach algebra $(B,B^+)$ over $\mathbb{Q}_p$ or $\mathbb{F}_p((t))$ which could be written as the following injective limit:
\begin{displaymath}
B=\varinjlim_{h} B_h,	
\end{displaymath}
where each $B_h$ is Banach (or more general Fr\'echet although we do not consider such generality). 	
\end{setting}

\begin{example}
There are many interesting models in our current context, for instance the period ring $B_e$ as considered by Berger in \cite{8Ber1} and the corresponding Robba rings in the full setting as in \cite{8KL2} (note that the latter is really indeed ind-Fr\'echet).	
\end{example}

\begin{definition} \mbox{\bf{(After Kedlaya-Liu \cite[Definition 2.4.1]{8KL2})}}
For any left $A\widehat{\otimes}B$-module $M$, we call it $m$-$B$-stably pseudocoherent if we have that it is $m$-$B$-pseudocoherent, complete for the natural topology and for any morphism $A\rightarrow A'$ which is the corresponding rational localization, the base change of $M$ to $A'\widehat{\otimes} B$ is complete for the natural LF topology as the corresponding left $A'\widehat{\otimes} B$-module.	As in \cite[Definition 2.4.1]{8KL2} we call that the corresponding left $A\widehat{\otimes}B$-module $M$ just $B$-stably pseudocoherent if we have that it is simply just $\infty$-$B$-stably pseudocoherent.
\end{definition}

\begin{definition}\mbox{\bf{(After Kedlaya-Liu \cite[Definition 2.4.4]{8KL2})}}
For any Banach left $A\widehat{\otimes}B$-module $M$, we call it is $m$-$B$-pseudoflat if for any right $A\widehat{\otimes}B$-module $M'$ $m$-$B$-stably pseudocoherent we have $\mathrm{Tor}_1^{A\widehat{\otimes}B}(M',M)=0$. For any Banach right $A\widehat{\otimes}B$-module $M$, we call it is $m$-$B$-pseudoflat if for any left $A\widehat{\otimes}B$-module $M'$ $m$-$B$-stably pseudocoherent we have $\mathrm{Tor}_1^{A\widehat{\otimes}B}(M,M')=0$.
\end{definition}

\begin{definition}\mbox{\bf{(After Kedlaya-Liu \cite[Definition 2.4.6]{8KL2})}} We consider the corresponding notion of the corresponding pro-projective module. We define over $A$ the corresponding $B$-pro-projective module $M$ to be a corresponding left $A\widehat{\otimes}B$-module such that one could find a filtered sequence of projectors such that in the sense of taking the prolimit we have the corresponding projector will converge to any chosen element in $M$. Here we assume the module is complete with respect to natural topology, and we assume the projectors are $A\widehat{\otimes}\varinjlim_h B_h$-linear and we assume that the corresponding image of the projectors are modules which are also finitely generated and projective.
	
\end{definition}

\begin{lemma}\mbox{\bf{(After Kedlaya-Liu \cite[Lemma 2.4.7]{8KL2})}}
Suppose we have over the space $\mathrm{Spa}(A,A^+)$ a corresponding $B$-pro-projective left module $M$. And suppose that we have a $2$-$B$-pseudoflat right module $C$ over $A$. And we assume that $C$ is complete for the LF topology. Then we have that the corresponding product $C{\otimes}_{A\widehat{\otimes}B}M$ is then complete under the corresponding natural topology. And moreover we have that in our situation:
\begin{displaymath}
\mathrm{Tor}_1^{A\widehat{\otimes}B}(C,M)=0.	
\end{displaymath}

\end{lemma}

\begin{proof}
This will be just a noncommutative and LF version of the corresponding  \cite[Lemma 2.4.7]{8KL2}. We adapt the corresponding argument in \cite[Lemma 2.4.7]{8KL2} to our situation. What we are going to consider in this case is then first consider the following presentation:
\[
\xymatrix@C+0pc@R+0pc{
M'  \ar[r]\ar[r]\ar[r] &(A\widehat{\otimes}\varinjlim_h B_h)^k \ar[r]\ar[r]\ar[r] & M\ar[r]\ar[r]\ar[r] &0,
}
\]
where the left module $M'$ is finitely presented. Then we consider the following commutative diagram:
\[
\xymatrix@C+0pc@R+3pc{
&C\otimes_{A\widehat{\otimes}\varinjlim_h B_h}M'  \ar[r]\ar[r]\ar[r]  \ar[d]\ar[d]\ar[d] &C\otimes_{A\widehat{\otimes}\varinjlim_h B_h}(A\widehat{\otimes}\varinjlim_h B_h)^k \ar[r]\ar[r]\ar[r] \ar[d]\ar[d]\ar[d] & C\otimes_{A\widehat{\otimes}\varinjlim_h B_h}M\ar[r]\ar[r]\ar[r] \ar[d]\ar[d]\ar[d] &0,\\
0 \ar[r]\ar[r]\ar[r] &C\widehat{\otimes}_{A\widehat{\otimes}\varinjlim_h B_h}M'  \ar[r]\ar[r]\ar[r] &C\widehat{\otimes}_{A\widehat{\otimes}B}(A\widehat{\otimes}\varinjlim_h B_h)^k \ar[r]\ar[r]\ar[r] & C\widehat{\otimes}_{A\widehat{\otimes}\varinjlim_h B_h}M\ar[r]\ar[r]\ar[r] &0.
}
\]\\
The first row is exact as in \cite[Lemma 2.4.7]{8KL2}, and we have the corresponding second row is also exact since we have that by hypothesis the corresponding module $M$ is $B$-pro-projective.  
Then as in \cite[Lemma 2.4.7]{8KL2} the results follow.	
\end{proof}

\begin{proposition} \mbox{\bf{(After Kedlaya-Liu \cite[Corollary 2.4.8]{8KL2})}}
Over the space \\$\mathrm{Spa}(A,A^+)$, we have any $B$-pro-projective left module is $2$-$B$-pseudoflat.
\end{proposition}

\begin{proof}
This is a direct consequence of the previous lemma.	
\end{proof}

\begin{lemma} \mbox{\bf{(After Kedlaya-Liu \cite[Corollary 2.4.9]{8KL2})}}
Keep the notation above, suppose we are working over the adic space $
\mathrm{Spa}(A,A^+)$. Suppose now the left $A\widehat{\otimes}B$-module is finitely generated, then we have that the following natural maps:
\begin{align}
A\{T\}\widehat{\otimes}B\otimes_{A\widehat{\otimes}B}M 	\rightarrow M\{T\},\\\
A\{T,T^{-1}\}\widehat{\otimes}B \otimes_{A\widehat{\otimes}B} M 	\rightarrow M\{T,T^{-1}\}	
\end{align}
are surjective. Suppose now the left $A\widehat{\otimes}B$-module is finitely presented, then we have that the following natural maps:
\begin{align}
A\{T\}\widehat{\otimes}B\otimes_{A\widehat{\otimes}B}M 	\rightarrow M\{T\},\\\
A\{T,T^{-1}\}\widehat{\otimes}B \otimes_{A\widehat{\otimes}B} M 	\rightarrow M\{T,T^{-1}\}	
\end{align}
are bijective. Here all the modules are assumed to be complete for the LF topology.

\end{lemma}

\begin{proof}
This is the corresponding consequence of the previous lemma.
\end{proof}

\begin{proposition} \mbox{\bf{(After Kedlaya-Liu \cite[Lemma 2.4.12]{8KL2})}}
We keep the corresponding notations in \cite[Lemma 2.4.10]{8KL2}. Then in our current sense we have that the corresponding morphism $A\rightarrow B_2$ is then $2$-$B$-pseudoflat. 
\end{proposition}

\begin{proof}
The corresponding proof could be made parallel to the corresponding proof of \cite[Lemma 2.4.12]{8KL2}. To be more precise we consider the following commutative diagram from the corresponding chosen exact sequence (which is just the analog of the corresponding one in \cite[Lemma 2.4.12]{8KL2}, namely as below we choose arbitrary short exact sequence with $P$ as a left $A\widehat{\otimes}\varinjlim_h B_h$-module which is $2$-$\varinjlim_h B_h$-stably pseudocoherent):
\[
\xymatrix@C+0pc@R+0pc{
0   \ar[r]\ar[r]\ar[r] &M \ar[r]\ar[r]\ar[r] &N \ar[r]\ar[r]\ar[r] & P \ar[r]\ar[r]\ar[r] &0,
}
\]
which then induces the following corresponding big commutative diagram:
\[\tiny
\xymatrix@C+0pc@R+4pc{
& &0 \ar[d]\ar[d]\ar[d] &0 \ar[d]\ar[d]\ar[d] &,\\
0   \ar[r]\ar[r]\ar[r]  &{A\{T\}\widehat{\otimes}\varinjlim_h B_h}\otimes_{A\widehat{\otimes}\varinjlim_h B_h} M \ar[r]\ar[r]\ar[r] \ar[d]^{1-fT}\ar[d]\ar[d] &{A\{T\}\widehat{\otimes}\varinjlim_h B_h}\otimes_{A\widehat{\otimes}\varinjlim_h B_h} N \ar[r]\ar[r]\ar[r] \ar[d]^{1-fT}\ar[d]\ar[d] &{A\{T\}\widehat{\otimes}\varinjlim_h B_h} \otimes_{A\widehat{\otimes}\varinjlim_h B_h}P \ar[r]\ar[r]\ar[r] \ar[d]^{1-fT}\ar[d]\ar[d] &0,\\
0   \ar[r]\ar[r]\ar[r] &{A\{T\}\widehat{\otimes}\varinjlim_h B_h}\otimes_{A\widehat{\otimes}\varinjlim_h B_h}M  \ar[r]\ar[r]\ar[r] \ar[d]\ar[d]\ar[d] &{A\{T\}\widehat{\otimes}\varinjlim_h B_h} \otimes_{A\widehat{\otimes}\varinjlim_h B_h}N\ar[r]\ar[r]\ar[r] \ar[d]\ar[d]\ar[d] & {A\{T\}\widehat{\otimes}\varinjlim_h B_h}\otimes_{A\widehat{\otimes}\varinjlim_h B_h}P \ar[r]\ar[r]\ar[r] \ar[d]\ar[d]\ar[d] &0,\\
0  \ar[r]^?\ar[r]\ar[r] &{B_2\widehat{\otimes}\varinjlim_h B_h}\otimes_{A\widehat{\otimes}\varinjlim_h B_h} M \ar[r]\ar[r]\ar[r] \ar[d]\ar[d]\ar[d] &{B_2\widehat{\otimes}\varinjlim_h B_h}\otimes_{A\widehat{\otimes}\varinjlim_h B_h}N \ar[r]\ar[r]\ar[r] \ar[d]\ar[d]\ar[d] & {B_2\widehat{\otimes}\varinjlim_h B_h}\otimes_{A\widehat{\otimes}\varinjlim_h B_h} P \ar[r]\ar[r]\ar[r] \ar[d]\ar[d]\ar[d] &0,\\
&0&0&0
}
\]
with the notations in \cite[Lemma 2.4.10]{8KL2} where we actually have the corresponding exactness along the horizontal direction at the corner around $M\otimes_{A\widehat{\otimes}\varinjlim_h B_h} {B_2\widehat{\otimes}\varinjlim_h B_h}$ marked with $?$, by diagram chasing.
\end{proof}

\begin{proposition} \mbox{\bf{(After Kedlaya-Liu \cite[Lemma 2.4.13]{8KL2})}} \label{8proposition4.9}
We keep the corresponding notations in \cite[Lemma 2.4.10]{8KL2}. Then in our current sense we have that the corresponding morphism $A\rightarrow B_1$ is then $2$-$B$-pseudoflat. And in our current sense we have that the corresponding morphism $A\rightarrow B_{1,2}$ is then $2$-$B$-pseudoflat. 
\end{proposition}

\begin{proof}
The statement for $A\rightarrow B_{1,2}$ could be proved by consider the composition $A\rightarrow B_2\rightarrow B_{12}$ where the $2$-$B$-pseudoflatness could be proved as in \cite[Lemma 2.4.13]{8KL2} by inverting the corresponding variable. For $A\rightarrow B_{1}$,
the corresponding proof could be made parallel to the corresponding proof of \cite[Lemma 2.4.13]{8KL2}. To be more precise we consider the following commutative diagram from the corresponding chosen exact sequence (which is just the analog of the corresponding one in \cite[Remark 2.4.5]{8KL2}, namely as below we choose arbitrary short exact sequence with $P$ as a left $A\widehat{\otimes}B$-module which is $2$-$B$-stably pseudocoherent):
\[
\xymatrix@C+0pc@R+0pc{
0   \ar[r]\ar[r]\ar[r] &M \ar[r]\ar[r]\ar[r] &N \ar[r]\ar[r]\ar[r] & P \ar[r]\ar[r]\ar[r] &0,
}
\]
which then induces the following corresponding big commutative diagram:
\[\tiny
\xymatrix@C-0.4pc@R+5pc{
& &0 \ar[d]\ar[d]\ar[d] &0 \ar[d]\ar[d]\ar[d] &,\\
0   \ar[r]\ar[r]\ar[r]  &M \ar[r]\ar[r]\ar[r] \ar[d]\ar[d]\ar[d] &N \ar[r]\ar[r]\ar[r] \ar[d]\ar[d]\ar[d] & P \ar[r]\ar[r]\ar[r] \ar[d]\ar[d]\ar[d] &0,\\
0   \ar[r]^?\ar[r]\ar[r] &({B_1\widehat{\otimes}B_\infty}\bigoplus {B_2\widehat{\otimes}B_\infty}) \otimes_{A\widehat{\otimes}B_\infty}M  \ar[r]\ar[r]\ar[r] \ar[d]\ar[d]\ar[d] &({B_1\widehat{\otimes}B_\infty}\bigoplus {B_2\widehat{\otimes}B_\infty})\otimes_{A\widehat{\otimes}B_\infty}N  \ar[r]\ar[r]\ar[r] \ar[d]\ar[d]\ar[d] &({B_1\widehat{\otimes}B_\infty}\bigoplus {B_2\widehat{\otimes}B_\infty})\otimes_{A\widehat{\otimes}B_\infty}P\ar[r]\ar[r]\ar[r] \ar[d]\ar[d]\ar[d] &0,\\
0  \ar[r]\ar[r]\ar[r] &{B_{12}\widehat{\otimes}B_\infty}\otimes_{A\widehat{\otimes}B_\infty}M \ar[r]\ar[r]\ar[r] \ar[d]\ar[d]\ar[d] &{B_{12}\widehat{\otimes}B_\infty}\otimes_{A\widehat{\otimes}B_\infty}N \ar[r]\ar[r]\ar[r] \ar[d]\ar[d]\ar[d] &{B_{12}\widehat{\otimes}B_\infty} \otimes_{A\widehat{\otimes}B_\infty} P \ar[r]\ar[r]\ar[r] \ar[d]\ar[d]\ar[d] &0,\\
&0&0&0
}
\]
with the notations in \cite[Lemma 2.4.10]{8KL2} where we actually have the corresponding exactness along the horizontal direction at the corner around
\begin{center}
 $M\otimes_{A\widehat{\otimes}\varinjlim_h B_h} ({B_1\widehat{\otimes}\varinjlim_h B_h}\bigoplus {B_2\widehat{\otimes}\varinjlim_h B_h})$
\end{center} 
marked with $?$, by diagram chasing. Here $B_\infty:=\varinjlim_h B_h$.
\end{proof}

\indent After these foundational results as in \cite{8KL2} we have the following proposition which is the corresponding noncommutative LF generalization of the corresponding result established in \cite[Theorem 2.4.15]{8KL2}.

\begin{proposition} \mbox{\bf{(After Kedlaya-Liu \cite[Theorem 2.4.15]{8KL2})}} \label{proposition4.10}
In our current context we have that for any rational localization $\mathrm{Spa}(A',A^{',+})\rightarrow \mathrm{Spa}(A,A^+)$ we have that along this base change the corresponding $B$-stably pseudocoherence is preserved.
	
\end{proposition}

\indent Then we have the corresponding Tate's acyclicity in the noncommutative deformed setting:

\begin{theorem}\mbox{\bf{(After Kedlaya-Liu \cite[Theorem 2.5.1]{8KL2})}} \label{theorem4.11} Now suppose we have in our corresponding \cref{setting4.1} a corresponding $B$-stably pseudocoherent module $M$. Then we consider the corresponding assignment such that for any $U\subset \mathrm{Spa}(A,A^+)$ we define $\widetilde{M}(U)$ as in the following:
\begin{align}
\widetilde{M}(U):=\varprojlim_{\mathrm{Spa}(S,S^+)\subset U,\mathrm{rational}} S\widehat{\otimes}B\otimes_{A\widehat{\otimes}B}M.	
\end{align}
Then we have that for any $\mathfrak{B}$ which is a rational covering of $U=\mathrm{Spa}(S,S^+)\subset \mathrm{Spa}(A,A^+)$ (certainly this $U$ is also assumed to be rational) we have that the vanishing of the following two cohomology groups:
\begin{align}
H^i(U,\widetilde{M}), \check{H}^i(U:\mathfrak{B},\widetilde{M})
\end{align}
for any $i>0$. When concentrating at the degree zero we have:
\begin{align}
H^0(U,\widetilde{M})=S\widehat{\otimes}B\otimes_{A\widehat{\otimes}B}M, \check{H}^0(U:\mathfrak{B},\widetilde{M})=S\widehat{\otimes}B\otimes_{A\widehat{\otimes}B}M.
\end{align}
	
\end{theorem}

\begin{proof}
By \cite[Propositions 2.4.20-2.4.21]{8KL1}, we could then finish proof as in \cite[Theorem 2.5.1]{8KL2} by using our previous \cref{proposition4.10} and \cref{8proposition4.9} as above.	
\end{proof}

\indent We now consider the corresponding noncommutative LF deformed version of Kiehl's glueing properties for stably pseudocoherent modules after \cite{8KL2}:

\begin{definition} \mbox{\bf{(After Kedlaya-Liu \cite[Definition 2.5.3]{8KL2})}}
Consider in our context (over $\mathrm{Spa}(A,A^+)$) the corresponding sheaves $\mathcal{O}\widehat{\otimes}B$, we will then define the corresponding pseudocoherent sheaves over $\mathrm{Spa}(A,A^+)$ to be those locally defined by attaching stably-pseudocoherent modules over the section. 
\end{definition}

\begin{lemma} \mbox{\bf{(After Kedlaya-Liu \cite[Lemma 2.5.4]{8KL2})}}\label{lemma4.13}
	Consider the corresponding notations in \cite[Lemma 2.4.10]{8KL2}, we have the corresponding morphism $A\rightarrow B_1\bigoplus B_2$. Then we have that this morphism is an descent morphism effective for the corresponding $B$-stably pseodocoherent Banach modules. 
\end{lemma}

\begin{proof}
See the proof of \cref{8lemma2.13}.
\end{proof}

\begin{theorem}\mbox{\bf{(After Kedlaya-Liu \cite[Theorem 2.5.5]{8KL2})}} \label{theorem4.14}
Taking global section will realize the equivalence between the following two categories: A. The category of all the pseudocoherent $\mathcal{O}\widehat{\otimes}B$-sheaves; B. The category of all the $B$-stably pseudocoherent modules over $A\widehat{\otimes}B$. 	
\end{theorem}

\begin{proof}
See \cite[Theorem 2.5.5]{8KL2}. We need to still apply \cite[Proposition 2.4.20]{8KL1}, as long as one considers instead in our situation \cref{lemma4.13} and \cref{theorem4.11}.\\	
\end{proof}

\subsection{Noncommutative Pseudocoherence in \'Etale Topology}

\indent We consider the extension of the corresponding discussion in the previous subsection to the \'etale topology. At this moment we consider the following same context on the geometric level as in \cite{8KL2}:

\begin{setting} \mbox{\bf{(After Kedlaya-Liu \cite[Hypothesis 2.5.8]{8KL2})}} \label{setting4.15}
As in the previous subsection we consider now the corresponding geometric setting namely an adic space $\mathrm{Spa}(A,A^+)$ where $A$ is assumed to be sheafy. Then we consider $B:=\varinjlim_h B_h$ as above. And we consider the corresponding \'etale site $\mathrm{Spa}(A,A^+)_\text{\'et}$. And by keeping the corresponding setting in the geometry as in \cite[Hypothesis 2.5.8]{8KL2}, we assume that there is a stable basis $\mathfrak{B}$ containing the space $\mathrm{Spa}(A,A^+)$ itself.	
\end{setting}

\begin{definition} \mbox{\bf{(After Kedlaya-Liu \cite[Definition 2.5.9]{8KL2})}}
For any left $A\widehat{\otimes}B$-module $M$, we call it $m$-$B$-\'etale stably pseudocoherent with respect to $\mathfrak{B}$ if we have that it is $m$-$B$-pseudocoherent, complete for the natural topology and for any rational localization $A\rightarrow A'$, the base change of $M$ to $A'\widehat{\otimes} B$ is complete for the natural LF topology as the corresponding left $A'\widehat{\otimes} B$-module.	As in \cite[Definition 2.4.1]{8KL2} we call that the corresponding left $A\widehat{\otimes}B$-module $M$ just $B$-stably pseudocoherent if we have that it is simply just $\infty$-$B$-stably pseudocoherent.
\end{definition}

\begin{definition}\mbox{\bf{(After Kedlaya-Liu \cite[Below Definition 2.5.9]{8KL2})}}
For any Banach left $A\widehat{\otimes}B$-module $M$, we call it is $m$-$B$-\'etale-pseudoflat if for any right $A\widehat{\otimes}B$-module $M'$ $m$-$B$-\'etale-stably pseudocoherent we have $\mathrm{Tor}_1^{A\widehat{\otimes}B}(M',M)=0$. For any Banach right $A\widehat{\otimes}B$-module $M$, we call it is $m$-$B$-\'etale-pseudoflat if for any left $A\widehat{\otimes}B$-module $M'$ $m$-$B$-\'etale-stably pseudocoherent we have $\mathrm{Tor}_1^{A\widehat{\otimes}B}(M,M')=0$.
\end{definition}

\indent The following proposition holds in our current setting.

\begin{proposition} \mbox{\bf{(After Kedlaya-Liu \cite[Lemma 2.5.10]{8KL2})}} \label{proposition4.18}
One can actually find another basis $\mathfrak{C}$ in $\mathfrak{B}$ such that any morphism in $\mathfrak{C}$ could be $2$-$B$-\'etale-pseudoflat with respect to either $\mathfrak{C}$ or $\mathfrak{B}$.	
\end{proposition}

\begin{proof}
We follow the proof of \cite[Lemma 2.5.10]{8KL2} in our current $B$-relative situation, however not that much needs to proof by relying on the proof of \cite[Lemma 2.5.10]{8KL2}. To be more precise the corresponding selection of the new basis $\mathfrak{C}$ comes from including all the morphism made up of some composition of rational localization and the corresponding finite \'etale ones. For these we have shown above the corresponding 2-$B$-\'etale pseudoflatness. Then the rest will be pure geometric for the analytic spaces which are just the same as the situation of \cite[Lemma 2.5.10]{8KL2}, therefore we just omit the corresponding argument, see \cite[Lemma 2.5.10]{8KL2}.	
\end{proof}

\indent Then we have the corresponding Tate's acyclicity in the noncommutative deformed setting in \'etale topology (here fix $\mathfrak{C}$ as above):

\begin{theorem}\mbox{\bf{(After Kedlaya-Liu \cite[Theorem 2.5.11]{8KL2})}} \label{theorem4.19}Now suppose we \\have in our corresponding \cref{setting4.1} a corresponding $B$-\'etale-stably pseudocoherent module $M$. Then we consider the corresponding assignment such that for any $U\subset \mathrm{Spa}(A,A^+)$ we define $\widetilde{M}(U)$ as in the following:
\begin{align}
\widetilde{M}(U):=\varprojlim_{\mathrm{Spa}(S,S^+)\subset U,\in \mathfrak{C}} S\widehat{\otimes}B\otimes_{A\widehat{\otimes}B}M.	
\end{align}
Then we have that for any $\mathfrak{B}$ which is a covering of $U=\mathrm{Spa}(S,S^+)\subset \mathrm{Spa}(A,A^+)$ (certainly this $U$ is also assumed to be in $\mathfrak{C}$, and we assume this covering is formed by using the corresponding members in $\mathfrak{C}$) we have that the vanishing of the following two cohomology groups:
\begin{align}
H^i(U,\widetilde{M}), \check{H}^i(U:\mathfrak{B},\widetilde{M})
\end{align}
for any $i>0$. When concentrating at the degree zero we have:
\begin{align}
H^0(U,\widetilde{M})=S\widehat{\otimes}B\otimes_{A\widehat{\otimes}B}M, \check{H}^0(U:\mathfrak{B},\widetilde{M})=S\widehat{\otimes}B\otimes_{A\widehat{\otimes}B}M.
\end{align}
	
\end{theorem}

\begin{proof}
By \cite[Proposition 8.2.21]{8KL1}, we could then finish proof as in \cite[Theorem 2.5.11]{8KL2} by using our previous \cref{proposition4.18} and the corresponding faithfully flat descent as in the situation of \cite[Theorem 2.5.11]{8KL2}.	
\end{proof}

\indent We now consider the corresponding noncommutative LF deformed version of Kiehl's glueing properties for stably pseudocoherent modules after \cite{8KL2}:

\begin{definition} \mbox{\bf{(After Kedlaya-Liu \cite[Definition 2.5.12]{8KL2})}}
Consider in our context (over $\mathrm{Spa}(A,A^+)_{\text{\'et}}$) the corresponding sheaves $\mathcal{O}_{\mathrm{Spa}(A,A^+)_{\text{\'et}}}\widehat{\otimes}B$, we will then define the corresponding pseudocoherent sheaves over $\mathrm{Spa}(A,A^+)_{\text{\'et}}$ to be those locally defined by attaching \'etale-stably-pseudocoherent modules over the section. 
\end{definition}

\begin{lemma} \mbox{\bf{(After Kedlaya-Liu \cite[Lemma 2.5.13]{8KL2})}}
	Consider the corresponding notations in \cite[Lemma 2.4.10]{8KL2}, we have the corresponding morphism $A\rightarrow B_1\bigoplus B_2$. Then we have that this morphism is an descent morphism effective for the corresponding $B$-\'etale-stably pseodocoherent Banach modules. 
\end{lemma}

\begin{proof}
See the proof of \cref{lemma2.21}.

\end{proof}

\begin{theorem}\mbox{\bf{(After Kedlaya-Liu \cite[Theorem 2.5.14]{8KL2})}} \label{theorem4.22}
Taking global section will realize the equivalence between the following two categories: A. The category of all the pseudocoherent $\mathcal{O}_{\mathrm{Spa}(A,A^+)_{\text{\'et}}}\widehat{\otimes}B$-sheaves; B. The category of all the $B$-\'etale-stably pseudocoherent modules over $A\widehat{\otimes}B$. 	
\end{theorem}

\begin{proof}
See \cite[Theorem 2.5.14]{8KL2}. We need to still apply \cite[Theorem 8.2.22]{8KL1}, as long as one considers instead in our situation, and \cite[Tag 03OD]{8SP}.\\	
\end{proof}

\subsection{Noncommutative Deformation over Quasi-Stein Spaces}

\indent In this subsection we consider the corresponding noncommutative deformation over the corresponding context in \cite[Chapter 2.6]{8KL2}.

\begin{setting}
Let $X$ be a corresponding quasi-Stein adic space over $\mathbb{Q}_p$ or $\mathbb{F}_p((t))$ in the sense of \cite[Definition 2.6.2]{8KL2}. Recall that what is happening is that $X$ could be written as the corresponding direct limit of affinoids $X:=\varinjlim_i X_i$. 	
\end{setting}

\begin{lemma} \mbox{\bf{(After Kedlaya-Liu \cite[Lemma 2.6.3]{8KL2})}}
We now consider the corresponding rings $A_i:=\mathcal{O}_{X_i}$ for all $i=0,1,...$, and in our current situation we consider the corresponding rings $A_i\widehat{\otimes}B$ (over $\mathbb{Q}_p$ or $\mathbb{F}_p((t))$) for all $i=0,1,...$. And in our situation we consider the corresponding modules $M_i$ over $A_i\widehat{\otimes}B$ for all $i=0,1,...$ with the same requirement as in \cite[Lemma 2.6.3]{8KL2} (namely those complete with respect to the natural topology which is LF topology in our current situation). Suppose that we have bounded surjective map from $f_i:A_{i}\widehat{\otimes}B\otimes_{A_{i+1}\widehat{\otimes}B} M_{i+1}\rightarrow M_i,i=0,1,...$. Then we have first the density of the corresponding image of $\varprojlim_i M_i$ in each $M_i$ for any $i=0,1,2,...$. And we have as well the corresponding vanishing of $R^1\varprojlim_i M_i$.
\end{lemma}

\begin{proof}
This is the $\varinjlim_h B_h$-relative version of the result in \cite[Lemma 2.6.3]{8KL2}. We need to then use the corresponding family of Banach norms parametrized by the index $h$. For the first statement we just choose sequence of Banach norms on all the corresponding modules for all $i=0,1,...$ and $h$ such that we have $\|f_i(x_{i+1})\|^h_i\leq 1/2\|x_{i+1}\|^h_{i+1}$ for any $x_{i+1}\in M_{i+1}$. Then for any $x_i\in M_i$ and any $\delta>0$, we now consider for any $j\geq 1$ and $\forall h$ the corresponding $x_{i+j}$ such that we have $\|x_{i+j}-f_{i+j+1}(x_{i+j+1})\|^h_{i+j}\leq \delta$. Then the sequence $x_{i+j+k},k=0,1,...$ will converge to some well-defined $y_{i+j}$ with in our situation the corresponding $y_{i}=f_{i}(y_{i+1})$. We then have $\|x_i-y_i\|^h_i\leq \delta$. This will prove the first statement. For the second statement as in \cite[Lemma 2.6.3]{8KL2} we form the product $M_0\times M_1\times M_2\times...$ and the consider the induced map $F$ from $M_{i+1}\rightarrow M_i$, and consider the corresponding cokernel of the map $1-F$ since this is just the corresponding limit we are considering. Then to show that the cokernel is zero we just look at the corresponding cokernel of the corresponding map on the corresponding completed direct summand which will project to the original one. But then we will have $\|f_i(v)\|^h_i\leq 1/2 \|v\|^h_i$, which produces an inverse to $1-F$ which will basically finish the proof for the second statement. 	
\end{proof}

\begin{proposition} \mbox{\bf{(After Kedlaya-Liu \cite[Lemma 2.6.4]{8KL2})}} In the same situation as above, suppose we have that the corresponding modules $M_i$ are basically $B$-stably pseudocoherent over the rings $A_i$ for all $i=0,1,...$. Now we consider the situation where $f_i:A_i\widehat{\otimes}B\otimes_{A_{i+1}\widehat{\otimes}B}M_{i+1}\rightarrow M_i$ is an isomorphism. Then the conclusion in our situation is then that the corresponding projection from $\varprojlim M_i$ to $M_i$ for each $i=0,1,2,...$.

\end{proposition}

\begin{proof}
This is a $\varinjlim_h B_h$-relative version of the \cite[Lemma 2.6.4]{8KL2}. We adapt the argument to our situation as in the following. First we choose some finite free covering $T$ of the limit $M$ such that we have for each $i$ the corresponding map $T_i\rightarrow M_i$ is surjective. Then we consider the index $j\geq i$ and set the kernel of the map from $T_j$ to $M_j$ to be $S_j$. By the direct analog of \cite[Lemma 2.5.6]{8KL2} we have that $A_i\widehat{\otimes}\varinjlim_h B_h\otimes S_j\overset{}{\rightarrow}S_i$ realizes the isomorphism, and we have that the corresponding surjectivity of the corresponding map from $\varprojlim_i S_i$ projecting to the $S_i$. Then one could finish the proof by 5-lemma to the following commutative diagram as in \cite[Lemma 2.6.4]{8KL2}:
\[ \tiny
\xymatrix@C+0pc@R+6pc{
 &(A_i\widehat{\otimes}\varinjlim_h B_h)\otimes \varprojlim_i S_i \ar[r]\ar[r]\ar[r] \ar[d]\ar[d]\ar[d] &(A_i\widehat{\otimes}\varinjlim_h B_h)\otimes \varprojlim_i F_i \ar[r]\ar[r]\ar[r] \ar[d]\ar[d]\ar[d] &(A_i\widehat{\otimes}\varinjlim_h B_h)\otimes \varprojlim_i M_i\ar[r]\ar[r]\ar[r]  \ar[d]\ar[d]\ar[d]&0,\\
0 \ar[r]\ar[r]\ar[r] &S_i   \ar[r]\ar[r]\ar[r] &F_i \ar[r]\ar[r]\ar[r] &M_i \ar[r]\ar[r]\ar[r] &0.\\
}
\]
 
\end{proof}

\begin{proposition}  \mbox{\bf{(After Kedlaya-Liu \cite[Theorem 2.6.5]{8KL2})}}
For any\\ quasi-compact adic affinoid space of $X$ which is denoted by $Y$, we have that the map $\mathcal{M}(X)\rightarrow \mathcal{M}(Y)$ is surjective for any $B$-stably pseudocoherent sheaf $\mathcal{M}$ over the sheaf $\mathcal{O}_X\widehat{\otimes}B$.	
\end{proposition}

\begin{proof}
This is just the corresponding corollary of the previous proposition.	
\end{proof}

\begin{proposition}  \mbox{\bf{(After Kedlaya-Liu \cite[Theorem 2.6.5]{8KL2})}}
We have that the stalk $\mathcal{M}_x$ is finitely generated over the stalk $\mathcal{O}_{X,x}$ for any $x\in X$ by $M(X)$, for any $B$-stably pseudocoherent sheaf $\mathcal{M}$ over the sheaf $\mathcal{O}_X\widehat{\otimes}B$.	
\end{proposition}

\begin{proof}
This is just the corresponding corollary of the proposition before the previous proposition.	
\end{proof}

\begin{proposition}  \mbox{\bf{(After Kedlaya-Liu \cite[Theorem 2.6.5]{8KL2})}}
For any\\ quasi-compact adic affinoid space of $X$ which is denoted by $Y$, we have that the corresponding vanishing of the corresponding sheaf cohomology groups $H^k(X,\mathcal{M})$ of $\mathcal{M}$ for higher $k>0$, for any $B$-stably pseudocoherent sheaf $\mathcal{M}$ over the sheaf $\mathcal{O}_X\widehat{\otimes}B$.	
\end{proposition}

\begin{proof}
See the proof of \cref{proposition2.28}.
\end{proof}

\begin{corollary}  \mbox{\bf{(After Kedlaya-Liu \cite[Corollary 2.6.6]{8KL2})}}
The corresponding functor from the corresponding $B$-deformed pseudocoherent sheaves over $X$ to the corresponding $B$-stably by taking the corresponding global section is an exact functor.
\end{corollary}

\begin{corollary}  \mbox{\bf{(After Kedlaya-Liu \cite[Corollary 2.6.8]{8KL2})}}
Consider a particular $\mathcal{O}_X\widehat{\otimes}B$-pseudocoherent sheaf $\mathcal{M}$ which is finite locally free throughout the whole space $X$. Then we have that the global section $\mathcal{M}(X)$ as $\mathcal{O}_X(X)\widehat{\otimes}B$ left module admits the corresponding structures of finite projective structure if and only if we have the corresponding global section is finitely generated.
\end{corollary}

\begin{proof}
As in \cite[Corollary 2.6.8]{8KL2} one could find some global splitting through the local splittings.
\end{proof}

\indent We now consider the following $B$-relative analog of \cite[Proposition 2.6.17]{8KL2}:

\begin{theorem} \mbox{\bf{(After Kedlaya-Liu \cite[Proposition 2.6.17]{8KL2})}} \label{theorem4.31}
Consider the following two statements for a particular $\mathcal{O}_X\widehat{\otimes}B$-pseudocoherent sheaf $\mathcal{M}$. First is that one can find finite many generators (the number is up to some uniform fixed integer $n\geq 0$) for each section of $\mathcal{M}(X_i)$ for each $i=1,2,...$. The second statement is that the global section $\mathcal{M}(X)$ is just finitely generated. Then in our situation the two statement is equivalent if we have that the corresponding space $X$ admits a $m$-uniformed covering in the exact same sense of \cite[Proposition 2.6.17]{8KL2}. 	
\end{theorem}

\begin{proof}
See \cref{theorem2.31}.
\end{proof}

%%\newpage

\newpage\section{Descent over Analytic Huber Pairs over $\mathbb{Z}_p$}

\subsection{Noncommutative Deformation over Analytic Huber Pairs in the Analytic Topology}

\indent We now study the corresponding glueing of stably pseudocoherent sheaves after our previous work \cite{8T3}. However since we currently are working in the corresponding framework of \cite{8Ked2}, we will consider the corresponding analytic Huber pair in \cite{8Ked2}.

\begin{setting}
Now we consider the corresponding Huber pair taking the general form of $(A,A^+)$ as in \cite[Definition 1.1.2]{8Ked2}, namely it is uniform analytic. We now assume we are going to work over some base $(V,V^+)=(\mathbb{Z}_p,\mathbb{Z}_p)$. And we will assume that we have another topological ring $Z$ which is assumed to satisfy the following condition: for any exact sequence of uniform analytic Huber rings
\[
\xymatrix@C+0pc@R+0pc{
0 \ar[r] \ar[r] \ar[r] &\Gamma_1 \ar[r] \ar[r] \ar[r] &\Gamma_2 \ar[r] \ar[r] \ar[r] &\Gamma_3 \ar[r] \ar[r] \ar[r] &0,
}
\] 
we have that the following is alway exact:
\[
\xymatrix@C+0pc@R+0pc{
0 \ar[r] \ar[r] \ar[r] &\Gamma_1\widehat{\otimes} Z \ar[r] \ar[r] \ar[r] &\Gamma_2\widehat{\otimes} Z \ar[r] \ar[r] \ar[r] &\Gamma_3\widehat{\otimes} Z \ar[r] \ar[r] \ar[r] &0.
}
\] 
This is achieved for instance when we have that the map $C\rightarrow Z$ splits in the category of all the topological modules. 
\end{setting}

\begin{setting}
We will maintain the corresponding assumption in \cite[Hypothesis 1.7.1]{8Ked2}. To be more precise we will have a map $(A,A^+)\rightarrow (B,B^+)$ which is a corresponding rational localization. And recall the corresponding complex in \cite{8Ked2}:
\[
\xymatrix@C+0pc@R+0pc{
0 \ar[r] \ar[r] \ar[r] &B \ar[r] \ar[r] \ar[r] &B\left\{\frac{f}{g}\right\}\bigoplus B\left\{\frac{g}{f}\right\} \ar[r] \ar[r] \ar[r] &B\left\{\frac{f}{g},\frac{g}{f}\right\} \ar[r] \ar[r] \ar[r] &0.
}
\]	
\end{setting}

\begin{setting}
We will need to consider the corresponding nice rational localizations in the sense of \cite[Definition 1.9.1]{8Ked2} where such localizations are defined to be those composites of the corresponding rational localizations in the Laurent or the balanced situation. For the topological modules, we will always assume that the modules are left modules.
\end{setting}

\begin{definition}\mbox{\bf{(After Kedlaya \cite[Definition 1.9.1]{8Ked2})}}
Over $A$, we define a corresponding $Z$-stably-pseudocoherent module $M$ to be a pseudocoherent module $M$ over $A\widehat{\otimes}Z$ which is complete with respect to the natural topology. And for any rational localization $A\rightarrow A'$ with respect to $A$, the completeness still holds. Similar we can define the corresponding $m$-$Z$-stably-pseudocoherent modules in the similar way where we just consider the corresponding $m$-pseudocoherent modules on the algebraic level. Over $A$, we define a corresponding $Z$-nice-stably-pseudocoherent module $M$ to be a pseudocoherent module $M$ over $A\widehat{\otimes}Z$ which is complete with respect to the natural topology. And for any nice rational localization $A\rightarrow A'$ with respect to $A$, the completeness still holds. Similar we can define the corresponding $m$-$Z$-nice-stably-pseudocoherent modules in the similar way where we just consider the corresponding $m$-pseudocoherent modules on the algebraic level.
\end{definition}

%%%%%%%%%%%%%%%%%%%%%%%%%%%%%%%%%%%%%%%%%!!!!!!!!!!

\begin{lemma} \mbox{\bf{(After Kedlaya \cite[Lemma 1.9.3]{8Ked2})}} \label{8lemma2.5}
Now assume that we are in the corresponding assumption of \cite[1.7.1]{8Ked2}. Also consider the corresponding complex in \cite[1.6.15.1]{8Ked2}:
\[
\xymatrix@C+0pc@R+0pc{
0 \ar[r] \ar[r] \ar[r] &B \ar[r] \ar[r] \ar[r] &B\left\{\frac{f}{g}\right\}\bigoplus B\left\{\frac{g}{f}\right\} \ar[r] \ar[r] \ar[r] &B\left\{\frac{f}{g},\frac{g}{f}\right\} \ar[r] \ar[r] \ar[r] &0.
}
\]	
And we assume now that this is exact. Now take $g$ to be $1-f$ or $1$. Here the corresponding elements $f,g$ come from $B$. Now over $B\widehat{\otimes}Z$ suppose we have a module $M$ which is assumed now to be finitely presented and complete with respect to the natural topology over $B\widehat{\otimes}Z$. Then we have that the group $\mathrm{Tor}_1(B\widehat{\otimes}Z\left\{\frac{g}{f}\right\},M)$ is zero. 	
\end{lemma}

\begin{proof}
This is a $Z$-relative version of the corresponding \cite[Lemma 1.9.3]{8Ked2}. We briefly recall the argument being adapted in our situation. As in \cite[Lemma 1.9.3]{8Ked2} we can reduce ourselves to the situation where $g$ is just $1$. Then regard the corresponding ring $B\widehat{\otimes}Z\left\{\frac{g}{f}\right\}$ as $B\widehat{\otimes}Z\left\{T\right\}/(1-fT)$. Then we consider the corresponding in the ring $B\widehat{\otimes}Z[[T]]$ the formal inverse of $1-fT$ namely the series $1+fT+f^2T^2+...$. This will force multiplication by this element to be injective. This just finishes the proof as in \cite[Lemma 1.9.3]{8Ked2}.
\end{proof}

\begin{lemma} \mbox{\bf{(After Kedlaya \cite[Lemma 1.9.4]{8Ked2})}} \label{lemma2.6}
Now assume that we are in the corresponding assumption of \cite[1.7.1]{8Ked2}. Also consider the corresponding complex in \cite[1.6.15.1]{8Ked2}:
\[
\xymatrix@C+0pc@R+0pc{
0 \ar[r] \ar[r] \ar[r] &B \ar[r] \ar[r] \ar[r] &B\left\{\frac{f}{g}\right\}\bigoplus B\left\{\frac{g}{f}\right\} \ar[r] \ar[r] \ar[r] &B\left\{\frac{f}{g},\frac{g}{f}\right\} \ar[r] \ar[r] \ar[r] &0.
}
\]	
And we assume now that this is exact. Now take $g$ to be $1-f$ or $1$. Here the corresponding element $f,g$ come from $B$. Now over $B\widehat{\otimes}Z$ suppose we have a module $M$ which is assumed now to be finitely presented and complete with respect to the natural topology over $B\widehat{\otimes}Z$, and furthermore the corresponding completeness is stable under any corresponding nice rational localization $B\rightarrow C$. Then we have that the group $\mathrm{Tor}_1(B\widehat{\otimes}Z\left\{\frac{f}{g}\right\},M)$ is zero and the group $\mathrm{Tor}_1(B\widehat{\otimes}Z\left\{\frac{f}{g},\frac{g}{f}\right\},M)$ is zero as well. 	
\end{lemma}

\begin{proof}
Choose a corresponding presentation for $M$ taking the form of:
\[
\xymatrix@C+0pc@R+0pc{
0 \ar[r] \ar[r] \ar[r] &X_1 \ar[r] \ar[r] \ar[r] &X_2\ar[r] \ar[r] \ar[r] &M \ar[r] \ar[r] \ar[r] &0.
}
\]
Here the module $X_2$ is a finite free left module. Then we have by the previous lemma the following exact sequence:
\[
\xymatrix@C+0pc@R+0pc{
0 \ar[r] \ar[r] \ar[r] &B\widehat{\otimes}Z\left\{\frac{g}{f}\right\}\otimes_{B\widehat{\otimes}Z}X_1 \ar[r] \ar[r] \ar[r] &B\widehat{\otimes}Z\left\{\frac{g}{f}\right\}\otimes_{B\widehat{\otimes}Z}X_2 \ar[r] \ar[r] \ar[r] &B\widehat{\otimes}Z\left\{\frac{g}{f}\right\}\otimes_{B\widehat{\otimes}Z}M \ar[r] \ar[r] \ar[r] &0.
}
\]	
Now apply the previous lemma again as in the proof of \cite[Lemma 1.9.4]{8Ked2} we have the following exact sequence:
\[\small
\xymatrix@C-0.4pc@R+0pc{
0 \ar[r] \ar[r] \ar[r] &B\widehat{\otimes}Z\left\{\frac{f}{g},\frac{g}{f}\right\}\otimes_{B\widehat{\otimes}Z}X_1 \ar[r] \ar[r] \ar[r] &B\widehat{\otimes}Z\left\{\frac{f}{g},\frac{g}{f}\right\}\otimes_{B\widehat{\otimes}Z}X_2 \ar[r] \ar[r] \ar[r] &B\widehat{\otimes}Z\left\{\frac{f}{g},\frac{g}{f}\right\}\otimes_{B\widehat{\otimes}Z}M \ar[r] \ar[r] \ar[r] &0.
}
\]	
Then to tackle the situation for $\frac{f}{g}$, we just further take the tensor product with the corresponding sequence as in the following:
\[
\xymatrix@C+0pc@R+0pc{
0 \ar[r] \ar[r] \ar[r] &B\widehat{\otimes}Z \ar[r] \ar[r] \ar[r] &B\widehat{\otimes}Z\left\{\frac{f}{g}\right\}\bigoplus B\widehat{\otimes}Z\left\{\frac{g}{f}\right\} \ar[r] \ar[r] \ar[r] &B\widehat{\otimes}Z\left\{\frac{f}{g},\frac{g}{f}\right\} \ar[r] \ar[r] \ar[r] &0.
}
\]
\end{proof}

\indent Then we have the following corollary which is the corresponding analog of \cite[Corollary 1.9.5]{8Ked2}:

\begin{corollary}\mbox{\bf{(After Kedlaya \cite[Corollary 1.9.5]{8Ked2})}}
Working over a uniform analytic Huber pair $(A,A^+)$ which is sheafy, suppose we consider a finitely generated $A\widehat{\otimes}Z$ module complete with respect to the natural topology. Then we have that for any nice rational localization $A\rightarrow B$ in the sense of \cite[Definition 1.9.1]{8Ked2} then tensoring with the corresponding $Z$ we have the vanishing of $\mathrm{Tor}_1(B,M)$.
	
\end{corollary}

\begin{proof}
By the previous lemmas.	
\end{proof}

\begin{corollary}\mbox{\bf{(After Kedlaya \cite[Corollary 1.9.6]{8Ked2})}}
Working over a uniform analytic Huber pair $(A,A^+)$ which is sheafy, suppose we consider a nice $Z$-stably pseudocoherent left $A\widehat{\otimes}Z$ module $M$. Then we have that for any rational localization $A\rightarrow B$ and any rational localization from $B$ to $C$ in the nice setting as in \cite[Definition 1.9.1]{8Ked2}. We have that base change along the map $B\rightarrow C$ (note that this is assumed to be nice) will preserve the corresponding stably pseudocoherence for $M$.

\end{corollary}

\begin{lemma} \mbox{\bf{(After Kedlaya \cite[Lemma 1.9.7]{8Ked2})}}
Now assume that we are in the corresponding assumption of \cite[1.7.1]{8Ked2}. Also consider the corresponding complex in \cite[1.6.15.1]{8Ked2}:
\[
\xymatrix@C+0pc@R+0pc{
0 \ar[r] \ar[r] \ar[r] &B \ar[r] \ar[r] \ar[r] &B\left\{\frac{f}{g}\right\}\bigoplus B\left\{\frac{g}{f}\right\} \ar[r] \ar[r] \ar[r] &B\left\{\frac{f}{g},\frac{g}{f}\right\} \ar[r] \ar[r] \ar[r] &0.
}
\]	
And we assume now that this is exact. Now take $g$ to be $1-f$ or $1$. Here the corresponding element $f,g$ come from $B$. Now over $B\widehat{\otimes}Z$ suppose we have a module $M$ which is assumed now to be finitely presented and complete with respect to the natural topology over $B\widehat{\otimes}Z$, and furthermore the corresponding completeness is stable under any corresponding nice rational localization $B\rightarrow C$. Then we have that over then $B\widehat{\otimes}Z$ in our situation tensoring with $M$ will preserve the corresponding exactness of the following:

\[
\xymatrix@C+0pc@R+0pc{
0 \ar[r] \ar[r] \ar[r] &B\widehat{\otimes}Z \ar[r] \ar[r] \ar[r] &B\widehat{\otimes}Z \left\{\frac{f}{g}\right\}\bigoplus B\widehat{\otimes}Z \left\{\frac{g}{f}\right\} \ar[r] \ar[r] \ar[r] &B\widehat{\otimes}Z \left\{\frac{f}{g},\frac{g}{f}\right\} \ar[r] \ar[r] \ar[r] &0.
}
\]

\end{lemma}

\begin{proof}
See \cite[Lemma 1.9.7]{8Ked2}.	
\end{proof}

\begin{theorem} \mbox{\bf{(After Kedlaya \cite[Corollary 1.9.8]{8Ked2})}}
Working over a uniform analytic Huber pair $(A,A^+)$ which is sheafy, suppose we consider a $Z$-nice-stably pseudocoherent left $A\widehat{\otimes}Z$ module $M$. Then we have that the presheaf $\widetilde{M}$ associated to the $M$ sheafified along the ring $A$ only is acyclic.
	
\end{theorem}

\begin{proof}
See \cite[Corollary 1.9.8]{8Ked2}.	
\end{proof}

\indent Now we establish the corresponding analogs of the corresponding results in \cite{8Ked2} which are needed in the corresponding descent of pseudocoherent sheaves.

\begin{lemma}\mbox{\bf{(After Kedlaya \cite[Lemma 1.9.9]{8Ked2})}}
Again as in \cite[Lemma 1.9.9]{8Ked2} we suppose we are in the situation of 1.7.1 of \cite{8Ked2} for the adic rings. Now one can actually find a neighbourhood around 0 of the corresponding ring $B\widehat{\otimes}Z\left\{\frac{f}{g},\frac{g}{f}\right\}$. Then one can find a corresponding decomposition for any invertible matrix $M$ where $M-1$ takes coefficients in this neighbourhood into $M_1M_2$ where invertible $M_1$ takes coefficients in $B\widehat{\otimes}Z\left\{\frac{f}{g}\right\}$ and invertible $M_2$ takes coefficients in $B\widehat{\otimes}Z\left\{\frac{g}{f}\right\}$. 
\end{lemma}

\begin{proof}
See \cite[Lemma 1.9.9]{8Ked2}.	
\end{proof}

%%%%%%%%%%%%%%%%%%%%%%%%%!!!!!!!!!!!!!!!!!!!!

\begin{lemma}\mbox{\bf{(After Kedlaya \cite[Lemma 1.9.10]{8Ked2})}} \label{lemma2.12}
Again as in \cite[Lemma 1.9.10]{8Ked2} we suppose we are in the situation of 1.7.1 of \cite{8Ked2} for the adic rings. Now for a glueing datum $M_1,M_2,M_{12}$ over the rings respectively $B\widehat{\otimes}Z\left\{\frac{f}{g}\right\}$, $B\widehat{\otimes}Z\left\{\frac{g}{f}\right\}$, $B\widehat{\otimes}Z\left\{\frac{f}{g},\frac{g}{f}\right\}$. We assume that they are finitely generated. We recall this means that we have:
\begin{align}
f_1: B\widehat{\otimes}Z\left\{\frac{f}{g},\frac{g}{f}\right\}\otimes_{B\widehat{\otimes}Z\left\{\frac{f}{g}\right\}}M_1\overset{\sim}{\rightarrow}M_{12},\\
f_2: B\widehat{\otimes}Z\left\{\frac{f}{g},\frac{g}{f}\right\}\otimes_{B\widehat{\otimes}Z\left\{\frac{f}{g}\right\}}M_2\overset{\sim}{\rightarrow}M_{12}.
\end{align}
Then we have that the corresponding strictly surjectivity of the map $M_1\bigoplus M_2\rightarrow M_{12}$ given by $f_1(x_1)-f_2(x_2)$. And we have that the corresponding equalizer $M$ will satisfy that $B\widehat{\otimes}Z\left\{\frac{f}{g}\right\} \otimes_{B\widehat{\otimes}Z}M\overset{}{\rightarrow}M_1$, $B\widehat{\otimes}Z\left\{\frac{g}{f}\right\}\otimes_{B\widehat{\otimes}Z}M\overset{}{\rightarrow}M_2$ are strictly surjective.	
\end{lemma}

\begin{proof}
This is a $Z$-relative version of the corresponding result in \cite[Lemma 1.9.10]{8Ked2}. We adapt the corresponding argument there to our situation. Choose basis of $M_1$ namely $e_1,...,e_n$, and choose basis of $M_2$ namely $h_1,...,h_n$. Now we set the corresponding matrix of the corresponding $f_2(w_i)$ under the basis $e_1,...,e_n$ to be $P$ while we set the corresponding matrix of the corresponding $f_1(v_i)$ under the basis $h_1,...,h_n$ to be $Q$. Then what is happening is that we can apply the previous lemma to our current situation to extract a neighbourhood such that $P(f^{-k}Q'-Q)$ (for some sufficiently $k$) takes the value in this neighbourhood and we have that the matrix $P(f^{-k}Q'-Q)+1$ will be decomposed into $K_1K_2^{-1}$ with the corresponding conditions on the entries as in the previous lemma. Then what we could do is as in \cite[Lemma 1.9.10]{8Ked2} to consider the setting $(x_1,x_2):=(\sum f^kX_{1ij}e_i,\sum Q'X_{2ij}h_i)$ which gives rise to that $f_1(x_1)-f_2(x_2)=0$ by our construction. Then as in \cite[Lemma 1.9.10]{8Ked2} the corresponding result in the first statement will follow once one considers as well \cite[Lemma 1.8.1]{8Ked2}. Then for the second statement, we consider application of the first statement as in \cite[Lemma 1.9.10]{8Ked2}. Briefly recalling, we consider the second part of statement (while the first part could be derived in the same fashion). To be more precise we consider any element $m_2\in M_2$ then send by $f_2$ to $M_{12}$, we will have by the first statement some $y_1\in M_1,y_2\in M_2$ such that $f_1(y_1)-f_2(y_2)=f_2(m_2)$. This will gives $(y_1,y_2+m_2)$ lives in the equalizer, which will proves the result as in \cite[Lemma 1.9.10]{8Ked2}.

\end{proof}

\begin{proposition} \mbox{\bf{(After Kedlaya \cite[Lemma 1.9.11]{8Ked2})}} \label{8lemma2.13}
Consider the same situation of the previous lemma. And assume that the corresponding \cite[1.6.15.1]{8Ked2} is exact. And we assume that the corresponding modules $M_1,M_2,M_{12}$ are $Z$-stably pseudocoherent over the rings respectively $B\widehat{\otimes}Z\left\{\frac{f}{g}\right\}$, $B\widehat{\otimes}Z\left\{\frac{g}{f}\right\}$, $B\widehat{\otimes}Z\left\{\frac{f}{g},\frac{g}{f}\right\}$. Then we have that actually the corresponding equalizer $M$ is just pseudocoherent in our current situation.	
\end{proposition}

%%%%%%%%%%%%%%%%%%%%%%%%%%%%%%%%%%%%%%%!!!!!!!!!!!!!!!!!!

\begin{proof}
This is a $Z$-relative version of \cite[Lemma 1.9.11]{8Ked2}. We adapt the corresponding argument to our situation. First we consider corresponding modules $M_1,M_2,M_{12}$, and choose finite free modules $L_1,L_2,L_{12}$ as in \cite[Lemma 1.9.11]{8Ked2} over $B\widehat{\otimes}Z\left\{\frac{f}{g}\right\}$, $B\widehat{\otimes}Z\left\{\frac{g}{f}\right\}$, $B\widehat{\otimes}Z\left\{\frac{f}{g},\frac{g}{f}\right\}$, then we can form the corresponding commutative diagram:
\[
\xymatrix@C+3pc@R+3pc{
 &0 \ar[d] \ar[d] \ar[d] &0 \ar[d] \ar[d] \ar[d]  &0 \ar[d] \ar[d] \ar[d]&\\
0 \ar[r] \ar[r] \ar[r] &F \ar[r] \ar[r] \ar[r] \ar[d] \ar[d] \ar[d] &F_1\bigoplus F_2 \ar[r] \ar[r] \ar[r] \ar[d] \ar[d] \ar[d] &F_{12} \ar[r]^? \ar[r] \ar[r] \ar[d] \ar[d] \ar[d] &0\\
0 \ar[r] \ar[r] \ar[r] &L \ar[r] \ar[r] \ar[r] \ar[d] \ar[d] \ar[d] &L_1\bigoplus L_2 \ar[r] \ar[r] \ar[r] \ar[d] \ar[d] \ar[d] &L_{12} \ar[r] \ar[r] \ar[r] \ar[d] \ar[d] \ar[d] &0\\
0 \ar[r] \ar[r] \ar[r] &M \ar[r] \ar[r] \ar[r] \ar[d]^? \ar[d] \ar[d] &M_1\bigoplus M_2\ar[r] \ar[r] \ar[r] \ar[d] \ar[d] \ar[d] &M_{12} \ar[r] \ar[r] \ar[r] \ar[d] \ar[d] \ar[d] &0\\
&0 &0  &0  &\\
}
\]	
where two arrows marked by $?$ are not a priori known to be exist to make the corresponding parts of the sequences exact. Here the left vertical arrows come from taking snake lemma (recall that one first chooses a corresponding finite free module $L$ mapping to $M$ which will base change to the corresponding module $L_1,L_2,L_{12}$, then $F$ will be taken as in \cite[Lemma 1.9.11]{8Ked2} to be the corresponding equalizer of the corresponding horizontal map, the idea is that the map $L\rightarrow M$ is just a map which is not a priori known to be surjective while eventually this is indeed the case). The corresponding modules $F_1,F_2,F_{12}$ are the corresponding kernels of the corresponding coverings by finite free modules. Then what is happening is that the corresponding first horizontal complex is actually within the same situation we are considering, therefore the previous lemma implies that the $?$-marked arrows exists and make the sequences exact. We then consider the following commutative diagram:
\[
\xymatrix@C+0pc@R+0pc{
 &B\widehat{\otimes}Z\left\{\frac{f}{g}\right\}\otimes F \ar[r] \ar[r] \ar[r] \ar[d] \ar[d] \ar[d] &B\widehat{\otimes}Z\left\{\frac{f}{g}\right\}\otimes L \ar[r] \ar[r] \ar[r] \ar[d] \ar[d] \ar[d] &B\widehat{\otimes}Z\left\{\frac{f}{g}\right\}\otimes M \ar[r] \ar[r] \ar[r] \ar[d] \ar[d] \ar[d] &0.\\
0 \ar[r] \ar[r] \ar[r] &F_1 \ar[r] \ar[r] \ar[r] &L_1 \ar[r] \ar[r] \ar[r] &M_1 \ar[r] \ar[r] \ar[r] &0.\\
}
\]
What we know is that the corresponding middle vertical arrow is isomorphism, while the left and right ones are surjective, which implies the right most one is an isomorphism. We then consider the following commutative diagram:
\[
\xymatrix@C+0pc@R+0pc{
 &B\widehat{\otimes}Z\left\{\frac{g}{f}\right\}\otimes F \ar[r] \ar[r] \ar[r] \ar[d] \ar[d] \ar[d] &B\widehat{\otimes}Z\left\{\frac{g}{f}\right\}\otimes L \ar[r] \ar[r] \ar[r] \ar[d] \ar[d] \ar[d] &B\widehat{\otimes}Z\left\{\frac{g}{f}\right\}\otimes M \ar[r] \ar[r] \ar[r] \ar[d] \ar[d] \ar[d] &0.\\
0 \ar[r] \ar[r] \ar[r] &F_2 \ar[r] \ar[r] \ar[r] &L_2 \ar[r] \ar[r] \ar[r] &M_2 \ar[r] \ar[r] \ar[r] &0.\\
}
\]
What we know is that the corresponding middle vertical arrow is isomorphism, while the left and right ones are surjective, which implies the right most one is an isomorphism. Finally to show that the module $M$ is pseudocoherent we consider the corresponding induction on the number of the finite projective modules participated in the resolution. $k=1$ will be obvious, while in general the induction step will achieved from $k$ to $k+1$ by considering as in \cite[Lemma 1.9.11]{8Ked2} the module $F$ and apply the corresponding construction in the same fashion.
\end{proof}

\
%%%%%%%%%%%%%%%%%%%%%%%%%%%%%%%%!!!!!!!!!!!!!!

\begin{proposition} \mbox{\bf{(After Kedlaya \cite[Lemma 1.9.13]{8Ked2})}}
Take the corresponding assumption in \cite[1.7.1]{8Ked2} and furthermore assume the sheafiness of $A$. Suppose we have that all the modules considered are $Z$-stably complete with respect to nice rational localization. Then consider $C_1$ the category of all such modules which are assumed to be pseudocoherent over $B\widehat{\otimes}Z\left\{\frac{f}{g},\frac{g}{f}\right\}$. Then consider $C_2$ the category of all such modules which are assumed to be pseudocoherent over $B\widehat{\otimes}Z\left\{\frac{f}{g}\right\}$. Then consider $C_3$ the category of all such modules which are assumed to be pseudocoherent over $B\widehat{\otimes}Z\left\{\frac{g}{f}\right\}$. Then consider $C_4$ the category of all such modules which are assumed to be pseudocoherent over $B\widehat{\otimes}Z$. Then we have a corresponding equalizer functor:
\begin{align}
C_2\times_{C_1}C_3\rightarrow C	
\end{align}
which is exact and fully faithful, and we have well-established inverse.
\end{proposition}

\begin{proof}
See \cite[Lemma 1.9.13]{8Ked2}.	
\end{proof}

\begin{proposition}  \mbox{\bf{(After Kedlaya \cite[Lemma 1.9.16]{8Ked2})}}
Take the corresponding assumption in \cite[1.7.1]{8Ked2} and furthermore assume the sheafiness of $A$. Suppose we have that all the modules considered are $Z$-stably complete with respect to nice rational localization. Let $C$ be the category of such $Z$-nice-stably pseudocoherent modules over $A$ and let $C'$ be the category of such $Z$-nice-stably pseudocoherent sheaves over $X$. If we consider the corresponding functor taking the form of the corresponding sheafification of $\widetilde{M}$	with respect to the corresponding nice rational subspaces. Then we have that in fact this realizes an exact equivalence between $C$ and $C'$.
\end{proposition}

\begin{proof}
See \cite[Lemma 1.9.16]{8Ked2}.	
\end{proof}

\begin{proposition}  \mbox{\bf{(After Kedlaya \cite[Lemma 1.9.17]{8Ked2})}}
Take the corresponding assumption in \cite[1.7.1]{8Ked2} and furthermore assume the sheafiness of $A$. Suppose we have that all the modules considered are $Z$-stably complete with respect to nice rational localization. Let $C$ be the category of such $Z$-nice-stably pseudocoherent modules over $A$ and let $C'$ be the category of such $Z$-nice-stably pseudocoherent modules over $B$. Then we have that in fact this realizes an exact functor from $C$ to $C'$.
\end{proposition}

\begin{proof}
See \cite[Lemma 1.9.17]{8Ked2}.	
\end{proof}

\begin{proposition}  \mbox{\bf{(After Kedlaya \cite[Corollary 1.9.18]{8Ked2})}}
Take the corresponding assumption in \cite[1.7.1]{8Ked2} and furthermore assume the sheafiness of $A$. Let $C$ be the category of such $Z$-stably pseudocoherent modules over $A$ and let $C'$ be the category of such $Z$-nice-stably pseudocoherent modules over $A$. Then we have that in fact this realizes an equality from $C\hookrightarrow C'$.
\end{proposition}

\begin{proof}
See \cite[Corollary 1.9.18]{8Ked2}.	
\end{proof}

\begin{theorem}\mbox{}\\
\mbox{\bf{(After Kedlaya \cite[Theorem 1.4.14, Theorem 1.4.16, Theorem 1.4.18]{8Ked2})}}\\
Take the corresponding assumption in \cite[1.7.1]{8Ked2} and furthermore assume the sheafiness of $A$. Then: I. We have the corresponding stability under rational localization for $Z$-stably pseudocoherent modules. II. We have the corresponding acyclicity of the corresponding presheaf $\widetilde{M}$ attached to any $Z$-stably pseudocoherent module. III. We have the corresponding glueing of the corresponding $Z$-stably pseudocoherent modules as in \cite[Theorem 1.4.18]{8Ked2}.
\end{theorem}

\begin{proof}
This is a very parallel $Z$-relative analog of \cite[Theorem 1.4.14, Theorem 1.4.16, Theorem 1.4.18]{8Ked2}, we refer the readers to the proof of \cite[Theorem 1.4.14, Theorem 1.4.16, Theorem 1.4.18]{8Ked2}. See \cite[After Corollary 1.9.18, the proof of Theorem 1.4.14, the proof of Theorem 1.4.16, the proof of Theorem 1.4.18]{8Ked2}.	
\end{proof}

%%\newpage

\newpage\section{Descent over Adic Banach Rings over $\mathbb{Z}_p$}

\subsection{Noncommutative Deformation over Adic Banach Rings in the Analytic Topology}

\indent We now study the corresponding glueing of stably pseudocoherent sheaves after our previous work \cite{8T3}. We now translate the corresponding discussion to the context which is more related to \cite{8KL2} over adic Banach rings. One can see that these are really parallel to the discussion over analytic Huber pair.

\begin{setting}
Now we consider the corresponding adic Banach ring taking the general form of $(A,A^+)$, namely it is uniform analytic. We now assume we are going to work over some base $(V,V^+)=(\mathbb{Z}_p,\mathbb{Z}_p)$. And we will assume that we have another Banach ring $Z$ over $V$ which is assumed to satisfy the following condition: for any exact sequence of uniform analytic adic Banach rings
\[
\xymatrix@C+0pc@R+0pc{
0 \ar[r] \ar[r] \ar[r] &\Gamma_1 \ar[r] \ar[r] \ar[r] &\Gamma_2 \ar[r] \ar[r] \ar[r] &\Gamma_3 \ar[r] \ar[r] \ar[r] &0,
}
\] 
we have that the following is alway exact:
\[
\xymatrix@C+0pc@R+0pc{
0 \ar[r] \ar[r] \ar[r] &\Gamma_1\widehat{\otimes} Z \ar[r] \ar[r] \ar[r] &\Gamma_2\widehat{\otimes} Z \ar[r] \ar[r] \ar[r] &\Gamma_3\widehat{\otimes} Z \ar[r] \ar[r] \ar[r] &0.
}
\] 
This is achieved for instance when we have that the map $C\rightarrow Z$ splits in the category of all the Banach modules.
\end{setting}

\begin{setting}
We will maintain the corresponding assumption in \cite[Hypothesis 1.7.1]{8Ked2}. To be more precise we will have a map $(A,A^+)\rightarrow (B,B^+)$ which is a corresponding rational localization. And recall the corresponding complex in \cite{8Ked2}:
\[
\xymatrix@C+0pc@R+0pc{
0 \ar[r] \ar[r] \ar[r] &B \ar[r] \ar[r] \ar[r] &B\left\{\frac{f}{g}\right\}\bigoplus B\left\{\frac{g}{f}\right\} \ar[r] \ar[r] \ar[r] &B\left\{\frac{f}{g},\frac{g}{f}\right\} \ar[r] \ar[r] \ar[r] &0.
}
\]	
\end{setting}

\begin{setting}
We will need to consider the corresponding nice rational localizations in the sense of \cite[Definition 1.9.1]{8Ked2} where such localizations are defined to be those composites of the corresponding rational localizations in the Laurent or the balanced situation. For the topological modules, we will always assume that the modules are left modules.
\end{setting}

\begin{definition}\mbox{\bf{(After Kedlaya \cite[Definition 1.9.1]{8Ked2})}}
Over $A$, we define a corresponding $Z$-stably-pseudocoherent module $M$ to be a pseudocoherent module $M$ over $A\widehat{\otimes}Z$ which is complete with respect to the natural topology. And for any rational localization $A\rightarrow A'$ with respect to $A$, the completeness still holds. Similar we can define the corresponding $m$-$Z$-stably-pseudocoherent modules in the similar way where we just consider the corresponding $m$-pseudocoherent modules on the algebraic level. Over $A$, we define a corresponding $Z$-nice-stably-pseudocoherent module $M$ to be a pseudocoherent module $M$ over $A\widehat{\otimes}Z$ which is complete with respect to the natural topology. And for any nice rational localization $A\rightarrow A'$ with respect to $A$, the completeness still holds. Similar we can define the corresponding $m$-$Z$-nice-stably-pseudocoherent modules in the similar way where we just consider the corresponding $m$-pseudocoherent modules on the algebraic level.
\end{definition}

\begin{lemma} \mbox{\bf{(After Kedlaya \cite[Lemma 1.9.3]{8Ked2})}}
Now assume that we are in the corresponding assumption of \cite[1.7.1]{8Ked2}. Also consider the corresponding complex in \cite[1.6.15.1]{8Ked2}:
\[
\xymatrix@C+0pc@R+0pc{
0 \ar[r] \ar[r] \ar[r] &B \ar[r] \ar[r] \ar[r] &B\left\{\frac{f}{g}\right\}\bigoplus B\left\{\frac{g}{f}\right\} \ar[r] \ar[r] \ar[r] &B\left\{\frac{f}{g},\frac{g}{f}\right\} \ar[r] \ar[r] \ar[r] &0.
}
\]	
And we assume now that this is exact. Now take $g$ to be $1-f$ or $1$. Here the corresponding element $f,g$ come from $B$. Now over $B\widehat{\otimes}Z$ suppose we have a module $M$ which is assumed now to be finitely presented and complete with respect to the natural topology over $B\widehat{\otimes}Z$. Then we have that the group $\mathrm{Tor}_1(B\widehat{\otimes}Z\left\{\frac{g}{f}\right\},M)$ is zero. 	
\end{lemma}

\begin{proof}
See \cref{8lemma2.5}.
\end{proof}

\begin{lemma} \mbox{\bf{(After Kedlaya \cite[Lemma 1.9.4]{8Ked2})}}
Now assume that we are in the corresponding parallel assumption of \cite[1.7.1]{8Ked2}. Also consider the corresponding complex in \cite[1.6.15.1]{8Ked2}:
\[
\xymatrix@C+0pc@R+0pc{
0 \ar[r] \ar[r] \ar[r] &B \ar[r] \ar[r] \ar[r] &B\left\{\frac{f}{g}\right\}\bigoplus B\left\{\frac{g}{f}\right\} \ar[r] \ar[r] \ar[r] &B\left\{\frac{f}{g},\frac{g}{f}\right\} \ar[r] \ar[r] \ar[r] &0.
}
\]	
And we assume now that this is exact. Now take $g$ to be $1-f$ or $1$. Here the corresponding element $f,g$ come from $B$. Now over $B\widehat{\otimes}Z$ suppose we have a module $M$ which is assumed now to be finitely presented and complete with respect to the natural topology over $B\widehat{\otimes}Z$, and furthermore the corresponding completeness is stable under any corresponding nice rational localization $B\rightarrow C$. Then we have that the group $\mathrm{Tor}_1(B\widehat{\otimes}Z\left\{\frac{f}{g}\right\},M)$ is zero and the group $\mathrm{Tor}_1(B\widehat{\otimes}Z\left\{\frac{f}{g},\frac{g}{f}\right\},M)$ is zero as well. 	
\end{lemma}

\begin{proof}
See \cref{lemma2.6}.
\end{proof}

\indent Then we have the following corollary which is the corresponding analog of \cite[Corollary 1.9.5]{8Ked2}:

\begin{corollary}\mbox{\bf{(After Kedlaya \cite[Corollary 1.9.5]{8Ked2})}}
Working over a uniform analytic adic Banach ring $(A,A^+)$ which is sheafy, suppose we consider a finitely generated $A\widehat{\otimes}Z$ module complete with respect to the natural topology. Then we have that for any nice rational localization $A\rightarrow B$ in the sense of \cite[Definition 1.9.1]{8Ked2} then tensoring with the corresponding $Z$ we have the vanishing of $\mathrm{Tor}_1(B,M)$.
	
\end{corollary}

\begin{proof}
By the previous lemmas.	
\end{proof}

\begin{corollary}\mbox{\bf{(After Kedlaya \cite[Corollary 1.9.6]{8Ked2})}}
Working over a uniform analytic adic Banach ring $(A,A^+)$ which is sheafy, suppose we consider a nice $Z$-stably pseudocoherent left $A\widehat{\otimes}Z$ module $M$. Then we have that for any rational localization $A\rightarrow B$ and any rational localization from $B$ to $C$ in the nice setting as in \cite[Definition 1.9.1]{8Ked2}. We have that base change along the map $B\rightarrow C$ (note that this is assumed to be nice) will preserve the corresponding stably pseudocoherence for $M$.

\end{corollary}

\begin{lemma} \mbox{\bf{(After Kedlaya \cite[Lemma 1.9.7]{8Ked2})}}
Now assume that we are in the corresponding parallel assumption of \cite[1.7.1]{8Ked2}. Also consider the corresponding complex in \cite[1.6.15.1]{8Ked2}:
\[
\xymatrix@C+0pc@R+0pc{
0 \ar[r] \ar[r] \ar[r] &B \ar[r] \ar[r] \ar[r] &B\left\{\frac{f}{g}\right\}\bigoplus B\left\{\frac{g}{f}\right\} \ar[r] \ar[r] \ar[r] &B\left\{\frac{f}{g},\frac{g}{f}\right\} \ar[r] \ar[r] \ar[r] &0.
}
\]	
And we assume now that this is exact. Now take $g$ to be $1-f$ or $1$. Here the corresponding element $f,g$ come from $B$. Now over $B\widehat{\otimes}Z$ suppose we have a module $M$ which is assumed now to be finitely presented and complete with respect to the natural topology over $B\widehat{\otimes}Z$, and furthermore the corresponding completeness is stable under any corresponding nice rational localization $B\rightarrow C$. Then we have that over then $B\widehat{\otimes}Z$ in our situation tensoring with $M$ will preserve the corresponding exactness of the following:

\[
\xymatrix@C+0pc@R+0pc{
0 \ar[r] \ar[r] \ar[r] &B\widehat{\otimes}Z \ar[r] \ar[r] \ar[r] &B\widehat{\otimes}Z \left\{\frac{f}{g}\right\}\bigoplus B\widehat{\otimes}Z \left\{\frac{g}{f}\right\} \ar[r] \ar[r] \ar[r] &B\widehat{\otimes}Z \left\{\frac{f}{g},\frac{g}{f}\right\} \ar[r] \ar[r] \ar[r] &0.
}
\]

\end{lemma}

\begin{proof}
See \cite[Lemma 1.9.7]{8Ked2}.	
\end{proof}

\begin{theorem} \mbox{\bf{(After Kedlaya \cite[Corollary 1.9.8]{8Ked2})}}
Working over a uniform analytic adic Banach ring $(A,A^+)$ which is sheafy, suppose we consider a $Z$-nice-stably pseudocoherent left $A\widehat{\otimes}Z$ module $M$. Then we have that the presheaf $\widetilde{M}$ associated to the $M$ sheafified along the ring $A$ only is acyclic.
	
\end{theorem}

\begin{proof}
See \cite[Corollary 1.9.8]{8Ked2}.	
\end{proof}

\indent Now we establish the corresponding analogs of the corresponding results in \cite{8Ked2} which are needed in the corresponding descent of pseudocoherent sheaves.

\begin{lemma}\mbox{\bf{(After Kedlaya \cite[Lemma 1.9.9]{8Ked2})}}
Again as in \cite[Lemma 1.9.9]{8Ked2} we suppose we are in the parallel situation of 1.7.1 of \cite[Lemma 1.9.9]{8Ked2} for the adic Banach rings. Now take a sufficiently small constant $\delta>0$. Then one can find a corresponding decomposition for any invertible matrix $M$ where $M-1$ takes norm $\|M-1\|\leq \delta$ into $M_1M_2$ where invertible $M_1$ takes coefficients in $B\widehat{\otimes}Z\left\{\frac{f}{g}\right\}$ and invertible $M_2$ takes coefficients in $B\widehat{\otimes}Z\left\{\frac{g}{f}\right\}$. 
\end{lemma}

\begin{proof}
See \cite[Lemma 1.9.9]{8Ked2}.	
\end{proof}

\begin{lemma}\mbox{\bf{(After Kedlaya \cite[Lemma 1.9.10]{8Ked2})}}
Again as in \cite[Lemma 1.9.10]{8Ked2} we suppose we are in the parallel situation of 1.7.1 of \cite[Lemma 1.9.10]{8Ked2} for the adic Banach rings. Now for a glueing datum $M_1,M_2,M_{12}$ over the rings respectively $B\widehat{\otimes}Z\left\{\frac{f}{g}\right\}$, $B\widehat{\otimes}Z\left\{\frac{g}{f}\right\}$, $B\widehat{\otimes}Z\left\{\frac{f}{g},\frac{g}{f}\right\}$. We assume that they are finitely generated. We recall this means that we have:
\begin{align}
f_1: B\widehat{\otimes}Z\left\{\frac{f}{g},\frac{g}{f}\right\}\otimes_{B\widehat{\otimes}Z\left\{\frac{f}{g}\right\}}M_1\overset{\sim}{\rightarrow}M_{12},\\
f_2: B\widehat{\otimes}Z\left\{\frac{f}{g},\frac{g}{f}\right\}\otimes_{B\widehat{\otimes}Z\left\{\frac{f}{g}\right\}}M_2\overset{\sim}{\rightarrow}M_{12}.
\end{align}
Then we have that the corresponding strictly surjectivity of the map $M_1\bigoplus M_2\rightarrow M_{12}$ given by $f_1(x_1)-f_2(x_2)$. And we have that the corresponding equalizer $M$ will satisfy that $B\widehat{\otimes}Z\left\{\frac{f}{g}\right\} \otimes_{B\widehat{\otimes}Z}M\overset{}{\rightarrow}M_1$, $B\widehat{\otimes}Z\left\{\frac{g}{f}\right\}\otimes_{B\widehat{\otimes}Z}M\overset{}{\rightarrow}M_2$ are strictly surjective.	
\end{lemma}

\begin{proof}
See \cref{lemma2.12}.

\end{proof}

\begin{proposition} \mbox{\bf{(After Kedlaya \cite[Lemma 1.9.11]{8Ked2})}}
Consider the same situation of the previous lemma. And assume that the corresponding \cite[1.6.15.1]{8Ked2} is exact. And we assume that the corresponding modules $M_1,M_2,M_{12}$ are $Z$-stably pseudocoherent over the rings respectively $B\widehat{\otimes}Z\left\{\frac{f}{g}\right\}$, $B\widehat{\otimes}Z\left\{\frac{g}{f}\right\}$, $B\widehat{\otimes}Z\left\{\frac{f}{g},\frac{g}{f}\right\}$. Then we have that actually the corresponding equalizer $M$ is just pseudocoherent in our current situation.	
\end{proposition}

\begin{proof}
This is a $Z$-relative version of \cite[Lemma 1.9.11]{8Ked2}. We adapt the corresponding argument to our situation. First we consider corresponding modules $M_1,M_2,M_{12}$, and choose finite free modules $L_1,L_2,L_{12}$ as in \cite[Lemma 1.9.11]{8Ked2} over $B\widehat{\otimes}Z\left\{\frac{f}{g}\right\}$, $B\widehat{\otimes}Z\left\{\frac{g}{f}\right\}$, $B\widehat{\otimes}Z\left\{\frac{f}{g},\frac{g}{f}\right\}$, then we can form the corresponding commutative diagram:
\[
\xymatrix@C+3pc@R+3pc{
 &0 \ar[d] \ar[d] \ar[d] &0 \ar[d] \ar[d] \ar[d]  &0 \ar[d] \ar[d] \ar[d]&\\
0 \ar[r] \ar[r] \ar[r] &F \ar[r] \ar[r] \ar[r] \ar[d] \ar[d] \ar[d] &F_1\bigoplus F_2 \ar[r] \ar[r] \ar[r] \ar[d] \ar[d] \ar[d] &F_{12} \ar[r]^? \ar[r] \ar[r] \ar[d] \ar[d] \ar[d] &0\\
0 \ar[r] \ar[r] \ar[r] &L \ar[r] \ar[r] \ar[r] \ar[d] \ar[d] \ar[d] &L_1\bigoplus L_2 \ar[r] \ar[r] \ar[r] \ar[d] \ar[d] \ar[d] &L_{12} \ar[r] \ar[r] \ar[r] \ar[d] \ar[d] \ar[d] &0\\
0 \ar[r] \ar[r] \ar[r] &M \ar[r] \ar[r] \ar[r] \ar[d]^? \ar[d] \ar[d] &M_1\bigoplus M_2\ar[r] \ar[r] \ar[r] \ar[d] \ar[d] \ar[d] &M_{12} \ar[r] \ar[r] \ar[r] \ar[d] \ar[d] \ar[d] &0\\
&0 &0  &0  &\\
}
\]	
where two arrows marked by $?$ are not a priori known to be exist to make the corresponding parts of the sequences exact. Here the left vertical arrows come from taking snake lemma (recall that one first chooses a corresponding finite free module $L$ mapping to $M$ which will base change to the corresponding module $L_1,L_2,L_{12}$, then $F$ will be taken as in \cite[Lemma 1.9.11]{8Ked2} to be the corresponding equalizer of the corresponding horizontal map, the idea is that the map $L\rightarrow M$ is just a map which is not a priori known to be surjective while eventually this is indeed the case). The corresponding modules $F_1,F_2,F_{12}$ are the corresponding kernels of the corresponding coverings by finite free modules. Then what is happening is that the corresponding first horizontal complex is actually within the same situation we are considering, therefore the previous lemma implies that the $?$-marked arrows exists and make the sequences exact. We then consider the following commutative diagram:
\[
\xymatrix@C+0pc@R+0pc{
 &B\widehat{\otimes}Z\left\{\frac{f}{g}\right\}\otimes F \ar[r] \ar[r] \ar[r] \ar[d] \ar[d] \ar[d] &B\widehat{\otimes}Z\left\{\frac{f}{g}\right\}\otimes L \ar[r] \ar[r] \ar[r] \ar[d] \ar[d] \ar[d] &B\widehat{\otimes}Z\left\{\frac{f}{g}\right\}\otimes M \ar[r] \ar[r] \ar[r] \ar[d] \ar[d] \ar[d] &0.\\
0 \ar[r] \ar[r] \ar[r] &F_1 \ar[r] \ar[r] \ar[r] &L_1 \ar[r] \ar[r] \ar[r] &M_1 \ar[r] \ar[r] \ar[r] &0.\\
}
\]
What we know is that the corresponding middle vertical arrow is isomorphism, while the left and right ones are surjective, which implies the right most one is an isomorphism. We then consider the following commutative diagram:
\[
\xymatrix@C+0pc@R+0pc{
 &B\widehat{\otimes}Z\left\{\frac{g}{f}\right\}\otimes F \ar[r] \ar[r] \ar[r] \ar[d] \ar[d] \ar[d] &B\widehat{\otimes}Z\left\{\frac{g}{f}\right\}\otimes L \ar[r] \ar[r] \ar[r] \ar[d] \ar[d] \ar[d] &B\widehat{\otimes}Z\left\{\frac{g}{f}\right\}\otimes M \ar[r] \ar[r] \ar[r] \ar[d] \ar[d] \ar[d] &0.\\
0 \ar[r] \ar[r] \ar[r] &F_2 \ar[r] \ar[r] \ar[r] &L_2 \ar[r] \ar[r] \ar[r] &M_2 \ar[r] \ar[r] \ar[r] &0.\\
}
\]
What we know is that the corresponding middle vertical arrow is isomorphism, while the left and right ones are surjective, which implies the right most one is an isomorphism. Finally to show that the module $M$ is pseudocoherent we consider the corresponding induction on the number of the finite projective modules participated in the resolution. $k=1$ will be obvious, while in general the induction step will achieved from $k$ to $k+1$ by considering as in \cite[Lemma 1.9.11]{8Ked2} the module $F$ and apply the corresponding construction in the same fashion.
\end{proof}

\

\begin{proposition} \mbox{\bf{(After Kedlaya \cite[Lemma 1.9.13]{8Ked2})}}
Take the corresponding parallel assumption in \cite[1.7.1]{8Ked2} and furthermore assume the sheafiness of $A$. Suppose we have that all the modules considered are $Z$-stably complete with respect to nice rational localization. Then consider $C_1$ the category of all such modules which are assumed to be pseudocoherent over $B\widehat{\otimes}Z\left\{\frac{f}{g},\frac{g}{f}\right\}$. Then consider $C_2$ the category of all such modules which are assumed to be pseudocoherent over $B\widehat{\otimes}Z\left\{\frac{f}{g}\right\}$. Then consider $C_3$ the category of all such modules which are assumed to be pseudocoherent over $B\widehat{\otimes}Z\left\{\frac{g}{f}\right\}$. Then consider $C_4$ the category of all such modules which are assumed to be pseudocoherent over $B\widehat{\otimes}Z$. Then we have a corresponding equalizer functor:
\begin{align}
C_2\times_{C_1}C_3\rightarrow C	
\end{align}
which is exact and fully faithful, and we have well-established inverse.
\end{proposition}

\begin{proof}
See \cite[Lemma 1.9.13]{8Ked2}.	
\end{proof}

\begin{proposition}  \mbox{\bf{(After Kedlaya \cite[Lemma 1.9.16]{8Ked2})}}
Take the corresponding parallel assumption in \cite[1.7.1]{8Ked2} and furthermore assume the sheafiness of the adic Banach $A$. Suppose we have that all the modules considered are $Z$-stably complete with respect to nice rational localization. Let $C$ be the category of such $Z$-nice-stably pseudocoherent modules over $A$ and let $C'$ be the category of such $Z$-nice-stably pseudocoherent sheaves over $X$. If we consider the corresponding functor taking the form of the corresponding sheafification of $\widetilde{M}$	with respect to the corresponding nice rational subspaces. Then we have that in fact this realizes an exact equivalence between $C$ and $C'$.
\end{proposition}

\begin{proof}
See \cite[Lemma 1.9.16]{8Ked2}.	
\end{proof}

\begin{proposition}  \mbox{\bf{(After Kedlaya \cite[Lemma 1.9.17]{8Ked2})}}
Take the corresponding parallel assumption in \cite[1.7.1]{8Ked2} and furthermore assume the sheafiness of the adic Banach ring $A$. Suppose we have that all the modules considered are $Z$-stably complete with respect to nice rational localization. Let $C$ be the category of such $Z$-nice-stably pseudocoherent modules over $A$ and let $C'$ be the category of such $Z$-nice-stably pseudocoherent modules over $B$. Then we have that in fact this realizes an exact functor from $C$ to $C'$.
\end{proposition}

\begin{proof}
See \cite[Lemma 1.9.17]{8Ked2}.	
\end{proof}

\begin{proposition}  \mbox{\bf{(After Kedlaya \cite[Corollary 1.9.18]{8Ked2})}}
Take the corresponding parallel assumption in \cite[1.7.1]{8Ked2} and furthermore assume the sheafiness of the adic Banach ring $A$. Let $C$ be the category of such $Z$-stably pseudocoherent modules over $A$ and let $C'$ be the category of such $Z$-nice-stably pseudocoherent modules over $A$. Then we have that in fact this realizes an equality from $C\subset C'$.
\end{proposition}

\begin{proof}
See \cite[Corollary 1.9.18]{8Ked2}.	
\end{proof}

\begin{theorem}\mbox{}\\
\mbox{\bf{(After Kedlaya \cite[Theorem 1.4.14, Theorem 1.4.16, Theorem 1.4.18]{8Ked2})}}\\
Take the corresponding parallel assumption in \cite[1.7.1]{8Ked2} and furthermore assume the sheafiness of the adic Baanch ring $A$. Then: I. we have the corresponding stability under rational localization for $Z$-stably pseudocoherent modules. II. We have the corresponding acyclicity of the corresponding presheaf $\widetilde{M}$ attached to any $Z$-stably pseudocoherent module. III. We have the corresponding glueing result as in \cite[Theorem 1.4.18]{8Ked2}.
\end{theorem}

\begin{proof}
This is a very parallel $Z$-relative analog of \cite[Theorem 1.4.14, Theorem 1.4.16, Theorem 1.4.18]{8Ked2}, we refer the readers to the proof of \cite[Theorem 1.4.14, Theorem 1.4.16, Theorem 1.4.18]{8Ked2}. See \cite[After Corollary 1.9.18, the proof of Theorem 1.4.14, the proof of Theorem 1.4.16, the proof of Theorem 1.4.18]{8Ked2}.	
\end{proof}

\

\subsection{Noncommutative Deformation over Quasi-Stein Spaces}

\indent We now consider the corresponding analytic quasi-Stein spaces over $\mathbb{Z}_p$, $Z$ will be an adic Banach ring as above:

\begin{setting} \mbox{\bf{(After Kedlaya-Liu \cite[Definition 2.6.2]{8KL2})}}
Fix an adic Banach space $X$ over $\mathbb{Z}_p$ which is assumed to be an injective limit of analytic adic affinoids $X=\varinjlim_i X_i$, where we assume that the corresponding transition map $\mathcal{O}_{X_{i+1}}\rightarrow \mathcal{O}_{X_{i}}$ takes dense images for all $i=1,2,...$.	
\end{setting}

\begin{lemma} \mbox{\bf{(After Kedlaya-Liu \cite[Lemma 2.6.3]{8KL2})}} \label{lemma6.24}
We now consider the corresponding rings $H_i:=\mathcal{O}_{X_i}$ for all $i=0,1,...$, and in our current situation we consider the corresponding rings $H_i\widehat{\otimes}Z$ (over $\mathbb{Z}_p$) for all $i=0,1,...$. And in our situation we consider the corresponding modules $M_i$ over $H_i\widehat{\otimes}Z$ for all $i=0,1,...$ with the same requirement as in \cite[Lemma 2.6.3]{8KL2} (namely those complete with respect to the natural topology). Suppose that we have bounded surjective map from $f_i:H_{i}\widehat{\otimes}Z\otimes_{H_{i+1}\widehat{\otimes}Z} M_{i+1}\rightarrow M_i,i=0,1,...$. Then we have first the density of the corresponding image of $\varprojlim_i M_i$ in each $M_i$ for any $i=0,1,2,...$. And we have as well the corresponding vanishing of $R^1\varprojlim_i M_i$.
\end{lemma}

\begin{proof}
This is the $Z$-relative version of the result in \cite[Lemma 2.6.3]{8KL2}. For the first statement we just choose sequence of Banach norms on all the corresponding modules for all $i=0,1,...$ such that we have $\|f_i(x_{i+1})\|_i\leq 1/2\|x_{i+1}\|_{i+1}$ for any $x_{i+1}\in M_{i+1}$. Then for any $x_i\in M_i$ and any $\delta>0$, we now consider for any $j\geq 1$ the corresponding $x_{i+j}$ such that we have $\|x_{i+j}-f_{i+j+1}(x_{i+j+1})\|_{i+j}\leq \delta$. Then the sequence $x_{i+j+k},k=0,1,...$ will converge to some well-defined $y_{i+j}$ with in our situation the corresponding $y_{i}=f_{i}(y_{i+1})$. We then have $\|x_i-y_i\|_i\leq \delta$. This will prove the first statement. For the second statement as in \cite[Lemma 2.6.3]{8KL2} we form the product $M_0\times M_1\times M_2\times...$ and the consider the induced map $F$ from $M_{i+1}\rightarrow M_i$, and consider the corresponding cokernel of the map $1-F$ since this is just the corresponding limit we are considering. Then to show that the cokernel is zero we just look at the corresponding cokernel of the corresponding map on the corresponding completed direct summand which will project to the original one. But then we will have $\|f_i(v)\|_i\leq 1/2 \|v\|_i$, which produces an inverse to $1-F$ which will basically finish the proof for the second statement. 	
\end{proof}

\begin{proposition} \mbox{\bf{(After Kedlaya-Liu \cite[Lemma 2.6.4]{8KL2})}} In the same situation as above, suppose we have that the corresponding modules $M_i$ are basically $Z$-stably pseudocoherent over the rings $H_i\widehat{\otimes}Z$ for all $i=0,1,...$. Now we consider the situation where $f_i:H_i\widehat{\otimes}Z\otimes_{H_{i+1}\widehat{\otimes}Z}M_{i+1}\rightarrow M_i$ is an isomorphism. Then the conclusion in our situation is then that the corresponding projection from $\varprojlim M_i$ to $M_i$ for each $i=0,1,2,...$ is an isomorphism.

\end{proposition}

\begin{proof}
This is a $Z$-relative and analytic version of the \cite[Lemma 2.6.4]{8KL2}. We adapt the argument to our situation as in the following. First we choose some finite free covering $T$ of the limit $M$ such that we have for each $i$ the corresponding map $T_i\rightarrow M_i$ is surjective. Then we consider the index $j\geq i$ and set the kernel of the map from $T_j$ to $M_j$ to be $S_j$. By the direct analog of \cite[Lemma 2.5.6]{8KL2} we have that $H_i\widehat{\otimes}Z\otimes S_j\overset{}{\rightarrow}S_i$ realizes the isomorphism, and we have that the corresponding surjectivity of the corresponding map from $\varprojlim_i S_i$ projecting to the $S_i$. Then one could finish the proof by 5-lemma to the following commutative diagram as in \cite[Lemma 2.6.4]{8KL2}:
\[ \tiny
\xymatrix@C+0pc@R+6pc{
 &(H_i\widehat{\otimes}Z)\otimes \varprojlim_i S_i \ar[r]\ar[r]\ar[r] \ar[d]\ar[d]\ar[d] &(H_i\widehat{\otimes}Z)\otimes \varprojlim_i F_i \ar[r]\ar[r]\ar[r] \ar[d]\ar[d]\ar[d] &(H_i\widehat{\otimes}Z)\otimes \varprojlim_i M_i\ar[r]\ar[r]\ar[r]  \ar[d]\ar[d]\ar[d]&0,\\
0 \ar[r]\ar[r]\ar[r] &S_i   \ar[r]\ar[r]\ar[r] &F_i \ar[r]\ar[r]\ar[r] &M_i \ar[r]\ar[r]\ar[r] &0.\\
}
\]
 
\end{proof}

\begin{proposition}  \mbox{\bf{(After Kedlaya-Liu \cite[Theorem 2.6.5]{8KL2})}}
For any\\ quasi-compact adic affinoid space of $X$ which is denoted by $Y$, we have that the map $\mathcal{M}(X)\rightarrow \mathcal{M}(Y)$ is surjective for any $Z$-stably pseudocoherent sheaf $\mathcal{M}$ over the sheaf $\mathcal{O}_X\widehat{\otimes}Z$.	
\end{proposition}

\begin{proof}
This is just the corresponding corollary of the previous proposition.	
\end{proof}

\begin{proposition}  \mbox{\bf{(After Kedlaya-Liu \cite[Theorem 2.6.5]{8KL2})}}
We have that the stalk $\mathcal{M}_x$ is generated over the stalk $\mathcal{O}_{X,x}$ for any $x\in X$ by $M(X)$, for any $Z$-stably pseudocoherent sheaf $\mathcal{M}$ over the sheaf $\mathcal{O}_X\widehat{\otimes}Z$.	
\end{proposition}

\begin{proof}
This is just the corresponding corollary of the proposition before the previous proposition.	
\end{proof}

\begin{proposition}  \mbox{\bf{(After Kedlaya-Liu \cite[Theorem 2.6.5]{8KL2})}} \label{proposition6.28}
For any\\ quasi-compact adic affinoid space of $X$ which is denoted by $Y$, we have that the corresponding vanishing of the corresponding sheaf cohomology groups $H^k(X,\mathcal{M})$ of $\mathcal{M}$ for higher $k>0$, for any $Z$-stably pseudocoherent sheaf $\mathcal{M}$ over the sheaf $\mathcal{O}_X\widehat{\otimes}Z$.	
\end{proposition}

\begin{proof}
We follow the idea of the proof of \cite[Theorem 2.6.5]{8KL2} by comparing this to the corresponding \v{C}ech cohomology with some covering $\mathfrak{X}=\{X_1,...,X_N,...\}$:
\begin{align}
\breve{H}^k(X,\mathfrak{X}=\{X_1,...,X_N,...\};\mathcal{M})=H^k(X,\mathcal{M}),\\
\breve{H}^k(X_i,\mathfrak{X}=\{X_1,...,X_i\};\mathcal{M})=H^k(X_i,\mathcal{M})=0,k\geq 1.
\end{align}
Here we have applied the corresponding \cite[Tag 01EW]{8SP}. Now we consider the situation where $k>1$:
\begin{align}
\breve{H}^{k-1}(X_{j+1},\mathfrak{X}=\{X_1,...,X_{j+1}\};\mathcal{M})	\rightarrow \breve{H}^{k-1}(X_j,\mathfrak{X}=\{X_1,...,X_{j}\};\mathcal{M})\rightarrow 0,
\end{align}
which induces the following isomorphism by \cite[2.6 Hilfssatz]{8Kie1}:
\begin{align}
\varprojlim_{j\rightarrow \infty} \breve{H}^{k-1}(X_j,\mathfrak{X}=\{X_1,...,X_{j}\};\mathcal{M})	\overset{\sim}{\rightarrow} \breve{H}^k(X,\mathfrak{X}=\{X_1,...,X_N,...\};\mathcal{M}). 
\end{align}
Then we have the corresponding results for the index $k>1$. For $k=1$, one can as in \cite[Theorem 2.6.5]{8KL2} relate this to the corresponding $R^1\varprojlim_i$, which will finishes the proof in the same fashion.
\end{proof}

\

\begin{corollary}  \mbox{\bf{(After Kedlaya-Liu \cite[Corollary 2.6.6]{8KL2})}}
The corresponding functor from the corresponding $Z$-deformed pseudocoherent sheaves over $X$ to the corresponding $Z$-stably by taking the corresponding global section is an exact functor.
\end{corollary}

\begin{corollary}  \mbox{\bf{(After Kedlaya-Liu \cite[Corollary 2.6.8]{8KL2})}}
Consider a particular $\mathcal{O}_X\widehat{\otimes}Z$-pseudocoherent sheaf $\mathcal{M}$ which is finite locally free throughout the whole space $X$. Then we have that the global section $\mathcal{M}(X)$ as $\mathcal{O}_X(X)\widehat{\otimes}Z$ left module admits the corresponding structures of finite projective structure if and only if we have the corresponding global section is finitely generated.
\end{corollary}

\begin{proof}
As in \cite[Corollary 2.6.8]{8KL2} one could find some global splitting through the local splittings.
\end{proof}

\

\indent We now consider the following $Z$-relative and analytic analog of \cite[Proposition 2.6.17]{8KL2}:

\begin{theorem} \mbox{\bf{(After Kedlaya-Liu \cite[Proposition 2.6.17]{8KL2})}} \label{theorem6.31}
Consider the following two statements for a particular $\mathcal{O}_X\widehat{\otimes}Z$-pseudocoherent sheaf $\mathcal{M}$. First is that one can find finite many generators (the number is up to some uniform fixed integer $n\geq 0$) for each section of $\mathcal{M}(X_i)$ for each $i=1,2,...$. The second statement is that the global section $\mathcal{M}(X)$ is just finitely generated. Then in our situation the two statement is equivalent if we have that the corresponding space $X$ admits an $m$-uniform covering in the exact same sense of \cite[Proposition 2.6.17]{8KL2}. 	
\end{theorem}

\begin{proof}
One direction is obvious, the other direct could be proved in the same way as in \cite[Proposition 2.6.17]{8KL2} where essentially the information on $X$ does not change at all. Namely as in \cite[Proposition 2.6.17]{8KL2} we could basically consider one single subcovering $\{Y_u\}_{u\in U}$ indexed by $U$ of the covering which is $m$-uniform. For any $u\in U$ we consider the smallest $i$ such that $X_i\bigcap Y_u$, and we denote this by $i(u)\geq 0$. Then we form:
\begin{align}
V_{u}:=X_{i(u)}\bigcup_{v\in U,v\leq u} Y_{v}.	
\end{align}
Then as in \cite[Proposition 2.6.17]{8KL2} we can find some $x_u,u\in U$ (which are denoted by the general form of $x_i$ in \cite[Proposition 2.6.17]{8KL2}) in the section $\mathcal{O}_X(V_{u})$ such that it restricts to some $0$ onto the space $V_u$ and it restricts to some element $a_1$ such that we have:
\begin{align}
a_1b_1+...+a_kb_k=1	
\end{align}
where $b_1,...,b_k$ are the corresponding topologically nilpotent elements to $Y_u$. By the corresponding density of the corresponding global section $\mathcal{O}_X(X)$ in $\mathcal{O}_X(V_{u})$ we could find some $x_i$ (with the same notation) in the global section $\mathcal{O}_X(X)$ such that it restricts to some topologically nilpotent element onto the space $V_u$ and it restricts to some element $a'_1$ such that we have:
\begin{align}
a'_1b'_1+...+a'_kb'_k=1	
\end{align}
where $b_1,...,b_k$ are the corresponding topologically nilpotent elements to $Y_u$, by simultaneously approximating the corresponding the coefficient $a_1,...,a_k$ and $b_1,...,b_k$ to achieve $a'_1,...,a'_k$ and $b'_1,...,b'_k$ such that we have:
\begin{align}
a'_1b'_1+...+a'_kb'_k=1.	
\end{align}
Then we can build up a corresponding generating set for the global section by the following approximation process just as in \cite[Proposition 2.6.17]{8KL2}. To be more precise what we have consider is to modify the corresponding generator for each $Y_u$ (for instance let them be denoted by $y_1,...,y_u$) to be $y_{u,1},...,y_{u,n}$ ($n$ is the corresponding uniform integer in the statement of the theorem) for all $u\in U$. It is achieved through induction, when we have $y_{u,1},...,y_{u,n}$ then we will set that to be $0$ otherwise in the situation where there exists some predecessor $u'$ of $u$ in the corresponding set $U$ we then set this to be $y_{u',1},...,y_{u',n}$. And we set as in \cite[Proposition 2.6.17]{8KL2}:
\begin{displaymath}
y_{u,j}:=	y_{u-1,j}+x_u^cy_u
\end{displaymath}
by lifting the corresponding power $c$ to be as large as possible. This will guarantee the convergence of:
\begin{align}
\lim_{u\rightarrow \infty}\{y_{u,1},...,y_{u,n}\},	
\end{align}
which gives the set of the corresponding global generators desired.
\end{proof}

%%\newpage

\newpage\section{Extensions and Applications}

\subsection{Deforming the $(\infty,1)$-Sheaves over Bambozzi-Kremnizer $\infty$-Analytic Spaces}

\indent Here we follow \cite{8BK1} to extend our discussion above partially to the corresponding $\infty$-analytic space attached to a Banach adic uniform algebra $A$ (required to be commutative) over $\mathbb{F}_p((t))$ or $\mathbb{Q}_p$ and we consider furthermore a Banach algebra $B$ over $\mathbb{F}_p((t))$ or $\mathbb{Q}_p$.  

\begin{remark}
We should definitely mention that one can apply the corresponding foundation from Clausen-Scholze \cite{8CS} to do this, namely the corresponding animated objects. But we will only apply the corresponding foundation in \cite{8BK1}. 
\end{remark}

\indent First we have the following $\infty$-version of Tate-\v{C}ech acyclicity in our context:

\begin{proposition}\mbox{\bf{(After Bambozzi-Kremnizer \cite[Theorem 4.15]{8BK1})}} 
The following totalization is strictly exact:
\[
\xymatrix@C+0pc@R+0pc{
\mathrm{Tot}(0  \ar[r]\ar[r]\ar[r] &A^h\widehat{\otimes} B \ar[r]\ar[r]\ar[r] & \prod_i A_i\widehat{\otimes} B \ar[r]\ar[r]\ar[r] &\prod_{i,j} A_i\widehat{\otimes} B \widehat{\otimes}^\mathbb{L} A_j\widehat{\otimes}B\ar[r]\ar[r]\ar[r] &...)
}
\]
for any homotopy Zariski covering $\{\mathrm{Spec}A_i\rightarrow \mathrm{Spec}A^h\}_{i\in I}$ in the sense of \cite[Theorem 4.15]{8BK1}. Here we have already attached the corresponding homotopical ring $A^h$ to $A$.
\end{proposition}

\begin{proof}
This is a direct consequence of the corresponding result in \cite[Theorem 4.15]{8BK1}.	
\end{proof}

\indent Along this idea one can then build up a sheaf of $\infty$-Banach ring over Bambozzi-Kremnizer's spectrum $\mathrm{Spa}^h(A)$, which we will denote it by $\mathcal{O}_{\mathrm{Spa}^h(A),B}$. Namely we have defined a new $\infty$-ringed space $(\mathrm{Spa}^h(A),\mathcal{O}_{\mathrm{Spa}^h(A),B})$.

\begin{definition} \mbox{\bf{(After Lurie \cite[Definition 2.9.1.1]{8Lu1})}} We will call a finite projective module spectrum $M$ over the ring $A^h\widehat{\otimes}B$ finite locally free, namely it will be basically sit in some finite coproduct of copies of $A^h\widehat{\otimes}B$ as a direct summand. Note that one can obviously construct a retract from some finite free module spectrum to such object, which is surjective namely with respect to the homotopy group $\pi_0$.
	
\end{definition}

\indent As in \cite[Proposition 7.2.4.20]{8Lu2}, we have that finite locally free sheaves are actually flat and pseudocoherent in the strong $\infty$-sense. This will cause the following to be true:

\begin{proposition} 
The following totalization is strictly exact:
\[
\xymatrix@C+0pc@R+0pc{
\mathrm{Tot}(0  \ar[r]\ar[r]\ar[r] &M \ar[r]\ar[r]\ar[r] & \prod_i A_i\widehat{\otimes} B \widehat{\otimes}^\mathbb{L} M  \ar[r]\ar[r]\ar[r] &\prod_{i,j} (A_i\widehat{\otimes} B \widehat{\otimes}^\mathbb{L} A_j\widehat{\otimes}B) \widehat{\otimes}^\mathbb{L} M\ar[r]\ar[r]\ar[r] &...)
}
\]
for any homotopy Zariski covering $\{\mathrm{Spec}A_i\rightarrow \mathrm{Spec}A^h\}_{i\in I}$ in the sense of \cite[Theorem 2.15, Theorem 4.15]{8BK1}. 
\end{proposition}

\begin{theorem} \mbox{\bf{(After Kedlaya-Liu \cite[Theorem 2.7.7]{8KL1})}}
Carrying the corresponding coefficient in noncommutative Banach ring $B$, we have that $\infty$-descent for locally free finitely generated module spectra with respect to simple Laurent covering. 	
\end{theorem}

\begin{proof} 
\indent Certainly along the idea of \cite[Theorem 2.7.7]{8KL1} one can actually try to consider the corresponding glueing of finite projective objects, but this will basically involve a derived or $\infty$ version of \cite[Proposition 2.4.20]{8KL1}. Instead we consider the following restricted situation which is also considered restrictively in \cite[Section 5]{8T3}. First we consider the corresponding glueing datum along $A$, namely suppose we have the corresponding short exact sequence of Banach adic uniform algebra over $\mathbb{Q}_p$:
\[
\xymatrix@C+0pc@R+0pc{
0   \ar[r]\ar[r]\ar[r] &\Pi \ar[r]\ar[r]\ar[r] &\Pi_1\bigoplus \Pi_2 \ar[r]\ar[r]\ar[r] & \Pi_{12} \ar[r]\ar[r]\ar[r] &0,
}
\]
which is assumed to satisfy conditions $(a)$ and $(b)$ of \cite[Definition 2.7.3]{8KL1}, namely this is strictly exact and the map $\Pi_1\rightarrow \Pi_{12}$ is of image that is dense. Now we consider the Bambozzi-Kremnizer spaces and rings:
\begin{align}
&\mathrm{Spa}^h(\Pi),\mathrm{Spa}^h(\Pi_1),\mathrm{Spa}^h(\Pi_2),\mathrm{Spa}^h(\Pi_{12}),\\
&\Pi^h,\Pi_1^h,\Pi_2^h,\Pi_{12}^h.  	
\end{align}
Now consider the following finite locally free bimodule spectra forming the corresponding glueing datum:
\begin{align}
M_1,M_2,M_{12}	
\end{align}
over:
\begin{align}
\Pi_1^h\widehat{\otimes}B,\Pi_2^h\widehat{\otimes}B,\Pi_{12}^h\widehat{\otimes}B.
\end{align}
Then we can regard these as certain sheaves of bimodules:
\begin{align}
\mathcal{M}_1,\mathcal{M}_2,\mathcal{M}_{12}	
\end{align}
over the $\infty$-sheaves:
\begin{align}
\mathcal{O}_{\Pi_1^h}\widehat{\otimes}B,\mathcal{O}_{\Pi_2^h}\widehat{\otimes}B,\mathcal{O}_{\Pi_{12}^h}\widehat{\otimes}B.
\end{align}
Then we can have the chance to take the global section along the space $\mathrm{Spa}^h(\Pi)$ which gives rise to a certain module spectrum $M$ over $\Pi^h$. However the corresponding $\pi_0(M)$ will then be finite projective over $\pi_0(\Pi)$ by \cite[Proposition 5.12]{8T3}. This finishes the corresponding glueing of vector bundle process in this specific situation.
\end{proof}

\begin{remark}
Again one can apply the whole machinery from \cite{8CS} to achieve this, as long as one would like to work with noncommutative analogs of the animated rings.	
\end{remark}

%\begin{theorem}   \mbox{\bf{(After Kedlaya-Liu \cite[Theorem 2.7.7]{8KL1})}}
%Taking the corresponding global section will establish a corresponding well-defined functor from the category of the corresponding locally finite projective sheaves over $\mathcal{O}_{\mathrm{Spa}^h(A),B}$ to the category of all the corresponding finite locally free modules over $A^h\widehat{\otimes}B$, which just establishes an equivalence.	
%\end{theorem}
%
%\begin{proof}
%See \cite[Theorem 2.7.7]{8KL1}.	\\
%\end{proof}

\

\indent Now we consider $\mathbb{Q}_p[[G]]$-adic objects in the situation over $\mathbb{Q}_p$. First we have the following $\infty$-version of Tate-\v{C}ech acyclicity in our context:

\begin{proposition}\mbox{\bf{(After Bambozzi-Kremnizer \cite[Theorem 4.15]{8BK1})}} 
The following totalization is strictly exact:
\[\tiny
\xymatrix@C+0pc@R+0pc{
\mathrm{Tot}(0  \ar[r]\ar[r]\ar[r] &A^h\widehat{\otimes} B[[G]]/E \ar[r]\ar[r]\ar[r] & \prod_i A_i\widehat{\otimes} B[[G]]/E \ar[r]\ar[r]\ar[r] &\prod_{i,j} A_i\widehat{\otimes} B[[G]]/E \widehat{\otimes} ^\mathbb{L} A_j\widehat{\otimes}B[[G]]/E \ar[r]\ar[r]\ar[r] &...),
}
\]
for any $E\subset G$, and for any homotopy Zariski covering $\{\mathrm{Spec}A_i\rightarrow \mathrm{Spec}A^h\}_{i\in I}$ in the sense of \cite[Theorem 2.15, Theorem 4.15]{8BK1}. Here we have already attached the corresponding homotopical ring $A^h$ to $A$.
\end{proposition}

\begin{proof}
This is a direct consequence of the corresponding result in \cite[Theorem 4.15]{8BK1}.	
\end{proof}

\indent Along this idea one can then build up a sheaf of $\infty$-Banach ring over Bambozzi-Kremnizer's spectrum $\mathrm{Spa}^h(A)$, which we will denote it by $\mathcal{O}_{\mathrm{Spa}^h(A),B[[G]]}$. Namely we have defined a new $\infty$-ringed space $(\mathrm{Spa}^h(A),\mathcal{O}_{\mathrm{Spa}^h(A),B[[G]]})$.

\begin{definition} \mbox{\bf{(After Lurie \cite[Definition 2.9.1.1]{8Lu1})}} We will call a finite projective module spectrum $M$ over the ring $A^h\widehat{\otimes}B[[G]]$ finite locally free, namely it will be basically sit in some finite coproduct of copies of $A^h\widehat{\otimes}B[[G]]$ as a direct summand. Note that one can obviously construct a retract from some finite free module spectrum to such object.
	
\end{definition}

\indent As in \cite[Proposition 7.2.4.20]{8Lu2}, we have that finite locally free sheaves are actually flat and pseudocoherent in the strong $\infty$-sense. This will cause the following to be true:

\begin{proposition} 
The following totalization is strictly exact:
\[\tiny
\xymatrix@C+0pc@R+0pc{
\mathrm{Tot}(0  \ar[r]\ar[r]\ar[r] &M \ar[r]\ar[r]\ar[r] & \prod_i A_i\widehat{\otimes} B[[G]]/E \widehat{\otimes}^\mathbb{L} M  \ar[r]\ar[r]\ar[r] &\prod_{i,j} (A_i\widehat{\otimes} B[[G]]/E \otimes^\mathbb{L} A_j\widehat{\otimes} B[[G]]/E) \otimes^\mathbb{L} M\ar[r]\ar[r]\ar[r] &...)
}
\]
for any homotopy Zariski covering $\{\mathrm{Spec}A_i\rightarrow \mathrm{Spec}A^h\}_{i\in I}$ in the sense of \cite[Theorem 4.15, Theorem 2.15]{8BK1}. 
\end{proposition}

%\begin{theorem} \mbox{\bf{(After Kedlaya-Liu \cite[Theorem 2.7.7]{8KL1})}}
%Taking the corresponding global section will establish a corresponding well-defined functor from the category of the corresponding locally finite projective sheaves over $\mathcal{O}_{\mathrm{Spa}^h(A),B[[G]]}$ to the category of all the corresponding finite locally free modules over $A^h\widehat{\otimes}B[[G]]$, which just establishes an equivalence.	
%\end{theorem}
%
%
%\begin{proof}
%See \cite[Theorem 2.7.7]{8KL1}.		
%\end{proof}

\

 \indent Now we consider the LF deformation, namely here $B$ could be written as $\varinjlim_h B_h$. First we have the following $\infty$-version of Tate-\v{C}ech acyclicity in our context:

\begin{proposition}\mbox{\bf{(After Bambozzi-Kremnizer \cite[Theorem 4.15]{8BK1})}} 
The following totalization is strictly exact:
\[
\xymatrix@C+0pc@R+0pc{
\mathrm{Tot}(0  \ar[r]\ar[r]\ar[r] &A^h\widehat{\otimes} B \ar[r]\ar[r]\ar[r] & \prod_i A_i\widehat{\otimes} B \ar[r]\ar[r]\ar[r] &\prod_{i,j} A_i\widehat{\otimes} B \widehat{\otimes}^\mathbb{L} A_j\widehat{\otimes}B\ar[r]\ar[r]\ar[r] &...)
}
\]
for any homotopy Zariski covering $\{\mathrm{Spec}A_i\rightarrow \mathrm{Spec}A^h\}_{i\in I}$ in the sense of \cite[Theorem 4.15]{8BK1}. Here we have already attached the corresponding homotopical ring $A^h$ to $A$.
\end{proposition}

\begin{proof}
This is a direct consequence of the corresponding result in \cite[Theorem 4.15]{8BK1}.	
\end{proof}

\indent Along this idea one can then build up a sheaf of $\infty$-Banach ring over Bambozzi-Kremnizer's spectrum $\mathrm{Spa}^h(A)$, which we will denote it by $\mathcal{O}_{\mathrm{Spa}^h(A),B}$. Namely we have defined a new $\infty$-ringed space $(\mathrm{Spa}^h(A),\mathcal{O}_{\mathrm{Spa}^h(A),B})$.

\begin{definition} \mbox{\bf{(After Lurie \cite[Definition 2.9.1.1]{8Lu1})}} We will call a finite projective module spectrum $M$ over the ring $A^h\widehat{\otimes}B$ finite locally free, namely it will be basically sit in some finite coproduct of copies of $A^h\widehat{\otimes}B$ as a direct summand. Note that one can obviously construct a retract from some finite free module spectrum to such object.
	
\end{definition}

\indent As in \cite[Proposition 7.2.4.20]{8Lu2}, we have that finite locally free sheaves are actually flat and pseudocoherent in the strong $\infty$-sense. This will cause the following to be true:

\begin{proposition} 
The following totalization is strictly exact:
\[
\xymatrix@C+0pc@R+0pc{
\mathrm{Tot}(0  \ar[r]\ar[r]\ar[r] &M \ar[r]\ar[r]\ar[r] & \prod_i A_i\widehat{\otimes} B \widehat{\otimes}^\mathbb{L} M  \ar[r]\ar[r]\ar[r] &\prod_{i,j} (A_i\widehat{\otimes} B \widehat{\otimes}^\mathbb{L} A_j\widehat{\otimes}B) \widehat{\otimes}^\mathbb{L} M\ar[r]\ar[r]\ar[r] &...)
}
\]
for any homotopy Zariski covering $\{\mathrm{Spec}A_i\rightarrow \mathrm{Spec}A^h\}_{i\in I}$ in the sense of \cite[Theorem 2.15, Theorem 4.15]{8BK1}. 
\end{proposition}

\begin{theorem} \mbox{\bf{(After Kedlaya-Liu \cite[Theorem 2.7.7]{8KL1})}}
Carrying the corresponding coefficient in the noncommutative Limit of Fr\'echet $B$, we have that $\infty$-descent for locally free finitely generated module spectra with respect to simple Laurent covering.	
\end{theorem}

\begin{proof} 
\indent Along the discussion in the Banach situation, we can also make the parallel discussion in the corresponding limit of Fr\'echet situation (LF). Certainly along the idea of \cite[Theorem 2.7.7]{8KL1} one can actually try to consider the corresponding glueing of finite projective objects, but this will basically involve a derived or $\infty$ version of \cite[Proposition 2.4.20]{8KL1}. Instead we consider the following restricted situation which is also considered restrictively in \cite[Section 5]{8T3}. First we consider the corresponding glueing datum along $A$, namely suppose we have the corresponding short exact sequence of Banach adic uniform algebra over $\mathbb{Q}_p$:
\[
\xymatrix@C+0pc@R+0pc{
0   \ar[r]\ar[r]\ar[r] &\Pi \ar[r]\ar[r]\ar[r] &\Pi_1\bigoplus \Pi_2 \ar[r]\ar[r]\ar[r] & \Pi_{12} \ar[r]\ar[r]\ar[r] &0,
}
\]
which is assumed to satisfy conditions $(a)$ and $(b)$ of \cite[Definition 2.7.3]{8KL1}, namely this is strictly exact and the map $\Pi_1\rightarrow \Pi_{12}$ is of image that is dense. Now we consider the Bambozzi-Kremnizer spaces and rings:
\begin{align}
&\mathrm{Spa}^h(\Pi),\mathrm{Spa}^h(\Pi_1),\mathrm{Spa}^h(\Pi_2),\mathrm{Spa}^h(\Pi_{12}),\\
&\Pi^h,\Pi_1^h,\Pi_2^h,\Pi_{12}^h.  	
\end{align}
Now consider the following finite locally free bimodule spectra forming the corresponding glueing datum:
\begin{align}
M_1,M_2,M_{12}	
\end{align}
over:
\begin{align}
\Pi_1^h\widehat{\otimes}\varinjlim_h B_h,\Pi_2^h\widehat{\otimes}\varinjlim_h B_h,\Pi_{12}^h\widehat{\otimes}\varinjlim_h B_h.
\end{align}
Then we can regard these as certain sheaves of bimodules:
\begin{align}
\mathcal{M}_1,\mathcal{M}_2,\mathcal{M}_{12}	
\end{align}
over the $\infty$-sheaves:
\begin{align}
\mathcal{O}_{\Pi_1^h}\widehat{\otimes}\varinjlim_h B_h,\mathcal{O}_{\Pi_2^h}\widehat{\otimes}\varinjlim_h B_h,\mathcal{O}_{\Pi_{12}^h}\widehat{\otimes}\varinjlim_h B_h.
\end{align}
Then we can have the chance to take the global section along the space $\mathrm{Spa}^h(\Pi)$ which gives rise to a certain module spectrum $M$ over $\Pi^h$. However the corresponding $\pi_0(M)$ will then be finite projective over $\pi_0(\Pi)$ by \cite[Proposition 5.12]{8T3}. This finishes the corresponding glueing of vector bundle process in this specific situation.
\end{proof}

\begin{remark}
Again one can apply the whole machinery from \cite{8CS} to achieve this, as long as one would like to work with noncommutative analogs of the animated rings.	
\end{remark}

\
%
%\begin{theorem}   \mbox{\bf{(After Kedlaya-Liu \cite[Theorem 2.7.7]{8KL1})}}
%Taking the corresponding global section will establish a corresponding well-defined functor from the category of the corresponding locally finite projective sheaves over $\mathcal{O}_{\mathrm{Spa}^h(A),B}$ to the category of all the corresponding finite locally free modules over $A^h\widehat{\otimes}B$, which just establishes an equivalence.	
%\end{theorem}
%
%\begin{proof}
%See \cite[Theorem 2.7.7]{8KL1}.	\\
%\end{proof}

\indent The above discussion could be made to be more general by looking at some derived rational localization $A\rightarrow D$ of $A$ as in \cite{8BK1} which we will denote by $D\in \infty-\mathbf{Comm}$.  For instance one can just take $D$ to be some derived rational localization taking the form such as $A\left<\frac{a_0}{a_0},...,\frac{a_k}{a_0}\right>^h$ (and consider the glueing over $\mathrm{Spec}A\left<\frac{a_0}{a_0},...,\frac{a_k}{a_0}\right>^h\in \infty-\mathbf{Comm}^\mathrm{op}$). We now assume the following assumption:

\begin{assumption}
Without any sheafiness condition on the ring $A$ in the corresponding classical sense, we assume the following totolization is strictly exact:
\[
\xymatrix@C+0pc@R+0pc{
\mathrm{Tot}(0  \ar[r]\ar[r]\ar[r] &D\widehat{\otimes} B \ar[r]\ar[r]\ar[r] & \prod_i D_i\widehat{\otimes} B \ar[r]\ar[r]\ar[r] &\prod_{i,j} D_i\widehat{\otimes} B \widehat{\otimes}^\mathbb{L} D_j\widehat{\otimes}B\ar[r]\ar[r]\ar[r] &...)
}
\]	
for any homotopy Zariski covering $\{\mathrm{Spec}D_i\rightarrow \mathrm{Spec}D\}_{i\in I}$ in the sense of \cite[Theorem 4.15]{8BK1}. Here recall the notation $\mathrm{Spec}D$ means the corresponding object in $\infty-\mathbf{Comm}^\mathrm{op}$. Here $B$ is just assumed to be Banach.
\end{assumption}

\begin{remark}
\indent By the corresponding Dold-Kan correspondence we could regard $D$ as an $(\infty,1)$-commutative object in the corresponding category of all the simplicial semi-normed monoids, namely an $E_\infty$-ring which carries the corresponding simplicial Banach structures. However then the corresponding ring $D\widehat{\otimes}B$ will then be a product of an $E_\infty$-ring with a discrete (in the sense of \cite[Chapter 7.2]{8Lu2}) $E_1$-ring (or if you wish a discrete $A_\infty$-ring).
\end{remark}

\begin{definition} \mbox{\bf{(After Lurie \cite[Definition 2.9.1.1]{8Lu1})}} We will call a finite projective module spectrum $M$ over the ring $D\widehat{\otimes}B$ finite locally free, namely it will be basically sit in some finite coproduct of copies of $D\widehat{\otimes}B$ as a direct summand. Note that one can obviously construct a retract from some finite free module spectrum to such object, which is surjective namely with respect to the homotopy group $\pi_0$.
	
\end{definition}

\begin{definition} \mbox{\bf{(After Lurie \cite[Definition 2.9.1.1]{8Lu1})}} We will call a sheaf $M$ over the $\infty$-ringed Grothendieck site:
\begin{align}
(\mathrm{Spec}D,\mathcal{O}_\mathrm{Spec}D\widehat{\otimes}B,\{\mathrm{Spec}D\left<\frac{f_0}{f_i},...,\frac{f_d}{f_i}\right>^h\}_{i=0,...,d})	
\end{align}
locally finite free if it is locally (over some derived rational subset $\mathrm{Spec}D_i\in \infty-\mathbf{Comm}^\mathrm{op}$ and $D_i:=D\left<\frac{f_0}{f_i},...,\frac{f_d}{f_i}\right>^h,i=1,...,d$) sitting in some coproduct of finite many copies of $D_i\widehat{\otimes}B$ as a direct summand. 
	
\end{definition}

\indent As in \cite[Proposition 7.2.4.20]{8Lu2}, we have that finite locally free sheaves are actually flat and pseudocoherent in the strong $\infty$-sense. This will cause the following to be true:

\begin{proposition} 
The following totalization is strictly exact:
\[
\xymatrix@C+0pc@R+0pc{
\mathrm{Tot}(0  \ar[r]\ar[r]\ar[r] &M \ar[r]\ar[r]\ar[r] & \prod_i D_i\widehat{\otimes} B \widehat{\otimes}^\mathbb{L} M  \ar[r]\ar[r]\ar[r] &\prod_{i,j} (D_i\widehat{\otimes} B \widehat{\otimes}^\mathbb{L} D_j\widehat{\otimes}B) \widehat{\otimes}^\mathbb{L} M\ar[r]\ar[r]\ar[r] &...)
}
\]
for any homotopy Zariski covering $\{\mathrm{Spec}D_i\rightarrow \mathrm{Spec}D\}_{i\in I}$ in the sense of \cite[Theorem 2.15, Theorem 4.15]{8BK1}. 
\end{proposition}

\indent Now we take a corresponding fiber sequence (under the corresponding Dold-Kan correspondence one can also consider the corresponding distinguished triangles 
\[
\xymatrix@C+0pc@R+0pc{
&D \ar[r]\ar[r]\ar[r] &\prod_{i=1,2} D_i \ar[r]\ar[r]\ar[r] & D_{1}\widehat{\otimes}^\mathbb{L} D_{2} \ar[r]\ar[r]\ar[r] & D[1],
}
\]
for the corresponding associated complexes):
\[
\xymatrix@C+0pc@R+0pc{
&D \ar[r]\ar[r]\ar[r] &\prod_{i=1,2} D_i \ar[r]\ar[r]\ar[r] & D_{12},
}
\]
and assume that the following:
\[
\xymatrix@C+0pc@R+0pc{
0   \ar[r]\ar[r]\ar[r] &\pi_0D \ar[r]\ar[r]\ar[r] &\pi_0D_1\bigoplus \pi_0D_2 \ar[r]\ar[r]\ar[r] & \pi_0D_{12} \ar[r]\ar[r]\ar[r] &0
}
\]
to satisfy \cite[Definition 2.7.3, condition (a), (b)]{8KL1}. Then we have the effective $\infty$-descentness of the corresponding morphism $D\widehat{\otimes}B\rightarrow D_1\widehat{\otimes}B\oplus D_2\widehat{\otimes}B$.

\

\begin{theorem} \mbox{\bf{(After Kedlaya-Liu \cite[Theorem 2.7.7]{8KL1})}}
Carrying the corresponding coefficient in noncommutative Banach ring $B$, we have that $\infty$-descent for locally free finitely generated module spectra with respect to binary covering on the corresponding $\pi_0$ in the situation described above.	
\end{theorem}

\begin{proof} 
\indent Certainly along the idea of \cite[Theorem 2.7.7]{8KL1} one can actually try to consider the corresponding glueing of finite projective objects in more general setting, but this will basically involve a derived or $\infty$ version of \cite[Proposition 2.4.20]{8KL1}. Instead we consider the following restricted situation which is also considered restrictively in \cite[Section 5]{8T3}. First we consider the corresponding glueing datum along $\pi_0D$, namely suppose we have the corresponding short exact sequence of Banach adic uniform algebra over $\mathbb{Q}_p$:
\[
\xymatrix@C+0pc@R+0pc{
0   \ar[r]\ar[r]\ar[r] &\pi_0D \ar[r]\ar[r]\ar[r] &\pi_0D_1\bigoplus \pi_0D_2 \ar[r]\ar[r]\ar[r] & \pi_0D_{12} \ar[r]\ar[r]\ar[r] &0
}
\]
which is assumed to satisfy conditions $(a)$ and $(b)$ of \cite[Definition 2.7.3]{8KL1}, namely this is strictly exact and the map $\Pi_1\rightarrow \Pi_{12}$ is of image that is dense. Now consider the following finite locally free bimodule spectra forming the corresponding glueing datum:
\begin{align}
M_1,M_2,M_{12}	
\end{align}
over:
\begin{align}
D_1\widehat{\otimes}B,D_2\widehat{\otimes}B,D_{12}\widehat{\otimes}B.
\end{align}
Then we can regard these as certain sheaves of bimodules:
\begin{align}
\mathcal{M}_1,\mathcal{M}_2,\mathcal{M}_{12}	
\end{align}
over the $\infty$-sheaves:
\begin{align}
\mathcal{O}_{\mathrm{Spec}D_1}\widehat{\otimes}B,\mathcal{O}_{\mathrm{Spec}D_2}\widehat{\otimes}B,\mathcal{O}_{\mathrm{Spec}D_{12}}\widehat{\otimes}B.
\end{align}
Then we can have the chance to take the global section along the space $\mathrm{Spec}D$ which gives rise to a certain module spectrum $M$ over $D$. However the corresponding $\pi_0(M)$ will then be finite projective over $\pi_0(\Pi)$ by \cite[Proposition 5.12]{8T3}. This finishes the corresponding glueing of vector bundle process in this specific situation.
\end{proof}

\

\subsection{Application to Equivariant Sheaves over Adic Fargues-Font
aine Curves}

\indent We now apply our study to the first specific situation around the corresponding mixed-type Hodge structures inspired essentially by the work \cite{8CKZ}, \cite{8PZ} and \cite{8Ked1}, which is also considered in \cite{8T3}. Let us first describe the corresponding situation:

\begin{setting}
We will use the notations in \cite{8T1}, \cite{8T2} and \cite{8T3}. We now consider the corresponding Fargues-Fontaine curves as our target to relativization. We will use the corresponding $\mathrm{FF}_R$ to denote the adic Fargues-Fontaine curve in our current situation (as in \cite[Theorem 4.6.1]{8KL2}) which is defined to be the corresponding suitable quotient of the following adic space:
\begin{align}
\bigcup_{0<r_1<r_2} \mathrm{Spa}(\widetilde{\Pi}^{[r_1,r_2]}_R,\widetilde{\Pi}^{[r_1,r_2],+}_R).	
\end{align}
Here $R$ is specified perfectoid Banach uniform algebra over $\mathbb{F}_{p^h}$.
\end{setting}

\begin{assumption}
In this section we assume that the corresponding base field will be just of mixed-characteristic.	
\end{assumption}

\indent In our current framework we look at the following deformation:

\begin{setting}
Over $\mathbb{Q}_p$ we consider the any Galois perfectoid tower $(H_\bullet,H_\bullet^+)$ (in the sense of \cite[Definition 5.1.1]{8KL2}) with Galois group $\mathrm{Gal}_H:=\mathrm{Gal}(\overline{H}_0/H_0)$, and we use the notation $\widetilde{H}_\infty$ to denote the corresponding perfectoid top level of the tower in mixed-characteristic situation, and we consider the corresponding ring $\overline{H}_\infty$ as its tilt. Now we consider the period rings $B_{\mathrm{dR},\overline{H}_\infty}^+$, $B_{\mathrm{dR},\overline{H}_\infty}$ and $B_{e,\overline{H}_\infty}$ (as constructed in \cite[Definition 8.6.5]{8KL2}). 	
\end{setting}

\indent We now try to be more specific, where we just consider the situation where the ring $\widetilde{H}_\infty$ is just the ring $\mathbb{C}_p$. We now consider some Banach algebra $B$ over $\mathbb{Q}_p$ as before. Now what is happening is then we consider the following ringed space:

\begin{align}
(\mathrm{FF}_R,\mathcal{O}_{\mathrm{FF}_R}\widehat{\otimes}B\widehat{\otimes}B_{\mathrm{dR},\overline{H}_\infty}^+,\mathcal{O}_{\mathrm{FF}_R}\widehat{\otimes}B\widehat{\otimes}B_{\mathrm{dR},\overline{H}_\infty},\mathcal{O}_{\mathrm{FF}_R}\widehat{\otimes}B\widehat{\otimes}B_{e,\overline{H}_\infty}).	
\end{align}

\indent Now we encode Berger's $B$-pair in some relativization after Kedlaya-Liu (\cite[Definition 9.3.11]{8KL1} and \cite[Definition 4.8.2]{8KL2}) and Kedlaya-Pottharst (\cite[Definition 2.17]{8KP}):

\begin{definition}
We now define a pseudocoherent $\mathrm{FF}_R$-$B$-$\mathrm{Gal}_H$ sheaf to be a triplet $(M^+_\mathrm{dR},M_\mathrm{dR},M_e)$ of $B\widehat{\otimes}?$-stably pseudocoherent sheaves such that the form a corresponding descent triplet with respect to the following diagram by natural base change:
\begin{align}
\mathcal{O}_{\mathrm{FF}_R}\widehat{\otimes}B\widehat{\otimes}B_{\mathrm{dR},\overline{H}_\infty}^+ \longrightarrow \mathcal{O}_{\mathrm{FF}_R}\widehat{\otimes}B\widehat{\otimes}B_{\mathrm{dR},\overline{H}_\infty}\longleftarrow \mathcal{O}_{\mathrm{FF}_R}\widehat{\otimes}B\widehat{\otimes}B_{e,\overline{H}_\infty},	
\end{align}
where $?=B_{\mathrm{dR},\overline{H}_\infty}^+,B_{\mathrm{dR},\overline{H}_\infty},B_{e,\overline{H}_\infty}$. For the rings:
\begin{align}
\mathcal{O}_{\mathrm{FF}_R}\widehat{\otimes}B\widehat{\otimes}B_{e,\overline{H}_\infty}	
\end{align}
the corresponding definition of pseudocoherent sheaves are defined to be locally associated to the $B\widehat{\otimes}B_{e,\overline{H}_\infty}$-stably pseudocoherent modules in the ind-Fr\'echet topology, namely the completeness with respect to the ind-Fr\'echet topology is then stable under the corresponding rational localization concentrated at the part for $\mathrm{FF}_R$. And we assume that the corresponding sheaves admit $\mathrm{Gal}_H$-equivariant action in the semilinear fashion.	
\end{definition}

\indent In the \'etale topology we consider the following ringed space:

\begin{align}
(\mathrm{FF}_R,\mathcal{O}_{\mathrm{FF}_R,\text{\'et}}\widehat{\otimes}B\widehat{\otimes}B_{\mathrm{dR},\overline{H}_\infty}^+,\mathcal{O}_{\mathrm{FF}_R,\text{\'et}}\widehat{\otimes}B\widehat{\otimes}B_{\mathrm{dR},\overline{H}_\infty},\mathcal{O}_{\mathrm{FF}_R,\text{\'et}}\widehat{\otimes}B\widehat{\otimes}B_{e,\overline{H}_\infty}).	
\end{align}

\indent Now we encode Berger's $B$-pair in some relativization after Kedlaya-Liu (\cite[Definition 9.3.11]{8KL1} and \cite[Definition 4.8.2]{8KL2}) and Kedlaya-Pottharst (\cite[Definition 2.17]{8KP}):

\begin{definition}
We now define a pseudocoherent $\mathrm{FF}_R$-$B$-$\mathrm{Gal}_H$ sheaf to be a triplet $(M^+_\mathrm{dR},M_\mathrm{dR},M_e)$ of $B\widehat{\otimes}?$-\'etale-stably pseudocoherent sheaves such that these form a corresponding descent triplet with respect to the following diagram by natural base change:
\begin{align}
\mathcal{O}_{\mathrm{FF}_R,\text{\'et}}\widehat{\otimes}B\widehat{\otimes}B_{\mathrm{dR},\overline{H}_\infty}^+ \longrightarrow \mathcal{O}_{\mathrm{FF}_R,\text{\'et}}\widehat{\otimes}B\widehat{\otimes}B_{\mathrm{dR},\overline{H}_\infty}\longleftarrow \mathcal{O}_{\mathrm{FF}_R,\text{\'et}}\widehat{\otimes}B\widehat{\otimes}B_{e,\overline{H}_\infty},	
\end{align}
where $?=B_{\mathrm{dR},\overline{H}_\infty}^+,B_{\mathrm{dR},\overline{H}_\infty},B_{e,\overline{H}_\infty}$. For the rings:
\begin{align}
\mathcal{O}_{\mathrm{FF}_R,\text{\'et}}\widehat{\otimes}B\widehat{\otimes}B_{e,\overline{H}_\infty}	
\end{align}
the corresponding definition of pseudocoherent sheaves are defined to be locally associated to the $B\widehat{\otimes}B_{e,\overline{H}_\infty}$-\'etale-stably pseudocoherent modules in the ind-Fr\'echet topology, namely the completeness with respect to the ind-Fr\'echet topology is then stable under the corresponding \'etale morphisms concentrated at the part for $\mathrm{FF}_R$. And we assume that the corresponding sheaves admit $\mathrm{Gal}_H$-equivariant action in the semilinear fashion.	
\end{definition}

\begin{definition}
We now define a projective $\mathrm{FF}_R$-$B$-$\mathrm{Gal}_H$ sheaf to be a triplet\\ $(M^+_\mathrm{dR},M_\mathrm{dR},M_e)$ of locally finite projective sheaves such that these form a corresponding descent triplet with respect to the following diagram by natural base change:
\begin{align}
\mathcal{O}_{\mathrm{FF}_R}\widehat{\otimes}B\widehat{\otimes}B_{\mathrm{dR},\overline{H}_\infty}^+ \longrightarrow \mathcal{O}_{\mathrm{FF}_R}\widehat{\otimes}B\widehat{\otimes}B_{\mathrm{dR},\overline{H}_\infty}\longleftarrow \mathcal{O}_{\mathrm{FF}_R}\widehat{\otimes}B\widehat{\otimes}B_{e,\overline{H}_\infty},	
\end{align}
where $?=B_{\mathrm{dR},\overline{H}_\infty}^+,B_{\mathrm{dR},\overline{H}_\infty},B_{e,\overline{H}_\infty}$. For the rings:
\begin{align}
\mathcal{O}_{\mathrm{FF}_R}\widehat{\otimes}B\widehat{\otimes}B_{e,\overline{H}_\infty}	
\end{align}
the corresponding definition of projective sheaves are defined to be locally associated to the finite projective modules in the ind-Fr\'echet topology, namely the completeness with respect to the ind-Fr\'echet topology is then stable under the corresponding rational localization concentrated at the part for $\mathrm{FF}_R$. And we assume that the corresponding sheaves admit $\mathrm{Gal}_H$-equivariant action in the semilinear fashion.	
\end{definition}

\indent In the \'etale topology we consider the following ringed space:

\begin{align}
(\mathrm{FF}_R,\mathcal{O}_{\mathrm{FF}_R,\text{\'et}}\widehat{\otimes}B\widehat{\otimes}B_{\mathrm{dR},\overline{H}_\infty}^+,\mathcal{O}_{\mathrm{FF}_R,\text{\'et}}\widehat{\otimes}B\widehat{\otimes}B_{\mathrm{dR},\overline{H}_\infty},\mathcal{O}_{\mathrm{FF}_R,\text{\'et}}\widehat{\otimes}B\widehat{\otimes}B_{e,\overline{H}_\infty}).	
\end{align}

\indent Now we encode Berger's $B$-pair in some relativization after Kedlaya-Liu (\cite[Definition 9.3.11]{8KL1} and \cite[Definition 4.8.2]{8KL2}) and Kedlaya-Pottharst (\cite[Definition 2.17]{8KP}):

\begin{definition}
We now define a projective $\mathrm{FF}_R$-$B$-$\mathrm{Gal}_H$ sheaf to be a triplet\\ $(M^+_\mathrm{dR},M_\mathrm{dR},M_e)$ of locally finite projective sheaves such that these form a corresponding descent triplet with respect to the following diagram by natural base change:
\begin{align}
\mathcal{O}_{\mathrm{FF}_R,\text{\'et}}\widehat{\otimes}B\widehat{\otimes}B_{\mathrm{dR},\overline{H}_\infty}^+ \longrightarrow \mathcal{O}_{\mathrm{FF}_R,\text{\'et}}\widehat{\otimes}B\widehat{\otimes}B_{\mathrm{dR},\overline{H}_\infty}\longleftarrow \mathcal{O}_{\mathrm{FF}_R,\text{\'et}}\widehat{\otimes}B\widehat{\otimes}B_{e,\overline{H}_\infty},	
\end{align}
where $?=B_{\mathrm{dR},\overline{H}_\infty}^+,B_{\mathrm{dR},\overline{H}_\infty},B_{e,\overline{H}_\infty}$. For the rings:
\begin{align}
\mathcal{O}_{\mathrm{FF}_R,\text{\'et}}\widehat{\otimes}B\widehat{\otimes}B_{e,\overline{H}_\infty}	
\end{align}
the corresponding definition of projective sheaves are defined to be locally associated to the finite projective modules in the ind-Fr\'echet topology, namely the completeness with respect to the ind-Fr\'echet topology is then stable under the corresponding \'etale morphisms concentrated at the part for $\mathrm{FF}_R$. And we assume that the corresponding sheaves admit $\mathrm{Gal}_H$-equivariant action in the semilinear fashion.	
\end{definition}

\begin{remark}
One can basically define the same $B$-$\mathrm{Gal}_H$-sheaves over any adic Banach affinoid algebra or space in the same fashion which we will not repeat again.
\end{remark}

\indent We now prove the following theorems:

\begin{theorem} \mbox{\bf{(After Kedlaya-Liu \cite[Theorem 4.6.1]{8KL2})}} \label{theorem5.18}
Consider the following categories:\\
I. The category of all the projective $\mathrm{FF}_R$-$B$-$\mathrm{Gal}_H$ sheaves over the site:
\begin{align}
(\mathrm{FF}_R,\mathcal{O}_{\mathrm{FF}_R}\widehat{\otimes}B\widehat{\otimes}B_{\mathrm{dR},\overline{H}_\infty}^+,\mathcal{O}_{\mathrm{FF}_R}\widehat{\otimes}B\widehat{\otimes}B_{\mathrm{dR},\overline{H}_\infty},\mathcal{O}_{\mathrm{FF}_R}\widehat{\otimes}B\widehat{\otimes}B_{e,\overline{H}_\infty}).	
\end{align}
II. The category of all the compatible families (with locally glueability and local cocycle condition) of all the projective $\{\mathrm{Spa}(\widetilde{\Pi}_R^{[s,r]},\widetilde{\Pi}_R^{[s,r],+})\}_{\{[s,r]\}}$-$B$-$\mathrm{Gal}_H$ sheaves over the sites:
\begin{align}
&(\mathrm{Spa}(\widetilde{\Pi}_R^{[s,r]},\widetilde{\Pi}_R^{[s,r],+}),\mathcal{O}_{\mathrm{Spa}(\widetilde{\Pi}_R^{[s,r]},\widetilde{\Pi}_R^{[s,r],+})}\widehat{\otimes}B\widehat{\otimes}B_{\mathrm{dR},\overline{H}_\infty}^+,\\
&\mathcal{O}_{\mathrm{Spa}(\widetilde{\Pi}_R^{[s,r]},\widetilde{\Pi}_R^{[s,r],+})}\widehat{\otimes}B\widehat{\otimes}B_{\mathrm{dR},\overline{H}_\infty},\mathcal{O}_{\mathrm{Spa}(\widetilde{\Pi}_R^{[s,r]},\widetilde{\Pi}_R^{[s,r],+})}\widehat{\otimes}B\widehat{\otimes}B_{e,\overline{H}_\infty}),
\end{align}	
for any $[s,r]\subset (0,\infty)$, carrying the corresponding semilinear action of the Frobenius coming from just the Robba ring part.\\
III. The category of all the compatible families (with locally glueability and local cocycle condition) of all the projective $\{\widetilde{\Pi}_R^{[s,r]}\}_{\{[s,r]\}}$-$B$-$\mathrm{Gal}_H$ modules over the rings:
\begin{align}
&(\mathcal{O}_{\mathrm{Spa}(\widetilde{\Pi}_R^{[s,r]},\widetilde{\Pi}_R^{[s,r],+})}\widehat{\otimes}B\widehat{\otimes}B_{\mathrm{dR},\overline{H}_\infty}^+,\\
&\mathcal{O}_{\mathrm{Spa}(\widetilde{\Pi}_R^{[s,r]},\widetilde{\Pi}_R^{[s,r],+})}\widehat{\otimes}B\widehat{\otimes}B_{\mathrm{dR},\overline{H}_\infty},\\
&\mathcal{O}_{\mathrm{Spa}(\widetilde{\Pi}_R^{[s,r]},\widetilde{\Pi}_R^{[s,r],+})}\widehat{\otimes}B\widehat{\otimes}B_{e,\overline{H}_\infty})|_{\mathrm{Spa}(\widetilde{\Pi}_R^{[s,r]},\widetilde{\Pi}_R^{[s,r],+})},
\end{align}	
for any $[s,r]\subset (0,\infty)$, carrying the corresponding semilinear action of the Frobenius coming from just the Robba ring part.\\
IV. The category of all the projective $\widetilde{\Pi}_R^{\infty}$-$B$-$\mathrm{Gal}_H$ modules over the rings:
\begin{align}
&(\mathcal{O}_{\mathrm{Spa}(\widetilde{\Pi}_R^{\infty},\widetilde{\Pi}_R^{\infty,+})}\widehat{\otimes}B\widehat{\otimes}B_{\mathrm{dR},\overline{H}_\infty}^+,\\
&\mathcal{O}_{\mathrm{Spa}(\widetilde{\Pi}_R^{\infty},\widetilde{\Pi}_R^{\infty,+})}\widehat{\otimes}B\widehat{\otimes}B_{\mathrm{dR},\overline{H}_\infty},\\
&\mathcal{O}_{\mathrm{Spa}(\widetilde{\Pi}_R^{\infty},\widetilde{\Pi}_R^{\infty,+})}\widehat{\otimes}B\widehat{\otimes}B_{e,\overline{H}_\infty})|_{\mathrm{Spa}(\widetilde{\Pi}_R^{\infty},\widetilde{\Pi}_R^{\infty,+})},
\end{align}	
carrying the corresponding semilinear action of the Frobenius coming from just the Robba ring part.\\
Then we have the corresponding categories above are equivalent.

\end{theorem}

\begin{proof}
The corresponding families of sheaves with respect to the intervals, the families of modules with respect to the intervals and the sheaves over $\mathrm{FF}_R$ could be related equivalently through the functor taking the corresponding local global section, by applying \cref{theorem4.14} and \cref{theorem3.15}. The corresponding families of sheaves with respect to the intervals, the families of modules with respect to the intervals and the modules over $\widetilde{\Pi}^\infty_R$ could be related equivalently through the functor taking the corresponding local projective limit global section, by applying \cref{theorem4.31} and \cref{theorem3.32} and by considering the corresponding Frobenius pullback to control the rank throughout. 	
\end{proof}

\begin{theorem} \mbox{\bf{(After Kedlaya-Liu \cite[Theorem 4.6.1]{8KL2})}} \label{theorem5.24}
Consider the following categories:\\
I. The category of all the pseudocoherent sheaves over the site:
\begin{align}
(\mathrm{FF}_R,\mathcal{O}_{\mathrm{FF}_R,\text{\'et}}\widehat{\otimes}B).	
\end{align}
II. The category of all the compatible families (with locally glueability and local cocycle condition) of all the pseudocoherent $\{\mathrm{Spa}(\widetilde{\Pi}_R^{[s,r]},\widetilde{\Pi}_R^{[s,r],+})\}_{\{[s,r]\}}$ sheaves over the sites:
\begin{align}
&(\mathrm{Spa}(\widetilde{\Pi}_R^{[s,r]},\widetilde{\Pi}_R^{[s,r],+}),\mathcal{O}_{\mathrm{Spa}(\widetilde{\Pi}_R^{[s,r]},\widetilde{\Pi}_R^{[s,r],+}),\text{\'et}}\widehat{\otimes}B),
\end{align}	
for any $[s,r]\subset (0,\infty)$, carrying the corresponding semilinear action of the Frobenius coming from just the Robba ring part.\\
III. The category of all the compatible families (with locally glueability and local cocycle condition) of all the stably-pseudocoherent $\{\widetilde{\Pi}_R^{[s,r]}\}_{\{[s,r]\}}$ modules over the rings:
\begin{align}
\mathcal{O}_{\mathrm{Spa}(\widetilde{\Pi}_R^{[s,r]},\widetilde{\Pi}_R^{[s,r],+}),\text{\'et}}\widehat{\otimes}B|_{\mathrm{Spa}(\widetilde{\Pi}_R^{[s,r]},\widetilde{\Pi}_R^{[s,r],+})},
\end{align}	
for any $[s,r]\subset (0,\infty)$, carrying the corresponding semilinear action of the Frobenius coming from just the Robba ring part.\\
IV. The category of all the stably-pseudocoherent $\widetilde{\Pi}_R^{\infty}$ modules over the rings:
\begin{align}
\mathcal{O}_{\mathrm{Spa}(\widetilde{\Pi}_R^{\infty},\widetilde{\Pi}_R^{\infty,+}),\text{\'et}}\widehat{\otimes}B|_{\mathrm{Spa}(\widetilde{\Pi}_R^{\infty},\widetilde{\Pi}_R^{\infty,+})},
\end{align}	
carrying the corresponding semilinear action of the Frobenius coming from just the Robba ring part.\\
Then we have the corresponding categories above are equivalent.

\end{theorem}

\begin{proof}
The corresponding families of sheaves with respect to the intervals, the families of modules with respect to the intervals and the sheaves over $\mathrm{FF}_R$ could be related equivalently through the functor taking the corresponding local global section, by applying \cref{theorem2.22}. The corresponding families of sheaves with respect to the intervals, the families of modules with respect to the intervals and the modules over $\widetilde{\Pi}^\infty_R$ could be related equivalently through the functor taking the corresponding local projective limit global section, by applying \cref{theorem2.31} and by considering the corresponding Frobenius pullback to control the rank throughout. 	
\end{proof}

\

This chapter is based on the following paper, where the author of this dissertation is the main author:
\begin{itemize}
\item Tong, Xin. "Period Rings with Big Coefficients and Application II." arXiv preprint arXiv:2101.03748 (2021). 
\end{itemize}

\newpage\chapter{Period Rings with Big Coefficients and Applications III}

\newpage\section{Introduction}

\subsection{Introduction and Main Results}

\indent In our previous work on the corresponding period rings with big coefficients and applications, we essentially studied many very general deformations of the corresponding pseudocoherent sheaves over adic Banach rings in the sense of \cite{9KL1} and \cite{9KL2}. We studied the corresponding glueing in the fashion of \cite{9KL1} and \cite{9KL2} carrying sufficiently large coefficients. Although we have some conditions on both the corresponding adic spaces we are considering and the corresponding coefficients, but the descent results on their own are already general enough to tackle some specific situations in the Hodge-Iwasawa theoretic consideration and noncommutative analytic geometry such as in \cite{9TX3} and \cite{9TX4}.\\

\indent First in our current scope of discussion we will first consider the corresponding extension of our paper \cite{9TX2} to the corresponding pro-\'etale topology setting. In the Tate setting, this will follow the corresponding unrelative situation in \cite{9KL2} while in the analytic situation the corresponding will following the corresponding unrelative situation in \cite{9Ked1} (in the context of \cref{chapter2} and \cref{9chapter3}):

\begin{theorem} \mbox{\bf{(See \cref{9proposition2.10}, \cref{proposition3.17})}} \label{9theorem9.1.1}
Over sheafy Tate adic Banach rings or analytic Huber rings, the corresponding descent of stably pseudocoherent modules (with respect to the underlying spaces) carrying the corresponding noncommutative coefficients holds over \'etale sites and pro-\'etale sites. Here for pro-\'etale sites we will consider only perfectoid rings. 	
\end{theorem}

\indent We then consider some application to the corresponding 'complete' noetherian situation (this will mean that we consider the noetherian deformation of noetherian rings) in the context of \cref{section6}:

\begin{theorem} \mbox{\bf{(See \cref{proposition6.9})}} \label{9theorem9.1.2} 
Over noetherian sheafy Tate adic Banach rings, the corresponding descent of finitely presented modules (with respect to the underlying spaces) carrying the corresponding noncommutative coefficients (such that the products of rings are noetherian as well) holds over \'etale sites.   	
\end{theorem}

\indent We also discussed the corresponding possible application to the descent in certain cases over the corresponding $\infty$-analytic stacks after Bambozzi-Ben-Bassat-Kremnizer \cite{9BBBK} and Bambozzi-Kremnizer \cite{9BK}, where the latter is really more related to the adic geometry we considered after \cite{9KL1} and \cite{9KL2}. Our approach certainly is definitely parallel to \cite{9KL1} and \cite{9KL2}, which is not actually parallel to more $\infty$-method considered by Ben-Bassat-Kremnizer in \cite{9BBK}.

\subsection{Future Work}

\indent We expect more interesting applications. For instance partial noetherian spaces could admit some desired descent carrying general noncommutative coefficients (certainly not need to be noetherian coefficients), such as different types of eigenvarieties or Shimura varieties.

%%\newpage

\newpage\section{Descent over Analytic Adic Banach Rings in the Rational Setting} \label{chapter2}

\subsection{Noncommutative Functional-Analytic Pseudocoherence ov
er Pro-\'Etale Topology}

\indent Now we consider the corresponding discussion of the corresponding glueing deformed pseudocoherent sheaves over \'etale topology which generalize the corresponding discussion in \cite{9Ked1}. 

\begin{setting}
Let $(W,W^+)$ be Tate adic Banach uniform ring which is defined over $\mathbb{Q}_p$. And we assume that the ring $V$ over $\mathbb{Q}_p$ is a Banach ring over the $\mathbb{Q}_p$. Assume that the ring $W$ is sheafy. 	
\end{setting}

\begin{remark}
The corresponding glueing of the corresponding pseudocoherent sheaves over the corresponding pro-\'etale site does not need the corresponding notions on the corresponding stability beyond the corresponding \'etale-stably pseudocoherence.	
\end{remark}

\begin{lemma} \mbox{\bf{(After Kedlaya-Liu \cite[Proposition 3.4.3]{9KL2})}}
Suppose that we are taking $(W,W^+)$ to be Fontaine perfectoid (note that we are considering the corresponding context of Tate adic Banach situation). Then we have that over the corresponding \'etale site, the corresponding group $H^i(\mathrm{Spa}(W,W^+)_\text{p\'et},\mathcal{O}\widehat{\otimes}V)$ vanishes for each $i>0$, while we have that the corresponding group $H^0(\mathrm{Spa}(W,W^+)_\text{p\'et},\mathcal{O}\widehat{\otimes}V)$ is just isomorphic $W\widehat{\otimes}V$.	
\end{lemma}

\begin{proof}
See \cite[Proposition 3.4.3]{9KL2}, where the corresponding \cite[Corollary 3.3.20]{9KL2} applies.	
\end{proof}

\indent Now we proceed to consider the corresponding pro-\'etale topology. First we recall the following result from \cite[Lemma 3.4.4]{9KL2} which does hold in our situation since we did not change much on the corresponding underlying adic spaces.

\begin{lemma}  \mbox{\bf{(Kedlaya-Liu \cite[Lemma 3.4.4]{9KL2})}}
Consider the ring $W$ as above which is furthermore assumed to be Fontaine perfectoid in the sense of \cite[Definition 3.3.1]{9KL2}, and consider any direct system of faithfully finite \'etale morphisms as:
\[
\xymatrix@C+0pc@R+0pc{
W_0 \ar[r] \ar[r] \ar[r] &W_1 \ar[r] \ar[r] \ar[r] &W_2  \ar[r] \ar[r] \ar[r] &...,   
}
\]	
where $A_0$ is just the corresponding base $A$. Then the corresponding completion of the infinite level could be decomposed as:
\begin{align}
W_0\oplus\widehat{\bigoplus}_{k=0}^\infty W_k/W_{k-1}.	
\end{align}
As in \cite[Lemma 3.4.4]{9KL2} we endow the corresponding completion mentioned above with the corresponding seminorm spectral. And for the corresponding quotient we endow with the corresponding quotient norm.
\end{lemma}

%%%%%%%%%%%%%%%%%%%%%%%%%%%%%!!!!!!!!!!!!!!!!

\indent Then we consider the corresponding deformation by the ring $V$ over $\mathbb{Q}_p$:

\begin{lemma}  \mbox{\bf{(After Kedlaya-Liu \cite[Lemma 3.4.4]{9KL2})}}
Consider the ring $A$ as above which is furthermore assumed to be Fontaine perfectoid in the sense of \cite[Definition 3.3.1]{9KL2}, and consider any direct system of faithfully finite \'etale morphisms as:
\[
\xymatrix@C+0pc@R+0pc{
W_{0,V} \ar[r] \ar[r] \ar[r] &W_{1,V} \ar[r] \ar[r] \ar[r] &W_{2,V}  \ar[r] \ar[r] \ar[r] &..., 
}  
\]	
where $W_0$ is just the corresponding base $W$. Then the corresponding completion of the infinite level could be decomposed as:
\begin{align}
W_{0,V}\oplus\widehat{\bigoplus}_{k=0}^\infty W_{k,V}/W_{k-1,V}.	
\end{align}
\end{lemma}

\

\begin{corollary}\mbox{\bf{(After Kedlaya-Liu \cite[Corollary 3.4.5]{9KL2})}} 
Starting with an adic Banach ring which is as in \cite[Corollary 3.4.5]{9KL2} assumed to be Fontaine perfectoid in the sense of \cite[Definition 3.3.1]{9KL2}. And as in the previous two lemmas we consider an admissible infinite direct system:
\[
\xymatrix@C+0pc@R+0pc{
W_{0,V} \ar[r] \ar[r] \ar[r] &W_{1,V} \ar[r] \ar[r] \ar[r] &W_{2,V}  \ar[r] \ar[r] \ar[r] &...,   
}
\]
whose infinite level will be assumed to take the corresponding spectral seminorm as in the \cite[Corollary 3.4.5]{9KL2}. Then carrying the corresponding coefficient $V$, we have in such situation the corresponding 2-pseudoflatness of the corresponding embedding map:
\begin{align}
W_{0,V}\rightarrow W_{\infty,V}.
\end{align}
Here $W_{\infty,V}$ denotes the corresponding completion of the limit $\varinjlim_{k\rightarrow\infty}W_{k,V}$.
\end{corollary}

\begin{proof}
See \cite[Corollary 3.4.5]{9KL2}.	
\end{proof}

\indent Now recall that from \cite{9KL2} since we do not have to modify the corresponding underlying spatial context, so we will also only have to consider the corresponding stability with respect to \'etale topology even although we are considering the corresponding profinite \'etale site in our current section. We first generalize the corresponding Tate acyclicity in \cite{9KL2} to the corresponding $V$-relative situation in our current situation:

%%%%%%%%%%%%%%%%%%%%%%%%%%%%%%!!!!!!!!!!!!!!!!!!!!

\begin{remark}
In the following we only consider the corresponding perfectoid ring $W$.	
\end{remark}

\begin{proposition}\mbox{\bf{(After Kedlaya-Liu \cite[Theorem 3.4.6]{9KL2})}} \label{proposition2.7} 
Now we consider the corresponding pro-\'etale site of $X$ which we denote it by $X_\text{p\'et}$. We assume we have a stable basis $\mathbb{H}$ of the corresponding perfectoid subdomains for this profinite \'etale site where each morphism therein will be \'etale pseudoflat. We consider as in \cite[Theorem 3.4.6]{9KL2} a corresponding module over $W\widehat{\otimes}V$ which is assumed to be \'etale-stably pseudocoherent. Then consider the corresponding presheaf $\widetilde{G}$ attached to this \'etale-stably pseudocoherent module with respect to in our situation the corresponding chosen well-defined basis $\mathbb{H}$, we have that the corresponding sheaf over some element $Y\in \mathbb{H}$ (that is to say over $\mathcal{O}_{Y_\text{p\'et}}\widehat{\otimes}V$) is acyclic and is acyclic with respect to some \v{C}eck covering coming from elements in $\mathbb{H}$. 
\end{proposition}

\begin{proof}
The corresponding proof could be made as in the corresponding nonrelative setting as in \cite[Theorem 3.4.6]{9KL2}. As in \cite[Theorem 3.4.6]{9KL2} the corresponding first step is to check the statement is true with respect to some specific basis consisting of all the perfectoid subdomains which are faithfully finite \'etale over some specific perfectoid subdomain $Y$ (along towers). Then the corresponding statement could be reduced to those with respect to covering which is simple Laurent rational and the corresponding covering in the sense of the previous corollary. Then actually in this first scope of consideration, the statement holds by previous corollary and \cite[Theorem 2.12]{9TX2}. Then in general for some general $Y$ we consider the corresponding covering of $Y$ consisting of elements satisfying the previous situation, then apply the previous argument to the corresponding basis in current situation satisfying the previous situation Then we consider the corresponding sheafification of the following chosen short exact sequence:
\[
\xymatrix@C+0pc@R+0pc{
0 \ar[r] \ar[r] \ar[r] &G'' \ar[r] \ar[r] \ar[r] &G'  \ar[r] \ar[r] \ar[r] &G \ar[r] \ar[r] \ar[r] &0 ,   
}
\]
with $G'$ finite projective. Then the idea is to consider the corresponding diagram factoring through the corresponding pro-\'etale sheafification of this short exact sequence down to the corresponding \'etale site. Then the five lemma will finish the corresponding proof as one considers the following commutative diagram:
\[
\xymatrix@C+0pc@R+3pc{
&G''_{V_Y} \ar[r] \ar[r] \ar[r] \ar[d] \ar[d] \ar[d] &G'_{V_Y}  \ar[r] \ar[r] \ar[r] \ar[d] \ar[d] \ar[d] &G_{V_Y} \ar[r] \ar[r] \ar[r] \ar[d] \ar[d] \ar[d] &0 ,  \\
0 \ar[r] \ar[r] \ar[r] &\Gamma(Y_\text{\'et},\widetilde{G''}) \ar[r] \ar[r] \ar[r] &\Gamma(Y_\text{\'et},\widetilde{G'})  \ar[r] \ar[r] \ar[r] &\Gamma(Y_\text{\'et},\widetilde{G}) \ar[r] \ar[r] \ar[r] &0.  \\ 
}
\]
Here the corresponding ring $V_Y$ is the adic Banach ring such that $Y=\mathrm{Spa}(V_Y,V_Y^+)$.
\end{proof}

\begin{definition}\mbox{\bf{(After Kedlaya-Liu \cite[Definition 3.4.7]{9KL2})}} 
Consider the pro-\'etale site of $X$ attached to the adic Banach ring $(W,W^+)$. We will define a sheaf of module $G$ over $\widehat{\mathcal{O}}_\text{p\'et}\widehat{\otimes}V$ to be $V$-pseudocoherent if locally we can define this as a sheaf attached to a $V$-\'etale-stably pseudocoherent module. As in \cite[Definition 3.4.7]{9KL2}, we do not have to consider the corresponding notion of $V$-profinite-\'etale-stably pseudocoherent module.	
\end{definition}

\begin{proposition}\mbox{\bf{(After Kedlaya-Liu \cite[Theorem 3.4.8]{9KL2})}} \label{9proposition2.10}
Taking the corresponding global section will realize the corresponding equivalence between the following two categories. The first one is the corresponding one of all the $V$-pseudocoherent sheaves over $\widehat{\mathcal{O}}_\text{p\'et}\widehat{\otimes}V$, while the second one is the corresponding one of all the $V$-\'etale-stably pseudocoherent modules over $W\widehat{\otimes}V$.
	
\end{proposition}

\begin{proof}
As in \cite[Theorem 3.4.8]{9KL2}, the corresponding \cite[Proposition 9.2.6]{9KL1} applies in the way that the corresponding conditions of \cite[Proposition 9.2.6]{9KL1} hold in our current situation.	
\end{proof}

\indent Obviously we have the following analog of \cite[Corollary 3.4.9]{9KL2}:

%%%%%%%%%%%%%%%%%%%%%%%%%%%%%%%%!!!!!!!!!!

\begin{corollary} \mbox{\bf{(After Kedlaya-Liu \cite[Corollary 3.4.9]{9KL2})}}
The following two categories are equivalent. The first is the corresponding category of all $V$-pseudocoherent sheaves over $\widehat{\mathcal{O}}_\text{p\'et}\widehat{\otimes}V$. The second is the corresponding category of all $V$-pseudocoherent sheaves over ${\mathcal{O}}_\text{\'et}\widehat{\otimes}V$. The corresponding functor is the corresponding pullback along the corresponding morphism of sites $X_{\text{p\'et}}\rightarrow X_{\text{\'et}}$. 	
\end{corollary}

%\
%
%\subsection{The Discussion for Exotic Topology}
%
%
%\indent Now we consider the corresponding exotic topology, namely the $v$-topology. Due to the corresponding nonflatness in this consideration, certainly we cannot expect the discussion around the corresponding pseudocoherent sheaves at all in our current situation.

%%\newpage

\newpage\section{Descent over Analytic Huber Rings in the Integral Setting} \label{9chapter3}

\subsection{Noncommutative Topological Pseudocoherence over \'Etale Topology}

\indent Now we consider the corresponding discussion of the corresponding glueing deformed pseudocoherent sheaves over \'etale topology which generalize the corresponding discussion in \cite{9Ked1}.  And we will consider the corresponding situation as assumed in the following setting:

\begin{setting}
Let $(W,W^+)$ be analytic Huber uniform pair which is defined over $\mathbb{Z}_p$. And we assume that the ring $V$ over $\mathbb{Z}_p$ is a topological ring (complete) and splitting over the $\mathbb{Z}_p$. We now fix a corresponding stable basis $\mathbb{H}$ for the corresponding \'etale site of the adic space $\mathrm{Spa}(W,W^+)$, locally consisting of the corresponding compositions of rational localizations and the corresponding finite \'etale morphisms. And as in \cite[Hypothesis 1.10.3]{9Ked1} we need to assume that the corresponding basis is made up of the corresponding adic spectrum of sheafy rings. Assume the corresponding sheafiness of the Huber ring $W$.
\end{setting}

\begin{definition} \mbox{\bf{(After Kedlaya-Liu \cite[Definition 2.5.9]{9KL2})}}
We define a $V$-\'etale stably pseudocoherent module over the corresponding Huber ring $W$ with respect to the corresponding basis $\mathbb{H}$ chosen above for the corresponding \'etale site of the analytic adic space $X$. We define a module over $W\widehat{\otimes}V$ is a $V$-\'etale stably pseudocoherent module if it is algebraically pseudocoherent (namely formed by the corresponding possibly infinite length resolution of finitely generated and projective modules) and at least complete with respect to the corresponding natural topology and also required to be also complete with respect to the natural topology along some base change with respect to any morphism in $\mathbb{H}$.
\end{definition}

\begin{definition} \mbox{\bf{(After Kedlaya-Liu \cite[Definition 2.5.9]{9KL2})}}
Along the previous definition we have the corresponding notion of $V$-\'etale-pseudoflat left (or right respectively) modules with respect to the corresponding chosen basis $\mathbb{H}$. A such module is defined to be a topological module $G$ over $W\widehat{\otimes}V$ complete with respect to the natural topology and for any right (or left respectively) $V$-\'etale stably pseudocoherent module, they will jointly give the $\mathrm{Tor}_1$ vanishing.
\end{definition}

\begin{lemma}\mbox{\bf{(After Kedlaya-Liu \cite[Lemma 2.5.10]{9KL2})}}
One could find another basis $\mathbb{H}'$ which is in our situation contained in the corresponding original basis $\mathbb{H}$ such that any morphism in $\mathbb{H}$ could be $V$-\'etale-pseudoflat with respect $\mathbb{H}'$ or just $\mathbb{H}$ itself.	
\end{lemma}

\begin{proof}
The derivation of such new basis is actually along the same way as in \cite[Lemma 2.5.10]{9KL2} since the corresponding compositions of rational localizations and the finite \'etales are actually satisfying the corresponding conditions in the statement by \cite[Theorem 2.12]{9TX2} over the analytic topology. In general one just show any general morphism will also decompose in the same way, which one can proves certainly as in \cite[Lemma 2.5.10]{9KL2}, where one instead consider in the current context the corresponding basis spreading result in \cite[Lemma 1.10.4]{9Ked1}. 
\end{proof}

%%%%%%%%%%%%%%%%%%%%%%%%%%%%%!!!!!!!!!!!!!!!!!!!!!

\begin{proposition} \mbox{\bf{(After Kedlaya-Liu \cite[Theorem 2.5.11]{9KL2})}}
Consider the site $X_\text{\'et}$ and consider the basis $\mathbb{H}$. Take any $V$-\'etale stably pseudocoherent module $M$ over $W\widehat{\otimes}V$. Consider the corresponding presheaf by taking the inverse limit throughout all the corresponding base change along morphisms in $\mathbb{H}$. Then we have the corresponding acyclicity of the presheaf over any element in $\mathbb{H}$ and any covering of this element by elements in $\mathbb{H}$.	
\end{proposition}

\begin{proof}
As in \cite[Theorem 2.5.11]{9KL2} apply the corresponding \cite[Proposition 8.2.21]{9KL1}.	
\end{proof}

\begin{proposition}\mbox{\bf{(After Kedlaya-Liu \cite[Lemma 2.5.13]{9KL2})}}
The corresponding glueing of $V$-\'etale-stably pseudocoherent modules holds in this current situation over $W$ along binary morphisms (namely along binary rational decomposition). 	
\end{proposition}

\begin{proof}
We adapt the argument of \cite[Lemma 2.5.13]{9KL2} to our current situation. Take the corresponding map to be $W\rightarrow W_1\bigoplus W_2$, and we denote the two spaces associated $W_1$ and $W_2$ by $Y_1$ and $Y_2$. Establish a corresponding descent datum of $V$-\'etale-stably pseudocoherent modules along this decomposition of $\mathrm{Spa}(W,W^+)$. One can realize this by acyclicity a sheaf $G$ over the whole space $\mathrm{Spa}(W,W^+)$ Then we have for each $k=1,2$ that the corresponding identification of $H^0(Y_{k,\text{\'et}},G)$ with $H^0(Y_k,G)$ and the vanishing of the $i$-th cohomology for higher $i$. By the corresponding results we have in the corresponding rational localization situation we vanishing of the corresponding $H^i(\mathrm{Spa}(W,W^+),G),i\geq 1$. Now consider a finite free covering $C$ of $G$ and take the kernel $K$ to form the corresponding exact sequence:
\[
\xymatrix@C+0pc@R+0pc{
0 \ar[r] \ar[r] \ar[r] &K \ar[r] \ar[r] \ar[r] &C  \ar[r] \ar[r] \ar[r] &G \ar[r] \ar[r] \ar[r] &0 ,   
}
\]
which will give rise to the following short exact sequence:
\[
\xymatrix@C+0pc@R+0pc{
0 \ar[r] \ar[r] \ar[r] &K(\mathrm{Spa}(W,W^+)) \ar[r] \ar[r] \ar[r] &C(\mathrm{Spa}(W,W^+))  \ar[r] \ar[r] \ar[r] &G(\mathrm{Spa}(W,W^+)) \ar[r] \ar[r] \ar[r] &0.   
}
\]
To finish we have also to show the corresponding resulting global section is actually $V$-\'etale-stably pseudocoherent. For this, we refer the readers to \cite[Lemma 2.22]{9TX2}.

\end{proof}

\begin{proposition} \mbox{\bf{(After Kedlaya-Liu \cite[Theorem 2.5.14]{9KL2})}}
Taking the corresponding global section will realize the equivalence between the following two categories. The first is the one of all the corresponding $V$-pseudocoherent sheaves over $\mathcal{O}_{\text{\'et}}\widehat{\otimes}V$. The second is the one of all the corresponding $V$-\'etale-stably pseudocoherent modules.	
\end{proposition}

\begin{proof}
This is by considering and applying the corresponding \cite[Lemma 1.10.4]{9Ked1} in our current context. The refinement comes from the acyclicity, the transitivity is straightforward, the corresponding binary rational localization situation is the previous proposition, while the corresponding last condition will basically comes from the corresponding f.p. descent \cite[Chapitre VIII]{9SGAI}.	
\end{proof}

\
\subsection{Noncommutative Topological Pseudocoherence over Pro-\'Etale Topology}

\begin{setting}
Let $(W,W^+)$ be analytic Huber ring which is defined over $\mathbb{Z}_p$. And we assume that the ring $V$ over $\mathbb{Z}_p$ is a complete topological ring over the $\mathbb{Z}_p$ which is completely exact. Assume that the ring $W$ is sheafy. We assume $p$ is a topologically nilpotent element.
\end{setting}

\begin{remark}
The corresponding glueing of the corresponding pseudocoherent sheaves over the corresponding pro-\'etale site does not need the corresponding notions on the corresponding stability beyond the corresponding \'etale-stably pseudocoherence.	
\end{remark}

\begin{lemma} \mbox{\bf{(After Kedlaya-Liu \cite[Proposition 3.4.3]{9KL2})}}
Suppose that we are taking $(W,W^+)$ to be perfectoid now in the sense of \cite[Definition 2.1.1]{9Ked1} (note that we are considering the corresponding context of analytic Huber pair situation). Then we have that over the corresponding \'etale site, the corresponding group $H^i(\mathrm{Spa}(W,W^+)_\text{p\'et},\mathcal{O}\widehat{\otimes}V)$ vanishes for each $i>0$, while we have that the corresponding group \\$H^0(\mathrm{Spa}(W,W^+)_\text{p\'et},\mathcal{O}\widehat{\otimes}V)$ is just isomorphic $W\widehat{\otimes}V$.	
\end{lemma}

\begin{proof}
See \cite[Proposition 3.4.3]{9KL2}.	
\end{proof}

\indent Now we proceed to consider the corresponding pro-\'etale topology. First we recall the following result from \cite[Lemma 3.4.4]{9KL2} which does hold in our situation since we did not change much on the corresponding underlying adic spaces.

\begin{lemma}  \mbox{\bf{(Kedlaya-Liu \cite[Lemma 3.4.4]{9KL2})}}
Consider the ring $W$ as above which is furthermore assumed to be perfectoid in the sense of instead \cite[Definition 2.1.1]{9Ked1}, and consider any direct system of faithfully finite \'etale morphisms as:
\[
\xymatrix@C+0pc@R+0pc{
W_0 \ar[r] \ar[r] \ar[r] &W_1 \ar[r] \ar[r] \ar[r] &W_2  \ar[r] \ar[r] \ar[r] &...,   
}
\]	
where $A_0$ is just the corresponding base $A$. Then the corresponding completion of the infinite level could be decomposed as:
\begin{align}
W_0\oplus\widehat{\bigoplus}_{k=0}^\infty W_k/W_{k-1}.	
\end{align}
As in \cite[Lemma 3.4.4]{9KL2} we endow the corresponding completion mentioned above with the corresponding seminorm spectral. And for the corresponding quotient we endow with the corresponding quotient norm.
\end{lemma}

\indent Then we consider the corresponding deformation by the ring $V$ over $\mathbb{Z}_p$:

\begin{lemma}  \mbox{\bf{(After Kedlaya-Liu \cite[Lemma 3.4.4]{9KL2})}}
Consider the ring $A$ as above which is furthermore assumed to be perfectoid in the sense of \cite[Definition 2.1.1]{9Ked1}, and consider any direct system of faithfully finite \'etale morphisms as:
\[
\xymatrix@C+0pc@R+0pc{
W_{0,V} \ar[r] \ar[r] \ar[r] &W_{1,V} \ar[r] \ar[r] \ar[r] &W_{2,V}  \ar[r] \ar[r] \ar[r] &..., 
}  
\]	
where $W_0$ is just the corresponding base $W$. Then the corresponding completion of the infinite level could be decomposed as:
\begin{align}
W_{0,V}\oplus\widehat{\bigoplus}_{k=0}^\infty W_{k,V}/W_{k-1,V}.	
\end{align}
\end{lemma}

\

\begin{corollary}\mbox{\bf{(After Kedlaya-Liu \cite[Corollary 3.4.5]{9KL2})}} 
Starting with an analytic Huber ring which is as in \cite[Corollary 3.4.5]{9KL2} assumed to be perfectoid in the sense of \cite[Definition 2.1.1]{9Ked1}. And as in the previous two lemmas we consider an admissible infinite direct system:
\[
\xymatrix@C+0pc@R+0pc{
W_{0,V} \ar[r] \ar[r] \ar[r] &W_{1,V} \ar[r] \ar[r] \ar[r] &W_{2,V}  \ar[r] \ar[r] \ar[r] &...,   
}
\]
whose infinite level will be assumed to take the corresponding spectral seminorm as in the \cite[Corollary 3.4.5]{9KL2}. Then carrying the corresponding coefficient $V$, we have in such situation the corresponding 2-pseudoflatness of the corresponding embedding map:
\begin{align}
W_{0,V}\rightarrow W_{\infty,V}.
\end{align}
Here $W_{\infty,V}$ denotes the corresponding completion of the limit $\varinjlim_{k\rightarrow\infty}W_{k,V}$.
\end{corollary}

\begin{proof}
See \cite[Corollary 3.4.5]{9KL2}.	
\end{proof}

\indent Now recall that from \cite{9KL2} since we do not have to modify the corresponding underlying spatial context, so we will also only have to consider the corresponding stability with respect to \'etale topology even although we are considering the corresponding profinite \'etale site in our current section. We first generalize the corresponding Tate acyclicity in \cite{9KL2} to the corresponding $V$-relative situation in our current situation:

\begin{remark}
In the following we only consider the corresponding perfectoid ring $W$.	
\end{remark}

\begin{proposition}\mbox{\bf{(After Kedlaya-Liu \cite[Theorem 3.4.6]{9KL2})}} 
Now we consider the corresponding pro-\'etale site of $X$ which we denote it by $X_\text{p\'et}$. Then we assume we have a stable basis $\mathbb{H}$ of the corresponding perfectoid subdomains for this profinite \'etale site where each morphism therein will be \'etale pseudoflat. Then we consider as in \cite[Theorem 3.4.6]{9KL2} a corresponding module over $W\widehat{\otimes}V$ which is assumed to be \'etale-stably pseudocoherent. Then consider the corresponding presheaf $\widetilde{G}$ attached to this \'etale-stably pseudocoherent with respect to in our situation the corresponding chosen well-defined basis $\mathbb{H}$, we have that the corresponding sheaf over some element $Y\in \mathbb{H}$ (that is to say over $\mathcal{O}_{Y_\text{p\'et}}\widehat{\otimes}V$) is acyclic and is acyclic with respect to some \v{C}eck covering coming from elements in $\mathbb{H}$. 
\end{proposition}

\begin{proof}
See \cref{proposition2.7}.
\end{proof}

\begin{definition}\mbox{\bf{(After Kedlaya-Liu \cite[Definition 3.4.7]{9KL2})}} 
Consider the pro-\'etale site of $X$ attached to the adic Banach ring $(W,W^+)$. We will define a sheaf of module $G$ over $\widehat{\mathcal{O}}_\text{p\'et}\widehat{\otimes}V$ to be $V$-pseudocoherent if locally we can define this as a sheaf attached to a $V$-\'etale-stably pseudocoherent module. As in \cite[Definition 3.4.7]{9KL2}, we do not have to consider the corresponding notion of $V$-profinite-\'etale-stably pseudocoherent module.	
\end{definition}

\begin{proposition}\mbox{\bf{(After Kedlaya-Liu \cite[Theorem 3.4.8]{9KL2})}} \label{proposition3.17}
Taking the corresponding global section will realize the corresponding equivalence between the following two categories. The first one is the corresponding one of all the $V$-pseudocoherent sheaves over $\widehat{\mathcal{O}}_\text{p\'et}\widehat{\otimes}V$, while the second one is the corresponding one of all the $V$-\'etale-stably pseudocoherent module over $W\widehat{\otimes}V$.
	
\end{proposition}

\begin{proof}
As in \cite[Theorem 3.4.8]{9KL2}, the corresponding analog of \cite[Proposition 9.2.6]{9KL1} applies in the way that the corresponding conditions of the analog of \cite[Proposition 9.2.6]{9KL1} hold in our current situation. Also see \cite[Theorem 2.9.9, Remark 2.9.10 and Lemma 1.10.4]{9Ked1}.
\end{proof}

\indent Obviously we have the following analog of \cite[Corollary 3.4.9]{9KL2}:

\begin{corollary} \mbox{\bf{(After Kedlaya-Liu \cite[Corollary 3.4.9]{9KL2})}}
The following two categories are equivalent. The first is the corresponding category of all $V$-pseudocoherent sheaves over $\widehat{\mathcal{O}}_\text{p\'et}\widehat{\otimes}V$. The second is the corresponding category of all $V$-pseudocoherent sheaves over ${\mathcal{O}}_\text{\'et}\widehat{\otimes}V$. The corresponding functor is the corresponding pullback along the corresponding morphism of sites $X_{\text{p\'et}}\rightarrow X_{\text{\'et}}$. 	
\end{corollary}

%%\newpage

%\newpage\section{Descent over adic Banach Rings in the Integral Setting}
%
%\subsection{Noncommutative Pseudocoherence over \'Etale Topology}
%
%\indent Now we consider the corresponding discussion of the corresponding glueing deformed pseudocoherent sheaves over \'etale topology which generalize the corresponding discussion in \cite{9Ked1}. This translates the corresponding discussion in the previous section to the setting of adic Banach rings, therefore we will skip the corresponding parallel argument in the proof.
%
%\begin{setting}
%Let $(A,A^+)$ is analytic uniform adic Banach ring which is defined over $\mathbb{Z}_p$. And we assume that the ring $Z$ over $\mathbb{Z}_p$ is a Banach ring and splitting over the $\mathbb{Z}_p$. 	
%\end{setting}
%
%
%
%
%
%\begin{definition}
%We defined a $Z$-\'etale stably pseudocoherent 	
%\end{definition}
%
%
%
%
%\subsection{Noncommutative Pseudocoherence over Pro-\'Etale Topology}
%
%\subsection{The Discussion for Exotic Topology}
%
%
%

%%%%%%%%%%%%%%%%%%%%%%%%%%%%%%%!!!!!!!!!!!!!

\newpage\section{Descent in General Setting for Analytic Huber Pai
rs}

\indent We now consider the situation where there is no condition on the base (except that we will assume certainly the corresponding sheafiness). We will discuss along \cite{9Ked1} the corresponding descent in analytic and \'etale topology, as well as the corresponding quasi-Stein spaces. The current situation will not assume we are working over $\mathbb{Z}_p$.

\begin{remark}
The corresponding discussion in this section and in the following section is actually an exercise proposed by Kedlaya included in \cite{9Ked1}, where the corresponding \'etale situation essentially is proposed by Kedlaya in \cite[Discussion after Lemma 1.10.4]{9Ked1} while we give the corresponding exposition in both \'etale and pro-\'etale situations closely after \cite{9KL1} and \cite{9KL2}.	
\end{remark}

\subsection{\'Etale Topology}

\begin{setting}
Let $(W,W^+)$ be analytic Huber uniform pair. We now fix a corresponding stable basis $\mathbb{H}$ for the corresponding \'etale site of the adic space $\mathrm{Spa}(W,W^+)$, locally consisting of the corresponding compositions of rational localizations and the corresponding finite \'etale morphisms. And as in \cite[Hypothesis 1.10.3]{9Ked1} we need to assume that the corresponding basis is made up of the corresponding adic spectrum of sheafy rings. Assume the corresponding sheafiness of the Huber ring $W$.
\end{setting}

\begin{remark}
One should treat our discussion in this section as essentially some detailization of some exercise predicted in \cite{9Ked1}.	
\end{remark}

\begin{definition} \mbox{\bf{(After Kedlaya-Liu \cite[Definition 2.5.9]{9KL2})}}
We define the \'etale stably pseudocoherent module over the corresponding Huber ring $W$ with respect to the corresponding basis $\mathbb{H}$ chosen above for the corresponding \'etale site of the analytic adic space $X$. We define a module over $W$ to be an \'etale stably pseudocoherent module if it is algebraically pseudocoherent (namely formed by the corresponding possibly infinite length resolution of finitely generated and projective modules) and at least complete with respect to the corresponding natural topology and also required to be also complete with respect to the natural topology along some base change with respect to any morphism in $\mathbb{H}$.
\end{definition}

\begin{definition} \mbox{\bf{(After Kedlaya-Liu \cite[Definition 2.5.9]{9KL2})}}
Along the previous definition we have the corresponding notion of \'etale-pseudoflat left (or right respectively) modules with respect to the corresponding chosen basis $\mathbb{H}$. A such module is defined to be a topological module $G$ over $W$ complete with respect to the natural topology and for any right (or left respectively) \'etale stably pseudocoherent module, they will jointly give the $\mathrm{Tor}_1$ vanishing.
\end{definition}

\begin{remark}
Although in this definition we considered the corresponding left or right modules, but essentially speaking this is no different from the usual situation since we are in the corresponding essential site theoretic situation.	
\end{remark}

\begin{lemma}\mbox{\bf{(After Kedlaya-Liu \cite[Lemma 2.5.10]{9KL2})}}\label{lemma4.4}
One could find another basis $\mathbb{H}'$ which is in our situation contained in the corresponding original basis $\mathbb{H}$ such that any morphism in $\mathbb{H}$ could be \'etale-pseudoflat with respect $\mathbb{H}'$ or just $\mathbb{H}$ itself.	
\end{lemma}

\begin{proof}
The derivation of such new basis is actually along the same way as in \cite[Lemma 2.5.10]{9KL2} since the corresponding compositions of rational localizations and the finite \'etales are actually satisfying the corresponding conditions in the statement by \cite[Theorem 2.12]{9TX2} over the analytic topology. In general one just show any general morphism will also decompose in the same way, which one can proves certainly as in \cite[Lemma 2.5.10]{9KL2}, where one instead consider in the current context the corresponding basis spreading result in \cite[Lemma 1.10.4]{9Ked1}. 
\end{proof}

\begin{proposition} \mbox{\bf{(After Kedlaya-Liu \cite[Theorem 2.5.11]{9KL2})}}
Consider the site $X_\text{\'et}$ and consider the basis $\mathbb{H}$. Take any \'etale stably pseudocoherent module $M$ over $W$. Consider the corresponding presheaf by taking the inverse limit throughout all the corresponding base change along morphisms in $\mathbb{H}$. Then we have the corresponding acyclicity of the presheaf over any element in $\mathbb{H}$ and any covering of this element by elements in $\mathbb{H}$.	
\end{proposition}

\begin{proof}
As in \cite[Theorem 2.5.11]{9KL2} apply the corresponding \cite[Proposition 8.2.21]{9KL1}.	
\end{proof}

\begin{proposition}\mbox{\bf{(After Kedlaya-Liu \cite[Lemma 2.5.13]{9KL2})}}
The corresponding glueing of \'etale-stably pseudocoherent modules holds in this current situation over $W$ along binary morphisms (namely along binary rational decomposition). 	
\end{proposition}

\begin{proof}
We adapt the argument of \cite[Lemma 2.5.13]{9KL2} to our current situation. Take the corresponding map to be $W\rightarrow W_1\bigoplus W_2$, and we denote the two spaces associated $W_1$ and $W_2$ by $Y_1$ and $Y_2$. Establish a corresponding descent datum of $V$-\'etale-stably pseudocoherent modules along this decomposition of $\mathrm{Spa}(W,W^+)$. One can realize this by acyclicity a sheaf $G$ over the whole space $\mathrm{Spa}(W,W^+)$ Then we have for each $k=1,2$ that the corresponding identification of $H^0(Y_{k,\text{\'et}},G)$ with $H^0(Y_k,G)$ and the vanishing of the $i$-th cohomology for higher $i$. By the corresponding results we have in the corresponding rational localization situation we have vanishing of the corresponding $H^i(\mathrm{Spa}(W,W^+),G),i\geq 1$. Now consider a finite free covering $C$ of $G$ and take the kernel $K$ to form the corresponding exact sequence:
\[
\xymatrix@C+0pc@R+0pc{
0 \ar[r] \ar[r] \ar[r] &K \ar[r] \ar[r] \ar[r] &C  \ar[r] \ar[r] \ar[r] &G \ar[r] \ar[r] \ar[r] &0 ,   
}
\]
which will give rise to the following short exact sequence:
\[
\xymatrix@C+0pc@R+0pc{
0 \ar[r] \ar[r] \ar[r] &K(\mathrm{Spa}(W,W^+)) \ar[r] \ar[r] \ar[r] &C(\mathrm{Spa}(W,W^+))  \ar[r] \ar[r] \ar[r] &G(\mathrm{Spa}(W,W^+)) \ar[r] \ar[r] \ar[r] &0.   
}
\]
To finish we have also to show the corresponding resulting global section is actually \'etale-stably pseudocoherent. For this, we refer the readers to \cite[Lemma 2.22]{9TX2}.

\end{proof}

\begin{proposition}\mbox{}\\
\mbox{\bf{(Kedlaya \cite[Below Lemma 1.10.4]{9Ked1})}}\\
\mbox{\bf{(After Kedlaya-Liu \cite[Theorem 2.5.14]{9KL2})}}
Taking the corresponding global section will realize the equivalence between the following two categories. The first is the one of all the corresponding pseudocoherent sheaves over $\mathcal{O}_{\text{\'et}}$. The second is the one of all the corresponding \'etale-stably pseudocoherent modules.	
\end{proposition}

\begin{proof}
This is by considering and applying the corresponding \cite[Lemma 1.10.4]{9Ked1} in our current context. The refinement comes from the acyclicity, the transitivity is straightforward, the corresponding binary rational localization situation is the previous proposition, while the corresponding last condition will basically come from the corresponding f.p. descent \cite[Chapitre VIII]{9SGAI}.	
\end{proof}

\
\subsection{Pro-\'etale Topology}

\begin{setting}
Let $(W,W^+)$ be an analytic Huber ring. Assume that the ring $W$ is sheafy. We assume $p$ is a topologically nilpotent element.	
\end{setting}

\begin{remark}
Although it seems that the current discussion is basically outside the corresponding deformed setting, but the foundation here may be generalized to the deformed setting by taking the corresponding completed tensor product over $\mathbb{F}_1$ as in \cite{9BBBK}.	
\end{remark}

\begin{remark}
The corresponding glueing of the corresponding pseudocoherent sheaves over the corresponding pro-\'etale site does not need the corresponding notions on the corresponding stability beyond the corresponding \'etale-stably pseudocoherence.	
\end{remark}

\begin{lemma} \mbox{\bf{(After Kedlaya-Liu \cite[Proposition 3.4.3]{9KL2})}}
Suppose that we are taking $(W,W^+)$ to be perfectoid now in the sense of \cite[Definition 2.1.1]{9Ked1} (note that we are considering the corresponding context of analytic Huber pair situation). Then we have that over the corresponding \'etale site, the corresponding group $H^i(\mathrm{Spa}(W,W^+)_\text{p\'et},\mathcal{O})$ vanishes for each $i>0$, while we have that the corresponding group $H^0(\mathrm{Spa}(W,W^+)_\text{p\'et},\mathcal{O})$ is just isomorphic to $W$.	
\end{lemma}

\begin{proof}
See \cite[Proposition 3.4.3]{9KL2}.	
\end{proof}

\indent Now we proceed to consider the corresponding pro-\'etale topology. First we recall the following result from \cite[Lemma 3.4.4]{9KL2} which does hold in our situation since we did not change much on the corresponding underlying adic spaces.\\

\begin{lemma}  \mbox{\bf{(Kedlaya-Liu \cite[Lemma 3.4.4]{9KL2})}}
Consider the ring $W$ as above which is furthermore assumed to be perfectoid in the sense of instead \cite[Definition 2.1.1]{9Ked1}, and consider any direct system of faithfully finite \'etale morphisms as:
\[
\xymatrix@C+0pc@R+0pc{
W_0 \ar[r] \ar[r] \ar[r] &W_1 \ar[r] \ar[r] \ar[r] &W_2  \ar[r] \ar[r] \ar[r] &...,   
}
\]	
where $W_0$ is just the corresponding base $W$. Then the corresponding completion of the infinite level could be decomposed as:
\begin{align}
W_0\oplus\widehat{\bigoplus}_{k=0}^\infty W_k/W_{k-1}.	
\end{align}
As in \cite[Lemma 3.4.4]{9KL2} we endow the corresponding completion mentioned above with the corresponding seminorm spectral. And for the corresponding quotient we endow with the corresponding quotient norm.
\end{lemma}

\

\begin{corollary}\mbox{\bf{(After Kedlaya-Liu \cite[Corollary 3.4.5]{9KL2})}} 
Starting with an analytic Huber ring which is as in \cite[Corollary 3.4.5]{9KL2} assumed to be perfectoid in the sense of \cite[Definition 2.1.1]{9Ked1}. And as in the previous two lemmas we consider an admissible infinite direct system:
\[
\xymatrix@C+0pc@R+0pc{
W_{0} \ar[r] \ar[r] \ar[r] &W_{1} \ar[r] \ar[r] \ar[r] &W_{2}  \ar[r] \ar[r] \ar[r] &...,   
}
\]
whose infinite level will be assumed to take the corresponding spectral seminorm as in the \cite[Corollary 3.4.5]{9KL2}. Then we have in such situation the corresponding 2-pseudoflatness of the corresponding embedding map:
\begin{align}
W_{0}\rightarrow W_{\infty}.
\end{align}
Here $W_{\infty}$ denotes the corresponding completion of the limit $\varinjlim_{k\rightarrow\infty}W_{k}$.
\end{corollary}

\begin{proof}
See \cite[Corollary 3.4.5]{9KL2}.	
\end{proof}

\indent Now recall that from \cite{9KL2} since we do not have to modify the corresponding underlying spatial context, so we will also only have to consider the corresponding stability with respect to \'etale topology even although we are considering the corresponding profinite \'etale site in our current section. We first generalize the corresponding Tate acyclicity in \cite{9KL2} to the corresponding our current situation.\\

\begin{remark}
In the following we only consider the base spaces which are perfectoid.	
\end{remark}

\begin{proposition}\mbox{\bf{(After Kedlaya-Liu \cite[Theorem 3.4.6]{9KL2})}} 
Now we consider the corresponding pro-\'etale site of $X$ which we denote it by $X_\text{p\'et}$. Then we assume we have a stable basis $\mathbb{H}$ of the corresponding perfectoid subdomains for this profinite \'etale site where each morphism therein will be \'etale pseudoflat. Then we consider as in \cite[Theorem 3.4.6]{9KL2} a corresponding module over $W$ which is assumed to be \'etale-stably pseudocoherent. Then consider the corresponding presheaf $\widetilde{G}$ attached to this \'etale-stably pseudocoherent with respect to in our situation the corresponding chosen well-defined basis $\mathbb{H}$, we have that the corresponding sheaf over some element $Y\in \mathbb{H}$ (that is to say over $\mathcal{O}_{Y_\text{p\'et}}$) is acyclic and is acyclic with respect to some \v{C}eck covering coming from elements in $\mathbb{H}$. 
\end{proposition}

\begin{proof}
See \cref{proposition2.7}.
\end{proof}

\begin{definition}\mbox{\bf{(After Kedlaya-Liu \cite[Definition 3.4.7]{9KL2})}} 
Consider the pro-\'etale site of $X$ attached to the analytic Huber pair $(W,W^+)$. We will define a sheaf of module $G$ over $\widehat{\mathcal{O}}_\text{p\'et}$ to be pseudocoherent if locally we can define this as a sheaf attached to an \'etale-stably pseudocoherent module. As in \cite[Definition 3.4.7]{9KL2}, we do not have to consider the corresponding notion of profinite-\'etale-stably pseudocoherent modules.	
\end{definition}

\begin{proposition}\mbox{}\\
\mbox{\bf{(Kedlaya \cite[Section 3.8]{9Ked1}, after Kedlaya-Liu \cite[Theorem 3.4.8]{9KL2})}} \label{proposition4.20}
Taking the corresponding global section will realize the corresponding equivalence between the following two categories. The first one is the corresponding one of all the pseudocoherent sheaves over $\widehat{\mathcal{O}}_\text{p\'et}$, while the second one is the corresponding one of all the \'etale-stably pseudocoherent modules over $W$.
	
\end{proposition}

\begin{proof}
As in \cite[Theorem 3.4.8]{9KL2}, the corresponding analog of \cite[Proposition 9.2.6]{9KL1} applies in the way that the corresponding conditions of the analog of \cite[Proposition 9.2.6]{9KL1} hold in our current situation. Also see \cite[Theorem 2.9.9, Remark 2.9.10 and Lemma 1.10.4]{9Ked1}.
\end{proof}

\indent Obviously we have the following analog of \cite[Corollary 3.4.9]{9KL2}:

\begin{corollary} \mbox{\bf{(After Kedlaya-Liu \cite[Corollary 3.4.9]{9KL2})}}
The following two categories are equivalent. The first is the corresponding category of all pseudocoherent sheaves over $\widehat{\mathcal{O}}_\text{p\'et}$. The second is the corresponding category of all pseudocoherent sheaves over ${\mathcal{O}}_\text{\'et}$. The corresponding functor is the corresponding pullback along the corresponding morphism of sites $X_{\text{p\'et}}\rightarrow X_{\text{\'et}}$. 	
\end{corollary}

\newpage\section{Descent in General Setting for Analytic Adic Banach Rings}

\indent We now translate the results in previous section to the corresponding analytic adic Banach context. Again everything is essentially as proposed in \cite[Below Lemma 1.10.4, also see Section 3.8]{9Ked1}.

\subsection{\'Etale Topology}

\begin{setting}
Let $(W,W^+)$ be analytic adic Banach ring. We now fix a corresponding stable basis $\mathbb{H}$ for the corresponding \'etale site of the adic space $\mathrm{Spa}(W,W^+)$, locally consisting of the corresponding compositions of rational localizations and the corresponding finite \'etale morphisms. And as in \cite[Hypothesis 1.10.3]{9Ked1} we need to assume that the corresponding basis is made up of the corresponding adic spectrum of sheafy rings. Assume the corresponding sheafiness of the Huber ring $W$.
\end{setting}

\begin{definition} \mbox{\bf{(After Kedlaya-Liu \cite[Definition 2.5.9]{9KL2})}}
We define the \'etale-stably pseudocoherent module over the corresponding adic Banach $W$ with respect to the corresponding basis $\mathbb{H}$ chosen above for the corresponding \'etale site of the analytic adic space $X$. We define a module over $W$ to be an \'etale stably pseudocoherent module if it is algebraically pseudocoherent (namely formed by the corresponding possibly infinite length resolution of finitely generated and projective modules) and at least complete with respect to the corresponding natural topology and also required to be also complete with respect to the natural topology along some base change with respect to any morphism in $\mathbb{H}$.
\end{definition}

%%%%%%%%%%%%%%%%%%%%%%%%%%%%%%%%!!!!!!!!!!!!

\begin{definition} \mbox{\bf{(After Kedlaya-Liu \cite[Definition 2.5.9]{9KL2})}}
Along the previous definition we have the corresponding notion of \'etale-pseudoflat left (or right respectively) modules with respect to the corresponding chosen basis $\mathbb{H}$. A such module is defined to be a topological module $G$ over $W$ complete with respect to the natural topology and for any right (or left respectively) \'etale stably pseudocoherent module, they will jointly give the $\mathrm{Tor}_1$ vanishing.
\end{definition}

\begin{lemma}\mbox{\bf{(After Kedlaya-Liu \cite[Lemma 2.5.10]{9KL2})}}
One could find another basis $\mathbb{H}'$ which is in our situation contained in the corresponding original basis $\mathbb{H}$ such that any morphism in $\mathbb{H}$ could be \'etale-pseudoflat with respect $\mathbb{H}'$ or just $\mathbb{H}$ itself.	
\end{lemma}

\begin{proof}
See \cref{lemma4.4}.
\end{proof}

\begin{proposition} \mbox{\bf{(After Kedlaya-Liu \cite[Theorem 2.5.11]{9KL2})}}
In our current analytic adic Banach situation, consider the site $X_\text{\'et}$ and consider the basis $\mathbb{H}$. Take any \'etale stably pseudocoherent module $M$ over $W$. Consider the corresponding presheaf by taking the inverse limit throughout all the corresponding base change along morphisms in $\mathbb{H}$. Then we have the corresponding acyclicity of the presheaf over any element in $\mathbb{H}$ and any covering of this element by elements in $\mathbb{H}$.	
\end{proposition}

\begin{proof}
As in \cite[Theorem 2.5.11]{9KL2} apply the corresponding \cite[Proposition 8.2.21]{9KL1}.	
\end{proof}

\begin{proposition}\mbox{\bf{(After Kedlaya-Liu \cite[Lemma 2.5.13]{9KL2})}}
In our current analytic adic Banach situation, the corresponding glueing of \'etale-stably pseudocoherent modules holds in this current situation over $W$ along binary morphisms (namely along binary rational decomposition). 	
\end{proposition}

\begin{proof}
We adapt the argument of \cite[Lemma 2.5.13]{9KL2} to our current situation over analytic adic Banach rings. Take the corresponding map to be $W\rightarrow W_1\bigoplus W_2$, and we denote the two spaces associated $W_1$ and $W_2$ by $Y_1$ and $Y_2$. Establish a corresponding descent datum of $V$-\'etale-stably pseudocoherent modules along this decomposition of $\mathrm{Spa}(W,W^+)$. One can realize this by acyclicity a sheaf $G$ over the whole space $\mathrm{Spa}(W,W^+)$ Then we have for each $k=1,2$ that the corresponding identification of $H^0(Y_{k,\text{\'et}},G)$ with $H^0(Y_k,G)$ and the vanishing of the $i$-th cohomology for higher $i$. By the corresponding results we have in the corresponding rational localization situation we vanishing of the corresponding $H^i(\mathrm{Spa}(W,W^+),G),i\geq 1$. Now consider a finite free covering $C$ of $G$ and take the kernel $K$ to form the corresponding exact sequence:
\[
\xymatrix@C+0pc@R+0pc{
0 \ar[r] \ar[r] \ar[r] &K \ar[r] \ar[r] \ar[r] &C  \ar[r] \ar[r] \ar[r] &G \ar[r] \ar[r] \ar[r] &0 ,   
}
\]
which will give rise to the following short exact sequence:
\[
\xymatrix@C+0pc@R+0pc{
0 \ar[r] \ar[r] \ar[r] &K(\mathrm{Spa}(W,W^+)) \ar[r] \ar[r] \ar[r] &C(\mathrm{Spa}(W,W^+))  \ar[r] \ar[r] \ar[r] &G(\mathrm{Spa}(W,W^+)) \ar[r] \ar[r] \ar[r] &0.   
}
\]
To finish we have also to show the corresponding resulting global section is actually \'etale-stably pseudocoherent. For this, we refer the readers to \cite[Lemma 2.22]{9TX2}.

\end{proof}

\begin{proposition} \mbox{\bf{(After Kedlaya-Liu \cite[Theorem 2.5.14]{9KL2})}}
In our current analytic adic Banach situation, taking the corresponding global section will realize the equivalence between the following two categories. The first is the one of all the corresponding pseudocoherent sheaves over $\mathcal{O}_{\text{\'et}}$. The second is the one of all the corresponding \'etale-stably pseudocoherent modules.	
\end{proposition}

\begin{proof}
This is by considering and applying the corresponding \cite[Lemma 1.10.4]{9Ked1} in our current context. The refinement comes from the acyclicity, the transitivity is straightforward, the corresponding binary rational localization situation is the previous proposition, while the corresponding last condition will basically comes from the corresponding f.p. descent \cite[Chapitre VIII]{9SGAI}.	
\end{proof}

\

\subsection{Pro-\'etale Topology}

\begin{setting}
Let $(W,W^+)$ be an analytic perfectoid adic Banach ring. Assume that the ring $W$ is sheafy.	We assume $p$ is a topologically nilpotent element.
\end{setting}

\begin{remark}
As in the corresponding analytic Huber ring situation, the corresponding glueing of the corresponding pseudocoherent sheaves over the corresponding pro-\'etale site does not need the corresponding notions on the corresponding stability beyond the corresponding \'etale-stably pseudocoherence.	
\end{remark}

\indent We translate the corresponding notions of perfectoid rings to the current context from \cite[Definition 2.1.1]{9Ked1}:

\begin{definition}\mbox{\bf{(Kedlaya \cite[Definition 2.1.1]{9Ked1})}}     \label{definition5.10}
Consider a general analytic adic Banach ring $(W,W^+)$, we will call it perfectoid if there exists a definition ideal $d\subset W^+$ such that $d^p$ contains $p$ and we have the corresponding Frobenius map from $W^+/d$ to $W^+/d^p$ is required to be surjective.	
\end{definition}

\begin{lemma} \mbox{\bf{(After Kedlaya-Liu \cite[Proposition 3.4.3]{9KL2})}}
Suppose that we are taking $(W,W^+)$ to be perfectoid now in the sense of \cref{definition5.10} (note that we are considering the corresponding context of analytic adic Banach situation). Then we have that over the corresponding \'etale site, the corresponding group $H^i(\mathrm{Spa}(W,W^+)_\text{p\'et},\mathcal{O})$ vanishes for each $i>0$, while we have that the corresponding group $H^0(\mathrm{Spa}(W,W^+)_\text{p\'et},\mathcal{O})$ is just isomorphic $W$.	
\end{lemma}

\begin{proof}
See \cite[Proposition 3.4.3]{9KL2}.	
\end{proof}

\indent Now we proceed to consider the corresponding pro-\'etale topology. First we recall the following result from \cite[Lemma 3.4.4]{9KL2}.

\begin{lemma}  \mbox{\bf{(Kedlaya-Liu \cite[Lemma 3.4.4]{9KL2})}}
Consider the ring $W$ as above which is furthermore assumed to be perfectoid in the sense of instead \cref{definition5.10}, and consider any direct system of faithfully finite \'etale morphisms as:
\[
\xymatrix@C+0pc@R+0pc{
W_0 \ar[r] \ar[r] \ar[r] &W_1 \ar[r] \ar[r] \ar[r] &W_2  \ar[r] \ar[r] \ar[r] &...,   
}
\]	
where $W_0$ is just the corresponding base $W$. Then the corresponding completion of the infinite level could be decomposed as:
\begin{align}
W_0\oplus\widehat{\bigoplus}_{k=0}^\infty W_k/W_{k-1}.	
\end{align}
As in \cite[Lemma 3.4.4]{9KL2} we endow the corresponding completion mentioned above with the corresponding seminorm spectral. And for the corresponding quotient we endow with the corresponding quotient norm.
\end{lemma}

\

\begin{corollary}\mbox{\bf{(After Kedlaya-Liu \cite[Corollary 3.4.5]{9KL2})}} 
Starting with an analytic adic Banach ring which is as in \cite[Corollary 3.4.5]{9KL2} assumed to be perfectoid in the sense of \cref{definition5.10}. And as in the previous two lemmas we consider an admissible infinite direct system:
\[
\xymatrix@C+0pc@R+0pc{
W_{0} \ar[r] \ar[r] \ar[r] &W_{1} \ar[r] \ar[r] \ar[r] &W_{2}  \ar[r] \ar[r] \ar[r] &...,   
}
\]
whose infinite level will be assumed to take the corresponding spectral seminorm as in the \cite[Corollary 3.4.5]{9KL2}. Then we have in such situation the corresponding 2-pseudoflatness of the corresponding embedding map:
\begin{align}
W_{0}\rightarrow W_{\infty}.
\end{align}
Here $W_{\infty}$ denotes the corresponding completion of the limit $\varinjlim_{k\rightarrow\infty}W_{k}$.
\end{corollary}

\begin{proof}
See \cite[Corollary 3.4.5]{9KL2}.	
\end{proof}

\indent Now recall that from \cite{9KL2} since we do not have modified the corresponding underlying spatial context, so we will also only have to consider the corresponding stability with respect to \'etale topology even although we are considering the corresponding profinite \'etale site in our current section. We first generalize the corresponding Tate acyclicity in \cite{9KL2} to the corresponding our current situation:

\begin{proposition}\mbox{\bf{(After Kedlaya-Liu \cite[Theorem 3.4.6]{9KL2})}} 
Now we consider the corresponding pro-\'etale site of $X$ which we denote it by $X_\text{p\'et}$. Then we assume we have a stable basis $\mathbb{H}$ of the corresponding perfectoid subdomains for this profinite \'etale site where each morphism therein will be \'etale pseudoflat. Then we consider as in \cite[Theorem 3.4.6]{9KL2} a corresponding module over $W$ which is assumed to be \'etale-stably pseudocoherent. Then consider the corresponding presheaf $\widetilde{G}$ attached to this \'etale-stably pseudocoherent with respect to in our situation the corresponding chosen well-defined basis $\mathbb{H}$, we have that the corresponding sheaf over some element $Y\in \mathbb{H}$ (that is to say over $\mathcal{O}_{Y_\text{p\'et}}$) is acyclic and is acyclic with respect to some \v{C}eck covering coming from elements in $\mathbb{H}$. 
\end{proposition}

\begin{proof}
See \cref{proposition2.7}.
\end{proof}

\begin{definition}\mbox{\bf{(After Kedlaya-Liu \cite[Definition 3.4.7]{9KL2})}} 
Consider the pro-\'etale site of $X$ attached to the adic Banach ring $(W,W^+)$. We will define a sheaf of module $G$ over $\widehat{\mathcal{O}}_\text{p\'et}$ to be pseudocoherent if locally we can define this as a sheaf attached to an \'etale-stably pseudocoherent module. As in \cite[Definition 3.4.7]{9KL2}, we do not have to consider the corresponding notion of profinite-\'etale-stably pseudocoherent module.	
\end{definition}

\begin{proposition}\mbox{\bf{(After Kedlaya-Liu \cite[Theorem 3.4.8]{9KL2})}}
Taking the corresponding global section will realize the corresponding equivalence between the following two categories. The first one is the corresponding one of all the pseudocoherent sheaves over $\widehat{\mathcal{O}}_\text{p\'et}$, while the second one is the corresponding one of all the \'etale-stably pseudocoherent modules over $W$.
	
\end{proposition}

\begin{proof}
As in \cite[Theorem 3.4.8]{9KL2}, the corresponding analog of \cite[Proposition 9.2.6]{9KL1} applies in the way that the corresponding conditions of the analog of \cite[Proposition 9.2.6]{9KL1} hold in our current situation. Also see \cite[Theorem 2.9.9, Remark 2.9.10 and Lemma 1.10.4]{9Ked1}. Certainly in this situation we need to translate the corresponding statemens in \cite[Theorem 2.9.9, Remark 2.9.10 and Lemma 1.10.4]{9Ked1} to the corresponding adic Banach setting.	
\end{proof}

\indent Obviously we have the following analog of \cite[Corollary 3.4.9]{9KL2}:

\begin{corollary} \mbox{\bf{(After Kedlaya-Liu \cite[Corollary 3.4.9]{9KL2})}}
The following two categories are equivalent. The first is the corresponding category of all pseudocoherent sheaves over $\widehat{\mathcal{O}}_\text{p\'et}$. The second is the corresponding category of all pseudocoherent sheaves over ${\mathcal{O}}_\text{\'et}$. The corresponding functor is the corresponding pullback along the corresponding morphism of sites $X_{\text{p\'et}}\rightarrow X_{\text{\'et}}$. 	
\end{corollary}

\
\subsection{Quasi-Stein Spaces}

We now generalize our discussion in \cite{9TX2} around pseudocoherent sheaves over quasi-Stein spaces to the general base situation without the corresponding restrictive assumption that we are working over $\mathbb{Z}_p$.

\begin{setting}
Assume that we are working over some quasi-Stein space $X$ which could be written as the corresponding direct limit of analytic adic Banach affinoids: $X:=\varinjlim_i  X_i$ such that $X_i$ could be written as the corresponding analytic affinoid $\mathrm{Spa}(W_i,W_i^+)$.	
\end{setting}

%%%%%%%%%%%%%%%%%%%%%%%%%%%%!!!!!!!!!!!!!!!

\begin{lemma}\mbox{\bf{(After Kedlaya-Liu \cite[Lemma 2.6.3]{9KL2})}}
Consider a compatible family of Banach modules over the projective system $\mathcal{O}_{X_i}(X_i)$. Suppose that the corresponding transition map $p_i:\mathcal{O}_{X_{i}}(X_{i})\widehat{\otimes}_{\mathcal{O}_{X_{i+1}}(X_{i+1})} M_{i+1}\overset{}{\rightarrow}M_i$ is surjective. Then we have: I. The corresponding global section of $M$ throughout the limit $\mathcal{O}(X)=\varprojlim_i \mathcal{O}_{X_{i}}(X_{i})$ is dense in each section $M_i$ for any $i\geq 0$; II. The corresponding vanishing of $R^1\varprojlim_{i\rightarrow \infty}$ holds in our situation. (Parallel statement in analytic Huber ring situation holds as well.)  	
\end{lemma}

\begin{proof}
This is basically the general version of the corresponding result in \cite[Lemma 2.6.3]{9KL2}. One could prove this in the same way. 
\end{proof}

\begin{lemma}\mbox{\bf{(After Kedlaya-Liu \cite[Corollary 2.6.4]{9KL2})}}
Consider a compatible family of Banach modules over the projective system $\mathcal{O}_{X_i}(X_i)$. Suppose that each member in the family is stably-pseudocoherent. And suppose that the corresponding transition map $p_i:\mathcal{O}_{X_{i}}(X_{i})\widehat{\otimes}_{\mathcal{O}_{X_{i+1}}(X_{i+1})} M_{i+1}\overset{}{\rightarrow}M_i$ is strictly isomorphism. Then we have in our situation the corresponding projection from the global section $M=\varprojlim_i M_i$ to each $M_i$ (for each $i\geq 0$) is then isomorphism in our current situation as in \cite[Corollary 2.6.4]{9KL2}. (Parallel statement in analytic Huber ring situation holds as well.)    	
\end{lemma}

\begin{proof}
This is basically the general version of the corresponding result in \cite[Corollary 2.6.4]{9KL2}. One could prove this in the same way. 
\end{proof}

\begin{proposition}\mbox{\bf{(After Kedlaya-Liu \cite[Theorem 2.6.5]{9KL2})}}
With the corresponding notations as above, we have that the global section of any pseudocoherent sheaf $M$ over $X$ will be dense in the section over any quasicompact. (Parallel statement in analytic Huber ring situation holds as well.)  	
\end{proposition}

\begin{proof}
See \cite[Theorem 2.6.5]{9KL2}.	
\end{proof}

\begin{proposition}\mbox{\bf{(After Kedlaya-Liu \cite[Theorem 2.6.5]{9KL2})}}
With the corresponding notations as above, we have that the global section of any pseudocoherent sheaf $M$ over $X$ will finitely generate each fiber $M_x$ over $\mathcal{O}_{X,x}$ for any $x\in X$. (Parallel statement in analytic Huber ring situation holds as well.)  
\end{proposition}

\begin{proof}
See \cite[Theorem 2.6.5]{9KL2}.	
\end{proof}

\begin{proposition}\mbox{\bf{(After Kedlaya-Liu \cite[Theorem 2.6.5]{9KL2})}}
With the corresponding notations as above, we have that the corresponding sheaf cohomology of any pseudocoherent sheaf $M$ over $X$ admits vanishing for degree bigger than 0. (Parallel statement in analytic Huber ring situation holds as well.)  
\end{proposition}

\begin{proof}
See \cite[Theorem 2.6.5]{9KL2}.	
\end{proof}

\begin{proposition} \mbox{\bf{(After Kedlaya-Liu \cite[Corollary 2.6.6]{9KL2})}}
For any pseudocoherent sheaf $M$ over the space $X$, we have in our current situation that the exact functor from the corresponding category of all the pseudocoherent sheaves to the corresponding one of all the corresponding $\mathcal{O}_X(X)$-modules through the corresponding taking the global section. (Parallel statement in analytic Huber ring situation holds as well.)  	
\end{proposition}

\begin{proposition}  \mbox{\bf{(After Kedlaya-Liu \cite[Corollary 2.6.8]{9KL2})}}
For any pseudocoherent sheaf $M$ admitting a structure of vector bundles over the space $X$, we have in our current situation that the corresponding sufficient and necessary condition for the global section to be finite projective is exactly that the global section of $M$ is finitely generated. (Parallel statement in analytic Huber ring situation holds as well.)  	
\end{proposition}

\begin{proof}
See \cite[Corollary 2.6.8]{9KL2}, it is straightforward to exact a global splitting from the local ones over quasi-compacts.	
\end{proof}

\begin{proposition}\mbox{\bf{(After Kedlaya-Liu \cite[Proposition 2.6.17]{9KL2})}} 
Assume that the space $X$ is basically $m$-uniform in the sense of \cite{9KL2}. Then we have the corresponding finiteness of the global section is equivalent to to the corresponding uniform finiteness through out all $X_i,i=0,1,2,...$. (Parallel statement in analytic Huber ring situation holds as well.)  	
\end{proposition}

\begin{proof}
See the proof given in \cite[Theorem 6.27]{9TX2}. 	
\end{proof}

\

%%%%%%%%%%%%%%%%%%%%%%%!!!!!!!!!!!!!!!!!!!!!!

\newpage\section{Applications} \label{section6}

\subsection{Application to Noetherian Cases}

\indent The corresponding noetherian situation is expected to be in some sense better than just being pseudoflat with respect to the corresponding rational localization.

\begin{setting}
We now work with the corresponding analytic adic Banach rings over $\mathbb{Q}_p$ or $\mathbb{F}_p((t))$ as the corresponding base spaces. The deformation will happen along some Banach ring over $\mathbb{Q}_p$ or $\mathbb{F}_p((t))$. And moreover we have to assume that the corresponding noetherianness preserves under the corresponding deformation over some ring $V$. Namely we have to assume that $W\widehat{\otimes}_* V$ is always noetherian, for instance this could be achieved when we have that $V$ over $\mathbb{Q}_p$ comes from truncations of distribution algebra over some $p$-adic Lie groups by considering some further $p$-adic microlocal analysis.
\end{setting}

\begin{lemma}
With the notations in the \cite[Lemma 2.4.10]{9KL2}, we have that the corresponding morphisms $W\widehat{\otimes}V\rightarrow B_1\widehat{\otimes}V$, $W\widehat{\otimes}V\rightarrow B_2\widehat{\otimes}V$ and $W\widehat{\otimes}V\rightarrow B_{12}\widehat{\otimes}V$ are 2-pseudoflat. And furthermore in our current context they are actually flat.
\end{lemma}

\begin{proof}
In our current situation, the corresponding 2-pseudoflatness is achieved since this is already true in the corresponding more general situation. However for finite presented modules (we do not have to consider the topological issues since we are in the noetherian situation) this is already flat which implies for finitely generated modules (we do not have to consider the topological issues since we are in the noetherian situation) this is already flat	as well. Then for any module which could be written as injective limit of finite ones, the result holds. 
\end{proof}

\begin{corollary}
The corresponding rational localization with respect to the adic Banach ring $(W,W^+)$ is flat (namely for all the corresponding finitely presented module over $W\widehat{\otimes}V$).	
\end{corollary}

\begin{proposition}
Then we consider the corresponding presheafification of any finitely generated module $M$ over $W\widehat{\otimes}V$, to be more precise over $\mathrm{Spa}(W,W^+)$ we will define the corresponding presheaf $\widetilde{M}$ by taking the inverse limit of the base changes of $M$ throughout all the rational localizations of $W$. Then we have that the Tate glueing property holds in our current situation for such presheaf $\widetilde{M}$.
\end{proposition}

\begin{proof}
This reduces to the previous lemma by using \cite[Proposition 2.4.20, Proposition 2.4.21]{9KL1}.	
\end{proof}

\indent The corresponding glueing finitely generated modules is also achieved in current noetherian situation. First we consider the corresponding result:

\begin{proposition}
The descent along the morphism $W\rightarrow B_1\bigoplus B_2$ is effective as long as one focuses on the category of all the finitely presented modules over $W\widehat{\otimes}V$. 	
\end{proposition}

\begin{proof}
This is the corresponding noetherian implication of the corresponding result in \cite[Lemma 2.14]{9TX2}.	
\end{proof}

\begin{proposition}
The corresponding functor of taking the corresponding global section will give rise to the corresponding equivalence between the corresponding $V$-coherent sheaves and the $V$ finitely presented modules.
\end{proposition}

\begin{proof}
This is also the corresponding noetherian implication of the corresponding result in \cite[Theorem 2.15]{9TX2}.	
\end{proof}

\indent Then fixing a stable basis $\mathbb{H}$ for the \'etale site of the space $\mathrm{Spa}(W,W^+)$, then we have the results over the \'etale site as well:

\begin{proposition}
We consider the corresponding presheafification of any finitely generated module $M$ over $W\widehat{\otimes}V$, to be more precise over $\mathrm{Spa}(W,W^+)$ we will define the corresponding presheaf $\widetilde{M}$ by taking the inverse limit of the base changes of $M$ throughout all the member in the basis $\mathbb{H}$. Then we have that the Tate glueing property holds in our current situation for such presheaf $\widetilde{M}$.
\end{proposition}

\begin{proposition}
The descent (in the \'etale topology) along the morphism $W\rightarrow B_1\bigoplus B_2$ is effective as long as one focuses on the category of all the finitely presented modules over $W\widehat{\otimes}V$. 	
\end{proposition}

\begin{proof}
This is the corresponding noetherian implication of the corresponding result in \cite[Lemma 2.22]{9TX2}.	
\end{proof}

\begin{proposition} \label{proposition6.9}
The corresponding functor (in the \'etale topology) of taking the corresponding global section will give rise to the corresponding equivalence between the corresponding $V$-coherent sheaves and the $V$ finitely presented modules.
\end{proposition}

\begin{proof}
This is also the corresponding noetherian implication of the corresponding result in \cite[Theorem 2.23]{9TX2}.	
\end{proof}

\

\subsection{Application to Descent over Noncommutative $\infty$-Analytic Prestacks after Bambozzi-Ben-Bassat-Kremnizer}

\indent We now contact with the corresponding derived analytic spaces from \cite{9BBBK}. Recall that we have the categories $\mathrm{Simp}(\mathrm{Ind}^m(\mathrm{BanSets}_{E}))$ and $\mathrm{Simp}(\mathrm{Ind}(\mathrm{BanSets}_{E}))$ of the corresponding simplicial sets in the corresponding inductive category of Banach sets over $E=\mathbb{Q}_p,\mathbb{F}_p((t))$.

\begin{theorem}\mbox{\bf{(Bambozzi-Ben-Bassat-Kremnizer)}}
The categories\\ $\mathrm{Simp}(\mathrm{Ind}^m(\mathrm{BanSets}_{E}))$ and $\mathrm{Simp}(\mathrm{Ind}(\mathrm{BanSets}_{E}))$ admit symmetric monoidal model categorical structure.
	
\end{theorem}

\indent Therefore based on this nice structure \cite{9BBBK} defined the corresponding $\infty$-categories of $\mathbb{E}_\infty$-rings:
\begin{align}
\mathrm{sComm}(\mathrm{Simp}(\mathrm{Ind}^m(\mathrm{BanSets}_{E})))\\
\mathrm{sComm}(\mathrm{Simp}(\mathrm{Ind}(\mathrm{BanSets}_{E}))).	
\end{align}

\indent We now consider the corresponding some object $A$ which is a corresponding $\infty$-locally convex ring in the categories above and we consider the corresponding object in the corresponding opposite category we call that $\mathrm{Spec}A$. And we consider the corresponding homotopy Zariski topology, which allows us to talk about the corresponding $\infty$-analytic stacks.

\indent Now recall that a connective $\mathbb{E}_1$-ring is called noetherian if we have that $\pi_0(A)$ is noetherian and $\pi_n(A)$ is basically finitely generated over $\pi_0(A)$.

\begin{remark}
The corresponding Koszul simplicial Banach ring considered in \cite{9BK} is actually concentrated in the nonpositive degrees\footnote{Thanks for Federico Bambozzi for reminding us of this.}, in the cohomological convention, namely $\mathrm{H}^{-n}(.)=0$ for $n<0$. But as in any abstract homotopy theory we consider the corresponding conventional transformation:
\begin{align}
\pi_n:=	\mathrm{H}^{-n},n\geq 0.
\end{align} 	
\end{remark}

\indent Now we consider an $\infty$-analytic stack $X$ as considered in \cite{9BBBK}. And we consider the corresponding ringed site attached to $X$ under the corresponding homotopy Zariski topology, denoted by $(X,\mathcal{O}_X)$.

\begin{example}
The examples which are very interesting to us are the corresponding $\infty$-Huber spectra constructed from any Banach rings over $E$ from \cite{9BK} by using Koszul complex. The corresponding classical sheafiness issue could be really forgotten.	
\end{example}

\indent Now we consider the corresponding $\mathbb{E}_1$-ring $A$ in the corresponding categories\\ $\mathrm{Simp}(\mathrm{Ind}^m(\mathrm{BanSets}_{E}))$ and $\mathrm{Simp}(\mathrm{Ind}(\mathrm{BanSets}_{E}))$. We make the following assumption:

\begin{assumption}
We assume that all the ring objects below are noetherian, as $\mathbb{E}_1$-rings.  	
\end{assumption}

\indent Suppose now we have four $\mathbb{E}_1$ Banach rings $A,A_1,A_2,A_{12}$\footnote{Namely all the homotopy groups have to be Banach instead of Bornological or ind-Banach.} in $\mathrm{Simp}(\mathrm{Ind}^m(\mathrm{BanSets}_{E}))$ and $\mathrm{Simp}(\mathrm{Ind}(\mathrm{BanSets}_{E}))$ respectively. And we assume that we have the following short strictly exact sequence in the sense of a glueing square:
\[
\xymatrix@C+0pc@R+0pc{
0 \ar[r] \ar[r] \ar[r] &\pi_0(A) \ar[r] \ar[r] \ar[r] &\pi_0(A_1)\oplus \pi_0(A_2)  \ar[r] \ar[r] \ar[r] &\pi_0(A_{12}) \ar[r] \ar[r] \ar[r] &0 ,   
}
\]
and we assume that the image of $\pi_0(A_i)$ in $\pi_0(A_{12}) $ is dense for each $i=1,2$. And we assume the corresponding morphism $\pi_0(A_i)\rightarrow \pi_0(A_{12})$ is flat for $i=1,2$. And we assume that $A,A_1,A_2,A_{12}$ form the corresponding derived glueing sequence.

\begin{example}
One can construct the following example. First consider any glueing sequence of Banach rings (not required to be sheafy or commutative) over $\mathbb{Q}_p$:
\[
\xymatrix@C+0pc@R+0pc{
0 \ar[r] \ar[r] \ar[r] &\Pi \ar[r] \ar[r] \ar[r] &\Pi_1\oplus \Pi_2  \ar[r] \ar[r] \ar[r] &\Pi_{12} \ar[r] \ar[r] \ar[r] &0.  
}
\]	
This is very natural in the corresponding situation for example where we deform from a corresponding nice short exact sequence of commutative sheafy rings of the same type, but in any rate we do not require the corresponding commutativity or the correponding sheafyness. Then take any derived global section of some $\mathrm{Spa}^h(R)$ from \cite{9BK}, denoted by $R^h$. Over $\mathbb{Q}_p$, we have the following short strictly exact sequence:
\[
\xymatrix@C+0pc@R+0pc{
0 \ar[r] \ar[r] \ar[r] &\Pi\widehat{\otimes}_{\mathbb{Q}_p} \pi_0(R^h) \ar[r] \ar[r] \ar[r] &\Pi_1\widehat{\otimes}_{\mathbb{Q}_p} \pi_0(R^h)\oplus \Pi_2\widehat{\otimes}_{\mathbb{Q}_p} \pi_0(R^h)  \ar[r] \ar[r] \ar[r] &\Pi_{12}\widehat{\otimes}_{\mathbb{Q}_p} \pi_0(R^h) \ar[r] \ar[r] \ar[r] &0 ,   
}
\]
but since the corresponding ring $\Pi,\Pi_1,\Pi_2,\Pi_{12}$ will actually realize the situation where $\Pi_*\widehat{\otimes} \pi_0(R^h)\overset{\sim}{\rightarrow}\pi_0(\Pi_*\widehat{\otimes}^\mathbb{L}R^h)$, we will have the desired situation as long as restrict to now the noetherian situation.
\end{example}

\begin{conjecture}
The map $A\rightarrow \prod_{i=1,2} A_i$ is an effective descent morphism with respect to the finitely presented left module spectra (meaning we have finitely presented $\pi_0$).	
\end{conjecture}

\begin{proposition}
The operation of taking equalizer along the corresponding map $\prod_{i=1,2} A_i\rightarrow A_{12}$ preserves the property of being finitely presented.	
\end{proposition}

\begin{proof}
The situation where we deform from the corresponding $\mathbb{E}_\infty$ by some noncommutative deformation could be implied by \cite[Lemma 2.14]{9TX2}. While in general this follows from \cite[Lemma 6.83]{9TX4}. Note that we achieve the finiteness for each $\pi_n(M)$ just by application of the results we know in the classical deformed situation, which is because the connecting homomorphism vanishes on each level.
\end{proof}

\indent What could happen in the commutative setting is actually also interesting. We have the following example in mind:

\begin{example}
Consider the context of the previous example, but assume that all the rings involved are Banach commutative. First we consider over $\mathbb{Q}_p$:
\[
\xymatrix@C+0pc@R+0pc{
0 \ar[r] \ar[r] \ar[r] &\Pi \ar[r] \ar[r] \ar[r] &\Pi_1\oplus \Pi_2  \ar[r] \ar[r] \ar[r] &\Pi_{12} \ar[r] \ar[r] \ar[r] &0,  
}
\]		
a strictly exact sequence of sheafy rings. We now deform this along some spectrum $\mathrm{Spa}^h(R)$ in \cite{9BK} which will produce a very interesting situation where we can actually obtain desired situation for the descent of finitely presented module spectra. And note that in commutative setting we could also have more geometric contexts to rely on as in \cite{9BK} and \cite{9BBBK}. This may have contact with the corresponding derived Galois deformation theory of \cite{9GV} for instance by considering the corresponding simplicial pro-Artin rings.
\end{example}

\begin{remark}
The corresponding machinery from \cite{9CS} should definitely reflect similar things here even in the $\mathbb{E}_1$-ring context. Also we would like to mention in the commutative setting that actually the descent for certain quasi-coherent modules are also considered extensively in \cite{9BBK}. Strikingly the ideas in \cite{9BBK} (although developed in a quite commutative setting) come from partially work from Kontsevich-Rosenberg \cite{9KR} namely essentially the noncommutative descent. Our feeling is that definitely the descent for certain quasi-coherent modules in \cite{9BBBK} could be established to noncommutative setting in some form both in the archimedean and nonarchimedean situations. 
\end{remark}

\

This chapter is based on the following paper, where the author of this dissertation is the main author:
\begin{itemize}
\item Tong, Xin. "Period Rings with Big Coefficients and Applications III." arXiv preprint arXiv:arXiv:2102.10692(2021).
\end{itemize}

\newpage

\bibliographystyle{ams}

\newpage\chapter{Topics on Geometric and Representation Theoretic Aspects of Period Rings I}

\newpage\section{Introduction}

\subsection{General Perfectoids and Preperfectoids}

\indent Our previous work on Hodge-Iwasawa theory aimed at the corresponding deformation of the corresponding Hodge structure after \cite{10KL1} and \cite{10KL2}. The corresponding application in our mind targets the corresponding noncommutative Iwasawa theory and the corresponding $p$-adic local systems with general Banach coefficients. The corresponding background of the foundations on the analytic geometry comes from essentially \cite{10KL1} where Banach analog of Fontaine's perfectoid spaces in the Tate case was defined and studied extensively. \\

\indent The corresponding analytic Huber analog of the main results of \cite{10KL1} was initiated in \cite{10Ked1} where the analytic Huber analog of perfectoids are defined, and the analytic Huber analog of the corresponding perfectoid correspondence is established and the corresponding descent for vector bundles and stably-pseudocoherent sheaves is established in the analytic topology.\\

\indent In our current development of the corresponding Hodge-Iwasawa theory we consider the corresponding generality at least parallel to \cite{10Ked1}. We now consider the corresponding analytic analog of our previous consideration which generalized the corresponding \cite{10KL1} and \cite{10KL2} to the corresponding big coefficient situation. Without deformation in our sense \cite{10KL1} and \cite{10KL2} have already considered the corresponding free of trivial norm cases in certain context and considerations. \\

\indent Certainly inverting $p$ takes back us to the Tate situation, but what we know is that this is at least not included in our previous consideration, since the relative Robba rings are not the same (even in our current situation many are still Tate as their owns). So we have to basically establish the corresponding parallel theory, in order to treat the situation where the base spaces or rings are not Tate. And more importantly the integral Robba rings are not Tate, which needs to be discussed in order to apply to the situation where the base spaces or rings are not Tate. \\

\indent With the corresponding notations in Chapter 5, we have:

\begin{theorem} \label{10theorem10.1.1}
Over $A$-relative preperfectoid Robba rings $\widetilde{\Pi}_{?,R,A}$\\ ($?=\empty,r,\infty,[r_1,r_2],\{[s,r]\}$, $0<r_1\leq r_2/p^{hk}$)  constructed over analytic base $(R,R^+)$, we have the equivalence among the categories of the Frobenius \'etale-stably-pseudocoherent modules, which could be further compared in equivalence to the the pseudocoherent sheaves over adic Fargues-Fontaine curves in both \'etale and pro-\'etale topology.	
\end{theorem}

\indent This is basically the corresponding commutative version where we encode the corresponding ring $A$ into the corresponding Huber spectrum in \cite{10KL1} and \cite{10KL2}. The corresponding noncommutative version is also true but we have to do this slightly different by deforming the structure sheaves.  With the corresponding notations again in Chapter 5, we have:

\begin{theorem} \label{10theorem10.1.2}
Over $B$-relative preperfectoid Robba rings $\widetilde{\Pi}_{?,R,B}$\\ ($?=\empty,r,\infty,[r_1,r_2],\{[s,r]\}$, $0<r_1\leq r_2/p^{hk}$) constructed over analytic base $(R,R^+)$, we have the equivalence among the categories of the Frobenius $B$-\'etale-stably-pseudocoherent modules, which could be further compared in equivalence to the the $B$-pseudocoherent sheaves over adic Fargues-Fontaine curves in \'etale topology.	
\end{theorem}

\indent Our ultimate consideration will after \cite{10GR}, where $(R,R^+)$ is not necessarily always assumed to be analytic, especially in the integral setting this is quite general. With the corresponding notations again in Chapter 6, we have:

\begin{theorem} \label{10theorem10.1.3}
Over $A$-relative preperfectoid Robba rings $\widetilde{\Pi}_{?,R,A}$ ($?=[r_1,r_2],\{[s,r]\}$, $0<r_1\leq r_2/p^{hk}$) constructed over analytic base $(R,R^+)$, we have the equivalence of the categories of the Frobenius \'etale-stably-pseudocoherent modules, which could be further compared in equivalence to the the pseudocoherent sheaves over adic Fargues-Fontaine curves in the \'etale or pro-\'etale topology.	
\end{theorem}

\begin{theorem} \label{10theorem10.1.4}
Over $B$-relative preperfectoid Robba rings $\widetilde{\Pi}_{?,R,B}$ ($?=[r_1,r_2],\{[s,r]\}$, $0<r_1\leq r_2/p^{hk}$) constructed over analytic base $(R,R^+)$, we have the equivalence of the categories of the Frobenius $B$-stably-pseudocoherent modules.	
\end{theorem}

\subsection{Some Consideration in the Future}

More thorough applications of \cite{10Lu2} will be expected in our study on the corresponding derived algebraic geometry of the corresponding period rings\footnote{The original motivation comes from the study of the algebraic pseudocoherent sheaves over algebraic Fargues-Fontaine curves originally in \cite{10KL2}, and also in \cite{10XT1} and \cite{10XT2} where the corresponding functional analytic information around locally convex vector spaces (after Bourbaki) is not considered.}. We only touched this to some very transparent extent which will be somehow reflecting some further consideration along the corresponding applications in our mind.\\

New robust developments from Clausen-Scholze \cite{10CS}, Bambozzi-Kremnizer \cite{10BK1} and Bambozzi-Ben-Bassat-Kremnizer \cite{10BBBK} may shed lights on the derived Galois (possibly in more general setting for fundamental groups) deformation theory (and some unobstructed crystalline substacks) such as in \cite{10GV} and possibly (we conjecture) the corresponding possible theory of some derived eigenvarieties and derived Selmer varieties in Kim's nonabelian Chabauty. We hope to give more hard thorough study and understanding around this. Especially in our situation we will come across many technical issues where these new developments could help.

\

%%\newpage

\newpage\section{Period Rings in General Setting with General Coefficients}

\indent In this section we start by defining Kedlaya-Liu's style period rings with coefficients in Tate adic Banach rings which are of finite type.

\begin{setting}
We consider the corresponding setting up which is as in the following. First we consider a field $E$ which is analytic nonarchimedean with normalized norm such that the corresponding uniformizer $\pi_E$ is of norm $1/p$. And we assume that the corresponding residue field of the field $E$ takes the form of $\mathbb{F}_{p^h}$. Then we fix a corresponding uniform adic Banach algebra $(R,R^+)$ over $\mathbb{F}_{p^h}$ which is not required to be contain topological nilpotent unit. However we need to assume that that this is analytic in the sense of \cite{10KL1} and \cite[Definition 1.1.2]{10Ked1}. We recall that this means the set of the corresponding topological nilpotents generates the unit ideal. And let $A$ be a general rigid analytic affinoid over $\mathbb{Q}_p$ or over $\mathbb{F}_p((t))$.
\end{setting}

\indent Then in this generality we define our large coefficient Robba rings following Kedlaya-Liu in the Tate algebra situation namely we have $A=\mathbb{Q}_p\{X_1,...,X_d\}$ or\\ $A=\mathbb{F}_p((t))\{X_1,...,X_d\}$:

%%%%%%%%%%%%%%%%%%%%%%%%%%%%%%%%%%!!!!!!!!!

\begin{definition}\mbox{\bf{(After Kedlaya-Liu \cite[Definition 4.1.1]{10KL2})}} Now consider the following constructions. First we consider the corresponding Witt vectors coming from the corresponding adic ring $(R,R^+)$. First we consider the corresponding generalized Witt vectors with respect to $(R,R^+)$ with the corresponding coefficients in the Tate algebra with the general notation $W(R^+)[[R]]$. The general form of any element in such deformed ring could be written as $\sum_{i\geq 0,i_1\geq 0,...,i_d\geq 0}\pi^i[\overline{y}_i]X_1^{i_1}...X_d^{i_d}$. Then we take the corresponding completion with respect to the following norm for some radius $t>0$:
\begin{align}
\|.\|_{t,A}(\sum_{i\geq 0,i_1\geq 0,...,i_d\geq 0}\pi^i[\overline{y}_i]X_1^{i_1}...X_d^{i_d}):= \max_{i\geq 0,i_1\geq 0,...,i_d\geq 0}p^{-i}\|.\|_R(\overline{y}_i)	
\end{align}
which will give us the corresponding ring $\widetilde{\Pi}_{\mathrm{int},t,R,A}$ such that we could put furthermore that:
\begin{align}
\widetilde{\Pi}_{\mathrm{int},R,A}:=\bigcup_{t>0} \widetilde{\Pi}_{\mathrm{int},t,R,A}.	
\end{align}
Then as in \cite[Definition 4.1.1]{10KL2}, we now put the ring $\widetilde{\Pi}_{\mathrm{bd},t,R,A}:=\widetilde{\Pi}_{\mathrm{int},t,R,A}[1/\pi]$ and we set:
\begin{align}
\widetilde{\Pi}_{\mathrm{bd},R,A}:=\bigcup_{t>0} \widetilde{\Pi}_{\mathrm{bd},t,R,A}.	
\end{align}
The corresponding Robba rings with respect to some intervals and some radius could be defined in the same way as in \cite[Definition 4.1.1]{10KL2}. To be more precise we consider the completion of the corresponding ring $W(R^+)[[R]][1/\pi]$ with respect to the following norm for some $t>0$ where $t$ lives in some prescribed interval $I=[s,r]$: 
\begin{align}
\|.\|_{t,A}(\sum_{i,i_1\geq 0,...,i_d\geq 0}\pi^i[\overline{y}_i]X_1^{i_1}...X_d^{i_d}):= \max_{i\geq 0,i_1\geq 0,...,i_d\geq 0}p^{-i}\|.\|_R(\overline{y}_i).	
\end{align}
This process will produce the corresponding Robba rings with respect to  the given interval $I=[s,r]$. Now for particular sorts of intervals $(0,r]$ we will have the corresponding Robba ring $\widetilde{\Pi}_{r,R,A}$ and we will have the corresponding Robba ring $\widetilde{\Pi}_{\infty,R,A}$	if the corresponding interval is taken to be $(0,\infty)$. Then in our situation we could just take the corresponding union throughout all the radius $r>0$ to define the corresponding full Robba ring taking the notation of $\widetilde{\Pi}_{R,A}$.
\end{definition}

\begin{remark}
The corresponding Robba rings $\widetilde{\Pi}_{\mathrm{bd},R,A}$, $\widetilde{\Pi}_{R,A}$, $\widetilde{\Pi}_{I,R,A}$, $\widetilde{\Pi}_{r,R,A}$, $\widetilde{\Pi}_{\infty,R,A}$ are actually themselves Tate adic Banach rings. However in many further application the non-Tateness of the ring $R$ will cause some reason for us to do the corresponding modification, which is considered on this level in fact in \cite{10KL1} in the context therein.	
\end{remark}

\begin{definition}
Then for any general affinoid algebra $A$ over the corresponding base analytic field, we just take the corresponding quotients of the corresponding rings defined in the previous definition over some Tate algebras in rigid analytic geometry, with the same notations though $A$ now is more general. Note that one can actually show that the definition does not depend on the corresponding choice of the corresponding presentations over $A$.
\end{definition}

\begin{remark}
Again in this situation more generally, the corresponding Robba rings $\widetilde{\Pi}_{\mathrm{bd},R,A}$, $\widetilde{\Pi}_{R,A}$, $\widetilde{\Pi}_{I,R,A}$, $\widetilde{\Pi}_{r,R,A}$, $\widetilde{\Pi}_{\infty,R,A}$ are actually themselves Tate adic Banach rings.	
\end{remark}

\indent Note that we can also as in the situation of \cite{10KL2} and \cite{10XT2} consider the corresponding the corresponding property checking of the corresponding period rings defined aboave. We collect the corresponding statements here while the the proof could be found in \cite{10XT2}:

\begin{lemma} \mbox{\bf{(After Kedlaya-Liu \cite[Lemma 5.2.6]{10KL2})}}
For any two radii $0<r_1<r_2$ we have the corresponding equality:
\begin{align}
\widetilde{\Pi}_{\mathrm{int},r_2,R,\mathbb{Q}_p\{T_1,...,T_d\}}\bigcap \widetilde{\Pi}_{[r_1,r_2],R,\mathbb{Q}_p\{T_1,...,T_d\}}	=\widetilde{\Pi}_{\mathrm{int},r_1,R,\mathbb{Q}_p\{T_1,...,T_d\}}.
\end{align}

\end{lemma}

\begin{proof}
See \cite[Lemma 5.2.6]{10KL2} and \cite[Proposition 2.13]{10XT2}.	
\end{proof}

\begin{lemma} \mbox{\bf{(After Kedlaya-Liu \cite[Lemma 5.2.6]{10KL2})}}
For any two radii $0<r_1<r_2$ we have the corresponding equality:
\begin{align}
\widetilde{\Pi}_{\mathrm{int},r_2,R,\mathbb{F}_p((t))\{T_1,...,T_d\}}\bigcap \widetilde{\Pi}_{[r_1,r_2],R,\mathbb{F}_p((t))\{T_1,...,T_d\}}	=\widetilde{\Pi}_{\mathrm{int},r_1,R,\mathbb{F}_p((t))\{T_1,...,T_d\}}.
\end{align}

\end{lemma}

\begin{proof}
See \cite[Lemma 5.2.6]{10KL2} and \cite[Proposition 2.13]{10XT2}.	
\end{proof}

\begin{lemma} \mbox{\bf{(After Kedlaya-Liu \cite[Lemma 5.2.6]{10KL2})}}
For general affinoid $A$ as above (over $\mathbb{Q}_p$ or $\mathbb{F}_p((t))$) and for any two radii $0<r_1<r_2$ we have the corresponding equality:
\begin{align}
\widetilde{\Pi}_{\mathrm{int},r_2,R,A}\bigcap \widetilde{\Pi}_{[r_1,r_2],R,A}	=\widetilde{\Pi}_{\mathrm{int},r_1,R,A}.
\end{align}

\end{lemma}

\begin{proof}
See \cite[Lemma 5.2.6]{10KL2} and \cite[Proposition 2.14]{10XT2}.	
\end{proof}

\begin{lemma} \mbox{\bf{(After Kedlaya-Liu \cite[Lemma 5.2.10]{10KL2})}}
For any four radii $0<r_1<r_2<r_3<r_4$ we have the corresponding equality:
\begin{align}
\widetilde{\Pi}_{[r_1,r_3],R,\mathbb{Q}_p\{T_1,...,T_d\}}\bigcap \widetilde{\Pi}_{[r_2,r_4],R,\mathbb{Q}_p\{T_1,...,T_d\}}	=\widetilde{\Pi}_{[r_1,r_4],R,\mathbb{Q}_p\{T_1,...,T_d\}}.
\end{align}

\end{lemma}

\begin{proof}
See \cite[Lemma 5.2.10]{10KL2} and \cite[Proposition 2.16]{10XT2}.	
\end{proof}

\begin{lemma} \mbox{\bf{(After Kedlaya-Liu \cite[Lemma 5.2.10]{10KL2})}}
For any four radii $0<r_1<r_2<r_3<r_4$ we have the corresponding equality:
\begin{align}
\widetilde{\Pi}_{[r_1,r_3],R,\mathbb{F}_p((t))\{T_1,...,T_d\}}\bigcap \widetilde{\Pi}_{[r_2,r_4],R,\mathbb{F}_p((t))\{T_1,...,T_d\}}	=\widetilde{\Pi}_{[r_1,r_4],R,\mathbb{F}_p((t))\{T_1,...,T_d\}}.
\end{align}

\end{lemma}

\begin{proof}
See \cite[Lemma 5.2.10]{10KL2} and \cite[Proposition 2.16]{10XT2}.	
\end{proof}

\begin{lemma} \mbox{\bf{(After Kedlaya-Liu \cite[Lemma 5.2.10]{10KL2})}}
For any four radii $0<r_1<r_2<r_3<r_4$ we have the corresponding equality:
\begin{align}
\widetilde{\Pi}_{[r_1,r_3],R,A}\bigcap \widetilde{\Pi}_{[r_2,r_4],R,A}	=\widetilde{\Pi}_{[r_1,r_4],R,A}.
\end{align}

\end{lemma}

\begin{proof}
See \cite[Lemma 5.2.10]{10KL2} and \cite[Proposition 2.17]{10XT2}.	
\end{proof}

%%%%%%%%%%%%%%%%%%%%%%%%%%%%%%%!!!!!!!!!!!!!!!!!!!

\

%%\newpage

\newpage\section{The Frobenius Modules and Frobenius Sheaves over the Period Rings} \label{10chapter3}
 
\indent We now consider the corresponding Frobenius actions and the corresponding action reflected on the corresponding Hodge-Iwasawa modules over general analytic $R$. As in the Tate situation we consider the corresponding 'arithmetic' Frobenius as well.

\begin{setting}
We are going to use the corresponding notation $F$ to denote the corresponding relative Frobenius induced from the corresponding $p^h$-power absolute Frobenius from $R$. This would mean that when we consider the corresponding Frobenius up to higher order.
\end{setting}

\indent Furthermore in our situation we have the corresponding sheaves of period rings (carrying the coefficient $A$) over $\mathrm{Spa}(R,R^+)$ in the corresponding pro-\'etale site:
\begin{align}
\widetilde{\Pi}_{\mathrm{bd},\mathrm{Spa}(R,R^+)_\text{pro\'et},A}, \widetilde{\Pi}_{\mathrm{Spa}(R,R^+)_\text{pro\'et},A}, \widetilde{\Pi}_{I,\mathrm{Spa}(R,R^+)_\text{pro\'et},A}, \widetilde{\Pi}_{r,\mathrm{Spa}(R,R^+)_\text{pro\'et},A}, \widetilde{\Pi}_{\infty,\mathrm{Spa}(R,R^+)_\text{pro\'et},A}	
\end{align}
which we will also briefly write as:
\begin{align}
\widetilde{\Pi}_{\mathrm{bd},*,A}, \widetilde{\Pi}_{*,A}, \widetilde{\Pi}_{I,*,A}, \widetilde{\Pi}_{r,*,A}, \widetilde{\Pi}_{\infty,*,A}	
\end{align}
when the corresponding topology is basically clear to the readers (one can also consider other kinds of topology in the same fashion such as the corresponding \'etale situation in our current context).

\begin{definition} \mbox{\bf{(After Kedlaya-Liu \cite[Definition 4.4.4]{10KL2})}} The corresponding Frobenius modules over the corresponding period sheaves will be defined to be finite projective modules carrying the semilinear action from $F^k$ (where we allow some power $k>0$) with some pullback isomorphic requirement which could be defined in the following way. Over the rings without intervals or radius this will mean that we have $F^{k*}M\overset{\sim}{\rightarrow}M$. While over the corresponding Robba rings with radius $r>0$ this will mean that we have $F^{k*}M\otimes_{*_{r/p^{hk}}} *_{r/p^{hk}}\overset{\sim}{\rightarrow}M\otimes_{*_r}*_{r/p^{hk}}$. While over the corresponding Robba rings with interval $[r_1,r_2]$ this will mean that we have $F^{k*}M\otimes_{*_{[r_1/p^{hk},r_2/p^{hk}]}} *_{[r_1,r_2/p^{hk}]}\overset{\sim}{\rightarrow}M\otimes_{*_{[r_1,r_2]}}*_{[r_1,r_2/p^{hk}]}$. Over the corresponding Robba rings without any intervals or radii we assume that they are base change from some modules over the corresponding Robba rings with some radii. 
	
\end{definition}

%%%%%%%%%%%%%%%%%%%%%%%%%%%%%%%%%%%%!!!!!!!!!!!!!!!!!!

\begin{definition} \mbox{\bf{(After Kedlaya-Liu \cite[Definition 4.4.4]{10KL2})}} The corresponding pseudocoherent or fpd Frobenius modules over the corresponding period sheaves will be defined to be pseudocoherent or fpd modules carrying the semilinear action from $F^k$ (where we allow some power $k>0$) with some pullback isomorphic requirement which could be defined in the following way. Over the rings without intervals or radius this will mean that we have $F^{k*}M\overset{\sim}{\rightarrow}M$. While over the corresponding Robba rings with radius $r>0$ this will mean that we have $F^{k*}M\otimes_{*_{r/p^{hk}}} *_{r/p^{hk}}\overset{\sim}{\rightarrow}M\otimes_{*_r}*_{r/p^{hk}}$. While over the corresponding Robba rings with interval $[r_1,r_2]$ this will mean that we have $F^{k*}M\otimes_{*_{[r_1/p^{hk},r_2/p^{hk}]}} *_{[r_1,r_2/p^{hk}]}\overset{\sim}{\rightarrow}M\otimes_{*_{[r_1,r_2]}}*_{[r_1,r_2/p^{hk}]}$. Over the corresponding Robba rings without any intervals or radii we assume that they are base changes from some modules over the corresponding Robba rings with some radii. And then all the modules are required to be complete with respect to the natural topology over any perfectoid subdomain within the corresponding pro-\'etale topology, and the corresponding modules over the corresponding Robba rings with respect to some radius $r>0$ will basically be required to be base change to any \'etale-stably pseudocoherent (over any specific chosen perfectoid subdomain) modules (note that this will include the corresponding glueing along the space $A$ with respect to the corresponding implicit \'etale topology induced over $A$). 
	
\end{definition}

\

\indent Then as in \cite[Definition 4.4.4]{10KL2} we consider the corresponding Frobenius modules over the period rings instead of period sheaves (and we do not either put topological conditions, which will be added later in specific consideration):

%%%%%%%%%%%%%%%%%%%%%%%%%%!!!!!!!!!!!!!

\begin{definition} \mbox{\bf{(After Kedlaya-Liu \cite[Definition 4.4.4]{10KL2})}} The corresponding Frobenius modules over the corresponding period rings will be defined to be finite projective modules carrying the semilinear action from $F^k$ (where we allow some power $k>0$) with some pullback isomorphic requirement which could be defined in the following way. Over the rings without intervals or radius this will mean that we have $F^{k*}M\overset{\sim}{\rightarrow}M$. While over the corresponding Robba rings with radius $r>0$ this will mean that we have $F^{k*}M\otimes_{*_{r/p^{hk}}} *_{r/p^{hk}}\overset{\sim}{\rightarrow}M\otimes_{*_r}*_{r/p^{hk}}$. While over the corresponding Robba rings with interval $[r_1,r_2]$ this will mean that we have $F^{k*}M\otimes_{*_{[r_1/p^{hk},r_2/p^{hk}]}} *_{[r_1,r_2/p^{hk}]}\overset{\sim}{\rightarrow}M\otimes_{*_{[r_1,r_2]}}*_{[r_1,r_2/p^{hk}]}$. Over the corresponding Robba rings without any intervals or radii we assume that they are base change from some modules over the corresponding Robba rings with some radii. 
	
\end{definition}

\begin{remark}
As in \cite[Definition 4.4.4]{10KL2} we did not out any topological conditions on the ring theoretic and algebraic representation theoretic objects defined above, but this will be definitely more precise in later development. Certainly this is more relevant with respect to the following definition.	
\end{remark}

\begin{definition} \mbox{\bf{(After Kedlaya-Liu \cite[Definition 4.4.4]{10KL2})}} The corresponding pseudocoherent or fpd Frobenius modules over the corresponding period rings will be defined to be pseudocoherent or fpd modules carrying the semilinear action from $F^k$ (where we allow some power $k>0$) with some pullback isomorphic requirement which could be defined in the following way. Over the rings without intervals or radius this will mean that we have $F^{k*}M\overset{\sim}{\rightarrow}M$. While over the corresponding Robba rings with radius $r>0$ this will mean that we have $F^{k*}M\otimes_{*_{r/p^{hk}}} *_{r/p^{hk}}\overset{\sim}{\rightarrow}M\otimes_{*_r}*_{r/p^{hk}}$. While over the corresponding Robba rings with interval $[r_1,r_2]$ this will mean that we have $F^{k*}M\otimes_{*_{[r_1/p^{hk},r_2/p^{hk}]}} *_{[r_1,r_2/p^{hk}]}\overset{\sim}{\rightarrow}M\otimes_{*_{[r_1,r_2]}}*_{[r_1,r_2/p^{hk}]}$. Over the corresponding Robba rings without any intervals or radii we assume that they are base change from some modules over the corresponding Robba rings with some radii. 
	
\end{definition}

\indent Then we define over our current analytic base the corresponding pseudocoherent and fpd sheaves carrying the corresponding Frobenius modules.

\begin{definition}\mbox{\bf{(After Kedlaya-Liu \cite[Definition 4.4.6]{10KL2})}}
Over the Robba rings with respect to the corresponding all finite intervals contained in $(0,\infty)$, we consider a family $(M_{[s,r]})_{[s,r]}$ of finite projective modules with respect to the corresponding intervals. Then we call this family a finite projective $F^k$-bundle over $\widetilde{\Pi}_{R,A}$ if each member in the family is assumed to be a corresponding finite projective $F^k$-modules.	
\end{definition}

\begin{definition}\mbox{\bf{(After Kedlaya-Liu \cite[Definition 4.4.6]{10KL2})}}
Over the Robba rings with respect to the corresponding all finite intervals contained in $(0,\infty)$, we consider a family $(M_{[s,r]})_{[s,r]}$ of fpd or pseudocoherent modules with respect to the corresponding intervals. Then we call this family a fpd or pseudocoherent $F^k$-sheaf over $\widetilde{\Pi}_{R,A}$ if each member in the family is assumed to be a corresponding stably or \'etale-stably fpd or pseudocoherent $F^k$-modules.	
\end{definition}

\begin{proposition} \mbox{\bf{(After Kedlaya-Liu \cite[Lemma 4.4.8]{10KL2})}}
Any pseudocoherent $F^k$-modules over the corresponding period rings above admits a surjective morphism from some finite projective $F^k$-modules over the same types of period rings.	
\end{proposition}

\begin{proof}
This is the corresponding relative and analytic version of the corresponding \cite[Lemma 4.4.8]{10KL2}. The proof is parallel, see \cite[Lemma 4.4.8]{10KL2}.	
\end{proof}

\

%%%%%%%%%%%%%%%%%%%%%%%%%%%%%%%%%%%!!!!!!!!!!!!!!!!

\newpage\section{Comparison in the \'Etale Setting}

\indent We now consider the corresponding parallel consideration to \cite[Section 4.5]{10KL2}. Essentially in our current context, this is really the corresponding purely analytic setting especially when we do not have the chance to invert $\pi$. Certainly we will also consider the corresponding deformed setting. The corresponding space we are going to work on is just the corresponding perfectoid spaces (in the corresponding analytic situation) taking the corresponding general form of $\mathrm{Spa}(R,R^+)$, and consider any topological ring $Z$ and the corresponding sheaf $\underline{Z}$. We will consider the corresponding finite projective, pseudocoherent and finite projective dimension local systems of $\underline{Z}$-modules. Note that again in this analytic setting the corresponding definitions of such modules are locally of the corresponding same types. This means in particular in the situation of finite projective situation the corresponding definition is locally the finite projective sheaves of modules over $\underline{Z}$ instead of the corresponding naive locally finite free ones especially in our large coefficient situation.

\begin{lemma} \mbox{\bf{(After Kedlaya-Liu \cite[Lemma 4.5.3]{10KL2})}}
The corresponding taking element to the corresponding element after the action of the corresponding operator $F-1$ will be to realize the following exact sequences as in \cite[Lemma 4.5.3]{10KL2}:
\[
\xymatrix@C+0pc@R+0pc{
0 \ar[r] \ar[r] \ar[r] &\underline{\mathcal{O}_{E_k}}  \ar[r] \ar[r] \ar[r] &\widetilde{\Omega}_{\mathrm{int},R} \ar[r] \ar[r] \ar[r] &\widetilde{\Omega}_{\mathrm{int},R} \ar[r] \ar[r] \ar[r] &0,\\
0 \ar[r] \ar[r] \ar[r] &\underline{{E}_{k}}  \ar[r] \ar[r] \ar[r] &\widetilde{\Omega}_{R} \ar[r] \ar[r] \ar[r] &\widetilde{\Omega}_{R} \ar[r] \ar[r] \ar[r] &0,\\
0 \ar[r] \ar[r] \ar[r] &\underline{\mathcal{O}_{E_k}}  \ar[r] \ar[r] \ar[r] &\widetilde{\Pi}_{\mathrm{int},R} \ar[r] \ar[r] \ar[r] &\widetilde{\Pi}_{\mathrm{int},R} \ar[r] \ar[r] \ar[r] &0,\\
0 \ar[r] \ar[r] \ar[r] &\underline{{E_k}}  \ar[r] \ar[r] \ar[r] &\widetilde{\Pi}_{R} \ar[r] \ar[r] \ar[r] &\widetilde{\Pi}_{R} \ar[r] \ar[r] \ar[r] &0.\\
}
\]	
\end{lemma}

\begin{proof}
See \cite[Lemma 4.5.3]{10KL2}.	
\end{proof}

\begin{proposition} \mbox{\bf{(After Kedlaya-Liu \cite[Lemma 4.5.4]{10KL2})}}
We have the corresponding strictly pseudocoherence of any corresponding finitely presented modules over $\widetilde{\Omega}_{\mathrm{int},R}$ or $\widetilde{\Pi}_{\mathrm{int},R}$. Therefore we do not have to consider the corresponding topological issue and the stability issue further in this certain special situation.	
\end{proposition}

\begin{proof}
See \cite[Lemma 4.5.4]{10KL2}.	
\end{proof}

\indent The following result is then achievable:

\begin{proposition}\mbox{\bf{(After Kedlaya-Liu \cite[Theorem 4.5.7]{10KL2})}}
Consider the following categories. The first one is the corresponding category of all the finite projective or pseudocoherent $\underline{\mathcal{O}_{E_k}}$-local systems. The second one is the corresponding category of all the finite projective or stably-pseudocoherent $F^k$-$\widetilde{\Omega}_{\mathrm{int},R}$-modules. The third one is the corresponding category of all the finite projective or stably-pseudocoherent $F^k$-$\widetilde{\Pi}_{\mathrm{int},R}$-modules. The fourth one is the corresponding category of all the finite projective or stably-pseudocoherent $F^k$-$\widetilde{\Omega}_{\mathrm{int},*}$-sheaves. The last one is the corresponding category of all the finite projective or stably-pseudocoherent $F^k$-$\widetilde{\Pi}_{\mathrm{int},*}$-sheaves. Here we consider the corresponding analytic topology. \\
\indent Then we have that these categories are actually equivalent. Similar parallel statement could then be made to fpd objects.	
\end{proposition}

\begin{remark}
The corresponding stably-pseudocoherent indication in the previous proposition is actually not serious as in \cite[Lemma 4.5.7]{10KL2}, the reason is that by the previous proposition, everything is already complete with respect to the natural topology.	
\end{remark}

\begin{proof}
This is just the corresponding analytic analog of the corresponding result in \cite[Lemma 4.5.7]{10KL2}. We just briefly mention the corresponding functors which realize the corresponding equivalence. The corresponding one from the first category to the corresponding sheaves is naturally just the corresponding base change. The corresponding functor from the sheaves to the corresponding modules is naturally the corresponding global section functor after the corresponding foundation in \cite[Theorem 1.4.2, Theorem 1.4.18]{10Ked1}. 
\end{proof}

\begin{remark}
Since we touched the corresponding foundations in \cite[Theorem 1.4.2, Theorem 1.4.18]{10Ked1} on the corresponding glueing of the vector bundles and the corresponding stably-pseudocoherent modules in the situation of analytic topology. However the corresponding descent in the \'etale, pro-\'etale and $v$-topology situations are actually parallel to the corresponding \cite[Theorem 2.5.5, Theorem 2.5.14, Theorem 3.4.8, Theorem 3.5.8]{10KL2}, which will basically upgrade naturally this proposition to the \'etale, pro-\'etale and $v$-topology situations, although we have not established the chance to present. 
\end{remark}

\indent Now carrying the corresponding coefficient $A$ we could have at least the corresponding fully faithfulness.

\begin{proposition}\mbox{\bf{(After Kedlaya-Liu \cite[Theorem 4.5.7]{10KL2})}}
Consider the following categories. The first one is the corresponding category of all the finite projective or pseudocoherent $\underline{\mathcal{O}_{{E_k},A}}$-local systems. The second one is the corresponding category of all the finite projective or stably-pseudocoherent $F^k$-$\widetilde{\Omega}_{\mathrm{int},R,A}$-modules. The third one is the corresponding category of all the finite projective or stably-pseudocoherent $F^k$-$\widetilde{\Pi}_{\mathrm{int},R,A}$-modules. The fourth one is the corresponding category of all the finite projective or stably-pseudocoherent $F^k$-$\widetilde{\Omega}_{\mathrm{int},*,A}$-sheaves. The last one is the corresponding category of all the finite projective or stably-pseudocoherent $F^k$-$\widetilde{\Pi}_{\mathrm{int},*,A}$-sheaves. Here we consider the corresponding analytic topology. \\
\indent Then we have that the first category could be embedded fully faithfully into the rest four categories and the corresponding modules could be compared to the sheaves under equivalences of the corresponding categories as long as we assume the corresponding sousperfectoidness of the ring $A$ which preserves the corresponding exactness of the exact sequences:
\[
\xymatrix@C+0pc@R+0pc{
0 \ar[r] \ar[r] \ar[r] &\underline{\mathcal{O}_{E_k}}  \ar[r] \ar[r] \ar[r] &\widetilde{\Omega}_{\mathrm{int},R} \ar[r] \ar[r] \ar[r] &\widetilde{\Omega}_{\mathrm{int},R} \ar[r] \ar[r] \ar[r] &0,\\
0 \ar[r] \ar[r] \ar[r] &\underline{{E}_{k}}  \ar[r] \ar[r] \ar[r] &\widetilde{\Omega}_{R} \ar[r] \ar[r] \ar[r] &\widetilde{\Omega}_{R} \ar[r] \ar[r] \ar[r] &0,\\
0 \ar[r] \ar[r] \ar[r] &\underline{\mathcal{O}_{E_k}}  \ar[r] \ar[r] \ar[r] &\widetilde{\Pi}_{\mathrm{int},R} \ar[r] \ar[r] \ar[r] &\widetilde{\Pi}_{\mathrm{int},R} \ar[r] \ar[r] \ar[r] &0,\\
0 \ar[r] \ar[r] \ar[r] &\underline{{E_k}}  \ar[r] \ar[r] \ar[r] &\widetilde{\Pi}_{R} \ar[r] \ar[r] \ar[r] &\widetilde{\Pi}_{R} \ar[r] \ar[r] \ar[r] &0.\\
}
\]	
Similar statement could be mede to the corresponding fpd objects.	
\end{proposition}

\begin{proof}
See \cite[Theorem 4.5.7]{10KL2}.	
\end{proof}

\

\newpage\section{Comparison in the Non-\'Etale Setting}

\indent We now consider the corresponding analog of the corresponding \cite[Chapter 4.6]{10KL2} over $\mathrm{Spa}(R,R^+)$ as above which is just assumed to be analytic.

%%%%%%%%%%%%%%%%%%%%%%%%%%%%%%%%%!!!!!!!!!!!!!!!

\begin{remark}
As we mentioned before, the corresponding analyticity of the ring $R$ will still not deduce from the analyticity of the corresponding rational Robba rings and sheaves.	
\end{remark}

\begin{theorem}\mbox{\bf{(After Kedlaya-Liu \cite[Theorem 4.6.1]{10KL2})}}
Assume the corresponding Robba rings with respect to closed intervals are sheafy. Consider the corresponding categories in the following:\\
A. The corresponding category of all the corresponding \'etale-stably-pseudocoherent sheaves over the adic Fargues-Fontaine curve (associated to $\{\widetilde{\Pi}_{I,R,A}\}_{\{I\subset (0,\infty)\}}$) in the corresponding \'etale topology;\\
B. The corresponding category of all the corresponding \'etale-stably-pseudocoherent $F^k$-sheaves over families of Robba rings $\{\widetilde{\Pi}_{I,R,A}\}_{\{I\subset (0,\infty)\}}$;\\
C. The corresponding category of all the corresponding \'etale-stably-pseudocoherent modules over any Robba rings $\widetilde{\Pi}_{I,R,A}$ such that the interval $I=[r_1,r_2]$ satisfies that $0\leq r_1 \leq r_2/p^{hk}$;\\
D. The corresponding category of all the corresponding pseudocoherent modules over any Robba rings $\widetilde{\Pi}_{R,A}$, but admitting models of strictly-pseudocoherent $F^k$-modules over some $\widetilde{\Pi}_{r,R,A}$ for some radius $r>0$ whose base changes to some $\widetilde{\Pi}_{[r_1,r_2],R,A}$ will produce corresponding $F^k$-modules which are \'etale-stably-pseudocoherent;\\
E. The corresponding category of all the corresponding strictly-pseudocoherent $F^k$-modules over some $\widetilde{\Pi}_{\infty,R,A}$ whose base changes to some $\widetilde{\Pi}_{[r_1,r_2],R,A}$ will provide corresponding $F^k$-modules which are \'etale-stably-pseudocoherent;\\
F. The corresponding category of all the corresponding \'etale-stably-pseudocoherent sheaves over the adic Fargues-Fontaine curve (associated to $\{\widetilde{\Pi}_{I,R,A}\}_{\{I\subset (0,\infty)\}}$) in the corresponding pro-\'etale topology.\\
Then we have that they are actually equivalent.

\end{theorem}

\begin{proof}
This is just parallel \cite[Theorem 4.6.1]{10KL2}, see \cite[Theorem 4.11]{10XT2}.	
\end{proof}

\begin{theorem}\mbox{\bf{(After Kedlaya-Liu \cite[Theorem 4.6.1]{10KL2})}}
Assume $A$ is sousperfectoid. Consider the corresponding categories of sheaves over $\mathrm{Spa}(R,R^+)_{\text{p\'et}}$ in the following:\\
%A. The corresponding category of all the corresponding stably-pseudocoherent sheaves over the adic Fargues-Fontaine curve (associated to $\{\widetilde{\Pi}_{I,R,A}\}_{\{I\subset (0,\infty)\}}$) in the corresponding \'etale topology;\\
%A. The corresponding category of all the corresponding \'etale-stably-pseudocoherent $F^k$-sheaves over families of Robba rings $\{\widetilde{\Pi}_{I,*,A}\}_{\{I\subset (0,\infty)\}}$;\\
A. The corresponding category of all the corresponding \'etale-stably-pseudocoherent modules over any Robba rings $\widetilde{\Pi}_{I,*,A}$ such that the interval $I=[r_1,r_2]$ satisfies that $0\leq r_1 \leq r_2/p^{hk}$;\\
B. The corresponding category of all the corresponding pseudocoherent modules over any Robba rings $\widetilde{\Pi}_{*,A}$, but admit models of strictly-pseudocoherent $F^k$-modules over some $\widetilde{\Pi}_{r,*,A}$ for some radius $r>0$ whose base changes to some $\widetilde{\Pi}_{[r_1,r_2],R,A}$ will provide corresponding $F^k$-modules which are \'etale-stably-pseudocoherent;\\
C. The corresponding category of all the corresponding \'etale-strictly-pseudocoherent $F^k$-modules over some $\widetilde{\Pi}_{\infty,*,A}$ whose base changes to some $\widetilde{\Pi}_{[r_1,r_2],*,A}$ will provide corresponding $F^k$-modules which are \'etale-stably-pseudocoherent.\\
%F. The corresponding category of all the corresponding stably-pseudocoherent sheaves over the adic Fargues-Fontaine curve (associated to $\{\widetilde{\Pi}_{I,R,A}\}_{\{I\subset (0,\infty)\}}$) in the corresponding pro-\'etale topology.\\
Then we have that they are actually equivalent.

\end{theorem}

\begin{proof}
This is parallel to \cite[Theorem 4.6.1]{10KL2}The proof is the same as the one where $*$ is just a ring $R$. See \cite[Theorem 4.11]{10XT2}.	
\end{proof}

\begin{theorem}\mbox{\bf{(After Kedlaya-Liu \cite[Corollary 4.6.2]{10KL2})}}
The corresponding analogs of the previous two theorems hold in the corresponding finite projective setting. Here we assume that $A$ is in the same hypotheses as above.	
\end{theorem}

\begin{proof}
See \cite[Corollary 4.6.2]{10KL2}.	
\end{proof}

%%%%%%%%%%%%%%%%%%%%%%%%%%%%!!!!!!!!!!!!!!!!!!!

\indent Now we drop the corresponding condition on the sheafiness on $A$ by using the corresponding derived spectrum $\mathrm{Spa}^h(A)$ from \cite{10BK1} where we have the corresponding $\infty$-sheaf $\mathcal{O}_{\mathrm{Spa}^h(A)}$, we now apply the corresponding construction to the Robba rings (in the rational setting over $\mathbb{Q}_p$) $*_{R,A}$ and denote the corresponding derived version by $*^h_{R,A}$. Here $A$ is assumed to be just Banach over $\mathbb{Q}_p$ or $\mathbb{F}_p((t))$, which is certainly in a very general situation.

\begin{remark}
One can apply the corresponding results of \cite{10CS} to achieve so as well, which we believe will be robust as well.	
\end{remark}

%%%%%%%%%%%%%%%%%%%%%%%%%%%%%!!!!!!!!!!!!!!!!

\begin{theorem}\mbox{\bf{(After Kedlaya-Liu \cite[Theorem 4.6.1]{10KL2})}}
Take a derived rational localization of $\widetilde{\Pi}_{I_1\bigcup I_2,R,A}$, which we denote it by $\widetilde{\Pi}^{h,*}_{I_1\bigcup I_2,R,A}$, where we choose two overlapped closed intervals $I_1$ and $I_2$. Now take the corresponding base changes of this localization along:
\begin{align}
\widetilde{\Pi}_{I_1\bigcup I_2,R,A}\rightarrow \widetilde{\Pi}_{I_1,R,A},\\
\widetilde{\Pi}_{I_1\bigcup I_2,R,A}\rightarrow \widetilde{\Pi}_{I_2,R,A},
\end{align}
which we will denote by $\widetilde{\Pi}^{h,*}_{I_1,R,A}$ and $\widetilde{\Pi}^{h,*}_{I_2,R,A}$. Consider the corresponding categories in the following:\\
A. The corresponding product of the category of all the corresponding finite projective $F^k$-modules over the Robba ring $\widetilde{\Pi}^{h,*}_{I_1,R,A}$ and the the category of all the corresponding finite projective $F^k$-modules over the Robba ring $\widetilde{\Pi}^{h,*}_{I_2,R,A}$, over the category of all the corresponding finite projective $F^k$-modules over $\widetilde{\Pi}^{h,*}_{I_1,R,A}\widehat{\otimes}^\mathbb{L}\widetilde{\Pi}^{h,*}_{I_2,R,A}$;\\
B. The corresponding category of all the corresponding finite projective modules over the Robba ring $\widetilde{\Pi}^{h,*}_{I_1\bigcup I_2,R,A}$.\\
Then we have that they are actually equivalent.

\end{theorem}

\begin{proof}
We consider a pair of overlapped intervals $I_1=[r_1,r_2],I_2=[s_1,s_2]$, and we have the corresponding $\infty$-Robba rings which form the desired situation, where we have the sequence which is basically exact up to higher homotopy:
\[
\xymatrix@C+2pc@R+3pc{
\widetilde{\Pi}^{h,*}_{I_1\bigcup I_2,R,A} \ar[r]\ar[r] \ar[r] \ar[r]  &\widetilde{\Pi}^{h,*}_{I_1,R,A}\bigoplus \widetilde{\Pi}^{h,*}_{I_2,R,A} \ar[r] \ar[r] \ar[r]\ar[r] &\widetilde{\Pi}^{h,*}_{I_1,R,A}\widehat{\otimes}^\mathbb{L}\widetilde{\Pi}^{h,*}_{I_2,R,A}.  
}
\] 
Now we consider the corresponding two finite projective module spectra (which admit retracts from finite free module spectra as in \cite[Proposition 7.2.2.7]{10Lu1}) $M_1,M_2,M_{12}$ over the rings in the middle and the rightmost positions, and assume that they form a corresponding glueing datum. Then we consider the corresponding presentation diagram for the module spectra:
\[
\xymatrix@C+2pc@R+3pc{
 &F \ar[r] \ar[r] \ar[r] \ar[d] \ar[d] \ar[d] &F_1\bigoplus F_2 \ar[d] \ar[d] \ar[d] \ar[r] \ar[r] \ar[r] &F_{12} \ar[d] \ar[d] \ar[d] \\
&G \ar[r] \ar[r] \ar[r] \ar[d] \ar[d] \ar[d] &G_1\bigoplus G_2 \ar[d]^{r_1\oplus r_2} \ar[d] \ar[d] \ar[r] \ar[r] \ar[r] &G_{12} \ar[d]^{r_{12}} \ar[d] \ar[d] \\
&M \ar[r] \ar[r] \ar[r] &M_1\bigoplus M_2 \ar[r] \ar[r] \ar[r] &M_{12} \\
}
\]  
by just taking the corresponding fibers, here $G,G_1,G_2,G_{12}$ are the finite free modules and we have the corresponding retracts $r_1,r_2,r_{12}$, but note that the corresponding map $r_1\oplus r_2$ might not be a priori a retract. However consider the following commutative diagram:
\[
\xymatrix@C+2pc@R+3pc{
& &0 \ar[d] \ar[d] \ar[d] &0 \ar[d] \ar[d] \ar[d] \\
 &\pi_0 F \ar[r] \ar[r] \ar[r] \ar[d] \ar[d] \ar[d] &\pi_0F_1\bigoplus \pi_0F_2 \ar[d] \ar[d] \ar[d] \ar[r] \ar[r] \ar[r] &\pi_0F_{12} \ar[d] \ar[d] \ar[d] \\
&\pi_0G \ar[r] \ar[r] \ar[r] \ar[d] \ar[d] \ar[d] &\pi_0G_1\bigoplus \pi_0G_2 \ar[d]^{\pi_0(r_1)\oplus \pi_0(r_2)} \ar[d] \ar[d] \ar[r] \ar[r] \ar[r] &\pi_0G_{12} \ar[d]^{r_{12}} \ar[d] \ar[d] \ar[r] \ar[r] \ar[r]&0 \\
&\pi_0M \ar[r] \ar[r] \ar[r] &\pi_0M_1\bigoplus \pi_0M_2\ar[d] \ar[d] \ar[d]  \ar[r] \ar[r] \ar[r] &\pi_0M_{12}\ar[d] \ar[d] \ar[d]  \ar[r] \ar[r] \ar[r]&0 \\
&&0&0.
}
\] 
We then consider the corresponding the argument of Kedlaya as in \cite[Proposition 5.11]{10XT3} which gives us modification on sections $\pi_0s_1$ and $\pi_0 s_2$ of $\pi_0 r_1$ and $\pi_0 r_2$ respectively such that we have (with the same notations) the new $\pi_0 s_1$ and $\pi_0 s_2$ give the corresponding new lifts $s_1$ and $s_2$ whose coproduct actually restricts to $M$ which produce a section for $M$ for the map $G\rightarrow M$, which proves the desired finite projectivity on $M$. From $A$ to $B$ we just take the corresponding projection while going back we consider the corresponding binary glueing to tackle the corresponding interval which cannot be reached by taking Frobenius translations but could be covered by two ones reachable by taking Frobenius translations. The lifts could be also modified by directly on the derived level.\\
\end{proof}

\indent In \cite{10GV}, Galatius and Venkatesh considered some derived Galois deformation theory. Now we make some discussion around this point by using the corresponding simplicial Banach rings (over $\mathbb{Q}_p$ or over $\mathbb{F}_p((t))$) in \cite{10BK1}, namely we use the notation $A^h$ to denote any local charts \footnote{Even one can take the more general simplicial Banach rings as in \cite{10BBBK}, especially one would like to focus on the corresponding analytification of derived Galois deformation rings as in \cite{10GV}, at least our feeling is that \cite{10BBBK} will allow one to take the corresponding derived adic generic fiber to produce some desired simplicial Banach rings to tackle derived Galois deformation problems in \cite{10GV}.} coming from the corresponding Bambozzi-Kremnizer spectrum in \cite{10BK1} attached to some Banach algebra over (over $\mathbb{Q}_p$ or over $\mathbb{F}_p((t))$). Over the analytic field specified we take the corresponding completed tensor products the period rings involved with $A^h$, which we will denote by $*_{R,A^h}$. Then we have the following similar results as above:

\begin{conjecture}\mbox{\bf{(After Kedlaya-Liu \cite[Theorem 4.6.1]{10KL2})}}
Consider the corresponding categories in the following:\\
A. The corresponding category of all the corresponding locally finite free sheaves over the adic Fargues-Fontaine curve (associated to $\{\widetilde{\Pi}_{I,R,A^h}\}_{\{I\subset (0,\infty)\}}$) in the corresponding homotopy Zariski topology from \cite{10BK1} and \cite{10BBBK} induced by Koszul derived rational localization;\\
B. The corresponding category of all the corresponding finite projective $F^k$-sheaves over families of Robba rings $\{\widetilde{\Pi}_{I,R,A^h}\}_{\{I\subset (0,\infty)\}}$;\\
C. The corresponding category of all the corresponding finite projective modules over any Robba rings $\widetilde{\Pi}_{I,R,A^h}$ such that the interval $I=[r_1,r_2]$ satisfies that $0\leq r_1 \leq r_2/p^{hk}$.\\
%D. The corresponding category of all the corresponding finite projective modules over any Robba rings $\widetilde{\Pi}^h_{R,A}$, but admit models of finite projective $F^k$-modules over some $\widetilde{\Pi}^h_{r,R,A}$ for some radius $r>0$ whose base changes to some $\widetilde{\Pi}^h_{[r_1,r_2],R,A}$ will provide a corresponding $F^k$-modules which is finite projective;\\
%E. The corresponding category of all the corresponding finite projective $F^k$-modules over some $\widetilde{\Pi}^h_{\infty,R,A}$ whose base changes to some $\widetilde{\Pi}^h_{[r_1,r_2],R,A}$ will provide a corresponding $F^k$-modules which is finite projective.\\
Then we have that they are actually equivalent.

\end{conjecture}

\begin{remark}
This conjecture and the one we made in the same fashion later namely \cref{conjecture6.9} could be formulated by using Clausen-Scholze analytic spaces in \cite{10CS}, which is in a different geometric setting, but we will believe the picture is quite similar. In fact once one replaces the geometric context we considered in order to use the context of \cite{10CS}, then this will be a direct consequence of the corresponding general descent result (\cite[Theorem 14.9, Remark 14.10]{10CS}) in the corresponding nonderived (but possibly highly nonsheafy) situation, which includes the corresponding descending of the corresponding general perfect complex as well. 	
\end{remark}

\begin{remark}
One can apply the corresponding results of \cite{10CS} to achieve so as well, which we believe will be robust as well.	
\end{remark}

\begin{theorem}\mbox{\bf{(After Kedlaya-Liu \cite[Theorem 4.6.1]{10KL2})}} \label{10theorem10.5.10}
Consider the corresponding categories in the following:\\
%A. The corresponding category of all the corresponding locally finite free sheaves over the adic Fargues-Fontaine curve (associated to $\{\widetilde{\Pi}_{I,S,A^h}\}_{\{I\subset (0,\infty)\}}$) in the corresponding homotopy Zariski topology from \cite{10BK1} and \cite{10BBBK} induced by Koszul derived rational localization;\\
A. The corresponding category of all the corresponding finite projective $F^k$-sheaves over families of Robba rings $\{\widetilde{\Pi}_{I,R,A^h}\}_{\{I\subset (0,\infty)\}}$;\\
B. The corresponding category of all the corresponding finite projective modules over any Robba rings $\widetilde{\Pi}_{I,R,A^h}$ such that the interval $I=[r_1,r_2]$ satisfies that $0\leq r_1 \leq r_2/p^{hk}$.\\
%D. The corresponding category of all the corresponding finite projective modules over any Robba rings $\widetilde{\Pi}^h_{R,A}$, but admit models of finite projective $F^k$-modules over some $\widetilde{\Pi}^h_{r,R,A}$ for some radius $r>0$ whose base changes to some $\widetilde{\Pi}^h_{[r_1,r_2],R,A}$ will provide a corresponding $F^k$-modules which is finite projective;\\
%E. The corresponding category of all the corresponding finite projective $F^k$-modules over some $\widetilde{\Pi}^h_{\infty,R,A}$ whose base changes to some $\widetilde{\Pi}^h_{[r_1,r_2],R,A}$ will provide a corresponding $F^k$-modules which is finite projective.\\
Then we have that they are actually equivalent.

\end{theorem}

\begin{proof}
Note that in our situation we do have the following nice short exact sequences:
\[
\xymatrix@C+0pc@R+0pc{
0\ar[r] \ar[r] \ar[r] &\pi_k \widetilde{\Pi}_{I_1\bigcup I_2,R,A^h} \ar[r] \ar[r] \ar[r] &\pi_k\widetilde{\Pi}_{I_1,R,A^h}\bigoplus \pi_k\widetilde{\Pi}_{I_2,R,A^h} \ar[r] \ar[r] \ar[r] &\pi_k\widetilde{\Pi}_{I_1\bigcap I_2,R,A^h} \ar[r] \ar[r] \ar[r] &0, k=0,1,....
}
\]
We consider a pair of overlapped intervals $I_1=[r_1,r_2],I_2=[s_1,s_2]$, and we have the corresponding $\infty$-Robba rings which form a corresponding glueing sequence in the sense \cite[Definition 2.7.3]{10KL1}:
\[
\xymatrix@C+0pc@R+0pc{
0\ar[r] \ar[r] \ar[r] &\pi_0\widetilde{\Pi}_{I_1\bigcup I_2,R,A^h} \ar[r] \ar[r] \ar[r] &\pi_0\widetilde{\Pi}_{I_1,R,A^h}\bigoplus \pi_0\widetilde{\Pi}_{I_2,R,A^h} \ar[r] \ar[r] \ar[r] &\pi_0\widetilde{\Pi}_{I_1\bigcap I_2,R,A^h} \ar[r] \ar[r] \ar[r] &0.
}
\] 	

Now we consider the corresponding two finite projective module spectra (which admit retracts from finite free module spectra as in \cite[Proposition 7.2.2.7]{10Lu1}) $M_1,M_2,M_{12}$ over the rings in the middle and the rightmost position, and assume that they form a corresponding glueing datum. Then we consider the corresponding presentation diagram for the module spectra:
\[
\xymatrix@C+2pc@R+3pc{
 &F \ar[r] \ar[r] \ar[r] \ar[d] \ar[d] \ar[d] &F_1\bigoplus F_2 \ar[d] \ar[d] \ar[d] \ar[r] \ar[r] \ar[r] &F_{12} \ar[d] \ar[d] \ar[d] \\
&G \ar[r] \ar[r] \ar[r] \ar[d] \ar[d] \ar[d] &G_1\bigoplus G_2 \ar[d]^{r_1\oplus r_2} \ar[d] \ar[d] \ar[r] \ar[r] \ar[r] &G_{12} \ar[d]^{r_{12}} \ar[d] \ar[d] \\
&M \ar[r] \ar[r] \ar[r] &M_1\bigoplus M_2 \ar[r] \ar[r] \ar[r] &M_{12} \\
}
\]  
by just taking the corresponding fibers, here $G,G_1,G_2,G_{12}$ are the finite free modules and we have the corresponding retracts $r_1,r_2,r_{12}$, but note that the corresponding map $r_1\oplus r_2$ might not be a priori a retract. However consider the following commutative diagram (which admits more exactness than in the previous theorem):
\[
\xymatrix@C+2pc@R+3pc{
&& &0 \ar[d] \ar[d] \ar[d] &0 \ar[d] \ar[d] \ar[d] \\
& &\pi_0 F \ar[r] \ar[r] \ar[r] \ar[d] \ar[d] \ar[d] &\pi_0F_1\bigoplus \pi_0F_2 \ar[d] \ar[d] \ar[d] \ar[r] \ar[r] \ar[r] &\pi_0F_{12} \ar[d] \ar[d] \ar[d] \\
&0 \ar[r] \ar[r] \ar[r] &\pi_0G \ar[r] \ar[r] \ar[r] \ar[d] \ar[d] \ar[d] &\pi_0G_1\bigoplus \pi_0G_2 \ar[d]^{\pi_0(r_1)\oplus \pi_0(r_2)} \ar[d] \ar[d] \ar[r] \ar[r] \ar[r] &\pi_0G_{12} \ar[d]^{r_{12}} \ar[d] \ar[d] \ar[r] \ar[r] \ar[r]&0 \\
&0 \ar[r] \ar[r] \ar[r] &\pi_0M \ar[r] \ar[r] \ar[r] &\pi_0M_1\bigoplus \pi_0M_2\ar[d] \ar[d] \ar[d]  \ar[r] \ar[r] \ar[r] &\pi_0M_{12}\ar[d] \ar[d] \ar[d]  \ar[r] \ar[r] \ar[r]&0 \\
&&&0&0.
}
\] 
We then consider the corresponding the argument of Kedlaya as in \cite[Proposition 5.11]{10XT3} which gives us modification on sections $\pi_0s_1$ and $\pi_0 s_2$ of $\pi_0 r_1$ and $\pi_0 r_2$ respectively such that we have (with the same notations) the new $\pi_0 s_1$ and $\pi_0 s_2$ give the corresponding new lifts $s_1$ and $s_2$ whose coproduct actually restricts to $M$ which produce a section for $M$ for the map $G\rightarrow M$, which proves the desired finite projectivity on $M$. From $A$ to $B$ we just take the corresponding projection while going back we consider the corresponding binary glueing to tackle the corresponding interval which cannot be reached by taking Frobenius translations but could be covered by two ones reachable by taking Frobenius translations. The lifts could be also modified by directly on the derived level. Or one can just show the flatness on the derived level directly by using derived Tor, which will directly proves that $M$ is finite projective as a module spectrum.
\end{proof}

\indent Now we even drop the hypothesis (by using the corresponding techniques in \cite{10XT3} and \cite{10XT4}) on the commutativity on $A$, namely now $A$ could be allowed to be strictly quotient coming from the free Tate algebras $\mathbb{Q}_p\left<Z_1,...,Z_d\right>$ and $\mathbb{F}_p((t))\left<Z_1,...,Z_d\right>$, which we will use some different notation $B$ to denote this. Note certainly that $A$ is some special case of such $B$, therefore now we are going to indeed generalize the discussion above to the noncommutative setting.

\begin{remark}
One can definitely consider more general coefficients in the Banach setting, but for Hodge-Iwasawa theory, we prefer to locate our discussion in the area where we could get the interesting rings by taking admissible and strictly quotient from free Tate algebras over analytic field, which certainly includes the corresponding situation of rigid analytic affinoids.	
\end{remark}

\begin{setting}
By the work \cite{10XT3} and \cite{10XT4} we can freely translate between the language of $B$-stably-pseudocoherent sheaves over analytic or \'etale topology and \'etale setting and the corresponding $B$-stably pseudocoherent modules or $B$-\'etale stably pseudocoherent modules, which makes the discussion in our context possible. What was happening indeed when we are in the corresponding sheafy tensor product situation is that one can direct read off such result in the next theorem by using the descent for the pseudocoherent modules with desired stability in \cite{10KL2} and \cite{10Ked1}.  	
\end{setting}

\begin{theorem}\mbox{\bf{(After Kedlaya-Liu \cite[Theorem 4.6.1]{10KL2})}}
Assume the corresponding Robba rings with respect to closed intervals are sheafy. Consider the corresponding categories in the following (all modules are left):\\
A. The corresponding category of all the corresponding $B$-\'etale-stably-pseudocoherent \\ 
sheaves over the adic Fargues-Fontaine curve (associated to $\{\widetilde{\Pi}_{I,R,B}\}_{\{I\subset (0,\infty)\}}$) in the corresponding \'etale topology;\\
B. The corresponding category of all the corresponding $B$-\'etale-stably-pseudocoherent $F^k$-sheaves over families of Robba rings $\{\widetilde{\Pi}_{I,R,B}\}_{\{I\subset (0,\infty)\}}$;\\
C. The corresponding category of all the corresponding $B$-\'etale-stably-pseudocoherent modules over any Robba rings $\widetilde{\Pi}_{I,R,B}$ such that the interval $I=[r_1,r_2]$ satisfies that $0\leq r_1 \leq r_2/p^{hk}$;\\
D. The corresponding category of all the corresponding $B$-pseudocoherent modules over any Robba rings $\widetilde{\Pi}_{R,B}$, but admit models of $B$-strictly-pseudocoherent $F^k$-modules over some $\widetilde{\Pi}_{r,R,B}$ for some radius $r>0$ whose base changes to some $\widetilde{\Pi}_{[r_1,r_2],R,B}$ will provide corresponding $F^k$-modules which are $B$-\'etale-stably-pseudocoherent;\\
E. The corresponding category of all the corresponding $B$-strictly-pseudocoherent $F^k$-modules over some $\widetilde{\Pi}_{\infty,R,B}$ whose base changes to some $\widetilde{\Pi}_{[r_1,r_2],R,B}$ will provide corresponding $F^k$-modules which are $B$-\'etale-stably-pseudocoherent.\\
%F. The corresponding category of all the corresponding $B$-stably-pseudocoherent sheaves over the adic Fargues-Fontaine curve (associated to $\{\widetilde{\Pi}_{I,R,A}\}_{\{I\subset (0,\infty)\}}$) in the corresponding pro-\'etale topology.\\
Then we have that they are actually equivalent.

\end{theorem}

%%%%%%%%%%%%%%%%%%%%%%%%%%%%%%%%!!!!!!!!!!!!!

\begin{proof}
This is just parallel \cite[Theorem 4.6.1]{10KL2}, see \cite[Theorem 4.11]{10XT2}. The corresponding modules could be regarded as sheaves over the deformed sites as in our previous work \cite{10XT3} and \cite{10XT4}.	
\end{proof}

\begin{theorem}\mbox{\bf{(After Kedlaya-Liu \cite[Theorem 4.6.1]{10KL2})}}
Consider the corresponding categories of sheaves over $\mathrm{Spa}(R,R^+)_{\text{\'et}}$ in the following (all the modules are left over the corresponding rings):\\
%A. The corresponding category of all the corresponding stably-pseudocoherent sheaves over the adic Fargues-Fontaine curve (associated to $\{\widetilde{\Pi}_{I,R,A}\}_{\{I\subset (0,\infty)\}}$) in the corresponding \'etale topology;\\
%A. The corresponding category of all the corresponding $B$-\'etale-stably-pseudocoherent $F^k$-sheaves over families of Robba rings $\{\widetilde{\Pi}_{I,*,B}\}_{\{I\subset (0,\infty)\}}$;\\
A. The corresponding category of all the corresponding $B$-\'etale-stably-pseudocoherent modules over any Robba rings $\widetilde{\Pi}_{I,*,B}$ such that the interval $I=[r_1,r_2]$ satisfies that $0\leq r_1 \leq r_2/p^{hk}$;\\
B. The corresponding category of all the corresponding $B$-pseudocoherent modules over any Robba rings $\widetilde{\Pi}_{*,B}$, but admit models of $B$-strictly-pseudocoherent $F^k$-modules over some $\widetilde{\Pi}_{r,*,B}$ for some radius $r>0$ whose base changes to some $\widetilde{\Pi}_{[r_1,r_2],R,B}$ will provide corresponding $F^k$-modules which are $B$-\'etale-stably-pseudocoherent;\\
C. The corresponding category of all the corresponding $B$-strictly-pseudocoherent $F^k$-modules over some $\widetilde{\Pi}_{\infty,*,B}$ whose base changes to some $\widetilde{\Pi}_{[r_1,r_2],*,B}$ will provide corresponding $F^k$-modules which are $B$-\'etale-stably-pseudocoherent.\\
%F. The corresponding category of all the corresponding stably-pseudocoherent sheaves over the adic Fargues-Fontaine curve (associated to $\{\widetilde{\Pi}_{I,R,A}\}_{\{I\subset (0,\infty)\}}$) in the corresponding pro-\'etale topology.\\
Then we have that they are actually equivalent.

\end{theorem}

\begin{proof}
This is parallel to \cite[Theorem 4.6.1]{10KL2}The proof is the same as the one where $*$ is just a ring $R$. See \cite[Theorem 4.11]{10XT2}.	
\end{proof}

\begin{theorem}\mbox{\bf{(After Kedlaya-Liu \cite[Corollary 4.6.2]{10KL2})}}
The corresponding analogs of the previous two theorems hold in the corresponding finite projective setting (one has to consider the bimodules). Here we assume that $B$ is in the same hypotheses as above.	
\end{theorem}

\begin{proof}
See \cite[Corollary 4.6.2]{10KL2} and \cite[Proposition 5.12]{10XT3}.	
\end{proof}

\

\newpage\section{Discussion on the Generality of Gabber-Ramero}

\subsection{General Period Rings}

\indent We can encode now the corresponding discussion in the previous section actually, but we choose to separately discuss the corresponding results in detail here. Certainly many results in \cite{10KL1} and \cite{10KL2}, and \cite{10XT1}, \cite{10XT2}, \cite{10XT3} and \cite{10XT4} rely on the corresponding topologically nilpotent units and systems of topologically nilpotents. Therefore we will discuss the corresponding admissible and reasonable generalization after \cite{10KL1}, \cite{10KL2}, \cite{10XT1}, \cite{10XT2}, \cite{10XT3}, \cite{10XT4} and \cite{10GR}.

\begin{setting}
Drop the condition on the analyticity on $(R,R^+)$ and we switch to the corresponding notation $(S,S^+)$.	
\end{setting}

%%%%%%%%%%%%%%%%%%%%%%%%!!!!!!!!!!!!!

After Kedlaya-Liu \cite[Definition 4.1.1]{10KL2}, we consider the following constructions. First we consider the corresponding Witt vectors coming from the corresponding adic ring $(S,S^+)$. First we consider the corresponding generalized Witt vectors with respect to $(S,S^+)$ with the corresponding coefficients in the Tate algebra with the general notation $W(S^+)[[S]]$. The general form of any element in such deformed ring could be written as $\sum_{i\geq 0,i_1\geq 0,...,i_d\geq 0}\pi^i[\overline{y}_i]X_1^{i_1}...X_d^{i_d}$. Then we take the corresponding completion with respect to the following norm for some radius $t>0$:
\begin{align}
\|.\|_{t,A}(\sum_{i\geq 0,i_1\geq 0,...,i_d\geq 0}\pi^i[\overline{y}_i]X_1^{i_1}...X_d^{i_d}):= \max_{i\geq 0,i_1\geq 0,...,i_d\geq 0}p^{-i}\|.\|_S(\overline{y}_i)	
\end{align}
which will give us the corresponding ring $\widetilde{\Pi}_{\mathrm{int},t,S,A}$ such that we could put furthermore that:
\begin{align}
\widetilde{\Pi}_{\mathrm{int},S,A}:=\bigcup_{t>0} \widetilde{\Pi}_{\mathrm{int},t,S,A}.	
\end{align}
Then as in \cite[Definition 4.1.1]{10KL2}, we now put the ring $\widetilde{\Pi}_{\mathrm{bd},t,S,A}:=\widetilde{\Pi}_{\mathrm{int},t,S,A}[1/\pi]$ and we set:
\begin{align}
\widetilde{\Pi}_{\mathrm{bd},S,A}:=\bigcup_{t>0} \widetilde{\Pi}_{\mathrm{bd},t,S,A}.	
\end{align}
The corresponding Robba rings with respect to some intervals and some radius could be defined in the same way as in \cite[Definition 4.1.1]{10KL2}. To be more precise we consider the completion of the corresponding ring $W(S^+)[[S]][1/\pi]$ with respect to the following norm for some $t>0$ where $t$ lives in some prescribed interval $I=[s,r]$: 
\begin{align}
\|.\|_{t,A}(\sum_{i,i_1\geq 0,...,i_d\geq 0}\pi^i[\overline{y}_i]X_1^{i_1}...X_d^{i_d}):= \max_{i\geq 0,i_1\geq 0,...,i_d\geq 0}p^{-i}\|.\|_S(\overline{y}_i).	
\end{align}
This process will produce the corresponding Robba rings with respect to  the given interval $I=[s,r]$. Now for particular sorts of intervals $(0,r]$ we will have the corresponding Robba ring $\widetilde{\Pi}_{r,S,A}$ and we will have the corresponding Robba ring $\widetilde{\Pi}_{\infty,S,A}$	if the corresponding interval is taken to be $(0,\infty)$. Then in our situation we could just take the corresponding union throughout all the radius $r>0$ to define the corresponding full Robba ring taking the notation of $\widetilde{\Pi}_{S,A}$.

\indent The corresponding Robba rings $\widetilde{\Pi}_{\mathrm{bd},S,A}$, $\widetilde{\Pi}_{S,A}$, $\widetilde{\Pi}_{I,S,A}$, $\widetilde{\Pi}_{r,S,A}$, $\widetilde{\Pi}_{\infty,S,A}$ are actually themselves Tate adic Banach rings. However in many further application the non-Tateness of the ring $S$ will cause some reason for us to do the corresponding modification.

\indent Then for any general affinoid algebra $A$ over the corresponding base analytic field, we just take the corresponding quotients of the corresponding rings defined in the previous definition over some Tate algebras in rigid analytic geometry, with the same notations though $A$ now is more general. Note that one can actually show that the definition does not depend on the corresponding choice of the corresponding presentations over $A$.

\indent Again in this situation more generally, the corresponding Robba rings $\widetilde{\Pi}_{\mathrm{bd},S,A}$, $\widetilde{\Pi}_{S,A}$, $\widetilde{\Pi}_{I,S,A}$, $\widetilde{\Pi}_{r,S,A}$, $\widetilde{\Pi}_{\infty,S,A}$ are actually themselves Tate adic Banach rings.

\begin{lemma} \mbox{\bf{(After Kedlaya-Liu \cite[Lemma 5.2.6]{10KL2})}}
For any two radii $0<r_1<r_2$ we have the corresponding equality:
\begin{align}
\widetilde{\Pi}_{\mathrm{int},r_2,S,\mathbb{Q}_p\{T_1,...,T_d\}}\bigcap \widetilde{\Pi}_{[r_1,r_2],S,\mathbb{Q}_p\{T_1,...,T_d\}}	=\widetilde{\Pi}_{\mathrm{int},r_1,S,\mathbb{Q}_p\{T_1,...,T_d\}}.
\end{align}

\end{lemma}

\begin{proof}
See \cite[Lemma 5.2.6]{10KL2} and \cite[Proposition 2.13]{10XT2}.	
\end{proof}

\begin{lemma} \mbox{\bf{(After Kedlaya-Liu \cite[Lemma 5.2.6]{10KL2})}}
For any two radii $0<r_1<r_2$ we have the corresponding equality:
\begin{align}
\widetilde{\Pi}_{\mathrm{int},r_2,S,\mathbb{F}_p((t))\{T_1,...,T_d\}}\bigcap \widetilde{\Pi}_{[r_1,r_2],S,\mathbb{F}_p((t))\{T_1,...,T_d\}}	=\widetilde{\Pi}_{\mathrm{int},r_1,S,\mathbb{F}_p((t))\{T_1,...,T_d\}}.
\end{align}

\end{lemma}

\begin{proof}
See \cite[Lemma 5.2.6]{10KL2} and \cite[Proposition 2.13]{10XT2}.	
\end{proof}

\begin{lemma} \mbox{\bf{(After Kedlaya-Liu \cite[Lemma 5.2.6]{10KL2})}}
For general affinoid $A$ as above (over $\mathbb{Q}_p$ or $\mathbb{F}_p((t))$) and for any two radii $0<r_1<r_2$ we have the corresponding equality:
\begin{align}
\widetilde{\Pi}_{\mathrm{int},r_2,S,A}\bigcap \widetilde{\Pi}_{[r_1,r_2],S,A}	=\widetilde{\Pi}_{\mathrm{int},r_1,S,A}.
\end{align}

\end{lemma}

\begin{proof}
See \cite[Lemma 5.2.6]{10KL2} and \cite[Proposition 2.14]{10XT2}.	
\end{proof}

\begin{lemma} \mbox{\bf{(After Kedlaya-Liu \cite[Lemma 5.2.10]{10KL2})}}
For any four radii $0<r_1<r_2<r_3<r_4$ we have the corresponding equality:
\begin{align}
\widetilde{\Pi}_{[r_1,r_3],S,\mathbb{Q}_p\{T_1,...,T_d\}}\bigcap \widetilde{\Pi}_{[r_2,r_4],S,\mathbb{Q}_p\{T_1,...,T_d\}}	=\widetilde{\Pi}_{[r_1,r_4],S,\mathbb{Q}_p\{T_1,...,T_d\}}.
\end{align}

\end{lemma}

\begin{proof}
See \cite[Lemma 5.2.10]{10KL2} and \cite[Proposition 2.16]{10XT2}.	
\end{proof}

\begin{lemma} \mbox{\bf{(After Kedlaya-Liu \cite[Lemma 5.2.10]{10KL2})}}
For any four radii $0<r_1<r_2<r_3<r_4$ we have the corresponding equality:
\begin{align}
\widetilde{\Pi}_{[r_1,r_3],S,\mathbb{F}_p((t))\{T_1,...,T_d\}}\bigcap \widetilde{\Pi}_{[r_2,r_4],S,\mathbb{F}_p((t))\{T_1,...,T_d\}}	=\widetilde{\Pi}_{[r_1,r_4],S,\mathbb{F}_p((t))\{T_1,...,T_d\}}.
\end{align}

\end{lemma}

\begin{proof}
See \cite[Lemma 5.2.10]{10KL2} and \cite[Proposition 2.16]{10XT2}.	
\end{proof}

\begin{lemma} \mbox{\bf{(After Kedlaya-Liu \cite[Lemma 5.2.10]{10KL2})}}
For any four radii $0<r_1<r_2<r_3<r_4$ we have the corresponding equality:
\begin{align}
\widetilde{\Pi}_{[r_1,r_3],S,A}\bigcap \widetilde{\Pi}_{[r_2,r_4],S,A}	=\widetilde{\Pi}_{[r_1,r_4],S,A}.
\end{align}

\end{lemma}

\begin{proof}
See \cite[Lemma 5.2.10]{10KL2} and \cite[Proposition 2.17]{10XT2}.	
\end{proof}

\begin{definition}
All the Frobenius finite projective, pseudocoherent and fpd modules under $F^k$ could be defined in the exact same way as in \cref{10chapter3}. We will not repeat the corresponding definition, but note that we change the notation for $R$ to be just $S$.	
\end{definition}

\begin{conjecture}\mbox{\bf{(After Kedlaya-Liu \cite[Theorem 4.6.1]{10KL2})}} \label{conjecture6.9}
Consider the corresponding categories in the following:\\
A. The corresponding category of all the corresponding locally finite free sheaves over the adic Fargues-Fontaine curve (associated to $\{\widetilde{\Pi}_{I,S,A^h}\}_{\{I\subset (0,\infty)\}}$) in the corresponding homotopy Zariski topology from \cite{10BK1} and \cite{10BBBK} induced by Koszul derived rational localization;\\
B. The corresponding category of all the corresponding finite projective $F^k$-sheaves over families of Robba rings $\{\widetilde{\Pi}_{I,S,A^h}\}_{\{I\subset (0,\infty)\}}$;\\
C. The corresponding category of all the corresponding finite projective modules over any Robba rings $\widetilde{\Pi}_{I,S,A^h}$ such that the interval $I=[r_1,r_2]$ satisfies that $0\leq r_1 \leq r_2/p^{hk}$.\\
%D. The corresponding category of all the corresponding finite projective modules over any Robba rings $\widetilde{\Pi}^h_{R,A}$, but admit models of finite projective $F^k$-modules over some $\widetilde{\Pi}^h_{r,R,A}$ for some radius $r>0$ whose base changes to some $\widetilde{\Pi}^h_{[r_1,r_2],R,A}$ will provide a corresponding $F^k$-modules which is finite projective;\\
%E. The corresponding category of all the corresponding finite projective $F^k$-modules over some $\widetilde{\Pi}^h_{\infty,R,A}$ whose base changes to some $\widetilde{\Pi}^h_{[r_1,r_2],R,A}$ will provide a corresponding $F^k$-modules which is finite projective.\\
Then we have that they are actually equivalent.

\end{conjecture}

\begin{remark}
One can apply the corresponding results of \cite{10CS} to achieve so as well, which we believe will be robust as well.	
\end{remark}

\begin{theorem}\mbox{\bf{(After Kedlaya-Liu \cite[Theorem 4.6.1]{10KL2})}}
Consider the corresponding categories in the following:\\
%A. The corresponding category of all the corresponding locally finite free sheaves over the adic Fargues-Fontaine curve (associated to $\{\widetilde{\Pi}_{I,S,A^h}\}_{\{I\subset (0,\infty)\}}$) in the corresponding homotopy Zariski topology from \cite{10BK1} and \cite{10BBBK} induced by Koszul derived rational localization;\\
A. The corresponding category of all the corresponding finite projective $F^k$-sheaves over families of Robba rings $\{\widetilde{\Pi}_{I,S,A^h}\}_{\{I\subset (0,\infty)\}}$;\\
B. The corresponding category of all the corresponding finite projective modules over any Robba rings $\widetilde{\Pi}_{I,S,A^h}$ such that the interval $I=[r_1,r_2]$ satisfies that $0\leq r_1 \leq r_2/p^{hk}$.\\
%D. The corresponding category of all the corresponding finite projective modules over any Robba rings $\widetilde{\Pi}^h_{R,A}$, but admit models of finite projective $F^k$-modules over some $\widetilde{\Pi}^h_{r,R,A}$ for some radius $r>0$ whose base changes to some $\widetilde{\Pi}^h_{[r_1,r_2],R,A}$ will provide a corresponding $F^k$-modules which is finite projective;\\
%E. The corresponding category of all the corresponding finite projective $F^k$-modules over some $\widetilde{\Pi}^h_{\infty,R,A}$ whose base changes to some $\widetilde{\Pi}^h_{[r_1,r_2],R,A}$ will provide a corresponding $F^k$-modules which is finite projective.\\
Then we have that they are actually equivalent.

\end{theorem}

\begin{proof}
As in the situation in the analytic setting we consider the following argument which will show the desired result. In our situation we do have the following nice short exact sequences:
\[
\xymatrix@C+0pc@R+0pc{
0\ar[r] \ar[r] \ar[r] &\pi_k \widetilde{\Pi}_{I_1\bigcup I_2,S,A^h} \ar[r] \ar[r] \ar[r] &\pi_k\widetilde{\Pi}_{I_1,S,A^h}\bigoplus \pi_k\widetilde{\Pi}_{I_2,S,A^h} \ar[r] \ar[r] \ar[r] &\pi_k\widetilde{\Pi}_{I_1\bigcap I_2,S,A^h} \ar[r] \ar[r] \ar[r] &0, k=0,1,....
}
\]
We consider a pair of overlapped intervals $I_1=[r_1,r_2],I_2=[s_1,s_2]$, and we have the corresponding $\infty$-Robba rings which form a corresponding glueing sequence in the sense \cite[Definition 2.7.3]{10KL1}:
\[
\xymatrix@C+0pc@R+0pc{
0\ar[r] \ar[r] \ar[r] &\pi_0\widetilde{\Pi}_{I_1\bigcup I_2,S,A^h} \ar[r] \ar[r] \ar[r] &\pi_0\widetilde{\Pi}_{I_1,S,A^h}\bigoplus \pi_0\widetilde{\Pi}_{I_2,S,A^h} \ar[r] \ar[r] \ar[r] &\pi_0\widetilde{\Pi}_{I_1\bigcap I_2,S,A^h} \ar[r] \ar[r] \ar[r] &0.
}
\]

Now we consider the corresponding two finite projective module spectra (which admits retracts from finite free module spectra as in \cite[Proposition 7.2.2.7]{10Lu1}) $M_1,M_2,M_{12}$ over the rings in the middle and the rightmost position, and assume that they form a corresponding glueing datum. Then we consider the corresponding presentation diagram for the module spectra:
\[
\xymatrix@C+2pc@R+3pc{
 &L \ar[r] \ar[r] \ar[r] \ar[d] \ar[d] \ar[d] &L_1\bigoplus L_2 \ar[d] \ar[d] \ar[d] \ar[r] \ar[r] \ar[r] &L_{12} \ar[d] \ar[d] \ar[d] \\
&N \ar[r] \ar[r] \ar[r] \ar[d] \ar[d] \ar[d] &N_1\bigoplus N_2 \ar[d]^{r_1\oplus r_2} \ar[d] \ar[d] \ar[r] \ar[r] \ar[r] &N_{12} \ar[d]^{r_{12}} \ar[d] \ar[d] \\
&M \ar[r] \ar[r] \ar[r] &M_1\bigoplus M_2 \ar[r] \ar[r] \ar[r] &M_{12} \\
}
\]  
by just taking the corresponding fibers, here $N,N_1,N_2,N_{12}$ are the finite free modules and we have the corresponding retracts $r_1,r_2,r_{12}$, but note that the corresponding map $r_1\oplus r_2$ might not be a priori a retract. However consider the following commutative diagram:
\[
\xymatrix@C+2pc@R+3pc{
&& &0 \ar[d] \ar[d] \ar[d] &0 \ar[d] \ar[d] \ar[d] \\
& &\pi_0 F \ar[r] \ar[r] \ar[r] \ar[d] \ar[d] \ar[d] &\pi_0F_1\bigoplus \pi_0F_2 \ar[d] \ar[d] \ar[d] \ar[r] \ar[r] \ar[r] &\pi_0F_{12} \ar[d] \ar[d] \ar[d] \\
&0 \ar[r] \ar[r] \ar[r] &\pi_0G \ar[r] \ar[r] \ar[r] \ar[d] \ar[d] \ar[d] &\pi_0G_1\bigoplus \pi_0G_2 \ar[d]^{\pi_0(r_1)\oplus \pi_0(r_2)} \ar[d] \ar[d] \ar[r] \ar[r] \ar[r] &\pi_0G_{12} \ar[d]^{r_{12}} \ar[d] \ar[d] \ar[r] \ar[r] \ar[r]&0 \\
&0 \ar[r] \ar[r] \ar[r] &\pi_0M \ar[r] \ar[r] \ar[r] &\pi_0M_1\bigoplus \pi_0M_2\ar[d] \ar[d] \ar[d]  \ar[r] \ar[r] \ar[r] &\pi_0M_{12}\ar[d] \ar[d] \ar[d]  \ar[r] \ar[r] \ar[r]&0 \\
&&&0&0.
}
\] 
We then consider the corresponding the argument of Kedlaya as in \cite[Proposition 5.11]{10XT3} which gives us modification on sections $\pi_0s_1$ and $\pi_0 s_2$ of $\pi_0 r_1$ and $\pi_0 r_2$ respectively such that we have (with the same notations) the new $\pi_0 s_1$ and $\pi_0 s_2$ give the coresponding new lifts $s_1$ and $s_2$ whose coproduct actually restricts to $M$ which produce a section for $M$ for the map $N\rightarrow M$, which proves the desired finite projectivity on $M$. From $A$ to $B$ we just take the corresponding projection while going back we consider the corresponding binary glueing to tackle the corresponding interval which cannot be reached by taking Frobenius translations but could be covered by two ones reachable by taking Frobenius translations.
One can actually derive this on the derived level directly by taking the corresponding lift of the difference of the base changes of $s_1$ and $s_2$ to the right in $\mathrm{Ext}^0(M_{12},L_{12})$ to $\mathrm{Ext}^0(M_1\bigoplus M_2,L_1\bigoplus L_2)$ to modify the corresponding sections $s_1$ and $s_2$ in order to restrict to $M$, which proves the result by the argument of Kedlaya as in \cite[Proposition 5.11]{10XT3}.
\end{proof}

%%%%%%%%%%%%%%%%%%%%%%%%%%%%%%%%%%%%!!!!!!!!!!!!!!

\indent Now again we even drop the hypothesis (by using the corresponding techniques in \cite{10XT3} and \cite{10XT4}) on the commutativity on $A$, namely now $A$ could be allowed to be strictly quotient coming from the free Tate algebras $\mathbb{Q}_p\left<Z_1,...,Z_d\right>$ and $\mathbb{F}_p((t))\left<Z_1,...,Z_d\right>$, which we will use some different notation $B$ to denote this. Note certainly that $A$ is some special case of such $B$, therefore now we are going to indeed generalize the discussion above to the noncommutative setting.

\begin{theorem}\mbox{\bf{(After Kedlaya-Liu \cite[Theorem 4.6.1]{10KL2})}}
Assume the corresponding Robba rings with respect to closed intervals are sheafy. Consider the corresponding categories in the following (all modules are left):\\
A. The corresponding category of all the corresponding finite locally projective sheaves over the adic Fargues-Fontaine curve (associated to $\{\widetilde{\Pi}_{I,S,B}\}_{\{I\subset (0,\infty)\}}$) in the corresponding \'etale topology;\\
B. The corresponding category of all the corresponding finite projective $F^k$-sheaves over families of Robba rings $\{\widetilde{\Pi}_{I,S,B}\}_{\{I\subset (0,\infty)\}}$;\\
C. The corresponding category of all the corresponding finite projective modules over any Robba rings $\widetilde{\Pi}_{I,S,B}$ such that the interval $I=[r_1,r_2]$ satisfies that $0\leq r_1 \leq r_2/p^{hk}$.\\
%D. The corresponding category of all the corresponding finite projective modules over any Robba rings $\widetilde{\Pi}_{R,B}$, but admit models of $B$-strictly-pseudocoherent $F^k$-modules over some $\widetilde{\Pi}_{r,R,B}$ for some radius $r>0$ whose base changes to some $\widetilde{\Pi}_{[r_1,r_2],R,B}$ will provide a corresponding $F^k$-modules which is $B$-\'etale-stably-pseudocoherent;\\
%E. The corresponding category of all the corresponding $B$-strictly-pseudocoherent $F^k$-modules over some $\widetilde{\Pi}_{\infty,R,B}$ whose base changes to some $\widetilde{\Pi}_{[r_1,r_2],R,B}$ will provide a corresponding $F^k$-modules which is $B$-\'etale-stably-pseudocoherent.\\
%F. The corresponding category of all the corresponding $B$-stably-pseudocoherent sheaves over the adic Fargues-Fontaine curve (associated to $\{\widetilde{\Pi}_{I,R,A}\}_{\{I\subset (0,\infty)\}}$) in the corresponding pro-\'etale topology.\\
Then we have that they are actually equivalent.

\end{theorem}

\begin{proof}
This is just parallel \cite[Theorem 4.6.1]{10KL2}, see \cite[Theorem 4.11]{10XT2}. The corresponding modules could be regarded as sheaves over the deformed sites as in our previous work \cite{10XT3} and \cite{10XT4}.	
\end{proof}

\begin{theorem}\mbox{\bf{(After Kedlaya-Liu \cite[Theorem 4.6.1]{10KL2})}}
Consider the corresponding categories of sheaves over $\mathrm{Spa}(S,S^+)_{\text{\'et}}$ in the following (all the modules are bimodules over the corresponding rings):\\
%A. The corresponding category of all the corresponding locally finite projecive sheaves over the adic Fargues-Fontaine curve (associated to $\{\widetilde{\Pi}_{I,R,A}\}_{\{I\subset (0,\infty)\}}$) in the corresponding \'etale topology;\\
A. The corresponding category of all the corresponding finite projective $F^k$-sheaves over families of Robba rings $\{\widetilde{\Pi}_{I,*,B}\}_{\{I\subset (0,\infty)\}}$;\\
B. The corresponding category of all the corresponding finite projective modules over any Robba rings $\widetilde{\Pi}_{I,*,B}$ such that the interval $I=[r_1,r_2]$ satisfies that $0\leq r_1 \leq r_2/p^{hk}$;\\
%C. The corresponding category of all the corresponding $B$-pseudocoherent modules over any Robba rings $\widetilde{\Pi}_{*,B}$, but admit models of $B$-strictly-pseudocoherent $F^k$-modules over some $\widetilde{\Pi}_{r,*,B}$ for some radius $r>0$ whose base changes to some $\widetilde{\Pi}_{[r_1,r_2],R,B}$ will provide a corresponding $F^k$-modules which is $B$-\'etale-stably-pseudocoherent;\\
%D. The corresponding category of all the corresponding $B$-strictly-pseudocoherent $F^k$-modules over some $\widetilde{\Pi}_{\infty,*,B}$ whose base changes to some $\widetilde{\Pi}_{[r_1,r_2],*,B}$ will provide a corresponding $F^k$-modules which is $B$-\'etale-stably-pseudocoherent.\\
%F. The corresponding category of all the corresponding stably-pseudocoherent sheaves over the adic Fargues-Fontaine curve (associated to $\{\widetilde{\Pi}_{I,R,A}\}_{\{I\subset (0,\infty)\}}$) in the corresponding pro-\'etale topology.\\
Then we have that they are actually equivalent.

\end{theorem}

\begin{proof}
Note that in this situation the corresponding space is not even analytic adic space, but the Robba rings over some certain perfectoid domains are analytic. This is parallel to \cite[Theorem 4.6.1]{10KL2}The proof is the same as the one where $*$ is just a ring $R$. See \cite[Theorem 4.11]{10XT2}.	
\end{proof}

\subsection{Derived Algebraic Geometry of $\infty$-Period Rings}

\indent We now consider the contact with the corresponding algebraic sheaves over the corresponding schemes associated to the corresponding period rings in \cite{10KL1}, \cite{10KL2} as in \cite{10XT1}, \cite{10XT2}, \cite{10XT3}, \cite{10XT4} along \cite{10KL1}.

\begin{remark}
Certainly we are now dropping the corresponding topology and functional analyticity. This is motivated by the corresponding discussion on the algebraic quasi-coherent sheaves over schematic Fargues-Fontaine curves in \cite{10KL2}.
\end{remark}

\begin{setting}
$A$ is assumed to be now general Banach (commutative at this moment) over the base analytic field, and keep $R$ as in the situation before in this section (in the most general setting we considered so far). Now we apply the whole machinery in \cite{10Lu1} to the corresponding $\infty$-Bambozzi-Kremnizer rings $*^h_{R,A}$ attached to the period rings $*_{R,A}$ carrying $A$. Then we regard $*^h_{R,A}$ as a corresponding $\mathbb{E}_1$-ring in \cite{10Lu1} \footnote{Note that the original convention of \cite{10BK1} is cohomological and vanishing in positive degree.}.
\end{setting}

%%%%%%%%%%%%%%%%%%%%%%%!!!!!!!!!!!!!!!!!!!!!

\begin{proposition}
The corresponding category of all the perfect, almost perfect (after \cite[Section 7.2.4]{10Lu1}, namely pseudocoherent) $F^k$-equivariant sheaves over the families of $\infty$-Robba \'etale $\infty$-Deligne-Mumford toposes $(\mathrm{Spec}\widetilde{\Pi}^h_{I,R,A},\mathcal{O})_{\{I \subset (0,\infty)\}}$ (as in \cite[Chapter 1.4]{10Lu2}) is equivalent to the category of all the perfect, almost perfect (after \cite[Section 7.2.4]{10Lu1}, namely pseudocoherent) $F^k$-equivariant sheaves over some $\infty$-Robba \'etale $\infty$-Deligne-Mumford topos $(\mathrm{Spec}\widetilde{\Pi}^h_{I,R,A},\mathcal{O})$ (as in \cite[Chapter 1.4]{10Lu2}) where $I=[r_1,r_2]$ ($0<r_1\leq r_2/p^{hk}$). Here we need to assume the corresponding \cref{conjecture6.22} holds.	
\end{proposition}

\begin{proof}
This is actually quite transparent that we just take the corresponding Frobenius translation to compare the corresponding sheaves.	
\end{proof}

\begin{proposition}
The corresponding category of all the locally finite free (after \cite[Section 7.2.2, 7.2.4]{10Lu1}, namely pseudocoherent and flat) $F^k$-equivariant sheaves over the families of $\infty$-Robba \'etale $\infty$-Deligne-Mumford toposes $(\mathrm{Spec}\widetilde{\Pi}^h_{I,R,A},\mathcal{O})_{\{I \subset (0,\infty)\}}$ (as in \cite[Chapter 1.4]{10Lu2}) is equivalent to the category of all the locally finite free (after \cite[Section 7.2.2, 7.2.4]{10Lu1}, namely pseudocoherent and flat) $F^k$-equivariant sheaves over some $\infty$-Robba \'etale $\infty$-Deligne-Mumford topos $(\mathrm{Spec}\widetilde{\Pi}^h_{I,R,A},\mathcal{O})$ (as in \cite[Chapter 1.4]{10Lu2}) where $I=[r_1,r_2]$ ($0<r_1\leq r_2/p^{hk}$).	
\end{proposition}

\begin{proof}
This is again actually quite transparent that we just take the corresponding Frobenius translation to compare the corresponding sheaves.	
\end{proof}

\begin{example}
Everything will be certainly more interesting when we maintain in the corresponding noetherian situation, where we do have the corresponding derived analytic consideration. But since our current goal is to study the corresponding derived algebraic geometry carrying relative $p$-adic Hodge structure we will not explicitly mention the corresponding analytic context.	
\end{example}

\indent Similarly as what we did in the corresponding derived analytic geometry above, we consider the corresponding contact with \cite{10GV} in the derived algebraic geometric setting. We keep the notation as above. Then we have the following results with the same arguments as in the above on derived deformation of Hodge structures (note that this is more related to the non-derived situation):

\begin{proposition}
The corresponding category of all the perfect, almost perfect (after \cite[Section 7.2.4]{10Lu1}, namely pseudocoherent) $F^k$-equivariant sheaves over the families of $\infty$-Robba \'etale $\infty$-Deligne-Mumford toposes $(\mathrm{Spec}\widetilde{\Pi}_{I,R,A^h},\mathcal{O})_{\{I \subset (0,\infty)\}}$ (as in \cite[Chapter 1.4]{10Lu2}) is equivalent to the category of all the perfect, almost perfect (after \cite[Section 7.2.4]{10Lu1}, namely pseudocoherent) $F^k$-equivariant sheaves over some $\infty$-Robba \'etale $\infty$-Deligne-Mumford topos $(\mathrm{Spec}\widetilde{\Pi}_{I,R,A^h},\mathcal{O})$ (as in \cite[Chapter 1.4]{10Lu2}) where $I=[r_1,r_2]$ ($0<r_1\leq r_2/p^{hk}$). Here we need to assume the corresponding \cref{conjecture6.22} holds.	
\end{proposition}

\begin{proposition}
The corresponding category of all the locally finite free (after \cite[Section 7.2.2, 7.2.4]{10Lu1}, namely pseudocoherent and flat) $F^k$-equivariant sheaves over the families of $\infty$-Robba \'etale $\infty$-Deligne-Mumford toposes $(\mathrm{Spec}\widetilde{\Pi}_{I,R,A^h},\mathcal{O})_{\{I \subset (0,\infty)\}}$ (as in \cite[Chapter 1.4]{10Lu2}) is equivalent to the category of all the locally finite free (after \cite[Section 7.2.2, 7.2.4]{10Lu1}, namely pseudocoherent and flat) $F^k$-equivariant sheaves over some $\infty$-Robba \'etale $\infty$-Deligne-Mumford topos $(\mathrm{Spec}\widetilde{\Pi}_{I,R,A^h},\mathcal{O})$ (as in \cite[Chapter 1.4]{10Lu2}) where $I=[r_1,r_2]$ ($0<r_1\leq r_2/p^{hk}$).	
\end{proposition}

\indent The following is expected to hold in full generality.

\begin{conjecture} \label{conjecture6.22}
In the two settings in this subsection above, we conjecture that the corresponding descent holds (along binary rational coverings) for perfect, almost perfect and finite projective modules (which will allow us to compare sheaves and the global sections) in the noetherian situation. 
\end{conjecture}

\begin{remark}
In fact we believe that this is a direct consequence (in the Huber ring setting) after we apply the results of \cite[Theorem 14.9, Remark 14.10]{10CS} since we are in the noetherian setting.
\end{remark}

\begin{remark}
The corresponding descent could happen when we have the noetherianness, where the corresponding rational localization is actually automatically flat, which will imply that one can certainly glue coherent sheaves in our situation, see \cite[Theorem 1.3.9]{10KL1}. Beyond the corresponding noetherianness it might be not safe to conjecture so. One can also use \cite{10CS} to tackle this in the corresponding context of even general $\infty$-quasi-coherent sheaves, and beyond the corresponding noetherian situation purely algebraic definition will be very hard to control (where we think we have to definitely replace the purely algebraic consideration with the corresponding consideration in \cite{10CS} or \cite{10BK1}). 
\end{remark}
 	
\

This chapter is based on the following paper, where the author of this dissertation is the main author:
\begin{itemize}
\item Tong, Xin. "Topics on Geometric and Representation Theoretic Aspects of Period Rings I." arXiv preprint arXiv:2102.10693 (2021).
\end{itemize}

\newpage\chapter{Topologization and Functional Analytification I: Intrinsic Morphisms of Commutative Algebras}

\newpage\section{Introduction}

\subsection{Main Consideration}

\indent Scholze's diamond \cite{11Sch2} is actually very general notion beyond the corresponding perfectoid spaces, partially because it contains the corresponding diamantine spaces after Hansen-Kedlaya \cite{11HK}. This point of view certainly gives the motivation for this notion from Hansen-Kedlaya. Therefore one could regard to some extent diamantine spaces as giving some (maybe better to say more ring theoretic) analogs of the corresponding Scholze's diamonds \cite{11Sch2} \cite{11FS}, moreover they behave as if they are perfectoid spaces \cite{11Sch3}, \cite{11KL1}, \cite{11KL2}, \cite{11Ked1} and \cite{11GR}. Similar discussion could be made to sousperfectoid space, which should be more 'perfectoid' generalization.\\

\indent Hansen-Kedlaya \cite{11HK} have given the definition of naive \'etale morphisms among any Tate Huber pairs namely these are locally composites of the rational localizations and finite \'etale morphisms (very importantly with strongly sheafy domains and targets). This is because one definitely believes that correct notion of \'etale morphism should admit such admissible decomposition and factorization. However what should be the correct intrinsic one has not been given in full detail yet. In the significant strongly sheafy situation, we are going to try to answer this question as proposed in \cite[Appendix 5]{11Ked1}. In this paper, we try to study the corresponding properties of the corresponding naive \'etale morphisms along the ideas of \cite[Appendix 5]{11Ked1} and \cite{11HK}. The goal is to accurately characterize the corresponding naive \'etale morphisms in some intrinsic way.\\

\indent Sheafiness plays a very crucial role in the discussion above. However suppose we do not have to worry about the sheafiness at all (in fact in some sense we really do not have to worry about this at all by the work of Clausen-Scholze \cite{11CS} and Bambozzi-Kremnizer \cite{11BK}), then one might want to believe that the robust definitions could be made even more robust. Therefore we investigate the corresponding morphisms of the corresponding ring objects where sheafiness could be replaced by $\infty$-sheafiness (namely sheafiness up to higher homotopy) after \cite{11BBBK} and \cite{11BK}. One should be able to consider Clausen-Scholze's foundation \cite{11CS} as well, however we will mainly focus on the $\infty$-locally convex objects in \cite{11BBBK} and \cite{11BK}, as in the corresponding schematic situation in \cite{11Lu1}, \cite{11Lu2}, \cite{11TV1} and \cite{11TV2}. We consider the corresponding interesting approaches through the corresponding formal and PD completions just as in the $\infty$-schematic situation in \cite{11R} which is very related to the corresponding Drinfeld's stacky construction \cite{11Dr1} and \cite{11Dr2} by using the \v{C}ech-Alexander complex on the corresponding crystalline cohomology and prismatic cohomology.\\

\indent The current list of definitions will be established for discrete $\mathbb{E}_\infty$-ring objects and $\mathbb{E}_\infty$-ring objects in suitable locally convex $\infty$-categories after after Bambozzi-Ben-Bassat-Kremnizer \cite{11BBBK}:\\

\indent D1. Localized intrinsic \'etale morphisms of open mapping Huber rings;\\
\indent D2. Localized intrinsic \'etale morphisms of open mapping adic Banach rings;\\ 
\indent D3. Localized intrinsic lisse morphisms of open mapping Huber rings;\\
\indent D4. Localized intrinsic lisse morphisms of open mapping adic Banach rings;\\
\indent D5. Localized intrinsic non-ramifi\'e morphisms of open mapping Huber rings;\\
\indent D6. Localized intrinsic non-ramifi\'e morphisms of open mapping adic Banach rings;\\
\indent $\infty$1. De Rham intrinsic \'etale-like morphisms of $\infty$-analytic functors;\\
\indent $\infty$2. De Rham intrinsic lisse-like morphisms of $\infty$-analytic functors;\\ 
\indent $\infty$3. De Rham intrinsic non-ramifi\'e-like morphisms of $\infty$-analytic functors;\\
\indent $\infty$4. PD (crystalline) intrinsic \'etale-like morphisms of $\infty$-analytic functors;\\
\indent $\infty$5. PD (crystalline) intrinsic lisse-like morphisms of $\infty$-analytic functors;\\ 
\indent $\infty$6. PD (crystalline) intrinsic non-ramifi\'e-like morphisms of $\infty$-analytic functors.\\

\indent Certainly for general locally convex spaces producing nice ring structures we really have to be very precise and accurate in any sorts of characterization. However we have not unfortunately achieve this due to some very subtle issues, mainly coming from the corresponding issues in very general functional analytification. That being otherwise all said, we still actually could literally talk about the desired definitions for simplicial noetherian Banach rings in certain situations.

\subsection{Further Consideration}

Our ultimate goal is certainly to study the corresponding geometric sites (\'etale, pro-\'etale, crystalline and prismatic \cite{11SGAIV},\cite{11Gro1},\cite{11Sch1},\cite{11Sch2},\cite{11Sch3},\cite{11KL1},\cite{11KL2},\\
\cite{11Ked1},\cite{11BS},\cite{11Dr1},\cite{11Dr2}) and the corresponding cohomologies (\'etale, pro-\'etale, crystalline and prismatic) for really general $\infty$-analytic spaces (possibly also noncommutative analogs of those in \cite{11KR1}) over $\mathbb{F}_1$ and try to apply to the locally noetherian situations, the strongly noetherian situations (such as in \cite{11G1}, \cite{11GL}), the strongly sheafy situations under the foundation of $\infty$-locally convex spaces (as in \cite{11HK}, \cite{11KL1} and more general situations), although our very beginning corresponding motivation for this article is an attempt to answer some questions in \cite[Appendix A5]{11Ked1}.

\
%%\newpage

\newpage\section{Affinoid Morphisms of Huber Rings}

\indent We start with the discussion on the corresponding intrinsic definition of \'etale morphisms.

\begin{setting}
We start with an analytic uniform Huber pair $(A,A^+)$. And we will consider the category of all such rings. We assume the corresponding completeness for the Huber pairs.
\end{setting}

\begin{definition}\mbox{\bf{(Hansen-Kedlaya \cite[Definition 5.1]{11HK})}}
We call a map of Huber rings $(A,A^+)\rightarrow (B,B^+)$ naive \'etale after \cite[Definition 5.1]{11HK} if it admit a factorization into rational localizations and finite \'etale morphisms. Here we assume $(A,A^+)$ is strong sheafy and we assume that $(B,B^+)$ is strongly sheafy. \footnote{Certainly one needs to be careful since we are now considering more general context than \cite{11HK} without assuming the corresponding Tateness.}	
\end{definition}

\begin{definition}\mbox{\bf{(Kedlaya \cite[Definition A5.2]{11Ked1})}}\label{11definition2.3}
Recall from \cite[Definition A5.2]{11Ked1}, we have the corresponding affinoid morphism from any strongly sheafy Huber ring $A$, namely a morphism $A\rightarrow B$, such that $B$ admits some surjective covering from $A\left<T_1,...,T_d\right>$ and through this map we have that $B$ is a stably-pseudocoherent sheaf over $A\left<T_1,...,T_d\right>$ and the corresponding ring $B$ is assumed to be sheafy \footnote{Cetainly one needs to be more careful since this is also slightly different from the original definition \cite[Definition A5.2]{11Ked1}, thanks Professor Kedlaya for telling me this should be better. We want to mention that this is a quite subtle point around the sheafiness (see \cite[Theorem 1.4.20]{11Ked1}), the point here is that we do not know the kernel of an affinoid morphism is closed or not, if it is closed then we could keep the knowledge that $B$ being stably-pseudocoherent is equivalent to $B$ being sheafy. However if this is not closed, then we do not have this sort of equivalence to our knowledge.}.

\end{definition}

\indent The belief (as proposed in \cite[Problem A5.3, Problem A5.4]{11Ked1}) is that somehow the corresponding affinoid morphisms in the definition should be directly used in the corresponding definitions of lisse morphisms and unramified morphisms, as well as certainly the \'etale morphisms. To investigate this kind of idea, we are going to first investigate the corresponding naive \'etale morphisms along this idea.\\

\begin{lemma}
Let $f_1:\Gamma_1\rightarrow \Gamma_2$ and $f_2:\Gamma_2\rightarrow \Gamma_3$ be two affinoid morphisms, then the composition $f_2\circ f_1$ is also affinoid.	
\end{lemma}

\begin{proof}
Straightforward.	
\end{proof}

%%%%%%%%%%%%%%%%%%%%%%%%%%%%%%%%%!!!!!!!!!!!!!!

\begin{lemma}\mbox{\bf{(Kedlaya)}} \label{11lemma2.5}
For any standard binary rational localization of $A$ with respect to $f,g\in A$, suppose we know that there are two surjective morphisms:
\begin{align}
s_1:A\left<\frac{f}{g}\right>\left<T_1,...,T_n\right>\rightarrow B\left<\frac{f}{g}\right>,\\
s_2:A\left<\frac{g}{f}\right>\left<T_1,...,T_{n}\right>\rightarrow B\left<\frac{g}{f}\right>.	
\end{align}
Then we have that there is a surjective morphism:
\begin{align}
s:A\left<T_1,...,T_{n'}\right>\rightarrow B.
\end{align}	
\end{lemma}

\begin{proof}
The following argument is due to Kedlaya \footnote[1]{11Thanks Professor Kedlaya for mentioning the similarity of this to the corresponding locality of morphisms of finite type as in Grothendieck's EGA I and II.}, we work out it for the convenience of the readers. First, we have the following short exact sequence:
\[
\xymatrix@C+0pc@R+0pc{
&0 \ar[r] \ar[r] \ar[r] &B \ar[r] \ar[r] \ar[r] &B\left<\frac{f}{g}\right>\bigoplus B\left<\frac{g}{f}\right> \ar[r] \ar[r] \ar[r] &B\left<\frac{f}{g},\frac{g}{f}\right> \ar[r] \ar[r] \ar[r] &0.
}
\]
Take any $b\in B$, and use the notation $(b_1,b_2)$ for the image in the middle. By the surjectivity of the maps $s_1,s_2$ we have that there exist some element $a_1\in A\left<\frac{f}{g}\right>\left<T_1,...,T_n\right>$ and some element $a_2\in A\left<\frac{f}{g}\right>\left<T_1,...,T_n\right>$ such that we have:
\begin{align}
s_1(a_1)=b_1,\\
s_2(a_2)=b_2.	
\end{align}
With more explicit expression we have the following:
\begin{align}
s_1(\sum_{i_1,...,i_n}\sum_{i} a_1^{i,i_1,...,i_n}u^iT_1^{i_1}...T_n^{i_n})=\sum_i b_1^iu^i,\\
s_2(\sum_{i_1,...,i_n}\sum_{i} a_2^{i,i_1,...,i_n}v^iT_1^{i_1}...T_n^{i_n})=\sum_i b_2^iv^i,	
\end{align}	
under the corresponding presentations up to liftings:
\begin{align}
B\left<\frac{f}{g}\right>=B\left<u\right>/(gu-f),\\
B\left<\frac{g}{f}\right>=B\left<v\right>/(fv-g).\\	
\end{align}
Then to finish we only have to take some finite sum in the summation to make approximation. We first claim that such finite sum approximation and modification will not change the corresponding surjectivity of the map $s_1$ and $s_2$. Namely for each $k=1,2$ the map $s_k$ will maintain surjective once we modify the image of $T_1,...,T_n$ infinitesimally around some neighbourhood $U$ of $0$, in other words it will maintain to be surjective even if we set $s_k(T_1),...,s_k(T_n)$ to be $x_1,...,x_n$ whenever $x_1-s_k(T_1),...,x_n-s_k(T_n)$ lives in the neighbourhood $U$, and moreover we have that the corresponding modification could be assumed to take $T_i$ to $x_i$ with $i=1,...,n$. By open mapping, we have that the corresponding lifts of the corresponding differences $x_1-s_k(T_1),...,x_n-s_k(T_n)$ could be made to be living in some arbitrarily chosen neighbourhood $V$ of $0$. Then we only have to consider the following map factoring through the corresponding map $s_k$:
\begin{align}
h: A_k\left<T_1,...,T_n\right>&\rightarrow A_k\left<T_1,...,T_n\right>\\
	T_i&\mapsto T_i+\mathrm{lifts~of}~x_i-s_k(T_k)
\end{align}
where $A_1$ is the ring $A\left<\frac{f}{g}\right>$ while we have $A_2$ is the ring $A\left<\frac{g}{f}\right>$, which basically proves the claim. Then this will indicate that one can find some joint finite subset $T:=\{T_1,...,T_{n'}\}$ for $B\left<\frac{f}{g}\right>$ and $B\left<\frac{f}{g}\right>$ such that the modified  
\begin{align}
s_1:A\left<\frac{f}{g}\right>\left<T_1,...,T_{n'}\right>\rightarrow B\left<\frac{f}{g}\right>,\\
s_2:A\left<\frac{g}{f}\right>\left<T_1,...,T_{n'}\right>\rightarrow B\left<\frac{g}{f}\right>,	
\end{align}
are basically surjective and they fit into the following commutative diagram:
\[\tiny
\xymatrix@C+0pc@R+3pc{
0 \ar[r] \ar[r] \ar[r] &A\left<T_1,...,T_{n'}\right> \ar[d] \ar[d] \ar[d] \ar[r] \ar[r] \ar[r] &A\left<\frac{f}{g}\right>\left<T_1,...,T_{n'}\right>\bigoplus A\left<\frac{g}{f}\right>\left<T_1,...,T_{n'}\right> \ar[d] \ar[d] \ar[d]\ar[r] \ar[r] \ar[r] &A\left<\frac{f}{g},\frac{g}{f}\right>\left<T_1,...,T_{n'}\right> \ar[d] \ar[d] \ar[d]\ar[r] \ar[r] \ar[r] &0\\
0 \ar[r] \ar[r] \ar[r] &B  \ar[r] \ar[r] \ar[r] &B\left<\frac{f}{g}\right>\bigoplus B\left<\frac{g}{f}\right> \ar[r] \ar[r] \ar[r] &B\left<\frac{f}{g},\frac{g}{f}\right> \ar[r] \ar[r] \ar[r] &0,
}
\]
where the middle and the rightmost vertical arrows are surjective. Then claim is then that the left vertical one is also surjective. The kernels $K_1\oplus K_2$ in the middle is mapped surjectively to the kernel $K_{12}$ of the rightmost vertical map. So the snake lemma will force the cokernel of the left vertical arrow to  be zero which shows the corresponding exactness at the corresponding location $?$ in the following commutative diagram:
\[\tiny
\xymatrix@C+0pc@R+3pc{
 &0 \ar[d] \ar[d] \ar[d]   &0 \ar[d] \ar[d] \ar[d]  &0 \ar[d] \ar[d] \ar[d] & \\
0 \ar[r] \ar[r] \ar[r]& K\ar[d] \ar[d] \ar[d]\ar[r] \ar[r] \ar[r]&K_1\bigoplus K_2 \ar[d] \ar[d] \ar[d]\ar[r] \ar[r] \ar[r]&K_{12} \ar[d] \ar[d] \ar[d]\ar[r] \ar[r] \ar[r]&\\
0 \ar[r] \ar[r] \ar[r] &A\left<T_1,...,T_{n'}\right> \ar[d] \ar[d] \ar[d] \ar[r] \ar[r] \ar[r] &A\left<\frac{f}{g}\right>\left<T_1,...,T_{n'}\right>\bigoplus A\left<\frac{g}{f}\right>\left<T_1,...,T_{n'}\right> \ar[d] \ar[d] \ar[d]\ar[r] \ar[r] \ar[r] &A\left<\frac{f}{g},\frac{g}{f}\right>\left<T_1,...,T_{n'}\right> \ar[d] \ar[d] \ar[d]\ar[r] \ar[r] \ar[r] &0\\
0 \ar[r] \ar[r] \ar[r] &B  \ar[d]^? \ar[d] \ar[d]\ar[r] \ar[r] \ar[r] &B\left<\frac{f}{g}\right>\bigoplus B\left<\frac{g}{f}\right> \ar[d] \ar[d] \ar[d]\ar[r] \ar[r] \ar[r] &B\left<\frac{f}{g},\frac{g}{f}\right> \ar[d] \ar[d] \ar[d]\ar[r] \ar[r] \ar[r] &0\\
 &0   &0  &0 &
}
\]
where $K_1,K_2,K_{12}$ are pseudocoherent, which implies that the corresponding module $K$ is also pseudocoherent.

\end{proof}

\begin{proposition} \label{11proposition2.5}
Any naive \'etale morphism is affinoid.
\end{proposition}

\begin{proof}
First let $f:(A,A^+)\rightarrow (B,B^+)$ be any naive \'etale morphism. Then locally this is basically composition of the corresponding rational localizations and finite \'etale maps. Locally rational localizations involved are actually affinoid, and locally the corresponding finite \'etale maps from strongly sheafy rings will have strongly sheafy target, which will imply that locally finite \'etale maps are affinoid. Then this could be globalized to force the global map $f$ to be affinoid. The properties of factoring through a surjection globally could be proved by glueing local ones through \cref{11lemma2.5}, by considering \cite[Proposition 2.4.20]{11KL1}. And globally the corresponding ring $B$ is stably-pseudocoherent over $A\left<T_1,...,T_n\right>$ for some $n$ since this is a local property. 
\end{proof}

\indent Therefore we have proved that the corresponding \'etale maps in the corresponding naive sense is actually affinoid in the above sense. Therefore it is now natural to try to find the corresponding properties which may completely characterize the corresponding naive \'etale morphisms which are affinoid.

\indent Certainly we may have the corresponding conjectures that all the naive \'etale morphisms will satisfy the corresponding properties of algebraically \'etale ones (such as in \cite[Chapitre 17]{11EGAIV4}, \cite[Tag 00U1]{11SP}). We now discuss the corresponding completed cotangent complex after Huber \cite[1.6.2]{11Hu1}. Recall for our current $B$ the corresponding completed differential $\Omega^1_{B/A,\mathrm{topo}}$ (see \cite[1.6.2]{11Hu1} for the construction for any $f$-adic rings in the noetherian setting). Therefore we consider the corresponding topological naive cotangent complex:
\begin{align}
\tau_{\leq 1}\mathbb{L}_{B/A,\mathrm{topo}}	
\end{align}
for any naive \'etale map $f:(A,A^+)\rightarrow (B,B^+)$. We now discuss the construction without the corresponding strongly noetherian requirement in our current situation. First we know that $B$ is of topologically finite type over $A$:
\begin{align}
B=A\left<X_1,...,X_n\right>_{T_1,...,T_n}/I.	
\end{align}
Then we could first define the topological free differentials:
\begin{align}
\Omega^1:=A\left<X_1,...,X_n\right>_{T_1,...,T_n}dX_1+...+A\left<X_1,...,X_n\right>_{T_1,...,T_n}dX_n.	
\end{align}
Then we have:
\begin{align}
\Omega^1_{B/A,\mathrm{topo}}:=	\Omega^1\slash (I\bigcup d(I))\Omega^1.
\end{align}

Here everything is assumed to be basically complete with respect to the corresponding natural topology. Namely we need to take the corresponding completion always with respect to the corresponding induced topology. Certainly here $\Omega^1$ is already complete due to the fact that it is finitely projective. Recall that a map $f:\Gamma_1\rightarrow \Gamma_2$ is called \'etale in the scheme theory if the naive cotangent complex (truncated and could be regarded as an $\infty$-module spectrum) is quasi-isomorphic to zero. The corresponding underlying complex reads:
\[
\xymatrix@C+0pc@R+0pc{
[I/I^2\ar[r] \ar[r] \ar[r] &\Omega^1_{\Gamma_2/\Gamma_1,\mathrm{topo}}].
}
\]

\indent In the situation where we consider $A\rightarrow B$ is affinoid, the corresponding ideal $I$ is actually stably-pseudocoherent over $A\left<X_1,...,X_n\right>$. It is peudocoherent by the corresponding two out of three property. The stability holds locally, so we have the case. And if the morphism if furthermore naive \'etale then we have $I/I^2$ is also stably-pseudocoherent, see \cref{lemma4.8}. 

\begin{remark}
Note that we are considering the very general and complicated non-noetherian situation, modules will need to be endowed with the natural topology and complete, although finite projective modules are complete automatically. This will have nontrivial things to do with the corresponding definition of $\Omega^1_{B/A,\mathrm{topo}}$.	
\end{remark}

\indent One can actually generalize the corresponding full cotangent complexes and derived de Rham complexes to this topological context following \cite{11III1}, \cite{11III2} and \cite{11B1}. First for the corresponding topological cotangent complex we consider the following definition (note that we have to assume the corresponding topologically finite type condition). We start with the corresponding algebraic ones for $B^h=A[X_1,...,X_n]_{T_1,...,T_n}/I$, under the topologization we have the corresponding derived cotangent complex:
\begin{align}
\mathbb{L}_{B^h/A,\mathrm{alg}},	
\end{align}
by taking the usual algebraic one. Then we take the corresponding completion with respect to the corresponding topologization which gives rise to the following topological one:
\begin{align}
\mathbb{L}_{B/A,\mathrm{topo}}.	
\end{align}

We define the corresponding de Rham complex in the following parallel way. What is happen is that consider the presentation $B^h=A[X_1,...,X_n]_{T_1,...,T_n}/I$ which gives rise to the corresponding algebraic de Rham complex:
\[
\xymatrix@C+0pc@R+0pc{
0\ar[r]\ar[r]\ar[r] &B^h \ar[r]\ar[r]\ar[r] & \Omega^1_{B^h/A,\mathrm{alg}}\ar[r]\ar[r]\ar[r] & \Omega^2_{B^h/A,\mathrm{alg}}  \ar[r]\ar[r]\ar[r] &...\ar[r]\ar[r]\ar[r] & \Omega^\bullet_{B^h/A,\mathrm{alg}} \ar[r]\ar[r]\ar[r] &..., 
}
\]
which will give rise to the corresponding topological one if we take the corresponding completion induced from the subset $T_1,...,T_n$:
\[
\xymatrix@C+0pc@R+0pc{
0\ar[r]\ar[r]\ar[r] &B \ar[r]\ar[r]\ar[r] & {\Omega}^1_{B^h/A,\mathrm{topo}}\ar[r]\ar[r]\ar[r] & {\Omega}^2_{B^h/A,\mathrm{topo}}  \ar[r]\ar[r]\ar[r] &...\ar[r]\ar[r]\ar[r] & {\Omega}^\bullet_{B^h/A,\mathrm{topo}} \ar[r]\ar[r]\ar[r] &.... 
}
\]	

From our construction for ${\Omega}^1_{B^h/A,\mathrm{topo}}$, one can actually define:
\begin{align}
{\Omega}^{\bullet,\mathrm{f}}_{B^h/A,\mathrm{topo}}:=\bigoplus_{i_1,...,i_\bullet\in \{1,...,n\}}A\left<X_1,...,X_n\right>_{T_1,...,T_n} dX_{i_1}\wedge dX_{i_2}\wedge...\wedge dX_{i_\bullet}
\end{align}
and then define:
\begin{align}
{\Omega}^{\bullet}_{B/A,\mathrm{topo}}:=\left(\bigoplus_{i_1,...,i_\bullet\in \{1,...,n\}}A\left<X_1,...,X_n\right>_{T_1,...,T_n} dX_{i_1}\wedge dX_{i_2}\wedge...\wedge dX_{i_\bullet}\right)\slash\\
\left((I\bigcup dI \bigcup d^\bullet I)\bigoplus_{i_1,...,i_\bullet\in \{1,...,n\}}A\left<X_1,...,X_n\right>_{T_1,...,T_n} dX_{i_1}\wedge dX_{i_2}\wedge...\wedge dX_{i_\bullet}\right),
\end{align}
after taking suitable completion when needed \footnote[1]{11As in \cite{11B1} and \cite{11GL} where one takes the corresponding derived $p$-completion out from the algebraic cotangent complex and the corresponding derived algebraic de Rham complex.}.

\

\newpage\section{Affinoid Morphisms of Banach Rings}

\indent We now consider the parallel situation of Banach rings.

\begin{setting}
We start with a uniform adic Banach ring $(A,A^+)$ in the general sense of \cite{11KL1} and \cite{11KL2} (without assumption on the topologically nilpotent units being existing, but we assume this is open mapping). And we will consider the category of all such rings. We assume the corresponding completeness as well.
\end{setting}

\begin{definition}\mbox{\bf{(Hansen-Kedlaya \cite[Definition 5.1]{11HK})}}
We call a map of adic Banach rings $(A,A^+)\rightarrow (B,B^+)$ naive \'etale after \cite[Definition 5.1]{11HK} if it admit a factorization into rational localizations and finite \'etale morphisms. Here we assume $(A,A^+)$ is strong sheafy and we assume that $(B,B^+)$ is sheafy.	
\end{definition}

\begin{remark}
Certainly this is in more general situation than the corresponding context of \cite{11HK}.	
\end{remark}

\begin{definition}\mbox{\bf{(Kedlaya \cite[Definition A5.2]{11Ked1})}}\label{11definition2.3}
Recall from \cite[Definition A5.2]{11Ked1}, we have the corresponding affinoid morphism from any strongly sheafy adic Banach ring $A$, namely a morphism $A\rightarrow B$, such that $B$ admits some surjective covering from $A\left<T_1,...,T_d\right>$ and through this map we have that $B$ is a stably-pseudocoherent sheaf over $A\left<T_1,...,T_d\right>$ and we assume that $(B,B^+)$ is strongly sheafy.

\end{definition}

\begin{remark}
Of course the corresponding foundation is not the same but parallel to such situation we are considering now, however it is definitely reasonable and parallel to follow \cite[Appendix A5]{11Ked1} to give the definition here.	
\end{remark}

\indent The belief (as proposed in \cite[Problem A5.3, Problem A5.4]{11Ked1}) is that somehow the corresponding affinoid morphisms in the definition should be directly used in the corresponding definitions of lisse morphisms and unramified morphisms, as well as certainly the \'etale morphisms. To investigate this kind of idea, we are going to first investigate the corresponding naive \'etale morphisms along this idea.\\

\begin{lemma}\mbox{\bf{(Kedlaya)}} \label{11lemma2.5}
For any standard binary rational localization of $A$ with respect to $f,g\in A$, suppose we know that there are two surjective morphisms:
\begin{align}
s_1:A\left<\frac{f}{g}\right>\left<T_1,...,T_n\right>\rightarrow B\left<\frac{f}{g}\right>,\\
s_2:A\left<\frac{g}{f}\right>\left<T_1,...,T_{n}\right>\rightarrow B\left<\frac{g}{f}\right>.	
\end{align}
Then we have that there is a surjective morphism:
\begin{align}
s:A\left<T_1,...,T_{n'}\right>\rightarrow B.
\end{align}	
\end{lemma}

\begin{proof}
The following argument is due to Kedlaya, we work out it for the convenience of the readers. And this is the corresponding Banach analog of corresponding parallel result in the Huber ring situation. First, we have the following short exact sequence:
\[
\xymatrix@C+0pc@R+0pc{
&0 \ar[r] \ar[r] \ar[r] &B \ar[r] \ar[r] \ar[r] &B\left<\frac{f}{g}\right>\bigoplus B\left<\frac{g}{f}\right> \ar[r] \ar[r] \ar[r] &B\left<\frac{f}{g},\frac{g}{f}\right> \ar[r] \ar[r] \ar[r] &0.
}
\]
Take any $b\in B$, and use the notation $(b_1,b_2)$ for the image in the middle. By the surjectivity of the maps $s_1,s_2$ we have that there exist some element $a_1\in A\left<\frac{f}{g}\right>\left<T_1,...,T_n\right>$ and some element $a_2\in A\left<\frac{f}{g}\right>\left<T_1,...,T_n\right>$ such that we have:
\begin{align}
s_1(a_1)=b_1,\\
s_2(a_2)=b_2.	
\end{align}
With more explicit expression we have the following:
\begin{align}
s_1(\sum_{i_1,...,i_n}\sum_{i} a_1^{i,i_1,...,i_n}u^iT_1^{i_1}...T_n^{i_n})=\sum_i b_1^iu^i,\\
s_2(\sum_{i_1,...,i_n}\sum_{i} a_2^{i,i_1,...,i_n}v^iT_1^{i_1}...T_n^{i_n})=\sum_i b_2^iv^i,	
\end{align}	
under the corresponding presentations up to liftings:
\begin{align}
B\left<\frac{f}{g}\right>=B\left<u\right>/(gu-f),\\
B\left<\frac{g}{f}\right>=B\left<v\right>/(fv-g).\\	
\end{align}
Then to finish we only have to take some finite sum in the summation to make approximation. We first claim that such finite sum approximation and modification will not change the corresponding surjectivity of the map $s_1$ and $s_2$. Namely for each $k=1,2$ the map $s_k$ will maintain surjective once we modify the image of $T_1,...,T_n$ infinitesimally around some neighbourhood $U$ of $0$, in other words it will maintain to be surjective even if we set $s_k(T_1),...,s_k(T_n)$ to be $x_1,...,x_n$ whenever $\|x_1-s_k(T_1)\|\leq \delta,...,\|x_n-s_k(T_n)\|\leq \delta$ for some prescribed constant $\delta<1$ and moreover we have that the corresponding modification could be assumed to take $T_i$ to $x_i$ with $i=1,...,n$. By open mapping in this current context, we have that the corresponding lifts of the corresponding differences $x_1-s_k(T_1),...,x_n-s_k(T_n)$ could be made to be living in some arbitrarily chosen neighbourhood $V$ of $0$ namely, we can find lifts $y_1,...,y_n$ of these differences such that
\begin{align}
\|y_1\|<1,...,\|y_n\|<1.	
\end{align}
Then we only have to consider the following map factors through the corresponding map $s_k$:
\begin{align}
h: A_k\left<T_1,...,T_n\right>&\rightarrow A_k\left<T_1,...,T_n\right>\\
	T_i&\mapsto T_i+\mathrm{lifts~of}~x_i-s_k(T_k)
\end{align}
where $A_1$ is the ring $A\left<\frac{f}{g}\right>$ while we have $A_2$ is the ring $A\left<\frac{g}{f}\right>$, which basically proves the claim. Then this will indicate that one can find some joint finite subset $T:=\{T_1,...,T_{n'}\}$ for $B\left<\frac{f}{g}\right>$ and $B\left<\frac{f}{g}\right>$ such that the modified  
\begin{align}
s_1:A\left<\frac{f}{g}\right>\left<T_1,...,T_{n'}\right>\rightarrow B\left<\frac{f}{g}\right>,\\
s_2:A\left<\frac{g}{f}\right>\left<T_1,...,T_{n'}\right>\rightarrow B\left<\frac{g}{f}\right>,	
\end{align}
are basically surjective and they fit into the following commutative diagram:
\[\tiny
\xymatrix@C+0pc@R+3pc{
0 \ar[r] \ar[r] \ar[r] &A\left<T_1,...,T_{n'}\right> \ar[d] \ar[d] \ar[d] \ar[r] \ar[r] \ar[r] &A\left<\frac{f}{g}\right>\left<T_1,...,T_{n'}\right>\bigoplus A\left<\frac{g}{f}\right>\left<T_1,...,T_{n'}\right> \ar[d] \ar[d] \ar[d]\ar[r] \ar[r] \ar[r] &A\left<\frac{f}{g},\frac{g}{f}\right>\left<T_1,...,T_{n'}\right> \ar[d] \ar[d] \ar[d]\ar[r] \ar[r] \ar[r] &0\\
0 \ar[r] \ar[r] \ar[r] &B  \ar[r] \ar[r] \ar[r] &B\left<\frac{f}{g}\right>\bigoplus B\left<\frac{g}{f}\right> \ar[r] \ar[r] \ar[r] &B\left<\frac{f}{g},\frac{g}{f}\right> \ar[r] \ar[r] \ar[r] &0,
}
\]
where the middle and the rightmost vertical arrows are surjective. Then claim is then that the left vertical one is also surjective. The kernels $K_1\oplus K_2$ in the middle is mapped surjectively to the kernel $K_{12}$ of the rightmost vertical map. So the snake lemma will force the cokernel of the left vertical arrow to  be zero which shows the corresponding exactness at the corresponding location $?$ in the following commutative diagram:
\[\tiny
\xymatrix@C+0pc@R+3pc{
 &0 \ar[d] \ar[d] \ar[d]   &0 \ar[d] \ar[d] \ar[d]  &0 \ar[d] \ar[d] \ar[d] & \\
0 \ar[r] \ar[r] \ar[r]& K\ar[d] \ar[d] \ar[d]\ar[r] \ar[r] \ar[r]&K_1\bigoplus K_2 \ar[d] \ar[d] \ar[d]\ar[r] \ar[r] \ar[r]&K_{12} \ar[d] \ar[d] \ar[d]\ar[r] \ar[r] \ar[r]&\\
0 \ar[r] \ar[r] \ar[r] &A\left<T_1,...,T_{n'}\right> \ar[d] \ar[d] \ar[d] \ar[r] \ar[r] \ar[r] &A\left<\frac{f}{g}\right>\left<T_1,...,T_{n'}\right>\bigoplus A\left<\frac{g}{f}\right>\left<T_1,...,T_{n'}\right> \ar[d] \ar[d] \ar[d]\ar[r] \ar[r] \ar[r] &A\left<\frac{f}{g},\frac{g}{f}\right>\left<T_1,...,T_{n'}\right> \ar[d] \ar[d] \ar[d]\ar[r] \ar[r] \ar[r] &0\\
0 \ar[r] \ar[r] \ar[r] &B  \ar[d]^? \ar[d] \ar[d]\ar[r] \ar[r] \ar[r] &B\left<\frac{f}{g}\right>\bigoplus B\left<\frac{g}{f}\right> \ar[d] \ar[d] \ar[d]\ar[r] \ar[r] \ar[r] &B\left<\frac{f}{g},\frac{g}{f}\right> \ar[d] \ar[d] \ar[d]\ar[r] \ar[r] \ar[r] &0\\
 &0   &0  &0 &
}
\]
where $K_1,K_2,K_{12}$ are pseudocoherent, which implies that the corresponding module $K$ is also pseudocoherent.

\end{proof}

\begin{lemma}
Let $f_1:\Gamma_1\rightarrow \Gamma_2$ and $f_2:\Gamma_2\rightarrow \Gamma_3$ be two affinoid morphisms, then the composition $f_2\circ f_1$ is also affinoid.	
\end{lemma}

\begin{proof}
Straightforward.	
\end{proof}

\begin{proposition} \label{11proposition2.5}
Any naive \'etale morphism is affinoid.
\end{proposition}

\begin{proof}
See \cref{11proposition2.5}.
\end{proof}

\indent Therefore as in the corresponding Huber pair situation, we have proved that the corresponding \'etale maps in the corresponding naive sense is actually affinoid in the above sense. Therefore it is now natural to try to find the corresponding properties which may completely characterize the corresponding naive \'etale morphisms which are affinoid.

\indent Certainly we may have the corresponding conjectures that all the naive \'etale morphisms will satisfy the corresponding properties of algebraically \'etale ones (such as in \cite[Chapitre 17]{11EGAIV4}, \cite[Tag 00U1]{11SP}). We now discuss the corresponding Banach completed cotangent complex after Huber \cite[1.6.2]{11Hu1}. Recall for our current $B$ the corresponding completed differential $\Omega^1_{B/A,\mathrm{topo}}$ (see \cite[1.6.2]{11Hu1} for the construction for any $f$-adic rings). Certainly in the context of adic Banach rings we have the parallel completed version of the differentials by taking Banach completion, which is in some trival way the current situation. Therefore we consider the corresponding topological naive cotangent complex:
\begin{align}
\tau_{\leq 1}\mathbb{L}_{B/A,\mathrm{topo}}	
\end{align}
for any naive \'etale map $f:(A,A^+)\rightarrow (B,B^+)$. We consider the construction without the corresponding strongly noetherian requirement in our current situation. First we know that $B$ is of topologically finite type over $A$:
\begin{align}
B=A\left<X_1,...,X_n\right>_{T_1,...,T_n}/I.	
\end{align}
Then we could first define the Banach free differentials:
\begin{align}
\Omega^1:=A\left<X_1,...,X_n\right>_{T_1,...,T_n}dX_1+...+A\left<X_1,...,X_n\right>_{T_1,...,T_n}dX_n.	
\end{align}
Then we have:
\begin{align}
\Omega^1_{B/A,\mathrm{topo}}:=	\Omega^1\slash (I\bigcup d(I))\Omega^1.
\end{align}

Here everything is assumed to be basically complete with respect to the corresponding natural topology. Namely we need to take the corresponding completion always with respect to the corresponding induced norms. Certainly here $\Omega^1$ is already complete due to the fact that it is finitely projective. Recall that a map $f:\Gamma_1\rightarrow \Gamma_2$ is called \'etale in the scheme theory if the naive cotangent complex (truncated and could be regarded as an $\infty$-module spectrum) is quasi-isomorphic to zero. The corresponding underlying complex reads:
\[
\xymatrix@C+0pc@R+0pc{
[I/I^2\ar[r] \ar[r] \ar[r] &\Omega^1_{\Gamma_2/\Gamma_1,\mathrm{topo}}].
}
\]

\indent In the situation where we consider $A\rightarrow B$ is affinoid, the corresponding ideal $I$ is actually stably-pseudocoherent over $A\left<T_1,...,T_n\right>$. It is peudocoherent by the corresponding two out of three property. The stability holds locally, so we have the case. And if the morphism if furthermore naive \'etale then we have $I/I^2$ is also stably-pseudocoherent, see \cref{lemma4.8}. 

\begin{remark}
Note that we are considering the very general and complicated non-noetherian situation, modules will need to be endowed with the natural topology coming from the Banach structures on the Banach rings and complete, although in this situation as well finite projective modules are complete automatically. This will have nontrivial things to do with the corresponding definition of $\Omega^1_{B/A,\mathrm{topo}}$.	
\end{remark}

\indent In the Banach world, one can actually generalize the corresponding full cotangent complexes and de Rham complex to this context. First for the corresponding topological cotangent complex we consider the following definition (note that we have to assume the corresponding topologically finite type condition). We start with the corresponding algebraic ones for $B^h=A[X_1,...,X_n]_{T_1,...,T_n}/I$, under the topologization we have the corresponding derived cotangent complex:
\begin{align}
\mathbb{L}_{B^h/A,\mathrm{alg}},	
\end{align}
by taking the usual algebraic one. Then we take the corresponding completion with respect to the corresponding topologization which gives rise to the following topological one:
\begin{align}
\mathbb{L}_{B/A,\mathrm{topo}}.	
\end{align}

We define the corresponding de Rham complex in the following way parallely. What is happen is that consider the presentation $B^h=A[X_1,...,X_n]_{T_1,...,T_n}/I$ which gives rise to the corresponding algebraic de Rham complex:
\[
\xymatrix@C+0pc@R+0pc{
0\ar[r]\ar[r]\ar[r] &B^h \ar[r]\ar[r]\ar[r] & \Omega^1_{B^h/A,\mathrm{alg}}\ar[r]\ar[r]\ar[r] & \Omega^2_{B^h/A,\mathrm{alg}}  \ar[r]\ar[r]\ar[r] &...\ar[r]\ar[r]\ar[r] & \Omega^\bullet_{B^h/A,\mathrm{alg}} \ar[r]\ar[r]\ar[r] &..., 
}
\]
which will give rise to the corresponding topological one if we take the corresponding completion under $(\|.\|,\mathrm{Ban})$ induced from the subset $T_1,...,T_n$:
\[
\xymatrix@C-0.4pc@R+0pc{
0\ar[r]\ar[r]\ar[r] &B^h_{\|.\|,\mathrm{Ban}} \ar[r]\ar[r]\ar[r] & \Omega^1_{B^h/A,\mathrm{alg},\|.\|,\mathrm{Ban}}\ar[r]\ar[r]\ar[r] & \Omega^2_{B^h/A,\mathrm{alg},\|.\|,\mathrm{Ban}}  \ar[r]\ar[r]\ar[r] &...\ar[r]\ar[r]\ar[r] & \Omega^\bullet_{B^h/A,\mathrm{alg},\|.\|,\mathrm{Ban}} \ar[r]\ar[r]\ar[r] &, \\
0\ar[r]\ar[r]\ar[r] &B \ar[r]\ar[r]\ar[r] & {\Omega}^1_{B^h/A,\mathrm{topo}}\ar[r]\ar[r]\ar[r] & {\Omega}^2_{B^h/A,\mathrm{topo}}  \ar[r]\ar[r]\ar[r] &...\ar[r]\ar[r]\ar[r] & {\Omega}^\bullet_{B^h/A,\mathrm{topo}} \ar[r]\ar[r]\ar[r] &. 
}
\]	
From our construction for ${\Omega}^1_{B^h/A,\mathrm{topo}}$, one can actually define:
\begin{align}
{\Omega}^{\bullet,\mathrm{f}}_{B^h/A,\mathrm{topo}}:=\bigoplus_{i_1,...,i_\bullet\in \{1,...,n\}}A\left<X_1,...,X_n\right>_{T_1,...,T_n} dX_{i_1}\wedge dX_{i_2}\wedge...\wedge dX_{i_\bullet}
\end{align}
and then define:
\begin{align}
{\Omega}^{\bullet}_{B/A,\mathrm{topo}}:=\left(\bigoplus_{i_1,...,i_\bullet\in \{1,...,n\}}A\left<X_1,...,X_n\right>_{T_1,...,T_n} dX_{i_1}\wedge dX_{i_2}\wedge...\wedge dX_{i_\bullet}\right)\slash\\
\left((I\bigcup dI \bigcup d^\bullet I)\bigoplus_{i_1,...,i_\bullet\in \{1,...,n\}}A\left<X_1,...,X_n\right>_{T_1,...,T_n} dX_{i_1}\wedge dX_{i_2}\wedge...\wedge dX_{i_\bullet}\right),
\end{align}
after taking suitable completion under Banach norms when needed. One can also follow the construction in \cite{11III1}, \cite{11III2}, \cite{11B1} and \cite{11GL} to first consider the corresponding polynomial resolution $P_\bullet$ for $B$, then consider the corresponding algebraic derived de Rham complex $\Omega^*_{P_\bullet/A,\mathrm{alg}}$, then take the corresponding Banach completion to produce the corresponding topological one $\Omega^*_{P_\bullet/A,\mathrm{topo}}$. Then as in \cite{11III1}, \cite{11III2}, \cite{11GL} and \cite{11B1} around analytic derived $p$-adic de Rham complex we can take the corresponding suitable derived filtered completion to get the complex $\widehat{\Omega}^{\bullet}_{B/A,\mathrm{topo}}$ (certainly we need to consider Banach version of some filtered derived category of simplicial Banach rings). In the pro-\'etale site theoretic setting for rigid spaces (namely the corresponding p-complete context) this recovers the corresponding construction in \cite{11GL}. Recall in \cite{11GL} the corresponding analytic derived $p$-adic de Rham complex is constructed by first define the integral version ${\Omega}^{\bullet}_{B_0/A_0,\mathrm{topo}}$, and then take the colimit throughout all such rings of definition, and invert $p$, and then take the filtered completion. 

\begin{remark}
Certainly after this previous discussion we can construct the corresponding Banach derived de Rham complex, Banach cotangent complex and Banach Andr\'e-Quillen homology for any morphism $A\rightarrow B$ of Banach rings admissible in our situation. Recall from \cite{11III1}, \cite{11III2}, \cite{11B1} we have the corresponding algebraic $p$-adic derived cotangent complex:
\begin{align}
\Omega^1_{A[B]^\bullet/A}\otimes_{A[B]^\bullet}B
\end{align}
where $A[B]^\bullet$ is just the corresponding standard cofibrant replacement (in the topology theoretic language) of $B/A$. Then we take the corresponding derived completion \footnote{The derived completion in our Banach situation is the derived Banach completion which for instance could happen by using the 'completion' functor from $\mathrm{Simp}(\mathrm{Ind}(\mathrm{NormSets}))$ to $\mathrm{Simp}(\mathrm{Ind}(\mathrm{BanachSets}))$ literally in \cite{11BBBK}.} under the induced norm to achieve the corresponding topological one:
\begin{align}
\mathbb{L}_{B/A,\mathrm{topo}}:=(\Omega^1_{A[B]^\bullet/A}\otimes_{A[B]^\bullet}B)^\wedge_{\|.\|}.
\end{align}
We have the corresponding algebraic $p$-adic derived de Rham complex:
\begin{align}
\Omega^\bullet_{A[B]^\bullet/A}
\end{align}
where $A[B]^\bullet$ is just the corresponding standard cofibrant replacement (in the topology theoretic language) of $B/A$. Then we take the corresponding completion under the induced norm to achieve the corresponding topological one:
\begin{align}
\mathbb{\mathrm{dR}}_{B/A,\mathrm{topo}}:=(\Omega^\bullet_{A[B]^\bullet/A})^\wedge_{\|.\|},
\end{align}
carrying certain filtration $\mathrm{Fil}^*_{\mathbb{\mathrm{dR}}_{B/A,\mathrm{topo}}}$, which allows one take the corresponding filtered completion to achieve the corresponding final object:
\begin{align}
\mathbb{\mathrm{dR}}^\wedge_{B/A,\mathrm{topo}}:=(\Omega^\bullet_{A[B]^\bullet/A})^\wedge_{\|.\|,\mathrm{Fil}^*_{\mathbb{\mathrm{dR}}_{B/A,\mathrm{topo}}}}.
\end{align}
\end{remark}

\

\newpage\section{Naive \'Etale Morphisms and Intrinsic \'Etale Morphisms}

\indent As discussed above we now study the corresponding properties of naive \'etale morphisms aiming at the corresponding characterization of the corresponding correct definitions of intrinsic \'etale morphisms. We now assume the corresponding analyticity of the adic rings. 

\begin{lemma}
Let $f:(A,A^+)\rightarrow (B,B^+)$ be any rational localization map. Then we have that $B$ is of topologically finite type.
\end{lemma}

\begin{proof}
Straightforward.	
\end{proof}

\begin{lemma}
Let $f:(A,A^+)\rightarrow (B,B^+)$ be any rational localization map. Then we have that $B$ is affinoid over $A$ in the sense \cref{11definition2.3}.
\end{lemma}

\begin{proof}
See the proof of \cref{11proposition2.5}.	
\end{proof}

\begin{lemma}
Let $f:(A,A^+)\rightarrow (B,B^+)$ be any rational localization map. Then we have
\begin{align}
\tau_{\leq 1}\mathbb{L}_{B/A,\mathrm{topo}}	
\end{align}
is quasi-isomorphic to zero.
\end{lemma}

\begin{proof}
It suffices to reduce to standard binary localization such as simple Laurent or balanced localization, where we give a proof in the case of simple Laurent one:
\begin{align}
A\rightarrow A\{T\}/(T-f)	
\end{align}
for some $f\in A$. Then we have that actually the corresponding topological cotangent complex will be the corresponding completion of the corresponding algebraic ones (see \cite[Proposition 1.6.3]{11Hu1}).	The corresponding quasi-isomorphism could be defined directly. One just considers the following algebraic differential map:
\begin{align}
I=(T-f) \rightarrow A[T]/(T-f) dT\\	
\end{align}
which is actually surjective since for any:
\begin{align}
\sum_{i\geq 0}a_iT^i+g(T)(T-f) dT \in A[T]/(T-f) dT	
\end{align}
one takes the corresponding integration of:
\begin{align}
\int_{f}^T	&\sum_{i\geq 0}a_iT^i+g(T)(T-f) dT \\
&= \int_{f}^T	\sum_{i\geq 0}a_iT^i+(\sum_{i\geq 0} g_iT^i) (T-f) dT\\
&= \int_{f}^T	\sum_{i\geq 0}a_iT^i+\int_{f}^T (\sum_{i\geq 0} g_iT^i)(T-f) dT \\
&= \int_{f}^T	\sum_{i\geq 0}a_iT^i+ \int_f^T \sum_{i\geq 0} g_iT^{i+1}- f\int_f^T \sum_{i\geq 0} g_iT^{i}\\
&=\sum_{i\geq 0}a_i\frac{1}{i+1}T^{i+1}|_{f}^T+\sum_{i\geq 0}g_i\frac{1}{i+2}T^{i+2}|_{f}^T-f \sum_{i\geq 0}g_i\frac{1}{i+1}T^{i+1}|_{f}^T  \\
&= (*)(T-f),
\end{align}
where we only have finite sums here since we are considering the corresponding topologized polynomial. For this algebraic map, we have that kernel is $(T-f)^2$, for instance consider:
\begin{align}
d((T-f)h(T))=0,	
\end{align}
we will have:
\begin{align}
h(T)+(T-f)h'(T)dT=0,	
\end{align}
which implies that the image of $(T-f)h(T)$ lives in the corresponding quotient $A[T]/(T-f)dT$. Therefore we have that the topologized (not complete yet) cotangent complex:
\[
\xymatrix@C+0pc@R+0pc{
[(T-f)/(T-f)^2\ar[r] \ar[r] \ar[r] &A[T]/(T-f)dT],
}
\]
which is quasi-isomorphic to zero. Then we take the corresponding completion with respect to the corresponding topology induced from $A$ we have the desired result.

\end{proof}

\begin{lemma}
Let $f:(A,A^+)\rightarrow (B,B^+)$ be any finite \'etale map. Then we have that $B$ is of topologically finite type.
\end{lemma}

\begin{proof}
Since we have that $B$ is affinoid over $A$ by \cref{11proposition2.5}.	
\end{proof}

\begin{lemma}
Let $f:(A,A^+)\rightarrow (B,B^+)$ be any finite \'etale map. Then we have that $B$ is affinoid over $A$ in the sense \cref{11definition2.3}.
\end{lemma}

\begin{proof}
See the proof of \cref{11proposition2.5}.	
\end{proof}

\begin{lemma}
Let $f:(A,A^+)\rightarrow (B,B^+)$ be any finite \'etale map. Then we have 
\begin{align}
\tau_{\leq 1}\mathbb{L}_{B/A,\mathrm{topo}}	
\end{align}
is quasi-isomorphic to zero.	
\end{lemma}

\begin{proof}
This is basically nontrivial due to the fact that we are discussing the corresponding topological cotangent complex. However, one takes the corresponding finite presentation $B=A[T_1,...,T_n]/(f_1,...,f_p)$ (note that in fact that we have that $B$ is finite over $A$) since we are considering a finite \'etale map, which will realize the corresponding desired algebraic cotangent complex:
\[
\xymatrix@C+0pc@R+0pc{
[I/I^2\ar[r] \ar[r] \ar[r] &\Omega^1_{B/A}].
}
\]
For the topological situation we have that $B=A\{T_1,...,T_n\}/(f_1,...,f_p)$ by taking the corresponding completion. Again note that this means that actually the corresponding $A$-algebra $B$ is still finite over $A$ (since that is the very assumption). Therefore we could have the chance to right $B$ as just $\oplus_{i} Ae_i$ this is basically inducing the same differential module $\bigoplus_j(\bigoplus_{i} Ae_i) dT_j$ in both the corresponding topological setting and the algebraic setting. Then one could get the corresponding desired topological cotangent complex which is quasi-isomorphic to zero.
\end{proof}

\indent Then we consider the corresponding local composition:

\begin{lemma}
Let $f:(A,A^+)\rightarrow (B,B^+)$ be any naive \'etale map. Then we have that $B$ is of topologically finite type.
\end{lemma}

\begin{proof}
Since we have that $B$ is affinoid over $A$ by \cref{11proposition2.5}.	
\end{proof}

\begin{lemma} \label{lemma4.8}
Let $f:(A,A^+)\rightarrow (B,B^+)$ be any naive \'etale morphism. Then we have 
\begin{align}
\tau_{\leq 1}\mathbb{L}_{B/A,\mathrm{topo}}	
\end{align}
is quasi-isomorphic to zero locally with respect to the corresponding rational localization.	
\end{lemma}

\begin{proof}
Locally we have that any naive \'etale morphism takes the corresponding truncated cotangent complex to be trivialized, by the corresponding composition properties of the cotangent complex \cite[Tag 08PN]{11SP}.   
\end{proof}

\indent In order to globalize the picture one has to work harder. First we have the following:

\begin{proposition}
Let $f:(A,A^+)\rightarrow (B,B^+)$ be any naive \'etale morphism. Then $f$ is affinoid, of topologically finite type.	
\end{proposition}

\begin{proof}
This is by \cref{11proposition2.5} for the affinoidness, which implies the corresponding second property. 
\end{proof}

\begin{definition}
We now define localized \footnote{However this could actually be globalized easily.} intrinsic \'etale morphism to be a morphism $f:(A,A^+)\rightarrow (B,B^+)$ which is affinoid with strongly sheafy target, and locally (with respect to the rational localization) the corresponding truncated topological cotangent complex is quasi-isomorphic to zero.	
\end{definition}

%\indent We now consider the corresponding properties of the localized intrinsic \'etale morphisms defined above.

\indent We now consider the corresponding intrinsic \'etale morphisms of the corresponding special adic spaces after \cite[Section 4, Section 7 and Section 11]{11HK}:

\begin{setting}
We now consider the three special adic spaces after \cite[Section 4, Section 7 and Section 11]{11HK}, they are the corresponding strongly sheafy adic spaces, the corresponding sousperfectoid adic spaces and the corresponding diamantine adic spaces. We will use the notations $T,S,D$ to denote them in general respectively.	
\end{setting}

\indent The corresponding categories of strongly sheafy adic spaces, sousperfectoid spaces and diamantine adic spaces are nice enough since at least we have well-defined notion of naive \'etale morphisms (which is certainly the correct one) and furthermore well-defined \'etale sites.

\begin{definition}
For strongly sheafy adic spaces, a morphism $T_1\rightarrow T_2$ is called localized intrinsic \'etale if locally on $T_1$ this is localized intrinsic \'etale, namely for any neighbourhood $U\subset T_1$ we have that the morphism $(\mathcal{O}_{T_2}(U'),\mathcal{O}^+_{T_2}(U'))\rightarrow (\mathcal{O}_{T_1}(U),\mathcal{O}^+_{T_1}(U))$ is localized intrinsic \'etale. 	
\end{definition}

\begin{definition}
For sousperfectoid adic spaces, a morphism $S_1\rightarrow S_2$ is called localized intrinsic \'etale if locally on $S_1$ this is localized intrinsic \'etale, namely for any neighbourhood $U\subset S_1$ we have that the morphism $(\mathcal{O}_{S_2}(U'),\mathcal{O}^+_{S_2}(U'))\rightarrow (\mathcal{O}_{S_1}(U),\mathcal{O}^+_{S_1}(U))$ is localized intrinsic \'etale. 	
\end{definition}

\begin{definition}
For diamantine adic spaces, a morphism $D_1\rightarrow D_2$ is called localized intrinsic \'etale if locally on $D_1$ this is localized intrinsic \'etale, namely for any neighbourhood $U\subset D_1$ we have that the morphism $(\mathcal{O}_{D_2}(U'),\mathcal{O}^+_{D_2}(U'))\rightarrow (\mathcal{O}_{D_1}(U),\mathcal{O}^+_{D_1}(U))$ is localized intrinsic \'etale. 	
\end{definition}

\

\newpage\section{Properties}

\indent We now study the corresponding properties of the corresponding localized intrinsic \'etale morphisms, following \cite[Chapitre 17]{11EGAIV4} and \cite[Especially Section 1.6, Section 1.7]{11Hu1}. We now assume the corresponding analyticity of the adic rings.

\begin{conjecture}
Any localized intrinsic \'etale morphism of strongly sheafy adic spaces is locally a composition of rational localization and finite \'etale morphism.	
\end{conjecture}

\indent Here is the special situation.

\begin{proposition}
As in \cite[Section 1.6, Section 1.7]{11Hu1}, namely in the strongly noetherian situation we have the conjecture holds.	
\end{proposition}

\begin{proof}
This is because in that setting our definition in the intrinsic setting coincides with the more algebraic one in \cite[Section 1.6, Section 1.7]{11Hu1}. And note that in this setting the affinoidness of the morphism reduces to just being admitting surjections from Tate algebra over the source.	
\end{proof}

%\begin{remark}
%If locally the presentation takes the form of $A\left<T\right>/(f)$, then we have that 
%\begin{align}
%\frac{d(f^2)}{dT}=c(T)f	
%\end{align}
%which implies that we have:
%\begin{align}
%2ff'=cf	
%\end{align}
%so we have that $f'=c(T)/2$ when $f\neq 0$. Therefore we have $f$ is just a linear function $xT+y$ times a general series $a(T)$.	
%\end{remark}

\indent If this is true then we have:

\begin{corollary}
Any localized intrinsic \'etale morphism of sousperfectoid adic spaces is locally a composition of rational localization and finite \'etale morphism. Any localized intrinsic \'etale morphism of diamantine adic spaces is locally a composition of rational localization and finite \'etale morphism.	
\end{corollary}

\begin{proposition}
Compositions of localized intrinsic \'etale morphisms of strongly sheafy adic spaces are again localized intrinsic \'etale morphism.	
\end{proposition}

\begin{proof}
Locally it is the corresponding compositions of topologically finite type morphism, and locally it is the corresponding compositions of the corresponding affinoid morphisms, and locally it is the corresponding compositions of morphisms giving rise to the quasi-isomorphic to zero truncated cotangent complex.	
\end{proof}

\begin{corollary}
Compositions of localized intrinsic \'etale morphisms of sousperfectoid adic spaces are again localized intrinsic \'etale morphism. Compositions of localized intrinsic \'etale morphisms of diamantine adic spaces are again localized intrinsic \'etale morphism.		
\end{corollary}

\begin{proposition}
The localized intrinsic \'etaleness of any morphism $T_1\rightarrow T_2$ of strongly sheafy adic rings is preserved under the base change along any morphism of $T_3\rightarrow T_2$.
\end{proposition}

\begin{proof}
The base change of any morphism of topologically finite type is again of topologically finite type. The affinoidness of morphism is also preserved under any base change morphism. Finally for the cotangent complex locally, we definitely have the corresponding result as well.	
\end{proof}

\begin{proposition}
The \'etale property of a morphism between strongly sheafy rings could be detected locally at each point.	
\end{proposition}

\begin{proof}
Straightforward.	
\end{proof}

\indent We now consider some functoriality issue in our current situation. Now we consider the corresponding localized intrinsic \'etale morphisms under the construction of Witt vectors. Now let:
\begin{align}
A\rightarrow B	
\end{align}
be a general morphism in positive characteristic. Therefore we can take the corresponding Witt vector construction:
\begin{align}
W(A^\flat)\rightarrow W(B^\flat),	
\end{align}
where we assume that $A^\flat\rightarrow B^\flat$ is localized intrinsic \'etale. Here we take the completion if needed along the corresponding Fontainisation.

\begin{proposition}
The map 
\begin{align}
W(A^\flat)\rightarrow W(B^\flat),	
\end{align}
is affinoid if the kernel is closed \footnote[1]{11This is again due to the very subtle point around the sheafiness such as in \cite[Theorem 1.4.20]{11Ked1}.}.	
\end{proposition}

\begin{proof}
We only need to check this locally, locally we have that there is a lifting:
\begin{align}
W(A^\flat)\{T_1,...\}\rightarrow W(B^\flat)\rightarrow 0	
\end{align}
from the corresponding surjection:
\begin{align}
A^\flat \{\overline{T}_1,...\}\rightarrow B^\flat \rightarrow 0.	
\end{align}
And what we have is that this map on the Witt vector level is also realizing the target as a stably-pseudocoherent module over the source since we have that the target is sheafy (\cite[Theorem 1.4.20]{11Ked1}). 
\end{proof}

\begin{proposition}
Same holds for the construction of integral Robba ring $\widetilde{\mathcal{R}}^r_*$ and Robba ring $\widetilde{\mathcal{R}}^{[s,r]}_*$ with respect to closed intervals as in \cite[Definition 5.1.1]{11KL1} and \cite[Definition 4.1.1]{11KL2}.	
\end{proposition}

\

\newpage\section{\'Etale-Like Morphisms of $\infty$-Banach Rings and the $\infty$-Analytic Stacks}

\subsection{Approach through De Rham Stacks}

\indent We now extend the corresponding discussion to the $\mathrm{E}_\infty$ objects in \cite[Remark 3.16]{11BBBK} by using the ideas as in \cite{11R}. Recall from \cite[Theorem 3.14]{11BBBK} we have the corresponding categories $\mathrm{Simp}\mathrm{Ind}^m(\mathrm{BanSets}_{H})$ and $\mathrm{Simp}\mathrm{Ind}(\mathrm{BanSets}_{H})$ which are the corresponding categories of the corresponding simplicial sets over the corresponding inductive categories of the corresponding Banach sets over some Banach ring $H$ \footnote{It is safer to assume the open mapping properties on homotopy groups.}.

\begin{theorem}\mbox{\bf{(Bambozzi-Ben-Bassat-Kremnizer)}} The corresponding categories\\
 $\mathrm{Simp}\mathrm{Ind}^m(\mathrm{BanSets}_{H})$ and $\mathrm{Simp}\mathrm{Ind}(\mathrm{BanSets}_{H})$ admit symmetric monoidal model categorical structure. Same holds for $\mathrm{Simp}\mathrm{Ind}^m(\mathrm{NrSets}_{H})$ and $\mathrm{Simp}\mathrm{Ind}(\mathrm{NrSets}_{H})$.	
\end{theorem}

\begin{corollary}\mbox{}\\
The corresponding categories $\mathrm{Simp}\mathrm{Ind}^m(\mathrm{BanSets}_{H})$ and $\mathrm{Simp}\mathrm{Ind}(\mathrm{BanSets}_{H})$ admit presentations as $(\infty,1)$-categories. Same holds for 
\begin{center}
$\mathrm{Simp}\mathrm{Ind}^m(\mathrm{NrSets}_{H})$
\end{center}
and
\begin{center}
$\mathrm{Simp}\mathrm{Ind}(\mathrm{NrSets}_{H})$.
\end{center}	
\end{corollary}

\indent Then recall from \cite[Remark 3.16]{11BBBK} we have the corresponding ring objects in the $\infty$-categories above:
\begin{align}
\mathrm{sComm}(\mathrm{Simp}\mathrm{Ind}^m(\mathrm{BanSets}_{H})),\\
\mathrm{sComm}(\mathrm{Simp}\mathrm{Ind}(\mathrm{BanSets}_{H})).	
\end{align}
and 
\begin{align}
\mathrm{sComm}(\mathrm{Simp}\mathrm{Ind}^m(\mathrm{NrSets}_{H})),\\
\mathrm{sComm}(\mathrm{Simp}\mathrm{Ind}(\mathrm{NrSets}_{H})).	
\end{align}

\indent Now we use general notation $A$ to denote any object in these categories, regarding as a general $\mathrm{E}_\infty$-ring. We consider the general morphism $A\rightarrow B$ in the first two categories in the following discussion.

\begin{definition}
For any general morphism $A\rightarrow B$, we call this affinoid if we have that that $\pi_0(B)$ is affinoid over $\pi_0(A)$, namely we have that there is a surjection map $\pi_0(A)\left<X_1,...,X_d\right>\rightarrow \pi_0(B)$. And moreover we assume that $\pi_0(B)\otimes_{\pi_0(A)}\pi_n(A)\overset{\sim}{\rightarrow}\pi_n(B)$, for any $n$. 	
\end{definition}

\begin{remark}
Kedlaya's theorem \cite[Theorem 1.4.20]{11Ked1} is actually expected to hold in more general setting, at least in the situation where the definition of the affinoid morphisms could be made independent from the corresponding stably-pseudocoherence for open mapping rings (note that we are working over analytic fields). However in the previous definition, we have been not really exact. To be really accurate in the characterization of some desired notion of the affinoidness we think that one has to add certain $\infty$-sheafiness (which certainly holds in \cite{11BK}). To be more precise for any general morphism $A\rightarrow B$ in \cite{11BK} (namely in current situation one considers the corresponding Banach algebras over the analytic fields in our situation), we call this affinoid if we have that that $\pi_0(B)$ is affinoid over $\pi_0(A)$, namely we have that there is a surjection map $\pi_0(A)\left<X_1,...,X_d\right>\rightarrow \pi_0(B)$. And moreover we assume that $\pi_0(B)\otimes_{\pi_0(A)}\pi_n(A)\overset{\sim}{\rightarrow}\pi_n(B)$, for any $n$. In this situation we have the nice sheafiness (up to higher homotopy). Again similar discussion could be made in the context of \cite{11CS}. Note that \cite[Theorem 1.4.20]{11Ked1} literally says that the sheafiness is equivalent (in some nice sense but in more flexible derived sense) to the stably-pseudocoherence.
\end{remark}

\begin{definition}
For any general morphism $A\rightarrow B$, we call this localized intrinsic \'etale if we have that that $\pi_0(B)$ is localized intrinsic \'etale over $\pi_0(A)$, namely we have that there is a surjection map $\pi_0(A)\left<X_1,...,X_d\right>\rightarrow \pi_0(B)$ and we have that locally the corresponding truncated topological cotangent complex is basically quasi-isomorphic to zero. And moreover we assume that $\pi_0(B)\otimes_{\pi_0(A)}\pi_n(A)\overset{\sim}{\rightarrow}\pi_n(B)$, for any $n$.	
\end{definition}

\indent We now use the corresponding $X=\mathrm{Spec}A$ to denote the corresponding $\infty$-stack in the opposite categories with respect to the ring $A$. We now define the corresponding de Rham stack attached to $X$ as in \cite[Remark 1.2]{11R}:

\begin{definition}
We now define:
\begin{align}
X_\mathrm{dR}(R):=\varinjlim_{I} X(\pi_0(R)/I)	
\end{align}
for any $R$ in 
\begin{align}
\mathrm{sComm}(\mathrm{Simp}\mathrm{Ind}^m(\mathrm{BanSets}_{H})),\\
\mathrm{sComm}(\mathrm{Simp}\mathrm{Ind}(\mathrm{BanSets}_{H})).	
\end{align}
And the injective limit is taking throughout all nilpotent ideals of $\pi_0(R)$. 	
\end{definition}

\begin{remark}
As in \cite[Definition 1.1, Remark 1.2]{11R}, one can actually define the corresponding de Rham and crystalline spaces for any functor from \\$\mathrm{sComm}(\mathrm{Simp}\mathrm{Ind}^m(\mathrm{BanSets}_{H})))	$ and $\mathrm{sComm}(\mathrm{Simp}\mathrm{Ind}(\mathrm{BanSets}_{H})))$ to $\underline{s\mathrm{Sets}} $. This means we do not have to consider $(\infty,1)$-sheaves satisfying certain $\infty$-descent conditions. 
\end{remark}

\begin{definition}
The corresponding formal completion of any morphism $X=\mathrm{Spec}B\rightarrow Y=\mathrm{Spec}A$:
\begin{align}
Y_{X,\mathrm{dR}}	
\end{align}
is defined to be:
\begin{align}
Y_{X,\mathrm{dR}}(R):=\varinjlim_{I} X(\pi_0(R)/I)\times_{Y(\pi_0(R)/I)}	Y(R),
\end{align}
for any $R$ in 
\begin{align}
\mathrm{sComm}(\mathrm{Simp}\mathrm{Ind}^m(\mathrm{BanSets}_{\mathbb{Q}_p})),\\
\mathrm{sComm}(\mathrm{Simp}\mathrm{Ind}(\mathrm{BanSets}_{\mathbb{Q}_p})).	
\end{align}
And the injective limit is taking throughout all nilpotent ideals of $\pi_0(R)$.	
\end{definition}

\begin{remark}
Certainly it is actually not clear how really we should deal with the corresponding ideals here, namely we are not for sure if we need to consider closed ideals. But for simplicial noetherian rings we really have some nice definitions, which will certainly be tangential to the corresponding Huber's original consideration.
\end{remark}

\begin{definition}
We now define the corresponding de Rham intrinsic \'etale morphism to be an affinoid morphism $X=\mathrm{Spec}B\rightarrow Y=\mathrm{Spec}A$ which satisfies the condition:
\begin{align}
\pi_0(X(R))\overset{\sim}{\rightarrow}	\pi_0(Y_{X,\mathrm{dR}}(R)),
\end{align}
for any $\mathrm{E}_\infty$-object $R$.	
\end{definition}

\begin{proposition}
Compositions of de Rham intrinsic \'etale morphisms are again $\mathrm{PD}$ intrinsic \'etale morphisms.
\end{proposition}

\begin{proof}
This is formal.	
\end{proof}

\subsection{Approach through Crystalline Stack and PD-morphisms}

\indent We now use the corresponding $X=\mathrm{Spec}A$ to denote the corresponding $\infty$-stack in the opposite categories with respect to the ring $A$. We now define the corresponding crystalline stack attached to $X$ as in \cite[Definition 1.1]{11R}:

\begin{definition}
We now define:
\begin{align}
X_\mathrm{crys}(R):=\varinjlim_{I,\gamma} X(\pi_0(R)/I)	
\end{align}
for any $R$ in 
\begin{align}
\mathrm{sComm}(\mathrm{Simp}\mathrm{Ind}^m(\mathrm{BanSets}_{H})),\\
\mathrm{sComm}(\mathrm{Simp}\mathrm{Ind}(\mathrm{BanSets}_{H})).	
\end{align}
And the injective limit is taking throughout all nilpotent ideals of $\pi_0(R)$ and the corresponding PD-structures. 	
\end{definition}

\begin{definition}
The corresponding PD completion of any morphism $X=\mathrm{Spec}B\rightarrow Y=\mathrm{Spec}A$:
\begin{align}
Y_{X,\mathrm{crys}}	
\end{align}
is defined to be:
\begin{align}
Y_{X,\mathrm{crys}}(R):=\varinjlim_{I,\gamma} X(\pi_0(R)/I)\otimes_{Y(\pi_0(R)/I)}	Y(R),
\end{align}
for any $R$ in 
\begin{align}
\mathrm{sComm}(\mathrm{Simp}\mathrm{Ind}^m(\mathrm{BanSets}_{H})),\\
\mathrm{sComm}(\mathrm{Simp}\mathrm{Ind}(\mathrm{BanSets}_{H})).	
\end{align}
And the injective limit is taking throughout all nilpotent ideals of $\pi_0(R)$ and all the corresponding PD structures.	
\end{definition}

\begin{remark}
Certainly it is actually not clear how really we should deal with the corresponding ideals here and the corresponding PD structures, namely we are not for sure if we need to consider closed ideals. But for simplicial noetherian rings we really have some nice definitions, which will certainly be tangential to the corresponding Huber's original consideration.
\end{remark}

\begin{definition}
We now define the corresponding PD intrinsic \'etale morphism to be an affinoid morphism $X=\mathrm{Spec}B\rightarrow Y=\mathrm{Spec}A$ which satisfies the condition:
\begin{align}
\pi_0(X(R))\overset{\sim}{\rightarrow}	\pi_0(Y_{X,\mathrm{crys}}(R)),
\end{align}
for any $\mathrm{E}_\infty$-object $R$.		
\end{definition}

\begin{proposition}
We have that any de Rham intrinsic \'etale morphism is a PD intrinsic \'etale morphism.	
\end{proposition}

\begin{proof}
This is formal.	
\end{proof}

\begin{proposition}
Compositions of $\mathrm{PD}$ intrinsic \'etale morphisms are again $\mathrm{PD}$ intrinsic \'etale morphisms.
\end{proposition}

\begin{proof}
This is formal.	
\end{proof}

\

\newpage\section{Lisse-Like and Non-Ramifi\'e-Like Morphisms of $\infty$-Banach Rings and the $\infty$-Analytic Stacks}

\subsection{Approach through De Rham Stacks}

\indent We now define the corresponding lisse-like morphisms along the idea in the previous section:

\begin{definition}
For any general morphism $A\rightarrow B$, we call this localized intrinsic lisse if we have that that $\pi_0(B)$ is localized intrinsic lisse over $\pi_0(A)$, namely we have that there is a surjection map $\pi_0(A)\left<X_1,...,X_d\right>\rightarrow \pi_0(B)$ and we have that locally the corresponding truncated topological cotangent complex is basically quasi-isomorphic to $\Omega_{\pi_0(B)/\pi_0(A)}[0]$. And moreover we assume that $\pi_0(B)\otimes_{\pi_0(A)}\pi_n(A)\overset{\sim}{\rightarrow}\pi_n(B)$, for any $n$.	
\end{definition}

\begin{definition}
We now define the corresponding de Rham intrinsic lisse morphism to be an affinoid morphism $X=\mathrm{Spec}B\rightarrow Y=\mathrm{Spec}A$ which satisfies the condition:
\begin{align}
\pi_0(X(R))\overset{}{\rightarrow}	\pi_0(Y_{X,\mathrm{dR}}(R))
\end{align}
being surjective, for any $\mathrm{E}_\infty$-object $R$.		
\end{definition}

\begin{definition}
For any general morphism $A\rightarrow B$, we call this localized intrinsic non-ramifi\'e if we have that that $\pi_0(B)$ is localized intrinsic non-ramifi\'e over $\pi_0(A)$, namely we have that there is a surjection map $\pi_0(A)\left<X_1,...,X_d\right>\rightarrow \pi_0(B)$ and we have that locally the corresponding truncated topological cotangent complex is basically quasi-isomorphic to $\Omega_{\pi_0(B)/\pi_0(A)}[0]$ which vanishes as well. And moreover we assume that $\pi_0(B)\otimes_{\pi_0(A)}\pi_n(A)\overset{\sim}{\rightarrow}\pi_n(B)$, for any $n$.	
\end{definition}

%\begin{definition}
%We now define the corresponding de Rham intrinsic non-ramifi\'e morphism to be an affinoid morphism $X=\mathrm{Spec}B\rightarrow Y=\mathrm{Spec}A$ which satisfies the condition:
%\begin{align}
%\pi_0(X)\overset{}{\rightarrow}	\pi_0(Y_{X,\mathrm{dR}})
%\end{align}
%being surjective.	
%\end{definition}

\subsection{Approach through Crystalline Stack and PD-morphisms}

\begin{definition}
We now define the corresponding PD intrinsic lisse morphism to be an affinoid morphism $X=\mathrm{Spec}B\rightarrow Y=\mathrm{Spec}A$ which satisfies the condition:
\begin{align}
\pi_0(X(R))\overset{}{\rightarrow}	\pi_0(Y_{X,\mathrm{crys}}(R))
\end{align}
being surjective, for any $\mathrm{E}_\infty$-object $R$\footnote{For non-ramifi\'e situation one considers injectivity.}.		
\end{definition}

\begin{proposition}
We have that any de Rham intrinsic lisse morphism is a PD intrinsic lisse morphism.	
\end{proposition}

\

\newpage\section{Perfectization and Fontainisation of $\infty$-Analytic Stacks Situation}

\subsection{Perfectization, Fontainisation and Crystalline Stacks}

\indent Now we consider the corresponding perfectoidization of $\infty$-analytic stacks after \cite{11R} and \cite{11Dr1} in the situation where $H$ is assumed to be of characteristic $p$. 

\begin{definition}
For any object in $\infty-\mathrm{Fun}(\mathrm{sComm}\mathrm{Simp}\mathrm{Ind}(\mathrm{Ban}_H),\underline{s\mathrm{Sets}})$, denoted by $X$, we define the corresponding \textit{perfectization} $X^{1/p^\infty}$ of $X$ to be the corresponding functor such that for any $R\in \mathrm{sComm}\mathrm{Simp}\mathrm{Ind}(\mathrm{Ban}_H)$ we have that $X^{1/p^\infty}(R):=X(R^\flat)$ where we define the corresponding tilting \textit{Fontainisation} $R^\flat$ of $R$ to be the corresponding derived completion of:
\begin{align}
\varprojlim\{ ...\overset{\mathrm{Fro}}{\longrightarrow}	R \overset{\mathrm{Fro}}{\longrightarrow} R \overset{\mathrm{Fro}}{\longrightarrow} R\}.
\end{align}

\end{definition}

\indent In the situation of the corresponding monomorphically inductive Banach sets, we have the parallel definition:

\begin{definition}
For any object in $\infty-\mathrm{Fun}(\mathrm{sComm}\mathrm{Simp}\mathrm{Ind}^m(\mathrm{BanSets}_H),\underline{s\mathrm{Sets}})$, denoted by $X$, we define the corresponding \textit{perfectization} $X^{1/p^\infty}$ of $X$ to be the corresponding functor such that for any ring $R\in \mathrm{sComm}\mathrm{Simp}\mathrm{Ind}^m(\mathrm{BanSets}_H)$ we have that $X^{1/p^\infty}(R):=X(R^\flat)$ where we define the corresponding tilting \textit{Fontainisation} $R^\flat$ of $R$ to be the corresponding derived completion of:
\begin{align}
\varprojlim\{ ...\overset{\mathrm{Fro}}{\longrightarrow}	R \overset{\mathrm{Fro}}{\longrightarrow} R \overset{\mathrm{Fro}}{\longrightarrow} R\}.
\end{align}

\end{definition}

\begin{remark}
This is very general notion beyond the corresponding $(\infty,1)$-sheaves satisfying certain descent with respect to the derived rational localizations or more general homotopy Zariski topology as in \cite{11BK} and \cite{11BBBK}. 	
\end{remark}

\indent Now we follow \cite[Proposition 5.3]{11R} and \cite[Section 1.1]{11Dr1} to give the following discussion around the corresponding Witt crystalline Stack:

\begin{definition}
We define the corresponding \textit{Witt Crystalline Stack} $X_W$ of any $X$ over $H/\mathbb{F}_p$ in 
\begin{center}
$\infty-\mathrm{Fun}(\mathrm{sComm}\mathrm{Simp}\mathrm{Ind}(\mathrm{Ban}_H),\underline{s\mathrm{Sets}})$
\end{center}
 or 
\begin{center} 
$\infty-\mathrm{Fun}(\mathrm{sComm}\mathrm{Simp}\mathrm{Ind}^m(\mathrm{Ban}_H),\underline{s\mathrm{Sets}})$
\end{center}
to be the functor $(W(\pi_0(X)^\flat)_{\pi_0(X),\mathrm{crys}})_p^\wedge$. And we define the corresponding pre-crystals to be sheaves of $\mathcal{O}$-modules over this functors when we have that $X$ is an $\infty$-analytic stack with reasonable topology.	
\end{definition}

\begin{example}
For instance if we have that $X=\mathrm{Spa}^h(R)$ coming from the corresponding Bambozzi-Kremnizer spectrum of any Banach ring over $\mathbb{F}_p((t))$ as constructed in \cite{11BK}. Then we have that the corresponding functor is $(W(\pi_0(X)^\flat)_{\pi_0(X),\mathrm{crys}})_{p}^\wedge$ is now admitting structures coming from the corresponding homotopy Zariski topology from $X$.	
\end{example}

\begin{example}
For instance if we have that $X=\mathrm{Spec}(\mathbb{F}_p[[t]])$ coming from the corresponding object in the opposite category of $\mathbb{F}_p[[t]]$. Then we have that the corresponding functor is $(W(\pi_0(X)^\flat)_{\pi_0(X),\mathrm{crys}})_{p}^\wedge$ is now admitting structures coming from the corresponding homotopy Zariski topology from $X$, is just the same as the corresponding one in the algebraic setting constructed in \cite[Proposition 5.3]{11R} and \cite[Section 1.1]{11Dr1}.	
\end{example}

\subsection{Perfectization, Fontainisation and Robba Stacks}

\indent Now we contact \cite{11KL1} and \cite{11KL2} to look at the corresponding Robba Stacks. Now take any $X$ to be any $\infty$-analytic stack which admits structures of simplicial complete Bornological rings or ind-Fr\'echet structures, namely we have the corresponding complete bornological topology or ind-Fr\'echet topology on $\pi_0(X)$. We work over $H/\mathbb{F}_p$ as well. 

\begin{example}
For instance we take that $X=\mathrm{Spa}^h(R)$ coming from the corresponding Bambozzi-Kremnizer spectrum of any Banach ring over $\mathbb{F}_p((t))$ as constructed in \cite{11BK}. 
\end{example}

\begin{definition}
For any $\infty$-analytic stack $X$ as above, we consider the corresponding Witt vector functor $W_n(\pi_0(X)^\flat)$ and then consider
\begin{center}
 $\varinjlim_{n\rightarrow \infty}W_n(\pi_0(X)^\flat)$, 
\end{center} 
namely $W(\pi_0(X)^\flat)$, then we consider the ring $W(\pi_0(X)^\flat)[1/p]$. Then we can take the corresponding completion with respect to the Gauss norm $\|.\|_{\pi_0(X),[s,r]}$ coming from the corresponding norm on $\pi_0(X)$ with respect to some interval $[s,r]\in (0,\infty)$ as in \cite[Definition 4.1.1]{11KL2}:
\begin{align}
\widetilde{\Pi}(X)_{[s,r]}:=(W(\pi_0(X)^\flat)[1/p])^\wedge_{\|.\|_{\pi_0(X),[s,r]},\text{Fr\'e}}.	
\end{align}
Then following \cite[Definition 4.1.1]{11KL2} we consider the following:
\begin{align}
\widetilde{\Pi}(X)_{r}:=\varprojlim_{s>0}(W(\pi_0(X)^\flat)[1/p])^\wedge_{\|.\|_{\pi_0(X),[s,r]},\text{Fr\'e}}	
\end{align}
and 
\begin{align}
\widetilde{\Pi}(X)_{}:=\varinjlim_{r>0}\varprojlim_{s>0}(W(\pi_0(X)^\flat)[1/p])^\wedge_{\|.\|_{\pi_0(X),[s,r]},\text{Fr\'e}}.	
\end{align}
We call these \textit{Robba stacks}.
\end{definition}
 
\begin{example}
In the situation where $X$ is some $\infty$-analytic stack carrying the corresponding sheaves of simplicial Banach rings (namely not in general bornological or ind-Fr\'echet) we have that the finite projective modules over the three Robba stacks (for by enough intervals) carrying semilinear Frobenius action which realizes the corresponding isomorphisms by Frobenius pullbacks are equivalent. In the noetherian setting we have the same holds for locally finite presented sheaves as well. This is the main results of \cite[Theorem 4.6.1]{11KL2}.
\end{example}

\

This chapter is based on the following paper, where the author of this dissertation is the main author:
\begin{itemize}
\item Tong, Xin. "Topologization and Functional Analytification I: Intrinsic Morphisms of Commutative Algebras." arXiv preprint arXiv:2102.10766 (2021).
\end{itemize}

\newpage\chapter{Topologization and Functional Analytification II: $\infty$-Categorical Motivic Constructions for Homotopical Contexts}

\newpage

%\chapter{Introduction}

\section{Introduction}

\indent In our previous paper \cite{12XT1} on the intrinsic morphisms of the corresponding topological, bornological and ind-Fr\'echet rings, we discussed some generalization of the corresponding contexts of \cite{12Huber1} and \cite{12Huber2}, along two directions. One is along some answer to a question of Kedlaya in \cite[Appendix 5]{12Ked1} namely the corresponding affinoid morphisms, and the other is along some stacky understanding of these such as in \cite{12Dr1}, \cite{12Dr2}, \cite{12R}.\\

\indent Under some certain foundation, \cite{12Pau1} actually studied many parallel considerations of our project. The parallel  considerations are the corresponding ind-Banach analytic spaces, the corresponding \'etale and pro-\'etale sites of ind-Banach analytic spaces, the corresponding de Rham stacks of such analytic spaces, global de Rham cohomology and comparison over ind-Banach spaces, and algebraic derived de Rham complex over the corresponding ind-Banach spaces. And the corresponding descent over the corresponding derived ind-Banach spaces. We will consider the corresponding topologization and functional analytification of the cohomology theories (topological derived de Rham after \cite{12B1}, \cite{12Bei},  \cite{12GL}, \cite{12Ill1}, \cite{12Ill2} topological logarithmic derived de Rham after Gabber such as in \cite{12B1} and \cite{12O}, certainly functional analytification construction after \cite{12KL1} and \cite{12KL2} and etc).  We would like to mention that our project is also largely inspired by the corresponding development in \cite{12CS1}, \cite{12CS2}, \cite{12FS}, \cite{12Sch2} where the corresponding $v$-stacky consideration is extensively developed. Note that $v$-stacks are very significant analytic spaces, especially in the corresponding geometrization of local Langlands correspondence in \cite{12FS}, where complicated derived categories are constructed through deep foundations in \cite{12Sch2} and \cite{12CS1}, \cite{12CS2} on condensed sets. One certainly believes that the corresponding foundations by using normed sets in \cite{12BBBK}, \cite{12BBK}, \cite{BBM}, \cite{12BK} and \cite{KKM} will also have potential applications to certainly motivic and functorial constructions in analytic geometry which are parallel to those given in \cite{M} by using the foundation from \cite{12CS1}, \cite{12CS2}. Note that the programs in \cite{12B1}, \cite{12BMS} and \cite{12BS} have already indicated that working with the corresponding correct simplicial spaces in the analytic setting is not only a pure generalization but also a very significant point around even the non simplicial analytic spaces especially in the corresponding singular situations.\\

\indent Then we follow \cite{12An1}, \cite{12An2}, \cite{12B1}, \cite{12Bei}, \cite{12BS}, \cite{12G1}, \cite{12GL}, \cite{12Ill1}, \cite{12Ill2}, \cite{12Qui}  to extend the corresponding discussion to certainly spaces, we mainly have focused on the corresponding rigid analytic spaces, pseudorigid spaces in some explicit way. And we then extend the corresponding discussion to the corresponding derived setting. The corresponding derived setting is literally following \cite{12An1}, \cite{12An2}, \cite{12B1}, \cite{12B2}, \cite{12Bei}, \cite{12BMS}, \cite{12BS}, \cite{12G1}, \cite{12GL}, \cite{12Ill1}, \cite{12Ill2}, \cite{12Qui} and the corresponding derived logarithmic setting is literally following \cite{12B1}, \cite{12Ko1} and \cite{12O}\footnote{The corresponding \cite{12Ko1} has already mentioned the corresponding derived logarithmic prismatic cohomology.}. \\

\indent We promote the construction from Bhatt, Illusie, Guo, Morrow, Scholze, Gabber \cite{12B1}, \cite{12G1}, \cite{12Ill1}, \cite{12Ill2}, \cite{12BMS} \cite{12O} in the context of topological derived de Rham complexes and topological derived logarithmic de Rham complexes, the construction from Nicolaus-Scholze \cite{12NS} in the context of derived THH, TP and TC on the level of $\mathbb{E}_1$-rings , the construction from Kedlaya-Liu in the context of the derived Robba rings and derived Frobenius sheaves \cite{12KL1} and \cite{12KL2}, the construction from Bhatt-Scholze, Koshikawa and in the context of derived prismatic cohomology and derived logarithmic prismatic cohomology \cite{12BS} and \cite{12Ko1} to the level of $\infty$-categorical functional analytic level after Lurie \cite{12Lu1}, \cite{12Lu2} in the stable $\infty$-cateogory of derived $J$-complete simplicial commutative algebras, Bambozzi-Ben-Bassat-Kremnizer \cite{12BBBK}, Ben-Bassat-Mukherjee \cite{BBM}, Bambozzi-Kremnizer \cite{12BK}, Clausen-Scholze \cite{12CS1} \cite{12CS2} and Kelly-Kremnizer-Mukherjee \cite{KKM} in the stable $\infty$-cateogory of simplicial functional analytic commutative algebras.\\

\indent Furthermore we promote the construction from Bhatt, Illusie, Guo, Morrow, Scholze, Gabber \cite{12B1}, \cite{12G1}, \cite{12Ill1}, \cite{12Ill2}, \cite{BMS2} \cite{12O} in the context of topological derived de Rham complexes and topological derived logarithmic de Rham complexes, the construction from Nicolaus-Scholze \cite{12NS} in the context of derived THH, TP and TC on the level of $\mathbb{E}_1$-rings, the construction from Kedlaya-Liu in the context of the derived Robba rings and derived Frobenius sheaves \cite{12KL1} and \cite{12KL2}, the construction from Bhatt-Scholze, Koshikawa and in the context of derived prismatic cohomology and derived logarithmic prismatic cohomology \cite{12BS} and \cite{12Ko1} not only to the level of $\infty$-categorical functional analytic level after \cite{12Lu1}, \cite{12Lu2}, but also to the corresponding $(\infty,1)$-ringed toposes level after Lurie \cite{12Lu1}, \cite{12Lu2} in the $\infty$-category of $\infty$-ringed toposes, Bambozzi-Ben-Bassat-Kremnizer \cite{12BBBK}, Ben-Bassat-Mukherjee \cite{BBM}, Bambozzi-Kremnizer \cite{12BK}, Clausen-Scholze \cite{12CS1} \cite{12CS2} and Kelly-Kremnizer-Mukherjee \cite{KKM} in the stable $\infty$-cateogory of $\infty$-functional analytic ringed toposes.\\

\indent We then work in the noncommutative setting after \cite{Kon1}, \cite{Ta}, \cite{KR1} and \cite{KR2}, with some philosophy rooted in some noncommtative motives and the corresponding nonabelian applications in noncommutative analytic geometry in the derived sense, and the noncommutative analogues of the corresponding Riemann hypothesis and the corresponding Tamagawa number conjectures, and so on. The issue is certainly that the usual Frobenius map looks somehow strange, which reminds us of the fact that actually we need to consider really large objects such as the corresponding Topological Hochschild Homologies and the corresponding nearby objects, in order to define the corresponding analogues of the corresponding prismatic cohomology through THH and nearby objects, the corresponding noncommutative 'perfectoidizations'. Here we choose to consider \cite{12NS} in order to apply the constructions to certain $\infty$-rings, which we will call them Fukaya-Kato analytifications from \cite{12FK} as a noncommutative analog of the constructions in \cite{BBM}.\\

\indent We first consider the following list of motivic constructions in this paper which we hope to establish in very general topological and Banach setting, where one can actually believe that the constructions are somehow parallel and in some sense equivalent even on the $\infty$-categorical level\footnote{We are actually talking about the spaces which might not be smooth.}:

\begin{setting}\mbox{}\\
A. Derived Topological de Rham complexes and Derived Topological Logarithmic de Rham complexes of simplicial derived $I$-complete rings and over pro-\'etale site, after \cite{12B1}, \cite{12B2}, \cite{12Bei}, \cite{12BMS}, \cite{12BS}, \cite{12G1}, \cite{12GL}, \cite{12Ill1}, \cite{12Ill2}, \cite{12O};\\
B. $\mathcal{O}\mathrm{B}_\mathrm{dR}$-sheaves and $\mathcal{O}\mathrm{B}_\mathrm{dR,log}$-sheaves over general analytic adic spaces\footnote{These are not necessarily $p$-adic Tate. But as long as $p$ is topological nilpotent we do have a basis consisting of perfectoid subdomains as in \cite[Theorem 2.9.9, Remark 2.9.10]{12Ked1}.} after \cite{12DLLZ2}, \cite{12Sch2};\\
C. $\varphi$-$\widetilde{C}_X$-sheaves, relative $B$-pairs over general analytic adic spaces after \cite{12KL1} and \cite{12KL2};\\
D. Derived Prismatic Cohomology and Derived Logarithmic Prismatic Cohomology, after \cite{12BS}, \cite{12Ko1}.\\
\end{setting}

\begin{remark}
One could also consider the corresponding derived $I$-complete version\footnote{In fact we are taking derived $I$-completion of the spectra instead of the $I$-completion of the spectra in $\infty$-category as in  \cite[Chapter 7.3]{12Lu2}.} of the corresponding left Kan extended THH and HH $\mathbb{E}_1$-ring spectra of derived $I$-topological $\mathbb{E}_1$-ring spectra as in \cite{12NS} and \cite{KKM}. In certain situation, this should be able to be compared to the constructions above.
\end{remark}

\newpage

\section{Preliminary}

\subsection{The Rings in Algebraic Topology}

We consider some very general topological rings. We will consider the corresponding Huber's adic rings, namely we allow the corresponding rings to contain some adic open subrings which have then the corresponding linear topology. And we consider more general simplicial topological rings which are more general adic in some derived sense.

\begin{setting}
We consider the corresponding category $C_{\infty,\mathrm{E}_\infty}$ of all the $\mathbb{E}_\infty$ rings\footnote{These are the corresponding ring objects in the algebraic topology, see \cite[Chapter 7]{12Lu1}. So we want to consider the corresponding interesting objects in classical algebraic topology such as in \cite{12MP} and \cite{12N}.}. We then consider the corresponding category $C_{\infty,\mathrm{E}_\infty,\mathrm{rat},I}$ of all the $\mathbb{E}_\infty$ rings which have subrings being derived $I$-adically complete\footnote{In fact we are taking derived $I$-completion of the spectra instead of the $I$-completion of the spectra in $\infty$-category as in  \cite[Chapter 7.3]{12Lu2}.}, where $I$ is finite generating set coming from elements in $\pi_0$. And we have the smaller category $C_{\infty,\mathrm{E}_\infty,\mathrm{int},I}$ of all the all the $\mathbb{E}_\infty$ rings being derived $I$-adically complete.	
\end{setting}

\begin{example}
We have many interesting ring spectra from classical algebraic topology. One can further takes the corresponding derived $p$-completion to achieve many interesting objects\footnote{In fact we are taking derived $I$-completion of the spectra instead of the $I$-completion of the spectra in $\infty$-category as in \cite[Chapter 7.3]{12Lu2}. This should be different in more general setting from the consideration in algebraic topology and more general $\infty$-category theory.}. For more discussion see \cite[Chapter 5-13]{12MP} and \cite{12N}.	
\end{example}

\begin{example}
As also discussed in \cite[Chapter 14-19]{12MP} on the level of just model categories or more general construction in \cite[Chapter 1.3]{12Lu1} one considers the stable homotopy category $D(R)$ of some ring spectra, namely the corresponding $\infty$-enhancement of the corresponding classical triangulated categories. One considers the corresponding derived complete objects in this $\infty$-category, then one will have some interesting $\mathbb{E}_\infty$-spectra. 	
\end{example}

\begin{example}
We then have the example that is the corresponding rigid analytic affinoids in \cite[Definition 4.1]{12Ta}. They have open subrings which are actually derived $p$-adically complete.	
\end{example}

\begin{example}
We then have the example that is the corresponding pseudorigid analytic affinoids as in \cite[Definition 3.1]{12Bel1},  \cite{12Bel2} and \cite[Definition 4.1]{12L}. They have open subrings which are actually derived $t$-adically complete for some specific pseudouniformizer $t$.
\end{example}

\begin{example}
We now fix a bounded morphism of simplicial adic rings $A\rightarrow B$ over $A^*$ where $A^*$ contains a corresponding ring of definition $A^*_0$ which is complete with respect to the $(p,I)$-topology and we assume that $A$ is adic and we assume that $(A^*_0,I)$ is a prism in \cite{12BS} namely we at least require that the corresponding $\delta$-structure on the corresponding ring will induce the map $\varphi(.):=.^p+p\delta(.)$ such that we have the situation where $p\in (I,\varphi(I))$. For $A$ or $B$ respectively we assume this contains a subring $A_0$ or $B_0$ (over $A_0^*$) respectively such that we have $A_0$ or $B_0$ respectively is derived complete with respect to the corresponding derived $(p,I)$-topology and we assume that $B=B_0[1/f,f\in I]$ (same for $A$). All the adic rings are assumed to be open mapping. We use the notation $d$ to denote a corresponding primitive element as in \cite[Section 2.3]{12BS} for $A^*$. 
\end{example}

\indent We also consider the corresponding simplicial prelog rings as in \cite[Chapter 6]{12B1}, and also consider the corresponding application of such algebraic derived constructions to some derived $I$-adic rings carrying the corresponding prelog structures such as in \cite[Chapter 2]{12DLLZ1}.

\begin{setting}
We consider the corresponding category $C_{\infty,\mathrm{E}_\infty,\mathrm{prelog}}$ of all the $\mathbb{E}_\infty$ rings carrying the corresponding prelog structures. We then consider the corresponding category $C_{\infty,\mathrm{E}_\infty,\mathrm{rat},I,\mathrm{prelog}}$ of all the $\mathbb{E}_\infty$ rings carrying prelog structures which have subrings being derived $I$-adically complete\footnote{In fact again we are taking derived $I$-completion of the spectra instead of the $I$-completion of the spectra in $\infty$-category as in  \cite[Chapter 7.3]{12Lu2}.}, where $I$ is finite generating set coming from elements in $\pi_0$. And we have the smaller category $C_{\infty,\mathrm{E}_\infty,\mathrm{int},I,\mathrm{prelog}}$ of all the all the $\mathbb{E}_\infty$ rings carrying prelog structures being derived $I$-adically complete.	
\end{setting}

\begin{example}
The first example is the corresponding rigid analytic affinoids in \cite[Definition 4.1]{12Ta}. They have open subrings which are actually derived $p$-adically complete. Then one adds the corresponding prelog structures from \cite[Chapter 2]{12DLLZ1}. 
\end{example}

\begin{example}
The second example is the corresponding pseudorigid analytic affinoids as in \cite[Definition 3.1]{12Bel1}, \cite{12Bel2} and \cite[Definition 4.1]{12L}. They have open subrings which are actually derived $t$-adically complete for some specific pseudouniformizer $t$. Then one adds the corresponding prelog structures from \cite[Chapter 2]{12DLLZ1}. 
\end{example}

\begin{example}
We now fix a bounded morphism of logarithmic simplicial adic rings $(A,M)\rightarrow (B,N)$ over $A^*$ where $A^*$ contains a corresponding ring of definition $A^*_0$ which is complete with respect to the $(p,I)$-topology and we assume that $A$ is adic and we assume that $(A^*_0,I)$ is a prism namely we at least require that the corresponding $\delta$-structure on the corresponding ring will induce the map $\varphi(.):=.^p+p\delta(.)$ such that we have the situation where $p\in (I,\varphi(I))$. For $(A,M)$ or $(B,N_0)$ respectively we assume this contains a subring $(A_0,M_0)$ or $(B_0,N_0)$ (over $A_0^*$) respectively such that we have $(A_0,M_0)$ or $(B_0,N_0)$ respectively is derived complete with respect to the corresponding derived $(p,I)$-topology and we assume that $B=B_0[1/f,f\in I]$ (same for $A$). All the adic rings are assumed to be open mapping. We use the notation $d$ to denote a corresponding primitive element as in \cite[Section 2.3]{12BS} for $A^*$. \\
\end{example}

\newpage

\section{Notations on $\infty$-Categories of $\infty$-Rings}

\begin{center}
\begin{tabularx}{\linewidth}{lX}
Notation & Description\\
\hline
$C_{\infty,\mathrm{E}_\infty}$ &$(\infty,1)$-category of $\mathbb{E}_\infty$ rings.\\
$C_{\infty,\mathrm{E}_\infty,\mathrm{rat},I}$ &$(\infty,1)$-category of $\mathbb{E}_\infty$ rings in the sense that the rings contain some derived rings rationally.\\
$C_{\infty,\mathrm{E}_\infty,\mathrm{int},I}$ &$(\infty,1)$-category of $\mathbb{E}_\infty$ rings which are derived adic rings.\\
$C_{\infty,\mathrm{E}_\infty,\text{prelog}}$ &$(\infty,1)$-category of $\mathbb{E}_\infty$ logarithmic rings.\\
$C_{\infty,\mathrm{E}_\infty,\mathrm{rat},I,\text{prelog}}$ &    $(\infty,1)$-category of $\mathbb{E}_\infty$ logarithmic rings in the sense that the rings contain some derived rings rationally.\\
$C_{\infty,\mathrm{E}_\infty,\mathrm{int},I,\text{prelog}}$ &    $(\infty,1)$-category of $\mathbb{E}_\infty$ logarithmic rings which are derived adic rings.\\

$\mathrm{Object}_{\mathrm{E}_\infty\mathrm{commutativealgebra},\mathrm{Simplicial}}(\mathrm{IndSNorm}_R)$& $(\infty,1)$-category of $\mathbb{E}_\infty$ commutative algebra objects in the $(\infty,1)$-category of $\mathbb{E}_\infty$ ind-seminormed modules over $R$. \\
$\mathrm{Object}_{\mathrm{E}_\infty\mathrm{commutativealgebra},\mathrm{Simplicial}}(\mathrm{Ind}^m\mathrm{SNorm}_R)$& $(\infty,1)$-category of $\mathbb{E}_\infty$ commutative algebra objects in the $(\infty,1)$-category of $\mathbb{E}_\infty$ monomorphic ind-seminormed modules over $R$. \\
$\mathrm{Object}_{\mathrm{E}_\infty\mathrm{commutativealgebra},\mathrm{Simplicial}}(\mathrm{IndNorm}_R)$& $(\infty,1)$-category of $\mathbb{E}_\infty$ commutative algebra objects in the $(\infty,1)$-category of $\mathbb{E}_\infty$ ind-normed modules over $R$. \\
$\mathrm{Object}_{\mathrm{E}_\infty\mathrm{commutativealgebra},\mathrm{Simplicial}}(\mathrm{Ind}^m\mathrm{Norm}_R)$& $(\infty,1)$-category of $\mathbb{E}_\infty$ commutative algebra objects in the $(\infty,1)$-category of $\mathbb{E}_\infty$ monomorphic ind-normed modules over $R$. \\
$\mathrm{Object}_{\mathrm{E}_\infty\mathrm{commutativealgebra},\mathrm{Simplicial}}(\mathrm{IndBan}_R)$& $(\infty,1)$-category of $\mathbb{E}_\infty$ commutative algebra objects in the $(\infty,1)$-category of $\mathbb{E}_\infty$ ind-Banach modules over $R$. \\
$\mathrm{Object}_{\mathrm{E}_\infty\mathrm{commutativealgebra},\mathrm{Simplicial}}(\mathrm{Ind}^m\mathrm{Ban}_R)$& $(\infty,1)$-category of $\mathbb{E}_\infty$ commutative algebra objects in the $(\infty,1)$-category of $\mathbb{E}_\infty$ monomorphic ind-Banach modules over $R$. \\
$\mathrm{Object}_{\mathrm{E}_\infty\mathrm{commutativealgebra},\mathrm{Simplicial}}(\mathrm{IndSNorm}_{\mathbb{F}_1})$& $(\infty,1)$-category of $\mathbb{E}_\infty$ commutative algebra objects in the $(\infty,1)$-category of $\mathbb{E}_\infty$ ind-seminormed sets over $\mathbb{F}_1$. \\
$\mathrm{Object}_{\mathrm{E}_\infty\mathrm{commutativealgebra},\mathrm{Simplicial}}(\mathrm{Ind}^m\mathrm{SNorm}_{\mathbb{F}_1})$& $(\infty,1)$-category of $\mathbb{E}_\infty$ commutative algebra objects in the $(\infty,1)$-category of $\mathbb{E}_\infty$ monomorphic ind-seminormed sets over $\mathbb{F}_1$. \\
$\mathrm{Object}_{\mathrm{E}_\infty\mathrm{commutativealgebra},\mathrm{Simplicial}}(\mathrm{IndNorm}_{\mathbb{F}_1})$& $(\infty,1)$-category of $\mathbb{E}_\infty$ commutative algebra objects in the $(\infty,1)$-category of $\mathbb{E}_\infty$ ind-normed sets over $\mathbb{F}_1$. \\
$\mathrm{Object}_{\mathrm{E}_\infty\mathrm{commutativealgebra},\mathrm{Simplicial}}(\mathrm{Ind}^m\mathrm{Norm}_{\mathbb{F}_1})$ & $(\infty,1)$-category of $\mathbb{E}_\infty$ commutative algebra objects in the $(\infty,1)$-category of $\mathbb{E}_\infty$ monomorphic ind-normed sets over $\mathbb{F}_1$. \\
$\mathrm{Object}_{\mathrm{E}_\infty\mathrm{commutativealgebra},\mathrm{Simplicial}}(\mathrm{IndBan}_{\mathbb{F}_1})$& $(\infty,1)$-category of $\mathbb{E}_\infty$ commutative algebra objects in the $(\infty,1)$-category of $\mathbb{E}_\infty$ ind-Banach sets over $\mathbb{F}_1$. \\
$\mathrm{Object}_{\mathrm{E}_\infty\mathrm{commutativealgebra},\mathrm{Simplicial}}(\mathrm{Ind}^m\mathrm{Ban}_{\mathbb{F}_1})$& $(\infty,1)$-category of $\mathbb{E}_\infty$ commutative algebra objects in the $(\infty,1)$-category of $\mathbb{E}_\infty$ monomorphic ind-Banach sets over $\mathbb{F}_1$. \\

 $\mathrm{Object}_{\mathrm{E}_\infty\mathrm{commutativealgebra},\mathrm{Simplicial}}(\mathrm{IndSNorm}_R)^\text{prelog}$& $(\infty,1)$-category of $\mathbb{E}_\infty$ commutative algebra objects in the $(\infty,1)$-category of $\mathbb{E}_\infty$ ind-seminormed modules over $R$ carrying the corresponding logarithmic structures. \\
$\mathrm{Object}_{\mathrm{E}_\infty\mathrm{commutativealgebra},\mathrm{Simplicial}}(\mathrm{Ind}^m\mathrm{SNorm}_R)^\text{prelog}$& $(\infty,1)$-category of $\mathbb{E}_\infty$ commutative algebra objects in the $(\infty,1)$-category of $\mathbb{E}_\infty$ monomorphic ind-seminormed modules over $R$ carrying the corresponding logarithmic structures. \\
$\mathrm{Object}_{\mathrm{E}_\infty\mathrm{commutativealgebra},\mathrm{Simplicial}}(\mathrm{IndNorm}_R)^\text{prelog}$& $(\infty,1)$-category of $\mathbb{E}_\infty$ commutative algebra objects in the $(\infty,1)$-category of $\mathbb{E}_\infty$ ind-normed modules over $R$ carrying the corresponding logarithmic structures. \\
$\mathrm{Object}_{\mathrm{E}_\infty\mathrm{commutativealgebra},\mathrm{Simplicial}}(\mathrm{Ind}^m\mathrm{Norm}_R)^\text{prelog}$& $(\infty,1)$-category of $\mathbb{E}_\infty$ commutative algebra objects in the $(\infty,1)$-category of $\mathbb{E}_\infty$ monomorphic ind-normed modules over $R$ carrying the corresponding logarithmic structures. \\
$\mathrm{Object}_{\mathrm{E}_\infty\mathrm{commutativealgebra},\mathrm{Simplicial}}(\mathrm{IndBan}_R)^\text{prelog}$& $(\infty,1)$-category of $\mathbb{E}_\infty$ commutative algebra objects in the $(\infty,1)$-category of $\mathbb{E}_\infty$ ind-Banach modules over $R$ carrying the corresponding logarithmic structures. \\
$\mathrm{Object}_{\mathrm{E}_\infty\mathrm{commutativealgebra},\mathrm{Simplicial}}(\mathrm{Ind}^m\mathrm{Ban}_R)^\text{prelog}$& $(\infty,1)$-category of $\mathbb{E}_\infty$ commutative algebra objects in the $(\infty,1)$-category of $\mathbb{E}_\infty$ monomorphic ind-Banach modules over $R$ carrying the corresponding logarithmic structures. \\
$\mathrm{Object}_{\mathrm{E}_\infty\mathrm{commutativealgebra},\mathrm{Simplicial}}(\mathrm{IndSNorm}_{\mathbb{F}_1})^\text{prelog}$& $(\infty,1)$-category of $\mathbb{E}_\infty$ commutative algebra objects in the $(\infty,1)$-category of $\mathbb{E}_\infty$ ind-seminormed sets over $\mathbb{F}_1$ carrying the corresponding logarithmic structures. \\
$\mathrm{Object}_{\mathrm{E}_\infty\mathrm{commutativealgebra},\mathrm{Simplicial}}(\mathrm{Ind}^m\mathrm{SNorm}_{\mathbb{F}_1})^\text{prelog}$& $(\infty,1)$-category of $\mathbb{E}_\infty$ commutative algebra objects in the $(\infty,1)$-category of $\mathbb{E}_\infty$ monomorphic ind-seminormed sets over $\mathbb{F}_1$ carrying the corresponding logarithmic structures. \\
$\mathrm{Object}_{\mathrm{E}_\infty\mathrm{commutativealgebra},\mathrm{Simplicial}}(\mathrm{IndNorm}_{\mathbb{F}_1})^\text{prelog}$& $(\infty,1)$-category of $\mathbb{E}_\infty$ commutative algebra objects in the $(\infty,1)$-category of $\mathbb{E}_\infty$ ind-normed sets over $\mathbb{F}_1$ carrying the corresponding logarithmic structures. \\
$\mathrm{Object}_{\mathrm{E}_\infty\mathrm{commutativealgebra},\mathrm{Simplicial}}(\mathrm{Ind}^m\mathrm{Norm}_{\mathbb{F}_1})^\text{prelog}$& $(\infty,1)$-category of $\mathbb{E}_\infty$ commutative algebra objects in the $(\infty,1)$-category of $\mathbb{E}_\infty$ monomorphic ind-normed sets over $\mathbb{F}_1$ carrying the corresponding logarithmic structures. \\
$\mathrm{Object}_{\mathrm{E}_\infty\mathrm{commutativealgebra},\mathrm{Simplicial}}(\mathrm{IndBan}_{\mathbb{F}_1})^\text{prelog}$& $(\infty,1)$-category of $\mathbb{E}_\infty$ commutative algebra objects in the $(\infty,1)$-category of $\mathbb{E}_\infty$ ind-Banach sets over $\mathbb{F}_1$ carrying the corresponding logarithmic structures. \\
$\mathrm{Object}_{\mathrm{E}_\infty\mathrm{commutativealgebra},\mathrm{Simplicial}}(\mathrm{Ind}^m\mathrm{Ban}_{\mathbb{F}_1})^\text{prelog}$& $(\infty,1)$-category of $\mathbb{E}_\infty$ commutative algebra objects in the $(\infty,1)$-category of $\mathbb{E}_\infty$ monomorphic ind-Banach sets over $\mathbb{F}_1$ carrying the corresponding logarithmic structures. \\

 $\mathrm{Object}_{\mathrm{E}_\infty\mathrm{commutativealgebra},\mathrm{Simplicial}}(\mathrm{IndSNorm}_R)^{\square}$ & $(\infty,1)$-category of $\mathbb{E}_\infty$ commutative algebra objects in the $(\infty,1)$-category of $\mathbb{E}_\infty$ ind-seminormed modules over $R$, generated by the corresponding formal series rings over $R$. Here $\square$ denotes \text{smoothformalseriesclosure}.  \\
$\mathrm{Object}_{\mathrm{E}_\infty\mathrm{commutativealgebra},\mathrm{Simplicial}}(\mathrm{Ind}^m\mathrm{SNorm}_R)^{\square}$& $(\infty,1)$-category of $\mathbb{E}_\infty$ commutative algebra objects in the $(\infty,1)$-category of $\mathbb{E}_\infty$ monomorphic ind-seminormed modules over $R$, generated by the corresponding formal series rings over $R$. \\
$\mathrm{Object}_{\mathrm{E}_\infty\mathrm{commutativealgebra},\mathrm{Simplicial}}(\mathrm{IndNorm}_R)^{\square}$& $(\infty,1)$-category of $\mathbb{E}_\infty$ commutative algebra objects in the $(\infty,1)$-category of $\mathbb{E}_\infty$ ind-normed modules over $R$, generated by the corresponding formal series rings over $R$, generated by the corresponding formal series rings over $R$. \\
$\mathrm{Object}_{\mathrm{E}_\infty\mathrm{commutativealgebra},\mathrm{Simplicial}}(\mathrm{Ind}^m\mathrm{Norm}_R)^{\square}$& $(\infty,1)$-category of $\mathbb{E}_\infty$ commutative algebra objects in the $(\infty,1)$-category of $\mathbb{E}_\infty$ monomorphic ind-normed modules over $R$, generated by the corresponding formal series rings over $R$. \\
$\mathrm{Object}_{\mathrm{E}_\infty\mathrm{commutativealgebra},\mathrm{Simplicial}}(\mathrm{IndBan}_R)^{\square}$& $(\infty,1)$-category of $\mathbb{E}_\infty$ commutative algebra objects in the $(\infty,1)$-category of $\mathbb{E}_\infty$ ind-Banach modules over $R$, generated by the corresponding formal series rings over $R$. \\
$\mathrm{Object}_{\mathrm{E}_\infty\mathrm{commutativealgebra},\mathrm{Simplicial}}(\mathrm{Ind}^m\mathrm{Ban}_R)^{\square}$& $(\infty,1)$-category of $\mathbb{E}_\infty$ commutative algebra objects in the $(\infty,1)$-category of $\mathbb{E}_\infty$ monomorphic ind-Banach modules over $R$, generated by the corresponding formal series rings over $R$. \\
$\mathrm{Object}_{\mathrm{E}_\infty\mathrm{commutativealgebra},\mathrm{Simplicial}}(\mathrm{IndSNorm}_{\mathbb{F}_1})^{\square}$& $(\infty,1)$-category of $\mathbb{E}_\infty$ commutative algebra objects in the $(\infty,1)$-category of $\mathbb{E}_\infty$ ind-seminormed sets over $\mathbb{F}_1$, generated by the corresponding formal series rings over $\mathbb{F}_1$. \\
$\mathrm{Object}_{\mathrm{E}_\infty\mathrm{commutativealgebra},\mathrm{Simplicial}}(\mathrm{Ind}^m\mathrm{SNorm}_{\mathbb{F}_1})^{\square}$& $(\infty,1)$-category of $\mathbb{E}_\infty$ commutative algebra objects in the $(\infty,1)$-category of $\mathbb{E}_\infty$ monomorphic ind-seminormed sets over $\mathbb{F}_1$, generated by the corresponding formal series rings over $\mathbb{F}_1$. \\
$\mathrm{Object}_{\mathrm{E}_\infty\mathrm{commutativealgebra},\mathrm{Simplicial}}(\mathrm{IndNorm}_{\mathbb{F}_1})^{\square}$& $(\infty,1)$-category of $\mathbb{E}_\infty$ commutative algebra objects in the $(\infty,1)$-category of $\mathbb{E}_\infty$ ind-normed sets over $\mathbb{F}_1$, generated by the corresponding formal series rings over $\mathbb{F}_1$. \\
$\mathrm{Object}_{\mathrm{E}_\infty\mathrm{commutativealgebra},\mathrm{Simplicial}}(\mathrm{Ind}^m\mathrm{Norm}_{\mathbb{F}_1})^{\square}$& $(\infty,1)$-category of $\mathbb{E}_\infty$ commutative algebra objects in the $(\infty,1)$-category of $\mathbb{E}_\infty$ monomorphic ind-normed sets over $\mathbb{F}_1$, generated by the corresponding formal series rings over $\mathbb{F}_1$. \\
$\mathrm{Object}_{\mathrm{E}_\infty\mathrm{commutativealgebra},\mathrm{Simplicial}}(\mathrm{IndBan}_{\mathbb{F}_1})^{\square}$& $(\infty,1)$-category of $\mathbb{E}_\infty$ commutative algebra objects in the $(\infty,1)$-category of $\mathbb{E}_\infty$ ind-Banach sets over $\mathbb{F}_1$, generated by the corresponding formal series rings over $\mathbb{F}_1$. \\
$\mathrm{Object}_{\mathrm{E}_\infty\mathrm{commutativealgebra},\mathrm{Simplicial}}(\mathrm{Ind}^m\mathrm{Ban}_{\mathbb{F}_1})^{\square}$& $(\infty,1)$-category of $\mathbb{E}_\infty$ commutative algebra objects in the $(\infty,1)$-category of $\mathbb{E}_\infty$ monomorphic ind-Banach sets over $\mathbb{F}_1$, generated by the corresponding formal series rings over $\mathbb{F}_1$. \\

$\mathrm{Object}_{\mathrm{E}_\infty\mathrm{commutativealgebra},\mathrm{Simplicial}}(\mathrm{IndSNorm}_R)^{\square,\text{prelog}}$& $(\infty,1)$-category of $\mathbb{E}_\infty$ commutative algebra objects in the $(\infty,1)$-category of $\mathbb{E}_\infty$ ind-seminormed modules over $R$ carrying the corresponding logarithmic structures, generated by the corresponding formal series rings over $R$. \\
$\mathrm{Object}_{\mathrm{E}_\infty\mathrm{commutativealgebra},\mathrm{Simplicial}}(\mathrm{Ind}^m\mathrm{SNorm}_R)^{\square,\text{prelog}}$& $(\infty,1)$-category of $\mathbb{E}_\infty$ commutative algebra objects in the $(\infty,1)$-category of $\mathbb{E}_\infty$ monomorphic ind-seminormed modules over $R$ carrying the corresponding logarithmic structures, generated by the corresponding formal series rings over $R$. \\
$\mathrm{Object}_{\mathrm{E}_\infty\mathrm{commutativealgebra},\mathrm{Simplicial}}(\mathrm{IndNorm}_R)^{\square,\text{prelog}}$& $(\infty,1)$-category of $\mathbb{E}_\infty$ commutative algebra objects in the $(\infty,1)$-category of $\mathbb{E}_\infty$ ind-normed modules over $R$ carrying the corresponding logarithmic structures, generated by the corresponding formal series rings over $R$. \\
$\mathrm{Object}_{\mathrm{E}_\infty\mathrm{commutativealgebra},\mathrm{Simplicial}}(\mathrm{Ind}^m\mathrm{Norm}_R)^{\square,\text{prelog}}$& $(\infty,1)$-category of $\mathbb{E}_\infty$ commutative algebra objects in the $(\infty,1)$-category of $\mathbb{E}_\infty$ monomorphic ind-normed modules over $R$ carrying the corresponding logarithmic structures, generated by the corresponding formal series rings over $R$. \\
$\mathrm{Object}_{\mathrm{E}_\infty\mathrm{commutativealgebra},\mathrm{Simplicial}}(\mathrm{IndBan}_R)^{\square,\text{prelog}}$& $(\infty,1)$-category of $\mathbb{E}_\infty$ commutative algebra objects in the $(\infty,1)$-category of $\mathbb{E}_\infty$ ind-Banach modules over $R$ carrying the corresponding logarithmic structures, generated by the corresponding formal series rings over $R$. \\
$\mathrm{Object}_{\mathrm{E}_\infty\mathrm{commutativealgebra},\mathrm{Simplicial}}(\mathrm{Ind}^m\mathrm{Ban}_R)^{\square,\text{prelog}}$& $(\infty,1)$-category of $\mathbb{E}_\infty$ commutative algebra objects in the $(\infty,1)$-category of $\mathbb{E}_\infty$ monomorphic ind-Banach modules over $R$ carrying the corresponding logarithmic structures, generated by the corresponding formal series rings over $R$. \\
$\mathrm{Object}_{\mathrm{E}_\infty\mathrm{commutativealgebra},\mathrm{Simplicial}}(\mathrm{IndSNorm}_{\mathbb{F}_1})^{\square,\text{prelog}}$& $(\infty,1)$-category of $\mathbb{E}_\infty$ commutative algebra objects in the $(\infty,1)$-category of $\mathbb{E}_\infty$ ind-seminormed sets over $\mathbb{F}_1$ carrying the corresponding logarithmic structures, generated by the corresponding formal series rings over ${\mathbb{F}_1}$. \\
$\mathrm{Object}_{\mathrm{E}_\infty\mathrm{commutativealgebra},\mathrm{Simplicial}}(\mathrm{Ind}^m\mathrm{SNorm}_{\mathbb{F}_1})^{\square,\text{prelog}}$& $(\infty,1)$-category of $\mathbb{E}_\infty$ commutative algebra objects in the $(\infty,1)$-category of $\mathbb{E}_\infty$ monomorphic ind-seminormed sets over $\mathbb{F}_1$ carrying the corresponding logarithmic structures, generated by the corresponding formal series rings over ${\mathbb{F}_1}$. \\
$\mathrm{Object}_{\mathrm{E}_\infty\mathrm{commutativealgebra},\mathrm{Simplicial}}(\mathrm{IndNorm}_{\mathbb{F}_1})^{\square,\text{prelog}}$& $(\infty,1)$-category of $\mathbb{E}_\infty$ commutative algebra objects in the $(\infty,1)$-category of $\mathbb{E}_\infty$ ind-normed sets over $\mathbb{F}_1$ carrying the corresponding logarithmic structures, generated by the corresponding formal series rings over ${\mathbb{F}_1}$. \\
$\mathrm{Object}_{\mathrm{E}_\infty\mathrm{commutativealgebra},\mathrm{Simplicial}}(\mathrm{Ind}^m\mathrm{Norm}_{\mathbb{F}_1})^{\square,\text{prelog}}$& $(\infty,1)$-category of $\mathbb{E}_\infty$ commutative algebra objects in the $(\infty,1)$-category of $\mathbb{E}_\infty$ monomorphic ind-normed sets over $\mathbb{F}_1$ carrying the corresponding logarithmic structures, generated by the corresponding formal series rings over ${\mathbb{F}_1}$. \\
$\mathrm{Object}_{\mathrm{E}_\infty\mathrm{commutativealgebra},\mathrm{Simplicial}}(\mathrm{IndBan}_{\mathbb{F}_1})^{\square,\text{prelog}}$& $(\infty,1)$-category of $\mathbb{E}_\infty$ commutative algebra objects in the $(\infty,1)$-category of $\mathbb{E}_\infty$ ind-Banach sets over $\mathbb{F}_1$ carrying the corresponding logarithmic structures, generated by the corresponding formal series rings over ${\mathbb{F}_1}$. \\
$\mathrm{Object}_{\mathrm{E}_\infty\mathrm{commutativealgebra},\mathrm{Simplicial}}(\mathrm{Ind}^m\mathrm{Ban}_{\mathbb{F}_1})^{\square,\text{prelog}}$& $(\infty,1)$-category of $\mathbb{E}_\infty$ commutative algebra objects in the $(\infty,1)$-category of $\mathbb{E}_\infty$ monomorphic ind-Banach sets over $\mathbb{F}_1$ carrying the corresponding logarithmic structures, generated by the corresponding formal series rings over ${\mathbb{F}_1}$. \\

\end{tabularx}
\end{center}

\begin{remark}
The generating process above is by adding all the colimits in the homotopy sense which are assumed to be sifted. $R$ will be Banach and carrying the corresponding $p$-adic topology.
\end{remark}

\newpage
\section{Notations on $\infty$-Categories of Commutative $\infty$-Ringed Toposes}

\indent We change the notations for $\infty$-ringed toposes slightly, therefore let us start from the rings.\\

\noindent Rings:\\

\noindent $\mathrm{Ind}^\text{smoothformalseriesclosure}\mathrm{Commutativealgebra}_{\mathrm{simplicial}}(\mathrm{Ind}\mathrm{Seminormed}_R)$: $(\infty,1)$-category of simplicial commutative algebra objects in the $\infty$-category of colimit completion of seminormed modules over a general Banach ring $R$, generated by the corresponding power series rings by taking colimits in the homotopical sense. \\ 
\noindent $\mathrm{Ind}^\text{smoothformalseriesclosure}\mathrm{Commutativealgebra}_{\mathrm{simplicial}}(\mathrm{Ind}^m\mathrm{Seminormed}_R)$: $(\infty,1)$-category of simplicial commutative algebra objects in the $\infty$-category of monomorphic colimit completion of seminormed modules over a general Banach ring $R$, generated by the corresponding power series rings by taking colimits in the homotopical sense. \\ 
\noindent $\mathrm{Ind}^\text{smoothformalseriesclosure}\mathrm{Commutativealgebra}_{\mathrm{simplicial}}(\mathrm{Ind}\mathrm{Normed}_R)$: $(\infty,1)$-category of simplicial commutative algebra objects in the $\infty$-category of colimit completion of normed modules over a general Banach ring $R$, generated by the corresponding power series rings by taking colimits in the homotopical sense. \\ 
\noindent $\mathrm{Ind}^\text{smoothformalseriesclosure}\mathrm{Commutativealgebra}_{\mathrm{simplicial}}(\mathrm{Ind}^m\mathrm{Normed}_R)$: $(\infty,1)$-category of simplicial commutative algebra objects in the $\infty$-category of monomorphic colimit completion of normed modules over a general Banach ring $R$, generated by the corresponding power series rings by taking colimits in the homotopical sense. \\ 
\noindent $\mathrm{Ind}^\text{smoothformalseriesclosure}\mathrm{Commutativealgebra}_{\mathrm{simplicial}}(\mathrm{Ind}\mathrm{Banach}_R)$: $(\infty,1)$-category of simplicial commutative algebra objects in the $\infty$-category of colimit completion of Banach modules over a general Banach ring $R$, generated by the corresponding power series rings by taking colimits in the homotopical sense. \\ 
\noindent $\mathrm{Ind}^\text{smoothformalseriesclosure}\mathrm{Commutativealgebra}_{\mathrm{simplicial}}(\mathrm{Ind}^m\mathrm{Banach}_R)$: $(\infty,1)$-category of simplicial commutative algebra objects in the $\infty$-category of monomorphic colimit completion of Banach modules over a general Banach ring $R$, generated by the corresponding power series rings by taking colimits in the homotopical sense. \\

\noindent $\mathrm{Ind}^\text{smoothformalseriesclosure}\mathrm{Commutativealgebra}_{\mathrm{simplicial}}(\mathrm{Ind}\mathrm{Seminormed}_{\mathbb{F}_1})$: $(\infty,1)$-category of simplicial commutative algebra objects in the $\infty$-category of colimit completion of seminormed sets over a general Banach ring $\mathbb{F}_1$, generated by the corresponding power series rings by taking colimits in the homotopical sense. \\ 
\noindent $\mathrm{Ind}^\text{smoothformalseriesclosure}\mathrm{Commutativealgebra}_{\mathrm{simplicial}}(\mathrm{Ind}^m\mathrm{Seminormed}_{\mathbb{F}_1})$: $(\infty,1)$-category of simplicial commutative algebra objects in the $\infty$-category of monomorphic colimit completion of seminormed sets over a general Banach ring $\mathbb{F}_1$, generated by the corresponding power series rings by taking colimits in the homotopical sense. \\ 
\noindent $\mathrm{Ind}^\text{smoothformalseriesclosure}\mathrm{Commutativealgebra}_{\mathrm{simplicial}}(\mathrm{Ind}\mathrm{Normed}_{\mathbb{F}_1})$: $(\infty,1)$-category of simplicial commutative algebra objects in the $\infty$-category of colimit completion of normed sets over a general Banach ring $\mathbb{F}_1$, generated by the corresponding power series rings by taking colimits in the homotopical sense. \\ 
\noindent $\mathrm{Ind}^\text{smoothformalseriesclosure}\mathrm{Commutativealgebra}_{\mathrm{simplicial}}(\mathrm{Ind}^m\mathrm{Normed}_{\mathbb{F}_1})$: $(\infty,1)$-category of simplicial commutative algebra objects in the $\infty$-category of monomorphic colimit completion of normed sets over a general Banach ring $\mathbb{F}_1$, generated by the corresponding power series rings by taking colimits in the homotopical sense. \\ 
\noindent $\mathrm{Ind}^\text{smoothformalseriesclosure}\mathrm{Commutativealgebra}_{\mathrm{simplicial}}(\mathrm{Ind}\mathrm{Banach}_{\mathbb{F}_1})$: $(\infty,1)$-category of simplicial commutative algebra objects in the $\infty$-category of colimit completion of Banach sets over a general Banach ring $\mathbb{F}_1$, generated by the corresponding power series rings by taking colimits in the homotopical sense. \\ 
\noindent $\mathrm{Ind}^\text{smoothformalseriesclosure}\mathrm{Commutativealgebra}_{\mathrm{simplicial}}(\mathrm{Ind}^m\mathrm{Banach}_{\mathbb{F}_1})$: $(\infty,1)$-category of simplicial commutative algebra objects in the $\infty$-category of monomorphic colimit completion of Banach sets over a general Banach ring $\mathbb{F}_1$, generated by the corresponding power series rings by taking colimits in the homotopical sense. \\

\noindent Prestacks:\\

\noindent $\infty-\mathrm{Prestack}_{\mathrm{Commutativealgebra}_{\mathrm{simplicial}}(\mathrm{Ind}\mathrm{Seminormed}_?)^\mathrm{opposite},\mathrm{Grotopology,homotopyepimorphism}}$: $\infty$-presheaves into $\infty$-groupoid over corresponding opposite category carrying the corresponding Grothendieck topology, and we will mainly consider the corresponding homotopy epimorphisms. $?=R,\mathbb{F}_1$.\\
\noindent $\infty-\mathrm{Prestack}_{\mathrm{Commutativealgebra}_{\mathrm{simplicial}}(\mathrm{Ind}^m\mathrm{Seminormed}_?)^\mathrm{opposite},\mathrm{Grotopology,homotopyepimorphism}}$: $\infty$-presheaves into $\infty$-groupoid over corresponding opposite category carrying the corresponding Grothendieck topology, and we will mainly consider the corresponding homotopy epimorphisms. $?=R,\mathbb{F}_1$.\\
\noindent $\infty-\mathrm{Prestack}_{\mathrm{Commutativealgebra}_{\mathrm{simplicial}}(\mathrm{Ind}\mathrm{Normed}_?)^\mathrm{opposite},\mathrm{Grotopology,homotopyepimorphism}}$: $\infty$-presheaves into $\infty$-groupoid over corresponding opposite category carrying the corresponding Grothendieck topology, and we will mainly consider the corresponding homotopy epimorphisms. $?=R,\mathbb{F}_1$.\\
\noindent $\infty-\mathrm{Prestack}_{\mathrm{Commutativealgebra}_{\mathrm{simplicial}}(\mathrm{Ind}^m\mathrm{Normed}_?)^\mathrm{opposite},\mathrm{Grotopology,homotopyepimorphism}}$: $\infty$-presheaves into $\infty$-groupoid over corresponding opposite category carrying the corresponding Grothendieck topology, and we will mainly consider the corresponding homotopy epimorphisms. $?=R,\mathbb{F}_1$.\\
\noindent $\infty-\mathrm{Prestack}_{\mathrm{Commutativealgebra}_{\mathrm{simplicial}}(\mathrm{Ind}\mathrm{Banach}_?)^\mathrm{opposite},\mathrm{Grotopology,homotopyepimorphism}}$: $\infty$-presheaves into $\infty$-groupoid over corresponding opposite category carrying the corresponding Grothendieck topology, and we will mainly consider the corresponding homotopy epimorphisms. $?=R,\mathbb{F}_1$.\\
\noindent $\infty-\mathrm{Prestack}_{\mathrm{Commutativealgebra}_{\mathrm{simplicial}}(\mathrm{Ind}^m\mathrm{Banach}_?)^\mathrm{opposite},\mathrm{Grotopology,homotopyepimorphism}}$: $\infty$-presheaves into $\infty$-groupoid over corresponding opposite category carrying the corresponding Grothendieck topology, and we will mainly consider the corresponding homotopy epimorphisms. $?=R,\mathbb{F}_1$.\\

\noindent Stacks:\\
 
 \noindent $\infty-\mathrm{Stack}_{\mathrm{Commutativealgebra}_{\mathrm{simplicial}}(\mathrm{Ind}\mathrm{Seminormed}_?)^\mathrm{opposite},\mathrm{Grotopology,homotopyepimorphism}}$: $\infty$-sheaves into $\infty$-groupoid over corresponding opposite category carrying the corresponding Grothendieck topology, and we will mainly consider the corresponding homotopy epimorphisms. $?=R,\mathbb{F}_1$.\\
\noindent $\infty-\mathrm{Stack}_{\mathrm{Commutativealgebra}_{\mathrm{simplicial}}(\mathrm{Ind}^m\mathrm{Seminormed}_?)^\mathrm{opposite},\mathrm{Grotopology,homotopyepimorphism}}$: $\infty$-sheaves into $\infty$-groupoid over corresponding opposite category carrying the corresponding Grothendieck topology, and we will mainly consider the corresponding homotopy epimorphisms. $?=R,\mathbb{F}_1$.\\
\noindent $\infty-\mathrm{Stack}_{\mathrm{Commutativealgebra}_{\mathrm{simplicial}}(\mathrm{Ind}\mathrm{Normed}_?)^\mathrm{opposite},\mathrm{Grotopology,homotopyepimorphism}}$: $\infty$-sheaves into $\infty$-groupoid over corresponding opposite category carrying the corresponding Grothendieck topology, and we will mainly consider the corresponding homotopy epimorphisms. $?=R,\mathbb{F}_1$.\\
\noindent $\infty-\mathrm{Stack}_{\mathrm{Commutativealgebra}_{\mathrm{simplicial}}(\mathrm{Ind}^m\mathrm{Normed}_?)^\mathrm{opposite},\mathrm{Grotopology,homotopyepimorphism}}$: $\infty$-sheaves into $\infty$-groupoid over corresponding opposite category carrying the corresponding Grothendieck topology, and we will mainly consider the corresponding homotopy epimorphisms. $?=R,\mathbb{F}_1$.\\
\noindent $\infty-\mathrm{Stack}_{\mathrm{Commutativealgebra}_{\mathrm{simplicial}}(\mathrm{Ind}\mathrm{Banach}_?)^\mathrm{opposite},\mathrm{Grotopology,homotopyepimorphism}}$: $\infty$-sheaves into $\infty$-groupoid over corresponding opposite category carrying the corresponding Grothendieck topology, and we will mainly consider the corresponding homotopy epimorphisms. $?=R,\mathbb{F}_1$.\\
\noindent $\infty-\mathrm{Stack}_{\mathrm{Commutativealgebra}_{\mathrm{simplicial}}(\mathrm{Ind}^m\mathrm{Banach}_?)^\mathrm{opposite},\mathrm{Grotopology,homotopyepimorphism}}$: $\infty$-sheaves into $\infty$-groupoid over corresponding opposite category carrying the corresponding Grothendieck topology, and we will mainly consider the corresponding homotopy epimorphisms. $?=R,\mathbb{F}_1$.\\

\noindent Ringed Toposes: \\
 
 \noindent $\infty-\mathrm{Toposes}^{\mathrm{ringed},\mathrm{commutativealgebra}_{\mathrm{simplicial}}(\mathrm{Ind}\mathrm{Seminormed}_?)}_{\mathrm{Commutativealgebra}_{\mathrm{simplicial}}(\mathrm{Ind}\mathrm{Seminormed}_?)^\mathrm{opposite},\mathrm{Grotopology,homotopyepimorphism}}$: $\infty$-sheaves into $\infty$-groupoid over corresponding opposite category carrying the corresponding Grothendieck topology, and we will mainly consider the corresponding homotopy epimorphisms. $?=R,\mathbb{F}_1$. And we assume the stack carries $\infty$-ringed toposes structure. \\
\noindent $\infty-\mathrm{Toposes}^{\mathrm{ringed},\mathrm{Commutativealgebra}_{\mathrm{simplicial}}(\mathrm{Ind}^m\mathrm{Seminormed}_?)}_{\mathrm{Commutativealgebra}_{\mathrm{simplicial}}(\mathrm{Ind}^m\mathrm{Seminormed}_?)^\mathrm{opposite},\mathrm{Grotopology,homotopyepimorphism}}$: $\infty$-sheaves into $\infty$-groupoid over corresponding opposite category carrying the corresponding Grothendieck topology, and we will mainly consider the corresponding homotopy epimorphisms. $?=R,\mathbb{F}_1$. And we assume the stack carries $\infty$-ringed toposes structure.\\
\noindent $\infty-\mathrm{Toposes}^{\mathrm{ringed},\mathrm{Commutativealgebra}_{\mathrm{simplicial}}(\mathrm{Ind}\mathrm{Normed}_?)}_{\mathrm{Commutativealgebra}_{\mathrm{simplicial}}(\mathrm{Ind}\mathrm{Normed}_?)^\mathrm{opposite},\mathrm{Grotopology,homotopyepimorphism}}$: $\infty$-sheaves into $\infty$-groupoid over corresponding opposite category carrying the corresponding Grothendieck topology, and we will mainly consider the corresponding homotopy epimorphisms. $?=R,\mathbb{F}_1$. And we assume the stack carries $\infty$-ringed toposes structure.\\
\noindent $\infty-\mathrm{Toposes}^{\mathrm{ringed},\mathrm{Commutativealgebra}_{\mathrm{simplicial}}(\mathrm{Ind}^m\mathrm{Normed}_?)}_{\mathrm{Commutativealgebra}_{\mathrm{simplicial}}(\mathrm{Ind}^m\mathrm{Normed}_?)^\mathrm{opposite},\mathrm{Grotopology,homotopyepimorphism}}$: $\infty$-sheaves into $\infty$-groupoid over corresponding opposite category carrying the corresponding Grothendieck topology, and we will mainly consider the corresponding homotopy epimorphisms. $?=R,\mathbb{F}_1$. And we assume the stack carries $\infty$-ringed toposes structure.\\
\noindent $\infty-\mathrm{Toposes}^{\mathrm{ringed},\mathrm{Commutativealgebra}_{\mathrm{simplicial}}(\mathrm{Ind}\mathrm{Banach}_?)}_{\mathrm{Commutativealgebra}_{\mathrm{simplicial}}(\mathrm{Ind}\mathrm{Banach}_?)^\mathrm{opposite},\mathrm{Grotopology,homotopyepimorphism}}$: $\infty$-sheaves into $\infty$-groupoid over corresponding opposite category carrying the corresponding Grothendieck topology, and we will mainly consider the corresponding homotopy epimorphisms. $?=R,\mathbb{F}_1$. And we assume the stack carries $\infty$-ringed toposes structure.\\
\noindent $\infty-\mathrm{Toposes}^{\mathrm{ringed},\mathrm{Commutativealgebra}_{\mathrm{simplicial}}(\mathrm{Ind}^m\mathrm{Banach}_?)}_{\mathrm{Commutativealgebra}_{\mathrm{simplicial}}(\mathrm{Ind}^m\mathrm{Banach}_?)^\mathrm{opposite},\mathrm{Grotopology,homotopyepimorphism}}$: $\infty$-sheaves into $\infty$-groupoid over corresponding opposite category carrying the corresponding Grothendieck topology, and we will mainly consider the corresponding homotopy epimorphisms. $?=R,\mathbb{F}_1$. And we assume the stack carries $\infty$-ringed toposes structure.\\

 \noindent $\mathrm{Proj}^\text{smoothformalseriesclosure}\infty-\mathrm{Toposes}^{\mathrm{ringed},\mathrm{commutativealgebra}_{\mathrm{simplicial}}(\mathrm{Ind}\mathrm{Seminormed}_?)}_{\mathrm{Commutativealgebra}_{\mathrm{simplicial}}(\mathrm{Ind}\mathrm{Seminormed}_?)^\mathrm{opposite},\mathrm{Grotopology,homotopyepimorphism}}$: $\infty$-sheaves into $\infty$-groupoid over corresponding opposite category carrying the corresponding Grothendieck topology, and we will mainly consider the corresponding homotopy epimorphisms. $?=R,\mathbb{F}_1$. And we assume the stack carries $\infty$-ringed toposes structure by specializing a ring object $\mathcal{O}$. \\
\noindent $\mathrm{Proj}^\text{smoothformalseriesclosure}\infty-\mathrm{Toposes}^{\mathrm{ringed},\mathrm{Commutativealgebra}_{\mathrm{simplicial}}(\mathrm{Ind}^m\mathrm{Seminormed}_?)}_{\mathrm{Commutativealgebra}_{\mathrm{simplicial}}(\mathrm{Ind}^m\mathrm{Seminormed}_?)^\mathrm{opposite},\mathrm{Grotopology,homotopyepimorphism}}$: $\infty$-sheaves into $\infty$-groupoid over corresponding opposite category carrying the corresponding Grothendieck topology, and we will mainly consider the corresponding homotopy epimorphisms. $?=R,\mathbb{F}_1$. And we assume the stack carries $\infty$-ringed toposes structure by specializing a ring object $\mathcal{O}$.\\
\noindent $\mathrm{Proj}^\text{smoothformalseriesclosure}\infty-\mathrm{Toposes}^{\mathrm{ringed},\mathrm{Commutativealgebra}_{\mathrm{simplicial}}(\mathrm{Ind}\mathrm{Normed}_?)}_{\mathrm{Commutativealgebra}_{\mathrm{simplicial}}(\mathrm{Ind}\mathrm{Normed}_?)^\mathrm{opposite},\mathrm{Grotopology,homotopyepimorphism}}$: $\infty$-sheaves into $\infty$-groupoid over corresponding opposite category carrying the corresponding Grothendieck topology, and we will mainly consider the corresponding homotopy epimorphisms. $?=R,\mathbb{F}_1$. And we assume the stack carries $\infty$-ringed toposes structure by specializing a ring object $\mathcal{O}$.\\
\noindent $\mathrm{Proj}^\text{smoothformalseriesclosure}\infty-\mathrm{Toposes}^{\mathrm{ringed},\mathrm{Commutativealgebra}_{\mathrm{simplicial}}(\mathrm{Ind}^m\mathrm{Normed}_?)}_{\mathrm{Commutativealgebra}_{\mathrm{simplicial}}(\mathrm{Ind}^m\mathrm{Normed}_?)^\mathrm{opposite},\mathrm{Grotopology,homotopyepimorphism}}$: $\infty$-sheaves into $\infty$-groupoid over corresponding opposite category carrying the corresponding Grothendieck topology, and we will mainly consider the corresponding homotopy epimorphisms. $?=R,\mathbb{F}_1$. And we assume the stack carries $\infty$-ringed toposes structure by specializing a ring object $\mathcal{O}$.\\
\noindent $\mathrm{Proj}^\text{smoothformalseriesclosure}\infty-\mathrm{Toposes}^{\mathrm{ringed},\mathrm{Commutativealgebra}_{\mathrm{simplicial}}(\mathrm{Ind}\mathrm{Banach}_?)}_{\mathrm{Commutativealgebra}_{\mathrm{simplicial}}(\mathrm{Ind}\mathrm{Banach}_?)^\mathrm{opposite},\mathrm{Grotopology,homotopyepimorphism}}$: $\infty$-sheaves into $\infty$-groupoid over corresponding opposite category carrying the corresponding Grothendieck topology, and we will mainly consider the corresponding homotopy epimorphisms. $?=R,\mathbb{F}_1$. And we assume the stack carries $\infty$-ringed toposes structure by specializing a ring object $\mathcal{O}$.\\
\noindent $\mathrm{Proj}^\text{smoothformalseriesclosure}\infty-\mathrm{Toposes}^{\mathrm{ringed},\mathrm{Commutativealgebra}_{\mathrm{simplicial}}(\mathrm{Ind}^m\mathrm{Banach}_?)}_{\mathrm{Commutativealgebra}_{\mathrm{simplicial}}(\mathrm{Ind}^m\mathrm{Banach}_?)^\mathrm{opposite},\mathrm{Grotopology,homotopyepimorphism}}$: $\infty$-sheaves into $\infty$-groupoid over corresponding opposite category carrying the corresponding Grothendieck topology, and we will mainly consider the corresponding homotopy epimorphisms. $?=R,\mathbb{F}_1$. And we assume the stack carries $\infty$-ringed toposes structure by specializing a ring object $\mathcal{O}$.\\

\begin{remark}
In $p$-adic Hodge theory, one usually will need to construct presheaves $M$ out of from the $\infty$-ring object $\mathcal{O}$.	Also one can define the corresponding pre-$\infty$-ringed Toposes, we will not continue provide the corresponding group of notations in the parallel way. 
\end{remark}

\noindent $\infty$-Quasicoherent Sheaves of Functional Analytic Modules over Ringed Toposes $\sharp=\mathrm{Seminormed},\mathrm{Normed},\mathrm{Banach}$: \\
 
 \noindent $\mathrm{Ind}\mathrm{\sharp Quasicoherent}_{\infty-\mathrm{Toposes}^{\mathrm{ringed},\mathrm{commutativealgebra}_{\mathrm{simplicial}}(\mathrm{Ind}\mathrm{Seminormed}_?)}_{\mathrm{Commutativealgebra}_{\mathrm{simplicial}}(\mathrm{Ind}\mathrm{Seminormed}_?)^\mathrm{opposite},\mathrm{Grotopology,homotopyepimorphism}}}$: Colimits completion of $\infty$-Quasicoherent Sheaves of Functional Analytic Modules over $\infty$-sheaves into $\infty$-groupoid over corresponding opposite category carrying the corresponding Grothendieck topology, and we will mainly consider the corresponding homotopy epimorphisms. $?=R,\mathbb{F}_1$. And we assume the stack carries $\infty$-ringed toposes structure. \\
\noindent $\mathrm{Ind}\mathrm{\sharp Quasicoherent}_{\infty-\mathrm{Toposes}^{\mathrm{ringed},\mathrm{Commutativealgebra}_{\mathrm{simplicial}}(\mathrm{Ind}^m\mathrm{Seminormed}_?)}_{\mathrm{Commutativealgebra}_{\mathrm{simplicial}}(\mathrm{Ind}^m\mathrm{Seminormed}_?)^\mathrm{opposite},\mathrm{Grotopology,homotopyepimorphism}}}$: Colimits completion of $\infty$-Quasicoherent Sheaves of Functional Analytic Modules over $\infty$-sheaves into $\infty$-groupoid over corresponding opposite category carrying the corresponding Grothendieck topology, and we will mainly consider the corresponding homotopy epimorphisms. $?=R,\mathbb{F}_1$. And we assume the stack carries $\infty$-ringed toposes structure.\\
\noindent $\mathrm{Ind}\mathrm{\sharp Quasicoherent}_{\infty-\mathrm{Toposes}^{\mathrm{ringed},\mathrm{Commutativealgebra}_{\mathrm{simplicial}}(\mathrm{Ind}\mathrm{Normed}_?)}_{\mathrm{Commutativealgebra}_{\mathrm{simplicial}}(\mathrm{Ind}\mathrm{Normed}_?)^\mathrm{opposite},\mathrm{Grotopology,homotopyepimorphism}}}$: Colimits completion of $\infty$-Quasicoherent Sheaves of Functional Analytic Modules over $\infty$-sheaves into $\infty$-groupoid over corresponding opposite category carrying the corresponding Grothendieck topology, and we will mainly consider the corresponding homotopy epimorphisms. $?=R,\mathbb{F}_1$. And we assume the stack carries $\infty$-ringed toposes structure.\\
\noindent $\mathrm{Ind}\mathrm{\sharp Quasicoherent}_{\infty-\mathrm{Toposes}^{\mathrm{ringed},\mathrm{Commutativealgebra}_{\mathrm{simplicial}}(\mathrm{Ind}^m\mathrm{Normed}_?)}_{\mathrm{Commutativealgebra}_{\mathrm{simplicial}}(\mathrm{Ind}^m\mathrm{Normed}_?)^\mathrm{opposite},\mathrm{Grotopology,homotopyepimorphism}}}$: Colimits completion of $\infty$-Quasicoherent Sheaves of Functional Analytic Modules over $\infty$-sheaves into $\infty$-groupoid over corresponding opposite category carrying the corresponding Grothendieck topology, and we will mainly consider the corresponding homotopy epimorphisms. $?=R,\mathbb{F}_1$. And we assume the stack carries $\infty$-ringed toposes structure.\\
\noindent $\mathrm{Ind}\mathrm{\sharp Quasicoherent}_{\infty-\mathrm{Toposes}^{\mathrm{ringed},\mathrm{Commutativealgebra}_{\mathrm{simplicial}}(\mathrm{Ind}\mathrm{Banach}_?)}_{\mathrm{Commutativealgebra}_{\mathrm{simplicial}}(\mathrm{Ind}\mathrm{Banach}_?)^\mathrm{opposite},\mathrm{Grotopology,homotopyepimorphism}}}$: Colimits completion of $\infty$-Quasicoherent Sheaves of Functional Analytic Modules over $\infty$-sheaves into $\infty$-groupoid over corresponding opposite category carrying the corresponding Grothendieck topology, and we will mainly consider the corresponding homotopy epimorphisms. $?=R,\mathbb{F}_1$. And we assume the stack carries $\infty$-ringed toposes structure.\\
\noindent $\mathrm{Ind}\mathrm{\sharp Quasicoherent}_{\infty-\mathrm{Toposes}^{\mathrm{ringed},\mathrm{Commutativealgebra}_{\mathrm{simplicial}}(\mathrm{Ind}^m\mathrm{Banach}_?)}_{\mathrm{Commutativealgebra}_{\mathrm{simplicial}}(\mathrm{Ind}^m\mathrm{Banach}_?)^\mathrm{opposite},\mathrm{Grotopology,homotopyepimorphism}}}$: Colimits completion of $\infty$-Quasicoherent Sheaves of Functional Analytic Modules over $\infty$-sheaves into $\infty$-groupoid over corresponding opposite category carrying the corresponding Grothendieck topology, and we will mainly consider the corresponding homotopy epimorphisms. $?=R,\mathbb{F}_1$. And we assume the stack carries $\infty$-ringed toposes structure.\\

 \noindent $\mathrm{Ind}\mathrm{\sharp Quasicoherent}_{\mathrm{Ind}^\text{smoothformalseriesclosure}\infty-\mathrm{Toposes}^{\mathrm{ringed},\mathrm{commutativealgebra}_{\mathrm{simplicial}}(\mathrm{Ind}\mathrm{Seminormed}_?)}_{\mathrm{Commutativealgebra}_{\mathrm{simplicial}}(\mathrm{Ind}\mathrm{Seminormed}_?)^\mathrm{opposite},\mathrm{Grotopology,homotopyepimorphism}}}$: Colimits completion of $\infty$-Quasicoherent Sheaves of Functional Analytic Modules over $\infty$-sheaves into $\infty$-groupoid over corresponding opposite category carrying the corresponding Grothendieck topology, and we will mainly consider the corresponding homotopy epimorphisms. $?=R,\mathbb{F}_1$. And we assume the stack carries $\infty$-ringed toposes structure by specializing a ring object $\mathcal{O}$. The difference is that we consider the corresponding $\mathcal{O}$ living in the colimit completion closure. \\
\noindent $\mathrm{Ind}\mathrm{\sharp Quasicoherent}_{\mathrm{Ind}^\text{smoothformalseriesclosure}\infty-\mathrm{Toposes}^{\mathrm{ringed},\mathrm{Commutativealgebra}_{\mathrm{simplicial}}(\mathrm{Ind}^m\mathrm{Seminormed}_?)}_{\mathrm{Commutativealgebra}_{\mathrm{simplicial}}(\mathrm{Ind}^m\mathrm{Seminormed}_?)^\mathrm{opposite},\mathrm{Grotopology,homotopyepimorphism}}}$: Colimits completion of $\infty$-Quasicoherent Sheaves of Functional Analytic Modules over $\infty$-sheaves into $\infty$-groupoid over corresponding opposite category carrying the corresponding Grothendieck topology, and we will mainly consider the corresponding homotopy epimorphisms. $?=R,\mathbb{F}_1$. And we assume the stack carries $\infty$-ringed toposes structure by specializing a ring object $\mathcal{O}$. The difference is that we consider the corresponding $\mathcal{O}$ living in the colimit completion closure.\\
\noindent $\mathrm{Ind}\mathrm{\sharp Quasicoherent}_{\mathrm{Ind}^\text{smoothformalseriesclosure}\infty-\mathrm{Toposes}^{\mathrm{ringed},\mathrm{Commutativealgebra}_{\mathrm{simplicial}}(\mathrm{Ind}\mathrm{Normed}_?)}_{\mathrm{Commutativealgebra}_{\mathrm{simplicial}}(\mathrm{Ind}\mathrm{Normed}_?)^\mathrm{opposite},\mathrm{Grotopology,homotopyepimorphism}}}$: Colimits completion of $\infty$-Quasicoherent Sheaves of Functional Analytic Modules over $\infty$-sheaves into $\infty$-groupoid over corresponding opposite category carrying the corresponding Grothendieck topology, and we will mainly consider the corresponding homotopy epimorphisms. $?=R,\mathbb{F}_1$. And we assume the stack carries $\infty$-ringed toposes structure by specializing a ring object $\mathcal{O}$. The difference is that we consider the corresponding $\mathcal{O}$ living in the colimit completion closure.\\
\noindent $\mathrm{Ind}\mathrm{\sharp Quasicoherent}_{\mathrm{Ind}^\text{smoothformalseriesclosure}\infty-\mathrm{Toposes}^{\mathrm{ringed},\mathrm{Commutativealgebra}_{\mathrm{simplicial}}(\mathrm{Ind}^m\mathrm{Normed}_?)}_{\mathrm{Commutativealgebra}_{\mathrm{simplicial}}(\mathrm{Ind}^m\mathrm{Normed}_?)^\mathrm{opposite},\mathrm{Grotopology,homotopyepimorphism}}}$: Colimits completion of $\infty$-Quasicoherent Sheaves of Functional Analytic Modules over $\infty$-sheaves into $\infty$-groupoid over corresponding opposite category carrying the corresponding Grothendieck topology, and we will mainly consider the corresponding homotopy epimorphisms. $?=R,\mathbb{F}_1$. And we assume the stack carries $\infty$-ringed toposes structure by specializing a ring object $\mathcal{O}$. The difference is that we consider the corresponding $\mathcal{O}$ living in the colimit completion closure.\\
\noindent $\mathrm{Ind}\mathrm{\sharp Quasicoherent}_{\mathrm{Ind}^\text{smoothformalseriesclosure}\infty-\mathrm{Toposes}^{\mathrm{ringed},\mathrm{Commutativealgebra}_{\mathrm{simplicial}}(\mathrm{Ind}\mathrm{Banach}_?)}_{\mathrm{Commutativealgebra}_{\mathrm{simplicial}}(\mathrm{Ind}\mathrm{Banach}_?)^\mathrm{opposite},\mathrm{Grotopology,homotopyepimorphism}}}$: Colimits completion of $\infty$-Quasicoherent Sheaves of Functional Analytic Modules over $\infty$-sheaves into $\infty$-groupoid over corresponding opposite category carrying the corresponding Grothendieck topology, and we will mainly consider the corresponding homotopy epimorphisms. $?=R,\mathbb{F}_1$. And we assume the stack carries $\infty$-ringed toposes structure by specializing a ring object $\mathcal{O}$. The difference is that we consider the corresponding $\mathcal{O}$ living in the colimit completion closure.\\
\noindent $\mathrm{Ind}\mathrm{\sharp Quasicoherent}_{\mathrm{Ind}^\text{smoothformalseriesclosure}\infty-\mathrm{Toposes}^{\mathrm{ringed},\mathrm{Commutativealgebra}_{\mathrm{simplicial}}(\mathrm{Ind}^m\mathrm{Banach}_?)}_{\mathrm{Commutativealgebra}_{\mathrm{simplicial}}(\mathrm{Ind}^m\mathrm{Banach}_?)^\mathrm{opposite},\mathrm{Grotopology,homotopyepimorphism}}}$: Colimits completion of $\infty$-Quasicoherent Sheaves of Functional Analytic Modules over $\infty$-sheaves into $\infty$-groupoid over corresponding opposite category carrying the corresponding Grothendieck topology, and we will mainly consider the corresponding homotopy epimorphisms. $?=R,\mathbb{F}_1$. And we assume the stack carries $\infty$-ringed toposes structure by specializing a ring object $\mathcal{O}$. The difference is that we consider the corresponding $\mathcal{O}$ living in the colimit completion closure.\\

\newpage
\section{Notations on $\infty$-Categories of Noncommutative $\infty$-Ringed Toposes}

\

\indent We change the notations for $\infty$-ringed noncommutative toposes slightly, therefore let us start from the rings. The corresponding generators we will choose in order to take the corresponding homotopy colimit completion and the corresponding homotopy limit completion are Fukaya-Kato rings in \cite{12FK}:
\begin{align}
R\left<Z_1,...,Z_n\right>,n\geq 1,\\
R\left[[Z_1,...,Z_n\right]],n\geq 1,	
\end{align}
with $Z_1,...,Z_n$ are noncommuting free variables. This would be the specific completions of the polynomials:
\begin{align}
R\left[Z_1,...,Z_n\right],n\geq 1.\\	
\end{align}

\noindent The corresponding analogs of analytification functors from Ben-Bassat-Mukherjee \cite[Section 4.2]{BBM} are given in the following. First we consider the $\infty$-category of $\mathbb{E}_1$-rings from \cite[Proposition 7.1.4.18, as well as the discussion above Proposition 7.1.4.18 on page 1225]{12Lu2} which we denote it by $\mathrm{Noncommutative}_{\mathbb{E}_1,\mathrm{Simplicial}}$, then we consider the corresponding category of all the polynomial rings with free variables over $R$, which we denote it by $\mathrm{Polynomial}^\mathrm{free}_R$, then we have the corresponding fully faithful embedding:
\begin{align}
\mathrm{Polynomial}^\mathrm{free}_R\rightarrow \mathrm{Noncommutative}_{\mathbb{E}_1,\mathrm{Simplicial}}.	
\end{align}
Then take the corresponding completion for each:
\begin{align}
R\left[Z_1,...,Z_n\right],n\geq 1,\\	
\end{align}
we have the Fukaya-Kato adic ring:
\begin{align}
R\left<Z_1,...,Z_n\right>,n\geq 1,\\
R\left[[Z_1,...,Z_n\right]],n\geq 1.	
\end{align}
As in \cite[Section 4.2]{BBM}, this will give the process what we call smooth formal series analytification by taking into account the corresponding homotopy colimit completion:\\
\begin{align}
\mathrm{Ind}^{\mathrm{smoothformalseriesclosure}}\mathrm{Polynomial}^\mathrm{free}_R\rightarrow \mathrm{Noncommutativealgebra}_{\mathrm{simplicial}}(\mathrm{Ind}\mathrm{Seminormed}_R),\\
\mathrm{Ind}^{\mathrm{smoothformalseriesclosure}}\mathrm{Polynomial}^\mathrm{free}_R\rightarrow \mathrm{Noncommutativealgebra}_{\mathrm{simplicial}}(\mathrm{Ind}\mathrm{Normed}_R),\\
\mathrm{Ind}^{\mathrm{smoothformalseriesclosure}}\mathrm{Polynomial}^\mathrm{free}_R\rightarrow \mathrm{Noncommutativealgebra}_{\mathrm{simplicial}}(\mathrm{Ind}\mathrm{Banach}_R).\\	
\end{align}
\begin{remark}
\indent The left actually spans all the $\mathbb{E}_1$-algebra in our setting by regarding the free variable polynomials as tensor algebras over $R$ as explained in \cite[Proposition 7.1.4.18, as well as the discussion above Proposition 7.1.4.18 on page 1225]{12Lu1} in analogy of \cite[Proposition 7.1.4.20, as well as the discussion above Proposition 7.1.4.20]{12Lu1} in the commutative situation.
\end{remark}

\

\newpage

\noindent Noncommutative Rings:\\

\noindent $\mathrm{Ind}^\text{smoothformalseriesclosure}\mathrm{Noncommutativealgebra}_{\mathrm{simplicial}}(\mathrm{Ind}\mathrm{Seminormed}_R)$: $(\infty,1)$-category of simplicial noncommutative algebra objects in the $\infty$-category of colimit completion of seminormed modules over a general Banach ring $R$, generated by the corresponding power series rings by taking colimits in the homotopical sense. \\ 
\noindent $\mathrm{Ind}^\text{smoothformalseriesclosure}\mathrm{Noncommutativealgebra}_{\mathrm{simplicial}}(\mathrm{Ind}^m\mathrm{Seminormed}_R)$: $(\infty,1)$-category of simplicial noncommutative algebra objects in the $\infty$-category of monomorphic colimit completion of seminormed modules over a general Banach ring $R$, generated by the corresponding power series rings by taking colimits in the homotopical sense. \\ 
\noindent $\mathrm{Ind}^\text{smoothformalseriesclosure}\mathrm{Noncommutativealgebra}_{\mathrm{simplicial}}(\mathrm{Ind}\mathrm{Normed}_R)$: $(\infty,1)$-category of simplicial noncommutative algebra objects in the $\infty$-category of colimit completion of normed modules over a general Banach ring $R$, generated by the corresponding power series rings by taking colimits in the homotopical sense. \\ 
\noindent $\mathrm{Ind}^\text{smoothformalseriesclosure}\mathrm{Noncommutativealgebra}_{\mathrm{simplicial}}(\mathrm{Ind}^m\mathrm{Normed}_R)$: $(\infty,1)$-category of simplicial noncommutative algebra objects in the $\infty$-category of monomorphic colimit completion of normed modules over a general Banach ring $R$, generated by the corresponding power series rings by taking colimits in the homotopical sense. \\ 
\noindent $\mathrm{Ind}^\text{smoothformalseriesclosure}\mathrm{Noncommutativealgebra}_{\mathrm{simplicial}}(\mathrm{Ind}\mathrm{Banach}_R)$: $(\infty,1)$-category of simplicial noncommutative algebra objects in the $\infty$-category of colimit completion of Banach modules over a general Banach ring $R$, generated by the corresponding power series rings by taking colimits in the homotopical sense. \\ 
\noindent $\mathrm{Ind}^\text{smoothformalseriesclosure}\mathrm{Noncommutativealgebra}_{\mathrm{simplicial}}(\mathrm{Ind}^m\mathrm{Banach}_R)$: $(\infty,1)$-category of simplicial noncommutative algebra objects in the $\infty$-category of monomorphic colimit completion of Banach modules over a general Banach ring $R$, generated by the corresponding power series rings by taking colimits in the homotopical sense. \\

\noindent $\mathrm{Ind}^\text{smoothformalseriesclosure}\mathrm{Noncommutativealgebra}_{\mathrm{simplicial}}(\mathrm{Ind}\mathrm{Seminormed}_{\mathbb{F}_1})$: $(\infty,1)$-category of simplicial noncommutative algebra objects in the $\infty$-category of colimit completion of seminormed sets over a general Banach ring $\mathbb{F}_1$, generated by the corresponding power series rings by taking colimits in the homotopical sense. \\ 
\noindent $\mathrm{Ind}^\text{smoothformalseriesclosure}\mathrm{Noncommutativealgebra}_{\mathrm{simplicial}}(\mathrm{Ind}^m\mathrm{Seminormed}_{\mathbb{F}_1})$: $(\infty,1)$-category of simplicial noncommutative algebra objects in the $\infty$-category of monomorphic colimit completion of seminormed sets over a general Banach ring $\mathbb{F}_1$, generated by the corresponding power series rings by taking colimits in the homotopical sense. \\ 
\noindent $\mathrm{Ind}^\text{smoothformalseriesclosure}\mathrm{Noncommutativealgebra}_{\mathrm{simplicial}}(\mathrm{Ind}\mathrm{Normed}_{\mathbb{F}_1})$: $(\infty,1)$-category of simplicial noncommutative algebra objects in the $\infty$-category of colimit completion of normed sets over a general Banach ring $\mathbb{F}_1$, generated by the corresponding power series rings by taking colimits in the homotopical sense. \\ 
\noindent $\mathrm{Ind}^\text{smoothformalseriesclosure}\mathrm{Noncommutativealgebra}_{\mathrm{simplicial}}(\mathrm{Ind}^m\mathrm{Normed}_{\mathbb{F}_1})$: $(\infty,1)$-category of simplicial noncommutative algebra objects in the $\infty$-category of monomorphic colimit completion of normed sets over a general Banach ring $\mathbb{F}_1$, generated by the corresponding power series rings by taking colimits in the homotopical sense. \\ 
\noindent $\mathrm{Ind}^\text{smoothformalseriesclosure}\mathrm{Noncommutativealgebra}_{\mathrm{simplicial}}(\mathrm{Ind}\mathrm{Banach}_{\mathbb{F}_1})$: $(\infty,1)$-category of simplicial noncommutative algebra objects in the $\infty$-category of colimit completion of Banach sets over a general Banach ring $\mathbb{F}_1$, generated by the corresponding power series rings by taking colimits in the homotopical sense. \\ 
\noindent $\mathrm{Ind}^\text{smoothformalseriesclosure}\mathrm{Noncommutativealgebra}_{\mathrm{simplicial}}(\mathrm{Ind}^m\mathrm{Banach}_{\mathbb{F}_1})$: $(\infty,1)$-category of simplicial noncommutative algebra objects in the $\infty$-category of monomorphic colimit completion of Banach sets over a general Banach ring $\mathbb{F}_1$, generated by the corresponding power series rings by taking colimits in the homotopical sense. \\

\noindent Prestacks:\\

\noindent $\infty-\mathrm{Prestack}_{\mathrm{Noncommutativealgebra}_{\mathrm{simplicial}}(\mathrm{Ind}\mathrm{Seminormed}_?)^\mathrm{opposite},\mathrm{Grotopology,homotopyepimorphism}}$: $\infty$-presheaves into $\infty$-groupoid over corresponding opposite category carrying the corresponding Grothendieck topology, and we will mainly consider the corresponding homotopy epimorphisms. $?=R,\mathbb{F}_1$.\\
\noindent $\infty-\mathrm{Prestack}_{\mathrm{Noncommutativealgebra}_{\mathrm{simplicial}}(\mathrm{Ind}^m\mathrm{Seminormed}_?)^\mathrm{opposite},\mathrm{Grotopology,homotopyepimorphism}}$: $\infty$-presheaves into $\infty$-groupoid over corresponding opposite category carrying the corresponding Grothendieck topology, and we will mainly consider the corresponding homotopy epimorphisms. $?=R,\mathbb{F}_1$.\\
\noindent $\infty-\mathrm{Prestack}_{\mathrm{Noncommutativealgebra}_{\mathrm{simplicial}}(\mathrm{Ind}\mathrm{Normed}_?)^\mathrm{opposite},\mathrm{Grotopology,homotopyepimorphism}}$: $\infty$-presheaves into $\infty$-groupoid over corresponding opposite category carrying the corresponding Grothendieck topology, and we will mainly consider the corresponding homotopy epimorphisms. $?=R,\mathbb{F}_1$.\\
\noindent $\infty-\mathrm{Prestack}_{\mathrm{Noncommutativealgebra}_{\mathrm{simplicial}}(\mathrm{Ind}^m\mathrm{Normed}_?)^\mathrm{opposite},\mathrm{Grotopology,homotopyepimorphism}}$: $\infty$-presheaves into $\infty$-groupoid over corresponding opposite category carrying the corresponding Grothendieck topology, and we will mainly consider the corresponding homotopy epimorphisms. $?=R,\mathbb{F}_1$.\\
\noindent $\infty-\mathrm{Prestack}_{\mathrm{Noncommutativealgebra}_{\mathrm{simplicial}}(\mathrm{Ind}\mathrm{Banach}_?)^\mathrm{opposite},\mathrm{Grotopology,homotopyepimorphism}}$: $\infty$-presheaves into $\infty$-groupoid over corresponding opposite category carrying the corresponding Grothendieck topology, and we will mainly consider the corresponding homotopy epimorphisms. $?=R,\mathbb{F}_1$.\\
\noindent $\infty-\mathrm{Prestack}_{\mathrm{Noncommutativealgebra}_{\mathrm{simplicial}}(\mathrm{Ind}^m\mathrm{Banach}_?)^\mathrm{opposite},\mathrm{Grotopology,homotopyepimorphism}}$: $\infty$-presheaves into $\infty$-groupoid over corresponding opposite category carrying the corresponding Grothendieck topology, and we will mainly consider the corresponding homotopy epimorphisms. $?=R,\mathbb{F}_1$.\\

\noindent Stacks:\\
 
 \noindent $\infty-\mathrm{Stack}_{\mathrm{Noncommutativealgebra}_{\mathrm{simplicial}}(\mathrm{Ind}\mathrm{Seminormed}_?)^\mathrm{opposite},\mathrm{Grotopology,homotopyepimorphism}}$: $\infty$-sheaves into $\infty$-groupoid over corresponding opposite category carrying the corresponding Grothendieck topology, and we will mainly consider the corresponding homotopy epimorphisms. $?=R,\mathbb{F}_1$.\\
\noindent $\infty-\mathrm{Stack}_{\mathrm{Noncommutativealgebra}_{\mathrm{simplicial}}(\mathrm{Ind}^m\mathrm{Seminormed}_?)^\mathrm{opposite},\mathrm{Grotopology,homotopyepimorphism}}$: $\infty$-sheaves into $\infty$-groupoid over corresponding opposite category carrying the corresponding Grothendieck topology, and we will mainly consider the corresponding homotopy epimorphisms. $?=R,\mathbb{F}_1$.\\
\noindent $\infty-\mathrm{Stack}_{\mathrm{Noncommutativealgebra}_{\mathrm{simplicial}}(\mathrm{Ind}\mathrm{Normed}_?)^\mathrm{opposite},\mathrm{Grotopology,homotopyepimorphism}}$: $\infty$-sheaves into $\infty$-groupoid over corresponding opposite category carrying the corresponding Grothendieck topology, and we will mainly consider the corresponding homotopy epimorphisms. $?=R,\mathbb{F}_1$.\\
\noindent $\infty-\mathrm{Stack}_{\mathrm{Noncommutativealgebra}_{\mathrm{simplicial}}(\mathrm{Ind}^m\mathrm{Normed}_?)^\mathrm{opposite},\mathrm{Grotopology,homotopyepimorphism}}$: $\infty$-sheaves into $\infty$-groupoid over corresponding opposite category carrying the corresponding Grothendieck topology, and we will mainly consider the corresponding homotopy epimorphisms. $?=R,\mathbb{F}_1$.\\
\noindent $\infty-\mathrm{Stack}_{\mathrm{Noncommutativealgebra}_{\mathrm{simplicial}}(\mathrm{Ind}\mathrm{Banach}_?)^\mathrm{opposite},\mathrm{Grotopology,homotopyepimorphism}}$: $\infty$-sheaves into $\infty$-groupoid over corresponding opposite category carrying the corresponding Grothendieck topology, and we will mainly consider the corresponding homotopy epimorphisms. $?=R,\mathbb{F}_1$.\\
\noindent $\infty-\mathrm{Stack}_{\mathrm{Noncommutativealgebra}_{\mathrm{simplicial}}(\mathrm{Ind}^m\mathrm{Banach}_?)^\mathrm{opposite},\mathrm{Grotopology,homotopyepimorphism}}$: $\infty$-sheaves into $\infty$-groupoid over corresponding opposite category carrying the corresponding Grothendieck topology, and we will mainly consider the corresponding homotopy epimorphisms. $?=R,\mathbb{F}_1$.\\

\indent The following is a noncommutative analog of \cite[Definition 5.6]{12BBBK}:

\begin{definition}
Here a homotopy epimorphism is defined to be such a morphism $A\rightarrow B$ such that $B\overline{\otimes}_A B^\mathrm{opp}\rightarrow B$ reflects isomorphism in the homotopy category.
\end{definition}

\noindent Ringed Toposes: \\
 
 \noindent $\infty-\mathrm{Toposes}^{\mathrm{ringed},\mathrm{Noncommutativealgebra}_{\mathrm{simplicial}}(\mathrm{Ind}\mathrm{Seminormed}_?)}_{\mathrm{Noncommutativealgebra}_{\mathrm{simplicial}}(\mathrm{Ind}\mathrm{Seminormed}_?)^\mathrm{opposite},\mathrm{Grotopology,homotopyepimorphism}}$: $\infty$-sheaves into $\infty$-groupoid over corresponding opposite category carrying the corresponding Grothendieck topology, and we will mainly consider the corresponding homotopy epimorphisms. $?=R,\mathbb{F}_1$. And we assume the stack carries $\infty$-ringed toposes structure. \\
\noindent $\infty-\mathrm{Toposes}^{\mathrm{ringed},\mathrm{Noncommutativealgebra}_{\mathrm{simplicial}}(\mathrm{Ind}^m\mathrm{Seminormed}_?)}_{\mathrm{Noncommutativealgebra}_{\mathrm{simplicial}}(\mathrm{Ind}^m\mathrm{Seminormed}_?)^\mathrm{opposite},\mathrm{Grotopology,homotopyepimorphism}}$: $\infty$-sheaves into $\infty$-groupoid over corresponding opposite category carrying the corresponding Grothendieck topology, and we will mainly consider the corresponding homotopy epimorphisms. $?=R,\mathbb{F}_1$. And we assume the stack carries $\infty$-ringed toposes structure.\\
\noindent $\infty-\mathrm{Toposes}^{\mathrm{ringed},\mathrm{Noncommutativealgebra}_{\mathrm{simplicial}}(\mathrm{Ind}\mathrm{Normed}_?)}_{\mathrm{Noncommutativealgebra}_{\mathrm{simplicial}}(\mathrm{Ind}\mathrm{Normed}_?)^\mathrm{opposite},\mathrm{Grotopology,homotopyepimorphism}}$: $\infty$-sheaves into $\infty$-groupoid over corresponding opposite category carrying the corresponding Grothendieck topology, and we will mainly consider the corresponding homotopy epimorphisms. $?=R,\mathbb{F}_1$. And we assume the stack carries $\infty$-ringed toposes structure.\\
\noindent $\infty-\mathrm{Toposes}^{\mathrm{ringed},\mathrm{Noncommutativealgebra}_{\mathrm{simplicial}}(\mathrm{Ind}^m\mathrm{Normed}_?)}_{\mathrm{Noncommutativealgebra}_{\mathrm{simplicial}}(\mathrm{Ind}^m\mathrm{Normed}_?)^\mathrm{opposite},\mathrm{Grotopology,homotopyepimorphism}}$: $\infty$-sheaves into $\infty$-groupoid over corresponding opposite category carrying the corresponding Grothendieck topology, and we will mainly consider the corresponding homotopy epimorphisms. $?=R,\mathbb{F}_1$. And we assume the stack carries $\infty$-ringed toposes structure.\\
\noindent $\infty-\mathrm{Toposes}^{\mathrm{ringed},\mathrm{Noncommutativealgebra}_{\mathrm{simplicial}}(\mathrm{Ind}\mathrm{Banach}_?)}_{\mathrm{Noncommutativealgebra}_{\mathrm{simplicial}}(\mathrm{Ind}\mathrm{Banach}_?)^\mathrm{opposite},\mathrm{Grotopology,homotopyepimorphism}}$: $\infty$-sheaves into $\infty$-groupoid over corresponding opposite category carrying the corresponding Grothendieck topology, and we will mainly consider the corresponding homotopy epimorphisms. $?=R,\mathbb{F}_1$. And we assume the stack carries $\infty$-ringed toposes structure.\\
\noindent $\infty-\mathrm{Toposes}^{\mathrm{ringed},\mathrm{Noncommutativealgebra}_{\mathrm{simplicial}}(\mathrm{Ind}^m\mathrm{Banach}_?)}_{\mathrm{Noncommutativealgebra}_{\mathrm{simplicial}}(\mathrm{Ind}^m\mathrm{Banach}_?)^\mathrm{opposite},\mathrm{Grotopology,homotopyepimorphism}}$: $\infty$-sheaves into $\infty$-groupoid over corresponding opposite category carrying the corresponding Grothendieck topology, and we will mainly consider the corresponding homotopy epimorphisms. $?=R,\mathbb{F}_1$. And we assume the stack carries $\infty$-ringed toposes structure.\\

 \noindent $\mathrm{Proj}^\text{smoothformalseriesclosure}\infty-\mathrm{Toposes}^{\mathrm{ringed},\mathrm{Noncommutativealgebra}_{\mathrm{simplicial}}(\mathrm{Ind}\mathrm{Seminormed}_?)}_{\mathrm{Noncommutativealgebra}_{\mathrm{simplicial}}(\mathrm{Ind}\mathrm{Seminormed}_?)^\mathrm{opposite},\mathrm{Grotopology,homotopyepimorphism}}$: $\infty$-sheaves into $\infty$-groupoid over corresponding opposite category carrying the corresponding Grothendieck topology, and we will mainly consider the corresponding homotopy epimorphisms. $?=R,\mathbb{F}_1$. And we assume the stack carries $\infty$-ringed toposes structure by specializing a ring object $\mathcal{O}$. \\
\noindent $\mathrm{Proj}^\text{smoothformalseriesclosure}\infty-\mathrm{Toposes}^{\mathrm{ringed},\mathrm{Noncommutativealgebra}_{\mathrm{simplicial}}(\mathrm{Ind}^m\mathrm{Seminormed}_?)}_{\mathrm{Noncommutativealgebra}_{\mathrm{simplicial}}(\mathrm{Ind}^m\mathrm{Seminormed}_?)^\mathrm{opposite},\mathrm{Grotopology,homotopyepimorphism}}$: $\infty$-sheaves into $\infty$-groupoid over corresponding opposite category carrying the corresponding Grothendieck topology, and we will mainly consider the corresponding homotopy epimorphisms. $?=R,\mathbb{F}_1$. And we assume the stack carries $\infty$-ringed toposes structure by specializing a ring object $\mathcal{O}$.\\
\noindent $\mathrm{Proj}^\text{smoothformalseriesclosure}\infty-\mathrm{Toposes}^{\mathrm{ringed},\mathrm{Noncommutativealgebra}_{\mathrm{simplicial}}(\mathrm{Ind}\mathrm{Normed}_?)}_{\mathrm{Noncommutativealgebra}_{\mathrm{simplicial}}(\mathrm{Ind}\mathrm{Normed}_?)^\mathrm{opposite},\mathrm{Grotopology,homotopyepimorphism}}$: $\infty$-sheaves into $\infty$-groupoid over corresponding opposite category carrying the corresponding Grothendieck topology, and we will mainly consider the corresponding homotopy epimorphisms. $?=R,\mathbb{F}_1$. And we assume the stack carries $\infty$-ringed toposes structure by specializing a ring object $\mathcal{O}$.\\
\noindent $\mathrm{Proj}^\text{smoothformalseriesclosure}\infty-\mathrm{Toposes}^{\mathrm{ringed},\mathrm{Noncommutativealgebra}_{\mathrm{simplicial}}(\mathrm{Ind}^m\mathrm{Normed}_?)}_{\mathrm{Noncommutativealgebra}_{\mathrm{simplicial}}(\mathrm{Ind}^m\mathrm{Normed}_?)^\mathrm{opposite},\mathrm{Grotopology,homotopyepimorphism}}$: $\infty$-sheaves into $\infty$-groupoid over corresponding opposite category carrying the corresponding Grothendieck topology, and we will mainly consider the corresponding homotopy epimorphisms. $?=R,\mathbb{F}_1$. And we assume the stack carries $\infty$-ringed toposes structure by specializing a ring object $\mathcal{O}$.\\
\noindent $\mathrm{Proj}^\text{smoothformalseriesclosure}\infty-\mathrm{Toposes}^{\mathrm{ringed},\mathrm{Noncommutativealgebra}_{\mathrm{simplicial}}(\mathrm{Ind}\mathrm{Banach}_?)}_{\mathrm{Noncommutativealgebra}_{\mathrm{simplicial}}(\mathrm{Ind}\mathrm{Banach}_?)^\mathrm{opposite},\mathrm{Grotopology,homotopyepimorphism}}$: $\infty$-sheaves into $\infty$-groupoid over corresponding opposite category carrying the corresponding Grothendieck topology, and we will mainly consider the corresponding homotopy epimorphisms. $?=R,\mathbb{F}_1$. And we assume the stack carries $\infty$-ringed toposes structure by specializing a ring object $\mathcal{O}$.\\
\noindent $\mathrm{Proj}^\text{smoothformalseriesclosure}\infty-\mathrm{Toposes}^{\mathrm{ringed},\mathrm{Noncommutativealgebra}_{\mathrm{simplicial}}(\mathrm{Ind}^m\mathrm{Banach}_?)}_{\mathrm{Noncommutativealgebra}_{\mathrm{simplicial}}(\mathrm{Ind}^m\mathrm{Banach}_?)^\mathrm{opposite},\mathrm{Grotopology,homotopyepimorphism}}$: $\infty$-sheaves into $\infty$-groupoid over corresponding opposite category carrying the corresponding Grothendieck topology, and we will mainly consider the corresponding homotopy epimorphisms. $?=R,\mathbb{F}_1$. And we assume the stack carries $\infty$-ringed toposes structure by specializing a ring object $\mathcal{O}$.\\

\noindent $\infty$-Quasicoherent Sheaves of Functional Analytic Modules over Ringed Toposes $\sharp=\mathrm{Seminormed},\mathrm{Normed},\mathrm{Banach}$: \\
 
 \noindent $\mathrm{Ind}\mathrm{\sharp Quasicoherent}_{\infty-\mathrm{Toposes}^{\mathrm{ringed},\mathrm{Noncommutativealgebra}_{\mathrm{simplicial}}(\mathrm{Ind}\mathrm{Seminormed}_?)}_{\mathrm{Noncommutativealgebra}_{\mathrm{simplicial}}(\mathrm{Ind}\mathrm{Seminormed}_?)^\mathrm{opposite},\mathrm{Grotopology,homotopyepimorphism}}}$: Colimits completion of $\infty$-Quasicoherent Sheaves of Functional Analytic Modules over $\infty$-sheaves into $\infty$-groupoid over corresponding opposite category carrying the corresponding Grothendieck topology, and we will mainly consider the corresponding homotopy epimorphisms. $?=R,\mathbb{F}_1$. And we assume the stack carries $\infty$-ringed toposes structure. \\
\noindent $\mathrm{Ind}\mathrm{\sharp Quasicoherent}_{\infty-\mathrm{Toposes}^{\mathrm{ringed},\mathrm{Noncommutativealgebra}_{\mathrm{simplicial}}(\mathrm{Ind}^m\mathrm{Seminormed}_?)}_{\mathrm{Noncommutativealgebra}_{\mathrm{simplicial}}(\mathrm{Ind}^m\mathrm{Seminormed}_?)^\mathrm{opposite},\mathrm{Grotopology,homotopyepimorphism}}}$: Colimits completion of $\infty$-Quasicoherent Sheaves of Functional Analytic Modules over $\infty$-sheaves into $\infty$-groupoid over corresponding opposite category carrying the corresponding Grothendieck topology, and we will mainly consider the corresponding homotopy epimorphisms. $?=R,\mathbb{F}_1$. And we assume the stack carries $\infty$-ringed toposes structure.\\
\noindent $\mathrm{Ind}\mathrm{\sharp Quasicoherent}_{\infty-\mathrm{Toposes}^{\mathrm{ringed},\mathrm{Noncommutativealgebra}_{\mathrm{simplicial}}(\mathrm{Ind}\mathrm{Normed}_?)}_{\mathrm{Noncommutativealgebra}_{\mathrm{simplicial}}(\mathrm{Ind}\mathrm{Normed}_?)^\mathrm{opposite},\mathrm{Grotopology,homotopyepimorphism}}}$: Colimits completion of $\infty$-Quasicoherent Sheaves of Functional Analytic Modules over $\infty$-sheaves into $\infty$-groupoid over corresponding opposite category carrying the corresponding Grothendieck topology, and we will mainly consider the corresponding homotopy epimorphisms. $?=R,\mathbb{F}_1$. And we assume the stack carries $\infty$-ringed toposes structure.\\
\noindent $\mathrm{Ind}\mathrm{\sharp Quasicoherent}_{\infty-\mathrm{Toposes}^{\mathrm{ringed},\mathrm{Noncommutativealgebra}_{\mathrm{simplicial}}(\mathrm{Ind}^m\mathrm{Normed}_?)}_{\mathrm{Noncommutativealgebra}_{\mathrm{simplicial}}(\mathrm{Ind}^m\mathrm{Normed}_?)^\mathrm{opposite},\mathrm{Grotopology,homotopyepimorphism}}}$: Colimits completion of $\infty$-Quasicoherent Sheaves of Functional Analytic Modules over $\infty$-sheaves into $\infty$-groupoid over corresponding opposite category carrying the corresponding Grothendieck topology, and we will mainly consider the corresponding homotopy epimorphisms. $?=R,\mathbb{F}_1$. And we assume the stack carries $\infty$-ringed toposes structure.\\
\noindent $\mathrm{Ind}\mathrm{\sharp Quasicoherent}_{\infty-\mathrm{Toposes}^{\mathrm{ringed},\mathrm{Noncommutativealgebra}_{\mathrm{simplicial}}(\mathrm{Ind}\mathrm{Banach}_?)}_{\mathrm{Noncommutativealgebra}_{\mathrm{simplicial}}(\mathrm{Ind}\mathrm{Banach}_?)^\mathrm{opposite},\mathrm{Grotopology,homotopyepimorphism}}}$: Colimits completion of $\infty$-Quasicoherent Sheaves of Functional Analytic Modules over $\infty$-sheaves into $\infty$-groupoid over corresponding opposite category carrying the corresponding Grothendieck topology, and we will mainly consider the corresponding homotopy epimorphisms. $?=R,\mathbb{F}_1$. And we assume the stack carries $\infty$-ringed toposes structure.\\
\noindent $\mathrm{Ind}\mathrm{\sharp Quasicoherent}_{\infty-\mathrm{Toposes}^{\mathrm{ringed},\mathrm{Noncommutativealgebra}_{\mathrm{simplicial}}(\mathrm{Ind}^m\mathrm{Banach}_?)}_{\mathrm{Noncommutativealgebra}_{\mathrm{simplicial}}(\mathrm{Ind}^m\mathrm{Banach}_?)^\mathrm{opposite},\mathrm{Grotopology,homotopyepimorphism}}}$: Colimits completion of $\infty$-Quasicoherent Sheaves of Functional Analytic Modules over $\infty$-sheaves into $\infty$-groupoid over corresponding opposite category carrying the corresponding Grothendieck topology, and we will mainly consider the corresponding homotopy epimorphisms. $?=R,\mathbb{F}_1$. And we assume the stack carries $\infty$-ringed toposes structure.\\

 \noindent $\mathrm{Ind}\mathrm{\sharp Quasicoherent}_{\mathrm{Ind}^\text{smoothformalseriesclosure}\infty-\mathrm{Toposes}^{\mathrm{ringed},\mathrm{Noncommutativealgebra}_{\mathrm{simplicial}}(\mathrm{Ind}\mathrm{Seminormed}_?)}_{\mathrm{Noncommutativealgebra}_{\mathrm{simplicial}}(\mathrm{Ind}\mathrm{Seminormed}_?)^\mathrm{opposite},\mathrm{Grotopology,homotopyepimorphism}}}$: Colimits completion of $\infty$-Quasicoherent Sheaves of Functional Analytic Modules over $\infty$-sheaves into $\infty$-groupoid over corresponding opposite category carrying the corresponding Grothendieck topology, and we will mainly consider the corresponding homotopy epimorphisms. $?=R,\mathbb{F}_1$. And we assume the stack carries $\infty$-ringed toposes structure by specializing a ring object $\mathcal{O}$. The difference is that we consider the corresponding $\mathcal{O}$ living in the colimit completion closure. \\
\noindent $\mathrm{Ind}\mathrm{\sharp Quasicoherent}_{\mathrm{Ind}^\text{smoothformalseriesclosure}\infty-\mathrm{Toposes}^{\mathrm{ringed},\mathrm{Noncommutativealgebra}_{\mathrm{simplicial}}(\mathrm{Ind}^m\mathrm{Seminormed}_?)}_{\mathrm{Noncommutativealgebra}_{\mathrm{simplicial}}(\mathrm{Ind}^m\mathrm{Seminormed}_?)^\mathrm{opposite},\mathrm{Grotopology,homotopyepimorphism}}}$: Colimits completion of $\infty$-Quasicoherent Sheaves of Functional Analytic Modules over $\infty$-sheaves into $\infty$-groupoid over corresponding opposite category carrying the corresponding Grothendieck topology, and we will mainly consider the corresponding homotopy epimorphisms. $?=R,\mathbb{F}_1$. And we assume the stack carries $\infty$-ringed toposes structure by specializing a ring object $\mathcal{O}$. The difference is that we consider the corresponding $\mathcal{O}$ living in the colimit completion closure.\\
\noindent $\mathrm{Ind}\mathrm{\sharp Quasicoherent}_{\mathrm{Ind}^\text{smoothformalseriesclosure}\infty-\mathrm{Toposes}^{\mathrm{ringed},\mathrm{Noncommutativealgebra}_{\mathrm{simplicial}}(\mathrm{Ind}\mathrm{Normed}_?)}_{\mathrm{Noncommutativealgebra}_{\mathrm{simplicial}}(\mathrm{Ind}\mathrm{Normed}_?)^\mathrm{opposite},\mathrm{Grotopology,homotopyepimorphism}}}$: Colimits completion of $\infty$-Quasicoherent Sheaves of Functional Analytic Modules over $\infty$-sheaves into $\infty$-groupoid over corresponding opposite category carrying the corresponding Grothendieck topology, and we will mainly consider the corresponding homotopy epimorphisms. $?=R,\mathbb{F}_1$. And we assume the stack carries $\infty$-ringed toposes structure by specializing a ring object $\mathcal{O}$. The difference is that we consider the corresponding $\mathcal{O}$ living in the colimit completion closure.\\
\noindent $\mathrm{Ind}\mathrm{\sharp Quasicoherent}_{\mathrm{Ind}^\text{smoothformalseriesclosure}\infty-\mathrm{Toposes}^{\mathrm{ringed},\mathrm{Noncommutativealgebra}_{\mathrm{simplicial}}(\mathrm{Ind}^m\mathrm{Normed}_?)}_{\mathrm{Noncommutativealgebra}_{\mathrm{simplicial}}(\mathrm{Ind}^m\mathrm{Normed}_?)^\mathrm{opposite},\mathrm{Grotopology,homotopyepimorphism}}}$: Colimits completion of $\infty$-Quasicoherent Sheaves of Functional Analytic Modules over $\infty$-sheaves into $\infty$-groupoid over corresponding opposite category carrying the corresponding Grothendieck topology, and we will mainly consider the corresponding homotopy epimorphisms. $?=R,\mathbb{F}_1$. And we assume the stack carries $\infty$-ringed toposes structure by specializing a ring object $\mathcal{O}$. The difference is that we consider the corresponding $\mathcal{O}$ living in the colimit completion closure.\\
\noindent $\mathrm{Ind}\mathrm{\sharp Quasicoherent}_{\mathrm{Ind}^\text{smoothformalseriesclosure}\infty-\mathrm{Toposes}^{\mathrm{ringed},\mathrm{Noncommutativealgebra}_{\mathrm{simplicial}}(\mathrm{Ind}\mathrm{Banach}_?)}_{\mathrm{Noncommutativealgebra}_{\mathrm{simplicial}}(\mathrm{Ind}\mathrm{Banach}_?)^\mathrm{opposite},\mathrm{Grotopology,homotopyepimorphism}}}$: Colimits completion of $\infty$-Quasicoherent Sheaves of Functional Analytic Modules over $\infty$-sheaves into $\infty$-groupoid over corresponding opposite category carrying the corresponding Grothendieck topology, and we will mainly consider the corresponding homotopy epimorphisms. $?=R,\mathbb{F}_1$. And we assume the stack carries $\infty$-ringed toposes structure by specializing a ring object $\mathcal{O}$. The difference is that we consider the corresponding $\mathcal{O}$ living in the colimit completion closure.\\
\noindent $\mathrm{Ind}\mathrm{\sharp Quasicoherent}_{\mathrm{Ind}^\text{smoothformalseriesclosure}\infty-\mathrm{Toposes}^{\mathrm{ringed},\mathrm{Noncommutativealgebra}_{\mathrm{simplicial}}(\mathrm{Ind}^m\mathrm{Banach}_?)}_{\mathrm{Noncommutativealgebra}_{\mathrm{simplicial}}(\mathrm{Ind}^m\mathrm{Banach}_?)^\mathrm{opposite},\mathrm{Grotopology,homotopyepimorphism}}}$: Colimits completion of $\infty$-Quasicoherent Sheaves of Functional Analytic Modules over $\infty$-sheaves into $\infty$-groupoid over corresponding opposite category carrying the corresponding Grothendieck topology, and we will mainly consider the corresponding homotopy epimorphisms. $?=R,\mathbb{F}_1$. And we assume the stack carries $\infty$-ringed toposes structure by specializing a ring object $\mathcal{O}$. The difference is that we consider the corresponding $\mathcal{O}$ living in the colimit completion closure.\\

Let us explain slightly what is happening here, the base $\infty$-rings are noncommutative in certain sense, which is not the same as in the foundation of \cite{12BBK}, \cite{BBM}, \cite{KKM}, \cite{12BK}. Certainly for instance one considers the corresponding $\infty$-category $\mathrm{Simpicial}(\mathrm{Ind}\mathrm{Banach}_{\mathbb{F}_1})$, then takes the corresponding fibrations over the corresponding noncommutative rings to achieve so.

\newpage

%\chapter{Topological Theory}

\section{Topological Andr\'e-Quillen Homology and Topological Derived de Rham Complexes}

\subsection{Derived $p$-Complete Derived de Rham Complex}

\indent We now first discuss the corresponding Banach version of  Andr\'e-Quillen Homology and the corresponding Banach version of  Derived de Rham complex parallel to \cite[Chapitre 3]{12An1}, \cite{12An2}, \cite[Chapter 2, Chapter 8]{12B1}, \cite[Chapter 1]{12Bei}, \cite[Chapter 5]{12G1}, \cite[Chapter 3, Chapter 4]{12GL}, \cite[Chapitre II, Chapitre III]{12Ill1}, \cite[Chapitre VIII]{12Ill2}, \cite[Section 4]{12Qui}. We would like to start from the corresponding context of \cite[Chapter 3, Chapter 4]{12GL}, and represent the construction for the convenience of the readers. We start from the corresponding construction of the algebraic $p$-adic derived de Rham complex for a map $A\rightarrow B$ of $p$-complete rings. This is the corresponding derived differential complex attached to the polynomial resolution of $B$:
\begin{align}
A[A[B]]...,	
\end{align}
which is now denoted by $\mathrm{Kan}_\mathrm{Left}\mathrm{deRham}^\text{degreenumber}_{B/A,,\mathrm{alg}}:=\mathrm{Kan}_\mathrm{Left}\mathrm{deRham}^\text{degreenumber}_{-/A,,\mathrm{alg}}(B)$ after taking the corresponding left Kan extension which will be the same for all the following constructions \footnote{We have already considered the corresponding left Kan extension to all the rings which are not concentrated at degree zero after \cite[Example 5.11, Example 5.12]{12BMS}, which is also discussed in \cite[Lecture 7]{12B2}.}. The corresponding cotangent complex associated is defined to be just:
\begin{align}
\mathbb{L}_{B/A,\mathrm{alg}}:=	\mathrm{deRham}^1_{A[B]^\text{degreenumber}/A,,\mathrm{alg}}\otimes_{A[B]^\text{degreenumber}} B.
\end{align}
The corresponding algebraic Andr\'e-Quillen homologies are defined to be:
\begin{align}
H_{\text{degreenumber},{\mathrm{AQ}},\mathrm{alg}}:=\pi_\text{degreenumber} (\mathbb{L}_{B/A,\mathrm{alg}}). 	
\end{align}

\indent The corresponding topological Andr\'e-Quillen complex is actually the completed version of the corresponding algebraic ones above by considering the corresponding certain $p$-completion over the simplicial module structure.

\indent Then we consider the corresponding derived algebraic de Rham complex which is just defined to be:
\begin{align}
\mathrm{Kan}_\mathrm{Left}\mathrm{deRham}^\text{degreenumber}_{B/A,\mathrm{alg}},\mathrm{Kan}_\mathrm{Left}\mathrm{Fil}^*_{\mathrm{deRham}^\text{degreenumber}_{B/A,\mathrm{alg}}}.	
\end{align}
We then take the corresponding Banach completion and we denote that by:
\begin{align}
\mathrm{Kan}_\mathrm{Left}\mathrm{deRham}^\text{degreenumber}_{B/A,\mathrm{topo}},\mathrm{Kan}_\mathrm{Left}\mathrm{Fil}^*_{\mathrm{deRham}^\text{degreenumber}_{B/A,\mathrm{topo}}}.	
\end{align}

\indent Then we need to take the corresponding Hodge-Filtered completion by using the corresponding filtration associated as above:
\begin{align}
\mathrm{Kan}_\mathrm{Left}\widehat{\mathrm{deRham}}^\text{degreenumber}_{B/A,\mathrm{topo}},\mathrm{Kan}_\mathrm{Left}\mathrm{Fil}^*_{\widehat{\mathrm{deRham}}^\text{degreenumber}_{B/A,\mathrm{topo}}}.	
\end{align}

This is basically the corresponding analytic and complete version the corresponding algebraic de Rham complex. Furthermore we allow large coefficients with rigid affinoid algebra $Z$ over $\mathbb{Q}_p$. Therefore we take the corresponding completed tensor product in the following. 

\begin{definition}
We define the following $Z$ deformed version of the corresponding complete version of the corresponding Andr\'e-Quillen homology and the corresponding complete version of the corresponding derived de Rham complex. We start from the corresponding construction of the algebraic $p$-adic derived de Rham complex for a map $A\rightarrow B$. Fix a pair of ring of definitions $A_0,B_0$ in $A,B$ respectively. Then this is the corresponding derived differential complex attached to the polynomial resolution of $B_0$:
\begin{align}
A_0[A_0[B_0]]...,	
\end{align}
which is now denoted by $\mathrm{Kan}_\mathrm{Left}\mathrm{deRham}^\text{degreenumber}_{B_0/A_0,\mathrm{alg}}$. The corresponding cotangent complex associated is defined to be just:
\begin{align}
\mathbb{L}_{B_0/A_0,\mathrm{alg}}:=	\mathrm{deRham}^1_{A_0[B_0]^\text{degreenumber}/A_0,\mathrm{alg}}\otimes_{A_0[B_0]^\text{degreenumber}} B_0.
\end{align}
The corresponding algebraic Andr\'e-Quillen homologies are defined to be:
\begin{align}
H_{\text{degreenumber},{\mathrm{AQ}},\mathrm{alg}}:=\pi_\text{degreenumber} (\mathbb{L}_{B_0/A_0,\mathrm{alg}}). 	
\end{align}
The corresponding topological Andr\'e-Quillen complex is actually the complete version of the corresponding algebraic ones above by considering the corresponding derived $p$-completion over the simplicial module structure:
\begin{align}
\mathbb{L}_{B_0/A_0,\mathrm{topo}}:=R\varprojlim_k	\mathrm{Kos}_{p^k}\left((\mathrm{deRham}^1_{A_0[B_0]^\text{degreenumber}/A_0,\mathrm{alg}}\otimes_{A_0[B_0]^\text{degreenumber}} B_0)\right).
\end{align}
Taking the product with $\mathcal{O}_Z$ we have the corresponding integral version of the topological Andr\'e-Quillen complex:
\begin{align}
\mathbb{L}_{B_0/A_0,\mathrm{topo},Z}:=\mathbb{L}_{B_0/A_0,\mathrm{topo}}\widehat{\otimes}_{\mathbb{Z}_p}\mathcal{O}_Z
\end{align}
Then we consider the corresponding derived algebraic de Rham complex which is just defined to be:
\begin{align}
\mathrm{Kan}_\mathrm{Left}\mathrm{deRham}^\text{degreenumber}_{B_0/A_0,\mathrm{alg}},\mathrm{Kan}_\mathrm{Left}\mathrm{Fil}^*_{\mathrm{deRham}^1_{B_0/A_0,\mathrm{alg}}}.	
\end{align}
We then take the corresponding derived $p$-completion and we denote that by:
\begin{align}
\mathrm{Kan}_\mathrm{Left}\mathrm{deRham}^\text{degreenumber}_{B_0/A_0,\mathrm{topo}}:=R\varprojlim_k\mathrm{Kos}_{p^k}\left(\mathrm{Kan}_\mathrm{Left}\mathrm{deRham}^\text{degreenumber}_{B_0/A_0,\mathrm{alg}},\mathrm{Kan}_\mathrm{Left}\mathrm{Fil}^*_{\mathrm{deRham}^\text{degreenumber}_{B_0/A_0,\mathrm{alg}}}\right),\\
\mathrm{Kan}_\mathrm{Left}\mathrm{Fil}^*_{\mathrm{deRham}^\text{degreenumber}_{B_0/A_0,\mathrm{topo}}}:=R\varprojlim_k\mathrm{Kos}_{p^k}\left(\mathrm{Kan}_\mathrm{Left}\mathrm{Fil}^*_{\mathrm{deRham}^\text{degreenumber}_{B_0/A_0,\mathrm{alg}}}\right).	
\end{align}
Before considering the corresponding integral version we just consider the corresponding product of these $\mathbb{E}_\infty$-rings with $\mathcal{O}_Z$ to get:
\begin{align}
\mathrm{Kan}_\mathrm{Left}\mathrm{deRham}^\text{degreenumber}_{B_0/A_0,\mathrm{topo},Z}:=\mathrm{Kan}_\mathrm{Left}\mathrm{deRham}^\text{degreenumber}_{B_0/A_0,\mathrm{topo}}\widehat{\otimes}^\mathbb{L}_{\mathbb{Z}_p}\mathcal{O}_Z,\\
\mathrm{Kan}_\mathrm{Left}\mathrm{Fil}^*_{\mathrm{deRham}^\text{degreenumber}_{B_0/A_0,\mathrm{topo}},Z}:=\mathrm{Kan}_\mathrm{Left}\mathrm{Fil}^*_{\mathrm{deRham}^\text{degreenumber}_{B_0/A_0,\mathrm{topo}}}\widehat{\otimes}^\mathbb{L}_{\mathbb{Z}_p}\mathcal{O}_Z.	
\end{align}
Then we consider the following construction for the map $A\rightarrow B$ by putting:
\begin{align}
\mathbb{L}_{B_0/A_0,\mathrm{topo},Z}:= \mathrm{Colim}_{A_0\rightarrow B_0}\mathbb{L}_{B_0/A_0,\mathrm{topo},Z}[1/p],\\
H_{\text{degreenumber},{\mathrm{AQ}},\mathrm{topo},Z}:=\pi_\text{degreenumber} (\mathbb{L}_{B/A,\mathrm{topo},Z}),	\\
\mathrm{Kan}_\mathrm{Left}\mathrm{deRham}^\text{degreenumber}_{B/A,\mathrm{topo},Z}:=\mathrm{Colim}_{A_0\rightarrow B_0}\mathrm{Kan}_\mathrm{Left}\mathrm{deRham}^\text{degreenumber}_{B_0/A_0,\mathrm{topo},Z}[1/p],\\
\mathrm{Kan}_\mathrm{Left}\mathrm{Fil}^*_{\mathrm{deRham}^\text{degreenumber}_{B/A,\mathrm{topo}},Z}:=\mathrm{Colim}_{A_0\rightarrow B_0}\mathrm{Kan}_\mathrm{Left}\mathrm{Fil}^*_{\mathrm{deRham}^1_{B_0/A_0,\mathrm{topo}},Z}[1/p].
\end{align}
Then we need to take the corresponding Hodge-Filtered completion by using the corresponding filtration associated as above to achieve the corresponding Hodge-complete objects in the corresponding filtered $\infty$-categories:
\begin{align}
\mathrm{Kan}_\mathrm{Left}{\mathrm{deRham}}^\text{degreenumber}_{B/A,\mathrm{topo,Hodge},Z},\mathrm{Kan}_\mathrm{Left}\mathrm{Fil}^*_{{\mathrm{deRham}}^\text{degreenumber}_{B/A,\mathrm{topo,Hodge}},Z}.	
\end{align} 	
\end{definition}

\begin{definition}
We define the corresponding finite projective filtered crystals to be the corresponding finite projective module spectra over the topological filtered $E_\infty$-ring $\mathrm{Kan}_\mathrm{Left}\mathrm{deRham}^\text{degreenumber}_{B/A,\mathrm{topo},Z}$ with the corresponding induced filtrations.	
\end{definition}

\begin{definition}
We define the corresponding almost perfect \footnote{This is the corresponding derived version of pseudocoherence from \cite{12Lu1}, \cite{12Lu2}.} filtered crystals to be the corresponding almost perfect module spectra over the topological filtered $E_\infty$-ring $\mathrm{Kan}_\mathrm{Left}\mathrm{deRham}^\text{degreenumber}_{B/A,\mathrm{topo},Z}$ with the corresponding induced filtrations.	
\end{definition}

\indent The following is derived from the main Poincar\'e Lemma from \cite[Theorem 1.2]{12GL} in the non-deformed situation. Consider a corresponding smooth rigid analytic space $X$ over $k/\mathbb{Q}_p$ (where $k$ is a corresponding unramified analytic field which is discretely-valued and the corresponding residue field is finite). Then we have the following:

\newpage

\begin{landscape}
\begin{proposition}
Consider the corresponding projective map $g:X_{\text{pro-\'etale}}\rightarrow X_{\text{\'et}}$ and the the projective map $f:X_{\text{pro-\'etale}}\rightarrow X$. Then we have the following two strictly exact long exact sequences 
\footnote{We should have the corresponding naturality taking into the following form:

\[\tiny
\xymatrix@C+0.4pc@R+0pc{
0\ar[r]\ar[r]\ar[r] \ar[r] &\mathrm{Kan}_\mathrm{Left}\mathrm{deRham}^\text{degreenumber}_{k[\widehat{\mathcal{O}}]^\text{degreenumber}/k,\mathrm{topo},Z}\ar[d]\ar[d]\ar[d] \ar[d] \ar[r]^\partial\ar[r]\ar[r] \ar[r] &\mathrm{Kan}_\mathrm{Left}\mathrm{deRham}^\text{degreenumber}_{X[\widehat{\mathcal{O}}]^\text{degreenumber}/X,\mathrm{topo},Z} \ar[d]\ar[d]\ar[d] \ar[d]\ar[r]^\partial\ar[r]\ar[r] \ar[r] &\mathrm{Kan}_\mathrm{Left}\mathrm{deRham}^\text{degreenumber}_{X[\widehat{\mathcal{O}}]^\text{degreenumber}/X,\mathrm{topo},Z}{\otimes}f^{-1} \mathrm{deRham}^1_{X,\mathrm{topo}} \ar[d]\ar[d]\ar[d] \ar[d]\ar[r]^\partial\ar[r]\ar[r] \ar[r]&\\
0\ar[r]\ar[r]\ar[r] \ar[r] &\mathrm{Kan}_\mathrm{Left}\mathrm{deRham}^\text{degreenumber}_{k[\widehat{\mathcal{O}}]^\text{degreenumber}/k,\mathrm{topo},Z} \ar[r]^\partial\ar[r]\ar[r] \ar[r] &\mathrm{Kan}_\mathrm{Left}\mathrm{deRham}^\text{degreenumber}_{X[\widehat{\mathcal{O}}]^\text{degreenumber}/X_{\text{\'et}},\mathrm{topo},Z} \ar[r]^\partial\ar[r]\ar[r] \ar[r] &\mathrm{Kan}_\mathrm{Left}\mathrm{deRham}^\text{degreenumber}_{X[\widehat{\mathcal{O}}]^\text{degreenumber}/X,\mathrm{topo},Z}{\otimes}f^{-1} \mathrm{deRham}^1_{X_{\text{\'et}},\mathrm{topo}} \ar[r]^\partial\ar[r]\ar[r] \ar[r]&.
}
\]
}:
\[\tiny
\xymatrix@C+0.4pc@R+0pc{
0\ar[r]\ar[r]\ar[r] \ar[r] &\mathrm{Kan}_\mathrm{Left}\mathrm{deRham}^\text{degreenumber}_{k[\widehat{\mathcal{O}}]^\text{degreenumber}/k,\mathrm{topo},Z} \ar[r]^\partial\ar[r]\ar[r] \ar[r] &\mathrm{Kan}_\mathrm{Left}\mathrm{deRham}^\text{degreenumber}_{X[\widehat{\mathcal{O}}]^\text{degreenumber}/X,\mathrm{topo},Z} \ar[r]^\partial\ar[r]\ar[r] \ar[r] &\mathrm{Kan}_\mathrm{Left}\mathrm{deRham}^\text{degreenumber}_{X[\widehat{\mathcal{O}}]^\text{degreenumber}/X,\mathrm{topo},Z}{\otimes}f^{-1} \mathrm{deRham}^1_{X,\mathrm{topo}} \ar[r]^\partial\ar[r]\ar[r] \ar[r]&...,
}
\]
and
\[\tiny
\xymatrix@C+0.4pc@R+0pc{
0\ar[r]\ar[r]\ar[r] \ar[r] &\mathrm{Kan}_\mathrm{Left}\mathrm{deRham}^\text{degreenumber}_{k[\widehat{\mathcal{O}}]^\text{degreenumber}/k,\mathrm{topo},Z} \ar[r]^\partial\ar[r]\ar[r] \ar[r] &\mathrm{Kan}_\mathrm{Left}\mathrm{deRham}^\text{degreenumber}_{X[\widehat{\mathcal{O}}]^\text{degreenumber}/X_{\text{\'et}},\mathrm{topo},Z} \ar[r]^\partial\ar[r]\ar[r] \ar[r] &\mathrm{Kan}_\mathrm{Left}\mathrm{deRham}^\text{degreenumber}_{X[\widehat{\mathcal{O}}]^\text{degreenumber}/X,\mathrm{topo},Z}{\otimes}f^{-1} \mathrm{deRham}^1_{X_{\text{\'et}},\mathrm{topo}} \ar[r]^\partial\ar[r]\ar[r] \ar[r] &....
}
\]
\end{proposition}

\begin{proof}
This is actually a direct consequence \cite[Theorem 1.2]{12GL}.	
\end{proof}

\begin{proposition}
Consider the corresponding projective map $g:X_{\text{pro-\'etale}}\rightarrow X_{\text{\'et}}$ and the the projective map $f:X_{\text{pro-\'etale}}\rightarrow X$. Let $M$ be a corresponding $Z$-projective differential crystal spectrum. Then we have the following two strictly exact long exact sequences \footnote{Here we drop the notation $\mathrm{Kan}_\mathrm{Left}$.}:
\[\tiny
\xymatrix@C+0.4pc@R+0pc{
M{\otimes}^\mathbb{L}(0\ar[r]\ar[r]\ar[r] \ar[r] &\mathrm{deRham}^\text{degreenumber}_{k[\widehat{\mathcal{O}}]^\text{degreenumber}/k,\mathrm{topo},Z} \ar[r]^\partial\ar[r]\ar[r] \ar[r] &\mathrm{deRham}^\text{degreenumber}_{X[\widehat{\mathcal{O}}]^\text{degreenumber}/X,\mathrm{topo},Z} \ar[r]^\partial\ar[r]\ar[r] \ar[r] &\mathrm{deRham}^\text{degreenumber}_{X[\widehat{\mathcal{O}}]^\text{degreenumber}/X,\mathrm{topo},Z}{\otimes}f^{-1} \mathrm{deRham}^1_{X,\mathrm{topo}} \ar[r]^\partial\ar[r]\ar[r] \ar[r] &...),
}
\]
and
\[\tiny
\xymatrix@C+0.4pc@R+0pc{
M{\otimes}^\mathbb{L}(0\ar[r]\ar[r]\ar[r] \ar[r] &\mathrm{deRham}^\text{degreenumber}_{k[\widehat{\mathcal{O}}]^\text{degreenumber}/k,\mathrm{topo},Z} \ar[r]^\partial\ar[r]\ar[r] \ar[r] &\mathrm{deRham}^\text{degreenumber}_{X[\widehat{\mathcal{O}}]^\text{degreenumber}/X_{\text{\'et}},\mathrm{topo},Z} \ar[r]^\partial\ar[r]\ar[r] \ar[r] &\mathrm{deRham}^\text{degreenumber}_{X[\widehat{\mathcal{O}}]^\text{degreenumber}/X,\mathrm{topo},Z}{\otimes}f^{-1} \mathrm{deRham}^1_{X_{\text{\'et}},\mathrm{topo}} \ar[r]^\partial\ar[r]\ar[r] \ar[r] &...).
}
\]
\end{proposition}

\begin{proof}
This is actually a direct consequence of the previous proposition due to the corresponding flatness.	
\end{proof}

\end{landscape}

\newpage

\subsection{Analytic Andr\'e-Quillen Homology and Analytic Derived de Rham complex of Pseudorigid Spaces}

%\subsection{The Definitions}

%%%%%%%%%%%%%%%%%%%%%%%%%%%%%%%-------------------------------

\indent We now consider the analytic Andr\'e-Quillen Homology and analytic Derived de Rham complex of pseudorigid space, which are very crucial in some development in \cite{12Bel1} in the arithmetic family. Therefore we just investigate the corresponding picture in the corresponding geometric family.

\begin{setting}
We consider now a corresponding morphism taking the corresponding form of $A\rightarrow B$ where $A$ is going to be a pseudorigid affinoid algebra over $\mathbb{Z}_p$, and $B$ is going to be a perfectoid chart of $A$ in the corresponding pro-\'etale site of the pseudorigid affinoid space attached to $A$. As in \cite[Definition 3.1, and below Definition 3.1]{12Bel1} in our situation $A$ is of topologically finite type over $\mathcal{O}_K[[t]]\left<\pi^a/t^b\right>[1/t]$, where $K$ is a discrete valued field containing $\mathbb{Q}_p$ and $(a,b)=1$ \footnote{Certainly one can also consider the characteristic $p$ situation.}.	
\end{setting}

\indent From \cite[Definition 3.1, and below Definition 3.1]{12Bel1} we have the following:

\begin{lemma} We have the following statements:\\
A. $A$ as above is Tate, complete over $\mathcal{O}_K$;\\
B. The ring $A$ has a ring of definition $A_0$ which is of $\mathcal{O}_K$-formally finite type;\\
C. The ring $A_0$ is of topologically finite type over $\mathcal{O}_K[[t]]\left<\pi^a/t^b\right>$.
\end{lemma}

\begin{proof}
It is very easy to see that the corresponding results hold in our current situation over $\mathcal{O}_K$.	
\end{proof}

As in the corresponding construction in the rigid situation in the previous section following \cite[Chapitre 3]{12An1}, \cite{12An2}, \cite[Chapter 2, Chapter 8]{12B1}, \cite[Chapter 1]{12Bei}, \cite[Chapter 5]{12G1}, \cite[Chapter 3, Chapter 4]{12GL}, \cite[Chapitre II, Chapitre III]{12Ill1}, \cite[Chapitre VIII]{12Ill2}, \cite[Section 4]{12Qui} we give the following definitions. First we consider the following new setting:

\begin{setting}
Now as in the above notion we will consider a general map of rings $A\rightarrow B$ over $\mathcal{O}_K[[t]]\left<\pi^a/t^b\right>[1/t]$ where we have the corresponding map of the associated ring of definitions $A_0\rightarrow B_0$ over $\mathcal{O}_K[[t]]\left<\pi^a/t^b\right>$ such that $A_0,B_0$ are basically $I$-adic (where $\mathcal{O}_K[[t]]\left<\pi^a/t^b\right>$ is $I$-adic).	
\end{setting}

 %%%%%%%%%%%%%%%%%%%%%%%%%%%%%%%%%%%%%%%-------------------------

\begin{definition}
We start from the corresponding construction of the algebraic $p$-adic derived de Rham complex for a map $A\rightarrow B$. Fix a pair of ring of definitions $A_0,B_0$ in $A,B$ respectively. Then this is the corresponding derived differential complex attached to the polynomial resolution of $B_0$:
\begin{align}
A_0[A_0[B_0]]...,	
\end{align}
which is now denoted by $\mathrm{Kan}_\mathrm{Left}\mathrm{deRham}^\text{degreenumber}_{B_0/A_0,\mathrm{alg}}$ after taking the corresponding left Kan extension. The corresponding cotangent complex associated is defined to be just:
\begin{align}
\mathbb{L}_{B_0/A_0,\mathrm{alg}}:=	\mathrm{deRham}^1_{A_0[B_0]^\text{degreenumber}/A_0,\mathrm{alg}}\otimes_{A_0[B_0]^\text{degreenumber}} B_0.
\end{align}
The corresponding algebraic Andr\'e-Quillen homologies are defined to be:
\begin{align}
H_{\text{degreenumber},{\mathrm{AQ}},\mathrm{alg}}:=\pi_\text{degreenumber} (\mathbb{L}_{B_0/A_0,\mathrm{alg}}). 	
\end{align}
The corresponding topological Andr\'e-Quillen complex is actually the complete version of the corresponding algebraic ones above by considering the corresponding derived $I$-completion over the simplicial module structure:
\begin{align}
\mathbb{L}_{B_0/A_0,\mathrm{topo}}:=R\varprojlim_I \mathrm{Kos}_{I}	\left((\mathrm{deRham}^1_{A_0[B_0]^\text{degreenumber}/A_0,\mathrm{alg}}\otimes_{A_0[B_0]^\text{degreenumber}} B_0)\right).
\end{align}
Then we consider the corresponding derived algebraic de Rham complex which is just defined to be:
\begin{align}
\mathrm{Kan}_\mathrm{Left}\mathrm{deRham}^\text{degreenumber}_{B_0/A_0,\mathrm{alg}},\mathrm{Kan}_\mathrm{Left}\mathrm{Fil}^*_{\mathrm{Kan}_\mathrm{Left}\mathrm{deRham}^1_{B_0/A_0,\mathrm{alg}}}.	
\end{align}
We then take the corresponding derived $I$-completion and we denote that by:
\begin{align}
\mathrm{Kan}_\mathrm{Left}\mathrm{deRham}^\text{degreenumber}_{B_0/A_0,\mathrm{topo}}:=R\varprojlim_I\mathrm{Kos}_I\left(\mathrm{Kan}_\mathrm{Left}\mathrm{deRham}^\text{degreenumber}_{B_0/A_0,\mathrm{alg}},\mathrm{Kan}_\mathrm{Left}\mathrm{Fil}^*_{\mathrm{Kan}_\mathrm{Left}\mathrm{deRham}^\text{degreenumber}_{B_0/A_0,\mathrm{alg}}}\right),\\
\mathrm{Fil}^*_{\mathrm{deRham}^\text{degreenumber}_{B_0/A_0,\mathrm{topo}}}:=R\varprojlim_I\mathrm{Kos}_I\left(\mathrm{Kan}_\mathrm{Left}\mathrm{Fil}^*_{\mathrm{Kan}_\mathrm{Left}\mathrm{deRham}^\text{degreenumber}_{B_0/A_0,\mathrm{alg}}}\right).	
\end{align}
Then we consider the following construction for the map $A\rightarrow B$ by putting:
\begin{align}
\mathbb{L}_{B_0/A_0,\mathrm{topo}}:= \mathrm{Colim}_{A_0\rightarrow B_0}\mathbb{L}_{B_0/A_0,\mathrm{topo}}[1/t],\\
H_{\text{degreenumber},{\mathrm{AQ}},\mathrm{topo}}:=\pi_\text{degreenumber} (\mathbb{L}_{B/A,\mathrm{topo}}),	\\
\mathrm{Kan}_\mathrm{Left}\mathrm{deRham}^\text{degreenumber}_{B/A,\mathrm{topo}}:=\mathrm{Colim}_{A_0\rightarrow B_0}\mathrm{Kan}_\mathrm{Left}\mathrm{deRham}^\text{degreenumber}_{B_0/A_0,\mathrm{topo}}[1/t],\\
\mathrm{Kan}_\mathrm{Left}\mathrm{Fil}^*_{\mathrm{deRham}^\text{degreenumber}_{B/A,\mathrm{topo}}}:=\mathrm{Colim}_{A_0\rightarrow B_0}\mathrm{Fil}^*_{\mathrm{Kan}_\mathrm{Left}\mathrm{deRham}^\text{degreenumber}_{B_0/A_0,\mathrm{topo}}}[1/t].
\end{align}
Then we need to take the corresponding Hodge-Filtered completion by using the corresponding filtration associated as above to achieve the corresponding Hodge-complete objects in the corresponding filtered $\infty$-categories:
\begin{align}
\mathrm{Kan}_\mathrm{Left}{\mathrm{deRham}}^\text{degreenumber}_{B/A,\mathrm{topo,Hodge}},\mathrm{Fil}^*_{\mathrm{Kan}_\mathrm{Left}{\mathrm{deRham}}^\text{degreenumber}_{B/A,\mathrm{topo,Hodge}}}.	
\end{align}
\end{definition}

\begin{definition}
We define the corresponding finite projective filtered crystals to be the corresponding finite projective module spectra over the topological filtered $E_\infty$-ring $\mathrm{Kan}_\mathrm{Left}\mathrm{deRham}^\text{degreenumber}_{B/A,\mathrm{topo}}$ with the corresponding induced filtrations.	
\end{definition}

\begin{definition}
We define the corresponding almost perfect \footnote{This is the corresponding derived version of pseudocoherence from \cite{12Lu1}, \cite{12Lu2}.} filtered crystals to be the corresponding almost perfect module spectra over the topological filtered $E_\infty$-ring $\mathrm{Kan}_\mathrm{Left}\mathrm{deRham}^\text{degreenumber}_{B /A,\mathrm{topo}}$ with the corresponding induced filtrations.	
\end{definition}

\indent We now consider the corresponding large coefficient local systems over pseudorigid spaces where $Z$ now is a topological algebra over $\mathbb{Z}_p$\footnote{It is better to assume that it is basically completely flat.}.

%%%%%%%%%%%%%%%%%%%%%%%%%%%%%%%%%%%-------------------------

\begin{definition}
We define the following $Z$ deformed version of the corresponding complete version of the corresponding Andr\'e-Quillen homology and the corresponding complete version of the corresponding derived de Rham complex. We start from the corresponding construction of the algebraic $p$-adic derived de Rham complex for a map $A\rightarrow B$. Fix a pair of ring of definitions $A_0,B_0$ in $A,B$ respectively. Then this is the corresponding derived differential complex attached to the polynomial resolution of $B_0$:
\begin{align}
A_0[A_0[B_0]]...,	
\end{align}
which is now denoted by $\mathrm{deRham}^\text{degreenumber}_{A_0[B_0]^\text{degreenumber}/A_0,\mathrm{alg}}$. The corresponding cotangent complex associated is defined to be just:
\begin{align}
\mathbb{L}_{B_0/A_0,\mathrm{alg}}:=	\mathrm{deRham}^1_{A_0[B_0]^\text{degreenumber}/A_0,\mathrm{alg}}\otimes_{A_0[B_0]^\text{degreenumber}} B_0.
\end{align}
The corresponding algebraic Andr\'e-Quillen homologies are defined to be:
\begin{align}
H_{\text{degreenumber},{\mathrm{AQ}},\mathrm{alg}}:=\pi_\text{degreenumber} (\mathbb{L}_{B_0/A_0,\mathrm{alg}}). 	
\end{align}
The corresponding topological Andr\'e-Quillen complex is actually the complete version of the corresponding algebraic ones above by considering the corresponding derived $I$-completion over the simplicial module structure:
\begin{align}
\mathbb{L}_{B_0/A_0,\mathrm{topo}}:=R\varprojlim_k	\mathrm{Kos}_I\left((\mathrm{deRham}^1_{A_0[B_0]^\text{degreenumber}/A_0,\mathrm{alg}}\otimes_{A_0[B_0]^\text{degreenumber}} B_0)\right).
\end{align}
Taking the product with $Z$ we have the corresponding integral version of the topological Andr\'e-Quillen complex\footnote{Here we did not take the corresponding derived completion, but in some situation this is achievable.}:
\begin{align}
\mathbb{L}_{B_0/A_0,\mathrm{topo},Z}:=\mathbb{L}_{B_0/A_0,\mathrm{topo}}{\otimes}^\mathbb{L}_{\mathbb{Z}_p}Z
\end{align}
Then we consider the corresponding derived algebraic de Rham complex which is just defined to be:
\begin{align}
\mathrm{Kan}_\mathrm{Left}\mathrm{deRham}^\text{degreenumber}_{B_0/A_0,\mathrm{alg}},\mathrm{Fil}^*_{\mathrm{Kan}_\mathrm{Left}\mathrm{deRham}^\text{degreenumber}_{B_0/A_0,\mathrm{alg}}}.	
\end{align}
We then take the corresponding derived $I$-completion and we denote that by:
\begin{align}
\mathrm{Kan}_\mathrm{Left}\mathrm{deRham}^\text{degreenumber}_{B_0/A_0,\mathrm{topo}}:=R\varprojlim_k\mathrm{Kos}_I\left(\mathrm{Kan}_\mathrm{Left}\mathrm{deRham}^\text{degreenumber}_{B_0/A_0,\mathrm{alg}},\mathrm{Fil}^*_{\mathrm{Kan}_\mathrm{Left}\mathrm{deRham}^\text{degreenumber}_{B_0/A_0,\mathrm{alg}}}\right),\\
\mathrm{Kan}_\mathrm{Left}\mathrm{Fil}^*_{\mathrm{deRham}^\text{degreenumber}_{B_0/A_0,\mathrm{topo}}}:=R\varprojlim_k\mathrm{Kos}_I\left(\mathrm{Kan}_\mathrm{Left}\mathrm{Fil}^*_{\mathrm{Kan}_\mathrm{Left}\mathrm{deRham}^\text{degreenumber}_{B_0/A_0,\mathrm{alg}}}\right).	
\end{align}
Before considering the corresponding integral version we just consider the corresponding product of these $\mathbb{E}_\infty$-rings with $Z$ to get\footnote{Here again we did not take the corresponding derived completion, but in some situation this is achievable, for instance when the ring $Z$ is endowed with derived $J$-complete topology.}:
\begin{align}
\mathrm{Kan}_\mathrm{Left}\mathrm{deRham}^\text{degreenumber}_{B_0/A_0,\mathrm{topo},Z}:=\mathrm{Kan}_\mathrm{Left}\mathrm{deRham}^\text{degreenumber}_{B_0/A_0,\mathrm{topo}}{\otimes}^\mathbb{L}_{\mathbb{Z}_p}Z,\\
\mathrm{Kan}_\mathrm{Left}\mathrm{Fil}^*_{\mathrm{deRham}^\text{degreenumber}_{B_0/A_0,\mathrm{topo}},Z}:=\mathrm{Fil}^*_{\mathrm{Kan}_\mathrm{Left}\mathrm{deRham}^\text{degreenumber}_{B_0/A_0,\mathrm{topo}}}{\otimes}^\mathbb{L}_{\mathbb{Z}_p}Z.	
\end{align}
Then we consider the following construction for the map $A\rightarrow B$ by putting:
\begin{align}
\mathbb{L}_{B_0/A_0,\mathrm{topo},Z}:= \mathrm{Colim}_{A_0\rightarrow B_0}\mathbb{L}_{B_0/A_0,\mathrm{topo},Z}[1/t],\\
H_{\text{degreenumber},{\mathrm{AQ}},\mathrm{topo},Z}:=\pi_\text{degreenumber} (\mathbb{L}_{B/A,\mathrm{topo},Z}),	\\
\mathrm{Kan}_\mathrm{Left}\mathrm{deRham}^\text{degreenumber}_{B/A,\mathrm{topo},Z}:=\mathrm{Colim}_{A_0\rightarrow B_0}\mathrm{Kan}_\mathrm{Left}\mathrm{deRham}^\text{degreenumber}_{B_0/A_0,\mathrm{topo},Z}[1/t],\\
\mathrm{Kan}_\mathrm{Left}\mathrm{Fil}^*_{\mathrm{Kan}_\mathrm{Left}\mathrm{deRham}^\text{degreenumber}_{B/A,\mathrm{topo}},Z}:=\mathrm{Colim}_{A_0\rightarrow B_0}\mathrm{Kan}_\mathrm{Left}\mathrm{Fil}^*_{\mathrm{Kan}_\mathrm{Left}\mathrm{deRham}^1_{B_0/A_0,\mathrm{topo}},Z}[1/t].
\end{align}
Then we need to take the corresponding Hodge-Filtered completion by using the corresponding filtration associated as above to achieve the corresponding Hodge-complete objects in the corresponding filtered $\infty$-categories:
\begin{align}
\mathrm{Kan}_\mathrm{Left}{\mathrm{deRham}}^\text{degreenumber}_{B/A,\mathrm{topo,Hodge},Z},\mathrm{Kan}_\mathrm{Left}\mathrm{Fil}^*_{\mathrm{Kan}_\mathrm{Left}{\mathrm{deRham}}^1_{B/A,\mathrm{topo,Hodge}},Z}.	
\end{align} 	
\end{definition}

\begin{example}
Now we construct the corresponding pseudorigid analog of the corresponding example in \cite[Example 4.7]{12GL}. Now we consider the corresponding the map:
\begin{align}
\mathbb{Z}_p[[u]]\left<\frac{p^a}{u^b}\right>[1/u]\left<X^\pm_1,X^\pm_2,...,X^\pm_d\right>\longrightarrow \mathbb{Z}_p[[u]]\left<\frac{p^a}{u^b}\right>[1/u]\left<X^{\pm/p^\infty}_1,X^{\pm/p^\infty}_2,...,X^{\pm/p^\infty}_d\right>	
\end{align}
which could be written as:	
\begin{align}
&\mathbb{Z}_p[[u]]\left<\frac{p^a}{u^b}\right>[1/u]\left<X^\pm_1,X^\pm_2,...,X^\pm_d\right>\longrightarrow \\
&\mathbb{Z}_p[[u]]\left<\frac{p^a}{u^b}\right>[1/u]\left<X^\pm_1,X^\pm_2,...,X^\pm_d\right>\left<Y^{\pm/p^\infty}_1,Y^{\pm/p^\infty}_2,...,Y^{\pm/p^\infty}_d\right>/(X_i-Y_i,i=1,...,d).	
\end{align}
So in our situation the corresponding Hodge complete topological derived de Rham complex will be just:
\begin{align}
\mathbb{Z}_p[[u]]\left<\frac{p^a}{u^b}\right>[1/u]\left<Y^{\pm/p^\infty}_1,Y^{\pm/p^\infty}_2,...,Y^{\pm/p^\infty}_d\right>[[Z_1,...,Z_d,Z_i=X_i-Y_i,i=1,...,d]].	
\end{align}

\end{example}

%\subsection{The Functoriality}

\indent We now follow the corresponding previous works \cite[Chapter 2, Chapter 8]{12B1}, \cite[Chapter 1]{12Bei}, \cite[Chapter 5]{12G1}, \cite[Chapter 3, Chapter 4]{12GL}, \cite{12Ill1}, \cite[Chapitre VIII]{12Ill2} to study the corresponding functoriality for some triple $A\rightarrow B \rightarrow C$ where we have first in the algebraic setting the following result. Here we first consider more general setting where we will consider those rings over $\mathcal{O}_K[[t]]\left<\frac{\pi^a}{t^b}\right>$ which is $t$-adically complete.

%%%%%%%%%%%%%%%%%%%%%%%%%%%%%%%%%%%%-----------------------

\begin{proposition}
The corresponding functoriality for the corresponding algebraic derived de Rham complex holds in our situation for the corresponding integral adic rings $A,B,C$ in the context of this section, namely we have the following corresponding commutative diagram:
\[
\xymatrix@C+0pc@R+3pc{
\mathrm{Kan}_\mathrm{Left}\mathrm{deRham}^\text{degreenumber}_{B/A,\mathrm{alg}} \ar[r] \ar[r] \ar[r]\ar[d] \ar[d] \ar[d] &\mathrm{Kan}_\mathrm{Left}\mathrm{deRham}^\text{degreenumber}_{B/A,\mathrm{alg}} \ar[d] \ar[d] \ar[d]\\
B \ar[r] \ar[r] \ar[r] &\mathrm{Kan}_\mathrm{Left}\mathrm{deRham}^\text{degreenumber}_{B/A,\mathrm{alg}}.
}
\]	
\end{proposition}

\begin{proof}
See \cite[Lemma 3.3]{12GL}.	
\end{proof}

\begin{proposition}
The corresponding functoriality for the corresponding algebraic derived de Rham complex holds in our situation for the corresponding integral adic rings $A,B,C$ in the context of this section, namely we have the following corresponding commutative diagram:
\[
\xymatrix@C+0pc@R+3pc{
\mathrm{Kan}_\mathrm{Left}\mathrm{deRham}^\text{degreenumber}_{B/A,\mathrm{topo,Hodge}} \ar[r] \ar[r] \ar[r]\ar[d] \ar[d] \ar[d] &\mathrm{Kan}_\mathrm{Left}\mathrm{deRham}^\text{degreenumber}_{B/A,\mathrm{topo,Hodge}} \ar[d] \ar[d] \ar[d]\\
B \ar[r] \ar[r] \ar[r] &\mathrm{Kan}_\mathrm{Left}\mathrm{deRham}^\text{degreenumber}_{B/A,\mathrm{topo,Hodge}}.
}
\]	
\end{proposition}

\begin{proof}
Just take the corresponding derived $I$-completion of the diagram in the previous proposition.	
\end{proof}

%\indent Now we discuss how to put the corresponding objects over the corresponding pro-\'etale site after \cite{12Sch1}, \cite{12Sch2}, \cite{12Ked1} and \cite{12GL}. The first consideration is we first put $A$ in to a family. Now let $A$ be smooth pseudorigid affinoid, we use the notation $X$ to denote a general smooth pseudorigid space. 

\indent The following conjectures are literally inspired by the rigid analytic situation in \cite[Theorem 1.2]{12GL}. We use the notation $*$ to denote the ring $\mathcal{O}_K[[t]]\left<\pi^a/t^b\right>[1/t]$. We assume the spaces in the following two conjectures are smooth pseudo-rigid analytic spaces over $\mathcal{O}_K[[t]]\left<\pi^a/t^b\right>[1/t]$.

\newpage

\begin{landscape}
\begin{conjecture}
The construction in this section could be carried over the pro-\'etale sites \footnote{For more on the foundations here for the pro-\'etale sites see \cite[Theorem 2.9.9, Remark 2.9.10]{12Ked1}}. Consider the corresponding projective map $g:X_{\text{pro-\'etale}}\rightarrow X_{\text{\'et}}$ and the the projective map $f:X_{\text{pro-\'etale}}\rightarrow X$. Then we have the following two strictly exact long exact sequences 
\footnote{We should have the corresponding naturality taking into the following form:
\[\tiny
\xymatrix@C+0.4pc@R+0pc{
0\ar[r]\ar[r]\ar[r] \ar[r] &\mathrm{deRham}^\text{degreenumber}_{*[\widehat{\mathcal{O}}]^\text{degreenumber}/*,\mathrm{topo},Z}\ar[d]\ar[d]\ar[d] \ar[d] \ar[r]^\partial\ar[r]\ar[r] \ar[r] &\mathrm{deRham}^\text{degreenumber}_{X[\widehat{\mathcal{O}}]^\text{degreenumber}/X,\mathrm{topo},Z} \ar[d]\ar[d]\ar[d] \ar[d]\ar[r]^\partial\ar[r]\ar[r] \ar[r] &\mathrm{deRham}^\text{degreenumber}_{X[\widehat{\mathcal{O}}]^\text{degreenumber}/X,\mathrm{topo},Z}{\otimes}f^{-1} \mathrm{deRham}^1_{X,\mathrm{topo}} \ar[d]\ar[d]\ar[d] \ar[d]\ar[r]^\partial\ar[r]\ar[r] \ar[r]&...\\
0\ar[r]\ar[r]\ar[r] \ar[r] &\mathrm{deRham}^\text{degreenumber}_{*[\widehat{\mathcal{O}}]^\text{degreenumber}/*,\mathrm{topo},Z} \ar[r]^\partial\ar[r]\ar[r] \ar[r] &\mathrm{deRham}^\text{degreenumber}_{X[\widehat{\mathcal{O}}]^\text{degreenumber}/X_{\text{\'et}},\mathrm{topo},Z} \ar[r]^\partial\ar[r]\ar[r] \ar[r] &\mathrm{deRham}^\text{degreenumber}_{X[\widehat{\mathcal{O}}]^\text{degreenumber}/X,\mathrm{topo},Z}{\otimes}f^{-1} \mathrm{deRham}^1_{X_{\text{\'et}},\mathrm{topo}} \ar[r]^\partial\ar[r]\ar[r] \ar[r] &....
}
\]
}:
\[\tiny
\xymatrix@C+0.4pc@R+0pc{
0\ar[r]\ar[r]\ar[r] \ar[r] &\mathrm{deRham}^\text{degreenumber}_{*[\widehat{\mathcal{O}}]^\text{degreenumber}/*,\mathrm{topo},Z} \ar[r]^\partial\ar[r]\ar[r] \ar[r] &\mathrm{deRham}^\text{degreenumber}_{X[\widehat{\mathcal{O}}]^\text{degreenumber}/X,\mathrm{topo},Z} \ar[r]^\partial\ar[r]\ar[r] \ar[r] &\mathrm{deRham}^\text{degreenumber}_{X[\widehat{\mathcal{O}}]^\text{degreenumber}/X,\mathrm{topo},Z}{\otimes}f^{-1} \mathrm{deRham}^1_{X,\mathrm{topo}} \ar[r]^\partial\ar[r]\ar[r] \ar[r]&...,
}
\]
and
\[\tiny
\xymatrix@C+0.4pc@R+0pc{
0\ar[r]\ar[r]\ar[r] \ar[r] &\mathrm{deRham}^\text{degreenumber}_{*[\widehat{\mathcal{O}}]^\text{degreenumber}/*,\mathrm{topo},Z} \ar[r]^\partial\ar[r]\ar[r] \ar[r] &\mathrm{deRham}^\text{degreenumber}_{X[\widehat{\mathcal{O}}]^\text{degreenumber}/X_{\text{\'et}},\mathrm{topo},Z} \ar[r]^\partial\ar[r]\ar[r] \ar[r] &\mathrm{deRham}^\text{degreenumber}_{X[\widehat{\mathcal{O}}]^\text{degreenumber}/X,\mathrm{topo},Z}{\otimes}f^{-1} \mathrm{deRham}^1_{X_{\text{\'et}},\mathrm{topo}} \ar[r]^\partial\ar[r]\ar[r] \ar[r]&....
}
\]
\end{conjecture}

\begin{conjecture}
Consider the corresponding projective map $g:X_{\text{pro-\'etale}}\rightarrow X_{\text{\'et}}$ and the the projective map $f:X_{\text{pro-\'etale}}\rightarrow X$. Let $M$ be a corresponding $Z$-projective differential crystal spectrum. Then we have the following two strictly exact long exact sequences:
\[\tiny
\xymatrix@C+0.4pc@R+0pc{
M{\otimes}^\mathbb{L}(0\ar[r]\ar[r]\ar[r] \ar[r] &\mathrm{deRham}^\text{degreenumber}_{\mathcal{O}_k[\widehat{\mathcal{O}}]^\text{degreenumber}/\mathcal{O}_k,\mathrm{topo},Z} \ar[r]^\partial\ar[r]\ar[r] \ar[r] &\mathrm{deRham}^\text{degreenumber}_{X[\widehat{\mathcal{O}}]^\text{degreenumber}/X,\mathrm{topo},Z} \ar[r]^\partial\ar[r]\ar[r] \ar[r] &\mathrm{deRham}^\text{degreenumber}_{X[\widehat{\mathcal{O}}]^\text{degreenumber}/X,\mathrm{topo},Z}{\otimes}f^{-1} \mathrm{deRham}^1_{X,\mathrm{topo}} \ar[r]^\partial\ar[r]\ar[r] \ar[r] &...),
}
\]
and
\[\tiny
\xymatrix@C+0.4pc@R+0pc{
M{\otimes}^\mathbb{L}(0\ar[r]\ar[r]\ar[r] \ar[r] &\mathrm{deRham}^\text{degreenumber}_{\mathcal{O}_k[\widehat{\mathcal{O}}]^\text{degreenumber}/\mathcal{O}_k,\mathrm{topo},Z} \ar[r]^\partial\ar[r]\ar[r] \ar[r] &\mathrm{deRham}^\text{degreenumber}_{X[\widehat{\mathcal{O}}]^\text{degreenumber}/X_{\text{\'et}},\mathrm{topo},Z} \ar[r]^\partial\ar[r]\ar[r] \ar[r] &\mathrm{deRham}^\text{degreenumber}_{X[\widehat{\mathcal{O}}]^\text{degreenumber}/X,\mathrm{topo},Z}{\otimes}f^{-1} \mathrm{deRham}^1_{X_{\text{\'et}},\mathrm{topo}} \ar[r]^\partial\ar[r]\ar[r] \ar[r]&...).
}
\]
\end{conjecture}

\end{landscape}

%%%!!!

\newpage

\subsection{Derived $(p,I)$-Complete Derived de Rham complex of Derived Adic Rings}

%\subsection{The Definitions}

\indent The corresponding construction of the derived de Rham complex could be defined for general derived spaces. Along our discussion in the situations of rigid analytic spaces and the corresponding pseudorigid spaces we now focus on the corresponding derived adic rings. 

\begin{setting}
We now fix a bounded morphism of simplicial adic rings $A\rightarrow B$ over $A^*$ where $A^*$ contains a corresponding ring of definition $A^*_0$ which is derived complete with respect to the $(p,I)$-topology and we assume that $A$ is adic and we assume that $(A^*_0,I)$ is a prism in \cite{12BS} namely we at least require that the corresponding $\delta$-structure on the corresponding ring will induce the map $\varphi(.):=.^p+p\delta(.)$ such that we have the situation where $p\in (I,\varphi(I))$. For $A$ or $B$ respectively we assume this contains a subring $A_0$ or $B_0$ (over $A_0^*$) respectively such that we have $A_0$ or $B_0$ respectively is derived complete with respect to the corresponding derived $(p,I)$-topology and we assume that $B=B_0[1/f,f\in I]$ (same for $A$). All the adic rings are assumed to be open mapping. We use the notation $d$ to denote a corresponding primitive element as in \cite[Section 2.3]{12BS} for $A^*$. We are going to assume that $p$ is a topologically nilpotent element.
\end{setting}

\indent As in the corresponding construction in the rigid and pseudocoherent situation in the previous section following \cite[Chapitre 3]{12An1}, \cite{12An2}, \cite[Chapter 2, Chapter 8]{12B1}, \cite[Chapter 1]{12Bei}, \cite[Chapter 5]{12G1}, \cite[Chapter 3, Chapter 4]{12GL}, \cite[Chapitre II, Chapitre III]{12Ill1}, \cite[Chapitre VIII]{12Ill2}, \cite[Section 4]{12Qui}, we have the following: 

\begin{definition}
We start from the corresponding construction of the algebraic $p$-adic derived de Rham complex for a map $A\rightarrow B$. Fix a pair of ring of definitions $A_0,B_0$ in $A,B$ respectively. Then this is the corresponding left Kan extended derived differential complex attached to $B_0/A_0$ which is now denoted by $\mathrm{Kan}_\mathrm{Left}\mathrm{deRham}^\text{degreenumber}_{{B_0}/A_0,\mathrm{alg}}$\footnote{Namely over $A_0$ we take the suitable corresponding left Kan extension, then apply to $B_0$.} from the corresponding $(p,I)$-complete commutative rings. The corresponding cotangent complex associated is defined to be just:
\begin{align}
\mathbb{L}_{B_0/A_0,\mathrm{alg}}:=	\mathrm{deRham}^1_{P^\text{degreenumber}_{B_0}/A_0,\mathrm{alg}}\otimes_{P^\text{degreenumber}_{B_0}} B_0.
\end{align}
The corresponding algebraic Andr\'e-Quillen homologies are defined to be:
\begin{align}
H_{\text{degreenumber},{\mathrm{AQ}},\mathrm{alg}}:=\pi_\text{degreenumber} (\mathbb{L}_{B_0/A_0,\mathrm{alg}}). 	
\end{align}
The corresponding topological Andr\'e-Quillen complex is actually the complete version of the corresponding algebraic ones above by considering the corresponding derived $(p,I)$-completion over the simplicial module structure:
\begin{align}
\mathbb{L}_{B_0/A_0,\mathrm{topo}}:=R\varprojlim \mathrm{Kos}_{(p,I)}	\left((\mathrm{deRham}^1_{P^\text{degreenumber}_{B_0}/A_0,\mathrm{alg}}\otimes_{P^\text{degreenumber}_{B_0}} B_0)\right).
\end{align}
Then we consider the corresponding derived algebraic de Rham complex which is just defined to be:
\begin{align}
\mathrm{Kan}_\mathrm{Left}\mathrm{deRham}^\text{degreenumber}_{-/A_0,\mathrm{alg}}(B_0),\mathrm{Fil}^*_{\mathrm{Kan}_\mathrm{Left}\mathrm{deRham}^\text{degreenumber}_{B_0/A_0,\mathrm{alg}}}.	
\end{align}
We then take the corresponding derived $(p,I)$-completion and we denote that by:
\begin{align}
\mathrm{Kan}_\mathrm{Left}\mathrm{deRham}^\text{degreenumber}_{B_0/A_0,\mathrm{topo}}:=R\varprojlim \mathrm{Kos}_{(p,I)}\left(\mathrm{Kan}_\mathrm{Left}\mathrm{deRham}^\text{degreenumber}_{B_0/A_0,\mathrm{alg}},\mathrm{Fil}^*_{\mathrm{Kan}_\mathrm{Left}\mathrm{deRham}^\text{degreenumber}_{B_0/A_0,\mathrm{alg}}}\right),\\
\mathrm{Kan}_\mathrm{Left}\mathrm{Fil}^*_{\mathrm{Kan}_\mathrm{Left}\mathrm{deRham}^\text{degreenumber}_{B_0/A_0,\mathrm{topo}}}:=R\varprojlim \mathrm{Kos}_{(p,I)}\left(\mathrm{Fil}^*_{\mathrm{Kan}_\mathrm{Left}\mathrm{deRham}^\text{degreenumber}_{B_0/A_0,\mathrm{alg}}}\right).	
\end{align}
Then in the situation that all the rings are classical adic rings we consider the following construction for the map $A\rightarrow B$ by putting:
\begin{align}
\mathbb{L}_{B/A,\mathrm{topo}}:= \mathrm{Colim}_{A_0\rightarrow B_0}\mathbb{L}_{B_0/A_0,\mathrm{topo}}[1/(d)],\\
H_{\text{degreenumber},{\mathrm{AQ}},\mathrm{topo}}:=\pi_\text{degreenumber} (\mathbb{L}_{B/A,\mathrm{topo}}),	\\
\mathrm{Kan}_\mathrm{Left}\mathrm{deRham}^\text{degreenumber}_{B/A,\mathrm{topo}}:=\mathrm{Colim}_{A_0\rightarrow B_0}\mathrm{Kan}_\mathrm{Left}\mathrm{deRham}^\text{degreenumber}_{B_0/A_0,\mathrm{topo}}[1/(d)],\\
\mathrm{Kan}_\mathrm{Left}\mathrm{Fil}^*_{\mathrm{Kan}_\mathrm{Left}\mathrm{deRham}^\text{degreenumber}_{B/A,\mathrm{topo}}}:=\mathrm{Colim}_{A_0\rightarrow B_0}\mathrm{Fil}^*_{\mathrm{Kan}_\mathrm{Left}\mathrm{deRham}^\text{degreenumber}_{B_0/A_0,\mathrm{topo}}}[1/(d)].
\end{align}
Then we need to take the corresponding Hodge-Filtered completion by using the corresponding filtration associated as above to achieve the corresponding Hodge-complete objects in the corresponding filtered $\infty$-categories:
\begin{align}
\mathrm{Kan}_\mathrm{Left}{\mathrm{deRham}}^\text{degreenumber}_{B/A,\mathrm{topo,Hodge}},\mathrm{Fil}^*_{\mathrm{Kan}_\mathrm{Left}{\mathrm{deRham}}^\text{degreenumber}_{B/A,\mathrm{topo,Hodge}}}.	
\end{align}
\end{definition}

\begin{definition}
We define the corresponding finite projective filtered crystals to be the corresponding finite projective module spectra over the topological filtered $E_\infty$-ring $\mathrm{Kan}_\mathrm{Left}\mathrm{deRham}^\text{degreenumber}_{B/A,\mathrm{topo}}$ with the corresponding induced filtrations.	
\end{definition}

\begin{definition}
We define the corresponding almost perfect \footnote{This is the corresponding derived version of pseudocoherence from \cite{12Lu1}, \cite{12Lu2}.} filtered crystals to be the corresponding almost perfect module spectra over the topological filtered $E_\infty$-ring $\mathrm{Kan}_\mathrm{Left}\mathrm{deRham}^\text{degreenumber}_{B/A,\mathrm{topo}}$ with the corresponding induced filtrations.	
\end{definition}

\indent We now consider the corresponding large coefficient local systems over pseudorigid spaces where $Z$ now is a simplicial topological algebra over $\mathbb{Z}_p$\footnote{It is better to assume that it is basically derived completely flat.}.

%%%!!!

\begin{definition}
We define the following $Z$ deformed version of the corresponding complete version of the corresponding Andr\'e-Quillen homology and the corresponding complete version of the corresponding derived de Rham complex. We start from the corresponding construction of the algebraic $p$-adic derived de Rham complex for a map $A\rightarrow B$. Fix a pair of ring of definitions $A_0,B_0$ in $A,B$ respectively. Then this is the corresponding left Kan extended derived differential complex attached to $B_0$ from the corresponding derived $(p,I)$-complete commutative rings which is now denoted by $\mathrm{Kan}_\mathrm{Left}\mathrm{deRham}^\text{degreenumber}_{B_0/A_0,\mathrm{alg}}$. The corresponding cotangent complex associated is defined to be just:
\begin{align}
\mathbb{L}_{B_0/A_0,\mathrm{alg}}:=	\mathrm{deRham}^1_{P^\text{degreenumber}_{B_0}/A_0,\mathrm{alg}}\otimes_{P^\text{degreenumber}_{B_0}} B_0.
\end{align}
The corresponding algebraic Andr\'e-Quillen homologies are defined to be:
\begin{align}
H_{\text{degreenumber},{\mathrm{AQ}},\mathrm{alg}}:=\pi_\text{degreenumber} (\mathbb{L}_{B_0/A_0,\mathrm{alg}}). 	
\end{align}
The corresponding topological Andr\'e-Quillen complex is actually the complete version of the corresponding algebraic ones above by considering the corresponding derived $(p,I)$-completion over the simplicial module structure:
\begin{align}
\mathbb{L}_{B_0/A_0,\mathrm{topo}}:=R\varprojlim \mathrm{Kos}_{(p,I)}	\left((\mathrm{deRham}^1_{P^\text{degreenumber}_{B_0}/A_0,\mathrm{alg}}\otimes_{P^\text{degreenumber}_{B_0}} B_0)\right).
\end{align}
Taking the product with $Z$ we have the corresponding integral version of the topological Andr\'e-Quillen complex:
\begin{align}
\mathbb{L}_{B_0/A_0,\mathrm{topo},Z}:=\mathbb{L}_{B_0/A_0,\mathrm{topo}}{\otimes}^\mathbb{L}_{\mathbb{Z}_p}Z
\end{align}
Then we consider the corresponding derived algebraic de Rham complex which is just defined to be:
\begin{align}
\mathrm{Kan}_\mathrm{Left}\mathrm{deRham}^\text{degreenumber}_{B_0/A_0,\mathrm{alg}},\mathrm{Fil}^*_{\mathrm{Kan}_\mathrm{Left}\mathrm{deRham}^\text{degreenumber}_{B_0/A_0,\mathrm{alg}}}.	
\end{align}
We then take the corresponding derived $(p,I)$-completion and we denote that by:
\begin{align}
\mathrm{Kan}_\mathrm{Left}\mathrm{deRham}^\text{degreenumber}_{B_0/A_0,\mathrm{topo}}:=R\varprojlim_k \mathrm{Kos}_{(p,I)}\left(\mathrm{Kan}_\mathrm{Left}\mathrm{deRham}^\text{degreenumber}_{B_0/A_0,\mathrm{alg}},\mathrm{Fil}^*_{\mathrm{deRham}^\text{degreenumber}_{B_0/A_0,\mathrm{alg}}}\right),\\
\mathrm{Fil}^*_{\mathrm{Kan}_\mathrm{Left}\mathrm{deRham}^\text{degreenumber}_{B_0/A_0,\mathrm{topo}}}:=R\varprojlim_k\left(\mathrm{Fil}^*_{\mathrm{Kan}_\mathrm{Left}\mathrm{deRham}^\text{degreenumber}_{B_0/A_0,\mathrm{alg}}}\right).	
\end{align}
Before considering the corresponding integral version we just consider the corresponding product of these $\mathbb{E}_\infty$-rings with $Z$ to get:
\begin{align}
\mathrm{Kan}_\mathrm{Left}\mathrm{deRham}^\text{degreenumber}_{B_0/A_0,\mathrm{topo},Z}:=\mathrm{Kan}_\mathrm{Left}\mathrm{deRham}^\text{degreenumber}_{B_0/A_0,\mathrm{topo}}{\otimes}^\mathbb{L}_{\mathbb{Z}_p}Z,\\
\mathrm{Fil}^*_{\mathrm{Kan}_\mathrm{Left}\mathrm{deRham}^\text{degreenumber}_{B_0/A_0,\mathrm{topo}},Z}:=\mathrm{Fil}^*_{\mathrm{Kan}_\mathrm{Left}\mathrm{deRham}^\text{degreenumber}_{B_0/A_0,\mathrm{topo}}}{\otimes}^\mathbb{L}_{\mathbb{Z}_p}Z.	
\end{align}
When we have that the corresponding ring $Z$ is also derived $(p,I)$-topologized and commutative, then we can further take the correspding derived $(p,I)$-completion to achieve the corresponding derived completed version:
\begin{align}
\mathrm{Kan}_\mathrm{Left}\mathrm{deRham}^\text{degreenumber}_{B_0/A_0,\mathrm{topo},Z}:=\mathrm{Kan}_\mathrm{Left}\mathrm{deRham}^\text{degreenumber}_{B_0/A_0,\mathrm{topo}}\widehat{\otimes}^\mathbb{L}_{\mathbb{Z}_p}Z,\\
\mathrm{Fil}^*_{\mathrm{Kan}_\mathrm{Left}\mathrm{deRham}^\text{degreenumber}_{B_0/A_0,\mathrm{topo}},Z}:=\mathrm{Fil}^*_{\mathrm{Kan}_\mathrm{Left}\mathrm{deRham}^\text{degreenumber}_{B_0/A_0,\mathrm{topo}}}\widehat{\otimes}^\mathbb{L}_{\mathbb{Z}_p}Z.	
\end{align}
Then in the situation that all the rings are classical adic rings we consider the following construction for the map $A\rightarrow B$ by putting:
\begin{align}
\mathbb{L}_{B/A,\mathrm{topo},Z}:= \mathrm{Colim}_{A_0\rightarrow B_0}\mathbb{L}_{B_0/A_0,\mathrm{topo},Z}[1/(d)],\\
H_{\text{degreenumber},{\mathrm{AQ}},\mathrm{topo},Z}:=\pi_\text{degreenumber} (\mathbb{L}_{B/A,\mathrm{topo},Z}),	\\
\mathrm{Kan}_\mathrm{Left}\mathrm{deRham}^\text{degreenumber}_{B/A,\mathrm{topo},Z}:=\mathrm{Colim}_{A_0\rightarrow B_0}\mathrm{Kan}_\mathrm{Left}\mathrm{deRham}^\text{degreenumber}_{B_0/A_0,\mathrm{topo},Z}[1/(d)],\\
\mathrm{Fil}^*_{\mathrm{Kan}_\mathrm{Left}\mathrm{deRham}^\text{degreenumber}_{B/A,\mathrm{topo}},Z}:=\mathrm{Colim}_{A_0\rightarrow B_0}\mathrm{Fil}^*_{\mathrm{Kan}_\mathrm{Left}\mathrm{deRham}^\text{degreenumber}_{B_0/A_0,\mathrm{topo}},Z}[1/(d)].
\end{align}
Then we need to take the corresponding Hodge-Filtered completion by using the corresponding filtration associated as above to achieve the corresponding Hodge-complete objects in the corresponding filtered $\infty$-categories:
\begin{align}
\mathrm{Kan}_\mathrm{Left}{\mathrm{deRham}}^\text{degreenumber}_{B/A,\mathrm{topo,Hodge},Z},\mathrm{Fil}^*_{\mathrm{Kan}_\mathrm{Left}{\mathrm{deRham}}^\text{degreenumber}_{B/A,\mathrm{topo,Hodge}},Z}.	
\end{align} 	
\end{definition}

%%%!!!

\newpage

\section{Topological Logarithmic Derived De Rham Complexes}

\subsection{Logarithmic Setting for Rigid Analytic Spaces}

%\subsection{The Logarithmic topological Derived de Rham Complex}

\indent In this section we are now going to follow Gabber's construction as in \cite[Chapter 8]{12O} and the extension by \cite[Chapter 5, Chapter 6, Chapter 7]{12B1} to consider the corresponding construction of the topological complete logarithmic cotangent complexes and the corresponding topological logarithmic derived de Rham complex following the corresponding construction for rigid spaces. Certainly the corresponding construction will also follow closely the corresponding \cite[Chapitre 3]{12An1}, \cite{12An2}, \cite[Chapter 2, Chapter 8]{12B1}, \cite[Chapter 1]{12Bei}, \cite[Chapter 5]{12G1}, \cite[Chapter 3, Chapter 4]{12GL}, \cite[Chapitre II, Chapitre III]{12Ill1}, \cite[Chapitre VIII]{12Ill2}, \cite[Section 4]{12Qui} as if we do not have the corresponding log structures, in particular the geometry context underlying is just as in \cite[Chapter 3, Chapter 4]{12GL}. For the geometric foundation of log adic spaces, see \cite{12DLLZ1}. We start from the corresponding construction of the algebraic $p$-adic logarithmic derived de Rham complex for a map $(A,M)\rightarrow (B,N)$ of $p$-complete rings carrying the corresponding log structures, here we use the corresponding notation $(*,?),*=A,B$ to denote the corresponding admissible log rings where $?$ represents the corresponding monoids in the consideration. This is the corresponding derived differential complex attached to the canonical resolution (which we will denote it by the same notation as in the non logarithmic setting) of $(B,N)$:
\begin{align}
(A,M)[(B,N)]^\text{degreenumber},	
\end{align}
which is now denoted by $\mathrm{Kan}_\mathrm{Left}\mathrm{deRham}^\text{degreenumber}_{(B,N)/(A,M),\mathrm{alg}}:=\mathrm{Kan}_\mathrm{Left}\mathrm{deRham}^\text{degreenumber}_{-/(A,M),\mathrm{alg}}((B,N))$ after taking the left Kan extension as in \cite[Chapter 6]{12B1} and applying to the ring $(B,N)$. The corresponding cotangent complex associated is defined to be just:
\begin{align}
\mathbb{L}_{(B,N)/(A,M),\mathrm{alg}}:=	\mathrm{Kan}_\mathrm{Left}\mathrm{deRham}^1_{(A,M)[(B,N)]^\text{degreenumber}/(A,M),\mathrm{alg}}\otimes_{(A,M)[(B,N)]^\text{degreenumber}} (B,N).
\end{align}
The corresponding algebraic Andr\'e-Quillen homologies are defined to be:
\begin{align}
H_{\text{degreenumber},{\mathrm{AQ}},\mathrm{alg}}:=\pi_\text{degreenumber} (\mathbb{L}_{(B,N)/(A,M),\mathrm{alg}}). 	
\end{align}

\indent The corresponding topological Andr\'e-Quillen complex is actually the complete version of the corresponding algebraic ones above by considering the corresponding certain $p$-completion over the simplicial module structure.

\indent Then we consider the corresponding derived algebraic de Rham complex which is just defined to be:
\begin{align}
\mathrm{Kan}_\mathrm{Left}\mathrm{deRham}^\text{degreenumber}_{(B,N)/(A,M),\mathrm{alg}},\mathrm{Fil}^*_{\mathrm{Kan}_\mathrm{Left}\mathrm{deRham}^1_{(B,N)/(A,M),\mathrm{alg}}}.	
\end{align}
We then take the corresponding Banach completion and we denote that by:
\begin{align}
\mathrm{Kan}_\mathrm{Left}\mathrm{deRham}^\text{degreenumber}_{(B,N)/(A,M),\mathrm{topo}},\mathrm{Fil}^*_{\mathrm{Kan}_\mathrm{Left}\mathrm{deRham}^1_{(B,N)/(A,M),\mathrm{topo}}}.	
\end{align}

\indent Then we need to take the corresponding Hodge-Filtered completion by using the corresponding filtration associated as above:
\begin{align}
\mathrm{Kan}_\mathrm{Left}\widehat{\mathrm{deRham}}^\text{degreenumber}_{(B,N)/(A,M),\mathrm{topo}},\mathrm{Fil}^*_{\mathrm{Kan}_\mathrm{Left}\widehat{\mathrm{deRham}}^1_{(B,N)/(A,M),\mathrm{topo}}}.	
\end{align}

This is basically the corresponding analytic and complete version the corresponding algebraic log de Rham complex. Furthermore we allow large coefficients with rigid affinoid algebra $Z$ over $\mathbb{Q}_p$. Therefore we take the corresponding completed tensor product in the following.

\begin{definition}
We define the following $Z$ deformed version of the corresponding complete version of the corresponding logarithmic Andr\'e-Quillen homology and the corresponding complete version of the corresponding logarithmic derived de Rham complex. We start from the corresponding construction of the algebraic $p$-adic logarithmic derived de Rham complex for a map $(A,M)\rightarrow (B,N)$. Fix a pair of ring of definitions $(A_0,M_0),(B_0,N_0)$ in $(A,M),(B,N)$ respectively. Then this is the corresponding derived differential complex attached to the canonical resolution of $(B_0,N_0)$:
\begin{align}
(A_0,M_0)[(B_0,N_0)]^\text{degreenumber},	
\end{align}
which is now denoted by $\mathrm{Kan}_\mathrm{Left}\mathrm{deRham}^\text{degreenumber}_{(A_0,M_0)[(B_0,N_0)]^\text{degreenumber}/(A_0,M_0),\mathrm{alg}}$ after taking the corresponding left Kan extension as in \cite[Chapter 6]{12B1}. The corresponding cotangent complex associated is defined to be just:
\begin{align}
\mathbb{L}_{B_0/(A_0,M_0),\mathrm{alg}}:=	\mathrm{Kan}_\mathrm{Left}\mathrm{deRham}^1_{(A_0,M_0)[B_0,N_0]^\text{degreenumber}/(A_0,M_0),\mathrm{alg}}\otimes_{(A_0,M_0)[B_0,N_0]^\text{degreenumber}} (B_0,N_0).
\end{align}
The corresponding algebraic logarithmic Andr\'e-Quillen homologies are defined to be:
\begin{align}
H_{\text{degreenumber},{\mathrm{AQ}},\mathrm{alg,log}}:=\pi_\text{degreenumber} (\mathbb{L}_{(B_0,N_0)/(A_0,M_0),\mathrm{alg}}). 	
\end{align}
The corresponding topological logarithmic Andr\'e-Quillen complex is actually the complete version of the corresponding algebraic ones above by considering the corresponding derived $p$-completion over the simplicial module structure:
%%%%%%%%%%%%%%%%%%%%%
\begin{align}
&\mathbb{L}_{(B_0,N_0)/(A_0,M_0),\mathrm{topo}}\\
&:=R\varprojlim_k\mathrm{Kos}_{p^k}	\left((\mathrm{Kan}_\mathrm{Left}\mathrm{deRham}^1_{(A_0,M_0)[(B_0,N_0)]^\text{degreenumber}/(A_0,M_0),\mathrm{alg}}\otimes_{(A_0,N_0)[(B_0,N_0)]^\text{degreenumber}} (B_0,N_0))\right).
\end{align}
Taking the product with $\mathcal{O}_Z$ we have the corresponding integral version of the topological logarithmic Andr\'e-Quillen complex:
\begin{align}
\mathbb{L}_{(B_0,N_0)/(A_0,M_0),\mathrm{topo},Z}:=\mathbb{L}_{(B_0,N_0)/(A_0,M_0),\mathrm{topo}}\widehat{\otimes}_{\mathbb{Z}_p}\mathcal{O}_Z
\end{align}
Then we consider the corresponding derived algebraic de Rham complex which is just defined to be:
\begin{align}
\mathrm{Kan}_\mathrm{Left}\mathrm{deRham}^\text{degreenumber}_{(B_0,N_0)/(A_0,M_0),\mathrm{alg}},\mathrm{Fil}^*_{\mathrm{Kan}_\mathrm{Left}\mathrm{deRham}^1_{(B_0,N_0)/(A_0,M_0),\mathrm{alg}}}.	
\end{align}
We then take the corresponding derived $p$-completion and we denote that by:
%%%%%%%%%%%%%%%%%%%%%%%
\begin{align}
&\mathrm{Kan}_\mathrm{Left}\mathrm{deRham}^\text{degreenumber}_{(B_0,N_0)/(A_0,M_0),\mathrm{topo}}\\
&:=R\varprojlim_k\mathrm{Kos}_{p^k}\left(\mathrm{Kan}_\mathrm{Left}\mathrm{deRham}^\text{degreenumber}_{(B_0,N_0)/(A_0,M_0),\mathrm{alg}},\mathrm{Fil}^*_{\mathrm{Kan}_\mathrm{Left}\mathrm{deRham}^\text{degreenumber}_{(B_0,N_0)/(A_0,M_0),\mathrm{alg}}}\right),\\
&\mathrm{Fil}^*_{\mathrm{Kan}_\mathrm{Left}\mathrm{deRham}^\text{degreenumber}_{(B_0,N_0)/(A_0,M_0),\mathrm{topo}}}:=R\varprojlim_k\mathrm{Kos}_{p^k}\left(\mathrm{Fil}^*_{\mathrm{Kan}_\mathrm{Left}\mathrm{deRham}^\text{degreenumber}_{(B_0,N_0)/(A_0,M_0),\mathrm{alg}}}\right).	
\end{align}
Before considering the corresponding integral version we just consider the corresponding product of these $\mathbb{E}_\infty$-rings with $\mathcal{O}_Z$ to get:
\begin{align}
\mathrm{Kan}_\mathrm{Left}\mathrm{deRham}^\text{degreenumber}_{(B_0,N_0)/(A_0,M_0),\mathrm{topo},Z}:=\mathrm{Kan}_\mathrm{Left}\mathrm{deRham}^\text{degreenumber}_{(B_0,N_0)/(A_0,M_0),\mathrm{topo}}\widehat{\otimes}^\mathbb{L}_{\mathbb{Z}_p}\mathcal{O}_Z,\\
\mathrm{Fil}^*_{\mathrm{Kan}_\mathrm{Left}\mathrm{deRham}^\text{degreenumber}_{(B_0,N_0)/(A_0,M_0),\mathrm{topo}},Z}:=\mathrm{Fil}^*_{\mathrm{Kan}_\mathrm{Left}\mathrm{deRham}^\text{degreenumber}_{(B_0,N_0)/(A_0,M_0),\mathrm{topo}}}\widehat{\otimes}^\mathbb{L}_{\mathbb{Z}_p}\mathcal{O}_Z.	
\end{align}
Then we consider the following construction for the map $(A,M)\rightarrow (B,N)$ by putting:

\begin{align}
&\mathbb{L}_{(B_0,N_0)/(A_0,M_0),\mathrm{topo},Z}:= \mathrm{Colim}_{(A_0,M_0)\rightarrow (B_0,N_0)}\mathbb{L}_{(B_0,N_0)/(A_0,M_0),\mathrm{topo},Z}[1/p],\\
&H_{\text{degreenumber},{\mathrm{AQ}},\mathrm{topo},Z}:=\pi_\text{degreenumber} (\mathbb{L}_{(B,N)/(A,M),\mathrm{topo},Z}),	\\
&\mathrm{Kan}_\mathrm{Left}\mathrm{deRham}^\text{degreenumber}_{(B,N)/(A,M),\mathrm{topo},Z}:=\mathrm{Colim}_{(A_0,M_0)\rightarrow (B_0,N_0)}\mathrm{Kan}_\mathrm{Left}\mathrm{deRham}^\text{degreenumber}_{(B_0,N_0)/(A_0,M_0),\mathrm{topo},Z}[1/p],\\
&\mathrm{Fil}^*_{\mathrm{Kan}_\mathrm{Left}\mathrm{deRham}^\text{degreenumber}_{(B,N)/(A,M),\mathrm{topo}},Z}:=\mathrm{Colim}\mathrm{Fil}^*_{\mathrm{Kan}_\mathrm{Left}\mathrm{deRham}^\text{degreenumber}_{(B_0,N_0)/(A_0,M_0),\mathrm{topo}},Z}[1/p].
\end{align}
Then we need to take the corresponding Hodge-Filtered completion by using the corresponding filtration associated as above to achieve the corresponding Hodge-complete objects in the corresponding filtered $\infty$-categories:
\begin{align}
{\mathrm{Kan}_\mathrm{Left}\mathrm{deRham}}^\text{degreenumber}_{(B,N)/(A,M),\mathrm{topo,Hodge},Z},\mathrm{Fil}^*_{{\mathrm{Kan}_\mathrm{Left}\mathrm{deRham}}^1_{(B,N)/(A,M),\mathrm{topo,Hodge}},Z}.	
\end{align} 	
\end{definition}

\begin{definition}
We define the corresponding finite projective filtered crystals to be the corresponding finite projective module spectra over the topological filtered $E_\infty$-ring $\mathrm{Kan}_\mathrm{Left}\mathrm{deRham}^\text{degreenumber}_{(B,N)/(A,M),\mathrm{topo},Z}$ with the corresponding induced filtrations.	
\end{definition}

\begin{definition}
We define the corresponding almost perfect \footnote{This is the corresponding derived version of pseudocoherence from \cite{12Lu1}, \cite{12Lu2}.} filtered crystals to be the corresponding almost perfect module spectra over the topological filtered $E_\infty$-ring $\mathrm{Kan}_\mathrm{Left}\mathrm{deRham}^\text{degreenumber}_{(B,N)/(A,N),\mathrm{topo},Z}$ with the corresponding induced filtrations.	
\end{definition}

\indent The following is inspired by the main Poincar\'e Lemma from \cite[Theorem 1.2]{12GL} in the non-deformed situation. Consider a corresponding smooth log rigid analytic space $X$ over $k/\mathbb{Q}_p$ (where $k$ is a corresponding unramified analytic field which is discretely-valued and the corresponding residue field is finite). Then we have the following:

\begin{conjecture}
Consider the corresponding projective map $g:X_{\text{Kummer-pro-\'etale}}\rightarrow X_{\text{Kummer-\'et}}$ and the the projective map $f:X_{\text{Kummer-pro-\'etale}}\rightarrow X$. Then we have the logarithmic versions of the strictly exact long exact sequences as in the Poincar\'e lemma in the rigid analytic situation. 
\end{conjecture}

%%%!!!

\newpage

\subsection{Logarithmic Setting for Pseudorigid Spaces}

%\subsection{The Logarithmic topological Derived de Rham Complex}

\indent We now consider the analytic logarithmic Andr\'e-Quillen Homology and analytic logarithmic Derived de Rham complex of pseudorigid space, which are very crucial in some development in \cite{12Bel1} in the arithmetic family. Therefore we just investigate the corresponding picture in the corresponding geometric family.

\begin{setting}
We consider now a corresponding morphism taking the corresponding form of $(A,M)\rightarrow (B,N)$ where $(A,M)$ is going to be a log pseudorigid affinoid algebra over $\mathbb{Z}_p$ (where we consider the corresponding foundation in \cite[Chapter 2]{12DLLZ1} for the corresponding log adic rings), and $(B,N)$ is going to be a Kummer perfectoid chart of $A$ in the corresponding Kummer pro-\'etale site of the log pseudorigid affinoid space attached to $A$ (where we consider the corresponding foundation in \cite[Chapter 5]{12DLLZ1} for the corresponding Kummer pro-\'etale sites). As in \cite[Definition 3.1, and below Definition 3.1]{12Bel1} in our situation $A$ is of topologically finite type over $\mathcal{O}_K[[t]]\left<\pi^a/t^b\right>[1/t]$, where $K$ is a discrete valued field containing $\mathbb{Q}_p$ and $(a,b)=1$ \footnote{Certainly one can also consider the characteristic $p$ situation.}.	
\end{setting}

\indent From \cite[Definition 3.1, and below Definition 3.1]{12Bel1} we have the following:

\begin{lemma} We have the following statements:\\
A. $A$ as above is Tate, complete over $\mathcal{O}_K$;\\
B. The ring $A$ has a ring of definition $A_0$ which is of $\mathcal{O}_K$-formally finite type;\\
C. The ring $A_0$ is of topologically finite type over $\mathcal{O}_K[[t]]\left<\pi^a/t^b\right>$.
\end{lemma}

\begin{proof}
Since we did not change the corresponding assumption on the corresponding ring theoretic consideration.
\end{proof}

As in the corresponding construction in the rigid situation in the previous section following \cite[Chapitre 3]{12An1}, \cite{12An2}, \cite[Chapter 5, Chapter 6, Chapter 7]{12B1}, \cite[Chapter 1]{12Bei}, \cite[Chapter 5]{12G1}, \cite[Chapter 3, Chapter 4]{12GL}, \cite[Chapitre II, Chapitre III]{12Ill1}, \cite[Chapitre VIII]{12Ill2}, \cite[Chapter 8]{12O}, \cite[Section 4]{12Qui}. We keep now the following setting:

\begin{setting}
Now as in the above notion we will consider a general map of logarithmic rings $(A,M)\rightarrow (B,N)$ over $\mathcal{O}_K[[t]]\left<\pi^a/t^b\right>[1/t]$ where we have the corresponding map of the associated ring of definitions $(A_0,M_0)\rightarrow (B_0,N_0)$ over $\mathcal{O}_K[[t]]\left<\pi^a/t^b\right>$ such that $A_0,B_0$ are basically $I$-adic (where $\mathcal{O}_K[[t]]\left<\pi^a/t^b\right>$ is $I$-adic).	
\end{setting}

\begin{definition}
We start from the corresponding construction of the algebraic $p$-adic log derived de Rham complex for a map $(A,M)\rightarrow (B,N)$. Fix a pair of ring of definitions $(A_0,M_0),(B_0,N_0)$ in $(A,M),(B,N)$ respectively. Then this is the corresponding derived differential complex attached to the canonical resolution of $(B_0,N_0)$:
\begin{align}
(A_0,M_0)[(B_0,N_0)]^\text{degreenumber},	
\end{align}
which is now denoted by $\mathrm{Kan}_\mathrm{Left}\mathrm{deRham}^\text{degreenumber}_{(B_0,N_0)/(A_0,M_0),\mathrm{alg}}$. This is after taking the left Kan extension as in \cite[Chapter 6]{12B1}. The corresponding logarithmic cotangent complex associated is defined to be just:
\begin{align}
\mathbb{L}_{(B_0,N_0)/(A_0,M_0),\mathrm{alg}}:=	\mathrm{Kan}_\mathrm{Left}\mathrm{deRham}^1_{(A_0,M_0)[(B_0,N_0)]^\text{degreenumber}/(A_0,M_0),\mathrm{alg}}\otimes_{(A_0,M_0)[(B_0,N_0)]^\text{degreenumber}} (B_0,N_0).
\end{align}
The corresponding algebraic logarithmic Andr\'e-Quillen homologies are defined to be:
\begin{align}
H_{\text{degreenumber},{\mathrm{AQ}},\mathrm{alg}}:=\pi_\text{degreenumber} (\mathbb{L}_{(B_0,N_0)/(A_0,M_0),\mathrm{alg}}). 	
\end{align}
The corresponding topological logarithmic Andr\'e-Quillen complex is actually the complete version of the corresponding algebraic ones above by considering the corresponding derived $I$-completion over the simplicial module structure:

\begin{align}
&\mathbb{L}_{(B_0,N_0)/(A_0,M_0),\mathrm{topo}}\\
&:=R\varprojlim_I\mathrm{Kos}_I	\left((\mathrm{Kan}_\mathrm{Left}\mathrm{deRham}^1_{(A_0,M_0)[(B_0,N_0)]^\text{degreenumber}/(A_0,M_0),\mathrm{alg}}\otimes_{(A_0,M_0)[(B_0,N_0)]^\text{degreenumber}} (B_0,N_0))\right).
\end{align}
Then we consider the corresponding logarithmic derived algebraic de Rham complex which is just defined to be:
\begin{align}
\mathrm{Kan}_\mathrm{Left}\mathrm{deRham}^\text{degreenumber}_{(B_0,N_0)/(A_0,M_0),\mathrm{alg}},\mathrm{Fil}^*_{\mathrm{Kan}_\mathrm{Left}\mathrm{deRham}^1_{(B_0,N_0)/(A_0,M_0),\mathrm{alg}}}.	
\end{align}
We then take the corresponding derived $I$-completion and we denote that by:
\begin{align}
&\mathrm{Kan}_\mathrm{Left}\mathrm{deRham}^\text{degreenumber}_{(B_0,N_0)/(A_0,M_0),\mathrm{topo}}:=\\
&R\varprojlim_I\mathrm{Kos}_I\left(\mathrm{Kan}_\mathrm{Left}\mathrm{deRham}^\text{degreenumber}_{(A_0,M_0)[(B_0,N_0)]^\text{degreenumber}/(A_0,M_0),\mathrm{alg}},\mathrm{Fil}^*_{\mathrm{Kan}_\mathrm{Left}\mathrm{deRham}^\text{degreenumber}_{(B_0,N_0)/(A_0,M_0),\mathrm{alg}}}\right),\\
&\mathrm{Fil}^*_{\mathrm{Kan}_\mathrm{Left}\mathrm{deRham}^\text{degreenumber}_{(B_0,N_0)/(A_0,M_0),\mathrm{topo}}}:=R\varprojlim_I\mathrm{Kos}_I\left(\mathrm{Fil}^*_{\mathrm{Kan}_\mathrm{Left}\mathrm{deRham}^\text{degreenumber}_{(B_0,N_0)/(A_0,M_0),\mathrm{alg}}}\right).	
\end{align}
Then we consider the following construction for the map $(A,M)\rightarrow (B,N)$ by putting:
\begin{align}
&\mathbb{L}_{(B_0,N_0)/(A_0,M_0),\mathrm{topo}}:= \mathrm{Colim}_{(A_0,M_0)\rightarrow (B_0,N_0)}\mathbb{L}_{(B_0,N_0)/(A_0,M_0),\mathrm{topo}}[1/t],\\
&H_{\text{degreenumber},{\mathrm{AQ}},\mathrm{topo}}:=\pi_\text{degreenumber} (\mathbb{L}_{(B,M)/(A,M),\mathrm{topo}}),	\\
&\mathrm{Kan}_\mathrm{Left}\mathrm{deRham}^\text{degreenumber}_{(B,N)/(A,M),\mathrm{topo}}:=\mathrm{Colim}_{(A_0,M_0)\rightarrow (B_0,N_0)}\mathrm{Kan}_\mathrm{Left}\mathrm{deRham}^\text{degreenumber}_{(B_0,N_0)/(A_0,M_0),\mathrm{topo}}[1/t],\\
&\mathrm{Fil}^*_{\mathrm{Kan}_\mathrm{Left}\mathrm{deRham}^\text{degreenumber}_{(B,N)/(A,M),\mathrm{topo}}}:=\mathrm{Colim}_{(A_0,M_0)\rightarrow (B_0,N_0)}\mathrm{Fil}^*_{\mathrm{Kan}_\mathrm{Left}\mathrm{deRham}^\text{degreenumber}_{(B_0,N_0)/(A_0,M_0),\mathrm{topo}}}[1/t].
\end{align}
Then we need to take the corresponding Hodge-Filtered completion by using the corresponding filtration associated as above to achieve the corresponding Hodge-complete objects in the corresponding filtered $\infty$-categories:
\begin{align}
{\mathrm{Kan}_\mathrm{Left}\mathrm{deRham}}^\text{degreenumber}_{(B,N)/(A,M),\mathrm{topo,Hodge}},\mathrm{Fil}^*_{{\mathrm{Kan}_\mathrm{Left}\mathrm{deRham}}^\text{degreenumber}_{(B,N)/(A,M),\mathrm{topo,Hodge}}}.	
\end{align}
\end{definition}

\begin{definition}
We define the corresponding finite projective filtered crystals to be the corresponding finite projective module spectra over the topological filtered $E_\infty$-ring $\mathrm{Kan}_\mathrm{Left}\mathrm{deRham}^\text{degreenumber}_{(B,N)/(A,M),\mathrm{topo}}$ with the corresponding induced filtrations.	
\end{definition}

\begin{definition}
We define the corresponding almost perfect \footnote{This is the corresponding derived version of pseudocoherence from \cite{12Lu1}, \cite{12Lu2}.} filtered crystals to be the corresponding almost perfect module spectra over the topological filtered $E_\infty$-ring $\mathrm{Kan}_\mathrm{Left}\mathrm{deRham}^\text{degreenumber}_{(B,N)/(A,M),\mathrm{topo}}$ with the corresponding induced filtrations.	
\end{definition}

\indent We now consider the corresponding large coefficient local systems over logarithmic pseudorigid spaces where $Z$ now is a topological algebra over $\mathbb{Z}_p$\footnote{It is better to assume that it is basically completely flat.}.

\begin{definition}
We define the following $Z$ deformed version of the corresponding complete version of the corresponding logarithmic Andr\'e-Quillen homology and the corresponding complete version of the corresponding logarithmic derived de Rham complex. We start from the corresponding construction of the logarithmic algebraic $p$-adic derived de Rham complex for a map $(A,M)\rightarrow (B,N)$. Fix a pair of ring of definitions $(A_0,M_0),(B_0,N_0)$ in $(A,M),(B,N)$ respectively. Then this is the corresponding derived differential complex attached to the canonical resolution of $(B_0,N_0)$:
\begin{align}
(A_0,M_0)[(B_0,N_0)]^\text{degreenumber},	
\end{align}
which is now denoted by $\mathrm{Kan}_\mathrm{Left}\mathrm{deRham}^\text{degreenumber}_{(B_0,N_0)/(A_0,M_0),\mathrm{alg}}$. The corresponding logarithmic cotangent complex associated is defined to be just:
\begin{align}
\mathbb{L}_{(B_0,N_0)/(A_0,M_0),\mathrm{alg}}:=	\mathrm{Kan}_\mathrm{Left}\mathrm{deRham}^1_{(A_0,M_0)[(B_0,N_0)]^\text{degreenumber}/(A_0,M_0),\mathrm{alg}}\otimes_{(A_0,M_0)[(B_0,N_0)]^\text{degreenumber}} (B_0,N_0).
\end{align}
The corresponding algebraic logarithmic Andr\'e-Quillen homologies are defined to be:
\begin{align}
H_{\text{degreenumber},{\mathrm{AQ}},\mathrm{alg}}:=\pi_\text{degreenumber} (\mathbb{L}_{(B_0,N_0)/(A_0,M_0),\mathrm{alg}}). 	
\end{align}
The corresponding topological logarithmic Andr\'e-Quillen complex is actually the complete version of the corresponding algebraic ones above by considering the corresponding derived $I$-completion over the simplicial module structure:
\begin{align}
&\mathbb{L}_{(B_0,N_0)/(A_0,M_0),\mathrm{topo}}\\
&:=R\varprojlim_k\mathrm{Kos}_I	\left((\mathrm{Kan}_\mathrm{Left}\mathrm{deRham}^1_{(A_0,M_0)[(B_0,N_0)]^\text{degreenumber}/(A_0,M_0),\mathrm{alg}}\otimes_{(A_0,M_0)[(B_0,N_0)]^\text{degreenumber}} (B_0,N_0))\right).
\end{align}
Taking the product with $Z$ we have the corresponding integral version of the topological logarithmic version of the corresponding Andr\'e-Quillen complex:
\begin{align}
\mathbb{L}_{(B_0,N_0)/(A_0,M_0),\mathrm{topo},Z}:=\mathbb{L}_{(B_0,N_0)/(A_0,M_0),\mathrm{topo}}{\otimes}_{\mathbb{Z}_p}\mathcal{O}_Z
\end{align}
Then we consider the corresponding logarithmic derived algebraic de Rham complex which is just defined to be:
\begin{align}
\mathrm{Kan}_\mathrm{Left}\mathrm{deRham}^\text{degreenumber}_{(B_0,N_0)/(A_0,M_0),\mathrm{alg}},\mathrm{Fil}^*_{\mathrm{Kan}_\mathrm{Left}\mathrm{deRham}^1_{(B_0,N_0)/(A_0,M_0),\mathrm{alg}}}.	
\end{align}
We then take the corresponding derived $I$-completion and we denote that by:
\begin{align}
&\mathrm{Kan}_\mathrm{Left}\mathrm{deRham}^\text{degreenumber}_{(B_0,N_0)/(A_0,M_0),\mathrm{topo}}:=\\
&R\varprojlim_k\mathrm{Kos}_I\left(\mathrm{Kan}_\mathrm{Left}\mathrm{deRham}^\text{degreenumber}_{(B_0,N_0)/(A_0,M_0),\mathrm{alg}},\mathrm{Fil}^*_{\mathrm{Kan}_\mathrm{Left}\mathrm{deRham}^\text{degreenumber}_{(B_0,N_0)/(A_0,M_0),\mathrm{alg}}}\right),\\
&\mathrm{Fil}^*_{\mathrm{Kan}_\mathrm{Left}\mathrm{deRham}^\text{degreenumber}_{(B_0,N_0)/(A_0,M_0),\mathrm{topo}}}:=R\varprojlim_k\mathrm{Kos}_I\left(\mathrm{Fil}^*_{\mathrm{Kan}_\mathrm{Left}\mathrm{deRham}^\text{degreenumber}_{(B_0,N_0)/(A_0,M_0),\mathrm{alg}}}\right).	
\end{align}
Before considering the corresponding integral version we just consider the corresponding product of these $\mathbb{E}_\infty$-rings with $Z$ to get:
\begin{align}
&\mathrm{Kan}_\mathrm{Left}\mathrm{deRham}^\text{degreenumber}_{(B_0,N_0)/(A_0,M_0),\mathrm{topo},Z}:=\mathrm{Kan}_\mathrm{Left}\mathrm{deRham}^\text{degreenumber}_{(B_0,N_0)/(A_0,M_0),\mathrm{topo}}{\otimes}^\mathbb{L}_{\mathbb{Z}_p}Z,\\
&\mathrm{Fil}^*_{\mathrm{Kan}_\mathrm{Left}\mathrm{deRham}^\text{degreenumber}_{(B_0,N_0)/(A_0,M_0),\mathrm{topo}},Z}:=\mathrm{Fil}^*_{\mathrm{Kan}_\mathrm{Left}\mathrm{deRham}^\text{degreenumber}_{(B_0,N_0)/(A_0,M_0),\mathrm{topo}}}{\otimes}^\mathbb{L}_{\mathbb{Z}_p}Z.	
\end{align}
Then we consider the following construction for the map $(A,M)\rightarrow (B,N)$ by putting:
\begin{align}
&\mathbb{L}_{(B_0,N_0)/(A_0,M_0),\mathrm{topo},Z}:= \mathrm{Colim}_{(A_0,M_0)\rightarrow (B_0,N_0)}\mathbb{L}_{(B_0,N_0)/(A_0,M_0),\mathrm{topo},Z}[1/t],\\
&H_{\text{degreenumber},{\mathrm{AQ}},\mathrm{topo},Z}:=\pi_\text{degreenumber} (\mathbb{L}_{(B,N)/(A,M),\mathrm{topo},Z}),	\\
&\mathrm{Kan}_\mathrm{Left}\mathrm{deRham}^\text{degreenumber}_{(B,N)/(A,M),\mathrm{topo},Z}:=\mathrm{Colim}_{(A_0,M_0)\rightarrow (B_0,N_0)}\mathrm{Kan}_\mathrm{Left}\mathrm{deRham}^\text{degreenumber}_{(B_0,N_0)/(A_0,M_0),\mathrm{topo},Z}[1/t],\\
&\mathrm{Fil}^*_{\mathrm{Kan}_\mathrm{Left}\mathrm{deRham}^\text{degreenumber}_{(B,N)/(A,M),\mathrm{topo}},Z}:=\mathrm{Colim}_{(A_0,M_0)\rightarrow (B_0,N_0)}\mathrm{Fil}^*_{\mathrm{Kan}_\mathrm{Left}\mathrm{deRham}^\text{degreenumber}_{(B_0,N_0)/(A_0,M_0),\mathrm{topo}},Z}[1/t].
\end{align}
Then we need to take the corresponding Hodge-Filtered completion by using the corresponding filtration associated as above to achieve the corresponding Hodge-complete objects in the corresponding filtered $\infty$-categories:
\begin{align}
{\mathrm{Kan}_\mathrm{Left}\mathrm{deRham}}^\text{degreenumber}_{(B,N)/(A,M),\mathrm{topo,Hodge},Z},\mathrm{Fil}^*_{{\mathrm{Kan}_\mathrm{Left}\mathrm{deRham}}^1_{(B,N)/(A,M),\mathrm{topo,Hodge}},Z}.	
\end{align} 	
\end{definition}

\begin{example}
Now we construct the corresponding pseudorigid analog of the corresponding example in \cite[Example 4.7]{12GL}. Now we consider the corresponding the map:
\begin{align}
\mathbb{Z}_p[[u]]\left<\frac{p^a}{u^b}\right>[1/u]\left<X^\pm_1,X^\pm_2,...,X^\pm_d,X_{d+1},...,X_e\right>\\
\longrightarrow \mathbb{Z}_p[[u]]\left<\frac{p^a}{u^b}\right>[1/u]\left<X^{\pm/p^\infty}_1,X^{\pm/p^\infty}_2,...,X^{\pm/p^\infty}_d,X^{1/p^\infty}_{d+1},...,X^{1/p^\infty}_e\right>	
\end{align}
which could be written as:	
\begin{align}
&\mathbb{Z}_p[[u]]\left<\frac{p^a}{u^b}\right>[1/u]\left<X^\pm_1,X^\pm_2,...,X^\pm_d,X_{d+1},...,X_e\right>\longrightarrow \\
&\mathbb{Z}_p[[u]]\left<\frac{p^a}{u^b}\right>[1/u]\left<X^\pm_1,X^\pm_2,...,X^\pm_d,X_{d+1},...,X_e\right>\left<Y^{\pm/p^\infty}_1,...,Y^{\pm/p^\infty}_d,Y^{1/p^\infty}_{d+1},...,Y^{1/p^\infty}_e\right>\\
&/(X_i-Y_i,i=1,...,e).		
\end{align}
So in our situation the corresponding Hodge complete topological derived de Rham complex will be just:
\begin{align}
\mathbb{Z}_p[[u]]\left<\frac{p^a}{u^b}\right>[1/u]&\left<Y^{\pm/p^\infty}_1,Y^{\pm/p^\infty}_2,...,Y^{\pm/p^\infty}_d,Y^{1/p^\infty}_{d+1},...,Y^{1/p^\infty}_e\right>\\
&[[Z_1,...,Z_e,Z_i=X_i-Y_i,i=1,...,e]].	
\end{align}

\end{example}

\

\indent The following conjectures are literally inspired by the rigid analytic situation in \cite[Theorem 1.2]{12GL}. We assume the spaces are smooth.

\begin{conjecture}
After the whole foundation of \cite{12DLLZ1} and \cite[Theorem 2.9.9, Remark 2.9.10]{12Ked1}, the construction in this section could be carried over the corresponding certain Kummer-pro-\'etale sites where logarithmic perfectoid subdomains form a corresponding basis of neighbourhood of the topology.	
\end{conjecture}

\indent Based on this conjecture one can conjecture the following:

\begin{conjecture}
Consider the corresponding projective map $g:X_{\text{Kummer-pro-\'etale}}\rightarrow X_{\text{Kummer-\'et}}$ and the the projective map $f:X_{\text{Kummer-pro-\'etale}}\rightarrow X$. Then we have the log versions of the strictly exact Poincar\'e long exact sequences as in the pseudorigid analytic situation.
\end{conjecture}

and

\begin{conjecture}
Consider the corresponding projective map $g:X_{\text{Kummer-pro-\'etale}}\rightarrow X_{\text{Kummer-\'et}}$ and the the projective map $f:X_{\text{Kummer-pro-\'etale}}\rightarrow X$. Let $M$ be a corresponding $Z$-projective differential crystal spectrum. Then we have the logarithmic versions of the strictly exact Poincar\'e long exact sequences for $M$ as in the pseudorigid analytic siatution.
\end{conjecture}

\newpage

\subsection{Derived $(p,I)$-Complete Logarithmic Derived de Rham complex}

\indent The corresponding construction of the derived de Rham complex could be defined for general derived spaces. Along our discussion in the situations of rigid analytic spaces and the corresponding pseudorigid spaces we now focus on the corresponding derived  rings. And we consider the corresponding logarithmic setting. 

\begin{setting}
We now fix a bounded morphism of log simplicial topological rings $(A,M)\rightarrow (B,N)$ over $A^*$ where $A^*$ contains a corresponding ring of definition $A^*_0$ which is derived complete with respect to the $(p,I)$-topology and we assume that $(A,M)$ is adic in the sense to be defined just below and we assume that $(A^*_0,I)$ is a prism namely we at least require that the corresponding $\delta$-structure on the corresponding ring will induce the map $\varphi(.):=.^p+p\delta(.)$ such that we have the situation where $p\in (I,\varphi(I))$. For $(A,M)$ or $(B,N)$ respectively we assume this contains a subring $(A_0,M)$ or $(B_0,N_0)$ (over $A^*_0$) respectively such that we have $A_0$ or $B_0$ respectively is derived complete with respect to the corresponding derived $(p,I)$-topology and we assume that $B=B_0[1/f,f\in I]$ (same for $A$) \footnote{One can also invert $p$ to consider the rationality with respect to $p$ as in the rigid analytic situation.}. All the adic rings are assumed to be open mapping. We use the notation $d$ to denote a corresponding primitive element as in \cite[Section 2.3]{12BS} for $A^*$. We are going to assume that $p$ is a topologically nilpotent element.
\end{setting}

\indent As in the corresponding non-logarithmic situation following \cite[Chapter 6]{12B1}\footnote{Here the corresponding left Kan extension happens along the embedding of the free prelog rings to the $\infty$-category of simplicial prelog rings. See \cite[Chapter 6]{12B1} for the construction.} we first consider the category $\mathrm{Alg}_{\mathrm{prelog},\mathrm{smooth},A_0^*}$ of all the smooth prelog algebras over $A_0^*$ and take the corresponding left Kan extension to the corresponding $\infty$-category $\mathrm{Alg}_{\infty,\mathrm{prelog},A_0^*}$ of all the prelog simplicial $A_0^*$-algebras. Then we may take the corresponding derived $(p,I)$-completion to get completed objects. Finally we could basically apply this to the relatively more specific rings in the previous setting. These are basically parallel to the situations where we work over $\mathbb{Z}_p$ and the pseudorigid disc\footnote{Certainly we could work in more general setting over any derived $J$-complete ring, such as in the situation where we do not require that the ring $A^*_0$ is a prism.}.

%%%!!!

\begin{remark}
The corresponding idea behind this is certainly inspired by \cite{12LL} for instance where the authors compared the corresponding prismatic cohomology of some suitable ring $R$ over $A^*_0/I$ for instance and the corresponding derived de Rham cohomology canonically attached to the corresponding morphism, and the corresponding Hodge filtration is compared to the Nygaard filtration. 	
\end{remark}

\indent As in the corresponding construction in the rigid and pseudocoherent situation in the previous section following \cite[Chapitre 3]{12An1}, \cite{12An2}, \cite[Chapter 5, Chapter 6, Chapter 7]{12B1}, \cite[Chapter 1]{12Bei}, \cite[Chapter 5]{12G1}, \cite[Chapter 3, Chapter 4]{12GL}, \cite[Chapitre II, Chapitre III]{12Ill1}, \cite[Chapitre VIII]{12Ill2}, \cite[Chapter 8]{12O}, \cite[Section 4]{12Qui} we give the following definition. 

\begin{definition}
We start from the corresponding construction of the algebraic $p$-adic logarithmic derived de Rham complex for a map $(A,M)\rightarrow (B,N)$. Fix a pair of ring of definitions $(A_0,M_0),(B_0,N_0)$ in $(A,M),(B,N)$ respectively. Then this is the corresponding derived differential complex attached to the cofibrant replacement of $(B_0,N_0)$:
\begin{align}
P^\text{degreenumber}_{(B_0,N_0)},	
\end{align}
which is now denoted by $\mathrm{Kan}_\mathrm{Left}\mathrm{deRham}^\text{degreenumber}_{(B_0,N_0)/(A_0,M_0),\mathrm{alg}}$ after taking the corresponding left Kan extension as in \cite[Chapter 6]{12B1}. The corresponding logarithmic cotangent complex associated is defined to be just:
\begin{align}
\mathbb{L}_{(B_0,N_0)/(A_0,M_0),\mathrm{alg}}:=	\mathrm{Kan}_\mathrm{Left}\mathrm{deRham}^1_{P^\text{degreenumber}_{(B_0,N_0)}/(A_0,M_0),\mathrm{alg}}\otimes_{P^\text{degreenumber}_{(B_0,N_0)}} (B_0,N_0).
\end{align}
The corresponding logarithmic algebraic Andr\'e-Quillen homologies are defined to be:
\begin{align}
H_{\text{degreenumber},{\mathrm{AQ}},\mathrm{alg}}:=\pi_\text{degreenumber} (\mathbb{L}_{(B_0,N_0)/(A_0,M_0),\mathrm{alg}}). 	
\end{align}
The corresponding topological Andr\'e-Quillen complex is actually the complete version of the corresponding algebraic ones above by considering the corresponding derived $(p,I)$-completion over the simplicial module structure:
\begin{align}
\mathbb{L}_{(B_0,N_0)/(A_0,M_0),\mathrm{topo}}:=R\varprojlim \mathrm{Kos}_{(p,I)}	\left((\mathrm{Kan}_\mathrm{Left}\mathrm{deRham}^1_{P^\text{degreenumber}_{(B_0,N_0)}/(A_0,M_0),\mathrm{alg}}\otimes_{P^\text{degreenumber}_{(B_0,N_0)}} (B_0,N_0))\right).
\end{align}
Then we consider the corresponding derived algebraic de Rham complex which is just defined to be:
\begin{align}
\mathrm{Kan}_\mathrm{Left}\mathrm{deRham}^\text{degreenumber}_{(B_0,N_0)/(A_0,M_0),\mathrm{alg}},\mathrm{Fil}^*_{\mathrm{Kan}_\mathrm{Left}\mathrm{deRham}^1_{(B_0,N_0)/(A_0,M_0),\mathrm{alg}}}.	
\end{align}
We then take the corresponding derived $(p,I)$-completion and we denote that by:
\begin{align}
\mathrm{Kan}_\mathrm{Left}&\mathrm{deRham}^\text{degreenumber}_{(B_0,N_0)/(A_0,M_0),\mathrm{topo}}\\
&:=R\varprojlim  \mathrm{Kos}_{(p,I)}\left(\mathrm{Kan}_\mathrm{Left}\mathrm{deRham}^\text{degreenumber}_{(B_0,N_0)/(A_0,M_0),\mathrm{alg}},\mathrm{Fil}^*_{\mathrm{Kan}_\mathrm{Left}\mathrm{deRham}^\text{degreenumber}_{(B_0,N_0)/(A_0,M_0),\mathrm{alg}}}\right),\\
&\mathrm{Fil}^*_{\mathrm{Kan}_\mathrm{Left}\mathrm{deRham}^\text{degreenumber}_{(B_0,N_0)/(A_0,M_0),\mathrm{topo}}}:=R\varprojlim \mathrm{Kos}_{(p,I)} \left(\mathrm{Fil}^*_{\mathrm{Kan}_\mathrm{Left}\mathrm{deRham}^\text{degreenumber}_{(B_0,N_0)/(A_0,M_0),\mathrm{alg}}}\right).	
\end{align}
Then in the situation that all the rings are classical logarithmic adic rings we consider the following construction for the map $A\rightarrow B$ by putting:
\begin{align}
\mathbb{L}_{(B,N)/(A,M),\mathrm{topo}}:= \mathrm{Colim}_{(A_0,M_0)\rightarrow (B_0,N_0)}\mathbb{L}_{(B_0,N_0)/(A_0,M_0),\mathrm{topo}}[1/(d)],\\
H_{\text{degreenumber},{\mathrm{AQ}},\mathrm{topo}}:=\pi_\text{degreenumber} (\mathbb{L}_{(B,N)/(A,M),\mathrm{topo}}),	\\
\mathrm{Kan}_\mathrm{Left}\mathrm{deRham}^\text{degreenumber}_{(B,N)/(A,M),\mathrm{topo}}:=\mathrm{Colim}_{(A_0,M_0)\rightarrow (B_0,N_0)}\mathrm{Kan}_\mathrm{Left}\mathrm{deRham}^\text{degreenumber}_{(B_0,N_0)/(A_0,M_0),\mathrm{topo}}[1/(d)],\\
\mathrm{Fil}^*_{\mathrm{Kan}_\mathrm{Left}\mathrm{deRham}^\text{degreenumber}_{(B,N)/(A,M),\mathrm{topo}}}:=\mathrm{Colim}_{(A_0,M_0)\rightarrow (B_0,N_0)}\mathrm{Fil}^*_{\mathrm{Kan}_\mathrm{Left}\mathrm{deRham}^\text{degreenumber}_{(B_0,N_0)/(A_0,M_0),\mathrm{topo}}}[1/(d)].
\end{align}
Then we need to take the corresponding Hodge-Filtered completion by using the corresponding filtration associated as above to achieve the corresponding Hodge-complete objects in the corresponding filtered $\infty$-categories:
\begin{align}
{\mathrm{Kan}_\mathrm{Left}\mathrm{deRham}}^\text{degreenumber}_{(B,N)/(A,M),\mathrm{topo,Hodge}},\mathrm{Fil}^*_{{\mathrm{Kan}_\mathrm{Left}\mathrm{deRham}}^\text{degreenumber}_{(B,N)/(A,M),\mathrm{topo,Hodge}}}.	
\end{align}
\end{definition}

\begin{definition}
We define the corresponding finite projective filtered crystals to be the corresponding finite projective module spectra over the topological filtered $E_\infty$-ring $\mathrm{Kan}_\mathrm{Left}\mathrm{deRham}^\text{degreenumber}_{(B,N)/{(A,M)},\mathrm{topo}}$ with the corresponding induced filtrations.	
\end{definition}

\begin{definition}
We define the corresponding almost perfect \footnote{This is the corresponding derived version of pseudocoherence from \cite{12Lu1}, \cite{12Lu2}.} filtered crystals to be the corresponding almost perfect module spectra over the topological filtered $E_\infty$-ring $\mathrm{Kan}_\mathrm{Left}\mathrm{deRham}^\text{degreenumber}_{(B,N)/(A,M),\mathrm{topo}}$ with the corresponding induced filtrations.	
\end{definition}

\indent We now consider the corresponding large coefficient local systems over pseudorigid spaces where $Z$ now is a simplicial topological algebra over $\mathbb{Z}_p$\footnote{It is better to assume that it is basically derived completely flat.}.

%%%!!!

\begin{definition}
We define the following $Z$ deformed version of the corresponding complete version of the corresponding logarithmic Andr\'e-Quillen homology and the corresponding complete version of the corresponding logarithmic derived de Rham complex. We start from the corresponding construction of the logarithmic algebraic derived de Rham complex for a map $(A,M)\rightarrow (B,N)$. Fix a pair of ring of definitions $(A_0,M_0),(B_0,N_0)$ in $(A,M),(B,N)$ respectively. Then this is the corresponding derived differential complex attached to the canonical resolution of $(B_0,N_0)$:
\begin{align}
(A_0,M_0)[(B_0,N_0)]^\text{degreenumber},	
\end{align}
which is now denoted by $\mathrm{Kan}_\mathrm{Left}\mathrm{deRham}^\text{degreenumber}_{(B_0,N_0)/(A_0,M_0),\mathrm{alg}}$. The corresponding logarithmic cotangent complex associated is defined to be just:
\begin{align}
\mathbb{L}_{(B_0,N_0)/(A_0,M_0),\mathrm{alg}}\\
&:=	\mathrm{Kan}_\mathrm{Left}\mathrm{deRham}^1_{(A_0,M_0)[(B_0,N_0)]^\text{degreenumber}/(A_0,M_0),\mathrm{alg}}\otimes_{(A_0,M_0)[(B_0,N_0)]^\text{degreenumber}} (B_0,N_0).
\end{align}
The corresponding algebraic logarithmic Andr\'e-Quillen homologies are defined to be:
\begin{align}
H_{\text{degreenumber},{\mathrm{AQ}},\mathrm{alg}}:=\pi_\text{degreenumber} (\mathbb{L}_{(B_0,N_0)/(A_0,M_0),\mathrm{alg}}). 	
\end{align}
The corresponding topological logarithmic Andr\'e-Quillen complex is actually the complete version of the corresponding algebraic ones above by considering the corresponding derived $(p,I)$-completion over the simplicial module structure:
\begin{align}
&\mathbb{L}_{(B_0,N_0)/(A_0,M_0),\mathrm{topo}}\\
&:=R\varprojlim_k\mathrm{Kos}_{(p,I)}	\left((\mathrm{Kan}_\mathrm{Left}\mathrm{deRham}^1_{(A_0,M_0)[(B_0,N_0)]^\text{degreenumber}/(A_0,M_0),\mathrm{alg}}\otimes_{(A_0,M_0)[(B_0,N_0)]^\text{degreenumber}} (B_0,N_0))\right).
\end{align}
Taking the product with $Z$ we have the corresponding integral version of the topological logarithmic version of the corresponding Andr\'e-Quillen complex:
\begin{align}
\mathbb{L}_{(B_0,N_0)/(A_0,M_0),\mathrm{topo},Z}:=\mathbb{L}_{(B_0,N_0)/(A_0,M_0),\mathrm{topo}}{\otimes}_{\mathbb{Z}_p}Z
\end{align}
Then we consider the corresponding logarithmic derived algebraic de Rham complex which is just defined to be:
\begin{align}
\mathrm{Kan}_\mathrm{Left}\mathrm{deRham}^\text{degreenumber}_{(B_0,N_0)/(A_0,M_0),\mathrm{alg}},\mathrm{Fil}^*_{\mathrm{Kan}_\mathrm{Left}\mathrm{deRham}^1_{(B_0,N_0)/(A_0,M_0),\mathrm{alg}}}.	
\end{align}
We then take the corresponding derived $(p,I)$-completion and we denote that by:
\begin{align}
\mathrm{Kan}_\mathrm{Left}&\mathrm{deRham}^\text{degreenumber}_{(B_0,N_0)/(A_0,M_0),\mathrm{topo}}:=\\
&R\varprojlim_k\mathrm{Kos}_{(p,I)}\left(\mathrm{Kan}_\mathrm{Left}\mathrm{deRham}^\text{degreenumber}_{(B_0,N_0)/(A_0,M_0),\mathrm{alg}},\mathrm{Fil}^*_{\mathrm{Kan}_\mathrm{Left}\mathrm{deRham}^\text{degreenumber}_{(B_0,N_0)/(A_0,M_0),\mathrm{alg}}}\right),\\
&\mathrm{Fil}^*_{\mathrm{Kan}_\mathrm{Left}\mathrm{deRham}^\text{degreenumber}_{(B_0,N_0)/(A_0,M_0),\mathrm{topo}}}:=R\varprojlim_k\mathrm{Kos}_{(p,I)}\left(\mathrm{Fil}^*_{\mathrm{Kan}_\mathrm{Left}\mathrm{deRham}^\text{degreenumber}_{(B_0,N_0)/(A_0,M_0),\mathrm{alg}}}\right).	
\end{align}
Before considering the corresponding integral version we just consider the corresponding product of these $\mathbb{E}_\infty$-rings with $Z$ to get:
\begin{align}
\mathrm{Kan}_\mathrm{Left}\mathrm{deRham}^\text{degreenumber}_{(B_0,N_0)/(A_0,M_0),\mathrm{topo},Z}:=\mathrm{Kan}_\mathrm{Left}\mathrm{deRham}^\text{degreenumber}_{(B_0,N_0)/(A_0,M_0),\mathrm{topo}}{\otimes}^\mathbb{L}_{\mathbb{Z}_p}Z,\\
\mathrm{Fil}^*_{\mathrm{Kan}_\mathrm{Left}\mathrm{deRham}^\text{degreenumber}_{(B_0,N_0)/(A_0,M_0),\mathrm{topo}},Z}:=\mathrm{Fil}^*_{\mathrm{Kan}_\mathrm{Left}\mathrm{deRham}^\text{degreenumber}_{(B_0,N_0)/(A_0,M_0),\mathrm{topo}}}{\otimes}^\mathbb{L}_{\mathbb{Z}_p}Z.	
\end{align}
When we have that the corresponding ring $Z$ is also derived $(p,I)$-topologized and commutative, then we can further take the correspding derived $(p,I)$-completion to achieve the corresponding derived completed version:
\begin{align}
\mathrm{Kan}_\mathrm{Left}\mathrm{deRham}^\text{degreenumber}_{(B_0,N_0)/(A_0,M_0),\mathrm{topo},Z}:=\mathrm{Kan}_\mathrm{Left}\mathrm{deRham}^\text{degreenumber}_{(B_0,N_0)/(A_0,M_0),\mathrm{topo}}\widehat{\otimes}^\mathbb{L}_{\mathbb{Z}_p}Z,\\
\mathrm{Fil}^*_{\mathrm{Kan}_\mathrm{Left}\mathrm{deRham}^\text{degreenumber}_{(B_0,N_0)/(A_0,M_0),\mathrm{topo}},Z}:=\mathrm{Fil}^*_{\mathrm{Kan}_\mathrm{Left}\mathrm{deRham}^\text{degreenumber}_{(B_0,N_0)/(A_0,M_0),\mathrm{topo}}}\widehat{\otimes}^\mathbb{L}_{\mathbb{Z}_p}Z.	
\end{align}
Then in the situation that all the rings are classical logarithmic adic rings we consider the following construction for the map $(A,M)\rightarrow (B,N)$ by putting:
\begin{align}
&\mathbb{L}_{(B_0,N_0)/(A_0,M_0),\mathrm{topo},Z}:= \mathrm{Colim}_{(A_0,M_0)\rightarrow (B_0,N_0)}\mathbb{L}_{(B_0,N_0)/(A_0,M_0),\mathrm{topo},Z}[1/(p,I)],\\
&H_{\text{degreenumber},{\mathrm{AQ}},\mathrm{topo},Z}:=\pi_\text{degreenumber} (\mathbb{L}_{(B,N)/(A,M),\mathrm{topo},Z}),	\\
&\mathrm{Kan}_\mathrm{Left}\mathrm{deRham}^\text{degreenumber}_{(B,N)/(A,M),\mathrm{topo},Z}:=\mathrm{Colim}_{(A_0,M_0)\rightarrow (B_0,N_0)}\mathrm{Kan}_\mathrm{Left}\mathrm{deRham}^\text{degreenumber}_{(B,N)/(A_0,M_0),\mathrm{topo},Z}[1/(p,I)],\\
&\mathrm{Fil}^*_{\mathrm{Kan}_\mathrm{Left}\mathrm{deRham}^\text{degreenumber}_{(B,N)/(A,M),\mathrm{topo}},Z}:=\mathrm{Colim}_{(A_0,M_0)\rightarrow (B_0,N_0)}\mathrm{Fil}^*_{\mathrm{Kan}_\mathrm{Left}\mathrm{deRham}^\text{degreenumber}_{(B_0,N_0)/(A_0,M_0),\mathrm{topo}},Z}[1/(p,I)].
\end{align}
Then we need to take the corresponding Hodge-Filtered completion by using the corresponding filtration associated as above to achieve the corresponding Hodge-complete objects in the corresponding filtered $\infty$-categories:
\begin{align}
{\mathrm{Kan}_\mathrm{Left}\mathrm{deRham}}^\text{degreenumber}_{(B,N)/(A,M),\mathrm{topo,Hodge},Z},\mathrm{Fil}^*_{{\mathrm{Kan}_\mathrm{Left}\mathrm{deRham}}^\text{degreenumber}_{(B,N)/(A,M),\mathrm{topo,Hodge}},Z}.	
\end{align} 	
\end{definition}

\newpage

\section{Robba Sheaves and Frobenius Sheaves}

\subsection{Pseudorigid Relative Toric Tower}

\indent In this section we discuss the corresponding Robba sheaves of \cite{12KL1} and \cite{12KL2} over general spaces. Here we consider the corresponding discussion for the pseudorigid situation which is certainly following the corresponding treatment in \cite[Chapter 8]{12KL2}. But we do not really understand if the corresponding generality could be achieved at this moment. So one definitely has to be very careful in the corresponding analysis \footnote{Pseudorigid spaces are actually along our generalization in our mind, which are the first kinds of spaces we would like to study in our project before considering more general topological or functional analytic spaces, namely those general spaces whose structure sheaves of simplicial rings carry topologies or norms such as in \cite{12BBBK} and \cite{12CS2}.}. These kinds of spaces are actually different from the rigid analytic space, although the theories are really related to each other, such as in \cite{12Bel1}, \cite{12Bel2} and \cite{12L}. Following the ideas in the rigid situation in \cite{12KL2}, we consider the towers in the smooth situation in the following:

\begin{setting}
We consider the following towers as in \cite[Chapter 5]{12KL2}. First we let $A_0$ be the ring $\mathcal{O}_K[[u]]\left<\frac{\pi^a}{u^b}\right>[1/u]\left<T_1^\pm,...,T_d^\pm\right>$ and we put $A_0^+$ to be the ring\\ $\mathcal{O}_K[[u]]\left<\frac{\pi^a}{u^b}\right>.\left<T_1^\pm,...,T_d^\pm\right>$. And for the higher level we have the following rings:
\begin{align}
&(A_n,A_n^+)\\
&:=(\mathcal{O}_K(\pi^{1/p^n})[[u]]\left<\frac{\pi^a}{u^b}\right>[1/u](u^{\pm 1/p^n})\left<T_1^{\pm 1/p^n},...,T_d^{\pm 1/p^n}\right>,\\
&\mathcal{O}_K(\pi^{1/p^n})[[u]]\left<\frac{\pi^a}{u^b}\right>[u^{1/p^n}]\left<T_1^{\pm 1/p^n},...,T_d^{\pm 1/p^n}\right>),
\end{align}
which implies that we have:
\[
\xymatrix@C+0pc@R+0pc{
A_0\ar[r] \ar[r] \ar[r] &A_1 \ar[r] \ar[r] \ar[r] &... \ar[r] \ar[r] \ar[r] &A_n \ar[r] \ar[r] \ar[r] &..., \forall n\geq 0. 
}
\]
\end{setting}

\begin{proposition}
The tower above is finite \'etale.	
\end{proposition}

\begin{proof}
Note that the corresponding pseudorigid space is actually formally of finite type over $\mathcal{O}_K$ namely we have that the ring at the zeroth level is just
\begin{center}
 $\mathcal{O}_K[[u,S_1,...,S_d]][1/u]\left<V_1^{\pm},...,V_e^{\pm}\right>$. 
\end{center} 
Therefore the corresponding tower under this sort of presentation will give:
\begin{align}
&(A_n,A_n^+)\\
&:=(\mathcal{O}_K[\pi^{1/p^n}][[u^{1/p^n},S_1^{1/p^n},...,S_d^{1/p^n}]][1/u^{1/p^n}]\left<V_1^{\pm 1/p^n},...,V_e^{\pm 1/p^n}\right>,\\
&\mathcal{O}_K[\pi^{1/p^n}][[u^{1/p^n},S_1^{1/p^n},...,S_d^{1/p^n}]]\left<V_1^{\pm 1/p^n},...,V_e^{\pm 1/p^n}\right>).
\end{align}	
This is certainly a finite \'etale tower.
\end{proof}

\begin{proposition}
The tower above is perfectoid.	
\end{proposition}

\begin{proof}
Note that the corresponding pseudorigid space is actually formally of finite type over $\mathcal{O}_K$ namely we have that the ring at the zeroth level is just
\begin{center}
 $\mathcal{O}_K[[u,S_1,...,S_d]][1/u]\left<V_1^{\pm},...,V_e^{\pm}\right>$.	
\end{center}
Therefore the corresponding tower under this sort of presentation will give:
\begin{align}
&(A_n,A_n^+)\\
&:=(\mathcal{O}_K[\pi^{1/p^n}][[u^{1/p^n},S_1^{1/p^n},...,S_d^{1/p^n}]][1/u^{1/p^n}]\left<V_1^{\pm 1/p^n},...,V_e^{\pm 1/p^n}\right>,\\
&\mathcal{O}_K[\pi^{1/p^n}][[u^{1/p^n},S_1^{1/p^n},...,S_d^{1/p^n}]]\left<V_1^{\pm 1/p^n},...,V_e^{\pm 1/p^n}\right>).
\end{align}	
The corresponding $\infty$-level is just:
\begin{align}
&A_\infty\\
&:=\mathcal{O}_K[\pi^{1/p^\infty}][[u^{1/p^\infty},S_1^{1/p^\infty},...,S_d^{1/p^\infty}]][1/u^{1/p^\infty}]\left<V_1^{\pm 1/p^\infty},...,V_e^{\pm 1/p^\infty}\right>^\wedge,\\
&A_\infty^+\\
&:=\mathcal{O}_K[\pi^{1/p^\infty}][[u^{1/p^\infty},S_1^{1/p^\infty},...,S_d^{1/p^\infty}]]\left<V_1^{\pm 1/p^\infty},...,V_e^{\pm 1/p^\infty}\right>^\wedge.
\end{align}
These are certainly the corresponding perfectoid rings in the sense of \cite[Definition 3.3.1]{12KL2}, also see \cite[Lemma 3.3.28]{12KL2}\footnote{As in \cite[Lemma 3.3.28]{12KL2} one takes the original pseudouniformizer, and goes to a suitably high level to achieve some root $u'$ of the function $X^{p^n}-Xu-u$, and evetually the corresponding element $u'^p$ will divide $p$.}.
\end{proof}

\begin{remark}
Here the corresponding rings are regarded as the corresponding topological rings instead of Banach rings, but we can certainly consider the corresponding Banach ring structure induced from the linear topology namely really the $\infty$-level of this tower is Fontaine perfectoid adic Banach ring which is Tate in the sense of \cite[Definition 3.1.1]{12KL2}.	
\end{remark}

\begin{proposition}
The tower is weakly decompleting.	
\end{proposition}

\begin{proof}
We give the proof where $K=\mathbb{Q}_p$, which is certainly carried over to more general situation. The corresponding tower is actually weakly decompleting once one considers the corresponding presentation as above for the corresponding ring of definition:
\begin{align}
&A_\infty\\
&:=\mathcal{O}_K[\pi^{1/p^\infty}][[u^{1/p^\infty},S_1^{1/p^\infty},...,S_d^{1/p^\infty}]][1/u^{1/p^\infty}]\left<V_1^{\pm 1/p^\infty},...,V_e^{\pm 1/p^\infty}\right>^\wedge,\\
&A_\infty^+\\
&:=\mathcal{O}_K[\pi^{1/p^\infty}][[u^{1/p^\infty},S_1^{1/p^\infty},...,S_d^{1/p^\infty}]]\left<V_1^{\pm 1/p^\infty},...,V_e^{\pm 1/p^\infty}\right>^\wedge.
\end{align}
This transfers to the positive characteristic situation under tilting:
\begin{align}
&R_\infty\\
&:={k}_K[\overline{\pi}^{1/p^\infty}][[\overline{u}^{1/p^\infty},\overline{S}_1^{1/p^\infty},...,\overline{S}_d^{1/p^\infty}]][1/\overline{u}^{1/p^\infty}]\left<\overline{V}_1^{\pm 1/p^\infty},...,\overline{V}_e^{\pm 1/p^\infty}\right>^\wedge,\\
&R_\infty^+\\
&:={k}_K[\overline{\pi}^{1/p^\infty}][[\overline{u}^{1/p^\infty},\overline{S}_1^{1/p^\infty},...,\overline{S}_d^{1/p^\infty}]]\left<\overline{V}_1^{\pm 1/p^\infty},...,\overline{V}_e^{\pm 1/p^\infty}\right>^\wedge.
\end{align}	
Then one can finish as in \cite[Lemma 7.1.2]{12KL2} by comparing this with $R_{(H_\text{degreenumber},H_\text{degreenumber}^+)}$.
\end{proof}

\newpage
\subsection{Pseudorigid Logarithmic Relative Toric Tower}

\begin{setting}
We consider the following towers as in \cite[Chapter 5]{12KL2}. First we let $A_0$ be the ring $\mathcal{O}_K[[u]]\left<\frac{\pi^a}{u^b}\right>[1/u]\left<T_1^\pm,...,T_d^\pm,T^-_{d+1},...,T^-_f\right>$ and we put $A_0^+$ to be the ring 
\begin{align}
\mathcal{O}_K[[u]]\left<\frac{\pi^a}{u^b}\right>.\left<T_1^\pm,...,T_d^\pm,T^-_{d+1},...,T^-_f\right>.	
\end{align}
And for the higher level we have the following rings:
\begin{align}
&(A_n,A_n^+)\\
&:=(\mathcal{O}_K(\pi^{1/p^n})[[u]]\left<\frac{\pi^a}{u^b}\right>[1/u](u^{\pm 1/p^n})\left<T_1^{\pm 1/p^n},...,T_d^{\pm 1/p^n},T^{-1/p^n}_{d+1},...,T^{-1/p^n}_f\right>,\\
&\mathcal{O}_K(\pi^{1/p^n})[[u]]\left<\frac{\pi^a}{u^b}\right>[u^{1/p^n}]\left<T_1^{\pm 1/p^n},...,T_d^{\pm 1/p^n},T^{-1/p^n}_{d+1},...,T^{-1/p^n}_f\right>),
\end{align}
which implies that we have:
\[
\xymatrix@C+0pc@R+0pc{
A_0\ar[r] \ar[r] \ar[r] &A_1 \ar[r] \ar[r] \ar[r] &... \ar[r] \ar[r] \ar[r] &A_n \ar[r] \ar[r] \ar[r] &..., \forall n\geq 0. 
}
\]
\end{setting}

\begin{proposition}
The tower above is finite \'etale.	
\end{proposition}

\begin{proof}
Note that the corresponding pseudorigid space is actually formally of finite type over $\mathcal{O}_K$ namely we have that the ring at the zeroth level is just
\begin{center}
 $\mathcal{O}_K[[u,S_1,...,S_d]][1/u]\left<V_1^{\pm},...,V_e^{\pm}\right>$. 
\end{center}
Therefore the corresponding tower under this sort of presentation will give:
\begin{align}
&(A_n,A_n^+)\\
&:=(\mathcal{O}_K[\pi^{1/p^n}][[u^{1/p^n},S_1^{1/p^n},...,S_d^{1/p^n}]][1/u^{1/p^n}]\left<V_1^{\pm 1/p^n},...,V_e^{\pm 1/p^n},T^{-1/p^n}_{e+1},...,T^{-1/p^n}_f\right>,\\
&\mathcal{O}_K[\pi^{1/p^n}][[u^{1/p^n},S_1^{1/p^n},...,S_d^{1/p^n}]]\left<V_1^{\pm 1/p^n},...,V_e^{\pm 1/p^n},T^{-1/p^n}_{d+1},...,T^{-1/p^n}_f\right>).
\end{align}	
This is certainly a finite \'etale tower.
\end{proof}

\begin{proposition}
The tower above is perfectoid.	
\end{proposition}

\begin{proof}
Note that the corresponding pseudorigid space is actually formally of finite type over $\mathcal{O}_K$ namely we have that the ring at the zeroth level is just
\begin{center}
 $\mathcal{O}_K[[u,S_1,...,S_d]][1/u]\left<V_1^{\pm},...,V_e^{\pm}\right>$. 
\end{center} 
Therefore the corresponding tower under this sort of presentation will give:
\begin{align}
&(A_n,A_n^+)\\
&:=(\mathcal{O}_K[\pi^{1/p^n}][[u^{1/p^n},S_1^{1/p^n},...,S_d^{1/p^n}]][1/u^{1/p^n}]\left<V_1^{\pm 1/p^n},...,V_e^{\pm 1/p^n},T^{-1/p^n}_{e+1},...,T^{-1/p^n}_f\right>,\\
&\mathcal{O}_K[\pi^{1/p^n}][[u^{1/p^n},S_1^{1/p^n},...,S_d^{1/p^n}]]\left<V_1^{\pm 1/p^n},...,V_e^{\pm 1/p^n},T^{-1/p^n}_{e+1},...,T^{-1/p^n}_f\right>).
\end{align}	
The corresponding $\infty$-level is just:
\begin{align}
&A_\infty:=\\
&\mathcal{O}_K[\pi^{1/p^\infty}][[u^{1/p^\infty},S_1^{1/p^\infty},...,S_d^{1/p^\infty}]][1/u^{1/p^\infty}]\left<V_1^{\pm 1/p^\infty},...,V_e^{\pm 1/p^\infty},T^{-1/p^\infty}_{e+1},...,T^{-1/p^\infty}_f\right>^\wedge,\\
&A_\infty^+\\
&:=\mathcal{O}_K[\pi^{1/p^\infty}][[u^{1/p^\infty},S_1^{1/p^\infty},...,S_d^{1/p^\infty}]]\left<V_1^{\pm 1/p^\infty},...,V_e^{\pm 1/p^\infty},T^{-1/p^\infty}_{e+1},...,T^{-1/p^\infty}_f\right>^\wedge.
\end{align}
These are certainly the corresponding perfectoid rings in the sense of \cite[Definition 3.3.1]{12KL2}, also see \cite[Lemma 3.3.28]{12KL2}\footnote{As in \cite[Lemma 3.3.28]{12KL2} one takes the original pseudouniformizer, and goes to a suitably high level to achieve some root $u'$ of the function $X^{p^n}-Xu-u$, and evetually the corresponding element $u'^p$ will divide $p$.}.  
\end{proof}

\begin{remark}
Here the corresponding rings are regarded as the corresponding topological rings instead of Banach rings, but we can certainly consider the corresponding Banach ring structure induced from the linear topology namely really the $\infty$-level of this tower is Fontaine perfectoid adic Banach ring which is Tate in the sense of \cite[Definition 3.1.1]{12KL2}.	
\end{remark}

\begin{proposition}
The tower is weakly decompleting.	
\end{proposition}

\begin{proof}
We give the proof where $K=\mathbb{Q}_p$, which is certainly carried over to more general situation. The corresponding tower is actually weakly decompleting once one considers the corresponding presentation as above for the corresponding ring of definition:
\begin{align}
&A_\infty:=\\
&\mathcal{O}_K[\pi^{1/p^\infty}][[u^{1/p^\infty},S_1^{1/p^\infty},...,S_d^{1/p^\infty}]][1/u^{1/p^\infty}]\left<V_1^{\pm 1/p^\infty},...,V_e^{\pm 1/p^\infty},T^{-1/p^\infty}_{e+1},...,T^{-1/p^\infty}_f\right>^\wedge,\\
&A_\infty^+\\
&:=\mathcal{O}_K[\pi^{1/p^\infty}][[u^{1/p^\infty},S_1^{1/p^\infty},...,S_d^{1/p^\infty}]]\left<V_1^{\pm 1/p^\infty},...,V_e^{\pm 1/p^\infty},T^{-1/p^\infty}_{e+1},...,T^{-1/p^\infty}_f\right>^\wedge.
\end{align}
This transfers to the positive characteristic situation under tilting:
\begin{align}
&R_\infty:=\\
&{k}_K[\overline{\pi}^{1/p^\infty}][[\overline{u}^{1/p^\infty},\overline{S}_1^{1/p^\infty},...,\overline{S}_d^{1/p^\infty}]][1/\overline{u}^{1/p^\infty}]\left<\overline{V}_1^{\pm 1/p^\infty},...,\overline{V}_e^{\pm 1/p^\infty},T^{-1/p^\infty}_{e+1},...,T^{-1/p^\infty}_f\right>^\wedge,\\
&R_\infty^+\\
&:={k}_K[\overline{\pi}^{1/p^\infty}][[\overline{u}^{1/p^\infty},\overline{S}_1^{1/p^\infty},...,\overline{S}_d^{1/p^\infty}]]\left<\overline{V}_1^{\pm 1/p^\infty},...,\overline{V}_e^{\pm 1/p^\infty},T^{-1/p^\infty}_{e+1},...,T^{-1/p^\infty}_f\right>^\wedge.
\end{align}	
Then one can finish as in \cite[Lemma 7.1.2]{12KL2} by comparing this with $R_{(H_\text{degreenumber},H_\text{degreenumber}^+)}$.
\end{proof}

\newpage
\subsection{Robba Sheaves over Pseudorigid Spaces}

\indent Now we consider Kedlaya-Liu's Robba sheaves, which are also some motivic and functorial construction. Our space is actually regarded over $\mathbb{Z}_p$, but note that the corresponding pseudorigid affinoid could be defined over arbitrary integral ring $\mathcal{O}_K$ of some analytic field $K$. Let $X$ be a pseudorigid space over $\mathcal{O}_K$, but we regard this as a corresponding Tate adic space over $\mathbb{Z}_p$. Now we consider the corresponding pro-\'etale site of $X$ which we denote it as $X_\text{pro\'et}$.

%%%!!!

\begin{definition}
We now apply the corresponding definitions of the Robba rings in \cite[Definition 4.1.1]{12KL2} with $E$ therein being just $\mathbb{Q}_p$ to any perfectoid subdomain $(P_\infty,P^+_\infty)$ of $X_\text{pro\'et}$. For each such perfectoid, we promote it to be a Banach ring $(P_\infty,P^+_\infty,\|.\|_{\infty})$, then we have the corresponding construction of the Robba rings in \cite[Definition 4.1.1]{12KL2} (we use the notation $\Pi$ to represent the notation $\mathcal{R}$):
\begin{align}
\widetilde{\Pi}^{[s,r]}_{(P_\infty,P^+_\infty,\|.\|_{\infty})},\\
\varprojlim_s \widetilde{\Pi}^{[s,r]}_{(P_\infty,P^+_\infty,\|.\|_{\infty})},\\
\varinjlim_r \varprojlim_s \widetilde{\Pi}^{[s,r]}_{(P_\infty,P^+_\infty,\|.\|_{\infty})}.	
\end{align}
Then one just organizes these to be certain presheaves over the site $X_\text{pro\'et}$:
\begin{align}
\widetilde{\Pi}^{[s,r]}_{X,\text{pro\'et}},\\
\varprojlim_s \widetilde{\Pi}^{[s,r]}_{X,\text{pro\'et}},\\
\varinjlim_r \varprojlim_s \widetilde{\Pi}^{[s,r]}_{X,\text{pro\'et}}.	
\end{align}	
However the corresponding construction is not canonical since the promotion to Banach rings locally is not functorial. But we do have the situation that they are actually sheaves due to the fact that we can regard them as sheaves over some preperfectoid spaces under the corresponding tilting of the total spaces, where we use the same notation to denote the sheaves.
\begin{align}
\widetilde{\Pi}^{[s,r]}_{X,\text{pro\'et}},\\
\varprojlim_s \widetilde{\Pi}^{[s,r]}_{X,\text{pro\'et}},\\
\varinjlim_r \varprojlim_s \widetilde{\Pi}^{[s,r]}_{X,\text{pro\'et}}.	
\end{align}	
The corresponding total spaces are defined as:
\begin{align}
\mathrm{Spectrumadic}_{\mathrm{total}}(\widetilde{\Pi}^{[s,r]}_{X,\text{pro\'et}},\widetilde{\Pi}^{[s,r],+}_{X,\text{pro\'et}}),\\
\mathrm{Spectrumadic}_{\mathrm{total}}(\varprojlim_s \widetilde{\Pi}^{[s,r]}_{X,\text{pro\'et}},\varprojlim_s \widetilde{\Pi}^{[s,r],+}_{X,\text{pro\'et}}),\\
\mathrm{Spectrumadic}_{\mathrm{total}}(\varinjlim_r \varprojlim_s \widetilde{\Pi}^{[s,r]}_{X,\text{pro\'et}},\varinjlim_r \varprojlim_s \widetilde{\Pi}^{[s,r],+}_{X,\text{pro\'et}}).	
\end{align}	
\end{definition}

\begin{definition}
As in \cite[Chapter 4.3 and Chapter 8]{12KL2}, consider the corresponding total space:
\begin{align}
\mathrm{Spectrumadic}_{\mathrm{total}}(\widetilde{\Pi}^{[s,r]}_{X,\text{pro\'et}},\widetilde{\Pi}^{[s,r],+}_{X,\text{pro\'et}}).	
\end{align}
We define the corresponding $\varphi$-sheaves (where $\varphi$ is Frobenius lifting from $p$-th power Frobenius coming from the characteristic $p$ rings encoded in the construction recalled above from \cite{12KL2}) to be sheaves locally attached to \'etale-stably pseudocoherent sheaves carrying the corresponding semilinear Frobenius morphisms realizing the isomorphisms under pullback. Then over 

\begin{align}
\mathrm{Spectrumadic}_{\mathrm{total}}(\varprojlim_s \widetilde{\Pi}^{[s,r]}_{X,\text{pro\'et}},\varprojlim_s \widetilde{\Pi}^{[s,r],+}_{X,\text{pro\'et}})	
\end{align}
we define the similar pseudocoherent sheaves (complete with respect to the natural topology) locally base change to some sheaves in the previous kind. Finally we define the similar pseudocoherent sheaves (complete with respect to the natural topology) over 
\begin{align}
\mathrm{Spectrumadic}_{\mathrm{total}}(\varinjlim_r \varprojlim_s \widetilde{\Pi}^{[s,r]}_{X,\text{pro\'et}},\varinjlim_r \varprojlim_s \widetilde{\Pi}^{[s,r],+}_{X,\text{pro\'et}}),
\end{align}
locally base change to some sheaves in the previous kind\footnote{Here it is certainly not expect to be the case where we could have uniform radius $r>0$.}.

\end{definition}

\newpage
\subsection{Robba Sheaves over $(p,I)$-Adic Spaces}

\indent Now we consider Kedlaya-Liu's Robba sheaves over more general adic spaces, which is also some motivic and functorial construction. Let $X$ be a Tate adic space over $\mathbb{Z}_p$ where $p$ is assumed to be topologically nilpotent. Now we consider the corresponding pro-\'etale site of $X$ which we denote it as $X_\text{pro\'et}$.

\begin{definition}
We now apply the corresponding definitions of the Robba rings in \cite[Definition 4.1.1]{12KL2} with $E$ therein being just $\mathbb{Q}_p$ to any perfectoid subdomain $(P_\infty,P^+_\infty)$ of $X_\text{pro\'et}$. For each such perfectoid, we promote it to be a Banach ring $(P_\infty,P^+_\infty,\|.\|_{\infty})$, then we have the corresponding construction of the Robba rings in \cite[Definition 4.1.1]{12KL2} (we use the notation $\Pi$ to represent the notation $\mathcal{R}$):
\begin{align}
\widetilde{\Pi}^{[s,r]}_{(P_\infty,P^+_\infty,\|.\|_{\infty})},\\
\varprojlim_s \widetilde{\Pi}^{[s,r]}_{(P_\infty,P^+_\infty,\|.\|_{\infty})},\\
\varinjlim_r \varprojlim_s \widetilde{\Pi}^{[s,r]}_{(P_\infty,P^+_\infty,\|.\|_{\infty})}.	
\end{align}
Then one just organizes these to be certain presheaves over the site $X_\text{pro\'et}$:
\begin{align}
\widetilde{\Pi}^{[s,r]}_{X,\text{pro\'et}},\\
\varprojlim_s \widetilde{\Pi}^{[s,r]}_{X,\text{pro\'et}},\\
\varinjlim_r \varprojlim_s \widetilde{\Pi}^{[s,r]}_{X,\text{pro\'et}}.	
\end{align}	
However the corresponding construction is not canonical since the promotion to Banach rings locally is not functorial. But we do have the situation that they are actually sheaves due to the fact that we can regard them as sheaves over some preperfectoid spaces under the corresponding tilting of the total spaces, where we use the same notation to denote the sheaves.
\begin{align}
\widetilde{\Pi}^{[s,r]}_{X,\text{pro\'et}},\\
\varprojlim_s \widetilde{\Pi}^{[s,r]}_{X,\text{pro\'et}},\\
\varinjlim_r \varprojlim_s \widetilde{\Pi}^{[s,r]}_{X,\text{pro\'et}}.	
\end{align}	
The corresponding total spaces are defined as:
\begin{align}
\mathrm{Spectrumadic}_{\mathrm{total}}(\widetilde{\Pi}^{[s,r]}_{X,\text{pro\'et}},\widetilde{\Pi}^{[s,r],+}_{X,\text{pro\'et}}),\\
\mathrm{Spectrumadic}_{\mathrm{total}}(\varprojlim_s \widetilde{\Pi}^{[s,r]}_{X,\text{pro\'et}},\varprojlim_s \widetilde{\Pi}^{[s,r],+}_{X,\text{pro\'et}}),\\
\mathrm{Spectrumadic}_{\mathrm{total}}(\varinjlim_r \varprojlim_s \widetilde{\Pi}^{[s,r]}_{X,\text{pro\'et}},\varinjlim_r \varprojlim_s \widetilde{\Pi}^{[s,r],+}_{X,\text{pro\'et}}).	
\end{align}	
\end{definition}

\begin{definition}
As in \cite[Chapter 4.3 and Chapter 8]{12KL2}, consider the corresponding total space:
\begin{align}
\mathrm{Spectrumadic}_{\mathrm{total}}(\widetilde{\Pi}^{[s,r]}_{X,\text{pro\'et}},\widetilde{\Pi}^{[s,r],+}_{X,\text{pro\'et}}).	
\end{align}
We define the corresponding $\varphi$-sheaves (where $\varphi$ is Frobenius lifting from $p$-th power Frobenius coming from the characteristic $p$ rings encoded in the construction recalled above from \cite{12KL2}) to be sheaves locally attached to \'etale-stably pseudocoherent sheaves carrying the corresponding semilinear Frobenius morphism realizing the isomorphism under pullback. Then over 

\begin{align}
\mathrm{Spectrumadic}_{\mathrm{total}}(\varprojlim_s \widetilde{\Pi}^{[s,r]}_{X,\text{pro\'et}},\varprojlim_s \widetilde{\Pi}^{[s,r],+}_{X,\text{pro\'et}})	
\end{align}
we define the similar pseudocoherent sheaves (complete with respect to the natural topology) locally base change to some sheaves in the previous kind. Finally we define the similar pseudocoherent sheaves (complete with respect to the natural topology) over 
\begin{align}
\mathrm{Spectrumadic}_{\mathrm{total}}(\varinjlim_r \varprojlim_s \widetilde{\Pi}^{[s,r]}_{X,\text{pro\'et}},\varinjlim_r \varprojlim_s \widetilde{\Pi}^{[s,r],+}_{X,\text{pro\'et}}),
\end{align}
locally base change to some sheaves in the previous kind\footnote{Here it is certainly not expected to be the case where we could have uniform radius $r>0$.}.

\end{definition}

\begin{remark}
The corresponding sheaves of the full Robba rings are well-defined since one could regard them as structure sheaves of some total spaces.	
\end{remark}

\newpage

\section{Derived $I$-Complete THH and HH of Derived $I$-Complete Rings}

\subsection{Derived $I$-Complete Objects}

\indent The corresponding THH and HH of derived $I$-complete rings are topological constructions which are very closely related to the corresponding relative $p$-adic motives. And certain derived $I$-complete versions are also relevant in some highly nontrivial way. We make some discussion closely after \cite[Section 2.2, Section 2.3]{12BMS}, \cite{12BS} and \cite[Chapter 3]{12NS}. We now consider the following rings:

\begin{setting}
We consider the derived $I$-complete $\mathbb{E}_1$-rings. For instance one can consider some derived $I$-complete $\mathbb{E}_1$-rings relative to the corresponding integral pseudorigid disc. In the following $I$ will be some two sided ideal in the $\pi_0$ of the base spectrum $R$. We regard all the ring spectra as being in the $\infty$-category of all ring spectra over $R$ and in the derived $\infty$-category $\mathbb{D}(R)$ of all the module spectra over $R$. We regard all the derived $I$-complete ring spectra as being in the $\infty$-category of all derived $I$-complete ring spectra over $R$ and in the derived $\infty$-category $\mathbb{D}_I(R)$ of all the derived $I$-complete module spectra over $R$. We assume that $I$ is central in the noncommutative setting. 	
\end{setting}

\begin{definition}
\indent For any such ring $R$ which is assumed to be $\mathbb{E}_\infty$, we will use the corresponding notations $L\mathrm{THH}(R)$ and $\mathrm{LHH}(R)$ to denote the corresponding left Kan extended THH or the corresponding HH functor from \cite[Chapter 3]{12NS}. And we will use the corresponding notations $L\mathrm{THH}(R)_I$ and $\mathrm{LHH}(R)_I$ to denote the corresponding left Kan extended THH or the corresponding HH functor after taking the corresponding derived $I$-completion as in \cite[Chapter 3]{12NS}.	
\end{definition}

\begin{definition}
\indent For any such ring $R$ which is assumed to be $\mathbb{E}_1$, we will use the corresponding notations $\mathrm{THH}(R)$ and $\mathrm{HH}(R)$ to denote the corresponding THH or the corresponding HH functor from \cite[Chapter 3]{12NS}. And we will use the corresponding notations $\mathrm{THH}(R)_I$ and $\mathrm{HH}(R)_I$ to denote the corresponding THH or the corresponding HH functor after taking the corresponding derived $I$-completion as in \cite[Chapter 3]{12NS} and \cite[Chapter 1 Notation]{12BS}.	
\end{definition}

\begin{example}
For instance one considers the adic rings in \cite[Section 1.4]{12FK}, one takes such a ring $R$ which is complete with respect to a corresponding two sided finitely generated ideal $I$, then the corresponding construction above namely $L\mathrm{THH}(R)_I$ and $L\mathrm{HH}(R)_I$ will also be basically topological versions of the spectra in \cite[Definition 1.3, Definition 1.4]{12ELS} in the corresponding derived nonocommutative deformation theory.
\end{example}

\indent As in \cite[Chapter 5]{12ELS} and more generally, one has the corresponding noncommutative version of the corresponding relative K\"ahler differential complexes $\mathrm{deRham}^\text{degreenumber}_{\mathrm{noncommutative}, B/A}$. For instance if we have the corresponding map $A\rightarrow B$ of adic rings in \cite[Section 1.4]{12FK} we could then take the corresponding derived $I$-completion from $\mathrm{deRham}^\text{degreenumber}_{\mathrm{noncommutative}, B/A}$ to achieve the corresponding topological version $\mathrm{deRham}^\text{degreenumber}_{\mathrm{noncommutative}, B/A,\mathrm{topo}}$.

\begin{remark}
As mentioned in the introduction, the corresponding topological derived $I$-adic version of the corresponding $L$THH and $L$HH spectra should be interesting to study, as over a quasisyntomic site of a quasisyntomic ring $Q$, \cite[Proposition 5.15]{12BMS} shows that the corresponding sheaf $\pi_0\mathrm{HC}^-(-/Q)_p$ will be closely related to the corresponding $p$-adic complete and Hodge-complete derived de Rham sheaf over the same site after taking the corresponding unfolding through the quasiregular semiperfectoids. We do not know if this relation could hold in some sense if we consider derived $I$-adic rings, but one might want to believe that in the situations we considered before the story should be in some sense easier to be established.	
\end{remark}

\newpage

%\chapter{Functional Analytic Theory}

\section{Functional Analytic Andr\'e-Quillen Homology and Topological Derived de Rham Complexes}

\subsection{The Construction from Kelly-Kremnizer-Mukherjee}

\indent We now work in the corresponding foundations from \cite{12BBBK}, \cite{12BBK}, \cite{BBM}, \cite{12BK} and \cite{KKM}, namely what we are going to consider will be literally the following $(\infty,1)$-categories and the associated constructions with some fixed Banach ring $R$ or $\mathbb{F}_1$\footnote{At this moment we do not work over more general rings such as $\mathbb{Z}$, but this is crucial in the corresponding globalization. That being said, working over $\mathbb{F}_1$ is in some sense more general even than these situations.}:
\begin{align}
\mathrm{Object}_{\mathrm{E}_\infty\mathrm{commutativealgebra},\mathrm{Simplicial}}(\mathrm{IndSNorm}_R),\\
\mathrm{Object}_{\mathrm{E}_\infty\mathrm{commutativealgebra},\mathrm{Simplicial}}(\mathrm{Ind}^m\mathrm{SNorm}_R),\\
\mathrm{Object}_{\mathrm{E}_\infty\mathrm{commutativealgebra},\mathrm{Simplicial}}(\mathrm{IndNorm}_R),\\
\mathrm{Object}_{\mathrm{E}_\infty\mathrm{commutativealgebra},\mathrm{Simplicial}}(\mathrm{Ind}^m\mathrm{Norm}_R),\\
\mathrm{Object}_{\mathrm{E}_\infty\mathrm{commutativealgebra},\mathrm{Simplicial}}(\mathrm{IndBan}_R),\\
\mathrm{Object}_{\mathrm{E}_\infty\mathrm{commutativealgebra},\mathrm{Simplicial}}(\mathrm{Ind}^m\mathrm{Ban}_R),
\end{align}
with
\begin{align}
\mathrm{Object}_{\mathrm{E}_\infty\mathrm{commutativealgebra},\mathrm{Simplicial}}(\mathrm{IndSNorm}_{\mathbb{F}_1}),\\
\mathrm{Object}_{\mathrm{E}_\infty\mathrm{commutativealgebra},\mathrm{Simplicial}}(\mathrm{Ind}^m\mathrm{SNorm}_{\mathbb{F}_1}),\\
\mathrm{Object}_{\mathrm{E}_\infty\mathrm{commutativealgebra},\mathrm{Simplicial}}(\mathrm{IndNorm}_{\mathbb{F}_1}),\\
\mathrm{Object}_{\mathrm{E}_\infty\mathrm{commutativealgebra},\mathrm{Simplicial}}(\mathrm{Ind}^m\mathrm{Norm}_{\mathbb{F}_1}),\\
\mathrm{Object}_{\mathrm{E}_\infty\mathrm{commutativealgebra},\mathrm{Simplicial}}(\mathrm{IndBan}_{\mathbb{F}_1}),\\
\mathrm{Object}_{\mathrm{E}_\infty\mathrm{commutativealgebra},\mathrm{Simplicial}}(\mathrm{Ind}^m\mathrm{Ban}_{\mathbb{F}_1}).
\end{align}
We then use the notation $\mathcal{X}$ to denote any of these. We now first discuss the corresponding topological version of  Andr\'e-Quillen Homology and the corresponding topological version of  Derived de Rham complex parallel to \cite[Chapitre 3]{12An1}, \cite{12An2}, \cite[Chapter 2, Chapter 8]{12B1}, \cite[Chapter 1]{12Bei}, \cite[Chapter 5]{12G1}, \cite[Chapter 3, Chapter 4]{12GL}, \cite[Chapitre II, Chapitre III]{12Ill1}, \cite[Chapitre VIII]{12Ill2}, \cite[Section 4]{12Qui}. We would like to start from the corresponding context of \cite[Section 5.2.1, Definition 5.5, Proposition 5.6, Definition 5.8, Definition 5.9]{KKM}\footnote{Also see \cite[Definition 4.4.7, Construction 4.4.10, Theorem 5.3.6, Definition 5.2.4, Corollary 5.3.9]{Ra}.}, and represent the construction for the convenience of the readers \footnote{Here we just present the homotopical contexts in Kelly-Kremnizer-Mukherjee for the general $(\infty,1)$-categories as above.}. Namely as in \cite[Section 5.2.1, Definition 5.5, Proposition 5.6]{KKM} we consider the corresponding cotangent complex associated to any pair object
\begin{center}
 $(A,B)\in \mathcal{X}\times \mathcal{X}_A$ 
\end{center} 
is defined to be (as in \cite[Section 5.2.1, Definition 5.5, Proposition 5.6]{KKM}) just:
\begin{align}
\mathbb{L}_{B/A,\textrm{KKM}}:=	\mathrm{deRham}^1_{A[B]^\text{degreenumber}/A,\textrm{KKM}}\otimes_{A[B]^\text{degreenumber}} B.
\end{align}
The corresponding topological Andr\'e-Quillen homologies (as in \cite[Section 5.2.1, Definition 5.5, Proposition 5.6]{KKM}) are defined to be:
\begin{align}
H_{\text{degreenumber},{\mathrm{AQ}},\textrm{KKM}}:=\pi_\text{degreenumber} (\mathbb{L}_{B/A,\textrm{KKM}}). 	
\end{align}
We then as in the definition \cite[Definition 5.8, Definition 5.9, Section 5.2.1]{KKM} have the corresponding different kinds of $(\infty,1)$-rings of de Rham complexes and we denote that by:
\begin{align}
\mathrm{Kan}_\mathrm{Left}\mathrm{deRham}^\text{degreenumber}_{B/A,\mathrm{functionalanalytic},\textrm{KKM}},\mathrm{Kan}_\mathrm{Left}\mathrm{Fil}^*_{\mathrm{deRham}^\text{degreenumber}_{B/A,\mathrm{functionalanalytic},\textrm{KKM}}}.	
\end{align}
Then we need to take the corresponding Hodge-Filtered completion by using the corresponding filtration associated as above:
\begin{align}
\mathrm{Kan}_\mathrm{Left}\widehat{\mathrm{deRham}}^\text{degreenumber}_{B/A,\mathrm{functionalanalytic},\mathrm{KKM}},\mathrm{Kan}_\mathrm{Left}\mathrm{Fil}^*_{\widehat{\mathrm{deRham}}^\text{degreenumber}_{B/A,\mathrm{functionalanalytic},\mathrm{KKM}}}.	
\end{align}

\newpage

\section{Functional Analytic Derived Prismatic Cohomology and Functional Analytic Derived Perfectoidizations}

\subsection{General Constructions}

\indent Lurie's book \cite[Section 5.5.8]{Lu3} illustrates in a very general way the corresponding constructions on how we could extend in a certain $(\infty,n)$-category the corresponding nonabelian derived functor from a smaller class of generators in some discrete sense to a large $(\infty,n)$-categorical closure by considering enough colimits being sifted. Therefore in our situation by using the foundation in \cite{12BBBK}, \cite{12BBK}, \cite{BBM}, \cite{12BK} and \cite{KKM} one can define and extend from Tate series rings, Stein series rings, formal series rings and dagger series rings (as those in \cite[Section 4.2]{BBM}) into very large $(\infty,1)$-categories in the following sense:
\begin{align}
\mathrm{Object}_{\mathrm{E}_\infty\mathrm{commutativealgebra},\mathrm{Simplicial}}(\mathrm{IndSNorm}_R),\\
\mathrm{Object}_{\mathrm{E}_\infty\mathrm{commutativealgebra},\mathrm{Simplicial}}(\mathrm{Ind}^m\mathrm{SNorm}_R),\\
\mathrm{Object}_{\mathrm{E}_\infty\mathrm{commutativealgebra},\mathrm{Simplicial}}(\mathrm{IndNorm}_R),\\
\mathrm{Object}_{\mathrm{E}_\infty\mathrm{commutativealgebra},\mathrm{Simplicial}}(\mathrm{Ind}^m\mathrm{Norm}_R),\\
\mathrm{Object}_{\mathrm{E}_\infty\mathrm{commutativealgebra},\mathrm{Simplicial}}(\mathrm{IndBan}_R),\\
\mathrm{Object}_{\mathrm{E}_\infty\mathrm{commutativealgebra},\mathrm{Simplicial}}(\mathrm{Ind}^m\mathrm{Ban}_R),
\end{align}
with
\begin{align}
\mathrm{Object}_{\mathrm{E}_\infty\mathrm{commutativealgebra},\mathrm{Simplicial}}(\mathrm{IndSNorm}_{\mathbb{F}_1}),\\
\mathrm{Object}_{\mathrm{E}_\infty\mathrm{commutativealgebra},\mathrm{Simplicial}}(\mathrm{Ind}^m\mathrm{SNorm}_{\mathbb{F}_1}),\\
\mathrm{Object}_{\mathrm{E}_\infty\mathrm{commutativealgebra},\mathrm{Simplicial}}(\mathrm{IndNorm}_{\mathbb{F}_1}),\\
\mathrm{Object}_{\mathrm{E}_\infty\mathrm{commutativealgebra},\mathrm{Simplicial}}(\mathrm{Ind}^m\mathrm{Norm}_{\mathbb{F}_1}),\\
\mathrm{Object}_{\mathrm{E}_\infty\mathrm{commutativealgebra},\mathrm{Simplicial}}(\mathrm{IndBan}_{\mathbb{F}_1}),\\
\mathrm{Object}_{\mathrm{E}_\infty\mathrm{commutativealgebra},\mathrm{Simplicial}}(\mathrm{Ind}^m\mathrm{Ban}_{\mathbb{F}_1}).
\end{align}

\begin{definition}
Since we are going to consider Bhatt-Scholze's derived prismatic functors \cite[Construction 7.6]{12BS}, so we take the formal series over a general prism $(A,I)$ (assume this to be bounded) such that $A/I$ is also Banach. Then we consider the corresponding formal series ring over $A/I$, then take the correponsding projection resolution compactly generated closure of them to get the sub $(\infty,1)$-categories of above $(\infty,1)$-categories, which we are going to denote them as:
\begin{align}
\mathrm{Object}_{\mathrm{E}_\infty\mathrm{commutativealgebra},\mathrm{Simplicial}}(\mathrm{IndSNorm}_{A/I})^{\text{smoothformalseriesclosure}},\\
\mathrm{Object}_{\mathrm{E}_\infty\mathrm{commutativealgebra},\mathrm{Simplicial}}(\mathrm{Ind}^m\mathrm{SNorm}_{A/I})^{\text{smoothformalseriesclosure}},\\
\mathrm{Object}_{\mathrm{E}_\infty\mathrm{commutativealgebra},\mathrm{Simplicial}}(\mathrm{IndNorm}_{A/I})^{\text{smoothformalseriesclosure}},\\
\mathrm{Object}_{\mathrm{E}_\infty\mathrm{commutativealgebra},\mathrm{Simplicial}}(\mathrm{Ind}^m\mathrm{Norm}_{A/I})^{\text{smoothformalseriesclosure}},\\
\mathrm{Object}_{\mathrm{E}_\infty\mathrm{commutativealgebra},\mathrm{Simplicial}}(\mathrm{IndBan}_{A/I})^{\text{smoothformalseriesclosure}},\\
\mathrm{Object}_{\mathrm{E}_\infty\mathrm{commutativealgebra},\mathrm{Simplicial}}(\mathrm{Ind}^m\mathrm{Ban}_{A/I})^{\text{smoothformalseriesclosure}},
\end{align}	
as in \cite[Section 4.2]{BBM} on the corresponding analytification defined from extension from formal series rings.\\
\end{definition}

\indent Consequently one applies this idea directly to Bhatt-Scholze's prismatic construction \cite{12BS}, one directly gets the corresponding functional analytic $(\infty,1)$-categorical prismatic cohomologies.

\newpage

\subsection{$\infty$-Categorical Functional Analytic Prismatic Cohomologies and $\infty$-Categorical Functional Analytic Preperfectoidizations}

\begin{definition}
\indent One can actually define the derived prismatic cohomologies through derived topological Hochschild cohomologies, derived topological period cohomologies and derived topological cyclic cohomologies as in \cite[Section 2.2, Section 2.3]{12BMS}, \cite[Theorem 1.13]{12BS}:
\begin{align}
	\mathrm{Kan}_{\mathrm{Left}}\mathrm{THH},\mathrm{Kan}_{\mathrm{Left}}\mathrm{TP},\mathrm{Kan}_{\mathrm{Left}}\mathrm{TC},
\end{align}
on the following $(\infty,1)$-compactly generated closures of the corresponding polynomials\footnote{Definitely, we need to put certain norms over in some relatively canonical way, as in \cite[Section 4.2]{BBM} one can basically consider rigid ones and dagger ones, and so on. We restrict to the \textit{formal} one. This means that we are going to consider the corresponding $p$-adic topology only such as the norm over $A/I\left<T_1,...,T_n\right>$, while rigid analytic situation is usually over a field such as in \cite{G2} over $B_\mathrm{dR}\left<T_1,...,T_n\right>$ for instance. One can then in the same way define the corresponding:
\begin{align}
	\mathrm{Kan}_{\mathrm{Left}}\mathrm{THH},\mathrm{Kan}_{\mathrm{Left}}\mathrm{TP},\mathrm{Kan}_{\mathrm{Left}}\mathrm{TC},
\end{align}
over 
\begin{align}
\mathrm{Object}_{\mathrm{E}_\infty\mathrm{commutativealgebra},\mathrm{Simplicial}}(\mathrm{IndSNorm}_{B_\mathrm{dR}})^{\mathrm{smoothformalseriesclosure}},\\
\mathrm{Object}_{\mathrm{E}_\infty\mathrm{commutativealgebra},\mathrm{Simplicial}}(\mathrm{Ind}^m\mathrm{SNorm}_{B_\mathrm{dR}})^{\mathrm{smoothformalseriesclosure}},\\
\mathrm{Object}_{\mathrm{E}_\infty\mathrm{commutativealgebra},\mathrm{Simplicial}}(\mathrm{IndNorm}_{B_\mathrm{dR}})^{\mathrm{smoothformalseriesclosure}},\\
\mathrm{Object}_{\mathrm{E}_\infty\mathrm{commutativealgebra},\mathrm{Simplicial}}(\mathrm{Ind}^m\mathrm{Norm}_{B_\mathrm{dR}})^{\mathrm{smoothformalseriesclosure}},\\
\mathrm{Object}_{\mathrm{E}_\infty\mathrm{commutativealgebra},\mathrm{Simplicial}}(\mathrm{IndBan}_{B_\mathrm{dR}})^{\mathrm{smoothformalseriesclosure}},\\
\mathrm{Object}_{\mathrm{E}_\infty\mathrm{commutativealgebra},\mathrm{Simplicial}}(\mathrm{Ind}^m\mathrm{Ban}_{B_\mathrm{dR}})^{\mathrm{smoothformalseriesclosure}},
\end{align}
by
\begin{align}
	&\mathrm{Kan}_{\mathrm{Left}}\mathrm{THH}_{\text{functionalanalytic,KKM},\text{BBM,formalanalytification}}(\mathcal{O}):=\\
	&(\underset{i}{\text{homotopycolimit}}_{\text{sifted},\text{derivedcategory}_{\infty}(B_\mathrm{dR}-\text{Module})}\mathrm{Kan}_{\mathrm{Left}}\mathrm{THH}_{\text{functionalanalytic,KKM}}(\mathcal{O}_i))_\text{BBM,formalanalytification},\\
	&\mathrm{Kan}_{\mathrm{Left}}\mathrm{TP}_{\text{functionalanalytic,KKM},\text{BBM,formalanalytification}}(\mathcal{O}):=\\
	&(\underset{i}{\text{homotopycolimit}}_{\text{sifted},\text{derivedcategory}_{\infty}(B_\mathrm{dR}-\text{Module})}\mathrm{Kan}_{\mathrm{Left}}\mathrm{TP}_{\text{functionalanalytic,KKM}}(\mathcal{O}_i))_\text{BBM,formalanalytification},\\
	&\mathrm{Kan}_{\mathrm{Left}}\mathrm{TC}_{\text{functionalanalytic,KKM},\text{BBM,formalanalytification}}(\mathcal{O}):=\\
	&(\underset{i}{\text{homotopycolimit}}_{\text{sifted},\text{derivedcategory}_{\infty}(B_\mathrm{dR}-\text{Module})}\mathrm{Kan}_{\mathrm{Left}}\mathrm{TC}_{\text{functionalanalytic,KKM}}(\mathcal{O}_i))_\text{BBM,formalanalytification},
\end{align}
where each $\mathcal{O}_i$ is given as some $B_\mathrm{dR}\left<T_1,...,T_n\right>$. See \cite{G2}, \cite{12G1}.} given over $A/I$ with a chosen prism $(A,I)$\footnote{In all the following, we assume this prism to be bounded and satisfy that $A/I$ is Banach.}:

\begin{align}
\mathrm{Object}_{\mathrm{E}_\infty\mathrm{commutativealgebra},\mathrm{Simplicial}}(\mathrm{IndSNorm}_{A/I})^{\mathrm{smoothformalseriesclosure}},\\
\mathrm{Object}_{\mathrm{E}_\infty\mathrm{commutativealgebra},\mathrm{Simplicial}}(\mathrm{Ind}^m\mathrm{SNorm}_{A/I})^{\mathrm{smoothformalseriesclosure}},\\
\mathrm{Object}_{\mathrm{E}_\infty\mathrm{commutativealgebra},\mathrm{Simplicial}}(\mathrm{IndNorm}_{A/I})^{\mathrm{smoothformalseriesclosure}},\\
\mathrm{Object}_{\mathrm{E}_\infty\mathrm{commutativealgebra},\mathrm{Simplicial}}(\mathrm{Ind}^m\mathrm{Norm}_{A/I})^{\mathrm{smoothformalseriesclosure}},\\
\mathrm{Object}_{\mathrm{E}_\infty\mathrm{commutativealgebra},\mathrm{Simplicial}}(\mathrm{IndBan}_{A/I})^{\mathrm{smoothformalseriesclosure}},\\
\mathrm{Object}_{\mathrm{E}_\infty\mathrm{commutativealgebra},\mathrm{Simplicial}}(\mathrm{Ind}^m\mathrm{Ban}_{A/I})^{\mathrm{smoothformalseriesclosure}}.
\end{align}
We call the corresponding functors are derived functional analytic Hochschild cohomologies, derived functional analytic period cohomologies and derived functional analytic cyclic cohomologies, which we are going to denote them as in the following:
\begin{align}
	&\mathrm{Kan}_{\mathrm{Left}}\mathrm{THH}_{\text{functionalanalytic,KKM},\text{BBM,formalanalytification}}(\mathcal{O}):=\\
	&(\underset{i}{\text{homotopycolimit}}_{\text{sifted},\text{derivedcategory}_{\infty}(A/I-\text{Module})}\mathrm{Kan}_{\mathrm{Left}}\mathrm{THH}_{\text{functionalanalytic,KKM}}(\mathcal{O}_i)\\
	&)_\text{BBM,formalanalytification},\\
	&\mathrm{Kan}_{\mathrm{Left}}\mathrm{TP}_{\text{functionalanalytic,KKM},\text{BBM,formalanalytification}}(\mathcal{O}):=\\
	&(\underset{i}{\text{homotopycolimit}}_{\text{sifted},\text{derivedcategory}_{\infty}(A/I-\text{Module})}\mathrm{Kan}_{\mathrm{Left}}\mathrm{TP}_{\text{functionalanalytic,KKM}}(\mathcal{O}_i)\\
	&)_\text{BBM,formalanalytification},\\
	&\mathrm{Kan}_{\mathrm{Left}}\mathrm{TC}_{\text{functionalanalytic,KKM},\text{BBM,formalanalytification}}(\mathcal{O}):=\\
	&(\underset{i}{\text{homotopycolimit}}_{\text{sifted},\text{derivedcategory}_{\infty}(A/I-\text{Module})}\mathrm{Kan}_{\mathrm{Left}}\mathrm{TC}_{\text{functionalanalytic,KKM}}(\mathcal{O}_i)\\
	&)_\text{BBM,formalanalytification},
\end{align}
by writing any object $\mathcal{O}$ as the corresponding colimit 
\begin{center}
$\underset{i}{\text{homotopycolimit}}_\text{sifted}\mathcal{O}_i$.
\end{center}
These are quite large $(\infty,1)$-commutative ring objects in the corresponding $(\infty,1)$-categories for $R=A/I$:

\begin{align}
\mathrm{Object}_{\mathrm{E}_\infty\mathrm{commutativealgebra},\mathrm{Simplicial}}(\mathrm{IndSNorm}_R),\\
\mathrm{Object}_{\mathrm{E}_\infty\mathrm{commutativealgebra},\mathrm{Simplicial}}(\mathrm{Ind}^m\mathrm{SNorm}_R),\\
\mathrm{Object}_{\mathrm{E}_\infty\mathrm{commutativealgebra},\mathrm{Simplicial}}(\mathrm{IndNorm}_R),\\
\mathrm{Object}_{\mathrm{E}_\infty\mathrm{commutativealgebra},\mathrm{Simplicial}}(\mathrm{Ind}^m\mathrm{Norm}_R),\\
\mathrm{Object}_{\mathrm{E}_\infty\mathrm{commutativealgebra},\mathrm{Simplicial}}(\mathrm{IndBan}_R),\\
\mathrm{Object}_{\mathrm{E}_\infty\mathrm{commutativealgebra},\mathrm{Simplicial}}(\mathrm{Ind}^m\mathrm{Ban}_R),
\end{align}
after taking the formal series ring left Kan extension analytification from \cite[Section 4.2]{BBM}, which is defined by taking the left Kan extension to all the $(\infty,1)$-ring objects in the $\infty$-derived category of all $A$-modules from formal series rings over $A$, into:
\begin{align}
\mathrm{Object}_{\mathrm{E}_\infty\mathrm{commutativealgebra},\mathrm{Simplicial}}(\mathrm{IndSNorm}_{\mathbb{F}_1})_A,\\
\mathrm{Object}_{\mathrm{E}_\infty\mathrm{commutativealgebra},\mathrm{Simplicial}}(\mathrm{Ind}^m\mathrm{SNorm}_{\mathbb{F}_1})_A,\\
\mathrm{Object}_{\mathrm{E}_\infty\mathrm{commutativealgebra},\mathrm{Simplicial}}(\mathrm{IndNorm}_{\mathbb{F}_1})_A,\\
\mathrm{Object}_{\mathrm{E}_\infty\mathrm{commutativealgebra},\mathrm{Simplicial}}(\mathrm{Ind}^m\mathrm{Norm}_{\mathbb{F}_1})_A,\\
\mathrm{Object}_{\mathrm{E}_\infty\mathrm{commutativealgebra},\mathrm{Simplicial}}(\mathrm{IndBan}_{\mathbb{F}_1})_A,\\
\mathrm{Object}_{\mathrm{E}_\infty\mathrm{commutativealgebra},\mathrm{Simplicial}}(\mathrm{Ind}^m\mathrm{Ban}_{\mathbb{F}_1})_A.
\end{align}
\end{definition}

\

\begin{remark}
One should actually do this in a more coherent way as in \cite{Ra}, \cite{12NS}, \cite{KKM}, by applying directly the corresponding construction to objects in the corresponding $(\infty,1)$-objects above, even the corresponding $A_\infty$-objects. The constructions here are the corresponding functional analytic counterpart of the corresponding condensed constructions in \cite{M} after \cite{12CS1} and \cite{12CS2}.
\end{remark}

\

\begin{definition}
Then we can in the same fashion consider the corresponding derived prismatic complexes \cite[Construction 7.6]{12BS}\footnote{One just applies \cite[Construction 7.6]{12BS} and then takes the left Kan extensions.} for the commutative algebras as in the above (for a given prism $(A,I)$):
\begin{align}
\mathrm{Kan}_{\mathrm{Left}}\Delta_{?/A},	
\end{align}
by the regular corresponding left Kan extension techniques on the following $(\infty,1)$-compactly generated closures of the corresponding polynomials given over $A/I$ with a chosen prism $(A,I)$:

\begin{align}
\mathrm{Object}_{\mathrm{E}_\infty\mathrm{commutativealgebra},\mathrm{Simplicial}}(\mathrm{IndSNorm}_{A/I})^{\mathrm{smoothformalseriesclosure}},\\
\mathrm{Object}_{\mathrm{E}_\infty\mathrm{commutativealgebra},\mathrm{Simplicial}}(\mathrm{Ind}^m\mathrm{SNorm}_{A/I})^{\mathrm{smoothformalseriesclosure}},\\
\mathrm{Object}_{\mathrm{E}_\infty\mathrm{commutativealgebra},\mathrm{Simplicial}}(\mathrm{IndNorm}_{A/I})^{\mathrm{smoothformalseriesclosure}},\\
\mathrm{Object}_{\mathrm{E}_\infty\mathrm{commutativealgebra},\mathrm{Simplicial}}(\mathrm{Ind}^m\mathrm{Norm}_{A/I})^{\mathrm{smoothformalseriesclosure}},\\
\mathrm{Object}_{\mathrm{E}_\infty\mathrm{commutativealgebra},\mathrm{Simplicial}}(\mathrm{IndBan}_{A/I})^{\mathrm{smoothformalseriesclosure}},\\
\mathrm{Object}_{\mathrm{E}_\infty\mathrm{commutativealgebra},\mathrm{Simplicial}}(\mathrm{Ind}^m\mathrm{Ban}_{A/I})^{\mathrm{smoothformalseriesclosure}}.
\end{align}
We call the corresponding functors functional analytic derived prismatic complexes which we are going to denote that as in the following:
\begin{align}
\mathrm{Kan}_{\mathrm{Left}}\Delta_{?/A,\text{functionalanalytic,KKM},\text{BBM,formalanalytification}}.	
\end{align}
This would mean the following definition:
\begin{align}
&\mathrm{Kan}_{\mathrm{Left}}\Delta_{?/A,\text{functionalanalytic,KKM},\text{BBM,formalanalytification}}(\mathcal{O})\\
&:=	((\underset{i}{\text{homotopycolimit}}_{\text{sifted},\text{derivedcategory}_{\infty}(A/I-\text{Module})}\mathrm{Kan}_{\mathrm{Left}}\Delta_{?/A,\text{functionalanalytic,KKM}}(\mathcal{O}_i))^\wedge\\
&)_\text{BBM,formalanalytification}
\end{align}
\footnote{Before the Ben-Bassat-Mukherjee $p$-adic formal analytification we take the corresponding derived $(p,I)$-completion.}by writing any object $\mathcal{O}$ as the corresponding colimit 
\begin{center}
$\underset{i}{\text{homotopycolimit}}_\text{sifted}\mathcal{O}_i$.
\end{center}
These are quite large $(\infty,1)$-commutative ring objects in the corresponding $(\infty,1)$-categories for $R=A/I$:
\begin{align}
\mathrm{Object}_{\mathrm{E}_\infty\mathrm{commutativealgebra},\mathrm{Simplicial}}(\mathrm{IndSNorm}_R),\\
\mathrm{Object}_{\mathrm{E}_\infty\mathrm{commutativealgebra},\mathrm{Simplicial}}(\mathrm{Ind}^m\mathrm{SNorm}_R),\\
\mathrm{Object}_{\mathrm{E}_\infty\mathrm{commutativealgebra},\mathrm{Simplicial}}(\mathrm{IndNorm}_R),\\
\mathrm{Object}_{\mathrm{E}_\infty\mathrm{commutativealgebra},\mathrm{Simplicial}}(\mathrm{Ind}^m\mathrm{Norm}_R),\\
\mathrm{Object}_{\mathrm{E}_\infty\mathrm{commutativealgebra},\mathrm{Simplicial}}(\mathrm{IndBan}_R),\\
\mathrm{Object}_{\mathrm{E}_\infty\mathrm{commutativealgebra},\mathrm{Simplicial}}(\mathrm{Ind}^m\mathrm{Ban}_R),\\
\end{align}
after taking the formal series ring left Kan extension analytification from \cite[Section 4.2]{BBM}, which is defined by taking the left Kan extension to all the $(\infty,1)$-ring objects in the $\infty$-derived category of all $A$-modules from formal series rings over $A$, into:
\begin{align}
\mathrm{Object}_{\mathrm{E}_\infty\mathrm{commutativealgebra},\mathrm{Simplicial}}(\mathrm{IndSNorm}_{\mathbb{F}_1})_A,\\
\mathrm{Object}_{\mathrm{E}_\infty\mathrm{commutativealgebra},\mathrm{Simplicial}}(\mathrm{Ind}^m\mathrm{SNorm}_{\mathbb{F}_1})_A,\\
\mathrm{Object}_{\mathrm{E}_\infty\mathrm{commutativealgebra},\mathrm{Simplicial}}(\mathrm{IndNorm}_{\mathbb{F}_1})_A,\\
\mathrm{Object}_{\mathrm{E}_\infty\mathrm{commutativealgebra},\mathrm{Simplicial}}(\mathrm{Ind}^m\mathrm{Norm}_{\mathbb{F}_1})_A,\\
\mathrm{Object}_{\mathrm{E}_\infty\mathrm{commutativealgebra},\mathrm{Simplicial}}(\mathrm{IndBan}_{\mathbb{F}_1})_A,\\
\mathrm{Object}_{\mathrm{E}_\infty\mathrm{commutativealgebra},\mathrm{Simplicial}}(\mathrm{Ind}^m\mathrm{Ban}_{\mathbb{F}_1})_A.
\end{align}
\end{definition}

\

\indent Then as in \cite[Definition 8.2]{12BS} we consider the corresponding perfectoidization in this analytic setting.

\begin{definition}
Let $(A,I)$ be a perfectoid prism, and we consider any $\mathrm{E}_\infty$-ring $\mathcal{O}$ in the following
\begin{align}
\mathrm{Object}_{\mathrm{E}_\infty\mathrm{commutativealgebra},\mathrm{Simplicial}}(\mathrm{IndSNorm}_{A/I})^{\mathrm{smoothformalseriesclosure}},\\
\mathrm{Object}_{\mathrm{E}_\infty\mathrm{commutativealgebra},\mathrm{Simplicial}}(\mathrm{Ind}^m\mathrm{SNorm}_{A/I})^{\mathrm{smoothformalseriesclosure}},\\
\mathrm{Object}_{\mathrm{E}_\infty\mathrm{commutativealgebra},\mathrm{Simplicial}}(\mathrm{IndNorm}_{A/I})^{\mathrm{smoothformalseriesclosure}},\\
\mathrm{Object}_{\mathrm{E}_\infty\mathrm{commutativealgebra},\mathrm{Simplicial}}(\mathrm{Ind}^m\mathrm{Norm}_{A/I})^{\mathrm{smoothformalseriesclosure}},\\
\mathrm{Object}_{\mathrm{E}_\infty\mathrm{commutativealgebra},\mathrm{Simplicial}}(\mathrm{IndBan}_{A/I})^{\mathrm{smoothformalseriesclosure}},\\
\mathrm{Object}_{\mathrm{E}_\infty\mathrm{commutativealgebra},\mathrm{Simplicial}}(\mathrm{Ind}^m\mathrm{Ban}_{A/I})^{\mathrm{smoothformalseriesclosure}}.
\end{align}
Then consider the derived prismatic object:
\begin{align}
\mathrm{Kan}_{\mathrm{Left}}\Delta_{?/A,\text{functionalanalytic,KKM},\text{BBM,formalanalytification}}(\mathcal{O}).
\end{align}	
Then as in \cite[Definition 8.2]{12BS} we have the following preperfectoidization:
\begin{align}
&(\mathcal{O})^{\text{preperfectoidization}}\\
&:=\mathrm{Colimit}(\mathrm{Kan}_{\mathrm{Left}}\Delta_{?/A,\text{functionalanalytic,KKM},\text{BBM,formalanalytification}}(\mathcal{O})\rightarrow \\
&\mathrm{Fro}_*\mathrm{Kan}_{\mathrm{Left}}\Delta_{?/A,\text{functionalanalytic,KKM},\text{BBM,formalanalytification}}(\mathcal{O})\\
&\rightarrow \mathrm{Fro}_* \mathrm{Fro}_*\mathrm{Kan}_{\mathrm{Left}}\Delta_{?/A,\text{functionalanalytic,KKM},\text{BBM,formalanalytification}}(\mathcal{O})\rightarrow...)^{\text{BBM,formalanalytification}},	
\end{align}
after taking the formal series ring left Kan extension analytification from \cite[Section 4.2]{BBM}, which is defined by taking the left Kan extension to all the $(\infty,1)$-ring object in the $\infty$-derived category of all $A$-modules from formal series rings over $A$. Then we define the corresponding perfectoidization:
\begin{align}
&(\mathcal{O})^{\text{perfectoidization}}\\
&:=\mathrm{Colimit}(\mathrm{Kan}_{\mathrm{Left}}\Delta_{?/A,\text{functionalanalytic,KKM},\text{BBM,formalanalytification}}(\mathcal{O})\longrightarrow \\
&\mathrm{Fro}_*\mathrm{Kan}_{\mathrm{Left}}\Delta_{?/A,\text{functionalanalytic,KKM},\text{BBM,formalanalytification}}(\mathcal{O})\\
&\longrightarrow \mathrm{Fro}_* \mathrm{Fro}_*\mathrm{Kan}_{\mathrm{Left}}\Delta_{?/A,\text{functionalanalytic,KKM},\text{BBM,formalanalytification}}(\mathcal{O})\longrightarrow...)^{\text{BBM,formalanalytification}}\times A/I.	
\end{align}
Furthermore one can take derived $(p,I)$-completion to achieve the derived $(p,I)$-completed versions:
\begin{align}
\mathcal{O}^\text{preperfectoidization,derivedcomplete}:=(\mathcal{O}^\text{preperfectoidization})^{\wedge},\\
\mathcal{O}^\text{perfectoidization,derivedcomplete}:=\mathcal{O}^\text{preperfectoidization,derivedcomplete}\times A/I.\\
\end{align}
These are large $(\infty,1)$-commutative algebra objects in the corresponding categories as in the above, attached to also large $(\infty,1)$-commutative algebra objects. When we apply this to the corresponding sub-$(\infty,1)$-categories of Banach perfectoid objects in \cite{BMS2}, \cite{GR}, \cite{12KL1}, \cite{12KL2}, \cite{12Ked1}, \cite{12Sch3},  we will recover the corresponding distinguished elemental deformation processes defined in \cite{BMS2}, \cite{GR}, \cite{12KL1}, \cite{12KL2}, \cite{12Ked1}, \cite{12Sch3}.
\end{definition}

\

\begin{remark}
One can then define such ring $\mathcal{O}$ to be \textit{preperfectoid} if we have the equivalence:
\begin{align}
\mathcal{O}^{\text{preperfectoidization}} \overset{\sim}{\longrightarrow}	\mathcal{O}.
\end{align}
One can then define such ring $\mathcal{O}$ to be \textit{perfectoid} if we have the equivalence:
\begin{align}
\mathcal{O}^{\text{preperfectoidization}}\times A/I \overset{\sim}{\longrightarrow}	\mathcal{O}.
\end{align}
	
\end{remark}

\newpage

\section{Functional Analytic Derived Logarithmic Prismatic Cohomology and Functional Analytic Derived Logarithmic Perfectoidizations}

\subsection{Functional Analytic Derived Logarithmic Prismatic Cohomology}

\indent We now extend the previous functional analytic discussion to logarithmic context, which is related to the logarithmic context for instance in \cite[Chapter 8]{12O}, \cite[Chapter 5, Chapter 6, Chapter 7]{12B1} and \cite{12Ko1}.

\begin{definition}
We define the logarithmic version of the $(\infty,1)$-categories we considered in the previous section (let $(B,*)$ be any Banach logarithmic ring or $\mathbb{F}_1$):
\begin{align}
\mathrm{Object}_{\mathrm{E}_\infty\mathrm{commutativealgebra},\mathrm{Simplicial}}(\mathrm{IndSNorm}_B)^{\text{prelog}},\\
\mathrm{Object}_{\mathrm{E}_\infty\mathrm{commutativealgebra},\mathrm{Simplicial}}(\mathrm{Ind}^m\mathrm{SNorm}_B)^{\text{prelog}},\\
\mathrm{Object}_{\mathrm{E}_\infty\mathrm{commutativealgebra},\mathrm{Simplicial}}(\mathrm{IndNorm}_B)^{\text{prelog}},\\
\mathrm{Object}_{\mathrm{E}_\infty\mathrm{commutativealgebra},\mathrm{Simplicial}}(\mathrm{Ind}^m\mathrm{Norm}_B)^{\text{prelog}},\\
\mathrm{Object}_{\mathrm{E}_\infty\mathrm{commutativealgebra},\mathrm{Simplicial}}(\mathrm{IndBan}_B)^{\text{prelog}},\\
\mathrm{Object}_{\mathrm{E}_\infty\mathrm{commutativealgebra},\mathrm{Simplicial}}(\mathrm{Ind}^m\mathrm{Ban}_B)^{\text{prelog}},
\end{align}
where the prelog object means a morphism $M\rightarrow ?$ where	$M$ is
in the following categories of $(\infty,1)$-monoid objects:
\begin{align}
\mathrm{Object}_{\mathrm{E}_\infty\mathrm{monoid},\mathrm{Simplicial}}(\mathrm{IndSNorm}_B),\\
\mathrm{Object}_{\mathrm{E}_\infty\mathrm{monoid},\mathrm{Simplicial}}(\mathrm{Ind}^m\mathrm{SNorm}_B),\\
\mathrm{Object}_{\mathrm{E}_\infty\mathrm{monoid},\mathrm{Simplicial}}(\mathrm{IndNorm}_B),\\
\mathrm{Object}_{\mathrm{E}_\infty\mathrm{monoid},\mathrm{Simplicial}}(\mathrm{Ind}^m\mathrm{Norm}_B),\\
\mathrm{Object}_{\mathrm{E}_\infty\mathrm{monoid},\mathrm{Simplicial}}(\mathrm{IndBan}_B),\\
\mathrm{Object}_{\mathrm{E}_\infty\mathrm{monoid},\mathrm{Simplicial}}(\mathrm{Ind}^m\mathrm{Ban}_B),\\
\end{align}
with $?$ being commutative algebra object in the previous section.
\end{definition}

\begin{assumption}
We are going to consider the possible colimits in the $(\infty,1)$-categories:
\begin{align}
\mathrm{Object}_{\mathrm{E}_\infty\mathrm{commutativealgebra},\mathrm{Simplicial}}(\mathrm{IndSNorm}_B)^{\text{prelog}},\\
\mathrm{Object}_{\mathrm{E}_\infty\mathrm{commutativealgebra},\mathrm{Simplicial}}(\mathrm{Ind}^m\mathrm{SNorm}_B)^{\text{prelog}},\\
\mathrm{Object}_{\mathrm{E}_\infty\mathrm{commutativealgebra},\mathrm{Simplicial}}(\mathrm{IndNorm}_B)^{\text{prelog}},\\
\mathrm{Object}_{\mathrm{E}_\infty\mathrm{commutativealgebra},\mathrm{Simplicial}}(\mathrm{Ind}^m\mathrm{Norm}_B)^{\text{prelog}},\\
\mathrm{Object}_{\mathrm{E}_\infty\mathrm{commutativealgebra},\mathrm{Simplicial}}(\mathrm{IndBan}_B)^{\text{prelog}},\\
\mathrm{Object}_{\mathrm{E}_\infty\mathrm{commutativealgebra},\mathrm{Simplicial}}(\mathrm{Ind}^m\mathrm{Ban}_B)^{\text{prelog}},
\end{align}
which means essentially that all the construction will only apply to the object $\mathcal{M}\rightarrow\mathcal{O}$ which could be written as the colimit of logarithmic formal series rings over $(B,*)$:
\begin{center}
$\underset{i}{\text{homotopycolimit}}_\text{sifted}(\mathcal{M}_i\rightarrow\mathcal{O}_i)$
\end{center}
in certain large enough $\infty$-category:
\begin{align}
\mathrm{Object}_{\mathrm{Simplicial}}(\mathrm{IndSets})\times \mathrm{Object}_{\mathrm{Simplicial}}(\mathrm{IndSets}).\\
\end{align}
We conjecture the corresponding homotopy colimits exist throughout in order to drop this technical assumption. The resulting $\infty$-categories are denoted to be:
\begin{align}
\mathrm{Object}_{\mathrm{E}_\infty\mathrm{commutativealgebra},\mathrm{Simplicial}}(\mathrm{IndSNorm}_{B})^{\mathrm{smoothformalseriesclosure},\text{prelog}},\\
\mathrm{Object}_{\mathrm{E}_\infty\mathrm{commutativealgebra},\mathrm{Simplicial}}(\mathrm{Ind}^m\mathrm{SNorm}_{B})^{\mathrm{smoothformalseriesclosure},\text{prelog}},\\
\mathrm{Object}_{\mathrm{E}_\infty\mathrm{commutativealgebra},\mathrm{Simplicial}}(\mathrm{IndNorm}_{B})^{\mathrm{smoothformalseriesclosure},\text{prelog}},\\
\mathrm{Object}_{\mathrm{E}_\infty\mathrm{commutativealgebra},\mathrm{Simplicial}}(\mathrm{Ind}^m\mathrm{Norm}_{B})^{\mathrm{smoothformalseriesclosure},\text{prelog}},\\
\mathrm{Object}_{\mathrm{E}_\infty\mathrm{commutativealgebra},\mathrm{Simplicial}}(\mathrm{IndBan}_{B})^{\mathrm{smoothformalseriesclosure},\text{prelog}},\\
\mathrm{Object}_{\mathrm{E}_\infty\mathrm{commutativealgebra},\mathrm{Simplicial}}(\mathrm{Ind}^m\mathrm{Ban}_{B})^{\mathrm{smoothformalseriesclosure},\text{prelog}}.
\end{align}	
Namely one considers all the logarithmic formal series rings, and then takes in 
\begin{align}
\mathrm{Object}_{\mathrm{Simplicial}}(\mathrm{IndSets})\times \mathrm{Object}_{\mathrm{Simplicial}}(\mathrm{IndSets})
\end{align}
the closure by colimits.\\
\end{assumption}

\

\begin{definition}
Then we can in the same fashion consider the corresponding derived prismatic complexes \cite[Construction 7.6]{12BS}, \cite[Definition 4.1]{12Ko1}\footnote{One just applies \cite[Construction 7.6]{12BS}, \cite[Definition 4.1]{12Ko1} and then takes the left Kan extensions. We note that the derived logarithmic prismatic cohomologies are considered in \cite[Just above Notation in Section 1]{12Ko1}.} for the commutative algebras as in the above (for a given log prism $(A,I,M)$):
\begin{align}
\mathrm{Kan}_{\mathrm{Left}}\Delta_{?/A},	
\end{align}
by the regular corresponding left Kan extension techniques on the following $(\infty,1)$-compactly generated closures of the corresponding polynomials given over $A/I$ with a chosen log prism $(A,I,M)$:

\begin{align}
\mathrm{Object}_{\mathrm{E}_\infty\mathrm{commutativealgebra},\mathrm{Simplicial}}(\mathrm{IndSNorm}_{A/I})^{\mathrm{smoothformalseriesclosure},\text{prelog}},\\
\mathrm{Object}_{\mathrm{E}_\infty\mathrm{commutativealgebra},\mathrm{Simplicial}}(\mathrm{Ind}^m\mathrm{SNorm}_{A/I})^{\mathrm{smoothformalseriesclosure},\text{prelog}},\\
\mathrm{Object}_{\mathrm{E}_\infty\mathrm{commutativealgebra},\mathrm{Simplicial}}(\mathrm{IndNorm}_{A/I})^{\mathrm{smoothformalseriesclosure},\text{prelog}},\\
\mathrm{Object}_{\mathrm{E}_\infty\mathrm{commutativealgebra},\mathrm{Simplicial}}(\mathrm{Ind}^m\mathrm{Norm}_{A/I})^{\mathrm{smoothformalseriesclosure},\text{prelog}},\\
\mathrm{Object}_{\mathrm{E}_\infty\mathrm{commutativealgebra},\mathrm{Simplicial}}(\mathrm{IndBan}_{A/I})^{\mathrm{smoothformalseriesclosure},\text{prelog}},\\
\mathrm{Object}_{\mathrm{E}_\infty\mathrm{commutativealgebra},\mathrm{Simplicial}}(\mathrm{Ind}^m\mathrm{Ban}_{A/I})^{\mathrm{smoothformalseriesclosure},\text{prelog}}.
\end{align}
We call the corresponding functors functional analytic derived logarithmic prismatic complexes which we are going to denote that as in the following:
\begin{align}
\mathrm{Kan}_{\mathrm{Left}}\Delta_{?/A,\text{functionalanalytic,logarithmic,KKM},\text{BBM,formalanalytification}}.	
\end{align}
This would mean the following definition{\footnote{Before the Ben-Bassat-Mukherjee $p$-adic formal analytification we take the corresponding derived $(p,I)$-completion.}}:
\begin{align}
\mathrm{Kan}_{\mathrm{Left}}&\Delta_{?/A,\text{functionalanalytic,logarithmic,KKM},\text{BBM,formalanalytification}}(\mathcal{O})\\
&:=	((\underset{i}{\text{homotopycolimit}}_{\text{sifted},\text{derivedcategory}_{\infty}(A/I-\text{Module})}\mathrm{Kan}_{\mathrm{Left}}\Delta_{?/A,\text{functionalanalytic,logarithmic,KKM}}\\
&(\mathcal{O}_i))^\wedge\\
&)_\text{BBM,formalanalytification}
\end{align}
by writing any object $\mathcal{O}$ as the corresponding colimit 
\begin{center}
$\underset{i}{\text{homotopycolimit}}_\text{sifted}\mathcal{O}_i$.
\end{center}
These are quite large $(\infty,1)$-commutative ring objects in the corresponding $(\infty,1)$-categories for $(R=A/I,*)$:

\begin{align}
\mathrm{Object}_{\mathrm{E}_\infty\mathrm{commutativealgebra},\mathrm{Simplicial}}(\mathrm{IndSNorm}_R)^{\text{prelog}},\\
\mathrm{Object}_{\mathrm{E}_\infty\mathrm{commutativealgebra},\mathrm{Simplicial}}(\mathrm{Ind}^m\mathrm{SNorm}_R)^{\text{prelog}},\\
\mathrm{Object}_{\mathrm{E}_\infty\mathrm{commutativealgebra},\mathrm{Simplicial}}(\mathrm{IndNorm}_R)^{\text{prelog}},\\
\mathrm{Object}_{\mathrm{E}_\infty\mathrm{commutativealgebra},\mathrm{Simplicial}}(\mathrm{Ind}^m\mathrm{Norm}_R)^{\text{prelog}},\\
\mathrm{Object}_{\mathrm{E}_\infty\mathrm{commutativealgebra},\mathrm{Simplicial}}(\mathrm{IndBan}_R)^{\text{prelog}},\\
\mathrm{Object}_{\mathrm{E}_\infty\mathrm{commutativealgebra},\mathrm{Simplicial}}(\mathrm{Ind}^m\mathrm{Ban}_R)^{\text{prelog}},\\
\end{align}
after taking the formal series ring left Kan extension analytification from \cite[Section 4.2]{BBM}, which is defined by taking the left Kan extension to all the $(\infty,1)$-ring objects in the $\infty$-derived category of all $A$-modules from formal series rings over $A$.
\end{definition}

\newpage

\subsection{Functional Analytic Derived Logarithmic Perfectoidizations}

\indent Then as in \cite[Definition 8.2]{12BS} we consider the corresponding perfectoidization in this analytic setting.

\begin{definition}
Let $(A,I,M)$ be a perfectoid logarithmic prism as in \cite[Definition 3.3]{12Ko1}, and we consider any $\mathrm{E}_\infty$-ring $\mathcal{O}$ in the following
\begin{align}
\mathrm{Object}_{\mathrm{E}_\infty\mathrm{commutativealgebra},\mathrm{Simplicial}}(\mathrm{IndSNorm}_{A/I})^{\mathrm{smoothformalseriesclosure},\text{prelog}},\\
\mathrm{Object}_{\mathrm{E}_\infty\mathrm{commutativealgebra},\mathrm{Simplicial}}(\mathrm{Ind}^m\mathrm{SNorm}_{A/I})^{\mathrm{smoothformalseriesclosure},\text{prelog}},\\
\mathrm{Object}_{\mathrm{E}_\infty\mathrm{commutativealgebra},\mathrm{Simplicial}}(\mathrm{IndNorm}_{A/I})^{\mathrm{smoothformalseriesclosure},\text{prelog}},\\
\mathrm{Object}_{\mathrm{E}_\infty\mathrm{commutativealgebra},\mathrm{Simplicial}}(\mathrm{Ind}^m\mathrm{Norm}_{A/I})^{\mathrm{smoothformalseriesclosure},\text{prelog}},\\
\mathrm{Object}_{\mathrm{E}_\infty\mathrm{commutativealgebra},\mathrm{Simplicial}}(\mathrm{IndBan}_{A/I})^{\mathrm{smoothformalseriesclosure},\text{prelog}},\\
\mathrm{Object}_{\mathrm{E}_\infty\mathrm{commutativealgebra},\mathrm{Simplicial}}(\mathrm{Ind}^m\mathrm{Ban}_{A/I})^{\mathrm{smoothformalseriesclosure},\text{prelog}}.
\end{align}
Then consider the derived prismatic object:
\begin{align}
\mathrm{Kan}_{\mathrm{Left}}&\Delta_{?/A,\text{functionalanalytic,logarithmic,KKM},\text{BBM,formalanalytification}}(\mathcal{O})\\
&:=	((\underset{i}{\text{homotopycolimit}}_{\text{sifted},\text{derivedcategory}_{\infty}(A/I-\text{Module})}\mathrm{Kan}_{\mathrm{Left}}\Delta_{?/A,\text{functionalanalytic,logarithmic,KKM}}\\
&(\mathcal{O}_i))^\wedge\\
&)_\text{BBM,formalanalytification}.
\end{align}
Then as in \cite[Definition 8.2]{12BS} we have the following preperfectoidization:
\begin{align}
&(\mathcal{O})^{\text{preperfectoidization}}\\
&:=\mathrm{Colimit}(\mathrm{Kan}_{\mathrm{Left}}\Delta_{?/A,\text{functionalanalytic,logarithmic,KKM},\text{BBM,formalanalytification}}(\mathcal{O})\longrightarrow \\
&\mathrm{Fro}_*\mathrm{Kan}_{\mathrm{Left}}\Delta_{?/A,\text{functionalanalytic,logarithmic,KKM},\text{BBM,formalanalytification}}(\mathcal{O})\\
&\longrightarrow \mathrm{Fro}_* \mathrm{Fro}_*\mathrm{Kan}_{\mathrm{Left}}\Delta_{?/A,\text{functionalanalytic,logarithmic,KKM},\text{BBM,formalanalytification}}(\mathcal{O})\longrightarrow...\\
&)^{\text{BBM,formalanalytification}}.	
\end{align}
\footnote{Again after taking the formal series ring left Kan extension analytification from \cite[Section 4.2]{BBM}, which is defined by taking the left Kan extension to all the $(\infty,1)$-ring objects in the $\infty$-derived category of all $A$-modules from formal series rings over $A$, into:
\begin{align}
\mathrm{Object}_{\mathrm{E}_\infty\mathrm{commutativealgebra},\mathrm{Simplicial}}(\mathrm{IndSNorm}_R)^{\text{prelog}},\\
\mathrm{Object}_{\mathrm{E}_\infty\mathrm{commutativealgebra},\mathrm{Simplicial}}(\mathrm{Ind}^m\mathrm{SNorm}_R)^{\text{prelog}},\\
\mathrm{Object}_{\mathrm{E}_\infty\mathrm{commutativealgebra},\mathrm{Simplicial}}(\mathrm{IndNorm}_R)^{\text{prelog}},\\
\mathrm{Object}_{\mathrm{E}_\infty\mathrm{commutativealgebra},\mathrm{Simplicial}}(\mathrm{Ind}^m\mathrm{Norm}_R)^{\text{prelog}},\\
\mathrm{Object}_{\mathrm{E}_\infty\mathrm{commutativealgebra},\mathrm{Simplicial}}(\mathrm{IndBan}_R)^{\text{prelog}},\\
\mathrm{Object}_{\mathrm{E}_\infty\mathrm{commutativealgebra},\mathrm{Simplicial}}(\mathrm{Ind}^m\mathrm{Ban}_R)^{\text{prelog}}.\\
\end{align}
}Then we define the corresponding perfectoidization:
\begin{align}
&(\mathcal{O})^{\text{perfectoidization}}\\
&:=\mathrm{Colimit}(\mathrm{Kan}_{\mathrm{Left}}\Delta_{?/A,\text{functionalanalytic,logarithmic,KKM},\text{BBM,formalanalytification}}(\mathcal{O})\longrightarrow \\
&\mathrm{Fro}_*\mathrm{Kan}_{\mathrm{Left}}\Delta_{?/A,\text{functionalanalytic,logarithmic,KKM},\text{BBM,formalanalytification}}(\mathcal{O})\\
&\longrightarrow \mathrm{Fro}_* \mathrm{Fro}_*\mathrm{Kan}_{\mathrm{Left}}\Delta_{?/A,\text{functionalanalytic,logarithmic,KKM},\text{BBM,formalanalytification}}(\mathcal{O})\longrightarrow\\
&...)^{\text{BBM,formalanalytification}}\times A/I.	
\end{align}
Furthermore one can take derived $(p,I)$-completion to achieve the derived $(p,I)$-completed versions:
\begin{align}
\mathcal{O}^\text{preperfectoidization,derivedcomplete}:=(\mathcal{O}^\text{preperfectoidization})^{\wedge},\\
\mathcal{O}^\text{perfectoidization,derivedcomplete}:=\mathcal{O}^\text{preperfectoidization,derivedcomplete}\times A/I.\\
\end{align}

\noindent These are large $(\infty,1)$-commutative algebra objects in the corresponding categories as in the above, attached to also large $(\infty,1)$-commutative algebra objects. When we apply this to the corresponding sub-$(\infty,1)$-categories of Banach perfectoid objects in \cite{BMS2}, \cite{12DLLZ1}, \cite{12DLLZ2}, \cite{GR}, \cite{12KL1}, \cite{12KL2}, \cite{12Ked1}, \cite{12Sch3} we will recover the corresponding distinguished elemental deformation processes defined in \cite{BMS2}, \cite{12DLLZ1}, \cite{12DLLZ2}, \cite{GR}, \cite{12KL1}, \cite{12KL2}, \cite{12Ked1}, \cite{12Sch3}.
\end{definition}

\

\begin{remark}
One can then define such ring $\mathcal{O}$ to be \textit{logarithmicpreperfectoid} if we have the equivalence:
\begin{align}
\mathcal{O}^{\text{preperfectoidization}} \overset{\sim}{\longrightarrow}	\mathcal{O}.
\end{align}
One can then define such ring $\mathcal{O}$ to be \textit{logarithmicperfectoid} if we have the equivalence:
\begin{align}
\mathcal{O}^{\text{preperfectoidization}}\times A/I \overset{\sim}{\longrightarrow}	\mathcal{O}.
\end{align}
	
\end{remark}

\newpage

\section{Functional Analytic Derived Prismatic Complexes for $(\infty,1)$-Analytic Stacks and the Preperfectoidizations}

\subsection{Functional Analytic Derived Prismatic Complexes for $(\infty,1)$-Analytic Stacks}

\indent We now promote the construction in the previous sections to the corresponding $(\infty,1)$-ringed toposes level after Lurie \cite{12Lu1}, \cite{12Lu2} and \cite{Lu3} in the $\infty$-category of $\infty$-ringed toposes, Bambozzi-Ben-Bassat-Kremnizer \cite{12BBBK}, Ben-Bassat-Mukherjee \cite{BBM}, Bambozzi-Kremnizer \cite{12BK}, Clausen-Scholze \cite{12CS1} \cite{12CS2} and Kelly-Kremnizer-Mukherjee \cite{KKM} in the $\infty$-cateogory of $\infty$-functional analytic ringed toposes.\\

\indent Now we consider the following $\infty$-categories of the corresponding $\infty$-analytic ringed toposes from Bambozzi-Ben-Bassat-Kremnizer \cite{12BBBK}:\\

\begin{align}
&\infty-\mathrm{Toposes}^{\mathrm{ringed},\mathrm{commutativealgebra}_{\mathrm{simplicial}}(\mathrm{Ind}\mathrm{Seminormed}_?)}_{\mathrm{Commutativealgebra}_{\mathrm{simplicial}}(\mathrm{Ind}\mathrm{Seminormed}_?)^\mathrm{opposite},\mathrm{Grothendiecktopology,homotopyepimorphism}}.\\
&\infty-\mathrm{Toposes}^{\mathrm{ringed},\mathrm{Commutativealgebra}_{\mathrm{simplicial}}(\mathrm{Ind}^m\mathrm{Seminormed}_?)}_{\mathrm{Commutativealgebra}_{\mathrm{simplicial}}(\mathrm{Ind}^m\mathrm{Seminormed}_?)^\mathrm{opposite},\mathrm{Grothendiecktopology,homotopyepimorphism}}.\\
&\infty-\mathrm{Toposes}^{\mathrm{ringed},\mathrm{Commutativealgebra}_{\mathrm{simplicial}}(\mathrm{Ind}\mathrm{Normed}_?)}_{\mathrm{Commutativealgebra}_{\mathrm{simplicial}}(\mathrm{Ind}\mathrm{Normed}_?)^\mathrm{opposite},\mathrm{Grothendiecktopology,homotopyepimorphism}}.\\
&\infty-\mathrm{Toposes}^{\mathrm{ringed},\mathrm{Commutativealgebra}_{\mathrm{simplicial}}(\mathrm{Ind}^m\mathrm{Normed}_?)}_{\mathrm{Commutativealgebra}_{\mathrm{simplicial}}(\mathrm{Ind}^m\mathrm{Normed}_?)^\mathrm{opposite},\mathrm{Grothendiecktopology,homotopyepimorphism}}.\\
&\infty-\mathrm{Toposes}^{\mathrm{ringed},\mathrm{Commutativealgebra}_{\mathrm{simplicial}}(\mathrm{Ind}\mathrm{Banach}_?)}_{\mathrm{Commutativealgebra}_{\mathrm{simplicial}}(\mathrm{Ind}\mathrm{Banach}_?)^\mathrm{opposite},\mathrm{Grothendiecktopology,homotopyepimorphism}}.\\
&\infty-\mathrm{Toposes}^{\mathrm{ringed},\mathrm{Commutativealgebra}_{\mathrm{simplicial}}(\mathrm{Ind}^m\mathrm{Banach}_?)}_{\mathrm{Commutativealgebra}_{\mathrm{simplicial}}(\mathrm{Ind}^m\mathrm{Banach}_?)^\mathrm{opposite},\mathrm{Grothendiecktopology,homotopyepimorphism}}.\\ 
&\mathrm{Proj}^\text{smoothformalseriesclosure}\infty-\mathrm{Toposes}^{\mathrm{ringed},\mathrm{commutativealgebra}_{\mathrm{simplicial}}(\mathrm{Ind}\mathrm{Seminormed}_?)}_{\mathrm{Commutativealgebra}_{\mathrm{simplicial}}(\mathrm{Ind}\mathrm{Seminormed}_?)^\mathrm{opposite},\mathrm{Grothendiecktopology,homotopyepimorphism}}. \\
&\mathrm{Proj}^\text{smoothformalseriesclosure}\infty-\mathrm{Toposes}^{\mathrm{ringed},\mathrm{Commutativealgebra}_{\mathrm{simplicial}}(\mathrm{Ind}^m\mathrm{Seminormed}_?)}_{\mathrm{Commutativealgebra}_{\mathrm{simplicial}}(\mathrm{Ind}^m\mathrm{Seminormed}_?)^\mathrm{opposite},\mathrm{Grothendiecktopology,homotopyepimorphism}}.\\
&\mathrm{Proj}^\text{smoothformalseriesclosure}\infty-\mathrm{Toposes}^{\mathrm{ringed},\mathrm{Commutativealgebra}_{\mathrm{simplicial}}(\mathrm{Ind}\mathrm{Normed}_?)}_{\mathrm{Commutativealgebra}_{\mathrm{simplicial}}(\mathrm{Ind}\mathrm{Normed}_?)^\mathrm{opposite},\mathrm{Grothendiecktopology,homotopyepimorphism}}.\\
&\mathrm{Proj}^\text{smoothformalseriesclosure}\infty-\mathrm{Toposes}^{\mathrm{ringed},\mathrm{Commutativealgebra}_{\mathrm{simplicial}}(\mathrm{Ind}^m\mathrm{Normed}_?)}_{\mathrm{Commutativealgebra}_{\mathrm{simplicial}}(\mathrm{Ind}^m\mathrm{Normed}_?)^\mathrm{opposite},\mathrm{Grothendiecktopology,homotopyepimorphism}}.\\
&\mathrm{Proj}^\text{smoothformalseriesclosure}\infty-\mathrm{Toposes}^{\mathrm{ringed},\mathrm{Commutativealgebra}_{\mathrm{simplicial}}(\mathrm{Ind}\mathrm{Banach}_?)}_{\mathrm{Commutativealgebra}_{\mathrm{simplicial}}(\mathrm{Ind}\mathrm{Banach}_?)^\mathrm{opposite},\mathrm{Grothendiecktopology,homotopyepimorphism}}.\\
&\mathrm{Proj}^\text{smoothformalseriesclosure}\infty-\mathrm{Toposes}^{\mathrm{ringed},\mathrm{Commutativealgebra}_{\mathrm{simplicial}}(\mathrm{Ind}^m\mathrm{Banach}_?)}_{\mathrm{Commutativealgebra}_{\mathrm{simplicial}}(\mathrm{Ind}^m\mathrm{Banach}_?)^\mathrm{opposite},\mathrm{Grothendiecktopology,homotopyepimorphism}}.\\ 
\end{align}

\begin{example}
\mbox{(Preadic nonsheafy spaces and their colimits)} The main example we would like consider comes from \cite{12BK} where one constructs the corresponding derived adic space $(\mathrm{Spectrumadic}_{\mathrm{BK}}(R),\mathcal{O}_{\mathrm{Spectrumadic}_{\mathrm{BK}}(R)})$ by using derived rational localization to reach some $\infty$-toposes carrying essentially $\infty$-sheaf of rings, from any Banach ring $R$\footnote{With some key essential assumptions but that is not so serious once one considers some foundation from Kedlaya in \cite{Ked2} on reified adic spaces as mentioned in the last section of \cite{12BK}.}. Then we apply this to the corresponding situation below. Let $(A,I)$ be a corresponding prism from Bhatt-Scholze, with the assumption that $A/I$ is Banach. Then what we do is consider all the formal series rings over $A/I$ then take the colimit completion in the homotopy sense, by embedding them through $\mathrm{Spectrumadic}_{\mathrm{BK}}$ into the $\infty$-toposes ringed. For instance for any such formal ring $F$, one could regard this as a inductive system of Banach rings:
\begin{align}
F=\underset{i}{\mathrm{homotopycolimit}} F_i,	
\end{align}
where we set:
\begin{align}
(\mathrm{Spectrumadic}_{\mathrm{BK}}(F),\mathcal{O}):=\underset{i}{\mathrm{homotopylimit}} (\mathrm{Spectrumadic}_{\mathrm{BK}}(F_i),\mathcal{O}_{\mathrm{Spectrumadic}_{\mathrm{BK}}(F_i)}).	
\end{align}
The resulting $\infty$-stacks generated are interesting to study. 	
\end{example}

\begin{definition}
\indent One can actually define the derived prismatic cohomology presheaves through derived topological Hochschild cohomology presheaves, derived topological period cohomology presheaves and derived topological cyclic cohomology presheaves as in \cite[Section 2.2, Section 2.3]{12BMS}, \cite[Theorem 1.13]{12BS}:
\begin{align}
	\mathrm{Kan}_{\mathrm{Left}}\mathrm{THH},\mathrm{Kan}_{\mathrm{Left}}\mathrm{TP},\mathrm{Kan}_{\mathrm{Left}}\mathrm{TC},
\end{align}
on the following $(\infty,1)$-compactly generated closures of the corresponding polynomials\footnote{Definitely, we need to put certain norms over in some relatively canonical way, as in \cite[Section 4.2]{BBM} one can basically consider rigid ones and dagger ones, and so on. We restrict to the \textit{formal} one.} given over $A/I$ with a chosen prism $(A,I)$\footnote{In all the following, we assume this prism to be bounded and satisfy that $A/I$ is Banach.}:
\begin{align}
\mathrm{Proj}^\text{smoothformalseriesclosure}\infty-\mathrm{Toposes}^{\mathrm{ringed},\mathrm{commutativealgebra}_{\mathrm{simplicial}}(\mathrm{Ind}\mathrm{Seminormed}_{A/I})}_{\mathrm{Commutativealgebra}_{\mathrm{simplicial}}(\mathrm{Ind}\mathrm{Seminormed}_{A/I})^\mathrm{opposite},\mathrm{Grothendiecktopology,homotopyepimorphism}}. \\
\mathrm{Proj}^\text{smoothformalseriesclosure}\infty-\mathrm{Toposes}^{\mathrm{ringed},\mathrm{Commutativealgebra}_{\mathrm{simplicial}}(\mathrm{Ind}^m\mathrm{Seminormed}_{A/I})}_{\mathrm{Commutativealgebra}_{\mathrm{simplicial}}(\mathrm{Ind}^m\mathrm{Seminormed}_{A/I})^\mathrm{opposite},\mathrm{Grothendiecktopology,homotopyepimorphism}}.\\
\mathrm{Proj}^\text{smoothformalseriesclosure}\infty-\mathrm{Toposes}^{\mathrm{ringed},\mathrm{Commutativealgebra}_{\mathrm{simplicial}}(\mathrm{Ind}\mathrm{Normed}_{A/I})}_{\mathrm{Commutativealgebra}_{\mathrm{simplicial}}(\mathrm{Ind}\mathrm{Normed}_{A/I})^\mathrm{opposite},\mathrm{Grothendiecktopology,homotopyepimorphism}}.\\
\mathrm{Proj}^\text{smoothformalseriesclosure}\infty-\mathrm{Toposes}^{\mathrm{ringed},\mathrm{Commutativealgebra}_{\mathrm{simplicial}}(\mathrm{Ind}^m\mathrm{Normed}_{A/I})}_{\mathrm{Commutativealgebra}_{\mathrm{simplicial}}(\mathrm{Ind}^m\mathrm{Normed}_{A/I})^\mathrm{opposite},\mathrm{Grothendiecktopology,homotopyepimorphism}}.\\
\mathrm{Proj}^\text{smoothformalseriesclosure}\infty-\mathrm{Toposes}^{\mathrm{ringed},\mathrm{Commutativealgebra}_{\mathrm{simplicial}}(\mathrm{Ind}\mathrm{Banach}_{A/I})}_{\mathrm{Commutativealgebra}_{\mathrm{simplicial}}(\mathrm{Ind}\mathrm{Banach}_{A/I})^\mathrm{opposite},\mathrm{Grothendiecktopology,homotopyepimorphism}}.\\
\mathrm{Proj}^\text{smoothformalseriesclosure}\infty-\mathrm{Toposes}^{\mathrm{ringed},\mathrm{Commutativealgebra}_{\mathrm{simplicial}}(\mathrm{Ind}^m\mathrm{Banach}_{A/I})}_{\mathrm{Commutativealgebra}_{\mathrm{simplicial}}(\mathrm{Ind}^m\mathrm{Banach}_{A/I})^\mathrm{opposite},\mathrm{Grothendiecktopology,homotopyepimorphism}}. 
\end{align}
We call the corresponding functors are derived functional analytic Hochschild cohomology presheaves, derived functional analytic period cohomology presheaves and derived functional analytic cyclic cohomology presheaves, which we are going to denote these presheaves as in the following for any $\infty$-ringed topos $(\mathbb{X},\mathcal{O})=\underset{i}{\text{homotopycolimit}}(\mathbb{X}_i,\mathcal{O}_i)$:
\begin{align}
	&\mathrm{Kan}_{\mathrm{Left}}\mathrm{THH}_{\text{functionalanalytic,KKM},\text{BBM,formalanalytification}}(\mathcal{O}):=\\
	&(\underset{i}{\text{homotopycolimit}}_{\text{sifted},\text{derivedcategory}_{\infty}(A/I-\text{Module})}\mathrm{Kan}_{\mathrm{Left}}\mathrm{THH}_{\text{functionalanalytic,KKM}}(\mathcal{O}_i)\\
	&)_\text{BBM,formalanalytification},\\
	&\mathrm{Kan}_{\mathrm{Left}}\mathrm{TP}_{\text{functionalanalytic,KKM},\text{BBM,formalanalytification}}(\mathcal{O}):=\\
	&(\underset{i}{\text{homotopycolimit}}_{\text{sifted},\text{derivedcategory}_{\infty}(A/I-\text{Module})}\mathrm{Kan}_{\mathrm{Left}}\mathrm{TP}_{\text{functionalanalytic,KKM}}(\mathcal{O}_i)\\
	&)_\text{BBM,formalanalytification},\\
	&\mathrm{Kan}_{\mathrm{Left}}\mathrm{TC}_{\text{functionalanalytic,KKM},\text{BBM,formalanalytification}}(\mathcal{O}):=\\
	&(\underset{i}{\text{homotopycolimit}}_{\text{sifted},\text{derivedcategory}_{\infty}(A/I-\text{Module})}\mathrm{Kan}_{\mathrm{Left}}\mathrm{TC}_{\text{functionalanalytic,KKM}}(\mathcal{O}_i)\\
	&)_\text{BBM,formalanalytification},
\end{align}
by writing any object $\mathcal{O}$ as the corresponding colimit 
\begin{center}
$\underset{i}{\text{homotopycolimit}}_\text{sifted}\mathcal{O}_i:=\underset{i}{\text{homotopycolimit}}_\text{sifted}\mathcal{O}_{\mathrm{Spectrumadic}_\mathrm{BK}(\mathrm{Formalseriesring}_i)}$.
\end{center}
These are quite large $(\infty,1)$-commutative ring objects in the corresponding $(\infty,1)$-categories for $?=A/I$ over $(\mathbb{X},\mathcal{O})$:

 \begin{align}
\mathrm{Ind}\mathrm{\sharp Quasicoherent}^{\text{presheaf}}_{\mathrm{Ind}^\text{smoothformalseriesclosure}\infty-\mathrm{Toposes}^{\mathrm{ringed},\mathrm{commutativealgebra}_{\mathrm{simplicial}}(\mathrm{Ind}\mathrm{Seminormed}_?)}_{\mathrm{Commutativealgebra}_{\mathrm{simplicial}}(\mathrm{Ind}\mathrm{Seminormed}_?)^\mathrm{opposite},\mathrm{Grotopology,homotopyepimorphism}}}. \\
\mathrm{Ind}\mathrm{\sharp Quasicoherent}^{\text{presheaf}}_{\mathrm{Ind}^\text{smoothformalseriesclosure}\infty-\mathrm{Toposes}^{\mathrm{ringed},\mathrm{Commutativealgebra}_{\mathrm{simplicial}}(\mathrm{Ind}^m\mathrm{Seminormed}_?)}_{\mathrm{Commutativealgebra}_{\mathrm{simplicial}}(\mathrm{Ind}^m\mathrm{Seminormed}_?)^\mathrm{opposite},\mathrm{Grotopology,homotopyepimorphism}}}.\\
\mathrm{Ind}\mathrm{\sharp Quasicoherent}^{\text{presheaf}}_{\mathrm{Ind}^\text{smoothformalseriesclosure}\infty-\mathrm{Toposes}^{\mathrm{ringed},\mathrm{Commutativealgebra}_{\mathrm{simplicial}}(\mathrm{Ind}\mathrm{Normed}_?)}_{\mathrm{Commutativealgebra}_{\mathrm{simplicial}}(\mathrm{Ind}\mathrm{Normed}_?)^\mathrm{opposite},\mathrm{Grotopology,homotopyepimorphism}}}.\\
\mathrm{Ind}\mathrm{\sharp Quasicoherent}^{\text{presheaf}}_{\mathrm{Ind}^\text{smoothformalseriesclosure}\infty-\mathrm{Toposes}^{\mathrm{ringed},\mathrm{Commutativealgebra}_{\mathrm{simplicial}}(\mathrm{Ind}^m\mathrm{Normed}_?)}_{\mathrm{Commutativealgebra}_{\mathrm{simplicial}}(\mathrm{Ind}^m\mathrm{Normed}_?)^\mathrm{opposite},\mathrm{Grotopology,homotopyepimorphism}}}.\\
\mathrm{Ind}\mathrm{\sharp Quasicoherent}^{\text{presheaf}}_{\mathrm{Ind}^\text{smoothformalseriesclosure}\infty-\mathrm{Toposes}^{\mathrm{ringed},\mathrm{Commutativealgebra}_{\mathrm{simplicial}}(\mathrm{Ind}\mathrm{Banach}_?)}_{\mathrm{Commutativealgebra}_{\mathrm{simplicial}}(\mathrm{Ind}\mathrm{Banach}_?)^\mathrm{opposite},\mathrm{Grotopology,homotopyepimorphism}}}.\\
\mathrm{Ind}\mathrm{\sharp Quasicoherent}^{\text{presheaf}}_{\mathrm{Ind}^\text{smoothformalseriesclosure}\infty-\mathrm{Toposes}^{\mathrm{ringed},\mathrm{Commutativealgebra}_{\mathrm{simplicial}}(\mathrm{Ind}^m\mathrm{Banach}_?)}_{\mathrm{Commutativealgebra}_{\mathrm{simplicial}}(\mathrm{Ind}^m\mathrm{Banach}_?)^\mathrm{opposite},\mathrm{Grotopology,homotopyepimorphism}}},\\ 
\end{align}
after taking the formal series ring left Kan extension analytification from \cite[Section 4.2]{BBM}, which is defined by taking the left Kan extension to all the $(\infty,1)$-ring objects in the $\infty$-derived category of all $A$-modules from formal series rings over $A$.
\end{definition}

\

\begin{definition}
Then we can in the same fashion consider the corresponding derived prismatic complex presheaves \cite[Construction 7.6]{12BS}\footnote{One just applies \cite[Construction 7.6]{12BS} and then takes the left Kan extensions.} for the commutative algebras as in the above (for a given prism $(A,I)$):
\begin{align}
\mathrm{Kan}_{\mathrm{Left}}\Delta_{?/A},	
\end{align}
by the regular corresponding left Kan extension techniques on the following $(\infty,1)$-compactly generated closures of the corresponding polynomials given over $A/I$ with a chosen prism $(A,I)$:\\

\begin{align}
\mathrm{Proj}^\text{smoothformalseriesclosure}\infty-\mathrm{Toposes}^{\mathrm{ringed},\mathrm{commutativealgebra}_{\mathrm{simplicial}}(\mathrm{Ind}\mathrm{Seminormed}_{A/I})}_{\mathrm{Commutativealgebra}_{\mathrm{simplicial}}(\mathrm{Ind}\mathrm{Seminormed}_{A/I})^\mathrm{opposite},\mathrm{Grothendiecktopology,homotopyepimorphism}}. \\
\mathrm{Proj}^\text{smoothformalseriesclosure}\infty-\mathrm{Toposes}^{\mathrm{ringed},\mathrm{Commutativealgebra}_{\mathrm{simplicial}}(\mathrm{Ind}^m\mathrm{Seminormed}_{A/I})}_{\mathrm{Commutativealgebra}_{\mathrm{simplicial}}(\mathrm{Ind}^m\mathrm{Seminormed}_{A/I})^\mathrm{opposite},\mathrm{Grothendiecktopology,homotopyepimorphism}}.\\
\mathrm{Proj}^\text{smoothformalseriesclosure}\infty-\mathrm{Toposes}^{\mathrm{ringed},\mathrm{Commutativealgebra}_{\mathrm{simplicial}}(\mathrm{Ind}\mathrm{Normed}_{A/I})}_{\mathrm{Commutativealgebra}_{\mathrm{simplicial}}(\mathrm{Ind}\mathrm{Normed}_{A/I})^\mathrm{opposite},\mathrm{Grothendiecktopology,homotopyepimorphism}}.\\
\mathrm{Proj}^\text{smoothformalseriesclosure}\infty-\mathrm{Toposes}^{\mathrm{ringed},\mathrm{Commutativealgebra}_{\mathrm{simplicial}}(\mathrm{Ind}^m\mathrm{Normed}_{A/I})}_{\mathrm{Commutativealgebra}_{\mathrm{simplicial}}(\mathrm{Ind}^m\mathrm{Normed}_{A/I})^\mathrm{opposite},\mathrm{Grothendiecktopology,homotopyepimorphism}}.\\
\mathrm{Proj}^\text{smoothformalseriesclosure}\infty-\mathrm{Toposes}^{\mathrm{ringed},\mathrm{Commutativealgebra}_{\mathrm{simplicial}}(\mathrm{Ind}\mathrm{Banach}_{A/I})}_{\mathrm{Commutativealgebra}_{\mathrm{simplicial}}(\mathrm{Ind}\mathrm{Banach}_{A/I})^\mathrm{opposite},\mathrm{Grothendiecktopology,homotopyepimorphism}}.\\
\mathrm{Proj}^\text{smoothformalseriesclosure}\infty-\mathrm{Toposes}^{\mathrm{ringed},\mathrm{Commutativealgebra}_{\mathrm{simplicial}}(\mathrm{Ind}^m\mathrm{Banach}_{A/I})}_{\mathrm{Commutativealgebra}_{\mathrm{simplicial}}(\mathrm{Ind}^m\mathrm{Banach}_{A/I})^\mathrm{opposite},\mathrm{Grothendiecktopology,homotopyepimorphism}}. 
\end{align}
We call the corresponding functors functional analytic derived prismatic complex presheaves which we are going to denote that as in the following:
\begin{align}
\mathrm{Kan}_{\mathrm{Left}}\Delta_{?/A,\text{functionalanalytic,KKM},\text{BBM,formalanalytification}}.	
\end{align}
This would mean the following definition{\footnote{Before the Ben-Bassat-Mukherjee $p$-adic formal analytification we take the corresponding derived $(p,I)$-completion.}}:
\begin{align}
&\mathrm{Kan}_{\mathrm{Left}}\Delta_{?/A,\text{functionalanalytic,KKM},\text{BBM,formalanalytification}}(\mathcal{O})\\
&:=	((\underset{i}{\text{homotopycolimit}}_{\text{sifted},\text{derivedcategory}_{\infty}(A/I-\text{Module})}\mathrm{Kan}_{\mathrm{Left}}\Delta_{?/A,\text{functionalanalytic,KKM}}(\mathcal{O}_i))^\wedge\\
&)_\text{BBM,formalanalytification}
\end{align}
by writing any object $\mathcal{O}$ as the corresponding colimit\footnote{Here we remind the readers of the corresponding foundation here, namely the presheaf $\mathcal{O}_i$ in fact takes the value in Koszul complex taking the following form:
\begin{align}
&\mathrm{Koszulcomplex}_{A/I\left<X_1,...,X_l\right>\left<T_1,...,T_m\right>}(a_1-T_1b_1,...,a_m-T_mb_m)\\
&= A/I\left<X_1,...,X_l\right>\left<T_1,...,T_m\right>/^\mathbb{L}(a_1-T_1b_1,...,a_m-T_mb_m).	
\end{align}
This is actually derived $p$-complete since each homotopy group is derived $p$-complete (from the corresponding Banach structure from $A/I\left<X_1,...,X_l\right>$ induced from the $p$-adic topology). Namely the definition of the presheaf:
\begin{align}
\mathrm{Kan}_{\mathrm{Left}}\Delta_{?/A,\text{functionalanalytic,KKM}}(\mathcal{O}_i)	
\end{align}
is directly the application of the derived prismatic functor from \cite[Construction 7.6]{12BS}.} 
\begin{center}
$\underset{i}{\text{homotopycolimit}}_\text{sifted}\mathcal{O}_i$.
\end{center}
These are quite large $(\infty,1)$-commutative ring objects in the corresponding $(\infty,1)$-categories for $R=A/I$:
\begin{align}
\mathrm{Ind}\mathrm{\sharp Quasicoherent}^{\text{presheaf}}_{\mathrm{Ind}^\text{smoothformalseriesclosure}\infty-\mathrm{Toposes}^{\mathrm{ringed},\mathrm{commutativealgebra}_{\mathrm{simplicial}}(\mathrm{Ind}\mathrm{Seminormed}_R)}_{\mathrm{Commutativealgebra}_{\mathrm{simplicial}}(\mathrm{Ind}\mathrm{Seminormed}_R)^\mathrm{opposite},\mathrm{Grotopology,homotopyepimorphism}}}. \\
\mathrm{Ind}\mathrm{\sharp Quasicoherent}^{\text{presheaf}}_{\mathrm{Ind}^\text{smoothformalseriesclosure}\infty-\mathrm{Toposes}^{\mathrm{ringed},\mathrm{Commutativealgebra}_{\mathrm{simplicial}}(\mathrm{Ind}^m\mathrm{Seminormed}_R)}_{\mathrm{Commutativealgebra}_{\mathrm{simplicial}}(\mathrm{Ind}^m\mathrm{Seminormed}_R)^\mathrm{opposite},\mathrm{Grotopology,homotopyepimorphism}}}.\\
\mathrm{Ind}\mathrm{\sharp Quasicoherent}^{\text{presheaf}}_{\mathrm{Ind}^\text{smoothformalseriesclosure}\infty-\mathrm{Toposes}^{\mathrm{ringed},\mathrm{Commutativealgebra}_{\mathrm{simplicial}}(\mathrm{Ind}\mathrm{Normed}_R)}_{\mathrm{Commutativealgebra}_{\mathrm{simplicial}}(\mathrm{Ind}\mathrm{Normed}_R)^\mathrm{opposite},\mathrm{Grotopology,homotopyepimorphism}}}.\\
\mathrm{Ind}\mathrm{\sharp Quasicoherent}^{\text{presheaf}}_{\mathrm{Ind}^\text{smoothformalseriesclosure}\infty-\mathrm{Toposes}^{\mathrm{ringed},\mathrm{Commutativealgebra}_{\mathrm{simplicial}}(\mathrm{Ind}^m\mathrm{Normed}_R)}_{\mathrm{Commutativealgebra}_{\mathrm{simplicial}}(\mathrm{Ind}^m\mathrm{Normed}_R)^\mathrm{opposite},\mathrm{Grotopology,homotopyepimorphism}}}.\\
\mathrm{Ind}\mathrm{\sharp Quasicoherent}^{\text{presheaf}}_{\mathrm{Ind}^\text{smoothformalseriesclosure}\infty-\mathrm{Toposes}^{\mathrm{ringed},\mathrm{Commutativealgebra}_{\mathrm{simplicial}}(\mathrm{Ind}\mathrm{Banach}_R)}_{\mathrm{Commutativealgebra}_{\mathrm{simplicial}}(\mathrm{Ind}\mathrm{Banach}_R)^\mathrm{opposite},\mathrm{Grotopology,homotopyepimorphism}}}.\\
\mathrm{Ind}\mathrm{\sharp Quasicoherent}^{\text{presheaf}}_{\mathrm{Ind}^\text{smoothformalseriesclosure}\infty-\mathrm{Toposes}^{\mathrm{ringed},\mathrm{Commutativealgebra}_{\mathrm{simplicial}}(\mathrm{Ind}^m\mathrm{Banach}_R)}_{\mathrm{Commutativealgebra}_{\mathrm{simplicial}}(\mathrm{Ind}^m\mathrm{Banach}_R)^\mathrm{opposite},\mathrm{Grotopology,homotopyepimorphism}}},\\ 
\end{align}
after taking the formal series ring left Kan extension analytification from \cite[Section 4.2]{BBM}, which is defined by taking the left Kan extension to all the $(\infty,1)$-ring objects in the $\infty$-derived category of all $A$-modules from formal series rings over $A$.\\
\end{definition}

\newpage

\subsection{Functional Analytic Derived Preperfectoidizations for $(\infty,1)$-Analytic Stacks}

\

\noindent Then as in \cite[Definition 8.2]{12BS} we consider the corresponding perfectoidization in this analytic setting.

\begin{definition}
Let $(A,I)$ be a perfectoid prism, and we consider any $\mathrm{E}_\infty$-ring $\mathcal{O}$ in the following
\begin{align}
\mathrm{Proj}^\text{smoothformalseriesclosure}\infty-\mathrm{Toposes}^{\mathrm{ringed},\mathrm{commutativealgebra}_{\mathrm{simplicial}}(\mathrm{Ind}\mathrm{Seminormed}_{A/I})}_{\mathrm{Commutativealgebra}_{\mathrm{simplicial}}(\mathrm{Ind}\mathrm{Seminormed}_{A/I})^\mathrm{opposite},\mathrm{Grothendiecktopology,homotopyepimorphism}}. \\
\mathrm{Proj}^\text{smoothformalseriesclosure}\infty-\mathrm{Toposes}^{\mathrm{ringed},\mathrm{Commutativealgebra}_{\mathrm{simplicial}}(\mathrm{Ind}^m\mathrm{Seminormed}_{A/I})}_{\mathrm{Commutativealgebra}_{\mathrm{simplicial}}(\mathrm{Ind}^m\mathrm{Seminormed}_{A/I})^\mathrm{opposite},\mathrm{Grothendiecktopology,homotopyepimorphism}}.\\
\mathrm{Proj}^\text{smoothformalseriesclosure}\infty-\mathrm{Toposes}^{\mathrm{ringed},\mathrm{Commutativealgebra}_{\mathrm{simplicial}}(\mathrm{Ind}\mathrm{Normed}_{A/I})}_{\mathrm{Commutativealgebra}_{\mathrm{simplicial}}(\mathrm{Ind}\mathrm{Normed}_{A/I})^\mathrm{opposite},\mathrm{Grothendiecktopology,homotopyepimorphism}}.\\
\mathrm{Proj}^\text{smoothformalseriesclosure}\infty-\mathrm{Toposes}^{\mathrm{ringed},\mathrm{Commutativealgebra}_{\mathrm{simplicial}}(\mathrm{Ind}^m\mathrm{Normed}_{A/I})}_{\mathrm{Commutativealgebra}_{\mathrm{simplicial}}(\mathrm{Ind}^m\mathrm{Normed}_{A/I})^\mathrm{opposite},\mathrm{Grothendiecktopology,homotopyepimorphism}}.\\
\mathrm{Proj}^\text{smoothformalseriesclosure}\infty-\mathrm{Toposes}^{\mathrm{ringed},\mathrm{Commutativealgebra}_{\mathrm{simplicial}}(\mathrm{Ind}\mathrm{Banach}_{A/I})}_{\mathrm{Commutativealgebra}_{\mathrm{simplicial}}(\mathrm{Ind}\mathrm{Banach}_{A/I})^\mathrm{opposite},\mathrm{Grothendiecktopology,homotopyepimorphism}}.\\
\mathrm{Proj}^\text{smoothformalseriesclosure}\infty-\mathrm{Toposes}^{\mathrm{ringed},\mathrm{Commutativealgebra}_{\mathrm{simplicial}}(\mathrm{Ind}^m\mathrm{Banach}_{A/I})}_{\mathrm{Commutativealgebra}_{\mathrm{simplicial}}(\mathrm{Ind}^m\mathrm{Banach}_{A/I})^\mathrm{opposite},\mathrm{Grothendiecktopology,homotopyepimorphism}}. 
\end{align}
Then consider the derived prismatic object:
\begin{align}
\mathrm{Kan}_{\mathrm{Left}}\Delta_{?/A,\text{functionalanalytic,KKM},\text{BBM,formalanalytification}}(\mathcal{O}).
\end{align}	
Then as in \cite[Definition 8.2]{12BS} we have the following preperfectoidization:
\begin{align}
&(\mathcal{O})^{\text{preperfectoidization}}\\
&:=\mathrm{Colimit}(\mathrm{Kan}_{\mathrm{Left}}\Delta_{?/A,\text{functionalanalytic,KKM},\text{BBM,formalanalytification}}(\mathcal{O})\rightarrow \\
&\mathrm{Fro}_*\mathrm{Kan}_{\mathrm{Left}}\Delta_{?/A,\text{functionalanalytic,KKM},\text{BBM,formalanalytification}}(\mathcal{O})\\
&\rightarrow \mathrm{Fro}_* \mathrm{Fro}_*\mathrm{Kan}_{\mathrm{Left}}\Delta_{?/A,\text{functionalanalytic,KKM},\text{BBM,formalanalytification}}(\mathcal{O})\rightarrow...)^{\text{BBM,formalanalytification}},	
\end{align}
after taking the formal series ring left Kan extension analytification from \cite[Section 4.2]{BBM}, which is defined by taking the left Kan extension to all the $(\infty,1)$-ring object in the $\infty$-derived category of all $A$-modules from formal series rings over $A$. Then we define the corresponding perfectoidization:
\begin{align}
&(\mathcal{O})^{\text{perfectoidization}}\\
&:=\mathrm{Colimit}(\mathrm{Kan}_{\mathrm{Left}}\Delta_{?/A,\text{functionalanalytic,KKM},\text{BBM,formalanalytification}}(\mathcal{O})\longrightarrow \\
&\mathrm{Fro}_*\mathrm{Kan}_{\mathrm{Left}}\Delta_{?/A,\text{functionalanalytic,KKM},\text{BBM,formalanalytification}}(\mathcal{O})\\
&\longrightarrow \mathrm{Fro}_* \mathrm{Fro}_*\mathrm{Kan}_{\mathrm{Left}}\Delta_{?/A,\text{functionalanalytic,KKM},\text{BBM,formalanalytification}}(\mathcal{O})\longrightarrow...)^{\text{BBM,formalanalytification}}\times A/I.	
\end{align}
Furthermore one can take derived $(p,I)$-completion to achieve the derived $(p,I)$-completed versions:
\begin{align}
\mathcal{O}^\text{preperfectoidization,derivedcomplete}:=(\mathcal{O}^\text{preperfectoidization})^{\wedge},\\
\mathcal{O}^\text{perfectoidization,derivedcomplete}:=\mathcal{O}^\text{preperfectoidization,derivedcomplete}\times A/I.\\
\end{align}
These are large $(\infty,1)$-commutative algebra objects in the corresponding categories as in the above, attached to also large $(\infty,1)$-commutative algebra objects. When we apply this to the corresponding sub-$(\infty,1)$-categories of Banach perfectoid objects in \cite{BMS2}, \cite{GR}, \cite{12KL1}, \cite{12KL2}, \cite{12Ked1}, \cite{12Sch3},  we will recover the corresponding distinguished elemental deformation processes defined in \cite{BMS2}, \cite{GR}, \cite{12KL1}, \cite{12KL2}, \cite{12Ked1}, \cite{12Sch3}. 
\end{definition}

\

\indent The $\infty$-presheaves in this section in the $\infty$-category:	
\begin{align}
\mathrm{Ind}\mathrm{\sharp Quasicoherent}^{\text{presheaf}}_{\mathrm{Ind}^\text{smoothformalseriesclosure}\infty-\mathrm{Toposes}^{\mathrm{ringed},\mathrm{commutativealgebra}_{\mathrm{simplicial}}(\mathrm{Ind}\mathrm{Seminormed}_R)}_{\mathrm{Commutativealgebra}_{\mathrm{simplicial}}(\mathrm{Ind}\mathrm{Seminormed}_R)^\mathrm{opposite},\mathrm{Grotopology,homotopyepimorphism}}}. \\
\mathrm{Ind}\mathrm{\sharp Quasicoherent}^{\text{presheaf}}_{\mathrm{Ind}^\text{smoothformalseriesclosure}\infty-\mathrm{Toposes}^{\mathrm{ringed},\mathrm{Commutativealgebra}_{\mathrm{simplicial}}(\mathrm{Ind}^m\mathrm{Seminormed}_R)}_{\mathrm{Commutativealgebra}_{\mathrm{simplicial}}(\mathrm{Ind}^m\mathrm{Seminormed}_R)^\mathrm{opposite},\mathrm{Grotopology,homotopyepimorphism}}}.\\
\mathrm{Ind}\mathrm{\sharp Quasicoherent}^{\text{presheaf}}_{\mathrm{Ind}^\text{smoothformalseriesclosure}\infty-\mathrm{Toposes}^{\mathrm{ringed},\mathrm{Commutativealgebra}_{\mathrm{simplicial}}(\mathrm{Ind}\mathrm{Normed}_R)}_{\mathrm{Commutativealgebra}_{\mathrm{simplicial}}(\mathrm{Ind}\mathrm{Normed}_R)^\mathrm{opposite},\mathrm{Grotopology,homotopyepimorphism}}}.\\
\mathrm{Ind}\mathrm{\sharp Quasicoherent}^{\text{presheaf}}_{\mathrm{Ind}^\text{smoothformalseriesclosure}\infty-\mathrm{Toposes}^{\mathrm{ringed},\mathrm{Commutativealgebra}_{\mathrm{simplicial}}(\mathrm{Ind}^m\mathrm{Normed}_R)}_{\mathrm{Commutativealgebra}_{\mathrm{simplicial}}(\mathrm{Ind}^m\mathrm{Normed}_R)^\mathrm{opposite},\mathrm{Grotopology,homotopyepimorphism}}}.\\
\mathrm{Ind}\mathrm{\sharp Quasicoherent}^{\text{presheaf}}_{\mathrm{Ind}^\text{smoothformalseriesclosure}\infty-\mathrm{Toposes}^{\mathrm{ringed},\mathrm{Commutativealgebra}_{\mathrm{simplicial}}(\mathrm{Ind}\mathrm{Banach}_R)}_{\mathrm{Commutativealgebra}_{\mathrm{simplicial}}(\mathrm{Ind}\mathrm{Banach}_R)^\mathrm{opposite},\mathrm{Grotopology,homotopyepimorphism}}}.\\
\mathrm{Ind}\mathrm{\sharp Quasicoherent}^{\text{presheaf}}_{\mathrm{Ind}^\text{smoothformalseriesclosure}\infty-\mathrm{Toposes}^{\mathrm{ringed},\mathrm{Commutativealgebra}_{\mathrm{simplicial}}(\mathrm{Ind}^m\mathrm{Banach}_R)}_{\mathrm{Commutativealgebra}_{\mathrm{simplicial}}(\mathrm{Ind}^m\mathrm{Banach}_R)^\mathrm{opposite},\mathrm{Grotopology,homotopyepimorphism}}},\\
\end{align}
are expected to be $\infty$-sheaves as long as one considers in the admissible situations the corresponding \v{C}ech $\infty$-descent for general seminormed modules as in \cite[Section 9.3]{KKM} and \cite{12BBK} such as Bambozzi-Kremnizer spaces in \cite{12BK}. Therefore we have:\\

\begin{proposition}
The motivic complex $\infty$-presheaf 
\begin{align}
\mathrm{Kan}_{\mathrm{Left}}\Delta_{?/A,\mathrm{functionalanalytic,KKM},\mathrm{BBM,formalanalytification}}(\mathcal{O}),\\
\end{align}
as well as the corresponding Hodge-Tate $\infty$-presheaf as in \cite{12BS}{\footnote{Before the Ben-Bassat-Mukherjee $p$-adic formal analytification we take the corresponding derived $(p,I)$-completion.}}:
\begin{align}
\mathrm{Kan}_{\mathrm{Left}}&\Delta_{?/A,\mathrm{functionalanalytic,KKM},\mathrm{BBM,formalanalytification}}(\mathcal{O})^{\mathrm{HodgeTate}}\\
&:=\mathrm{Kan}_{\mathrm{Left}}\Delta_{?/A,\mathrm{functionalanalytic,KKM},\mathrm{BBM,formalanalytification}}(\mathcal{O})\times A/I\\
&=((\underset{i}{\mathrm{homotopycolimit}}_{\mathrm{sifted},\mathrm{derivedcategory}_{\infty}(A/I-\mathrm{Module})}\mathrm{Kan}_{\mathrm{Left}}\overline{\Delta}_{?/A,\mathrm{functionalanalytic,KKM}}(\mathcal{O}_i))^\wedge\\
&)_\mathrm{BBM,formalanalytification}	
\end{align}
and the preperfectoidization $\infty$-presheaves:
\begin{align}
&(\mathcal{O})^{\mathrm{preperfectoidization}},\\
&(\mathcal{O})^{\mathrm{perfectoidization}},\\
&(\mathcal{O})^{\mathrm{preperfectoidization,derivedcompleted}},\\
&(\mathcal{O})^{\mathrm{perfectoidization,derivedcompleted}}
\end{align}
are $\infty$-sheaves over Bambozzi-Kremnizer topos
\begin{align}
(\mathrm{Spectrumadic}_{\mathrm{BK}}(F),\mathcal{O}),	
\end{align}
attached to any colimit of formal series rings $F$ over $A/I$ in \cite{12BK}.\\	
\end{proposition}

\begin{remark}
One can then define such ring $\mathcal{O}$ to be \textit{preperfectoid} if we have the equivalence:
\begin{align}
\mathcal{O}^{\text{preperfectoidization}} \overset{\sim}{\longrightarrow}	\mathcal{O}.
\end{align}
One can then define such ring $\mathcal{O}$ to be \textit{perfectoid} if we have the equivalence:
\begin{align}
\mathcal{O}^{\text{preperfectoidization}}\times A/I \overset{\sim}{\longrightarrow}	\mathcal{O}.
\end{align}
	
\end{remark}

\newpage

\subsection{Functional Analytic Derived de Rham Complexes for $(\infty,1)$-Analytic Stacks and de Rham Preperfectoidizations}

\indent As in \cite{12LL} we have the comparison between the derived prismatic cohomology and the corresponding derived de Rham cohomology in some very well-defined way which respects the corresponding filtrations, we can then in our situation take the corresponding definition of some derived de Rham complex as the one side of the comparison from \cite{12LL}\footnote{Our goal here is actually study the corresponding derived de Rham period rings and the corresponding applications in $p$-adic Hodge theory extending work of \cite{12DLLZ1}, \cite{12DLLZ2}, \cite{12Sch2} in the motivation from \cite{12GL}, when applyting the construction to derived $p$-adic formal stacks and derived logarithmic $p$-adic formal stacks.}. To be more precise after \cite[Chapitre 3]{12An1}, \cite{12An2}, \cite[Chapter 2, Chapter 8]{12B1}, \cite[Chapter 1]{12Bei}, \cite[Chapter 5]{12G1}, \cite[Chapter 3, Chapter 4]{12GL}, \cite[Chapitre II, Chapitre III]{12Ill1}, \cite[Chapitre VIII]{12Ill2}, \cite[Section 4]{12Qui}, and \cite[Example 5.11, Example 5.12]{BMS2} we define the corresponding:

\begin{definition}
Then we can in the same fashion consider the corresponding derived de Rham complex presheaves for the commutative algebras as in the above (for a given prism $(A,I)$):
\begin{align}
\mathrm{Kan}_{\mathrm{Left}}\mathrm{deRham}_{?/A},	
\end{align}
\footnote{After taking derived $p$-completion.}by the regular corresponding left Kan extension techniques on the following $(\infty,1)$-compactly generated closures of the corresponding polynomials given over $A/I$ with a chosen prism $(A,I)$:\\

\begin{align}
\mathrm{Proj}^\text{smoothformalseriesclosure}\infty-\mathrm{Toposes}^{\mathrm{ringed},\mathrm{commutativealgebra}_{\mathrm{simplicial}}(\mathrm{Ind}\mathrm{Seminormed}_{A/I})}_{\mathrm{Commutativealgebra}_{\mathrm{simplicial}}(\mathrm{Ind}\mathrm{Seminormed}_{A/I})^\mathrm{opposite},\mathrm{Grothendiecktopology,homotopyepimorphism}}. \\
\mathrm{Proj}^\text{smoothformalseriesclosure}\infty-\mathrm{Toposes}^{\mathrm{ringed},\mathrm{Commutativealgebra}_{\mathrm{simplicial}}(\mathrm{Ind}^m\mathrm{Seminormed}_{A/I})}_{\mathrm{Commutativealgebra}_{\mathrm{simplicial}}(\mathrm{Ind}^m\mathrm{Seminormed}_{A/I})^\mathrm{opposite},\mathrm{Grothendiecktopology,homotopyepimorphism}}.\\
\mathrm{Proj}^\text{smoothformalseriesclosure}\infty-\mathrm{Toposes}^{\mathrm{ringed},\mathrm{Commutativealgebra}_{\mathrm{simplicial}}(\mathrm{Ind}\mathrm{Normed}_{A/I})}_{\mathrm{Commutativealgebra}_{\mathrm{simplicial}}(\mathrm{Ind}\mathrm{Normed}_{A/I})^\mathrm{opposite},\mathrm{Grothendiecktopology,homotopyepimorphism}}.\\
\mathrm{Proj}^\text{smoothformalseriesclosure}\infty-\mathrm{Toposes}^{\mathrm{ringed},\mathrm{Commutativealgebra}_{\mathrm{simplicial}}(\mathrm{Ind}^m\mathrm{Normed}_{A/I})}_{\mathrm{Commutativealgebra}_{\mathrm{simplicial}}(\mathrm{Ind}^m\mathrm{Normed}_{A/I})^\mathrm{opposite},\mathrm{Grothendiecktopology,homotopyepimorphism}}.\\
\mathrm{Proj}^\text{smoothformalseriesclosure}\infty-\mathrm{Toposes}^{\mathrm{ringed},\mathrm{Commutativealgebra}_{\mathrm{simplicial}}(\mathrm{Ind}\mathrm{Banach}_{A/I})}_{\mathrm{Commutativealgebra}_{\mathrm{simplicial}}(\mathrm{Ind}\mathrm{Banach}_{A/I})^\mathrm{opposite},\mathrm{Grothendiecktopology,homotopyepimorphism}}.\\
\mathrm{Proj}^\text{smoothformalseriesclosure}\infty-\mathrm{Toposes}^{\mathrm{ringed},\mathrm{Commutativealgebra}_{\mathrm{simplicial}}(\mathrm{Ind}^m\mathrm{Banach}_{A/I})}_{\mathrm{Commutativealgebra}_{\mathrm{simplicial}}(\mathrm{Ind}^m\mathrm{Banach}_{A/I})^\mathrm{opposite},\mathrm{Grothendiecktopology,homotopyepimorphism}}. 
\end{align}
We call the corresponding functors functional analytic derived de Rham complex presheaves which we are going to denote that as in the following:
\begin{align}
\mathrm{Kan}_{\mathrm{Left}}\mathrm{deRham}_{?/A,\text{functionalanalytic,KKM},\text{BBM,formalanalytification}}.	
\end{align}
This would mean the following definition{\footnote{Before the Ben-Bassat-Mukherjee $p$-adic formal analytification we take the corresponding derived $(p,I)$-completion.}}:
\begin{align}
\mathrm{Kan}_{\mathrm{Left}}&\mathrm{deRham}_{?/A,\text{functionalanalytic,KKM},\text{BBM,formalanalytification}}(\mathcal{O})\\
&:=	((\underset{i}{\text{homotopycolimit}}_{\text{sifted},\text{derivedcategory}_{\infty}(A/I-\text{Module})}\mathrm{Kan}_{\mathrm{Left}}\mathrm{deRham}_{?/A,\text{functionalanalytic,KKM}}\\
&(\mathcal{O}_i))^\wedge\\
&)_\text{BBM,formalanalytification}
\end{align}
by writing any object $\mathcal{O}$ as the corresponding colimit 
\begin{center}
$\underset{i}{\text{homotopycolimit}}_\text{sifted}\mathcal{O}_i$.
\end{center}
These are quite large $(\infty,1)$-commutative ring objects in the corresponding $(\infty,1)$-categories for $R=A/I$:
\begin{align}
\mathrm{Ind}\mathrm{\sharp Quasicoherent}^{\text{presheaf}}_{\mathrm{Ind}^\text{smoothformalseriesclosure}\infty-\mathrm{Toposes}^{\mathrm{ringed},\mathrm{commutativealgebra}_{\mathrm{simplicial}}(\mathrm{Ind}\mathrm{Seminormed}_R)}_{\mathrm{Commutativealgebra}_{\mathrm{simplicial}}(\mathrm{Ind}\mathrm{Seminormed}_R)^\mathrm{opposite},\mathrm{Grotopology,homotopyepimorphism}}}. \\
\mathrm{Ind}\mathrm{\sharp Quasicoherent}^{\text{presheaf}}_{\mathrm{Ind}^\text{smoothformalseriesclosure}\infty-\mathrm{Toposes}^{\mathrm{ringed},\mathrm{Commutativealgebra}_{\mathrm{simplicial}}(\mathrm{Ind}^m\mathrm{Seminormed}_R)}_{\mathrm{Commutativealgebra}_{\mathrm{simplicial}}(\mathrm{Ind}^m\mathrm{Seminormed}_R)^\mathrm{opposite},\mathrm{Grotopology,homotopyepimorphism}}}.\\
\mathrm{Ind}\mathrm{\sharp Quasicoherent}^{\text{presheaf}}_{\mathrm{Ind}^\text{smoothformalseriesclosure}\infty-\mathrm{Toposes}^{\mathrm{ringed},\mathrm{Commutativealgebra}_{\mathrm{simplicial}}(\mathrm{Ind}\mathrm{Normed}_R)}_{\mathrm{Commutativealgebra}_{\mathrm{simplicial}}(\mathrm{Ind}\mathrm{Normed}_R)^\mathrm{opposite},\mathrm{Grotopology,homotopyepimorphism}}}.\\
\mathrm{Ind}\mathrm{\sharp Quasicoherent}^{\text{presheaf}}_{\mathrm{Ind}^\text{smoothformalseriesclosure}\infty-\mathrm{Toposes}^{\mathrm{ringed},\mathrm{Commutativealgebra}_{\mathrm{simplicial}}(\mathrm{Ind}^m\mathrm{Normed}_R)}_{\mathrm{Commutativealgebra}_{\mathrm{simplicial}}(\mathrm{Ind}^m\mathrm{Normed}_R)^\mathrm{opposite},\mathrm{Grotopology,homotopyepimorphism}}}.\\
\mathrm{Ind}\mathrm{\sharp Quasicoherent}^{\text{presheaf}}_{\mathrm{Ind}^\text{smoothformalseriesclosure}\infty-\mathrm{Toposes}^{\mathrm{ringed},\mathrm{Commutativealgebra}_{\mathrm{simplicial}}(\mathrm{Ind}\mathrm{Banach}_R)}_{\mathrm{Commutativealgebra}_{\mathrm{simplicial}}(\mathrm{Ind}\mathrm{Banach}_R)^\mathrm{opposite},\mathrm{Grotopology,homotopyepimorphism}}}.\\
\mathrm{Ind}\mathrm{\sharp Quasicoherent}^{\text{presheaf}}_{\mathrm{Ind}^\text{smoothformalseriesclosure}\infty-\mathrm{Toposes}^{\mathrm{ringed},\mathrm{Commutativealgebra}_{\mathrm{simplicial}}(\mathrm{Ind}^m\mathrm{Banach}_R)}_{\mathrm{Commutativealgebra}_{\mathrm{simplicial}}(\mathrm{Ind}^m\mathrm{Banach}_R)^\mathrm{opposite},\mathrm{Grotopology,homotopyepimorphism}}},\\ 
\end{align}
after taking the formal series ring left Kan extension analytification from \cite[Section 4.2]{BBM}, which is defined by taking the left Kan extension to all the $(\infty,1)$-ring objects in the $\infty$-derived category of all $A$-modules from formal series rings over $A$.\\
\end{definition}

\noindent Now we apply this construction to any $A/I$ formal scheme $(X,\mathcal{O}_X)$. But in order to use prismatic technology to reach the corresponding period derived de Rham sheaves, we need to consider more, this would be some definition for 'perfectoidizations':

\begin{definition}
Let $(A,I)$ be a perfectoid prism, and we consider any $\mathrm{E}_\infty$-ring $\mathcal{O}$ in the following
\begin{align}
\mathrm{Proj}^\text{smoothformalseriesclosure}\infty-\mathrm{Toposes}^{\mathrm{ringed},\mathrm{commutativealgebra}_{\mathrm{simplicial}}(\mathrm{Ind}\mathrm{Seminormed}_{A/I})}_{\mathrm{Commutativealgebra}_{\mathrm{simplicial}}(\mathrm{Ind}\mathrm{Seminormed}_{A/I})^\mathrm{opposite},\mathrm{Grothendiecktopology,homotopyepimorphism}}. \\
\mathrm{Proj}^\text{smoothformalseriesclosure}\infty-\mathrm{Toposes}^{\mathrm{ringed},\mathrm{Commutativealgebra}_{\mathrm{simplicial}}(\mathrm{Ind}^m\mathrm{Seminormed}_{A/I})}_{\mathrm{Commutativealgebra}_{\mathrm{simplicial}}(\mathrm{Ind}^m\mathrm{Seminormed}_{A/I})^\mathrm{opposite},\mathrm{Grothendiecktopology,homotopyepimorphism}}.\\
\mathrm{Proj}^\text{smoothformalseriesclosure}\infty-\mathrm{Toposes}^{\mathrm{ringed},\mathrm{Commutativealgebra}_{\mathrm{simplicial}}(\mathrm{Ind}\mathrm{Normed}_{A/I})}_{\mathrm{Commutativealgebra}_{\mathrm{simplicial}}(\mathrm{Ind}\mathrm{Normed}_{A/I})^\mathrm{opposite},\mathrm{Grothendiecktopology,homotopyepimorphism}}.\\
\mathrm{Proj}^\text{smoothformalseriesclosure}\infty-\mathrm{Toposes}^{\mathrm{ringed},\mathrm{Commutativealgebra}_{\mathrm{simplicial}}(\mathrm{Ind}^m\mathrm{Normed}_{A/I})}_{\mathrm{Commutativealgebra}_{\mathrm{simplicial}}(\mathrm{Ind}^m\mathrm{Normed}_{A/I})^\mathrm{opposite},\mathrm{Grothendiecktopology,homotopyepimorphism}}.\\
\mathrm{Proj}^\text{smoothformalseriesclosure}\infty-\mathrm{Toposes}^{\mathrm{ringed},\mathrm{Commutativealgebra}_{\mathrm{simplicial}}(\mathrm{Ind}\mathrm{Banach}_{A/I})}_{\mathrm{Commutativealgebra}_{\mathrm{simplicial}}(\mathrm{Ind}\mathrm{Banach}_{A/I})^\mathrm{opposite},\mathrm{Grothendiecktopology,homotopyepimorphism}}.\\
\mathrm{Proj}^\text{smoothformalseriesclosure}\infty-\mathrm{Toposes}^{\mathrm{ringed},\mathrm{Commutativealgebra}_{\mathrm{simplicial}}(\mathrm{Ind}^m\mathrm{Banach}_{A/I})}_{\mathrm{Commutativealgebra}_{\mathrm{simplicial}}(\mathrm{Ind}^m\mathrm{Banach}_{A/I})^\mathrm{opposite},\mathrm{Grothendiecktopology,homotopyepimorphism}}. 
\end{align}
Then consider the derived prismatic object:
\begin{align}
\mathrm{Kan}_{\mathrm{Left}}\mathrm{deRham}_{?/A,\text{functionalanalytic,KKM},\text{BBM,formalanalytification}}(\mathcal{O}).
\end{align}	
Then as in \cite[Definition 8.2]{12BS} we have the following de Rham  preperfectoidization:
\begin{align}
&(\mathcal{O})^{\text{deRham,preperfectoidization}}\\
&:=\mathrm{Colimit}(\mathrm{Kan}_{\mathrm{Left}}\mathrm{deRham}_{?/A,\text{functionalanalytic,KKM},\text{BBM,formalanalytification}}(\mathcal{O})\rightarrow \\
&\mathrm{Fro}_*\mathrm{Kan}_{\mathrm{Left}}\mathrm{deRham}_{?/A,\text{functionalanalytic,KKM},\text{BBM,formalanalytification}}(\mathcal{O})\\
&\rightarrow \mathrm{Fro}_* \mathrm{Fro}_*\mathrm{Kan}_{\mathrm{Left}}\mathrm{deRham}_{?/A,\text{functionalanalytic,KKM},\text{BBM,formalanalytification}}(\mathcal{O})\rightarrow...\\
&)^{\text{BBM,formalanalytification}},	
\end{align}
after taking the formal series ring left Kan extension analytification from \cite[Section 4.2]{BBM}, which is defined by taking the left Kan extension to all the $(\infty,1)$-ring object in the $\infty$-derived category of all $A/I$-modules from formal series rings over $A/I$. 
Furthermore one can take derived $(p,I)$-completion to achieve the derived $(p,I)$-completed versions.

\end{definition}

\newpage
\subsection{Functional Analytic Derived Prismatic Complexes for Inductive $(\infty,1)$-Analytic Stacks}

\indent We now promote the construction in the previous sections to the corresponding inductive systems of $(\infty,1)$-ringed toposes level after Lurie \cite{12Lu1}, \cite{12Lu2} and \cite{Lu3} in the $\infty$-category of $\infty$-ringed toposes, Bambozzi-Ben-Bassat-Kremnizer \cite{12BBBK}, Ben-Bassat-Mukherjee \cite{BBM}, Bambozzi-Kremnizer \cite{12BK}, Clausen-Scholze \cite{12CS1} \cite{12CS2} and Kelly-Kremnizer-Mukherjee \cite{KKM} in the $\infty$-cateogory of $\infty$-functional analytic ringed toposes.\\

\indent Now we consider the following $\infty$-categories of the corresponding $\infty$-analytic ringed toposes from Bambozzi-Ben-Bassat-Kremnizer \cite{12BBBK}:\\

\begin{align}
&\infty-\mathrm{Toposes}^{\mathrm{ringed},\mathrm{commutativealgebra}_{\mathrm{simplicial}}(\mathrm{Ind}\mathrm{Seminormed}_?)}_{\mathrm{Commutativealgebra}_{\mathrm{simplicial}}(\mathrm{Ind}\mathrm{Seminormed}_?)^\mathrm{opposite},\mathrm{Grothendiecktopology,homotopyepimorphism}}.\\
&\infty-\mathrm{Toposes}^{\mathrm{ringed},\mathrm{Commutativealgebra}_{\mathrm{simplicial}}(\mathrm{Ind}^m\mathrm{Seminormed}_?)}_{\mathrm{Commutativealgebra}_{\mathrm{simplicial}}(\mathrm{Ind}^m\mathrm{Seminormed}_?)^\mathrm{opposite},\mathrm{Grothendiecktopology,homotopyepimorphism}}.\\
&\infty-\mathrm{Toposes}^{\mathrm{ringed},\mathrm{Commutativealgebra}_{\mathrm{simplicial}}(\mathrm{Ind}\mathrm{Normed}_?)}_{\mathrm{Commutativealgebra}_{\mathrm{simplicial}}(\mathrm{Ind}\mathrm{Normed}_?)^\mathrm{opposite},\mathrm{Grothendiecktopology,homotopyepimorphism}}.\\
&\infty-\mathrm{Toposes}^{\mathrm{ringed},\mathrm{Commutativealgebra}_{\mathrm{simplicial}}(\mathrm{Ind}^m\mathrm{Normed}_?)}_{\mathrm{Commutativealgebra}_{\mathrm{simplicial}}(\mathrm{Ind}^m\mathrm{Normed}_?)^\mathrm{opposite},\mathrm{Grothendiecktopology,homotopyepimorphism}}.\\
&\infty-\mathrm{Toposes}^{\mathrm{ringed},\mathrm{Commutativealgebra}_{\mathrm{simplicial}}(\mathrm{Ind}\mathrm{Banach}_?)}_{\mathrm{Commutativealgebra}_{\mathrm{simplicial}}(\mathrm{Ind}\mathrm{Banach}_?)^\mathrm{opposite},\mathrm{Grothendiecktopology,homotopyepimorphism}}.\\
&\infty-\mathrm{Toposes}^{\mathrm{ringed},\mathrm{Commutativealgebra}_{\mathrm{simplicial}}(\mathrm{Ind}^m\mathrm{Banach}_?)}_{\mathrm{Commutativealgebra}_{\mathrm{simplicial}}(\mathrm{Ind}^m\mathrm{Banach}_?)^\mathrm{opposite},\mathrm{Grothendiecktopology,homotopyepimorphism}}.\\ 
&\mathrm{Ind}^\text{smoothformalseriesclosure}\infty-\mathrm{Toposes}^{\mathrm{ringed},\mathrm{commutativealgebra}_{\mathrm{simplicial}}(\mathrm{Ind}\mathrm{Seminormed}_?)}_{\mathrm{Commutativealgebra}_{\mathrm{simplicial}}(\mathrm{Ind}\mathrm{Seminormed}_?)^\mathrm{opposite},\mathrm{Grothendiecktopology,homotopyepimorphism}}. \\
&\mathrm{Ind}^\text{smoothformalseriesclosure}\infty-\mathrm{Toposes}^{\mathrm{ringed},\mathrm{Commutativealgebra}_{\mathrm{simplicial}}(\mathrm{Ind}^m\mathrm{Seminormed}_?)}_{\mathrm{Commutativealgebra}_{\mathrm{simplicial}}(\mathrm{Ind}^m\mathrm{Seminormed}_?)^\mathrm{opposite},\mathrm{Grothendiecktopology,homotopyepimorphism}}.\\
&\mathrm{Ind}^\text{smoothformalseriesclosure}\infty-\mathrm{Toposes}^{\mathrm{ringed},\mathrm{Commutativealgebra}_{\mathrm{simplicial}}(\mathrm{Ind}\mathrm{Normed}_?)}_{\mathrm{Commutativealgebra}_{\mathrm{simplicial}}(\mathrm{Ind}\mathrm{Normed}_?)^\mathrm{opposite},\mathrm{Grothendiecktopology,homotopyepimorphism}}.\\
&\mathrm{Ind}^\text{smoothformalseriesclosure}\infty-\mathrm{Toposes}^{\mathrm{ringed},\mathrm{Commutativealgebra}_{\mathrm{simplicial}}(\mathrm{Ind}^m\mathrm{Normed}_?)}_{\mathrm{Commutativealgebra}_{\mathrm{simplicial}}(\mathrm{Ind}^m\mathrm{Normed}_?)^\mathrm{opposite},\mathrm{Grothendiecktopology,homotopyepimorphism}}.\\
&\mathrm{Ind}^\text{smoothformalseriesclosure}\infty-\mathrm{Toposes}^{\mathrm{ringed},\mathrm{Commutativealgebra}_{\mathrm{simplicial}}(\mathrm{Ind}\mathrm{Banach}_?)}_{\mathrm{Commutativealgebra}_{\mathrm{simplicial}}(\mathrm{Ind}\mathrm{Banach}_?)^\mathrm{opposite},\mathrm{Grothendiecktopology,homotopyepimorphism}}.\\
&\mathrm{Ind}^\text{smoothformalseriesclosure}\infty-\mathrm{Toposes}^{\mathrm{ringed},\mathrm{Commutativealgebra}_{\mathrm{simplicial}}(\mathrm{Ind}^m\mathrm{Banach}_?)}_{\mathrm{Commutativealgebra}_{\mathrm{simplicial}}(\mathrm{Ind}^m\mathrm{Banach}_?)^\mathrm{opposite},\mathrm{Grothendiecktopology,homotopyepimorphism}}.\\ 
\end{align}

\begin{definition}
\indent One can actually define the derived prismatic cohomology presheaves through derived topological Hochschild cohomology presheaves, derived topological period cohomology presheaves and derived topological cyclic cohomology presheaves as in \cite[Section 2.2, Section 2.3]{12BMS}, \cite[Theorem 1.13]{12BS}:
\begin{align}
	\mathrm{Kan}_{\mathrm{Left}}\mathrm{THH},\mathrm{Kan}_{\mathrm{Left}}\mathrm{TP},\mathrm{Kan}_{\mathrm{Left}}\mathrm{TC},
\end{align}
on the following $(\infty,1)$-compactly generated closures of the corresponding polynomials\footnote{Definitely, we need to put certain norms over in some relatively canonical way, as in \cite[Section 4.2]{BBM} one can basically consider rigid ones and dagger ones, and so on. We restrict to the \textit{formal} one.} given over $A/I$ with a chosen prism $(A,I)$\footnote{In all the following, we assume this prism to be bounded and satisfy that $A/I$ is Banach.}:
\begin{align}
\mathrm{Ind}^\text{smoothformalseriesclosure}\infty-\mathrm{Toposes}^{\mathrm{ringed},\mathrm{commutativealgebra}_{\mathrm{simplicial}}(\mathrm{Ind}\mathrm{Seminormed}_{A/I})}_{\mathrm{Commutativealgebra}_{\mathrm{simplicial}}(\mathrm{Ind}\mathrm{Seminormed}_{A/I})^\mathrm{opposite},\mathrm{Grothendiecktopology,homotopyepimorphism}}. \\
\mathrm{Ind}^\text{smoothformalseriesclosure}\infty-\mathrm{Toposes}^{\mathrm{ringed},\mathrm{Commutativealgebra}_{\mathrm{simplicial}}(\mathrm{Ind}^m\mathrm{Seminormed}_{A/I})}_{\mathrm{Commutativealgebra}_{\mathrm{simplicial}}(\mathrm{Ind}^m\mathrm{Seminormed}_{A/I})^\mathrm{opposite},\mathrm{Grothendiecktopology,homotopyepimorphism}}.\\
\mathrm{Ind}^\text{smoothformalseriesclosure}\infty-\mathrm{Toposes}^{\mathrm{ringed},\mathrm{Commutativealgebra}_{\mathrm{simplicial}}(\mathrm{Ind}\mathrm{Normed}_{A/I})}_{\mathrm{Commutativealgebra}_{\mathrm{simplicial}}(\mathrm{Ind}\mathrm{Normed}_{A/I})^\mathrm{opposite},\mathrm{Grothendiecktopology,homotopyepimorphism}}.\\
\mathrm{Ind}^\text{smoothformalseriesclosure}\infty-\mathrm{Toposes}^{\mathrm{ringed},\mathrm{Commutativealgebra}_{\mathrm{simplicial}}(\mathrm{Ind}^m\mathrm{Normed}_{A/I})}_{\mathrm{Commutativealgebra}_{\mathrm{simplicial}}(\mathrm{Ind}^m\mathrm{Normed}_{A/I})^\mathrm{opposite},\mathrm{Grothendiecktopology,homotopyepimorphism}}.\\
\mathrm{Ind}^\text{smoothformalseriesclosure}\infty-\mathrm{Toposes}^{\mathrm{ringed},\mathrm{Commutativealgebra}_{\mathrm{simplicial}}(\mathrm{Ind}\mathrm{Banach}_{A/I})}_{\mathrm{Commutativealgebra}_{\mathrm{simplicial}}(\mathrm{Ind}\mathrm{Banach}_{A/I})^\mathrm{opposite},\mathrm{Grothendiecktopology,homotopyepimorphism}}.\\
\mathrm{Ind}^\text{smoothformalseriesclosure}\infty-\mathrm{Toposes}^{\mathrm{ringed},\mathrm{Commutativealgebra}_{\mathrm{simplicial}}(\mathrm{Ind}^m\mathrm{Banach}_{A/I})}_{\mathrm{Commutativealgebra}_{\mathrm{simplicial}}(\mathrm{Ind}^m\mathrm{Banach}_{A/I})^\mathrm{opposite},\mathrm{Grothendiecktopology,homotopyepimorphism}}. 
\end{align}
We call the corresponding functors are derived functional analytic Hochschild cohomology presheaves, derived functional analytic period cohomology presheaves and derived functional analytic cyclic cohomology presheaves, which we are going to denote these presheaves as in the following for any $\infty$-ringed topos $(\mathbb{X},\mathcal{O})=\underset{i}{\text{homotopycolimit}}(\mathbb{X}_i,\mathcal{O}_i)$:
\begin{align}
	&\mathrm{Kan}_{\mathrm{Left}}\mathrm{THH}_{\text{functionalanalytic,KKM},\text{BBM,formalanalytification}}(\mathcal{O}):=\\
	&(\underset{i}{\text{homotopylimit}}_{\text{sifted},\text{derivedcategory}_{\infty}(A/I-\text{Module})}\mathrm{Kan}_{\mathrm{Left}}\mathrm{THH}_{\text{functionalanalytic,KKM}}(\mathcal{O}_i)\\
	&)_\text{BBM,formalanalytification},\\
	&\mathrm{Kan}_{\mathrm{Left}}\mathrm{TP}_{\text{functionalanalytic,KKM},\text{BBM,formalanalytification}}(\mathcal{O}):=\\
	&(\underset{i}{\text{homotopylimit}}_{\text{sifted},\text{derivedcategory}_{\infty}(A/I-\text{Module})}\mathrm{Kan}_{\mathrm{Left}}\mathrm{TP}_{\text{functionalanalytic,KKM}}(\mathcal{O}_i)\\
	&)_\text{BBM,formalanalytification},\\
	&\mathrm{Kan}_{\mathrm{Left}}\mathrm{TC}_{\text{functionalanalytic,KKM},\text{BBM,formalanalytification}}(\mathcal{O}):=\\
	&(\underset{i}{\text{homotopylimit}}_{\text{sifted},\text{derivedcategory}_{\infty}(A/I-\text{Module})}\mathrm{Kan}_{\mathrm{Left}}\mathrm{TC}_{\text{functionalanalytic,KKM}}(\mathcal{O}_i)\\
	&)_\text{BBM,formalanalytification},
\end{align}
by writing any object $\mathcal{O}$ as the corresponding limit\footnote{Note that we are working with ind-$\infty$-ringed $\infty$-stacks.} 
\begin{center}
$\underset{i}{\text{homotopylimit}}_\text{sifted}\mathcal{O}_i:=\underset{i}{\text{homotopylimit}}_\text{sifted}\mathcal{O}_{\mathrm{Spectrumadic}_\mathrm{BK}(\mathrm{Formalseriesring}_i)}$.
\end{center}
These are quite large $(\infty,1)$-commutative ring objects in the corresponding $(\infty,1)$-categories for $?=A/I$ over $(\mathbb{X},\mathcal{O})$:

 \begin{align}
\mathrm{Ind}\mathrm{\sharp Quasicoherent}^{\text{presheaf}}_{\mathrm{Ind}^\text{smoothformalseriesclosure}\infty-\mathrm{Toposes}^{\mathrm{ringed},\mathrm{commutativealgebra}_{\mathrm{simplicial}}(\mathrm{Ind}\mathrm{Seminormed}_?)}_{\mathrm{Commutativealgebra}_{\mathrm{simplicial}}(\mathrm{Ind}\mathrm{Seminormed}_?)^\mathrm{opposite},\mathrm{Grotopology,homotopyepimorphism}}}. \\
\mathrm{Ind}\mathrm{\sharp Quasicoherent}^{\text{presheaf}}_{\mathrm{Ind}^\text{smoothformalseriesclosure}\infty-\mathrm{Toposes}^{\mathrm{ringed},\mathrm{Commutativealgebra}_{\mathrm{simplicial}}(\mathrm{Ind}^m\mathrm{Seminormed}_?)}_{\mathrm{Commutativealgebra}_{\mathrm{simplicial}}(\mathrm{Ind}^m\mathrm{Seminormed}_?)^\mathrm{opposite},\mathrm{Grotopology,homotopyepimorphism}}}.\\
\mathrm{Ind}\mathrm{\sharp Quasicoherent}^{\text{presheaf}}_{\mathrm{Ind}^\text{smoothformalseriesclosure}\infty-\mathrm{Toposes}^{\mathrm{ringed},\mathrm{Commutativealgebra}_{\mathrm{simplicial}}(\mathrm{Ind}\mathrm{Normed}_?)}_{\mathrm{Commutativealgebra}_{\mathrm{simplicial}}(\mathrm{Ind}\mathrm{Normed}_?)^\mathrm{opposite},\mathrm{Grotopology,homotopyepimorphism}}}.\\
\mathrm{Ind}\mathrm{\sharp Quasicoherent}^{\text{presheaf}}_{\mathrm{Ind}^\text{smoothformalseriesclosure}\infty-\mathrm{Toposes}^{\mathrm{ringed},\mathrm{Commutativealgebra}_{\mathrm{simplicial}}(\mathrm{Ind}^m\mathrm{Normed}_?)}_{\mathrm{Commutativealgebra}_{\mathrm{simplicial}}(\mathrm{Ind}^m\mathrm{Normed}_?)^\mathrm{opposite},\mathrm{Grotopology,homotopyepimorphism}}}.\\
\mathrm{Ind}\mathrm{\sharp Quasicoherent}^{\text{presheaf}}_{\mathrm{Ind}^\text{smoothformalseriesclosure}\infty-\mathrm{Toposes}^{\mathrm{ringed},\mathrm{Commutativealgebra}_{\mathrm{simplicial}}(\mathrm{Ind}\mathrm{Banach}_?)}_{\mathrm{Commutativealgebra}_{\mathrm{simplicial}}(\mathrm{Ind}\mathrm{Banach}_?)^\mathrm{opposite},\mathrm{Grotopology,homotopyepimorphism}}}.\\
\mathrm{Ind}\mathrm{\sharp Quasicoherent}^{\text{presheaf}}_{\mathrm{Ind}^\text{smoothformalseriesclosure}\infty-\mathrm{Toposes}^{\mathrm{ringed},\mathrm{Commutativealgebra}_{\mathrm{simplicial}}(\mathrm{Ind}^m\mathrm{Banach}_?)}_{\mathrm{Commutativealgebra}_{\mathrm{simplicial}}(\mathrm{Ind}^m\mathrm{Banach}_?)^\mathrm{opposite},\mathrm{Grotopology,homotopyepimorphism}}},\\ 
\end{align}
after taking the formal series ring left Kan extension analytification from \cite[Section 4.2]{BBM}, which is defined by taking the left Kan extension to all the $(\infty,1)$-ring objects in the $\infty$-derived category of all $A$-modules from formal series rings over $A$.
\end{definition}

\

\begin{definition}
Then we can in the same fashion consider the corresponding derived prismatic complex presheaves \cite[Construction 7.6]{12BS}\footnote{One just applies \cite[Construction 7.6]{12BS} and then takes the left Kan extensions.} for the commutative algebras as in the above (for a given prism $(A,I)$):
\begin{align}
\mathrm{Kan}_{\mathrm{Left}}\Delta_{?/A},	
\end{align}
by the regular corresponding left Kan extension techniques on the following $(\infty,1)$-compactly generated closures of the corresponding polynomials given over $A/I$ with a chosen prism $(A,I)$:\\

\begin{align}
\mathrm{Ind}^\text{smoothformalseriesclosure}\infty-\mathrm{Toposes}^{\mathrm{ringed},\mathrm{commutativealgebra}_{\mathrm{simplicial}}(\mathrm{Ind}\mathrm{Seminormed}_{A/I})}_{\mathrm{Commutativealgebra}_{\mathrm{simplicial}}(\mathrm{Ind}\mathrm{Seminormed}_{A/I})^\mathrm{opposite},\mathrm{Grothendiecktopology,homotopyepimorphism}}. \\
\mathrm{Ind}^\text{smoothformalseriesclosure}\infty-\mathrm{Toposes}^{\mathrm{ringed},\mathrm{Commutativealgebra}_{\mathrm{simplicial}}(\mathrm{Ind}^m\mathrm{Seminormed}_{A/I})}_{\mathrm{Commutativealgebra}_{\mathrm{simplicial}}(\mathrm{Ind}^m\mathrm{Seminormed}_{A/I})^\mathrm{opposite},\mathrm{Grothendiecktopology,homotopyepimorphism}}.\\
\mathrm{Ind}^\text{smoothformalseriesclosure}\infty-\mathrm{Toposes}^{\mathrm{ringed},\mathrm{Commutativealgebra}_{\mathrm{simplicial}}(\mathrm{Ind}\mathrm{Normed}_{A/I})}_{\mathrm{Commutativealgebra}_{\mathrm{simplicial}}(\mathrm{Ind}\mathrm{Normed}_{A/I})^\mathrm{opposite},\mathrm{Grothendiecktopology,homotopyepimorphism}}.\\
\mathrm{Ind}^\text{smoothformalseriesclosure}\infty-\mathrm{Toposes}^{\mathrm{ringed},\mathrm{Commutativealgebra}_{\mathrm{simplicial}}(\mathrm{Ind}^m\mathrm{Normed}_{A/I})}_{\mathrm{Commutativealgebra}_{\mathrm{simplicial}}(\mathrm{Ind}^m\mathrm{Normed}_{A/I})^\mathrm{opposite},\mathrm{Grothendiecktopology,homotopyepimorphism}}.\\
\mathrm{Ind}^\text{smoothformalseriesclosure}\infty-\mathrm{Toposes}^{\mathrm{ringed},\mathrm{Commutativealgebra}_{\mathrm{simplicial}}(\mathrm{Ind}\mathrm{Banach}_{A/I})}_{\mathrm{Commutativealgebra}_{\mathrm{simplicial}}(\mathrm{Ind}\mathrm{Banach}_{A/I})^\mathrm{opposite},\mathrm{Grothendiecktopology,homotopyepimorphism}}.\\
\mathrm{Ind}^\text{smoothformalseriesclosure}\infty-\mathrm{Toposes}^{\mathrm{ringed},\mathrm{Commutativealgebra}_{\mathrm{simplicial}}(\mathrm{Ind}^m\mathrm{Banach}_{A/I})}_{\mathrm{Commutativealgebra}_{\mathrm{simplicial}}(\mathrm{Ind}^m\mathrm{Banach}_{A/I})^\mathrm{opposite},\mathrm{Grothendiecktopology,homotopyepimorphism}}. 
\end{align}
We call the corresponding functors functional analytic derived prismatic complex presheaves which we are going to denote that as in the following:
\begin{align}
\mathrm{Kan}_{\mathrm{Left}}\Delta_{?/A,\text{functionalanalytic,KKM},\text{BBM,formalanalytification}}.	
\end{align}
This would mean the following definition{\footnote{Before the Ben-Bassat-Mukherjee $p$-adic formal analytification we take the corresponding derived $(p,I)$-completion.}}:
\begin{align}
&\mathrm{Kan}_{\mathrm{Left}}\Delta_{?/A,\text{functionalanalytic,KKM},\text{BBM,formalanalytification}}(\mathcal{O})\\
&:=	((\underset{i}{\text{homotopylimit}}_{\text{sifted},\text{derivedcategory}_{\infty}(A/I-\text{Module})}\mathrm{Kan}_{\mathrm{Left}}\Delta_{?/A,\text{functionalanalytic,KKM}}(\mathcal{O}_i))^\wedge\\
&)_\text{BBM,formalanalytification}
\end{align}
by writing any object $\mathcal{O}$ as the corresponding limit\footnote{Here we remind the readers of the corresponding foundation here, namely the presheaf $\mathcal{O}_i$ in fact takes the value in Koszul complex taking the following form:
\begin{align}
&\mathrm{Koszulcomplex}_{A/I\left<X_1,...,X_l\right>\left<T_1,...,T_m\right>}(a_1-T_1b_1,...,a_m-T_mb_m)\\
&= A/I\left<X_1,...,X_l\right>\left<T_1,...,T_m\right>/^\mathbb{L}(a_1-T_1b_1,...,a_m-T_mb_m).	
\end{align}
This is actually derived $p$-complete since each homotopy group is derived $p$-complete (from the corresponding Banach structure from $A/I\left<X_1,...,X_l\right>$ induced from the $p$-adic topology). Namely the definition of the presheaf:
\begin{align}
\mathrm{Kan}_{\mathrm{Left}}\Delta_{?/A,\text{functionalanalytic,KKM}}(\mathcal{O}_i)	
\end{align}
is directly the application of the derived prismatic functor from \cite[Construction 7.6]{12BS}.} 
\begin{center}
$\underset{i}{\text{homotopylimit}}_\text{sifted}\mathcal{O}_i$.
\end{center}
These are quite large $(\infty,1)$-commutative ring objects in the corresponding $(\infty,1)$-categories for $R=A/I$:
\begin{align}
\mathrm{Ind}\mathrm{\sharp Quasicoherent}^{\text{presheaf}}_{\mathrm{Ind}^\text{smoothformalseriesclosure}\infty-\mathrm{Toposes}^{\mathrm{ringed},\mathrm{commutativealgebra}_{\mathrm{simplicial}}(\mathrm{Ind}\mathrm{Seminormed}_R)}_{\mathrm{Commutativealgebra}_{\mathrm{simplicial}}(\mathrm{Ind}\mathrm{Seminormed}_R)^\mathrm{opposite},\mathrm{Grotopology,homotopyepimorphism}}}. \\
\mathrm{Ind}\mathrm{\sharp Quasicoherent}^{\text{presheaf}}_{\mathrm{Ind}^\text{smoothformalseriesclosure}\infty-\mathrm{Toposes}^{\mathrm{ringed},\mathrm{Commutativealgebra}_{\mathrm{simplicial}}(\mathrm{Ind}^m\mathrm{Seminormed}_R)}_{\mathrm{Commutativealgebra}_{\mathrm{simplicial}}(\mathrm{Ind}^m\mathrm{Seminormed}_R)^\mathrm{opposite},\mathrm{Grotopology,homotopyepimorphism}}}.\\
\mathrm{Ind}\mathrm{\sharp Quasicoherent}^{\text{presheaf}}_{\mathrm{Ind}^\text{smoothformalseriesclosure}\infty-\mathrm{Toposes}^{\mathrm{ringed},\mathrm{Commutativealgebra}_{\mathrm{simplicial}}(\mathrm{Ind}\mathrm{Normed}_R)}_{\mathrm{Commutativealgebra}_{\mathrm{simplicial}}(\mathrm{Ind}\mathrm{Normed}_R)^\mathrm{opposite},\mathrm{Grotopology,homotopyepimorphism}}}.\\
\mathrm{Ind}\mathrm{\sharp Quasicoherent}^{\text{presheaf}}_{\mathrm{Ind}^\text{smoothformalseriesclosure}\infty-\mathrm{Toposes}^{\mathrm{ringed},\mathrm{Commutativealgebra}_{\mathrm{simplicial}}(\mathrm{Ind}^m\mathrm{Normed}_R)}_{\mathrm{Commutativealgebra}_{\mathrm{simplicial}}(\mathrm{Ind}^m\mathrm{Normed}_R)^\mathrm{opposite},\mathrm{Grotopology,homotopyepimorphism}}}.\\
\mathrm{Ind}\mathrm{\sharp Quasicoherent}^{\text{presheaf}}_{\mathrm{Ind}^\text{smoothformalseriesclosure}\infty-\mathrm{Toposes}^{\mathrm{ringed},\mathrm{Commutativealgebra}_{\mathrm{simplicial}}(\mathrm{Ind}\mathrm{Banach}_R)}_{\mathrm{Commutativealgebra}_{\mathrm{simplicial}}(\mathrm{Ind}\mathrm{Banach}_R)^\mathrm{opposite},\mathrm{Grotopology,homotopyepimorphism}}}.\\
\mathrm{Ind}\mathrm{\sharp Quasicoherent}^{\text{presheaf}}_{\mathrm{Ind}^\text{smoothformalseriesclosure}\infty-\mathrm{Toposes}^{\mathrm{ringed},\mathrm{Commutativealgebra}_{\mathrm{simplicial}}(\mathrm{Ind}^m\mathrm{Banach}_R)}_{\mathrm{Commutativealgebra}_{\mathrm{simplicial}}(\mathrm{Ind}^m\mathrm{Banach}_R)^\mathrm{opposite},\mathrm{Grotopology,homotopyepimorphism}}},\\ 
\end{align}
after taking the formal series ring left Kan extension analytification from \cite[Section 4.2]{BBM}, which is defined by taking the left Kan extension to all the $(\infty,1)$-ring objects in the $\infty$-derived category of all $A$-modules from formal series rings over $A$.\\
\end{definition}

\newpage

\subsection{Functional Analytic Derived Preperfectoidizations for Inductive $(\infty,1)$-Analytic Stacks}

\

\noindent Then as in \cite[Definition 8.2]{12BS} we consider the corresponding perfectoidization in this analytic setting.

\begin{definition}
Let $(A,I)$ be a perfectoid prism, and we consider any $\mathrm{E}_\infty$-ring $\mathcal{O}$ in the following
\begin{align}
\mathrm{Ind}^\text{smoothformalseriesclosure}\infty-\mathrm{Toposes}^{\mathrm{ringed},\mathrm{commutativealgebra}_{\mathrm{simplicial}}(\mathrm{Ind}\mathrm{Seminormed}_{A/I})}_{\mathrm{Commutativealgebra}_{\mathrm{simplicial}}(\mathrm{Ind}\mathrm{Seminormed}_{A/I})^\mathrm{opposite},\mathrm{Grothendiecktopology,homotopyepimorphism}}. \\
\mathrm{Ind}^\text{smoothformalseriesclosure}\infty-\mathrm{Toposes}^{\mathrm{ringed},\mathrm{Commutativealgebra}_{\mathrm{simplicial}}(\mathrm{Ind}^m\mathrm{Seminormed}_{A/I})}_{\mathrm{Commutativealgebra}_{\mathrm{simplicial}}(\mathrm{Ind}^m\mathrm{Seminormed}_{A/I})^\mathrm{opposite},\mathrm{Grothendiecktopology,homotopyepimorphism}}.\\
\mathrm{Ind}^\text{smoothformalseriesclosure}\infty-\mathrm{Toposes}^{\mathrm{ringed},\mathrm{Commutativealgebra}_{\mathrm{simplicial}}(\mathrm{Ind}\mathrm{Normed}_{A/I})}_{\mathrm{Commutativealgebra}_{\mathrm{simplicial}}(\mathrm{Ind}\mathrm{Normed}_{A/I})^\mathrm{opposite},\mathrm{Grothendiecktopology,homotopyepimorphism}}.\\
\mathrm{Ind}^\text{smoothformalseriesclosure}\infty-\mathrm{Toposes}^{\mathrm{ringed},\mathrm{Commutativealgebra}_{\mathrm{simplicial}}(\mathrm{Ind}^m\mathrm{Normed}_{A/I})}_{\mathrm{Commutativealgebra}_{\mathrm{simplicial}}(\mathrm{Ind}^m\mathrm{Normed}_{A/I})^\mathrm{opposite},\mathrm{Grothendiecktopology,homotopyepimorphism}}.\\
\mathrm{Ind}^\text{smoothformalseriesclosure}\infty-\mathrm{Toposes}^{\mathrm{ringed},\mathrm{Commutativealgebra}_{\mathrm{simplicial}}(\mathrm{Ind}\mathrm{Banach}_{A/I})}_{\mathrm{Commutativealgebra}_{\mathrm{simplicial}}(\mathrm{Ind}\mathrm{Banach}_{A/I})^\mathrm{opposite},\mathrm{Grothendiecktopology,homotopyepimorphism}}.\\
\mathrm{Ind}^\text{smoothformalseriesclosure}\infty-\mathrm{Toposes}^{\mathrm{ringed},\mathrm{Commutativealgebra}_{\mathrm{simplicial}}(\mathrm{Ind}^m\mathrm{Banach}_{A/I})}_{\mathrm{Commutativealgebra}_{\mathrm{simplicial}}(\mathrm{Ind}^m\mathrm{Banach}_{A/I})^\mathrm{opposite},\mathrm{Grothendiecktopology,homotopyepimorphism}}. 
\end{align}
Then consider the derived prismatic object:
\begin{align}
\mathrm{Kan}_{\mathrm{Left}}\Delta_{?/A,\text{functionalanalytic,KKM},\text{BBM,formalanalytification}}(\mathcal{O}).
\end{align}	
Then as in \cite[Definition 8.2]{12BS} we have the following preperfectoidization:
\begin{align}
&(\mathcal{O})^{\text{preperfectoidization}}\\
&:=\mathrm{Colimit}(\mathrm{Kan}_{\mathrm{Left}}\Delta_{?/A,\text{functionalanalytic,KKM},\text{BBM,formalanalytification}}(\mathcal{O})\rightarrow \\
&\mathrm{Fro}_*\mathrm{Kan}_{\mathrm{Left}}\Delta_{?/A,\text{functionalanalytic,KKM},\text{BBM,formalanalytification}}(\mathcal{O})\\
&\rightarrow \mathrm{Fro}_* \mathrm{Fro}_*\mathrm{Kan}_{\mathrm{Left}}\Delta_{?/A,\text{functionalanalytic,KKM},\text{BBM,formalanalytification}}(\mathcal{O})\rightarrow...)^{\text{BBM,formalanalytification}},	
\end{align}
after taking the formal series ring left Kan extension analytification from \cite[Section 4.2]{BBM}, which is defined by taking the left Kan extension to all the $(\infty,1)$-ring object in the $\infty$-derived category of all $A$-modules from formal series rings over $A$. Then we define the corresponding perfectoidization:
\begin{align}
&(\mathcal{O})^{\text{perfectoidization}}\\
&:=\mathrm{Colimit}(\mathrm{Kan}_{\mathrm{Left}}\Delta_{?/A,\text{functionalanalytic,KKM},\text{BBM,formalanalytification}}(\mathcal{O})\longrightarrow \\
&\mathrm{Fro}_*\mathrm{Kan}_{\mathrm{Left}}\Delta_{?/A,\text{functionalanalytic,KKM},\text{BBM,formalanalytification}}(\mathcal{O})\\
&\longrightarrow \mathrm{Fro}_* \mathrm{Fro}_*\mathrm{Kan}_{\mathrm{Left}}\Delta_{?/A,\text{functionalanalytic,KKM},\text{BBM,formalanalytification}}(\mathcal{O})\longrightarrow...)^{\text{BBM,formalanalytification}}\times A/I.	
\end{align}
Furthermore one can take derived $(p,I)$-completion to achieve the derived $(p,I)$-completed versions:
\begin{align}
\mathcal{O}^\text{preperfectoidization,derivedcomplete}:=(\mathcal{O}^\text{preperfectoidization})^{\wedge},\\
\mathcal{O}^\text{perfectoidization,derivedcomplete}:=\mathcal{O}^\text{preperfectoidization,derivedcomplete}\times A/I.\\
\end{align}
These are large $(\infty,1)$-commutative algebra objects in the corresponding categories as in the above, attached to also large $(\infty,1)$-commutative algebra objects. When we apply this to the corresponding sub-$(\infty,1)$-categories of Banach perfectoid objects in \cite{BMS2}, \cite{GR}, \cite{12KL1}, \cite{12KL2}, \cite{12Ked1}, \cite{12Sch3},  we will recover the corresponding distinguished elemental deformation processes defined in \cite{BMS2}, \cite{GR}, \cite{12KL1}, \cite{12KL2}, \cite{12Ked1}, \cite{12Sch3}. 
\end{definition}

\begin{remark}
One can then define such ring $\mathcal{O}$ to be \textit{preperfectoid} if we have the equivalence:
\begin{align}
\mathcal{O}^{\text{preperfectoidization}} \overset{\sim}{\longrightarrow}	\mathcal{O}.
\end{align}
One can then define such ring $\mathcal{O}$ to be \textit{perfectoid} if we have the equivalence:
\begin{align}
\mathcal{O}^{\text{preperfectoidization}}\times A/I \overset{\sim}{\longrightarrow}	\mathcal{O}.
\end{align}
	
\end{remark}

\newpage

\subsection{Functional Analytic Derived de Rham Complexes for Inductive $(\infty,1)$-Analytic Stacks and de Rham Preperfectoidizations}

\indent As in \cite{12LL} we have the comparison between the derived prismatic cohomology and the corresponding derived de Rham cohomology in some very well-defined way which respects the corresponding filtrations, we can then in our situation take the corresponding definition of some derived de Rham complex as the one side of the comparison from \cite{12LL}\footnote{Our goal here is actually study the corresponding derived de Rham period rings and the corresponding applications in $p$-adic Hodge theory extending work of \cite{12DLLZ1}, \cite{12DLLZ2}, \cite{12Sch2} in the motivation from \cite{12GL}, when applyting the construction to derived $p$-adic formal stacks and derived logarithmic $p$-adic formal stacks.}. To be more precise after \cite[Chapitre 3]{12An1}, \cite{12An2}, \cite[Chapter 2, Chapter 8]{12B1}, \cite[Chapter 1]{12Bei}, \cite[Chapter 5]{12G1}, \cite[Chapter 3, Chapter 4]{12GL}, \cite[Chapitre II, Chapitre III]{12Ill1}, \cite[Chapitre VIII]{12Ill2}, \cite[Section 4]{12Qui}, and \cite[Example 5.11, Example 5.12]{BMS2} we define the corresponding:

\begin{definition}
Then we can in the same fashion consider the corresponding derived de Rham complex presheaves for the commutative algebras as in the above (for a given prism $(A,I)$):
\begin{align}
\mathrm{Kan}_{\mathrm{Left}}\mathrm{deRham}_{?/A},	
\end{align}
\footnote{After taking derived $p$-completion.}by the regular corresponding left Kan extension techniques on the following $(\infty,1)$-compactly generated closures of the corresponding polynomials given over $A/I$ with a chosen prism $(A,I)$:\\

\begin{align}
\mathrm{Ind}^\text{smoothformalseriesclosure}\infty-\mathrm{Toposes}^{\mathrm{ringed},\mathrm{commutativealgebra}_{\mathrm{simplicial}}(\mathrm{Ind}\mathrm{Seminormed}_{A/I})}_{\mathrm{Commutativealgebra}_{\mathrm{simplicial}}(\mathrm{Ind}\mathrm{Seminormed}_{A/I})^\mathrm{opposite},\mathrm{Grothendiecktopology,homotopyepimorphism}}. \\
\mathrm{Ind}^\text{smoothformalseriesclosure}\infty-\mathrm{Toposes}^{\mathrm{ringed},\mathrm{Commutativealgebra}_{\mathrm{simplicial}}(\mathrm{Ind}^m\mathrm{Seminormed}_{A/I})}_{\mathrm{Commutativealgebra}_{\mathrm{simplicial}}(\mathrm{Ind}^m\mathrm{Seminormed}_{A/I})^\mathrm{opposite},\mathrm{Grothendiecktopology,homotopyepimorphism}}.\\
\mathrm{Ind}^\text{smoothformalseriesclosure}\infty-\mathrm{Toposes}^{\mathrm{ringed},\mathrm{Commutativealgebra}_{\mathrm{simplicial}}(\mathrm{Ind}\mathrm{Normed}_{A/I})}_{\mathrm{Commutativealgebra}_{\mathrm{simplicial}}(\mathrm{Ind}\mathrm{Normed}_{A/I})^\mathrm{opposite},\mathrm{Grothendiecktopology,homotopyepimorphism}}.\\
\mathrm{Ind}^\text{smoothformalseriesclosure}\infty-\mathrm{Toposes}^{\mathrm{ringed},\mathrm{Commutativealgebra}_{\mathrm{simplicial}}(\mathrm{Ind}^m\mathrm{Normed}_{A/I})}_{\mathrm{Commutativealgebra}_{\mathrm{simplicial}}(\mathrm{Ind}^m\mathrm{Normed}_{A/I})^\mathrm{opposite},\mathrm{Grothendiecktopology,homotopyepimorphism}}.\\
\mathrm{Ind}^\text{smoothformalseriesclosure}\infty-\mathrm{Toposes}^{\mathrm{ringed},\mathrm{Commutativealgebra}_{\mathrm{simplicial}}(\mathrm{Ind}\mathrm{Banach}_{A/I})}_{\mathrm{Commutativealgebra}_{\mathrm{simplicial}}(\mathrm{Ind}\mathrm{Banach}_{A/I})^\mathrm{opposite},\mathrm{Grothendiecktopology,homotopyepimorphism}}.\\
\mathrm{Ind}^\text{smoothformalseriesclosure}\infty-\mathrm{Toposes}^{\mathrm{ringed},\mathrm{Commutativealgebra}_{\mathrm{simplicial}}(\mathrm{Ind}^m\mathrm{Banach}_{A/I})}_{\mathrm{Commutativealgebra}_{\mathrm{simplicial}}(\mathrm{Ind}^m\mathrm{Banach}_{A/I})^\mathrm{opposite},\mathrm{Grothendiecktopology,homotopyepimorphism}}. 
\end{align}
We call the corresponding functors functional analytic derived de Rham complex presheaves which we are going to denote that as in the following:
\begin{align}
\mathrm{Kan}_{\mathrm{Left}}\mathrm{deRham}_{?/A,\text{functionalanalytic,KKM},\text{BBM,formalanalytification}}.	
\end{align}
This would mean the following definition{\footnote{Before the Ben-Bassat-Mukherjee $p$-adic formal analytification we take the corresponding derived $(p,I)$-completion.}}:
\begin{align}
\mathrm{Kan}_{\mathrm{Left}}&\mathrm{deRham}_{?/A,\text{functionalanalytic,KKM},\text{BBM,formalanalytification}}(\mathcal{O})\\
&:=	((\underset{i}{\text{homotopylimit}}_{\text{sifted},\text{derivedcategory}_{\infty}(A/I-\text{Module})}\mathrm{Kan}_{\mathrm{Left}}\mathrm{deRham}_{?/A,\text{functionalanalytic,KKM}}\\
&(\mathcal{O}_i))^\wedge\\
&)_\text{BBM,formalanalytification}
\end{align}
by writing any object $\mathcal{O}$ as the corresponding colimit 
\begin{center}
$\underset{i}{\text{homotopylimit}}_\text{sifted}\mathcal{O}_i$.
\end{center}
These are quite large $(\infty,1)$-commutative ring objects in the corresponding $(\infty,1)$-categories for $R=A/I$:
\begin{align}
\mathrm{Ind}\mathrm{\sharp Quasicoherent}^{\text{presheaf}}_{\mathrm{Ind}^\text{smoothformalseriesclosure}\infty-\mathrm{Toposes}^{\mathrm{ringed},\mathrm{commutativealgebra}_{\mathrm{simplicial}}(\mathrm{Ind}\mathrm{Seminormed}_R)}_{\mathrm{Commutativealgebra}_{\mathrm{simplicial}}(\mathrm{Ind}\mathrm{Seminormed}_R)^\mathrm{opposite},\mathrm{Grotopology,homotopyepimorphism}}}. \\
\mathrm{Ind}\mathrm{\sharp Quasicoherent}^{\text{presheaf}}_{\mathrm{Ind}^\text{smoothformalseriesclosure}\infty-\mathrm{Toposes}^{\mathrm{ringed},\mathrm{Commutativealgebra}_{\mathrm{simplicial}}(\mathrm{Ind}^m\mathrm{Seminormed}_R)}_{\mathrm{Commutativealgebra}_{\mathrm{simplicial}}(\mathrm{Ind}^m\mathrm{Seminormed}_R)^\mathrm{opposite},\mathrm{Grotopology,homotopyepimorphism}}}.\\
\mathrm{Ind}\mathrm{\sharp Quasicoherent}^{\text{presheaf}}_{\mathrm{Ind}^\text{smoothformalseriesclosure}\infty-\mathrm{Toposes}^{\mathrm{ringed},\mathrm{Commutativealgebra}_{\mathrm{simplicial}}(\mathrm{Ind}\mathrm{Normed}_R)}_{\mathrm{Commutativealgebra}_{\mathrm{simplicial}}(\mathrm{Ind}\mathrm{Normed}_R)^\mathrm{opposite},\mathrm{Grotopology,homotopyepimorphism}}}.\\
\mathrm{Ind}\mathrm{\sharp Quasicoherent}^{\text{presheaf}}_{\mathrm{Ind}^\text{smoothformalseriesclosure}\infty-\mathrm{Toposes}^{\mathrm{ringed},\mathrm{Commutativealgebra}_{\mathrm{simplicial}}(\mathrm{Ind}^m\mathrm{Normed}_R)}_{\mathrm{Commutativealgebra}_{\mathrm{simplicial}}(\mathrm{Ind}^m\mathrm{Normed}_R)^\mathrm{opposite},\mathrm{Grotopology,homotopyepimorphism}}}.\\
\mathrm{Ind}\mathrm{\sharp Quasicoherent}^{\text{presheaf}}_{\mathrm{Ind}^\text{smoothformalseriesclosure}\infty-\mathrm{Toposes}^{\mathrm{ringed},\mathrm{Commutativealgebra}_{\mathrm{simplicial}}(\mathrm{Ind}\mathrm{Banach}_R)}_{\mathrm{Commutativealgebra}_{\mathrm{simplicial}}(\mathrm{Ind}\mathrm{Banach}_R)^\mathrm{opposite},\mathrm{Grotopology,homotopyepimorphism}}}.\\
\mathrm{Ind}\mathrm{\sharp Quasicoherent}^{\text{presheaf}}_{\mathrm{Ind}^\text{smoothformalseriesclosure}\infty-\mathrm{Toposes}^{\mathrm{ringed},\mathrm{Commutativealgebra}_{\mathrm{simplicial}}(\mathrm{Ind}^m\mathrm{Banach}_R)}_{\mathrm{Commutativealgebra}_{\mathrm{simplicial}}(\mathrm{Ind}^m\mathrm{Banach}_R)^\mathrm{opposite},\mathrm{Grotopology,homotopyepimorphism}}},\\ 
\end{align}
after taking the formal series ring left Kan extension analytification from \cite[Section 4.2]{BBM}, which is defined by taking the left Kan extension to all the $(\infty,1)$-ring objects in the $\infty$-derived category of all $A$-modules from formal series rings over $A$.\\
\end{definition}

\noindent Now we apply this construction to any $A/I$ formal scheme $(X,\mathcal{O}_X)$. But in order to use prismatic technology to reach the corresponding period derived de Rham sheaves, we need to consider more, this would be some definition for 'perfectoidizations':

\begin{definition}
Let $(A,I)$ be a perfectoid prism, and we consider any $\mathrm{E}_\infty$-ring $\mathcal{O}$ in the following
\begin{align}
\mathrm{Ind}^\text{smoothformalseriesclosure}\infty-\mathrm{Toposes}^{\mathrm{ringed},\mathrm{commutativealgebra}_{\mathrm{simplicial}}(\mathrm{Ind}\mathrm{Seminormed}_{A/I})}_{\mathrm{Commutativealgebra}_{\mathrm{simplicial}}(\mathrm{Ind}\mathrm{Seminormed}_{A/I})^\mathrm{opposite},\mathrm{Grothendiecktopology,homotopyepimorphism}}. \\
\mathrm{Ind}^\text{smoothformalseriesclosure}\infty-\mathrm{Toposes}^{\mathrm{ringed},\mathrm{Commutativealgebra}_{\mathrm{simplicial}}(\mathrm{Ind}^m\mathrm{Seminormed}_{A/I})}_{\mathrm{Commutativealgebra}_{\mathrm{simplicial}}(\mathrm{Ind}^m\mathrm{Seminormed}_{A/I})^\mathrm{opposite},\mathrm{Grothendiecktopology,homotopyepimorphism}}.\\
\mathrm{Ind}^\text{smoothformalseriesclosure}\infty-\mathrm{Toposes}^{\mathrm{ringed},\mathrm{Commutativealgebra}_{\mathrm{simplicial}}(\mathrm{Ind}\mathrm{Normed}_{A/I})}_{\mathrm{Commutativealgebra}_{\mathrm{simplicial}}(\mathrm{Ind}\mathrm{Normed}_{A/I})^\mathrm{opposite},\mathrm{Grothendiecktopology,homotopyepimorphism}}.\\
\mathrm{Ind}^\text{smoothformalseriesclosure}\infty-\mathrm{Toposes}^{\mathrm{ringed},\mathrm{Commutativealgebra}_{\mathrm{simplicial}}(\mathrm{Ind}^m\mathrm{Normed}_{A/I})}_{\mathrm{Commutativealgebra}_{\mathrm{simplicial}}(\mathrm{Ind}^m\mathrm{Normed}_{A/I})^\mathrm{opposite},\mathrm{Grothendiecktopology,homotopyepimorphism}}.\\
\mathrm{Ind}^\text{smoothformalseriesclosure}\infty-\mathrm{Toposes}^{\mathrm{ringed},\mathrm{Commutativealgebra}_{\mathrm{simplicial}}(\mathrm{Ind}\mathrm{Banach}_{A/I})}_{\mathrm{Commutativealgebra}_{\mathrm{simplicial}}(\mathrm{Ind}\mathrm{Banach}_{A/I})^\mathrm{opposite},\mathrm{Grothendiecktopology,homotopyepimorphism}}.\\
\mathrm{Ind}^\text{smoothformalseriesclosure}\infty-\mathrm{Toposes}^{\mathrm{ringed},\mathrm{Commutativealgebra}_{\mathrm{simplicial}}(\mathrm{Ind}^m\mathrm{Banach}_{A/I})}_{\mathrm{Commutativealgebra}_{\mathrm{simplicial}}(\mathrm{Ind}^m\mathrm{Banach}_{A/I})^\mathrm{opposite},\mathrm{Grothendiecktopology,homotopyepimorphism}}. 
\end{align}
Then consider the derived prismatic object:
\begin{align}
\mathrm{Kan}_{\mathrm{Left}}\mathrm{deRham}_{?/A,\text{functionalanalytic,KKM},\text{BBM,formalanalytification}}(\mathcal{O}).
\end{align}	
Then as in \cite[Definition 8.2]{12BS} we have the following de Rham  preperfectoidization:
\begin{align}
&(\mathcal{O})^{\text{deRham,preperfectoidization}}\\
&:=\mathrm{Colimit}(\mathrm{Kan}_{\mathrm{Left}}\mathrm{deRham}_{?/A,\text{functionalanalytic,KKM},\text{BBM,formalanalytification}}(\mathcal{O})\rightarrow \\
&\mathrm{Fro}_*\mathrm{Kan}_{\mathrm{Left}}\mathrm{deRham}_{?/A,\text{functionalanalytic,KKM},\text{BBM,formalanalytification}}(\mathcal{O})\\
&\rightarrow \mathrm{Fro}_* \mathrm{Fro}_*\mathrm{Kan}_{\mathrm{Left}}\mathrm{deRham}_{?/A,\text{functionalanalytic,KKM},\text{BBM,formalanalytification}}(\mathcal{O})\rightarrow...\\
&)^{\text{BBM,formalanalytification}},	
\end{align}
after taking the formal series ring left Kan extension analytification from \cite[Section 4.2]{BBM}, which is defined by taking the left Kan extension to all the $(\infty,1)$-ring object in the $\infty$-derived category of all $A/I$-modules from formal series rings over $A/I$. 
Furthermore one can take derived $(p,I)$-completion to achieve the derived $(p,I)$-completed versions.

\end{definition}

\newpage

\section{Functional Analytic Noncommutative Motives, Noncommutative Primatic Cohomologies and the Preperfectoidizations}

\indent We now work in the noncommutative setting after \cite{Kon1}, \cite{Ta}, \cite{KR1} and \cite{KR2}, with some philosophy rooted in some noncommtative motives and the corresponding nonabelian applications in noncommutative analytic geometry in the derived sense, and the noncommutative analogs of the corresponding Riemann hypothesis and the corresponding Tamagawa number conjectures. The issue is certainly that the usual Frobenius map looks strange, which tells us of the fact that actually we need to consider really large objects such as the corresponding Topological Hochschild Homologies and the corresponding nearby objects. Here we choose to consider \cite{12NS} in order to apply the constructions to certain $\infty$-rings, which we will call them Fukaya-Kato analytifications from \cite{12FK}. \\

\indent Now we consider the following $\infty$-categories of the corresponding $\infty$-analytic ringed toposes from Bambozzi-Ben-Bassat-Kremnizer \cite{12BBBK} in some paralle way:\\

\begin{align}
&\infty-\mathrm{Toposes}^{\mathrm{ringed},\mathrm{Noncommutativealgebra}_{\mathrm{simplicial}}(\mathrm{Ind}\mathrm{Seminormed}_?)}_{\mathrm{Noncommutativealgebra}_{\mathrm{simplicial}}(\mathrm{Ind}\mathrm{Seminormed}_?)^\mathrm{opposite},\mathrm{Grothendiecktopology,homotopyepimorphism}}.\\
&\infty-\mathrm{Toposes}^{\mathrm{ringed},\mathrm{Noncommutativealgebra}_{\mathrm{simplicial}}(\mathrm{Ind}^m\mathrm{Seminormed}_?)}_{\mathrm{Noncommutativealgebra}_{\mathrm{simplicial}}(\mathrm{Ind}^m\mathrm{Seminormed}_?)^\mathrm{opposite},\mathrm{Grothendiecktopology,homotopyepimorphism}}.\\
&\infty-\mathrm{Toposes}^{\mathrm{ringed},\mathrm{Noncommutativealgebra}_{\mathrm{simplicial}}(\mathrm{Ind}\mathrm{Normed}_?)}_{\mathrm{Noncommutativealgebra}_{\mathrm{simplicial}}(\mathrm{Ind}\mathrm{Normed}_?)^\mathrm{opposite},\mathrm{Grothendiecktopology,homotopyepimorphism}}.\\
&\infty-\mathrm{Toposes}^{\mathrm{ringed},\mathrm{Noncommutativealgebra}_{\mathrm{simplicial}}(\mathrm{Ind}^m\mathrm{Normed}_?)}_{\mathrm{Noncommutativealgebra}_{\mathrm{simplicial}}(\mathrm{Ind}^m\mathrm{Normed}_?)^\mathrm{opposite},\mathrm{Grothendiecktopology,homotopyepimorphism}}.\\
&\infty-\mathrm{Toposes}^{\mathrm{ringed},\mathrm{Noncommutativealgebra}_{\mathrm{simplicial}}(\mathrm{Ind}\mathrm{Banach}_?)}_{\mathrm{Noncommutativealgebra}_{\mathrm{simplicial}}(\mathrm{Ind}\mathrm{Banach}_?)^\mathrm{opposite},\mathrm{Grothendiecktopology,homotopyepimorphism}}.\\
&\infty-\mathrm{Toposes}^{\mathrm{ringed},\mathrm{Noncommutativealgebra}_{\mathrm{simplicial}}(\mathrm{Ind}^m\mathrm{Banach}_?)}_{\mathrm{Noncommutativealgebra}_{\mathrm{simplicial}}(\mathrm{Ind}^m\mathrm{Banach}_?)^\mathrm{opposite},\mathrm{Grothendiecktopology,homotopyepimorphism}}.\\ 
&\mathrm{Proj}^\text{smoothformalseriesclosure}\infty-\mathrm{Toposes}^{\mathrm{ringed},\mathrm{Noncommutativealgebra}_{\mathrm{simplicial}}(\mathrm{Ind}\mathrm{Seminormed}_?)}_{\mathrm{Noncommutativealgebra}_{\mathrm{simplicial}}(\mathrm{Ind}\mathrm{Seminormed}_?)^\mathrm{opposite},\mathrm{Grothendiecktopology,homotopyepimorphism}}. \\
&\mathrm{Proj}^\text{smoothformalseriesclosure}\infty-\mathrm{Toposes}^{\mathrm{ringed},\mathrm{Noncommutativealgebra}_{\mathrm{simplicial}}(\mathrm{Ind}^m\mathrm{Seminormed}_?)}_{\mathrm{Noncommutativealgebra}_{\mathrm{simplicial}}(\mathrm{Ind}^m\mathrm{Seminormed}_?)^\mathrm{opposite},\mathrm{Grothendiecktopology,homotopyepimorphism}}.\\
&\mathrm{Proj}^\text{smoothformalseriesclosure}\infty-\mathrm{Toposes}^{\mathrm{ringed},\mathrm{Noncommutativealgebra}_{\mathrm{simplicial}}(\mathrm{Ind}\mathrm{Normed}_?)}_{\mathrm{Noncommutativealgebra}_{\mathrm{simplicial}}(\mathrm{Ind}\mathrm{Normed}_?)^\mathrm{opposite},\mathrm{Grothendiecktopology,homotopyepimorphism}}.\\
&\mathrm{Proj}^\text{smoothformalseriesclosure}\infty-\mathrm{Toposes}^{\mathrm{ringed},\mathrm{Noncommutativealgebra}_{\mathrm{simplicial}}(\mathrm{Ind}^m\mathrm{Normed}_?)}_{\mathrm{Noncommutativealgebra}_{\mathrm{simplicial}}(\mathrm{Ind}^m\mathrm{Normed}_?)^\mathrm{opposite},\mathrm{Grothendiecktopology,homotopyepimorphism}}.\\
&\mathrm{Proj}^\text{smoothformalseriesclosure}\infty-\mathrm{Toposes}^{\mathrm{ringed},\mathrm{Noncommutativealgebra}_{\mathrm{simplicial}}(\mathrm{Ind}\mathrm{Banach}_?)}_{\mathrm{Noncommutativealgebra}_{\mathrm{simplicial}}(\mathrm{Ind}\mathrm{Banach}_?)^\mathrm{opposite},\mathrm{Grothendiecktopology,homotopyepimorphism}}.\\
&\mathrm{Proj}^\text{smoothformalseriesclosure}\infty-\mathrm{Toposes}^{\mathrm{ringed},\mathrm{Noncommutativealgebra}_{\mathrm{simplicial}}(\mathrm{Ind}^m\mathrm{Banach}_?)}_{\mathrm{Noncommutativealgebra}_{\mathrm{simplicial}}(\mathrm{Ind}^m\mathrm{Banach}_?)^\mathrm{opposite},\mathrm{Grothendiecktopology,homotopyepimorphism}}.\\ 
\end{align}

\indent Starting from the $\infty$-rings, we have:\\

\begin{definition}
\indent One can actually define the derived prismatic cohomologies through derived topological Hochschild cohomologies, derived topological period cohomologies and derived topological cyclic cohomologies as in \cite[Section 2.2, Section 2.3]{12BMS}, \cite[Theorem 1.13]{12BS}:
\begin{align}
	\mathrm{Kan}_{\mathrm{Left}}\mathrm{THH}^\mathrm{noncommutative},\mathrm{Kan}_{\mathrm{Left}}\mathrm{TP}^\mathrm{noncommutative},\mathrm{Kan}_{\mathrm{Left}}\mathrm{TC}^\mathrm{noncommutative},
\end{align}
\footnote{One has the corresponding $p$-completed versions as well.}on the following $(\infty,1)$-compactly generated closures of the corresponding polynomials\footnote{Definitely, we need to put certain norms over in some relatively canonical way, as in \cite[Section 4.2]{BBM} one can basically consider rigid ones and dagger ones, and so on. In this noncommutative setting we do not actually need to fix the type of the analytification, as in commutative setting since we are going to apply directly construction from Bhatt-Scholze and Koshikawa \cite{12BS} and \cite{12Ko1}, though the derived characterization of the prismatic cohomology might not need to be restricted to the corresponding $p$-adic formal scheme situations.} given over $A/I$ with a chosen prism $(A,I)$\footnote{In all the following, we assume this prism to be bounded and satisfy that $A/I$ is Banach.}:

\begin{align}
\mathrm{Object}_{\mathrm{E}_\infty\mathrm{Noncommutativealgebra},\mathrm{Simplicial}}(\mathrm{IndSNorm}_{A/I})^{\mathrm{smoothformalseriesclosure}},\\
\mathrm{Object}_{\mathrm{E}_\infty\mathrm{Noncommutativealgebra},\mathrm{Simplicial}}(\mathrm{Ind}^m\mathrm{SNorm}_{A/I})^{\mathrm{smoothformalseriesclosure}},\\
\mathrm{Object}_{\mathrm{E}_\infty\mathrm{Noncommutativealgebra},\mathrm{Simplicial}}(\mathrm{IndNorm}_{A/I})^{\mathrm{smoothformalseriesclosure}},\\
\mathrm{Object}_{\mathrm{E}_\infty\mathrm{Noncommutativealgebra},\mathrm{Simplicial}}(\mathrm{Ind}^m\mathrm{Norm}_{A/I})^{\mathrm{smoothformalseriesclosure}},\\
\mathrm{Object}_{\mathrm{E}_\infty\mathrm{Noncommutativealgebra},\mathrm{Simplicial}}(\mathrm{IndBan}_{A/I})^{\mathrm{smoothformalseriesclosure}},\\
\mathrm{Object}_{\mathrm{E}_\infty\mathrm{Noncommutativealgebra},\mathrm{Simplicial}}(\mathrm{Ind}^m\mathrm{Ban}_{A/I})^{\mathrm{smoothformalseriesclosure}}.
\end{align}
We call the corresponding functors are derived functional analytic Hochschild cohomologies, derived functional analytic period cohomologies and derived functional analytic cyclic cohomologies, which we are going to denote them as in the following:
\begin{align}
	&\mathrm{Kan}_{\mathrm{Left}}\mathrm{THH}^\mathrm{noncommutative}_{\text{functionalanalytic,KKM},\text{BBM,formalanalytification,nc}}(\mathcal{O}):=\\
	&(\underset{i}{\text{homotopycolimit}}_{\text{sifted},\text{derivedcategory}_{\infty}(A/I-\text{Module})}\mathrm{Kan}_{\mathrm{Left}}\mathrm{THH}^\mathrm{noncommutative}_{\text{functionalanalytic,KKM}}(\mathcal{O}_i)\\
	&)_\text{BBM,formalanalytification,nc},\\
	&\mathrm{Kan}_{\mathrm{Left}}\mathrm{TP}^\mathrm{noncommutative}_{\text{functionalanalytic,KKM},\text{BBM,formalanalytification,nc}}(\mathcal{O}):=\\
	&(\underset{i}{\text{homotopycolimit}}_{\text{sifted},\text{derivedcategory}_{\infty}(A/I-\text{Module})}\mathrm{Kan}_{\mathrm{Left}}\mathrm{TP}^\mathrm{noncommutative}_{\text{functionalanalytic,KKM}}(\mathcal{O}_i)\\
	&)_\text{BBM,formalanalytification,nc},\\
	&\mathrm{Kan}_{\mathrm{Left}}\mathrm{TC}^\mathrm{noncommutative}_{\text{functionalanalytic,KKM},\text{BBM,formalanalytification,nc}}(\mathcal{O}):=\\
	&(\underset{i}{\text{homotopycolimit}}_{\text{sifted},\text{derivedcategory}_{\infty}(A/I-\text{Module})}\mathrm{Kan}_{\mathrm{Left}}\mathrm{TC}^\mathrm{noncommutative}_{\text{functionalanalytic,KKM}}(\mathcal{O}_i)\\
	&)_\text{BBM,formalanalytification,nc},
\end{align}
by writing any object $\mathcal{O}$ as the corresponding colimit 
\begin{center}
$\underset{i}{\text{homotopycolimit}}_\text{sifted}\mathcal{O}_i$.
\end{center}
These are quite large $(\infty,1)$-commutative ring objects in the corresponding $(\infty,1)$-categories for $R=A/I$:

\begin{align}
\mathrm{Object}_{\mathrm{E}_\infty\mathrm{Noncommutativealgebra},\mathrm{Simplicial}}(\mathrm{IndSNorm}_R),\\
\mathrm{Object}_{\mathrm{E}_\infty\mathrm{Noncommutativealgebra},\mathrm{Simplicial}}(\mathrm{Ind}^m\mathrm{SNorm}_R),\\
\mathrm{Object}_{\mathrm{E}_\infty\mathrm{Noncommutativealgebra},\mathrm{Simplicial}}(\mathrm{IndNorm}_R),\\
\mathrm{Object}_{\mathrm{E}_\infty\mathrm{Noncommutativealgebra},\mathrm{Simplicial}}(\mathrm{Ind}^m\mathrm{Norm}_R),\\
\mathrm{Object}_{\mathrm{E}_\infty\mathrm{Noncommutativealgebra},\mathrm{Simplicial}}(\mathrm{IndBan}_R),\\
\mathrm{Object}_{\mathrm{E}_\infty\mathrm{Noncommutativealgebra},\mathrm{Simplicial}}(\mathrm{Ind}^m\mathrm{Ban}_R),
\end{align}
after taking the formal series ring left Kan extension analytification from \cite[Section 4.2]{BBM}, which is defined by taking the left Kan extension to all the $(\infty,1)$-ring objects in the $\infty$-derived category of all $A$-modules from formal series rings over $A$, into:
\begin{align}
\mathrm{Object}_{\mathrm{E}_\infty\mathrm{Noncommutativealgebra},\mathrm{Simplicial}}(\mathrm{IndSNorm}_{\mathbb{F}_1})_A,\\
\mathrm{Object}_{\mathrm{E}_\infty\mathrm{Noncommutativealgebra},\mathrm{Simplicial}}(\mathrm{Ind}^m\mathrm{SNorm}_{\mathbb{F}_1})_A,\\
\mathrm{Object}_{\mathrm{E}_\infty\mathrm{Noncommutativealgebra},\mathrm{Simplicial}}(\mathrm{IndNorm}_{\mathbb{F}_1})_A,\\
\mathrm{Object}_{\mathrm{E}_\infty\mathrm{Noncommutativealgebra},\mathrm{Simplicial}}(\mathrm{Ind}^m\mathrm{Norm}_{\mathbb{F}_1})_A,\\
\mathrm{Object}_{\mathrm{E}_\infty\mathrm{Noncommutativealgebra},\mathrm{Simplicial}}(\mathrm{IndBan}_{\mathbb{F}_1})_A,\\
\mathrm{Object}_{\mathrm{E}_\infty\mathrm{Noncommutativealgebra},\mathrm{Simplicial}}(\mathrm{Ind}^m\mathrm{Ban}_{\mathbb{F}_1})_A.
\end{align}
\end{definition}

\

\begin{definition}
Then we can in the same fashion consider the corresponding derived prismatic complexes \cite[Construction 7.6]{12BS}\footnote{One just applies \cite[Construction 7.6]{12BS} and then takes the left Kan extensions.} for the commutative algebras as in the above (for a given prism $(A,I)$):
\begin{align}
\mathrm{Kan}_{\mathrm{Left}}\Delta_{?/A}:=\mathrm{Kan}_{\mathrm{Left}}\mathrm{TP}^\mathrm{noncommutative,pcomplete}(?/A),	
\end{align}
by the regular corresponding left Kan extension techniques on the following $(\infty,1)$-compactly generated closures of the corresponding polynomials given over $A/I$ with a chosen prism $(A,I)$:

\begin{align}
\mathrm{Object}_{\mathrm{E}_\infty\mathrm{Noncommutativealgebra},\mathrm{Simplicial}}(\mathrm{IndSNorm}_{A/I})^{\mathrm{smoothformalseriesclosure}},\\
\mathrm{Object}_{\mathrm{E}_\infty\mathrm{Noncommutativealgebra},\mathrm{Simplicial}}(\mathrm{Ind}^m\mathrm{SNorm}_{A/I})^{\mathrm{smoothformalseriesclosure}},\\
\mathrm{Object}_{\mathrm{E}_\infty\mathrm{Noncommutativealgebra},\mathrm{Simplicial}}(\mathrm{IndNorm}_{A/I})^{\mathrm{smoothformalseriesclosure}},\\
\mathrm{Object}_{\mathrm{E}_\infty\mathrm{Noncommutativealgebra},\mathrm{Simplicial}}(\mathrm{Ind}^m\mathrm{Norm}_{A/I})^{\mathrm{smoothformalseriesclosure}},\\
\mathrm{Object}_{\mathrm{E}_\infty\mathrm{Noncommutativealgebra},\mathrm{Simplicial}}(\mathrm{IndBan}_{A/I})^{\mathrm{smoothformalseriesclosure}},\\
\mathrm{Object}_{\mathrm{E}_\infty\mathrm{Noncommutativealgebra},\mathrm{Simplicial}}(\mathrm{Ind}^m\mathrm{Ban}_{A/I})^{\mathrm{smoothformalseriesclosure}}.
\end{align}
We call the corresponding functors functional analytic derived prismatic complexes which we are going to denote that as in the following:
\begin{align}
\mathrm{Kan}_{\mathrm{Left}}&\Delta_{?/A,\text{functionalanalytic,KKM},\text{BBM,formalanalytification,nc}}\\
&:=\mathrm{Kan}_{\mathrm{Left}}\mathrm{TP}^\mathrm{noncommutative,pcomplete}(?/A)_{\text{functionalanalytic,KKM},\text{BBM,formalanalytification,nc}}.	
\end{align}
This would mean the following definition{\footnote{Before the Ben-Bassat-Mukherjee $p$-adic formal analytification we take the corresponding $p$-completion.}}:
\begin{align}
&\mathrm{Kan}_{\mathrm{Left}}\Delta_{?/A,\text{functionalanalytic,KKM},\text{BBM,formalanalytification,nc}}(\mathcal{O})\\
&:=	((\underset{i}{\text{homotopycolimit}}_{\text{sifted},\text{derivedcategory}_{\infty}(A/I-\text{Module})}\mathrm{Kan}_{\mathrm{Left}}\Delta_{?/A,\text{functionalanalytic,KKM}}(\mathcal{O}_i))^\wedge\\
&)_\text{BBM,formalanalytification,nc}\\
&=((\underset{i}{\text{homotopycolimit}}_{\text{sifted},\text{derivedcategory}_{\infty}(A/I-\text{Module})}\mathrm{Kan}_{\mathrm{Left}}\mathrm{TP}^\mathrm{noncommutative,pcomplete}\\
&(?/A)_{\text{functionalanalytic,KKM}}(\mathcal{O}_i))^\wedge)_\text{BBM,formalanalytification,nc}\\
\end{align}
by writing any object $\mathcal{O}$ as the corresponding colimit 
\begin{center}
$\underset{i}{\text{homotopycolimit}}_\text{sifted}\mathcal{O}_i$.
\end{center}
These are quite large $(\infty,1)$-commutative ring objects in the corresponding $(\infty,1)$-categories for $R=A/I$:
\begin{align}
\mathrm{Object}_{\mathrm{E}_\infty\mathrm{Noncommutativealgebra},\mathrm{Simplicial}}(\mathrm{IndSNorm}_R),\\
\mathrm{Object}_{\mathrm{E}_\infty\mathrm{Noncommutativealgebra},\mathrm{Simplicial}}(\mathrm{Ind}^m\mathrm{SNorm}_R),\\
\mathrm{Object}_{\mathrm{E}_\infty\mathrm{Noncommutativealgebra},\mathrm{Simplicial}}(\mathrm{IndNorm}_R),\\
\mathrm{Object}_{\mathrm{E}_\infty\mathrm{Noncommutativealgebra},\mathrm{Simplicial}}(\mathrm{Ind}^m\mathrm{Norm}_R),\\
\mathrm{Object}_{\mathrm{E}_\infty\mathrm{Noncommutativealgebra},\mathrm{Simplicial}}(\mathrm{IndBan}_R),\\
\mathrm{Object}_{\mathrm{E}_\infty\mathrm{Noncommutativealgebra},\mathrm{Simplicial}}(\mathrm{Ind}^m\mathrm{Ban}_R),\\
\end{align}
after taking the formal series ring left Kan extension analytification from \cite[Section 4.2]{BBM}, which is defined by taking the left Kan extension to all the $(\infty,1)$-ring objects in the $\infty$-derived category of all $A$-modules from formal series rings over $A$, into:
\begin{align}
\mathrm{Object}_{\mathrm{E}_\infty\mathrm{Noncommutativealgebra},\mathrm{Simplicial}}(\mathrm{IndSNorm}_{\mathbb{F}_1})_A,\\
\mathrm{Object}_{\mathrm{E}_\infty\mathrm{Noncommutativealgebra},\mathrm{Simplicial}}(\mathrm{Ind}^m\mathrm{SNorm}_{\mathbb{F}_1})_A,\\
\mathrm{Object}_{\mathrm{E}_\infty\mathrm{Noncommutativealgebra},\mathrm{Simplicial}}(\mathrm{IndNorm}_{\mathbb{F}_1})_A,\\
\mathrm{Object}_{\mathrm{E}_\infty\mathrm{Noncommutativealgebra},\mathrm{Simplicial}}(\mathrm{Ind}^m\mathrm{Norm}_{\mathbb{F}_1})_A,\\
\mathrm{Object}_{\mathrm{E}_\infty\mathrm{Noncommutativealgebra},\mathrm{Simplicial}}(\mathrm{IndBan}_{\mathbb{F}_1})_A,\\
\mathrm{Object}_{\mathrm{E}_\infty\mathrm{Noncommutativealgebra},\mathrm{Simplicial}}(\mathrm{Ind}^m\mathrm{Ban}_{\mathbb{F}_1})_A.
\end{align}
\end{definition}

\

\indent Then as in \cite[Definition 8.2]{12BS} we consider the corresponding 'perfectoidization' in this analytic setting.

\begin{definition}
Let $(A,I)$ be a perfectoid prism, and we consider any $\mathrm{E}_\infty$-ring $\mathcal{O}$ in the following
\begin{align}
\mathrm{Object}_{\mathrm{E}_\infty\mathrm{Noncommutativealgebra},\mathrm{Simplicial}}(\mathrm{IndSNorm}_{A/I})^{\mathrm{smoothformalseriesclosure}},\\
\mathrm{Object}_{\mathrm{E}_\infty\mathrm{Noncommutativealgebra},\mathrm{Simplicial}}(\mathrm{Ind}^m\mathrm{SNorm}_{A/I})^{\mathrm{smoothformalseriesclosure}},\\
\mathrm{Object}_{\mathrm{E}_\infty\mathrm{Noncommutativealgebra},\mathrm{Simplicial}}(\mathrm{IndNorm}_{A/I})^{\mathrm{smoothformalseriesclosure}},\\
\mathrm{Object}_{\mathrm{E}_\infty\mathrm{Noncommutativealgebra},\mathrm{Simplicial}}(\mathrm{Ind}^m\mathrm{Norm}_{A/I})^{\mathrm{smoothformalseriesclosure}},\\
\mathrm{Object}_{\mathrm{E}_\infty\mathrm{Noncommutativealgebra},\mathrm{Simplicial}}(\mathrm{IndBan}_{A/I})^{\mathrm{smoothformalseriesclosure}},\\
\mathrm{Object}_{\mathrm{E}_\infty\mathrm{Noncommutativealgebra},\mathrm{Simplicial}}(\mathrm{Ind}^m\mathrm{Ban}_{A/I})^{\mathrm{smoothformalseriesclosure}}.
\end{align}
Then consider the derived prismatic object:
\begin{align}
\mathrm{Kan}_{\mathrm{Left}}&\Delta_{?/A,\text{functionalanalytic,KKM},\text{BBM,formalanalytification,nc}}(\mathcal{O})\\
&:=\mathrm{Kan}_{\mathrm{Left}}\mathrm{TP}^\mathrm{noncommutative,pcomplete}(?/A)_{\text{functionalanalytic,KKM},\text{BBM,formalanalytification,nc}}(\mathcal{O}).
\end{align}	
Then as in \cite[Definition 8.2]{12BS} we have the following 'preperfectoidization':
\begin{align}
&(\mathcal{O})^{\text{preperfectoidization}}\\
&:=\mathrm{Colimit}(\mathrm{Kan}_{\mathrm{Left}}\Delta_{?/A,\text{functionalanalytic,KKM},\text{BBM,formalanalytification,nc}}(\mathcal{O})\rightarrow \\
&\phi_*\mathrm{Kan}_{\mathrm{Left}}\Delta_{?/A,\text{functionalanalytic,KKM},\text{BBM,formalanalytification,nc}}(\mathcal{O})\\
&\rightarrow \phi_* \phi_*\mathrm{Kan}_{\mathrm{Left}}\Delta_{?/A,\text{functionalanalytic,KKM},\text{BBM,formalanalytification,nc}}(\mathcal{O})\rightarrow...)^{\text{BBM,formalanalytification,nc}},	
\end{align}
after taking the formal series ring left Kan extension analytification from \cite[Section 4.2]{BBM}, which is defined by taking the left Kan extension to all the $(\infty,1)$-ring object in the $\infty$-derived category of all $A$-modules from formal series rings over $A$. Then we define the corresponding 'perfectoidization':
\begin{align}
&(\mathcal{O})^{\text{perfectoidization}}\\
&:=\mathrm{Colimit}(\mathrm{Kan}_{\mathrm{Left}}\Delta_{?/A,\text{functionalanalytic,KKM},\text{BBM,formalanalytification,nc}}(\mathcal{O})\longrightarrow \\
&\phi_*\mathrm{Kan}_{\mathrm{Left}}\Delta_{?/A,\text{functionalanalytic,KKM},\text{BBM,formalanalytification,nc}}(\mathcal{O})\\
&\longrightarrow \phi_* \phi_*\mathrm{Kan}_{\mathrm{Left}}\Delta_{?/A,\text{functionalanalytic,KKM},\text{BBM,formalanalytification,nc}}(\mathcal{O})\longrightarrow...\\
&)^{\text{BBM,formalanalytification,nc}}\times A/I.	
\end{align}
Furthermore one can take derived $(p,I)$-completion to achieve the derived $(p,I)$-completed versions:
\begin{align}
\mathcal{O}^\text{preperfectoidization,derivedcomplete}:=(\mathcal{O}^\text{preperfectoidization})^{\wedge},\\
\mathcal{O}^\text{perfectoidization,derivedcomplete}:=\mathcal{O}^\text{preperfectoidization,derivedcomplete}\times A/I.\\
\end{align}
These are large $(\infty,1)$-noncommutative algebra objects in the corresponding categories as in the above, attached to also large $(\infty,1)$-noncommutative algebra objects. When we apply this to the corresponding sub-$(\infty,1)$-categories of Banach perfectoid objects as in \cite{BMS2}, \cite{GR}, \cite{12KL1}, \cite{12KL2}, \cite{12Ked1}, \cite{12Sch3},  we will recover the corresponding noncommutative analogues of the distinguished elemental deformation processes defined in \cite{BMS2}, \cite{GR}, \cite{12KL1}, \cite{12KL2}, \cite{12Ked1}, \cite{12Sch3}.
\end{definition}

\

\begin{remark}
One can then define such ring $\mathcal{O}$ to be \textit{preperfectoid} if we have the equivalence:
\begin{align}
\mathcal{O}^{\text{preperfectoidization}} \overset{\sim}{\longrightarrow}	\mathcal{O}.
\end{align}
One can then define such ring $\mathcal{O}$ to be \textit{perfectoid} if we have the equivalence:
\begin{align}
\mathcal{O}^{\text{preperfectoidization}}\times A/I \overset{\sim}{\longrightarrow}	\mathcal{O}.
\end{align}
	
\end{remark}

\

\begin{definition}
\indent One can actually define the derived prismatic cohomology presheaves through derived topological Hochschild cohomology presheaves, derived topological period cohomology presheaves and derived topological cyclic cohomology presheaves as in \cite[Section 2.2, Section 2.3]{12BMS}, \cite[Theorem 1.13]{12BS}:
\begin{align}
	\mathrm{Kan}_{\mathrm{Left}}\mathrm{THH}^\mathrm{noncommutative},\mathrm{Kan}_{\mathrm{Left}}\mathrm{TP}^\mathrm{noncommutative},\mathrm{Kan}_{\mathrm{Left}}\mathrm{TC}^\mathrm{noncommutative},
\end{align}
on the following $(\infty,1)$-compactly generated closures of the corresponding polynomials\footnote{Definitely, we need to put certain norms over in some relatively canonical way, as in \cite[Section 4.2]{BBM} one can basically consider rigid ones and dagger ones, and so on. Again in this noncommutative setting we do not actually need to fix the type of the analytification, as in commutative setting since we are going to apply directly construction from Bhatt-Scholze and Koshikawa \cite{12BS} and \cite{12Ko1}, though the derived characterization of the prismatic cohomology might not need to be restricted to the corresponding $p$-adic formal scheme situations.} given over $A/I$ with a chosen prism $(A,I)$\footnote{In all the following, we assume this prism to be bounded and satisfy that $A/I$ is Banach.}:
\begin{align}
\mathrm{Proj}^\text{smoothformalseriesclosure}\infty-\mathrm{Toposes}^{\mathrm{ringed},\mathrm{Noncommutativealgebra}_{\mathrm{simplicial}}(\mathrm{Ind}\mathrm{Seminormed}_{A/I})}_{\mathrm{Noncommutativealgebra}_{\mathrm{simplicial}}(\mathrm{Ind}\mathrm{Seminormed}_{A/I})^\mathrm{opposite},\mathrm{Grothendiecktopology,homotopyepimorphism}}. \\
\mathrm{Proj}^\text{smoothformalseriesclosure}\infty-\mathrm{Toposes}^{\mathrm{ringed},\mathrm{Noncommutativealgebra}_{\mathrm{simplicial}}(\mathrm{Ind}^m\mathrm{Seminormed}_{A/I})}_{\mathrm{Noncommutativealgebra}_{\mathrm{simplicial}}(\mathrm{Ind}^m\mathrm{Seminormed}_{A/I})^\mathrm{opposite},\mathrm{Grothendiecktopology,homotopyepimorphism}}.\\
\mathrm{Proj}^\text{smoothformalseriesclosure}\infty-\mathrm{Toposes}^{\mathrm{ringed},\mathrm{Noncommutativealgebra}_{\mathrm{simplicial}}(\mathrm{Ind}\mathrm{Normed}_{A/I})}_{\mathrm{Noncommutativealgebra}_{\mathrm{simplicial}}(\mathrm{Ind}\mathrm{Normed}_{A/I})^\mathrm{opposite},\mathrm{Grothendiecktopology,homotopyepimorphism}}.\\
\mathrm{Proj}^\text{smoothformalseriesclosure}\infty-\mathrm{Toposes}^{\mathrm{ringed},\mathrm{Noncommutativealgebra}_{\mathrm{simplicial}}(\mathrm{Ind}^m\mathrm{Normed}_{A/I})}_{\mathrm{Noncommutativealgebra}_{\mathrm{simplicial}}(\mathrm{Ind}^m\mathrm{Normed}_{A/I})^\mathrm{opposite},\mathrm{Grothendiecktopology,homotopyepimorphism}}.\\
\mathrm{Proj}^\text{smoothformalseriesclosure}\infty-\mathrm{Toposes}^{\mathrm{ringed},\mathrm{Noncommutativealgebra}_{\mathrm{simplicial}}(\mathrm{Ind}\mathrm{Banach}_{A/I})}_{\mathrm{Noncommutativealgebra}_{\mathrm{simplicial}}(\mathrm{Ind}\mathrm{Banach}_{A/I})^\mathrm{opposite},\mathrm{Grothendiecktopology,homotopyepimorphism}}.\\
\mathrm{Proj}^\text{smoothformalseriesclosure}\infty-\mathrm{Toposes}^{\mathrm{ringed},\mathrm{Noncommutativealgebra}_{\mathrm{simplicial}}(\mathrm{Ind}^m\mathrm{Banach}_{A/I})}_{\mathrm{Noncommutativealgebra}_{\mathrm{simplicial}}(\mathrm{Ind}^m\mathrm{Banach}_{A/I})^\mathrm{opposite},\mathrm{Grothendiecktopology,homotopyepimorphism}}. 
\end{align}
We call the corresponding functors are derived functional analytic Hochschild cohomology presheaves, derived functional analytic period cohomology presheaves and derived functional analytic cyclic cohomology presheaves, which we are going to denote these presheaves as in the following for any $\infty$-ringed topos $(\mathbb{X},\mathcal{O})=\underset{i}{\text{homotopycolimit}}(\mathbb{X}_i,\mathcal{O}_i)$:
\begin{align}
	&\mathrm{Kan}_{\mathrm{Left}}\mathrm{THH}^\mathrm{noncommutative}_{\text{functionalanalytic,KKM},\text{BBM,formalanalytification,nc}}(\mathcal{O}):=\\
	&(\underset{i}{\text{homotopycolimit}}_{\text{sifted},\text{derivedcategory}_{\infty}(A/I-\text{Module})}\mathrm{Kan}_{\mathrm{Left}}\mathrm{THH}^\mathrm{noncommutative}_{\text{functionalanalytic,KKM}}(\mathcal{O}_i)\\
	&)_\text{BBM,formalanalytification,nc},\\
	&\mathrm{Kan}_{\mathrm{Left}}\mathrm{TP}^\mathrm{noncommutative}_{\text{functionalanalytic,KKM},\text{BBM,formalanalytification,nc}}(\mathcal{O}):=\\
	&(\underset{i}{\text{homotopycolimit}}_{\text{sifted},\text{derivedcategory}_{\infty}(A/I-\text{Module})}\mathrm{Kan}_{\mathrm{Left}}\mathrm{TP}^\mathrm{noncommutative}_{\text{functionalanalytic,KKM}}(\mathcal{O}_i)\\
	&)_\text{BBM,formalanalytification,nc},\\
	&\mathrm{Kan}_{\mathrm{Left}}\mathrm{TC}^\mathrm{noncommutative}_{\text{functionalanalytic,KKM},\text{BBM,formalanalytification,nc}}(\mathcal{O}):=\\
	&(\underset{i}{\text{homotopycolimit}}_{\text{sifted},\text{derivedcategory}_{\infty}(A/I-\text{Module})}\mathrm{Kan}_{\mathrm{Left}}\mathrm{TC}^\mathrm{noncommutative}_{\text{functionalanalytic,KKM}}(\mathcal{O}_i)\\
	&)_\text{BBM,formalanalytification,nc},
\end{align}
by writing any object $\mathcal{O}$ as the corresponding colimit\footnote{Here we assume that in the following the presheaves $\mathcal{O}_i$ are taking values in the derived $p$-completed $\mathbb{E}_1$-algebras over $A/I$, which then are generated by those $p$-adic formal series ring $A/I\left<Z_1,...,Z_n\right>$, $n=1,2,...$ by the corresponding derived colimit completion with free variables $Z_1,...,Z_n$, $n=1,2,...$.} 
\begin{center}
$\underset{i}{\text{homotopycolimit}}_\text{sifted}\mathcal{O}_i$.
\end{center}
These are quite large $(\infty,1)$-commutative ring objects in the corresponding $(\infty,1)$-categories for $?=A/I$ over $(\mathbb{X},\mathcal{O})$:

 \begin{align}
\mathrm{Ind}\mathrm{\sharp Quasicoherent}^{\text{presheaf}}_{\mathrm{Ind}^\text{smoothformalseriesclosure}\infty-\mathrm{Toposes}^{\mathrm{ringed},\mathrm{Noncommutativealgebra}_{\mathrm{simplicial}}(\mathrm{Ind}\mathrm{Seminormed}_?)}_{\mathrm{Noncommutativealgebra}_{\mathrm{simplicial}}(\mathrm{Ind}\mathrm{Seminormed}_?)^\mathrm{opposite},\mathrm{Grotopology,homotopyepimorphism}}}. \\
\mathrm{Ind}\mathrm{\sharp Quasicoherent}^{\text{presheaf}}_{\mathrm{Ind}^\text{smoothformalseriesclosure}\infty-\mathrm{Toposes}^{\mathrm{ringed},\mathrm{Noncommutativealgebra}_{\mathrm{simplicial}}(\mathrm{Ind}^m\mathrm{Seminormed}_?)}_{\mathrm{Noncommutativealgebra}_{\mathrm{simplicial}}(\mathrm{Ind}^m\mathrm{Seminormed}_?)^\mathrm{opposite},\mathrm{Grotopology,homotopyepimorphism}}}.\\
\mathrm{Ind}\mathrm{\sharp Quasicoherent}^{\text{presheaf}}_{\mathrm{Ind}^\text{smoothformalseriesclosure}\infty-\mathrm{Toposes}^{\mathrm{ringed},\mathrm{Noncommutativealgebra}_{\mathrm{simplicial}}(\mathrm{Ind}\mathrm{Normed}_?)}_{\mathrm{Noncommutativealgebra}_{\mathrm{simplicial}}(\mathrm{Ind}\mathrm{Normed}_?)^\mathrm{opposite},\mathrm{Grotopology,homotopyepimorphism}}}.\\
\mathrm{Ind}\mathrm{\sharp Quasicoherent}^{\text{presheaf}}_{\mathrm{Ind}^\text{smoothformalseriesclosure}\infty-\mathrm{Toposes}^{\mathrm{ringed},\mathrm{Noncommutativealgebra}_{\mathrm{simplicial}}(\mathrm{Ind}^m\mathrm{Normed}_?)}_{\mathrm{Noncommutativealgebra}_{\mathrm{simplicial}}(\mathrm{Ind}^m\mathrm{Normed}_?)^\mathrm{opposite},\mathrm{Grotopology,homotopyepimorphism}}}.\\
\mathrm{Ind}\mathrm{\sharp Quasicoherent}^{\text{presheaf}}_{\mathrm{Ind}^\text{smoothformalseriesclosure}\infty-\mathrm{Toposes}^{\mathrm{ringed},\mathrm{Noncommutativealgebra}_{\mathrm{simplicial}}(\mathrm{Ind}\mathrm{Banach}_?)}_{\mathrm{Noncommutativealgebra}_{\mathrm{simplicial}}(\mathrm{Ind}\mathrm{Banach}_?)^\mathrm{opposite},\mathrm{Grotopology,homotopyepimorphism}}}.\\
\mathrm{Ind}\mathrm{\sharp Quasicoherent}^{\text{presheaf}}_{\mathrm{Ind}^\text{smoothformalseriesclosure}\infty-\mathrm{Toposes}^{\mathrm{ringed},\mathrm{Noncommutativealgebra}_{\mathrm{simplicial}}(\mathrm{Ind}^m\mathrm{Banach}_?)}_{\mathrm{Noncommutativealgebra}_{\mathrm{simplicial}}(\mathrm{Ind}^m\mathrm{Banach}_?)^\mathrm{opposite},\mathrm{Grotopology,homotopyepimorphism}}},\\ 
\end{align}
after taking the formal series ring left Kan extension analytification from \cite[Section 4.2]{BBM}, which is defined by taking the left Kan extension to all the $(\infty,1)$-ring objects in the $\infty$-derived category of all $A$-modules from formal series rings over $A$.
\end{definition}

\begin{assumption}\mbox{(Technical Assumption)} 
Here we assume that in the following the presheaves $\mathcal{O}_i$ are taking values in the derived $p$-completed $\mathbb{E}_1$-algebras over $A/I$, which then are generated by those $p$-adic formal series ring $A/I\left<Z_1,...,Z_n\right>$, $n=1,2,...$ by the corresponding derived colimit completion with free variables $Z_1,...,Z_n$, $n=1,2,...$. However this assumption does not really matter so significantly in this noncommutative situation since we are just taking the direct definition through $\text{TP}$.  	
\end{assumption}

\

\begin{definition}
Then we can in the same fashion consider the corresponding derived noncommutative prismatic complex presheaves \cite[Construction 7.6]{12BS}\footnote{One just applies \cite[Construction 7.6]{12BS} and then takes the left Kan extensions.} for the commutative algebras as in the above (for a given prism $(A,I)$):
\begin{align}
\mathrm{Kan}_{\mathrm{Left}}\Delta_{?/A}:=\mathrm{Kan}_{\mathrm{Left}}\mathrm{TP}^\mathrm{noncommutative,pcomplete}(?/A),	
\end{align}
\footnote{The motivation for this definition comes from the corresponding commutative picture, namely the corresponding topological characterization of the corresponding prismatic cohomology and the corresponding completed version by using the corresponding Nygaard filtrations. See \cite[Theorem 1.13]{12BS}.}by the regular corresponding left Kan extension techniques on the following $(\infty,1)$-compactly generated closures of the corresponding polynomials given over $A/I$ with a chosen prism $(A,I)$:\\

\begin{align}
\mathrm{Proj}^\text{smoothformalseriesclosure}\infty-\mathrm{Toposes}^{\mathrm{ringed},\mathrm{Noncommutativealgebra}_{\mathrm{simplicial}}(\mathrm{Ind}\mathrm{Seminormed}_{A/I})}_{\mathrm{Noncommutativealgebra}_{\mathrm{simplicial}}(\mathrm{Ind}\mathrm{Seminormed}_{A/I})^\mathrm{opposite},\mathrm{Grothendiecktopology,homotopyepimorphism}}. \\
\mathrm{Proj}^\text{smoothformalseriesclosure}\infty-\mathrm{Toposes}^{\mathrm{ringed},\mathrm{Noncommutativealgebra}_{\mathrm{simplicial}}(\mathrm{Ind}^m\mathrm{Seminormed}_{A/I})}_{\mathrm{Noncommutativealgebra}_{\mathrm{simplicial}}(\mathrm{Ind}^m\mathrm{Seminormed}_{A/I})^\mathrm{opposite},\mathrm{Grothendiecktopology,homotopyepimorphism}}.\\
\mathrm{Proj}^\text{smoothformalseriesclosure}\infty-\mathrm{Toposes}^{\mathrm{ringed},\mathrm{Noncommutativealgebra}_{\mathrm{simplicial}}(\mathrm{Ind}\mathrm{Normed}_{A/I})}_{\mathrm{Noncommutativealgebra}_{\mathrm{simplicial}}(\mathrm{Ind}\mathrm{Normed}_{A/I})^\mathrm{opposite},\mathrm{Grothendiecktopology,homotopyepimorphism}}.\\
\mathrm{Proj}^\text{smoothformalseriesclosure}\infty-\mathrm{Toposes}^{\mathrm{ringed},\mathrm{Noncommutativealgebra}_{\mathrm{simplicial}}(\mathrm{Ind}^m\mathrm{Normed}_{A/I})}_{\mathrm{Noncommutativealgebra}_{\mathrm{simplicial}}(\mathrm{Ind}^m\mathrm{Normed}_{A/I})^\mathrm{opposite},\mathrm{Grothendiecktopology,homotopyepimorphism}}.\\
\mathrm{Proj}^\text{smoothformalseriesclosure}\infty-\mathrm{Toposes}^{\mathrm{ringed},\mathrm{Noncommutativealgebra}_{\mathrm{simplicial}}(\mathrm{Ind}\mathrm{Banach}_{A/I})}_{\mathrm{Noncommutativealgebra}_{\mathrm{simplicial}}(\mathrm{Ind}\mathrm{Banach}_{A/I})^\mathrm{opposite},\mathrm{Grothendiecktopology,homotopyepimorphism}}.\\
\mathrm{Proj}^\text{smoothformalseriesclosure}\infty-\mathrm{Toposes}^{\mathrm{ringed},\mathrm{Noncommutativealgebra}_{\mathrm{simplicial}}(\mathrm{Ind}^m\mathrm{Banach}_{A/I})}_{\mathrm{Noncommutativealgebra}_{\mathrm{simplicial}}(\mathrm{Ind}^m\mathrm{Banach}_{A/I})^\mathrm{opposite},\mathrm{Grothendiecktopology,homotopyepimorphism}}. 
\end{align}
We call the corresponding functors functional analytic derived prismatic complex presheaves which we are going to denote that as in the following:
\begin{align}
\mathrm{Kan}_{\mathrm{Left}}&\Delta_{?/A,\text{functionalanalytic,KKM},\text{BBM,formalanalytification,nc}}\\
&:=\mathrm{Kan}_{\mathrm{Left}}\mathrm{TP}^\mathrm{noncommutative,pcomplete}(?/A)_{\text{functionalanalytic,KKM},\text{BBM,formalanalytification,nc}}.	
\end{align}
This would mean the following definition{\footnote{Before the Ben-Bassat-Mukherjee $p$-adic formal analytification we take the corresponding $p$-completion.}}:
\begin{align}
\mathrm{Kan}_{\mathrm{Left}}&\Delta_{?/A,\text{functionalanalytic,KKM},\text{BBM,formalanalytification,nc}}(\mathcal{O})\\
&:=	((\underset{i}{\text{homotopycolimit}}_{\text{sifted},\text{derivedcategory}_{\infty}(A/I-\text{Module})}\mathrm{Kan}_{\mathrm{Left}}\Delta_{?/A,\text{functionalanalytic,KKM}}(\mathcal{O}_i))^\wedge\\
&)_\text{BBM,formalanalytification,nc}
\end{align}
by writing any object $\mathcal{O}$ as the corresponding colimit 
\begin{center}
$\underset{i}{\text{homotopycolimit}}_\text{sifted}\mathcal{O}_i$.
\end{center}
These are quite large $(\infty,1)$-commutative ring objects in the corresponding $(\infty,1)$-categories for $R=A/I$:
\begin{align}
\mathrm{Ind}\mathrm{\sharp Quasicoherent}^{\text{presheaf}}_{\mathrm{Ind}^\text{smoothformalseriesclosure}\infty-\mathrm{Toposes}^{\mathrm{ringed},\mathrm{Noncommutativealgebra}_{\mathrm{simplicial}}(\mathrm{Ind}\mathrm{Seminormed}_R)}_{\mathrm{Noncommutativealgebra}_{\mathrm{simplicial}}(\mathrm{Ind}\mathrm{Seminormed}_R)^\mathrm{opposite},\mathrm{Grotopology,homotopyepimorphism}}}. \\
\mathrm{Ind}\mathrm{\sharp Quasicoherent}^{\text{presheaf}}_{\mathrm{Ind}^\text{smoothformalseriesclosure}\infty-\mathrm{Toposes}^{\mathrm{ringed},\mathrm{Noncommutativealgebra}_{\mathrm{simplicial}}(\mathrm{Ind}^m\mathrm{Seminormed}_R)}_{\mathrm{Noncommutativealgebra}_{\mathrm{simplicial}}(\mathrm{Ind}^m\mathrm{Seminormed}_R)^\mathrm{opposite},\mathrm{Grotopology,homotopyepimorphism}}}.\\
\mathrm{Ind}\mathrm{\sharp Quasicoherent}^{\text{presheaf}}_{\mathrm{Ind}^\text{smoothformalseriesclosure}\infty-\mathrm{Toposes}^{\mathrm{ringed},\mathrm{Noncommutativealgebra}_{\mathrm{simplicial}}(\mathrm{Ind}\mathrm{Normed}_R)}_{\mathrm{Noncommutativealgebra}_{\mathrm{simplicial}}(\mathrm{Ind}\mathrm{Normed}_R)^\mathrm{opposite},\mathrm{Grotopology,homotopyepimorphism}}}.\\
\mathrm{Ind}\mathrm{\sharp Quasicoherent}^{\text{presheaf}}_{\mathrm{Ind}^\text{smoothformalseriesclosure}\infty-\mathrm{Toposes}^{\mathrm{ringed},\mathrm{Noncommutativealgebra}_{\mathrm{simplicial}}(\mathrm{Ind}^m\mathrm{Normed}_R)}_{\mathrm{Noncommutativealgebra}_{\mathrm{simplicial}}(\mathrm{Ind}^m\mathrm{Normed}_R)^\mathrm{opposite},\mathrm{Grotopology,homotopyepimorphism}}}.\\
\mathrm{Ind}\mathrm{\sharp Quasicoherent}^{\text{presheaf}}_{\mathrm{Ind}^\text{smoothformalseriesclosure}\infty-\mathrm{Toposes}^{\mathrm{ringed},\mathrm{Noncommutativealgebra}_{\mathrm{simplicial}}(\mathrm{Ind}\mathrm{Banach}_R)}_{\mathrm{Noncommutativealgebra}_{\mathrm{simplicial}}(\mathrm{Ind}\mathrm{Banach}_R)^\mathrm{opposite},\mathrm{Grotopology,homotopyepimorphism}}}.\\
\mathrm{Ind}\mathrm{\sharp Quasicoherent}^{\text{presheaf}}_{\mathrm{Ind}^\text{smoothformalseriesclosure}\infty-\mathrm{Toposes}^{\mathrm{ringed},\mathrm{Noncommutativealgebra}_{\mathrm{simplicial}}(\mathrm{Ind}^m\mathrm{Banach}_R)}_{\mathrm{Noncommutativealgebra}_{\mathrm{simplicial}}(\mathrm{Ind}^m\mathrm{Banach}_R)^\mathrm{opposite},\mathrm{Grotopology,homotopyepimorphism}}},\\ 
\end{align}
after taking the formal series ring left Kan extension analytification from \cite[Section 4.2]{BBM}, which is defined by taking the left Kan extension to all the $(\infty,1)$-ring objects in the $\infty$-derived category of all $A$-modules from formal series rings over $A$.\\
\end{definition}

\

\indent Then as in \cite[Definition 8.2]{12BS} we consider the corresponding noncommutative perfectoidization in this analytic setting.

\begin{definition}
Let $(A,I)$ be a perfectoid prism, and we consider any $\mathrm{E}_\infty$-ring $\mathcal{O}$ in the following
\begin{align}
\mathrm{Proj}^\text{smoothformalseriesclosure}\infty-\mathrm{Toposes}^{\mathrm{ringed},\mathrm{Noncommutativealgebra}_{\mathrm{simplicial}}(\mathrm{Ind}\mathrm{Seminormed}_{A/I})}_{\mathrm{Noncommutativealgebra}_{\mathrm{simplicial}}(\mathrm{Ind}\mathrm{Seminormed}_{A/I})^\mathrm{opposite},\mathrm{Grothendiecktopology,homotopyepimorphism}}. \\
\mathrm{Proj}^\text{smoothformalseriesclosure}\infty-\mathrm{Toposes}^{\mathrm{ringed},\mathrm{Noncommutativealgebra}_{\mathrm{simplicial}}(\mathrm{Ind}^m\mathrm{Seminormed}_{A/I})}_{\mathrm{Noncommutativealgebra}_{\mathrm{simplicial}}(\mathrm{Ind}^m\mathrm{Seminormed}_{A/I})^\mathrm{opposite},\mathrm{Grothendiecktopology,homotopyepimorphism}}.\\
\mathrm{Proj}^\text{smoothformalseriesclosure}\infty-\mathrm{Toposes}^{\mathrm{ringed},\mathrm{Noncommutativealgebra}_{\mathrm{simplicial}}(\mathrm{Ind}\mathrm{Normed}_{A/I})}_{\mathrm{Noncommutativealgebra}_{\mathrm{simplicial}}(\mathrm{Ind}\mathrm{Normed}_{A/I})^\mathrm{opposite},\mathrm{Grothendiecktopology,homotopyepimorphism}}.\\
\mathrm{Proj}^\text{smoothformalseriesclosure}\infty-\mathrm{Toposes}^{\mathrm{ringed},\mathrm{Noncommutativealgebra}_{\mathrm{simplicial}}(\mathrm{Ind}^m\mathrm{Normed}_{A/I})}_{\mathrm{Noncommutativealgebra}_{\mathrm{simplicial}}(\mathrm{Ind}^m\mathrm{Normed}_{A/I})^\mathrm{opposite},\mathrm{Grothendiecktopology,homotopyepimorphism}}.\\
\mathrm{Proj}^\text{smoothformalseriesclosure}\infty-\mathrm{Toposes}^{\mathrm{ringed},\mathrm{Noncommutativealgebra}_{\mathrm{simplicial}}(\mathrm{Ind}\mathrm{Banach}_{A/I})}_{\mathrm{Noncommutativealgebra}_{\mathrm{simplicial}}(\mathrm{Ind}\mathrm{Banach}_{A/I})^\mathrm{opposite},\mathrm{Grothendiecktopology,homotopyepimorphism}}.\\
\mathrm{Proj}^\text{smoothformalseriesclosure}\infty-\mathrm{Toposes}^{\mathrm{ringed},\mathrm{Noncommutativealgebra}_{\mathrm{simplicial}}(\mathrm{Ind}^m\mathrm{Banach}_{A/I})}_{\mathrm{Noncommutativealgebra}_{\mathrm{simplicial}}(\mathrm{Ind}^m\mathrm{Banach}_{A/I})^\mathrm{opposite},\mathrm{Grothendiecktopology,homotopyepimorphism}}. 
\end{align}
Then consider the derived prismatic object:
\begin{align}
\mathrm{Kan}_{\mathrm{Left}}\Delta_{?/A,\text{functionalanalytic,KKM},\text{BBM,formalanalytification,nc}}(\mathcal{O}).
\end{align}	
Then as in \cite[Definition 8.2]{12BS} we have the following preperfectoidization:
\begin{align}
&(\mathcal{O})^{\text{preperfectoidization}}\\
&:=\mathrm{Colimit}(\mathrm{Kan}_{\mathrm{Left}}\Delta_{?/A,\text{functionalanalytic,KKM},\text{BBM,formalanalytification,nc}}(\mathcal{O})\rightarrow \\
&\phi_*\mathrm{Kan}_{\mathrm{Left}}\Delta_{?/A,\text{functionalanalytic,KKM},\text{BBM,formalanalytification,nc}}(\mathcal{O})\\
&\rightarrow \phi_* \phi_*\mathrm{Kan}_{\mathrm{Left}}\Delta_{?/A,\text{functionalanalytic,KKM},\text{BBM,formalanalytification,nc}}(\mathcal{O})\rightarrow...)^{\text{BBM,formalanalytification,nc}},	
\end{align}
after taking the formal series ring left Kan extension analytification from \cite[Section 4.2]{BBM}, which is defined by taking the left Kan extension to all the $(\infty,1)$-ring object in the $\infty$-derived category of all $A$-modules from formal series rings over $A$. Then we define the corresponding perfectoidization:
\begin{align}
&(\mathcal{O})^{\text{perfectoidization}}\\
&:=\mathrm{Colimit}(\mathrm{Kan}_{\mathrm{Left}}\Delta_{?/A,\text{functionalanalytic,KKM},\text{BBM,formalanalytification,nc}}(\mathcal{O})\longrightarrow \\
&\phi_*\mathrm{Kan}_{\mathrm{Left}}\Delta_{?/A,\text{functionalanalytic,KKM},\text{BBM,formalanalytification,nc}}(\mathcal{O})\\
&\longrightarrow \phi_* \phi_*\mathrm{Kan}_{\mathrm{Left}}\Delta_{?/A,\text{functionalanalytic,KKM},\text{BBM,formalanalytification,nc}}(\mathcal{O})\longrightarrow...)^{\text{BBM,formalanalytification,nc}}\times A/I.	
\end{align}
Furthermore one can take derived $(p,I)$-completion to achieve the derived $(p,I)$-completed versions:
\begin{align}
\mathcal{O}^\text{preperfectoidization,derivedcomplete}:=(\mathcal{O}^\text{preperfectoidization})^{\wedge},\\
\mathcal{O}^\text{perfectoidization,derivedcomplete}:=\mathcal{O}^\text{preperfectoidization,derivedcomplete}\times A/I.\\
\end{align}
These are large $(\infty,1)$-commutative algebra objects in the corresponding categories as in the above, attached to also large $(\infty,1)$-commutative algebra objects. When we apply this to the corresponding sub-$(\infty,1)$-categories of Banach perfectoid objects in \cite{BMS2}, \cite{GR}, \cite{12KL1}, \cite{12KL2}, \cite{12Ked1}, \cite{12Sch3},  we will have the corresponding noncommutative analogues of the distinguished elemental deformation processes defined in \cite{BMS2}, \cite{GR}, \cite{12KL1}, \cite{12KL2}, \cite{12Ked1}, \cite{12Sch3}. 
\end{definition}

\begin{remark}
One can then define such ring $\mathcal{O}$ to be \textit{preperfectoid} if we have the equivalence:
\begin{align}
\mathcal{O}^{\text{preperfectoidization}} \overset{\sim}{\longrightarrow}	\mathcal{O}.
\end{align}
One can then define such ring $\mathcal{O}$ to be \textit{perfectoid} if we have the equivalence:
\begin{align}
\mathcal{O}^{\text{preperfectoidization}}\times A/I \overset{\sim}{\longrightarrow}	\mathcal{O}.
\end{align}
	
\end{remark}

\end{document}